\batchmode
\makeatletter
\def\input@path{{/tmp/book_jyo/}}
\makeatother
\documentclass[american,paperwidth=210mm,  paperheight=297mm,  top=10mm,  bottom=25mm,  inner=10mm,  outer=10mm,  marginparsep=10mm,  marginparwidth=10mm]{scrbook}
\usepackage[T1]{fontenc}
\usepackage[utf8]{inputenc}
\usepackage{colortbl}
\usepackage{babel}
\usepackage{float}
\usepackage{mathrsfs}
\usepackage{varwidth}
\usepackage{amsmath}
\usepackage{amssymb}
\usepackage{stmaryrd}
\usepackage{stackrel}
\usepackage{splitidx}
\makeindex
\usepackage{graphicx}
\usepackage[a4paper]{geometry}
\usepackage{esint}
\PassOptionsToPackage{normalem}{ulem}
\usepackage{ulem}
\newindex{idx}
\newindex{fam}
\usepackage[pdftex,pdfusetitle,
 bookmarks=true,bookmarksnumbered=false,bookmarksopen=false,
 breaklinks=false,pdfborder={0 0 0},pdfborderstyle={},backref=false,colorlinks=false]
 {hyperref}

\makeatletter

\let\SF@@footnote\footnote
\def\footnote{\ifx\protect\@typeset@protect
    \expandafter\SF@@footnote
  \else
    \expandafter\SF@gobble@opt
  \fi
}
\expandafter\def\csname SF@gobble@opt \endcsname{\@ifnextchar[%]
  \SF@gobble@twobracket
  \@gobble
}
\edef\SF@gobble@opt{\noexpand\protect
  \expandafter\noexpand\csname SF@gobble@opt \endcsname}
\def\SF@gobble@twobracket[#1]#2{}
\providecommand{\tabularnewline}{\\}
\newenvironment{cellvarwidth}[1][t]
    {\begin{varwidth}[#1]{\linewidth}}
    {\@finalstrut\@arstrutbox\end{varwidth}}

\usepackage{amsmath}
\usepackage{tikz}
\usepackage{pgfplots}

\usepackage[T1]{fontenc}
\usepackage{lmodern}
\usepackage[svgnames]{xcolor}
\usepackage{chngcntr}
\usepackage{amssymb}
\usepackage{aliascnt}
\usepackage{tikz}
\usepackage{lipsum}\usetikzlibrary{shapes.misc,calc}
\usepackage[hooks]{tcolorbox}
\usepackage{framed}
\tcbuselibrary{skins,theorems,breakable}

\definecolor{hellgrau}{RGB}{215,215,215}
\definecolor{blau}{RGB}{17,94,140}
\definecolor{orange}{RGB}{229,94,30}
\definecolor{green}{RGB}{17,140,130}
\makeatletter

\tcbset{
  mytheorem/.code args={#1#2#3#4}{
    \refstepcounter{#2}\label{#4}
    \pgfkeysalso{title={\setbox\z@=\hbox{#1~\csname the#2\endcsname\ }\hangindent\wd\z@\hangafter=1 \mbox{#1~\csname the#2\endcsname\ } #3 }}},
}
\newcommand{\refpar}[1]{
	\left(\ref{#1}\right)
}

\newcommand{\salg}[0]{
	$\sigma-$algebra
}
\newcommand{\salgs}[0]{
	$\sigma-$algebras
}
\newcommand{\preds}[0]{
	$\left(\Omega,\mathscr{A},P\right)$
}
\newcommand{\mtcbmaketheorem}[5]{
  \newtcolorbox{#1}[3][]{#3,mytheorem={#2}{#4}{##2}{#5:##3},##1}
}
\newcommand{\mtcbmakeexo}[5]{
  \newtcolorbox{#1}[3][]{
		before={\markright{Exercices and Solutions}},%#2 \csname the#4\endcsname}},
		#3,
		mytheorem={#2}{#4}{##2}{#5:##3},		
		##1
	}
}

\tcbset{
  myproof/.code args={#1#2#3#4}{
	\label{#4}
    \pgfkeysalso{title={\setbox\z@=\hbox{#1~\csname the#2\endcsname\ }\hangindent\wd\z@\hangafter=1 \mbox{#1~\csname the#2\endcsname\ } #3}}},
}

\newcommand{\mtcbmakeproof}[5]{
  \newtcolorbox{#1}[3][]{#3,myproof={#2}{#4}{##2}{#5:##3},##1}
}

\tcbset{
  myrmk/.code args={#1#2#3#4}{
	\label{#4}
    \pgfkeysalso{title={\setbox\z@=\hbox{#1~ }\hangindent\wd\z@\hangafter=1 \mbox{#1~} #3 }}},
}

\newcommand{\mtcbmakermk}[5]{
  \newtcolorbox{#1}[3][]{#3,myrmk={#2}{#4}{##2}{#5:##3},##1}
}

\makeatother

\newaliascnt{lemm}{defi}
\counterwithin{defi}{chapter}
\counterwithin{lemm}{chapter}

\counterwithin{exo}{chapter}
\counterwithin{solexo}{chapter}

\tcbset{
satzstyle/.style={
	standard jigsaw,
    enhanced,
	skin first = enhanced,
    skin middle = enhanced,
    skin last = enhanced,
    breakable,
    freelance,
    boxrule=0pt,
    width=\linewidth,
	colbacktitle=hellgrau,
	frame hidden,
	overlay unbroken={%
    \path[fill=#1]
        ([yshift=-7.5pt]frame.north west) -- ([xshift=7.5pt]frame.north west) --
        (frame.north east) -- (frame.north east|-interior.north east) --
        (frame.north west|-interior.north west) -- cycle;
    },
    interior titled code={
    \path[fill=#1!20!white]
        (frame.west|-interior.north west) -- (frame.east|-interior.north east) --   
        (frame.east|-interior.south east) -- (frame.west|-interior.south west) -- cycle;
    \path[draw=white, line width=1pt] ([xshift=-1pt]frame.west|-interior.north west) -- ([xshift=1pt]frame.east|-interior.north east);
    },
    interior code={
    \path[fill=#1!20!white]
        (frame.west|-interior.north west) -- (frame.east|-interior.north east) --   
        (frame.east|-interior.south east) -- (frame.west|-interior.south west) -- cycle;
    \path[draw=white, line width=1pt] ([xshift=-1pt]frame.west|-interior.north west) -- ([xshift=1pt]frame.east|-interior.north east);
    },
    overlay broken = {
        \path[fill=#1]
        ([yshift=-7.5pt]frame.north west) -- ([xshift=7.5pt]frame.north west) --
        (frame.north east) -- (frame.north east|-interior.north east) --
        (frame.north west|-interior.north west) -- cycle;
		interior titled code={
		\path[fill=#1!20!white]
        	(frame.west|-interior.north west) -- (frame.east|-interior.north east) --   
        	(frame.east|-interior.south east) -- (frame.west|-interior.south west) -- cycle;
    	\path[draw=white, line width=1pt] ([xshift=-1pt]frame.west|-interior.north west) -- ([xshift=1pt]frame.east|-interior.north east);
		},
		interior code={    
		\path[fill=#1!20!white]
        	(frame.west|-interior.north west) -- (frame.east|-interior.north east) --   
        	(frame.east|-interior.south east) -- (frame.west|-interior.south west) -- cycle;
    	\path[draw=white, line width=1pt] ([xshift=-1pt]frame.west|-interior.north west) -- ([xshift=1pt]frame.east|-interior.north east);    
		}
	},	
    fonttitle=\bfseries\sffamily,
    toprule at break=0pt,
}
}

\mtcbmaketheorem{definition}{Definition}{satzstyle=orange}{defi}{df}
\mtcbmaketheorem{theorem}{Theorem}{satzstyle=green}{defi}{th}
\mtcbmaketheorem{corollary}{Corollary}{satzstyle=green!80!white}{defi}{co}
\mtcbmaketheorem{proposition}{Proposition}{satzstyle=green!80!white}{defi}{pr}
\mtcbmaketheorem{lemma}{Lemma}{satzstyle=green!80!white}{defi}{lm}
\mtcbmaketheorem{example}{Example}{satzstyle=blau}{defi}{ex}
\mtcbmakermk{examplen}{Example}{satzstyle=blau}{defi}{ex}
\mtcbmaketheorem{counterexample}{Counterexample}{satzstyle=blau}{defi}{ex}
\mtcbmakermk{solutioncounterexample}{Solution}{satzstyle=blau!50!white}{defi}{rm}
\mtcbmaketheorem{examples}{Examples}{satzstyle=blau}{defi}{exs}
\mtcbmakermk{solutionexample}{Solution}{satzstyle=blau!50!white}{defi}{rm}
\mtcbmakeexo{exercise}{Exercise}{satzstyle=blau}{exo}{exo}
\mtcbmakeexo{solution}{Solution}{satzstyle=blau!80!white}{solexo}{sol}
\mtcbmakeproof{proof}{Proof}{satzstyle=gray}{defi}{pr}
\mtcbmakermk{remark}{Remark}{satzstyle=brown!70!white}{defi}{rm}
\mtcbmakermk{remarks}{Remarks}{satzstyle=brown!70!white}{defi}{rm}
\mtcbmakermk{reminders}{Reminders}{satzstyle=brown!70!white}{defi}{remind}
\mtcbmakermk{reminder}{Reminder}{satzstyle=brown!70!white}{defi}{remind}
\mtcbmakermk{denotation}{Notation}{satzstyle=gray!70!white}{defi}{de}
\mtcbmakermk{denotations}{Notations}{satzstyle=gray!70!white}{defi}{de}
\mtcbmakermk{particular}{Particular Case}{satzstyle=gray!70!white}{defi}{pa}
\mtcbmakermk{application}{Application}{satzstyle=gray!70!white}{defi}{lm}
\mtcbmakermk{complem}{Complement}{satzstyle=gray!70!white}{}{cp}
\mtcbmakermk{summary}{}{satzstyle=blau!70!white}{}{}
\mtcbmakermk{objective}{Objectives}{satzstyle=blau!70!white}{}{}
\usetikzlibrary{decorations.pathreplacing}
\usetikzlibrary{arrows.meta}
\usepackage{pdfpages}

\usetikzlibrary{decorations.pathreplacing}
\usetikzlibrary{arrows.meta}
\usepackage{pgfplots}
\usepackage{multicol}
\pgfplotsset{compat=1.18}
\usepgfplotslibrary{fillbetween}
\tikzset{every node/.style={font=\sffamily}}
\newcommand{\boxeq}[1]{
  \begin{tcolorbox}[
	standard jigsaw,
	colframe=black!90!white,
	colback=white!95!black,
	opacityback=0.5,
	opacityframe=0.5
	]
	  #1
	\end{tcolorbox}
}
\usepackage[automark]{scrlayer-scrpage}
\clearpairofpagestyles

\ihead{\pagemark}
\ohead{\headmark}

\automark[section]{chapter}
\newpairofpagestyles{noheader}{
  \clearpairofpagestyles
}

\newcommand{\mindex}[1]{%
  \begingroup
  \index{#1}%
  \endgroup
}
\makeatother

\begin{document}
\pagenumbering{gobble}

\pagestyle{noheader}

\includepdf{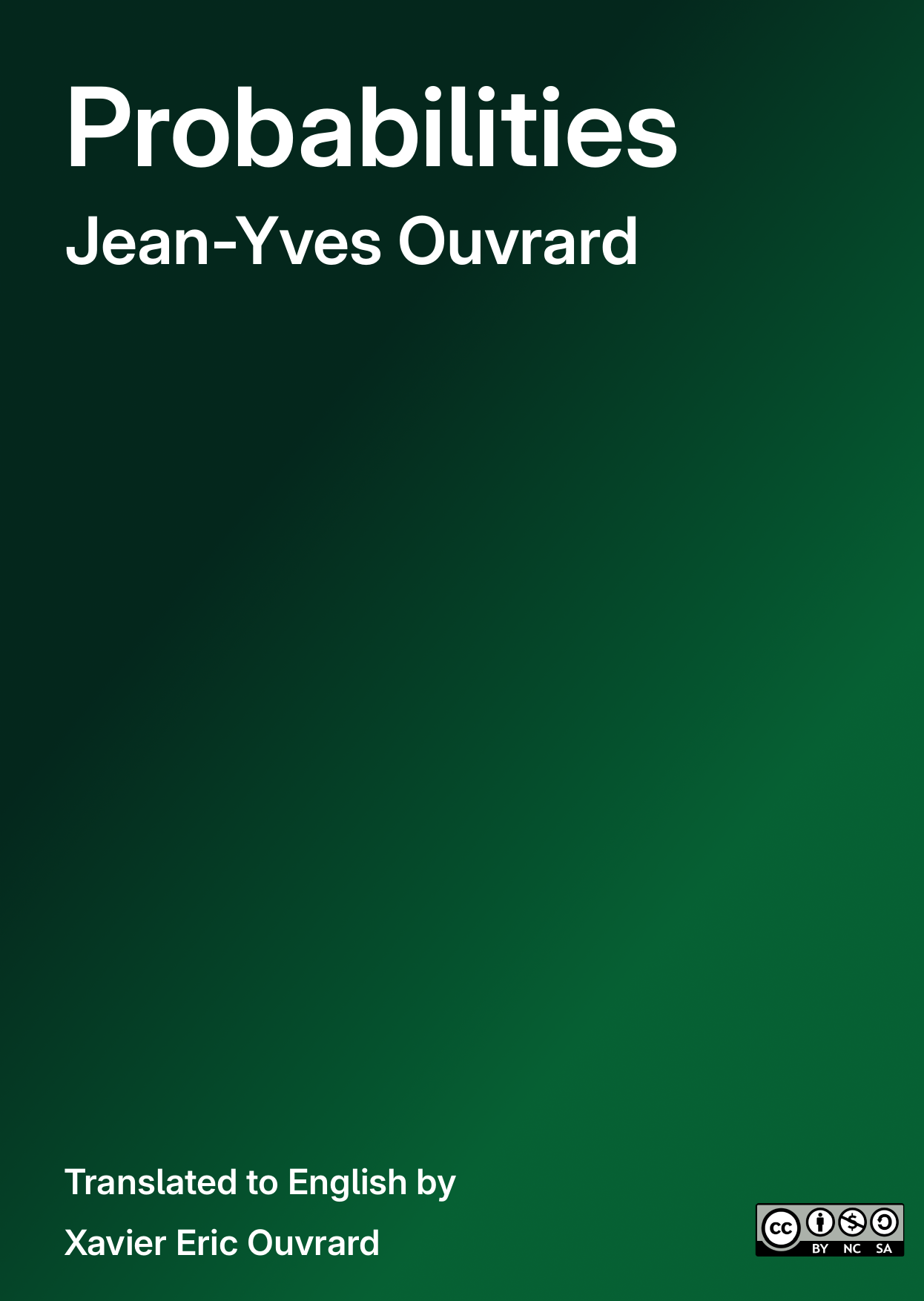}

\vspace*{\fill}

ISBN: 979-10-415-8994-4

EAN: 9791041589944

\vspace{0.5cm}

Legal deposit at Bibliothèque National de France currently running

\vspace{0.5cm}

This work is published under the Creative Commons Attribution--NonCommercial--ShareAlike
4.0 International License (\href{https://creativecommons.org/licenses/by-nc-sa/4.0/}{CC BY-NC-SA 4.0}).

The digital content may be copied, shared, and adapted under the terms
of this license, provided that it is distributed free of charge and
without any form of financial compensation.

\pagenumbering{arabic}\pagestyle{scrheadings}

\tableofcontents{}

\pagestyle{noheader}

\chapter*{Preface\thispagestyle{empty}\pagestyle{empty}}

\textit{Ad memoriam, Jean-Yves Ouvrard (1945-2024)}

\begin{figure}[H]
\centering
\begin{center}\includegraphics[scale=0.2]{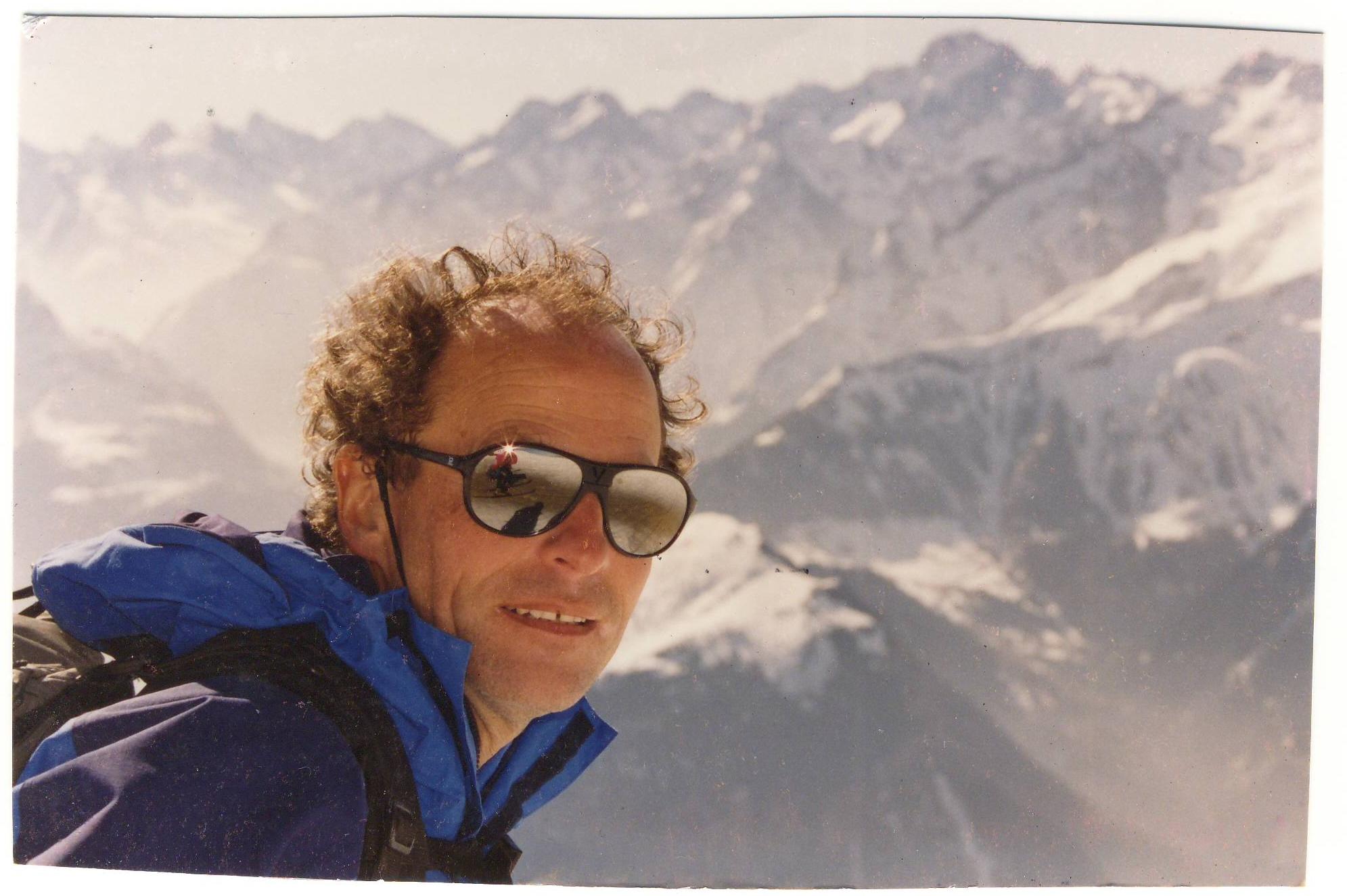}\end{center}
\end{figure}

This book is the English translation of the book titled \textit{Probabilités}
authored by Professor Jean-Yves Ouvrard and published by Cassini Edition
in French. The translation was undertaken by his elder son Dr. Xavier
Ouvrard, during the final days of his father's life. It is intended
as an \textit{ad memoriam} work, to make his father's contributions
to probability theory accessible to a wider audience. 

In translating this book, the translator took the liberty of making
corrections and adaptations where he felt they would improve clarity
and readability. Particular care was taken to ensure that the English
version is comfortable to read on modern digital platforms. 

The translator is deeply grateful to Editions Cassini for allowing
this work to be translated and published freely in English. Part 1
of the original French book remains available through their catalog.
He wants also to give his warmest thanks to Prof. Christophe Leuridant
(University of Grenoble-Alpes, France), former colleague of my father,
who helped him in the first steps of this project, to Prof. Jordan
Stoyanov, who had also known my father during his stay at Grenoble
University, for all his comments and greetings. The translator is
also thankful to Prof. Robert Dalang (EPFL) for his help in making
this work published on Arxiv. Last, the translator is deeply grateful
to his wife, Cécile Brunet-Ouvrard, for her support during all this
translation project, that spanned over 18 months.

If you find this book helpful, the translator would be truly grateful
if you considered making a donation to the French association Leucémie
Espoir (Leukemia Hope):

\href{https://www.leucemie-espoir.org/}{https://www.leucemie-espoir.org/}.

Alternatively, if you are able, you might consider donating blood
to support leukemia patients who truly rely on transfusions to extend
their life.

The English content of this book is released under a \textbf{CC-BY-SA-NC
license}. Anyone who wishes to adapt or use it to educate others about
probability is warmly encouraged to do so.

Any feedback is welcome to help improve this work further.

\subsection*{A note on the Author}

Professor Jean-Yves Ouvrard was an Associate Professor at the Université
Joseph Fourier in Grenoble (Isère, France), where he spent most of
his academic career. He held a \textit{Doctorat d'État} (State PhD)
in Mathematics and specialized in Probability Theory. Professor Ouvrard
was a regular jury member for the French national Agrégation in Mathematics,
a prestigious competitive examination for teacher certification. He
was also deeply involved in preparing students for this examination
at Université Joseph Fourier at Grenoble, contributing actively to
mathematics education.

\subsection*{A note on the Translator}

Dr. Xavier Eric Ouvrard holds a PhD in Science with a specialization
in Computer Science from the University of Geneva. He currently works
at Ecole Polytechnique Fédérale de Lausanne (EPFL) in a research center
focused on sustainable computing. He previously worked at CERN, where
he completed his doctoral research on theoretical developments to
model complex co-occurrence networks and tools for visually navigating
and querying textual data introducing hyper-bag-graphs, an extension
of hypergraphs to multisets family. 

He taught as Agrégé professor in mathematics for over 20 years in
an international secondary school. His academic journey began with
a Master's degree in Process Engineering and a qualification as a
Papermaking Engineer from Institut National Polytechnique de Grenoble
(INPG) (Grenoble, France). He is definitely grateful to his father
for all his support throughout his journey to science.

\chapter*{Introduction\thispagestyle{empty}\pagestyle{empty}}
\begin{quote}
\begin{flushright}Nothing frightens me more that the certainties
that prescribe good law, the normality of the social contract, justice...

Behind the words, we forget those who speak---their looks or their
lack of looks, their smiles or their pretense---signifying everything
left unsaid beneath the ``right-thinking'' word.

Yvon Chaix\\
Rio Theatre, Grenoble, 1995\end{flushright}
\end{quote}
The author has long served as a jury member for the French competitive
examinations in Mathematics teaching, known as the prestigious \textit{Agrégation
de mathématiques}. He contributed to both the external examination
(from 1987 to 1990) and to the internal examination (from 1990 to
1995). He prepared the students for these exams at the University
Joseph Fourier in Grenoble until his retirement. This book is rooted
in that experience. It was written with the intent to help candidates
prepare for various French competitive examinations for secondary-level
mathematics teaching.

The content is divided into two parts, each aligned with the curricula
of the relevant competitive examinations at the time the French editions
of these two books were written.
\begin{itemize}
\item \textbf{Part \ref{part:Introduction-to-Probability}} is aimed at
all candidates and is sufficient to access the French internal Agrégation
and the CAPES examinations. It focuses primarily on discrete probabilities,
with an emphasis on introducing the reader to probabilistic modelling.
A chapter on continuous random variables with densities is also included
although it avoids delving into measure theory. This first part is
also suitable for undergraduate mathematics students.
\item \textbf{Part \ref{part:Deepening-Probability-Theory}} targets students
preparing for the external Agrégation in mathematics. It covers all
topics in that curriculum. Of course, some of the chapters, and particularly
the ones covering martingales and Markov chains, first concern the
candidates choosing at the oral examination the option ``Probabilities
and Statistics''. In this second part, an appendix provides a summary
of measure theory, offering the foundational results necessary to
understand the material. This part is especially valuable for students
pursuing a Master's degree in Mathematics.
\end{itemize}
Practice is essential for any student or competitive examination candidate.
Each chapter concludes with a series of fully worked exercises, covering
all key concepts.

In the original French edition, readers were encouraged to use the
material to craft their own lessons and avoid giving formulaic or
stereotyped presentations at the examination.

The original content in French was published in two separated books
at \href{https://store.cassini.fr/gb/}{Éditions Cassini}.

\subsubsection*{A Word of the Author (from the French Edition)}

\textit{``I would like to deeply thank the Cassini Éditions. By making
accessible these books to an audience motivated by the preparation
of competitive examinations and a genuine curiosity for learning and
reflection, they have enabled me to offer support---which I hope
fruitful---to anyone patient enough to follow along.}

\textit{I am especially grateful to André Bellaïche, with whom I had
many long and fruitful discussions during the development of this
work. André Bellaïche authored the Appendix in Section \ref{sec:Appendix:-the-Riemann}
on the Riemann integral over $\mathbb{R}^{n},$ which brought greater
rigor to Chapter \ref{chap:Random-Variables-with}---a rigor difficult
to achieve when ``continuous'' random variables are treated only
through the Riemann integral.}

\textit{Finally, I would like to thank the reviewers of this work.
Their comments helped refine the manuscript and guided it toward its
final form of this book. I hope that readers find this book both rewarding
and enjoyable.''}

\part{Introduction to Probability Theory}\label{part:Introduction-to-Probability}

\pagestyle{scrheadings}

\chapter{Random Phenomena and Probabilistic Models}\label{chap:Random-Phenomena}
\begin{quote}
\begin{flushright}\textit{Most men have, like plants, hidden properties
that chance discovers.}

\textbf{\sindex[fam]{La Rochefoucauld, François}\href{https://en.wikipedia.org/wiki/Fran\%25C3\%25A7ois_de_La_Rochefoucauld_(writer)}{La Rochefoucauld}}{\bfseries\footnote{François VI, second duke of \textbf{\href{https://en.wikipedia.org/wiki/Fran\%25C3\%25A7ois_de_La_Rochefoucauld_(writer)}{La Rochefoucauld}}
(1613-1680) is a French writer, moralist, memorialist and did a military
career. He is famous for his Memoirs, Maximes and letters.}} (1613-1680)\end{flushright}
\end{quote}
\begin{objective}{}{}

Chapter \ref{chap:Random-Phenomena} aims at introducing \textbf{Probability
Theory}.
\begin{itemize}
\item Section \ref{sec:Random-Experiment} introduces the fundamental concept
of a\textbf{ random experiment}, laying the groundwork for probability
theory. 
\item Section \ref{sec:The-Algebra-of-Events} develops the basic vocabulary
of probability, including: 
\begin{itemize}
\item \textbf{Sample space}, 
\item \textbf{Random experiment}, 
\item \textbf{Outcome}
\item And, \textbf{algebra of events}.
\end{itemize}
\item Section \ref{sec:Discrete-Probabilized-Spaces} introduces the notion
of \textbf{$\sigma-$algebra} followed definitions of:
\begin{itemize}
\item \textbf{Probabilizable space},
\item \textbf{Probability}
\item \textbf{Probabilized space}. \\
This section also explores initial properties of probability and presents
\textbf{\href{https://en.wikipedia.org/wiki/Henri_Poincar\%25C3\%25A9}{Poincaré}\sindex[fam]{Poincaré, Henri}
formula} for finite and countable unions of events.
\end{itemize}
\item Section \ref{sec:Discrete-Probabilized-Spaces} focuses on \textbf{discrete
probabilizable space} introducing:
\begin{itemize}
\item The concept of \textbf{germ of a probability} law in a discrete setting
\item Common\textbf{ discrete probability laws} key for modelling random
phenomena.
\end{itemize}
\item Section \ref{sec:Random-Variables} concludes the Chapter by defining:
\begin{itemize}
\item \textbf{Random variables}
\item The \textbf{probability law} \textbf{induced} by a random variable.
\end{itemize}
\end{itemize}
\end{objective}

\begin{figure}[t]
\begin{center}\includegraphics[width=0.4\textwidth]{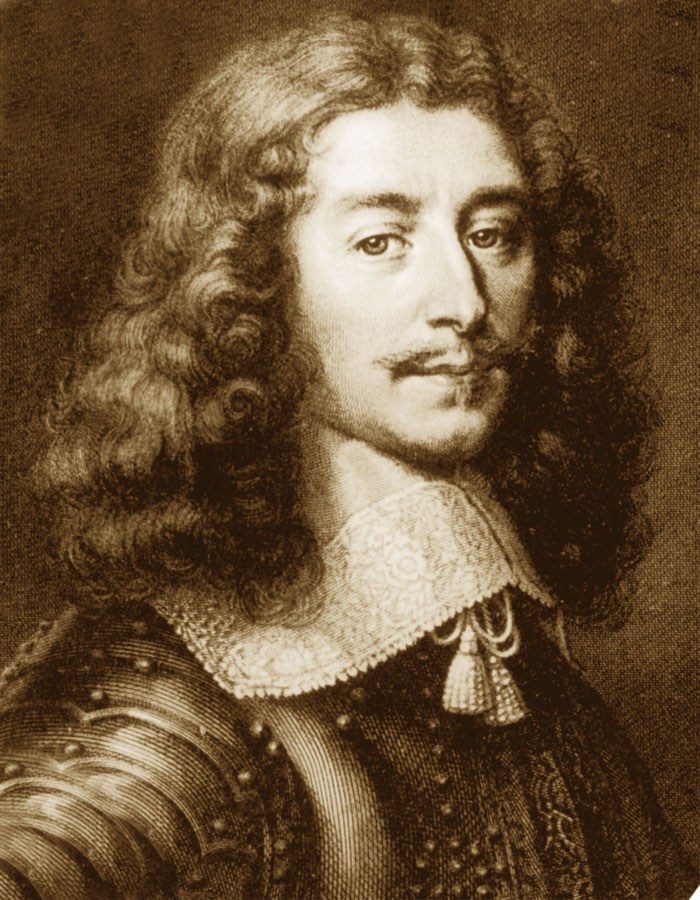}

{\tiny Par \href{https://commons.wikimedia.org/w/index.php?curid=80761264}{PdeBardon}
--- Travail personnel, CC BY-SA 4.0}\end{center}

\caption{\textbf{\protect\href{https://en.wikipedia.org/wiki/Fran\%25C3\%25A7ois_de_La_Rochefoucauld_(writer)}{François de La Rochefoucauld}}
(1613 - 1680)}\sindex[fam]{La Rochefoucauld, François}
\end{figure}

\section*{Introduction}

Probability theory is the branch of mathematics that focuses on the
study of random phenomena---phenomena influenced by chance.

Its origins trace back to the analysis of gambling, which introduced
foundational concepts such as event probability and expected gain.
Over time, probability theory has grown into a mathematical discipline
comparable to Geometry, Algebra or Analysis.

Today, probability theory finds applications across numerous fields.
In Physics and Biology, randomness and chance play a pivotal role---for
instance, in the vast diversity of traits among species---where statistical
methods are essential for analysis. In Economics and Technology understanding
and managing probabilities is key to controlling outcomes and navigating
uncertainty. Finally, probability theory underpins advancements in
machine learning and artificial intelligence, serving as a cornerstone
in these rapidly evolving domains.

\section{The Notion of a Random Experiment}\label{sec:Random-Experiment}

To express the concept of chance in mathematical terms, it is essential
to focus on specific circumstan-ces---such as those occurring repeatedly
in gambling, like rolling dice or dealing cards. These controlled
and often repeatable scenarios are referred to as random experiments.

A \textbf{random experiment}\index{random experiment}---or \textbf{\index{randomized experiment}randomized
experiment} in American English---refers to an experiment that can
be repeated---at least theoretically---under identical conditions,
but whose outcomes vary unpredictably from trial to trial.

Typical examples of random experiments include:
\begin{itemize}
\item Throwing two dice;
\item Distributing cards in a game of bridge, where the 52 cards are dealt
among four players;
\item Observing the broken pieces of an object, such as a shattered glass;
\item Examining the genetic traits of an individual based on the genetic
characteristics of their parents;
\item Observing the decay of a radioactive atomic nucleus;
\item Waiting for a bus at a specific stop after 6 p.m.
\end{itemize}
When referring to identical conditions in a random experiment, it
means: ``identical to the extent that the observer can ensure such
identical conditions.'' At first glance, one might think that chance
arises from the gap between idealized and real-world conditions. In
this view, an experiment has both predictable and unpredictable components,
with chance serving as the source of the unpredictable. This is the
common interpretation of chance.

Chance often reflects our ignorance of certain conditions within the
experiment. For instance, choosing between two boxes---one containing
a reward and the other empty---can be seen as a form of drawing lots.
In this case, what appears random to one observer may be entirely
known to another. Thus, randomness arises when the number of influencing
factors becomes so vast that it is impossible to account for them
all.

However, more intriguing aspects of chance emerge in some phenomena
governed by \textit{determinist} principles. For example, in solid
mechanics, the laws of motion are well understood, and knowledge of
initial conditions theoretically allows precise prediction of a system
future state. Yet, consider the act of throwing a die: could the thrower
predict the result even if they knew all the die characteristics,
from its geometric shape and mechanical structure to the distribution
of its mass and the uniformity of its engraved digits?

Any minor precision in the knowledge of the initial conditions is
amplified during the motion. While some systems limit this amplification,
others---such as a die rolling accross a table---exhibit extraordinarily
rapid magnification of uncertainties, making precise prediction impossible.
This is how randomness emerges in many deterministic systems\footnote{This perspective, articulated by \href{https://en.wikipedia.org/wiki/Henri_Poincar\%25C3\%25A9}{Henri Poincaré}\sindex[fam]{Poincaré, Henri},
was underappreciated during his lifetime but regained prominence with
the advent of chaos theory. A compelling discussion of this idea can
be found in Ivar Ekeland book, Au Hasard: La Chance, La Science et
Le Monde \cite{ekeland1991hasard} (English Translation: The broken
dice, and other mathematical tales of chance, \cite{ekeland1996broken}).}. There is, however, a domain where chance is not merely the result
of limited knowledge. At the atomic and subatomic levels studied by
quantum mechanics, chance exists fundamentally and cannot be eliminated,
even in principle.

\begin{figure}[t]
\begin{center}\includegraphics[width=0.4\textwidth]{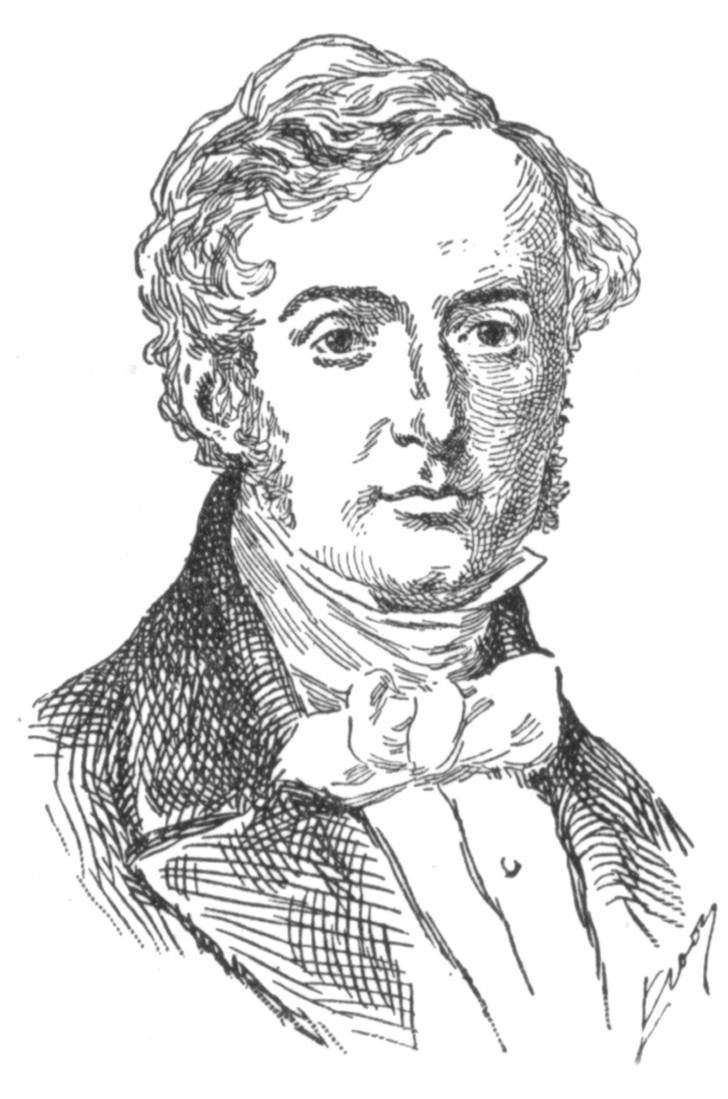}

{\tiny Credits: Nouveau Larousse Illustré XIXth century. Public Domain}\end{center}

\caption{\textbf{\protect\href{https://en.wikipedia.org/w/index.php?title=Robert_Brown_(botanist,_born_1773)}{Robert Brown}}
(1773 - 1858)}\sindex[fam]{Brown, Robert}
\end{figure}

Below is a list of random phenomena where chance plays a role and
which are the focus on ongoing scientific and technological research
(the corresponding field is mentioned in \textit{italics}):
\begin{itemize}
\item The occurrence of defective item in a manufacturing process\textit{---reliability
theory, quality control.}
\item The effects of chemical fertilizers on cereals growth---\textit{agricultural
science, statistics}. Variability is inherent to living systems, requiring
probabilistic or statistical approaches.
\item Incoming calls at a call center or the usage patterns of computers
in a network\textit{---queuing theory, network theory.}
\item Precise arrival time of buses at a station\textit{---queuing theory,
network theory};
\item Brownian motion: the random movement of particles in a a liquid, first
observed by Scottish botanist \textbf{\href{https://en.wikipedia.org/w/index.php?title=Robert_Brown_(botanist,_born_1773)}{Robert Brown}}\sindex[fam]{Brown, Robert}
in 1827. Its mathematical study has sparked ongoing research with
applications in fields like telecommunications, medicine, insurance,
and financial market theory---\textit{stochastic processes}.
\item Uncertainty in the trajectory of a spacecraft: various perturbations
cause deviation from the spacecraft theoretical path, introducing
randomness---\textit{filtering theory, stochastic control}.
\end{itemize}
To study a phenomenon influenced by chance---as defined earlier---a
practitioner seeks to isolate a random experiment and construct a
probabilistic model. This model enables them to give some predictions
by calculating the probabilities of specific outcomes or events. However,
any conclusions drawn are only meaningful within the framework of
the model! 

To extend their analysis, the practitioner must validate the model,
by comparing its predictions with data collected during the experiment.
This validation process is the domain of \textit{statistics}. If the
model passes this test successfully, it can then be used to make reliable
predictions or provide explanations for a phenomenon under investigation.

The first step in mathematically modelling a random experiment is
to define the set of all possible outcomes---also called realizations
or results---of the experiment. This set is typically denoted by
the letter $\Omega$ and is referred to as the \textbf{sample space}\mindex{sample!space}---also
known as the \textbf{sample description space}\mindex{sample!description space},
\textbf{possibility space\index{possibility space}}, \textbf{outcome
space}\mindex{outcome!space}, or \textbf{universe}\index{universe}.\footnote{Translator's note (Tr.N. for short): While the term universe is less
common in English literature, it is preferred here due to its clarity,
self-contained nature, and alignment with the author's original preference
in French.}

For example:
\begin{itemize}
\item When two standard dice are thrown, the possible results are pair of
integers between 1 and 6. The sample space can be
\[
\Omega=\left\{ 1,2,...,6\right\} \times\left\{ 1,2,...,6\right\} ,
\]
which contains 36 elements. 
\item In the card game named bridge, where 52 cards are distributed among
4 players, the sample space $\Omega$ can consist of all possible
ways to distribute the 52 cards among the four players. This set has
a finite number of elements, specifically 
\[
\left(\begin{array}{c}
52\\
13\,13\,13\,13
\end{array}\right)=53'644'737'765'488'792'839'237'440'000.
\]
\item For the radioactive decay of a nucleus, if the focus is on its lifetime,
the sample space can be $\left]0\,;\,+\infty\right[,$ the set of
strictly positive real numbers.
\end{itemize}
A same experiment can be described in multiple ways. For instance,
if two players throw two dice sequentially with the objective to have
the higher sum, the sample space could be represented as either 
\[
\left\{ 1,2,...,6\right\} ^{4}
\]
 or 
\[
\left\{ 2,3,...,12\right\} \times\left\{ 2,3,...,12\right\} .
\]

\textbf{The choice of representation depends on the objective of the
study}, which influences how realizations are defined.

Finally, we observe that any random experiment can be translated into
the selection of a random point $\omega$ within a carefully chosen
universe $\Omega.$

\section{The Algebra of Events}\label{sec:The-Algebra-of-Events}

A \textbf{random event}\index{random event} is an occurrence linked
to a \textbf{random experiment}\index{random experiment}. It may
or may not happen, and its realization---or non realization----depends
entirely on the \textbf{outcome\index{outcome}} $\omega$ of the
experiment. Conceptually, an event can be thought of as a property
of the result $\omega,$ that can either hold true or not, or also
as a statement about $\omega$ that may be true or false.

Mathematically, a \textbf{random event} is represented by specifying
the \textbf{set of outcomes} $\omega$ from $\Omega$ for which this
event occurs.

For example, consider the case of rolling two dice. The event ``the
sum of the numbers on the top faces is less than or equal to 4''
corresponds to the subset of outcomes
\[
\left\{ \left(1,1\right),\left(1,2\right),\left(2,1\right),\left(1,3\right),\left(2,2\right),\left(3,1\right)\right\} 
\]
from the universe
\[
\Omega=\left\{ 1,2,...,6\right\} \times\left\{ 1,2,...,6\right\} .
\]

In probability theory, it is standard to identify an event directly
with the subset of $\Omega$ it represents. This subset is often denoted
by the same symbol as the event itself. This identification is fundamental
in probability theory as the concepts and logical operations defined
for events correspond directly to the notions and set operations in
set theory. This correspondence is summarized in Table \ref{Tab: probabilist set terminology correspondance}
and translated from \cite{neveu1964bases}.

\begin{table}
\begin{center}%
\begin{tabular}{|c|c|}
\hline 
\textbf{Probabilist terminology} &
\textbf{Set terminology}\tabularnewline
\hline 
Universe $\Omega$ &
Set $\Omega$\tabularnewline
\hline 
Possible outcome $\omega$ &
$\omega,$ element of $\Omega$\tabularnewline
\hline 
Event $A$ &
$A,$ subset of $\Omega$\tabularnewline
\hline 
$A$ is realised by the outcome $\omega$  &
$\omega\in A$\tabularnewline
\hline 
The event $A$ implies the event $B$ &
$A\subseteq B$\tabularnewline
\hline 
Event $A$ or event $B$ &
$A\cup B$ ($A$ union $B$)\tabularnewline
\hline 
Event $A$ and event $B$ &
$A\cap B$ ($A$ inter $B$)\tabularnewline
\hline 
Event $A$ does not happen &
$A^{c}$ (complement of $A$)\tabularnewline
\hline 
Impossible event &
$\emptyset$\tabularnewline
\hline 
Sure event &
$\Omega$\tabularnewline
\hline 
Incompatible events $A$ and $B$ &
$A\cap B=\emptyset$\tabularnewline
\hline 
\end{tabular}\end{center}

\caption{Correspondence between probabilist terminology and the set terminology,
translated from \cite{neveu1964bases}}
\label{Tab: probabilist set terminology correspondance}
\end{table}

It is important to note that not every subset of $\Omega$ is necessarily
associated with an event. This distinction is explored in greater
depth in Part \ref{part:Deepening-Probability-Theory}. For now, consider
the case of rolling two identical dice. In this scenario, the outcome
$\left(1,2\right)$ is indistinguishable from $\left(2,1\right)$,
making it impossible to observe $\omega=\left(1,2\right)$ as a distinct
result. Consequently, the subset $\left\{ \left(1,2\right)\right\} $
is not considered an event.

In general, we are interested in a collection $\mathcal{L}$ of properties
of the outcome, which it is reasonnable to assume stable under the
following operations:
\begin{itemize}
\item \textbf{Negation:} if $L\in\mathcal{L},$ then its negation $\lnot L$---the
event that occurs when $L$ does not occur---must also belong to
$\mathcal{L}.$
\item \textbf{Conjunction:} if $L_{1}\in\mathcal{L}$ and $L_{2}\in\mathcal{L},$
then the event where both $L_{1}$ and $L_{2}$ occur, denoted as
$L_{1}\land L_{2}$ must also be in $\mathcal{L}.$
\item \textbf{Disjunction:} if $L_{1}\in\mathcal{L}$ and $L_{2}\in\mathcal{L},$
then the event where at least one of $L_{1}$ or (non exclusive) $L_{2}$
occurs, denoted as $L_{1}\lor L_{2}$ must also be in $\mathcal{L}.$
\end{itemize}
The stability---or closeness from a set point of view---of those
operations ensure that the collection $\mathcal{L}$ satisfies the
basic logical operations necessary for a coherent probabilistic framework.

If $\mathcal{A}$ denotes the family of events---mathematically corresponding
to each of the properties of the collection $\mathcal{L}$---this
is equivalent to suppose that $\mathcal{A}$ is stable---closed---under
the complement, the intersection and the union. A family of subsets
satisfying these properties is called an algebra. Thus, $\mathcal{A}$
is referred as the \textbf{algebra of events}\mindex{algebra!events}.
For deeper mathematical reasons, such as handling limits, we require
$\mathcal{A}$ to be stable under countable union. This leads to a
more general structure known as a \textbf{$\sigma-$algebra}\mindex{sigma-algebra @ \salg},
which will be discussed in the next subsection.

\section{$\sigma-$Algebra and Probability Axioms. First Properties.}\label{sec:Algebra-and-Probability}

We now introduce the mathematical framework necessary for probability
theory. Before defining what a probability is, we first establish
the definitions of $\sigma-$algebras and probabilizable spaces, as
they serve as fundamental prerequisites to probability definition.

\begin{definition}{$\sigma-$ algebra}{sigmaalg} A family $\mathcal{A}$
of subsets of a set $\Omega$ is called a \textbf{$\sigma-$algebra}\mindex{sigma-algebra @ \salg}---or
a \mindex{sigma-field @ $\sigma-$field}\textbf{$\sigma-$field}
or more rarely in English a \textbf{tribe}\index{tribe}---on $\Omega,$
if it satisfies the following three axioms:

(i) \textbf{Containment of the Universal Set} 

\boxeq{
\[
\Omega\in\mathcal{A}.
\]
}

(ii) \textbf{Stability under Complementation} 

\boxeq{
\[
\text{If}\,A\in\mathcal{A},\,\text{then}\,A^{c}\in\mathcal{A}.
\]
}

(iii) \textbf{Stability under Countable Unions} 

For every sequence $\left(A_{i}\right)_{i\in\mathbb{N}}$ of elements
of $\mathcal{A},$ 

\boxeq{
\[
\underset{i\in\mathbb{N}}{\bigcup}A_{i}\in\mathcal{A}.
\]
}

\end{definition}

\begin{example}{}{}1. The power set $\mathcal{P}\left(\Omega\right),$
consisting of all subsets of $\Omega,$ is a $\sigma-$algebra on
$\Omega.$

2. For every $A\in\mathcal{P}\left(\Omega\right),$ the family $\left\{ A,A^{c},\Omega,\emptyset\right\} $
forms a $\sigma-$algebra on $\Omega,$ called the \textbf{$\sigma-$algebra
generated by the event} $A.$\mindex{sigma-algebra @ \salg ! generated by an event}

3. The $\sigma-$algebra $\left\{ \emptyset,\Omega\right\} $ is called
the \mindex{sigma-algebra @ \salg! trivial}\textbf{trivial $\sigma-$algebra}
on $\Omega.$ It is the $\sigma-$algebra with the smallest possible
cardinality that can be generated.

\end{example}

\begin{proposition}{Immediate Properties of a $\sigma-$ algebra}{immedpropsigalg}Let
$\mathcal{A}$ be a $\sigma-$algebra on $\Omega.$ The following
properties hold:

(i) \textbf{Containment of the Empty Set}\boxeq{
\[
\emptyset\in\mathcal{A}.
\]
}

(ii) \textbf{Stability under Finite Unions}

For every finite sequence $\left(A_{i}\right)_{1\leqslant i\leqslant n}$
of elements in $\mathcal{A},$\boxeq{
\[
\underset{i=1}{\bigcup^{n}}A_{i}\in\mathcal{A}.
\]
}

(iii) \textbf{Stability under Finite Intersections} 

For every finite sequence $\left(A_{i}\right)_{1\leqslant i\leqslant n}$
of elements in $\mathcal{A},$ we have\boxeq{
\[
\underset{i=1}{\bigcap^{n}}A_{i}\in\mathcal{A}.
\]
}

(iv) \textbf{Stability under Countable Intersections}

For every sequence $\left(A_{i}\right)_{n\in\mathbb{N}}$ of elements
in $\mathcal{A},$ we have\boxeq{
\[
\underset{n\in\mathbb{N}}{\bigcap}A_{i}\in\mathcal{A}.
\]
}

(v) \textbf{Stability under Set Difference}

For every $A,B\in\mathcal{A}$, we have\footnotemark\boxeq{
\[
A\backslash B\in\mathcal{A}.
\]
}

\end{proposition}

\footnotetext{The notation $A\backslash B$ denotes the set of elements
that belong to $A$ but not to $B.$ Formally, we have 
\[
A\backslash B=A\cap B^{c}.
\]
}

\begin{proof}{}{}

Let $\mathcal{A}$ be a $\sigma-$algebra on $\Omega.$ We now prove
the stated properties.

(i) \textbf{Containment of the Empty Set}

Since $\Omega\in\mathcal{A}$ and $\mathcal{A}$ is a $\sigma-$algebra,
it follows from the stability under complementation that
\[
\Omega^{c}=\emptyset\in\mathcal{A}.
\]

(ii) \textbf{Stability under Finite Union}

Consider a finite family $\left(A_{i}\right)_{i\in\left\llbracket 1,n\right\rrbracket }$\footnotemark
of elements of $\mathcal{A}$.

Define the family $\left(B_{i}\right)_{i\in\mathbb{N}}$, such that:
\begin{itemize}
\item For every $i\in\left\llbracket 1,n\right\rrbracket ,$ let $B_{i}=A_{i}.$
\item For every $i\geqslant n+1,$ let $B_{i}=\emptyset.$
\end{itemize}
Since for every $i\in\mathbb{N},$ $B_{i}\in\mathcal{A},$ the third
axiom of a $\sigma-$algebra ensures that
\[
\underset{i\in\mathbb{N}}{\bigcup}B_{i}\in\mathcal{A}
\]

Moreover, by construction
\[
\underset{i\in\mathbb{N}}{\bigcup}B_{i}=\underset{i=1}{\bigcup^{n}}A_{i}.
\]
Thus, the result follows.

(iii) \textbf{Stability under Finite Intersections}

Consider a finite famility of elements $\left(A_{i}\right)_{i\in\left\llbracket 1,n\right\rrbracket }$
in $\mathcal{A}$. A classical approach is to use the complements,
\[
\left(\underset{i=1}{\bigcap^{n}}A_{i}\right)^{c}=\underset{i=1}{\bigcup^{n}}A^{c}_{i}.
\]

Since each $A_{i}\in\mathcal{A},$ the second axiom of a $\sigma-$algebra
ensures that $A^{c}_{i}\in\mathcal{A}$ for every $i\in\left\llbracket 1,n\right\rrbracket .$

By the stability under finite union, using (ii),
\[
\underset{i=1}{\bigcup^{n}}A^{c}_{i}\in\mathcal{A}
\]
Taking the complement again, and using the second axiom of a $\sigma-$algebra
once more, we conclude that
\[
\left[\left(\underset{i=1}{\bigcap^{n}}A_{i}\right)^{c}\right]^{c}=\underset{i=1}{\bigcap^{n}}A_{i}\in\mathcal{A}.
\]

(iv) \textbf{Stability under Countable Intersections}

The same reasoning than in (iii) applies by considering an infinite
sequence and taking complements.

(v) \textbf{Stability under Set Difference}

Since $A\in\mathcal{A}$ and $B^{c}\in\mathcal{A}$---by stability
under complementation---if follows that
\[
A\backslash B=A\cap B^{c}\in\mathcal{A}.
\]

\end{proof}

\footnotetext{The notation $\left\llbracket a,b\right\rrbracket $
where $a,b\in\mathbb{N}$ and $a\leqslant b$ designates the set $\left\{ i:\,i\in\mathbb{N}\land a\leqslant i\leqslant b\right\} .$}

\begin{definition}{Probabilizable Space}{probabilizablesp}

Let $\Omega$ be a set and $\mathcal{A}$ a $\sigma-$algebra on $\Omega.$

The pair $\left(\Omega,\mathcal{A}\right)$ is called a \index{probabilizable space}\textbf{probabilizable
space}\footnotemark.

When a probabilizable space $\left(\Omega,\mathcal{A}\right)$ is
given, the $\sigma-$algebra $\mathcal{A}$ is referred to as the
\mindex{sigma-algebra @ \salg ! event}\mindex{event!sigma-algebra @ \salg}\textbf{event
$\sigma-$algebra}.

\end{definition}

\footnotetext{In measure theory, a probabilizable space is nothing
more than a \textbf{measurable space}\index{measurable space}.}

The modelling of a random phenomenon and its associated properties
or events begins with the selection of a \textbf{probabilizable space,}
which encapsulates the set of all possible outcomes and the collection
of properties or events under study. This space provides a \textit{qualitative}
description of the phenomenon and the properties being analyzed. 

We now introduce the concept of \textbf{probability\index{probability}}
on a probabilizable space, which will allow for a \textit{quantitative}
assessment of the phenomenon under consideration. For now, we will
keep the definition at an intuitive level without delving into more
precise technical details.

\begin{figure}[t]
\begin{center}\includegraphics[width=0.4\textwidth]{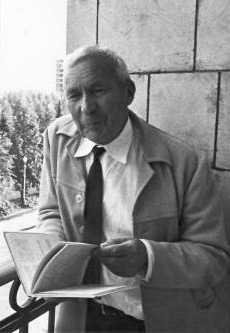}

{\tiny Credits: \href{https://opc.mfo.de/detail?photoID=7493}{Konrad Jacobs}
under CC-BY-SA}\end{center}

\caption{\textbf{\protect\href{https://en.wikipedia.org/wiki/Andrey_Kolmogorov}{Andrey Kolmogorov}}
(1903 - 1987)}\sindex[fam]{Kolmogorov, Andrey}
\end{figure}

\begin{definition}{Probability. Probabilized Space}{probabilizedspace}Let
$\left(\Omega,\mathcal{A}\right)$ be a \textbf{probabilizable space}.
A \textbf{probability\index{probability}} $P$ on this space is a
function $P:\mathcal{A}\rightarrow\mathbb{R}^{+}$ that satisfies
the following two axioms:

(i) \textbf{Normalization}

\boxeq{
\[
P\left(\Omega\right)=1.
\]
}

(ii) \textbf{$\boldsymbol{\sigma-}$additivity}

The function $P$ is \textbf{$\sigma-$additive}\mindex{sigma-additive @ $\sigma-$additive}:
for every sequence $\left(A_{i}\right)_{i\in\mathbb{N}}$ of pairwise
disjoint sets in $\mathcal{A},$

\boxeq{
\[
P\left(\underset{i\in\mathbb{N}}{\biguplus}A_{i}\right)=\sum\limits^{+\infty}_{i=0}P\left(A_{i}\right).
\]
}

The triple $\left(\Omega,\mathcal{A},P\right)$ is called a \textbf{\mindex{probability!space}probability
space}---or sometimes a \textbf{probabilized space}\index{probabilized space}\footnotemark.

\end{definition}

\footnotetext{\textbf{\href{https://en.wikipedia.org/wiki/Andrey_Kolmogorov}{Andrey Kolmogorov}}\sindex[fam]{Kolmogorov, Andrey}
(1903 -1987), a Russian mathematician, introduced this axiomatic approach
in 1929 in his work \textquotedbl General Theory of Measure and the
Calculus of Probability\textquotedbl , where he began formalizing
probability using measure theory. He later expanded and fully developed
this framework in 1933, in his seminal monograph \textquotedbl Grundbegriffe
der Wahrscheinlichkeitsrechnung\textquotedbl{} (Foundations of the
Theory of Probability). This 1933 publication established the modern
axioms of probability and rigorously connected probability theory
with \textbf{measure theory}. He worked also on Markov and stationary
processes. Andreï Kolmogorov was designated Professor at Moscow University
in 1931 and in 1933 was already directing its Mathematics Institute.
\label{fn:Kolmogorov}}

The term \index{law of probability}\mindex{probability ! law}\textbf{law
of probability} is often used interchangeably with probability in
some contexts.

The following proposition presents the elementary properties of probabilized
spaces.

\begin{proposition}{First Properties of Probabilities}{probafirstprop}

Let $\left(\Omega,\mathcal{A},P\right)$ be a probabilized space. 

The following properties hold:

(i) \textbf{Probability of the Empty Set}

\boxeq{
\[
P\left(\emptyset\right)=0.
\]
}

(ii) \textbf{Finite Additivity}

\boxeq{For every finite sequence $\left(A_{i}\right)_{1\leqslant i\leqslant n}$
of pairwise disjoint sets of $\mathcal{A}$\footnotemark\footnotemark,

\[
P\left(\biguplus^{n}_{i=1}A_{i}\right)=\sum\limits^{n}_{i=1}P\left(A_{i}\right).
\]
}

(iii) \textbf{Monotonicity and Additivity for Difference}

\boxeq{

If $A,B\in\mathcal{A}$ and $A\subset B,$ then
\[
P\left(A\right)\leqslant P\left(B\right)
\]

and
\begin{equation}
P\left(B\backslash A\right)=P\left(B\right)-P\left(A\right).
\end{equation}
}

(iv) \textbf{Probability Bounds}

\boxeq{For every $A\in\mathcal{A},$ 
\[
P\left(A\right)\in\left[0;1\right].
\]
}

(v) \textbf{Complement Rule}

\boxeq{For every $A\in\mathcal{A},$
\begin{equation}
P\left(A^{c}\right)=1-P\left(A\right).
\end{equation}
}

(vi) \textbf{Continuity from Below---for nondecreasing Sequences}

For every nondecreasing sequence $\left(A_{i}\right)_{i\in\mathbb{N}}$
in $\mathcal{A}$---i.e. such that for every $i\in\mathbb{N},$ $A_{i}\subset A_{i+1}$---we
have\footnotemark\boxeq{
\[
P\left(\bigcup\limits_{i\in\mathbb{N}}A_{i}\right)=\lim\limits_{i\rightarrow+\infty}\uparrow P\left(A_{i}\right).
\]
}

(vii) \textbf{Continuity from Above---for nonincreasing Sequences}

For every nonincreasing sequence $\left(A_{i}\right)_{i\in\mathbb{N}}$
in $\mathcal{A}$---i.e. such that for every $i\in\mathbb{N},$ $A_{i}A_{i+1}$---we
have\boxeq{
\[
P\left(\bigcap_{i\in\mathbb{N}}A_{i}\right)=\lim\limits_{i\rightarrow+\infty}\downarrow P\left(A_{i}\right).
\]
}

(viii) \textbf{Inclusion-Exclusion Principle for Two Sets}

\boxeq{For every $A,B\in\mathcal{A},$

\begin{equation}
P\left(A\cup B\right)=P\left(A\right)+P\left(B\right)-P\left(A\cap B\right).\label{eq:P_AunionB}
\end{equation}
}

\end{proposition}

\addtocounter{footnote}{-2}%
\footnotetext{We say that $P$ is a set function \index{finitely additive}\mindex{additive ! finitely}\textbf{finitely
additive}}

\stepcounter{footnote}
\footnotetext{$A\uplus B$ denotes the union of $A$ and $B,$ while mentioning
at the same time that the sets are disjoint. A similar notation is
used for sequence of pairwise disjoint sets. Some authors use $A+B$
or $A\sqcup B.$}

\stepcounter{footnote}
\footnotetext{We denote $\lim_{n\to+\infty}\uparrow u_{n}$ to designate the limit
$\lim_{n\to+\infty}u_{n}$ while recalling that the sequence $\left(u_{n}\right)$
is non-{}-decreasing. Similarly, mutatis mutandis, for $\lim_{n\to+\infty}\downarrow u_{n}.$}

\begin{proof}{}{}

(i) \textbf{Probability of the Empty Set}

Consider the family $\left(B_{i}\right)_{i\in\mathbb{N}}$ defined
as follows:
\begin{itemize}
\item $B_{0}=\Omega$
\item $B_{i}=\emptyset$ for every $i>0.$ 
\end{itemize}
Since this family consists of pairwise disjoint sets and satisfies
$\biguplus_{i\in\mathbb{N}}B_{i}=\Omega,$ the $\sigma-$additivity
of $P$ implies
\[
P\left(\Omega\right)=1=\sum\limits_{i\in\mathbb{N}}P\left(B_{i}\right)=P\left(B_{0}\right)+\sum\limits_{i>0}P\left(B_{i}\right).
\]
Since $P\left(B_{0}\right)=P\left(\Omega\right)=1,$ it follows that
$\sum\limits_{i>0}P\left(B_{i}\right)=0,$ and hence $P\left(\emptyset\right)=0,$
as $P$ has nonnegative values.

(ii)\textbf{ Finite Additivity}

Let $\left(A_{i}\right)_{0\leqslant i\leqslant n}$ be a finite family
of pairwise disjoint elements of $\mathcal{A}.$

Define the family $\left(B_{i}\right)_{i\in\mathbb{N}}$ by setting
\begin{itemize}
\item $B_{i}=A_{i},$ for $0\leqslant i\leqslant n;$ 
\item $B_{i}=\emptyset,$ for every $i>n.$ 
\end{itemize}
Hence, $\biguplus\limits_{i\in\mathbb{N}}B_{i}=\biguplus\limits^{n}_{i=0}A_{i}$
and the family $\left(B_{i}\right)_{i\in\mathbb{N}}$ has its elements
in $\mathcal{A}$ pairwise disjoint. Then, since $P$ is $\sigma-$additive,
\[
P\left(\biguplus^{n}_{i=0}A_{i}\right)=P\left(\biguplus_{i\in\mathbb{N}}B_{i}\right)=\sum\limits^{+\infty}_{i=0}P\left(B_{i}\right)
\]

Using Property (i), the sum simplifies to $\sum\limits^{n}_{i=0}P\left(A_{i}\right),$
proving the result.

(iii)\textbf{ Monotonicity and Additivity for Difference}

Since $B$ can be decomposed as $B=A\uplus\left(B\backslash A\right)$\footnotemark,
the $\sigma-$additivity of $P$ gives 
\[
P\left(B\right)=P\left(A\right)+P\left(B\backslash A\right).
\]
 that is $P\left(B\backslash A\right)=P\left(B\right)-P\left(A\right).$ 

Since $P$ is nonnegative, it follows that $P\left(A\right)\leqslant P\left(B\right).$

(iv) \textbf{Probability Bounds}

Since, for every $A\in\mathcal{A},$ 
\[
\emptyset\subset A\subset\Omega
\]
applying (iii) to these inclusions gives 
\[
0=P\left(\emptyset\right)\leqslant P\left(A\right)\leqslant P\left(\Omega\right)=1.
\]

(v) \textbf{Complementation}

For every $A\in\mathcal{A},$
\[
\Omega=A\uplus A^{c}.
\]
 The $\sigma-$additivity of $P$ implies 
\[
P\left(A\right)+P\left(A^{c}\right)=P\left(\Omega\right)=1.
\]
Rearranging gives $P\left(A^{c}\right)=1-P\left(A\right).$

(vi) \textbf{Continuity from Below---for nondecreasing Sequences:}

Let $\left(A_{i}\right)_{i\in\mathbb{N}}$ be a family of nondecreasing
subset in $\mathcal{A}.$ 

For every $p\in\mathbb{N},$ decomposing $A=\bigcup\limits_{i\in\mathbb{N}}A_{i}$
as 
\[
A=A_{p}\uplus\left[\biguplus\limits_{i\geqslant p}\left(A_{i+1}\backslash A_{i}\right)\right].
\]

the $\sigma-$additivity of $P$ gives the relation
\[
P\left(A\right)=P\left(A_{p}\right)+\sum\limits_{i\geqslant p}P\left(A_{i+1}\backslash A_{i}\right).
\]

It means that, taking $p=0,$ $\sum_{i\geqslant0}P\left(A_{i+1}\backslash A_{i}\right)$
is finite and equal to $P\left(A\right)-P\left(A_{0}\right).$ 

Hence the series, of general term $P\left(A_{i+1}\backslash A_{i}\right)$
for $i\in\mathbb{N}$ converges and its remainder $\sum_{i\geqslant p}P\left(A_{i+1}\backslash A_{i}\right)$
tends to zero as $p$ tends to $+\infty.$ 

Thus, $P\left(A_{p}\right)$ tends to $P\left(A\right)$ when $p$
tends to $+\infty,$ i.e. 
\[
\lim\limits_{p\rightarrow+\infty}P\left(A_{p}\right)=P\left(\bigcup\limits_{i\in\mathbb{N}}A_{i}\right).
\]

(vii) \textbf{Continuity from Above---for Decreasing Sequences}

Let $\left(A_{i}\right)_{i\in\mathbb{N}}$ being a decreasing sequence
of elements of $\mathcal{A}.$ We construct, for every $i\in\mathbb{N},$
$B_{i}=A^{c}_{i}.$ By applying (vi), it holds
\[
\lim\limits_{p\rightarrow+\infty}P\left(B_{p}\right)=P\left(\bigcup\limits_{i\in\mathbb{N}}B_{i}\right)
\]

Taking complements and using (v), we conclude the proof.

(viii) \textbf{Inclusion-Exclusion Principle for Two Sets}

Using the decomposition
\[
A\cup B=A\uplus\left(B\backslash\left(A\cap B\right)\right),
\]
and applying the $\sigma-$additivity for disjoint union, we obtain
\[
P\left(A\cup B\right)=P\left(A\right)+P\left(B\backslash\left(A\cap B\right)\right).
\]
Since (iii) gives $P\left(B\backslash\left(A\cap B\right)\right)=P\left(B\right)-P\left(A\cap B\right),$
the result follows.

\end{proof}

\footnotetext{$\uplus$ designates the union of two sets, considering
that they are disjoint.}

The formula $\refpar{eq:P_AunionB}$ is generalized to $n$ subsets
of $\mathcal{A},$ providing an inclusion-exclusion principle for
probabilities, systematically accounting for overcounting in unions.

\begin{figure}[t]
\begin{center}\includegraphics[width=0.4\textwidth]{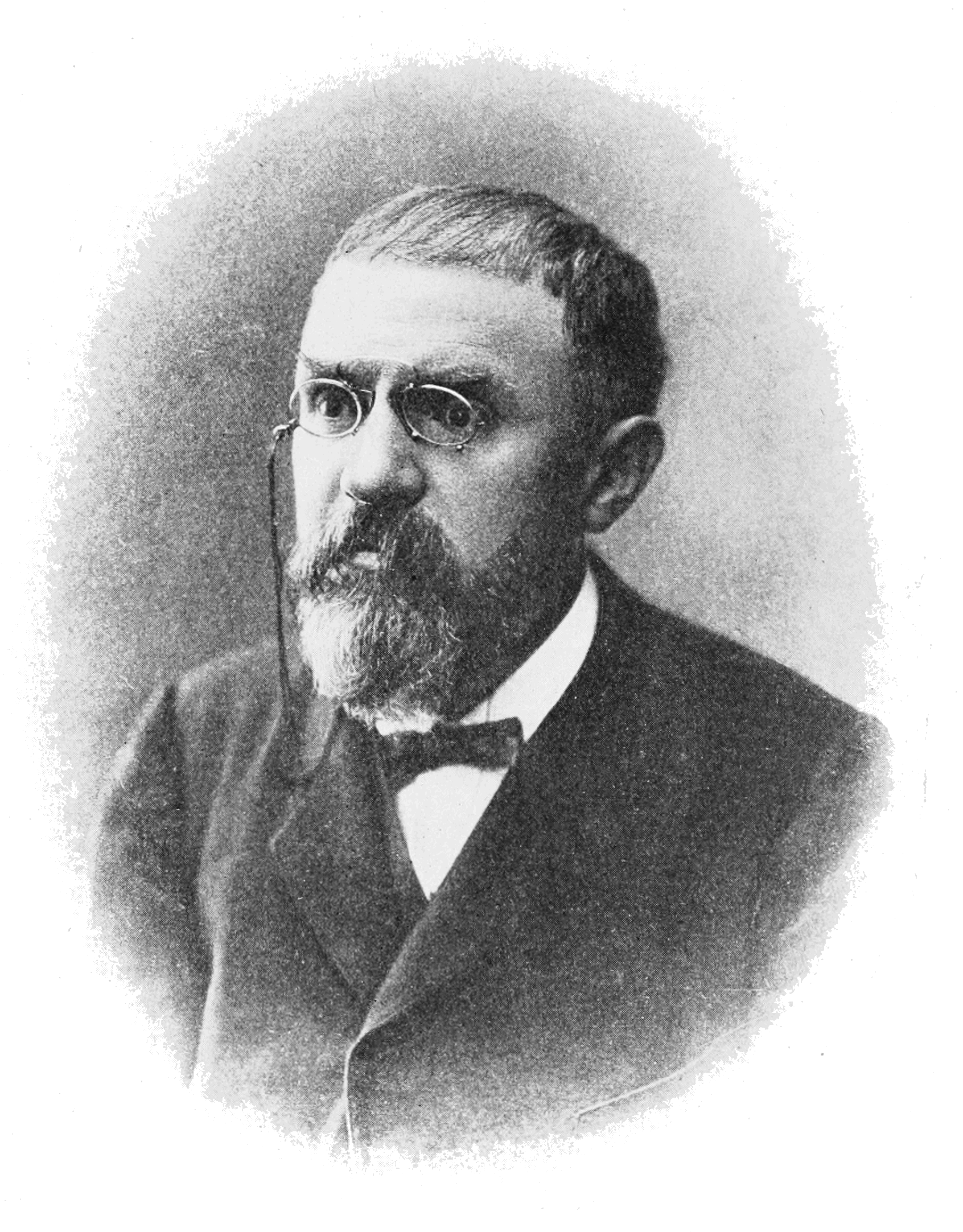}

{\tiny Credits: Popular Science Monthly Volume 82, Public Domain}\end{center}

\caption{\textbf{\protect\href{https://en.wikipedia.org/w/index.php?title=Henri_Poincar\%C3\%A9}{Henri Poincaré}}
(1854 - 1912)}
\end{figure}

\begin{proposition}{Poincaré Formula}{poincareform} For every integer
$n\geqslant2,$ and every finite sequence $\left(A_{i}\right)_{1\leqslant i\leqslant n}$
of elements of $\mathcal{A},$

\boxeq{
\begin{multline*}
P\left(\bigcup\limits^{n}_{i=1}A_{i}\right)=\sum\limits^{n}_{i=1}P\left(A_{i}\right)-\sum\limits_{1\leqslant i<j\leqslant n}P\left(A_{i}\cap A_{j}\right)+\sum\limits_{1\leqslant i<j<k\leqslant n}P\left(A_{i}\cap A_{j}\cap A_{k}\right)\\
-\dots+\left(-1\right)^{n-1}P\left(A_{1}\cap A_{2}\cap\dots\cap A_{n}\right).
\end{multline*}
}

This formula is named the \textbf{\index{Poincaré formula}Poincaré
formula}\footnotemark.

More generally, let $I$ be a finite set of cardinal at least 2 and
let $\left(A_{i}\right)_{i\in I}$ be a family of elements of $\mathcal{A}.$
Then denoting $\left|J\right|$ the cardinal of the set $J\subset I,$
the Poincaré\sindex[fam]{Poincaré, Henri}\index{Poincaré formula}
formula can be expressed as

\boxeq{
\begin{multline*}
P\left(\bigcup\limits_{i\in I}A_{i}\right)=\sum\limits_{\substack{J\subset I\\
\left|J\right|=1
}
}P\left(\bigcap\limits_{j\in J}A_{j}\right)-\sum\limits_{\substack{J\subset I\\
\left|J\right|=2
}
}P\left(\bigcap\limits_{j\in J}A_{j}\right)+\sum\limits_{\substack{J\subset I\\
\left|J\right|=3
}
}P\left(\bigcap\limits_{j\in J}A_{j}\right)\\
-\dots+\left(-1\right)^{\left|I\right|-1}P\left(\bigcap\limits_{j\in I}A_{j}\right).
\end{multline*}
}

\end{proposition}

\footnotetext{\textbf{\href{https://en.wikipedia.org/w/index.php?title=Henri_Poincar\%C3\%A9}{Henri Poincaré}}\sindex[fam]{Poincaré, Henri}
(1854 - 1912) was a French mathematician, and he is known as The Last
Universalist since he excelled in all fields of mathematics, and was
certainly one of the last to be able to do so. He was also theoretical
physicist, engineer and philosopher of science. He particularly studied
the problem of the three bodies, by introducing key concepts that
led to the theory of chaos.}

\begin{proof}{}{}

This proposition is proved by induction on $n.$

\textbf{Initialization step}

From the equation $\refpar{eq:P_AunionB}$, this formula is true at
the order 2. 

\textbf{Induction step}

We suppose the formula to be true at the order $n.$ We seek to prove
it for $n+1.$

Observing that
\[
\bigcup\limits^{n+1}_{i=1}A_{i}=\left(\bigcup\limits^{n}_{i=1}A_{i}\right)\cup A_{n+1}.
\]
and since each term belongs to $\mathcal{A},$ applying the formula
at the order 2 gives
\[
P\left(\bigcup\limits^{n+1}_{i=1}A_{i}\right)=P\left(\bigcup\limits^{n}_{i=1}A_{i}\right)+P\left(A_{n+1}\right)-P\left(\left(\bigcup\limits^{n}_{i=1}A_{i}\right)\cap A_{n+1}\right)
\]

Since
\[
\left(\bigcup\limits^{n}_{i=1}A_{i}\right)\cap A_{n+1}=\bigcup\limits^{n}_{i=1}A_{i}\cap A_{n+1},
\]
and noting that each $A_{i}\cap A_{n+1}$ belongs to $\mathcal{A},$
we apply the induction hypothesis twice to the right-hand side member
and rearrange accordingly the terms. This completes the proof.

\end{proof}

\section{Discrete Probabilized Spaces}\label{sec:Discrete-Probabilized-Spaces}

\subsection{Definition}

The notion of countable set plays a fundamental role in this book.
For clarity, we begin by summarizing key definitions and well-known
results concerning countable sets.

A \textbf{countable set}\index{countable set}\mindex{set ! countable},
as defined in this book, is a set that can be placed in bijective
correspondence with a subset of $\mathbb{N}.$ Consequently, finite
sets are also countable\footnote{In some texts, a set is called countable only if it can be placed
in bijective correspondence with the entirety of $\mathbb{N},$ i.e.
sets that are countably infinite sets. In such cases, when we talk
of countable sets, those texts consider finite sets or countably infinite
sets.}.

A subset of a countable set is countable. The image of a countable
set under any function remains countable. Moreover, if $A_{1},A_{2},\dots,A_{n}$
are countable sets, then both their union $A_{1}\cup A_{2}\cup...\cup A_{n}$
and their Cartesian product $A_{1}\times A_{2}\times\dots\times A_{n}$
are countable. Similarly, if $\left(A_{i}\right)_{i\in\mathbb{N}}$
is an infinite sequence of countable sets, then the union $\bigcup_{i\in\mathbb{N}}A_{i}$
is countable.

The sets $\mathbb{N},$ $\mathbb{Z},$ $\mathbb{Q}$ are countable,
while $\mathbb{R}$ and non-trivial intervals of $\mathbb{R}$ are
uncountable.

Now consider a \textbf{probabilizable space} $\left(\Omega,\mathcal{A}\right)$
where $\Omega$ is countable, and assume that for every $\omega\in\Omega,$
the singleton $\left\{ \omega\right\} $ belongs to $\mathcal{A}.$
It follows that every subset $A$ of $\Omega$ belongs to the $\sigma-$algebra
$\mathcal{A}.$ Indeed, we can express 

\boxeq{
\[
A=\biguplus\limits_{\omega\in A}\left\{ \omega\right\} 
\]
}as a countable union of elements of the $\sigma-$algebra $\mathcal{A}$.
Consequently, $\mathcal{A}=\mathcal{P}\left(\Omega\right).$ From
this observation, we can enounce the following definition.

\begin{definition}{Discrete Probabilizable Space}{discretprobasp}

A \textbf{discrete probabilizable space\mindex{discrete!probabilizable space}\mindex{probabilizable space ! discrete}}
is a probabilizable space $\left(\Omega,\mathcal{A}\right)$ with
$\Omega$ countable and $\mathcal{A}=\mathcal{P}\left(\Omega\right).$ 

\end{definition}

A probability on a discrete probabilizable space $\Omega$\footnote{Often, when talking about a discrete probabilizable space $\left(\Omega,\mathcal{P}\left(\Omega\right)\right),$
the $\sigma-$algebra $\mathcal{P}\left(\Omega\right)$ is left implicit.} is fully determined once the probabilities of the elementary events
$\left\{ \omega\right\} $ are known. Since any subset $A\subset\Omega$
can be expressed as a disjoint union 
\[
A=\biguplus\limits_{\omega\in A}\left\{ \omega\right\} ,
\]
 it follows that\boxeq{
\begin{equation}
P\left(A\right)=\sum_{\omega\in A}P\left(\left\{ \omega\right\} \right)\label{eq:sum_elementary_proba}
\end{equation}
}for every subset $A$ of $\Omega.$

In the next chapter, we will rigourously define the sum appearing
in the second part of the equation $\refpar{eq:sum_elementary_proba}$.
For now, we proceed with a ``formal'' reasonning, which remains
straightforward in the case of a finite\footnote{However, this property does not hold in general. In Chapter \ref{chap:Random-Variables-with},
we will encounter probability laws---particularly those with densities---defined
on $\Omega=\mathbb{R}$---for which, for every $x\in\mathbb{R},$
$P\left(\left\{ x\right\} \right)=0.$ Nonetheless, in such cases,
the value of $P\left(\left[a,b\right]\right)$ for $a\in\mathbb{R},\,b\in\mathbb{R}$
and $a\neq b$ cannot be deduced from the values of $P\left(\left\{ x\right\} \right).$} $\Omega.$

\subsection{Germ of a Probability Law in a Discrete Space}

Define the function $g:\Omega\rightarrow\mathbb{R}^{+}$ by setting
\[
g\left(\omega\right)=P\left(\left\{ \omega\right\} \right).
\]
Since $P$ is a probability, we have for every $\omega\in\Omega,$
$g\left(\omega\right)\geqslant0$ and $\sum_{\omega\in\Omega}g\left(\omega\right)=1.$

Conversely, if $g$ is a nonnegative function defined on $\Omega,$
such that$\sum\limits_{\omega\in\Omega}g\left(\omega\right)=1,$ then
there exists a unique probability $P$ on $\Omega$ such that $P\left(\left\{ \omega\right\} \right)=g\left(\omega\right).$ 

This is formalized in the following lemma.

\begin{lemma}{Probability Law from a Function}{probalaw}

Let $\Omega$ be a countable set.

(i) \textbf{Probability Law from a Function}

Let $g:\Omega\rightarrow\mathbb{R}^{+}$ be a function such that
\[
\sum\limits_{\omega\in\Omega}g\left(\omega\right)=1.
\]

\boxeq{For $A\in\mathcal{P}\left(\Omega\right),$ we set
\begin{equation}
P\left(A\right)=\sum\limits_{\omega\in A}g\left(\omega\right).\label{eq:probability_from_germ}
\end{equation}
}

The so defined function $P$ yields a probability on the discrete
probabilizable space $\left(\Omega,\mathcal{P}\left(\Omega\right)\right).$

(ii) \textbf{Germ Induced by a Probability}

Any probability on $\left(\Omega,\mathcal{P}\left(\Omega\right)\right)$
arises in this way, with the function $g\left(\omega\right)=P\left(\left\{ \omega\right\} \right).$

\end{lemma}

\begin{proof}{}{}

(i) \textbf{Probability Law from a Function}

We immediately obtain $P\left(\emptyset\right)=0.$ 

Furthermore, let $\left(A_{i}\right)_{i\in\mathbb{N}}$ be a sequence
of pairwise disjoint elements of $\mathcal{A}.$ By the definition
of $P,$ we have
\begin{align*}
P\left(\biguplus\limits_{i\in\mathbb{N}}A_{i}\right) & =\sum\limits_{\omega\in\biguplus\limits_{i\in\mathbb{N}}A_{i}}g\left(\omega\right).
\end{align*}

Since the sets $A_{i}$ are pairwise disjoint, we can rewrite this
sum as

\[
\sum\limits_{\omega\in\biguplus\limits_{i\in\mathbb{N}}A_{i}}g\left(\omega\right)=\sum\limits_{i\in\mathbb{N}}\left(\sum\limits_{\omega\in A_{i}}g\left(\omega\right)\right)=\sum\limits_{i\in\mathbb{N}}P\left(A_{i}\right).
\]

Thus, $P$ satisfies the countable additivity.

(ii) \textbf{Germ Induced by a Probability}

This is straightforward from the discussion preceding the lemma.

\end{proof}

\begin{remark}{}{}

In this proof, we have used the equality
\[
\sum\limits_{i\in\mathbb{N}}\left(\sum\limits_{\omega\in A_{i}}g\left(\omega\right)\right)=\sum\limits_{i\in\mathbb{N}}P\left(A_{i}\right).
\]

This follows from an associativity property, which will be established
in the next chapter in the context of summable families under the
name ``\textit{summation per packet property}''.

\end{remark}

\begin{definition}{Germ of a Probability}{probagerm}A function $g$
from the set $\Omega$ to $\mathbb{R}^{+}$ satisfying 
\[
\sum\limits_{\omega\in\Omega}g\left(\omega\right)=1
\]
is called the \textbf{germ\index{germ of a probability}\mindex{probability ! germ}\footnotemark}
\textbf{of the probability $P,$} as defined by the equation $\refpar{eq:probability_from_germ},$
that is for every $A\in\mathcal{P}\left(\Omega\right),$ $P\left(A\right)=\sum\limits_{\omega\in A}g\left(\omega\right).$

\end{definition}

\footnotetext{The designation ``germ'' is not widely used but appears
in some foreign texts. In this book, we adopt it systematically to
avoid circumlocutions to name this function or the more cumbersome
expression ``function of probability''.}

\subsection{Some Discrete Probability Laws}

Probabilities defined on a discrete probabilizable space are often
referred to as\textbf{ \mindex{probability ! law ! discrete}\mindex{discrete ! probability law}discrete
probability laws}.

Below, we present a series of classical examples of discrete probabilized
spaces that are somewhat basic but fundamental models, often used
as references for constructing and analyzing more complex probabilistic
frameworks.

Later, we will explore the circumstances in which these models naturally
arise.

\subsubsection{Uniform Law (or Uniform Probability) on a Finite Set}\label{subsec:Uniform-law-(or}

\begin{definition}{Uniform Law}{uniflaw} The \textbf{uniform law
}\mindex{law ! uniform}\index{uniform law}on a finite set $\Omega$
is the probability that assigns the same value to each elementary
event $\left\{ \omega\right\} .$ The existence of such a probability
$P$ is guaranteed by the fact that this probability is associated
to a germ $g,$ which is both a constant function on $\Omega$ and
satisfies the condition $\sum\limits_{\omega\in\Omega}g\left(\omega\right)=1.$

Since $g$ is constant, it follows that\footnotemark, for every $\omega\in\Omega,$
$g\left(\omega\right)=\dfrac{1}{\left|\Omega\right|}.$

Thus, for every subset $A$ of $\Omega,$ the probability of $A$
is given by\boxeq{
\[
P\left(A\right)=\dfrac{\left|A\right|}{\left|\Omega\right|}.
\]
}

\end{definition}

\footnotetext{$\left|.\right|$ denotes the number of elements in
a set.}

\begin{remarks}{Two Important Observations}{}
\begin{enumerate}
\item When referring to a random choice from a finite set, it is often implicitly
assumed that the selection follows the \textbf{uniform law}, meaning
each element of the finite set is chosen with equal probability.
\item For an infinite set, however, defining an uniform probability in the
same manner is not possible.
\end{enumerate}
\end{remarks}

\begin{example}{Example of Modelling Using the Uniform Probability}{example_uniform}

A fair die is rolled $n$ times and we seek the probability to obtain
the number 6 exactly $k$ times, where $k\leqslant n.$ A realization
of this experiment is a sequence of $n$ integers, each between 1
and 6. Thus, the sample space is 
\[
\Omega=\left\{ 1,...,6\right\} ^{n}.
\]
 The event ``obtaining the number 6 exactly $k$ times'' is represented
by the subset:
\[
A_{k}=\left\{ \left(x_{1},x_{2},...,x_{n}\right):\,x_{i}=6\text{ for exactly }k\text{ indices}\right\} \subseteq\Omega.
\]
Since the die is fair, a natural way to model this random experiment\footnotemark
is to attach to the discrete probabilizable space $\left(\Omega,\mathcal{P}\left(\Omega\right)\right)$
the uniform probability $P$ assigning the same importance to each
possible outcome.

There are\footnotemark $\binom{n}{k}$ way to choose the $k$ positions
where $x_{i}=6,$ and for the $n-k$ positions remaining there are
$5^{n-k}$ possibilities to assign any of the numbers 1 to 5 in each
of the position, 6 being excluded. Thus,
\[
\left|A_{k}\right|=\binom{n}{k}5^{n-k}.
\]
 Furthermore, $\left|\Omega\right|=6^{n}.$ Hence, the probability
of the event $A_{k}$ is
\[
P\left(A_{k}\right)=\binom{n}{k}\dfrac{5^{n-k}}{6^{n}}
\]

This can also be written
\[
P\left(A_{k}\right)=\binom{n}{k}\left(\dfrac{1}{6}\right)^{k}\left(\dfrac{5}{6}\right)^{n-k}.
\]

This result describes the probability of obtaining exactly $k$ sixes
in $n$ rolls of a fair six-sided die.

\end{example}

\addtocounter{footnote}{-1}

\footnotetext{Another way to justify this choice will be done in
Chapter \ref{chap:Independence}.}

\stepcounter{footnote}

\footnotetext{$\binom{n}{k}$ designates the number of combination
of $n$ objects taken $k$ by $k,$ that is the number of subsets
with $k$ elements in a set with $n$ elements.}

\subsubsection{Geometric Modelling on $\boldsymbol{\mathbb{N}}$ and on $\boldsymbol{\mathbb{N}^{*}}$}

\begin{definition}{Geometric Law on $\mathbb{N}$}{geomlaw}

The outcome space is $\Omega=\mathbb{N}$.

Let $p$ be a real number in $\left]0,1\right[.$

The \textbf{geometric law on $\boldsymbol{\mathbb{N}}$}\index{geometric law on ensuremath{mathbb{N}}@geometric law on $\mathbb{N}$}\mindex{law ! geometric ! on $\mathbb{N}$}
with parameter $p$---denoted by $\mathcal{G}_{\mathbb{N}}\left(p\right)$---is
the probability law induced by the germ $g$ defined as

\boxeq{
\begin{align*}
\forall n\in\mathbb{N},\,\,\,\, & g(n)=q^{n}p
\end{align*}
}

where $q=1-p.$

\end{definition}

\begin{definition}{Geometric Law on $\mathbb{N}^\ast$}{geomlawN}

The outcome space is now $\Omega=\mathbb{N}^{*}$.

Let $p$ be a real number in $\left]0,1\right[.$

The \textbf{geometric law on $\boldsymbol{\mathbb{N}^{*}}$}\index{geometric law on ensuremath{mathbb{N}}^{*}@geometric law on $\mathbb{N}^{*}$}\mindex{law ! geometric ! on $\mathbb{N}^\ast$}
with parameter $p$---denoted by $\mathcal{G}_{\mathbb{N}^{*}}\left(p\right)$---is
the probability law induced by the germ $g$ defined as

\boxeq{
\begin{align*}
\forall n\in\mathbb{N}^{\ast},\,\,\,\, & g(n)=q^{n-1}p
\end{align*}
}

where $q=1-p.$

\end{definition}

\begin{remark}{}{}The geometric law on $\mathbb{N}^{\ast}$ reports
the probability that making $n$ independent trials the first occurrence
of the first success is at the $n$-th step, while the geometric law
on $\mathbb{N}$ reports the number of failures before a first success.

\end{remark}

\subsubsection{Binomial Model}

\begin{definition}{Binomial Law}{binomiallaw}

Let $n$ be a positive integer, that is $n\in\mathbb{N^{*}},$ and
let $p$ be a real number in the open interval $\left]0,1\right[$\footnotemark.

The outcome space is $\Omega=\left\llbracket 0,n\right\rrbracket .$

The \textbf{binomial law}\mindex{binomial!law}\mindex{law ! binomial}---also
known as the \textbf{binomial probability}\mindex{binomial!probability}---with
parameters $n$ and $p$---denoted by $\mathcal{B}\left(n,p\right)$---is
the probability law induced by the germ $g$ defined as\boxeq{
\[
\forall k\in\left\llbracket 0,n\right\rrbracket ,\,\,\,\,g\left(k\right)=\binom{n}{k}p^{k}q^{n-k},
\]
}where $q=1-p.$ 

\end{definition}

\footnotetext{Tr.N. Technically, even if it appears as a degenerated
case, $p$ can be taken in the closed interval $p\in\left[0,1\right].$}

\begin{figure}[t]
\begin{center}\includegraphics[width=0.4\textwidth]{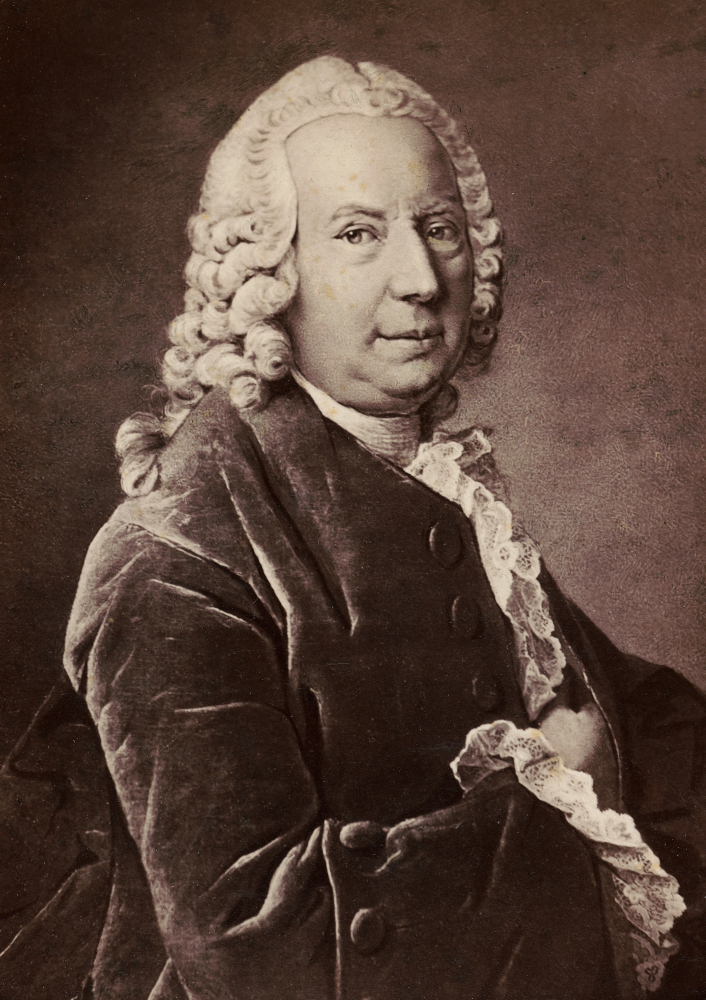}

{\tiny Credits: Collection of the \href{https://commons.wikimedia.org/wiki/User:ETH-Bibliothek}{ETH Bibliotek}
Public Domain}\end{center}

\caption{\textbf{\protect\href{https://en.wikipedia.org/wiki/Daniel_Bernoulli}{Daniel Bernoulli}}
(1700 - 1782) }\sindex[fam]{Bernoulli, Daniel}
\end{figure}

\begin{remarks}{}{}

1. In the example given in Subsection \ref{subsec:Uniform-law-(or},
where we compute the probability of obtaining exactly $k$ occurrences
of the number 6 in $n$ throws of a fair six-sided classical die,
we may adopt an alternative model if we are only interested in the
result, and not on the way to obtain it. In this case, the set of
realizations becomes
\[
\Omega^{\prime}=\left\llbracket 0,n\right\rrbracket ,
\]
and we define a probability law $P^{\prime}$ on $\Omega^{\prime}$
such that 
\[
P^{\prime}\left(\left\{ k\right\} \right)=P\left(A_{k}\right).
\]

Thus, the probability law $P^{\prime}$ is a binomial law, 
\[
P^{\prime}=\mathcal{B}\left(n,\dfrac{1}{6}\right).
\]

2. More generally, consider an experiment with only two possible outcomes,
referred to as a \textbf{Bernoulli\footnotemark trial}\index{Bernoulli trial}\mindex{trial ! index}:
``success'' with probability $p$ and ``failure'' with probability
$q=1-p.$ When repeating this Bernoulli trial $n$ times independently,\footnotemark
the probability law governing the event ``obtain exactly $k$ successes
in $n$ trials'' follows the binomial law $\mathcal{B}\left(n,p\right).$

\end{remarks}

\addtocounter{footnote}{-1}

\footnotetext{\textbf{\href{https://en.wikipedia.org/wiki/Daniel_Bernoulli}{Daniel Bernoulli}}\sindex[fam]{Bernoulli, Daniel}
(1700 - 1782) is a Swiss-French mathematician and physicist. His work
focused on fluid mechanics, probability and statistics.}

\stepcounter{footnote}

\footnotetext{Tr.N. The exact meaning of independence will be given
in Chapter \ref{chap:Independence}.}

\subsubsection{Poisson Model}

\begin{figure}[t]
\begin{center}\includegraphics[width=0.4\textwidth]{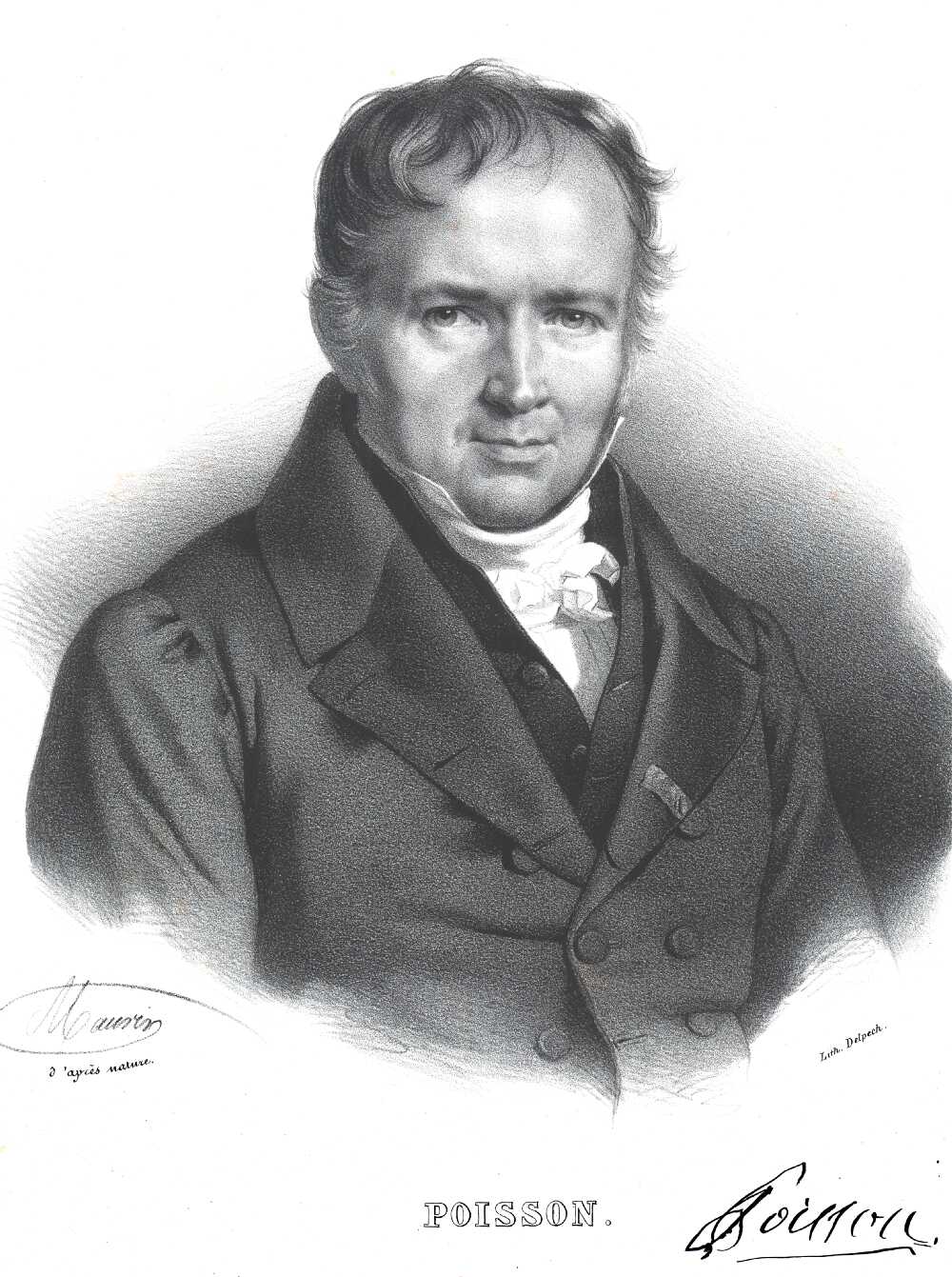}

{\tiny Credits: \href{http://www.sil.si.edu/digitalcollections/hst/scientific-identity/CF/display_results.cfm?alpha_sort=P}{SIL}
Public Domain}\end{center}

\caption{\textbf{\protect\href{https://en.wikipedia.org/wiki/Sim\%C3\%A9on_Denis_Poisson}{Siméon Denis Poisson}}
(1781-1840)}\sindex[fam]{Poisson, Siméon Denis}
\end{figure}

\begin{definition}{Poisson Law}{poissonlaw}

The sample space is $\Omega=\mathbb{N}.$

Let $\lambda$ be a positive real number.

The \textbf{Poisson\footnotemark law}\index{Poisson law}\mindex{law ! Poisson}
with parameter $\lambda$---denoted by $\mathcal{P}\left(\lambda\right)$---is
the probability induced by the germ defined as\boxeq{
\[
\forall n\in\mathbb{N},\,\,\,\,g\left(n\right)=\textrm{e}^{-\lambda}\dfrac{\lambda^{n}}{n!}.
\]
}

\end{definition}

\footnotetext{Named after \textbf{\href{https://en.wikipedia.org/wiki/Sim\%C3\%A9on_Denis_Poisson}{Siméon Denis Poisson}}\sindex[fam]{Poisson, Siméon Denis}
(1781-1840), a professor of mathematics at the Ecole Polytechnique
(1806) and to the Paris Faculty of Sciences (1809), whose work contributed
to mathematical analysis, probability theory, celestial mechanics
and mathematical physics.}

\subsubsection{Hypergeometric Model}

Let $U_{1}$ and $U_{2}$ be two disjoint subsets of a set $U$ such
that $U=U_{1}\uplus U_{2}.$ Denote the cardinalities by $\left|U\right|=r,$
$\left|U_{1}\right|=r_{1}$ and, $\left|U_{2}\right|=r_{2},$ where
$r_{2}=r-r_{1}.$ Let $n$ be an integer such that $1\leqslant n<r.$ 

We randomly select $n$ elements of $U$ and study the probability
of obtaining exactly $k_{1}$ elements from $U_{1}$---which implies
selecting exactly $k_{2}=n-k_{1}$ elements from $U_{2}.$

A realization of this experiment is a subset $A$ of $U$ containing
exactly $n$ elements. The outcome space $\Omega$ is thus given by
\[
\Omega=\left\{ A\in\mathcal{P}\left(U\right):\,\left|A\right|=n\right\} .
\]
Its cardinality is $\left|\Omega\right|=\binom{r}{n}.$

The studied event is then represented by the subset $A_{k_{1}}$ of
$\Omega$ defined as
\[
A_{k_{1}}=\left\{ A\in\Omega:\,\left|A\cap U_{1}\right|=k_{1}\right\} .
\]

The probabilizable space $\left(\Omega,\text{\ensuremath{\mathcal{P}}\ensuremath{\left(\Omega\right)}}\right)$
is then equiped with the uniform probability $P.$

The set $A_{k_{1}}$ is empty if and only if $r_{1}<k_{1}<n$, or
if $0<k_{1}<n-\left(r-r_{1}\right).$ 

Otherwise, if
\begin{equation}
\text{max}\left(0,n-\left(r-r_{1}\right)\right)\leqslant k_{1}\leqslant\text{min}\left(n,r_{1}\right),\label{eq:hypergeometric_val_k_1}
\end{equation}
then its cardinality is given by 
\[
\left|A_{k_{1}}\right|=\binom{r_{1}}{k_{1}}\binom{r-r_{1}}{n-k_{1}}.
\]

Hence, the probability of the event is\boxeq{
\begin{equation}
P\left(A_{k_{1}}\right)=\dfrac{\binom{r_{1}}{k_{1}}\binom{r-r_{1}}{n-k_{1}}}{\binom{r}{n}}\label{eq:hypergeometric probability eq}
\end{equation}
}

\begin{definition}{Hypergeometric Law}{hypergeomlaw}The \textbf{hypergeometric
law\index{hypergeometric law}\mindex{law ! hypergeometric}} is
the law of probability associated with the germ $g$ defined on the
set of integers $k_{1}$ satisfying:
\begin{equation}
\text{max}\left(0,n-\left(r-r_{1}\right)\right)\leqslant k_{1}\leqslant\text{min}\left(n,r_{1}\right),\label{eq:hypergeometric_val_k_1-1}
\end{equation}
such that\boxeq{
\[
g\left(k_{1}\right)=\dfrac{\binom{r_{1}}{k_{1}}\binom{r-r_{1}}{n-k_{1}}}{\binom{r}{n}}.
\]
}

\end{definition}

\begin{example}{Hypergeometric Law Modelling}{}

A lake contains $r$ fishes, among which $r_{1}$ are infected with
a disease. We randomly select $n$ fishes from the lake. The probability
of obtaining exactly $k_{1}$ diseased fishes follows the hypergeometric
model, where:
\begin{itemize}
\item $U$ represents the total fish population,
\item $U_{1}$ represents the diseased fishes,
\item $U_{2}$ represents the healthy fishes.
\end{itemize}
This probability is given by the hypergeometric probability formula
\begin{equation}
P\left(A_{k_{1}}\right)=\dfrac{\binom{r_{1}}{k_{1}}\binom{r-r_{1}}{n-k_{1}}}{\binom{r}{n}}.\label{eq:hypergeometric probability eq-1}
\end{equation}

\end{example}

\begin{remark}{}{}

This model forms the foundation of \textbf{polling theory}, where
a randomly selected sample from a population is analyzed to estimate
the characteristics of the entire population.

\end{remark}

\section{Random Variables}\label{sec:Random-Variables}

When studying a phenomenon, numerical or vectorial values associated
with it are often analyzed. In broad terms, a \textbf{random variable\index{random variable}}
is a numerical quantity associated with a random experiment, whose
value depends solely on the outcome $\omega$ of the experiment. Mathematically,
it is simply a function defined on a set $\Omega.$

\begin{example}{Coin Tossing and a Random Variable}{}{}

Consider a game in which a fair coin is tossed $n$ times. The number
of tails observed depends on the sequence of tosses. To model this
game:
\begin{itemize}
\item Assign 0 to heads;
\item Assign 1 to tails.
\end{itemize}
Thus, an outcome of an experiment is an $n-$tuple $\omega=\left(\omega_{1},...,\omega_{n}\right),$where
each $\omega_{i}$ takes values in $\left\{ 0,1\right\} .$ The sample
space is then $\Omega=\left\{ 0,1\right\} ^{n}.$

Since the coin is fair and the game is honest, we model the experiment
using the uniform probability $P$ on the probabilizable space $\left(\Omega,\mathcal{P}\left(\Omega\right)\right).$

Now, define $X$ the function on $\Omega$ over $\mathbb{N}$ by 
\[
X\left(\omega\right)=\sum\limits^{n}_{i=1}\omega_{i}.
\]

The function $X$ is called a \textbf{real random variable}\index{real random variable},
representing the number of tails obtained.

For an integer $k,$ consider the pre-image $X^{-1}\left(\left\{ k\right\} \right)$
of the singleton $\left\{ k\right\} $
\[
X^{-1}\left(\left\{ k\right\} \right)=\left\{ \omega\in\Omega:\,X\left(\omega\right)=k\right\} .
\]
This set corresponds to the event: ``exactly $k$ tails are obtained
in $n$ tosses''.

The standard notation for the event $X^{-1}\left(\left\{ k\right\} \right)$
is $\left\{ X=k\right\} $ or simply $\left(X=k\right).$

Instead of writing $P\left(\left\{ X=k\right\} \right),$ we use a
more concise notation $P\left(X=k\right).$

It follows that\boxeq{
\[
P\left(X=k\right)=\dfrac{\binom{n}{k}}{2^{n}}.
\]
}The function $\left\{ \begin{array}{c}
\left\llbracket 0,n\right\rrbracket \rightarrow\left[0,1\right]\\
k\mapsto P\left(X=k\right)
\end{array}\right.$ is the \textbf{germ} of the \textbf{binomial probability law} $\mathcal{B}\left(n,\dfrac{1}{2}\right).$

\end{example}

After this introductory example, we now present the general definition
of a random variable.

\begin{definition}{Random Variable}{randomvar}\label{def: random variable}

Let $\left(\Omega,\mathcal{A}\right)$ and $\left(E,\mathcal{E}\right)$
be two probabilizable spaces. A function $X$ from $\Omega$ to $E$
is said to be a \textbf{random variable\index{random variable}} with
values in $E$ if, for every $A\in\mathcal{E}$ the set $X^{-1}\left(A\right)$
belongs to $\mathcal{A},$ where we define the inverse image of $A$
under $X$ as\boxeq{
\[
X^{-1}\left(A\right)=\left\{ \omega\in\Omega:\,X\left(\omega\right)\in A\right\} .
\]
}

\end{definition}

\begin{remark}{}{}

Suppose we have chosen the probabilizable space $\left(\Omega,\mathcal{A}\right)$
to model the random experiment. Given a random variable $X,$ the
inverse image of a set $A$ under $X$ is the event consisting of
all outcomes $\omega$ such that $X\left(\omega\right)$ is in $A.$
This event is often denoted succinctly as $\left(X\in A\right),$
and we refer to it as ``the event that $X$ belongs to $A.$''

\end{remark}

The following proposition serves a prerequisite for defining discrete
random variables in Definition \ref{def: discrete random variable}.

\begin{proposition}{Conditions for a Function to be a Random Variable}{randomvarcond}Let
$X$ be a function from $\Omega$ to $E$ such that $X\left(\Omega\right)$
is countable. Suppose the $\sigma-$algebra $\mathcal{E}$ on $E$
satisfies such, that for every $x\in E,$ the singleton $\left\{ x\right\} $
belongs to $\mathcal{E}.$ Then, for $X$ to be a random variable,
it is necessary and sufficient that for every $x\in E,$ the preimage
$X^{-1}\left(\left\{ x\right\} \right)$ belongs to $\mathcal{A}.$

\end{proposition}

\begin{proof}{}{}

The necessity of the condition follows immediately from the definition
of a random variable.

Conversely, assume that for every $x\in E,$ we have $X^{-1}\left(\left\{ x\right\} \right)\in\mathcal{A}.$

Let $A\in\mathcal{E}.$ Then, we can express the preimage of $A$
under $X$ as follows:
\begin{align*}
X^{-1}\left(A\right) & =X^{-1}\left(X\left(\Omega\right)\cap A\right)\\
 & =\biguplus\limits_{x\in X\left(\Omega\right)\cap A}X^{-1}\left(\left\{ x\right\} \right).
\end{align*}

Since $X\left(\Omega\right)\cap A$ is countable, we have $X^{-1}\left(A\right)\in\mathcal{A}.$
Thus, $X$ is a random variable.

\end{proof}

\begin{definition}{Discrete Random Variable}{discreterandomvar}\label{def: discrete random variable}Let
$\left(\Omega,\mathcal{A}\right)$ be a probabilizable space, and
let $E$ be a set. A mapping $X$ from $\Omega$ to $E$ is said to
be a \textbf{\mindex{discrete ! random variable}\mindex{random variable ! discrete}discrete
random variable} if the two following conditions hold:

\boxeq{(i) The set $X\left(\Omega\right)$ of values taken by $X$
is countable.

(ii) For every $x\in E,$ we have $X^{-1}\left(\left\{ x\right\} \right)\in\mathcal{A}.$}

\end{definition}

\begin{remarks}{}{}

1. If we equip $E$ with a $\sigma-$algebra $\mathcal{E},$ then
every discrete random variable taking values in $E$ is also a random
variable with values in $\left(E,\mathcal{E}\right),$ in the sense
of Definition \ref{def: random variable}, if and only if the $\sigma-$algebra
$\mathcal{E}$ contains all singletons of $E.$ 

2. Initially, we restrict our study to discrete random variables,
where $E$ will be one of $\mathbb{N},$ $\mathbb{Z},$ $\mathbb{N}^{n}$
or $\mathbb{Z}^{n},$ with $\mathcal{E}$ taken as the power set $\mathcal{P}\left(E\right).$

3. If $\mathcal{A}$ is the $\sigma-$algebra consisting of all subsets
of $\Omega$---a common situation for discrete random variables---then
any mapping from $\Omega$ to a space $E$ equipped with an arbitrary
$\sigma-$algebra $\mathcal{E}$ is a random variable from $\left(\Omega,\mathcal{A}\right)$
to $\left(E,\mathcal{E}\right).$ We will later see that if $\Omega$
is uncountable, we will need to refine the choice of $\mathcal{A},$
restricting to a smaller $\sigma-$algebra than the full power set
of $\Omega.$ 

\end{remarks}

\begin{proposition}{Probability Law of a Random Variable}{randomvarprobalaw}

Let $X$ be a random variable defined on a probabilized space $\left(\Omega,\mathcal{A},P\right),$
taking values in a probabilizable space $\left(E,\mathcal{E}\right).$
The mapping
\[
P_{X}:\left\{ \begin{array}{c}
\mathcal{E}\rightarrow\left[0,1\right]\\
A\mapsto P\left[X^{-1}\left(A\right)\right]
\end{array}\right.
\]
defines a probability on the probabilizable space $\left(E,\mathcal{E}\right),$
called the \textbf{probability law}\mindex{probability!law} of the
random variable $X.$ 

\end{proposition}

\begin{proof}{}{}

If $\left(A_{i}\right)_{i\in\mathbb{N}}$ is a sequence of pairwise
disjoint sets in $\mathcal{E},$ then we have
\[
X^{-1}\left(\biguplus_{i\in\mathbb{N}}A_{i}\right)=\biguplus_{i\in\mathbb{N}}X^{-1}\left(A_{i}\right).
\]

Thus, by applying the $\sigma-$additivity of $P,$ we obtain
\begin{align*}
P_{X}\left(\biguplus_{i\in\mathbb{N}}\left(A_{i}\right)\right) & =P\left(X^{-1}\left(\biguplus_{i\in\mathbb{N}}\left(A_{i}\right)\right)\right)\\
 & =P\left(\biguplus_{i\in\mathbb{N}}X^{-1}\left(A_{i}\right)\right)\\
 & =\sum_{i\in\mathbb{N}}P\left(X^{-1}\left(A_{i}\right)\right)\\
 & =\sum_{i\in\mathbb{N}}P_{X}\left(A_{i}\right).
\end{align*}

Furthermore, since $X^{-1}\left(E\right)=\Omega,$ it follows that
$P_{X}\left(E\right)=1.$

\end{proof}

\begin{remarks}{}{}

1. The notion of the \textbf{law of a random variable} is fundamental,
as it describes the probabilistic behaviour of the random variable
and serves as the primary means through which a statistician can extract
probabilistic information. While the probabilized space $\left(\Omega,\mathcal{A},P\right)$
is essential for probability theory, it is, in most cases, not directly
observable in practical applications. Often, making its existence
explicit is unnecessary; knowing that it exists suffices. It is important
to note that the law of a random variable is a probability measure
on the space of values that the random variable takes.

2. If $X$ is a discrete random variable taking values in $\left(E,\mathcal{E}\right),$
its \textbf{law} is completely \textbf{determined by the probability
germ} $x\mapsto P\left(X=x\right)$ on $\left(E,\mathcal{E}\right).$ 

Indeed, for every $A\in\mathcal{E},$
\begin{align*}
P\left[X^{-1}\left(A\right)\right] & =P\left[\biguplus\limits_{x\in X\left(\Omega\right)\cap A}X^{-1}\left(\left\{ x\right\} \right)\right]\\
 & =\sum_{x\in X\left(\Omega\right)\cap A}P\left[X^{-1}\left(\left\{ x\right\} \right)\right]\\
 & =\sum_{x\in A}P\left(X=x\right).
\end{align*}

The second equality holds trivially when $X\left(\Omega\right)$ is
finite. If $X\left(\Omega\right)$ is infinite, then an additional
justification is required, involving the concept of summable families.
The systematic study of summable families will be addressed in the
next chapter.

\end{remarks}

\section*{Exercises}

\addcontentsline{toc}{section}{Exercises}

\begin{exercise}{Unique Faces in Dice Rolls}{ex1.1}

What is the probability that, when rolling six fair and distinguishable
dice, all six faces show a different number?

\end{exercise}

\begin{exercise}{Coin Tossing and Events}{ex1.2}

We toss $n$ coins, with $n\geqslant2.$

1. What is the probability that the outcome contains both heads and
tails (event $A$)? 

2. What is the probability that the outcome consists entirely of heads
and contains at most one tail (event $B$)?

3. Construct a probabilistic model and compute the probability of
$A\cap B.$ Compare it to the product of the probabilities of $A$
and $B.$

\end{exercise}

\begin{exercise}{Birthday Problem and Repeated Outcomes}{ex1.3}

1. A box contains $M$ tokens numbered from 1 to $M.$ We draw $n$
tokens successively, placing back the drawn token and mixing well
the tokens in the box each time. What is the probability that no token
is drawn more than once?

2. A classroom contains $n$ students. Assuming that birthdays are
uniformly distributed across the 365 days of the year, and ignoring
leap years, what is the probability that at least two students share
the same birthday, given that the number of students $n\leqslant365?$

\end{exercise}

\begin{exercise}{The Lottery Problem}{ex1.4}

A lottery prints $M$ vouchers numbered from 1 to $M.$ 

Without loss of generality, suppose that the first $n$ vouchers $\left(2n\leqslant M\right)$
are the winning ones---though buyers do not know this. What is the
probability that a buyer who purchases $n$ tickets gets at least
one winning ticket?

\end{exercise}

\begin{exercise}{An Application of Poincaré Formula}{ex1.5}

Let $\left(\Omega,\mathcal{A},P\right)$ be a probabilized space,
and let $\left\{ A_{i}\right\} _{1\leqslant i\leqslant n}$ be a finite
family of events, where each $A_{i}$ belongs to $\mathcal{A}.$ For
every non-empty subset $J$ of $I=\left\{ 1,\dots,n\right\} ,$ define
\[
\widehat{A_{J}}=\bigcap\limits_{j\in J}A_{j}
\]
and
\[
B^{J}=\widehat{A_{J}}\cap\left(\bigcap\limits_{j\in I\backslash J}A^{c}_{j}\right).
\]

For $1\leqslant m\leqslant n,$
\[
B_{m}=\biguplus\limits_{\substack{J\in\mathcal{P}\left(I\right)\\
\left|J\right|=m
}
}B^{J}.
\]

Finally, let
\[
B_{0}=\bigcap\limits^{n}_{j=1}A^{c}_{j}.
\]

Preliminary question: Show that for two distinct non-empty subsets
$J$ and $J^{\prime}$ of $I,$ the sets $B^{J}$ and $B^{J^{\prime}}$
are disjoint.

A point $\omega\in\Omega$ belongs to $B_{m}$ $\left(1\leqslant m\leqslant n\right)$
if and only if it belongs to exactly $m$ events $A_{i}.$ $\omega$
belongs to $B_{0}$ if and only if it does not belong to any of the
events $A_{i}.$

Define
\[
S_{0}=1
\]

and, for $1\leqslant r\leqslant n,$
\[
S_{r}=\sum\limits_{\substack{J\in\mathcal{P}\left(I\right)\\
\left|J\right|=r
}
}P\left(\widehat{A_{J}}\right).
\]

1. Prove that:

(a) For every $A\in\mathcal{A},$
\begin{align}
P\left(B_{0}\cap A\right) & =P\left(A\right)-\sum\limits^{n}_{j=1}P\left(A_{j}\cap A\right)+\sum\limits_{1\leqslant j_{1}<j_{2}\leqslant n}P\left(A_{j_{1}}\cap A_{j_{2}}\cap A\right)\nonumber \\
 & -\sum\limits_{1\leqslant j_{1}<j_{2}<j_{3}\leqslant n}P\left(A_{j_{1}}\cap A_{j_{2}}\cap A_{j_{3}}\cap A\right)\nonumber \\
 & +\dots+\left(-1\right)^{n}P\left(A_{1}\cap A_{2}\cap\cdots\cap A_{n}\cap A\right);\label{eq:ex5qu1a}
\end{align}

(b) For $1\leqslant m\leqslant n,$
\begin{equation}
P\left(B_{m}\right)=S_{m}-\binom{m+1}{m}S_{m+1}+\cdots+\left(-1\right)^{n-m}\binom{n}{m}S_{n}.\label{eq:ex5qu1b}
\end{equation}

2. Example. During a social event, $n$ people each write their name
on a card and place it in a large hat---assuming no duplicate names.
At the end of the event, after thoroughly mixing the cards, each person
draws one card at random. 

a. What is the probability that nobody draws their own name?

b. What is the probability that exactly $m$ $\left(1\leqslant m\leqslant n\right)$
persons draw a card with their own name?

\end{exercise}

\begin{exercise}{Example of Random Variable: the Head and Tail Game}{1.6}

A player tosses a fair coin $N$ times $\left(N\geqslant2\right).$
We are interested in studying the random variable representing the
number of the trial in which ``tail'' appears for the first time.

Construct a probabilistic model and determine the probability law
of this random variable.

\end{exercise}

\section*{Solution of Exercises}

\addcontentsline{toc}{section}{Solution of Exercises}

\begin{solution}{}{}

We choose as the outcome space the set of realizations $\Omega=\left\{ 1,2,...,6\right\} ^{6},$
equipped with the $\sigma-$algebra $\mathcal{P}\left(\Omega\right).$
As probability $P$ on $\left(\Omega,\mathcal{A}\right),$ we consider
the uniform probability corresponding to the assumption that the dice
are fair. The event: ``all faces show an upper face with a different
number'' is represented by the subset

\[
A=\left\{ \left(\sigma\left(1\right),\sigma\left(2\right),\dots,\sigma\left(6\right)\right):\,\sigma\in\mathcal{S}_{6}\right\} ,
\]
where $\mathcal{S}_{6}$ corresponds to the set of all permutations
of the integers from 1 to 6.

We have $\left|\Omega\right|=6^{6}$ and $\left|A\right|=6!,$ which
yields
\[
P\left(A\right)=\dfrac{\left|A\right|}{\left|\Omega\right|}=\dfrac{6!}{6^{6}}.
\]

Hence, \boxeq{
\[
P\left(A\right)=\dfrac{720}{46'656}\approx0.015.
\]
}

\end{solution}

\begin{solution}{}{}

We take as outcome space the set $\Omega=\left\{ H,T\right\} ^{n},$
equipped with the $\sigma-$algebra $\mathcal{P}\left(\Omega\right).$
As probability $P$ on $\left(\Omega,\mathcal{A}\right),$ we consider
the uniform probability.

\textbf{1. Probability that the outcome contains both heads and tails
(event A)}

The complementary event to the event of having both heads and tails
is the event where there are only heads or only tails. This is represented
by the event
\[
A^{c}=\left\{ H\right\} ^{n}\uplus\left\{ T\right\} ^{n}.
\]

We have $\left|\Omega\right|=2^{n}$ and $\left|A^{c}\right|=2.$

Therefore,

\boxeq{
\begin{align*}
P\left(A\right) & =1-P\left(A^{c}\right)=1-\dfrac{1}{2^{n-1}}.
\end{align*}
}

\textbf{2. Probability of having at most one tail (event B)}

The event of having only heads is denoted by $B_{0}$ and the event
of having exactly one tail is denoted by $B_{1}.$

We have $\left|B_{0}\right|=1$ and $\left|B_{1}\right|=n.$

Since $B=B_{0}\biguplus B_{1},$ we obtain\boxeq{
\begin{align*}
P\left(B\right) & =P\left(B_{0}\right)+P\left(B_{1}\right)=\dfrac{n+1}{2^{n}}.
\end{align*}
}

\textbf{3. Probability of $A\cap B$ and comparison to $P\left(A\right)\times P\left(B\right).$}

The event $A\cap B$ is the event where both event $A$ and $B$ occur
simultaneously. 

Since $A\cap B=B_{1},$
\[
P\left(A\cap B\right)=\dfrac{n}{2^{n}}.
\]

Additionally,
\[
P\left(A\right)\times P\left(B\right)=\dfrac{\left(2^{n-1}-1\right)\left(n+1\right)}{2^{2n-1}}.
\]

We now determine, if there exists, the value(s) of $n$ such that
\[
P\left(A\cap B\right)=P\left(A\right)\times P\left(B\right)
\]

This leads to the equation
\[
\left(2^{n-1}-1\right)\left(n+1\right)=n2^{n-1}
\]
which is equivalent to:
\[
2^{n-1}=n+1
\]

The above equation holds only for $n=3.$

Tr/N: Indeed, consider $f:\mathbb{R}^{+}\mapsto\mathbb{R},$ defined
for $x\in\mathbb{R}^{+}$ by
\[
f(x)=2^{x-1}-x-1.
\]
This function is continuous and derivable, and we have
\[
f^{\prime}\left(x\right)=\left(\ln2\right)e^{\left(x-1\right)\ln2}-1.
\]

$f'\left(x_{0}\right)>0$ for $x_{0}>1-\dfrac{\ln\left(\ln2\right)}{\ln2}\approx1.529.$

We have $f\left(0\right)=-\dfrac{1}{2},$ $f$ strictly decreasing
on $\left[0,x_{0}\right],$ thus $f\left(x_{0}\right)<f\left(0\right)<0.$
Additionally $\lim_{x\to+\infty}f\left(x\right)=+\infty.$ $f$ being
a continuous and strictly increasing function on $\left[x_{0};+\infty\right[,$
thus there exists a unique value by the intermediate value theorem
such that $f\left(x_{1}\right)=0$ on $\left[x_{0};+\infty\right[.$
We have $f\left(3\right)=0.$ Thus $x_{1}=3.$

\boxeq{Hence, for $n\neq3,$ 
\[
P\left(A\cap B\right)\neq P\left(A\right)\times P\left(B\right).
\]
}

\end{solution}

\begin{solution}{}{}

\textbf{1. Probability that no token is drawn more than once}

Each outcome is an $n-$tuple of integers between 1 and $M.$ Also,
we take as the outcome space the set 
\[
\Omega=\left\{ 1,2,...,M\right\} ^{n}.
\]

The event under study can be written as
\[
A=\left\{ \left(\omega_{1},\omega_{2},...,\omega_{n}\right):\,\omega_{i}\neq\omega_{j}\,\text{whenever}\,i\neq j\right\} .
\]

We have
\[
\left|\Omega\right|=M^{n},
\]

and\footnotemark
\[
\left|A\right|=P^{n}_{M}=\dfrac{M!}{n!}.
\]

It follows that
\begin{align*}
P\left(A\right) & =\dfrac{\left|A\right|}{\left|\Omega\right|}\\
 & =\dfrac{M!}{n!}\times\dfrac{1}{M^{n}}\\
 & =\dfrac{M\left(M-1\right)\times\dots\times\left(M-n+1\right)}{M^{n}}\\
 & =\dfrac{M-1}{M}\times\dots\times\dfrac{M-\left(n-1\right)}{M}
\end{align*}
Hence,\boxeq{
\[
P\left(A\right)=\prod^{n-1}_{i=1}\left(1-\dfrac{i}{M}\right).
\]
}

\textbf{2. Birthday Problem}

Replacing the word ``token'' with ``student'', the numbered tokens
correspond to the 365 days of the year. We study the complementary
event, denoted $B,$ to the event $\overline{B}:$ ``no student share
the same birthday''. The event $\overline{B}$ corresponds to $M=365.$ 

Hence, the probability of the event $B$ is\boxeq{
\[
P_{n}\left(B\right)=1-\dfrac{365!}{365^{n}\times\left(365-n\right)!}.
\]
}

Below are some probability values depending on $n.$

\begin{center}%
\begin{tabular}{|c|c|c|c|c|c|c|}
\hline 
$n$ &
4 &
16 &
22 &
23 &
40 &
64\tabularnewline
\hline 
$P_{n}\left(B\right)$ &
0,016 &
0,284 &
0.476 &
0.507 &
0.891 &
0.997\tabularnewline
\hline 
\end{tabular}\end{center}

\end{solution}

\footnotetext{$P^{n}_{r}$ corresponds to the number of permutations
of $n$ distinct objects taken $r$ at a time without repetition.}

\begin{solution}{}{}

Each outcome is a $n-$uple of distinct integers between $1$ and
$M.$ We take as outcome space the set
\[
\Omega=\left\{ A\in\mathcal{P}\left(\left\{ 1,2,...,M\right\} \right):\,\left|A\right|=n\right\} .
\]

Let $G$ be the event: ``the client received at least one winning
ticket''. 

Its complementary event is
\[
G^{c}=\left\{ A\in\mathcal{P}\left(\left\{ n+1,...,M\right\} \right):\,\left|A\right|=n\right\} .
\]

We have
\[
\begin{array}{ccccc}
\left|\Omega\right|=\binom{M}{n} &  & \text{et} &  & \left|G^{c}\right|=\binom{M-n}{n}.\end{array}
\]

Therefore, the probability of the event $G$ is
\begin{align*}
P\left(G\right) & =1-P\left(G^{c}\right)=1-\dfrac{\binom{M-n}{n}}{\binom{M}{n}}
\end{align*}

This simplifies to\boxeq{
\[
P\left(G\right)=1-\prod^{n-1}_{i=0}\left(1-\dfrac{n}{M-i}\right).
\]
}

\end{solution}

\begin{solution}{}{}

\textbf{Preliminary question}

Given that $J$ and $J'$ are distinct subsets of $I,$ there exists
$j_{0}\in J$ such that $j_{0}\notin J'.$
\begin{itemize}
\item Since $j_{\ensuremath{0}}\in J$, it follows that
\[
B^{J}=A_{j_{0}}\cap\left(\bigcap\limits_{j\in J\backslash\left\{ j_{0}\right\} }A_{j}\right)\cap\left(\bigcap\limits_{j\in I\backslash J}A^{c}_{j}\right).
\]
\item Since $j_{0}\notin J',$ we have
\[
B^{J'}=\left(\bigcap\limits_{j\in J^{\prime}}A_{j}\right)\cap A^{c}_{j_{0}}\cap\left(\bigcap\limits_{j\in\left(I\backslash J^{\prime}\right)\backslash\left\{ j_{0}\right\} }A^{c}_{j}\right).
\]
\end{itemize}
Hence, $B^{J}\cap B^{J^{\prime}}=\emptyset$ since $A_{j_{0}}$ and
$A^{c}_{j_{0}}$ cannot occur simultaneously.

\textbf{1. (a) Computation of $P\left(B_{0}\cap A\right)$}

For every event $A\in\mathcal{A},$
\begin{align*}
A & =\Omega\cap A=\left(B_{0}\uplus B^{c}_{0}\right)\cap A=\left(B_{0}\cap A\right)\uplus\left(B^{c}_{0}\cap A\right),
\end{align*}
which implies
\begin{equation}
P\left(B_{0}\cap A\right)=P\left(A\right)-P\left(B^{c}_{0}\cap A\right).\label{eq:1.12}
\end{equation}

Using the complement
\[
B^{c}_{0}=\left(\bigcap\limits^{n}_{j=1}A^{c}_{j}\right)^{c}=\bigcup\limits^{n}_{j=1}A_{j},
\]
 we have, by sustituting in the equation $\refpar{eq:1.12},$
\begin{align*}
P\left(B_{0}\cap A\right) & =P\left(A\right)-P\left(\left(\bigcup\limits^{n}_{j=1}A_{j}\right)\cap A\right)=P\left(A\right)-P\left(\bigcup\limits^{n}_{j=1}\left(A_{j}\cap A\right)\right).
\end{align*}

Applying the Poincaré formula, we obtain the result.

\textbf{(b) Computation of $P\left(B_{m}\right)$}

Let $J$ be a non-empty subset of $I$ with $m$ elements, and let
$i_{1},i_{2},\dots,i_{n-m}$ be a nondecreasing enumeration of $I\backslash J.$ 

Using the previous formula with
\[
\begin{array}{ccccc}
A=\widehat{A_{J}} &  & \text{and} &  & B_{0}=\bigcap\limits^{n-m}_{h=1}A^{c}_{i_{h}},\end{array}
\]
and remarking that
\begin{align*}
B^{J} & =\widehat{A_{J}}\cap\left(\bigcap_{j\in I\backslash J}A^{c}_{j}\right)=\widehat{A_{J}}\cap\left(\bigcap^{n-m}_{h=1}A^{c}_{i_{h}}\right),
\end{align*}
we have
\begin{multline*}
P\left(B^{J}\right)=P\left(\widehat{A_{J}}\right)-\sum\limits^{n-m}_{h=1}P\left(A_{i_{h}}\cap\widehat{A_{J}}\right)+\sum\limits_{1\leqslant j_{1}<j_{2}\leqslant n-m}P\left(A_{j_{1}}\cap A_{j_{2}}\cap\widehat{A_{J}}\right)\\
-\sum\limits_{1\leqslant j_{1}<j_{2}<j_{3}\leqslant n-m}P\left(A_{j_{1}}\cap A_{j_{2}}\cap A_{j_{3}}\cap\widehat{A_{J}}\right)\\
+\dots+\left(-1\right)^{n-m}P\left(A_{1}\cap A_{2}\cap\cdots\cap A_{n}\cap\widehat{A_{J}}\right).
\end{multline*}
which can be written 
\begin{multline}
P\left(B^{J}\right)=P\left(\widehat{A_{J}}\right)-\sum\limits_{\substack{J^{\prime}J\\
\left|J^{\prime}\right|=m+1
}
}P\left(\widehat{A_{J^{\prime}}}\right)+\sum\limits_{\substack{J^{\prime}J\\
\left|J^{\prime}\right|=m+2
}
}P\left(\widehat{A_{J^{\prime}}}\right)+\dots+\left(-1\right)^{n-m}P\left(\widehat{A_{I}}\right).\label{eq:ex_app_poincare_1.b}
\end{multline}

But, since
\[
P\left(B_{m}\right)=\sum\limits_{\left|J\right|=m}P\left(B^{J}\right),
\]
using the equation $\refpar{eq:ex_app_poincare_1.b}$ yields
\begin{multline*}
P\left(B_{m}\right)=\sum\limits_{\left|J\right|=m}P\left(\widehat{A_{J}}\right)-\sum\limits_{\left|J\right|=m}\sum\limits_{\substack{J^{\prime}J\\
\left|J^{\prime}\right|=m+1
}
}P\left(\widehat{A_{J^{\prime}}}\right)+\sum\limits_{\left|J\right|=m}\sum\limits_{\substack{J^{\prime}J\\
\left|J^{\prime}\right|=m+2
}
}P\left(\widehat{A_{J^{\prime}}}\right)\\
+\dots+\left(-1\right)^{n-m}\sum\limits_{\left|J\right|=m}P\left(\widehat{A_{I}}\right).
\end{multline*}

By associativity of the sum,
\begin{align*}
\sum\limits_{\left|J\right|=m}\sum\limits_{\substack{J^{\prime}J\\
\left|J^{\prime}\right|=m+1
}
}P\left(\widehat{A_{J^{\prime}}}\right) & =\sum\limits_{\left|J^{\prime}\right|=m+1}\sum\limits_{\substack{J\subset J^{\prime}\\
\left|J^{\prime}\right|=m
}
}P\left(\widehat{A_{J^{\prime}}}\right)\\
 & =\sum\limits_{\left|J^{\prime}\right|=m+1}\left(\begin{array}{c}
m+1\\
m
\end{array}\right)P\left(\widehat{A_{J^{\prime}}}\right)\\
 & =\left(\begin{array}{c}
m+1\\
m
\end{array}\right)S_{m+1}.
\end{align*}

Proceeding by the same way with other sums, we have the enounced formula.

\textbf{2.} The names of the individuals are numbered from 1 to $n,$
and we define $I=\left\{ 1,\dots,n\right\} .$ The outcome space $\Omega$
is taken as the set of bijections from $I$ onto itself, and the $\sigma-$algebra
$\mathcal{A}$ is taken as the power set of $\Omega.$ The probabilized
space is then given by $\left(\Omega,\mathcal{A},P\right)$ where
$P$ is the uniform probability.

We define 
\[
A_{i}=\left\{ f\in\Omega:\,f\left(i\right)=i\right\} 
\]
as the set of bijections having $i$ as a fixed point.

\textbf{(a) Probability that nobody draws their own name}

The subset 
\[
B_{0}=\bigcap\limits^{n}_{i=1}A^{c}_{i}
\]
 corresponds to the event ``no individual draws their own name.''
Our goal is to compute $P\left(B_{0}\right).$ We use the equality
$\refpar{eq:ex5qu1a}$ with $A=\Omega.$ We obtain
\[
P\left(A_{i}\right)=\dfrac{\left|A_{i}\right|}{\left|\Omega\right|}=\dfrac{\left(n-1\right)!}{n!}.
\]
For $i\neq j,$
\[
P\left(A_{i}\cap A_{j}\right)=\dfrac{\left|A_{i}\cap A_{j}\right|}{\left|\Omega\right|}=\dfrac{\left(n-2\right)!}{n!},
\]
and so on...

Applying the equality $\refpar{eq:ex5qu1a}$ and simplifying, we obtain\boxeq{
\[
P\left(B_{0}\right)=\dfrac{1}{2!}-\dfrac{1}{3!}+\dfrac{1}{4!}+\dots+\left(-1\right)^{n}\dfrac{1}{n!}.
\]
}

As $n$ tends to infinity, we observe that $P\left(B_{0}\right)$
tends to $\textrm{e}^{-1}\approx0.368.$

\textbf{(b) Probability that exactly $m$ persons draw a card with
their name}

The event $B_{m}$ corresponds to the event where ``exactly $m$
individuals draw a card with their own name.'' We seek to compute
$P\left(B_{m}\right).$ We use the equality $\refpar{eq:ex5qu1b}.$

Let $J$ be a subset of $I$ with cardinality $r.$ Then
\[
P\left(\widehat{A_{J}}\right)=\dfrac{\left(n-r\right)!}{n!}
\]

Thus, we compute
\[
S_{r}=\left(\begin{array}{c}
n\\
r
\end{array}\right)\dfrac{\left(n-r\right)!}{n!}=\dfrac{1}{r!}.
\]
Substituting into the equality $\refpar{eq:ex5qu1b}$ and simplifying,
we obtain\boxeq{
\[
P\left(B_{m}\right)=\dfrac{1}{m!}\left(\dfrac{1}{2!}-\dfrac{1}{3!}+\dfrac{1}{4!}-\dots+\left(-1\right)^{n-m}\dfrac{1}{\left(n-m\right)!}\right).
\]
}

As $n$ tends to infinity, we observe that: $P\left(B_{m}\right)$
tends to $\dfrac{\textrm{e}^{-1}}{m!}.$ 

Additionally, for a fixed $n,$ the number of bijections from $I$
onto itself with exactly $m$ fixed points is $n!P\left(B_{m}\right).$

\end{solution}

\begin{solution}{}{}

We encode ``tails'' as 1 and ``heads'' as 0. An outcome of this
game is thus a sequence of length $N$ consisting of 0s or 1s. Since
the coin is fair, we define the probabilized space $\left(\Omega,\mathcal{A},P\right),$
where $\Omega=\left\{ 0,1\right\} ^{N},$ $\mathcal{A}$ is the power
set of $\Omega,$ and $P$ is the uniform probability on the probabilizable
space $\left(\Omega,\mathcal{A}\right).$

For $\omega\in\Omega,$ we write $\omega=\left(\omega_{1},...,\omega_{N}\right)$
and define the $j-$th projection $X_{j}$ of $\omega$ on $\Omega$
as $X_{j}\left(\omega\right)=\omega_{j}.$

We then define the random variable $T$ by
\[
\forall\omega\in\Omega,\,\,\,\,T\left(\omega\right)=\inf\left\{ j\in\mathbb{N}^{\ast}:\,X_{j}\left(\omega\right)=1\right\} ,
\]
by conveining here that $\inf\emptyset=N+1.$

The random variable $T$ takes value in the set of integers with values
between 1 and $N+1;$ the integer $T\left(\omega\right)$ represents
the number of draws until the first appearance of ``tails'' in the
sequence $\omega.$ The law of $T$ is fully determined by the probabilities
$P\left(T=j\right)$ for $1\leqslant j\leqslant N+1.$

We have
\[
\left(T=1\right)=\left(X_{1}=1\right)
\]
 and
\[
\left(X_{1}=1\right)=\left\{ \left(1,\omega_{2},\dots,\omega_{N}\right):\,\left(\omega_{2},\dots,\omega_{N}\right)\in\left\{ 0,1\right\} ^{N-1}\right\} ,
\]

Thus, since $\left|\left(X_{1}=1\right)\right|=2^{N-1}$ and $\left|\Omega\right|=2^{N}$
and, that $P$ is the uniform probability, we obtain
\[
P\left(T=1\right)=\dfrac{\left|\left(X_{1}=1\right)\right|}{\left|\Omega\right|}=\dfrac{1}{2}.
\]

For $2\leqslant j\leqslant N,$
\[
\left(T=j\right)=\left[\bigcap^{j-1}_{i=1}\left(X_{i}=0\right)\right]\cap\left(X_{j}=1\right).
\]
 This last can also be described as
\[
\left(T=j\right)=\left\{ \left(0,...,0,1,\omega_{j+1},\dots,\omega_{N}\right):\,\left(\omega_{j+1},\dots,\omega_{N}\right)\in\left\{ 0,1\right\} ^{N-j}\right\} .
\]

Thus, since $\left|\left(T=j\right)\right|=2^{N-j},$\boxeq{
\[
P\left(T=j\right)=\dfrac{1}{2^{j}}.
\]
}

Finally, since $\left(T=N+1\right)=\left\{ \left(0,0,\dots,0\right)\right\} ,$
$\left|\left(T=N+1\right)\right|=1,$ and \boxeq{
\[
P\left(T=N+1\right)=\dfrac{1}{2^{N}}.
\]
}

As a verification, we check that
\[
\sum\limits^{N+1}_{j=1}P\left(T=j\right)=1.
\]

\end{solution}

\chapter{Real Number Summable Families}\label{chap:Real-Number-Summable}

\begin{objective}{}{}

Chapter \ref{chap:Real-Number-Summable} explores \textbf{summable
families of real numbers}. It highlights the versatility and utility
of using summable families---especially the advantage of using arbitrary
indexing sets---and emphasizes that for nonnegative real numbers
the sum is always well-defined, even when the family is non-summable.
\begin{itemize}
\item Section \ref{sec:Sum-of-a-Family} begins with a review on the \textbf{extended
real line $\overline{\mathbb{R}},$} then introduces the \textbf{sum
of families of nonnegative real numbers}. It defines \textbf{summability}
and provides conditions under which a family is summable, including
connections between summable families and classical series.
\item Section \ref{sec:Arithmetic-in-rbar.} addresses \textbf{arithmetic
operations in $\overline{\mathbb{R}},$} then defines the \textbf{sum
of a family of nonnegative elements in $\overline{\mathbb{R}}.$}
It presents fundamental properties, including \textbf{linearity},
\textbf{order preservation}, and \textbf{sum of disjoint index parts}.
It concludes with the important notion of \textbf{packet summation}.
\item Section \ref{sec:Sum-of-any-sign} extends the discussion to families
of real numbers of any sign. It establishes that \textbf{summability}
in this context corresponds with \textbf{absolute summability}, and
presents a theorem connecting \textbf{summable families} with \textbf{absolute
convergent series}. This section explores the \textbf{set additivity
function} associated to the sum of two summable families, examines
the \textbf{linearity} of the summability and introduces a theorem
ensuring \textbf{packet summation} for real-valued families. It also
presents a version of \textbf{Fubini theorem} for families, facilitating
the manipulation of \textbf{double-indexed family} sums. Finally,
it addresses \textbf{products} and \textbf{factorization} of terms
within summable families.
\end{itemize}
\end{objective}

\section*{Introduction}

Infinite sums---i.e. sums of an infinite number of real numbers---are
frequently used in Probability Theory. One of the simplest way to
handle this concept---as illustrated by the formula $\sum_{i\in I}x_{i}$---is
through the notion of a \textbf{summable family}\index{summable family}\mindex{family ! summable},
which is the focus of this Chapter. 

The definition of the sum of a family of real numbers differs from
that of a \textbf{convergent series}\index{convergent series}\mindex{series ! convergent}---something
that the reader is supposed to be familiar with---since no specific
condition, particularly regarding the \textbf{countability} of the
index set $I$, is required. Moreover, the definition of $\sum_{i\in I}x_{i}$
does not depend on any particular ordering of the set $I.$ There
is also a practical distinction: under certain conditions, $\sum_{i\in I}x_{i}$
can be treated as a \textbf{finite sum}\mindex{sum ! finite}\index{finite sum},
using methods such as \textbf{packet summation}---as outlined in
Theorem \ref{Theorem Packet-summation any real}.

However, using \textbf{summable families} instead of \textbf{convergent
series} does not lead to significantly new results. In practice, the
index set is still \textbf{countable}, and the concept of summing
a family of real numbers is essentially equivalent to summing the
terms of an absolutely convergent series.

\section{Sum of a Family of Nonnegative Real Numbers}\label{sec:Sum-of-a-Family}

\subsection{The Extended Real Line $\overline{\mathbb{R}}$}

Recall that any \textbf{non-empty} and \textbf{upper-bounded} subset
$A$ of $\mathbb{R}$---that is, a subset with an \textbf{upper bound}---has
a \textbf{least upper bound}\index{least upper bound}, also know
as the \textbf{supremum}\index{supremum}, denoted by $\text{sup}A.$
By definition, the real number $\text{sup}A$ is the \textbf{smallest
value among all the upper bounds} of $A.$ The existence of $\text{sup}A$
is guaranteed by a theorem\footnote{Tr.N.: This theorem states that all non empty and upper-bounded subset
of $\mathbb{R}$ admits a supremum. }. 

If the subset $A$ is not upper-bounded, we set
\[
\text{sup}A=+\infty.
\]

Similarly, the \textbf{\index{greatest lower bound}greatest lower
bound}, or \textbf{infinimum}\index{infinimum}, of $A$, denoted
by $\text{inf}A$, is the largest value among all lower bounds of
$A.$

In this framework, the real line $\mathbb{R}$ can be \textbf{usefully
extended} by adding to $\mathbb{R}$ two elements, $+\infty$ and
$-\infty.$ The resulting set, 
\[
\overline{\mathbb{R}}=\mathbb{R}\cup\left\{ +\infty,-\infty\right\} ,
\]
 is known as the \textbf{\index{extended real line}\mindex{real line! extended}extended
real line}. The standard order relation on $\mathbb{R}$ is naturally
extended to $\overline{\mathbb{R}}$ by stipulating that, for every
$x\in\mathbb{R},$ 
\[
\begin{array}{ccccc}
-\infty\leqslant x &  & \text{and} &  & x\leqslant+\infty.\end{array}
\]

This extension simplifies certain analysis compared to working solely
in $\mathbb{R}.$ In $\overline{\mathbb{R}},$ every subset---including
the empty set $\emptyset$---has both a supremum and an infinum.
For instance, if $+\infty$ is an element of $A,$ or if $A$ is a
subset of $\mathbb{R}$ that is\textbf{ not upper-bounded}, then $\text{sup}A=+\infty.$
Indeed in the latest, $+\infty$ is an upper bound of $A,$ and no
smaller element of $\overline{\mathbb{R}}$ can serve as an upper
bound of $A.$ The introduction of $\overline{\mathbb{R}}$ thus gives
subtantive meaning to the relation $\text{sup}A=+\infty,$ which was
previously treated as a mere notational convention.

Occasionally, and somewhat informally, all elements of $\overline{\mathbb{R}}$
are referred to as ``\textbf{real numbers}'', while the elements
of $\mathbb{R}$ specifically are called ``\textbf{\index{finite real numbers}\mindex{real number! finite}finite
real numbers}''. Under this convention, \textbf{real random variables}
are understood to take values in $\overline{\mathbb{R}},$ and \textbf{real
functions} or \textbf{numerical functions} are those with values in
$\overline{\mathbb{R}}.$ \textbf{Finite real random variables} and
\textbf{finite numerical functions} are restricted to values in $\mathbb{R}.$
This language, while an abuse of terminology, is always clarified
in context.

\subsection{Summable Families of Nonnegative Real Numbers}

We have already encountered the notion of a family of elements in
a set $E,$ indexed by a set $I.$ In general, this concept is distinct
from that of a mapping from $I$ to $E.$ The difference lies only
in notation and perspective: when we speak of a family $\left(x_{i}\right)_{i\in I},$
the focus is on the elements $x_{i},$ which we imagine as being ``located''
or ``enumerated'' by the elements of $I.$ For instance, a \textbf{sequence}
is simply a \textbf{family indexed by $\mathbb{N}$}---or by a \textbf{finite
interval of integers}.

The sum of a \textbf{family of nonnegative real numbers}---where
``nonnegative'' means values greater than or equal to zero---is
naturally defined as an extension of the notion of the sum for finite
families.

\begin{definition}{Sum of a Family of Nonnegative Real Numbers}{}\label{definition sum of non-negative real numbers}

Let $\left(x_{i}\right)_{i\in I}$ be a family of nonnegative real
numbers.

The \textbf{\mindex{family ! sum}sum\index{sum of the family} of
the family $\left(x_{i}\right)_{i\in I}$} is defined as the element
of $\left[0;+\infty\right]$\footnotemark. The sum is denoted as\boxeq{
\[
\sum_{i\in I}x_{i}
\]
}and is given by\boxeq{
\[
\sum_{i\in I}x_{i}=\text{sup}\left\{ \sum_{i\in F}x_{i}:\,F\,\text{finite part of }I\right\} .
\]
}

\end{definition}

\footnotetext{The set $\left[0;+\infty\right]$ can be written interchangeably
as $\overline{\mathbb{R}}^{+},$ representing the set of nonnegative
elements of $\overline{\mathbb{R}}$---that is, $\mathbb{R}^{+}\cup\left\{ +\infty\right\} .$}

Thus, the sum of a family of \textbf{nonnegative real numbers} is
always \textbf{well-defined}. It is either a finite real number or
$+\infty.$

\begin{definition}{Summable and Non-Summable Families}{}

A family $\left(x_{i}\right)_{i\in I}$ of nonnegative real numbers
is called \textbf{summable}\index{summable}\mindex{family ! summable}
if its sum $\sum_{i\in I}x_{i}$ is finite.

If $\sum_{i\in I}x_{i}=+\infty,$ then the family $\left(x_{i}\right)_{i\in I}$
is said to be \textbf{non-summable}\index{non-summable}\mindex{family ! non-summable},
even though the sum is still assigned a value.

\end{definition}

\subsection{Comparison with Nonnegative Term Series}

\begin{proposition}{Condition of Summability}{}\label{proposition:countable indices}

Let $\left(x_{i}\right)_{i\in I}$ be a summable family of nonnegative
real numbers.

Then, the set of indices $i$ for which $x_{i}>0$ is countable.

\end{proposition}

In other words, $x_{i}=0$ for any but a countable set of indices.

\begin{proof}{}{}

Let $J$ be the set of indices $i\in I$ such that $x_{i}>0.$ Define
the set $J$ as the union of the sets $J_{n}$ for $n\geqslant1,$
where
\[
J_{n}=\left\{ i\in I:\,x_{i}\geqslant\dfrac{1}{n}\right\} .
\]

We now demonstrate that each $J_{n}$ is finite. Suppose, for contradiction,
that the set $J_{n}$ is infinite. Then, for every $k,$ there exists
a subset $K$ of $J_{n}$ with $k$ elements, yielding
\[
\sum_{i\in I}x_{i}\geqslant\sum_{i\in K}x_{i}\geqslant\dfrac{k}{n}.
\]

Since this holds for every $k\geqslant1,$ it follows that $\sum_{i\in I}x_{i}=+\infty,$
contradicting the assumption that the family $\left(x_{i}\right)_{i\in I}$
is summable.

Therefore, each $J_{n}$ must be finite, and since $J$ is a countable
union of finite sets, it is countable.

\end{proof}

\begin{proposition}{Family Sum and Enumeration}{fam_sum_enum}\label{Prop:countable enumeration}

Let $J$ be a countable infinite set, and let $\left(j_{k}\right)_{k\in\mathbb{N}}$
be an enumeration\footnotemark of $J.$ Then\boxeq{
\[
\sum_{j\in J}x_{j}=\sum^{+\infty}_{k=0}x_{j_{k}}.
\]
}

\end{proposition}

\footnotetext{i.e. a bijection from $\mathbb{N}$ onto $J$}

As usual, the notation $\sum^{+\infty}_{k=0}x_{j_{k}}$ in the second
part of the equality refers to the limits of the partial sums $\sum^{N}_{k=0}x_{j_{k}}$
as $N$ tends to $+\infty.$ Since the terms $x_{j_{k}}$ are nonnegative,
the sequence of partial sums is nondecreasing, ensuring the existence
of the limit. 
\begin{itemize}
\item Either the limit of the partial sums is finite: we say in this case
that ``the series with general term $x_{j_{k}}$ converges'', or
equivalently, we say that ``the series $\sum x_{j_{k}}$ converges'';
\item Or the limit of the partial sums is infinite: in this case we say
that ``the series diverges''.
\end{itemize}
As a corollary of Proposition \ref{Prop:countable enumeration}, for
the family of nonnegative terms $\left(x_{j}\right)_{j\in J}$ to
be summable, it is both necessary and sufficient that the series $\sum x_{j_{k}}$
converges.

\begin{proof}{}{}

For every integer $N,$ define the finite subset $J_{N}=\left\{ j_{0},j_{1},...,j_{N}\right\} $
of $J.$ Then
\[
\sum^{N}_{k=0}x_{j_{k}}=\sum_{j\in J_{N}}x_{j}\leqslant\sum_{j\in J}x_{j}.
\]

Taking the limit as $N$ tends to $+\infty,$ we obtain
\begin{equation}
\sum^{+\infty}_{k=0}x_{j_{k}}\leqslant\sum_{j\in J}x_{j}.\label{eq:majorant_of_the_sum}
\end{equation}

If $\sum_{j\in J}x_{j}<+\infty,$ then the inequality $\refpar{eq:majorant_of_the_sum}$
remains valid upon taking the limit in $\mathbb{R}.$ If $\sum_{j\in J}x_{j}=+\infty,$
the inequality $\refpar{eq:majorant_of_the_sum}$ holds trivially.

Conversely, let $F$ be a finite subset of $J.$ There exists a unique
$N,$ such that $F\subset J_{N}$ and such that for $k<N,$ $F$ is
not a subset of $J_{k}$ yielding
\[
\sum_{j\in F}x_{j}\leqslant\sum_{j\in J_{N}}x_{j}=\sum^{N}_{k=0}x_{j_{k}}.
\]

Thus,
\begin{equation}
\sum_{j\in F}x_{j}\leqslant\sum^{+\infty}_{k=0}x_{j_{k}}.
\end{equation}

This inequality, combined with $\refpar{eq:majorant_of_the_sum},$
established the desired result.

\end{proof}

\section{Arithmetic in $\overline{\mathbb{R}}$. Sum of a Family of $\overline{\mathbb{R}}$
Elements}\label{sec:Arithmetic-in-rbar.}

\subsection{Arithmetic in $\overline{\mathbb{R}}$ }

To clearly state the main properties of summation for a family of
nonnegative real numbers---especially the key property of summation
by packets---see Proposition \ref{prop:packet-summation non-negative reals}---we
extend our framework to include elements of a family that take the
value $+\infty.$ This necessitates defining the sum of $+\infty$
and a nonnegative real number.

We extend the basic arithmetic operations (addition, subtraction,
and multiplication) to include cases where one of the terms or factors
is infinite, with one important exception: the value of the expression
$\left(+\infty\right)+\left(-\infty\right)$ is not defined, nor are
those reducible to it, such as $\left(+\infty\right)-\left(+\infty\right)$
and $\left(-\infty\right)-\left(-\infty\right).$

In all other cases, the extension of operations follows naturally.
We define for every finite real number $x,$

\boxeq{
\begin{align*}
x+\left(+\infty\right) & =+\infty,\\
x+\left(-\infty\right) & =-\infty.
\end{align*}
}

We also define

\boxeq{
\begin{align*}
+\infty+\left(+\infty\right) & =+\infty,\\
-\infty+\left(-\infty\right) & =-\infty.
\end{align*}
}

Subtraction is defined as

\boxeq{
\[
x-y=x+\left(-y\right),
\]
}

and we set

\boxeq{
\[
-\left(+\infty\right)=-\infty.
\]
}

Multiplication---defined in every case for $x\ne0$---follows

\boxeq{
\begin{align*}
x\times\left(+\infty\right) & =\left\langle \text{sign}\left(x\right)\right\rangle \infty,\\
x\times\left(-\infty\right) & =\left\langle -\text{sign}\left(x\right)\right\rangle \infty.
\end{align*}
}

Here, $\left\langle \text{sign}\left(x\right)\right\rangle $ denotes
the sign of $x,$ which is $+$ if $x>0$ and $-$ if $x<0.$

And\footnote{Tr.N. It is important to note that the convention $0\times\left(\pm\infty\right)=0$
is typical of extended arithmetic, even though it may seem unusual
at first glance. The reasoning behind this definition is that, in
the extended real number framework, multiplying zero by any infinite
quantity is typically defined to result in zero. However, in other
contexts, such as standard real analysis, this operation leads to
indeterminate forms and requires a more nuanced approach, often resolved
through limit processes. In non-standard analysis or other specialized
mathematical frameworks, the treatment of $0\times\infty$ may vary
based on the underlying theory.}

\boxeq{
\[
0\times\left(\pm\infty\right)=0.
\]
}

As usual, we write $\infty$ instead of $+\infty,$ but it is understood
that the only operations that are undefined are those that can be
reduced to $\infty-\infty.$

Lastly, we define, as it is natural\boxeq{
\[
\left|+\infty\right|=\left|-\infty\right|=+\infty.
\]
}

\subsection{Sum of a Family of Nonnegative Elements of $\overline{\mathbb{R}}.$}

\begin{definition}{Sum of a Family of Nonnegative Elements of $\overline{\mathbb{R}}$}{}The
notion of sum of a family of nonnegative real numbers---possibly
infinites---is defined as in the case of nonnegative finite real
numbers by

\boxeq{
\[
\sum\limits_{i\in I}x_{i}=\text{sup}\left\{ \sum\limits_{i\in F}x_{i}:\,F\text{ finite subset of }I\right\} .
\]
}

The family $\left(x_{i}\right)_{i\in I}$ is said to be \textbf{summable\index{summable}}\mindex{family ! summable}
if $\sum\limits_{i\in I}x_{i}<+\infty.$

\end{definition}

The following result is straightforward.

\begin{proposition}{Finity of Elements of a Summable Family}{}

If the family $\left(x_{i}\right)_{i\in I}$ is summable, then for
every $i\in I,$ $x_{i}<+\infty.$

\end{proposition}

\begin{proof}{}{}

Assume there exists an index $i_{0}\in I$ such that $x_{i_{0}}=+\infty.$ 

Then
\[
\sum\limits_{i\in I}x_{i}\geqslant\sum\limits_{i\in\left\{ i_{0}\right\} }x_{i}=x_{i_{0}}=+\infty,
\]
 which contradicts the assumption of summability.

\end{proof}

\subsection{Properties of Summable Families of Nonnegative Real Numbers}

\begin{proposition}{Properties of Summable Families of Nonnegative Real Numbers}{}\label{Proposition:Properties-of-summable-non-neg-reals}Let
$\left(x_{i}\right)_{i\in I}$ and $\left(y_{i}\right)_{i\in I}$
be two families of nonnegative real numbers (possibly infinite).

(i) If for every $i\in I,$ $x_{i}\leqslant y_{i},$ then

\boxeq{
\[
\sum\limits_{i\in I}x_{i}\leqslant\sum\limits_{i\in I}y_{i}.
\]
}

(ii) If $a$ is a nonnegative element of $\overline{\mathbb{R}},$
then

\boxeq{
\[
a\sum\limits_{i\in I}x_{i}=\sum\limits_{i\in I}ax_{i}.
\]
}

(iii) We have

\boxeq{
\[
\sum\limits_{i\in I}\left(x_{i}+y_{i}\right)=\sum\limits_{i\in I}x_{i}+\sum\limits_{i\in I}y_{i}.
\]
}

\end{proposition}

\begin{proof}{}{}

\textbf{General principle:} To prove an inequality of the type: $\text{sup}A\leqslant\beta,$
it suffices to show that $\beta$ is an upper bound of $A,$ as $\text{sup}A$
is the least of all upper bounds of $A.$ That is, we must verify
that $\alpha\leqslant\beta$ for every $\alpha\in A.$

In particular, to establish the inequality
\[
\sum\limits_{i\in I}x_{i}\leqslant\beta
\]
where the terms $x_{i}$ are nonnegative or equal to $+\infty$, it
is enough to show that for every finite subset $F$ of $I,$ we have
\[
\sum\limits_{i\in F}x_{i}\leqslant\beta.
\]

(i) For every finite subset $F$ of $I,$ since $x_{i}$ and $y_{i}$
are both nonnegative and that $x_{i}\leqslant y_{i}$, we have
\[
\sum\limits_{i\in F}x_{i}\leqslant\sum\limits_{i\in F}y_{i}\leqslant\sum\limits_{i\in I}y_{i}.
\]

By applying the general principle stated above, we obtain the desired
formula.

(ii) By definition of the multiplication in $\overline{\mathbb{R}}^{+},$
for every finite subset $F$ of $I,$ we have
\[
\sum\limits_{i\in F}ax_{i}=a\sum\limits_{i\in F}x_{i}.
\]
By definition of the sum of a nonnegative family, $\sum\limits_{i\in I}x_{i}$
is the supremum of the set constituted of all the sums on finite subsets
of $I,$ $\sum\limits_{i\in F}x_{i},$ thus
\[
\sum\limits_{i\in F}x_{i}\leqslant\sum\limits_{i\in I}x_{i}.
\]
Since $a$ is nonnegative, we have
\[
\sum\limits_{i\in F}ax_{i}=a\sum\limits_{i\in F}x_{i}\leqslant a\sum\limits_{i\in I}x_{i}.
\]
Applying the general principle, we conclude
\[
\sum\limits_{i\in I}ax_{i}\leqslant a\sum\limits_{i\in I}x_{i}.
\]

Following the same reasoning in the reverse inequality, we obtain
equality.

(iii) It is well known that ``the supremum of a sum is not always
the sum of the suprema.'' Again, to establish equality, we proceed
via a double inequality.

For every finite subset $F$ in $I,$ we have

\begin{align*}
\sum\limits_{i\in F}\left(x_{i}+y_{i}\right) & =\sum\limits_{i\in F}x_{i}+\sum\limits_{i\in F}y_{i}\\
 & \leqslant\sum\limits_{i\in I}x_{i}+\sum\limits_{i\in I}y_{i}.
\end{align*}
Thus, applying the general principle
\begin{equation}
\sum\limits_{i\in I}\left(x_{i}+y_{i}\right)\leqslant\sum\limits_{i\in I}x_{i}+\sum\limits_{i\in I}y_{i}.\label{eq:sum_(x_i+y_i)leqslant}
\end{equation}

For the reverse inequality, we proceed in two steps. 

Let $F$ and $G$ be two finite subsets of $I.$ Then
\[
\sum\limits_{i\in F}x_{i}+\sum\limits_{i\in G}y_{i}\leqslant\sum\limits_{i\in F\cup G}\left(x_{i}+y_{i}\right)\leqslant\sum\limits_{i\in I}\left(x_{i}+y_{i}\right),
\]
and thus
\[
\sum\limits_{i\in F}x_{i}+\sum\limits_{i\in G}y_{i}\leqslant\sum\limits_{i\in I}\left(x_{i}+y_{i}\right).
\]

Taking the supremum over all finite subsets $G$ of $I,$ we obtain
\[
\sum\limits_{i\in F}x_{i}+\sum\limits_{i\in I}y_{i}\leqslant\sum\limits_{i\in I}\left(x_{i}+y_{i}\right).
\]

Applying the same argument with respect to $F,$ we deduce
\begin{equation}
\sum\limits_{i\in I}x_{i}+\sum\limits_{i\in I}y_{i}\leqslant\sum\limits_{i\in I}\left(x_{i}+y_{i}\right).\label{eq:sum_x_i+sum_y_ileqslant}
\end{equation}

Combining Inequalities $\refpar{eq:sum_(x_i+y_i)leqslant}$ and $\refpar{eq:sum_x_i+sum_y_ileqslant}$,
we conclude the desired result.

\end{proof}

We start by recalling the definition of an indicator function.

\begin{definition}{Indicator Function}{}

Let $X$ be a set and $A$ a subset of $X.$ 

The \textbf{indicator function\index{indicator function}\mindex{function ! indicator}}
relatively to $A$ in $X$---also called the \textbf{characteristic
function}\index{characteristic function}\textbf{\mindex{function ! characteristic}}---is
a function defined on $X$ denoted $\boldsymbol{1}_{A}$ defined from
$X$ onto $\left\{ 0,1\right\} $ such that for every $x\in X,$\boxeq{
\[
\boldsymbol{1}_{A}\left(x\right)=\begin{cases}
1, & \text{if\,}x\in A,\\
0, & \text{otherwise.}
\end{cases}
\]
}

\end{definition}

The use of an indicator function $\boldsymbol{1}_{A}$ allows one
to express the sum relatively to a subset $A$ of $I$ as a sum on
$I,$ which proves to be useful in certain proofs.

\begin{proposition}{Sum of a Family Using An Indicator Function}{sum_with_indicator}

Let $\left(x_{i}\right)_{i\in I}$ be a family of nonnegative real
numbers---possibly infinite---and let $A$ be any subset of $I.$
Then\boxeq{
\begin{equation}
\sum\limits_{i\in A}x_{i}=\sum\limits_{i\in I}\boldsymbol{1}_{A}\left(i\right)x_{i}.\label{eq:indicator_function_and_summable}
\end{equation}
}

\end{proposition}

\begin{proof}{}{}

Let $F$ be a finite subset of $A.$ Since $\boldsymbol{1}_{A}\left(i\right)=1,$
for every $i\in F,$
\[
\sum\limits_{i\in F}x_{i}=\sum\limits_{i\in F}\boldsymbol{1}_{A}\left(i\right)x_{i}.
\]

\begin{itemize}
\item Since $F$ is also a finite subset of $I,$ it follows that
\[
\sum\limits_{i\in F}x_{i}\leqslant\sum\limits_{i\in I}\boldsymbol{1}_{A}\left(i\right)x_{i}.
\]
Taking the supremum over all finite subsets $F$ of $A,$ we deduce
\begin{equation}
\sum\limits_{i\in A}x_{i}\leqslant\sum\limits_{i\in I}\boldsymbol{1}_{A}\left(i\right)x_{i}.\label{eq:one_way_indicator}
\end{equation}
\item Conversely, let $G$ be a finite subset of $I.$ Then
\[
\sum\limits_{i\in G}\boldsymbol{1}_{A}\left(i\right)x_{i}=\sum\limits_{i\in G\cap A}x_{i}.
\]
Since $G\cap A$ is a finite subset of $A,$
\[
\sum_{i\in G\cap A}x_{i}\leqslant\sum_{i\in A}x_{i}.
\]
Thus
\[
\sum\limits_{i\in G}\boldsymbol{1}_{A}\left(i\right)x_{i}\leqslant\sum\limits_{i\in A}x_{i}.
\]
The preceding inequality being verified for every finite subset $G$
of $I,$ taking the supremum over all finite subsets $G$ of $I,$
we obtain
\begin{equation}
\sum\limits_{i\in I}\boldsymbol{1}_{A}\left(i\right)x_{i}\leqslant\sum\limits_{i\in A}x_{i}.\label{eq:other_way_indicator}
\end{equation}
\item Since we have established both inequalities $\refpar{eq:one_way_indicator}$
and $\refpar{eq:other_way_indicator}$, we conclude with the desired
equality $\refpar{eq:indicator_function_and_summable}.$
\end{itemize}
\end{proof}

\begin{proposition}{Comparison of Sums. Sum of Disjoint Parts}{}\label{Proposition inequality sum}

Let $\left(x_{i}\right)_{i\in I}$ be a family of nonnegative real
numbers (possibly infinite). Let $A$ and $B$ be two subsets of $I.$

(i) If $A\subset B,$ then\boxeq{
\[
\sum\limits_{i\in A}x_{i}\leqslant\sum\limits_{i\in B}x_{i}.
\]
}

(ii) If $A$ and $B$ are disjoints, then\boxeq{
\[
\sum\limits_{i\in A\uplus B}x_{i}=\sum\limits_{i\in A}x_{i}+\sum\limits_{i\in B}x_{i}.
\]
}

\end{proposition}

Property (ii) states that the mapping $A\mapsto\sum\limits_{i\in A}x_{i}$
is a \textbf{\index{set-additive function}\mindex{function ! set-additive}set-additive
function}. This property extends naturally to the case of a finite
collection of pairwise disjoint subsets of I.

\begin{proof}{}{}(i) Since every finite subset $F$ of $A$ is also
a finite subset of $B,$
\[
\sum\limits_{i\in F}x_{i}\leqslant\sum\limits_{i\in B}x_{i}.
\]

Taking the supremum over all such $F$ finite subset of $A$ yields
the enounced result.

(ii) Since $A$ and $B$ are disjoint, we have
\[
\boldsymbol{1}_{A\uplus B}=\boldsymbol{1}_{A}+\boldsymbol{1}_{B}.
\]
Applying the previous result in Proposition $\ref{pr:sum_with_indicator},$
we obtain
\[
\sum\limits_{i\in A\uplus B}x_{i}=\sum\limits_{i\in I}\boldsymbol{1}_{A\uplus B}\left(i\right)x_{i}=\sum\limits_{i\in I}\boldsymbol{1}_{A}\left(i\right)x_{i}+\sum\limits_{i\in I}\boldsymbol{1}_{B}\left(i\right)x_{i}=\sum\limits_{i\in A}x_{i}+\sum\limits_{i\in B}x_{i}.
\]

\end{proof}

The following proposition is evident in the case of finite sets; in
this context, it can be viewed as a generalization of the associativity
property of finite sums. Along with the corresponding proposition
in the case of real-numbers of any sign, it plays a fundamental role
in the developments that follow in this book.

It is worth noting that no analogous result can be easily formulated
for series. This is the primary motivation for introducing the concept
of a summable family. 

\begin{proposition}{Packet Summation}{}\textbf{\label{prop:packet-summation non-negative reals}}Let
$\left(x_{i}\right)_{i\in I}$ be a family of nonnegative real numbers
(possibly infinite) and let $\left(A_{j}\right)_{j\in J}$ be a partition
of $I.$ Then, we have\boxeq{
\[
\sum\limits_{i\in I}x_{i}=\sum\limits_{j\in J}\left(\sum\limits_{i\in A_{j}}x_{i}\right).
\]
}

This method of summation is called the \textbf{\index{packet summation}packet
summation}.

\end{proposition}

\begin{proof}{}{}For every finite subset $F$ of $I,$
\[
F=\biguplus\limits_{j\in J}\left(A_{j}\cap F\right).
\]

Since the set $F$ is finite, it does not intersect more than a finite
number of sets $A_{j};$ that is, there exists a finite subset $J_{F}$
of $J$ such that
\[
F=\biguplus\limits_{j\in J_{F}}\left(A_{j}\cap F\right).
\]

Using the associativity of finite sums, we obtain
\[
\sum\limits_{i\in F}x_{i}=\sum\limits_{j\in J_{F}}\left(\sum\limits_{i\in A_{j}\cap F}x_{i}\right).
\]

Since $\left(A_{j}\cap F\right)\subset A_{j},$ it follows that
\[
\sum\limits_{i\in F}x_{i}\leqslant\sum\limits_{j\in J_{F}}\left(\sum\limits_{i\in A_{j}}x_{i}\right)
\]

By replacing in the previous inequality $J_{F}$ with $J,$ and then
$F$ with $I$---or, alternatively by appealing to well-known properties
of the supremum, we obtain
\[
\sum\limits_{i\in I}x_{i}\leqslant\sum\limits_{j\in J}\left(\sum\limits_{i\in A_{j}}x_{i}\right).
\]

Conversely, for every finite subset $J_{0}$ of $J,$ the simple additivity---Proposition
\ref{Proposition inequality sum} (ii)---implies
\[
\sum\limits_{j\in J_{0}}\left(\sum\limits_{i\in A_{j}}x_{i}\right)=\sum\limits_{i\in\biguplus\limits_{j\in J_{0}}A_{j}}x_{i}\leqslant\sum\limits_{i\in I}x_{i}
\]

Taking the supremum over all such finite subsets $J_{0}$ yields
\[
\sum\limits_{j\in J}\left(\sum\limits_{i\in A_{j}}x_{i}\right)\leqslant\sum\limits_{i\in I}x_{i}.
\]

Since both inequalities hold, we conclude
\[
\sum\limits_{i\in I}x_{i}=\sum\limits_{j\in J}\left(\sum\limits_{i\in A_{j}}x_{i}\right).
\]

\end{proof}

\section{Sum of a Family of Real Numbers of Any Sign}\label{sec:Sum-of-any-sign}

We now consider families of real numbers that may take any sign. The
notion of a summable family and the corresponding sum are defined
in a straightforward manner by separating the positive and negative
parts. However, in contrast to the case of nonnegative real numbers,
we can no longer allow infinite values, as this could lead to the
ill-defined operation $\infty-\infty.$

\subsection{Summable Families of Real Numbers of Any Sign}

To begin, we recall the standard decomposition of a real number $x.$

Given any $x\in\mathbb{R},$ we define\boxeq{
\[
\begin{array}{ccccc}
x^{+}=\begin{cases}
x, & \text{if\,}x\geqslant0,\\
0, & \text{if\,}x<0,
\end{cases} & \,\,\,\, & \text{and} & \,\,\,\, & x^{-}=\begin{cases}
0, & \text{if\,}x\geqslant0,\\
\left|x\right|, & \text{if\,}x<0.
\end{cases}\end{array}
\]
}

The following fundamental identities are then satisfied
\[
\begin{array}{ccc}
x=x^{+}-x^{-}, &  & \left|x\right|=x^{+}+x^{-}.\end{array}
\]

The quantity $x^{+}$ is called the \textbf{nonnegative part\index{nonnegative part}\mindex{part ! nonnegative}}
of $x,$ while $x^{-}$ is referred to as its \textbf{nonpositive
part}\index{nonpositive part}\textbf{\mindex{part ! nonpositive}}.
This decomposition is often used to express a function $f$ as the
difference of two nonnegative functions $f^{+}$ and $f^{-},$ where
for every $x$ in the domain of $f,$ 
\[
\begin{array}{ccccc}
f^{+}\left(x\right)=\left(f\left(x\right)\right)^{+} & \,\,\,\, & \text{and} & \,\,\,\, & f^{-}\left(x\right)=\left(f\left(x\right)\right)^{-}.\end{array}
\]

These can also be expressed as
\[
\begin{array}{ccccc}
f^{+}=\text{sup}\left(f,0\right) & \,\,\,\, & \text{and} & \,\,\,\, & f^{-}=\text{sup}\left(-f,0\right).\end{array}
\]

We now extend these notions to families of real numbers. Specifically,
any family $\left(x_{i}\right)_{i\in I}$ of real numbers can be expressed
as the difference of two families of nonnegative real numbers. This
allows us to extend the concept of a summable family to include families
of real numbers of arbitrary sign.

\begin{definition}{Sum of Family of Real Numbers of Any Sign}{}A
family of real numbers $\left(x_{i}\right)_{i\in I}$ is said to be
\textbf{summable\index{summable}\mindex{family ! summable}} if
\[
\begin{array}{ccccc}
\sum\limits_{i\in I}x^{+}_{i}<+\infty & \,\,\,\, & \text{and} & \,\,\,\, & \sum\limits_{i\in I}x^{-}_{i}<+\infty.\end{array}
\]
 If $\left(x_{i}\right)_{i\in I}$ is summable, then its \textbf{sum}\index{sum}\textbf{\mindex{family ! sum}},
denoted $\sum_{i\in I}x_{i},$ is the real number defined as\boxeq{
\[
\sum_{i\in I}x_{i}=\sum_{i\in I}x^{+}_{i}-\sum_{i\in I}x^{-}_{i}.
\]
}

\end{definition}

Unlike families of nonnegative real numbers, the sum of a family $\left(x_{i}\right)_{i\in I}$
is not always well-defined. However, in the case $\left(x_{i}\right)_{i\in I}$
is a family of nonnegative real numbers---as in Definition \ref{definition sum of non-negative real numbers}---this
new definition coincides with the previously introduced notion of
sum. 

\begin{lemma}{Summability of the Absolute Value of Elements of a Summable Family}{}\label{lemma: summability of absolute value}A
family $\left(x_{i}\right)_{i\in I}$ is summable if and only if the
family $\left(\left|x_{i}\right|\right)_{i\in I}$ is summable, or,
equivalently,
\[
\sum_{i\in I}\left|x_{i}\right|<+\infty.
\]

\end{lemma}

\begin{proof}{}{}

Since $\left|x_{i}\right|=x^{+}_{i}+x^{-}_{i},$
\[
\sum_{i\in I}\left|x_{i}\right|=\sum_{i\in I}x^{+}_{i}+\sum_{i\in I}x^{-}_{i}.
\]

If $\left(x_{i}\right)_{i\in I}$ is summable, then
\[
\begin{array}{ccc}
\sum_{i\in I}x^{+}_{i}<+\infty, & \,\,\,\, & \sum_{i\in I}x^{-}_{i}<+\infty.\end{array}
\]

Summing these inequalities gives
\[
\sum_{i\in I}\left|x_{i}\right|<+\infty.
\]

Conversely, if $\sum_{i\in I}\left|x_{i}\right|<+\infty,$ then necessarily
we have
\[
\begin{array}{ccccc}
\sum_{i\in I}x^{+}_{i}<+\infty & \,\,\,\, & \text{and} & \,\,\,\, & \sum_{i\in I}x^{-}_{i}<+\infty.\end{array}
\]

Thus, the family $\left(x_{i}\right)_{i\in I}$ is summable.

\end{proof}

Hence, there is no distinction between the notions of \textbf{summability}
and \textbf{absolute summability}.

\subsection{Summable Families and Series}

The following theorem establishes the link between the concept of
a \textbf{summable family} and that of an \textbf{absolutely convergent
series}.

\begin{theorem}{Summable Family and Absolutely Convergent Series}{}For
a family of real numbers $\left(x_{i}\right)_{i\in I}$ to be summable,
it is necessary that the set $J=\left\{ i\in I:\,x_{i}\neq0\right\} $
is \textbf{countable}.

Suppose this condition is satisfied, and assume that $J$ is infinite.
Let $\left(j_{n}\right)_{n\in\mathbb{N}}$ be any enumeration of $J.$
Then, the family $\left(x_{i}\right)_{i\in I}$ is summable\textbf{
if and only if} the series $\sum^{+\infty}_{n=0}x_{j_{n}}$ is absolutely
convergent. 

Moreover, if this last condition holds, \boxeq{
\begin{equation}
\sum_{i\in I}x_{i}=\sum^{+\infty}_{n=0}x_{j_{n}}.\label{eq:link_summable_family_abs_conv_series}
\end{equation}
}

\end{theorem}

\begin{proof}{}{}Apart from Equality $\refpar{eq:link_summable_family_abs_conv_series},$
the previous assertions result immediately from Lemma \ref{lemma: summability of absolute value}
and Propositions \ref{proposition:countable indices} and \ref{Prop:countable enumeration}.

To establish Equality $\refpar{eq:link_summable_family_abs_conv_series},$
we use Proposition \ref{proposition:countable indices}. 

We have on the one hand
\[
\begin{array}{ccccc}
\sum\limits_{i\in I}x^{+}_{i}=\sum\limits^{+\infty}_{n=0}x^{+}_{j_{n}} & \,\,\,\, & \text{and} & \,\,\,\, & \sum\limits_{i\in I}x^{-}_{i}=\sum\limits^{+\infty}_{n=0}x^{-}_{j_{n}},\end{array}
\]
and on the other hand, by definition,
\[
\sum\limits_{i\in I}x_{i}=\sum\limits_{i\in I}x^{+}_{i}-\sum\limits_{i\in I}x^{-}_{i}.
\]

Since each term satisfies $x_{j_{n}}=x^{+}_{j_{n}}-x^{-}_{j_{n}}$
and, since we have already established the convergence of the two
nonnegative series on the right-hand side, it follows that
\[
\sum\limits^{+\infty}_{n=0}x_{j_{n}}=\sum\limits^{+\infty}_{n=0}x^{+}_{j_{n}}-\sum\limits^{+\infty}_{n=0}x^{-}_{j_{n}}
\]
 which converges, proving the desired equality.

\end{proof}

\begin{remark}{}{}

As an immediate corollary of the preceding theorem, we recover the
well-known result concerning the \textbf{absolute convergence of a
series} $\sum\limits^{+\infty}_{n=0}x_{n}.$ In particular, when a
series is \textbf{absolutely convergent}, its sum is well-defined
and \textbf{does not depend on the order of its terms}. 

Throughout this book, whenever a series is absolutely convergent,
we will designate indifferrently its sum by using either the \textbf{series
notation} $\sum\limits^{+\infty}_{n=0}x_{n}$ or the ``summable family''
notation $\sum\limits_{n\in\mathbb{N}}x_{n}$ interchangeably. The
same convention will be applied to series with \textbf{nonnegative
terms}, even when their sum is infinite.

Furthermore, the $\sigma-$additivity property of probability can
be immediately reformulated in terms of summable families. If $I$
is countable and if $\left(A_{i}\right)_{i\in I}$ is a family of
pairwise disjoint events, then\boxeq{
\[
P\left(\biguplus\limits_{i\in I}A_{i}\right)=\sum\limits_{i\in I}P\left(A_{i}\right).
\]
}

\end{remark}

\subsection{Summable Families Properties}

\begin{lemma}{Restriction of a Summable Family. Set Additive Function}{set_add_function}\label{Lemma:Summable-families-restriction}Let
$\left(x_{i}\right)_{i\in I}$ be a summable family of real numbers.

(i) For every subset $A$ of $I,$ the restricted family to the subset
$A$, $\left(x_{i}\right)_{i\in A}$, remains summable.

(ii) The mapping $A\mapsto\sum_{i\in A}x_{i}$ defines a \textbf{set-additive
function}.

\end{lemma}

\begin{proof}{}{}(i) Since $A\subset I,$ this follows directly from
the inequality 
\[
\sum\limits_{i\in A}\left|x_{i}\right|\leqslant\sum\limits_{i\in I}\left|x_{i}\right|
\]
---Proposition \ref{Proposition inequality sum}---and Lemma \ref{lemma: summability of absolute value}.

(Tr.N) Lemma \ref{lemma: summability of absolute value} ensures that
as $\left(x_{i}\right)_{i\in I}$ is a summable family of real numbers,
$\sum\limits_{i\in I}\left|x_{i}\right|<+\infty$ and by consequence
$\sum\limits_{i\in A}\left|x_{i}\right|<+\infty$and, the fact that
$\left(x_{i}\right)_{i\in A}$ remains summable.

(ii) The result follows from the identity
\[
\sum_{i\in A}x_{i}=\sum_{i\in A}x^{+}_{i}-\sum_{i\in A}x^{-}_{i}
\]
along with the analogous statement for the family of nonnegative real
numbers---Proposition \ref{Proposition inequality sum} (ii). 

(Tr.N) More precisely, let $A$ and $B$ be two disjoint sets, we
have
\begin{align*}
\sum_{i\in A\uplus B}x_{i} & =\sum_{i\in A\uplus B}x^{+}_{i}-\sum_{i\in A\uplus B}x^{-}_{i}\\
 & =\sum_{i\in A}x^{+}_{i}+\sum_{i\in B}x^{+}_{i}-\left(\sum_{i\in A}x^{-}_{i}+\sum_{i\in B}x^{-}_{i}\right)\\
 & =\sum_{i\in A}x^{+}_{i}-\sum_{i\in A}x^{-}_{i}+\sum_{i\in B}x^{+}_{i}-\sum_{i\in B}x^{-}_{i}\\
 & =\sum_{i\in A}x_{i}+\sum_{i\in B}x_{i}.
\end{align*}

\end{proof}

\begin{remark}{}{}The analogue of the assertion (i) of Lemma $\ref{lm:set_add_function}$
does not hold for convergent series that are not absolutely convergent.
For instance, consider the alternating harmonic series:
\[
1-\dfrac{1}{2}+\dfrac{1}{3}+\dots+\dfrac{\left(-1\right)^{n}}{n}+\dots.
\]

For families of real numbers of any sign, the linearity of summation
is not as straightforward as in the case of nonnegative real-valued
families. That is, unlike for nonnegative real-valued summable families,
the identity
\[
\sum\limits_{i\in I}\left(x_{i}+y_{i}\right)=\sum\limits_{i\in I}x_{i}+\sum\limits_{i\in I}y_{i}
\]
is not automatically guaranteed. 

To better understand why, note that writing $x_{i}=x^{+}_{i}-x^{-}_{i}$
amounts to splitting the index set $I$ in two subsets: the indices
such that $x_{i}>0$ and those such that $x_{i}<0$---disregarding
indices where $x_{i}=0.$ However, there is no reason why this decomposition
should align in the same way for the three families $\left(x_{i}\right)_{i\in I},$
$\left(y_{i}\right)_{i\in I}$ and $\left(x_{i}+y_{i}\right)_{i\in I},$
leading to potential inconsistencies.

\end{remark}

\begin{proposition}{Sum of Two Summable Families}{}Let $\left(x_{i}\right)_{i\in I}$
and $\left(y_{i}\right)_{i\in I}$ be two summable families of real
numbers. Then, the family $\left(x_{i}+y_{i}\right)_{i\in I}$ is
also summable, and we have\boxeq{
\begin{equation}
\sum\limits_{i\in I}\left(x_{i}+y_{i}\right)=\sum\limits_{i\in I}x_{i}+\sum\limits_{i\in I}y_{i}.\label{eq:sum_of_two_summable_families}
\end{equation}
}

\end{proposition}

\begin{proof}{}{}
\begin{itemize}
\item \textbf{Summability of $\left(x_{i}+y_{i}\right)_{i\in I}$}
\end{itemize}
We observe that
\[
\sum\limits_{i\in I}\left|x_{i}+y_{i}\right|\leqslant\sum\limits_{i\in I}\left(\left|x_{i}\right|+\left|y_{i}\right|\right)=\sum\limits_{i\in I}\left|x_{i}\right|+\sum\limits_{i\in I}\left|y_{i}\right|.
\]

(Tr.N) Since both $\left(x_{i}\right)_{i\in I}$ and $\left(y_{i}\right)_{i\in I},$
using Lemma \ref{lemma: summability of absolute value},
\[
\begin{array}{ccccc}
\sum\limits_{i\in I}\left|x_{i}\right|<+\infty & \,\,\,\, & \text{and} & \,\,\,\, & \sum\limits_{i\in I}\left|y_{i}\right|<+\infty,\end{array}
\]
which implies that
\[
\sum\limits_{i\in I}\left|x_{i}+y_{i}\right|<+\infty
\]
and $\left(x_{i}+y_{i}\right)_{i\in I}$ to be summable.
\begin{itemize}
\item \textbf{Computation of $\sum\limits_{i\in I}\left(x_{i}+y_{i}\right)$}
\end{itemize}
To establish the desired equality $\refpar{eq:sum_of_two_summable_families},$
we partition the set $I$ into six subsets:
\begin{itemize}
\item $A=\left\{ i\in I:\,x_{i}\geqslant0\text{\,and\,}y_{i}\geqslant0\right\} $
\item $B=\left\{ i\in I:\,x_{i}<0\text{\,and\,}y_{i}<0\right\} $
\item $C=\left\{ i\in I:\,x_{i}\geqslant0\text{,\,}y_{i}<0\text{\,and\,}x_{i}+y_{i}\geqslant0\right\} $
\item $D=\left\{ i\in I:\,x_{i}\geqslant0\text{,\,}y_{i}<0\text{\,and\,}x_{i}+y_{i}<0\right\} $
\item $E=\left\{ i\in I:\,x_{i}<0\text{,\,}y_{i}\geqslant0\text{\,and\,}x_{i}+y_{i}\geqslant0\right\} $
\item $F=\left\{ i\in I:\,x_{i}<0\text{,\,}y_{i}\geqslant0\text{\,and\,}x_{i}+y_{i}<0\right\} $
\end{itemize}
By Proposition \ref{Proposition:Properties-of-summable-non-neg-reals}---additivity
of summable families with nonnegative terms---we deduce
\[
\sum\limits_{i\in A}\left(x_{i}+y_{i}\right)=\sum\limits_{i\in A}x_{i}+\sum\limits_{i\in A}y_{i}.
\]

A similar equality holds for $B$, replacing $x_{i}$ by $-x_{i},$
and $y_{i}$ by $-y_{i}.$

To prove
\[
\sum\limits_{i\in C}\left(x_{i}+y_{i}\right)=\sum\limits_{i\in C}x_{i}+\sum\limits_{i\in C}y_{i}
\]
we note that since the family $\left(-y_{i}\right)_{i\in C}$ is summable,
it suffices to show
\[
\sum\limits_{i\in C}\left(x_{i}+y_{i}\right)+\sum\limits_{i\in C}\left(-y_{i}\right)=\sum\limits_{i\in C}x_{i}.
\]

Since for $i\in C,$ we have $x_{i}+y_{i}\geqslant0$ and $-y_{i}\geqslant0,$
the result follows once again from Proposition \ref{Proposition:Properties-of-summable-non-neg-reals}. 

Applying the same reasoning to subsets $D,$ $E$ and $F,$ we establish
the analogous equalities. Finally, we conclude using Lemma \ref{Lemma:Summable-families-restriction}
(ii).

\end{proof}

\begin{corollary}{Summability is Linear}{}\label{corol: linearity summability}

Let $\left(x_{i}\right)_{i\in I}$ and $\left(y_{i}\right)_{i\in I}$
be two summable families of real numbers. 

For every real numbers $\lambda$ and $\mu,$ the family $\left(\lambda x_{i}+\mu y_{i}\right)_{i\in I}$
is summable, and\boxeq{
\[
\sum\limits_{i\in I}\left(\lambda x_{i}+\mu y_{i}\right)=\lambda\sum\limits_{i\in I}x_{i}+\mu\sum\limits_{i\in I}y_{i}.
\]
}

\end{corollary}

\begin{proof}{}{}This follows immediately from the previous proposition.

\end{proof}

We now turn to the question of how the \textbf{packet summation property}
extends to families of real numbers of arbitrary sign.

As mentioned in the introduction of this chapter, the existence of
this property is a fundamental justification for using summable families
rather than absolutely convergent series. This notion enables us to
handle the sums $\sum_{i\in I}x_{i}$---and, in particular, double
sums, triple sums and higher dimension sums---in a way that is entirely
analogous to finite sums, ensuring the same level of ease and consistency.

\begin{theorem}{Packet Summation for Families of Real Numbers}{}\label{Theorem Packet-summation any real}Let
$\left(x_{i}\right)_{i\in I}$ be a family of real numbers, and let
$\left(A_{j}\right)_{j\in J}$ be a partition of $I.$

Then, for the family $\left(x_{i}\right)_{i\in I}$ to be summable,
it is \textbf{necessary and sufficient} that:
\begin{itemize}
\item Each restricted family $\left(x_{i}\right)_{i\in A_{j}}$ is summable,
\item And that the family $\left(\sum\limits_{i\in A_{j}}\left|x_{i}\right|\right)_{j\in J}$is
summable.
\end{itemize}
Moreover, if the family $\left(x_{i}\right)_{i\in I}$ is summable,
then\boxeq{
\begin{equation}
\sum\limits_{i\in I}x_{i}=\sum\limits_{j\in J}\left(\sum\limits_{i\in A_{j}}x_{i}\right).\label{eq:packet_summation_reals_any_sign}
\end{equation}
}

\end{theorem}

\begin{remark}{}{}It is worth noting that for a family $\left(x_{i}\right)_{i\in I}$
to be summable, it is not sufficient that each restricted family $\left(x_{i}\right)_{i\in A_{j}}$
to be summable and that the family of the sums $\left(\sum_{i\in A_{j}}x_{i}\right)_{j\in J}$---without
absolute values---is summable. 

A counter-example can be easily constructed by considering the series
\[
\left(1-1\right)+(1-1)+(1-1)+\cdots.
\]

\end{remark}

\begin{proof}{}{}
\begin{itemize}
\item \textbf{Necessary and sufficient condition of summability via a partition
of $I$}
\end{itemize}
By Proposition \ref{prop:packet-summation non-negative reals}---packet
summation for families of nonnegative terms---,
\begin{equation}
\sum\limits_{i\in I}\left|x_{i}\right|=\sum\limits_{j\in J}\left(\sum\limits_{i\in A_{j}}\left|x_{i}\right|\right),\label{eq:summation_per_packet_of_absolute_value}
\end{equation}
which establishes the first assertion.
\begin{itemize}
\item \textbf{Calculation of the sum using packets}
\end{itemize}
Now suppose that both sides of $\refpar{eq:summation_per_packet_of_absolute_value}$
are finite. In this case, we obtain
\begin{align*}
\sum\limits_{i\in I}x_{i} & =\sum\limits_{i\in I}x^{+}_{i}-\sum\limits_{i\in I}x^{-}_{i}\\
 & =\sum\limits_{j\in J}\left(\sum\limits_{i\in A_{j}}x^{+}_{i}\right)-\sum\limits_{j\in J}\left(\sum\limits_{i\in A_{j}}x^{-}_{i}\right)
\end{align*}
again by Proposition \ref{prop:packet-summation non-negative reals}. 

This proposition also guarantees that the families
\[
\begin{array}{ccccc}
\left(\sum\limits_{i\in A_{j}}x^{+}_{i}\right)_{j\in J} &  & \text{and} &  & \left(\sum\limits_{i\in A_{j}}x^{-}_{i}\right)_{j\in J}\end{array}
\]
are summable. 

By Corollary \ref{corol: linearity summability}, their difference
is also summable, which allows us to write
\begin{equation}
\sum\limits_{i\in I}x_{i}=\sum\limits_{j\in J}\left(\sum\limits_{i\in A_{j}}x^{+}_{i}-\sum\limits_{i\in A_{j}}x^{-}_{i}\right).\label{eq:sum_as_a_sum_of_differences}
\end{equation}

By Proposition \ref{prop:packet-summation non-negative reals}, for
each index $j$, the families $\left(x^{+}_{i}\right)_{i\in A_{j}}$
and $\left(x^{-}_{i}\right)_{i\in A_{j}}$ are both summable on $A_{j}.$
The family $\left(x_{i}\right)_{i\in A_{j}},$ which is obtained by
subtracting them, is thus also summable and, we can write
\[
\sum\limits_{i\in A_{j}}x^{+}_{i}-\sum\limits_{i\in A_{j}}x^{-}_{i}=\sum\limits_{i\in A_{j}}x_{i}.
\]

Substituting this into $\refpar{eq:sum_as_a_sum_of_differences},$
we deduce $\refpar{eq:packet_summation_reals_any_sign}.$

\end{proof}

\begin{figure}[t]
\begin{center}\includegraphics[width=0.4\textwidth]{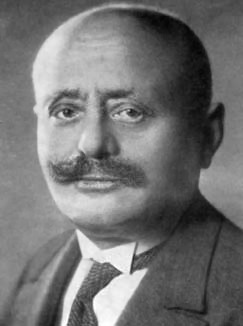}

{\tiny Credits: Unknown Public domain}\end{center}

\caption{\textbf{\protect\href{https://en.wikipedia.org/wiki/Guido_Fubini}{Guido Fubini}}
(1879 - 1943)}\sindex[fam]{Fubini, Guido}
\end{figure}

We now present the Fubini\footnote{\textbf{\href{https://en.wikipedia.org/wiki/Guido_Fubini}{Guido Fubini}}\sindex[fam]{Fubini, Guido}
(1879 - 1943) was an Italian mathematician. His main contributions
are known as the Fubini theorem and the Fubini-Study metric.} theorem for families.

\begin{corollary}{Fubini Theorem for Families}{}

(i) Let $\left(x_{ij}\right)_{\left(i,j\right)\in I\times J}$ be
a \textbf{double family\index{double family}\mindex{family ! double}}
of \textbf{nonnegative} real numbers, possibly infinite. Then\boxeq{
\begin{equation}
\sum\limits_{\left(i,j\right)\in I\times J}x_{ij}=\sum\limits_{i\in I}\left(\sum\limits_{j\in J}x_{ij}\right)=\sum\limits_{j\in J}\left(\sum\limits_{i\in I}x_{ij}\right).\label{eq:double family sum of non-neg reals}
\end{equation}
}

(ii) Let $\left(x_{ij}\right)_{\left(i,j\right)\in I\times J}$ be
a double family of real numbers (finite) of \textbf{arbitrary sign}. 

Then, if one of the sums
\[
\begin{array}{ccccc}
\sum\limits_{\left(i,j\right)\in I\times J}\left|x_{ij}\right|, & \,\,\,\, & \sum\limits_{i\in I}\left(\sum\limits_{j\in J}\left|x_{ij}\right|\right) & \,\,\,\,\text{and,}\,\,\,\, & \sum\limits_{j\in J}\left(\sum\limits_{i\in I}\left|x_{ij}\right|\right)\end{array}
\]
is finite, the two others are also finite, and\boxeq{
\begin{equation}
\sum\limits_{\left(i,j\right)\in I\times J}x_{ij}=\sum\limits_{i\in I}\left(\sum\limits_{j\in J}x_{ij}\right)=\sum\limits_{j\in J}\left(\sum\limits_{i\in I}x_{ij}\right).\label{eq:double_sum_of_reals}
\end{equation}
}---this implies that all families considered in $\refpar{eq:double_sum_of_reals}$
are summable.

\end{corollary}

\begin{proof}{}{}This follows directly from the previous proposition
by noting that
\[
I\times J=\biguplus\limits_{i\in I}\left\{ i\right\} \times J=\biguplus\limits_{j\in J}I\times\left\{ j\right\} .
\]

\end{proof}

\begin{corollary}{Product of Summable Families}{}\label{Corollary: product of families}(i)
Let $\left(x_{i}\right)_{i\in I}$ and $\left(y_{j}\right)_{j\in J}$
be two families of \textbf{nonnegative} real numbers, possibly infinite.
Then\boxeq{
\begin{equation}
\sum\limits_{\left(i,j\right)\in I\times J}x_{i}y_{j}=\left(\sum\limits_{i\in I}x_{i}\right)\left(\sum\limits_{j\in J}x_{j}\right).\label{eq:factorization}
\end{equation}
}

(ii) Let $\left(x_{i}\right)_{i\in I}$ and $\left(y_{j}\right)_{j\in J}$
be two summable families of real numbers of \textbf{arbitrary sign}.
Then, the double family $\left(x_{i}y_{j}\right)_{\left(i,j\right)\in I\times J}$
is summable and, the equality\boxeq{
\begin{equation}
\sum\limits_{\left(i,j\right)\in I\times J}x_{i}y_{j}=\left(\sum\limits_{i\in I}x_{i}\right)\left(\sum\limits_{j\in J}x_{j}\right)\label{eq:factorization-1}
\end{equation}
}remains valid.

\end{corollary}

\begin{proof}{}{}Applying the previous corollary directly yields
the result.

\end{proof}

In the remainder of this book, we often use a generalization of Corollary
\ref{Corollary: product of families} to the case of an arbitrary
finite number of families.

\begin{corollary}{Factorization of Terms of Arbitrary Finite Number of Families}{}(i)
Let $I_{1},\dots,I_{k}$ be sets, and let $\left(x^{\left(1\right)}_{i_{1}}\right)_{i_{1}\in I_{1}},\dots,\left(x^{\left(k\right)}_{i_{k}}\right)_{i_{k}\in I_{k}}$
be families of nonnegative real numbers, respectively indexed by $I_{1},\dots,I_{k}.$
Then\boxeq{
\begin{equation}
\sum\limits_{\left(i_{1},\dots,i_{k}\right)\in I_{1}\times\dots\times I_{k}}x^{(1)}_{i_{1}}\dots x^{(k)}_{i_{k}}=\left(\sum\limits_{i_{1}\in I_{1}}x^{(1)}_{i_{1}}\right)\times\dots\times\left(\sum\limits_{i_{k}\in I_{k}}x^{(k)}_{i_{k}}\right).\label{eq:factorization of k terms sum}
\end{equation}
}

(ii) Now, suppose that $\left(x^{\left(1\right)}_{i_{1}}\right)_{i_{1}\in I_{1}},\dots,\left(x^{\left(k\right)}_{i_{k}}\right)_{i_{k}\in I_{k}}$are
summable families of real numbers of arbitrary sign. Then, the family
\[
\left(x^{(1)}_{i_{1}}\dots x^{(k)}_{i_{k}}\right)_{\left(i_{1},\dots,i_{k}\right)\in I_{1}\times\dots\times I_{k}}
\]
is also summable, and the equality \boxeq{
\begin{equation}
\sum\limits_{\left(i_{1},\dots,i_{k}\right)\in I_{1}\times\dots\times I_{k}}x^{(1)}_{i_{1}}\dots x^{(k)}_{i_{k}}=\left(\sum\limits_{i_{1}\in I_{1}}x^{(1)}_{i_{1}}\right)\times\dots\times\left(\sum\limits_{i_{k}\in I_{k}}x^{(k)}_{i_{k}}\right).\label{eq:factorization of k terms sum-1}
\end{equation}
}remains valid.

\end{corollary}

\begin{proof}{}{}The proof follows by induction from the previous
corollary.

\end{proof}

\section*{Exercices}

\addcontentsline{toc}{section}{Exercises}

\begin{exercise}{}{}

1. Prove that the family $\left(x_{mn}\right)_{\left(m,n\right)\in\mathbb{N}^{\ast^{2}}}$
indexed by $\mathbb{N}^{\ast^{2}},$ and defined by
\[
\forall\left(m,n\right)\in\mathbb{N}^{\ast^{2}},\,\,\,\,x_{mn}=\dfrac{\left(-1\right)^{mn}}{m^{2}n^{2}}
\]
is summable. 

2. Compute its sum $\Sigma.$

\textit{Hint: Recall that
\[
\sum^{+\infty}_{n=1}\dfrac{1}{n^{2}}=\dfrac{\pi^{2}}{6},
\]
a result that can be obtained using the Fourier series theory. }

\end{exercise}

\begin{exercise}{Summable Families and Probability Germs}{}

Let $\alpha$ and $\beta$ be two real numbers strictly between 0
and 1, and let $g$ be the function with nonnegative values defined
on $\mathbb{N}^{2}$ by
\[
\forall\left(i,j\right)\in\mathbb{N}^{2},\,\,\,\,g\left(i,j\right)=\alpha\beta\left(1-\alpha\right)^{i}\left(1-\beta\right)^{j}.
\]

1. \textbf{Verification of Probability Germ}

Show that $g$ is a probability germ on the probabilizable space $\left(\mathbb{N}^{2},\mathcal{P}\left(\mathbb{N}^{2}\right)\right).$ 

Let $P$ be the probability associated to $g$ on the probabilizable
space $\left(\mathbb{N}^{2},\mathcal{P}\left(\mathbb{N}^{2}\right)\right).$

2. \textbf{Law of Random Variables $X$ and $Y$}

Let $X$ and $Y$ be discrete random variables defined on the probabilized
space $\left(\mathbb{N}^{2},\mathcal{P}\left(\mathbb{N}^{2}\right),P\right)$
defined by
\[
\forall\left(i,j\right)\in\mathbb{N}^{2},\,\,\,\,X\left(i,j\right)=i\,\,\,\,\text{and,\,\,\,\,}Y\left(i,j\right)=j.
\]

Determine the laws followed by the random variables $X$ and $Y.$
Are those laws identifiable to known laws? If so, which ones?

3. \textbf{Probability Computation}

Compute the following probabilities:

(a) $P\left(X=Y\right),$

(b) $P\left(X>Y\right).$

4. \textbf{Summability and Sum of the Family $Zg$}

Define the discrete random variable $Z$ on the probabilized space
$\left(\mathbb{N}^{2},\mathcal{P}\left(\mathbb{N}^{2}\right),P\right)$
by
\[
\forall\left(i,j\right)\in\mathbb{N}^{2},\,\,\,X\left(i,j\right)=\begin{cases}
1, & \text{if }i\text{\,and\,}j\,\text{are even},\\
-1, & \text{if }i\text{\,and\,}j\,\text{are odd},\\
0, & \text{otherwise}.
\end{cases}
\]

(a) \textbf{Summability of $Zg$ over $\mathbb{N}^{2}$}

Prove that the family $Zg$ defined on $\mathbb{N}^{2}$ by
\[
\forall\left(i,j\right)\in\mathbb{N}^{2},\,\,\,\,\left(Zg\right)\left(i,j\right)=Z\left(i,j\right)g\left(i,j\right)
\]
is summable over $\mathbb{N}^{2}.$ 

Compute its sum $\sum\limits_{\left(i,j\right)\in\mathbb{N}^{2}}Zg\left(i,j\right).$

(b) \textbf{Summability of $Zg$ over the Diagonal of $\mathbb{N}^{2}$}

Let $D$ denote the diagonal set $\left\{ \left(i,j\right)\in\mathbb{N}^{2}:\,i=j\right\} .$

Verify that the family $Zg$ is summable on $D,$ and determine the
value of the sum $\sum\limits_{\left(i,j\right)\in D}Zg\left(i,j\right).$

\end{exercise}

\section*{Solutions of Exercises}

\addcontentsline{toc}{section}{Solutions of Exercises}

\begin{solution}{}{}

1. \textbf{Summability of the family $\left(x_{mn}\right)_{\left(m,n\right)\in\mathbb{N}^{\ast^{2}}}$}

We have
\[
\sum\limits_{\left(m,n\right)\in\mathbb{N}^{\ast^{2}}}\left|x_{mn}\right|=\sum\limits_{\left(m,n\right)\in\mathbb{N}^{\ast^{2}}}\dfrac{1}{m^{2}n^{2}}.
\]

Applying the factorization property $\refpar{eq:factorization},$
a particular case of Fubini theorem for nonnegative families, we obtain
\[
\sum\limits_{\left(m,n\right)\in\mathbb{N}^{\ast^{2}}}\left|x_{mn}\right|=\left(\sum\limits_{n\in\mathbb{N}^{\ast}}\dfrac{1}{n^{2}}\right)^{2}.
\]

As$\sum\limits_{n\in\mathbb{N}^{\ast}}\dfrac{1}{n^{2}}=\dfrac{\pi^{2}}{6}<+\infty,$
it follows that
\[
\sum\limits_{\left(m,n\right)\in\mathbb{N}^{\ast^{2}}}\left|x_{mn}\right|<+\infty.
\]

Thus, the family $\left(x_{mn}\right)$ is summable on $\mathbb{N}^{\ast2}.$ 

\textbf{2. Computation of the sum of the family $\left(x_{mn}\right)_{\left(m,n\right)\in\mathbb{N}^{\ast^{2}}}$}

Now, we compute its sum
\begin{equation}
\Sigma=\sum\limits_{\substack{\left(m,n\right)\in\mathbb{N}^{\ast^{2}}\\
\text{s.t.\,}mn\text{\,even}
}
}\dfrac{1}{m^{2}n^{2}}-\sum\limits_{\substack{\left(m,n\right)\in\mathbb{N}^{\ast^{2}}\\
\text{s.t.\,}mn\text{\,odd}
}
}\dfrac{1}{m^{2}n^{2}}.\label{eq:sum_ex1_even-odd}
\end{equation}

We write
\[
\begin{array}{ccccc}
E=\left\{ 2k,k\in\mathbb{N}^{*}\right\}  & \,\,\,\, & \text{and} & \,\,\,\, & O=\left\{ 2k+1,k\in\mathbb{N}\right\} .\end{array}
\]

Then
\[
\left\{ \left(m,n\right)\in\mathbb{N}^{\ast^{2}}:\,mn\text{ is even}\right\} =\left(E\times E\right)\uplus\left(E\times O\right)\uplus\left(O\times E\right)
\]
and
\[
\left\{ \left(m,n\right)\in\mathbb{N}^{\ast^{2}}:\,mn\text{ is odd}\right\} =O\times O.
\]

By addivity over these sets,
\[
\sum\limits_{\substack{\left(m,n\right)\in\mathbb{N}^{\ast^{2}}\\
\text{s.t.\,}mn\text{\,even}
}
}\dfrac{1}{m^{2}n^{2}}=\sum\limits_{\substack{\left(m,n\right)\in E\times E}
}\dfrac{1}{m^{2}n^{2}}+\sum\limits_{\substack{\left(m,n\right)\in E\times O}
}\dfrac{1}{m^{2}n^{2}}+\sum\limits_{\substack{\left(m,n\right)\in O\times E}
}\dfrac{1}{m^{2}n^{2}},
\]
and
\[
\sum\limits_{\substack{\left(m,n\right)\in\mathbb{N}^{\ast^{2}}\\
\text{s.t.\,}mn\text{\,odd}
}
}\dfrac{1}{m^{2}n^{2}}=\sum\limits_{\substack{\left(m,n\right)\in O\times O}
}\dfrac{1}{m^{2}n^{2}}.
\]

Define
\[
\begin{array}{ccccc}
U=\sum\limits_{n\in E}\dfrac{1}{n^{2}} & \,\,\,\, & \text{and} & \,\,\,\, & V=\sum\limits_{n\in O}\dfrac{1}{n^{2}}\end{array}.
\]

Using Fubini theorem for nonnegative families, we obtain
\[
\sum\limits_{\substack{\left(m,n\right)\in E\times E}
}\dfrac{1}{m^{2}n^{2}}=\left(\sum\limits_{m\in E}\dfrac{1}{m^{2}}\right)\left(\sum\limits_{n\in E}\dfrac{1}{n^{2}}\right)=U^{2},
\]
and similarly
\[
\sum\limits_{\substack{\left(m,n\right)\in E\times O}
}\dfrac{1}{m^{2}n^{2}}=\left(\sum\limits_{m\in E}\dfrac{1}{m^{2}}\right)\left(\sum\limits_{n\in O}\dfrac{1}{n^{2}}\right)=UV,
\]
\[
\sum\limits_{\substack{\left(m,n\right)\in O\times E}
}\dfrac{1}{m^{2}n^{2}}=\left(\sum\limits_{m\in O}\dfrac{1}{m^{2}}\right)\left(\sum\limits_{n\in E}\dfrac{1}{n^{2}}\right)=VU,
\]
\[
\sum\limits_{\substack{\left(m,n\right)\in O\times O}
}\dfrac{1}{m^{2}n^{2}}=\left(\sum\limits_{m\in O}\dfrac{1}{m^{2}}\right)\left(\sum\limits_{n\in O}\dfrac{1}{n^{2}}\right)=V^{2}.
\]
Thus
\[
\Sigma=U^{2}+2UV-V^{2}.
\]

Since, we have by additivity
\[
S=U+V,
\]
and as
\[
U=\sum\limits_{k\in\mathbb{N}^{\ast}}\dfrac{1}{\left(2k\right)^{2}}=\dfrac{1}{4}\sum\limits_{k\in\mathbb{N}^{\ast}}\dfrac{1}{k^{2}}=\dfrac{S}{4},
\]
we deduce
\[
V=S-U=\dfrac{3S}{4},
\]
Finally, using $S=\dfrac{\pi^{2}}{6},$ we obtain\boxeq{
\[
\Sigma=\dfrac{S^{2}}{16}+\dfrac{3S^{2}}{8}-\dfrac{9S^{2}}{16}=-\dfrac{S^{2}}{8}=-\dfrac{\pi^{4}}{288}\approx-0.338.
\]
}

\end{solution}

\begin{solution}{}{}

\textbf{1. Verification that $g$ is a probability germ}

Since $g$ is a nonnegative function with separate variables---i.e.
$g\left(i,j\right)=g_{1}\left(i\right)g_{2}\left(j\right)$---we
compute
\begin{align*}
\sum\limits_{\left(i,j\right)\in\mathbb{N}^{2}}g\left(i,j\right) & =\alpha\beta\left(\sum\limits_{i\in\mathbb{N}}\left(1-\alpha\right)^{i}\right)\left(\sum\limits_{j\in\mathbb{N}}\left(1-\beta\right)^{j}\right).\\
 & =\alpha\beta\left(\sum\limits^{+\infty}_{i=0}\left(1-\alpha\right)^{i}\right)\left(\sum\limits^{+\infty}_{j=0}\left(1-\beta\right)^{j}\right)\\
 & =\alpha\beta\underset{n\rightarrow+\infty}{\lim}\dfrac{1-\left(1-\alpha\right)^{n+1}}{1-\left(1-\alpha\right)}\underset{n\rightarrow+\infty}{\lim}\dfrac{1-\left(1-\beta\right)^{n+1}}{1-\left(1-\beta\right)}\\
 & =\alpha\beta\dfrac{1}{1-\left(1-\alpha\right)}\dfrac{1}{1-\left(1-\beta\right)}\\
 & =1,
\end{align*}
using the fact that geometrical series of reason strictly lower than
1 in absolute value converge and Proposition \ref{Prop:countable enumeration}. 

Since $g$ is nonnegative and sums to 1 over $\mathbb{N}^{2},$ it
defines a probability germ on the probabilizable space $\left(\mathbb{N}^{2},\mathcal{P}\left(\mathbb{N}^{2}\right)\right).$

\textbf{2. Laws of $X$ and $Y$}

For every $i\in\mathbb{N},$
\[
\left(X=i\right)=\left\{ i\right\} \times\mathbb{N}.
\]

Applying the definition of the probability $P$ associated with $g,$
and applying Proposition \ref{Prop:countable enumeration} on countable
additivity, we obtain
\begin{align*}
P\left(X=i\right) & =\sum\limits_{\left(j,k\right)\in\left\{ i\right\} \times\mathbb{N}}g\left(j,k\right)\\
 & =\alpha\left(1-\alpha\right)^{i}\beta\sum\limits_{k\in\mathbb{N}}\left(1-\beta\right)^{k}\\
 & =\alpha\left(1-\alpha\right)^{i}\beta\underset{n\rightarrow+\infty}{\lim}\dfrac{1-\left(1-\beta\right)^{n+1}}{1-\left(1-\beta\right)}\\
 & =\alpha\left(1-\alpha\right)^{i}\beta\dfrac{1}{1-(1-\beta)}\\
 & =\alpha\left(1-\alpha\right)^{i}.
\end{align*}

This shows that the law of the random variable $X$ is the geometric
law of parameter $\alpha$ on $\mathbb{N}.$ Similarly, the law of
the random variable $Y$ is the geometric law of parameter $\beta$
on $\mathbb{N}.$

\textbf{3. Computation of probabilities}

The method is to partition the events.

\textbf{(a) $P\left(X=Y\right)$}

Partition the event
\[
\left(X=Y\right)=\biguplus\limits_{i\in\mathbb{N}}\left(\left(X=i\right)\cap\left(Y=i\right)\right)=\left\{ \left(i,i\right):\,i\in\mathbb{N}\right\} .
\]
Thus,
\begin{align*}
P\left(X=Y\right) & =\sum\limits_{i\in\mathbb{N}}g\left(i,i\right)=\alpha\beta\sum\limits_{i\in\mathbb{N}}\left(1-\alpha\right)^{i}\left(1-\beta\right)^{i}.
\end{align*}
Since $0<\left(1-\alpha\right)<1$ and $0<\left(1-\beta\right)<1,$
we have $0<\left(1-\alpha\right)\left(1-\beta\right)<1,$ and, since
the previous sum is the sum of a geometric series, and by applying
Proposition \ref{Prop:countable enumeration}, we have\boxeq{
\[
P\left(X=Y\right)=\alpha\beta\underset{n\rightarrow+\infty}{\lim}\dfrac{1-\left[\left(1-\alpha\right)\left(1-\beta\right)\right]^{n+1}}{1-\left(1-\alpha\right)\left(1-\beta\right)}=\dfrac{\alpha\beta}{1-\left(1-\alpha\right)\left(1-\beta\right)}.
\]
}

\textbf{(b) }$P\left(X>Y\right)$

Partitioning the event, it holds
\[
\left(X>Y\right)=\biguplus\limits_{\substack{\left(i,j\right)\in\mathbb{N}^{2}\\
i>j
}
}\left(\left(X=i\right)\cap\left(Y=j\right)\right)=\left\{ \left(i,j\right)\in\mathbb{N}^{2}:\,i>j\right\} ,
\]
and, consequently\boxeq{
\[
P\left(X>Y\right)=\sum\limits_{\substack{\left(i,j\right)\in\mathbb{N}^{2}\\
i>j
}
}\left[\alpha\beta\left(1-\alpha\right)^{i}\left(1-\beta\right)^{j}\right].
\]
}As $g$ is nonnegative, we can apply the Fubini Theorem and Proposition
\ref{Prop:countable enumeration}, and write successively
\begin{align*}
P\left(X>Y\right) & =\alpha\beta\sum\limits_{j\in\mathbb{N}}\left[\left(1-\beta\right)^{j}\sum\limits_{\substack{i\in\mathbb{N}\\
i>j
}
}\left(1-\alpha\right)^{i}\right]\\
 & =\alpha\beta\sum\limits^{+\infty}_{j=0}\left[\left(1-\beta\right)^{j}\sum\limits^{+\infty}_{i=j+1}\left(1-\alpha\right)^{i}\right]\\
 & =\alpha\beta\sum\limits^{+\infty}_{j=0}\left[\left(1-\beta\right)^{j}\underset{n\rightarrow+\infty}{\lim}\sum\limits^{n}_{i=j+1}\left(1-\alpha\right)^{i}\right]\\
 & =\alpha\beta\sum\limits^{+\infty}_{j=0}\left[\left(1-\beta\right)^{j}\underset{n\rightarrow+\infty}{\lim}\left(1-\alpha\right)^{j+1}\dfrac{1-\left(1-\alpha\right)^{n+1}}{1-\left(1-\alpha\right)}\right]\\
 & =\alpha\beta\sum\limits^{+\infty}_{j=0}\left[\left(1-\beta\right)^{j}\dfrac{\left(1-\alpha\right)^{j+1}}{\alpha}\right]\\
 & =\beta\left(1-\alpha\right)\underset{n\to+\infty}{\lim}\sum\limits^{n}_{j=0}\left[\left(1-\beta\right)\left(1-\alpha\right)\right]^{j}\\
 & =\dfrac{\beta\left(1-\alpha\right)}{1-\left(1-\alpha\right)\left(1-\beta\right)}.
\end{align*}

\textbf{4. Summability of $Zg$ and computation of sums}

We apply in this question the fundamental properties of summable families.

(a) \textbf{Summability of $Zg$ over $\mathbb{N}^{2}:$}

Since $\left|Z\right|\leqslant1,$ we have $\left|Zg\right|\leqslant\left|g\right|=g.$
Since the family $g$ is summable on $\mathbb{N}^{2},$ it follows
that the family $Zg$ is also summable.

We compute its sum
\[
\sum\limits_{\left(i,j\right)\in\mathbb{N}^{2}}Zg\left(i,j\right)=\sum\limits_{\left(i,j\right)\in\left(2\mathbb{N}\right)^{2}}g\left(i,j\right)-\sum\limits_{\left(i,j\right)\in\left(2\mathbb{N}+1\right)^{2}}g\left(i,j\right).
\]

Using the Fubini theorem, we have (writing $i=2k$ and $j=2l$)
\begin{align*}
\sum\limits_{\left(i,j\right)\in\left(2\mathbb{N}\right)^{2}}g\left(i,j\right) & =\sum\limits_{\left(k,l\right)\in\mathbb{N}^{2}}\left[\alpha\beta\left(1-\alpha\right)^{2k}\left(1-\beta\right)^{2l}\right]\\
 & =\alpha\beta\left(\sum\limits_{k\in\mathbb{N}}\left(1-\alpha\right)^{2k}\right)\left(\sum\limits_{k\in\mathbb{N}}\left(1-\beta\right)^{2k}\right)\\
 & =\alpha\beta\dfrac{1}{1-\left(1-\alpha\right)^{2}}\dfrac{1}{1-\left(1-\beta\right)^{2}}\\
 & =\dfrac{1}{\left(2-\alpha\right)\left(2-\beta\right)}.
\end{align*}

Following the same method, we have
\begin{align*}
\sum\limits_{\left(i,j\right)\in\left(2\mathbb{N}+1\right)^{2}}g\left(i,j\right) & =\sum\limits_{\left(k,l\right)\in\mathbb{N}^{2}}\left[\alpha\beta\left(1-\alpha\right)^{2k+1}\left(1-\beta\right)^{2l+1}\right]\\
 & =\left(1-\alpha\right)\left(1-\beta\right)\sum\limits_{\left(i,j\right)\in\left(2\mathbb{N}\right)^{2}}g\left(i,j\right),
\end{align*}
thus
\begin{align*}
\sum\limits_{\left(i,j\right)\in\mathbb{N}^{2}}Zg\left(i,j\right) & =\left[1-\left(1-\alpha\right)\left(1-\beta\right)\right]\sum\limits_{\left(i,j\right)\in\left(2\mathbb{N}\right)^{2}}g\left(i,j\right)\\
 & =\dfrac{\alpha+\beta-\alpha\beta}{\left(2-\alpha\right)\left(2-\beta\right)}.
\end{align*}

(b) \textbf{Summability of $Zg$ over $D$}

Since the family $Zg$ is summable on $\mathbb{N}^{2},$ it is also
summable on $D.$

We compute its sum
\begin{align*}
\sum\limits_{\left(i,j\right)\in D}Zg\left(i,j\right) & =\sum\limits_{\left(i,j\right)\in D\cap\left(2\mathbb{N}\right)^{2}}g\left(i,j\right)-\sum\limits_{\left(i,j\right)\in D\cap\left(2\mathbb{N}+1\right)^{2}}g\left(i,j\right)\\
 & =\sum\limits_{i\in\mathbb{N}}\left[\alpha\beta\left(1-\alpha\right)^{2i}\left(1-\beta\right)^{2i}\right]-\sum\limits_{i\in\mathbb{N}}\left[\alpha\beta\left(1-\alpha\right)^{2i+1}\left(1-\beta\right)^{2i+1}\right]\\
 & =\alpha\beta\left[1-\left(1-\alpha\right)\left(1-\beta\right)\right]\dfrac{1}{1-\left(1-\alpha\right)^{2}\left(1-\beta\right)^{2}}\\
 & =\dfrac{\alpha\beta}{1+\left(1-\alpha\right)\left(1-\beta\right)}.
\end{align*}

\end{solution}

\chapter{Independence}\label{chap:Independence}

\begin{objective}{}{}

Chapter \ref{chap:Independence} begins with an introductory example
to present the concept of independence.
\begin{itemize}
\item Section \ref{sec:Independence-of-Events} provides the formal definition
of \textbf{independence} between two events, as well as for family
of events---both \textbf{mutually and pairwise}. The \textbf{complementation}
of events and its impact on independence is also discussed. The concept
of independence is then extended to random variables, including a
method to \textbf{verify independence via $n$-tuples}. Finally, the
\textbf{independence of functions of independent random variables}
is examined covering both univariate and multivariate functions.
\item Section \ref{sec:Law-of-the-Sum-of-Inpendent-RV} explores the law
governing the \textbf{sum of independent random variables}, introducing
the \textbf{convolution product induced by the sum}. A special result
on the \textbf{convolution of two Poisson laws} is also presented.
\item Section \ref{sec:Independence-and-Cartesian} analyzes independence
in the context of \textbf{Cartesian products}. It introduces the \textbf{product
probability}, and concludes on the existence of a $\sigma-$algebra
where variables are independent given a Cartesian product of probabilizable
spaces.
\item Section \ref{sec:Geometric-and-Binomial} begins with the definition
of a \textbf{Bernoulli random variable}. It then examines the case
of \textbf{geometric and binomial laws}, and concludes with the study
of the \textbf{convolution of two binomial laws}.
\end{itemize}
\end{objective}

\section*{Introduction}

The modelling of successive rolls of a same die naturally leads us,
when considering events within the set of two rolls, to consider as
the fundamental space, or realization space, the Cartesian product
$\Omega=\Omega_{1}\times\Omega_{2}$ where $\Omega_{1}=\Omega_{2}=\left\{ 1,\dots,6\right\} .$
We then equip this space with the uniform probability as the probability
$P$ on the probabilizable space $\left(\Omega,\mathcal{P}\left(\Omega\right)\right).$ 

An event $A$ that depends only on the outcome of the first roll is
represented by a subset of $\Omega$ of the form 
\[
A=C_{1}\times\Omega_{2}.
\]
Similarly, an event $B$ that depends only on the second roll is represented
by a subset of $\Omega$ of the form 
\[
B=\Omega_{1}\times C_{2}.
\]
Intuitively, these events are independent, as the outcome of the first
roll should not influence the outcome of the second roll. Let us now
compute the probability of the simultaneous occurrence of these events.
We have
\[
A\cap B=\left(C_{1}\times\Omega_{2}\right)\cap\left(\Omega_{1}\times C_{2}\right)=C_{1}\times C_{2}.
\]

The additivity of a probability induces
\begin{align*}
P\left(C_{1}\times C_{2}\right) & =\sum\limits_{\omega\in C_{1}}P\left(\left\{ \omega\right\} \times C_{2}\right)\\
 & =\sum\limits_{\omega\in C_{1}}\dfrac{\left|C_{2}\right|}{\left|\Omega_{1}\right|\left|\Omega_{2}\right|}\\
 & =\dfrac{\left|C_{1}\right|\left|C_{2}\right|}{\left|\Omega_{1}\right|\left|\Omega_{2}\right|}.
\end{align*}

In particular, for $C_{2}=\Omega_{2},$ we get
\[
P\left(C_{1}\times\Omega_{2}\right)=\dfrac{\left|C_{1}\right|}{\left|\Omega_{1}\right|}.
\]

Similarly, with $C_{1}=\Omega_{1},$ we find
\[
P\left(\Omega_{1}\times C_{2}\right)=\dfrac{\left|C_{2}\right|}{\left|\Omega_{2}\right|}.
\]

Hence, we have shown that\boxeq{
\begin{equation}
P\left(A\cap B\right)=P\left(A\right)P\left(B\right).\label{eq:P(AandB)_independent}
\end{equation}
}

This relation is taken as definition of the \textbf{independence of
two events\mindex{independence!two events}\mindex{event ! independence}}
$A$ and $B$. In the example, it is simply a direct translation of
the physical independence of the two die rolls.

Another way to describe this example is by introducing the random
variables $X_{1}$ and $X_{2}$ defined on $\Omega,$ taking values
in $\Omega_{1}$ and $\Omega_{2}$ respectively, and given by
\[
\begin{array}{cccccc}
\forall\left(\omega_{1},\omega_{2}\right)\in\Omega, & \,\,\,\, & X_{1}\left(\omega_{1},\omega_{2}\right)=\omega_{1} & \,\,\,\, & \text{and} & \,\,\,\,\end{array}X_{2}\left(\omega_{1},\omega_{2}\right)=\omega_{2}.
\]
From a set-theoretic perspective, these variables correspond to the
natural projections of the Cartesian product onto its components.
From a probabilistic perspective, $X_{1}$ and $X_{2}$ represent
the outcome of the first and second rolls, respectively. The event
$A$ can now be written as
\[
A=\left(X_{1}\in C_{1}\right).
\]

Similarly,
\[
B=\left(X_{2}\in C_{2}\right).
\]

Hence,
\[
A\cap B=\left(X_{1}\in C_{1}\right)\cap\left(X_{2}\in C_{2}\right).
\]

To designate this event which is a conjunction or a simultaneous realization
of events $\left(X_{1}\in C_{1}\right)$ and $\left(X_{2}\in C_{2}\right),$
the common usage denotes it as
\[
\left(X_{1}\in C_{1},X_{2}\in C_{2}\right).
\]

Using this notation, the independence relation in the formula $\refpar{eq:P(AandB)_independent}$
can be rewritten as\boxeq{
\begin{equation}
P\left(X_{1}\in C_{1},X_{2}\in C_{2}\right)=P\left(X_{1}\in C_{1}\right)P\left(X_{2}\in C_{2}\right).\label{eq:independence_random_var}
\end{equation}
}

This relation serves as the definition of the \textbf{independence
of two random variables\mindex{independence!two random variables}}.
More generally, let $X_{1}$ and $X_{2}$ be two random variables
defined on the same probabilized space $\left(\Omega,\mathcal{A},P\right),$
taking values in the probabilizable spaces $\left(\Omega_{1},\mathcal{A}_{1}\right)$
and $\left(\Omega_{2},\mathcal{A}_{2}\right)$ respectively. The variables
$X_{1}$ and $X_{2}$ are said to be \textbf{independent\index{independent}\mindex{random variables ! independent}}
if the relation $\refpar{eq:independence_random_var}$ holds for every
$C_{1}\in\mathcal{A}_{1}$ and $C_{2}\in\mathcal{A}_{2}.$

It is worth noting that this independence condition in $\refpar{eq:independence_random_var}$
can also be expressed in term of random variable laws\footnote{Tr.N. Refer to Proposition $\ref{pr:randomvarprobalaw}$ for the definition
of the law of a random variable.} as
\begin{equation}
P_{\left(X_{1},X_{2}\right)}\left(C_{1}\times C_{2}\right)=P_{X_{1}}\left(C_{1}\right)P_{X_{2}}\left(C_{2}\right).\label{eq:random var laws independence}
\end{equation}

\begin{remark}{}{}Returning to our previous example, we observe that
\[
P_{X_{1}}\left(C_{1}\right)=P\left(X_{1}\in C_{1}\right)=P\left(C_{1}\times\Omega_{1}\right)
\]

Thus,
\[
\begin{array}{cc}
\forall C_{i}\in\mathcal{P}\left(\Omega_{i}\right), & \,\,\,\,P_{X_{i}}\left(C_{i}\right)=\dfrac{\left|C_{i}\right|}{\left|\Omega_{_{i}}\right|},\end{array}\text{with }i=1,2.
\]

This confirms that the random variables $X_{1}$ and $X_{2}$ follow
the uniform law on $\Omega_{1}$ and $\Omega_{2}$ respectively.

\end{remark}

As the following example illustrates, stating that the probability
$P$ is uniform on $\Omega$ conveys more information (and ``more''
is, as will be seen in Proposition $\ref{pr:product_probability}$,
underlies the notion of independence) than merely stating that each
random variable $X_{i}$ follows a uniform law on $\Omega_{i}.$ Indeed,
in this example, we construct a probability $P$ on $\Omega$ that
is not uniform, yet each random variable $X_{i}$ follows a uniform
law on $\Omega_{i}.$

Define the subsets $\Lambda_{1}$ and $\Lambda_{2}$ of $\Omega$
as
\[
\Lambda_{1}=\left\{ \left(0,1\right),\left(1,2\right),\left(2,0\right)\right\} 
\]
and
\[
\Lambda_{2}=\left\{ \left(0,2\right),\left(2,1\right),\left(1,0\right)\right\} .
\]

Let $\epsilon$ be a real number such that
\[
0<\epsilon<\dfrac{1}{36}.
\]

Define the probability germ $g$ on $\Omega$ by
\[
\begin{array}{ccc}
\forall\left(i,j\right)\in\Omega, & \,\,\,\, & g\left(i,j\right)=\begin{cases}
\frac{1}{36}+\epsilon, & \text{if }\left(i,j\right)\in\Lambda_{1},\\
\frac{1}{36}-\epsilon, & \text{if }\left(i,j\right)\in\Lambda_{2},\\
\frac{1}{36}, & \text{otherwise.}
\end{cases}\end{array}
\]

Let $P$ be the (non uniform) probability on $\left(\Omega,\mathcal{P}\left(\Omega\right)\right)$
associated to $g.$

For every $i\in\left\{ 1,2,\dots,6\right\} ,$
\[
\left(X_{1}=i\right)=\biguplus\limits^{6}_{j=1}\left\{ \left(i,j\right)\right\} ,
\]
and thus, as it is easily verified
\[
P\left(X_{1}=i\right)=\sum^{6}_{j=1}P\left(\left\{ \left(i,j\right)\right\} \right)=\sum^{6}_{j=1}g\left(i,j\right)=\dfrac{1}{6}.
\]

It means that $X_{1}$ follows the uniform law on $\Omega_{1},$ and
similarly, $X_{2}$ follows the uniform law on $\Omega_{2}.$ However,
the random variables $X_{1}$ and $X_{2}$ are not independent in
the sense defined by the relation $\refpar{eq:independence_random_var}$
since, for instance, we have
\[
P\left(X_{1}=1,X_{2}=2\right)=g\left(1,2\right)=\dfrac{1}{36}+\epsilon,
\]
which differs from product
\[
P\left(X_{1}=1\right)P\left(X_{2}=2\right)=\dfrac{1}{36}.
\]

\section{Independence of Events and Random Variables}\label{sec:Independence-of-Events}

We now study the notion of independence in a more general setting.
In some cases, however, we must restrict ourselves to the case of
discrete random variables\footnote{The general study will be done in the second part of this book, within
the framework of probability theory based on measure theory.}.

\subsection{Independent Events}

\begin{definition}{Independence of Two Events and of a Family of Events}{}

Let $\left(\Omega,\mathcal{A},P\right)$ be a probabilized space.

(i) Two events $A\in\mathcal{A}$ and $B\in\mathcal{A}$ are said
to be \textbf{independent}\index{independent}\mindex{events ! independent}
if\boxeq{
\[
P\left(A\cap B\right)=P\left(A\right)P\left(B\right).
\]
}

(ii) Let $\left(A_{i}\right)_{i\in I}$ be a family of events. These
events are said to be \textbf{independent}, sometimes\index{independent}\footnotemark
we also say \textbf{independent as a whole}\index{independent as a whole}\mindex{events ! independent ! as a whole}\mindex{independent ! as a whole}
, if, for every non-empty finite subset $J$ of $I,$ we have\boxeq{
\[
P\left(\bigcap\limits_{j\in J}A_{j}\right)=\prod\limits_{j\in J}P\left(A_{j}\right).
\]
}

In particular, some events $A_{1},A_{2},\dots,A_{n}$ are mutually
independent, if for any finite sequence $\left(j_{k}\right)_{k\geqslant2}$
such that $1\leqslant j_{1}<j_{2}<\dots<j_{k}\leqslant n,$ the following
holds\boxeq{
\[
P\left(A_{j_{1}}\cap A_{j_{2}}\cap\dots\cap A_{j_{k}}\right)=P\left(A_{j_{1}}\right)P\left(A_{j_{2}}\right)\dots P\left(A_{j_{k}}\right).
\]
}

(iii) Let $\left(A_{i}\right)_{i\in I}$ be a family of events. These
events are said to be \textbf{pairwise independent}\mindex{independent ! pairwise}\textbf{\mindex{events ! family ! pairwise independent}}
if for any pair of distinct event of this family the events are independent.

\end{definition}

\footnotetext{Sometimes also called \textbf{mutually independent}\index{mutually independent}\mindex{independent ! mutually}\textbf{\mindex{events ! independent ! mutually}}.}

\begin{remarks}{}{}1. If $I=\left\{ 1,2\right\} ,$ then the definition
reduces to the independence of two events.

2. The independence of the events $A_{1},A_{2},\dots,A_{n}$ imposes
$2^{n}-n-1$ conditions to be filled---corresponding to the number
of subsets of $I=\left\{ 1,2,\dots,n\right\} $ with at least two
elements. 

3. If $n$ events are \textbf{mutually independent}\index{mutually independent},
then they are also \textbf{pairwise independent}\index{pairwise independent},
but the converse is not necessarily true, as shown in the next example.

\end{remarks}

\begin{example}{}{}Consider the probabilized space $\left(\Omega,\mathcal{A},P\right),$
where $\Omega$ is a set with four elements, 
\[
\Omega=\left\{ \omega_{1},\omega_{2},\omega_{3},\omega_{4}\right\} .
\]

Take the $\sigma-$algebra of all subsets of $\Omega$ and equip it
with $P,$ the uniform probability on $\Omega.$ Consider the events
\[
\begin{array}{ccccc}
A=\left\{ \omega_{1},\omega_{2}\right\} , & \,\,\,\, & B=\left\{ \omega_{1},\omega_{3}\right\} , &  & C=\left\{ \omega_{1},\omega_{4}\right\} .\end{array}
\]

Show that these event are pairwise independent, but not mutually.

\end{example}

\begin{solutionexample}{}{}

Compute
\[
P\left(A\cap B\right)=P\left(A\cap C\right)=P\left(B\cap C\right)=\dfrac{1}{4}
\]
and
\[
P\left(A\right)P\left(B\right)=P\left(A\right)P\left(C\right)=P\left(B\right)P\left(C\right)=\dfrac{1}{4}.
\]

Thus, the events $A,$ $B$ and $C$ are pairwise independent. However,
\[
\begin{array}{ccccc}
P\left(A\cap B\cap C\right)=\dfrac{1}{4}, & \,\,\,\, & \text{whereas} & \,\,\, & P\left(A\right)P\left(B\right)P\left(C\right)=\dfrac{1}{8}\end{array}.
\]

The events $A,$ $B$ and $C$ are then \textbf{not mutually independent}.

\end{solutionexample}

4. It is possible to have
\[
P\left(A\cap B\cap C\right)=P\left(A\right)P\left(B\right)P\left(C\right)
\]
without the events $A,$ $B$ and $C$ being pairwise independent.

\begin{example}{}{}

Consider the probabilized space $\left(\Omega,\mathcal{P}\left(\Omega\right),P\right),$
where
\[
\Omega=\left\{ 1,2,\dots,6\right\} ^{2}
\]
equipped with the $\sigma-$algebra of all subsets and where the probability
$P$ is the uniform probability. 

Consider the three events
\[
\begin{array}{ccccc}
A=\Omega\times\left\{ 1,2,5\right\} , & \,\,\,\, & B=\Omega\times\left\{ 4,5,6\right\} , & \,\,\,\, & C=\left\{ \left(i,j\right)\in\Omega:\,i+j=9\right\} .\end{array}
\]

Show that these events are three-way independent, but are neither
independent mutually nor pairwise.

\end{example}{}{}

\begin{solutionexample}{}{}

The event $C$ can also be rewritten explicitly as
\[
C=\left\{ \left(3,6\right),\left(4,5\right),\left(5,4\right),\left(6,3\right)\right\} .
\]

Thus, we compute
\[
P\left(C\right)=\dfrac{4}{36}.
\]

Furthermore, 
\[
\begin{array}{ccc}
P\left(A\cap B\right)=\dfrac{6}{36} & \,\text{and,}\, & P\left(A\right)P\left(B\right)=\dfrac{1}{2}\times\dfrac{1}{2}=\dfrac{1}{4},\end{array}
\]
which shows that\boxeq{
\[
P\left(A\cap B\right)\neq P\left(A\right)P\left(B\right).
\]
}

Similarly,
\[
\begin{array}{ccc}
P\left(A\cap C\right)=\dfrac{1}{36}, & \,\text{and}\, & P\left(A\right)P\left(C\right)=\dfrac{1}{2}\times\dfrac{1}{9}=\dfrac{1}{18}.\end{array}
\]

Thus,\boxeq{
\[
P\left(A\cap C\right)\neq P\left(A\right)P\left(C\right).
\]
}

Last,
\[
\begin{array}{ccc}
P\left(B\cap C\right)=\dfrac{3}{36}, & \,\text{and}\, & P\left(A\right)P\left(C\right)=\dfrac{1}{2}\times\dfrac{1}{9}=\dfrac{1}{18}.\end{array}
\]

Hence,\boxeq{
\[
P\left(B\cap C\right)\neq P\left(B\right)P\left(C\right).
\]
}

Nonetheless,
\[
\begin{array}{ccc}
P\left(A\cap B\cap C\right)=\dfrac{1}{36} & \,\text{and}\, & P\left(A\right)P\left(B\right)P\left(C\right)=\dfrac{1}{36}.\end{array}
\]

Then\boxeq{
\[
P\left(A\cap B\cap C\right)=P\left(A\right)P\left(B\right)P\left(C\right).
\]
}

This shows that the events $A,$$B$ and $C$ are \textbf{three-way
independent} but are neither\textbf{ independent mutually}, nor\textbf{
pairwise independent.}

\end{solutionexample}

\begin{remark}{}{}Tr.N.: To dig further the notion of dependency,
the interested readed may refer to the article by J. Stoyanov \cite{stoyanov1998global}
on the global dependency measure for sets of random elements, where
the author considers all possible combinations of random elements,
defining a hierarchical property of independence. The amount of dependency
contained in the whole set of random elements is then considered by
defining a dependency measure taking values in the interval $\left[0,1\right].$
This is then applied to the Italian problem corresponding to the existence
of a probability model and a set of random variables for any given
independence / dependence structure.

\end{remark}

\subsection{Independent Events and Complementation}

Saying that two events $A$ and $B$ are independent means that the
occurrence or the non-occurrence of $B$ has no influence on the probability
of $A$ to occur. It is then natural to consider wether this remains
true if we replace $B$ by its complement $B^{c}.$ The following
proposition confirms this intuition.

\begin{proposition}{Independence and Complementation}{}\label{Proposition 3.2 independence}Let
$\left(\Omega,\mathcal{A},P\right)$ be a probabilized space.

(i) If $A$ and $B$ are two independent events, the events $A$ and
$B^{c},$ $A^{c}$ and $B^{c}$ and last, $A^{c}$ and $B$ are also
independent.

(ii) This property expends to families of independent events as follows: 

Let $\left(A_{i}\right)_{i\in I}$ be a family of independent events
indexed by an arbitrary set $I.$ Suppose that the set $I$ is partitioned
in two disjoint subsets $I_{1}$ and $I_{2}.$

Define the family $\left(B_{i}\right)_{i\in I}$ by
\[
B_{i}=\begin{cases}
A_{i}, & \text{if }i\in I_{1},\\
A^{c}_{i}, & \text{if }i\in I_{2.}
\end{cases}
\]

The events $\left(B_{i}\right)_{i\in I}$ are mutually independent.

\end{proposition}

\begin{proof}{}{}(i) We prove the results for $A$ and $B^{c},$
the other results can be deduced by taking complements. 

By definition of probability, as $A=\left(A\cap B^{c}\right)\uplus\left(A\cap B\right),$
\begin{equation}
P\left(A\cap B^{c}\right)=P\left(A\right)-P\left(A\cap B\right).\label{eq:prob_A_inter_B^C}
\end{equation}

Since $A$ and $B$ are independent,
\[
P\left(A\cap B\right)=P\left(A\right)P\left(B\right).
\]

Substituting this into the previous equation,
\[
P\left(A\cap B^{c}\right)=P\left(A\right)-P\left(A\right)P\left(B\right).
\]

Factoring out $P\left(A\right)$ and using the fact that $P\left(B^{c}\right)=1-P\left(B\right),$
yields
\[
P\left(A\cap B^{c}\right)=P\left(A\right)P\left(B^{c}\right).
\]

This proves the independence of $A$ and $B^{c}.$ The independence
of the other pairs follows by applying the same reasoning.

(ii) We have to prove that, for any non-empty finite subset $J$ of
$I,$ the following holds
\begin{equation}
P\left(\bigcap\limits_{j\in J}B_{j}\right)=\prod\limits_{j\in J}P\left(B_{j}\right).\label{eq:product_of_j_independent_event}
\end{equation}

Define
\[
\begin{array}{ccccc}
J_{1}=I_{1}\cap J & \,\,\,\, & \text{and} & \,\,\,\, & J_{2}=I_{2}\cap J.\end{array}
\]

We assume that both sets are non-empty---proofs in other cases are
treated by a similar manner. Then
\begin{align*}
\bigcap\limits_{j\in J}B_{j} & =\left(\bigcap\limits_{j\in J_{1}}A_{j}\right)\cap\left(\bigcap\limits_{j\in J_{2}}A^{c}_{j}\right)\\
 & =\left(\bigcap\limits_{j\in J_{1}}A_{j}\right)\cap\left(\bigcup_{j\in J_{2}}A_{j}\right)^{c}.
\end{align*}

Using $\refpar{eq:prob_A_inter_B^C},$ 
\begin{equation}
P\left(\bigcap\limits_{j\in J}B_{j}\right)=P\left(\bigcap\limits_{j\in J_{1}}A_{j}\right)-P\left(\left(\bigcap\limits_{j\in J_{1}}A_{j}\right)\cap\left(\bigcup_{j\in J_{2}}A_{j}\right)\right).\label{eq:P_inter_B_j}
\end{equation}

Since
\[
\left(\bigcap\limits_{j\in J_{1}}A_{j}\right)\cap\left(\bigcup_{j\in J_{2}}A_{j}\right)=\bigcup_{k\in J_{2}}\left(A_{k}\cap\left(\bigcap_{j\in J_{1}}A_{j}\right)\right),
\]
we have, by applying the Poincaré formula to this last union,
\begin{multline*}
P\left(\bigcup_{k\in J_{2}}\left(A_{k}\cap\left(\bigcap_{j\in J_{1}}A_{j}\right)\right)\right)=\sum_{\substack{K\subset J_{2}\\
\left|K\right|=1
}
}P\left(\bigcap_{k\in K}A_{k}\cap\left(\bigcap_{j\in J_{1}}A_{j}\right)\right)\,\,\,\,\,\,\,\,\,\,\,\,\,\,\,\,\,\,\,\,\,\,\,\,\,\,\,\,\,\,\,\,\,\,\,\,\,\,\,\,\\
-\sum_{\substack{K\subset J_{2}\\
\left|K\right|=2
}
}P\left(\bigcap_{k\in K}A_{k}\cap\left(\bigcap_{j\in J_{1}}A_{j}\right)\right)+\dots+\left(-1\right)^{\left|J_{2}\right|-1}P\left(\bigcap_{k\in K}A_{k}\cap\left(\bigcap_{j\in J_{1}}A_{j}\right)\right).
\end{multline*}

Which after insertion in $\refpar{eq:P_inter_B_j}$ and simplifications
gives
\begin{multline*}
P\left(\bigcap\limits_{j\in J}B_{j}\right)=P\left(\bigcap\limits_{j\in J_{1}}A_{j}\right)-\sum\limits_{\substack{K\subset J_{2}\\
\left|K\right|=1
}
}P\left(\bigcap\limits_{j\in K\cup J_{1}}A_{j}\right)+\sum\limits_{\substack{K\subset J_{2}\\
\left|K\right|=2
}
}P\left(\bigcap\limits_{j\in K\cup J_{1}}A_{j}\right)\\
-\dots+\left(-1\right)^{\left|J_{2}\right|}P\left(\bigcap\limits_{j\in J}A_{j}\right).
\end{multline*}

Since the $\left(A_{i}\right)_{i\in I}$ are mutually independent,
\begin{multline*}
P\left(\bigcap\limits_{j\in J}B_{j}\right)=\prod_{j\in J_{1}}P\left(A_{j}\right)-\sum\limits_{\substack{K\subset J_{2}\\
\left|K\right|=1
}
}\prod_{j\in K\cup J_{1}}P\left(A_{j}\right)+\sum\limits_{\substack{K\subset J_{2}\\
\left|K\right|=2
}
}\prod_{j\in K\cup J_{1}}P\left(A_{j}\right)\\
-\dots+\left(-1\right)^{\left|J_{2}\right|}\prod_{j\in J}P\left(A_{j}\right)
\end{multline*}

By factoring,
\begin{align*}
P\left(\bigcap\limits_{j\in J}B_{j}\right) & =\left[\prod_{j\in J_{1}}P\left(A_{j}\right)\right]\\
 & \times\left[1-\sum\limits_{\substack{j\in J_{2}}
}P\left(A_{j}\right)+\sum\limits_{\substack{K\subset J_{2}\\
\left|K\right|=2
}
}\prod_{j\in K}P\left(A_{j}\right)-\dots+\left(-1\right)^{\left|J_{2}\right|}\prod_{j\in J_{2}}P\left(A_{j}\right)\right].
\end{align*}

Recognizing the second factor, as the development of a product,
\begin{align*}
P\left(\bigcap\limits_{j\in J}B_{j}\right) & =\left[\prod_{j\in J_{1}}P\left(A_{j}\right)\right]\times\left[\prod_{j\in J_{2}}\left(1-P\left(A_{j}\right)\right)\right].
\end{align*}

Since
\[
P\left(B_{j}\right)=\begin{cases}
P\left(A_{j}\right), & \text{if }j\in J_{1},\\
1-P\left(A_{j}\right), & \text{if }j\in J_{2},
\end{cases}
\]
we conclude that
\[
P\left(\bigcap_{j\in J}B_{j}\right)=\prod_{j\in J}P\left(B_{j}\right).
\]

This proves that the family $\left(B_{i}\right)_{i\in I}$ is independent.

\end{proof}

\subsection{Independent Random Variables}

We now define the notion of \textbf{independent random variables}\index{independent random variables}\mindex{random variables ! independent}.

\begin{definition}{Independence of Random Variables}{}All considered
random variables are defined on the same probabilized space $\left(\Omega,\mathcal{A},P\right).$

(i) Let $X_{1}$ and $X_{2}$ be two random variables taking values
in probabilizable spaces $\left(E_{1},\mathcal{E}_{1}\right)$ and
$\left(E_{2},\mathcal{E}_{2}\right),$ respectively.

The random variables $X_{1}$ and $X_{2}$ are \index{independent}\mindex{random variables ! independent}\textbf{independent}
if, for every $A_{1}\in\mathcal{E}_{1}$ and, for every $A_{2}\in\mathcal{E}_{2},$
the events $\left(X_{1}\in A_{1}\right)$ and $\left(X_{2}\in A_{2}\right)$
are independent.

(ii) This notion extends to any family of random variables as follows:

Let $\left(X_{i}\right)_{i\in I}$ be a family of random variables
taking values in the respective probabilizable spaces $\left(E_{i},\mathcal{E}_{i}\right)$
where $i\in I.$

The random variables $X_{i}$ are \index{independent}\mindex{random variables ! family ! independent}\textbf{independent}
if, for every family of sets $\left(A_{i}\right)_{i\in I}$ such that
for every $i\in I,$ $A_{i}\in\mathcal{E}_{i},$ the events $\left(X_{i}\in A_{i}\right),\,i\in I,$
are independent.

\end{definition}

\begin{proposition}{Necessary and Sufficient Condition for a Family of Random Variables to Be Independent}{}\label{Prop:independence random variables}Let
$I$ be a non-empty finite set. Let $\left(X_{i}\right)_{i\in I}$
be a family of random variables taking values in the respective probabilizable
spaces $\left(E_{i},\mathcal{E}_{i}\right),$ where $i\in I.$

The random variables $X_{i}$ are independent if and only if, for
every family of sets $\left(A_{i}\right)_{i\in I}$ with $A_{i}\in\mathcal{E}_{i}$
for every $i\in I,$ the following relation holds\boxeq{
\begin{equation}
P\left[\bigcap_{i\in I}\left(X_{i}\in A_{i}\right)\right]=\prod\limits_{i\in I}P\left(X_{i}\in A_{i}\right).\label{eq:independence_random_variables_I}
\end{equation}
}

\end{proposition}

\begin{proof}{}{}The necessity of the condition is straightforward.

To prove sufficiency, assume that the relation $\refpar{eq:independence_random_variables_I}$
holds for the family $\left(A_{i}\right)_{i\in I}.$

Let $\left(B_{i}\right)_{i\in I}$ be any family of sets with $B_{i}\in\mathcal{E}_{i}$
for every $i\in I,$ and let $J$ be any subset of $I.$

Define a new family $\left(A_{i}\right)_{i\in I}$ by
\[
A_{i}=\begin{cases}
B_{i}, & \text{if }i\in J,\\
E_{i}, & \text{otherwise.}
\end{cases}
\]

Since for every $i\notin J,$
\[
\left(X_{i}\in A_{i}\right)=\left(X_{i}\in E_{i}\right)=\Omega,
\]
it follows that
\[
P\left(X_{i}\in A_{i}\right)=1.
\]

Thus, applying the relation $\refpar{eq:independence_random_variables_I},$
we can write
\begin{align*}
P\left[\bigcap\limits_{j\in J}\left(X_{j}\in B_{j}\right)\right] & =P\left[\bigcap\limits_{i\in I}\left(X_{i}\in A_{i}\right)\right]\\
 & =\prod\limits_{i\in I}P\left(X_{i}\in A_{i}\right)\\
 & =\prod\limits_{j\in J}P\left(X_{j}\in B_{j}\right).
\end{align*}

Since $J$ is an arbitrary subset of $I,$ this proves that the events
$\left(X_{i}\in B_{i}\right),\,i\in I,$ are independent; therefore,
the random variables $X_{i},\,i\in I,$ are independent.

\end{proof}

In the case where the random variables $X_{i},\,i\in I$ are discrete,
we can subtitute in Proposition \ref{Prop:independence random variables}
the events $\bigcap\limits_{i\in I}\left(X_{i}\in A_{i}\right)$ by
the events $\bigcap\limits_{i\in I}\left(X_{i}=x_{i}\right).$ This
is the content of Proposition \ref{Prop:independent random var and n uplet}
below, stated for the sake of simplicity, in the case $I=\left\{ 1,\dots,n\right\} .$ 

We introduce the notation\boxeq{
\[
\left(X_{1}=x_{1},\dots,X_{n}=x_{n}\right)
\]
}to represent the event 
\[
\bigcap^{n}_{i=1}\left(X_{i}=x_{i}\right).
\]
 This notation will be used throughout the remainder of this work.

\begin{proposition}{Checking Random Variables Independence via n-tuples}{}\label{Prop:independent random var and n uplet}

Let $\left(X_{i}\right)_{i\in\left\{ 1,\dots,n\right\} }$ be a family
of discrete random variables taking values respectively in the sets
$E_{i},\,i\in\left\{ 1,\dots,n\right\} .$

The random variables $X_{i}$ are \textbf{independent} if and only
if, for every $n-$tuples $\left(x_{1},\dots,x_{n}\right)\in X_{1}\left(\Omega\right)\times\dots\times X_{n}\left(\Omega\right),$\boxeq{
\begin{equation}
P\left(X_{1}=x_{1},\dots,X_{n}=x_{n}\right)=P\left(X_{1}=x_{1}\right)\dots P\left(X_{n}=x_{n}\right).\label{eq:discrete_random_var_indep}
\end{equation}
}

\end{proposition}

\begin{proof}{}{}

\textbf{Necessity of the Condition}

Let $x_{1}\in X_{1}\left(\Omega\right),\dots,x_{n}\in X_{n}\left(\Omega\right).$ 

Since the random variables $X_{1},\dots,X_{n}$ are independent, applying
the relation $\refpar{eq:independence_random_variables_I}$ to the
sets $A_{i}=\left\{ x_{i}\right\} $ for $i\in\left\{ 1,\dots,n\right\} $
gives the relation $\refpar{eq:discrete_random_var_indep}.$

\textbf{Sufficiency of the Condition}

Assume that the relation $\refpar{eq:discrete_random_var_indep}$
holds for every $n-$uple $\left(x_{1},\dots,x_{n}\right)\in X_{1}\left(\Omega\right)\times\dots\times X_{n}\left(\Omega\right).$

Let $A_{i}\in\mathcal{E}_{i}$ be any set for $i\in\left\{ 1,\dots,n\right\} .$
Then
\[
\left(X_{i}\in A_{i}\right)=\biguplus\limits_{x_{i}\in X_{i}\left(\Omega\right)\cap A_{i}}\left(X_{i}=x_{i}\right).
\]

Using the distributivity in union and intersection operations, we
get
\[
P\left(\bigcap\limits_{i\in I}\left(X_{i}\in A_{i}\right)\right)=P\left[\biguplus\limits_{\left(x_{1},\dots,x_{n}\right)\in\left(X_{1}\left(\Omega\right)\cap A_{1}\right)\times\dots\times\left(X_{n}\left(\Omega\right)\cap A_{n}\right)}\left(\bigcap\limits^{n}_{i=1}\left(X_{i}=x_{i}\right)\right)\right].
\]

Since the random variables $X_{i}$ are discrete, the sets $X_{i}\left(\Omega\right)\cap A_{i}$
are countable. By $\sigma-$additivity of $P,$ we obtain 
\[
P\left(\bigcap\limits_{i\in I}\left(X_{i}\in A_{i}\right)\right)=\sum_{\left(x_{1},\dots,x_{n}\right)\in\left(X_{1}\left(\Omega\right)\cap A_{1}\right)\times\dots\times\left(X_{n}\left(\Omega\right)\cap A_{n}\right)}P\left[\bigcap\limits^{n}_{i=1}\left(X_{i}=x_{i}\right)\right].
\]

Using the relation $\refpar{eq:discrete_random_var_indep},$ we obtain
\[
P\left(\bigcap\limits_{i\in I}\left(X_{i}\in A_{i}\right)\right)=\sum_{\left(x_{1},\dots,x_{n}\right)\in\left(X_{1}\left(\Omega\right)\cap A_{1}\right)\times\dots\times\left(X_{n}\left(\Omega\right)\cap A_{n}\right)}\prod^{n}_{i=1}P\left(X_{i}=x_{i}\right).
\]

Applying the Fubini property in the case of direct product of families
\[
P\left(\bigcap\limits_{i\in I}\left(X_{i}\in A_{i}\right)\right)=\prod^{n}_{i=1}\sum_{x_{i}\in X_{i}\left(\Omega\right)\cap A_{i}}P\left(X_{i}=x_{i}\right),
\]
which can be reformulated as
\[
P\left(\bigcap\limits_{i\in I}\left(X_{i}\in A_{i}\right)\right)=\prod^{n}_{i=1}P\left(X_{i}\in A_{i}\right).
\]

Thus, the random variables $X_{1},X_{2},\dots,X_{n}$ are independent.

\end{proof}

\begin{remark}{}{}A useful exercise in order to assimilate this proof,
which, despite its quite abstract nature, is of common usage, is to
verify it explicitly for two random variables taking values in $\mathbb{N}.$

\end{remark}

\subsection{Functions of Random Variables (Discrete Case)}

Let $X$ be a random variable defined on a probabilized space $\left(\Omega,\mathcal{A},P\right)$
taking values in a probalizable space $\left(E,\mathcal{E}\right).$
Let $f$ be a random variable defined on $\left(E,\mathcal{E}\right),$
taking values in a space $\left(F,\mathcal{F}\right).$ Then, the
composed function $f\circ X$ is a random variable defined from $\left(\Omega,\mathcal{A},P\right)$
onto $\left(F,\mathcal{F}\right).$

Indeed, for every $B\in\mathcal{F},$ the set
\[
\left(f\circ X\right)^{-1}\left(B\right)=X^{-1}\left(f^{-1}\left(B\right)\right)
\]
belongs to $\mathcal{A}.$ 

Instead of writing $f\circ X,$ this random variable is often denoted
as $f\left(X\right):$ it is referred to as the \textbf{function $f$
of the random variable\mindex{function!random variable}\mindex{random variable ! function}}
$X.$

If $X$ is discrete---that is, if the set of values of $X$ is countable---and
if, for every $x\in E,$ the set $X^{-1}\left(\left\{ x\right\} \right)$
belongs to $\mathcal{A},$ then $f\circ X$ is a random variable,
whatever the function $f:E\rightarrow F$ is. For instance, if $X$
is a discrete real-valued random variable, the functions $\sin X$
and $X^{2},$ defined on $\Omega$ are also discrete random variables.

Still in the discrete case\footnote{To avoid complex technical difficulties, the interested reader may
refer to Chapter \ref{chap:PartIIChap9} of Part \ref{part:Deepening-Probability-Theory}
for a more general treatment.}, we can consider functions of multiple random variables. 

Let $X_{1},\dots,X_{n}$ be discrete random variables on $\left(\Omega,\mathcal{A},P\right),$
taking values in $\left(E_{1},\mathcal{P}\left(E_{1}\right)\right),\dots,\left(E_{n},\mathcal{P}\left(E_{n}\right)\right),$
respectively. 

Let $F$ be any set, and consider a function
\[
f:E_{1}\times\dots\times E_{n}\rightarrow F.
\]

Then the composition of the function
\[
\omega\mapsto\left(X_{1}\left(\omega\right),\dots,X_{n}\left(\omega\right)\right)
\]
mapping $\Omega$ to $E_{1}\times\dots\times E_{n},$ followed by
the function $f,$ is denoted $f\left(X_{1},\dots,X_{n}\right),$
and is a discrete random variable taking values in $F.$

The following result states that ``functions of independent random
variables are independent'': this result is of constant usage, as
it often allows us to establish the independence of random variables
without explicit calculation.

\begin{proposition}{Independence of Functions of Independent Random Variables}{}Let
$\left(X_{i}\right)_{i\in I}$ be any family of random variables defined
on the probabilized space $\left(\Omega,\mathcal{A},P\right)$ taking
values respectively in probabilizable spaces $\left(E_{i},\mathcal{E}_{i}\right),\,i\in I.$
Suppose that the random variables $X_{i}$ are \textbf{independent}. 

For each $i\in I,$ let $f_{i}$ be a random variable defined on the
probabilizable space $\left(E_{i},\mathcal{E}_{i}\right)$ taking
values in the probabilizable space $\left(F_{i},\mathcal{F}_{i}\right).$

Then the random variables $f_{i}\left(X_{i}\right)$ are also independent.

\end{proposition}

\begin{proof}{}{}

It suffices to remark that, for every $i\in I$ and for every $B_{i}\in\mathcal{F}_{i},$
we have
\[
\left(f_{i}\left(X_{i}\right)\in B_{i}\right)=\left(f_{i}\circ X_{i}\right)^{-1}\left(B_{i}\right)=X^{-1}_{i}\left(f^{-1}_{i}\left(B_{i}\right)\right)=\left(X_{i}\in f^{-1}_{i}\left(B_{i}\right)\right),
\]
and that
\[
f^{-1}_{i}\left(B_{i}\right)\in\mathcal{E}_{i}.
\]

The random variables $X_{i}$ being independent, it ensures that the
random variables $f_{i}\left(X_{i}\right)$ are independent.

\end{proof}

\begin{corollary}{Independence of Multivariate Functions of Independent Variables}{indep_multivar_functions_indep}

Let $X_{1},\dots,X_{n}$ be $n$ discrete random variables defined
on the same probabilized space $\left(\Omega,\mathcal{A},P\right),$
with respective values in $E_{1},\dots,E_{n}.$ 

Let
\[
\begin{array}{ccccc}
f_{1}:E_{1}\times\dots\times E_{k}\rightarrow F_{1} & \,\,\,\, & \text{and} & \,\,\,\, & f_{2}:E_{k+1}\times\dots\times E_{n}\rightarrow F_{2}\end{array}
\]
be any two functions. 

If the random variables $X_{1},\dots,X_{n}$ are independent, then
the random variables $f_{1}\left(X_{1},\dots,X_{k}\right)$ and $f_{2}\left(X_{k+1},\dots,X_{n}\right)$
are independent.

\end{corollary}

\begin{proof}{}{}

This follows immediately from the previous proposition.

\end{proof}

\begin{example}{}{}

We can assert without any calculation, that if $X_{1},X_{2,},X_{3}$
and $X_{4}$ are four independent random variables taking values in
$\mathbb{Z},$ then the random variables
\begin{itemize}
\item $X_{1}+X_{2}$ and $X_{3}+X_{4}$ are independent, 
\item $X_{1}+X_{4}$ and $X_{2}+X_{3}$ are independent---it suffices to
use the corollary by reindexing the $X_{i}$---
\item Or any variables, such as $\left(-1\right)^{X_{2}}\sin\left(\pi X_{3}\right)$
and $X^{5}_{1}+X^{1000}_{4}$ are also independent.
\end{itemize}
\end{example}

\section{Law of the Sum of Independent Random Variables}\label{sec:Law-of-the-Sum-of-Inpendent-RV}

Very often there is the need to determine the law of a sum of random
variables. When the random variables $X_{1}$ and $X_{2}$ are independent,
the law of their sum can be calculated from the laws of $X_{1}$ and
$X_{2}.$ We present the calculation method, known as the \textbf{\index{convolution}convolution}
method, in the case of discrete random variables. 

\begin{proposition}{Law of the Sum of Two Discrete Independent Random  Variables}{}Let
$X_{1}$ and $X_{2}$ be two \textbf{discrete} and \textbf{independent}
random variables defined on the probabilized space $\left(\Omega,\mathcal{A},P\right)$,
taking values in the same probabilizable space $\left(E,\mathcal{E}\right),$
where $E$ is a subset of an abelian group, stable under addition\footnotemark.

Then, the random variable $X_{1}+X_{2}$ is discrete and its law is
given by either of the two following relations\boxeq{
\begin{equation}
\forall x\in E,\,\,\,\,P\left(X_{1}+X_{2}=x\right)=\sum\limits_{x_{1}\in E}P\left(X_{1}=x_{1}\right)P\left(X_{2}=x-x_{1}\right)\label{eq:prob_of_X_1+X_2_rel_1}
\end{equation}
}or equivalently,\boxeq{
\begin{equation}
\forall x\in E,\,\,\,\,P\left(X_{1}+X_{2}=x\right)=\sum\limits_{x_{2}\in E}P\left(X_{1}=x-x_{2}\right)P\left(X_{2}=x_{2}\right)
\end{equation}
}

\end{proposition}

\footnotetext{Typical examples of $E$ include groups such as $\mathbb{Z},\mathbb{Z}^{d},\mathbb{R},$
and so on..., or other common cases like $\mathbb{N}$ and $\mathbb{R}^{+}.$}

\begin{proof}{}{}We assert previously, without prooving it, that
any function $f\left(X_{1},\dots,X_{n}\right)$ of discrete random
variables is also a discrete random variable. We now give a proof
in the case of $X_{1}+X_{2}.$

The set $\left(X_{1}+X_{2}\right)\left(\Omega\right)$ is countable
because it is the image of the set $X_{1}\left(\Omega\right)\times X_{2}\left(\Omega\right),$
which is countable as the Cartesian product of two countable sets,
under the mapping $\left(x_{1},x_{2}\right)\mapsto x_{1}+x_{2}$.

Moreover, for every $x\in E,$
\begin{equation}
\left(X_{1}+X_{2}=x\right)=\biguplus\limits_{x_{1}\in X_{1}\left(\Omega\right)}\left(X_{1}=x_{1}\right)\cap\left(X_{2}=x-x_{1}\right),\label{eq:X_1_plus_X_2_equal_x}
\end{equation}
which shows that the set $\left(X_{1}+X_{2}=x\right)$ effectively
belongs to $\mathcal{A}.$

Having established this, we obtain the relation $\refpar{eq:prob_of_X_1+X_2_rel_1}$
immediately from $\refpar{eq:X_1_plus_X_2_equal_x}$ by using the
$\sigma-$additivity of $P$ and the independence of the random variables
$X_{1}$ and $X_{2}.$

Similarly, we prove the second relation.

\end{proof}

\begin{remark}{}{}The relation $\refpar{eq:prob_of_X_1+X_2_rel_1}$
can also be expressed in terms of random variables probability laws
as follows\boxeq{
\begin{equation}
\begin{array}{ccc}
\forall x\in E, & \,\,\,\, & P_{X_{1}+X_{2}}\left(\left\{ x\right\} \right)=\sum\limits_{x_{1}\in E}P_{X_{1}}\left(\left\{ x_{1}\right\} \right)P_{X_{2}}\left(\left\{ x-x_{1}\right\} \right).\end{array}\label{eq:proba_law_of_X_1+X_2}
\end{equation}
}

\end{remark}

\begin{definition}{Convolution of Two Probabilities}{} 

Let $X_{1}$ and $X_{2}$ be two discrete and independent random variables
defined on the probabilized space $\left(\Omega,\mathcal{A},P\right)$,
taking values in the same probabilizable space $\left(E,\mathcal{E}\right),$
where $E$ is a subset of an abelian group, stable under addition\footnotemark,
with respective discrete laws $P_{X_{1}}$ and $P_{X_{2}}.$

The product of convolution of $P_{X_{1}}$ and $P_{X_{2}}$ or more
simply the \textbf{\index{convolution}convolution} of the probabilities
$P_{X_{1}}$ and $P_{X_{2}},$ denoted as $P_{X_{1}}\ast P_{X_{2}},$
is defined, for every $x\in E,$ by

\boxeq{
\[
\left(P_{X_{1}}\ast P_{X_{2}}\right)\left(\left\{ x\right\} \right)=\sum\limits_{x_{1}\in E}P_{X_{1}}\left(\left\{ x_{1}\right\} \right)P_{X_{2}}\left(\left\{ x-x_{1}\right\} \right).
\]
}

\end{definition}

With this notation, we then have\boxeq{
\[
P_{X_{1}+X_{2}}=P_{X_{1}}\ast P_{X_{2}}.
\]
}The convolution is a commutative operation between probability laws.

We now provide two examples for computing the law of the sum of independent
random variables.

\begin{example}{Triangular Law}{}

Let $n\geqslant1$ be an integer. Let $X$ and $Y$ be two independent
random variables, defined on the probabilized space $\left(\Omega,\mathcal{A},P\right),$
both following the uniform law on the set $\left\llbracket 0,n\right\rrbracket $.

Study the law of the sum $Z=X+Y$ and determine a property of symmetry.

\end{example}

\begin{solutionexample}{}{}The random variables $X$ and $Y$ take
values in $\mathbb{N},$ which is a stable subset by addition in the
$\mathbb{Z}$ group. Thus, for every $k\in\mathbb{N},$
\[
P\left(Z=k\right)=\sum\limits_{j\in\mathbb{N}}P\left(X=j\right)P\left(Y=k-j\right).
\]

Since
\[
P\left(X=j\right)=P\left(Y=j\right)=\begin{cases}
\dfrac{1}{n+1}, & \text{if }j\in\left\llbracket 0,n\right\rrbracket ,\\
0, & \text{otherwise,}
\end{cases}
\]
it follows that:
\begin{itemize}
\item If $0\leqslant k\leqslant n,$ then
\[
P\left(Z=k\right)=\sum\limits^{k}_{j=0}\left(\dfrac{1}{n+1}\right)^{2}=\dfrac{k+1}{\left(n+1\right)^{2}}.
\]
\item If $n<k\leqslant2n,$ then
\[
P\left(Z=k\right)=\sum\limits^{n}_{j=k-n}\left(\dfrac{1}{n+1}\right)^{2}=\dfrac{2n-k+1}{\left(n+1\right)^{2}}.
\]
\end{itemize}
Thus,
\[
\sum\limits^{2n}_{j=0}P\left(Z=k\right)=1,
\]
ensuring that
\[
P_{Z}\left(\left\{ 0,1,2,\dots,2n\right\} \right)=1.
\]

If $k$ and $l$ are symmetric with respect to $n,$ i.e. if $l=2n-k,$
then
\[
P\left(Z=l\right)=\dfrac{2n-\left(2n-l\right)+1}{\left(n+1\right)^{2}}
\]
and thus
\[
P\left(Z=l\right)=P\left(Z=k\right)=\dfrac{k+1}{\left(n+1\right)^{2}}.
\]
We then say that the law of $Z$ is symmetric with respect to $n.$

In Figure \ref{Fig:Triangular_law}, we give an example of the triangular
law for $n=7.$

\end{solutionexample}

\begin{figure}
\begin{center}\begin{tikzpicture}
	\definecolor{orange}{HTML}{FFA200}
	\definecolor{purple}{HTML}{9500FF}
  \begin{axis}[
		axis x line=center,
		axis y line=center,
		xtick = {0,1,2,3,4,5,6,7,8,9,10,11,12,13,14},
		ytick = {0,0.05,0.1},
		xmin=-0.5,
		ymax=14.9,
		ymin=-0.02,
		ymax=0.15,
		xlabel={$k$},
		ylabel={$P\left(Z=k\right)$},
every axis x label/.style={
    at={(ticklabel* cs:1.02)},
    anchor=west,
},
every axis y label/.style={
    at={(ticklabel* cs:1.02)},
    anchor=south,
},
	]
\addplot[white, samples=80, smooth, thick, domain=-0.5:14.5 ]
      {0.15};
	\pgfplotsinvokeforeach {0,1,2,3,4,5,6,7} {
    	\draw[orange, line width=3] (axis cs: #1,0)
      -- (axis cs: #1,{(#1+1)/64});
  }
\pgfplotsinvokeforeach {8,9,10,11,12,13,14} {
    	\draw[orange, line width=3] (axis cs: #1,0)
      -- (axis cs: #1,{(14-#1+1)/64});
  }
  \end{axis}

\end{tikzpicture}\end{center}

\caption{Triangular law, case $n=7$}
\label{Fig:Triangular_law}
\end{figure}

\begin{example}{Sum of Two Independent Random Variables Following the Poisson Law}{sum_rv_poisson_law} 

Let $X_{1}$ and $X_{2}$ be two independent random variables, defined
on the probabilized space $\left(\Omega,\mathcal{A},P\right),$ following
the Poisson laws with respective parameters $\lambda_{1}$ and $\lambda_{2},$
both positive.

Study the law of $Z=X_{1}+X_{2}.$\end{example}

\begin{solutionexample}{}{}Since the random variables $X_{1}$ and
$X_{2}$ take values in $\mathbb{N},$ the same holds for the random
variable $Z.$ 

The convolution relation gives, since $X_{1}$ and $X_{2}$ are independent
\[
\begin{array}{ccc}
\forall k\in\mathbb{N}, & \,\,\,\, & P\left(Z=k\right)=\sum\limits_{j\in\mathbb{N}}P\left(X_{1}=j\right)P\left(X_{2}=k-j\right).\end{array}
\]

Since, for $i=1,2,$ and for every $j\in\mathbb{N},$
\[
P\left(X_{i}=j\right)=\exp\left(-\lambda_{i}\right)\dfrac{\lambda^{j}_{i}}{j!},
\]
we obtain that, for every $k\in\mathbb{N},$
\begin{align*}
P\left(Z=k\right) & =\sum^{k}_{j=0}\left(\exp\left(-\lambda_{1}\right)\dfrac{\lambda^{j}_{1}}{j!}\right)\left(\exp\left(-\lambda_{2}\right)\dfrac{\lambda^{k-j}_{2}}{\left(k-j\right)!}\right)\\
 & =\dfrac{\exp\left(-\left(\lambda_{1}+\lambda_{2}\right)\right)}{k!}\sum^{k}_{j=0}\binom{k}{j}\lambda^{j}_{1}\lambda^{k-j}_{2}.
\end{align*}
Hence,

\boxeq{
\[
P\left(Z=k\right)=\dfrac{\exp\left(-\left(\lambda_{1}+\lambda_{2}\right)\right)}{k!}\left(\lambda_{1}+\lambda_{2}\right)^{k}.
\]
}

\end{solutionexample}

Thus, we established the following proposition.

\begin{proposition}{Convolution of Two Poisson Laws}{}The convolution
of two Poisson laws $\mathcal{P}\left(\lambda_{1}\right)$ and $\mathcal{P}\left(\lambda_{2}\right)$
with nonnegative parameters $\lambda_{1}$ and $\lambda_{2},$ is
a Poisson law with parameter $\lambda_{1}+\lambda_{2}.$ 

\end{proposition}

We say that the family of the Poisson law is \textbf{\index{stable under convolution}\mindex{convolution ! stable}stable}
\textbf{under convolution}.

\section{Independence and Cartesian Products: Construction of a Model}\label{sec:Independence-and-Cartesian}

Returning to the introductive example, that is the modelling of the
consecutive draws from the same die, we naturally chose as the fundamental
space, or outcome space, the Cartesian product 
\[
\Omega=\Omega_{1}\times\Omega_{2},
\]
where $\Omega_{1}=\Omega_{2}=\left\{ 1,2,\dots,6\right\} ,$ and we
endowed it with the uniform probability $P$ on the probabilizable
space $\left(\Omega,\mathcal{P}\left(\Omega\right)\right).$ We then
showed that the projection random variables $X_{1}$ and $X_{2},$
which separately describe all information available from each die
roll, are independent random variables following the uniform law on
$\left\{ 1,2,\dots,6\right\} .$ 

However, from the perspective of how the game is played, it is actually
the two last properties---independence and uniform law---that are
``natural.'' If we adopted this viewpoint in constructing the model,
we would have been lead to set the following ``inverse'' problem: 

\begin{leftbar}

``\textit{Construct a probabilized space $\left(\Omega,\mathcal{A},P\right)$
on which we can define two independent random variables both following
the uniform law on $\left\{ 1,2,\dots,6\right\} .$}'' 

\end{leftbar}

As we are going to see it in a more general setting, a solution to
this problem is precisely the probabilized space we have qualified
as ``natural''.

\subsubsection*{Problem statement}

We ask ourselves the following problem.

\begin{leftbar}

Given a family $\left(E_{i},\mathcal{E}_{i}\right),\,i\in I$ of probabilizable
spaces, along a law of given probability $P_{i}$ on each of those
spaces, could we construct a probabilized space $\left(\Omega,\mathcal{A},P\right)$
and a family of random variables $\left(X_{i}\right)_{i\in I}$ defined
on this space, taking values respectively in $\left(E_{i},\mathcal{E}_{i}\right),$
such that:
\begin{itemize}
\item The random variables $X_{i}$ are independent, and
\item For every $i\in I,$ the random variable $X_{i}$ follows the law
$P_{i}$?
\end{itemize}
\end{leftbar}

This is a fundamental problem for the construction of models. Here,
we solve it within a restricted framework:
\begin{itemize}
\item The index set $I$ is finite.
\item The sets $E_{i}$ are countable. 
\end{itemize}
In this framework, a solution is provided by the following proposition.

\begin{proposition}{Product Probability of Probabilities}{product_probability}Let
$I$ be $\left\llbracket 1,n\right\rrbracket .$ For each $i\in I,$
let $E_{i}$ be a countable set, and let $P_{i}$ be a probability
on the probabilizable space $\left(E_{i},\mathcal{P}\left(E_{i}\right)\right),$
generated by the probability germ $p_{i.}$ 

Consider the Cartesian product 
\[
\Omega=\prod^{n}_{i=1}E_{i},
\]
 and let $X_{i}$ be the projection of $\Omega$ onto $E_{i}.$ 

Define the function $p$ on $\Omega$ by the relation\boxeq{
\[
\begin{array}{ccc}
\forall\left(\omega_{1},\omega_{2},\dots,\omega_{n}\right)\in\Omega, & \,\,\,\, & p\left(\omega_{1},\omega_{2},\dots,\omega_{n}\right)=\prod\limits^{n}_{i=1}p_{i}\left(\omega_{i}\right)\end{array}.
\]
}

Then, $p$ is the germ of a probability $P$ on $\left(\Omega,\mathcal{P}\left(\Omega\right)\right),$
called the \index{product probability}\textbf{product probability}
of the probabilities $P_{i}.$ 

The probabilized space $\left(\Omega,\mathcal{P}\left(\Omega\right),P\right)$
satisfies the following properties:
\begin{itemize}
\item The random variables $X_{i}$ are independent\footnotemark.
\item The random variables $X_{i}$ follow the respective laws $P_{i}.$
\end{itemize}
\end{proposition}

\footnotetext{We also say that the random variables are independent
for the probability $P.$}

\begin{proof}{}{}

By the Fubini property,
\[
\sum\limits_{\left(\omega_{1},\omega_{2},\dots,\omega_{n}\right)\in\Omega}p\left(\omega_{1},\omega_{2},\dots,\omega_{n}\right)=\prod\limits^{n}_{i=1}\left(\sum\limits_{\omega_{i}\in E_{i}}p_{i}\left(\omega_{i}\right)\right)=1.
\]

This shows that the nonnegative function $p$ defines a probability
germ.

Moreover, for every fixed $i_{0}\in I$ and, for every $\omega_{i_{0}}\in E_{i_{0}},$
\[
\left(X_{i_{0}}=\omega_{i_{0}}\right)=\biguplus\limits_{\substack{\left(x_{1},\dots,x_{n}\right)\in\Omega\\
x_{i_{0}}=\omega_{i_{0}}
}
}\left\{ \left(x_{1},\dots,x_{n}\right)\right\} .
\]

Thus, again by the $\sigma-$additivity of $P$ and by the Fubini
property,
\begin{align*}
P\left(X_{i_{0}}=\omega_{i_{0}}\right) & =\sum_{\substack{\left(x_{1},\dots,x_{n}\right)\in\Omega\\
x_{i_{0}}=\omega_{i_{0}}
}
}p\left(x_{1},\dots,x_{n}\right)\\
 & =p_{i_{0}}\left(\omega_{i_{0}}\right)\prod\limits_{\substack{i\in I\\
i\neq i_{0}
}
}\left(\sum\limits_{x_{i}\in E_{i}}p_{i}\left(x_{i}\right)\right).
\end{align*}

Thus, it shows that, for every $\omega_{i_{0}}\in E_{i_{0}},$
\[
P\left(X_{i_{0}}=\omega_{i_{0}}\right)=p_{i_{0}}\left(\omega_{i_{0}}\right),
\]
that is the law of the random variable $X_{i_{0}}$ is $p_{i_{0}}.$ 

Finally, by definition of $P,$
\[
\begin{array}{ccc}
\forall\left(\omega_{1},\omega_{2},\dots,\omega_{n}\right)\in\Omega, & \,\,\,\, & P\left[\bigcap\limits^{n}_{i=1}\left(X_{i}=\omega_{i}\right)\right]=\end{array}\prod\limits^{n}_{i=1}P\left(X_{i}=\omega_{i}\right).
\]

This proves the independence of the random variables $X_{i}.$

\end{proof}

If we wish to model a coin-tossing game where we do not predetermine
the number of trials---for instance, if we stop to play after having
ten times a ``tail''---we have to take as outcome space the set
\[
\left\{ 0,1\right\} ^{\mathbb{N}^{\ast}},
\]
of infinite sequences of 0's and 1's---where 1 represents ``tails''
and 0 represents ``heads''. This set is uncountable---it is more
or less put in bijection with the real number interval $\left[0,1\right]$
via the diadic development 
\[
\left(x_{n}\right)\mapsto\sum^{+\infty}_{n=1}\dfrac{x_{n}}{2^{n}}.
\]
 Thus, we immediately move beyond the framework of discrete probabilized
spaces. We will revisit this question later in Part \ref{part:Deepening-Probability-Theory}
Chapter \ref{chap:PartIIChap10}, as solving this modelling problem
requires measure theory. 

Meanwhile, we admit the following proposition which generalizes to
the case where $I=\mathbb{N}$ in the following manner

\begin{proposition}{Existence of $\sigma-$algebra Where Variables are Independent}{}\label{Generalized proposition for independent random variables}For
each $i\in\mathbb{N},$ let $E_{i}$ be a countable set and $P_{i}$
be a probability on the probabilizable space $\left(E_{i},\mathcal{P}\left(E_{i}\right)\right),$
generated by the probability germ $p_{i}.$

Consider the Cartesian product 
\[
\Omega=\prod_{i\in\mathbb{N}}E_{i},
\]
and, for each $i\in\mathbb{N},$ let $X_{i}$ be the projection of
$\Omega$ on the $i-$th factor $E_{i}.$ 

Then, there exists a $\sigma-$algebra $\mathcal{A}$ on $\Omega$
and a probability $P$ on the probabilizable space $\left(\Omega,\mathcal{A}\right)$
such that:
\begin{itemize}
\item The random variables $X_{i}$ are independent;
\item The random variables $X_{i}$ follow the respective laws $P_{i}.$
\end{itemize}
\end{proposition}

\section{Geometric and Binomial Models}\label{sec:Geometric-and-Binomial}

Let $\left(\Omega,\mathcal{A},P\right)$ be a probabilized space,
and let $\left(A_{n}\right)_{n\in\mathbb{N}}$ be a sequence of \textbf{independent
events, each occuring with the same probability} $p$ with $p\in\left]0,1\right[.$ 

We define the random variables $N$ and $N^{\prime}$ taking values
in $\overline{\mathbb{N}}$ by
\begin{equation}
\begin{array}{ccc}
\forall\omega\in\Omega, & \,\,\,\, & \begin{cases}
N\left(\omega\right)=\inf\left\{ n\in\mathbb{N}:\,\omega\in A_{n}\right\} ,\\
N^{\prime}\left(\omega\right)=\inf\left\{ n\in\mathbb{N}^{\ast}:\,\omega\in A_{n}\right\} ,
\end{cases}\end{array}\label{eq:random_var_N_and_Nprime}
\end{equation}
with the convention that $\inf\emptyset=+\infty.$

For every integer $n\geqslant1,$ we define the random variable $S_{n}$
taking values in $\mathbb{N}$ by
\begin{equation}
S_{n}=\sum\limits^{n}_{i=1}\boldsymbol{1}_{A_{i}},\label{eq:S_n_sum_indicator}
\end{equation}
where $\boldsymbol{1}_{A_{i}}$ is the \textbf{\index{indicator function}indicator
function} of the event $A_{i}$ relatively to $\Omega.$

The data of $\left(\Omega,\mathcal{A},P\right)$ and of $\left(A_{n}\right)_{n\in\mathbb{N}}$
provide a model of the following situation: indefinitely the ``same''
random experiment is repeated---for instance, we can roll an infinite
time the same die, or rolling a new die each time, with exactly the
same properties as the first, the model being the same---and, we
are interested each time to the ``same'' outcome realization.

Then, the previously defined random variables have the following meaning:
\begin{itemize}
\item The integer $N\left(\omega\right)$ represents the \textbf{index of
the first experiment where the event occurs} assuming experiments
are numbered starting from $0.$
\item The integer $N^{\prime}\left(\omega\right)$ represents the \textbf{index
of the first experiment where the event is realized}, supposing the
experiments are numbered from 1.
\item The integer $S_{n}\left(\omega\right)$ represents \textbf{the number
of times the event occurs} during the experiments numbered from 1
to $n.$
\end{itemize}
We now study the laws of the random variables $N,$ $N^{\prime}$
and $S_{n}$ for $n\geqslant1.$ 

First, we start by a definition.

\begin{definition}{Bernoulli Variable}{}We call \index{variable of Bernoulli}\textbf{variable
of Bernoulli\footnotemark} all random variable $X$ with integer
values whom law is determined by
\[
\begin{array}{ccccc}
P\left(X=1\right)=p & \,\,\,\, & \text{and} & \,\,\,\, & P\left(X=0\right)=1-p,\end{array}
\]
---assuming generally that $0<p<1$\footnotemark---and we classicaly
set: $q=1-p$. This is called a \textbf{Bernoulli variable\index{Bernoulli variable}}
with parameter $p.$

\end{definition}

\addtocounter{footnote}{-1}

\footnotetext{In memory to Jacques Bernoulli (Bâle 1654-1705), first
of a dynasty of mathematicians and physicians. His work is dedicated
in particular to infinitesimal calculus, to differential equations
and to the isoperimeter problem as well as probability calculus.}

\stepcounter{footnote}

\footnotetext{It can be simpler to take as definition of a Bernoulli
variable, any random variable with values in $\left\{ 0,1\right\} .$
In this case, its law is a Bernoulli law possibly degenerated.}

\begin{remark}{}{}A Bernoulli variable is simply a random variable
taking values 0 or 1 with probability 1. 
\begin{itemize}
\item The indicator function $\boldsymbol{1}_{A}$ of an event $A$ is a
Bernoulli variable.
\item Conversely, a Bernoulli variable $X$ does not differ from the indicator
function of an event---the event $\left(X=1\right)$---except possibly
on a set of probability zero.
\end{itemize}
The random variable $S_{n}$ defined in $\refpar{eq:S_n_sum_indicator}$
is a sum of $n$ independent random variables of Bernoulli.

\end{remark}

\begin{proposition}{Law of $N$ and $N^\prime$}{}

Let $\left(\Omega,\mathcal{A},P\right)$ be a probabilized space,
and let $\left(A_{n}\right)_{n\in\mathbb{N}}$ be a sequence of \textbf{independent
events, each occuring with the same probability} $p$ with $p\in\left]0,1\right[.$ 
\begin{itemize}
\item The random variable $N,$ such that for every $\omega\in\Omega,$
$N\left(\omega\right)=\inf\left\{ n\in\mathbb{N}:\,\omega\in A_{n}\right\} ,$
follows a \textbf{geometric law on $\mathbb{N},$} denoted $\mathcal{G}_{\mathbb{N}}\left(p\right).$
\item The random variable $N^{\prime}$ such that that for every $\omega\in\Omega,$
$N^{\prime}\left(\omega\right)=\inf\left\{ n\in\mathbb{N}^{\ast}:\,\omega\in A_{n}\right\} ,$
follows a \textbf{geometric law} on $\mathbb{N}^{\ast},$ denoted
$\mathcal{G}_{\mathbb{N}^{\ast}}\left(p\right).$
\item Moreover, the law of $N^{\prime}$ is identical to the law of the
random variable $N+1.$
\end{itemize}
\end{proposition}

\begin{proof}{}{}

We have
\[
\left(N=0\right)=A_{0},
\]
and
\[
\begin{array}{ccc}
\forall k\geqslant1, & \,\,\,\, & \left(N=k\right)=\left(\bigcap\limits^{k-1}_{j=0}A^{c}_{j}\right)\cap A_{k}.\end{array}
\]

Since the events $A_{i},i\in\mathbb{N}$ are independent, the events
$A^{c}_{j},\,0\leqslant j\leqslant k-1,$ and $A_{k}$ are also independent.

Thus,

\[
\begin{array}{ccc}
\forall k\in\mathbb{N}, & \,\,\,\, & P\left(N=k\right)=p\left(1-p\right)^{k}.\end{array}
\]

Furthermore,
\[
P\left(N=+\infty\right)=1-P\left(N\in\mathbb{N}\right)=1-\sum\limits^{+\infty}_{k=0}p\left(1-p\right)^{k},
\]
and consequently $P\left(N=+\infty\right)=0.$ 

We effectively show that $P_{N}=\mathcal{G}_{\mathbb{N}}\left(p\right).$

Similarly, for the computation of the law of the random variable $N^{\prime}$
we observe that

\[
\left(N^{\prime}=1\right)=A_{1},
\]
and that
\[
\begin{array}{ccc}
\forall k\geqslant2, & \,\,\,\, & \left(N^{\prime}=k\right)=\left(\bigcap\limits^{k-1}_{j=1}A^{c}_{j}\right)\cap A_{k}.\end{array}
\]

Following the same computation, we obtain the law of $N^{\prime}.$

Finally, we study the law $P_{N^{\prime\prime}}$ of the random variable
$N^{\prime\prime}=N+1.$ We have, for every $k\in\mathbb{N}^{\ast},$
\[
P\left(N^{\prime\prime}=k\right)=P\left(N=k-1\right)=p\left(1-p\right)^{k-1},
\]

This shows that $P_{N^{\prime\prime}}=\mathcal{G}_{\mathbb{N}^{\ast}}\left(p\right),$
completing the proof.

\end{proof}

In Figure \ref{Fig:Geometric_laws}, we present the geometric laws
on $\mathbb{N}$ and $\mathbb{N}^{\ast}$ for $p=0.4.$

\begin{figure}
\begin{center}%
\begin{tabular}{cc}
\begin{tikzpicture}
	\definecolor{orange}{HTML}{FFA200}
	\definecolor{purple}{HTML}{9500FF}
  \begin{axis}[
		axis x line=center,
		axis y line=center,
		xtick = {0,1,2,3,4,5,6,7,8,9,10},
		ytick = {0,0.1,0.2,0.3,0.4},
		ymin=-0.02,
		ymax=0.45
	]
\addplot[white, samples=80, smooth, thick, domain=-0.5:10.5 ]{0.42};
\pgfplotsinvokeforeach {0,1,2,3,4,5,6,7,8,9,10} {
    	\draw[green!80!black, line width=3] (axis cs: #1,0)
      -- (axis cs: #1,{(0.4*0.6^#1)});
  }
  \end{axis}

\end{tikzpicture} &
\begin{tikzpicture}
	\definecolor{orange}{HTML}{FFA200}
	\definecolor{purple}{HTML}{9500FF}
  \begin{axis}[
		axis x line=center,
		axis y line=center,
		xtick = {0,1,2,3,4,5,6,7,8,9,10},
		ytick = {0,0.1,0.2,0.3,0.4},
		ymin=-0.02,
		ymax=0.45
	]
\addplot[white, samples=80, smooth, thick, domain=-0.5:10.5 ]{0.42};
\pgfplotsinvokeforeach {0,1,2,3,4,5,6,7,8,9,10} {
    	\draw[green!80!black, line width=3] (axis cs: #1+1,0)
      -- (axis cs: #1+1,{(0.4*0.6^#1)});
  }
  \end{axis}

\end{tikzpicture}\tabularnewline
Geometric law on $\mathbb{N}$ with $p=0.4$ &
Geometric law on $\mathbb{N}^{\ast}$ with $p=0.4$\tabularnewline
\end{tabular}\end{center}

\caption{Geometric laws on $\mathbb{N}$ and $\mathbb{N}^{\ast}$}
\label{Fig:Geometric_laws}
\end{figure}

\begin{proposition}{Law of $S_n$}{}

Let $\left(\Omega,\mathcal{A},P\right)$ be a probabilized space,
and let $\left(A_{n}\right)_{n\in\mathbb{N}}$ be a sequence of \textbf{independent
events, each occuring with the same probability} $p$ with $p\in\left]0,1\right[.$ 

The random variable $S_{n}$ defined for every integer $n\geqslant1,$
by
\[
S_{n}=\sum\limits^{n}_{i=1}\boldsymbol{1}_{A_{i}},
\]
where $\boldsymbol{1}_{A_{i}}$ is the \textbf{indicator function}
of the event $A_{i}$ relatively to $\Omega,$ follows a \textbf{binomial
law} $\mathcal{B}\left(n,p\right).$

\end{proposition}

\begin{proof}{}{}

First, we note that
\[
S_{n}\left(\Omega\right)\subset\left\llbracket 0,n\right\rrbracket .
\]

For every $k\in\left\llbracket 0,n\right\rrbracket ,$ we have
\[
\left(S_{n}=k\right)=\biguplus_{\substack{I\in\mathcal{P}\left(\left\llbracket 1,n\right\rrbracket \right)\\
\left|I\right|=k
}
}\left[\left(\bigcap\limits_{i\in I}A_{i}\right)\cap\left(\bigcap\limits_{i\in\left\{ 1,2,\dots,n\right\} \backslash I}A^{c}_{i}\right)\right]
\]

Since this is a disjoint union, using the additivity of $P,$ the
independence of the events $A_{i}$---using Proposition \ref{Proposition 3.2 independence}---and,
the fact that all events $A_{i}$ are of same probability $p$---and
their complement of probability $1-p$---it comes
\[
P\left(S_{n}=k\right)=\sum\limits_{\substack{I\in\mathcal{P}\left(\left\llbracket 1,n\right\rrbracket \right)\\
\left|I\right|=k
}
}p^{k}\left(1-p\right)^{n-k}.
\]

We have thus shown that, for every $k\in\left\llbracket 0,n\right\rrbracket ,$
\[
P\left(S_{n}=k\right)=\binom{n}{k}p^{k}\left(1-p\right)^{n-k},
\]
which corresponds to the fact that $S_{n}$ follows the binomial law
$\mathcal{B}\left(n,p\right).$

\end{proof}

\begin{figure}
\begin{center}%
\begin{tabular}{cc}
\begin{tikzpicture}[
    declare function={binom(\k,\n,\p)=\n!/(\k!*(\n-\k)!)*\p^\k*(1-\p)^(\n-\k);}
]
	\definecolor{orange}{HTML}{FFA200}
	\definecolor{purple}{HTML}{9500FF}
  \begin{axis}[
		axis x line=center,
		axis y line=center,
		xtick = {0,1,2,3,4,5,6,7,8,9,10},
		ytick = {0,0.1,0.2,0.3,0.4},
		ymin=-0.02,
		ymax=0.32
	]
\addplot[white, samples=80, smooth, thick, domain=-0.5:10.5 ]{0.31};
\pgfplotsinvokeforeach {0,1,2,3,4,5,6,7,8,9,10} {
    	\draw[green!80!black, line width=3] (axis cs: #1,0)
      -- (axis cs: #1,{binom(#1,10,0.2)});
  }
  \end{axis}

\end{tikzpicture} &
\begin{tikzpicture}[
    declare function={binom(\k,\n,\p)=\n!/(\k!*(\n-\k)!)*\p^\k*(1-\p)^(\n-\k);}
]
	\definecolor{orange}{HTML}{FFA200}
	\definecolor{purple}{HTML}{9500FF}
  \begin{axis}[
		axis x line=center,
		axis y line=center,
		xtick = {0,1,2,3,4,5,6,7,8,9,10},
		ytick = {0,0.1,0.2,0.3,0.4},
		ymin=-0.02,
		ymax=0.32
	]
\addplot[white, samples=80, smooth, thick, domain=-0.5:10.5 ]{0.32};
\pgfplotsinvokeforeach {0,1,2,3,4,5,6,7,8,9,10} {
    	\draw[green!80!black, line width=3] (axis cs: #1,0)
      -- (axis cs: #1,{binom(#1,10,0.6)});
  }
  \end{axis}

\end{tikzpicture}\tabularnewline
$\mathcal{B}\left(10,0.2\right)$ &
$\mathcal{B}\left(10,0.6\right)$\tabularnewline
\end{tabular}\end{center}

\caption{Examples of Binomial Laws}
\label{Fig:Binomial_law_examples}
\end{figure}

In Figure \ref{Fig:Binomial_law_examples}, we present two examples
of binomial laws $\mathcal{B}\left(n,p\right)$ for the same value
of $n$ and two different values of $p.$

\begin{proposition}{Convolution of Two Binomial Laws}{}

The \textbf{convolution of two binomial laws\index{convolution of two binomial laws}\mindex{binomial ! law ! convolution}\mindex{law ! binomial ! convolution}}
$\mathcal{B}\left(n_{1},p\right)$ and $\mathcal{B}\left(n_{2},p\right)$
is the binomial law $\mathcal{B}\left(n_{1}+n_{2},p\right).$

\end{proposition}

\begin{remark}{}{}

In other words, the family of binomial laws with the second parameter
$p$ is \textbf{\index{stable under convolution}\mindex{convolution ! stable}stable
under convolution}.

\end{remark}

\begin{proof}{}{}

We consider the random variable
\[
S_{n_{1}+n_{2}}=S_{n_{1}}+\sum\limits^{n_{1}+n_{2}}_{i=n_{1}+1}\boldsymbol{1}_{A_{i}}.
\]

The random variables $S_{n_{1}}$ and $\sum\limits^{n_{1}+n_{2}}_{i=n_{1}+1}\boldsymbol{1}_{A_{i}}$
are independent and follow the respective laws $\mathcal{B}\left(n_{1},p\right)$
and $\mathcal{B}\left(n_{2},p\right).$ 

By definition of convolution, the convolution of the laws $\mathcal{B}\left(n_{1},p\right)$
and $\mathcal{B}\left(n_{2},p\right)$ is the law of $S_{n_{1}+n_{2}}.$

Nonetheless, as
\[
S_{n_{1}+n_{2}}=\sum\limits^{n_{1}+n_{2}}_{i=1}\boldsymbol{1}_{A_{i}},
\]
it comes from the previous proposition that the law of $S_{n_{1}+n_{2}}$
is the binomial law $\mathcal{B}\left(n_{1}+n_{2},p\right),$ which
shows the result.

\end{proof}

\begin{example}{A Game}{}

A player repeatedly flips a biased coin---where $p$ denotes the
probability to obtain ``tail''---until the first tail occurs. If
this happens at the $k-$th flip, they roll a fair die $k$ times.
They win if they roll exactly one 6.

Find the probability that the player wins.

\end{example}

\begin{solutionexample}{}{} To model this game, we consider a probabilized
space $\left(\Omega,\mathcal{A},P_{p}\right)$ on which the following
random variables are defined:
\begin{itemize}
\item For each integer $k\in\mathbb{N}^{\ast},$ let $X_{k}$ be a random
variable following the binomial law $\mathcal{B}\left(k,\dfrac{1}{6}\right),$
representing the number of $6$ obtained during $k$ throwing of a
die. 
\item Let $R$ be a random variable following the geometric law of parameter
$p$ on $\mathbb{N}^{\ast}.$ This random variable represents the
rank of the first throw where tail appears. 
\item We assume that the random variables $R$ and $X_{k},\,k\in\mathbb{N}^{\ast}$
are independent---Proposition \ref{Generalized proposition for independent random variables}
ensures the existence of such a model. 
\end{itemize}
The event $G$ ``the player wins'' can be written:
\[
G=\biguplus\limits^{+\infty}_{k=1}\left(R=k,X_{k}=1\right).
\]

By the $\sigma-$additivity of $P_{p}$ and the independence of the
events $\left(R=k\right)$ and $\left(X_{k}=1\right),$ we obtain,
setting $q=1-p:$
\begin{align*}
P_{p}\left(G\right) & =\sum\limits^{+\infty}_{k=1}P_{p}\left(R=k\right)P_{p}\left(X_{k}=1\right)\\
 & =\sum\limits^{+\infty}_{k=1}pq^{k-1}\binom{k}{1}\dfrac{1}{6}\left(\dfrac{5}{6}\right)^{k-1}\\
 & =\dfrac{p}{6}\sum\limits^{+\infty}_{k=1}k\left(\dfrac{5q}{6}\right)^{k-1}.
\end{align*}

We observe that $0<\dfrac{5q}{6}<1.$ 

Moreover, classical results on the differentiation of power series
on its---real---domain of convergence allow to write that, for every
$x$ such that $0<x<1,$
\begin{align*}
\sum\limits^{+\infty}_{k=1}kx^{k-1} & =\dfrac{\text{d}}{\text{d}x}\left(\sum\limits^{+\infty}_{k=1}x^{k}\right)=\dfrac{\text{d}}{\text{d}x}\left(\dfrac{1}{1-x}\right)=\dfrac{1}{\left(1-x\right)^{2}}.
\end{align*}

It follows, after simplification, that\boxeq{
\[
P_{p}\left(G\right)=\dfrac{6p}{\left(1+5p\right)^{2}}.
\]
}

We can notice that
\[
\max\left\{ P_{p}\left(G\right):\,p\in\left]0,1\right[\right\} =0.3.
\]

This maximum is reached for $p=\dfrac{1}{5}:$ which is strictly greater
than $\dfrac{1}{6}.$

\end{solutionexample}

\section*{Exercises}

\addcontentsline{toc}{section}{Exercises}

Probability Theory can provide alternative methods for proving results
in various domain of mathematics. Here is an example in number theory.

\begin{exercise}{Euler Totient Function via Probability}{exercise3.1}

Let $\Omega=\left\llbracket 1,n\right\rrbracket ,$ where $n\geqslant2$
is a non prime integer. We consider the probabilized space $\left(\Omega,\mathcal{P}\left(\Omega\right),P\right),$
where $P$ is the uniform probability.

If $d$ divides $n,$ define the event
\[
A_{d}=\left\{ kd:\,k\in\Omega\,\,\,\,\text{and\,\,\,\,}kd\in\Omega\right\} .
\]

1. What is the probability of $A_{d}?$

2. Let $p_{1}<p_{2}<\dots<p_{r}$ be the sequence of prime divisors
of $n,$ ranked in increasing order.

(a) Prove that $\left(A_{p_{i}}\right)_{1\leqslant i\leqslant r}$
is a family of independent events.

(b) Deduce the cardinal $\varphi\left(n\right)$ of the set $A$ of
integers smaller or equal to $n$ that are coprime to $n.$ The function
$\varphi\left(n\right)$ is known as Euler totient function. 

\end{exercise}

\begin{exercise}{ Independence and Intuition}{exercise3.2}

Let $\left(\Omega,\mathcal{A},P\right)$ be a probabilized space on
which are defined the independent random variables $U$ and $V,$
each taking values in $\left\{ -1,1\right\} $ and of same law defined
by the relations
\[
\begin{array}{ccc}
P_{U}\left(-1\right)=\dfrac{1}{3} & \,\,\,\,\text{and}\,\,\,\, & P_{U}\left(1\right)=\dfrac{2}{3}.\end{array}
\]

Let $X$ and $Y$ be the random variables defined by
\[
\begin{array}{ccc}
X=U & \,\,\,\, & Y=\text{sign}\left(U\right)V.\end{array}
\]

1. What is the law of the random variable $\left(X,Y\right)?$ Are
the random variables $X$ and $Y$ independent?

2. Are the random variables $X^{2}$ and $Y^{2}$ independent? 

\end{exercise}

\begin{exercise}{Sum, Poisson Law and Lack of Independence}{exercise3.3} 

\textbf{This exercise illustrates how the sum of two random variables
following the Poisson law can also follow a Poisson law without having
those random variables independent.}

Let denote $\Omega$ the cartesian product $\mathbb{N}\times\mathbb{N}$
and let $\Lambda_{1}$ and $\Lambda_{2}$ be the subsets of $\Omega$
defined as
\[
\begin{array}{ccc}
\Lambda_{1}=\left\{ \left(0,1\right),\left(1,2\right),\left(2,0\right)\right\}  & \,\,\,\, & \Lambda_{2}=\left\{ \left(0,2\right),\left(2,1\right),\left(1,0\right)\right\} .\end{array}
\]

Let $q$ and $r$ be two positive real numbers. Let for every $i\in\mathbb{N},$
$q_{i}$ and $r_{i}$ be two real numbers such that
\[
\begin{array}{ccc}
q_{i}=\exp\left(-q\right)\dfrac{q^{i}}{i!} & \,\,\,\,\text{and}\,\,\,\, & r_{i}=\exp\left(-r\right)\dfrac{r^{i}}{i!}.\end{array}
\]

Let $\epsilon$ be a real number such that
\[
0<\epsilon<\min\left\{ q_{i}r_{j}:\,\left(i,j\right)\in\Lambda_{1}\cup\Lambda_{2}\right\} .
\]

Define the function $\mu:\Omega\to\mathbb{R}$ by
\[
\begin{array}{ccc}
\forall\left(i,j\right)\in\Omega, & \,\,\,\, & \mu\left(i,j\right)=\begin{cases}
q_{i}r_{j}+\epsilon, & \text{if\,}\left(i,j\right)\in\Lambda_{1},\\
q_{i}r_{j}-\epsilon, & \text{if\,}\left(i,j\right)\in\Lambda_{2},\\
q_{i}r_{j}, & \text{otherwise.}
\end{cases}\end{array}
\]

Let $P$ be the function defined on the probabilizable space $\left(\Omega,\mathcal{P}\left(\Omega\right)\right)$
by
\[
\begin{array}{ccc}
\forall A\in\mathcal{P}\left(\Omega\right), & \,\,\,\, & P\left(A\right)=\sum\limits_{\left(i,j\right)\in A}\mu\left(i,j\right).\end{array}
\]

1. Prove that $P$ is a probability on $\left(\Omega,\mathcal{P}\left(\Omega\right)\right).$

We denote $X_{1}$ and $X_{2}$ the canonical projections of $\Omega$
on $\mathbb{N},$ defined by setting for every $\left(i,j\right)\in\Omega$,
\[
\begin{array}{ccc}
X_{1}\left(i,j\right)=i & \,\,\,\, & X_{2}\left(i,j\right)=j\end{array}.
\]

2. Compute the law of the sum of random variables $X_{1}$ and $X_{2}$
and identify it with a known law.

3. What are the laws of the random variables $X_{1}$ and $X_{2}?$

4. Are the random variables $X_{1}$ and $X_{2}$ independent?

\end{exercise}

\begin{exercise}{Negative Binomial Law}{exercise3.4}

Let $\left(X_{n}\right)_{n\in\mathbb{N}^{\ast}}$ be a sequence of
independent random variables, defined on a probabilized space $\left(\Omega,\mathcal{A},P\right)$
of same geometric law on $\mathbb{N}^{\ast}$ with parameter $p.$

Define, for every $n\in\mathbb{N}^{\ast},$ the random variable $S_{n}$
by
\[
S_{n}=\sum\limits^{n}_{j=1}X_{j}.
\]

1. Compute the law of the random variable $S_{2}.$

2. Use induction to find the law of the random variable $S_{n}.$

\end{exercise}

\begin{exercise}{Around Geometric Law}{exercise3.5} 

Let $X$ and $Y$ be two independent random variables defined on the
probabilized space $\left(\Omega,\mathcal{A},P\right)$ of same geometric
law on $\mathbb{N}$ with parameter $p.$

1. The following probability computations are relevant for modelling,
for instance, the time until success for two players acting simultaneously
under the same conditions.

(a) Compute $P\left(Y\geqslant X\right).$ Analyze the special case
where $p=\dfrac{1}{2}.$

(b) Compute $P\left(Y=X\right).$ Analyze the case where $p=\dfrac{1}{2}.$

(c) Prove that $P\left(Y>X\right)=P\left(X>Y\right),$ and use this
to recover the probability $P\left(Y\geqslant X\right).$

2. Define two new random variables $U$ and $V$ by
\[
\begin{array}{ccc}
U=\max\left(X,Y\right) & \,\,\,\, & V=\min\left(X,Y\right).\end{array}
\]

(a) Compute, for every $\left(u,v\right)\in\mathbb{N}^{2},$ the probability
$P\left(U\leqslant u,V\geqslant v\right).$

(b) Deduce the laws of the random variables $U$ and $V.$

(c) Identify the law of the random variable $V.$

3. Define the random variable $W=U-V.$ 

Determine its law.

\end{exercise}

\begin{exercise}{Binomial Law, Poisson Law and Stock Management}{exercise3.6} 

The number of customers entering a store in one selling day is supposed
to follow the Poisson law\footnotemark with parameter $\lambda.$
Each customer independently purchases at most one article $A,$ with
probability $p.$ The initial stock of item $A$ at the opening of
the store is $s,$ with $s\geqslant1.$

We aim to compute the law of the total number of bought items purchased
during the day and, the probability that there is no stock outage.

We model this situation by setting a random variable $N$ representing
the number of customers entering the store in a day. We define a sequence
of random variables $\left(X_{n}\right)_{n\in\mathbb{N}^{\ast}},$
all defined on the same probabilized space $\left(\Omega,\mathcal{A},P\right),$
representing the buying decision of the $n-$th client, where:
\begin{itemize}
\item $X_{n}=1$ if the $n$-th customer buys an item, 
\item $X_{n}=0$ otherwise. 
\end{itemize}
We assume all the random variables to be independent. We also assume
that the law of $N$ is a Poisson law $\mathcal{P}\left(\lambda\right)$
and that, for every $n\in\mathbb{N}^{\ast},$ the random variable
$X_{n}$ follows a Bernoulli law with parameter $p.$ The total number
of requested items in one day is a random variable $T$ defined by
\[
T=\boldsymbol{1}_{\left(N\geqslant1\right)}\sum^{N}_{j=1}X_{j}.
\]

Determine the law of $T$ and compute the probability $P\left(T\leqslant s\right).$

\end{exercise}

\footnotetext{The choice of this hypothesis can be justified by Poisson
theorem---Theorem $\ref{th:Poisson-theorem}$.}

\section*{Solutions of Exercises}

\addcontentsline{toc}{section}{Solutions of Exercises}

\begin{solution}{}{solexercise3.1}

\textbf{1. Probability of $A_{d}$}

Let $k_{n}$ be the integer such that
\[
n=k_{n}d.
\]

We have
\[
\left|A_{d}\right|=k_{n}
\]

and thus, since $P$ is the uniform probability\boxeq{
\[
P\left(A_{d}\right)=\dfrac{k_{n}}{n}=\dfrac{1}{d}.
\]
}

\textbf{2.} We will first prove the independence of the events $A_{p_{i}},$
then use the fact that the complements of independent events are also
independent.

\textbf{(a) $\left(A_{p_{i}}\right)_{1\leqslant i\leqslant r}$ is
a family of independent events}

Since the integers $p_{i}$ and $p_{j}$ are primes, and are distinct
when $i\neq j,$ any common multiple of $p_{i}$ and $p_{j}$ is a
multiple of $p_{i}p_{j}.$ Therefore, for $i\neq j,$
\[
A_{p_{i}}\cap A_{p_{j}}=\left\{ kp_{i}p_{j}:\,k\in\Omega,\,kp_{i}p_{j}\in\Omega\right\} =A_{p_{i}p_{j}}.
\]

In particular, since the integer $n$ is a multiple of $p_{i}p_{j},$
\[
P\left(A_{p_{i}}\cap A_{p_{j}}\right)=P\left(A_{p_{i}p_{j}}\right)=\dfrac{1}{p_{i}p_{j}}=P\left(A_{p_{i}}\right)P\left(A_{p_{j}}\right).
\]

By similar reasoning, for every non-empty subset $J$ of $\left\llbracket 1,r\right\rrbracket ,$
we obtain\boxeq{
\[
P\left(\bigcap\limits_{i\in J}A_{p_{i}}\right)=\prod\limits_{i\in J}P\left(A_{p_{i}}\right),
\]
}which shows that the events $A_{p_{i}},$ $1\leqslant i\leqslant r,$
are independent.

\textbf{(b) Computation of $\varphi\left(n\right)$}

Since the events $A_{p_{i}},$ $1\leqslant i\leqslant r,$ are independent,
their complement $A^{c}_{p_{i}}$ are also independent. Moreover,
since $\bigcap^{r}_{i=1}A^{c}_{i}=A,$
\[
\dfrac{\varphi\left(n\right)}{n}=\prod^{r}_{i=1}\left(1-\dfrac{1}{p_{i}}\right),
\]
so that\boxeq{
\[
\varphi\left(n\right)=n\prod^{r}_{i=1}\left(1-\dfrac{1}{p_{i}}\right).
\]
}

\end{solution}

\begin{solution}{}{solexercise3.2}

\textbf{1. Law of the random variable $\left(X,Y\right).$ Independence
of $X$ and $Y$}

For every $x\in\left\{ -1,1\right\} $ and $y\in\left\{ -1,1\right\} ,$
\[
\left(X=x,Y=y\right)=\left(U=x,V=\text{sign}\left(x\right)y\right)
\]
and thus, by independence of the random variables $U$ and $V,$
\[
P\left(X=x,Y=y\right)=P\left(U=x\right)P\left(V=\text{sign}\left(x\right)y\right).
\]

Hence, we compute
\[
\begin{array}{ccc}
P\left(X=1,Y=1\right)=\dfrac{4}{9}, & \,\,\,\, & P\left(X=1,Y=-1\right)=\dfrac{2}{9},\\
P\left(X=-1,Y=1\right)=\dfrac{1}{9}, & \,\,\,\, & P\left(X=-1,Y=-1\right)=\dfrac{2}{9}.
\end{array}
\]

Also
\[
P\left(X=1\right)=\sum\limits_{y\in\left\{ -1,1\right\} }P\left(X=1,Y=y\right)=\dfrac{4}{9}+\dfrac{2}{9}=\dfrac{2}{3}.
\]

Similarly,
\[
P\left(Y=1\right)=\sum\limits_{x\in\left\{ -1,1\right\} }P\left(X=x,Y=1\right)=\dfrac{4}{9}+\dfrac{1}{9}=\dfrac{5}{9}.
\]

It follows
\[
P\left(X=1\right)P\left(Y=1\right)=\dfrac{2}{3}\times\dfrac{5}{9}=\dfrac{10}{27}.
\]

Thus,\boxeq{
\[
P\left(X=1\right)P\left(Y=1\right)\neq P\left(X=1,Y=1\right),
\]
}which shows that the events $X=1$ and $Y=1$ are not independent,
which implies that the random variables $X$ and $Y$ are also not
independent.

\textbf{2. Independence of $X^{2}$ and $Y^{2}$}

Nonetheless
\[
\begin{array}{ccccc}
X^{2}=U^{2} & \,\,\,\, & \text{and} & \,\,\,\, & Y^{2}=V^{2},\end{array}
\]
which implies that the random variables $X^{2}$ and $Y^{2}$ are
independent since the variables $U^{2}$ and $V^{2}$ are independent
as function of random variables that are independent.

\end{solution}

\begin{solution}{}{solexercise3}

\textbf{1. $P$ is a probability on $\left(\Omega,\mathcal{P}\left(\Omega\right)\right)$}

By the choice of $\epsilon,$ for every $\left(i,j\right)\in\Omega,$
we have $\mu\left(i,j\right)>0.$ For conciseness, let us denote:
$\mu_{ij}=\mu\left(i,j\right).$ By the additivity of nonnegative
families,
\[
\sum\limits_{\left(i,j\right)\in\Omega}\mu_{ij}=\sum\limits_{\left(i,j\right)\in\Lambda_{1}}\mu_{ij}+\sum\limits_{\left(i,j\right)\in\Lambda_{2}}\mu_{ij}+\sum\limits_{\left(i,j\right)\in\left(\Lambda_{1}\cup\Lambda_{2}\right)^{c}}\mu_{ij}.
\]

That is
\[
\sum\limits_{\left(i,j\right)\in\Omega}\mu_{ij}=\left(\sum\limits_{\left(i,j\right)\in\Lambda_{1}}q_{i}r_{j}\right)+3\epsilon+\left(\sum\limits_{\left(i,j\right)\in\Lambda_{2}}q_{i}r_{j}\right)-3\epsilon+\sum\limits_{\left(i,j\right)\in\left(\Lambda_{1}\cup\Lambda_{2}\right)^{c}}q_{i}r_{j}.
\]

Thus, it follows
\[
\sum\limits_{\left(i,j\right)\in\Omega}\mu_{ij}=\sum\limits_{\left(i,j\right)\in\mathbb{N}\times\mathbb{N}}q_{i}r_{j}.
\]

By the Fubini theorem,
\[
\sum\limits_{\left(i,j\right)\in\Omega}\mu_{ij}=\left(\sum\limits_{i\in\mathbb{N}}q_{i}\right)\left(\sum\limits_{j\in\mathbb{N}}r_{j}\right).
\]

Since
\[
\begin{array}{ccc}
\sum\limits_{i\in\mathbb{N}}q_{i}=1 & \,\,\,\,\text{and\,\,\,\,} & \sum\limits_{j\in\mathbb{N}}r_{j}=1,\end{array}
\]
it follows that \boxeq{
\[
\sum\limits_{\left(i,j\right)\in\Omega}\mu_{ij}=1.
\]
}

The function $\mu$ is then indeed a probability germ: it generates
the probability $P.$

\textbf{2. Law of the sum of random variables $X_{1}$ and $X_{2}$
and identification.}

Let $k\in\mathbb{N}.$ Then
\[
\left(X_{1}+X_{2}=k\right)=\biguplus^{k}_{i=0}\left(X_{1}=i,X_{2}=k-i\right).
\]

Since the random variable $\left(X_{1},X_{2}\right)$ is the identity
map on $\Omega,$ its law is given by the probability $P$ itself.
Consequently, we have
\begin{align*}
P\left(X_{1}+X_{2}=k\right) & =\sum^{k}_{i=0}\mu_{i,k-i}\\
 & =\sum_{i\in\Lambda_{1,k}}\mu_{i,k-i}+\sum_{i\in\Lambda_{2,k}}\mu_{i,k-i}+\sum_{i\in\Lambda_{3,k}}\mu_{i,k-i},
\end{align*}
where the sets $\Lambda_{1,k},$ $\Lambda_{2,k}$ and $\Lambda_{3,k}$
are defined by
\begin{align*}
\Lambda_{1,k} & =\left\{ i\in\mathbb{N}:\,0\leqslant i\leqslant k\,\,\,\,\text{and}\,\,\,\,\left(i,k-i\right)\in\Lambda_{1}\right\} \\
\Lambda_{2,k} & =\left\{ i\in\mathbb{N}:\,0\leqslant i\leqslant k\,\,\,\,\text{and}\,\,\,\,\left(i,k-i\right)\in\Lambda_{2}\right\} \\
\Lambda_{3,k} & =\left\{ i\in\mathbb{N}:\,0\leqslant i\leqslant k\right\} \backslash\left(\Lambda_{1,k}\cup\Lambda_{2,k}\right).
\end{align*}

Explicitly,
\[
\Lambda_{1,k}=\begin{cases}
\left\{ 0\right\} , & \text{if\,}k=1,\\
\left\{ 2\right\} , & \text{if\,}k=2,\\
\left\{ 1\right\} , & \text{if}\,k=3,\\
\emptyset, & \text{otherwise,}
\end{cases}
\]
and
\[
\Lambda_{2,k}=\begin{cases}
\left\{ 1\right\} , & \text{if\,}k=1,\\
\left\{ 0\right\} , & \text{if\,}k=2,\\
\left\{ 2\right\} , & \text{if}\,k=3,\\
\emptyset, & \text{otherwise,}
\end{cases}
\]

So
\[
\left|\Lambda_{1,k}\right|=\left|\Lambda_{2,k}\right|=\begin{cases}
1, & \text{if}\,1\leqslant k\leqslant3,\\
0, & \text{otherwise.}
\end{cases}
\]

Thus,
\[
P\left(X_{1}+X_{2}=k\right)=\left(\sum\limits_{i\in\Lambda_{1,k}}q_{i}r_{k-i}\right)+\epsilon\left|\Lambda_{1,k}\right|+\left(\sum\limits_{i\in\Lambda_{2,k}}q_{i}r_{k-i}\right)-\epsilon\left|\Lambda_{2,k}\right|+\sum\limits_{i\in\Lambda_{3,k}}q_{i}r_{k-i}
\]
which gives, since $q_{i}=\exp\left(-q\right)\dfrac{q^{i}}{i!}$ and
$r_{k-i}=\exp\left(-r\right)\dfrac{r^{k-i}}{\left(k-i\right)!},$
\begin{align*}
P\left(X_{1}+X_{2}=k\right) & =\sum\limits^{k}_{i=0}q_{i}r_{k-i}\\
 & =\exp\left(-\left(q+r\right)\right)\sum\limits^{k}_{i=0}\dfrac{q^{i}r^{k-i}}{i!\left(k-i\right)!}\\
 & =\exp\left(-\left(q+r\right)\right)\dfrac{1}{k!}\sum\limits^{k}_{i=0}\left(\begin{array}{c}
k\\
i
\end{array}\right)q^{i}r^{k-i}
\end{align*}
Hence,\boxeq{
\[
P\left(X_{1}+X_{2}=k\right)=\exp\left(-\left(q+r\right)\right)\dfrac{\left(q+r\right)^{k}}{k!}.
\]
}

The law of the random variable $X_{1}+X_{2}$ is thus a Poisson law
of parameter $q+r.$

\textbf{3. Laws of the random variables $X_{1}$ and $X_{2}$}

For every $k\in\mathbb{N},$
\[
\left(X_{1}=k\right)=\biguplus\limits_{j\in\mathbb{N}}\left\{ \left(k,j\right)\right\} 
\]
and thus,
\begin{equation}
P\left(X_{1}=k\right)=\sum\limits_{j\in\mathbb{N}}\mu_{k,j}.\label{eq:P(X_1=00003Dk)}
\end{equation}

We compute this formula with different values of $k:$

(a) If $k=0,$ then
\begin{align*}
P\left(X_{1}=0\right) & =\mu_{0,0}+\mu_{0,1}+\mu_{0,2}+\sum\limits^{+\infty}_{j=3}\mu_{0,j}\\
 & =q_{0}r_{0}+\left(q_{0}r_{1}+\epsilon\right)+\left(q_{0}r_{2}-\epsilon\right)+\sum\limits^{+\infty}_{j=3}q_{0}r_{j}\\
 & =q_{0}\left(\sum\limits^{+\infty}_{j=0}r_{j}\right)\\
 & =q_{0}.
\end{align*}

(b) If $k=1,$ then
\begin{align*}
P\left(X_{1}=1\right) & =\mu_{1,0}+\mu_{1,2}+\mu_{1,1}+\sum\limits^{+\infty}_{j=3}\mu_{1,j}\\
 & =\left(q_{1}r_{0}-\epsilon\right)+\left(q_{1}r_{2}+\epsilon\right)+q_{1}r_{1}+\sum\limits^{+\infty}_{j=3}q_{1}r_{j}\\
 & =q_{1}\left(\sum\limits^{+\infty}_{j=0}r_{j}\right)\\
 & =q_{1}.
\end{align*}

(c) If $k=2,$ then
\begin{align*}
P\left(X_{1}=2\right) & =\mu_{2,0}+\mu_{2,1}+\sum\limits^{+\infty}_{j=2}\mu_{2,j}\\
 & =\left(q_{2}r_{0}-\epsilon\right)+\left(q_{2}r_{1}+\epsilon\right)+\sum\limits^{+\infty}_{j=2}q_{2}r_{j}\\
 & =q_{2}\left(\sum\limits^{+\infty}_{j=0}r_{j}\right)\\
 & =q_{2}.
\end{align*}

(d) If $k\geqslant3,$ then
\[
P\left(X_{1}=k\right)=\sum\limits^{+\infty}_{j=0}q_{k}r_{j}=q_{k}.
\]

In summary, we have, for every $k\in\mathbb{N},$\boxeq{
\[
P\left(X_{1}=k\right)=q_{k}.
\]
}

Similarly, we could show that, for every $k\in\mathbb{N},$\boxeq{
\[
P\left(X_{2}=k\right)=r_{k}.
\]
}

We show that the random variables $X_{1}$ and $X_{2}$ are Poisson
laws of respective parameters $q$ and $r.$

\textbf{4. Independence of the random variables $X_{1}$ and $X_{2}$}

We saw, that for every $\left(i,j\right)\in\mathbb{N}\times\mathbb{N},$
\[
P\left(X_{1}=i,X_{2}=j\right)=\mu_{ij}.
\]

In particular,
\[
P\left(X_{1}=1,X_{2}=0\right)=\mu_{1,0}
\]
and thus\boxeq{
\[
P\left(X_{1}=1,X_{2}=0\right)\neq P\left(X_{1}=1\right)P\left(X_{2}=0\right)
\]
}which shows that the random variables $X_{1}$ and $X_{2}$ are
not independent.

\end{solution}

\begin{solution}{}{solexercise3.4}

\textbf{1. Computation of the law of the random variable $S_{2}.$}

The random variables $X_{n}$ take values in $\mathbb{N}^{\ast}$
and we have $S_{2}\left(\Omega\right)\subset\mathbb{N}^{\ast}.$ For
$k\in\mathbb{N}^{*},$ using the \textbf{convolution method}, we get
\[
P\left(S_{2}=k\right)=\sum^{k}_{j=1}P\left(X_{1}=j\right)P\left(X_{2}=k-j\right).
\]

In particular, since each $X_{n}\geqslant1,$ $P\left(S_{2}=1\right)=0.$

Setting $q=1-p$, for $k\geqslant2,$ it holds
\begin{align*}
P\left(S_{2}=k\right) & =\sum^{k-1}_{j=1}pq^{j-1}pq^{k-j-1}=\left(k-1\right)p^{2}q^{k-2}.
\end{align*}
Hence, for $k\geqslant2,$\boxeq{
\[
P\left(S_{2}=k\right)=\binom{k-1}{1}p^{2}q^{k-2}.
\]
}

\textbf{2. Law of the random variable $S_{n}.$}

Suppose that the law of the random variable $S_{n}$ is given by
\begin{equation}
P\left(S_{n}=k\right)=\begin{cases}
0, & \text{if }1\leqslant k\leqslant n-1,\\
\binom{k-1}{n-1}p^{n}q^{k-n}, & \text{if }k\geqslant n.
\end{cases}\label{eq:induction S_n}
\end{equation}

Since $S_{n+1}=S_{n}+X_{n+1}$ and that the random variables $S_{n}$
and $X_{n+1}$ are independent, the convolution method gives
\[
P\left(S_{n+1}=k\right)=\sum^{k}_{j=1}P\left(S_{n}=j\right)P\left(X_{n+1}=k-j\right).
\]
From the induction hypothesis, we deduce:
\begin{itemize}
\item For $1\leqslant k\leqslant n,$
\[
P\left(S_{n+1}=k\right)=0.
\]
\item For $k\geqslant n+1,$
\begin{align*}
P\left(S_{n+1}=k\right) & =\sum^{k-1}_{j=n}\binom{j-1}{n-1}p^{n}q^{j-n}pq^{k-j-1}\\
 & =\left(\sum^{k-1}_{j=n}\binom{j-1}{n-1}\right)p^{n+1}q^{k-\left(n+1\right)}.
\end{align*}
\end{itemize}
It remains to calculate the coefficient of $p^{n+1}q^{k-\left(n+1\right)};$
to achieve it, we consider the polynom $Q\left(x\right)$ defined
by
\[
Q\left(x\right)=\sum^{k-1}_{n=1}\left(\sum^{k-1}_{j=n}\binom{j-1}{n-1}\right)x^{n}.
\]

By sum permutation and by translation of indices, we successively
obtain
\begin{align*}
Q\left(x\right) & =\sum^{k-1}_{j=1}\left(\sum^{j}_{n=1}\binom{j-1}{n-1}x^{n}\right)\\
 & =\sum^{k-1}_{j=1}\left(\sum^{j-1}_{p=0}\binom{j-1}{p}x^{p+1}\right)\\
 & =\sum^{k-1}_{j=1}x\left(1+x\right)^{j-1}\\
 & =x\dfrac{1-\left(1+x\right)^{k-1}}{1-\left(1+x\right)}\\
 & =\left(1+x\right)^{k-1}-1\\
 & =\sum^{k-1}_{j=1}\binom{k-1}{j}x^{j}.
\end{align*}

By comparing powers of $x$ one by one in this last with the definition
of $Q\left(x\right),$
\[
\sum^{k-1}_{j=n}\binom{j-1}{n-1}=\binom{k-1}{n}.
\]

It proves that the law of $S_{n}$ is given by Relation \ref{eq:induction S_n}.
This law is also called the \textbf{\index{negative binomial law}\mindex{law ! negative binomial}negative
binomial law} with parameters $n$ and $p,$ also called the Pascal-Montmort
law.

\end{solution}

\begin{solution}{}{solexercise3.5}

\textbf{1. (a) Computation of $P\left(Y\geqslant X\right).$ Special
case where $p=\dfrac{1}{2}.$}

We have
\[
\left(Y\geqslant X\right)=\biguplus\limits_{\substack{\left(i,j\right)\in\mathbb{N}^{2}\\
i\geqslant j
}
}\left(Y=i,X=j\right)
\]
and, by $\sigma-$additivity of $P,$
\[
P\left(Y\geqslant X\right)=\sum_{\substack{\left(i,j\right)\in\mathbb{N}^{2}\\
i\geqslant j
}
}P\left(Y=i,X=j\right).
\]

The random variables $X$ and $Y$ being independent and of same geometric
law on $\mathbb{N}$ with parameter $p,$ we can successively write,
by setting as usual $q=1-p,$
\begin{align*}
P\left(Y\geqslant X\right) & =\sum\limits_{\substack{\left(i,j\right)\in\mathbb{N}^{2}\\
i\geqslant j
}
}pq^{i}pq^{j}\\
 & =p^{2}\sum^{+\infty}_{i=0}\left[q^{i}\left(\sum^{i}_{j=0}q^{j}\right)\right]\\
 & =p^{2}\sum^{+\infty}_{i=0}\left[q^{i}\dfrac{1-q^{i+1}}{1-q}\right]\\
 & =\sum^{+\infty}_{i=0}pq^{i}-pq\sum^{+\infty}_{i=0}q^{2i}\\
 & =1-\dfrac{pq}{1-q^{2}}\\
 & =\dfrac{1}{1+q}.
\end{align*}
In particular, when $p=\dfrac{1}{2},$ we obtain\boxeq{
\[
P\left(Y\geqslant X\right)=\dfrac{2}{3}.
\]
}

\textbf{(b) Computation of $P\left(Y=X\right).$ Special case where
$p=\dfrac{1}{2}.$}

We have
\[
\left(Y=X\right)=\biguplus_{i\in\mathbb{N}}\left(Y=i,X=i\right)
\]
and, by $\sigma-$additivity of $P,$
\[
P\left(Y=X\right)=\sum_{i\in\mathbb{N}}P\left(Y=i,X=i\right).
\]

The random variables $X$ and $Y$ being independent, of same geometric
law on $\mathbb{N}$ of parameter $p,$ we can write
\begin{align*}
P\left(Y=X\right) & =\sum_{i\in\mathbb{N}}p^{2}q^{2i}=\dfrac{p^{2}}{1-q^{2}}=\dfrac{p}{1+q}.
\end{align*}

In particular, when $p=\dfrac{1}{2},$ we obtain \boxeq{
\[
P\left(Y=X\right)=\dfrac{1}{3}.
\]
}

\textbf{(c) $P\left(Y>X\right)=P\left(X>Y\right),$ recovery of $P\left(Y\geqslant X\right).$}

Since
\[
P\left(Y>X\right)=\sum_{\substack{\left(i,j\right)\in\mathbb{N}^{2}\\
i>j
}
}P\left(Y=i\right)P\left(X=j\right)
\]
 and, since the random variables $X$ and $Y$ are of same law, it
holds
\begin{align*}
P\left(Y>X\right) & =\sum_{\substack{\left(i,j\right)\in\mathbb{N}^{2}\\
i>j
}
}P\left(X=i\right)P\left(Y=j\right)\\
 & =P\left(X>Y\right).
\end{align*}

Hence, since
\[
P\left(X<Y\right)+P\left(Y=X\right)+P\left(Y>X\right)=1
\]
we have
\[
P\left(Y>X\right)=\dfrac{1}{2}\left[1-P\left(Y=X\right)\right].
\]

It follows
\[
P\left(Y>X\right)=\dfrac{q}{1+q}.
\]

In particular, when $p=\dfrac{1}{2},$
\[
P\left(Y>X\right)=P\left(Y<X\right)=P\left(Y=X\right)=\dfrac{1}{3}.
\]

We find back the probability $P\left(Y\geqslant X\right)$ by the
following manner
\begin{align*}
P\left(Y\geqslant X\right) & =P\left(Y>X\right)+P\left(Y=X\right)\\
 & =\dfrac{q}{1+q}+\dfrac{p}{1+q}\\
 & =\dfrac{1}{1+q}.
\end{align*}

\textbf{2. (a) Computation, for every $\left(u,v\right)\in\mathbb{N}^{2},$
of the probability $P\left(U\leqslant u,V\geqslant v\right).$}

If $0\leqslant u<v,$ then
\[
\left(U\leqslant u\right)\cap\left(V\geqslant v\right)=\emptyset,
\]
and thus
\[
P\left(U\leqslant u,V\geqslant v\right)=0.
\]

If $0\leqslant v\leqslant u,$ then
\[
\left(U\leqslant u\right)\cap\left(V\geqslant v\right)=\left(v\leqslant X\leqslant u,v\leqslant Y\leqslant u\right).
\]

The random variables $X$ and $Y$ are independent, of same geometric
law on $\mathbb{N}$ with parameter $p,$ we can write
\begin{align*}
P\left(U\leqslant u,V\geqslant v\right) & =\left[P\left(v\leqslant X\leqslant u\right)\right]^{2}=\left(\sum^{u}_{j=v}pq^{j}\right)^{2},
\end{align*}
hence, after simplifications,\boxeq{
\begin{equation}
P\left(U\leqslant u,V\geqslant v\right)=\begin{cases}
\left(q^{v}-q^{u+1}\right)^{2}, & \text{if }0\leqslant v\leqslant u,\\
0, & \text{otherwise.}
\end{cases}\label{eq:ex5_P(max_min)}
\end{equation}
}

\textbf{(b) Laws of the random variables $U$ and $V$}

We have $\left(V\geqslant0\right)=\Omega.$

Hence,
\[
P\left(U\leqslant u\right)=\begin{cases}
\left(1-q^{u+1}\right)^{2}, & \text{if\,}u\geqslant0,\\
0, & \text{otherwise.}
\end{cases}
\]
But, since 
\[
P\left(U=u\right)=P\left(U\leqslant u\right)-P\left(U\leqslant u-1\right),
\]
it follows, after computation, that\boxeq{
\[
P\left(U=u\right)=\begin{cases}
p\left[2q^{u}-q^{2u}\left(1+q\right)\right], & \text{if\,}u\geqslant0,\\
0, & \text{otherwise.}
\end{cases}
\]
}

\textbf{(c) Law of the random variable $V.$}

The sequence of sets $\left(U\leqslant u\right)_{u\in\mathbb{N}}$
is nondecreasing---for the inclusion---and their union is $\Omega.$
We can thus write that
\[
P\left(V\geqslant v\right)=\lim_{u\rightarrow+\infty}P\left(U\leqslant u,V\geqslant v\right).
\]

Thus,
\[
P\left(V\geqslant v\right)=\begin{cases}
q^{2v}, & \text{if }v\geqslant0,\\
0, & \text{otherwise.}
\end{cases}
\]

But, since
\[
P\left(V=v\right)=P\left(V\geqslant v\right)-P\left(V\geqslant v+1\right),
\]
it follows, after computation,\boxeq{
\[
P\left(V=v\right)=\begin{cases}
q^{2v}\left[1-q^{2}\right], & \text{if }v\geqslant0,\\
0, & \text{otherwise.}
\end{cases}
\]
}which shows that the law of $V$ is a geometric law on $\mathbb{N}$
with parameter $1-q^{2}.$

\textbf{3. Law of $W=U-V.$ }

We can write
\[
P\left(U\leqslant u,V=v\right)=P\left(U\leqslant u,V\geqslant v\right)-P\left(U\leqslant u,V\geqslant v+1\right).
\]

Using $\refpar{eq:ex5_P(max_min)}$ and after simplifications, it
holds:
\begin{itemize}
\item If $u\geqslant v+1$ and $v\geqslant0,$ then
\begin{equation}
P\left(U\leqslant u,V=v\right)=pq^{v}\left[q^{v}\left(1+q\right)-2q^{u+1}\right].\label{eq:ex5_3_first_case}
\end{equation}
\item If $u=v\geqslant0,$ then
\begin{align*}
P\left(U\leqslant u,V=v\right) & =P\left(V=U=u\right)=\left[q^{u}-q^{u+1}\right]^{2}=p^{2}q^{2u}.
\end{align*}
\end{itemize}
We can remark, in this last case, that the relation $\refpar{eq:ex5_3_first_case}$
is still true.

Then, for every $w\in\mathbb{N},$ we have
\[
\left(W\leqslant w\right)=\biguplus_{v\in\mathbb{N}}\left(U\leqslant w+v,V=v\right),
\]
it comes by using the relation $\refpar{eq:ex5_3_first_case}$ and
after computation,
\[
\begin{array}{ccc}
\forall w\in\mathbb{N}, & \,\,\,\, & P\left(W\leqslant w\right)=1-\dfrac{2}{1+q}q^{w+1}.\end{array}
\]

We then deduce the probability $P\left(W\leqslant w\right):$
\begin{itemize}
\item If $w\geqslant1,$ then
\begin{align*}
P\left(W=w\right) & =P\left(W\leqslant w\right)-P\left(W\leqslant w-1\right)=\dfrac{2p}{1+q}q^{w}.
\end{align*}
\item If $w=0,$ then
\begin{align*}
P\left(W=0\right) & =P\left(X=Y\right)=\dfrac{p}{1+q}.
\end{align*}
\end{itemize}
\end{solution}

\begin{solution}{}{solexercise3.6}

\textbf{Law of $T$}

We partition
\[
\left(T=0\right)=\left(N=0\right)\biguplus\left[\biguplus_{n\geqslant1}\left(N=n,\sum^{n}_{j=1}X_{j}=0\right)\right].
\]

Denoting $q=1-p,$ by $\sigma-$additivity of $P$ and independence
of the random variables $N$ and $X_{n},$ it holds
\[
P\left(T=0\right)=P\left(N=0\right)+\sum_{n\geqslant1}\exp\left(-\lambda\right)\dfrac{\lambda^{n}}{n!}q^{n},
\]
which gives, after computation
\begin{align*}
P\left(T=0\right) & =\exp\left(-\lambda\right)+\sum_{n\geqslant1}\exp\left(-\lambda\right)\dfrac{\lambda^{n}}{n!}q^{n}=\exp\left(-\lambda p\right).
\end{align*}

Moreover, for every $k\geqslant1,$ we can write
\[
\left(T=k\right)=\biguplus_{n\geqslant1}\left(N=n,\sum^{n}_{j=1}X_{j}=k\right).
\]

Since the law of $\sum^{n}_{j=1}X_{j}$ is the binomial law $\mathcal{B}\left(n,p\right),$
it follows
\begin{align*}
P\left(T=k\right) & =\sum_{n\geqslant k}\exp\left(-\lambda\right)\dfrac{\lambda^{n}}{n!}\binom{n}{k}p^{k}q^{n-k}\\
 & =\exp\left(-\lambda\right)\dfrac{\left(p\lambda\right)^{k}}{k!}\sum_{n\geqslant k}\dfrac{\left(\lambda q\right)^{n-k}}{\left(n-k\right)!}\\
 & =\exp\left(-\lambda p\right)\dfrac{\left(p\lambda\right)^{k}}{k!}.
\end{align*}

We just proved that the law of $T$ is a Poisson law with parameter
$\lambda p.$

\textbf{Computation of the probability $P\left(T\leqslant s\right)$}

Hence, the probability that there is no stock outage is\boxeq{
\[
P\left(T\leqslant s\right)=\exp\left(-\lambda p\right)\sum^{s}_{k=0}\dfrac{\left(p\lambda\right)^{k}}{k!}.
\]
}

\end{solution}

\chapter{Probabilities and Conditional Laws}\label{chap:Probabilities-and-Conditional}

\begin{objective}{}{}

Chapter \ref{chap:Probabilities-and-Conditional} introduces the concepts
of conditional probabilities and conditional laws.
\begin{itemize}
\item Section \ref{sec:Conditional-Probabilities} focuses on conditional
probability. It starts by defining the concept, then introduces the
chain rule for conditional probabilities and the law of total probability.
The section ends with a discussion on probability of causes and states
Bayes theorem.
\item Section \ref{sec:Conditional-Laws} explores conditional laws.
\item Section \ref{sec:Evolutive-phenomena-modeling} addresses how conditional
probabilities can be used to model an evolutive phenomena. This leads
to the definition of a Markov sequence of random variables.
\end{itemize}
\end{objective}

In the previous Chapter, we studied the notion of independent events.
In this Chapter, we are going to see how to translate mathematically
the influence of one event on another. The notion of conditional probability
allows, in some extent, to take into account the information carried
by the realization of an event on the possible realization of another. 

\section{Conditional Probabilities}\label{sec:Conditional-Probabilities}

In what follows, $P$ is any probability on a probabilizable space
$\left(\Omega,\mathcal{A}\right).$

\subsection{Definitions. Total Probability Formula.}

\begin{definition}{Conditional Probability}{}

Let $B\in\mathcal{A}$ be an event such that $P\left(B\right)>0.$ 

Let $A\in\mathcal{A}$ be another event. 

The \textbf{probability of $A$ conditioned by $B$}\mindex{probability!of $A$ conditioned by $B$}---or
\textbf{the conditional probability of $A$ given $B$}\index{conditional probability of $A$ given $B$}---is
the real number, denoted $P\left(A\mid B\right)$ or $P^{B}\left(A\right),$
and defined by\boxeq{
\begin{equation}
P\left(A\mid B\right)=P^{B}\left(A\right)=\dfrac{P\left(A\cap B\right)}{P\left(B\right)}.\label{eq:Conditional probability def}
\end{equation}
}

\end{definition}

\begin{proposition}{A Conditional Probability is a Probability}{} 

Let $B\in\mathcal{A}$ be an event such that $P\left(B\right)>0.$ 

The function $P^{B}$ on $\mathcal{A}$ to $\mathbb{R}^{+}$ which
maps every $A\in\mathcal{A}$ to $P^{B}\left(A\right)$ is a probability
on the probabilizable space $\left(\Omega,\mathcal{A}\right).$ It
is refered to as the \textbf{\mindex{probability!conditioned to an event}probability
conditioned\footnotemark with respect to $B$}---or more commonly,
the \textbf{conditional probability with respect to $B$.}

\end{proposition}

\footnotetext{It is often mentioned as the \textbf{\index{conditional probability}\mindex{probability ! conditional}conditional
probability $P^{B}.$}}

\begin{proof}{}{} 

We check that $P^{B}:\mathcal{A}\to\mathbb{R}^{+}$ behaves like a
probability.
\begin{itemize}
\item \textbf{Non-negativity} is immediate, as
\[
P^{B}\left(A\right)\geqslant0.
\]
\item \textbf{Normalization} happens, since
\[
P^{B}\left(\Omega\right)=\dfrac{P\left(\Omega\cap B\right)}{P\left(B\right)}=1.
\]
\item The mapping $P^{B}$ is \textbf{$\sigma-$additive}.
\end{itemize}
Indeed, for every sequence $\left(A_{i}\right)_{i\in\mathbb{N}}$
of events, disjoints two by two,
\[
\left(\biguplus_{i\in\mathbb{N}}A_{i}\right)\cap B=\biguplus_{i\in\mathbb{N}}\left(A_{i}\cap B\right),
\]
 and thus, using the definition of $P^{B}$ and the $\sigma-$additivity
of $P,$
\begin{align*}
P^{B}\left(\biguplus_{i\in\mathbb{N}}A_{i}\right) & =\dfrac{P\left(\left(\biguplus_{i\in\mathbb{N}}A_{i}\right)\cap B\right)}{P\left(B\right)}=\dfrac{P\left(\biguplus_{i\in\mathbb{N}}\left(A_{i}\cap B\right)\right)}{P\left(B\right)}=\dfrac{\sum_{i\in\mathbb{N}}P\left(A_{i}\cap B\right)}{P\left(B\right)}.
\end{align*}
Hence,
\[
P^{B}\left(\biguplus_{i\in\mathbb{N}}A_{i}\right)=\sum_{i\in\mathbb{N}}P^{B}\left(A_{i}\right).
\]

\end{proof}

\begin{remarks}{}{}

1. The importance of this proposition resides in the fact that all
the already proven properties---and thus that will be further proven---for
any probability are true for the \textbf{conditional probability}
$P^{B}.$ 

2. Let $B\in\mathcal{A}$ be an event such that $P\left(B\right)>0$
and any other event $A.$ For the events $A$ and $B$ to be independent,
it must and it suffices that
\[
P^{B}\left(A\right)=P\left(A\right).
\]

That is the knowledge of $B$ does not affect the likelihood of $A.$
Independence can now be seen as a special case of conditional probability
where ``nothing changes'' under the additional knowledge of $B.$

\end{remarks}

The following example illustrates how a conditional probability gives
an idea of the information given by the occurrence of one event about
the occurrence of another.

\begin{example}{Two Rolls of a Die}{} A die is rolled two times.
Let $A$ be the event: ``a 6 is obtained on the first roll'', and
let $B_{k},\,2\leqslant k\leqslant12,$ be the event: ``the sum of
the two integers obtained is $k$''. 

The two die rolls are modelled by the probabilized space $\left(\Omega,\mathcal{P}\left(\Omega\right),P\right)$
where $\Omega=\left\{ 1,2,\dots,6\right\} ^{2}$ and $P$ is the uniform
probability. We then have:
\begin{itemize}
\item $A=\left\{ 6\right\} \times\left\{ 1,2,\dots,6\right\} $
\item $B_{k}=\left\{ \left(\omega_{1},\omega_{2}\right)\in\Omega:\,\omega_{1}+\omega_{2}=k\right\} .$
\end{itemize}
Compute:
\begin{itemize}
\item The probability of $A,$
\item The probability of $A$ knowing $B_{12},$
\item The probability of $A$ knowing $B_{11}.$
\end{itemize}
\end{example}

\begin{solutionexample}{}{}

We observe that:
\[
\begin{array}{ccc}
B_{12}\subset A & \,\,\,\,\text{and\,\,\,\,} & B_{11}=\left\{ \left(5,6\right),\left(6,5\right)\right\} \end{array}.
\]

Therefore:
\[
\begin{array}{ccccc}
P\left(A\right)=\dfrac{1}{6}, & \,\,\,\, & P\left(A\mid B_{12}\right)=1 & \text{\,and,\,} & P\left(A\mid B_{11}\right)=\dfrac{1}{2}.\end{array}
\]

\end{solutionexample}

The following proposition is deduced directly from the definition.

\begin{proposition}{Probability of Intersection and Conditional Probability}{}\label{Prop P inter in function cond}Let
$B\in\mathcal{A}$ be an event such that $P\left(B\right)>0.$ Let
$A\in\mathcal{A}$ be another event.\boxeq{
\begin{equation}
P\left(A\cap B\right)=P\left(A\mid B\right)P\left(B\right).\label{eq:p_intersection}
\end{equation}
}

\end{proposition}

It is often under the form of Proposition \ref{Prop P inter in function cond}
that the definition of $P\left(A\mid B\right)$ is used. Indeed, in
several problems involving probabilistic reasoning, the information
provided---whether as part of experimental data or through a problem
analysis---naturally lends itself to interpretation in a straightforward
manner in term of conditional probabilities.

Proposition \ref{Prop P inter in function cond} admits the following
generalization, particularly useful in case of sequences of dependent
events.

\begin{proposition}{Chain Rule of Conditional Probabilities}{} 

Let $\left(A_{i}\right)_{1\leqslant i\leqslant n}$ be a finite sequence
of events such that $P\left(\bigcap_{1\leqslant i\leqslant n-1}A_{i}\right)>0.$

Then, the probability of their joint occurrence is given by the product
of successive conditional probabilities:\boxeq{
\[
P\left(\bigcap^{n}_{i=1}A_{i}\right)=P\left(A_{1}\right)P\left(A_{2}\mid A_{1}\right)P\left(A_{3}\mid A_{1}\cap A_{2}\right)\dots P\left(A_{n}\mid A_{1}\cap\dots\cap A_{n-1}\right).
\]
}

\end{proposition}

\begin{proof}{}{}We begin by noting that all the conditional probabilities
introduced are well-defined. Indeed, for every $j\in\left\{ 1,2,\dots,n-1\right\} ,$
\[
P\left(\bigcap^{j}_{i=1}A_{i}\right)\geqslant P\left(\bigcap^{n-1}_{i=1}A_{i}\right)>0.
\]

The desired formula then follows directly by induction.

\end{proof}

\begin{example}{}{}

We successively draw 4 cards from a standard deck of 52 cards. 

What is the probability that all 4 cards drawn are aces?

\end{example}

\begin{solutionexample}{}{}Let $A_{i}$ denote the event ``the $i-$th
card drawn is an ace'' for $i=1,2,3,4.$ We are interested in computing
the probability 
\[
P\left(A_{1}\cap A_{2}\cap A_{3}\cap A_{4}\right).
\]
 Among the probabilities we can compute, apart $P\left(A_{1}\right)=\dfrac{4}{52},$
are the conditional probabilities:
\begin{itemize}
\item $P\left(A_{2}\mid A_{1}\right)=\dfrac{3}{51},$
\item $P\left(A_{3}\mid A_{1}\cap A_{2}\right)=\dfrac{2}{50},$
\item $P\left(A_{4}\mid A_{1}\cap A_{2}\cap A_{3}\right)=\dfrac{1}{49}.$
\end{itemize}
Therefore, the probability of drawing 4 aces in a row is
\begin{align*}
P\left(A_{1}\cap A_{2}\cap A_{3}\cap A_{4}\right) & =P\left(A_{1}\right)P\left(A_{2}\mid A_{1}\right)P\left(A_{3}\mid A_{1}\cap A_{2}\right)P\left(A_{4}\mid A_{1}\cap A_{2}\cap A_{3}\right)\\
 & =\dfrac{4\times3\times2\times1}{52\times51\times50\times49}\\
 & \approx3.69\times10^{-6}.
\end{align*}

We can observe that we did not take care in this example to construct
a space $\Omega$ where the $A_{i}$ could be considered as subsets.
This omission is intentional and will be justified later in Sub-section
\ref{sec:Evolutive-phenomena-modeling}, where the motivation of Proposition
$\left(\right.$\ref{pr:evol_phenomenon_modelling}$\left.\right)$
is.

\end{solutionexample}

\begin{definition}{Complete System of Constituents}{} Let $\left(A_{i}\right)_{i\in I}$
be a countable family of events pairwise disjoint\footnotemark---that
is no two events in the family occurs simultaneously. Suppose further
that
\begin{equation}
P\left(\biguplus_{i\in I}A_{i}\right)=1.
\end{equation}

Such a family is called a \textbf{\index{complete system of constituents}\mindex{constituent ! complete system}complete
system of constituents}.

\end{definition}

\footnotetext{We often say that those events are pairwise incompatible.}

Let $N$ be the complement of $\biguplus_{i\in I}A_{i}$ in the outcome
space $\Omega.$ Then $P\left(N\right)=0.$ The events $\left(A_{i}\right)_{i\in I},$
along with $N$ constitutes a partition of $\Omega.$ In probabilistic
terms, we interpret this as: with probability exactly 1, one---and
only one---of the events $A_{i}$ occurs.

\begin{example}{}{}We repeatedly toss a coin until the first occurrence
of a tail. 

Let $A_{i}$ be the event ``tail appears for the first time on the
i-th toss''. The events $A_{i},i\in\mathbb{N}^{\ast}$ are pairwise
incompatible. There is also a theoretical possibility that none of
the $A_{i}$ occurs---i.e. the coin indefinitely many times lands
heads---but the probability of this event is 0.

Therefore, the events $A_{i},i\in\mathbb{N}^{\ast}$ form a complete
system of constituents.

\end{example}

\begin{theorem}{Total Probability Formula}{}\label{Theo-Total-probability}

Let $\left(A_{i}\right)_{i\in I}$ be a complete system of constituents
such that for every $i\in I,\,P\left(A_{i}\right)>0.$ Then, for every
event $A\in\mathcal{A},$ \boxeq{
\begin{equation}
P\left(A\right)=\sum_{i\in I}P\left(A\mid A_{i}\right)P\left(A_{i}\right).\label{eq:total probability formula}
\end{equation}
}

\end{theorem}

\begin{proof}{}{}Let $N$ be the complement of $\biguplus_{i\in I}A_{i}.$
\[
A=\left[\biguplus_{i\in I}\left(A\cap A_{i}\right)\right]\biguplus\left(A\cap N\right).
\]

Since
\[
P\left(A\cap N\right)\leqslant P\left(N\right)=0,
\]
 we obtain immediately, by using the $\sigma-$additivity of $P,$
\begin{align*}
P\left(A\right) & =P\left[\biguplus_{i\in I}\left(A\cap A_{i}\right)\right]=\sum_{i\in I}P\left(A\cap A_{i}\right).
\end{align*}

Using the definition of the conditional probabilities $P\left(A\cap A_{i}\right)=P\left(A\mid A_{i}\right)P\left(A_{i}\right),$
we conclude
\begin{equation}
P\left(A\right)=\sum_{i\in I}P\left(A\mid A_{i}\right)P\left(A_{i}\right).\label{eq:total probability formula-1}
\end{equation}

\end{proof}

A special case of system of constituents is the system $\left(B,B^{c}\right),$
where $B$ is an event such that: $0<P\left(B\right)<1.$ In this
case, the total probability formula becomes\boxeq{
\[
P\left(A\right)=P\left(A\mid B\right)P\left(B\right)+P\left(A\mid B^{c}\right)P\left(B^{c}\right).
\]
}

\begin{example}{Urn and Balls}{example4.4}\label{example A-poll-with}

Suppose we have two urns $U_{1}$ and $U_{2}.$ Each urn $U_{i},$
with $i=1,2$ contains $w_{i}$ white balls and $b_{i}$ blue balls.
We randomly select one of the urns and draw one ball from it. 

What is the probability to draw a blue ball?

\end{example}

\begin{solutionexample}{}{}On Example \ref{ex:example4.9}, we will
revisit more deeply this example to build the model associated to
this experiment. For now, assume that a probabilized space $\left(\Omega,\mathcal{A},P\right)$
has been constructed such that, denoting $B$ the event ``we draw
a blue ball'' and $U_{i}$ the drawing takes place in the urn $i,\,i=1,2,$
we have for $i=1,2,$
\[
P\left(U_{i}\right)=\dfrac{1}{2}\,\,\,\,\text{and}\,\,\,\,P\left(B\mid U_{i}\right)=\dfrac{b_{i}}{b_{i}+w_{i}}.
\]

The events $U_{1}$ and $U_{2}$ constitute a complete system of constituents
and the probability $P\left(B\right)$ to draw a blue ball is given
by
\begin{align*}
P\left(B\right) & =P\left(B\mid U_{1}\right)P\left(U_{1}\right)+P\left(B\mid U_{2}\right)P\left(U_{2}\right).
\end{align*}

Hence,\boxeq{
\[
P(B)=\dfrac{1}{2}\left(\dfrac{b_{1}}{b_{1}+w_{1}}+\dfrac{b_{2}}{b_{2}+w_{2}}\right).
\]
}

\end{solutionexample}

\begin{figure}[t]
\begin{center}\includegraphics[width=0.4\textwidth]{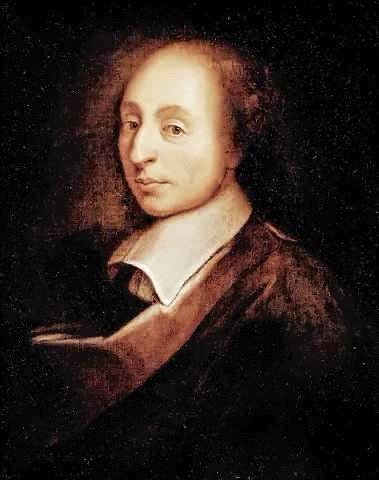}

{\tiny Public Domain}\end{center}

\caption{\textbf{\protect\href{https://en.wikipedia.org/wiki/Blaise_Pascal}{Blaise Pascal}}
(1623-1662)}\sindex[fam]{Pascal, Blaise}
\end{figure}

\begin{example}{The Second Problem of the Knight of Méré}{}

The origins of probability theory trace back to 1654, when \textbf{\href{https://en.wikipedia.org/wiki/Blaise_Pascal}{Blaise Pascal}}\sindex[fam]{Pascal, Blaise}\footnotemark
solved two problems posed by the Knight of Méré. One of them is the
following: 

\begin{leftbar}

``Two players are engaged in a game of chance played over several
rounds. The first player to win three rounds wins the entire stake.
Suppose now the game is interrupted when the first player needs only
one more round to win and the second needs two additional rounds to
win. How should the stake be divided to be fairly distributed?''

\end{leftbar}

\end{example}

\footnotetext{Blaise Pascal (1623-1662) scholar, thinker and writer.
Before 1654, his mathematical work are essentially of geometric nature.
In 1654, the Méré knight introduced him to Fermat. The writing correspondence
between Pascal and Fermat is partly at the origins of the probability
calculus.}

\begin{solutionexample}{}{}

To compute how the stake should be divided between players---which
Pascal explains that it has to be equal to what we would named today
the probability of the player to win---Pascal proposes a reasoning
on what we now call \textbf{conditional probabilities} and the \textbf{total
probability formula}.

Here is the way of how he argued: 
\begin{itemize}
\item Suppose one more round is played, the first player has one chance
over two of winning immediately. Thus, they should receive half the
stake. 
\item If the first player loses this round, then the players are now tied,
and from that point onward, each has the same probability to win.
So, each should receive half of the remaining half of the stake. 
\end{itemize}
Hence to fairly distribute the stake in this case, the first player
has to receive three quarters of the stake and the second player the
remaining.

\end{solutionexample}

\subsection{\textquotedblleft Probability of Causes\textquotedblright{} and Bayes
Formula}

In a two-phase experiment like the urn and ball in Example \ref{example A-poll-with},
conditional probability can be interpreted intuitively.

A priori, i.e. before that the choice of the urn is done, the probability
to draw a blue ball is $P\left(B\right).$ Suppose now that the first
phase of the experiment is done, i.e. the choice at random of an urn,
and suppose, for instance, that $U_{1}$ is chosen. The updated probability
of drawing a blue ball linked to this knowledge is 
\[
P\left(B\left|U_{1}\right.\right)=\dfrac{b_{1}}{b_{1}+w_{1}}.
\]

What matter is not the phases order of the experiment. Even if the
entire experiment is already performed, and we have learned afterwards
that the chosen urn was the first one---i.e. we just have a partial
information on the experiment---we still have to ensure to take into
account this last information. The probability to draw a blue ball
is $P\left(B\left|U_{1}\right.\right)$ as we have gained in knowledge.
Any decision that has to be taken---like gambling---depends on the
probability of the blue ball to appear in the first urn, hence it
convenes to use the conditional probability $P\left(B\left|U_{1}\right.\right)$
instead of the probability $P\left(B\right).$ It is why $P\left(U_{1}\right)$
is seen as a \textbf{\index{prior probability}}\mindex{probability ! prior}\textbf{
prior probability}: $U_{1}$ is observed first, then we observe $B.$ 

Now suppose that on the contrary, we know the color of the drawn ball---let
say a blue ball---but that we ignore in which urn it has been drawn.
It is still an incomplete information. Then the relevant quantity
is in this case,
\[
P\left(U_{1}\left|B\right.\right)
\]
which is the \textbf{posterior probability}\index{posterior probability}\mindex{probability ! posterior}.
It represents the probability to have drawn in urn $U_{1}$, knowing
that a blue ball was observed. 

This is often called \textbf{probability of causes}\mindex{probability!of causes},
because we attempt to infer the hidden cause---the urn---from an
observed effect---the color of the ball. To be clearer, if for instance,
the first urn contains a lot of blue balls, while the second has only
a few of them, the fact that a blue ball has been drawn let suppose
that it has most likely been drawn from the first urn when the drawing
was done. The number $P\left(U_{1}\mid B\right)$ measures the value
of the likelihood that we should give to this hypothesis. Mathematically,
however, there is no difference in nature between $P\left(B\mid U_{1}\right)$
and $P\left(U_{1}\mid B\right):$ they are just conditional probabilities!

\begin{proposition}{Link Between Prior and Posterior Probability}{}\label{Prop: link between prior and posterior}Let
$A_{1},A_{2}\in\mathcal{A}$ be two events of nonzero probability.
Then\boxeq{
\begin{equation}
P\left(A_{1}\mid A_{2}\right)=P\left(A_{2}\mid A_{1}\right)\dfrac{P\left(A_{1}\right)}{P\left(A_{2}\right)}.\label{eq:prior_posterior_bayes_at_2}
\end{equation}
}

\end{proposition}

\begin{proof}{}{}By definition of a conditional property,
\begin{align*}
P\left(A_{1}\mid A_{2}\right) & =\dfrac{P\left(A_{1}\cap A_{2}\right)}{P\left(A_{2}\right)}=\dfrac{P\left(A_{2}\mid A_{1}\right)P\left(A_{1}\right)}{P\left(A_{2}\right)}.
\end{align*}

\end{proof}

\begin{remark}{}{}

Tr.N.: In the previous formula $\refpar{eq:prior_posterior_bayes_at_2}$,
as it has many applications, often people refers to:
\begin{itemize}
\item $P\left(A_{1}\right)$ as the \textbf{prior probability\index{prior probability}}\mindex{probability ! prior}
of event $A_{1}$ occuring before observing event $A_{2}.$
\item $P\left(A_{2}\mid A_{1}\right)$ as the \textbf{likelihood\index{likelihood}}
as it is the probability of observing $A_{2}$ knowing that $A_{1}$
has been observed.
\item $P\left(A_{2}\right)$ as the \textbf{evidence probability\index{evidence probability}}
\mindex{probability ! evidence}as it is the probability of $A_{2}$
to occur.
\item $P\left(A_{1}\mid A_{2}\right)$ as the \textbf{posterior probability\index{posterior probability}\mindex{probability ! posterior}}
of event $A_{1}$ occuring knowing that $A_{2}$ has been observed.
\end{itemize}
This formula allowing the computation of the posterior probability
is intensively used in many domains, particularly in medical applications,
risk assessment, economics and machine learning with the different
flavours of Bayesian classifiers. As already mentioned, this naming
is a convention as it all depends on the way the $A_{1}$ and $A_{2}$
are observed!

\end{remark}

\begin{theorem}{Bayes Theorem}{}

Let $\left(A_{i}\right)_{i\in I}$ be a complete system of constituents
such that for every $i\in I,$ $P\left(A_{i}\right)>0.$ 

Then for every event $A\in\mathcal{A}$ of nonzero probability and
for every $i\in I,$\boxeq{
\begin{equation}
P\left(A_{i}\mid A\right)=\dfrac{P\left(A\mid A_{i}\right)P\left(A_{i}\right)}{\sum_{j\in I}P\left(A\mid A_{j}\right)P\left(A_{j}\right)}.\label{eq:Bayes formula}
\end{equation}
}

\end{theorem}

\begin{proof}{}{}From Proposition \ref{Prop: link between prior and posterior},
we have for every $i\in I,$
\[
P\left(A_{i}\mid A\right)=\dfrac{P\left(A\mid A_{i}\right)P\left(A_{i}\right)}{P\left(A\right)}.
\]

It is then enough to apply Theorem \ref{Theo-Total-probability}.

\end{proof}

\begin{example}{}{}Let us return to Example \ref{example A-poll-with}. 

What is the posterior probability that the draw was made in $U_{1}$
knowing that a blue ball was drawn?

\end{example}

\begin{solutionexample}{}{}We apply Bayes theorem since the events
$U_{1}$ and $U_{2}$ are a complete system of constituents. For $i=1,2,$
Bayes formula gives
\begin{align*}
P\left(U_{i}\mid B\right) & =\dfrac{P\left(U_{i}\right)P\left(B\mid U_{i}\right)}{P\left(U_{1}\right)P\left(B\mid U_{1}\right)+P\left(U_{2}\right)P\left(B\mid U_{2}\right)}
\end{align*}

Since $P\left(U_{1}\right)=P\left(U_{2}\right),$ we have for $i=1,2,$\boxeq{
\[
P\left(U_{i}\mid B\right)=\dfrac{\dfrac{b_{i}}{b_{i}+w_{i}}}{\dfrac{b_{1}}{b_{1}+w_{1}}+\dfrac{b_{2}}{b_{2}+w_{2}}}.
\]
}

\end{solutionexample}

\section{Conditional Laws}\label{sec:Conditional-Laws}

When studying a random phenomena, we are often given a partial information
in the form of a random variable $X$ taking values in a probabilizable
space $\left(E,\mathcal{E}\right).$ Suppose that $X$ is a discrete
random variable with values in a set $E,$ and let $x$ be a value
in $E$ such that $P\left(X=x\right)>0.$ 

Let $Y$ be another random variable, taking values in a probabililizable
space $\left(F,\mathcal{F}\right).$ By definition, the \textbf{law
of $Y$ conditional}{\bfseries\footnote{We also say: the ``law of $Y$ knowing that $X=x.$''}}\textbf{
to the event}\index{law of ensuremath{Y} conditional to the event $left(X=xright)@law of $Y$ conditional to the event $\left(X=x\right)$}
$\left(X=x\right)$ is the probability on $\left(F,\mathcal{F}\right)$
defined by the function
\[
B\mapsto P^{\left(X=x\right)}\left(Y\in B\right)\equiv P\left(Y\in B\mid X=x\right).
\]
This last is written $P^{\left(X=x\right)}_{Y}.$

If $Y$ is also a discrete random variable, the family of probabilities
$P^{\left(X=x\right)}_{Y},$ where $x$ describes the set of events
$E$ such that $P\left(X=x\right)>0$ together with the law $P_{X}$
of $X,$ completely determines the law of the random variable $\left(X,Y\right).$

Indeed, since $\left(X,Y\right)$ is also a discrete random variable,
its law is fully specified by the fact that probabilities $P\left(X=x,Y=y\right)$
are given for every $x\in E$ and $y\in F.$ 

But if $x\in E$ and $y\in F,$ there are two cases:
\begin{itemize}
\item If $P\left(X=x\right)=0,$ then
\[
P\left(X=x,Y=y\right)\leqslant P\left(X=x\right)=0
\]
which induces that
\[
P\left(X=x,Y=y\right)=0.
\]
\item If $P\left(X=x\right)>0,$ then we can compute
\begin{align*}
P\left(X=x,Y=y\right) & =P\left(Y=y\mid X=x\right)P\left(X=x\right)\\
 & =P^{\left(X=x\right)}\left(Y=y\right)P\left(X=x\right).
\end{align*}
\end{itemize}
Thus, in all cases, the joint probability $P\left(X=x,Y=y\right)$
can always be determined from the conditional laws $P^{\left(X=x\right)}$
and the marginal law of $X.$

\begin{example}{}{}

Let $X$ and $Y$ be two discrete random variables taking values in
$\mathbb{N},$ defined on a probabilized space $\left(\Omega,\mathcal{A},P\right).$
Assume that $X$ follows a Poisson law $\mathcal{P}\left(\lambda\right),$
where $\lambda>0.$ Furthermore, suppose that for every integer $n>0,$
the conditional law of $Y$ given $\left(X=n\right)$ is a binomial
law $\mathcal{B}\left(n,p\right),$ and that if $X=0,$ the random
variable $Y$ takes the value 0 with probability 1.

Our goal is:

1. To determine the law of $Y.$

2. For each $k\in\mathbb{N},$ to dermine the law of $X$ conditional
on the event $\left(Y=k\right).$

\end{example}

\begin{solutionexample}{}{}

\textbf{1. Determination of the law of $Y.$}

Let $q=1-p.$ For every $\left(n,k\right)\in\mathbb{N}^{2},$
\[
P\left(X=n,Y=k\right)=P^{\left(X=n\right)}\left(Y=k\right)P\left(X=n\right).
\]

From the hypothesis, we deduce that
\[
P\left(X=n,Y=k\right)=\begin{cases}
0, & \text{if }0\leqslant n\leqslant k,\\
\text{e}^{-\lambda}\dfrac{\lambda^{n}}{n!}\binom{n}{k}p^{k}q^{n-k}, & \text{if\,}0\leqslant k\leqslant n.
\end{cases}
\]

Since $\biguplus_{n\in\mathbb{N}}\left(X=n\right)=\Omega,$
\[
P\left(Y=k\right)=\sum_{n\in\mathbb{N}}P\left(X=n,Y=k\right).
\]

Hence, it holds:
\begin{itemize}
\item If $k<0,$ 
\[
P\left(Y=k\right)=0.
\]
\item If $k\geqslant0,$
\[
P\left(Y=k\right)=\sum^{+\infty}_{n=k}\text{e}^{-\lambda}\dfrac{\lambda^{n}}{n!}\left(\begin{array}{c}
n\\
k
\end{array}\right)p^{k}q^{n-k}.
\]
\end{itemize}
By factorization of terms that are independent of $n,$
\[
P\left(Y=k\right)=\text{e}^{-\lambda}\dfrac{\left(\lambda p\right)^{k}}{k!}\sum^{+\infty}_{n=k}\dfrac{1}{\left(n-k\right)!}\left(\lambda q\right)^{n-k},
\]

Substituting $n-k$ by $m$ in the summation, gives
\[
P\left(Y=k\right)=\text{e}^{-\lambda}\dfrac{\left(\lambda p\right)^{k}}{k!}\sum^{+\infty}_{m=0}\dfrac{1}{m!}\left(\lambda q\right)^{m},
\]

Which yields after reduction\boxeq{
\[
P\left(Y=k\right)=\text{e}^{-\lambda p}\dfrac{\left(\lambda p\right)^{k}}{k!}.
\]
}

The law of $Y$ is thus the Poisson law of parameter $\lambda p.$

\textbf{2. Determination of the law of $X$ conditional on the event
$Y=k.$}

By definition of a conditional probability, for every $n\in\mathbb{N}$
and for every $k\in\mathbb{N},$
\[
P^{\left(Y=k\right)}\left(X=n\right)=\dfrac{P\left(X=n,Y=k\right)}{P\left(Y=k\right)}.
\]

We then deduce using the previous question:
\begin{itemize}
\item If $n<k,$ then
\[
P^{\left(Y=k\right)}\left(X=n\right)=0.
\]
\item If $0\leqslant k\leqslant n,$ then
\begin{align*}
P^{\left(Y=k\right)}\left(X=n\right) & =\dfrac{\text{e}^{-\lambda}\dfrac{\lambda^{n}}{n!}\left(\begin{array}{c}
n\\
k
\end{array}\right)p^{k}q^{n-k}}{\text{e}^{-\lambda p}\dfrac{\left(\lambda p\right)^{k}}{k!}}=\text{e}^{-\lambda q}\dfrac{\left(\lambda q\right)^{n-k}}{\left(n-k\right)!}.
\end{align*}
\end{itemize}
\boxeq{For any fixed integer $k,$ we conclude that the conditional
law of $X$ knowing $\left(Y=k\right)$ is the \textbf{Poisson law
with parameter $\lambda q$ shifted by $k.$}}

\end{solutionexample}

\begin{remark}{}{}A concrete example of such modelling can be illustrated
using a crossroad. Let $X$ be the random variable representing the
number of cars arriving during one hour interval. It is assumed that
$X$ follows a Poisson law $\mathcal{P}\left(\lambda\right)$---a
common assumption in this kind of scenario\footnotemark.

Cars can only turn one of the two directions $A$ or $B.$ Let $Y$
be the random variable representing the number of cars that take direction
$A$ during that time interval. Each car independently choose direction
$A$ with probability $p.$ Therefore, if we suppose that $n$ cars
arrive at the crossroad during the time slot, the conditional law
of $Y$ follows a binomial law $\mathcal{B}\left(n,p\right).$

The conditional probability $P^{\left(Y=k\right)}\left(X=n\right)$
represents the probability that $n$ vehicles arrived at the crossroads
knowing that $k$ of them took direction $A.$ This is thus an example
of a ``probability of causes.''

\end{remark}

\footnotetext{See the footnote of Exercice \ref{exo:exercise3.6}.}

\section{Evolutive Phenomenon Modelling}\label{sec:Evolutive-phenomena-modeling}

The construction of a probabilistic model for a concrete random experiment
generally begins with the definition of a probabilized space $\left(\Omega,\mathcal{A},P\right).$
We then define the events of interest---such as the subsets of $\Omega$---as
well as certain random variables---such as functions defined on $\Omega$---representing
the quantities related to the experiment. However, explicit knowledge
of $\Omega$ is rarely required. It is often just needed to know which
events are independent, which random variable $X$ follows a given
law, and so forth. The space $\Omega$ thus takes a backseat.

In everyday modelling problems, the space $\Omega$ is seldom determined
by the situation under consideration---unless in a few cases related
to games of chance. What is typically accessible are some imposed
constaints on the data: for instance, that the events $A_{1},\dots,A_{n}$
are equally probable---often due to symmetry considerations---or
that the random variables $X$ and $Y$ are independent, or that the
random variable $X$ follows a particular law. It is then the mathematician's
task to show the existence of a space $\Omega,$ of events $A_{1},A_{2},\dots,A_{n}$
defined as subsets of this space, and of random variables defined
as functions on $\Omega,$ all satisfying the required properties.

In Chapter \ref{chap:Independence}, we showed that, given some discrete
probability laws $P_{1},\dots,P_{n}$ on some countable sets $E_{1},\dots,E_{n},$
there exists a probabilized space $\left(\Omega,\mathcal{A},P\right)$
and some random variables $X_{1},\dots,X_{n}$\footnote{Recall that we set $\Omega=E_{1}\times\dots\times E_{n}$ and defined
$X_{i}$ as the $i-$th projection $\text{pr}_{i}:\Omega\to E_{i}.$} defined on $\left(\Omega,\mathcal{A},P\right),$ with values respectively
in $E_{1},\dots,E_{n}$ and such that:

(a) $X_{i}$ follows the law $P_{i}\,\left(i=1,\dots,n\right);$

(b) $X_{1},\dots,X_{n}$ are independent.

In many random experiments composed of successive steps---such as
those served as example in this chapter---what is immediately accessible
in practice are the \textbf{conditional probabilities}.

\begin{proposition}{Evolutive Phenomenon Modelling}{evol_phenomenon_modelling}

Let $E_{1},\dots,E_{n}$ be countable sets. Suppose we are given:
\begin{itemize}
\item A germ of probability law on $E_{1},$ denoted $p_{1}.$
\item For each $x_{1}\in E_{1},$ a germ of probability law on $E_{2},$
denoted $p^{x_{1}}_{2}.$
\item ...
\item For each $\left(x_{1},\dots,x_{n-1}\right)\in E_{1}\times\dots\times E_{n-1},$
a germ of probability law on $E_{n},$ denoted $p^{x_{1},\dots,x_{n-1}}_{n}.$
\end{itemize}
Then, there exists a probabilized space $\left(\Omega,\mathcal{A},P\right)$
and some random variables $X_{1},\dots,X_{n}$ defined on this space,
taking values in $E_{1},\dots,E_{n}$ respectively, such that for
every $x_{1}\in E_{1},\dots,x_{n}\in E_{n},$

\boxeq{
\begin{equation}
\begin{cases}
P\left(X_{1}=x_{1}\right)=p_{1}\left(x_{1}\right)\\
P\left(X_{2}=x_{2}\mid X_{1}=x_{1}\right)=p^{x_{1}}_{2}\left(x_{2}\right)\\
\dots\\
P\left(X_{n}=x_{n}\mid X_{1}=x_{1},\dots,X_{n-1}=x_{n-1}\right)=p^{x_{1},\dots,x_{n-1}}_{n}\left(x_{n}\right).
\end{cases}\label{eq:cond_prob_ev_phenomenon_modelling}
\end{equation}
}

For instance, we can take as $\Omega=E_{1}\times\dots\times E_{n},$
$\mathcal{A}=\mathcal{P}\left(\Omega\right),$ and define $P$ from
its germ $p$ by setting, for $\omega=\left(x_{1},\dots,x_{n}\right),$\boxeq{
\begin{equation}
p\left(\omega\right)=p_{1}\left(x_{1}\right)p^{x_{1}}_{2}\left(x_{2}\right)\dots p^{x_{1},\dots,x_{n-1}}_{n}\left(x_{n}\right),\label{eq:P_from_germ_ev_phen_modelling}
\end{equation}
} and take for $X_{i}\,\left(1\leqslant i\leqslant n\right)$ the
$i-th$ projection defined by $X_{i}\left(x_{1},\dots,x_{n}\right)=x_{i}.$

\end{proposition}

\begin{proof}{}{}

It suffices to verify that the proposed solution given in the previous
paragraph remains valid.

By the Fubini theorem for families with nonnegative terms, we can
write
\[
\sum_{\omega\in\Omega}p\left(\omega\right)=\sum_{x_{1}\in E_{1}}p_{1}\left(x_{1}\right)\left(\sum_{x_{2}\in E_{2}}p^{x_{1}}_{2}\left(x_{2}\right)\dots\left(\sum_{x_{n}\in E_{n}}p^{x_{1},\dots,x_{n-1}}_{n}\left(x_{n}\right)\right)\right).
\]

By successively performing the summations from right to left, we verify
that
\[
\sum_{\omega\in\Omega}p\left(\omega\right)=1
\]
 so that $p$ indeed defines a probability germ. The relations stated
in $\refpar{eq:cond_prob_ev_phenomenon_modelling}$ then follow directly
from the chain rule for conditionnal probabilities. 

\end{proof}

\begin{remarks}{}{}

1. If $p^{x_{1},\dots,x_{i-1}}_{i}\left(x_{i}\right)$ depends only
on $x_{i}$---that is, if for each $i=1,\dots,n$ and for every $\left(x_{1},\dots,x_{i-1}\right)\in E_{1}\times\dots\times E_{i-1},$
the germ $p^{x_{1},\dots,x_{i-1}}_{i}$ is a germ $g_{i}$ on $E_{i}$
independent of $x_{1},\dots,x_{i-1}$---then the probability $P$
is the \textbf{product probability} of the probabilities associated
with the germs $g_{i}.$

In this case, the random variables $X_{1},\dots,X_{n}$ satisfying
the relations in $\refpar{eq:cond_prob_ev_phenomenon_modelling}$
are necessarily independent, since the relations
\[
P\left(X_{i}=x_{i}\mid X_{1}=x_{1},\dots,X_{i-1}=x_{i-1}\right)=P\left(X_{i}=x_{i}\right)
\]
imply, by staightforward induction, that
\[
P\left(X_{1}=x_{1},\dots,X_{i-1}=x_{i-1},X_{i}=x_{i}\right)=P\left(X_{1}=x_{1}\right)\dots P\left(X_{i-1}=x_{i-1}\right)P\left(X_{i}=x_{i}\right).
\]

Thus we recover the construction presented in Section \ref{sec:Independence-and-Cartesian}.

2. If $p^{x_{1},\dots,x_{i-1}}_{i}\left(x_{i}\right)$ only depends
on $x_{i}$ and $x_{i-1},$ meaning that for each $i=2,\dots,n,$
the germ $p^{x_{1},\dots,x_{i-1}}_{i}$ is a germ on $E_{i}$ depending
solely on $x_{i-1},$ then\boxeq{
\begin{equation}
P\left(X_{i}=x_{i}\mid X_{1}=x_{1},\dots,X_{i-1}=x_{i-1}\right)=P\left(X_{i}=x_{i}\mid X_{i-1}=x_{i-1}\right).\label{eq:markovian_condition}
\end{equation}
}

A sequence of random variables $\left(X_{i}\right)_{1\leqslant i\leqslant n}$
satisfying the condition $\refpar{eq:markovian_condition}$ is called
a \index{Markovian sequence}\textbf{Markovian}\footnotemark\textbf{
sequence}. 

\end{remarks}

\footnotetext{In reference to the russian mathematician \textbf{Andreï
Markov\sindex[fam]{Markov, Andreï}} (1856-1922), a student of P.
Tchebichev, author of many foundational works in probability theory.}

\begin{figure}[t]
\begin{center}\includegraphics[width=0.4\textwidth]{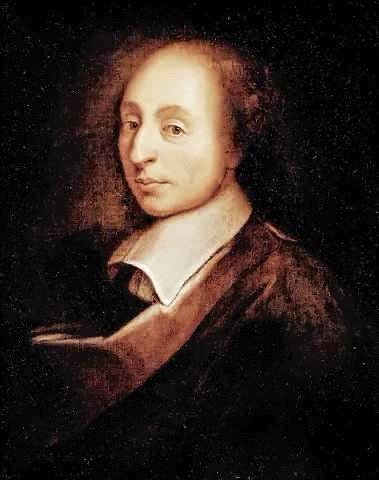}

{\tiny Public Domain}\end{center}

\caption{\textbf{\protect\href{https://en.wikipedia.org/wiki/Blaise_Pascal}{Blaise Pascal}}
(1623-1662)}\sindex[fam]{Pascal, Blaise}
\end{figure}

\begin{example}{Markovian Sequence}{}

Let $\left(U_{i}\right)_{1\leqslant i\leqslant n}$ be a sequence
of independent random variables taking values in $\mathbb{Z}.$ 

Define, for $i=1,\dots,n,$ the random variables 
\[
X_{i}=U_{1}+\dots+U_{i}.
\]

Prove that the sequence $\left(X_{i}\right)_{1\leqslant i\leqslant n}$
is a Markovian sequence taking values in $\mathbb{Z}.$ 

\end{example}

\begin{solutionexample}{}{}

For $i\geqslant2,$ we can recursively express $X_{i}$ as
\[
X_{i}=X_{i-1}+U_{i}.
\]
Then, using the previous notations, for every $x_{1},\dots,x_{i}\in\mathbb{Z},$
\begin{multline*}
P\left(X_{i}=x_{i}\mid X_{1}=x_{1},\dots,X_{i-1}=x_{i-1}\right)\\
\,\,\,\,\,\,\,\,\,\,\,\,\,=P\left(U_{i}=x_{i}-x_{i-1}\mid X_{1}=x_{1},\dots,X_{i-1}=x_{i-1}\right).
\end{multline*}

Since the random variable $\left(X_{1},\dots,X_{i-1}\right)$ is only
function of $U_{1},\dots,U_{i-1},$ which are all independent of $U_{i},$
it follows that $U_{i}$ is independent of $\left(X_{1},\dots,X_{i-1}\right).$
Therefore, 
\[
P\left(X_{i}=x_{i}\mid X_{1}=x_{1},\dots,X_{i-1}=x_{i-1}\right)=P\left(U_{i}=x_{i}-x_{i-1}\right).
\]

Moreover, by the same reasoning, \boxeq{
\[
P\left(X_{i}=x_{i}\mid X_{i-1}=x_{i-1}\right)=P\left(U_{i}=x_{i}-x_{i-1}\right).
\]
}which proves the result.

\end{solutionexample}

\begin{example}{Follow-up on the Urn Example}{example4.9}

Returning to Example \ref{ex:example4.4}, construct the associated
probabilized space, as previously announced.

\end{example}

\begin{solutionexample}{}{}

We start by reformulating the modelling problem in terms of random
variables.

Let $E_{1}=\left\{ 1,2\right\} $ be the set representing the choice
of the urn, and $E_{2}=\left\{ b,w\right\} $ the choice of colors
(blue or white). Our goal is to construct a model consisting, on the
one hand, of a probabilized space $\left(\Omega,\mathcal{A},P\right),$
and on the other, of two random variables $X_{1}:\Omega\to E_{1}$
and $X_{2}:\Omega\to E_{2}$ satisfying the following relations for
each $i=1,2,$
\[
\begin{cases}
P\left(X_{1}=i\right)=\dfrac{1}{2},\\
P\left(X_{2}=b\mid X_{1}=i\right)=\dfrac{b_{i}}{b_{i}+w_{i}}, & P\left(X_{2}=w\mid X_{1}=i\right)=\dfrac{w_{i}}{b_{i}+w_{i}}.
\end{cases}
\]

Proposition $\left(\ref{pr:evol_phenomenon_modelling}\right)$ indicates
that we can set the universe as $\Omega=E_{1}\times E_{2},$ the $\sigma-$algebra
as $\mathcal{A}=\mathcal{P}\left(\Omega\right),$ and a probability
$P$ by specifying for each $i=1,2,$
\[
\begin{cases}
P\left(\left\{ \left(i,b\right)\right\} \right)=\dfrac{1}{2}\dfrac{b_{i}}{b_{i}+w_{i}},\\
P\left(\left\{ \left(i,w\right)\right\} \right)=\dfrac{1}{2}\dfrac{w_{i}}{b_{i}+w_{i}}.
\end{cases}
\]

The random variables $X_{1}$ and $X_{2}$ are then simply the projection
of the Cartesian product $\Omega=E_{1}\times E_{2}$ onto their respective
factor spaces.

\end{solutionexample}

\section*{Exercises}

\addcontentsline{toc}{section}{Exercises}

\begin{exercise}{Bayes Formula and the Law of Total Probability. Heads and Tails Game}{exercise4.1}

A player plays a game of heads and tails using two fair and balanced
coins as follows: 
\begin{itemize}
\item Step 1: They toss both coins simultaneously a first time. 
\begin{itemize}
\item If no tails appear, their gain is zero and the game ends. 
\end{itemize}
\item Step 2: If at least one tail appears, then they toss both coins again
simultaneously as many times as the number of tails obtained in the
first toss.
\begin{itemize}
\item Their gain is equal to the total number of tails obtained during this
second series of tosses.
\end{itemize}
\end{itemize}
1. What is the probability that the player's gain is not zero? 

2. What is the probability that the player obtain two tails on the
first toss, given that exactly one tail was obtained in the second
series?

\end{exercise}

\begin{exercise}{Urns Model}{exercice4.2}

Two urns, $A$ and $B,$ contain the following:
\begin{itemize}
\item Urn $A:$ 6 blue balls and 5 green balls.
\item Urn $B:$ 4 blue balls and 8 green balls.
\end{itemize}
Two balls are randomly transferred from urn $B$ to urn $A.$ Then
a single ball is drawn at random from urn $A.$

1. What is the probability that the drawn ball is blue?

2. Given that the drawn ball is blue, what is the probability that
at least one of the transferred balls was blue?

\end{exercise}

\begin{exercise}{A Simplified Demographic Model}{exercise4.3}

In a given population, the probability that a household to have exactly
$k$ children is defined as follows, for a real number $a$ in the
range 0 to 1 :
\begin{itemize}
\item $p_{0}=p_{1}=a.$
\item For every $k\geqslant2,$ $p_{k}=\left(1-2a\right)2^{-\left(k-1\right)}.$
\end{itemize}
Furthermore, we assume (for simplicity) that each child is independtly
assigned the label ``boy'' or ``girl''\footnotemark with equal
probability $\dfrac{1}{2}.$

1. What is the probability that a household with exactly two boys
is only composed of two children?

2. What is the probability that a household has two girls, given that
it has two boys?

\end{exercise}

\footnotetext{Tr.N: we suppose to be in a country where only two
genders are recognized, which is no more the case in more advanced
countries.}

\begin{exercise}{Chain Rule of Conditional Probabilities}{exercise4.4}

An urn contains $n$ red balls and $n$ green balls. Pairs of balls
are drawn sequentially, without replacement, until the urn is empty.

What is the probability that each drawn pair consists of one red and
one green ball?

\end{exercise}

\begin{exercise}{Chain Rule of Conditional Probabilities and Time Before Success. Non-positive Binomial Law}{exercise4.5}

Let $\left(\Omega,\mathcal{A},P\right)$ be a probabilized space,
and let $\left(A_{n}\right)_{n\in\mathbb{N}^{\ast}}$ be a sequence
of independent events, each occuring with probability $p.$ Let $q=1-p.$
With the convention $\text{inf}\emptyset=+\infty,$ define a sequence
of random variables $\left(T_{n}\right)_{n\in\mathbb{N}^{\ast}}$
by setting, for every $\omega\in\Omega,$
\begin{align*}
T_{1}\left(\omega\right) & =\text{inf\ensuremath{\left(j\geqslant1:\,\omega\in A_{j}\right)},}\\
T_{2}\left(\omega\right) & =\text{inf\ensuremath{\left(j\geqslant T_{1}\left(\omega\right):\,\omega\in A_{j}\right)},}\\
\dots\\
T_{n+1}\left(\omega\right) & =\text{inf}\left(j>T_{n}\left(\omega\right):\,\omega\in A_{j}\right).
\end{align*}

1. (a) For every $k\in\mathbb{N}^{\ast},$ determine the law of the
random variable $T_{k}.$

(b) For every increasing sequence of integers $n_{1},n_{2},\dots,n_{k+1},$
compute the conditional probability
\[
P\left(T_{k+1}-T_{k}=n_{k+1}-n_{k}\left|T_{1}=n_{1},\dots,T_{k}=n_{k}\right.\right).
\]

(c) Deduce that, if $n_{1}<n_{2}<\dots<n_{k+1},$ then
\[
P\left(T_{k+1}=n_{k+1}\left|T_{1}=n_{1},\dots,T_{k}=n_{k}\right.\right)=P\left(T_{k+1}=n_{k+1}\left|T_{k}=n_{k}\right.\right).
\]

(d) Determine the value of the joint probability $P\left(T_{1}=n_{1},\dots,T_{k}=n_{k}\right).$

2. Prove that the random variables $T_{1},T_{2}-T_{1},\dots,T_{k}-T_{k-1}$
are independent, determine their law and recover the fact that the
$k$-fold convolution of the geometric law $\mathcal{G}_{\mathbb{N}^{\ast}}\left(p\right)$---i.e.
the law of the sum of $k$ independent random variables of same law
$\mathcal{G}_{\mathbb{N}^{\ast}}\left(p\right)$---is the negative
binomial law $\mathcal{B}^{-}\left(k,p\right).$

3. Compute the probability
\[
P\left(T_{1}=3,T_{5}=9,T_{7}=17\right).
\]
 Give numerical values rounded to three significant digits for the
following values of $p=$0.3, 0.5 and 0.8.

\end{exercise}

\begin{example}{Game of Heads and Tails}{}

Let $A_{n}$ be the event ``a tail appears on the $n$-th coin toss.''
Then the random variable $T_{k}$ represents the number of the toss
on which the $k$-th tail appears in a sequence of independent coin
toss.

More generally, this framework applies to any sequence of independent
experiments where a particular property may or may not occur and where
the probability of the outcome $p$ remains constant accross trials.
In this context, $A_{n}$ is the event ``the property is realized
at the $n$-th trial,'' and the random variable $T_{k}$ represents
the number of the trial on which the outcome occurs for the $k$-th
time. 

\end{example}

\begin{exercise}{Probability on a Product Space and Total Probability Formula}{exercise4.6}

A game involves three cards:
\begin{itemize}
\item The first card has two red faces,
\item The second card has two green faces,
\item The last one has one red face and one green face.
\end{itemize}
A card is selected at random, and then one of its two faces is shown
at random to a spectator. The spectator must then bet on the color
of the hidden face.

What is the optimal strategy for the spectator to maximize their chances
of winning?

\end{exercise}

\begin{exercise}{Random Variables of Geometric Law, Conditioning and Independence}{exercise4.7}

All random variables introduced are defined on the same probabilizable
space $\left(\Omega,\mathcal{A}\right).$ Let $P$ be a probability
on that space.

Let $q$ and $r$ be two real numbers strictly bounded by 0 and 1.
Consider two independent discrete random variables $U$ and $V,$
each following a geometric law on $\mathbb{N}^{\ast},$ with respective
parameters $1-q$ and $1-r.$ Those laws are respectively denoted
$\mathcal{G}_{\mathbb{N}^{\ast}}\left(1-q\right)$ and $\mathcal{G}_{\mathbb{N}^{\ast}}\left(1-r\right).$

1. The following probability computations are useful for instance
in comparing the time before success of two players playing simultaneously
under conditions that potentially differ.

(a) Compute $P\left(U<V\right).$

(b) What is the value of this probability in the case where $q=r,$
and then specifically for $q=r=\dfrac{1}{2}?$

2. We study the law of $U$ under the conditional probability $P^{\left(U<V\right)}.$

(a) For every $k\in\mathbb{N}^{\ast},$ compute the conditional probability
\[
P\left(U=k\left|U<V\right.\right).
\]

Identify the conditional probability of $U$ given ``$U<V$'', i.e.
the probability on $\mathbb{N}^{\ast}$ generated by the germ: 
\[
k\mapsto P\left(U=k\left|U<V\right.\right).
\]

(b) Determine, for every $k\in\mathbb{N},$ the conditional probability
\[
P\left(U>k\left|U<V\right.\right).
\]

3. Let $X_{1}$ and $X_{2}$ be two independent discrete random variables
each following the same geometric law on $\mathbb{N}^{\ast}$ with
parameter $q.$ 

Define the sets
\[
A_{1}=\left\{ \omega\in\Omega:\,X_{1}\left(\omega\right)<X_{2}\left(\omega\right)\right\} \,\,\,\,\,A_{2}=\left\{ \omega\in\Omega:\,X_{2}\left(\omega\right)<X_{1}\left(\omega\right)\right\} 
\]
and
\[
H=\left\{ \omega\in\Omega:\,X_{1}\left(\omega\right)\neq X_{2}\left(\omega\right)\right\} .
\]

Consider the random variables
\[
X=\min\left(X_{1},X_{2}\right)\,\,\,\,\,\text{and}\,\,\,\,\,J=\boldsymbol{1}_{A_{1}}+2\boldsymbol{1}_{A_{2}}.
\]

Let $P^{H}$ denote the conditional probability on the probabilizable
space $\left(\Omega,\mathcal{A}\right)$ defined, for every $A\in\mathcal{A},$
by 
\[
P^{H}\left(A\right)=P\left(A\mid H\right).
\]

(a) Compute, for $k\in\mathbb{N}^{\ast},$ the probability
\[
P\left(H\cap\left(J=1\right)\cap\left(X>k\right)\right).
\]

(b) Justify the equality
\[
P\left(A_{1}\right)=P\left(A_{2}\right)
\]
and, deduce for every $k\in\mathbb{N},$ the value of the probability
\[
P^{H}\left(\left(J=1\right)\cap\left(X=k\right)\right).
\]

(c) What is the value of $P^{H}\left(\left(J=1\right)\right)?$

4. We study the independence property of the random variables $X$
and $J$ under the probability $P^{H}.$

(a) Note that we have the equality
\[
P\left(H\cap\left(J=1\right)\cap\left(X>k\right)\right)=P\left(H\cap\left(J=2\right)\cap\left(X>k\right)\right)
\]
and deduce from previous results, the value of the probability $P^{H}\left(X>k\right).$

(b) Prove that the random variables $X$ and $J$ are independent
on the probabilized space $\left(\Omega,\mathcal{A},P^{H}\right).$

\begin{examplen}{Heads and Tails Game with Two Players}{}

Two players play a game of heads and tails using the same coin, taking
turns to toss the coin alternately. The random variable $X_{1}$---respectively
$X_{2}$---corresponds to the number of the toss on which the first---respectively
the second---player obtains a tail for the first time. We assume
that tosses are counted independently for the two players.

The random variable $X$ corresponds to the number of the toss on
which one of the players obtains their first tail. 

Define a random variable $J$ as follows:
\begin{itemize}
\item $J=1,$ if the first player gets a tail before the second
\item $J=2,$ if the second player gets a tail before the first,
\item $J=0,$ if both players succeed to get a tail at the same time.
\end{itemize}
We can show that conditionally to the fact the players do not get
a tail with the same number of tosses, the random variables $X$ and
$J$ are independent.

\end{examplen}

\end{exercise}

\section*{Solutions of Exercises}

\addcontentsline{toc}{section}{Solutions of Exercises}

\begin{solution}{}{solexercise4.1}

\textbf{1. Probability of a nonzero gain}

We suppose defined a probabilized space $\left(\Omega,\mathcal{A},P\right)$
such that, $X_{1}$ designating the random variable giving the number
of tails obtained at the first step, its law is given by
\[
\begin{array}{ccccc}
P\left(X_{1}=0\right)=\dfrac{1}{4}, &  & P\left(X_{1}=1\right)=\dfrac{1}{2}, &  & P\left(X_{1}=2\right)=\dfrac{1}{4}\end{array}
\]
and, $X_{2}$ denoting the random variable representing the player's
gain, it holds
\[
\left(X_{1}=0\right)\subset\left(X_{2}=0\right).
\]

Furthermore, the conditional law of $X_{2}$ given $X_{1}$ is---partially---specified
as
\begin{align*}
P\left(X_{2}=0\mid X_{1}=1\right) & =\dfrac{1}{4},\\
P\left(X_{2}=1\mid X_{1}=1\right) & =\dfrac{1}{2},\\
P\left(X_{2}=2\mid X_{1}=1\right) & =\dfrac{1}{4},\\
P\left(X_{2}=0\mid X_{1}=2\right) & =\dfrac{1}{16},\\
P\left(X_{2}=1\mid X_{1}=2\right) & =\dfrac{1}{4}.
\end{align*}

Now observe that
\[
\left(X_{2}=0\right)=\left(X_{1}=0\right)\biguplus\left[\left(X_{1}>0\right)\cap\left(X_{2}=0\right)\right],
\]

By the law of total probability,
\begin{align*}
P\left(X_{2}=0\right) & =P\left(X_{1}=0\right)+P\left(X_{1}=1\right)P\left(X_{2}=0\mid X_{1}=1\right)\\
 & \qquad\qquad\qquad+P\left(X_{1}=2\right)P\left(X_{2}=0\mid X_{1}=2\right)\\
 & =\dfrac{25}{64}.
\end{align*}

Thus, since we look for the probability to have a nonzero gain,\boxeq{
\begin{align*}
P\left(X_{2}>0\right) & =1-P\left(X_{2}=0\right)=\dfrac{39}{64}.
\end{align*}
}

\textbf{2. Probability to obtain two tails at the first toss given
that only one tail is obtained at the second step.}

The probability to obtain two tails at the first throwing knowing
they obtained only one tail at the second step is $P\left(X_{1}=2\mid X_{2}=1\right).$ 

By using Bayes theorem,
\begin{align*}
P\left(X_{1}=2\mid X_{2}=1\right) & =\dfrac{P\left(X_{1}=2\right)P\left(X_{2}=1\mid X_{1}=2\right)}{\sum^{2}_{j=0}P\left(X_{1}=j\right)P\left(X_{2}=1\mid X_{1}=j\right)}\\
 & =\dfrac{\dfrac{1}{4}\times\dfrac{1}{4}}{\dfrac{1}{2}\times\dfrac{1}{2}+\dfrac{1}{4}\times\dfrac{1}{4}}
\end{align*}
Hence,\boxeq{
\[
P\left(X_{1}=2\mid X_{2}=1\right)=\dfrac{1}{5}.
\]
}

\begin{remark}{}{}

To delve deeper into this exercise, one could attempt to construct
a probabilized space $\left(\Omega,\mathcal{A},P\right)$ using the
method discussed in Section \ref{sec:Evolutive-phenomena-modeling}.
This would require a complete definition of the conditional law system
of $X_{2}$ knowing $X_{1},$ which in the current setting has only
been partially used.

\end{remark}

\end{solution}

\begin{solution}{}{solexercise4.2}

We suppose to be constructed a probabilized space $\left(\Omega,\mathcal{A},P\right)$
that models successive uniform drawings.

\textbf{1. Probability that the drawn ball is blue}

Let $X_{1}$ denote the random variable representing the number of
blue balls transferred from urn $B$ to urn $A.$ Let $E$ be the
event ``a blue ball is drawn from urn $A.$''

The probabilities for each value of $X_{1}$ are
\[
\begin{array}{c}
P\left(X_{1}=2\right)=\dfrac{\left(\begin{array}{c}
8\\
0
\end{array}\right)\left(\begin{array}{c}
4\\
2
\end{array}\right)}{\left(\begin{array}{c}
12\\
2
\end{array}\right)}=\dfrac{1}{11},\\
P\left(X_{1}=1\right)=\dfrac{\left(\begin{array}{c}
8\\
1
\end{array}\right)\left(\begin{array}{c}
4\\
1
\end{array}\right)}{\left(\begin{array}{c}
12\\
2
\end{array}\right)}=\dfrac{16}{33},\\
P\left(X_{1}=0\right)=\dfrac{\left(\begin{array}{c}
8\\
2
\end{array}\right)\left(\begin{array}{c}
4\\
0
\end{array}\right)}{\left(\begin{array}{c}
12\\
2
\end{array}\right)}=\dfrac{14}{33}.
\end{array}
\]

Next, as we assumed uniform drawings from urn $A,$ the conditional
probabilities of drawing a blue ball are:
\begin{align*}
P\left(E\mid X_{1}=2\right) & =\dfrac{8}{13},\\
P\left(E\mid X_{1}=1\right) & =\dfrac{7}{13},\\
P\left(E\mid X_{1}=0\right) & =\dfrac{6}{13}.
\end{align*}

Hence, using the total probability formula,
\[
P\left(E\right)=\sum^{2}_{j=0}P\left(E\mid X_{1}=j\right)P\left(X_{1}=j\right).
\]

Thus, all computation done,\boxeq{
\[
P\left(E\right)=\dfrac{220}{429}.
\]
}

\textbf{2. Probability that at least one of the transferred ball to
be blue knowing that the drawn ball is blue.}

We want to compute $P\left(X_{1}\geqslant1\mid E\right).$

By Bayes theorem,
\begin{align*}
P\left(X_{1}=2\mid E\right) & =\dfrac{P\left(E\mid X_{1}=2\right)P\left(X_{1}=2\right)}{P\left(E\right)}=\dfrac{\dfrac{8}{13}\times\dfrac{1}{11}}{\dfrac{220}{429}}=\dfrac{6}{55},
\end{align*}

and
\begin{align*}
P\left(X_{1}=1\mid E\right) & =\dfrac{P\left(E\mid X_{1}=1\right)P\left(X_{1}=1\right)}{P\left(E\right)}=\dfrac{\dfrac{7}{13}\times\dfrac{16}{33}}{\dfrac{220}{429}}=\dfrac{28}{55}.
\end{align*}

Therefore, by $\sigma-$additivity of a conditional probability
\begin{align*}
P\left(X_{1}\geqslant1\mid E\right) & =P\left(X_{1}=1\mid E\right)+P\left(X_{1}=2\mid E\right).
\end{align*}
Hence,
\[
P\left(X_{1}\geqslant1\mid E\right)=\dfrac{34}{55}.
\]

\end{solution}

\begin{solution}{}{solexercise4.3}

We assume the probabilized space $\left(\Omega,\mathcal{A},P\right)$
to be constructed.

For every integer $n,$ we define
\begin{itemize}
\item $C_{n}:$ the event ``the household consists of $n$ children''.
\item $B_{n}:$ the event ``the household consists of $n$ boys''.
\item $G_{n}:$ the event ``the household consists of $n$ girls''
\end{itemize}
We are given that $P\left(C_{n}\right)=p_{n}.$

\textbf{1. Probability that a household with exactly two boys has
exactly two children.}

We want to compute the conditional probability $P\left(C_{2}\mid B_{2}\right).$ 

Using the Bayes theorem,
\begin{equation}
P\left(C_{2}\mid B_{2}\right)=\dfrac{P\left(C_{2}\right)P\left(B_{2}\mid C_{2}\right)}{P\left(B_{2}\right)},\label{eq:Bayes_formula_household}
\end{equation}
and,---Tr.N since there are $\binom{n}{2}$ ways to choose 2 boys
in a household of $n$ children and that there are $2^{n}$ ways of
ordering a $n-$uple of births likely gendered---
\[
P\left(B_{2}\mid C_{n}\right)=\binom{n}{2}\dfrac{1}{2^{n}}.
\]

From the law of total probability,
\[
P\left(B_{2}\right)=\sum^{+\infty}_{n=0}P\left(C_{n}\right)P\left(B_{2}\mid C_{n}\right).
\]

Clearly,
\[
P\left(B_{2}\mid C_{0}\right)=P\left(B_{2}\mid C_{1}\right)=0.
\]

So the sum reduces to
\begin{align*}
P\left(B_{2}\right) & =\sum^{+\infty}_{n=2}p_{n}\binom{n}{2}\dfrac{1}{2^{n}}=\dfrac{1-2a}{4^{2}}\sum^{+\infty}_{n=2}n\left(n-1\right)\dfrac{1}{4^{n-2}}.
\end{align*}

It remains to compute this sum. From the term-by-term differentiation
theorem for power series, it holds
\begin{align*}
\sum^{+\infty}_{n=2}n\left(n-1\right)x^{n-2} & =\dfrac{\text{d}^{2}}{\text{d}x^{2}}\left(\sum^{+\infty}_{n=0}x^{n}\right)=\dfrac{\text{d}^{2}}{\text{d}x^{2}}\left(\dfrac{1}{1-x}\right)=\dfrac{2}{\left(1-x\right)^{3}}.
\end{align*}

Hence, after simplifications,
\[
P\left(B_{2}\right)=\dfrac{8\left(1-2a\right)}{27}.
\]

Then, substituting in the equation $\refpar{eq:Bayes_formula_household},$
\[
P\left(C_{2}\mid B_{2}\right)=\dfrac{p_{2}\binom{2}{2}\dfrac{1}{2^{2}}}{\dfrac{8\left(1-2a\right)}{27}}.
\]

Thus, after simplification,\boxeq{
\[
P\left(C_{2}\mid B_{2}\right)=\dfrac{27}{64}\simeq0.421.
\]
}

\textbf{2. Probability that a household with two boys also has two
girls}

We now compute the conditional probability $P\left(G_{2}\mid B_{2}\right).$ 

By the total probability theorem,
\[
P\left(G_{2}\cap B_{2}\right)=\sum^{+\infty}_{n=0}P\left(C_{n}\right)P\left(G_{2}\cap B_{2}\mid C_{n}\right).
\]

Since all the conditional probabilities are equal to zero but the
one for $n=4,$ it holds
\[
P\left(G_{2}\cap B_{2}\right)=P\left(C_{4}\right)P\left(G_{2}\cap B_{2}\mid C_{4}\right).
\]

Thus, we obtain numerically
\begin{align*}
P\left(G_{2}\cap B_{2}\right) & =\left(1-2a\right)\times\dfrac{1}{2^{3}}\times\binom{4}{2}\times\dfrac{1}{2^{4}}=\dfrac{3\left(1-2a\right)}{64}.
\end{align*}

By the definition of the conditional probability,
\[
P\left(G_{2}\mid B_{2}\right)=\dfrac{P\left(G_{2}\cap B_{2}\right)}{P\left(G_{2}\right)}
\]
we obtain\boxeq{
\[
P\left(G_{2}\mid B_{2}\right)=\dfrac{81}{512}\simeq0.518.
\]
}

\end{solution}

\begin{solution}{}{solexercise4.4}

We suppose constructed a probabilized space $\left(\Omega,\mathcal{A},P\right)$
where all drawings are performed uniformly at random at each step.
For each $j\in\left\llbracket 1,n\right\rrbracket ,$ denote $E_{j}$
the event: ``at the j-th drawing, a bawl of each color is obtained.'' 

We have for every $j\in\left\llbracket 1,n-1\right\rrbracket $---Tr.N:
since there are $n-j$ red balls and $n-j$ green balls remaining
in the urn, i.e. $\left(n-j\right)^{2}$ ways of chosing a pair of
distinct color balls in the urn, while there are $\binom{2\left(n-j\right)}{2}$
possible pairs of balls that can be drawn.---
\[
P\left(E_{j+1}\left|\bigcap^{j}_{i=1}E_{i}\right.\right)=\dfrac{\left(n-j\right)^{2}}{\binom{2\left(n-j\right)}{2}}
\]

Additionally,
\[
P\left(E_{1}\right)=\dfrac{n^{2}}{\binom{2n}{2}}.
\]

Using the chain rule of conditional probabilities,
\[
P\left(\bigcap^{n}_{i=1}E_{i}\right)=P\left(E_{1}\right)\prod^{n-1}_{j=1}P\left(E_{j+1}\left|\bigcap^{j}_{i=1}E_{i}\right.\right).
\]

Thus, after simplifications,\boxeq{
\[
P\left(\bigcap^{n}_{i=1}E_{i}\right)=\dfrac{2^{n}\left(n!\right)^{2}}{\left(2n\right)!}.
\]
}

\end{solution}

\begin{solution}{}{solexercise4.5}

Observe that we can write
\[
T_{n+1}=\inf\left\{ j>T_{n}:\,\boldsymbol{1}_{A_{j}}=1\right\} .
\]

\textbf{1. (a) Law of the random variable $T_{k}$}

We already know that the law of $T_{1}$ is the geometric law on $\mathbb{N}^{\ast}$
with parameter $p.$

Now for $l\geqslant k\geqslant2,$
\[
\left(T_{k}=l\right)=\left(\sum^{l-1}_{j=1}\boldsymbol{1}_{A_{j}}=k-1\right)\cap\left(\boldsymbol{1}_{A_{l}}=1\right).
\]

The events $A_{n}$ being independent and with same probability $p,$
the law of the random variable $\sum^{l-1}_{j=1}\boldsymbol{1}_{A_{j}}$
is the binomial law $\mathcal{B}\left(l-1,p\right)$ and the events
$\left(\sum^{l-1}_{j=1}\boldsymbol{1}_{A_{j}}=k-1\right)$ and $\left(\boldsymbol{1}_{A_{k}}=1\right)$
are independent.

It holds\boxeq{
\[
P\left(T_{k}=l\right)=\binom{l-1}{k-1}p^{k}q^{l-k}.
\]
}

The law of the random variable $T_{k}$ is the \textbf{negative binomial
law\index{negative binomial law}} with parameters $k$ and $p,$
denoted $\mathcal{B}^{-}\left(k,p\right).$ 

\textbf{(b) Conditional probability $P\left(T_{k+1}-T_{k}=n_{k+1}-n_{k}\,\left|\,T_{1}=n_{1},\dots,T_{k}=n_{k}\right.\right).$}

We write
\begin{multline*}
\bigcap^{k}_{j=1}\left(T_{j}=n_{j}\right)\cap\left(T_{k+1}-T_{k}=n_{k+1}-n_{k}\right)\\
=\left[\bigcap^{k}_{j=1}\left(T_{j}=n_{j}\right)\right]\cap\left(\sum^{n_{k+1}-1}_{j=n_{k}+1}\boldsymbol{1}_{A_{j}}=0\right)\cap\left(\boldsymbol{1}_{A_{n_{k+1}}}=1\right).
\end{multline*}

By independence of the events within the second member of this last
equality,
\begin{multline*}
P\left(T_{k+1}-T_{k}=n_{k+1}-n_{k}\,\left|\,T_{1}=n_{1},\dots,T_{k}=n_{k}\right.\right)\\
=P\left[\left(\sum^{n_{k+1}-1}_{j=n_{k}+1}\boldsymbol{1}_{A_{j}}=0\right)\cap\left(\boldsymbol{1}_{A_{n_{k+1}}}=1\right)\right],
\end{multline*}

Thus\boxeq{
\begin{equation}
P\left(T_{k+1}-T_{k}=n_{k+1}-n_{k}\,\left|\,T_{1}=n_{1},\dots,T_{k}=n_{k}\right.\right)=pq^{n_{k+1}-n_{k}-1}.\label{eq:ex4.5rel4.9}
\end{equation}
}

\textbf{(c) Proof of $P\left(T_{k+1}=n_{k+1}\,\left|\,T_{1}=n_{1},\dots,T_{k}=n_{k}\right.\right)=P\left(T_{k+1}=n_{k+1}\,\left|\,T_{k}=n_{k}\right.\right)$}

A similar computation proves that

\[
P\left(T_{k+1}-T_{k}=n_{k+1}-n_{k}\,\left|\,T_{k}=n_{k}\right.\right)=pq^{n_{k+1}-n_{k}-1}.
\]

Since
\begin{multline*}
P\left(T_{k+1}=n_{k+1}\,\left|\,T_{1}=n_{1},\dots,T_{k}=n_{k}\right.\right)\\
=P\left(T_{k+1}-T_{k}=n_{k+1}-n_{k}\,\left|\,T_{1}=n_{1},\dots,T_{k}=n_{k}\right.\right)
\end{multline*}
and
\[
P\left(T_{k+1}=n_{k+1}\,\left|\,T_{k}=n_{k}\right.\right)=P\left(T_{k+1}-T_{k}=n_{k+1}-n_{k}\,\left|\,T_{k}=n_{k}\right.\right),
\]
we can deduce that\boxeq{
\[
P\left(T_{k+1}=n_{k+1}\,\left|\,T_{1}=n_{1},\dots,T_{k}=n_{k}\right.\right)=P\left(T_{k+1}=n_{k+1}\,\left|\,T_{k}=n_{k}\right.\right).
\]
}

\begin{remark}{}{}

This last relation expresses the Markovian character of the sequence
of random variables $T_{k},$ i.e. the probabilistic evolution of
the sequence at the time $k+1$ depends only on the current state
$T_{k}$ at the instant $k,$ and not on the full past, and the states
taken before. The conditional probability is given by
\[
P\left(T_{k+1}=n_{k+1}\,\left|\,T_{k}=n_{k}\right.\right)=pq^{n_{k+1}-n_{k}-1}.
\]

\end{remark}

\textbf{(d) Value of the joint probability $P\left(T_{1}=n_{1},\dots,T_{k}=n_{k}\right)$}

For $n_{1}<n_{2}<\dots<n_{k},$ we deduce from the chain rule of conditional
probabilities
\begin{multline*}
P\left(T_{1}=n_{1},\dots,T_{k}=n_{k}\right)=P\left(T_{1}=n_{1}\right)P\left(T_{2}=n_{2}\,\left|\,T_{1}=n_{1}\right.\right)\\
\dots P\left[\left(T_{k}=n_{k}\right)\,\left|\,\bigcap^{k-1}_{j=1}\left(T_{j}=n_{j}\right)\right.\right].
\end{multline*}

Hence, using the previous question
\[
P\left(T_{1}=n_{1},\dots,T_{k}=n_{k}\right)=pq^{n_{1}-1}pq^{n_{2}-n_{1}-1}\dots pq^{n_{k}-n_{k-1}-1},
\]
which gives, after simplifications\boxeq{
\[
P\left(T_{1}=n_{1},\dots,T_{k}=n_{k}\right)=p^{k}q^{n_{k}-k}.
\]
}

\textbf{2. Independence of the increments $T_{k+1}-T_{k}$}

From the equation $\refpar{eq:ex4.5rel4.9},$ it holds that, for every
$m\geqslant1,$
\[
P\left(T_{k+1}-T_{k}=m\,\left|\,T_{1}=n_{1},\dots,T_{k}=n_{k}\right.\right)=pq^{m-1}.
\]

Thus, for every $k\geqslant1,$ the random variables $T_{k+1}-T_{k}$
and $\left(T_{1},T_{2},\dots,T_{k}\right)$ are independent. The random
variables $T_{1},T_{2}-T_{1},\dots,T_{k}-T_{k-1}$ are also independent.
Those variables follow the geometric law on $\mathbb{N}^{\ast},$
denoted $\mathcal{G}_{\mathbb{N}^{\ast}}\left(p\right).$ 

Since
\[
T_{k}=T_{1}+\sum^{k-1}_{j=1}\left(T_{j+1}-T_{j}\right)
\]
it follows from the definition of the convolution of laws---law of
sum of independent random variables---that the $k$-th convolution
of the geometric law on $\mathbb{N}^{\ast}$is the negative binomial
law $\mathcal{B}^{-}\left(k,p\right).$

\textbf{3. Computation of $P\left(T_{1}=3,T_{5}=9,T_{7}=17\right).$}

By independence of the random variables $T_{1},$ $T_{5}-T_{1},$
$T_{7}-T_{5},$ it follows that
\begin{align*}
P\left(T_{1}=3,T_{5}=9,T_{7}=17\right) & =P\left(T_{1}=3,T_{5}-T_{1}=6,T_{7}-T_{5}=8\right)\\
 & =P\left(T_{1}=3\right)P\left(T_{5}-T_{1}=6\right)P\left(T_{7}-T_{5}=8\right).
\end{align*}

As we have
\[
T_{5}-T_{1}=\sum^{5}_{j=2}\left(T_{j}-T_{j-1}\right),
\]
the law of $T_{5}-T_{1}$ is the negative binomial law $\mathcal{B}^{-}\left(4,p\right).$

Similarly, since
\[
T_{7}-T_{5}=\sum^{7}_{j=6}\left(T_{j}-T_{j-1}\right),
\]
the law of $T_{7}-T_{5}$ is the negative binomial law $\mathcal{B}^{-}\left(2,p\right).$

It follows
\[
P\left(T_{1}=3,T_{5}=9,T_{7}=17\right)=pq^{3-1}\left(\binom{6-1}{4-1}p^{4}q^{2}\right)\left(\binom{8-1}{2-1}p^{2}q^{6}\right),
\]
which is equivalent to\boxeq{
\[
P\left(T_{1}=3,T_{5}=9,T_{7}=17\right)=70p^{7}q^{10}.
\]
}

Approximated values of those probabilities for the $p$ values given
are gathered in Table \ref{Tab Approximated values}.

\end{solution}

\begin{table}
\begin{center}%
\begin{tabular}{|c|c|c|c|}
\hline 
$p$ &
0.3 &
0.5 &
0.8\tabularnewline
\hline 
$P\left(T_{1}=3,T_{5}=9,T_{7}=17\right)$ &
$4.32\times10^{-4}$ &
$5.34\times10^{-4}$ &
$1.5\times10^{-6}$\tabularnewline
\hline 
\end{tabular}\end{center}

\caption{Approximated values of $P\left(T_{1}=3,T_{5}=9,T_{7}=17\right)$ for
different values of $p.$}
\label{Tab Approximated values}
\end{table}

\begin{solution}{}{solexercice4.6}

We model this experiment in two steps, in the following manner
\begin{itemize}
\item Step 1 (card selection): to model the first step, i.e. by choosing
the card, we choose the probabilized space $\left(\Omega_{1},\mathcal{P}\left(\Omega_{1}\right),P_{1}\right)$
where
\[
\Omega_{1}=\left\{ RR,RG,GG\right\} ,
\]
and $P_{1}$ the uniform probability.
\item Step 2 (face revelation): then, we define the set $\Omega_{2}=\left\{ R,G\right\} $
giving the outcome corresponding to the colour of the visible face.
\end{itemize}
We define the following probability germs
\begin{itemize}
\item The germ $p_{1}$ which generates the uniform probability on $\Omega_{1}$
and which is defined for every $\omega_{1}\in\Omega_{1}$ by $p_{1}\left(\omega_{1}\right)=\dfrac{1}{3}.$
\item The family of germs on $\Omega_{2}$ defined by
\[
\begin{array}{ccccc}
p^{RR}_{2}\left(R\right)=1 &  & \text{and} &  & p^{RR}_{2}\left(G\right)=0,\\
p^{GG}_{2}\left(R\right)=0 &  & \text{and} &  & p^{GG}_{2}\left(G\right)=1,\\
p^{RG}_{2}\left(R\right)=\dfrac{1}{2} &  & \text{and} &  & p^{RG}_{2}\left(G\right)=\dfrac{1}{2}.
\end{array}
\]
\end{itemize}
This family of probability germs defines a unique probability $P$
on the probabilizable space $\left(\Omega_{1}\times\Omega_{2},\mathcal{P}\left(\Omega_{1}\times\Omega_{2}\right)\right)$
by the relation
\[
\forall\left(\omega_{1},\omega_{2}\right)\in\Omega_{1}\times\Omega_{2},\,\,\,\,\,P\left(\left\{ \omega_{1},\omega_{2}\right\} \right)=p_{1}\left(\omega_{1}\right)p^{\omega_{1}}_{2}\left(\omega_{2}\right).
\]

Let $H_{R}$ designates the event ``the hidden face is red'' and
$H_{G}$ the event ``the hidden face is green''; denote $X_{1}$
and $X_{2}$ the canonical projections from $\Omega_{1}\times\Omega_{2}$
onto $\Omega_{1}$ and $\Omega_{2}$ respectively. 

We first compute the conditional probability $P\left(H_{R}\left|X_{2}=R\right.\right).$

The event $H_{R}$ can be partitioned in the following manner
\[
H_{R}=\left(X_{1}=RR\right)\uplus\left[\left(X_{1}=RG\right)\cap\left(X_{2}=G\right)\right],
\]
 which gives, by $\sigma-$additivity of a conditional probability
\[
P\left(H_{R}\left|X_{2}=R\right.\right)=P\left(\left(X_{1}=RR\right)\left|X_{2}=R\right.\right)+P\left(\left(X_{1}=RG\right)\cap\left(X_{2}=G\right)\left|X_{2}=R\right.\right).
\]

But,
\[
P\left(\left(X_{1}=RG\right)\cap\left(X_{2}=G\right)\left|X_{2}=R\right.\right)=0
\]
 and using the Bayes formula,
\[
P\left(\left(X_{1}=RR\right)\left|X_{2}=R\right.\right)=P\left(X_{2}=R\left|\left(X_{1}=RR\right)\right.\right)\times\dfrac{P\left(X_{1}=RR\right)}{P\left(X_{2}=R\right)}.
\]

If follows that
\[
P\left(H_{R}\left|X_{2}=R\right.\right)=\dfrac{1}{3P\left(X_{2}=R\right)}.
\]

Additionally, from the total probability theorem,
\begin{align*}
P\left(X_{2}=R\right) & =P\left(X_{2}=R\left|\left(X_{1}=RR\right)\right.\right)P\left(X_{1}=RR\right)\\
 & \qquad\qquad+P\left(X_{2}=R\left|\left(X_{1}=GG\right)\right.\right)P\left(X_{1}=GG\right)\\
 & \qquad\qquad+P\left(X_{2}=R\left|\left(X_{1}=RG\right)\right.\right)P\left(X_{1}=RG\right).
\end{align*}

Thus,
\[
P\left(X_{2}=R\right)=\dfrac{1}{3}+\left(\dfrac{1}{2}\times\dfrac{1}{3}\right)=\dfrac{1}{2}.
\]

We then obtain
\[
P\left(H_{R}\left|X_{2}=R\right.\right)=\dfrac{2}{3}
\]
and consequently
\[
P\left(H_{G}\left|X_{2}=R\right.\right)=\dfrac{1}{3}.
\]

Exchanging the role of $R$ and $G,$ we obtain similarly\boxeq{
\[
\begin{array}{ccccc}
P\left(H_{G}\left|X_{2}=G\right.\right)=\dfrac{2}{3} & \,\,\,\, & \text{and} & \,\,\,\,\, & P\left(H_{R}\left|X_{2}=G\right.\right)=\dfrac{1}{3}.\end{array}
\]
}

The best strategy for the spectator is thus to gamble on the color
of the exposed face.

\end{solution}

\begin{solution}{}{solexercise4.7}

1. We partition the studied event.

\textbf{(a) Computation of $P\left(U<V\right)$}

We use the disjoint union
\[
\left(U<V\right)=\biguplus_{\substack{1\leqslant u<v\\
\left(u,v\right)\in\mathbb{N}^{\ast2}
}
}\left[\left(U=u\right)\cap\left(V=v\right)\right].
\]

Hence, by $\sigma-$additivity of $P$ and independence of the random
variables $U$ and $V,$
\[
P\left(U<V\right)=\sum_{\substack{1\leqslant u<v\\
\left(u,v\right)\in\mathbb{N}^{\ast2}
}
}\left(1-q\right)q^{u-1}\left(1-r\right)r^{v-1}.
\]

Since the family under the sum sign is nonnegative, the Fubini property
induces
\begin{align*}
P\left(U<V\right) & =\left(1-q\right)\left(1-r\right)\sum^{+\infty}_{u=1}q^{u-1}\sum^{+\infty}_{v=u+1}r^{v-1}\\
 & =\left(1-q\right)\left(1-r\right)\sum^{+\infty}_{u=1}\dfrac{r^{u}}{1-r}q^{u-1}.
\end{align*}
Hence,\boxeq{
\[
P\left(U<V\right)=\dfrac{\left(1-q\right)r}{1-rq}.
\]
}

\textbf{(b) Special cases $q=r$ and} $q=r=\dfrac{1}{2}.$
\begin{itemize}
\item If $q=r,$ then $P\left(U<V\right)=\dfrac{q}{1+q}.$
\item If $q=r=\dfrac{1}{2},$ then $P\left(U<V\right)=\dfrac{1}{3}.$
\end{itemize}
\textbf{2. Conditional law under $\left(U<V\right)$}

The event $\left(U<V\right)$ is of nonzero probability.

\textbf{(a) Computation of $P\left(U=k\left|U<V\right.\right).$}

Since
\[
\left(U=k\right)\cap\left(U<V\right)=\left(U=k\right)\cap\left(k<V\right),
\]

by the Bayes formula, it follows that
\[
P\left(U=k\,\left|\,U<V\right.\right)=\dfrac{P\left(\left(U=k\right)\cap\left(k<V\right)\right)}{P\left(U<V\right)}.
\]

By independence of $U$ and $V,$ it comes
\[
P\left(U=k\,\left|\,U<V\right.\right)=\dfrac{P\left(U=k\right)P\left(k<V\right)}{P\left(U<V\right)}.
\]

Since
\[
\left(k<V\right)=\biguplus_{\substack{k<v}
}\left(V=v\right)
\]
we obtain
\begin{equation}
P\left(k<V\right)=\sum^{+\infty}_{v=k+1}\left(1-r\right)r^{v-1}=r^{k}.\label{eq:P(k<V)}
\end{equation}

Thus, using the first question
\begin{align*}
P\left(U=k\,\left|\,U<V\right.\right) & =\dfrac{\left(1-q\right)q^{k-1}r^{k}}{\dfrac{\left(1-q\right)r}{1-rq}}\\
 & =\left(1-rq\right)\left(rq\right)^{k-1}.
\end{align*}

It shows that the conditional law generated by the germ $P\left(U=\cdotp\,\left|\,U<V\right.\right)$
is the geometric law on $\mathbb{N}^{\ast}$ of parameter $1-rq,$
denoted $\mathcal{G}_{\mathbb{N}^{\ast}}\left(1-rq\right).$

\textbf{(b) Computation of $P\left(U>k\,\left|\,U<V\right.\right).$}

From the equation $\refpar{eq:P(k<V)},$ we have, for every $k\in\mathbb{N},$
\begin{align*}
P\left(U>k\,\left|\,U<V\right.\right) & =\sum^{+\infty}_{j=k+1}P\left(U=j\,\mid\,U<V\right)\\
 & =\sum^{+\infty}_{j=k+1}\left(1-rq\right)\left(rq\right)^{j-1}\\
 & =\left(1-rq\right)\left(rq\right)^{k}\dfrac{1}{1-rq}
\end{align*}
Hence,\boxeq{
\[
P\left(U>k\,\mid\,U<V\right)=\left(rq\right)^{k}.
\]
}

\textbf{3. Properties of $J$ and $X$ under conditioning on $H$}

This question states some properties of random variables $J$ and
$X$ under the conditioning by the event $H.$

\textbf{(a) Computation of $P\left(H\cap\left(J=1\right)\cap\left(X>k\right)\right).$}

Since $\left(J=1\right)\subset H$ and that, on the set $\left(J=1\right),$
$X=X_{1},$ and
\[
P\left(H\cap\left(J=1\right)\cap\left(X>k\right)\right)=P\left(\left(J=1\right)\cap\left(X_{1}>k\right)\right).
\]
 Hence,
\begin{equation}
P\left(H\cap\left(J=1\right)\cap\left(X>k\right)\right)=P\left(X_{1}>k\,\left|\,X_{1}<X_{2}\right.\right)P\left(X_{1}<X_{2}\right).\label{eq:P_HcapJ_is_1capX>k}
\end{equation}

Using the results of the first two questions in this particular case,
it holds
\[
\begin{array}{ccccc}
P\left(X_{1}>k\,\left|\,X_{1}<X_{2}\right.\right)=\left(1-q\right)^{2k} & \,\,\,\, & \text{and} & \,\,\,\, & P\left(X_{1}<X_{2}\right)=\dfrac{1-q}{2-q}\end{array}
\]
and thus\boxeq{
\[
P\left(H\cap\left(J=1\right)\cap\left(X>k\right)\right)=\dfrac{\left(1-q\right)^{2k+1}}{2-q}.
\]
}

\textbf{(b) Computation of $P^{H}\left(\left(J=1\right)\cap\left(X=k\right)\right)$}

The random variables $X_{1}$ and $X_{2}$ are independent and of
same law; we then write
\begin{align*}
P\left(A_{1}\right) & =\sum_{x_{1}<x_{2}}P_{X_{1}}\left(\left\{ x_{1}\right\} \right)P_{X_{2}}\left(\left\{ x_{2}\right\} \right)\\
 & =\sum_{x_{1}<x_{2}}P_{X_{2}}\left(\left\{ x_{1}\right\} \right)P_{X_{1}}\left(\left\{ x_{2}\right\} \right),
\end{align*}
which shows the equality
\[
P\left(A_{1}\right)=P\left(A_{2}\right).
\]

Since $H=A_{1}\uplus A_{2},$ we deduce that we have
\[
P\left(H\right)=2P\left(A_{1}\right).
\]

Thus, using the equation $\refpar{eq:P_HcapJ_is_1capX>k}$ and the
Bayes formula,
\begin{align*}
P^{H}\left(\left(J=1\right)\cap\left(X>k\right)\right) & =\dfrac{P\left(H\cap\left(J=1\right)\cap\left(X>k\right)\right)}{P\left(H\right)}\\
 & =\dfrac{P\left(X_{1}>k\,\left|\,X_{1}<X_{2}\right.\right)P\left(X_{1}<X_{2}\right)}{P\left(H\right)}\\
 & =\dfrac{1}{2}P\left(X_{1}>k\,\left|\,X_{1}<X_{2}\right.\right).
\end{align*}
Thus,\boxeq{
\[
P^{H}\left(\left(J=1\right)\cap\left(X>k\right)\right)=\dfrac{1}{2}\left(1-q\right)^{2k}.
\]
}

\textbf{(c) Computation of $P^{H}\left(J=1\right)$}

The definition of a conditional probability gives
\[
P^{H}\left(J=1\right)=\dfrac{P\left(\left(J=1\right)\cap H\right)}{P\left(H\right)}=\dfrac{P\left(A_{1}\right)}{P\left(H\right)}.
\]

Thus\boxeq{
\begin{equation}
P^{H}\left(J=1\right)=\dfrac{1}{2}.\label{eq:P_H_Jis1}
\end{equation}
}

4. We study the laws of the random variables $X$ and $J$ under the
probability $P^{H}$ and conclude to the independence under $P^{H}$
of these random variables.

\textbf{(a) Computation of $P^{H}\left(X>k\right)$}

We have the following set equalities
\[
H\cap\left(J=1\right)\cap\left(X>k\right)=\left(X_{1}>k\right)\cap\left(X_{1}<X_{2}\right)
\]
and
\[
H\cap\left(J=2\right)\cap\left(X>k\right)=\left(X_{2}>k\right)\cap\left(X_{2}<X_{1}\right).
\]

The random variables $X_{1}$ and $X_{2}$ being independent and of
same law, an analoguous computation to the one in the Question 3.b.
ensures that all sets are of same probability.

Additionally, we note that
\[
H\cap\left(J=0\right)=\emptyset.
\]

We thus have
\begin{align*}
P^{H}\left(X>k\right) & =P^{H}\left(\left(J=1\right)\cap\left(X>k\right)\right)+P^{H}\left(\left(J=2\right)\cap\left(X>k\right)\right)\\
 & =2P^{H}\left(\left(J=1\right)\cap\left(X>k\right)\right)
\end{align*}
Hence,\boxeq{
\[
P^{H}\left(X>k\right)=\left(1-q\right)^{2k}.
\]
}

\textbf{(b) X and J are independent under $P^{H}$}

Since $P^{H}\left(J=0\right)=0,$ it follows from the equation $\refpar{eq:P_H_Jis1}$
that
\[
P^{H}\left(J=2\right)=\dfrac{1}{2}.
\]

By gathering the results from the previous questions, we have then
shown for $j=0,1,2$ and for every $k\in\mathbb{N},$ that we have\boxeq{
\[
P^{H}\left(\left(J=j\right)\cap\left(X>k\right)\right)=P^{H}\left(J=j\right)P^{H}\left(X>k\right),
\]
}and thus, that \textbf{under the probability $P^{H},$ the random
variables $X$ and $J$ are independent.}

\end{solution}

\chapter{Moments of a Discrete Random Variable}\label{chap:Moments-of-a}

\begin{objective}{}{}

Chapter \ref{chap:Moments-of-a} introduces the concept of moments
for discrete random variables.
\begin{itemize}
\item Section \ref{sec:Mean-or-Mathematical} begins by defining the mean,
or mathematical expectation, of a random variable. It outlines basic
properties of expectation, then extends the concept to function of
a random variable. The section ends with some expectations for classical
discrete laws.
\item Section \ref{sec:Higher-Order-Moments} covers the definition of higher
order moments. It introduces Hölder and Minkowski inequalities, and
apply them to higher moment expectation of two random variables. It
places particular emphasis on the properties of second-order moments,
leading to the definition of variance and standard deviation. Leibniz
formula connecting variance and expectation is then presented. The
section continues with the definition and properties of covariance
of two random variables. Some variance of classical discrete laws
are then given, before moving to the Markov and Chebyshev inequalities.
Last, the correlation coefficient is defined and the linear regression
problem is exposed.
\item Section \ref{sec:Generative-Functions} introduces the concept of
generating function of a random variable which characterizes the law
of this random variable. Generating functions for classical laws are
then presented, along with their relationship to moments. The Chapter
concludes by addressing the sum of a random number of random variables.
\end{itemize}
\end{objective}

\section*{Introduction}

Moments of a random variable are numerical parameters that provide
insight into the law of this random variable---though, they do not
completely determine it in general. The most commonly used are the
mean---or mathematical expectation---and the variance.

In this Chapter, we focus only on the moments of discrete, real-valued
random variables. The summation process involved is the summation
over families of real numbers. In the next Chapter, we will define
the same notions for random variables with density; formally, it will
suffice to change the process of summation, the sum $\Sigma$ being
replaced by the integral $\int.$ A general treatment that unifies
both cases relies on measure theory and integration, which is the
subject of Part \ref{part:Deepening-Probability-Theory} of this book.

In this chapter, unless otherwise stated, all random variables are
defined on a probabilized space $\left(\Omega,\mathcal{A},P\right)$
and are discrete. Moreover they take real values, either finite or
infinite. In practice, the variables studied are finite with probability
1---that is as we will see in Chapter \ref{chap:Approximation-of-laws.}
almost surely finite.

\section{Mean or Mathematical Expectation}\label{sec:Mean-or-Mathematical}

\subsection{Definition}

\begin{definition}{Mean. Mathematical Expectation}{}

Let $X$ be a discrete random variable taking real values. 

If the family of real numbers $\left(xP\left(X=x\right)\right)_{x\in\mathbb{R}}$
is summable, i.e. if the sum $\sum_{x\in\mathbb{R}}\left|x\right|P\left(X=x\right)$
is finite, then we consider the number
\[
\sum_{x\in\mathbb{R}}xP\left(X=x\right)
\]
denoted, following the typographic needs, $\mathbb{E}\left(X\right)$
or, occasionally, $\mathbb{E}X$---Tr.N.: For clarity, we prefer
the first notation, contrary to the original French version of this
book which privilieges the second notation. $\mathbb{E}\left(X\right)$
is called \index{mean}\textbf{mean}, or \index{mathematical expectation}\textbf{mathematical
expectation\index{expectation}} of the random variable $X.$ 

We say in the case where the sum is finite that the random variable
$X$ admits a \textbf{mean\index{mean}} or a \textbf{mathematical
expectation}\index{mathematical expectation}. 

In summary, when the sum is finite\boxeq{
\begin{equation}
\mathbb{E}\left(X\right)=\sum_{x\in\mathbb{R}}xP\left(X=x\right).\label{eq:expectation-def}
\end{equation}
}

\end{definition}

\begin{remarks}{}{}

\textbf{1. Law-Based Definition}

The notion of expectation is defined here without refering to a probabilized
space $\left(\Omega,\mathcal{A},P\right)$ and appears as a concept
linked directly to the law of the random variable.

\textbf{2. Restriction of the Expectation to the Nonzero Probability
Set}

Since $X$ is supposed discrete, the set of its values $X\left(\Omega\right)$
is a countable set; the set $\text{val}\left(X\right)$ of its taken
values with a nonzero probability is thus also countable.

Of course, in the definiton of $\mathbb{E}\left(X\right)$ we can
restrict the summation to the set $\text{val}\left(X\right)$ constituted
of the real numbers $x$ such that $P\left(X=x\right)>0.$ We can
then write\boxeq{
\[
\mathbb{E}\left(X\right)=\sum_{x\in\text{val}\left(X\right)}xP\left(X=x\right).
\]
}

This is this formula which serves for the computation of the expectation
from the law of $X.$ The explicit computation is done by noting that,
depending on the finiteness of the set of values of $X,$ the definition
of $\mathbb{E}\left(X\right)$ takes one of the following form:
\begin{itemize}
\item If $\text{val}\left(X\right)$ is finite and can be enumerated as
a set of the form $\left\{ x_{1},x_{2},\dots,x_{n}\right\} ,$ then
$X$ has an expectation, and we have\boxeq{
\begin{equation}
\mathbb{E}\left(X\right)=\sum^{n}_{i=1}x_{i}P\left(X=x_{i}\right).\label{eq:mean finite case}
\end{equation}
}
\item If $\text{val}\left(X\right)$ is infinite and can be enumerated as
$\left\{ x_{n}:\,n\in\mathbb{N}\right\} ,$ then $X$ has an expectation
if and only if $\sum^{+\infty}_{i=0}\left|x_{i}\right|P\left(X=x_{i}\right)<+\infty$---that
is, the series is absolutely convergent---and in that case\boxeq{
\begin{equation}
\mathbb{E}\left(X\right)=\sum^{+\infty}_{i=0}x_{i}P\left(X=x_{i}\right).\label{eq:mean_infinite_case}
\end{equation}
}
\end{itemize}
The existence of the expectation and of its value do not depend on
the order of the chosen enumeration of the values of $X$---see Proposition
$\ref{pr:fam_sum_enum}$ and the remarks that follow.

\textbf{3. Condition for a Random Variable to Have a expectation}

Last, we can observe that for a random variable $X\geqslant0,$ the
sum $\sum_{x\in\text{val}\left(X\right)}xP\left(X=x\right)$ has always
a meaning, even if it takes the value $+\infty.$ We write in this
case $\mathbb{E}\left(X\right)=+\infty.$ With this convention, the
condition for $X$ to have an expectation is written $\mathbb{E}\left(\left|X\right|\right)<+\infty.$

In practical terms, the expectation $\mathbb{E}\left(X\right)$ is
simply the \textbf{weighted mean\index{weighted mean}\mindex{mean ! weighted}}
of the values $x$ taken by the random variable $X,$ with each value
weighted by the probability that $X$ takes that value. 

\end{remarks}

\subsection{First Properties of the Expectation}

The following properties of the expectation are often used.

\begin{proposition}{Common Expectation Properties}{}

\textbf{1. Expectation of a Constant}

Let $a$ be a real number. 

If the random variable is such that $P\left(X=a\right)=1$---i.e.
$X$ is a constant with probability 1---then its expectation equals
to $a.$ 

In particular, we have\boxeq{
\[
\mathbb{E}\left(a\right)=a.
\]
}

\textbf{2. Bounded Discrete Random Variable}

Any bounded---in the sense there exists a nonnegative real number
$M$ such that for every $\omega\in\Omega,$ $\left|X\left(\omega\right)\right|\leqslant M$---discrete
random variable $X$ admits an expectation.

\textbf{3. Non-Negativity}

If $X\geqslant0,$ then $\mathbb{E}\left(X\right)\geqslant0.$ 

If $X\geqslant0$ and $\mathbb{E}\left(X\right)=0,$ then $X=0$ with
probability 1.

\textbf{4. Monotonicity and Absolute Value Inequality}

If $X$ and $Y$ admit an expectation and verify $X\leqslant Y,$
then
\[
\mathbb{E}\left(X\right)\leqslant\mathbb{E}\left(Y\right).
\]

Also, we have the inequality\boxeq{
\[
\left|\mathbb{E}\left(X\right)\right|\leqslant\mathbb{E}\left(\left|X\right|\right).
\]
}

\end{proposition}

\begin{proof}

These properties result from similar properties for summable families.

Tr/N: For instance, we give the proof for 1. and 3.

1. $P\left(X=a\right)=1$ implies that $P\left(X\neq a\right)=0.$ 

Thus,
\[
\sum\limits_{x\in\mathbb{R}}xP\left(X=x\right)=aP\left(X=a\right)=a
\]
 and the result.

3. We have $\sum_{x>0}xP\left(X=x\right)=0,$ thus $xP\left(X=x\right)=0$
for every real number $x>0.$

Since $x>0,$ this implies $P\left(X=x\right)=0$ for every real number
$x>0.$

As $X$ is discrete,
\[
P\left(X>0\right)=\sum_{x>0}xP\left(X=x\right)=0,
\]
which implies $P\left(X=0\right)=1.$

\end{proof}

\begin{example}{Expectation of a Random Variable Following a Uniform Law}{}

If $X$ is a discrete random variable taking real values, of law the
uniform probability\footnotemark over the finite set $\left\{ x_{1},\dots,x_{n}\right\} ,$
then the expectation $\mathbb{E}\left(X\right)$ of $X$ is simply
the arithmetic mean of the values $x$ taken by $X,$\boxeq{
\[
\mathbb{E}\left(X\right)=\dfrac{1}{n}\sum^{n}_{i=1}x_{i}.
\]
}

\end{example}

\footnotetext{We say that the random variable $X$ is \mindex{random variable!uniformly distributed}\index{uniformly distributed}\textbf{uniformly
distributed} over the set $\left\{ x_{1},\dots,x_{n}\right\} .$}

\begin{example}{Expectation of an Indicator Function. Expectation of a Bernoulli Law}{expectation_bernoulli_law}

Let $A\in\mathcal{A}$ be an event. Its indicator function $\boldsymbol{1}_{A}$
is a discrete real-valued random variable and whose law is determined
by the relations
\[
P\left(\boldsymbol{1}_{A}=1\right)=P\left(A\right),\,\,\,\,P\left(\boldsymbol{1}_{A}=0\right)=1-P\left(A\right).
\]

Its expectation is then
\[
\mathbb{E}\left(\boldsymbol{1}_{A}\right)=P\left(A\right).
\]

Equivalently, if $X$ is a Bernouilli random variable with parameter
$p,$ then its expectation is
\[
\mathbb{E}\left(X\right)=p.
\]

\end{example}

\subsection{Expectation of a Function of a Random Variable}

\begin{theorem}{Transfer Theorem}{transfer_theorem}

Let $X$ be a discrete random variable taking values in the set $E.$
Let $f$ be a function from $E$ to $\mathbb{R}.$ 

The composition $Y=f\circ X$---very often denoted $f\left(X\right)$---is
a real-valued discrete random variable.

$Y$ admits an expectation if and only if the sum
\[
\sum_{x\in E}\left|f\left(x\right)\right|P\left(X=x\right)
\]
is finite. 

In this case, the expectation of $f\left(X\right)$ exists and is
given by\boxeq{
\begin{equation}
\mathbb{E}\left(f\left(X\right)\right)=\sum_{x\in E}f\left(x\right)P\left(X=x\right)\label{eq:transfer_theorem}
\end{equation}
}

\end{theorem}

Thanks to this Theorem, we can compute $\mathbb{E}\left(f\left(X\right)\right),$
by performing the summation on the set of values taken by $X,$ rather
than over the set of values taken by $f\left(X\right).$ This justifies
the name transfer theorem. Tr/N: This theorem is also known as the
LOTUS theorem, that is the ``law of the unconscious statistician'',
since it is often mistakenly taken as the definition of the expectation
of $f\left(X\right)$ while its mean is given by $\refpar{eq:expectation-def}$
with $Y=f\left(X\right).$

In practice, it is then enough to know the law of $X$ to ensure the
existence of the expectation of $f\left(X\right)$ and compute explicitly
its value. 

\begin{proof}{}{}

Let $y\in\mathbb{R}.$ The event $\left\{ f\left(X\right)=y\right\} $
can be expressed as the disjoint union of the events $\left\{ X=x\right\} $
where $x$ describes $f^{-1}\left(y\right).$ Since $X$ is a discrete
random variable, the set of values taken by $X$ is countable. 

Hence, the set of inverse images $x$ of $y$ for which the event
$\left\{ X=x\right\} $ is not empty constitutes a countable subset
of $f^{-1}\left(y\right),$ which we denote $G_{y}.$ Invoking the
$\sigma-$additivity of $P,$ we thus have
\[
P\left(f\left(X\right)=y\right)=\sum_{x\in G_{y}}P\left(X=x\right).
\]

But since $P\left(X=x\right)=0$ for $x\in f^{-1}\left(y\right)\backslash G_{y},$
we can write
\[
P\left(f\left(X\right)=y\right)=\sum_{x\in f^{-1}\left(y\right)}P\left(X=x\right).
\]

We can then write
\begin{align*}
\left|y\right|P\left(f\left(X\right)=y\right) & =\left|y\right|\sum_{x\in f^{-1}\left(y\right)}P\left(X=x\right)\\
 & =\sum_{y\in\mathbb{R}}\left|y\right|P\left(f\left(X\right)=y\right).
\end{align*}

The finiteness of the sum
\[
\sum_{x\in E}\left|f\left(x\right)\right|P\left(X=x\right)
\]
is equivalent to the one of the sum
\[
\sum_{y\in\mathbb{R}}\left|y\right|P\left(f\left(X\right)=y\right)
\]
i.e. to the existence of the expectation of $f\left(X\right).$

When this condition is realized, the packet summation is licit, and
we then have
\begin{align*}
\mathbb{E}\left(f\left(X\right)\right) & =\sum_{y\in\mathbb{R}}yP\left(f\left(X\right)=y\right)\\
 & =\sum_{y\in\mathbb{R}}\left(\sum_{x\in f^{-1}\left(y\right)}f\left(x\right)P\left(X=x\right)\right)\\
 & =\sum_{x\in E}f\left(x\right)P\left(X=x\right).
\end{align*}

\end{proof}

In the case where $\Omega$ is countable, we can take $E=\Omega,$
and let $X$ be the identity function from $\Omega$ to itself. Applying
Theorem $\ref{th:transfer_theorem}$ yields the following statement.

\begin{corollary}{Nessecary and Sufficient Condition for the Existence of The Expectation of a Random Variable}{}

Let $Y$ be a discrete random variable taking real values. 

$Y$ admits an expectation if and only if the family $\left(Y\left(\omega\right)P\left(\left\{ \omega\right\} \right)\right)_{\omega\in\Omega}$
is summable. In this case, the expectation of $Y$ is given by\boxeq{
\[
\mathbb{E}\left(Y\right)=\sum_{\omega\in\Omega}Y\left(\omega\right)P\left(\left\{ \omega\right\} \right).
\]
}

\end{corollary}

\begin{denotation}{Set of Discrete Random Variables with Real Values having an expectation}{}

We denote\footnotemark $\mathcal{L}^{1}_{d}\left(\Omega,\mathcal{A},P\right)$\mindex{L1d @ $\mathcal{L}^{1}_{d}\left(\Omega,\mathcal{A},P\right)$}
the \textbf{set of real-valued discrete\footnotemark random variables
defined on the probabilized space $\left(\Omega,\mathcal{A},P\right)$
that admit an expectation}: that is, all random variables $X$ that
are such that $\mathbb{E}\left(\left|X\right|\right)<+\infty.$

\end{denotation}

\addtocounter{footnote}{-1}

\footnotetext{This is an ``unofficial'' denotation. We also draw
the attention to the readers of Part II, on the fact that, unlike
the $\mathscr{L}^{1}$ classical spaces as they are defined there,
these sets lack of significant topological properties, as in particular,
they are not complete spaces.}

\stepcounter{footnote}

\footnotetext{Tr.N: Giving the $d$ to the notation.}

\begin{proposition}{$\mathcal{L}_d^1\left(\Omega,\mathcal{A},P\right)$  Vectorial Space}{L1_vect_space}

The set $\mathcal{L}^{1}_{d}\left(\Omega,\mathcal{A},P\right)$ is
a vector space and the function $X\mapsto\mathbb{E}\left(X\right)$
is a linear form on $\mathcal{L}^{1}_{d}\left(\Omega,\mathcal{A},P\right).$

\end{proposition}

\begin{proof}{}{}

Let $X$ and $Y$ be two elements of $\mathcal{L}^{1}_{d}\left(\Omega,\mathcal{A},P\right)$
and $a$ and $b$ be any two real numbers. The linear combination
$aX+bY$ is a real-valued discrete random variable. The random variable
$Z=\left(X,Y\right)$ is discrete and taking values in $E=X\left(\Omega\right)\times Y\left(\Omega\right).$

Let $\pi_{1}$ and $\pi_{2}$ be the canonical projections, i.e. the
functions defined for every $\left(x,y\right)\in E$ by $\pi_{1}\left(x,y\right)=x$
and $\pi_{2}\left(x,y\right)=y.$ The random variable $aX+bY$ depends
on the random variable $Z;$ it is denoted
\[
aX+bY=a\pi_{1}\left(Z\right)+b\pi_{2}\left(Z\right).
\]

By Theorem $\ref{th:transfer_theorem},$ the families $\left\{ \pi_{i}\left(z\right)P\left(Z=z\right)\right\} _{z\in E},$
$i=1,2$ are summable. Similarly, by the triangular inequality, the
family $\left\{ \left(a\pi_{1}\left(z\right)+b\pi_{2}\left(z\right)\right)P\left(Z=z\right)\right\} _{z\in E}$
is summable, which shows that the random variable $aX+bY$ admits
an expectation.

Always using Theorem $\ref{th:transfer_theorem}$, we have
\begin{align*}
\mathbb{E}\left(aX+bY\right) & =\sum_{z\in E}\left(a\pi_{1}\left(z\right)+b\pi_{2}\left(z\right)\right)P\left(Z=z\right)\\
 & =a\sum_{z\in E}\pi_{1}\left(z\right)P\left(Z=z\right)+b\sum_{z\in E}\pi_{2}\left(z\right)P\left(Z=z\right),
\end{align*}
which shows that
\[
\mathbb{E}\left(aX+bY\right)=a\mathbb{E}\left(X\right)+b\mathbb{E}\left(Y\right).
\]

\end{proof}

\begin{remark}{}{}

Theorem $\ref{th:transfer_theorem}$ still applies without modification,
with the same proof, when $f$ is a function taking values in $\overline{\mathbb{R}},$
and where $Y=f\left(X\right)$ can take infinite values---however,
it results from the hypotheses that if $Y=\pm\infty,$ this is with
a zero probability.

Proposition $\ref{pr:L1_vect_space}$ extends in the same way, using
the usual convention $0\times\infty=0,$ unless the computation of
$aX\left(\omega\right)+bY\left(\omega\right)$ leads to the expression
$\infty-\infty.$ As previously, if it happens, then it is with a
zero probability.

\end{remark}

\subsection{Expectations of Classical Discrete Laws}

We always refer to the same probabilized space $\left(\Omega,\mathcal{A},P\right).$

\begin{example}{Binomial Law}{expectation_binomial_law}

Let $\left(A_{i}\right)_{1\leqslant i\leqslant n}$ be a family of
$n$ independent events, each with the same probability $p,$ where
$0<p<1.$

Recall that the law of the random variable $S_{n}=\sum^{n}_{i=1}\boldsymbol{1}_{A_{i}}$
is the binomial law $\mathcal{B}\left(n,p\right).$ By the linearity
of expectation, we have
\begin{align*}
\mathbb{E}\left(S_{n}\right) & =\sum^{n}_{i=1}\mathbb{E}\left(\boldsymbol{1}_{A_{i}}\right)=\sum^{n}_{i=1}P\left(A_{i}\right)=np.
\end{align*}

\end{example}

\begin{example}{Geometric Laws}{expectation_geometric_laws}

Let $\left(A_{i}\right)_{1\leqslant i\leqslant n}$ be a family of
$n$ independent events, each with the same probability $p,$ where
$0<p<1.$
\begin{itemize}
\item Recall that the random variable $N,$ defined by
\[
\forall\omega\in\Omega,\,\,\,\,N\left(\omega\right)=\inf\left(n\in\mathbb{N}:\,\omega\in A_{n}\right)
\]
with the convention $\inf\emptyset=+\infty,$ follows the geometric
law $\mathcal{G}_{\mathbb{N}}\left(p\right).$\\
Denoting $q=1-p,$ and using the equation $\refpar{eq:mean_infinite_case},$
we have
\begin{align*}
\mathbb{E}\left(N\right) & =\sum^{+\infty}_{n=0}nP\left(N=n\right)=\sum^{+\infty}_{n=0}npq^{n}=\dfrac{p}{q}.
\end{align*}
\item Recall also that the random variable $N^{\prime},$ defined by 
\[
\forall\omega\in\Omega,\,\,\,\,N^{\prime}\left(\omega\right)=\inf\left(n\in\mathbb{N^{\ast}}:\,\omega\in A_{n}\right)
\]
with the same convention $\inf\emptyset=+\infty,$ follows the geometric
law $\mathcal{G}_{\mathbb{N}^{\ast}}\left(p\right).$ Moreover, defining
$N^{\prime\prime}=1+N,$ the law of $N^{\prime\prime}$ is the same
as the one of $N^{\prime}.$\\
Hence, we obtain 
\begin{align*}
\mathbb{E}\left(N^{\prime}\right) & =\mathbb{E}\left(N^{\prime\prime}\right)=1+\mathbb{E}\left(N\right)=\dfrac{1}{p}.
\end{align*}
\end{itemize}
\end{example}

\begin{example}{Poisson Law}{expectation_poisson_law}

If $X$ is a random variable following the Poisson law $\mathcal{P}\left(\lambda\right),$
with $\lambda>0,$ then by the equality $\refpar{eq:mean_infinite_case},$
\begin{align*}
\mathbb{E}\left(X\right) & =\sum^{+\infty}_{n=0}nP\left(X=n\right)=\sum^{+\infty}_{n=0}n\exp\left(-\lambda\right)\dfrac{\lambda^{n}}{n!}=\lambda.
\end{align*}

\end{example}

\begin{remark}{}{}

If we know that a random variable follows a geometric law---respectively
a Poisson law---then knowing its expectation is sufficient to fully
determine that law. This is not the case for the binomial law, as
we will see later.

\end{remark}

\section{Higher Order Moments}\label{sec:Higher-Order-Moments}

\subsection{Any Order Moments}

\begin{definition}{$p$-th Order Moment}{p_order_moment_definition}

Let $p$ be a real number greater than or equal to 1. We denote by
$\mathcal{L}^{p}_{d}\left(\Omega,\mathcal{A},P\right)$ the set of
discrete real-valued random variables $X$ such that $\mathbb{E}\left(\left|X\right|^{p}\right)<\infty.$

If $p$ is an integer at least 1, and if $X\in\mathcal{L}^{p}_{d}\left(\Omega,\mathcal{A},P\right),$
the real number $\mathbb{E}\left(X^{p}\right)$ is called the \textbf{$p-$th
order moment\mindex{p-th order moment@$p-$th order moment}} of $X.$

\end{definition}

We are going to prove the Hölder and Minkowski inequalities in the
context of summable families, from which we will deduce properties
of the set $\mathcal{L}^{p}_{d}\left(\Omega,\mathcal{A},P\right).$

\begin{definition}{Conjugate Real Numbers}{conjugate_real_numbers}

Two positive real numbers $p$ and $q$ are said to be \textbf{conjugate}\index{conjugate},
if they satisfy the relation
\[
\dfrac{1}{p}+\dfrac{1}{q}=1.
\]

\end{definition}

\begin{remark}{}{}

This condition implies that both $p$ and $q$ must be strictly greater
than $1.$

\end{remark}

\begin{lemma}{}{lemma_pre_holder}

Let $p$ and $q$ be two conjugate real numbers. For every two nonnegative
real numbers $a$ and $b,$
\begin{equation}
ab\leqslant\dfrac{a^{p}}{p}+\dfrac{b^{q}}{q}.\label{eq:ab_lemma}
\end{equation}

\end{lemma}

\begin{proof}{}{}

The function $x\mapsto-\ln x$ is convex on $\left]0,+\infty\right[.$
Therefore for every $x,y>0,$
\[
-\ln\left(\dfrac{x}{p}+\dfrac{y}{q}\right)\leqslant-\dfrac{1}{p}\ln x-\dfrac{1}{q}\ln y,
\]

which implies
\[
\ln\left(\dfrac{x}{p}+\dfrac{y}{q}\right)\geqslant\ln x^{\frac{1}{p}}+\ln y^{\frac{1}{q}}
\]

Hence,
\[
x^{\frac{1}{p}}y^{\frac{1}{q}}\leqslant\dfrac{x}{p}+\dfrac{y}{q}.
\]

It then suffices to substitute $x$ by $a^{p}$ and $y$ by $b^{q}$
to obtain the inequality $\refpar{eq:ab_lemma}$ for $a>0$ and $b>0.$ 

If $a=0$ or $b=0,$ the inequality $\refpar{eq:ab_lemma}$ holds
trivially---provided that the expression $0^{p}$ and $0^{q}$ are
given a meaning.

\end{proof}

\begin{figure}[t]
\begin{center}\includegraphics[width=0.4\textwidth]{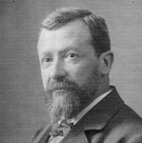}

{\tiny\href{By\%20Archiv\%20Universit\%C3\%A4t\%20Leipzig\%20-\%20Archiv\%20Universit\%C3\%A4t\%20Leipzig,\%20Public\%20Domain,\%20https://commons.wikimedia.org/w/index.php?curid=9110118}{Universität Leipzig}
Public Domain}\end{center}

\caption{\textbf{\protect\href{https://en.wikipedia.org/wiki/Otto_H\%25C3\%25B6lder}{Otto Hölder}}
(1859 - 1937)}\sindex[fam]{Hölder, Otto}
\end{figure}

From this lemma, we can deduce the \textbf{Hölder}{\bfseries\footnote{\href{https://en.wikipedia.org/wiki/Otto_H\%25C3\%25B6lder}{Otto Hölder}
(1859 - 1937) was a German mathematician, known for the Hölder inequality,
in fact proben earlier by Leonard Rogers. He produced different theorems
in geometry and different area of analysis.\label{fn:Holder}}}\textbf{\sindex[fam]{Hölder, Otto} inequality}\index{Hölder inequality}.

\begin{proposition}{Hölder inequality}{holder_inequality}

Let $p$ and $q$ be two conjugate real numbers.

Let $\left(x_{i}\right)_{i\in I},$ $\left(y_{i}\right)_{i\in I}$
and $\left(\rho_{i}\right)_{i\in I}$ be three families of nonnegative
real numbers indexed by the same set $I.$

Then, we have\boxeq{
\begin{equation}
\sum_{i\in I}x_{i}y_{i}\rho_{i}\leqslant\left(\sum_{i\in I}x^{p}_{i}\rho_{i}\right)^{\frac{1}{p}}\left(\sum_{i\in I}x^{q}_{i}\rho_{i}\right)^{\frac{1}{q}}.\label{eq:Holder_inequality}
\end{equation}
}

\end{proposition}

\begin{proof}{}{}

It suffices to prove the inequality $\refpar{eq:Holder_inequality}$
under the assumption that the two terms on the right-hand side are
finite. If one of them is zero---say the first one---then, for every
$i\in I,$ we have $x^{p}_{i}\rho_{i}=0$ and, thus $x_{i}\rho_{i}=0.$ 

We then have
\[
\sum_{i\in I}x_{i}y_{i}\rho_{i}=0.
\]

If both terms on the right-hand side are nonzero, by Lemma $\ref{lm:lemma_pre_holder},$
for each $i\in I,$ we have
\[
\dfrac{x_{i}}{\left(\sum_{i\in I}x^{p}_{i}\rho_{i}\right)^{\frac{1}{p}}}\times\dfrac{y_{i}}{\left(\sum_{i\in I}x^{q}_{i}\rho_{i}\right)^{\frac{1}{q}}}\leqslant\dfrac{1}{p}\dfrac{x^{p}_{i}}{\sum_{i\in I}x^{p}_{i}\rho_{i}}+\dfrac{1}{q}\dfrac{y^{q}_{i}}{\sum_{i\in I}x^{q}_{i}\rho_{i}}.
\]
Multiplying by $\rho_{i}$ which is nonnegative and summing over $I,$
we obtain
\[
\dfrac{\sum_{i\in I}x_{i}y_{i}\rho_{i}}{\left(\sum_{i\in I}x^{p}_{i}\rho_{i}\right)^{\frac{1}{p}}\left(\sum_{i\in I}x^{q}_{i}\rho_{i}\right)^{\frac{1}{q}}}\leqslant\dfrac{1}{p}+\dfrac{1}{q}=1
\]
which proves the inequality $\refpar{eq:Holder_inequality}$.

\end{proof}

\begin{figure}[t]
\begin{center}\includegraphics[width=0.4\textwidth]{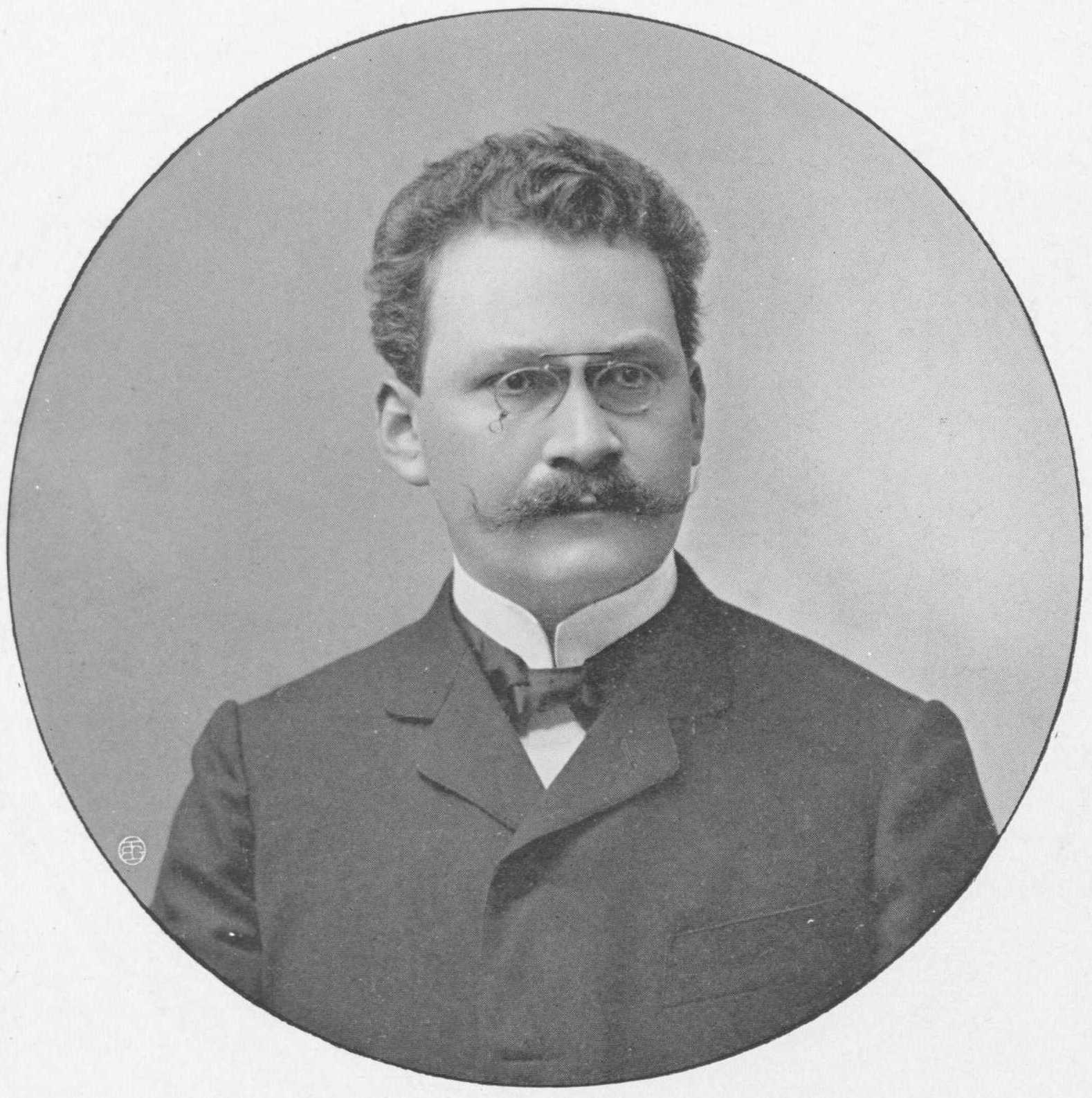}

{\tiny Public Domain}\end{center}

\caption{\textbf{\protect\href{https://en.wikipedia.org/wiki/Hermann_Minkowski}{Hermann Minkowski}}
(1864 - 1909)}\sindex[fam]{Minkowski, Hermann}
\end{figure}

The \textbf{\sindex[fam]{Minkowski, Hermann}Minkowski}{\bfseries\footnote{\href{https://en.wikipedia.org/wiki/Hermann_Minkowski}{Hermann Minkowski}
(1864 - 1909) was a German, Polish, Lithuanian-German or Russian mathematician.
He created and developed the geometry of numbers and convex geometry.
He applied many geometrical methods to solve various physics problem,
including theory of relativity, and some problems in number theory.
He introduced the four-dimensional space, named Minkowski spacetime,
which is part of the ground work for Albert Einstein theory of relativity.\label{fn:Minkowski}}}\textbf{ inequality\index{Minkowski inequality}} can now be deduced.

\begin{proposition}{Minkowski inequality}{minkowski_inequality}

Under the same hypotheses than in Proposition $\ref{pr:holder_inequality}$,
we have

\boxeq{
\begin{equation}
\left(\sum_{i\in I}\left(x_{i}+y_{i}\right)^{p}\rho_{i}\right)^{\frac{1}{p}}\leqslant\left(\sum_{i\in I}x^{p}_{i}\rho_{i}\right)^{\frac{1}{p}}+\left(\sum_{i\in I}y^{p}_{i}\rho_{i}\right)^{\frac{1}{p}}.\label{eq:minkowski ineq}
\end{equation}
}

\end{proposition}

\begin{proof}{}{}As for $p>1,$ the additivity of the sum induces

\[
\sum_{i\in I}\left(x_{i}+y_{i}\right)^{p}\rho_{i}=\sum_{i\in I}\left(x_{i}+y_{i}\right)^{p-1}x_{i}\rho_{i}+\sum_{i\in I}\left(x_{i}+y_{i}\right)^{p-1}y_{i}\rho_{i}.
\]

We now apply the inequality $\refpar{eq:Holder_inequality}$ to each
of the factors on the right part.

Let $q$ be the conjugate exponent of $p,$ that is 
\[
\dfrac{1}{p}+\dfrac{1}{q}=1.
\]

Applying the inequality $\refpar{eq:Holder_inequality},$ we have
\[
\sum_{i\in I}\left(x_{i}+y_{i}\right)^{p-1}x_{i}\rho_{i}\leqslant\left(\sum_{i\in I}\left(x_{i}+y_{i}\right)^{q\left(p-1\right)}\rho_{i}\right)^{\frac{1}{q}}\left(\sum_{i\in I}x^{p}_{i}\rho_{i}\right)^{\frac{1}{p}}
\]
\[
\sum_{i\in I}\left(x_{i}+y_{i}\right)^{p-1}y_{i}\rho_{i}\leqslant\left(\sum_{i\in I}\left(x_{i}+y_{i}\right)^{q\left(p-1\right)}\rho_{i}\right)^{\frac{1}{q}}\left(\sum_{i\in I}y^{p}_{i}\rho_{i}\right)^{\frac{1}{p}}
\]

Adding these two last inequalities and factoring
\[
\sum_{i\in I}\left(x_{i}+y_{i}\right)^{p}\rho_{i}\leqslant\left(\sum_{i\in I}\left(x_{i}+y_{i}\right)^{q\left(p-1\right)}\rho_{i}\right)^{\frac{1}{q}}\left[\left(\sum_{i\in I}x^{p}_{i}\rho_{i}\right)^{\frac{1}{p}}+\left(\sum_{i\in I}y^{p}_{i}\rho_{i}\right)^{\frac{1}{p}}\right]
\]

It suffices then to note that 
\[
q\left(p-1\right)=p
\]
and by dividing by the factor on the right-hand side the inequality
and using the fact that $\dfrac{1}{p}=1-\dfrac{1}{q},$ we obtain
the desired inequality.

\end{proof}

\begin{figure}[t]
\begin{center}\includegraphics[width=0.4\textwidth]{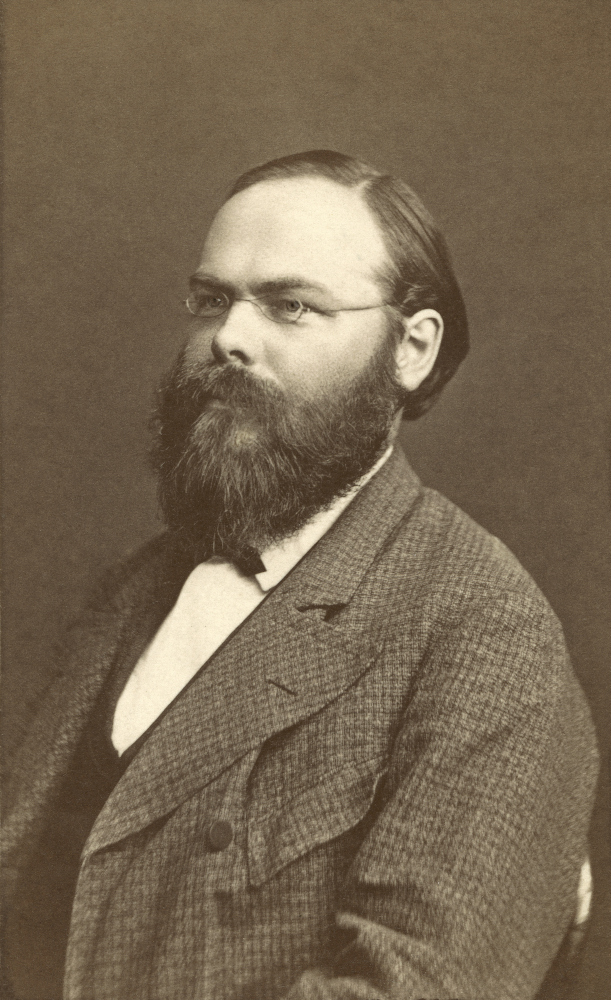}

{\tiny Collection of the \href{https://commons.wikimedia.org/wiki/User:ETH-Bibliothek}{ETH-Bibliothek}
Public Domain}\end{center}

\caption{\textbf{\protect\href{https://en.wikipedia.org/wiki/Hermann_Schwarz}{Hermann Schwarz}}
(1843 - 1921)}\sindex[fam]{Schwarz, Hermann}
\end{figure}

\begin{proposition}{Schwartz Inequality}{}

(i) If $p\geqslant1,$ $\mathcal{L}^{p}_{d}\left(\Omega,\mathcal{A},P\right)$
is a vector space.

(ii) Let $p$ and $q$ be two conjugate real numbers. Let $X\in\mathcal{L}^{p}_{d}\left(\Omega,\mathcal{A},P\right)$
and $Y\in\mathcal{L}^{p}_{d}\left(\Omega,\mathcal{A},P\right).$ Then
$XY\in\mathcal{L}^{1}_{d}\left(\Omega,\mathcal{A},P\right),$ and
the following inequality holds\boxeq{
\begin{equation}
\mathbb{E}\left(\left|XY\right|\right)\leqslant\left(\mathbb{E}\left(\left|X\right|^{p}\right)\right)^{\frac{1}{p}}\left(\mathbb{E}\left(\left|Y\right|^{q}\right)\right)^{\frac{1}{q}}.\label{eq:Schwarz_ineq_p}
\end{equation}
}

With $p=q=2$ in $\refpar{eq:Schwarz_ineq_p}$ and by observing that
$\left|\mathbb{E}\left(XY\right)\right|\leqslant\mathbb{E}\left(\left|XY\right|\right),$
we obtain the so-called \index{Schwarz inequality}\textbf{\sindex[fam]{Schwarz, Hermann}Schwarz}\footnotemark\textbf{
inequality}\boxeq{
\begin{equation}
\left|\mathbb{E}\left(XY\right)\right|\leqslant\left(\mathbb{E}\left(X^{2}\right)\right)^{\frac{1}{2}}\left(\mathbb{E}\left(Y^{2}\right)\right)^{\frac{1}{2}}.\label{eq:Schwarz ineq p-1}
\end{equation}
}

(iii) Let $\alpha$ and $\beta$ be two integers, such that: $1\leqslant\alpha\leqslant\beta.$
Then, we have the set inclusion 
\[
\mathcal{L}^{\beta}_{d}\left(\Omega,\mathcal{A},P\right)\subset\mathcal{L}^{\alpha}_{d}\left(\Omega,\mathcal{A},P\right)
\]
 and
\begin{equation}
\left(\mathbb{E}\left(\left|X^{\alpha}\right|\right)\right)^{\frac{1}{\alpha}}\leqslant\left(\mathbb{E}\left(\left|X^{\beta}\right|\right)\right)^{\frac{1}{\beta}}.\label{eq:moment_order_inequality}
\end{equation}

In particular, a discrete real-valued random variable that has a moment
of order $p\geqslant1,$ also admits any moment of order $p^{\prime}$
such that $1\leqslant p^{\prime}\leqslant p.$

\end{proposition}

\footnotetext{\href{https://en.wikipedia.org/wiki/Hermann_Schwarz}{Hermann Schwarz}
(1843 - 1921) (not to be confused with \href{https://en.wikipedia.org/wiki/Laurent_Schwartz}{Laurent Schwarz},
a French mathematician) was a German mathematician known for his work
in complex analysis, and its Cauchy-Schwarz inequality.\label{fn:schwartz}}

\begin{proof}{}{}

\textbf{(i) $\mathcal{L}^{p}_{d}\left(\Omega,\mathcal{A},P\right)$
is a vector space}

Let $X$ and $Y$ be two elements of $\mathcal{L}^{p}_{d}\left(\Omega,\mathcal{A},P\right)$
and let $a$ and $b$ be any real numbers. 

The function $aX+bY$ is a discrete real-valued random variable. 

The random variable $Z=\left(X,Y\right)$ is discrete and takes values
in $E=\text{val}\left(X\right)\times\text{val}\left(Y\right).$ 

Let $\pi_{1}$ and $\pi_{2}$ be the canonical projections on the
Cartesian product $\text{val}\left(X\right)\times\text{val}\left(Y\right),$
i.e. the functions defined on $E$ by the relations $\pi_{1}\left(x,y\right)=x$
and $\pi_{2}\left(x,y\right)=y$ with $\left(x,y\right)\in E.$

The random variable $aX+bY$ can be expressed as a function of the
random variable $Z:$
\[
aX+bY=a\pi_{1}\left(Z\right)+b\pi_{2}\left(Z\right).
\]

By Theorem $\ref{th:transfer_theorem}$, the families $\left(\left|\pi_{i}\left(z\right)\right|^{p}P\left(Z=z\right)\right)_{z\in E},\,i=1,2,$
are summable. Applying Minkowski inequality stated in Proposition
$\ref{pr:minkowski_inequality}$ yields
\begin{multline*}
\left[\sum_{z\in E}\left|a\pi_{1}\left(z\right)+b\pi_{2}\left(z\right)\right|^{p}P\left(Z=z\right)\right]^{\frac{1}{p}}\\
\leqslant\left|a\right|\left[\sum_{z\in E}\left|\pi_{1}\left(z\right)\right|^{p}P\left(Z=z\right)\right]^{\frac{1}{p}}+\left|b\right|\left[\sum_{z\in E}\left|\pi_{2}\left(z\right)\right|^{p}P\left(Z=z\right)\right]^{\frac{1}{p}}.
\end{multline*}
Since the right-hand side of the inequality is finite, the family
\[
\left(\left|a\pi_{1}\left(z\right)+b\pi_{2}\left(z\right)\right|^{p}P\left(Z=z\right)\right)_{z\in E}
\]
is summable. By the transfer theorem---Theorem $\ref{th:transfer_theorem}$---it
follows that the random variable $aX+bY$ admits a moment of order
$p.$ Hence, $\mathcal{L}^{p}_{d}\left(\Omega,\mathcal{A},P\right)$
is closed under linear combinations and is therefore a vector space.

\textbf{(ii) Schwartz inequality}

Let $X$ and $Y$ be two elements of $\mathcal{L}^{p}_{d}\left(\Omega,\mathcal{A},P\right).$
Using the same notations, it results from the Hölder inequality---Proposition
$\ref{pr:holder_inequality}$---that
\[
\sum_{z\in E}\left|\pi_{1}\left(z\right)\pi_{2}\left(z\right)\right|P\left(Z=z\right)\leqslant\left[\sum_{z\in E}\left|\pi_{1}\left(z\right)\right|^{p}P\left(Z=z\right)\right]^{\frac{1}{p}}\left[\sum_{z\in E}\left|\pi_{2}\left(z\right)\right|^{q}P\left(Z=z\right)\right]^{\frac{1}{q}}.
\]

The right-hand side of this inequality is finite by hypothesis, so
the left-hand side is finite as well. By Theorem $\ref{th:transfer_theorem}$,
the random variable $\left|XY\right|\equiv\left|\pi_{1}\left(Z\right)\pi_{2}\left(Z\right)\right|$
admits an expectation, and confirms that the inequality $\refpar{eq:Schwarz_ineq_p}$
holds.

\textbf{(iii) $\mathcal{L}^{\beta}_{d}\left(\Omega,\mathcal{A},P\right)\subset\mathcal{L}^{\alpha}_{d}\left(\Omega,\mathcal{A},P\right)$
and $\left(\mathbb{E}\left(\left|X^{\alpha}\right|\right)\right)^{\frac{1}{\alpha}}\leqslant\left(\mathbb{E}\left(\left|X^{\beta}\right|\right)\right)^{\frac{1}{\beta}}.$}

Let $X\in\mathcal{L}^{\beta}_{d}\left(\Omega,\mathcal{A},P\right),$
with $\beta>\alpha>0.$ 

Choose $\gamma$ such that
\[
\dfrac{1}{\dfrac{\beta}{\alpha}}+\dfrac{1}{\gamma}=1,
\]
i.e. such that $\gamma=\dfrac{\beta}{\beta-\alpha}.$

Applying the Hölder inequality from Proposition $\ref{pr:holder_inequality}$
to $\left|x\right|^{\alpha}\cdot1,$ yields
\[
\sum_{x\in\mathbb{R}}\left[\left|x\right|^{\alpha}\cdot1\right]P\left(X=x\right)\leqslant\left(\sum_{x\in\mathbb{R}}\left(\left|x\right|^{\alpha}\right)^{\frac{\beta}{\alpha}}P\left(X=x\right)\right)^{\frac{\alpha}{\beta}}\left(\sum_{x\in\mathbb{R}}1^{\gamma}P\left(X=x\right)\right)^{\gamma}.
\]

Since $\sum_{x\in\mathbb{R}}P\left(X=x\right)=1$ and $X\in\mathcal{L}^{\beta}_{d}\left(\Omega,\mathcal{A},P\right),$
it holds
\[
\sum_{x\in\mathbb{R}}\left|x\right|^{\alpha}P\left(X=x\right)\leqslant\left(\sum_{x\in\mathbb{R}}\left|x\right|^{\beta}P\left(X=x\right)\right)^{\frac{\alpha}{\beta}}<+\infty.
\]

From the transfer theorem---Theorem $\ref{th:transfer_theorem}$---we
deduce that $X\in\mathcal{L}^{\alpha}_{d}\left(\Omega,\mathcal{A},P\right)$
and that the inequality $\refpar{eq:moment_order_inequality}$ on
the moments is proved.

\end{proof}

\subsection{Moments of Order Two}

We now undertake a focused study of second-order moments, which play
a fundamental role in numerous theorems and problems in both probability
theory and statistics.

Independently of the previous Subsection, we present a new proposition
that relates the first- and second-order moments of a discrete numerical
random variable---see Definition $\ref{df:p_order_moment_definition}$.

\begin{proposition}{Schwarz Inequality}{schwarz_inequality}

Let $X\in\mathcal{L}^{2}_{d}\left(\Omega,\mathcal{A},P\right).$ 

Then $X\in\mathcal{L}^{1}_{d}\left(\Omega,\mathcal{A},P\right),$
and\boxeq{
\begin{equation}
\mathbb{E}\left(\left|X\right|\right)\leqslant\left(\mathbb{E}\left(X^{2}\right)\right)^{\frac{1}{2}}.\label{eq:mean_order_2_ineq}
\end{equation}
}

Moreover, if $X$ and $Y$ are discrete random variables admitting
second-order moments---that is, $X$ and $Y$ belong to $\mathcal{L}^{2}_{d}\left(\Omega,\mathcal{A},P\right)$---then
the Schwarz inequality holds\boxeq{
\begin{equation}
\left|\mathbb{E}\left(XY\right)\right|\leqslant\left(\mathbb{E}\left(X^{2}\right)\right)^{\frac{1}{2}}\left(\mathbb{E}\left(Y^{2}\right)\right)^{\frac{1}{2}}.\label{eq:Schwarz_ineq_p-1-1}
\end{equation}
}

\end{proposition}

\begin{proof}{}{}

We first observe that, for every $x\in\overline{\mathbb{R}},$ the
inequality $\left|x\right|\leqslant1+x^{2}$ holds. Thus
\[
\left|X\right|\leqslant1+X^{2}.
\]

Since the random variable $X^{2}$ belongs, by hypothesis, to $\mathcal{L}^{1}_{d}\left(\Omega,\mathcal{A},P\right),$
it is as well the same for $X.$ 

Now, define the second degree polynom $Q$ with real number coefficients
by
\begin{align*}
Q\left(\lambda\right) & =\mathbb{E}\left(\left[\left|X\right|-\lambda\mathbb{E}\left(\left|X\right|\right)\right]^{2}\right)\\
 & =\lambda^{2}\left(\mathbb{E}\left(\left|X\right|\right)\right)^{2}-2\lambda\left(\mathbb{E}\left(X\right)\right)^{2}+\mathbb{E}\left(X^{2}\right).
\end{align*}

Since $Q\left(\lambda\right)\geqslant0$ for every $\lambda\in\mathbb{R},$
the reduced discriminant $\Delta^{\prime}$ is non-positive, which
can be written
\[
\Delta^{\prime}=\left(\mathbb{E}\left(\left|X\right|\right)\right)^{2}\left[\mathbb{E}\left(\left|X\right|\right)^{2}-\mathbb{E}\left(X^{2}\right)\right]\leqslant0.
\]

This shows the inequality $\refpar{eq:mean_order_2_ineq}.$

The same reasoning applied to the polynomial $S$ defined by
\[
S\left(\lambda\right)=\mathbb{E}\left(\left[X+\lambda Y\right]^{2}\right)
\]
allows to show the inequality $\refpar{eq:Schwarz_ineq_p-1-1},$ after
having noticed that, since $X$ and $Y$ are discrete real random
variables, $X+\lambda Y$ is also defined.

\end{proof}

\begin{remark}{}{}Under the same hypotheses, by applying the inequality
$\refpar{eq:Schwarz_ineq_p-1-1}$ to the random variables $\left|X\right|$
and $\left|Y\right|,$ we obtain the following inequality\boxeq{
\begin{equation}
\mathbb{E}\left(\left|XY\right|\right)\leqslant\left(\mathbb{E}\left(X^{2}\right)\right)^{\frac{1}{2}}\left(\mathbb{E}\left(Y^{2}\right)\right)^{\frac{1}{2}}.\label{eq:Schwarz ineq p-1-1-1}
\end{equation}
}

\end{remark}

\begin{definition}{Centered Random Variable. Variance. Standard Deviation. Reduced Centered Random Variable}{}

(i) Let $X\in\mathcal{L}^{1}_{d}\left(\Omega,\mathcal{A},P\right).$
The discrete numerical random variable \boxeq{
\[
\mathring{X}=X-\mathbb{E}\left(X\right)
\]
}is called \textbf{the centered random variable\index{centered random variable}\mindex{random variable ! centered}}
associated with $X.$

It satisfies
\[
\mathbb{E}\left(\mathring{X}\right)=0.
\]

(ii) Let $X\in\mathcal{L}^{2}_{d}\left(\Omega,\mathcal{A},P\right).$
Then the random variable $X$ admits an expectation. The real number
$\mathbb{E}\left(\left[X-\mathbb{E}\left(X\right)\right]^{2}\right)$
is called the \textbf{variance}\footnotemark\textbf{\index{variance}}
of the random variable $X$ and denoted $\text{var}\left(X\right).$
Its positive root square is called the \textbf{standard-deviation\index{standard-deviation}}
of $X$ and denoted $\sigma_{X}.$

We thus have\boxeq{
\[
\text{var}\left(X\right)=\sigma^{2}_{X}.
\]
}

(iii) If $X\in\mathcal{L}^{2}_{d}\left(\Omega,\mathcal{A},P\right)$
and $\sigma_{X}\neq0,$ the random variable \boxeq{
\[
\tilde{X}=\dfrac{X-\mathbb{E}\left(X\right)}{\sigma_{X}}
\]
} is called the \textbf{reduced centered random variable}\index{reduced centered random variable}\mindex{random variable ! centered reduced}
associated with $X.$ 

It satisfies
\[
\begin{array}{ccccc}
\mathbb{E}\left(\tilde{X}\right)=0 & \,\,\,\, & \text{and} & \,\,\,\, & \sigma_{\tilde{X}}=1.\end{array}
\]

\end{definition}

\footnotetext{The variance of $X$ is thus the second-order moment
of the centered variable associated to $X$ sometimes referred as
the ``centered moment of order two of $X$'', though this terminology
is rarely used.}

\begin{remarks}{}{}By the transfer theorem---Theorem $\ref{th:transfer_theorem}$---if
$X\in\mathcal{L}^{2}_{d}\left(\Omega,\mathcal{A},P\right),$ then
the variance of $X$ is given by\boxeq{
\begin{equation}
\sigma^{2}_{X}=\sum_{x\in\text{val}\left(X\right)}\left(x-\mathbb{E}\left(X\right)\right)^{2}P\left(X=x\right).\label{eq:variance_from_law}
\end{equation}
}

This formula provides the basis for computing the variance from the
law of $X.$ Since the random variable is discrete, the set $X\left(\Omega\right)$
is countable. Depending on whether it is finite or infinite, this
last formula is often explicited under one of the following forms:
\begin{itemize}
\item If $X\left(\Omega\right)$ is finite, say $X\left(\Omega\right)=\left\{ x_{1},x_{2},\dots,x_{n}\right\} ,$
then\boxeq{
\begin{equation}
\sigma^{2}_{X}=\sum^{n}_{i=1}\left(x_{i}-\mathbb{E}\left(X\right)\right)^{2}P\left(X=x_{i}\right).\label{eq:variance_with_finite_discrete_var}
\end{equation}
}
\item If $X\left(\Omega\right)$ is infinite, and can be enumerated as $\left\{ x_{n}:\,n\in\mathbb{N}\right\} ,$
then\boxeq{
\begin{equation}
\sigma^{2}_{X}=\sum^{+\infty}_{i=0}\left(x_{i}-\mathbb{E}\left(X\right)\right)^{2}P\left(X=x_{i}\right).\label{eq:variance_with_infinite_discrete_var-1}
\end{equation}
}
\end{itemize}
\end{remarks}

The following proposition allows in general an easier computation
of the variance.

\begin{proposition}{Variance Computation from Expectations}{}

Let $X\in\mathcal{L}^{2}_{d}\left(\Omega,\mathcal{A},P\right).$ Then:

(i) The variance of $X$ satisfies\boxeq{
\begin{equation}
\sigma^{2}_{X}=\mathbb{E}\left(X^{2}\right)-\left(\mathbb{E}\left(X\right)\right)^{2}.\label{eq:variance_from_expectations}
\end{equation}
}

(ii) For every real numbers $a$ and $b,$ we have\boxeq{
\begin{equation}
\sigma^{2}_{aX+b}=a^{2}\sigma^{2}_{X}.\label{eq:translation_of_variance}
\end{equation}
}

\end{proposition}

\begin{proof}{}{}

(i) We compute the variance by expanding the square and applying the
linearity of the expectation.
\begin{align*}
\sigma^{2}_{X} & =\mathbb{E}\left(\left[X-\mathbb{E}\left(X\right)\right]^{2}\right)\\
 & =\mathbb{E}\left(X^{2}-2X\mathbb{E}\left(X\right)+\left(\mathbb{E}\left(X\right)\right)^{2}\right)\\
 & =\mathbb{E}\left(X^{2}\right)-\left(\mathbb{E}\left(X\right)\right)^{2}.
\end{align*}

(ii) For every $a,b\in\mathbb{R},$ we have
\[
\mathbb{E}\left(aX+b\right)=a\mathbb{E}\left(X\right)+b.
\]
Thus, the variance is
\begin{align*}
\sigma^{2}_{aX+b} & =\mathbb{E}\left(\left[aX+b-\mathbb{E}\left(aX+b\right)\right]^{2}\right)\\
 & =\mathbb{E}\left(\left[aX+b-\left(a\mathbb{E}\left(X\right)+b\right)\right]^{2}\right)\\
 & =a^{2}\mathbb{E}\left(\left[X-\mathbb{E}\left(X\right)\right]^{2}\right)\\
 & =a^{2}\sigma^{2}_{X}.
\end{align*}

\end{proof}

\begin{remark}{}{}

The formula $\refpar{eq:variance_from_expectations}$ is sometimes
referred to as the \textbf{Leibniz formula\index{Leibniz formula}}
and is commonly used to compute the variance of a random variable.

The formula $\refpar{eq:translation_of_variance}$ expresses the invariance
of the standard deviation under translation and its homogeneity with
respect to scaling.

In summary, changes in origin and scale are reflected through the
formula $\refpar{eq:translation_of_variance}.$

\end{remark}

\subsection{Covariance of Two Random Variables}

\begin{definition}{Covariance of Two Random Variables}{}

Let $X$ and $Y$ be two random variables belonging to $\mathcal{L}^{2}_{d}\left(\Omega,\mathcal{A},P\right).$
By the Schwarz inequality, the random variable $\left(X-\mathbb{E}\left(X\right)\right)\left(Y-\mathbb{E}\left(Y\right)\right)$
admits an expectation, which is called \textbf{covariance\index{covariance}}
of $X$ and $Y,$ and is denoted $\text{cov}\left(X,Y\right),$ given
by\boxeq{
\[
\text{cov}\left(X,Y\right)=\mathbb{E}\left(\left(X-\mathbb{E}\left(X\right)\right)\left(Y-\mathbb{E}\left(Y\right)\right)\right).
\]
}

\end{definition}

The next proposition generally provides a more convenient way to compute
the covariance of two random variables. It also helps to compute the
\index{variance of a sum}\textbf{variance of a sum\mindex{sum ! variance}}
of random variables.

\begin{proposition}{Covariance Computation by Expectation. Variance of a Sum}{}

(i) If $X$ and $Y$ belong to $\mathcal{L}^{2}_{d}\left(\Omega,\mathcal{A},P\right),$
then\boxeq{
\begin{equation}
\text{cov}\left(X,Y\right)=\mathbb{E}\left(XY\right)-\mathbb{E}\left(X\right)\mathbb{E}\left(Y\right).\label{eq:covariance_expectation}
\end{equation}
}

(ii) If $X$ and $Y$ belong to $\mathcal{L}^{2}_{d}\left(\Omega,\mathcal{A},P\right),$
then the variance of their sum satisfies\boxeq{
\begin{equation}
\sigma^{2}_{X+Y}=\sigma^{2}_{X}+\sigma^{2}_{Y}+2\text{cov}\left(X,Y\right).\label{eq:variance_sum_and_cov}
\end{equation}
}

\end{proposition}

\begin{proof}{}{}

(i) We apply the linearity of the expectation after expanding the
product,---Tr.N: i.e. 
\begin{align*}
\text{cov}\left(X,Y\right) & =\mathbb{E}\left(\left(X-\mathbb{E}\left(X\right)\right)\left(Y-\mathbb{E}\left(Y\right)\right)\right)\\
 & =\mathbb{E}\left(XY-X\mathbb{E}\left(Y\right)-\mathbb{E}\left(X\right)Y+\mathbb{E}\left(X\right)\mathbb{E}\left(Y\right)\right)\\
 & =\mathbb{E}\left(XY\right)-\mathbb{E}\left(X\mathbb{E}\left(Y\right)\right)-\mathbb{E}\left(\mathbb{E}\left(X\right)Y\right)+\mathbb{E}\left(X\right)\mathbb{E}\left(Y\right)\\
 & =\mathbb{E}\left(XY\right)-\mathbb{E}\left(X\right)\mathbb{E}\left(Y\right)-\mathbb{E}\left(X\right)\mathbb{E}\left(Y\right)+\mathbb{E}\left(X\right)\mathbb{E}\left(Y\right)\\
 & =\mathbb{E}\left(XY\right)-\mathbb{E}\left(X\right)\mathbb{E}\left(Y\right).
\end{align*}

(ii) First observe that the random variable $X+Y$ is well defined.
Then
\begin{align*}
\sigma^{2}_{X+Y} & =\mathbb{E}\left(\left[X+Y-\mathbb{E}\left(X+Y\right)\right]^{2}\right)\\
 & =\mathbb{E}\left(\left[\left(X-\mathbb{E}\left(X\right)\right)+\left(Y-\mathbb{E}\left(Y\right)\right)\right]^{2}\right)\\
 & =\mathbb{E}\left(\left[X-\mathbb{E}\left(X\right)\right]^{2}\right)+2\mathbb{E}\left(\left[X-\mathbb{E}\left(X\right)\right]\left[Y-\mathbb{E}\left(Y\right)\right]\right)+\mathbb{E}\left(\left[Y-\mathbb{E}\left(Y\right)\right]^{2}\right)\\
 & =\sigma^{2}_{X}+2\text{cov\ensuremath{\left(X,Y\right)}}+\sigma^{2}_{Y}.
\end{align*}

\end{proof}

The following proposition is frequently used.

\begin{proposition}{}{}

Let $X$ and $Y$ be two independent random variables belonging to
$\mathcal{L}^{1}_{d}\left(\Omega,\mathcal{A},P\right),$ the discrete
random variable $XY$ admits an expectation, and\boxeq{
\begin{equation}
\mathbb{E}\left(XY\right)=\mathbb{E}\left(X\right)\mathbb{E}\left(Y\right).\label{eq:expectectation_product}
\end{equation}
}

\end{proposition}

\begin{proof}{}{}

The product $XY$ is a discrete numerical random variable. The pair
$Z=\left(X,Y\right)$ is also a discrete random variable. Using the
same notation as in Proposition $\ref{pr:L1_vect_space},$ we have
\[
XY=\pi_{1}\left(Z\right)\pi_{2}\left(Z\right).
\]

Since the random variables $X$ and $Y$ are independent, we have,
for every $z=\left(x,y\right)\in X\left(\Omega\right)\times Y\left(\Omega\right),$
\begin{equation}
P\left(Z=z\right)=P\left(X=x\right)P\left(Y=y\right).\label{eq:probability_couple}
\end{equation}

Now, using the Fubini property for nonnegative families, we compute
\begin{align*}
\sum_{x\in X\left(\Omega\right)\times Y\left(\Omega\right)}\left|\pi_{1}\left(Z\right)\pi_{2}\left(Z\right)\right|P\left(Z=z\right) & =\sum_{\left(x,y\right)\in X\left(\Omega\right)\times Y\left(\Omega\right)}\left|xy\right|P\left(X=x\right)P\left(Y=y\right)\\
 & =\left(\sum_{x\in X\left(\Omega\right)}\left|x\right|P\left(X=x\right)\right)\\
 & \,\,\,\,\,\,\,\,\,\,\,\,\,\,\,\,\,\,\,\,\,\,\,\,\,\,\,\,\,\,\times\left(\sum_{y\in Y\left(\Omega\right)}\left|y\right|P\left(Y=y\right)\right)\\
 & =\mathbb{E}\left(\left|X\right|\right)\mathbb{E}\left(\left|Y\right|\right)<+\infty.
\end{align*}

It follows that the random variable $\left|XY\right|$ admits an expectation.
By a similar computation---now valid without the absolute values---we
obtain the formula $\refpar{eq:expectectation_product}.$

\end{proof}

\begin{corollary}{Covariance and Variance of Sum of Independent Variables}{covariance_independent}

Let $X$ and $Y$ be independent random variables belonging to $\mathcal{L}^{2}_{d}\left(\Omega,\mathcal{A},P\right).$ 

Then\boxeq{
\[
\text{cov}\left(X,Y\right)=0\,\,\,\,\text{and\,\,\,\,}\sigma^{2}_{X+Y}=\sigma^{2}_{X}+\sigma^{2}_{Y}.
\]
}

\end{corollary}

\begin{proof}{}{} The previous proposition shows that 
\[
\text{cov}\left(X,Y\right)=0.
\]

To obtain the second equality, it suffices to substitute into $\refpar{eq:variance_sum_and_cov}.$

\end{proof}

\begin{remarks}{}{}

1. This corollary is frequently used. In particular, it implies that
if $X_{1},\dots,X_{n}$ are independent random variables with same
law and finite variance, then
\[
\sigma_{X_{1}+\dots+X_{n}}=\sigma_{X_{1}}\sqrt{n}.
\]

2. As the following example shows, the condition $\text{cov}\left(X,Y\right)=0$
does not necessarily imply that $X$ and $Y$ are independent.

\end{remarks}

\begin{example}{Non-Independent Variables with Null Covariance}{}

Let $U$ and $V$ be independent random variables belonging to $\mathcal{L}^{2}_{d}\left(\Omega,\mathcal{A},P\right).$

Let $a,b,c,d$ be four nonzero real numbers.

Define the discrete random variables with real values $X$ and $Y$
by
\[
\left\{ \begin{array}{c}
X=aU+bV\\
Y=cU+dV.
\end{array}\right.
\]

1. Compute $\text{cov}\left(X,Y\right).$

2. Now suppose that $U$ and $V$ have the same law given by
\[
\begin{array}{ccc}
P\left(U=1\right)=p & \,\,\,\, & P\left(U=0\right)=1-p.\end{array}
\]

Additionally, suppose that $ac+bd=0$ with $a+b\neq0$ and $c+d\neq0.$

Show that $\text{cov}\left(X,Y\right)=0,$ but that $X$ and $Y$
are not independent.

\end{example}

\begin{solutionexample}{}{}

1. Since $U$ and $V$ are independent random variables, we have
\begin{align*}
\text{cov}\left(X,Y\right) & =\mathbb{E}\left(\left(\mathring{\overbrace{aU+bV}}\right)\left(\mathring{\overbrace{cU+dV}}\right)\right)\\
 & =\mathbb{E}\left(\left(a\mathring{U}+b\mathring{V}\right)\left(c\mathring{U}+d\mathring{V}\right)\right)\\
 & =ac\sigma^{2}_{U}+\left(ad+bc\right)\text{cov}\left(U,V\right)+bd\sigma^{2}_{V}\\
 & =ac\sigma^{2}_{U}+bd\sigma^{2}_{V}.
\end{align*}

2. As $U$ and $V$ follow the same law, we have $\sigma^{2}_{U}=\sigma^{2}_{V}.$

Since $ac+bd=0,$ $\text{cov}\left(X,Y\right)=0.$

Nevertheless, since $a+b\neq0$ and $c+d\neq0,$ $X$ and $Y$ are
not independent, as shown by the following equalities between events
\[
\left(X=0\right)=\left(U=0,V=0\right)=\left(Y=0\right).
\]

Thus, on the one hand
\begin{align*}
P\left(X=0\right) & =P\left(Y=0\right)\\
 & =P\left(U=0,V=0\right)\\
 & =P\left(U=0\right)P\left(V=0\right)\\
 & =\left(1-p\right)^{2},
\end{align*}
and on the other hand
\begin{align*}
P\left(X=0,Y=0\right) & =P\left(U=0,V=0\right)\\
 & =\left(1-p\right)^{2}.
\end{align*}

Thus for $0<p<1,$
\[
P\left(X=0,Y=0\right)\neq P\left(X=0\right)P\left(Y=0\right).
\]

Hence, the random variables $X$ and $Y$ are not independent in this
case!

\end{solutionexample}

\subsection{Variances of Classical Discrete Laws}

We continue with the notations introduced in Examples $\ref{ex:expectation_bernoulli_law}$,
$\ref{ex:expectation_binomial_law}$, $\ref{ex:expectation_geometric_laws}$,
$\ref{ex:expectation_poisson_law}$.

\subsubsection*{Bernoulli Law}

Since $\boldsymbol{1}^{2}_{A}=\boldsymbol{1}_{A},$
\begin{align*}
\sigma^{2}_{\boldsymbol{1}_{A}} & =\mathbb{E}\left(\boldsymbol{1}^{2}_{A}\right)-\left[\mathbb{E}\left(\boldsymbol{1}_{A}\right)\right]^{2}=p-p^{2}
\end{align*}

Hence\boxeq{
\[
\sigma^{2}_{\boldsymbol{1}_{A}}=p\left(1-p\right)
\]
}

\subsubsection*{Binomial Law}

The family of random variables $\left(\boldsymbol{1}_{A_{i}}\right)_{1\leqslant i\leqslant n}$
is a family of independent random variables, with the same law, and
therefore the same variance. From Corollary $\ref{co:covariance_independent}$,
we obtain
\begin{align*}
\sigma^{2}_{S_{n}} & =\sum^{n}_{i=1}\sigma^{2}_{\boldsymbol{1}_{A_{i}}}=n\sigma^{2}_{\boldsymbol{1}_{A_{1}}}.
\end{align*}

Thus,\boxeq{
\[
\sigma^{2}_{S_{n}}=np\left(1-p\right).
\]
}

\subsubsection*{Geometric Law on $\mathbb{N}$}

Let the random variable $N$ follow the law $\mathcal{G}_{\mathbb{N}}\left(p\right).$
We have
\[
\mathbb{E}\left(N^{2}\right)=\mathbb{E}\left(N\left(1+N\right)\right)+\mathbb{E}\left(N\right).
\]

Using the transfer theorem
\[
\mathbb{E}\left(N\left(N-1\right)\right)=\sum^{+\infty}_{n=0}n\left(n-1\right)pq^{n}.
\]

Applying the differentiation theorem for power series---an argument
already given---we obtain for every $x$ such that $0\leqslant x<1$
\[
\sum^{+\infty}_{n=2}n\left(n-1\right)x^{n-2}=\dfrac{\text{d}^{2}}{\text{d}x^{2}}\left(\sum^{+\infty}_{n=0}x^{n}\right).
\]

It follows that
\begin{align*}
\mathbb{E}\left(N^{2}\right) & =\dfrac{2pq^{2}}{\left(1-q\right)^{3}}+\dfrac{q}{p}=\dfrac{q}{p}\left(2\dfrac{q}{p}+1\right).
\end{align*}

Therefore, using the equality $\refpar{eq:variance_from_expectations}$
and simplifying\boxeq{
\[
\sigma^{2}_{N}=\dfrac{q}{p^{2}}.
\]
}

\subsubsection*{Geometric Law on $\mathbb{N}^{\ast}$}

The random variables $N^{\prime}$ and $1+N$ follow the same geometric
law $\mathcal{G}_{\mathbb{N^{\ast}}}\left(p\right),$ and thus have
the same variance. Since the variance is invariant under translation,
it follows that

\boxeq{
\[
\sigma^{2}_{N^{\prime}}=\dfrac{q}{p^{2}}.
\]
}

\subsubsection*{Poisson Law}

Let $X$ be a random variable following the Poisson law $\mathcal{P}\left(\lambda\right).$ 

We have
\[
\mathbb{E}\left(X^{2}\right)=\mathbb{E}\left(X\left(X-1\right)\right)+\mathbb{E}\left(X\right),
\]
and, by the transfer theorem,
\[
\mathbb{E}\left(X\left(X-1\right)\right)=\sum^{+\infty}_{n=0}n\left(n-1\right)\exp\left(-\lambda\right)\dfrac{\lambda^{n}}{n!}.
\]

By reindexing the sum,
\[
\sum^{+\infty}_{n=2}n\left(n-1\right)\dfrac{\lambda^{n}}{n!}=\lambda^{2}\exp\lambda.
\]

Thus,
\[
\mathbb{E}\left(X^{2}\right)=\lambda^{2}+\lambda,
\]

Using the equality $\refpar{eq:variance_from_expectations},$ we obtain\boxeq{
\[
\sigma^{2}_{X}=\lambda.
\]
}

\begin{remark}{}{}

For a random variable following a Poisson law, the parameter of the
law represents both the expectation and the variance.

\end{remark}

\subsection{Markov and Chebyshev Inequalities}

The classical Chebyshev inequality provides an upper-bound of the
probability that a random variable deviates ``too far'' from its
expectation. It is considered a rough estimate compared to more accurate
results that are available in some particular cases---see Chapter
\ref{chap:Approximation-of-laws.}, Example $\ref{ex:example7.1}$
and Exercise $\ref{exo:exercise7.2}$. This inequality is mainly of
theoretical importance, and appears, among others, in the proof of
the weak law of large numbers. 

We first give the \textbf{\index{Markov inequality}Markov inequality}
that is useful to prove the Chebyshev inequality.

\begin{proposition}{Markov Inequality}{markov_inequality}

If $X\in\mathcal{L}^{1}_{d}\left(\Omega,\mathcal{A},P\right),$ then,
for every $\epsilon>0,$\boxeq{
\begin{equation}
P\left(\left|X\right|>\epsilon\right)\leqslant\dfrac{\mathbb{E}\left(\left|X\right|\right)}{\epsilon}.
\end{equation}
}

\end{proposition}

\begin{proof}{}{}

Consider the set $D=\left\{ x\in X\left(\Omega\right):\,\left|x\right|>\epsilon\right\} .$

We then have the following successive lower bounds:
\begin{align*}
\mathbb{E}\left(\left|X\right|\right) & =\sum_{x\in\mathbb{R}}\left|x\right|P\left(X=x\right)\\
 & \geqslant\sum_{x\in D}\left|x\right|P\left(X=x\right)\\
 & \geqslant\epsilon\sum_{x\in D}P\left(X=x\right).
\end{align*}

Since $X\left(\Omega\right)$ is countable, it follows that:
\begin{align*}
\mathbb{E}\left(\left|X\right|\right) & \geqslant\epsilon P\left(X\in D\right)\\
 & =\epsilon P\left(X^{-1}\left(D\right)\right),
\end{align*}
 which proves the inequality.

\end{proof}

We deduce the \textbf{\index{Chebyshev inequality}Chebyshev inequality}.

\begin{proposition}{Chebyshev Inequality}{chebychev_inequality}

If $X\in\mathcal{L}^{2}_{d}\left(\Omega,\mathcal{A},P\right),$ then,
for every $\epsilon>0,$\boxeq{
\begin{equation}
P\left(\left|X-\mathbb{E}\left(X\right)\right|>\epsilon\right)\leqslant\dfrac{\sigma^{2}_{X}}{\epsilon^{2}}.\label{eq:Chebychev_inequality}
\end{equation}
}

\end{proposition}

The Chebyshev inequality is sometimes used in the following equivalent
form:\boxeq{
\[
P\left(\left|X-\mathbb{E}\left(X\right)\right|>\epsilon\sigma_{X}\right)\leqslant\dfrac{1}{\epsilon^{2}}.
\]
}

\begin{proof}{}{}

It suffices to apply the Markov inequality to the random variable
$\left(X-\mathbb{E}\left(X\right)\right)^{2}$ and noting that
\[
\left\{ \left(X-\mathbb{E}\left(X\right)\right)^{2}>\epsilon^{2}\right\} =\left\{ \left|X-\mathbb{E}\left(X\right)\right|>\epsilon\right\} .
\]

The definition of the variance then yields the desired result.

\end{proof}

\begin{remarks}{}{}

1. A direct proof of this inequality can be given as follows.

Consider the set
\[
D=\left\{ x\in X\left(\Omega\right):\,\left|x-\mathbb{E}\left(X\right)\right|>\epsilon\right\} .
\]

Then we have the following sequence of lower bounds
\begin{align*}
\sigma^{2}_{X} & =\sum_{x\in\mathbb{R}}\left[x-\mathbb{E}\left(X\right)\right]^{2}P\left(X=x\right)\\
 & \geqslant\sum_{x\in D}\left[x-\mathbb{E}\left(X\right)\right]^{2}P\left(X=x\right)\\
 & \geqslant\sum_{x\in D}\epsilon^{2}P\left(X=x\right).
\end{align*}
Since $X\left(\Omega\right)$ is countable,
\[
\sigma^{2}_{X}\geqslant\epsilon^{2}P\left(X\in D\right)=\epsilon^{2}P\left[X^{-1}\left(D\right)\right],
\]
which proves the inequality.

2. The smaller $\sigma_{X}$ is, the more concentrated the random
variable $X$ is around its expectation. The extreme case of this
concentration occurs when $\sigma_{X}=0.$ It is clear that if a random
variable $X$ is equal to a constant with probability 1, then its
expectation is that constant, and $\sigma_{X}=0.$ The Chebyshev inequality
gives the converse: if $\sigma_{X}=0,$ we have for every integer
$n,$ 
\[
P\left(\left|X-\mathbb{E}\left(X\right)\right|\geqslant\dfrac{1}{n}\right)=0.
\]
 Hence
\[
P\left(\left|X-\mathbb{E}\left(X\right)\right|>0\right)=\lim_{n\to+\infty}\downarrow P\left(\left|X-\mathbb{E}\left(X\right)\right|\geqslant\dfrac{1}{n}\right)=0
\]
 and so
\[
P\left(X=\mathbb{E}\left(X\right)\right)=1.
\]
 This can also be shown by a direct argument.

\end{remarks}

\subsection{Correlation coefficient. Linear regression}

As we will see, the \textbf{correlation coefficient} of two random
variables allows us to measure a certain degree of relationships between
them.

\begin{definition}{Correlation Coefficient}{correlation_coefficient}

Let $X$ and $Y$ belong to $\mathcal{L}^{2}_{d}\left(\Omega,\mathcal{A},P\right),$
and assume that $\sigma_{X}\neq0$ and $\sigma_{Y}\neq0.$

We define the \textbf{correlation coefficient\index{correlation coefficient}}
of $X$ and $Y$ as the real number denoted $\rho_{X,Y}$ and defined
by\boxeq{
\[
\rho_{X,Y}=\dfrac{\text{cov}\left(X,Y\right)}{\sigma_{X}\sigma_{Y}}.
\]
}

\end{definition}

\begin{proposition}{}{}

Let $X$ and $Y$ be in $\mathcal{L}^{2}_{d}\left(\Omega,\mathcal{A},P\right)$
and such that $\sigma_{X}\neq0$ and $\sigma_{Y}\neq0$ and having
as correlation coefficient $\rho_{X,Y}.$

(i) We have\boxeq{
\[
\left|\rho_{X,Y}\right|\leqslant1.
\]
}

(ii) Equality $\left|\rho_{X,Y}\right|=1$ holds if and only if there
exists three nonzero real numbers $a,b,c$ such that
\[
P\left(aX+bY+c=0\right)=1.
\]

\end{proposition}

\begin{proof}{}{}

(i) The Schwarz inequality gives
\begin{align*}
\left|\text{cov}\left(X,Y\right)\right| & =\left|\mathbb{E}\left(\mathring{X}\mathring{Y}\right)\right|\leqslant\left(\mathbb{E}\left(\mathring{X}^{2}\right)\right)^{\frac{1}{2}}\left(\mathbb{E}\left(\mathring{Y}^{2}\right)\right)^{\frac{1}{2}}=\sigma_{X}\sigma_{Y},
\end{align*}
which proves the stated inequality.

(ii) (a) Suppose $\left|\rho_{X,Y}\right|=1.$ Then the second-degree
polynomial in $\lambda,$ $\mathbb{E}\left(\left(\mathring{X}+\lambda\mathring{Y}\right)^{2}\right),$
admits a double root $\lambda_{0}.$ 

Therefore, its discriminant is equal to zero, and
\[
\mathbb{E}\left(\left(\mathring{X}+\lambda\mathring{Y}\right)^{2}\right)=0.
\]

By the remark following Proposition $\ref{pr:chebychev_inequality},$
it follows that
\[
P\left(\mathring{X}+\lambda_{0}\mathring{Y}=0\right)=1.
\]
 (b) Conversely, supposing there exists three real numbers $a,b,c$
not all zero, such that
\begin{equation}
P\left(aX+bY+c=0\right)=1.\label{eq:probability_affine_relation}
\end{equation}

If $c\neq$0, then both $a$ and $b$ are different from 0. Indeed,
if for instance, $a=0,$ then $P\left(bY+c=0\right)=1,$ and thus
$\sigma^{2}_{bY+c}=\sigma^{2}_{0}=0,$ which implies $b^{2}\sigma^{2}_{Y}=0$
and thus $b=0,$ which leads to $P\left(c=0\right)=1,$ a contradiction.

In this case, we can rewrite the equation as
\begin{equation}
P\left(X=\alpha Y+\beta\right)=1\label{eq:probability_reduced_affine_relation}
\end{equation}
where $\alpha\neq0.$

If $c=0,$ then we have $a\neq0$ or $b\neq0.$ If for instance $a\neq0,$
the equality $\refpar{eq:probability_reduced_affine_relation}$ is
still satisfied with $\beta=0$---with $b\neq0,$ the computation
is similar.

In the two cases, we thus have
\begin{align*}
\text{cov\ensuremath{\left(X,Y\right)}} & =\mathbb{E}\left(\left(\alpha\mathring{Y}\right)\mathring{Y}\right)=\alpha\sigma^{2}_{Y}
\end{align*}
and
\[
\sigma^{2}_{X}=\sigma^{2}_{\alpha Y+\beta}=\alpha^{2}\sigma^{2}_{Y}.
\]

Thus, the correlation coefficient is
\[
\rho_{X,Y}=\dfrac{\alpha\sigma^{2}_{Y}}{\left|\alpha\right|\sigma^{2}_{Y}}.
\]
which gives $\left|\rho_{X,Y}\right|=1.$

\end{proof}

\begin{remark}{}{}

We have used the following result: if $X\in\mathcal{L}^{1}_{d}\left(\Omega,\mathcal{A},P\right),$
and if $Y$ is a discrete random variable such that $P\left(X=Y\right)=1,$
then $Y\in\mathcal{L}^{1}_{d}\left(\Omega,\mathcal{A},P\right)$ and
$\mathbb{E}\left(X\right)=\mathbb{E}\left(Y\right).$

Indeed, consider
\[
D=\left\{ \left(x,y\right)\in\mathbb{R}^{2}:\,x=y\right\} .
\]
 If $Z=\left(X,Y\right),$ then $P_{Z}\left(D\right)=1.$

Still using the notations of Proposition $\ref{pr:L1_vect_space},$
we have
\begin{align*}
\mathbb{E}\left(\left|X\right|\right) & =\mathbb{E}\left(\left|\pi_{1}\left(Z\right)\right|\right)=\sum_{z\in\mathbb{R}^{2}}\left|\pi_{1}\left(Z\right)\right|P\left(Z=z\right).
\end{align*}

Since $P\left(Z\in D\right)=1$ and $\pi_{1}\left(z\right)=\pi_{2}\left(z\right)$
on $D,$ it follows that
\begin{align*}
\mathbb{E}\left(\left|X\right|\right) & =\sum_{z\in D}\left|\pi_{1}\left(z\right)\right|P\left(Z=z\right)\\
 & =\sum_{z\in D}\left|\pi_{2}\left(z\right)\right|P\left(Z=z\right)\\
 & =\sum_{z\in\mathbb{R}^{2}}\left|\pi_{2}\left(z\right)\right|P\left(Z=z\right),
\end{align*}
which proves that $Y\in\mathcal{L}^{1}_{d}\left(\Omega,\mathcal{A},P\right).$
The same computation without the absolute values is then valid and
yields the equality
\[
\mathbb{E}\left(X\right)=\mathbb{E}\left(Y\right).
\]

\end{remark}

Proposition $\ref{df:correlation_coefficient}$ states that if $\rho_{X,Y}=\pm1,$
then one of the two variables $X$ and $Y$ is a linear function of
the other. Additionally, if $X$ and $Y$ are independent, then $\rho_{X,Y}=0$---although
the converse is false. Hence, if $\rho_{X,Y}\neq0,$ we can conclude
that $X$ and $Y$ are not independent. 

The correlation coefficient is primarily used in Statistics for series
of empirical observations---which can be considered as random variables
taking a finite number of values and following a uniform law over
this set of values. Generally, we interpret a value of $\left|\rho_{X,Y}\right|$
close to 1 as indicating that $Y$ is approximately a linear function
of $X.$

\subsubsection*{The problem of the linear regression}\label{subsec:The-problem-of-linear-regression}

\begin{leftbar}

Given the random variables $X$ and $Y$ belonging to $\mathcal{L}^{2}_{d}\left(\Omega,\mathcal{A},P\right),$
we look for the best approximation of $Y$ as linear function of $X$
in the \index{sense of least squares}\textbf{sense of least squares}.
That is, an element $\left(\widehat{a},\widehat{b}_{\widehat{a}}\right)\in\mathbb{R}^{2}$
of the following set corresponding to the minimization problem\boxeq{
\[
\left\{ \left(x,y\right)\in\mathbb{R}^{2}:\phi\left(x,y\right)=\left\{ \inf\left\{ \phi\left(a,b\right):\,\left(a,b\right)\in\mathbb{R}^{2}\right\} \right\} \right\} ,
\]
where
\[
\phi\left(a,b\right)=\mathbb{E}\left(\left(Y-\left(aX+b\right)\right)^{2}\right).
\]
}

This problem is also known as the \textbf{\index{linear regression problem}linear
regression problem}.

\end{leftbar}

We can rewrite $\phi\left(a,b\right)$ as
\begin{align*}
\phi\left(a,b\right) & =\mathbb{E}\left(\left(\mathring{Y}-a\mathring{X}+\left(\mathbb{E}\left(Y\right)-a\mathbb{E}\left(X\right)-b\right)\right)^{2}\right)\\
 & =\mathbb{E}\left(\left(\mathring{Y}-a\mathring{X}\right)^{2}\right)+\left(\mathbb{E}\left(Y\right)-a\mathbb{E}\left(X\right)-b\right)^{2}.
\end{align*}

Given any fixed $a,$ this expression is minimized when $\hat{b}_{a}=\mathbb{E}\left(Y\right)-a\mathbb{E}\left(X\right),$
that is when 
\[
\mathbb{E}\left(Y\right)-a\mathbb{E}\left(X\right)-b=0
\]

We now still have to minimize in $a$ the polynomial
\begin{align*}
f(a)=\phi\left(a,\hat{b}_{a}\right) & =\mathbb{E}\left(\left(\mathring{Y}-a\mathring{X}\right)^{2}\right)\\
 & =\sigma^{2}_{Y}-2a\text{cov}\left(X,Y\right)+a^{2}\sigma^{2}_{X}.
\end{align*}

Taking the derivative, we have
\[
f^{\prime}\left(a\right)=2a\sigma^{2}_{X}-2\text{cov}\left(X,Y\right).
\]

We find
\[
f'\left(a\right)\geqslant0\Leftrightarrow a\geqslant\dfrac{\text{cov\ensuremath{\left(X,Y\right)}}}{\sigma^{2}_{X}}.
\]

Let
\[
\hat{a}=\dfrac{\text{cov\ensuremath{\left(X,Y\right)}}}{\sigma^{2}_{X}}
\]
then $f$ reaches a minimum at $\hat{a}.$ The solution to the linear
regression problem is thus the pair $\left(\hat{a},\hat{b}_{\hat{a}}\right)$
given by
\[
\begin{cases}
\hat{a}=\rho_{X,Y}\dfrac{\sigma_{Y}}{\sigma_{X}}\\
\hat{b}_{\hat{a}}=\mathbb{E}\left(Y\right)-\hat{a}\mathbb{E}\left(X\right)=\mathbb{E}\left(Y\right)-\mathbb{E}\left(X\right)\cdot\rho_{X,Y}\dfrac{\sigma_{Y}}{\sigma_{X}}.
\end{cases}
\]

Thus, the line $D$ with equation\boxeq{
\[
\left(y-\mathbb{E}\left(Y\right)\right)-\rho_{X,Y}\dfrac{\sigma_{Y}}{\sigma_{X}}\left(x-\mathbb{E}\left(X\right)\right)=0
\]
}is called the \textbf{linear regression line of $Y$ in $X.$\index{linear regression line of $Y$ in $X$}} 

The best approximation $\widetilde{Y}$ of $Y$ as a linear function
of $X,$ in the sense of the least squares, is\boxeq{
\[
\widetilde{Y}=\mathbb{E}\left(Y\right)-\rho_{X,Y}\dfrac{\sigma_{Y}}{\sigma_{X}}\left(X-\mathbb{E}\left(X\right)\right).
\]
}

Additionally, $P\left(\left(X,Y\right)\in D\right)=1$ if and only
if $\phi\left(\hat{a},\hat{b_{\hat{a}}}\right)=0.$

\subsubsection*{Particular case}

If the random variable follows the uniform law on the set of $n$
points in the plane
\[
\left\{ \left(x_{i},y_{i}\right):\,1\leqslant i\leqslant n\right\} ,
\]
then\boxeq{
\[
\phi\left(a,b\right)=\dfrac{1}{n}\sum^{n}_{i=1}\left(y_{i}-\left(ax_{i}+b\right)\right)^{2}.
\]
}

In this way, we recover the \textbf{least squares approximation line}
used in physics. It is let as an exercise to find the equation of
this line.

\section{Generating Functions}\label{sec:Generative-Functions}

A generating function of a random variable taking values in $\mathbb{N}$
characterizes the law of that random variable. Generating functions
are then a convenient computational tool. They are especially useful
to compute the laws of sums of independent random variables, as well
as for studying the convergence in law, as it will be seen in the
second part of this book.

In this section, unless otherwise specified, the random variables
are assumed to be defined on a probabilized space $\left(\Omega,\mathcal{A},P\right)$
and to take values in $\mathbb{N}.$

\subsection{Definition}

\begin{lemma}{}{}

Let $X$ be a random variable defined on the probabilized space $\left(\Omega,\mathcal{A},P\right)$
and taking values in $\mathbb{N}.$ 

For every $s\in\left[-1,1\right],$ the random variable $s^{X}$ belongs
to $\mathcal{L}^{1}_{d}\left(\Omega,\mathcal{A},P\right).$

\end{lemma}

\begin{proof}{}{}

We observe that $\left|s^{X}\right|\leqslant1.$

Since the constant random variable equals to $1$ admits an expectation,
it follows that $s^{X}$ also has an expectation.

We denote, when this quantity exists,
\[
G_{X}\left(s\right)=\mathbb{E}\left(s^{X}\right).
\]

\end{proof}

\begin{definition}{Generating Function}{}

The function $G_{X}$ is called the \textbf{\index{generating function}generating
function} of the random variable $X.$

\end{definition}

\begin{proposition}{Generating Function: Properties and Law Caracterization}{}

Let $X$ be a random variable defined on the probabilized space $\left(\Omega,\mathcal{A},P\right),$
taking values in $\mathbb{N},$ with generating function $G_{X},$
and following the law defined for every $n\in\mathbb{N}$ by 
\[
P\left(X=n\right)=p_{n}
\]
where $0\leqslant p_{n}\leqslant1,\,n\in\mathbb{N}$ and such that
$\sum_{n\in\mathbb{N}}p_{n}=1.$

(i) The domain of definition of $G_{X}$ contains the interval $\left[-1,1\right].$
We have\boxeq{
\[
\forall s\in\left[-1,1\right],\,\,\,\,\left|G_{X}\left(s\right)\right|\leqslant1\,\,\,\,\text{and}\,\,\,\,G_{X}\left(1\right)=1.
\]
}

(ii) For every $s\in\left[-1,1\right],$\boxeq{
\[
G_{X}\left(s\right)=\sum^{+\infty}_{n=0}p_{n}s^{n}.
\]
}

(iii) The function $G_{X}$ is continuous on $\left[-1,1\right],$
and $C^{\infty}$ on $\left]-1,1\right[.$

(iv) The generating function $G_{X}$---and even its restriction
to the interval $\left[-1,1\right]$---\textbf{\index{characterizes the law}characterizes
the law} of $X;$ specifically\boxeq{
\[
\forall n\in\mathbb{N},\,\,\,\,P\left(X=n\right)=\dfrac{G^{\left(n\right)}_{X}\left(0\right)}{n!}.
\]
}

\end{proposition}

\begin{proof}{}{}

(i) From the previous lemma, we know that $G_{X}\left(s\right)=\mathbb{E}\left(s^{X}\right)$
is defined for every $s\in\left[-1,1\right].$ To prove the inequality,
observe that, for every $s\in\left[-1,1\right],$
\[
\left|G_{X}\left(s\right)\right|\leqslant\mathbb{E}\left(s^{X}\right)\leqslant\mathbb{E}\left(1\right)=1.
\]

(ii) This follows directly from the transfer theorem.

(iii) It suffices to note that $G_{X}$ is a power series with a radius
of convergence at least 1---which is also the case for the last point.

\end{proof}

\subsection{Generating Functions of Classical Laws with Values in $\mathbb{N}$}

We keep the notations previously introduced. Let us fix an arbitrary
$s\in\left[-1,1\right].$

\subsubsection{Generating Function of the Binomial Law}

\begin{proposition}{}{}

If $P_{X}=\mathcal{B}\left(n,p\right),$ then the generating function
of $X$ is\boxeq{
\[
G_{X}\left(s\right)=\left(ps+q\right)^{n}
\]
}where $q=1-p.$

\end{proposition}

\begin{proof}{}{}

Indeed, by expansion via the binomial theorem,
\[
G_{X}\left(s\right)=\sum^{n}_{k=0}\left(\begin{array}{c}
n\\
k
\end{array}\right)p^{k}q^{n-k}s^{k}.
\]

\end{proof}

\subsubsection{Generating Function of the Poisson Law}

\begin{proposition}{}{}

If $P_{X}=\mathcal{P}\left(\lambda\right),$ the Poisson law with
parameter $\lambda>0,$ then the generating function of $X$ is\boxeq{
\[
G_{X}\left(s\right)=\exp\left(\lambda\left(s-1\right)\right).
\]
}

\end{proposition}

\begin{proof}{}{}

Indeed, we compute
\begin{align*}
G_{X}\left(s\right) & =\sum^{+\infty}_{k=0}\exp\left(-\lambda\right)\dfrac{\lambda^{k}}{k!}s^{k}\\
 & =\exp\left(-\lambda\right)\sum^{+\infty}_{k=0}\dfrac{\left(\lambda s\right)^{k}}{k!}\\
 & =\exp\left(-\lambda\right)\exp\left(\lambda s\right)\\
 & =\exp\left(\lambda\left(s-1\right)\right).
\end{align*}

\end{proof}

\subsubsection{Generating Function of the Geometric Law on $\mathbb{N}$}

\begin{proposition}{}{}

If $P_{X}=\mathcal{G}_{\mathbb{N}}\left(p\right),$ geometric law
with parameter $p\in\left]0,1\right[,$ then\boxeq{
\[
G_{X}\left(s\right)=\dfrac{p}{1-qs}.
\]
}

\end{proposition}

\begin{proof}{}{}

Indeed, for $\left|s\right|\leqslant1,$
\begin{align*}
G_{X}\left(s\right) & =\sum^{+\infty}_{k=0}pq^{k}s^{k}=p\sum^{+\infty}_{k=0}\left(qs\right)^{k}=\dfrac{p}{1-qs}.
\end{align*}

\end{proof}

\subsubsection{Generating Function of the Geometric Law on $\mathbb{N}^{\ast}$}

\begin{proposition}{}{}

If $P_{X}=\mathcal{G}_{\mathbb{N}^{\ast}}\left(p\right),$ geometric
law with parameter $p\in\left]0,1\right[,$ then\boxeq{
\[
G_{X}\left(s\right)=\dfrac{ps}{1-qs}.
\]
}

\end{proposition}

\begin{proof}{}{}

Indeed,
\begin{align*}
G_{X}\left(s\right) & =\sum^{+\infty}_{k=1}pq^{k-1}s^{k}=ps\sum^{+\infty}_{k=1}\left(qs\right)^{k-1}=\dfrac{ps}{1-qs}.
\end{align*}

\end{proof}

\subsubsection{Generating Function of the Negative Binomial Law}

\begin{proposition}{Generating Function of the Sum of Two Random Variables}{generating_sum}

Let $X$ and $Y$ be two \textbf{independent} random variables defined
on a probabilized space $\left(\Omega,\mathcal{A},P\right)$ taking
values in $\mathbb{N},$ with generating functions $G_{X}$ and $G_{Y},$
respectively. 

Then for every $s\in\left[-1,1\right],$ the generating function of
$X+Y$ is given by\boxeq{
\[
G_{X+Y}\left(s\right)=G_{X}\left(s\right)G_{Y}\left(s\right).
\]
}

\end{proposition}

\begin{proof}{}{}

Since $X$ and $Y$ are independent random variables, the random variables
$s^{X}$ and $s^{Y}$ are also independent for every $s\in\left[-1,1\right].$ 

Therefore, 
\[
G_{X+Y}\left(s\right)=\mathbb{E}\left(s^{X}s^{Y}\right)=\mathbb{E}\left(s^{X}\right)\mathbb{E}\left(s^{Y}\right)=G_{X}\left(s\right)G_{Y}\left(s\right)
\]

\end{proof}

\begin{corollary}{Generating Function of the Negative Binomial Law}{}

Let $X$ be a random variable defined on the probabilized space $\left(\Omega,\mathcal{A},P\right),$
taking values in $\mathbb{N},$ and following the negative binomial
law $\mathcal{B}^{-}\left(n,p\right).$ 

Then, the generating function of $X,$ defined on $s\in\left[-1;1\right],$
is the function given by\boxeq{
\[
G_{X}\left(s\right)=\left(\dfrac{ps}{1-qs}\right)^{n},
\]
}with $q=1-p.$

\end{corollary}

\begin{proof}{}{}

Let $\left(X_{i}\right)_{1\leqslant i\leqslant n}$ be a family of
independent random variables, each following the geometric law $\mathcal{G}_{\mathbb{N}^{\ast}}\left(p\right).$ 

We have previously seen that the random variable $S_{n}=\sum^{n}_{i=1}X_{i}$
follows the negative binomial law $\mathcal{B}^{-}\left(n,p\right).$ 

By applying the previous proposition, extended to the case of $n$
independent random variables, we conclude.

\end{proof}

\subsection{Generating Function and Moments}

Since the generating function of a random variable determines its
law, it is natural that it also encodes its \textbf{moments}\index{moment},
when they exist.

We denote, for every $r\in\mathbb{N}^{\ast},$ $G^{\left(r\right)}_{X}\left(1^{-}\right)$
the $r-$th left-hand derivative of $G_{X}$ at 1, when it exists.

\begin{proposition}{Generating Function and Moments}{generating_moments}

Let $X$ be a random variable defined on a probabilized space $\left(\Omega,\mathcal{A},P\right),$
taking values on $\mathbb{N}.$

Then, for $X$ to admit a moment of order $r\in\mathbb{N}^{\ast},$
it is necessary and sufficient that its generating function $G_{X}$
is $r$ times differentiable from the left at $1.$ In that case,
we have\boxeq{
\begin{equation}
G^{\left(r\right)}_{X}\left(1^{-}\right)=\sum^{+\infty}_{k=r}k\left(k-1\right)\dots\left(k-r+1\right)p_{k}.\label{eq:r-th derivative G_X in 1-}
\end{equation}
}

This can also be written as
\begin{equation}
\mathbb{E}\left(X\left(X-1\right)\dots\left(X-r+1\right)\right)=G^{\left(r\right)}_{X}\left(1^{-}\right).\label{eq:Expectation order r and r-th derivative}
\end{equation}

In particular, for $r=1,$ we have
\begin{equation}
\mathbb{E}\left(X\right)=G^{\prime}_{X}\left(1^{-}\right).\label{eq:Expectation and 1st derivative G_X}
\end{equation}

\end{proposition}

\begin{proof}{}{}

We prove the results only for $r=1.$ The general case is left as
an exercise.

For every $s\in\left]0;1\right[$ and all $n\in\mathbb{N^{*}},$ we
have
\[
G_{X}\left(1\right)-G_{X}\left(s\right)=\sum^{+\infty}_{n=1}p_{n}\left(1-s^{n}\right).
\]

Since
\[
1-s^{n}=\left(1-s\right)\sum^{n-1}_{j=0}s^{j},
\]
we obtain
\begin{equation}
\dfrac{G_{X}\left(1\right)-G_{X}\left(s\right)}{1-s}=\sum^{+\infty}_{n=1}p_{n}\left(\sum^{n-1}_{j=0}s^{j}\right).\label{eq:variation_toll_G_X}
\end{equation}

Now, for every $n\in\mathbb{N}^{\ast}$ and every $s\in\left[0;1\right],$
\[
0\leqslant p_{n}\left(\sum^{n-1}_{j=0}s^{j}\right)\leqslant np_{n}.
\]

\begin{itemize}
\item If $X$ admits an expectation, then the series with general term $p_{n}\left(\sum^{n-1}_{j=0}s^{j}\right)$
is normally, thus uniformly convergent on $\left[0;1\right].$ It
follows that, by the equality $\refpar{eq:variation_toll_G_X},$ the
left-hand limit at 1 of $\dfrac{G_{X}\left(1\right)-G_{X}\left(s\right)}{1-s}$
exists and is equal to $\sum^{+\infty}_{n=1}np_{n}\equiv\mathbb{E}\left(X\right).$
\item Conversely, suppose that the left-hand limit at 1 of $\dfrac{G_{X}\left(1\right)-G_{X}\left(s\right)}{1-s}$
exists. Then, for every $N\in\mathbb{N}^{\ast},$ by the equality
$\refpar{eq:variation_toll_G_X},$
\begin{align*}
0 & \leqslant\sum^{N}_{n=1}np_{n}=\lim_{s\to1^{-}}\sum^{N}_{n=1}p_{n}\left(\sum^{n-1}_{j=0}s^{j}\right)\\
 & \leqslant\lim_{s\to1^{-}}\sum^{+\infty}_{n=1}p_{n}\left(\sum^{n-1}_{j=0}s^{j}\right)=G^{\prime}_{X}\left(1^{-}\right).
\end{align*}
This implies that the series with general nonnegative term $np_{n}$
is convergent and that the expectation of $X$ exists.
\end{itemize}
\end{proof}

\subsection{Sum of an Arbitrary Number of Random Variables}

We already encountered examples of situations where we need to study
the law of the sum of an arbitrary number of independent random variables.
The generating functions provide a powerful tool to tackle such problems.

\begin{proposition}{Sum of an Arbitrary Number of Random Variables}{sum_arbitrary}

Let $\left(X_{n}\right)_{n\in\mathbb{N}^{\ast}}$ be a sequence of
random variables of same non degenerated law, taking values in $\mathbb{N},$
and let $T$ be a random variable taking values in $\mathbb{N}^{\ast}$
such that $P\left(T=1\right)\neq1.$ 

Suppose the random variables $T,X_{n},n\in\mathbb{N}^{\ast},$ to
be independent. 

Define, for every $n\in\mathbb{N}^{\ast},$ the random variable 
\[
S_{n}=\sum^{n}_{j=1}X_{j}.
\]

Define the random variable $S$ for every $\omega\in\Omega$ by
\[
S\left(\omega\right)=S_{T\left(\omega\right)}\left(\omega\right).
\]

Let $G_{T}$ and $G_{X_{1}}$ denote the generating functions of $T$
and $X_{1}$ respectively. 

Then the generating function $G_{S}$ of $S$ yields by\boxeq{
\[
G_{S}=G_{T}\circ G_{X_{1}}.
\]
}

\end{proposition}

\begin{proof}{}{}

Since $T$ takes values in $\mathbb{N}^{\ast},$ we have, for every
$s\in\left[-1;1\right],$ 
\begin{align*}
G_{S}\left(s\right) & =\sum^{+\infty}_{k=0}s^{k}P\left(S=k\right)\\
 & =\sum^{+\infty}_{k=0}s^{k}P\left(\left(S=k\right)\cap\left(\biguplus_{n\in\mathbb{N}^{\ast}}\left(T=n\right)\right)\right)\\
 & =\sum_{k\in\mathbb{N}}s^{k}\left(\sum_{n\in\mathbb{N}^{\ast}}P\left(S=k,T=n\right)\right).
\end{align*}

We note that
\[
\sum_{k\in\mathbb{N}}\left|s\right|^{k}\left(\sum_{n\in\mathbb{N}^{\ast}}P\left(S=k,T=n\right)\right)\leqslant\sum_{k\in\mathbb{N}}\left(\sum_{n\in\mathbb{N}^{\ast}}P\left(S=k,T=n\right)\right)=1.
\]
so the family $\left(s^{k}P\left(S=k,T=n\right)\right)_{\left(k,n\right)\in\mathbb{N}\times\mathbb{N}^{\ast}}$
is summable.

By the Fubini theorem and the independence of the random variables
$S_{n}$ and $T,$ we can write
\begin{align*}
G_{S}\left(s\right) & =\sum_{n\in\mathbb{N}^{\ast}}\left(\sum_{k\in\mathbb{N}}s^{k}P\left(S_{n}=k,T=n\right)\right)\\
 & =\sum_{n\in\mathbb{N}^{\ast}}\left(\sum_{k\in\mathbb{N}}s^{k}P\left(S_{n}=k\right)\right)P\left(T=n\right)\\
 & =\sum_{n\in\mathbb{N}^{\ast}}G_{S_{n}}\left(s\right)P\left(T=n\right).
\end{align*}

Since the random variables $X_{n}$ are independent and follow the
same law, and thus have the same generating function, we deduce from
Proposition $\ref{pr:generating_sum}$ that
\begin{align*}
G_{S}\left(s\right) & =\sum^{+\infty}_{n=1}\left[G_{X_{1}}\left(s\right)\right]^{n}P\left(T=n\right)=G_{T}\left(G_{X_{1}}\left(s\right)\right).
\end{align*}

This proves the stated result.

\end{proof}

\begin{figure}[t]
\begin{center}\includegraphics[width=0.4\textwidth]{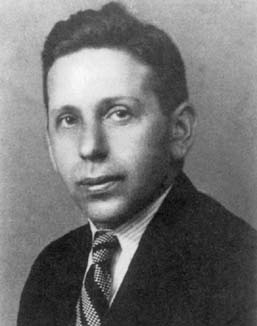}

{\tiny Credits: Public Domain}\end{center}

\caption{\textbf{\protect\href{https://en.wikipedia.org/wiki/Abraham_Wald}{Abraham Wald}}
(1902-1950)}\sindex[fam]{Wald, Abraham}
\end{figure}

In the same context, we give some relations, named \textbf{Wald}{\bfseries\footnote{\textbf{\sindex[fam]{Wald, Abraham}\href{https://en.wikipedia.org/wiki/Abraham_Wald}{Abraham Wald}}
(1902-1950) was a Hungarian mathematician and statistician. He worked
in different fiels such as decision theory, geometry and econometrics.
He also founded the field of sequential analysis. He obtained his
Ph.D. in mathematics in 1931 from the University of Vienna, but fled
Austria in 1938, follwing the discrimination against Jews after Nazis
annexed Austria. Wald emigrated to the USA where he worked on econometrics
reseach. He worked at Columbia University, until his death in a plane
crash in India.}}\textbf{\sindex[fam]{Wald, Abraham} identities}\index{Wald identities},
on the expectation and variance of $S.$

\begin{corollary}{Wald Identities}{}

Under the same hypotheses as in Proposition $\ref{pr:sum_arbitrary}$:

(i) If $X_{1}$ and $T$ have an expectation, then $S$ also admits
an expectation, given by\boxeq{
\begin{equation}
\mathbb{E}\left(S\right)=\mathbb{E}\left(X_{1}\right)\mathbb{E}\left(T\right).\label{eq:expectation}
\end{equation}
}

(ii) If $X_{1}$ and $T$ admit a second-order moment, then $S$ also
admits a second-order moment, and its variance is given by\boxeq{
\begin{equation}
\sigma^{2}_{S}=\sigma^{2}_{X_{1}}\mathbb{E}\left(T\right)+\left(\mathbb{E}\left(X_{1}\right)\right)^{2}\sigma^{2}_{T}.\label{eq:variance}
\end{equation}
}

\end{corollary}

\begin{proof}{}{}

(i) We cannot, in general, differentiate a composed function under
the left-hand limit (or the right one) without justification, as there
is no general theorem for that---one can look for counter-examples.
Therefore, even though the generating functions $G_{T}$ and $G_{X_{1}}$
admit a left-hand derivative in $1$---Proposition $\ref{pr:generating_moments}$---we
cannot directly assert $\refpar{eq:expectation}$ by writing directly
\[
G^{\prime}_{S}\left(1^{-}\right)=G^{\prime}_{T}\left(1^{-}\right)G^{\prime}_{X_{i}}\left(1^{-}\right).
\]

However, this equality can be justified as follow, since
\[
\begin{array}{ccccc}
G_{X_{1}}\left(1\right)=1 & \,\,\,\, & \text{and} & \,\,\,\, & G^{\prime}_{X_{1}}\left(1^{-}\right)>0\end{array}
\]
---this last, as $G^{\prime}_{X_{1}}\left(1^{-}\right)=\mathbb{E}\left(X_{1}\right)$
and $X_{1}$ is non-degenerated---there exists $s_{0}\in\left[0;1\right[$
such that, for every $s\in\left[s_{0};1\right[,$ we have $G_{X_{1}}\left(s\right)<1.$

For such $s\in\left[s_{0};1\right[,$ consider
\begin{align*}
\dfrac{G_{S}\left(1\right)-G_{S}\left(s\right)}{1-s} & =\dfrac{1-G_{T}\left(G_{X_{1}}\left(s\right)\right)}{1-s}\\
 & =\dfrac{1-G_{T}\left(G_{X_{1}}\left(s\right)\right)}{1-G_{X_{1}}\left(s\right)}\times\dfrac{1-G_{X_{1}}\left(s\right)}{1-s}.
\end{align*}

Since
\[
\lim_{s\to1}G_{X_{1}}\left(s\right)=1,
\]
we deduce that the limit of $\dfrac{G_{S}\left(1\right)-G_{S}\left(s\right)}{1-s}$
when $s$ tends to 1 exists, and consequentely that the left-hand
derivative in 1 of $G_{S}$ exists, and it is given by
\[
G^{\prime}_{S}\left(1^{-}\right)\equiv\lim_{s\to1}\dfrac{G_{S}\left(1\right)-G_{S}\left(s\right)}{1-s}=G^{\prime}_{T}\left(1^{-}\right)G^{\prime}_{X_{i}}\left(1^{-}\right).
\]

Proposition $\ref{pr:generating_moments}$ ensures that $S$ admits
an expectation given by the equality $\refpar{eq:expectation}.$

(ii) For every $s\in\left]-1;1\right[,$
\[
G^{\prime\prime}_{S}\left(s\right)=G^{\prime\prime}_{T}\left(G_{X_{1}}\left(s\right)\right)\left[G^{\prime\prime}_{X_{1}}\left(s\right)\right]^{2}+G^{\prime}_{T}\left(G_{X_{1}}\left(s\right)\right)G^{\prime\prime}_{X_{1}}\left(s\right).
\]

By Proposition $\ref{pr:generating_moments},$ we deduce, that, since
$\lim_{s\to1}G_{X_{1}}\left(s\right)=1,$ the left-hand limit of $G^{\prime\prime}_{S}\left(s\right)$
when $s$ tends to 1 exists and is given by
\begin{equation}
\lim_{s\to1}G^{\prime\prime}_{S}\left(s\right)=\mathbb{E}\left(T\left(T-1\right)\right)\left(\mathbb{E}\left(X_{1}\right)\right)^{2}+\mathbb{E}\left(T\right)\mathbb{E}\left(X_{1}\left(X_{1}-1\right)\right)\label{eq:lim_G_second}
\end{equation}

We now prove the existence of the expectation of $S\left(S-1\right).$
By applying Proposition $\ref{pr:generating_moments}$, for every
$s\in\left]-1;1\right[,$
\[
G^{\prime\prime}_{S}\left(s\right)=\sum^{+\infty}_{n=2}n\left(n-1\right)s^{n-2}P\left(S=n\right).
\]

Since $G^{\prime\prime}_{S}\left(s\right)$ being bounded in the neighborhood
of 1, there exists $c>0$ and $s_{1}\in\left[0;1\right[$ such that
\[
\forall s\in\left[s_{1};1\right[,\,\forall N>2,\,\,\,\,\,0\leqslant\sum^{N}_{n=2}n\left(n-1\right)s^{n-2}P\left(S=n\right)\leqslant c.
\]

By continuity, we then have, for every $N>2,$
\[
0\leqslant\sum^{N}_{n=2}n\left(n-1\right)P\left(S=n\right)\leqslant c
\]
which shows that the expectation of $S\left(S-1\right)$ exists and
that
\[
\mathbb{E}\left(S\left(S-1\right)\right)=\lim_{s\to1^{-}}G^{\prime\prime}_{S}\left(s\right).
\]

So we have, using the equality $\refpar{eq:lim_G_second},$
\[
\left(S\left(S-1\right)\right)=\mathbb{E}\left(T\left(T-1\right)\right)\left(\mathbb{E}\left(X_{1}\right)\right)^{2}+\mathbb{E}\left(T\right)\mathbb{E}\left(X_{1}\left(X_{1}-1\right)\right).
\]

An easy algebraical computation, taking into account of the relation
\[
\sigma^{2}_{S}=\mathbb{E}\left(S\left(S-1\right)\right)+\mathbb{E}\left(S\right)-\left(\mathbb{E}\left(S\right)\right)^{2}
\]
gives
\begin{equation}
\sigma^{2}_{S}=\sigma^{2}_{X_{1}}\mathbb{E}\left(T\right)+\left(\mathbb{E}\left(X_{1}\right)\right)^{2}\sigma^{2}_{T}.\label{eq:variance-1}
\end{equation}

\end{proof}

\begin{remark}{}{}

Wald moments can, of course, be proved by directly reasoning on the
moments, whithout using the generating functions---which is equivalent
to studying the law of $S.$ It is left as an exercise!

\end{remark}

\section*{Exercises}

\addcontentsline{toc}{section}{Exercises}

\begin{exercise}{}{exercise5.1}

A gambler plays heads and tails against a bank, using a coin that
is necessarily unfair---let $p$ be the probability of obtaining
tails in one toss. Their initial stake is $a\in\mathbb{N}^{\ast}.$
Until they obtain heads, the gambler loses what was at stake and decides
to bet $k>1$ times the previous stake in the next round . When the
gambler obtains tails, they win $k$ times their last stake and stop
gambling.

What is the gambler expected gain?

\end{exercise}

\begin{exercise}{Expectation and Moments of Higher Order}{exercise5.2}

Let $X$ and $Y$ be two independent discrete random variables defined
on a probabilized space $\left(\Omega,\mathcal{A},P\right),$ each
following the same uniform law on the set $\left\{ i\in\mathbb{N}:\,0\leqslant i\leqslant n\right\} .$
Define the random variables $Z$ and $T$ as
\[
\begin{array}{ccccc}
Z=\left|X-Y\right| & \,\,\,\, & \text{and} & \,\,\,\, & T=\min\left(X,Y\right).\end{array}
\]

1. This question concerns the transfer theorem and the linearity of
expectation.

(a) Justify the existence of moments of any order for $Z$ and $T.$

(b) Compute $\mathbb{E}\left(Z\right)$ without determining the law
of $Z,$ and find an asymptotic equivalent as $n$ tends to $+\infty.$

(c) Deduce $\mathbb{E}\left(T\right)$ from the previous question
and give an asymptotic equivalent as $n$ tends to $+\infty.$

\textit{Hint: Recall the equality, for every $\left(a,b\right)\in\mathbb{R}^{2},$
\[
\left|a-b\right|=a+b-2\min\left(a,b\right).
\]
}

2. Let $U$ be a random variable taking values in $\mathbb{N},$ such
that there exists $K\in\mathbb{N}^{\ast}$ with $0\leqslant U\leqslant K.$

(a) Express $\sum^{K}_{j=1}P\left(U\geqslant j\right)$ as a function
of $\mathbb{E}\left(U\right).$

(b) Similarly, compute $\sum^{K}_{j=1}j^{2}P\left(U\geqslant j\right)$
as a function of $\mathbb{E}\left(U\right),$ $\mathbb{E}\left(U^{2}\right)$
and $\mathbb{E}\left(U^{3}\right).$

3. When studying the law of the minimum of random variables, it is
often convenient to compute the probability that it is greater than
or equal to a given number.

(a) Compute, for every $j\in\mathbb{N},$ the probability $P\left(T\geqslant j\right).$

(b) Using question 2.a., find the value of $\mathbb{E}\left(T\right).$

4. Hint: Use the property of the linearity of expectation as much
as possible.

(a) Compute $\mathbb{E}\left(Z^{2}\right)$ as a function of the variance
$\sigma^{2}_{X}$ of the random variable $X.$

(b) What is the value of the variance $\sigma^{2}_{Z}$ of the random
variable $Z?$

\end{exercise}

\begin{exercise}{Law of a Maximum, Expectation and Higher Order Moments}{exercise5.3}

All the random variables are defined on a probabilized space $\left(\Omega,\mathcal{A},P\right).$

1. Let $U$ be a discrete random variable taking values in $\mathbb{N}$
admitting a second-order moment.

a. Prove the formula
\[
\mathbb{E}\left(U\right)=\sum_{j\geqslant1}P\left(U\geqslant j\right).
\]

b. Express the same sum
\[
\sum_{j\geqslant1}P\left(U\geqslant j\right)
\]
as a function of $\mathbb{E}\left(U^{2}\right)$ and $\mathbb{E}\left(U\right).$

2. Let $\left(X_{n}\right)_{n\in\mathbb{N}^{\ast}}$ be a sequence
of independent discrete random variables in $\mathbb{N},$ all following
the same law $\mu.$ 

Denote, for every $k\in\mathbb{N},$
\[
p_{k}=P\left(X_{n}=k\right).
\]

Let, for every $k\in\mathbb{N},$ $M_{n}$ be the random variable
such that
\[
M_{n}=\max_{1\leqslant i\leqslant n}\left(X_{i}\right).
\]

Compute, for every $k\in\mathbb{N},$ the probability $P\left(M_{n}\leqslant k\right)$
as a function of $F_{k}$ and $n.$

3. Suppose that $\mu$ is the uniform law on the set of integers $\left\llbracket 1,K\right\rrbracket $
where $K$ is an integer strictly greater than 1.

(a) Compute the probability $P\left(M_{n}=k\right)$ for every $k\in\left\llbracket 1,K\right\rrbracket .$

(b) Three fair dice are rolled. 

What is the probability that the maximum of the digit obtained is
4?

4. Now suppose $\mu$ is the geometric law on $\mathbb{N}^{\ast}$
with parameter $p\in\left]0;1\right[$ and let $q=1-p.$

(a) Compute the expectation $\mathbb{E}\left(M_{n}\right).$

(b) Three gamblers toss to heads and tails with a fair coin, and stop
as soon as they obtain tails. The random variable $M_{3}$ is the
number of toss that have been done by one or more players to obtain
tail at last.

Compute $\mathbb{E}\left(M_{3}\right).$ Provide the exact value,
then compute an approximate value rounded to the nearest hundredth.

\end{exercise}

\begin{exercise}{Expectation and Independence}{exercise5.4}

Let $X$ and Y be two independent discrete random variables defined
on the same probabilized space $\left(\Omega,\mathcal{A},P\right),$
both following the same geometric law on $\mathbb{N}^{\ast}$ with
parameter $p,\,p\in\left]0;1\right[.$ We denote $q=1-p.$ 

Define the random variables $T,$ $Z,$ and $G$ by
\[
\begin{array}{ccccc}
T=\min\left(X,Y\right), & \,\,\,\, & Z=\left|X-Y\right|, & \,\,\,\,\text{and}\,\,\,\, & G=\dfrac{Z}{T}.\end{array}
\]

1. The goal of this question is to study the law of the minimum of
two random variables.

(a) Compute, for every $x\in\mathbb{N}^{\ast},$ the probability $P\left(X\geqslant x\right).$

(b) Compute, for every $t\in\mathbb{N}^{\ast},$ the probability $P\left(T\geqslant t\right),$
and determine the law followed by $T.$

2. In this question, we compute expectations and study the law of
the pair $\left(T,Z\right).$

(a) Compute the expectations $\mathbb{E}\left(X\right)$ and $\mathbb{E}\left(\dfrac{1}{X}\right).$

(b) Compute, for every $\left(t,z\right)\in\mathbb{N}^{\ast}\times\mathbb{N},$
the probability $P\left(T\geqslant t,Z=z\right).$ Analyze the case
$z=0$ separately.

(c) Deduce the law of $Z.$

3. Show that the random variables $T$ and $Z$ are independent. That
is, for every $\left(t,z\right)\in\mathbb{N}^{\ast}\times\mathbb{N},$
we have
\[
P\left(T=t,Z=z\right)=P\left(T=t\right)P\left(Z=z\right).
\]

Note that the joint probability $P\left(T\geqslant t,Z\geqslant z\right)$
can be written in the form $f\left(t\right)g\left(z\right),$ where
$f$ and $g$ are functions defined on $\mathbb{N}^{\ast}$ and $\mathbb{N}$
respectively.

4. What is the value of the expectation $\mathbb{E}\left(G\right)?$

\end{exercise}

\begin{exercise}{Moments and Independence}{exercise5.5}

Let $U$ and $V$ be two independent discrete random variables on
a probabilized space $\left(\Omega,\mathcal{A},P\right)$ taking values
in $\mathbb{Z}$ and admitting second-order moments. Furthermore,
suppose $U$ is centered.

Define two discrete random variables $X$ and $Y$ by
\[
\begin{array}{ccc}
X=\left(-1\right)^{V}U & \,\,\,\, & Y=V.\end{array}
\]

1. In this part, we study questions of independence and correlation.

(a) Justify the existence of the expectation of $X$ and compute it.

(b) Justify the existence of the expectation of $XY$ and compute
it.

(c) What is the covariance of $X$ and $Y?$

(d) Are the random variables $X^{2}$ and $Y^{2}$ independent?

2. In this question, assume the law of $U$ is given by
\[
\begin{array}{ccccc}
P\left(U=-2\right)=\dfrac{1}{3} & \,\,\,\, & \text{and} & \,\,\,\, & P\left(U=1\right)=\dfrac{2}{3}\end{array},
\]
and the one of $V$ by
\[
\begin{array}{ccccc}
P\left(V=1\right)=\dfrac{1}{2} & \,\,\,\, & \text{and} & \,\,\,\, & P\left(V=2\right)=\dfrac{1}{2}\end{array}.
\]

(a) Compute $\mathbb{E}\left(X^{3}\right)$ and $\mathbb{E}\left(U^{3}\right).$

(b) Compute $\mathbb{E}\left(\boldsymbol{1}_{\left(V=1\right)}X^{3}\right).$

(c) Are the random variables $X$ and $Y$ independent?

3. In this question, we suppose that $U$ follows a law that is symmetric---i.e.
$U$ and $-U$ follow the same law---and that $P\left(U=0\right)=0.$

Let $f$ and $g$ be any two real-valued bounded functions defined
on $\mathbb{Z}.$

(a) Justify the existence of the expectation of the discrete random
variable $f\left(X\right)g\left(Y\right).$

(b) Show that
\[
\mathbb{E}\left(f\left(X\right)g\left(Y\right)\right)=\begin{cases}
\mathbb{E}\left(f\left(U\right)\right)\mathbb{E}\left(f\left(V\right)\right), & \text{if }f\text{ is even,}\\
0, & \text{if }f\text{ is odd.}
\end{cases}
\]

(c) Deduce from the previous question that, for every real-valued
bounded functions $f$ and $g$ defined on $\mathbb{Z},$ we have
\[
\mathbb{E}\left(f\left(X\right)g\left(Y\right)\right)=\mathbb{E}\left(f\left(X\right)\right)\mathbb{E}\left(g\left(Y\right)\right).
\]

(d) What can be said about the independence of the random variables
$X$ and $Y?$

\end{exercise}

\begin{exercise}{Chebyshev Inequality and Bernstein Polynomial}{exercise5.6}

Let $f$ be a continuous real-valued function on the closed interval
$\left[0;1\right].$

For every $n\in\mathbb{N}^{\ast},$ we denote $B_{n}$ the \textbf{\index{Bernstein polynomial}Bernstein\sindex[fam]{Bernstein, Sergei}
polynomial} defined by
\[
B_{n}\left(x\right)=\sum^{n}_{k=0}f\left(\dfrac{k}{n}\right)\binom{n}{k}x^{k}\left(1-x\right)^{n-k},
\]
with the convention
\[
0^{0}=1.
\]

Now, consider a probabilized space $\left(\Omega,\mathcal{A},P\right),$
and for each $x\in\left]0;1\right[$ a sequence of independent random
variables, following the same Bernoulli law with parameter $x.$

Define
\[
S_{n}=\sum^{n}_{k=1}X_{k}.
\]

1. Determine the expectation $\mathbb{E}\left(f\left(\dfrac{S_{n}}{n}\right)\right).$

2. For each $\epsilon>0,$ consider the real $\delta\left(\epsilon\right)$
defined by
\[
\delta\left(\epsilon\right)=\left\{ \left|f\left(x\right)-f\left(y\right)\right|\,:\,x,y\in\left[0,1\right]\,\text{and}\,\left|x-y\right|\leqslant\epsilon\right\} .
\]

a. Prove that $\delta\left(\epsilon\right)$ tends to $0$ with $\epsilon.$

b. Show that
\[
_{x\in\left[0,1\right]}\left|B_{n}\left(x\right)-f\left(x\right)\right|\leqslant\delta\left(\epsilon\right)+\dfrac{2\left\Vert f\right\Vert _{\infty}}{n\epsilon^{2}}.
\]

Deduce that the sequence of polynomials $B_{n}$ converges uniformly
to $f$ on $\left[0;1\right].$

\end{exercise}

\begin{exercise}{Expectation, Generating Function, Chebyshev Inequality}{exercise5.7}

Let $\left(X_{n}\right)_{n\in\mathbb{N}^{\ast}}$ be a sequence of
independent random variables defined on the probabilized space $\left(\Omega,\mathcal{A},P\right),$
all following the same geometric law on $\mathbb{N}^{\ast}$ with
parameter $p\in\left]0;1\right[.$

For $n\in\mathbb{N}^{\ast},$ define
\[
S_{n}=\sum^{n}_{k=1}X_{k}.
\]

1. Compute the expectation $\mathbb{E}\left(\dfrac{S_{n}}{n}\right)$
and the variance $\sigma^{2}_{\frac{S_{n}}{n}}$ of the random variable
$\dfrac{S_{n}}{n}.$

2. Compute the generating function $G_{S_{n}}$ of the random variable
$S_{n}$ for every $t\in\left[-1,1\right].$

3. Justify that $G_{S_{n}}$ can be expanded as a power series on
the interval $\left]-1;1\right[,$ and deduce for every $k\in\mathbb{N}^{\ast},$
the probability $P\left(S_{n}=k\right).$

Let $f$ be a real-valued function, uniformly continuous and bounded
on the half-line $\left[1;+\infty\right[.$ For each $n\in\mathbb{N}^{\ast},$
denote $B_{n}$ the function defined on $\left]0;1\right]$ by
\[
\forall x\in\left]0;1\right],\,\,\,\,B_{n}\left(x\right)=x^{n}\sum^{+\infty}_{k=0}f\left(1+\dfrac{k}{n}\right)\left(\begin{array}{c}
k+n-1\\
n-1
\end{array}\right)\left(1-x\right)^{k}.
\]

4. Compute the expectation $\mathbb{E}\left(f\left(\dfrac{S_{n}}{n}\right)\right)$
as a function of $B_{n}\left(p\right),$ thereby justifying a posteriori
its existence.

5. For every $\epsilon>0,$ let $\delta\left(\epsilon\right)$ be
the real number defined by
\[
\delta\left(\epsilon\right)=\left\{ \left|f\left(x\right)-f\left(y\right)\right|:\,\left|x\right|\geqslant1,\,\left|y\right|\geqslant1,\,\left|x-y\right|\leqslant\epsilon\right\} .
\]

(a) Prove that
\[
\left|\mathbb{E}\left(f\left(\dfrac{S_{n}}{n}\right)\right)-f\left(\dfrac{1}{p}\right)\right|\leqslant\delta\left(\epsilon\right)+2\left\Vert f\right\Vert _{\infty}P\left(\left|\dfrac{S_{n}}{n}-\dfrac{1}{p}\right|>\epsilon\right).
\]

(b) Deduce that, for every $a>1,$ we have
\[
\lim_{n\to+\infty}{}_{x\in\left[1,a\right]}\left|B_{n}\left(\dfrac{1}{x}\right)-f\left(x\right)\right|=0.
\]

\end{exercise}

\begin{exercise}{Generating Function and Moments}{exercise5.8}

The number $N$ of customers entering a mall during a day is a random
variable following the Poisson law with parameter $\lambda>0.$

The respective probabilities that a customer buys zero, one or two
articles of a brand $A$ are $\dfrac{1}{6},\dfrac{1}{2}$ and $\dfrac{1}{3},$
respectively. The total number of articles bought in a day is modelled
by a random variable, denoted $S.$ We study the law of $S$ thanks
to its generating function.

We refine the way the problem is modelled as follows: let $N$ be
a random variable on a probabilized space $\left(\Omega,\mathcal{A},P\right),$
following the Poisson law with parameter $\lambda>0.$ Let $\left(X_{n}\right)_{n\in\mathbb{N}^{\ast}}$
be a sequence of independent random variables, all following the same
law, where $X_{n}$ represents the number of articles bought by the
$n$-th client. The law is defined by
\[
P\left(X_{1}=0\right)=\dfrac{1}{6},\,\,\,\,P\left(X_{1}=1\right)=\dfrac{1}{2},\,\,\,\,P\left(X_{1}=2\right)=\dfrac{1}{3}.
\]

Moreover, suppose that $N$ and the $X_{n},\,n\in\mathbb{N}^{\ast},$
form a family of independent random variables. Finally, define the
random variable $S$ as
\[
S=\boldsymbol{1}_{\left(N\geqslant1\right)}\sum^{N}_{j=1}X_{j},
\]
with the convention
\[
\sum^{0}_{j=1}X_{j}=0.
\]

1. Compute the generating function $G_{S}$ of the random variable
$S$ for every $t\in\left[-1;1\right].$

2. Deduce the probability $P\left(S=3\right)$ and compute it numerically
for $\lambda=6.$

3. Justify the existence of the expectation and variance of the random
variable $S,$ and compute their respective values, denoting $\mathbb{E}\left(S\right)$
the expectation of $S$ and $\sigma^{2}_{S}$ the variance of $S.$
Compute them numerically in the case $\lambda=6.$

\end{exercise}

\begin{exercise}{Trinomial Law, Binomial Law, Generating Function and Independence. Poisson Law Caracterization.}{exercise5.9}

All the introduced random variables are defined on the same probabilized
space $\left(\Omega,\mathcal{A},P\right).$

Let $p,q,r$ be three positive real numbers such that $p+q+r=1.$

For every integer $n\geqslant1,$ consider the random variable $Y_{n}=\left(U_{n},V_{n}\right)$
with values in $\mathbb{N}^{2},$ and following the trinomial law
defined for every $\left(k,l\right)\in\mathbb{N}^{2}$ such that $0\leqslant k+l\leqslant n,$
by
\[
P\left(Y_{n}=\left(k,l\right)\right)=\dfrac{n!}{k!l!\left(n-\left(k+l\right)\right)!}p^{k}q^{l}r^{n-\left(k+l\right)}.
\]

Define $Y_{0}=\left(0,0\right).$ With the convention $0!=1,$ we
observe that $Y_{0}$ still satisfies the above formula.

1. Show that the integer-valued random variables $U_{n}$ and $V_{n}$
follow binomial laws. Determine their parameters.

2. Are the random variables $U_{n}$ and $V_{n}$ independent?

3. Prove the formula
\begin{equation}
\begin{array}{ccc}
\forall x,y\in\mathbb{R} & \,\,\,\, & \sum\limits^{n}_{k=1}k\left(\begin{array}{c}
n\\
k
\end{array}\right)x^{k-1}y^{n-k}=n\left(x+y\right)^{n-1},\end{array}\label{eq:formula_ex59_3}
\end{equation}
and deduce from it the expectation $\mathbb{E}\left(U_{n},V_{n}\right).$

4. Compute the covariance of the random variables $U_{n}$ and $V_{n}$
and the variance of the random variable $U_{n}+V_{n}.$

Let $N$ be a random variable taking values in $\mathbb{N},$ and
denote $a_{n}=P\left(N=n\right),$ for every $n\in\mathbb{N}.$ Suppose
that the family of random variables $\left\{ Y_{n},\,n\in\mathbb{N},\,N\right\} $
are independent. We define the random variables $U$ and $V$ for
every $\omega\in\Omega,$ by
\[
\begin{array}{ccccc}
U\left(\omega\right)=U_{N\left(\omega\right)}\left(\omega\right) & \,\,\,\, & \text{and} & \,\,\,\, & V\left(\omega\right)=V_{N\left(\omega\right)}\left(\omega\right)\end{array}.
\]

We denote $Y=\left(U,V\right).$

5. Suppose that $N$ follows a Poisson law with parameter $\lambda>0.$

Compute, for every $\left(k,l\right)\in\mathbb{N}^{2},$ the probability
$P\left(Y=\left(k,l\right)\right).$

Hint: Use the complete system of constituents $\left\{ \left(N=n\right)\right\} _{n\in\mathbb{N}}.$

Deduce from this that the random variables $U$ and $V$ are independent
and identify their laws.

6. Still assuming that $N$ follows a Poisson law with parameter $\lambda>0,$
we now aim to prove the independence of the random variables $U$
and $V,$ using generating functions.

The generating function of the random variable $Y$ is defined on
$\left[0;1\right]^{2}$ for every $\left(a,b\right)\in\left[0;1\right]^{2}$
by
\[
G_{Y}\left(a,b\right)=\mathbb{E}\left(a^{U}b^{V}\right).
\]

Compute the generating function of $Y,$ and deduce from it the generating
functions of the random variables $U$ and $V.$ From these results,
recover their laws. Use this to conclude that the random variables
$U$ and $V$ are independent.

7. Now suppose that $N$ has an arbitrary law, but that the random
variables $U$ and $V$ are independent. Prove, using the generating
function of $Y,$ that the random variable $N$ must follow a Poisson
law.

\end{exercise}

\section*{Solutions of Exercises}

\addcontentsline{toc}{section}{Solutions of Exercises}

\begin{solution}{}{solexercise5.1}

We consider the probabilized space $\left(\Omega,\mathcal{A},P\right).$
Let $\left(X_{n}\right)_{n\in\mathbb{N}^{\ast}}$ be a sequence of
independent random variables, each following the same Bernoulli law
$\mathcal{B}\left(1,p\right),$ where the event $\left(X_{n}=1\right)$
corresponds to the gambler obtaining a tail on the $n$-th toss.

Let $T$ be the random variable representing the index of the first
toss that results in a tail
\[
T=\inf\left(n\geqslant1:\,X_{n}=1\right),
\]
with the usual convention that $\inf\emptyset=+\infty.$ 

The random variable $G,$ representing the gambler algebraic gain,
is then defined as follows:
\begin{itemize}
\item If $\left(T=1\right),$ then $G=-a+ka.$
\item If $\left(T=n\right),$ with $n\geqslant2,$ then
\begin{align*}
G & =-a\left(1+k+\dots+k^{n-1}\right)+k^{n}a\\
 & =-a\dfrac{1-k^{n}}{1-k}+k^{n}a\\
 & =k^{n}a\left(\dfrac{2-k}{1-k}\right)-\dfrac{a}{1-k}.
\end{align*}
\end{itemize}
This last formula is also valid for $n=1,$ and we have $P\left(T<+\infty\right)=1.$

Define, for every $n\in\mathbb{N}^{\ast},$
\[
u_{n}=k^{n}a\left(\dfrac{2-k}{1-k}\right)-\dfrac{a}{1-k}.
\]

Let $\varphi$ be the function from $\mathbb{N}^{\ast}\cup\left\{ +\infty\right\} $
to $\mathbb{R},$ defined for every $t\in\mathbb{N}^{\ast}\cup\left\{ +\infty\right\} $
by
\[
\varphi\left(t\right)=\sum_{n\in\mathbb{N}^{\ast}}\boldsymbol{1}_{\left\{ n\right\} }\left(t\right)u_{n}.
\]

The gain can then be written as $G=\varphi\left(T\right).$

By the transfer theorem,
\[
\mathbb{E}\left(G\right)=\sum_{t\in\overline{\mathbb{N}^{\ast}}}\varphi\left(t\right)P\left(T=t\right).
\]

Since $P\left(T=+\infty\right)=0,$ it follows that
\[
\mathbb{E}\left(G\right)=\sum_{n\in\mathbb{N}^{\ast}}u_{n}P\left(T=n\right).
\]

Since $T$ follows a geometric law on $\mathbb{N}^{\ast}$ with parameter
$p,$ and letting $q=1-p,$ we get
\[
\mathbb{E}\left(G\right)=\sum_{n\in\mathbb{N}^{\ast}}u_{n}pq^{n-1}.
\]

It follows that
\begin{itemize}
\item If $qk\geqslant1$ and $k\neq2,$ then $\mathbb{E}\left(G\right)=+\infty,$
\item If $k=2,$ then $\mathbb{E}\left(G\right)=a,$
\item If $qk<1,$ recognizing a sum of a geometric law, then
\begin{align*}
\mathbb{E}\left(G\right) & =a\left(\dfrac{2-k}{1-k}\right)\sum^{+\infty}_{n=1}k^{n}pq^{n-1}-\dfrac{a}{1-k}\\
 & =a\left(\dfrac{2-k}{1-k}\right)\left(\dfrac{pk}{1-qk}\right)-\dfrac{a}{1-k}\\
 & =a\dfrac{pk-1}{1-qk}.
\end{align*}
\end{itemize}
\end{solution}

\begin{solution}{}{solexercise5.2}

\textbf{1. }Any nonnegative random variable dominated by a random
variable which admits an expectation also has an expectation.

\textbf{(a) Existence of the moments of all orders for $Z$ and $T$}

Since $X$ and $Y$ are nonnegative,
\[
\begin{array}{ccccc}
0\leqslant Z\leqslant X+Y & \,\,\,\, & \text{and} & \,\,\,\, & 0\leqslant T\leqslant X.\end{array}
\]

As $X$ and $Y$ admits an expectation, the same holds for $Z$ and
$T.$

As a remark, note that with probability 1, these random variables
are bounded and thus admit expectations.

\textbf{(b) Computation of $\mathbb{E}\left(Z\right)$}

Using successively the transfer theorem and the independence of $X$
and $Y,$ we obtain, defining $E=\left\llbracket 0,n\right\rrbracket ,$
\begin{align*}
\mathbb{E}\left(Z\right) & =\sum_{\left(x,y\right)\in E^{2}}\left|x-y\right|P\left(X=x,Y=y\right)\\
 & =\sum_{\left(x,y\right)\in E^{2}}\left|x-y\right|P\left(X=x\right)P\left(Y=y\right)\\
 & =\dfrac{1}{\left(n+1\right)^{2}}\sum_{\left(x,y\right)\in E^{2}}\left|x-y\right|\\
 & =\dfrac{2}{\left(n+1\right)^{2}}\sum^{n}_{i=1}\left(\sum^{i-1}_{j=0}\left(i-j\right)\right)\\
 & =\dfrac{2}{\left(n+1\right)^{2}}\sum^{n}_{i=1}\left(i^{2}-\sum^{i-1}_{j=0}j\right)\\
 & =\dfrac{2}{\left(n+1\right)^{2}}\sum^{n}_{i=1}\left(i^{2}-\dfrac{\left(i-1\right)i}{2}\right)\\
 & =\dfrac{1}{\left(n+1\right)^{2}}\sum^{n}_{i=1}\left(i^{2}+i\right).
\end{align*}

Recalling that
\[
\sum^{n}_{i=1}i^{2}=\dfrac{n\left(n+1\right)\left(2n+1\right)}{6},
\]
we find after simplification
\[
\mathbb{E}\left(Z\right)=\dfrac{n\left(n+2\right)}{3\left(n+1\right)}.
\]

Hence, when $n$ tends to the infinity, we have
\[
\mathbb{E}\left(Z\right)\sim\dfrac{n}{3}.
\]

\textbf{(c) Computation of $\mathbb{E}\left(T\right)$ and asymptotic
behavior as $n$ tends to $+\infty.$}

We have
\[
T=\dfrac{1}{2}\left(X+Y-\left|X-Y\right|\right),
\]
and therefore
\[
\mathbb{E}\left(T\right)=\dfrac{1}{2}\left(\mathbb{E}\left(X\right)+\mathbb{E}\left(Y\right)-\mathbb{E}\left(Z\right)\right).
\]

Since
\[
\mathbb{E}\left(X\right)=\mathbb{E}\left(Y\right)=\sum^{n}_{j=0}\dfrac{j}{n+1}=\dfrac{n}{2},
\]
we get
\[
\mathbb{E}\left(T\right)=\dfrac{n\left(2n+1\right)}{6\left(n+1\right)}.
\]
As $n$ tends to infinity
\[
\mathbb{E}\left(T\right)\sim\dfrac{n}{3}.
\]

\textbf{2. (a) Expression of $\sum^{K}_{j=1}P\left(U\geqslant j\right)$
in terms of $\mathbb{E}\left(U\right).$}
\begin{align*}
\sum^{K}_{j=1}P\left(U\geqslant j\right) & =\sum^{K}_{j=1}\mathbb{E}\left(\boldsymbol{1}_{\left(U\geqslant j\right)}\right)=\mathbb{E}\left(\sum^{K}_{j=1}\boldsymbol{1}_{\left(U\geqslant j\right)}\right)=\mathbb{E}\left(\sum^{\min\left(K,U\right)}_{j=1}1\right).
\end{align*}

Since $0\leqslant U\leqslant K,$ this simplifies to
\[
\sum^{K}_{j=1}P\left(U\geqslant j\right)=\mathbb{E}\left(U\right).
\]

This formula connecting the expectation and the probabilities $P\left(U\geqslant j\right)$
are in fact very general and useful.

\textbf{(b) Computation of $\sum^{K}_{j=1}j^{2}P\left(U\geqslant j\right)$
in terms of $\mathbb{E}\left(U\right),$ $\mathbb{E}\left(U^{2}\right)$
and $\mathbb{E}\left(U^{3}\right).$}

Similarly,
\begin{align*}
\sum^{K}_{j=1}j^{2}P\left(U\geqslant j\right) & =\sum^{K}_{j=1}j^{2}\mathbb{E}\left(\boldsymbol{1}_{\left(U\geqslant j\right)}\right)\\
 & =\mathbb{E}\left(\sum^{K}_{j=1}j^{2}\boldsymbol{1}_{\left(U\geqslant j\right)}\right)\\
 & =\mathbb{E}\left(\sum^{\min\left(K,U\right)}_{j=1}j^{2}\right)\\
 & =\mathbb{E}\left(\sum^{U}_{j=1}j^{2}\right)\\
 & =\mathbb{E}\left(\dfrac{U\left(U+1\right)\left(2U+1\right)}{6}\right).
\end{align*}

Thus,
\[
\sum^{K}_{j=1}j^{2}P\left(U\geqslant j\right)=\dfrac{1}{3}\mathbb{E}\left(U^{3}\right)+\dfrac{1}{2}\mathbb{E}\left(U^{2}\right)+\dfrac{1}{6}\mathbb{E}\left(U\right).
\]

\textbf{3. }Saying that the minimum of two numbers is greater than
or equal to a third, is equivalent to saying that both numbers are
greater or equal to the third.

\textbf{(a) Computation of the probability $P\left(T\geqslant j\right)$
for every $j\in\mathbb{N}.$ }

Using the independence of $X$ and $Y,$ we have, for every $j\in\mathbb{N},$
\[
P\left(T\geqslant j\right)=P\left(X\geqslant j,Y\geqslant j\right)=P\left(X\geqslant j\right)P\left(Y\geqslant j\right).
\]

So
\[
P\left(T\geqslant j\right)=\begin{cases}
0, & \text{if }j>n,\\
\left(\sum^{n}_{k=j}P\left(X=k\right)\right)^{2}, & \text{if }0\leqslant j\leqslant n,
\end{cases}
\]
which is equivalent to

\[
P\left(T\geqslant j\right)=\begin{cases}
\left(\dfrac{n-j+1}{n+1}\right)^{2}, & \text{if }0\leqslant j\leqslant n,\\
0, & \text{if }j>n.
\end{cases}
\]

\textbf{(b) Value of $\mathbb{E}\left(T\right).$}

From the previous results, since $0\leqslant T\leqslant n,$
\begin{align*}
\mathbb{E}\left(T\right) & =\sum^{n}_{j=1}P\left(T\geqslant j\right)=\sum^{n}_{j=1}\left(\dfrac{n-j+1}{n+1}\right)^{2}=\dfrac{n\left(2n+1\right)}{6\left(n+1\right)}.
\end{align*}

\textbf{4. }We must make optimal use of the Leibniz formula.

\textbf{(a) Computation of $\mathbb{E}\left(Z^{2}\right)$ in terms
of the variance $\sigma^{2}_{X}$ of the random variable $X.$}

We have
\[
\mathbb{E}\left(Z^{2}\right)=\mathbb{E}\left(\left(X-Y\right)^{2}\right)=\mathbb{E}\left(X^{2}\right)+\mathbb{E}\left(Y^{2}\right)-2\mathbb{E}\left(XY\right).
\]

Since $X$ and $Y$ follow the same law, and have in particular the
same moment, and are independent
\begin{align*}
\mathbb{E}\left(Z^{2}\right) & =2\mathbb{E}\left(X^{2}\right)-2\left(\mathbb{E}\left(X\right)\right)^{2}\\
 & =2\sigma^{2}_{X}.
\end{align*}

\textbf{(b) Computation of the variance $\sigma^{2}_{Z}$ of the random
variable $Z$}

The variance of $Z$ is given by
\[
\sigma^{2}_{Z}=\mathbb{E}\left(Z^{2}\right)-\left(\mathbb{E}\left(Z\right)\right)^{2},
\]
so
\begin{equation}
\sigma^{2}_{Z}=2\sigma^{2}_{X}-\left(\mathbb{E}\left(Z\right)\right)^{2}.\label{eq:sigma_z_ex2_4.b}
\end{equation}

Since
\[
\mathbb{E}\left(X^{2}\right)=\sum^{n}_{j=0}\dfrac{j^{2}}{n+1}=\dfrac{n\left(2n+1\right)}{6}
\]
and
\begin{align*}
\sigma^{2}_{X} & =\mathbb{E}\left(X^{2}\right)-\left(\mathbb{E}\left(X\right)\right)^{2}\\
 & =\dfrac{n\left(2n+1\right)}{6}-\left(\dfrac{1}{n+1}\times\dfrac{n\left(n+1\right)}{2}\right)^{2}\\
 & =\dfrac{n\left(n+2\right)}{12}
\end{align*}

By substituting into the equality $\refpar{eq:sigma_z_ex2_4.b},$
it yields
\[
\sigma^{2}_{Z}=\dfrac{n\left(n+2\right)}{6}-\dfrac{n^{2}\left(n+2\right)^{2}}{9\left(n+1\right)^{2}}.
\]

After simplification, we obtain
\[
\sigma^{2}_{Z}=\dfrac{n\left(n+2\right)\left(n^{2}+2n+3\right)}{18\left(n+1\right)^{2}}.
\]

We can note that, as $n$ tends to $+\infty,$
\[
\begin{array}{ccccc}
\sigma^{2}_{X}\sim\dfrac{n^{2}}{12} & \,\,\,\, & \text{and} & \,\,\,\, & \sigma^{2}_{Z}\sim\dfrac{n^{2}}{18}.\end{array}
\]

\end{solution}

\begin{solution}{}{solexercise5.3}

\textbf{1.} The question is similar to the previous exercise. 

Here, we partition the set $\left(U\geqslant j\right)$ into an infinite
sequence of sets and we use the Fubini property for nonnegative families.

\textbf{a. Proof that $\mathbb{E}\left(U\right)=\sum_{j\geqslant1}P\left(U\geqslant j\right).$}

Since
\[
\left(U\geqslant j\right)=\biguplus_{k\geqslant j}\left(U=k\right),
\]
it follows from the $\sigma-$additivity of $P,$ the positivity of
the terms, and the Fubini property that we have successively
\begin{align*}
\sum_{j\geqslant1}P\left(U\geqslant j\right) & =\sum_{j\geqslant1}P\left(\biguplus_{k\geqslant j}\left(U=k\right)\right)\\
 & =\sum_{j\geqslant1}\left(\sum_{k\geqslant j}P\left(U=k\right)\right)\\
 & =\sum_{k\geqslant1}\left(\sum^{k}_{j=1}P\left(U=k\right)\right)\\
 & =\sum_{k\geqslant1}kP\left(U=k\right)\\
 & =\mathbb{E}\left(U\right).
\end{align*}

\textbf{b. Computation of the sum, $\sum_{j\geqslant1}P\left(U\geqslant j\right)$}

Similarly,
\begin{align*}
\sum_{j\geqslant1}jP\left(U\geqslant j\right) & =\sum_{j\geqslant1}jP\left(\biguplus_{k\geqslant j}\left(U=k\right)\right)\\
 & =\sum_{j\geqslant1}\left(\sum_{k\geqslant j}jP\left(U=k\right)\right)\\
 & =\sum_{k\geqslant1}\left(\sum^{k}_{j=1}jP\left(U=k\right)\right)\\
 & =\sum_{k\geqslant1}\dfrac{k\left(k+1\right)}{2}P\left(U=k\right)\\
 & =\dfrac{1}{2}\left[\mathbb{E}\left(U^{2}\right)+\mathbb{E}\left(U\right)\right].
\end{align*}

\textbf{2. Computation, for every $k\in\mathbb{N},$ of the probability
$P\left(M_{n}\leqslant k\right)$ in terms of $F_{k}$ and $n.$}

Maximum and cumulative distribution functions are concepts relative
to order structure. It is therefore natural to access to the law of
the maximum of random variables in terms of the cumulative distribution
function.

We have
\[
\left(M_{n}\leqslant k\right)=\bigcap^{n}_{i=1}\left(X_{i}\leqslant k\right).
\]

The random variables $X_{i}$ are independent and follow the same
law, so for every $k\in\mathbb{N},$
\[
P\left(M_{n}\leqslant k\right)=\prod^{n}_{i=1}P\left(X_{i}\leqslant k\right)=\left(F_{k}\right)^{n}.
\]

\textbf{3. (a) Computation of $P\left(M_{n}=k\right)$ for every $k\in\left\llbracket 1,K\right\rrbracket .$}

If $\mu$ is the uniform law on $\left\llbracket 1,K\right\rrbracket ,$
then for every such $k\in\left\llbracket 1,K\right\rrbracket ,$ 
\[
p_{k}=\dfrac{1}{K}.
\]

Moreover, $F_{0}=P\left(X_{n}=0\right)=0.$

For every $k\in\left\llbracket 1,K\right\rrbracket ,$ it follows
\[
F_{k}=P\left(X_{n}=k\right)=\dfrac{k}{K},
\]
and
\begin{align*}
P\left(M_{n}=k\right) & =P\left(M_{n}\leqslant k\right)-P\left(M_{n}\leqslant k-1\right)\\
 & =\left(F_{k}\right)^{n}-\left(F_{k-1}\right)^{n}\\
 & =\dfrac{1}{K^{n}}\left(k^{n}-\left(k-1\right)^{n}\right)
\end{align*}

\textbf{(b) Computation of the probability that the maximum digit
obtained is 4 when rolling three fair dice.}

Since the dice are fair, $\mu$ is the uniform law on $\left\{ 1,2,\dots,6\right\} .$ 

We compute $P\left(M_{3}=4\right):$
\[
P\left(M_{3}=4\right)=\dfrac{1}{6^{3}}\left(4^{3}-3^{3}\right)=\dfrac{37}{216}\approx0.171.
\]

\textbf{4.} We compute the expectation of the maximum of independent
random variables following a geometric law and illustrate it with
a game involving three players.

\textbf{(a) Compute the expectation $\mathbb{E}\left(M_{n}\right).$}

Since $\mu$ is the geometric law on $\mathbb{N}^{\ast}$ with parameter
$p,$
\[
F_{k}=\sum^{k}_{j=1}pq^{j-1}=p\dfrac{1-q^{k}}{1-q}=1-q^{k}.
\]

Since $P\left(M_{n}\leqslant k\right)=\left(F_{k}\right)^{n},$ then
\[
P\left(M_{n}\leqslant k\right)=\left(1-q^{k}\right)^{n}.
\]

Therefore
\[
P\left(M_{n}\geqslant1\right)=1,
\]
and, for $k\geqslant2,$
\begin{align*}
P\left(M_{n}\geqslant k\right) & =1-P\left(M_{n}\leqslant k-1\right)\\
 & =1-\left(1-q^{k-1}\right)^{n}\\
 & =1-\sum^{n}_{i=0}\left(\begin{array}{c}
n\\
i
\end{array}\right)\left(-1\right)^{n-i}q^{\left(n-i\right)\left(k-1\right)}.
\end{align*}

Reindexing the sum using $j=n-i,$ and noting $\left(\begin{array}{c}
n\\
n-j
\end{array}\right)=\left(\begin{array}{c}
n\\
j
\end{array}\right),$ we obtain
\begin{align*}
P\left(M_{n}\geqslant k\right) & =1-\sum^{n}_{j=0}\left(\begin{array}{c}
n\\
j
\end{array}\right)\left(-1\right)^{j}q^{j\left(k-1\right)}\\
 & =\sum^{n}_{j=1}\left(\begin{array}{c}
n\\
j
\end{array}\right)\left(-1\right)^{j+1}q^{j\left(k-1\right)}.
\end{align*}

Thus, using the first question, we have
\begin{align*}
\mathbb{E}\left(M_{n}\right) & =1+\sum_{k\geqslant2}\left(\sum^{n}_{j=1}\left(\begin{array}{c}
n\\
j
\end{array}\right)\left(-1\right)^{j+1}q^{j\left(k-1\right)}\right)\\
 & =1+\sum^{n}_{j=1}\left(\begin{array}{c}
n\\
j
\end{array}\right)\left(-1\right)^{j+1}\sum_{k\geqslant2}q^{j\left(k-1\right)}\\
 & =1+\sum^{n}_{j=1}\left(\begin{array}{c}
n\\
j
\end{array}\right)\left(-1\right)^{j+1}\dfrac{q^{j}}{1-q^{j}}.
\end{align*}

\textbf{(b) Computation of $\mathbb{E}\left(M_{3}\right).$}

In this case
\[
\begin{array}{ccccc}
p=q=\dfrac{1}{2} & \,\,\,\, & \text{and} & \,\,\,\, & n=3,\end{array}
\]
which gives
\begin{align*}
\mathbb{E}\left(M_{3}\right) & =1+\left(\begin{array}{c}
3\\
1
\end{array}\right)\dfrac{\dfrac{1}{2}}{1-\dfrac{1}{2}}-\left(\begin{array}{c}
3\\
2
\end{array}\right)\dfrac{\dfrac{1}{4}}{1-\dfrac{1}{4}}+\left(\begin{array}{c}
3\\
3
\end{array}\right)\dfrac{\dfrac{1}{8}}{1-\dfrac{1}{8}}=\dfrac{22}{7}\approx3.14.
\end{align*}

\end{solution}

\begin{solution}{}{solexercise5.4}

\textbf{1.} We already note that the law of a discrete random variable
$X$ taking values in $\mathbb{N}$ is determined by the probabilities
$P\left(X\geqslant x\right).$

\textbf{(a) Computation of $P\left(X\geqslant x\right)$}

For every $x\in\mathbb{N}^{\ast},$ since $0<q<1$ and $q=1-p,$
\begin{align*}
P\left(X\geqslant x\right) & =\sum^{+\infty}_{k=x}pq^{k-1}=p\dfrac{q^{x-1}}{1-q}=q^{x-1}.
\end{align*}

Since $Y$ follows the same law, we have
\[
P\left(Y\geqslant y\right)=q^{y-1}.
\]

\textbf{(b) Computation of $P\left(T\geqslant t\right),$ and law
followed by $T$}

The random variables $X$ and $Y$ are independent. Thus, for every
$x\in\mathbb{N}^{\ast},$ we have
\begin{align*}
P\left(T\geqslant t\right) & =P\left(X\geqslant t,Y\geqslant t\right)=P\left(X\geqslant t\right)P\left(Y\geqslant t\right)=q^{2\left(t-1\right)}.
\end{align*}

Also
\begin{align*}
P\left(T=t\right) & =P\left(T\geqslant t\right)-P\left(T\geqslant t+1\right)=q^{2\left(t-1\right)}\left(1-q^{2}\right).
\end{align*}

Thus, $T$ follows the geometric law on $\mathbb{N}^{\ast}$ with
parameter $1-q^{2}.$

\textbf{2. (a) Computation of the expectations $\mathbb{E}\left(X\right)$
and $\mathbb{E}\left(\dfrac{1}{X}\right)$}

Since $X$ follows the geometric law on $\mathbb{N}^{\ast}$ with
parameter $p,$ we recall that
\[
\mathbb{E}\left(X\right)=\dfrac{1}{p}.
\]

By the transfer theorem, 
\begin{align*}
\mathbb{E}\left(\dfrac{1}{X}\right) & =\sum^{+\infty}_{k=1}\dfrac{1}{k}P\left(X=k\right)=\sum^{+\infty}_{k=1}\dfrac{1}{k}pq^{k-1}=p\sum^{+\infty}_{k=0}\dfrac{1}{k+1}q^{k}.
\end{align*}

Since, for $0<x<1,$
\[
\sum^{+\infty}_{k=0}\dfrac{x^{k+1}}{k+1}=-\ln\left(1-x\right).
\]
we obtain
\[
\mathbb{E}\left(\dfrac{1}{X}\right)=-\dfrac{p}{q}\ln p.
\]

\textbf{(b) Computation of $P\left(T\geqslant t,Z=z\right)$ }

For every $\left(t,z\right)\in\mathbb{N}^{\ast}\times\mathbb{N},$
\[
P\left(T\geqslant t,Z=z\right)=P\left(X\geqslant t,Y\geqslant t,\left|X-Y\right|=z\right).
\]

Denoting
\[
D_{t,z}=\left\{ \left(x,y\right)\in\mathbb{N}^{\ast}\times\mathbb{N}^{\ast}:\,x\geqslant t,\,y\geqslant t,\,\left|x-y\right|=z\right\} ,
\]
it follows
\[
P\left(T\geqslant t,Z=z\right)=\sum_{\left(x,y\right)\in D_{t,z}}p^{2}q^{x+y+2}.
\]

We have, in the case where $z=0,$
\[
D_{t,0}=\left\{ \left(x,y\right)\in\mathbb{N}^{\ast}\times\mathbb{N}^{\ast}:\,x\geqslant t,\,x=y\right\} ,
\]
which gives
\[
P\left(T\geqslant t,Z=0\right)=\sum_{x\geqslant t}p^{2}q^{2\left(x-1\right)}.
\]
Thus, after computation, for every $t\in\mathbb{N}^{\ast},$
\[
P\left(T\geqslant t,Z=0\right)=\dfrac{p}{1+q}q^{2\left(t-1\right)}.
\]

For every $z\in\mathbb{N}^{\ast},$
\begin{align*}
D_{t,z} & =\left\{ \left(x,y\right)\in\mathbb{N}^{\ast}\times\mathbb{N}^{\ast}:\,x\geqslant t,\,y=x+z\right\} \\
 & \qquad\qquad\biguplus\left\{ \left(x,y\right)\in\mathbb{N}^{\ast}\times\mathbb{N}^{\ast}:\,y\geqslant t,\,x=y+z\right\} ,
\end{align*}
which gives
\begin{align*}
P\left(T\geqslant t,Z=z\right) & =2p^{2}\sum_{y\geqslant t}q^{2\left(y-1\right)+z}=2p^{2}\dfrac{q^{2\left(t-1\right)+z}}{1-q^{2}}.
\end{align*}

Thus, for every $\left(t,z\right)\in\mathbb{N}^{*}\times\mathbb{N},$
\[
P\left(T\geqslant t,Z=z\right)=\dfrac{2p}{1+q}q^{2\left(t-1\right)+z}.
\]

\textbf{(c) Law of $Z$}

For every $z\in\mathbb{N}^{\ast},$
\[
P\left(Z=z\right)=P\left(T\geqslant1,Z=z\right)=\dfrac{2p}{1+q}q^{z},
\]
 and
\[
P\left(Z=0\right)=P\left(T\geqslant1,Z=0\right)=\dfrac{p}{1+q}.
\]

\textbf{3. $T$ and $Z$ are independent}

From the previous question, there exists nonnegative functions $f$
and $g$ such that for every $\left(t,z\right)\in\mathbb{N}^{\ast}\times\mathbb{N},$
\[
P\left(T\geqslant t,Z=z\right)=f\left(t\right)g\left(z\right).
\]

Thus,
\begin{align*}
P\left(T\geqslant t\right) & =\sum_{z\in\mathbb{N}}P\left(T\geqslant t,Z=z\right)=f\left(t\right)\sum_{z\in\mathbb{N}}g\left(z\right).
\end{align*}

By denoting 
\[
K=\sum_{z\in\mathbb{N}}g\left(z\right),
\]
 we have for every $t\in\mathbb{N}^{\ast},$
\[
P\left(T\geqslant t\right)=Kf\left(t\right).
\]

Taking $t=1,$ we find
\[
K=\dfrac{1}{f\left(1\right)}.
\]

Moreover, for every $z\in\mathbb{N},$
\begin{align*}
P\left(Z=z\right) & =P\left(T\geqslant1,Z=z\right)=f\left(1\right)g\left(z\right).
\end{align*}

Combining these results, we obtain for every $\left(t,z\right)\in\mathbb{N}^{\ast}\times\mathbb{N},$
\[
P\left(T\geqslant t\right)P\left(Z=z\right)=f\left(t\right)g\left(z\right).
\]

Since for every $\left(t,z\right)\in\mathbb{N}^{\ast}\times\mathbb{N},$
\[
P\left(T=t,Z=z\right)=P\left(T\geqslant t,Z=z\right)-P\left(T\geqslant t+1,Z=z\right)
\]
it follows
\begin{align*}
P\left(T=t,Z=z\right) & =f\left(t\right)g\left(z\right)-f\left(t+1\right)g\left(z\right)\\
 & =\left[f\left(t\right)-f\left(t+1\right)\right]g\left(z\right)\\
 & =\dfrac{1}{K}\left[P\left(T\geqslant t\right)-P\left(T\geqslant t+1\right)\right]g\left(z\right)\\
 & =P\left(T=t\right)P\left(Z=z\right).
\end{align*}

This proves that the random variables $T$ and $Z$ are independent.

\textbf{4. Computation of $\mathbb{E}\left(G\right)$}

Since $T$ and $Z$ are independent, $Z$ and $\dfrac{1}{T}$ are
independent, and we have
\[
\mathbb{E}\left(G\right)=\mathbb{E}\left(Z\right)\mathbb{E}\left(\dfrac{1}{T}\right).
\]

The random variable $T$ follows the geometric law on $\mathbb{N}^{\ast}$
with parameter $1-q^{2}.$ By using Question 2, we have
\begin{align*}
\mathbb{E}\left(\dfrac{1}{T}\right) & =-\dfrac{1-q^{2}}{q^{2}}\ln\left(1-q^{2}\right)=-\dfrac{\left(1+q\right)p}{q^{2}}\ln\left(p\left(1+q\right)\right).
\end{align*}

Moreover, from the previous question,
\begin{align*}
\mathbb{E}\left(Z\right) & =\sum_{z\geqslant1}zP\left(Z=z\right)=\dfrac{2p}{1+q}\sum_{z\geqslant1}zq^{z}=\dfrac{2q}{1+q}\sum_{z\geqslant1}zpq^{z-1}.
\end{align*}

This last sum can be identified to the expectation of a random variable
following the geometric law on $\mathbb{N}^{\ast}$ with parameter
$p,$ so
\[
\mathbb{E}\left(Z\right)=\dfrac{2q}{\left(1+q\right)p}.
\]

Thus
\begin{align*}
\mathbb{E}\left(G\right) & =\mathbb{E}\left(Z\right)\mathbb{E}\left(\dfrac{1}{T}\right)=-\dfrac{2}{q}\ln\left(p\left(1+q\right)\right).
\end{align*}

\end{solution}

\begin{solution}{}{solexercise5.5}

\textbf{1. }It is natural to expect that, in general, the random variables
$X$ and $Y$ are not independent.

\textbf{(a) Existence and computation of the expectation of $X$}

Since $X=\left(-1\right)^{V}U,$ we have $\left|X\right|\leqslant\left|U\right|.$ 

As the random variable $U$ admits an expectation, so does $X.$ 

Since $U$ and $V$ are independent random variables, $\left(-1\right)^{V}$
and $U$ are independent as well.

Thus
\[
\mathbb{E}\left(X\right)\equiv\mathbb{E}\left(\left(-1\right)^{V}U\right)=\mathbb{E}\left(\left(-1\right)^{V}\right)\mathbb{E}\left(U\right).
\]

Since $U$ is centered, $\mathbb{E}\left(U\right)=0,$ so $\mathbb{E}\left(X\right)=0.$

\textbf{(b) Existence and computation of the expectation of $XY$}

As the random variables $U$ and $V$ have a second order moment,
the product $UV$ also admits an expectation. Conjugated with the
fact that
\[
\left|XY\right|=\left|\left(-1\right)^{V}UV\right|\leqslant\left|UV\right|,
\]
the random variable $XY$ also has an expectation.

Since $U$ and $V$ are independent random variables, so are $\left(-1\right)^{V}U$
and $V,$ and
\[
\mathbb{E}\left(XY\right)=\mathbb{E}\left(\left(-1\right)^{V}V\right)\mathbb{E}\left(U\right),
\]
and thus
\[
\mathbb{E}\left(XY\right)=0.
\]

\textbf{(c) Covariance of $X$ and $Y$}
\[
\text{cov}\left(X,Y\right)=\mathbb{E}\left(XY\right)-\mathbb{E}\left(X\right)\mathbb{E}\left(Y\right)=0.
\]

\textbf{(d) Independence of $X^{2}$ and $Y^{2}$}

We have
\[
\begin{array}{ccccc}
X^{2}=U^{2} & \,\,\,\, & \text{and} & \,\,\,\, & Y^{2}=V^{2}.\end{array}
\]

Since $U$ and $V$ are independent random variables, so are $U^{2}$
and $V^{2},$ which implies that the random variables $X^{2}$ and
$Y^{2}$ are independent.

\textbf{2.} We now want to show that $X$ and $Y$ are not independent
random variables in this case.

\textbf{(a) Computation of $\mathbb{E}\left(X^{3}\right)$ and $\mathbb{E}\left(U^{3}\right)$}

Since
\[
\mathbb{E}\left(U\right)=\left(-2\right)\times\dfrac{1}{3}+1\times\dfrac{2}{3}=0,
\]
we observe that U is well centered.

We have
\begin{align*}
\mathbb{E}\left(X^{3}\right) & =\mathbb{E}\left(\boldsymbol{1}_{\left(V=1\right)}X^{3}\right)+\mathbb{E}\left(\boldsymbol{1}_{\left(V=2\right)}X^{3}\right)\\
 & =\mathbb{E}\left(\boldsymbol{1}_{\left(V=1\right)}\left(-1\right)^{3}U^{3}\right)+\mathbb{E}\left(\boldsymbol{1}_{\left(V=2\right)}\left(-1\right)^{6}U^{3}\right)
\end{align*}

Since $U$ and $V$ are independent,
\begin{align*}
\mathbb{E}\left(X^{3}\right) & =-\mathbb{E}\left(\boldsymbol{1}_{\left(V=1\right)}\right)\mathbb{E}\left(U^{3}\right)+\mathbb{E}\left(\boldsymbol{1}_{\left(V=2\right)}\right)\mathbb{E}\left(U^{3}\right)\\
 & =\mathbb{E}\left(U^{3}\right)\left[P\left(V=2\right)-P\left(V=1\right)\right]\\
 & =0.
\end{align*}

We compute
\[
\mathbb{E}\left(U^{3}\right)=\left(-2\right)^{3}\dfrac{1}{3}+1^{3}\dfrac{2}{3}=-2.
\]

\textbf{(b) Computation of $\mathbb{E}\left(\boldsymbol{1}_{\left(V=1\right)}X^{3}\right)$}

We have
\[
\mathbb{E}\left(\boldsymbol{1}_{\left(V=1\right)}X^{3}\right)=\mathbb{E}\left(\boldsymbol{1}_{\left(V=1\right)}\left(-1\right)^{3}U^{3}\right),
\]
thus, by independence of the random variables $U$ and $V,$
\begin{align*}
\mathbb{E}\left(\boldsymbol{1}_{\left(V=1\right)}X^{3}\right) & =-\mathbb{E}\left(\boldsymbol{1}_{\left(V=1\right)}\right)\mathbb{E}\left(U^{3}\right)\\
 & =-P\left(V=1\right)\mathbb{E}\left(U^{3}\right)\\
 & =1.
\end{align*}

\textbf{(c) Independence of $X$ and $Y$}

We just show that, on the one hand
\[
\mathbb{E}\left(\boldsymbol{1}_{\left(Y=1\right)}X^{3}\right)=\mathbb{E}\left(\boldsymbol{1}_{\left(V=1\right)}X^{3}\right)=1,
\]
and on the other hand
\[
\mathbb{E}\left(\boldsymbol{1}_{\left(Y=1\right)}\right)\mathbb{E}\left(X^{3}\right)=P\left(V=1\right)\mathbb{E}\left(X^{3}\right)=0.
\]
Hence,
\[
\mathbb{E}\left(\boldsymbol{1}_{\left(Y=1\right)}X^{3}\right)\neq\mathbb{E}\left(\boldsymbol{1}_{\left(Y=1\right)}\right)\mathbb{E}\left(X^{3}\right).
\]

The random variables $X$ and $Y$ are not independent, even if $X^{2}$
and $Y^{2}$ are.

\textbf{3.} We aim to show in this question, that in this case the
random variables $X$ and $Y$ are independent.

\textbf{(a) Existence of the expectation of $f\left(X\right)g\left(Y\right)$}

Denoting
\[
\left\Vert f\right\Vert _{\infty}={}_{x\in\mathbb{Z}}\left|f\left(x\right)\right|,
\]
we have the inequality
\[
\left|f\left(X\right)g\left(Y\right)\right|\leqslant\left\Vert f\right\Vert _{\infty}\left\Vert g\right\Vert _{\infty}.
\]

Therefore, the random variable $f\left(X\right)g\left(Y\right)$ is
bounded and admits an expectation.

\textbf{(b) Proof that} 
\[
\mathbb{E}\left(f\left(X\right)g\left(Y\right)\right)=\begin{cases}
\mathbb{E}\left(f\left(U\right)\right)\mathbb{E}\left(f\left(V\right)\right), & \text{if }f\text{ is even,}\\
0, & \text{if }f\text{ is odd.}
\end{cases}
\]

By the transfer theorem,
\begin{align*}
\mathbb{E}\left(f\left(X\right)g\left(Y\right)\right) & =\mathbb{E}\left(f\left(\left(-1\right)^{V}U\right)g\left(V\right)\right)\\
 & =\sum_{\left(u,v\right)\in\mathbb{Z}^{2}}f\left(\left(-1\right)^{v}u\right)g\left(v\right)P\left(U=u,V=v\right).
\end{align*}

Since the random variables $U$ and $V$ are independent,
\[
\mathbb{E}\left(f\left(X\right)g\left(Y\right)\right)=\sum_{\left(u,v\right)\in\mathbb{Z}^{2}}f\left(\left(-1\right)^{v}u\right)g\left(v\right)P\left(U=u\right)P\left(V=v\right).
\]

Because the functions $f$ and $g$ are bounded,
\[
\sum_{\left(u,v\right)\in\mathbb{Z}^{2}}\left|f\left(\left(-1\right)^{v}u\right)g\left(v\right)\right|P\left(U=u\right)P\left(V=v\right)<+\infty.
\]

By the Fubini theorem,
\begin{equation}
\mathbb{E}\left[f\left(X\right)g\left(Y\right)\right]=\sum_{v\in\mathbb{Z}}g\left(v\right)\left[\sum_{u\in\mathbb{Z}}f\left(\left(-1\right)^{v}u\right)P\left(U=u\right)\right]P\left(V=v\right)\label{eq:esp_fXgX}
\end{equation}

\begin{itemize}
\item If $f$ is even, it follows, by the transfer theorem,
\[
\mathbb{E}\left[f\left(X\right)g\left(Y\right)\right]=\sum_{v\in\mathbb{Z}}g\left(v\right)\left[\mathbb{E}\left(f\left(U\right)\right)\right]P\left(V=v\right).
\]
\end{itemize}
Still by applying the transfer theorem,
\[
\mathbb{E}\left[f\left(X\right)g\left(Y\right)\right]=\mathbb{E}\left(f\left(U\right)\right)\mathbb{E}\left(g\left(V\right)\right).
\]

\begin{itemize}
\item If $f$ is odd, it results from the equality $\refpar{eq:esp_fXgX}$
that
\begin{multline}
\mathbb{E}\left[f\left(X\right)g\left(Y\right)\right]=\sum_{v\in2\mathbb{Z}}g\left(v\right)\left[\mathbb{E}\left(f\left(U\right)\right)\right]P\left(V=v\right)\\
-\sum_{v\in2\mathbb{Z}+1}g\left(v\right)\left[\mathbb{E}\left(f\left(U\right)\right)\right]P\left(V=v\right).\label{eq:E_fX_gY}
\end{multline}
\end{itemize}
Since the family $\left\{ f\left(u\right)P\left(U=u\right)\right\} _{u\in\mathbb{Z}}$
is summable, it results from the transfer theorem that
\begin{align*}
\mathbb{E}\left(f\left(U\right)\right) & =\sum_{u\in\mathbb{N}}f\left(u\right)P\left(U=u\right)+\sum_{u\in-\mathbb{N}^{\ast}}f\left(u\right)P\left(U=u\right)\\
 & =\sum_{u\in\mathbb{N}}f\left(u\right)P\left(U=u\right)+\sum_{u\in\mathbb{N}^{\ast}}f\left(-u\right)P\left(U=-u\right).
\end{align*}

Since the law $P_{U}$ is symmetric and $f$ is odd, we obtain
\begin{align*}
\mathbb{E}\left(f\left(U\right)\right) & =\sum_{u\in\mathbb{N}}f\left(u\right)P\left(U=u\right)-\sum_{u\in\mathbb{N}^{\ast}}f\left(u\right)P\left(U=u\right)\\
 & =f\left(0\right)P\left(U=0\right).
\end{align*}

Since $P\left(U=0\right)=0,$ it follows that $\mathbb{E}\left(f\left(U\right)\right)=0$
and applying to the equality $\refpar{eq:E_fX_gY},$ it comes
\[
\mathbb{E}\left[f\left(X\right)g\left(Y\right)\right]=0.
\]

\textbf{(c) Proof of the identity $\mathbb{E}\left(f\left(X\right)g\left(Y\right)\right)=\mathbb{E}\left(f\left(X\right)\right)\mathbb{E}\left(g\left(Y\right)\right)$}

Any function $f$ can be decomposed as the sum of an even function
$f_{P}$ and an odd function $f_{I}.$ 

From the above
\begin{align*}
\mathbb{E}\left(f\left(X\right)g\left(Y\right)\right) & =\mathbb{E}\left[\left(f_{P}\left(X\right)+f_{I}\left(X\right)\right)g\left(Y\right)\right]\\
 & =\mathbb{E}\left[f_{P}\left(X\right)g\left(Y\right)\right]+\mathbb{E}\left[f_{I}\left(X\right)g\left(Y\right)\right]\\
 & =\mathbb{E}\left(f_{P}\left(U\right)\right)\mathbb{E}\left(g\left(V\right)\right).
\end{align*}

Taking $g=1,$we get
\[
\mathbb{E}\left(f\left(X\right)\right)=\mathbb{E}\left(f_{P}\left(U\right)\right).
\]

So
\[
\mathbb{E}\left(f\left(X\right)g\left(Y\right)\right)=\mathbb{E}\left(f\left(X\right)\right)\mathbb{E}\left(g\left(Y\right)\right).
\]

\textbf{(d) Independence of the random variables $X$ and $Y$}

Taking $f=\boldsymbol{1}_{A}$ and $g=\boldsymbol{1}_{B}$ for every
subset $A$ and $B$ of $\mathbb{Z},$ we obtain \textbf{the independence
of $X$ and $Y.$}

\end{solution}

\begin{figure}[t]
\begin{center}(Cliché en attente d'autorisation)

\copyright Renault Communication, cliché C 24654 - Ce cliché est
actuellement consultable sur: \href{https://journals.openedition.org/artefact/6711}{openedition.org}\end{center}

\caption{\protect\href{https://en.wikipedia.org/wiki/Pierre_B\%25C3\%25A9zier}{Pierre Bézier}
(1910-1999)}\sindex[fam]{Bézier, Pierre}
\end{figure}

\begin{solution}{}{solexercise5.6}

\textbf{1. Computation of $\mathbb{E}\left(f\left(\dfrac{S_{n}}{n}\right)\right)$}

We know that the law followed by $S_{n}$ is the binomial law $\mathcal{B}\left(n,x\right).$

By applying the transfer theorem, we have
\[
\mathbb{E}\left(f\left(\dfrac{S_{n}}{n}\right)\right)=\sum^{n}_{k=0}f\left(\dfrac{k}{n}\right)\left(\begin{array}{c}
n\\
k
\end{array}\right)x^{k}\left(1-x\right)^{n-k}.
\]

\textbf{2.} The aim is to establish a constructive proof of the Weierstrass
approximation theorem in the case where the functions are defined
on the interval $\left[0;1\right].$ This approximation method leads
to \textbf{Bézier}\footnotemark\sindex[fam]{Bézier, Pierre} \textbf{polynomials\index{Bézier polynomials}}
and the Bézier procedure used for linear modelling and curve generation
of arbitrary shape passing through given points.

\textbf{a. Proof that $\delta\left(\epsilon\right)$ tends to $0$
with $\epsilon$}

The function $f$ is continuous on the compact interval $\left[0;1\right]:$
it implies that $f$ is uniformly continuous. Therefore, $\delta\left(\epsilon\right)$
tends to 0 with $\epsilon.$

\textbf{b. Proof that $_{x\in\left[0,1\right]}\left|B_{n}\left(x\right)-f\left(x\right)\right|\leqslant\delta\left(\epsilon\right)+\dfrac{2\left\Vert f\right\Vert _{\infty}}{n\epsilon^{2}}.$
Uniform convergence of $B_{n}$ to $f$ on $\left[0;1\right]$}

For every $x\in\left]0;1\right[,$
\[
\left|B_{n}\left(x\right)-f\left(x\right)\right|=\left|\mathbb{E}\left(f\left(\dfrac{S_{n}}{n}\right)\right)-f\left(x\right)\right|,
\]
so
\begin{multline*}
\left|B_{n}\left(x\right)-f\left(x\right)\right|\\
\,\,\,\,\leqslant\mathbb{E}\left(\boldsymbol{1}_{\left(\left|\dfrac{S_{n}}{n}-x\right|\leqslant\epsilon\right)}\left|f\left(\dfrac{S_{n}}{n}\right)-f\left(x\right)\right|\right)+\mathbb{E}\left(\boldsymbol{1}_{\left(\left|\dfrac{S_{n}}{n}-x\right|>\epsilon\right)}\left|f\left(\dfrac{S_{n}}{n}\right)-f\left(x\right)\right|\right).
\end{multline*}

Thus
\[
\left|B_{n}\left(x\right)-f\left(x\right)\right|\leqslant\delta\left(\epsilon\right)+2\left\Vert f\right\Vert _{\infty}P\left(\left|\dfrac{S_{n}}{n}-x\right|>\epsilon\right).
\]

By the Chebyshev inequality,
\[
P\left(\left|\dfrac{S_{n}}{n}-\mathbb{E}\left(\dfrac{S_{n}}{n}\right)\right|>\epsilon\right)\leqslant\dfrac{1}{\epsilon^{2}}\sigma^{2}_{\frac{S_{n}}{n}}.
\]

On the one hand, we have
\[
\mathbb{E}\left(\dfrac{S_{n}}{n}\right)=x,
\]
and on the other hand, since the random variables $X_{k}$ are independent
and follow the same law, they have same variance. 

Hence,
\[
\sigma^{2}_{\frac{S_{n}}{n}}=\dfrac{1}{n^{2}}\sigma^{2}_{S_{n}}=\dfrac{1}{n}\sigma^{2}_{X_{1}}.
\]

Since $\sigma^{2}_{X_{1}}=x\left(1-x\right)\leqslant1,$ we obtain,
for every $x\in\left]0;1\right[,$
\[
\left|B_{n}\left(x\right)-f\left(x\right)\right|\leqslant\delta\left(\epsilon\right)+\dfrac{2\left\Vert f\right\Vert _{\infty}}{n\epsilon^{2}}.
\]

The function $B_{n}-f$ is countinuous on $\left[0;1\right],$ so
\[
_{x\in\left[0;1\right]}\left|B_{n}\left(x\right)-f\left(x\right)\right|\leqslant\delta\left(\epsilon\right)+\dfrac{2\left\Vert f\right\Vert _{\infty}}{n\epsilon^{2}}.
\]

It follows that, for every $\epsilon>0,$
\[
0\leqslant\limsup_{n\to+\infty}{}_{x\in\left[0;1\right]}\left|B_{n}\left(x\right)-f\left(x\right)\right|\leqslant\delta\left(\epsilon\right),
\]
which, taking into account the question (a), shows that
\[
\limsup_{n\to+\infty}{}_{x\in\left[0;1\right]}\left|B_{n}\left(x\right)-f\left(x\right)\right|=0.
\]

Therefore,
\[
\lim_{n\to+\infty}{}_{x\in\left[0;1\right]}\left|B_{n}\left(x\right)-f\left(x\right)\right|=0,
\]
which proves that the sequence of Berstein polynomials converges uniformly
on $\left[0;1\right]$ to the function $f.$

\end{solution}

\footnotetext{\href{https://en.wikipedia.org/wiki/Pierre_B\%25C3\%25A9zier}{Pierre Bézier}
(1910-1999) was a French engineer at Renault. He is one of the founders
of solid, geometric and physical modelling field. He his well-known
for the Bézier curve and surface, that he patented and disseminated.}

\begin{solution}{}{solexercise5.7}

\textbf{1. Computation of $\mathbb{E}\left(\dfrac{S_{n}}{n}\right)$
and $\sigma^{2}_{\frac{S_{n}}{n}}$}

By linearity of the expectation,
\[
\mathbb{E}\left(\dfrac{S_{n}}{n}\right)=\dfrac{1}{n}\sum^{n}_{j=1}\mathbb{E}\left(X_{j}\right).
\]

Since the random variables $X_{j}$ follow the same geometric law
on $\mathbb{N}^{\ast}$ with parameter $p$---and therefore have
the same expectation---we obtain
\[
\mathbb{E}\left(\dfrac{S_{n}}{n}\right)=\dfrac{1}{p}.
\]

Moreover, since these random variables are independent and follow
the same law, with the same variance, we have
\[
\sigma^{2}_{\frac{S_{n}}{n}}=\dfrac{1}{n^{2}}\sigma^{2}_{S_{n}}=\dfrac{1}{n}\sigma^{2}_{X_{1}}.
\]

Thus,
\[
\sigma^{2}_{\frac{S_{n}}{n}}=\dfrac{1-p}{np^{2}}.
\]

\textbf{2. Computation of $G_{S_{n}}$}

Recall that the generating function of $X_{1}$ is, for every $t\in\left[-1;1\right],$
\[
G_{X_{1}}\left(t\right)=\dfrac{pt}{1-t\left(1-p\right)}.
\]

The random variables $X_{j}$ are independent and follow the same
law, and thus share the same generating function. It follows that,
for every $t\in\left[-1;1\right],$
\begin{equation}
G_{S_{n}}\left(t\right)=\prod^{n}_{j=1}G_{X_{j}}\left(t\right)=\left[\dfrac{pt}{1-t\left(1-p\right)}\right]^{n}.\label{eq:G_S_n}
\end{equation}

Furthermore, by the transfer theorem, for every $t\in\left[-1;1\right],$
\[
G_{S_{n}}\left(t\right)=\sum^{+\infty}_{k=1}t^{k}P\left(S_{n}=k\right).
\]

Since
\[
\sum^{+\infty}_{k=1}P\left(S_{n}=k\right)=1,
\]
the function $G_{S_{n}}$ can be developed into a power series on
$\left]-1;1\right[,$ which is also evident from the equality $\refpar{eq:G_S_n}.$ 

The same equality allows us to write, for every $t\in\left]-1;1\right[,$
\[
G_{S_{n}}\left(t\right)=\left(pt\right)^{n}\left[1+\sum^{+\infty}_{k=1}\dfrac{\left(-n\right)\left(-n-1\right)\dots\left(-n-k+1\right)}{k!}\left[t\left(1-p\right)\right]^{k}\right],
\]
which simplifies to
\[
G_{S_{n}}\left(t\right)=p^{n}\sum^{+\infty}_{k=0}\left(\begin{array}{c}
n+k-1\\
n-1
\end{array}\right)\left(1-p\right)^{k}t^{k+n},
\]
or also, after the change of index $l=k+n,$
\[
G_{S_{n}}\left(t\right)=\sum^{+\infty}_{l=n}\left(\begin{array}{c}
l-1\\
n-1
\end{array}\right)p^{n}\left(1-p\right)^{l-n}t^{l}.
\]

\textbf{3. Justification for $G_{S_{n}}$ being expandable as a power
series on $\left]-1;1\right[.$ Computation of $P\left(S_{n}=k\right)$}

The uniqueness of power series expansion gives
\[
P\left(S_{n}=k\right)=\begin{cases}
\left(\begin{array}{c}
k-1\\
n-1
\end{array}\right)p^{n}\left(1-p\right)^{k-n}, & \text{if }k\ensuremath{\geqslant}n,\\
0, & \text{otherwise.}
\end{cases}
\]

Thus, the law of $S_{n}$ is the negative binomial law $\mathcal{B}^{-}\left(n,p\right).$

\textbf{4. Computation of $\mathbb{E}\left(f\left(\dfrac{S_{n}}{n}\right)\right)$
as a function of $B_{n}\left(p\right).$ Justification of the expectation}

The random variable $f\left(\dfrac{S_{n}}{n}\right)$ is bounded and
admits an expectation. By the transfer theorem,
\begin{align*}
\mathbb{E}\left(f\left(\dfrac{S_{n}}{n}\right)\right) & =\sum^{+\infty}_{k=1}f\left(\dfrac{k}{n}\right)P\left(S_{n}=k\right)=\sum^{+\infty}_{k=n}f\left(\dfrac{k}{n}\right)\binom{k-1}{n-1}p^{n}\left(1-p\right)^{k-n}.
\end{align*}

Making the change of indices $l=k-n,$
\begin{align*}
\mathbb{E}\left(f\left(\dfrac{S_{n}}{n}\right)\right) & =p^{n}\sum^{+\infty}_{l=0}f\left(1+\dfrac{l}{n}\right)\binom{l+n-1}{n-1}\left(1-p\right)^{l}.
\end{align*}

We thus obtain the existence of $B_{n}\left(p\right)$ for every $p\in\left]0;1\right[.$
Additionally, $B_{n}\left(1\right)=0$ and
\[
\mathbb{E}\left(f\left(\dfrac{S_{n}}{n}\right)\right)=B_{n}\left(p\right).
\]

\textbf{5.} \textbf{(a) Proof that $\left|\mathbb{E}\left(f\left(\dfrac{S_{n}}{n}\right)\right)-f\left(\dfrac{1}{p}\right)\right|\leqslant\delta\left(\epsilon\right)+2\left\Vert f\right\Vert _{\infty}P\left(\left|\dfrac{S_{n}}{n}-\dfrac{1}{p}\right|>\epsilon\right)$}

We have
\begin{multline*}
\left|\mathbb{E}\left(f\left(\dfrac{S_{n}}{n}\right)\right)-f\left(\dfrac{1}{p}\right)\right|\leqslant\mathbb{E}\left(\boldsymbol{1}_{\left(\left|\dfrac{S_{n}}{n}-\dfrac{1}{p}\right|\leqslant\epsilon\right)}\left|f\left(\dfrac{S_{n}}{n}\right)-f\left(\dfrac{1}{p}\right)\right|\right)\\
+\mathbb{E}\left(\boldsymbol{1}_{\left(\left|\dfrac{S_{n}}{n}-\dfrac{1}{p}\right|>\epsilon\right)}\left|f\left(\dfrac{S_{n}}{n}\right)-f\left(\dfrac{1}{p}\right)\right|\right).
\end{multline*}

By upper-bounding each term---note that $\dfrac{S_{n}}{n}\geqslant1$
and $\dfrac{1}{p}>1$---we obtain the inequality
\begin{equation}
\left|\mathbb{E}\left(f\left(\dfrac{S_{n}}{n}\right)\right)-f\left(\dfrac{1}{p}\right)\right|\leqslant\delta\left(\epsilon\right)+2\left\Vert f\right\Vert _{\infty}P\left(\left|\dfrac{S_{n}}{n}-\dfrac{1}{p}\right|>\epsilon\right).\label{eq:expectation-f}
\end{equation}

\textbf{(b) Proof that $\lim_{n\to+\infty}{}_{x\in\left[1,a\right]}\left|B_{n}\left(\dfrac{1}{x}\right)-f\left(x\right)\right|=0$}

By the Chebyshev inequality,
\[
P\left(\left|\dfrac{S_{n}}{n}-\dfrac{1}{p}\right|>\epsilon\right)\leqslant\dfrac{\sigma^{2}_{\frac{S_{n}}{n}}}{\epsilon^{2}}.
\]

Using the result from Question 1, it follows that
\[
P\left(\left|\dfrac{S_{n}}{n}-\dfrac{1}{p}\right|>\epsilon\right)\leqslant\dfrac{1-p}{np^{2}\epsilon^{2}}.
\]

Then, using the equality $\refpar{eq:expectation-f}$ and the third
question, we have, for every $a>1,$
\[
_{x\in\left[1,a\right]}\left|B_{n}\left(\dfrac{1}{x}\right)-f\left(x\right)\right|={}_{p\in\left[1/a,1\right]}\left|B_{n}\left(p\right)-f\left(\dfrac{1}{p}\right)\right|\leqslant\delta\left(\epsilon\right)+2\left\Vert f\right\Vert _{\infty}\dfrac{a^{2}}{n\epsilon^{2}}.
\]

The uniform continuity of $f$ can be expressed as
\[
\lim_{\epsilon\to0}\delta\left(\epsilon\right)=0.
\]

Although we could conclude as in Exercise $\ref{exo:exercise5.6},$
we provide another argument to conclude here.

For every $\eta>0,$ choose $\epsilon>0$ such that $\delta\left(\epsilon\right)\leqslant\dfrac{\eta}{2},$
and then $N$ such that
\[
2\left\Vert f\right\Vert _{\infty}\dfrac{a^{2}}{N\epsilon^{2}}\leqslant\dfrac{\eta}{2}.
\]

It follows that, for every $n\geqslant N,$
\[
_{x\in\left[1,a\right]}\left|B_{n}\left(\dfrac{1}{x}\right)-f\left(x\right)\right|\leqslant\eta,
\]
which proves that
\[
\lim_{n\to+\infty}{}_{x\in\left[1,a\right]}\left|B_{n}\left(\dfrac{1}{x}\right)-f\left(x\right)\right|=0.
\]

Contrary to the Bernstein polynomial result in the previous exercise,
this approximation result has, to the best of our knowledge, no known
practical applications.

\end{solution}

\begin{solution}{}{solexercise5.8}

\textbf{1. Computation of $G_{S}$ }

Since $S$ takes its values in $\mathbb{N},$ the transfer theorem
allows to write, for every $t\in\left[-1;1\right],$
\[
G_{S}\left(t\right)=\sum^{+\infty}_{k=0}t^{k}P\left(S=k\right).
\]

This can be rewritten by inserting the complete system of constituents
formed by the sets $\left(N=n\right)$ as 
\begin{align*}
G_{S}\left(t\right) & =\sum^{+\infty}_{k=0}t^{k}P\left(\left(S=k\right)\cap\left(\biguplus_{n\in\mathbb{N}}\left(N=n\right)\right)\right)\\
 & =\sum_{k\in\mathbb{N}}t^{k}\left(\sum_{n\in\mathbb{N}}P\left(S=k,N=n\right)\right).
\end{align*}

Since
\[
\sum_{k\in\mathbb{N}}\left|t\right|^{k}\left(\sum_{n\in\mathbb{N}}P\left(S=k,N=n\right)\right)\leqslant\sum_{k\in\mathbb{N}}\left(\sum_{n\in\mathbb{N}}P\left(S=k,N=n\right)\right)=1,
\]
the family $\left\{ t^{k}P\left(S=k,N=n\right)\right\} _{\left(k,n\right)\in\mathbb{N}\times\mathbb{N}}$
is summable, and by the Fubini theorem,
\[
G_{S}\left(t\right)=\sum_{n\in\mathbb{N}}\left(\sum_{k\in\mathbb{N}}t^{k}P\left(S=k,N=n\right)\right).
\]

Since $\left(N=0\right)\subset\left(S=0\right),$
\[
G_{S}\left(t\right)=P\left(N=0\right)+\sum_{n\in\mathbb{N}^{\ast}}\left(\sum_{k\in\mathbb{N}}t^{k}P\left(S=k,N=n\right)\right).
\]

Denoting, for $n\geqslant1,$
\[
S_{n}=\sum^{n}_{j=1}X_{j},
\]
we have
\[
\left(S=k\right)\cap\left(N=n\right)=\left(S_{n}=k\right)\cap\left(N=n\right)
\]
and thus
\[
G_{S}\left(t\right)=P\left(N=0\right)+\sum_{n\in\mathbb{N}^{\ast}}\left(\sum_{k\in\mathbb{N}}t^{k}P\left(S_{n}=k,N=n\right)\right).
\]

Since $S_{n}$ and $N$ are independent, it follows that
\begin{align*}
G_{S}\left(t\right) & =P\left(N=0\right)+\sum_{n\in\mathbb{N}^{\ast}}\left(\sum_{k\in\mathbb{N}}t^{k}P\left(S_{n}=k\right)\right)P\left(N=n\right)\\
 & =P\left(N=0\right)+\sum_{n\in\mathbb{N}^{\ast}}G_{S_{n}}\left(t\right)P\left(N=n\right).
\end{align*}

The random variables $X_{n}$ are independent and follow the same
law, so they share the same generating function. It then follows from
Proposition $\ref{pr:generating_sum}$ that
\begin{align*}
G_{S}\left(t\right) & =P\left(N=0\right)+\sum^{+\infty}_{n=1}\left(G_{X_{1}}\left(t\right)\right)^{n}P\left(N=n\right).
\end{align*}

Hence,
\[
G_{S}\left(t\right)=G_{N}\left(G_{X_{1}}\left(t\right)\right).
\]

We recall that
\[
G_{N}\left(t\right)=\exp\left(\lambda\left(t-1\right)\right).
\]

Additionally, we have
\[
G_{X_{1}}\left(t\right)=\mathbb{E}\left(t^{X_{1}}\right)=\dfrac{1}{6}+\dfrac{1}{2}t+\dfrac{1}{3}t^{2},
\]

Therefore, for every $t\in\left[-1;1\right],$
\[
G_{S}\left(t\right)=\exp\left(\lambda\left(-\dfrac{5}{6}+\dfrac{1}{2}t+\dfrac{1}{3}t^{2}\right)\right).
\]

\textbf{2. Computation of $P\left(S=3\right),$ numerical value in
the case $\lambda=6$}

Since the power series expansion of the function $G_{S}$ on $\left]-1;1\right[$
is unique, we have
\[
P\left(S=3\right)=\dfrac{G^{\prime\prime\prime}_{S}\left(0\right)}{3!}.
\]

Let
\[
f\left(t\right)=\lambda\left(-\dfrac{5}{6}+\dfrac{1}{2}t+\dfrac{1}{3}t^{2}\right),
\]

A straightforward computation yields
\[
G^{\prime\prime\prime}_{S}\left(t\right)=\left[f^{\prime\prime\prime}\left(t\right)+3f^{\prime}\left(t\right)f^{\prime\prime}\left(t\right)+\left(f^{\prime}\left(t\right)\right)^{3}\right]\exp\left(f\left(t\right)\right)
\]
with
\[
\begin{array}{ccccc}
f^{\prime}\left(t\right)=\lambda\left(\dfrac{1}{2}+\dfrac{2}{3}t\right), & \,\,\,\, & f^{\prime\prime}\left(t\right)=\dfrac{2}{3}\lambda, & \,\,\,\, & f^{\prime\prime\prime}\left(t\right)=0.\end{array}
\]

It follows that
\[
P\left(S=3\right)=\dfrac{\lambda^{2}}{6}\left(1+\dfrac{\lambda}{8}\right)\exp\left(-\dfrac{5\lambda}{6}\right).
\]

In the case $\lambda=6,$ we get
\[
P\left(S=3\right)=\dfrac{21}{2}\exp\left(-5\right)\approx0.07.
\]

\textbf{3. Existence and computation of $\mathbb{E}\left(S\right)$
and $\sigma^{2}_{S},$ numerical value for $\lambda=6$}

We have
\[
0\leqslant S\leqslant2N
\]
and since the random variable $N$ admits higher order moments, so
does $S.$ The Wald formula still applies here, as the generating
function of $S$ is obtained by composition of $N$ and $X_{1}.$
Thus
\[
\mathbb{E}\left(S\right)=G^{\prime}_{S}\left(1\right)=\mathbb{E}\left(N\right)\mathbb{E}\left(X_{1}\right).
\]

We compute
\[
\mathbb{E}\left(X_{1}\right)=\dfrac{1}{2}+2\times\dfrac{1}{3}=\dfrac{7}{6},
\]
so
\[
\mathbb{E}\left(S\right)=\dfrac{7}{6}\lambda.
\]

In the case $\lambda=6,$ we get $\mathbb{E}\left(S\right)=7.$

Similarly,
\[
\sigma^{2}_{S}=\sigma^{2}_{X_{1}}\mathbb{E}\left(N\right)+\left(\mathbb{E}\left(X_{1}\right)\right)^{2}\sigma^{2}_{N}.
\]

We have
\[
\begin{array}{ccccc}
\mathbb{E}\left(X^{2}_{1}\right)=\dfrac{1}{2}+4\times\dfrac{1}{3}=\dfrac{11}{6} & \,\,\,\, & \text{and} & \,\,\,\, & \sigma^{2}_{X_{1}}=\dfrac{11}{6}-\left(\dfrac{7}{6}\right)^{2}=\dfrac{17}{36}.\end{array}
\]

Since $\sigma^{2}_{N}=\lambda,$ we conclude
\[
\sigma^{2}_{S}=\dfrac{11}{6}\lambda.
\]

For $\lambda=6,$ this gives $\sigma^{2}_{S}=11.$

\end{solution}

\begin{solution}{}{solexercise5.9}

\textbf{1. Law of $U_{n}$ and $V_{n}$}

The family of events $\left(V_{n}=l\right),\,l\in\mathbb{N},$ forms
a complete system of constituents. Thus, we have
\[
P\left(U_{n}=k\right)=\sum_{l\in\mathbb{N}}P\left(U_{n}=k,V_{n}=l\right).
\]

Then
\[
P\left(U_{n}=k\right)=\sum^{n-k}_{l=0}\dfrac{n!}{k!l!\left[n-\left(k+l\right)\right]!}p^{k}q^{l}r^{n-\left(k+l\right)},
\]
which gives
\begin{align*}
P\left(U_{n}=k\right) & =\dfrac{p^{k}n!}{k!\left(n-k\right)!}\sum^{n-k}_{l=0}\binom{n-k}{l}q^{l}r^{\left(n-k\right)-l}=\binom{n}{k}p^{k}\left(q+r\right)^{n-k}.
\end{align*}

Since $q+r=1-p,$ it comes
\[
P\left(U_{n}=k\right)=\binom{n}{k}p^{k}\left(1-p\right)^{n-k},
\]
that is the law followed by $U_{n}$ is the binomial law $\mathcal{B}\left(n,p\right).$
A similar proof shows that the law of $V_{n}$ corresponds to the
binomial law $\mathcal{B}\left(n,q\right).$

\textbf{2. Independence of $U_{n}$ and $V_{n}$}

We have
\[
P\left(U_{n}=0,V_{n}=0\right)=r^{n}
\]
and
\[
P\left(U_{n}=0\right)P\left(V_{n}=0\right)=\left(1-p\right)^{n}\left(1-q\right)^{n}.
\]

Since $r^{n}\neq\left(1-p\right)^{n}\left(1-q\right)^{n},$ it yields
\[
P\left(U_{n}=0,V_{n}=0\right)\neq P\left(U_{n}=0\right)P\left(V_{n}=0\right),
\]
and the random variables $U_{n}$ and $V_{n}$ are not independent.

\textbf{3. Proof of$\begin{array}{ccc}
\forall x,y\in\mathbb{R}, & \,\,\,\, & \sum\limits^{n}_{k=1}k\binom{n}{k}x^{k-1}y^{n-k}=n\left(x+y\right)^{n-1}.\end{array}$ Computation of $\mathbb{E}\left(U_{n},V_{n}\right)$}

Differentiating with respect to $x$ both sides of the binomial identity
\[
\left(x+y\right)^{n}=\sum^{n}_{k=0}\binom{n}{k}x^{k}y^{n-k},
\]
we immediately obtain the equality $\refpar{eq:formula_ex59_3}.$

Using the transfer theorem
\begin{align*}
\mathbb{E}\left(U_{n}V_{n}\right) & =\sum_{\left(k,l\right)\in\mathbb{N}^{2}}klP\left(U_{n}=k,V_{n}=l\right)\\
 & =\sum_{\substack{0\leqslant k+l\leqslant n\\
k,l\geqslant1
}
}kl\dfrac{n!}{k!l!\left(n-\left(k+l\right)\right)!}p^{k}q^{l}r^{n-\left(k+l\right)}\\
 & =\sum^{n}_{k=1}kp^{k}\dfrac{n!}{k!\left(n-k\right)!}\left[\sum^{n-k}_{l=1}l\binom{n-k}{l}q^{l}r^{\left(n-k\right)-l}\right].
\end{align*}

By applying the equality $\refpar{eq:formula_ex59_3}$, it follows
\begin{align*}
\mathbb{E}\left(U_{n}V_{n}\right) & =\sum^{n}_{k=1}k\binom{n}{k}p^{k}\left[\left(n-k\right)q\left(q+r\right)^{\left(n-k\right)-1}\right]\\
 & =q\sum^{n}_{k=1}k\left(n-k\right)\binom{n}{k}p^{k}\left(1-p\right)^{\left(n-k\right)-1}.
\end{align*}

Since, by simple computation we have
\[
k\left(n-k\right)\binom{n}{k}=nk\binom{n-1}{k},
\]
 we find
\[
\mathbb{E}\left(U_{n}V_{n}\right)=qn\sum^{n-1}_{k=1}k\binom{n-1}{k}p^{k}\left(1-p\right)^{\left(n-1\right)-k}.
\]

Applying the equality $\refpar{eq:formula_ex59_3}$, it holds that
\[
\mathbb{E}\left(U_{n}V_{n}\right)=pqn\left(n-1\right).
\]

\textbf{4. Computation of the Covariance of $U_{n}$ and $V_{n}.$
Variance of $U_{n}+V_{n}$}

The expectations of $U_{n}$ and $V_{n}$ are respectively $np$ and
$nq.$ 

The covariance of $U_{n}$ and $V_{n}$ is computed by
\begin{align*}
\text{cov}\left(U_{n},V_{n}\right) & =\mathbb{E}\left(U_{n}V_{n}\right)-\mathbb{E}\left(U_{n}\right)\mathbb{E}\left(V_{n}\right)\\
 & =pqn\left(n-1\right)-n^{2}pq\\
 & =-npq.
\end{align*}

The variance of $U_{n}+V_{n}$ is
\begin{align*}
\sigma^{2}_{U_{n}+V_{n}} & =\sigma^{2}_{U_{n}}+\sigma^{2}_{V_{n}}+2\text{cov}\left(U_{n},V_{n}\right)\\
 & =np\left(1-p\right)+nq\left(1-q\right)-2npq\\
 & =nr\left(1-r\right).
\end{align*}

\textbf{5. Computation of $P\left(Y=\left(k,l\right)\right).$ Independence
of $U$ and $V$}

The family of events $\left(N=n\right),\,n\in\mathbb{N},$ forms a
complete system of constituents, so
\begin{align*}
P\left(Y=\left(k,l\right)\right) & =\sum_{n\in\mathbb{N}}P\left(Y=\left(k,l\right),N=n\right)\\
 & =\sum_{n\in\mathbb{N}}P\left(Y_{n}=\left(k,l\right),N=n\right).
\end{align*}

Taking into account the independence of the random variables $Y_{n}$
and $N,$
\[
P\left(Y=\left(k,l\right)\right)=\sum_{n\in\mathbb{N}}P\left(Y_{n}=\left(k,l\right)\right)P\left(N=n\right).
\]

Substituting the probabilities by their values, we obtain
\begin{align}
P\left(Y=\left(k,l\right)\right) & =\sum^{+\infty}_{n=l+k}\dfrac{n!}{k!l!\left[n-\left(k+l\right)\right]!}p^{k}q^{l}r^{n-\left(k+l\right)}\text{e}^{-\lambda}\dfrac{\lambda^{n}}{n!}\nonumber \\
 & =\dfrac{\left(\lambda p\right)^{k}}{k!}\dfrac{\left(\lambda q\right)^{l}}{l!}\text{e}^{-\lambda}\sum^{+\infty}_{n=l+k}\dfrac{\left(r\lambda\right)^{n-\left(k+l\right)}}{\left[n-\left(k+l\right)\right]!}\\
 & =\dfrac{\left(\lambda p\right)^{k}}{k!}\dfrac{\left(\lambda q\right)^{l}}{l!}\text{e}^{-\lambda\left(1-r\right)}\nonumber \\
 & =\left[\text{e}^{-\lambda p}\dfrac{\left(\lambda p\right)^{k}}{k!}\right]\left[\text{e}^{-\lambda q}\dfrac{\left(\lambda q\right)^{l}}{l!}\right].\label{eq:5.41}
\end{align}

We then obtain the laws of the marginales $U$ and $V.$ We have for
every $k\in\mathbb{N},$
\begin{align*}
P\left(U=k\right) & =\sum_{l\in\mathbb{N}}P\left(Y=\left(k,l\right)\right)\\
 & =\text{e}^{-\lambda p}\dfrac{\left(\lambda p\right)^{k}}{k!}.
\end{align*}
and for every $l\in\mathbb{N},$ we have
\begin{align*}
P\left(V=l\right) & =\sum_{k\in\mathbb{N}}P\left(Y=\left(k,l\right)\right)\\
 & =\text{e}^{-\lambda q}\dfrac{\left(\lambda q\right)^{l}}{l!}.
\end{align*}

Thus, the random variables $U$ and $V$ follow the Poisson law with
respective parameters $\lambda p$ and $\lambda q.$

The equality $\refpar{eq:5.41}$ shows the independence of the random
variables $U$ and $V.$

\textbf{6. Computation of the generating function of $Y,$ $U$ and
$V.$ Law and independence of $U$ and $V.$} 

By the transfer theorem, we can write, for every $\left(a,b\right)\in\left[0;1\right]^{2},$
the generating function of $Y$ under the form:
\begin{align*}
G_{Y}\left(a,b\right) & =\mathbb{E}\left(a^{U}b^{V}\right)\\
 & =\sum_{\left(k,l\right)\in\mathbb{N}^{2}}a^{k}b^{l}P\left(Y=\left(k,l\right)\right).
\end{align*}

The family of events $\left(N=n\right),n\in\mathbb{N},$ is a complete
system of constituents. 

Thus,
\begin{align*}
G_{Y}\left(a,b\right) & =\sum_{\left(k,l\right)\in\mathbb{N}^{2}}a^{k}b^{l}\left[\sum_{n\in\mathbb{N}}P\left(Y=\left(k,l\right),N=n\right)\right]\\
 & =\sum_{\left(k,l\right)\in\mathbb{N}^{2}}a^{k}b^{l}\left[\sum_{n\in\mathbb{N}}P\left(Y_{n}=\left(k,l\right),N=n\right)\right].
\end{align*}

Since the random variables $Y_{n}$ and $N$ are independent, we obtain
\[
G_{Y}\left(a,b\right)=\sum_{\left(k,l\right)\in\mathbb{N}^{2}}a^{k}b^{l}\left[\sum_{n\in\mathbb{N}}P\left(Y_{n}=\left(k,l\right)\right)P\left(N=n\right)\right].
\]

The terms being nonnegative, we also have
\[
G_{Y}\left(a,b\right)=\sum_{n\in\mathbb{N}}\left[\sum_{\left(k,l\right)\in\mathbb{N}^{2}}a^{k}b^{l}P\left(Y_{n}=\left(k,l\right)\right)\right]P\left(N=n\right),
\]
 that is
\begin{equation}
G_{Y}\left(a,b\right)=\sum_{n\in\mathbb{N}}G_{Y_{n}}\left(a,b\right)P\left(N=n\right),\label{eq:G_Y(a_b)}
\end{equation}
where $G_{Y_{n}}\left(a,b\right)$ is the generating function of $Y_{n}.$

The transfer theorem yields
\begin{align*}
G_{Y_{n}}\left(a,b\right) & =\sum_{\left(k,l\right)\in\mathbb{N}^{2}}a^{k}b^{l}P\left(Y_{n}=\left(k,l\right)\right)\\
 & =\sum_{\substack{0\leqslant k+l\leqslant n\\
k,l\geqslant1
}
}a^{k}b^{l}\dfrac{n!}{k!l!\left[n-\left(k+l\right)\right]!}p^{k}q^{l}r^{n-\left(k+l\right)}.
\end{align*}

By the trinome formula and taking into account that $p+q+r=1,$ we
obtain
\begin{equation}
G_{Y_{n}}\left(a,b\right)=\left(pa+qb+r\right)^{n}.\label{eq:G_Y_n}
\end{equation}

When $N$ follows the Poisson law $\mathcal{P}\left(\lambda\right),$
it yields
\begin{align*}
G_{Y}\left(a,b\right) & =\sum^{+\infty}_{n=0}\left[\left(pa+qb+r\right)^{n}\text{e}^{-\lambda}\dfrac{\lambda^{n}}{n!}\right]\\
 & =\text{e}^{\lambda\left(pa+qb+r-1\right)}\\
 & =\text{e}^{\lambda p\left(a-1\right)}\text{e}^{\lambda q\left(b-1\right)}.
\end{align*}

We have to remark that, in general, we obtain the generating functions
of the marginals $U$ and $V$ by the relations:
\begin{itemize}
\item For every $a\in\left[0,1\right],$
\[
G_{U}\left(a\right)=\mathbb{E}\left(a^{U}\right)=G_{Y}\left(a,1\right).
\]
\item For every $b\in\left[0,1\right],$
\[
G_{V}\left(b\right)=\mathbb{E}\left(b^{V}\right)=G_{Y}\left(1,b\right).
\]
\end{itemize}
Thus,
\begin{itemize}
\item For every $a\in\left[0,1\right],$
\[
G_{U}\left(a\right)=\text{e}^{\lambda p\left(a-1\right)}.
\]
\item For every $b\in\left[0,1\right],$
\[
G_{V}\left(b\right)=\text{e}^{\lambda q\left(b-1\right)}.
\]
\end{itemize}
We have just shown that, for every $\left(a,b\right)\in\left[0,1\right]^{2},$
\begin{equation}
G_{Y}\left(a,b\right)=G_{U}\left(a\right)G_{V}\left(b\right).\label{eq:G_Y_is_G_UG_V}
\end{equation}

Thanks to the transfer theorem, we can now write, that for every $\left(a,b\right)\in\left[0;1\right]^{2},$
\[
\sum_{\left(k,l\right)\in\mathbb{N}^{2}}a^{k}b^{l}P\left(Y=\left(k,l\right)\right)=\sum_{k\in\mathbb{N}}a^{k}P\left(U=k\right)\sum_{l\in\mathbb{N}}b^{l}P\left(V=l\right).
\]

Since all the terms are nonnegative, applying the Fubini theorem,
\[
\sum_{\left(k,l\right)\in\mathbb{N}^{2}}a^{k}b^{l}P\left(Y=\left(k,l\right)\right)=\sum_{\left(k,l\right)\in\mathbb{N}^{2}}a^{k}b^{l}P\left(U=k\right)P\left(V=l\right)<+\infty.
\]

Thus, for every $\left(k,l\right)\in\mathbb{N}^{2},$
\[
P\left(Y=\left(k,l\right)\right)=P\left(U=k,V=l\right)=P\left(U=k\right)P\left(V=l\right)
\]
which shows again the independence of the variables $U$ and $V$.

\begin{remark}{}{}

This method is general. The equality $\refpar{eq:G_Y_is_G_UG_V}$
is a necessary and sufficient condition for independence of $U$ and
$V,$ marginals of the random variable $Y.$

\end{remark}

\textbf{7. Law of $N.$}

Equalities $\refpar{eq:G_Y(a_b)}$ and $\refpar{eq:G_Y_n}$ yield,
for every $\left(a,b\right)\in\left[0;1\right]^{2},$
\begin{align*}
G_{Y}\left(a,b\right) & =\sum_{n\in\mathbb{N}}G_{Y_{n}}\left(a,b\right)a_{n}\\
 & =\sum_{n\in\mathbb{N}}\left(pa+qb+r\right)^{n}a_{n}\\
 & =G_{N}\left(pa+qb+r\right),
\end{align*}
where $G_{N}$ is the generating function of $N.$ 

The hypothesis that $U$ and $V$ are independent implies, by the
computation method of the marginal generating functions from the previous
question, that
\[
G_{Y}\left(a,b\right)=G_{Y}\left(a,1\right)G_{Y}\left(1,b\right),
\]
which gives the relation, for every $\left((a,b\right)\in\left[0;1\right]^{2},$
\[
G_{N}\left(p\left(a-1\right)+q\left(b-1\right)+1\right)=G_{N}\left(p\left(a-1\right)+1\right)G_{N}\left(q\left(b-1\right)+1\right).
\]

We denote, for every $a\in\left[0;1\right],$ 
\[
g\left(a\right)=G_{N}\left(a+1\right).
\]

Then, for every $\left(a,b\right)\in\left[0,1\right]^{2},$
\[
g\left(a+b\right)=g\left(a\right)g\left(b\right).
\]

The function $g$ is continuous and satisfies the Cauchy functional
equation, so $g$ must be an exponential. 

Since $g$ is greater or equal to 1, there exists $\lambda>0$ such
that for every $a\in\left[0,1\right],$ we have $g\left(a\right)=\text{e}^{\lambda a},$
which yields
\[
G_{N}\left(a\right)=\text{e}^{\lambda\left(a-1\right)}.
\]

The generating function determining the law of the random variable,
it shows that the law followed by $N$ is the Poisson law $\mathcal{P}\left(\lambda\right).$

\begin{remark}{}{}

In summary, for the random variable $N$ to follow the Poisson law,
it is necessary and sufficient that the random variables $U$ and
$V$ are independent.

\end{remark}

\end{solution}

\chapter{Random Variables with Density}\label{chap:Random-Variables-with}

\begin{objective}{}{}

Chapter \ref{chap:Random-Variables-with} aims to introduce random
variables with densities in an elementary way, without addressing
the possible existence issues, that will be dealt with in Part \ref{part:Deepening-Probability-Theory}.
\begin{itemize}
\item Section \ref{sec:Probabilities-on} begins by defining the Borel $\sigma-$algebra
on $\mathbb{R}^{n}$ and introduces the concept of a density for a
probability on $\mathbb{R}^{n}.$ It concludes by examining the special
case of $\mathbb{R}.$ Several classical probability laws on $\mathbb{R}$
are then presented: the uniform law on a closed interval, the exponential
law, the Cauchy law, the normal law and the chi-squared law. The section
ends with examples of classical probability laws on $\mathbb{R}^{2}:$
the uniform law on a closed rectangle, on a closed disk, and the standard
normal law.
\item Section \ref{sec:Law-of-a} addresses the density of a random variable
taking values in $\mathbb{R}^{n}$ and its cumulative distribution
function, then explores their properties and their relationship to
probabilities. The marginal law of a random variable is then defined
in the special case $n=2.$ The law of a function of a random variable
is examined in the case of a monotonic function. 
\item Section \ref{sec:Expectation-and-Variance} begins by defining the
mathematical expectation of a real-valued random variable, followed
by its variance. Classical properties are then presented. The covariance
of two real-valued random variables is defined, as well as the correlation
coefficient.
\item The previously introduced concepts are applied in Section \ref{sec:Laplace-Gauss-Law-in}
to the case of the two-dimensional Laplace-Gauss law.
\item Section \ref{sec:Independence-of-Two} discusses the independence
of two real-valued random variables by examining criteria for independence.
\item In Section \ref{sec:Sum-of-Independent}, the sum of two independent
real-valued random variables is considered, and the convolution of
two densities is introduced.
\item Section \ref{sec:Conditional-Densities} focuses on conditional densities,
providing a formal definition followed by a proposition adapting Bayes
theorem to the case of densities.
\item Section \ref{sec:Appendix:-the-Riemann} is an appendix that provides
a rigorous framework for the Riemann integral used throughout the
chapter.
\end{itemize}
\end{objective}

The extension of the theory developed in previous chapters to non-countable
space $\Omega$ raises serious mathematical difficulties.

For instance, if we take $\Omega=\left[0,1\right]\times\left[0,1\right],$
there does not exist any probability on $\left(\Omega,\mathcal{P}\left(\Omega\right)\right)$
that:
\begin{itemize}
\item On the one hand, assigns the same probability to two subsets of $\left[0,1\right]\times\left[0,1\right]$
that are translations of each other;
\item And on the other hand, is such that $P\left(\left[a,b\right]\times\left[c,d\right]\right)$
is equal to the area of the rectangle $\left[a,b\right]\times\left[c,d\right].$
\end{itemize}
Nevertheless, these are basic properties one would expect from a uniform
probability law.

However, a probability $P$ satisfying these properties can be defined
on a smaller $\sigma-$algebra than $\mathcal{P}\left(\Omega\right),$
namely the $\sigma-$algebra of Borel sets. The notion of a Borel
set belongs to measure theory, which provides the proper framework
for an in-depth study of probability calculus. That, however, is not
the framework we use in this part.

The aim of this chapter is to introduce the notion of a random variable
with a density in an elementary way, without addressing existence
issues. The point of view is therefore mostly descriptive\footnote{This is the perspective adopted in French national examinations (Mathematics
CAPES and Internal Agrégation) for becoming a mathematics teacher
in secondary schools.}.

\section{Probabilities on $\mathbb{R}^{n}$}\label{sec:Probabilities-on}

\subsection{Density of a Probability on $\mathbb{R}^{n}$}

We assert, without proof, the existence of the \textbf{\sindex[fam]{Borel, Emile}\mindex{Borel sigma-algebra @ Borel $\sigma-$algebra}\mindex{sigma-algebra@$\sigma-$algebra ! Borel}Borel}
\textbf{$\sigma-$algebra}.

\begin{proposition}{Borel $\sigma-$algebra}{}

There exists a smaller $\sigma-$algebra---smaller in the sense of
inclusion---on $\mathbb{R}^{n}$ containing the family of open rectangles
$\prod^{n}_{i=1}\left]a_{i},b_{i}\right[.$ This $\sigma-$algebra
is called the \textbf{\sindex[fam]{Borel, Emile}\mindex{Borel sigma-algebra@Borel $\sigma-$algebra}Borel}
\textbf{$\sigma-$algebra} of $\mathbb{R}^{n}$ and denoted $\mathcal{B}_{\mathbb{R}^{n}}.$
Its elements are called the \textbf{Borel subsets\index{Borel subsets}}
of $\mathbb{R}^{n}.$ Every open set is Borel, any close set is Borel,
and every half-open rectangle $\prod^{n}_{i=1}\left]a_{i},b_{i}\right]$
is also Borel.

\end{proposition}

\begin{figure}[t]
\begin{center}\includegraphics[width=0.4\textwidth]{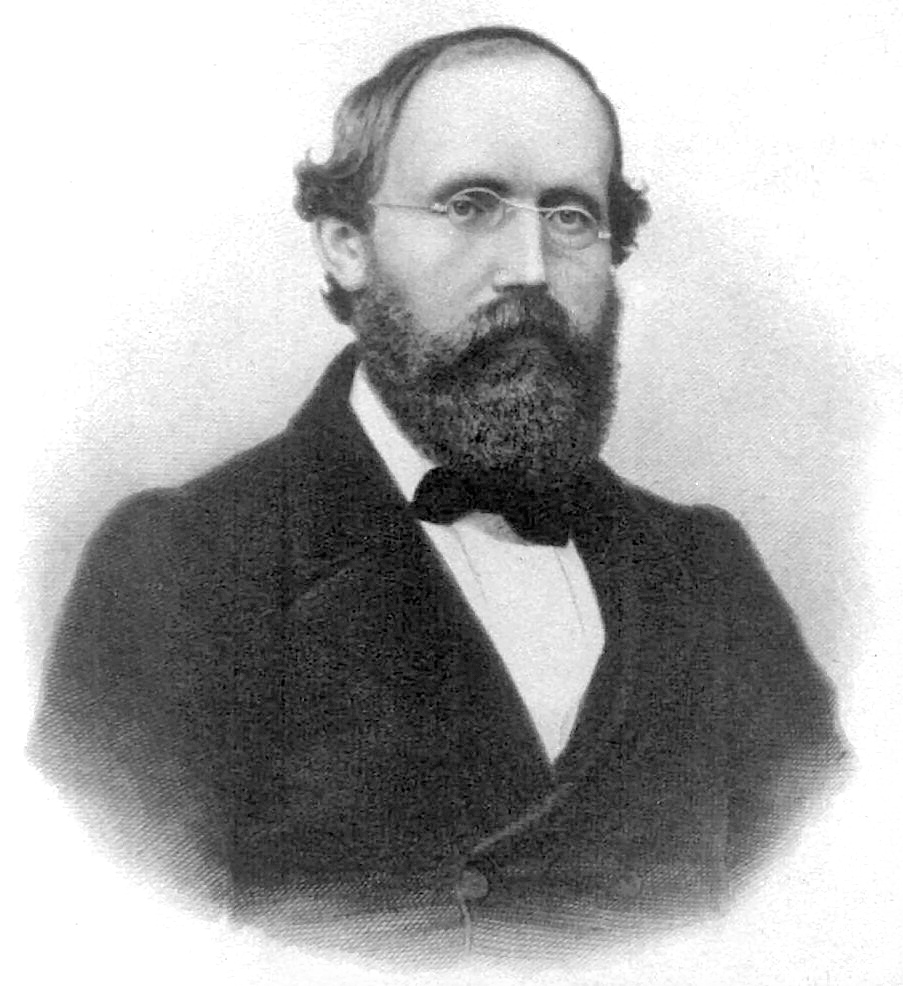}

{\tiny Unknown author---\href{http://www.sil.si.edu/digitalcollections/hst/scientific-identity/explore.htm}{http://www.sil.si.edu/digitalcollections/hst/scientific-identity/explore.htm}---Public
domain}\end{center}

\caption{\textbf{\protect\href{https://en.wikipedia.org/wiki/Bernhard_Riemann}{George Friedrich Bernhard Riemann}}
(1826 - 1866)}\sindex[fam]{Riemann, Friedrich}
\end{figure}

Let $f$ be a function from $\mathbb{R}^{n}$ to $\mathbb{R}^{+}$
\sindex[fam]{Riemann, Friedrich}\textbf{Riemann}{\bfseries\footnote{\textbf{\href{https://en.wikipedia.org/wiki/Bernhard_Riemann}{George Friedrich Bernhard Riemann}\sindex[fam]{Riemann, Friedrich}}
(1826 - 1866) was a German mathematician. He made profound contributions
to analysis, number theory, and differential geometry. Particularly,
he gave a rigorous formulation of the integral, known now as the Riemann
integral, and he is also known for his work on Fourier series. He
introduced the Riemann surfaces in complex analysis. The Riemann hypothesis,
published in 1859, marks the foundation of analytic number theory.
He made pioneer contributions to differential geometry and built the
foundations of the mathematics of general relativity. He is considered
as one of the greatest mathematicians of all time. }}\textbf{-integrable\index{Riemann-integrable}} on $\mathbb{R}^{n}$
and such that $\intop_{\mathbb{R}^{n}}f\left(x\right)\textrm{d}x=1.$
We assert, and we accept it without proof, that there exists a unique
probability $P$ on the probabilizable space $\left(\mathbb{R}^{n},\mathcal{B}_{\mathbb{R}^{n}}\right)$
such that, for every half-open rectangles $A$ in $\mathbb{R}^{n}$
of the form $\prod^{n}_{i=1}\left]a_{i},b_{i}\right],$\boxeq{
\begin{equation}
P\left(A\right)=\intop_{A}f(x)\textrm{d}x.\label{eq:P(A)asinton_A}
\end{equation}
}

The function $f$ is called the \index{density}\textbf{density} of
the probability $P.$ We also say that the probability $P$ has density
$f.$

\begin{particular}{}{}

If $n=1,$ and if $f$ is a function defined on $\mathbb{R},$ taking
nonnegative real values, continuous except at a finite number of points,
and such that
\[
\intop^{+\infty}_{-\infty}f\left(x\right)\textrm{d}x=1,
\]
then there exists a unique probability $P$ on the probabilizable
space $\left(\mathbb{R},\mathcal{B}_{\mathbb{R}}\right)$ such that,
for every real numbers $a$ and $b$ with $a<b,$

\boxeq{
\[
P\left(\left]a,b\right]\right)=\intop^{b}_{a}f\left(x\right)\textrm{d}x.
\]
}

We admit that this equality $\refpar{eq:P(A)asinton_A}$ remains valid
when $A$ is a finite union of intervals on $\mathbb{R},$ whether
they are closed, open, or semi-open, and whether they are bounded
or not.

\end{particular}

\subsection{Classical Examples of Probability Laws on $\mathbb{R}$}

We present classical examples of probability laws defined on the probabilizable
space $\left(\mathbb{R},\mathcal{B}_{\mathbb{R}}\right)$ by a density
function $f.$

\subsubsection{Uniform Law on the Interval $\left[a,b\right]$}

The \textbf{uniform law on the interval $\left[a,b\right]$\index{uniform law}\mindex{law ! uniform}}
is denoted $\mathcal{U}\left(\left[a,b\right]\right).$ 

Its density $f$ is defined for every $x\in\mathbb{R}$ by

\boxeq{
\[
f\left(x\right)=\dfrac{1}{b-a}\boldsymbol{1}_{\left[a,b\right]}\left(x\right),
\]
}

where $a$ and $b$ are two real numbers such that $a<b.$

The uniform law assigns equal probability to any two sub-intervals
of the same length within the interval $\left[a,b\right].$

\subsubsection{Exponential Law with Parameter $p>0$}

The \textbf{exponential law with parameter $p>0$\index{exponential law}\mindex{law ! exponential}}
is denoted $\exp\left(p\right).$ Its density $f$ is defined for
every $x\in\mathbb{R}$ by

\boxeq{
\[
f\left(x\right)=\boldsymbol{1}_{\mathbb{R}^{+}}\left(x\right)p\text{e}^{-px}.
\]
}

It is worth to mentioning that this probability assigns zero probability
to any interval contained in $\mathbb{R}^{-}.$ This probability law
is often used to model waiting times.

\begin{figure}
\begin{center}\begin{tikzpicture}[
    declare function={explaw(\t,\p)=\p*exp(-\p*\t);}
]
	\definecolor{orange}{HTML}{FFA200}
	\definecolor{purple}{HTML}{9500FF}
  \begin{axis}[
		axis x line=center,
		axis y line=center,
		xtick = \empty,
		ytick = \empty,
		ymin=-0.02,
		ymax=0.5
	]
\addplot[green, samples=80, smooth, thick, domain=-3:0 ]{0};
\addplot[green, samples=80, smooth, thick, domain=0:10.5 ]{explaw(x,0.4)};
\draw[green,thick] (0,0)
      -- (0,0.4);
\node at (-0.5,0.4) {$p$};

  \end{axis}

\end{tikzpicture}\end{center}

\caption{Exponential law with parameter $p$}
\end{figure}

\subsubsection{Cauchy Law}

\begin{figure}[t]
\begin{center}\includegraphics[width=0.4\textwidth]{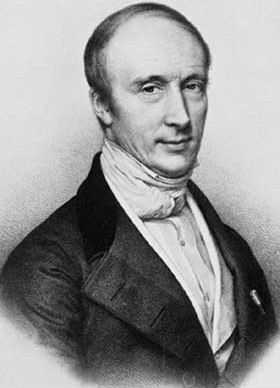}

{\tiny Library of Congress Prints and Photographs --- Public Domain}\end{center}

\caption{\textbf{\protect\href{https://en.wikipedia.org/wiki/Augustin-Louis_Cauchy}{Baron Augustin-Louis Cauchy}}
(1789 - 1857)}\sindex[fam]{Cauchy, Augustin-Louis}
\end{figure}

The density of the \textbf{\sindex[fam]{Cauchy, Augustin-Louis}Cauchy}{\bfseries\footnote{\textbf{\href{https://en.wikipedia.org/wiki/Augustin-Louis_Cauchy}{Baron Augustin-Louis Cauchy}\sindex[fam]{Cauchy, Augustin-Louis}}
(1789-1857) was a French mathematician, engineer, and physicist. By
stating and proving rigorously key theorems of calculus, he was creating
the field of real analysis, and pioneering the fiel of complex analysis.
He wrote more than eight hundred research articles on various topics
in both mathematics and mathematical physics, in particular in continuum
physics.}}\textbf{ law\index{Cauchy law}\mindex{law ! Cauchy}} is defined
for every $x\in\mathbb{R}$ by

\boxeq{
\[
f\left(x\right)=\dfrac{1}{\pi}\dfrac{1}{1+x^{2}}.
\]
}

\begin{figure}
\begin{center}\begin{tikzpicture}
  \begin{axis}[axis lines=center]
    \addplot[blue, samples=80, smooth, thick, domain=-4.25:4.25 ]
      {1/pi*1/(1+x^2)};
  \end{axis}
\end{tikzpicture}\end{center}

\caption{Cauchy law}
\end{figure}

\subsubsection{Laplace-Gauss Law, or Normal, with Parameters $m\in\mathbb{R}$ and
$\sigma^{2}>0$}

\begin{figure}[t]
\begin{center}\includegraphics[width=0.4\textwidth]{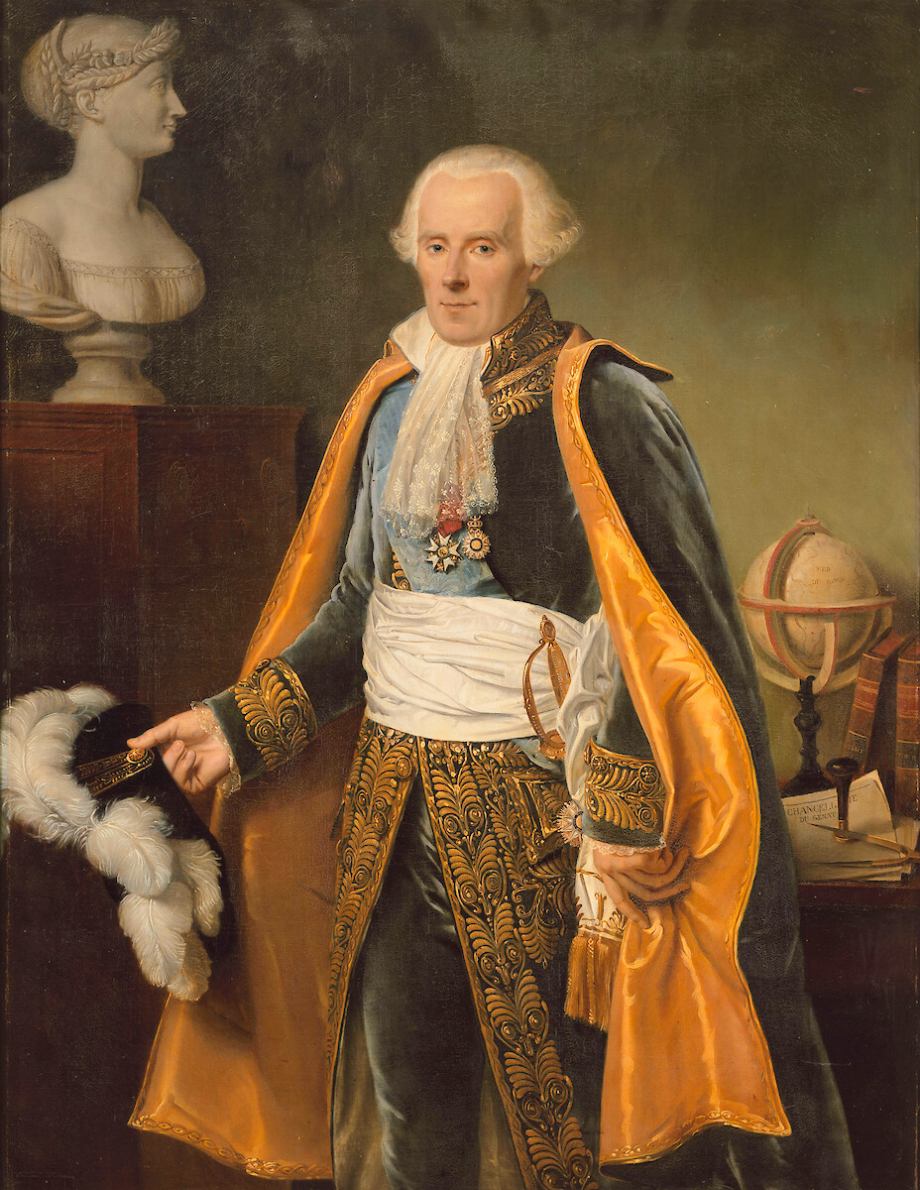}

{\tiny Paulin Guérin --- http://www.photo.rmn.fr/--- Public Domain}\end{center}

\caption{\textbf{\protect\href{https://en.wikipedia.org/wiki/Pierre-Simon_Laplace}{Pierre-Simon, Marquis de Laplace}}
(1749 - 1827)}\sindex[fam]{Laplace, Pierre-Simon}
\end{figure}
\begin{figure}[t]
\begin{center}\includegraphics[width=0.4\textwidth]{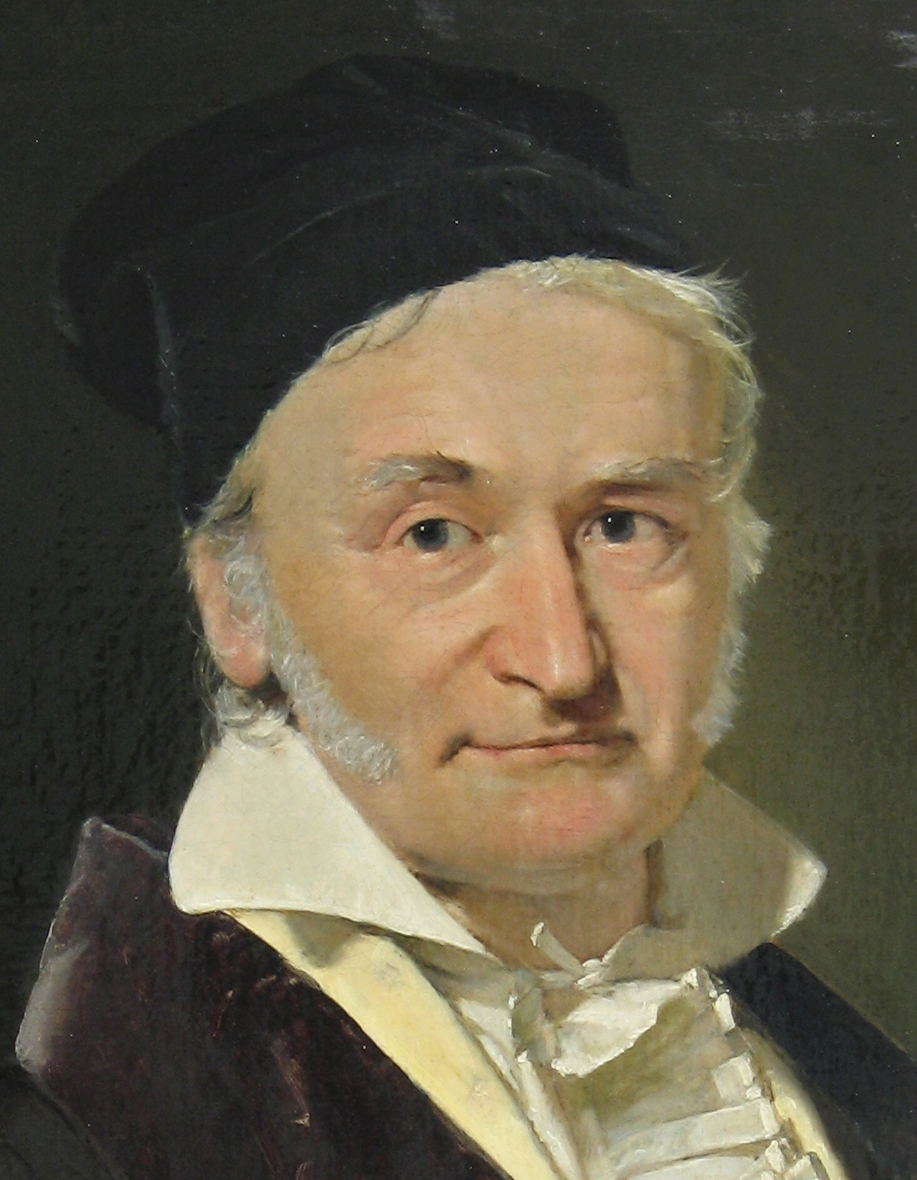}

{\tiny Gottlieb Biermann / D’après Christian Albrecht Jensen--- Public
Domain}\end{center}

\caption{\textbf{\protect\href{https://en.wikipedia.org/wiki/Carl_Friedrich_Gauss}{Carl Friedrich Gauss}}
(1777 - 1855)}\sindex[fam]{Gauss, Carl Friedrich}
\end{figure}

The \textbf{\index{Laplace-Gauss law}\mindex{law ! Laplace-Gauss}Laplace}{\bfseries\footnote{\textbf{\href{https://en.wikipedia.org/wiki/Pierre-Simon_Laplace}{Pierre-Simon, Marquis de Laplace}\sindex[fam]{Laplace, Pierre-Simon}}
(1749 - 1827) was a French mathematician and contributed also in physics,
astronomy, engineering and statistics. He formulated the Laplace's
equation and introduced the Laplace transform; particularly he extended
in his Celestial Mechanics Treaty the work of his predecessors, contributing
to the birth of mathematical astronomy.}}\textbf{\sindex[fam]{Laplace, Pierre-Simon}-Gauss}{\bfseries\footnote{\textbf{\href{https://en.wikipedia.org/wiki/Carl_Friedrich_Gauss}{Carl Friedrich Gauss}\sindex[fam]{Gauss, Carl Friedrich}}
(1777 - 1855) was a German mathematician, astronomer, geodesist, and
physicist. He did several contributions in mathematics and science,
in particular in number theory, algebra, analysis, geometry, statistics
and probability. He was professor of astronomy from 1807 until his
death.}}\textbf{\sindex[fam]{Gauss, Carl Friedrich} law}, or \textbf{normal
law}\index{normal law}\textbf{\mindex{law ! normal}}, with parameters
$m\in\mathbb{R}$ and $\sigma^{2}>0$ is denoted $\mathcal{N}\left(m,\sigma^{2}\right).$
Its density function\footnote{The classical formula $\intop^{+\infty}_{-\infty}e^{-x^{2}/2}\textrm{d}x=\sqrt{2\pi}$
allows, by variable change, to ensure that $f$ is a well-defined
density.} $f$ is defined for every $x\in\mathbb{R}$ by

\boxeq{
\[
f\left(x\right)=\dfrac{1}{\sigma\sqrt{2\pi}}\exp\left(-\dfrac{\left(x-m\right)^{2}}{2\sigma^{2}}\right).
\]
}

Its graph is the famous bell-shaped curbe of Gauss. It has two inflection
points of abscissas $x_{1}=m-\sigma$ and $x_{2}=m+\sigma,$ and 
\[
f\left(x_{1}\right)=f\left(x_{2}\right)=\dfrac{1}{\sigma\sqrt{2\pi e}},
\]
which shows that the higher and narrower the peak is, the smaller
$\sigma$ must be. This law appears very frequently in modelling,
due to the central limit theorem, which will be stated later.

\begin{figure}
\begin{center}\begin{tikzpicture}
	\definecolor{orange}{HTML}{FFA200}
	\definecolor{purple}{HTML}{9500FF}
  \begin{axis}[
		axis x line=center,
		axis y line=center,
		xtick = \empty,
		ytick =\empty,
		ymin=-0.05,
		ymax=0.45
	]
 \addplot[green, samples=80, smooth, thick, domain=-1.75:5.75 ]
      {1/(sqrt(pi*2))*exp(-(x-2)^2/2)};
	\pgfplotsinvokeforeach {1,3} {
    	\draw[purple, thick] (axis cs: #1,0)
      -- (axis cs: #1,{(1/sqrt(2*pi))*exp((-1/2)*(#1-2)^2)});
  }
\pgfplotsinvokeforeach {2} {
    	\draw[orange, thick] (axis cs: #1,0)
      -- (axis cs: #1,{(1/sqrt(2*pi))*exp((-1/2)*(#1-2)^2)});
  }

    	\draw[orange,thick] (axis cs: 0,{(1/sqrt(2*pi)))})
      -- (axis cs: 2,{(1/sqrt(2*pi))});
  \draw[purple,thick] (axis cs: 0,{(1/sqrt(2*pi))*exp((-1/2))})
      -- (axis cs: 3,{(1/sqrt(2*pi))*exp((-1/2))});

	\node[purple] at (axis cs: 1,-0.03) {$m-\sigma$};
\node[orange] at (axis cs: 2,-0.03) {$m$};
\node[purple] at (axis cs: 3,-0.03) {$m+\sigma$};
\node[purple] at (axis cs: -0.7,0.24) {$\frac{1}{\sigma\sqrt{2\pi\textrm{e}}}$};
\node[orange] at (axis cs: -0.7,0.39) {$\frac{1}{\sigma\sqrt{2\pi}}$};
  \end{axis}

\end{tikzpicture}\end{center}

\caption{Normal law}
\end{figure}

\begin{figure}
\begin{center}\begin{tikzpicture}
  \begin{axis}[
		axis lines*=center
	]
    \addplot[blue, samples=80, smooth, thick, domain=-6.25:6.25 ]
      {1/(sqrt(pi*2)*0.3)*exp(-x^2/(2*0.3^2)};
\addplot[green, samples=80, smooth, thick, domain=-6.25:6.25 ]
      {1/sqrt(pi*2)*exp(-x^2/2};
\addplot[red, samples=80, smooth, thick, domain=-6.25:6.25 ]
      {1/(sqrt(pi*2)*3)*exp(-x^2/(2*3^2))};
  \end{axis}
\end{tikzpicture}    \end{center}

\caption{Different values of $\sigma$}
\end{figure}

\subsubsection{Chi-Squared Law with $n-$Degrees of Freedom}

The \textbf{Chi-Squared law with $n-$degrees of freedom\index{Chi-squared law}\mindex{law ! Chi-squared}}
is denoted $\chi^{2}_{n}.$\footnote{The greek letter $\chi$ is transcripted as chi, but is pronouced
ki.} Its density $f$ is defined for every $x\in\mathbb{R}$ by

\boxeq{
\[
f\left(x\right)=\boldsymbol{1}_{\mathbb{R}^{+}}\left(x\right)\dfrac{1}{K_{n}}\exp\left(-\dfrac{x}{2}\right)x^{\frac{n}{2}-1},
\]
}

where, for each integer $p\geqslant1,$ we denote
\[
\begin{array}{ccccc}
K_{2p}=2^{p}\left(p-1\right)! & \,\,\,\, & \text{and} & \,\,\,\, & K_{2p+1}=\dfrac{\left(2p-1\right)!}{2^{p-1}\left(p-1\right)!}\sqrt{2\pi}.\end{array}
\]

We can show that the Chi-Squared law is the law followed by a random
variable of the form $X^{2}_{1}+\dots+X^{2}_{n},$ where $X_{1},\dots,X_{n}$
are independent random variables, each following the standard normal
law $\mathcal{N}\left(0,1\right).$ This explains the expression ``$n-$degrees
of freedom''. Exercise $\ref{exo:exercise6.5}$ gives a proof\footnote{Once this chapter has been read, it can be a good exercise to directly
derive the formula for the density of the law $\chi^{2}_{n}$ in the
case $n=2,$and $n=3.$ For $n=2,$ this requires knowledge of integration
in polar coordinates; for $n=3,$ spherical coordinates must be used.} of this fact in the cases $n=1$ and $n=2.$ The Chi-Squared law
is fundamental in statistics.

\subsection{Classical Example of Probability Laws on $\mathbb{R}^{2}$}

\subsubsection{Uniform Law on the Rectangle $\left[a,b\right]\times\left[c,d\right]$}

The \textbf{uniform law on the rectangle $\left[a,b\right]\times\left[c,d\right]$\mindex{uniform law ! rectangle}\mindex{law ! uniform ! rectangle}}
is denoted $\mathcal{U}\left(\left[a,b\right]\times\left[c,d\right]\right),$
where $a,b,c$ and $d$ are real numbers such that $a<b$ and $c<d.$ 

Its density $f$ is defined for every $x\in\mathbb{R}^{2}$ by

\boxeq{
\[
f\left(x\right)=\dfrac{1}{\left(b-a\right)\left(d-c\right)}\boldsymbol{1}_{\left[a,b\right]\times\left[c,d\right]}\left(x\right).
\]
}

\subsubsection{Uniform Law on the Disk $D\left(O,r\right)$}

The \textbf{uniform law on the disk $D\left(O,r\right)$\mindex{uniform law ! disk}\mindex{law ! uniform ! disk}}
of center $O\left(0,0\right)$ and of radius $r>0$ is denoted $\mathcal{U}\left(D\left(O,r\right)\right).$

Its density $f$ is defined for every $x\in\mathbb{R}^{2}$ by

\boxeq{
\[
f\left(x\right)=\dfrac{1}{\pi r^{2}}\boldsymbol{1}_{D\left(0,r\right)}\left(x\right).
\]
}

More generally, the \textbf{uniform law on the subset $A$\mindex{uniform law ! subset}\mindex{law ! uniform ! subset}}
of the plane $\mathbb{R}^{2}$ that has a well-defined, nonzero area
is defined for every $x\in\mathbb{R}^{2}$ by

\boxeq{
\[
f\left(x\right)=\dfrac{1}{\text{area}\left(A\right)}\boldsymbol{1}_{A}\left(x\right).
\]
}

\subsubsection{Gaussian Law, or Normal, Centered Reduced on $\mathbb{R}^{2}$}

The \textbf{\sindex[fam]{Gauss, Carl Friedrich}Gaussian law}\mindex{law ! Gaussian}\mindex{Gaussian!law}---also
called \textbf{normal law}\index{normal law}---centered and reduced
on $\mathbb{R}^{2}$ is the probability law of density $f$ defined
for every $x=\left(x_{1},x_{2}\right)\in\mathbb{R}^{2}$ by

\boxeq{
\[
f\left(x\right)=\dfrac{1}{2\pi}\exp\left(-\dfrac{x^{2}_{1}+x^{2}_{2}}{2}\right).
\]
}

\section{Law of a Random Variable taking values in $\mathbb{R}^{n}$}\label{sec:Law-of-a}

\subsection{Density of a Random Variable. Cumulative Distribution Function}

In what follows, $n$ is an integer greater than or equal to 1, and
$X$ is a random variable defined on the probabilized space $\left(\Omega,\mathcal{A},P\right),$
taking values in $\left(\mathbb{R}^{n},\mathcal{B}_{\mathbb{R}^{n}}\right).$
We recall that the law followed by $X,$ denoted $P_{X},$ is a probability
defined on the probabilizable space $\left(\mathbb{R}^{n},\mathcal{B}_{\mathbb{R}^{n}}\right)$
such that for every $A\in\mathcal{B}_{\mathbb{R}^{n}}$ by\boxeq{
\[
P_{X}\left(A\right)=P\left(X\in A\right).
\]
}

\begin{definition}{Density of a Random Variable}{}

If there exists a function $f_{X},$ defined on $\mathbb{R}^{n},$
taking nonnegative values and Riemann-integrable, such that
\[
\intop_{\mathbb{R}^{n}}f_{X}\left(x\right)\text{d}x=1,
\]
 and such that for every rectangle $A$ in $\mathbb{R}^{n},$\boxeq{
\[
P_{X}\left(A\right)=\intop_{A}f_{X}\left(x\right)\text{d}x,
\]
}then this function $f_{X}$ is called the \mindex{density!random variable}\textbf{\mindex{random variable ! density}density
of the random variable $X.$} This density fully determines the law
followed by $X,$ and corresponds to the density of the probability
$P_{X}.$

\end{definition}

\begin{definition}{Cumulative Distribution Function of a Random Variable}{cum_dist_function}

If $n=1,$ the function $F_{X}$ from $\mathbb{R}$ to $\left[0,1\right],$
defined for every $x\in\mathbb{R}$ by\\
\boxeq{
\[
F_{X}\left(x\right)=P\left(X\leqslant x\right)
\]
}is called the \index{cumulative distribution function}\textbf{cumulative
distribution function} of the random variable $X.$

\end{definition}

\begin{remark}{}{}

If $n=1,$ and if the random variable $X$ admits a density $f_{X},$
it follows from the definitions that for every real numbers $a$ and
$b$ such that $a<b,$

\boxeq{
\begin{equation}
P\left(a\leqslant X\leqslant b\right)=F_{X}\left(b\right)-F_{X}\left(a\right)=\intop^{b}_{a}f_{X}\left(x\right)\text{d}x.\label{eq:P(a<=00003DX<=00003Db)repartition}
\end{equation}
}

Intuitively, we can write
\[
P\left(x<X\leqslant x+\text{d}x\right)=f_{X}\left(x\right)\text{d}x
\]
where $\text{d}x$ is considered as ``infinitesimally small''. Also
this explicit formulation is not mathematically rigourous, it is often
used in practice by physicists and engineers. This expression is equivalent,
for a regular function $f_{X},$ to performing a first order Taylor
expansion of $F_{X},$ and neglecting the remainder. 

\end{remark}

We now give the main properties of the cumulative distribution functions.

\begin{proposition}{Properties of a Cumulative Distribution Function}{}

Let $X$ be a random variable taking real values, with cumulative
distribution function $F_{X}.$

(i) The cumulative distribution function $F_{X}$ fully determines
the law followed by the random variable $X.$

(ii) For every real numbers $a$ and $b$ such that $a<b,$\boxeq{
\begin{equation}
F_{X}\left(b\right)-F_{X}\left(a\right)=P\left(a\leqslant X\leqslant b\right),\label{eq:F_X_and_P}
\end{equation}
}and, for every $x\in\mathbb{R},$\boxeq{
\begin{equation}
\begin{array}{ccccc}
\lim_{x\to+\infty}F_{X}\left(x\right)=1 & \,\,\,\, & \text{and} & \,\,\,\, & \lim_{x\to-\infty}F_{X}\left(x\right)=0.\end{array}\label{eq:lim+inf-infrepfct}
\end{equation}
}

(iii) The cumulative distribution function $F_{X}$ is a nondecreasing
function, right-continuous and admits a left-hand limit at every point.
Moreover, for every $x\in\mathbb{R},$\boxeq{
\begin{equation}
P\left(X=x\right)=F_{X}\left(x\right)-\lim_{y\nearrow x}F_{X}\left(y\right),\label{eq:P_Xisx_function_of_F_X}
\end{equation}
}that is $P\left(X=x\right)$ corresponds to the jump of $F_{X}$
at $x.$

(iv) If the random variable $X$ admits a density $f_{X},$ then the
cumulative distribution function $F_{X}$ is continuous at every point,
and for every $x\in\mathbb{R},$
\[
P\left(X=x\right)=0.
\]

Moreover, for every $x\in\mathbb{R},$

\boxeq{
\begin{equation}
P\left(X\leqslant x\right)=F_{X}\left(x\right)=\intop^{x}_{-\infty}f_{X}\left(u\right)\text{d}u\label{eq:rep_function_as_int}
\end{equation}
}

and the cumulative distribution function $F_{X}$ is differentiable
at every point where the density $f_{X}$ is continuous.

\end{proposition}

\begin{proof}{}{}

(i) This point comes from measure theory, and is therefore taken as
given.

(ii) If $a<b,$ we have, since $\left(X\leqslant a\right)\subset\left(X\leqslant b\right),$
\begin{align*}
F_{X}\left(b\right)-F_{X}\left(a\right) & =P\left(X\leqslant b\right)-P\left(X\leqslant a\right)=P\left(a<X\leqslant b\right).
\end{align*}

Moreover, for every real sequence $\left(x_{n}\right)_{n\in\mathbb{N}}$
that is nondecreasing and converges to $+\infty,$
\[
\bigcup_{n\in\mathbb{N}}\left(X\leqslant x_{n}\right)=\Omega.
\]
The sequence of sets $\left(X\leqslant x_{n}\right)$ being nondecreasing,
it follows that
\[
\lim_{n\to+\infty}P\left(X\leqslant x_{n}\right)=P\left(\Omega\right)=1.
\]

Similarly, for every sequence $\left(x_{n}\right)_{n\in\mathbb{N}}$
that is decreasing and converges to $-\infty,$
\[
\bigcap_{n\in\mathbb{N}}\left(X\leqslant x_{n}\right)=\emptyset.
\]

The sequence of sets $\left(X\leqslant x_{n}\right)$ being decreasing,
\[
\lim_{n\to+\infty}P\left(X\leqslant x_{n}\right)=P\left(\emptyset\right)=0.
\]

Thus, the relationships $\refpar{eq:lim+inf-infrepfct}$ are established.

(iii) It follows from $\refpar{eq:F_X_and_P}$ that the function $F_{X}$
is nondecreasing. Moreover, for every sequence $\left(b_{n}\right)_{n\in\mathbb{N}}$
that decreases and converges to $b,$
\[
F_{X}\left(b_{n}\right)-F_{X}\left(b\right)=P\left(X\in\left]b,b_{n}\right]\right).
\]

The sequence of sets $\left(X\in\left]b,b_{n}\right]\right)$ is decreasing
and has empty intersection. It follows that
\[
\lim_{n\to+\infty}F_{X}\left(b_{n}\right)=F_{X}\left(b\right),
\]
which shows the right-continuity.

Since the function $F_{X}$ is nondecreasing and bounded, it admits
a left-hand limit at every point\footnotemark. Finally, since for
every strictly increasing sequence $\left(x_{n}\right)_{n\in\mathbb{N}}$
that converges to $x,$
\[
\left(X=x\right)=\bigcap_{n\in\mathbb{N}}\left(x_{n}<X\leqslant x\right),
\]
we obtain
\begin{align*}
P\left(X=x\right) & =\lim_{n\to+\infty}P\left(x_{n}<X\leqslant x\right)\\
 & =\lim_{n\to+\infty}\left(F_{X}\left(x\right)-F_{X}\left(x_{n}\right)\right),
\end{align*}
which proves the relation $\refpar{eq:P_Xisx_function_of_F_X}.$

(iv) If $X$ admits a density $f_{X},$ then from the relation $\refpar{eq:P_Xisx_function_of_F_X},$
for every $x\in\mathbb{R},$
\begin{align*}
P\left(X=x\right) & =\lim_{n\to+\infty}\left(F_{X}\left(x\right)-F_{X}\left(x-\dfrac{1}{n}\right)\right),\\
 & =\lim_{n\to+\infty}\intop^{x}_{x-\frac{1}{n}}f_{X}\left(u\right)\text{d}u\\
 & =0.
\end{align*}

Moreover, the relation $\refpar{eq:rep_function_as_int}$ can be otained
from $\refpar{eq:F_X_and_P}$ and $\refpar{eq:lim+inf-infrepfct}.$
It implies the differentiability property of the cumulative distribution
function.

\end{proof}

\footnotetext{The set of its points of discontinuity is at most countable.}

\begin{remarks}{}{}

A random variable $X$ that satisfies, for every $x\in\mathbb{R},$
the condition $P\left(X=x\right)=0$ is said to follow a \index{diffusive law}\textbf{\mindex{law ! diffusive}diffusive
law}. This property does not induce that the random variable admits
a density.

In the context of the remark following Definition $\ref{df:cum_dist_function},$
it results from this proposition that the formula $\refpar{eq:conditional_density}$
leads to the equality
\[
P\left(a\leqslant X\leqslant b\right)=\intop^{b}_{a}f_{X}\left(x\right)\textrm{d}x.
\]

\end{remarks}

\subsection{Marginals of a Random Variable taking Values in $\mathbb{R}^{2}$}

\begin{definition}{Marginals}{}

If $X=\left(X_{1},X_{2}\right)$ is a random variable taking values
in $\mathbb{R}^{2},$ the random variables $X_{1}$ and $X_{2}$ are
called the \textbf{marginals\index{marginal of a random variable}\mindex{random variable!marginal}}
of the random variable $X.$

\end{definition}

\begin{remark}{}{}

Of course, this definition extends to $\mathbb{R}^{n}$ and to any
projection of the random variable $X$ on a subspace generated by
vectors of the canonical basis.

For simplicity, we describe the law of a marginal only in the case
where the random variable takes values in $\mathbb{R}^{2};$ the extension
to $\mathbb{R}^{n}$ raises only a minor complication in writing.

\end{remark}

For what follows, we say that a function defined on $\mathbb{R}^{2}$
is \textbf{regular}\index{regular}, if the Fubini theorem\footnote{Some sufficient conditions of regularity are given in Section \ref{sec:Appendix:-the-Riemann}.
This terminology is not standard and is used only within the framework
of the results in this chapter.} can be applied to it on $\mathbb{R}^{2}.$ 

\begin{proposition}{Density of Marginals}{density_marginals}

If the random variable $X=\left(X_{1},X_{2}\right)$ taking values
in $\mathbb{R}^{2},$ admits a regular density $f_{X},$ then the
marginals $X_{1}$ and $X_{2}$ admit as densities $f_{X_{1}}$ and
$f_{X_{2}}$ respectively, given by:
\begin{itemize}
\item For every $x_{1}\in\mathbb{R},$\\
\boxeq{
\[
f_{X_{1}}\left(x_{1}\right)=\intop^{+\infty}_{-\infty}f_{X}\left(x_{1},x_{2}\right)\text{d}x_{2}.
\]
}
\item For every $x_{2}\in\mathbb{R},$\\
\boxeq{
\[
f_{X_{2}}\left(x_{2}\right)=\intop^{+\infty}_{-\infty}f_{X}\left(x_{1},x_{2}\right)\text{d}x_{1}.
\]
}
\end{itemize}
\end{proposition}

\begin{proof}{}{}

Since for every $x_{1}\in\mathbb{R},$ 
\[
\left(X_{1}\leqslant x_{1}\right)=\bigcup_{n\in\mathbb{N}}\left[\left(X_{1}\leqslant x_{1}\right)\cap\left(X_{2}\leqslant n\right)\right]
\]
and since the sequence of sets inside the brackets is nondecreasing,
\begin{align*}
P\left(X_{1}\leqslant x_{1}\right) & =\lim_{n\to+\infty}P\left(\left(X_{1}\leqslant x_{1}\right)\cap\left(X_{2}\leqslant n\right)\right)\\
 & =\lim_{n\to+\infty}\intop_{\left]-\infty,x_{1}\right]\times\left]-\infty,n\right]}f_{X}\left(u_{1},u_{2}\right)\text{d}u_{1}\text{d}u_{2}\\
 & =\lim_{n\to+\infty}\intop_{\left]-\infty,x_{1}\right]\times\left]-\infty,+\infty\right]}f_{X}\left(u_{1},u_{2}\right)\text{d}u_{1}\text{d}u_{2}.
\end{align*}

By applying the Fubini theorem\footnotemark,
\[
P\left(X_{1}\leqslant x_{1}\right)=\intop^{x_{1}}_{-\infty}\left[\intop^{+\infty}_{-\infty}f_{X}\left(u_{1},u_{2}\right)\text{d}u_{2}\right]\text{d}u_{1}.
\]

This proves the existence of the density $f_{X_{1}}$ and gives its
expression.

\end{proof}

\footnotetext{See the Appendix on the Rieman integral in $\mathbb{R}^{n}$
in Section \ref{sec:Appendix:-the-Riemann}. }

\begin{example}{Marginals of a Reduced Centered Normal Law on $\mathbb{R}^2$}{}

Let $X$ be a random variable following the centered and reduced normal
law on $\mathbb{R}^{2}.$ Determine the laws followed by the marginals
$X_{1}$ and $X_{2}$ of $X.$

\end{example}

\begin{solutionexample}{}{}

We have for every $x_{1}\in\mathbb{R},$
\[
f_{X_{1}}=\intop^{+\infty}_{-\infty}\dfrac{1}{2\pi}\text{e}^{-\frac{x^{2}_{1}+x^{2}_{2}}{2}}\text{d}x_{2}.
\]

Additionally,
\[
\intop^{+\infty}_{-\infty}\dfrac{1}{\sqrt{2\pi}}\text{e}^{-\frac{x^{2}_{2}}{2}}\text{d}x_{2}=1.
\]

It follows that, for every $x_{1}\in\mathbb{R},$\boxeq{
\[
f_{x_{1}}\left(x_{1}\right)=\dfrac{1}{\sqrt{2\pi}}\text{e}^{-\frac{x^{2}_{1}}{2}}.
\]
}

Thus, the marginals $X_{1}$ and $X_{2}$ follow the normal law with
parameters 0 and 1, denoted $\mathcal{N}\left(0,1\right).$

\end{solutionexample}

\subsection{Law of a Random Variable Function}

It is useful to know how to study the law followed by a function of
a random variable. Still for the same reasons, we consider in this
Chapter only the case where this function is monotonic.

\begin{proposition}{Density of a Composed Function}{}

Let $X$ be a real-valued random variable, admitting a density $f_{X}$
and let $g$ be a real-valued function defined on $\mathbb{R},$ strictly
monotonic and differentiable. Then the real-valued random variable
$Y=g\left(X\right)$ admits a density $f_{Y}$ given by

\boxeq{
\[
f_{Y}\left(y\right)=\begin{cases}
f_{X}\left(g^{-1}\left(y\right)\right)\left|\left(g^{-1}\right)^{\prime}\left(y\right)\right|, & \text{if }y\in g\left(\mathbb{R}\right),\\
0, & \text{otherwise.}
\end{cases}
\]
}

\end{proposition}

\begin{proof}{}{}

The function $g$ is defined on $\mathbb{R},$ strictly monotonic
and differentiable, and is therefore bijective from $\mathbb{R}$
onto the interval $g\left(\mathbb{R}\right).$ 

Let us suppose that $g$ is strictly decreasing; its derivative is
then negative, and never zero at any point.

Then, for every real $y\in g\left(\mathbb{R}\right),$ we have---taking
care to reverse the inequality due to the fact that $g$ is strictly
decreasing
\begin{align*}
P\left(Y\leqslant y\right) & =P\left(X\geqslant g^{-1}\left(y\right)\right)=P\left(X>g^{-1}\left(y\right)\right)=\intop^{+\infty}_{g^{-1}\left(y\right)}f_{X}\left(x\right)\text{d}x,
\end{align*}
the second equality holding because $X$ admits a density.

By making the variable change defined by $v=g\left(x\right),$we obtain
\begin{align*}
P\left(Y\leqslant y\right) & =\intop^{-\infty}_{y}f_{X}\left(g^{-1}\left(v\right)\right)\left(g^{-1}\right)^{\prime}\left(v\right)\text{d}v=\intop^{y}_{-\infty}f_{X}\left(g^{-1}\left(v\right)\right)\left|\left(g^{-1}\right)^{\prime}\left(v\right)\right|\text{d}v.
\end{align*}

If additionally $g\left(\mathbb{R}\right)=\left]\alpha,\beta\right[,$
then for every $y<\alpha,$
\[
P\left(Y\leqslant y\right)=P\left(g\left(X\right)\leqslant y\right)=0,
\]
and for every $y>\beta,$
\[
P\left(Y\leqslant y\right)=P\left(g\left(X\right)\leqslant y\right)=1.
\]

Thus, the proposition is proven in the case of a strictly decreasing
function. If $g$ is strictly increasing, the proof is similar and
is let to the interested reader as an exercise.

\end{proof}

\section{Expectation and Variance of a Real Random Variable}\label{sec:Expectation-and-Variance}

\subsection{Mathematical Expectation}

\begin{definition}{Average or Mathematical Expectation}{}

Let $X$ be a real-valued random variable defined on the probabilized
space $\left(\Omega,\mathcal{A},P\right),$ and admitting a density
$f_{X}.$

If the function $x\mapsto\left|x\right|f_{X}\left(x\right)$ is Riemann-integrable
over $\mathbb{R},$ we say that $X$ admits an expectation.

The \textbf{mean\index{mean}} or \textbf{mathematical expectation\index{expectation}\index{mathematical expectation}}
of $X,$ denoted $\mathbb{E}\left(X\right)$ or $\mathbb{E}X,$ is
defined by

\boxeq{
\[
\mathbb{E}\left(X\right)=\intop_{\mathbb{R}}xf_{X}\left(x\right)\text{d}x.
\]
}

\end{definition}

We observe the analogy between this definition and that of the expectation
of a discrete random variable.

If now $\phi$ is a real-valued function defined on $\mathbb{R},$
such that $Y=\phi\left(X\right)$ is a random variable, we can compute
the expectation of $\phi\left(X\right)$ using a formula similar analogous
to the discrete case---see the formula $\refpar{eq:transfer_theorem}$
known as the transfer theorem---that is
\begin{equation}
\mathbb{E}\left(\phi\left(X\right)\right)=\intop_{\mathbb{R}}\phi\left(x\right)f_{X}\left(x\right)\text{d}x,\label{eq:transfer_theorem_density_var}
\end{equation}
under the obvious condition that the function $\left|\phi\right|f_{X}$
is Riemann-integrable. 

Unfortunately, in this continuous setting, such a formula is not only
difficult to prove rigorously, but even to state properly. Its left-hand
side only makes sense if $Y=\phi\left(X\right)$ admits a density---which
is the case if $\phi$ is a differentiable and strictly monotonic
function\footnote{More generally, $Y=\phi\left(X\right)$ has a density if $\mathbb{R}$
can be decomposed in an union of finite number of intervals on which
$\phi$ is differentiable and strictly monotonic. One can refer to
Exercise $\ref{exo:exercise6.5}$ for a proof in the case $\phi\left(X\right)=X^{2}.$}. Alternately, $Y=\phi\left(X\right)$ is a discrete random variable---which
occurs for instance, if $\phi$ is a constant, or if $\phi\left(x\right)=\left\lfloor x\right\rfloor $
(integer part of $x$).

But there are other possibilities:
\begin{itemize}
\item If $\phi\left(x\right)=0$ for $x<0,$ and $\phi\left(x\right)=1+x$
for $x\geqslant0,$
\item The law of $\phi\left(X\right)$ may include a discrete part at 0
and a continuous part elsewhere,
\item And many other more complex cases.
\end{itemize}
We admit that, even when $\phi\left(X\right)$ is neither discrete
nor have a density, it is still possible to give meaning to the expression
$\mathbb{E}\left(\phi\left(X\right)\right)$---provided that $\left|\phi\right|f_{X}$
is integrable---and that the formula $\refpar{eq:transfer_theorem_density_var}$
may be used to compute its value in all such cases.

We state, without proof, the following theorem for a random variable
$X=\left(X_{1},\dots,X_{n}\right)$ taking values in $\mathbb{R}^{n}.$

\begin{theorem}{Transfer Theorem for a Random Variable with Density in $\mathbb{R}^n$}{transfer_r_n}

Let $X$ be a random variable defined on a probabilizable space $\left(\Omega,\mathcal{A},P\right),$
taking values in $\mathbb{R}^{n},$ and admitting a density $f_{X}.$

Let $\phi$ be a real-valued function defined on $\mathbb{R}^{n}.$

If the function $\left|\phi\right|f_{X}$ is Riemann-integrable on
$\mathbb{R}^{n},$ then the random variable $\phi\left(X\right)$
admits an expectation, given by the formula

\boxeq{
\[
\mathbb{E}\left(\phi\left(X\right)\right)=\intop_{\mathbb{R}^{n}}\phi\left(x\right)f_{X}\left(x\right)\text{d}x.
\]
}

\end{theorem}

\begin{proposition}{Linearity of the Mathematical Expectation}{linerity_mean}

(i) Let $X_{1}$ and $X_{2}$ be two real random variables such that
the random variable $\left(X_{1},X_{2}\right)$ admits a regular\footnotemark
density.

If $X_{1}$ and $X_{2}$ admit an expectation, then the random variable
$\lambda_{1}X_{1}+\lambda_{2}X_{2}$ also admits an expectation for
every real numbers $\lambda_{1}$ and $\lambda_{2},$ and

\boxeq{
\[
\mathbb{E}\left(\lambda_{1}X_{1}+\lambda_{2}X_{2}\right)=\lambda_{1}\mathbb{E}\left(X_{1}\right)+\lambda_{2}\mathbb{E}\left(X_{2}\right).
\]
}

(ii) Let $X$ be a real random variable with density $f_{X}$ and
suppose $X$ admits an expectation.

Then, for every real numbers $a$ and $b,$ the random variable $aX+b$
admits an expectation, and

\boxeq{
\begin{equation}
\mathbb{E}\left(aX+b\right)=a\mathbb{E}\left(X\right)+b.\label{eq:e_a_X_pl_b}
\end{equation}
}

\end{proposition}

\footnotetext{See previous note in Proposition $\ref{pr:density_marginals}$}

\begin{proof}{}{}

(i) Consider the function $\phi$ defined for every $\left(x_{1},x_{2}\right)\in\mathbb{R}^{2},$
by
\[
\phi\left(x_{1},x_{2}\right)=\left|\lambda_{1}x_{1}+\lambda_{2}x_{2}\right|.
\]

By the triangular inequality, for every $\left(x_{1},x_{2}\right)\in\mathbb{R}^{2},$
\[
0\leqslant\phi\left(x_{1},x_{2}\right)f_{\left(X_{1},X_{2}\right)}\left(x_{1},x_{2}\right)\leqslant\left(\left|\lambda_{1}\right|\left|x_{1}\right|+\left|\lambda_{2}\right|\left|x_{2}\right|\right)f_{\left(X_{1},X_{2}\right)}\left(x_{1},x_{2}\right).
\]

We want to show that the right-hand side of this inequality defines
a regular function. By Proposition $\ref{pr:density_marginals}$,
the function $x_{2}\mapsto\left|x_{1}\right|f_{\left(X_{1},X_{2}\right)}\left(x_{1},x_{2}\right)$
is Riemann-integrable, with integral $\left|x_{1}\right|f_{X_{1}}\left(x_{1}\right).$
Moreover, since $X_{1}$ admits an expectation, the function $x_{1}\mapsto\left|x_{1}\right|f_{X_{1}}\left(x_{1}\right)$
is also Riemann-integrable. By exchanging the roles of $x_{1}$ and
$x_{2},$ we obtain the announced regularity. The Fubini theorem then
allows us to write
\begin{multline*}
\intop_{\mathbb{R}^{2}}\left(\left|\lambda_{1}\right|\left|x_{1}\right|+\left|\lambda_{2}\right|\left|x_{2}\right|\right)f_{\left(X_{1},X_{2}\right)}\left(x_{1},x_{2}\right)\text{d}x_{1}\text{d}x_{2}=\left|\lambda_{1}\right|\mathbb{E}\left(X_{1}\right)+\left|\lambda_{2}\right|\mathbb{E}\left(X_{2}\right)<+\infty.
\end{multline*}

This shows that the random variable $\lambda_{1}X_{1}+\lambda_{2}X_{2}$
admits an expectation. We can then write, using the Fubini theorem
\begin{align*}
\mathbb{E}\left(\lambda_{1}X_{1}+\lambda_{2}X_{2}\right) & =\intop_{\mathbb{R}^{2}}\left(\lambda_{1}x_{1}+\lambda_{2}x_{2}\right)f_{\left(X_{1},X_{2}\right)}\left(x_{1},x_{2}\right)\text{d}x_{1}\text{d}x_{2}\\
 & =\lambda_{1}\intop^{+\infty}_{-\infty}x_{1}\left(\intop^{+\infty}_{-\infty}f_{\left(X_{1},X_{2}\right)}\left(x_{1},x_{2}\right)\text{d}x_{2}\right)\text{d}x_{1}\\
 & \qquad\qquad+\lambda_{2}\intop^{+\infty}_{-\infty}x_{2}\left(\intop^{+\infty}_{-\infty}f_{\left(X_{1},X_{2}\right)}\left(x_{1},x_{2}\right)\text{d}x_{1}\right)\text{d}x_{2},
\end{align*}
which ensures the result.

(ii) Taking for $\phi$ the function defined for every $x\in\mathbb{R}$
by
\[
\phi\left(x\right)=ax+b,
\]
we have
\[
\left|\phi\left(x\right)\right|f_{X}\left(x\right)\leqslant\left(\left|a\right|\left|x\right|+\left|b\right|\right)f_{X}\left(x\right),
\]
which shows, since both the functions $f_{X}$ and $x\mapsto\left|x\right|f_{X}\left(x\right)$
are Riemann-integrable, that the function $\left|\phi\right|f_{X}$
is Riemann-integrable. Therefore, the random variable $aX+b$ admits
an expectation. 

Moreover,
\[
\mathbb{E}\left(aX+b\right)=\intop^{+\infty}_{-\infty}\left(ax+b\right)f_{X}\left(x\right)\text{d}x,
\]
and applying the linearity of the integral with the fact that
\[
\intop^{+\infty}_{-\infty}f_{X}\left(x\right)\textrm{d}x=1
\]
we obtain the result.

\end{proof}

\begin{example}{Normal Law Expectation}{}

If the random variable $X$ follows the normal law $\mathcal{N}\left(m,\sigma^{2}\right),$
then\boxeq{
\[
\mathbb{E}\left(X\right)=m.
\]
}

\end{example}

\begin{solutionexample}{}{}

Indeed, we have
\[
\mathbb{E}\left(X\right)=\intop^{+\infty}_{-\infty}x\dfrac{1}{\sigma\sqrt{2\pi}}\exp\left(-\dfrac{\left(x-m\right)^{2}}{2\sigma^{2}}\right)\textrm{d}x.
\]

Using the change of variable
\[
y=\dfrac{x-m}{\sigma},
\]
we obtain
\begin{align*}
\mathbb{E}\left(X\right) & =\intop^{+\infty}_{-\infty}\left(\sigma y+m\right)\dfrac{1}{\sqrt{2\pi}}\exp\left(-\dfrac{y^{2}}{2}\right)\text{d}y\\
 & =\sigma\intop^{+\infty}_{-\infty}y\dfrac{1}{\sqrt{2\pi}}\exp\left(-\dfrac{y^{2}}{2}\right)\text{d}y+m\intop^{+\infty}_{-\infty}\dfrac{1}{\sqrt{2\pi}}\exp\left(-\dfrac{y^{2}}{2}\right)\text{d}y.
\end{align*}

The first integral is equal to $0,$ since it is the integral of an
odd function over a symmetric interval and the second integral equals
1, since it is the integral of a density. Therefore,
\[
\mathbb{E}\left(X\right)=m.
\]

\end{solutionexample}

\subsection{Moments of Order 2. Variance.}

\begin{proposition}{Existence of the Expectation Under the Existence of the Expectation of the Square of a Random Variable}{}

Let $X$ be a real-valued random variable with density $f_{X}.$ If
$X^{2}$ admits an expectation, then $X$ admits an expectation.

\end{proposition}

\begin{proof}{}{}

We have for every $x\in\mathbb{R},$ 
\[
\left|x\right|\leqslant x^{2}+1.
\]

The function $x\mapsto\left|x\right|f_{X}\left(x\right)$ is then
Riemann-integrable on $\mathbb{R}.$

\end{proof}

\begin{definition}{Second-Order Moment and Variance of a Random Variable with Density}{}

If $X$ is a real random variable with density $f_{X}$ such that
$X^{2}$ admits an expectation, the nonnegative real number $\mathbb{E}\left(X^{2}\right)$
is called the \textbf{\index{second-order moment}\mindex{moment ! second-order}second-order
moment} of $X,$ and the nonnegative real number $\mathbb{E}\left(\left(X-\mathbb{E}\left(X\right)\right)^{2}\right)$
is called the \textbf{variance\index{variance}} of $X$ and denoted
$\sigma^{2}_{X}.$

\end{definition}

\begin{proposition}{Properties of the Variance of a Random Variable with Density}{}

If $X$ is a real random variable with density $f_{X}$ such that
$X^{2}$ admits an expectation, then the variance of $X$ verifies:

(i)

\boxeq{
\begin{equation}
\sigma^{2}_{X}=\intop^{+\infty}_{-\infty}\left(x-\mathbb{E}\left(X\right)\right)^{2}f_{X}\left(x\right)\textrm{d}x.\label{eq:varX^2}
\end{equation}
}

(ii)

\boxeq{
\[
\sigma^{2}_{X}=\mathbb{E}\left(X^{2}\right)-\left(\mathbb{E}\left(X\right)\right)^{2}.
\]
}

(iii) For every real numbers $a$ and $b,$

\boxeq{
\[
\sigma^{2}_{aX+b}=a^{2}\sigma^{2}_{X}.
\]
}

\end{proposition}

\begin{proof}{}{}

(i) It is enough to apply Theorem $\ref{th:transfer_r_n}$ with, as
function $\phi,$ the function defined for every $x\in\mathbb{R},$
by
\[
\phi\left(x\right)=\left(x-\mathbb{E}\left(X\right)\right)^{2}.
\]

(ii) By developing the square and using the linearity of the integral
in the equality $\refpar{eq:varX^2},$ we obtain successively
\begin{align*}
\sigma^{2}_{X} & =\intop^{+\infty}_{-\infty}\left(x^{2}-2x\mathbb{E}\left(X\right)+\left(\mathbb{E}\left(X\right)\right)^{2}\right)f_{X}\left(x\right)\textrm{d}x\\
 & =\intop^{+\infty}_{-\infty}x^{2}f_{X}\left(x\right)\textrm{d}x-2\mathbb{E}\left(X\right)\intop^{+\infty}_{-\infty}xf_{X}\left(x\right)\textrm{d}x+\left(\mathbb{E}\left(X\right)\right)^{2}\intop^{+\infty}_{-\infty}xf_{X}\left(x\right)\textrm{d}x\\
 & =\mathbb{E}\left(X^{2}\right)-2\left(\mathbb{E}\left(X\right)\right)^{2}+\left(\mathbb{E}\left(X\right)\right)^{2}\\
 & =\mathbb{E}\left(X^{2}\right)-\left(\mathbb{E}\left(X\right)\right)^{2}.
\end{align*}

(iii) From the equality
\[
\left(aX+b\right)-\mathbb{E}\left(aX+b\right)=a\left(X-\mathbb{E}\left(X\right)\right)
\]
it follows that
\begin{align*}
\sigma^{2}_{aX+b} & =\mathbb{E}\left(\left[a\left(X-\mathbb{E}\left(X\right)\right)\right]^{2}\right)=a^{2}\mathbb{E}\left(\left(X-\mathbb{E}\left(X\right)\right)^{2}\right)=a^{2}\sigma^{2}_{X}.
\end{align*}

\end{proof}

We now give an example of the computation of the variance of a random
variable.

\begin{example}{}{Normal Law Variance}

If a random variable $X$ follows the normal law $\mathcal{N}\left(m,\sigma^{2}\right),$
show that its variance $\sigma^{2}_{X}$ is equal to $\sigma^{2}.$

\end{example}

\begin{solutionexample}{}{}

The expectation of the random variable $X$ following the normal law
$\mathcal{N}\left(m,\sigma^{2}\right)$ is $\mathbb{E}\left(X\right)=m.$

We have
\[
\sigma^{2}_{X}=\mathbb{E}\left(\left(X-m\right)^{2}\right).
\]

Using the equality $\refpar{eq:varX^2},$ we obtain
\[
\sigma^{2}_{X}=\intop^{+\infty}_{-\infty}\left(x-m\right)^{2}\dfrac{1}{\sigma\sqrt{2\pi}}\text{e}^{-\frac{\left(x-m\right)^{2}}{2\sigma^{2}}}\textrm{d}x.
\]

Using the change of variable defined by $y=\dfrac{x-m}{\sigma},$
the above equality becomes
\begin{align*}
\sigma^{2}_{X} & =\intop^{+\infty}_{-\infty}\left(\sigma y\right)^{2}\dfrac{1}{\sqrt{2\pi}}\text{e}^{-\frac{y^{2}}{2}}\text{d}y\\
 & =\sigma^{2}\intop^{+\infty}_{-\infty}y^{2}\dfrac{1}{\sqrt{2\pi}}\text{e}^{-\frac{y^{2}}{2}}\text{d}y.
\end{align*}

Using an integration by parts, justified by the convergence of the
two integrals, it holds that
\begin{align*}
\intop^{+\infty}_{-\infty}y^{2}\text{e}^{-\frac{y^{2}}{2}}\text{d}y & =\left[-y\text{e}^{-\frac{y^{2}}{2}}\right]^{+\infty}_{-\infty}+\intop^{+\infty}_{-\infty}\text{e}^{-\frac{y^{2}}{2}}\text{d}y\\
 & =\sqrt{2}\intop^{+\infty}_{-\infty}\text{e}^{-u^{2}}\text{d}u\\
 & =\sqrt{2}\left[\sqrt{\pi}\text{erf}\left(z\right)\right]^{+\infty}_{0}\\
 & =\sqrt{2\pi}.
\end{align*}

Hence,
\[
\sigma^{2}_{X}=\sigma^{2}.
\]

\end{solutionexample}

\begin{definition}{Centered Random Variable. Centered Reduced Random Variable.}{}

Suppose that $X$ is a real random variable of density $f_{X}$ and
that $X$ admits a moment of order two. 
\begin{itemize}
\item The random variable $X-\mathbb{E}\left(X\right)$ is called a \textbf{centered
random variable}\index{centered random variable}\mindex{random variable!centered}:
its expectation is zero.
\item The random variable $\dfrac{X-\mathbb{E}\left(X\right)}{\sigma_{X}}$
is called the \textbf{centered reduced random variable\index{centered reduced random variable}\mindex{random variable!centered!reduced}}
associated with $X:$ its expectation is zero and its standard deviation
is equal to 1.
\end{itemize}
\end{definition}

\begin{figure}[t]
\begin{center}\includegraphics[width=0.4\textwidth]{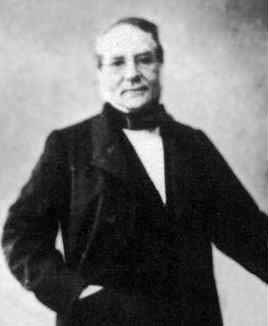}

{\tiny Credits: Public Domain}\end{center}

\caption{\textbf{\protect\href{https://en.wikipedia.org/wiki/Ir\%25C3\%25A9n\%25C3\%25A9e-Jules_Bienaym\%25C3\%25A9}{Irénée-Jules Bienaymé}}
(1796 - 1878)}\sindex[fam]{Bienaymé, Irénée-Jules}
\end{figure}

\begin{figure}[t]
\begin{center}\includegraphics[width=0.4\textwidth]{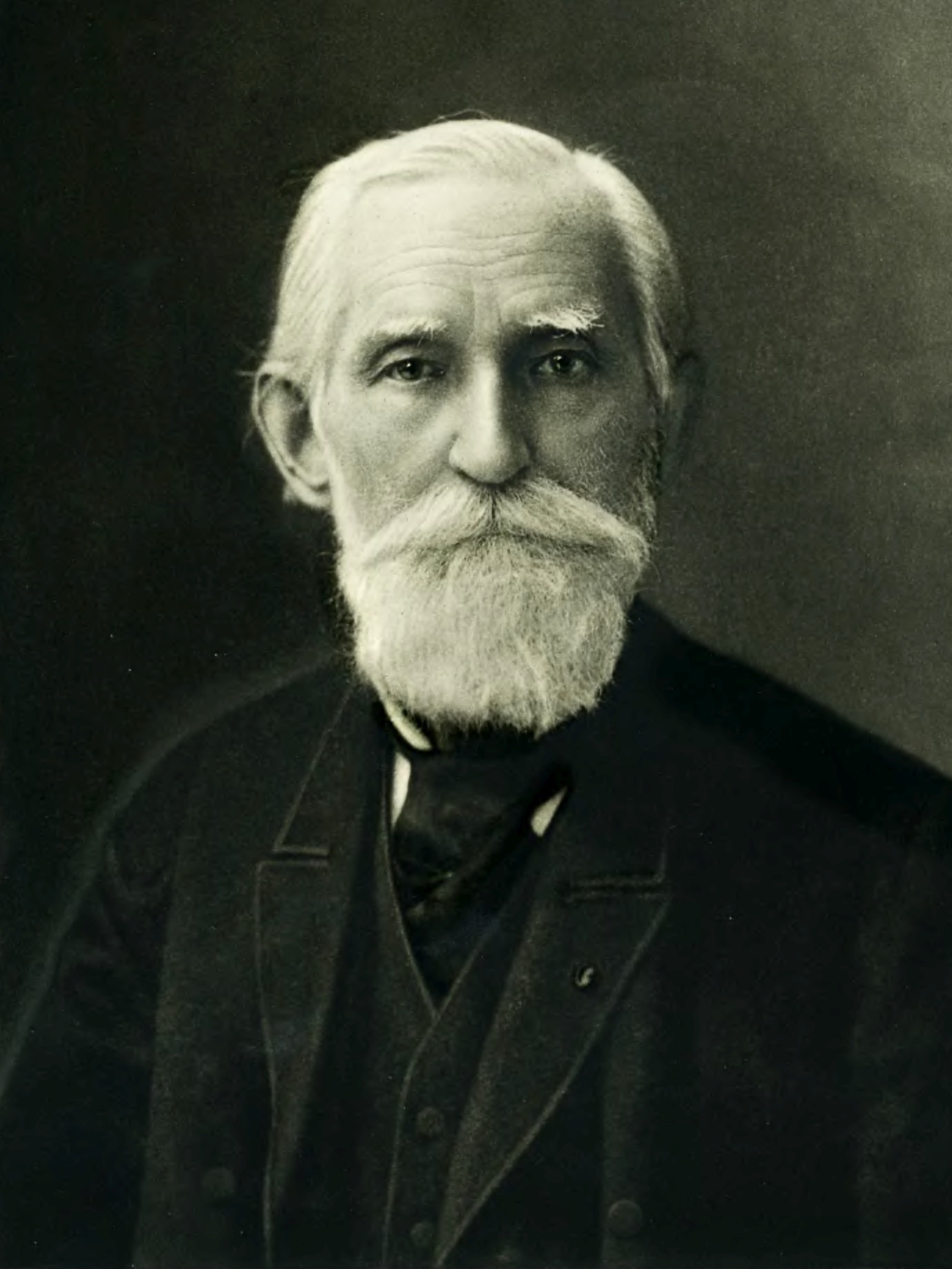}

{\tiny Credits: Public Domain}\end{center}

\caption{\textbf{\protect\href{https://en.wikipedia.org/wiki/Pafnuty_Chebyshev}{Pafnuty Chebyshev}}
(1821-1894)}\sindex[fam]{Chebyshev, Pafnuty}
\end{figure}

The computation of expectation and variance for classical laws will
be tackled in Exercise $\ref{exo:exercise6.7}$.

We now present the \textbf{Bienaymé}{\bfseries\footnote{\textbf{\href{https://en.wikipedia.org/wiki/Ir\%25C3\%25A9n\%25C3\%25A9e-Jules_Bienaym\%25C3\%25A9}{Irénée-Jules Bienaymé}}
(1796 - 1878) was a French statistician. He contributed to the fields
of probability and statistics, with different applications in finance,
demography and social sciences.}}\textbf{\sindex[fam]{Bienaymé, Irénée-Jules}-Chebyshev}{\bfseries\footnote{\textbf{\href{https://en.wikipedia.org/wiki/Pafnuty_Chebyshev}{Pafnuty Chebyshev}}
(1821 - 1894) was a Russian mathematician. He is known for several
fundamental contributions to the fields of probability, statistics,
mechanics and number theory. Among them the Chebyshev inequality used
to prove the weak law of large numbers.}}\textbf{\sindex[fam]{Chebyshev, Pafnuty} inequality\index{Bienaymé-Chebyshev inequality}}
for a random variable with density. Although coarse, it provides an
idea of how the values of $X$ are distributed around its mean. This
inequality is particularly useful for proving the convergence in probability.

\begin{proposition}{Bienaymé-Chebyshev Inequality}{}Let $X$ be
a real-valued random variable with density $f_{X},$ admitting a moment
of order two.

For every real number $\epsilon>0,$

\boxeq{
\[
P\left(\left|X-\mathbb{E}\left(X\right)\right|>\epsilon\right)\leqslant\dfrac{\sigma^{2}_{X}}{\epsilon^{2}}.
\]
}

\end{proposition}

\begin{proof}{}{}

For every real number $\epsilon>0,$ let $m=\mathbb{E}\left(X\right).$
Then
\begin{align*}
\sigma^{2}_{X} & =\intop^{+\infty}_{-\infty}\left(x-m\right)^{2}f_{X}\left(x\right)\textrm{d}x\\
 & \geqslant\intop^{m-\epsilon}_{-\infty}\left(x-m\right)^{2}f_{X}\left(x\right)\textrm{d}x+\intop^{+\infty}_{m+\epsilon}\left(x-m\right)^{2}f_{X}\left(x\right)\textrm{d}x\\
 & \geqslant\epsilon^{2}\left(\intop^{m-\epsilon}_{-\infty}f_{X}\left(x\right)\textrm{d}x+\intop^{+\infty}_{m+\epsilon}f_{X}\left(x\right)\textrm{d}x\right)\\
 & \geqslant\epsilon^{2}P\left(\left|X-\mathbb{E}\left(X\right)\right|>\epsilon\right).
\end{align*}

This yields the inequality.

\end{proof}

\begin{remark}{}{}

The Bienaymé-Chebyshev inequality is often used in the following form,
for $\rho>0$

\boxeq{
\[
P\left(\left|X-\mathbb{E}\left(X\right)\right|>\rho\sigma_{X}\right)\leqslant\dfrac{1}{\rho^{2}}.
\]
}

\end{remark}

\subsection{Covariance and Correlation Coefficient}

\begin{proposition}{Expectation of the Product of Marginals of a Random Variable in $\mathbb{R}^2$ with Density}{}

If the random variable $X=\left(X_{1},X_{2}\right)$ taking values
in $\mathbb{R}^{2}$ admits a regular density $f_{X},$ and if $X_{1}$
and $X_{2}$ have second-order moments, then the product random variable
$X_{1}X_{2}$ admits an expectation.

\end{proposition}

\begin{proof}{}{}

We have, for every $\left(x_{1},x_{2}\right)\in\mathbb{R}^{2},$
\[
\left|x_{1}x_{2}\right|\leqslant\dfrac{1}{2}\left(x^{2}_{1}+x^{2}_{2}\right).
\]

Since $f_{X}$ is a density, it is nonnegative, so
\[
\left|x_{1}x_{2}\right|f_{X}\left(x_{1},x_{2}\right)\leqslant\dfrac{1}{2}\left(x^{2}_{1}+x^{2}_{2}\right)f_{X}\left(x_{1},x_{2}\right).
\]

By Proposition $\ref{pr:density_marginals}$, the marginals admit
densities. A reasoning similar to the one used to show the assertion
(i) of Proposition $\ref{pr:linerity_mean}$ ensures that we can apply
the Fubini theorem and write
\begin{align*}
\intop_{\mathbb{R}^{2}}x^{2}_{1}f_{X}\left(x_{1},x_{2}\right)\text{d}x_{1}\text{d}x_{2} & =\intop^{+\infty}_{-\infty}x^{2}_{1}\left(\intop^{+\infty}_{-\infty}f_{X}\left(x_{1},x_{2}\right)\text{d}x_{2}\right)\text{d}x_{1}\\
 & \leqslant\mathbb{E}\left(X^{2}_{1}\right)\\
 & <+\infty,
\end{align*}
and similarly
\begin{align*}
\intop_{\mathbb{R}^{2}}x^{2}_{1}f_{X}\left(x_{1},x_{2}\right)\text{d}x_{1}\text{d}x_{2} & =\intop^{+\infty}_{-\infty}x^{2}_{2}\left(\intop^{+\infty}_{-\infty}f_{X}\left(x_{1},x_{2}\right)\text{d}x_{1}\right)\text{d}x_{2}\\
 & \leqslant\mathbb{E}\left(X^{2}_{2}\right)\\
 & <+\infty.
\end{align*}

It thus follows that the function
\[
\left(x_{1},x_{2}\right)\mapsto\left|x_{1}x_{2}\right|f_{X}\left(x_{1},x_{2}\right)
\]
is Riemann-integrable on $\mathbb{R}^{2}$ and then the random variable
$X_{1}X_{2}$ admits an expectation.

\end{proof}

Just as in the discrete case, the \textbf{covariance} and the \textbf{correlation
coefficient} of two random variables help measure the link between
them.

\begin{definition}{}{Covariance. Correlation Coefficient}

Let $X=\left(X_{1},X_{2}\right)$ be the random variable taking values
in $\mathbb{R}^{2},$ and suppose $X$ admits a density $f_{X},$
with both $X_{1}$ and $X_{2}$ admitting second order moments.
\begin{itemize}
\item The \textbf{covariance\index{covariance}} of $X_{1}$ and $X_{2},$
denoted $\text{cov}\left(X_{1},X_{2}\right),$ is the real number
defined by\\
\boxeq{
\[
\text{cov}\left(X_{1},X_{2}\right)=\mathbb{E}\left(\left(X_{1}-\mathbb{\mathbb{E}}\left(X_{1}\right)\right)\left(X_{2}-\mathbb{\mathbb{E}}\left(X_{2}\right)\right)\right).
\]
}
\item The \textbf{correlation coefficient\index{correlation coefficient}}
of $X_{1}$ and $X_{2},$ denoted $\rho_{X_{1},X_{2}}$ is then defined
by\\
\boxeq{
\[
\rho_{X_{1},X_{2}}=\dfrac{\text{cov}\left(X_{1},X_{2}\right)}{\sigma_{X_{1}}\sigma_{X_{2}}}.
\]
}
\end{itemize}
\end{definition}

\begin{remark}{}{}

Under these assumptions,$\begin{array}{c}
\sigma_{X_{1}}\neq0\end{array}$ and $\sigma_{X_{2}}\neq0.$

\end{remark}

\begin{proposition}{Correlation Coefficient Values}{}

Let $X=\left(X_{1},X_{2}\right)$ be a random variable taking values
in $\mathbb{R}^{2}.$ Suppose $X$ admits a density $f_{X},$ with
both $X_{1}$ and $X_{2}$ admitting moments of order two. Then the
correlation coefficient $\rho_{X_{1},X_{2}}$ satisfies\boxeq{
\[
\left|\rho_{X_{1},X_{2}}\right|\leqslant1.
\]
}

\end{proposition}

\begin{proof}{}{}

After ensuring that the inequality of Cauchy-Schwarz is still valid,
we proceed as in the discrete case. 

\end{proof}

\section{Laplace-Gauss Law in Two Dimensions}\label{sec:Laplace-Gauss-Law-in}

We apply these concepts to the study of the Laplace-Gauss law in two
dimensions.

Recall that an arbitrary quadratic form $q$ on $\mathbb{R}^{2}$
is written\boxeq{
\[
q\left(x_{1},x_{2}\right)=ax^{2}_{1}+2bx_{1}x_{2}+cx^{2}_{2},
\]
} and that, if $a\neq0$, its square decomposition is written\boxeq{
\[
q\left(x_{1},x_{2}\right)=a\left(x_{1}+\dfrac{b}{a}x_{2}\right)^{2}+x^{2}_{2}\left(c-\dfrac{b^{2}}{a}\right).
\]
}

We deduce from this that, for $q$ to be nonnegative definite, it
is necessary and sufficient that\boxeq{
\[
\begin{array}{ccccc}
a>0, & \,\,\,\, & c>0, & \,\,\,\, & c-\dfrac{b^{2}}{a}>0.\end{array}
\]
}

We assume this condition to be satisfied here after. It then follows
using the square decomposition, that
\[
I\equiv\intop_{\mathbb{R}^{2}}\exp\left(-\dfrac{1}{2}q\left(x_{1},x_{2}\right)\right)\text{d}x_{1}\text{d}x_{2}<+\infty.
\]

\begin{definition}{Laplace-Gauss---or Normal---Centered Random Variable in $\mathbb{R}^2$}{}

A random variable $X=\left(X_{1},X_{2}\right)$ taking values in $\mathbb{R}^{2}$
is called a \textbf{Lapace-Gauss random variable}\index{Lapace-Gauss random variable}\mindex{random variable!Laplace-Gauss}---or
\textbf{normal random variable}\index{normal random variable}\mindex{random variable!normal}---centered
in two-dimensions if it admits a density $f_{X}$ of the form

\boxeq{
\[
f_{X}\left(x_{1},x_{2}\right)=K\text{e}^{-\frac{1}{2}q\left(x_{1},x_{2}\right)},
\]
where
\[
K=\dfrac{1}{\intop_{\mathbb{R}^{2}}\text{e}^{-\frac{1}{2}q\left(x_{1},x_{2}\right)}\text{d}x_{1}\text{d}x_{2}}
\]
}

and $q$ is a nonnegative definite quadratic form on $\mathbb{R}^{2}.$

\end{definition}

We now compute, as an example, the moments of order one and two of
such a random variable of Laplace-Gauss in two dimensions and express
its density with the help of these moments.

\subsubsection*{Computation of $K$}

We denote $I=\intop_{\mathbb{R}^{2}}\text{e}^{-\frac{1}{2}q\left(x_{1},x_{2}\right)}\text{d}x_{1}\text{d}x_{2}.$
We use the above square decomposition of $q\left(x_{1},x_{2}\right).$

Let $m$ be a real number and $\sigma>0.$ 

Recall that
\begin{equation}
\intop^{+\infty}_{-\infty}\dfrac{1}{\sigma\sqrt{2\pi}}\text{e}^{-\frac{\left(x-m\right)^{2}}{2\sigma^{2}}}\textrm{d}x=1,\label{eq:full_density_normal}
\end{equation}
By applying the Fubini theorem and using the equality $\refpar{eq:full_density_normal},$
as well as the square decomposition of a quadratic form, we obtain
\begin{align*}
I & =\intop_{\mathbb{R}^{2}}\text{e}^{-\frac{1}{2}q\left(x_{1},x_{2}\right)}\text{d}x_{1}\text{d}x_{2}\\
 & =\intop^{+\infty}_{-\infty}\text{e}^{-\frac{1}{2}x^{2}_{2}\left(c-\frac{b^{2}}{a}\right)}\left(\intop^{+\infty}_{-\infty}\text{e}^{-\frac{a}{2}\left(x_{1}+\frac{b}{a}x_{2}\right)^{2}}\text{d}x_{1}\right)\text{d}x_{2}\\
 & =\dfrac{\sqrt{2\pi}}{\sqrt{a}}\intop^{+\infty}_{-\infty}\text{e}^{-\frac{1}{2}x^{2}_{2}\left(c-\frac{b^{2}}{a}\right)}\text{d}x_{2}\\
 & =\sqrt{\dfrac{2\pi}{a\left(c-\dfrac{b^{2}}{a}\right)}.}
\end{align*}

Thus,\\
\boxeq{
\[
K=\sqrt{\dfrac{ac-b^{2}}{2\pi}}.
\]
}

\subsubsection*{Computation of the Expectation $\mathbb{E}\left(X_{1}\right)$}

We have
\[
\mathbb{E}\left(X_{1}\right)=K\intop_{\mathbb{R}^{2}}x_{1}\text{e}^{-\frac{1}{2}q\left(x_{1},x_{2}\right)}\text{d}x_{1}\text{d}x_{2}.
\]

We define $F$ the function from $\mathbb{R}^{2}$ to $\mathbb{R},$
by
\[
F\left(x_{1},x_{2}\right)=x_{1}\text{e}^{-\frac{1}{2}q\left(x_{1},x_{2}\right)}.
\]

We observe that for every $\left(x_{1},x_{2}\right)\in\left(\mathbb{R}^{+}\right)^{2},$
\[
F\left(-x_{1},-x_{2}\right)=-F\left(x_{1},x_{2}\right),
\]
 and
\[
F\left(-x_{1},x_{2}\right)=-F\left(x_{1},-x_{2}\right).
\]

It follows that\\
\boxeq{
\[
\mathbb{E}\left(X_{1}\right)=0.
\]
}

By a similar argument,\\
\boxeq{
\[
\mathbb{E}\left(X_{2}\right)=0.
\]
}

\subsubsection*{Computation of $\mathbb{E}\left(X_{1}X_{2}\right)$}

We have
\[
\mathbb{E}\left(X_{1}X_{2}\right)=K\intop_{\mathbb{R}^{2}}x_{1}x_{2}\text{e}^{-\frac{1}{2}q\left(x_{1},x_{2}\right)}\text{d}x_{1}\text{d}x_{2}.
\]

The reduction of the quadratic form $q$ suggests the following change
of variables on $\mathbb{R}^{2},$
\[
\left\{ \substack{u=\sqrt{a}\left(x_{1}+\dfrac{b}{a}x_{2}\right)\\
v=x_{2}\sqrt{c-\dfrac{b^{2}}{a}}
}
\right.\Leftrightarrow\left\{ \substack{x_{1}=\dfrac{u}{\sqrt{a}}-\dfrac{bv}{a\sqrt{c-\dfrac{b^{2}}{a}}}\\
x_{2}=\dfrac{v}{\sqrt{c-\dfrac{b^{2}}{a}}}
}
\right..
\]

The Jacobian of this diffeomorphism is
\begin{align*}
J & =\left|\begin{array}{cc}
\dfrac{\partial x_{1}\left(u,v\right)}{\partial u} & \dfrac{\partial x_{2}\left(u,v\right)}{\partial v}\\
\dfrac{\partial x_{2}\left(u,v\right)}{\partial u} & \dfrac{\partial x_{2}\left(u,v\right)}{\partial v}
\end{array}\right|=\dfrac{1}{\sqrt{ac-b^{2}}}.
\end{align*}

Thus,
\[
\mathbb{E}\left(X_{1}X_{2}\right)=\dfrac{1}{2\pi}\intop_{\mathbb{R}^{2}}\left(\dfrac{u}{\sqrt{a}}-\dfrac{bv}{a\sqrt{c-\dfrac{b^{2}}{a}}}\right)\left(\dfrac{v}{\sqrt{c-\dfrac{b^{2}}{a}}}\right)\text{e}^{-\frac{u^{2}+v^{2}}{2}}\text{d}u\text{d}v.
\]

By using a symmetry argument as before
\[
\mathbb{E}\left(X_{1}X_{2}\right)=\dfrac{1}{2\pi}\dfrac{-b}{ac-b^{2}}\intop_{\mathbb{R}^{2}}v^{2}\text{e}^{-\frac{u^{2}+v^{2}}{2}}\text{d}u\text{d}v.
\]

Since
\[
\begin{array}{ccccc}
\intop^{+\infty}_{-\infty}\text{e}^{-\frac{x^{2}}{2}}\text{d}x=\sqrt{2\pi} & \,\,\,\, & \text{and} & \,\,\,\, & \intop^{+\infty}_{-\infty}x^{2}\text{e}^{-\frac{x^{2}}{2}}\text{d}x=\sqrt{2\pi},\end{array}
\]
and applying the Fubini theorem, we obtain, omitting the details,\\
\boxeq{
\[
\text{cov}\left(X_{1},X_{2}\right)=\mathbb{E}\left(X_{1}X_{2}\right)=\dfrac{-b}{ac-b^{2}}.
\]
}

\subsubsection*{Computation of $\mathbb{E}\left(X^{2}_{1}\right)$}

We have
\[
\mathbb{E}\left(X^{2}_{1}\right)=K\intop_{\mathbb{R}^{2}}x^{2}_{1}\text{e}^{-\frac{1}{2}q\left(x_{1},x_{2}\right)}\text{d}x_{1}\text{d}x_{2}.
\]

Using the same change of variables and the same reasoning as before,
we can write
\[
\mathbb{E}\left(X^{2}_{1}\right)=\dfrac{1}{2\pi}\intop_{\mathbb{R}^{2}}\left(\dfrac{u}{\sqrt{a}}-\dfrac{bv}{a\sqrt{c-\dfrac{b^{2}}{a}}}\right)^{2}\text{e}^{-\frac{u^{2}+v^{2}}{2}}\text{d}u\text{d}v.
\]

Expanding and simplifying the integrand yields
\begin{alignat*}{1}
\mathbb{E}\left(X^{2}_{1}\right) & =\dfrac{1}{a\sqrt{2\pi}}\intop_{\mathbb{R}^{2}}u^{2}\text{e}^{-\frac{u^{2}}{2}}\text{d}u+\dfrac{1}{\sqrt{2\pi}}\dfrac{b^{2}}{a\left(ac-b^{2}\right)}\intop_{\mathbb{R}^{2}}v^{2}\text{e}^{-\frac{v^{2}}{2}}\text{d}v\\
 & =\dfrac{1}{a}+\dfrac{b^{2}}{a\left(ac-b^{2}\right)}.
\end{alignat*}
Hence,\\
\boxeq{
\[
\mathbb{E}\left(X^{2}_{1}\right)=\dfrac{c}{ac-b^{2}}
\]
}

We therefore find the variance of $X_{1},$\\
\boxeq{
\[
\sigma^{2}_{X_{1}}\equiv\mathbb{E}\left(X^{2}_{1}\right)=\dfrac{c}{ac-b^{2}}.
\]
}

Similarly, the variance of $X_{2}$ is\\
\boxeq{
\[
\sigma^{2}_{X_{2}}\equiv\mathbb{E}\left(X^{2}_{2}\right)=\dfrac{a}{ac-b^{2}}.
\]
}

The correlation coefficient of $X_{1}$ and $X_{2}$ is thus given
by\\
\boxeq{
\[
\rho_{X_{1},X_{2}}\equiv\dfrac{\text{cov\ensuremath{\left(X_{1},X_{2}\right)}}}{\sigma_{X_{1}}\sigma_{X_{2}}}=-\dfrac{b}{\sqrt{ac}}.
\]
}

\subsubsection*{Density of a Centered Laplace-Gauss Random Variable in Two Dimensions}

We observe that
\[
1-\rho^{2}_{X_{1},X_{2}}=\dfrac{ac-b^{2}}{ac}.
\]

A straightforward algebraic computation leads to the following expression
for the quadratic form $q$ written in terms of the variances and
the correlation coefficient of $X_{1}$ and $X_{2}.$ 

For every $\left(x_{1},x_{2}\right)\in\mathbb{R}^{2},$
\[
q\left(x_{1},x_{2}\right)=\dfrac{1}{1-\rho^{2}_{X_{1},X_{2}}}\left(\left(\dfrac{x_{1}}{\sigma_{X_{1}}}\right)^{2}-2\rho_{X_{1},X_{2}}\dfrac{x_{1}}{\sigma_{X_{1}}}\dfrac{x_{2}}{\sigma_{X_{2}}}+\left(\dfrac{x_{2}}{\sigma_{X_{2}}}\right)^{2}\right).
\]

Moreover,
\[
K=\dfrac{1}{2\pi\sigma_{X_{1}}\sigma_{X_{2}}\sqrt{1-\rho^{2}_{X_{1},X_{2}}}}
\]
which gives the expression for \textbf{the density of a two-dimensional
centered Laplace-Gauss random variable}.

For every $\left(x_{1},x_{2}\right)\in\mathbb{R}^{2},$\\
\boxeq{
\[
f_{X}\left(x_{1},x_{2}\right)=\dfrac{\exp\left[-\dfrac{1}{2}\dfrac{1}{1-\rho^{2}_{X_{1},X_{2}}}\left(\left(\dfrac{x_{1}}{\sigma_{X_{1}}}\right)^{2}-2\rho_{X_{1},X_{2}}\dfrac{x_{1}}{\sigma_{X_{1}}}\dfrac{x_{2}}{\sigma_{X_{2}}}+\left(\dfrac{x_{2}}{\sigma_{X_{2}}}\right)^{2}\right)\right]}{2\pi\sigma_{X_{1}}\sigma_{X_{2}}\sqrt{1-\rho^{2}_{X_{1},X_{2}}}}
\]
}

\section{Independence of Two Real Random Variables}\label{sec:Independence-of-Two}

In this section, $X_{1}$ and $X_{2}$ are two real-valued random
variables.

We denote $X=\left(X_{1},X_{2}\right)$ the vector random variable
taking values in $\mathbb{R}^{2}$ associated with the two random
variables $X_{1}$ and $X_{2}.$ We study the criterion for independence
of the random variables $X_{1}$ and $X_{2}.$

\begin{proposition}{Independence Criterion of Two Real Random Variables}{}

(i) The random variables $X_{1}$ and $X_{2}$ are independent, if
and only if, for every $\left(x_{1},x_{2}\right)\in\mathbb{R}^{2},$\\
\boxeq{
\begin{equation}
P\left(X_{1}\leqslant x_{1},X_{2}\leqslant x_{2}\right)=P\left(X_{1}\leqslant x_{1}\right)P\left(X_{2}\leqslant x_{2}\right).\label{eq:independence_criterion}
\end{equation}
}

(ii) If the random variables $X_{1}$ and $X_{2}$ are independent
and admit as densities $f_{X_{1}}$ and $f_{X_{2}},$ then the random
variable $X=\left(X_{1},X_{2}\right)$ also admits a density $f_{X},$
function defined on $\mathbb{R}^{2}$ by\\
\boxeq{
\[
f_{X}\left(x_{1},x_{2}\right)=f_{X_{1}}\left(x_{1}\right)f_{X_{2}}\left(x_{2}\right).
\]
}

(iii) Conversely, if the random variable $X=\left(X_{1},X_{2}\right)$
admits a density $f_{X}$ that satisfies, for every $\left(x_{1},x_{2}\right)\in\mathbb{R}^{2},$
\[
f_{X}\left(x_{1},x_{2}\right)=g\left(x_{1}\right)h\left(x_{2}\right),
\]
where $g$ and $h$ are some nonnegative functions that are Riemann-integrable
over $\mathbb{R},$ then the random variables $X_{1}$ and $X_{2}$
are independent. Furthermore $g$ and $h$ are, up to a nonnegative
multiplicative coefficient, the respective densities of $X_{1}$ and
$X_{2}.$ 

\end{proposition}

\begin{proof}{}{}

(i) The direct implication follows directly from the definition of
the independence of random variables. 

The converse is accepted without proof, as it relies on results from
measure theory.

(ii) Under this hypothesis, the relation $\refpar{eq:independence_criterion}$
is written, for every real number $a,b,c,d$ with $a<b$ and $c<d,$
\[
P\left(a<X_{1}\leqslant b,c<X_{2}\leqslant d\right)=\left(\intop^{b}_{a}f_{X_{1}}\left(u_{1}\right)\text{d}u_{1}\right)\left(\intop^{d}_{c}f_{X_{2}}\left(u_{2}\right)\text{d}u_{2}\right).
\]

By applying the Fubini theorem and noting that the function 
\[
\left(x_{1},x_{2}\right)\mapsto f_{X_{1}}\left(x_{1}\right)f_{X_{2}}\left(x_{2}\right)
\]
 is regular, we obtain
\[
P\left(a<X_{1}\leqslant b,c<X_{2}\leqslant d\right)=\intop^{b}_{a}\intop^{d}_{c}f_{X_{1}}\left(u_{1}\right)f_{X_{2}}\left(u_{2}\right)\text{d}u_{1}\text{d}u_{2}.
\]

This last equality shows the result.

(iii) The random variable $X=\left(X_{1},X_{2}\right)$ admits a density
$f_{X},$ satisfying for every $\left(x_{1},x_{2}\right)\in\mathbb{R}^{2},$
\[
f_{X}\left(x_{1},x_{2}\right)=g\left(x_{1}\right)h\left(x_{2}\right),
\]
which by hypothesis is regular. 

Then the marginals $X_{1}$ and $X_{2}$ also admit densities, given
respectively on $\mathbb{R}$ by
\[
f_{X_{1}}\left(x_{1}\right)=\intop^{+\infty}_{-\infty}f_{X}\left(x_{1},x_{2}\right)\text{d}x_{2},
\]
and
\[
f_{X_{2}}\left(x_{2}\right)=\intop^{+\infty}_{-\infty}f_{X}\left(x_{1},x_{2}\right)\text{d}x_{1}.
\]

From the hypothesis, we deduce that for every $x_{1}\in\mathbb{R},$
\[
f_{X_{1}}\left(x_{1}\right)=g\left(x_{1}\right)\intop^{+\infty}_{-\infty}h\left(x_{2}\right)\text{d}x_{2},
\]
and for every $x_{2}\in\mathbb{R},$
\[
f_{X_{2}}\left(x_{2}\right)=h\left(x_{2}\right)\intop^{+\infty}_{-\infty}g\left(x_{1}\right)\text{d}x_{1}.
\]

Since we have on the one hand
\[
\intop^{+\infty}_{-\infty}f_{X_{1}}\left(x_{1}\right)\text{d}x_{1}=1
\]
and, on the other hand
\begin{align*}
\intop^{+\infty}_{-\infty}f_{X_{1}}\left(x_{1}\right)\text{d}x_{1} & =\intop^{+\infty}_{-\infty}g\left(x_{1}\right)\left(\intop^{+\infty}_{-\infty}h\left(x_{2}\right)\text{d}x_{2}\right)\text{d}x_{1}\\
 & =\left(\intop^{+\infty}_{-\infty}h\left(x_{2}\right)\text{d}x_{2}\right)\left(\intop^{+\infty}_{-\infty}g\left(x_{1}\right)\text{d}x_{1}\right),
\end{align*}
we have
\[
1=\left(\intop^{+\infty}_{-\infty}h\left(x_{2}\right)\text{d}x_{2}\right)\left(\intop^{+\infty}_{-\infty}g\left(x_{1}\right)\text{d}x_{1}\right).
\]

Using this last relation, we deduce, that for every $\left(x_{1},x_{2}\right)\in\mathbb{R}^{2},$
\begin{align*}
f_{X}\left(x_{1},x_{2}\right) & =g\left(x_{1}\right)h\left(x_{2}\right)\\
 & =g\left(x_{1}\right)h\left(x_{2}\right)\times1\\
 & =g\left(x_{1}\right)h\left(x_{2}\right)\left(\intop^{+\infty}_{-\infty}h\left(x_{2}\right)\text{d}x_{2}\right)\left(\intop^{+\infty}_{-\infty}g\left(x_{1}\right)\text{d}x_{1}\right)\\
 & =\left(g\left(x_{1}\right)\intop^{+\infty}_{-\infty}h\left(x_{2}\right)\text{d}x_{2}\right)\left(h\left(x_{2}\right)\intop^{+\infty}_{-\infty}g\left(x_{1}\right)\text{d}x_{1}\right)\\
 & =f_{X_{1}}\left(x_{1}\right)f_{X_{2}}\left(x_{2}\right).
\end{align*}

We then have, for every $\left(x_{1},x_{2}\right)\in\mathbb{R}^{2},$
\begin{align*}
P\left(X_{1}\leqslant x_{1},X_{2}\leqslant x_{2}\right) & =\intop^{x_{1}}_{-\infty}\intop^{x_{2}}_{-\infty}f_{X}\left(u_{1},u_{2}\right)\text{d}u_{1}\text{d}u_{2}\\
 & =\intop^{x_{1}}_{-\infty}\intop^{x_{2}}_{-\infty}f_{X_{1}}\left(u_{1}\right)f_{X_{2}}\left(u_{2}\right)\text{d}u_{1}\text{d}u_{2}\\
 & =\left(\intop^{x_{1}}_{-\infty}f_{X_{1}}\left(u_{1}\right)\text{d}u_{1}\right)\left(\intop^{x_{2}}_{-\infty}f_{X_{2}}\left(u_{2}\right)\text{d}u_{2}\right)\\
 & =P\left(X_{1}\leqslant x_{1}\right)P\left(X_{2}\leqslant x_{2}\right).
\end{align*}

This final equality proves the independence of the random variables
$X_{1}$ and $X_{2}.$

\end{proof}

\begin{corollary}{Expectation and Variance of Independent Random Variables with Density}{}

Let $X_{1}$ and $X_{2}$ be two independent real-valued random variables
with densities.

(i) If $\phi_{1}$ and $\phi_{2}$ are two real-valued functions such
that $\phi_{1}\left(X_{1}\right)$ and $\phi_{2}\left(X_{2}\right)$
admit an expectation, then the product random variable $\phi_{1}\left(X_{1}\right)\phi_{2}\left(X_{2}\right)$
admits also an expectation, and\\
\boxeq{
\begin{equation}
\mathbb{E}\left(\phi_{1}\left(X_{1}\right)\phi_{2}\left(X_{2}\right)\right)=\mathbb{E}\left(\phi_{1}\left(X_{1}\right)\right)\mathbb{E}\left(\phi_{2}\left(X_{2}\right)\right).\label{eq:mean_product_independent}
\end{equation}
}

(ii) Moreover, if the random variables $X_{1}$ and $X_{2}$ both
admit a variance, then\\
\boxeq{
\[
\text{cov}\left(X_{1},X_{2}\right)=0
\]
}and, consequently,\\
\boxeq{
\begin{equation}
\sigma^{2}_{X_{1}+X_{2}}=\sigma^{2}_{X_{1}}+\sigma^{2}_{X_{2}}.\label{eq:covariance_sum}
\end{equation}
}

\end{corollary}

\begin{remark}{}{}

The converse of (ii) is false: it is possible for two random variables
$X_{1}$ and $X_{2}$ to satisfy
\[
\text{cov}\left(X_{1},X_{2}\right)=0
\]
without being independent. 

\end{remark}

\begin{proof}{}{}

\textbf{(i) Expectation of $\phi_{1}\left(X_{1}\right)\phi_{2}\left(X_{2}\right)$}

The random variable $X=\left(X_{1},X_{2}\right)$ admits a density
$f_{X}$ direct product of $f_{X_{1}}$ and $f_{X_{2}}.$

Since the random variables $\phi_{1}\left(X_{1}\right)$ and $\phi_{2}\left(X_{2}\right)$
admit an expectation, we have for each $i=1,2,$
\[
\intop^{+\infty}_{-\infty}\left|\phi_{i}\left(u_{i}\right)\right|f_{X_{i}}\left(u_{i}\right)\text{d}u_{i}<+\infty
\]

Thus, the product of this integrals is also finite\boxeq{
\[
\left(\intop^{+\infty}_{-\infty}\left|\phi_{1}\left(u_{1}\right)\right|f_{X_{1}}\left(u_{1}\right)\text{d}u_{1}\right)\left(\intop^{+\infty}_{-\infty}\left|\phi_{2}\left(u_{2}\right)\right|f_{X_{2}}\left(u_{2}\right)\text{d}u_{2}\right)<+\infty
\]
}

By the Fubini theorem, this implies
\begin{align*}
 & \left(\intop^{+\infty}_{-\infty}\left|\phi_{1}\left(u_{1}\right)\right|f_{X_{1}}\left(u_{1}\right)\text{d}u_{1}\right)\left(\intop^{+\infty}_{-\infty}\left|\phi_{2}\left(u_{2}\right)\right|f_{X_{2}}\left(u_{2}\right)\text{d}u_{2}\right)\\
 & \qquad\qquad=\intop^{+\infty}_{-\infty}\intop^{+\infty}_{-\infty}\left|\phi_{1}\left(u_{1}\right)\right|\left|\phi_{2}\left(u_{2}\right)\right|f_{X_{1}}\left(u_{1}\right)f_{X_{2}}\left(u_{2}\right)\text{d}u_{1}\text{d}u_{2}<+\infty
\end{align*}
and consequently the random variable $\left|\phi_{1}\left(X_{1}\right)\phi_{2}\left(X_{2}\right)\right|$
admits an expectation.

Applying the Fubini theorem without absolute values yields the relation
$\refpar{eq:mean_product_independent}.$

\textbf{(ii) Covariance of $X_{1}$ and $X_{2}$} 

In particular, we obtain using (i),
\[
\mathbb{E}\left(X_{1}X_{2}\right)=\mathbb{E}\left(X_{1}\right)\mathbb{E}\left(X_{2}\right),
\]
which implies
\[
\text{cov}\left(X_{1},X_{2}\right)=0.
\]

We then obtain the equality $\refpar{eq:covariance_sum}$ by expanding
the square $\left(\mathring{X_{1}}+\mathring{X_{2}}\right)^{2}$ and
by using the linearity of the expectation. 

\end{proof}

\section{Sum of Independent Real Random Variables}\label{sec:Sum-of-Independent}

We study the law of the sum of two independent random variables that
admit a density.

\begin{proposition}{Density of a Sum of Random Variables}{}

Let $X_{1}$ and $X_{2}$ be two independent real-valued random variables
with densities $f_{X_{1}}$ and $f_{X_{2}}$ respectively. Then the
random variable $X=X_{1}+X_{2}$ also admits a density $f_{X}$ given
by the \textbf{\index{convolution} convolution} of the functions
$f_{X_{1}}$ and $f_{X_{2}}.$ That is, the function defined, for
every $x\in\mathbb{R},$ by\\
\boxeq{
\begin{equation}
f_{X}\left(x\right)=\intop^{+\infty}_{-\infty}f_{X_{1}}\left(x_{1}\right)f_{X_{2}}\left(x-x_{1}\right)\text{d}x_{1}=\intop^{+\infty}_{-\infty}f_{X_{1}}\left(x-x_{2}\right)f_{X_{2}}\left(x_{2}\right)\text{d}x_{2}.\label{eq:f_X_x}
\end{equation}
}

\end{proposition}

\begin{proof}{}{}

The random variable $\left(X_{1},X_{2}\right)$ admits a density $f_{\left(X_{1},X_{2}\right)}$
which is the direct product of the densities $f_{X_{1}}$ and $f_{X_{2}}.$
Therefore, for every $x\in\mathbb{R},$
\[
P\left(X_{1}+X_{2}\leqslant x\right)=\intop^{+\infty}_{-\infty}\intop^{+\infty}_{-\infty}\boldsymbol{1}_{D_{x}}\left(x_{1},x_{2}\right)f_{X_{1}}\left(x_{1}\right)f_{X_{2}}\left(x_{2}\right)\text{d}x_{1}\text{d}x_{2},
\]
where $D_{x}$ is the subset of $\mathbb{R}^{2}$ defined by
\[
D_{x}=\left\{ \left(x_{1},x_{2}\right)\in\mathbb{R}^{2}:\,x_{1}+x_{2}\leqslant x\right\} .
\]

Now, perform the change of variables on $\mathbb{R}^{2}$ onto itself,
which has Jacobian 1, defined by
\[
\left\{ \begin{array}{c}
u=x_{1}\\
v=x_{1}+x_{2}
\end{array}\right.\Leftrightarrow\left\{ \begin{array}{c}
x_{1}=u\\
x_{2}=v-u.
\end{array}\right.
\]
It follows, that for every $x\in\mathbb{R},$
\[
P\left(X_{1}+X_{2}\leqslant x\right)=\intop^{+\infty}_{-\infty}\intop^{+\infty}_{-\infty}\boldsymbol{1}_{\left]-\infty,x\right]}\left(v\right)f_{X_{1}}\left(u\right)f_{X_{2}}\left(v-u\right)\text{d}u\text{d}v.
\]

Applying the Fubini theorem, we have, for every $x\in\mathbb{R},$
\[
P\left(X_{1}+X_{2}\leqslant x\right)=\intop^{x}_{-\infty}\left(\intop^{+\infty}_{-\infty}f_{X_{1}}\left(u\right)f_{X_{2}}\left(v-u\right)\text{d}u\right)\text{d}v,
\]
which proves the result.

\end{proof}

\begin{example}{}{}

Let $X_{1}$ and $X_{2}$ be two independent real-valued random variables,
each following the same exponential law with parameter $p>0.$

Determine the law of the sum $X=X_{1}+X_{2}$ by finding its density
$f_{X}$---the random variable $X$ can, for instance, represent
the sum of two independent waiting times, each following the same
exponential law.

\end{example}

\begin{solutionexample}{}{}

We recall that for every $x\in\mathbb{R},$
\[
f_{X_{1}}\left(x\right)=f_{X_{2}}\left(x\right)=\boldsymbol{1}_{\mathbb{R}^{+}}\left(x\right)p\text{e}^{-px}.
\]

The random variable $X$ admits a density given, for every $x\in\mathbb{R},$
by
\[
f_{X}\left(x\right)=\intop^{+\infty}_{-\infty}\boldsymbol{1}_{\mathbb{R}^{+}}\left(x_{1}\right)p\text{e}^{-px_{1}}\boldsymbol{1}_{\mathbb{R}^{+}}\left(x-x_{1}\right)p\text{e}^{-p\left(x-x_{1}\right)}\text{d}x_{1}.
\]
Thus,
\[
f_{X}\left(x\right)=p^{2}\boldsymbol{1}_{\mathbb{R}^{+}}\left(x\right)\text{e}^{-px}\intop^{x}_{0}\text{d}x_{1}.
\]

Therefore, the random variable $X$ admits a density $f_{X}$ defined,
for every $x\in\mathbb{R},$ by\boxeq{
\[
f_{X}\left(x\right)=p^{2}\boldsymbol{1}_{\mathbb{R}^{+}}\left(x\right)x\text{e}^{-px}.
\]
}

\end{solutionexample}

\section{Conditional Densities}\label{sec:Conditional-Densities}

Here, we consider only random variables $\left(X_{1},X_{2}\right)$
taking values in $\mathbb{R}^{2}$ and admitting a density $f_{\left(X_{1},X_{2}\right)}.$
In this case, we know that for every $x_{1}\in\mathbb{R},$
\[
P\left(X_{1}=x_{1}\right)=0,
\]
since $X_{1}$ admits a density. It is therefore impossible to define,
as it was done in the discrete case, a conditional probability such
as $P\left(X_{2}\leqslant x_{2}\left|X_{1}=x_{1}\right.\right).$
Nonetheless, by analogy we define a notion of \textbf{\index{conditional density}conditional
density}\footnote{The general theory of conditional probabilities will be presented
in the second part.}.

\begin{definition}{Conditional Density of a Random Variable}{}

Let $\left(X_{1},X_{2}\right)$ be a random variable taking values
in $\mathbb{R}^{2}$ and admitting a density $f_{\left(X_{1},X_{2}\right)}$
that is regular\footnotemark. 

We define, for every $x_{1}\in\mathbb{R}$ such that $f_{X_{1}}\left(x_{1}\right)\neq0,$
the \textbf{conditional density of $X_{2}$ knowing\index{conditional density of $X_{2}$ knowing $X_{1}=x_{1}$}
$X_{1}=x_{1},$ }denoted $f^{X_{1}=x_{1}}_{X_{2}},$ as the function
defined on $\mathbb{R},$ for every $x_{2}\in\mathbb{R},$ by\\
\boxeq{
\begin{equation}
f^{X_{1}=x_{1}}_{X_{2}}\left(x_{2}\right)=\dfrac{f_{\left(X_{1},X_{2}\right)}\left(x_{1},x_{2}\right)}{f_{X_{1}}\left(x_{1}\right)}.\label{eq:conditional_density}
\end{equation}
}

\end{definition}

\footnotetext{See the footnote of Proposition $\ref{pr:density_marginals}$}

\begin{proposition}{}{}

Under the assumptions of the previous definition, and if the support
of $f_{X_{1}}$---that is the closure of the set $\left\{ x_{1}\in\mathbb{R}:\,f_{X_{1}}\left(x_{1}\right)\neq0\right\} $---is
an interval $I,$ then we have, for every real numbers $a_{1},\,a_{2},\,b_{1},\,b_{2}$
such that $a_{1}<b_{1},$ $\left[a_{1},\,b_{1}\right]\subset I,$
and $a_{2}<b_{2},$\\
\boxeq{
\begin{equation}
P\left(a_{1}<X_{1}\leqslant b_{1},a_{2}<X_{2}\leqslant b_{2}\right)=\intop^{b_{1}}_{a_{1}}\left(\intop^{b_{2}}_{a_{2}}f^{X_{1}=x_{1}}_{X_{2}}\left(x_{2}\right)\text{d}x_{2}\right)f_{X_{1}}\left(x_{1}\right)\text{d}x_{1}.\label{eq:proba_intersection_dens}
\end{equation}
}

\end{proposition}

\begin{proof}{}{}

Since, for every $x_{1}\in\left[a_{1},b_{1}\right],$ we have $f_{X_{1}}\left(x_{1}\right)\neq0,$
it follows from the Fubini theorem and from the relation $\refpar{eq:conditional_density}$
that 
\begin{align*}
P\left(a_{1}<X_{1}\leqslant b_{1},a_{2}<X_{2}\leqslant b_{2}\right) & =\intop^{b_{1}}_{a_{1}}\left(\intop^{b_{2}}_{a_{2}}f_{\left(X_{1},X_{2}\right)}\left(x_{1},x_{2}\right)\text{d}x_{2}\right)\text{d}x_{1}\\
 & =\intop^{b_{1}}_{a_{1}}\left(\intop^{b_{2}}_{a_{2}}f^{X_{1}=x_{1}}_{X_{2}}\left(x_{2}\right)\text{d}x_{2}\right)f_{X_{1}}\left(x_{1}\right)\text{d}x_{1}.
\end{align*}

\end{proof}

Conversely, suppose that $X_{1}$ admits a density $f_{X_{1}}$ whose
support is an interval $I,$ and that there exists a \textbf{family
of densities\mindex{family!densities}} $\left(f^{x}\left(\cdot\right)\right)_{x\in I}$
such that, for every real number $a_{1},a_{2},b_{1},b_{2}$ with $a_{1}<b_{1},$
$\left[a_{1},b_{1}\right]\subset I,$ and $a_{2}<b_{2},$ the following
equality holds
\[
P\left(a_{1}<X_{1}\leqslant b_{1},a_{2}<X_{2}\leqslant b_{2}\right)=\intop^{b_{1}}_{a_{1}}\left(\intop^{b_{2}}_{a_{2}}f^{x_{1}}\left(x_{2}\right)\text{d}x_{2}\right)f_{X_{1}}\left(x_{1}\right)\text{d}x_{1},
\]
then the random variable $\left(X_{1},X_{2}\right)$ admits a density
$f_{\left(X_{1},X_{2}\right)}$ whose support is contained in $I\times\mathbb{R},$
and for every $\left(x_{1},x_{2}\right)\in I\times\mathbb{R},$\boxeq{
\[
f_{\left(X_{1},X_{2}\right)}\left(x_{1},x_{2}\right)=f_{X_{1}}\left(x_{1}\right)f^{x_{1}}\left(x_{2}\right).
\]
}

In this case, the random variable $X_{2}$ admits a \textbf{\index{conditional density knowing}conditional
density given} $X_{1}=x_{1},$ which is exactly $f^{x_{1}}.$

\begin{example}{}{}

Let $X=\left(X_{1},X_{2}\right)$ be a random variable following the
centered Laplace-Gauss law in two dimensions. We recall that its density
$f_{X}$ is defined for every $\left(x_{1},x_{2}\right)\in\mathbb{R}^{2}$
by
\[
f_{X}\left(x_{1},x_{2}\right)=K\exp\left(-\dfrac{1}{2}q\left(x_{1},x_{2}\right)\right),
\]
where the quadratic form $q$ is defined, for every $\left(x_{1},x_{2}\right)\in\mathbb{R}^{2},$
by
\[
q\left(x_{1},x_{2}\right)=\dfrac{1}{1-\rho^{2}}\left(\left(\dfrac{x_{1}}{\sigma_{1}}\right)^{2}-2\rho\dfrac{x_{1}}{\sigma_{1}}\dfrac{x_{2}}{\sigma_{2}}+\left(\dfrac{x_{2}}{\sigma_{2}}\right)^{2}\right).
\]

The coefficient $K$ is given by
\[
K=\dfrac{1}{2\pi\sigma_{1}\sigma_{2}\sqrt{1-\rho^{2}}}.
\]

The parameters are such that $\sigma_{1}>0,$ $\sigma_{2}>0$ and
$\left|\rho\right|<1.$ 

The quadratic form $q$ can also be written, for every $\left(x_{1},x_{2}\right)\in\mathbb{R}^{2},$
as
\[
q\left(x_{1},x_{2}\right)=\dfrac{1}{1-\rho^{2}}\left(\left(\dfrac{x_{2}}{\sigma_{2}}-\rho\dfrac{x_{1}}{\sigma_{1}}\right)^{2}+\left(1-\rho^{2}\right)\left(\dfrac{x_{1}}{\sigma_{1}}\right)^{2}\right).
\]

It follows that the marginal $X_{1}$ admits a density $f_{X_{1}}$
given, for every $x_{1}\in\mathbb{R},$ by
\begin{align*}
f_{X_{1}}\left(x_{1}\right) & =K\exp\left(-\dfrac{1}{2}\left(\dfrac{x_{1}}{\sigma_{1}}\right)^{2}\right)\intop^{+\infty}_{-\infty}\exp\left(-\dfrac{1}{2\left(1-\rho^{2}\right)\sigma^{2}_{2}}\left(x_{2}-\rho\dfrac{\sigma_{2}}{\sigma_{1}}x_{1}\right)^{2}\right)\text{d}x_{1}\\
 & =\dfrac{1}{\sqrt{2\pi}\sigma_{1}}\exp\left(-\dfrac{1}{2}\left(\dfrac{x_{1}}{\sigma_{1}}\right)^{2}\right).
\end{align*}

It follows that the conditional density of $X_{2}$ given $X_{1}=x_{1}$
for every $x_{2}\in\mathbb{R},$ takes the form \boxeq{
\[
f^{X_{1}=x_{1}}_{X_{2}}=\dfrac{1}{\sqrt{2\pi}\sigma_{2}\sqrt{1-\rho^{2}}}\exp\left(-\dfrac{1}{2\left(1-\rho^{2}\right)\sigma^{2}_{2}}\left(x_{2}-\rho\dfrac{\sigma_{2}}{\sigma_{1}}x_{1}\right)^{2}\right),
\]
}that is, for each $x_{1}\in\mathbb{R},$ the function $f^{X_{1}=x_{1}}_{X_{2}}$
is the density of a \textbf{Gaussian law\mindex{Gaussian!law}} on
$\mathbb{R},$ namely
\[
\mathcal{N}\left(\rho\dfrac{\sigma_{2}}{\sigma_{1}}x_{1},\left(1-\rho^{2}\right)\sigma^{2}_{2}\right).
\]

Conversely, suppose that the random variable $X=\left(X_{1},X_{2}\right)$
is such that the law followed by $X_{1}$ is the \textbf{\mindex{Gaussian!law}Gaussian
law} $\mathcal{N}\left(0,\sigma^{2}_{1}\right)$ and, that the random
variable $X_{2}$ admits a conditional density given $X_{1}=x_{1},$
defined for every $x_{2}\in\mathbb{R},$ by
\[
f^{x_{1}}\left(x_{2}\right)=\dfrac{1}{\sqrt{2\pi}a}\exp\left(-\dfrac{1}{2a^{2}}\left(x_{2}-bx_{1}\right)^{2}\right),
\]
where $a>0.$ Then the random variable $X=\left(X_{1},X_{2}\right)$
follows a \textbf{\index{Laplace-Gauss law}centered Laplace-Gauss
law} in two dimensions.

\end{example}

\section{Appendix: the Riemann Integral in $\mathbb{R}^{n}$}\label{sec:Appendix:-the-Riemann}

Tr. Note: The Author of this Appendix is André Bellaïche.

We present here, without proof, the main definitions and the most
important properties of the \textbf{Riemann integral\index{Riemann integral}}
in $\mathbb{R}^{n}.$

\subsection{Definition of the Riemann Integral in $\mathbb{R}^{n}$}

\subsubsection{Integral of a Function on a Closed Rectangle.}

Let 
\[
P=\left[a_{1},b_{1}\right]\times\cdots\times\left[a_{n},b_{n}\right]
\]
be a bounded, closed rectangle in $\mathbb{R}^{n},$ that is a compact
rectangle in $\mathbb{R}^{n}.$

The concept of a Riemann-integrable function on $P,$ and the corresponding
integral, are defined exactely as in dimension 1, with intervals replaced
by rectangles and the length of an interval replaced by the volume
of a rectangle.

We begin by defining the integral of a step function. 

\begin{definition}{Step Function}{}

A \textbf{step function\index{step function}} on a rectangle $P$
is, by definition, a function of the form\\
\boxeq{
\begin{equation}
f=\sum^{k}_{i=1}\lambda_{i}\boldsymbol{1}_{A_{i}},\label{eq:step_function}
\end{equation}
}where the $A_{i}$ are bounded closed rectangles contained in $P.$

\end{definition}

\begin{remark}{}{}

Note that indicator functions of open rectangles, semi-open, etc...
are step functions.

\end{remark}

\begin{definition}{Integral of a Step Function}{}

The \textbf{integral of a step function}\index{integral of a step function}\mindex{step function ! integral},
as defined in $\refpar{eq:step_function},$ is given by\\
\boxeq{
\begin{equation}
\intop_{P}f\left(x\right)\text{d}x=\sum^{k}_{i=1}\lambda_{i}\text{vol}A_{i}.\label{eq:integral_step_function}
\end{equation}
}

\end{definition}

Although the representation of a step function under the form $\refpar{eq:step_function}$
is not unique, the sum on the right-hand side of the equality $\refpar{eq:integral_step_function}$
does not depend on the specific representation chosen.

\begin{definition}{Function Integrable in the Sense of Riemann}{}

We say that a bounded function $f,$ defined on $P,$ is \textbf{integrable
in the sense of Riemann}\index{integrable in the sense of Riemann},
or \textbf{Riemann-integrable}\index{Riemann-integrable}, on $P$
if, for every $\epsilon>0,$ there exists some step functions $\mathbb{\phi}$
and $\psi,$ defined on $P,$ such that
\[
\phi\leqslant f\leqslant\psi
\]
and
\[
\intop_{P}\left(\psi\left(x\right)-\phi\left(x\right)\right)\text{d}x<\epsilon.
\]

If $f$ is Riemann-integrable on $P,$ the numbers
\[
\begin{array}{ccccc}
\underset{\phi\leqslant f}{}\intop_{P}\phi\left(x\right)\textrm{d}x & \,\,\,\, & \text{and} & \,\,\,\, & \underset{\psi\geqslant f}{\text{inf}}\intop_{P}\psi\left(x\right)\textrm{d}x,\end{array}
\]
where $\phi$ and $\psi$ are step functions, are finite and equal.

This value is by definition the \textbf{integral of $f$ on $P.$}
It is denoted
\[
\begin{array}{ccccc}
\intop_{P}f\left(x\right)\textrm{d}x & \,\,\,\, & \text{or} & \,\,\,\, & \intop_{P}f\left(x_{1},\dots,x_{n}\right)\text{d}x_{1}\dots\text{d}x_{n}.\end{array}
\]

\end{definition}

\subsubsection{Properties}

The basic properties of the integral are the same than in dimension
1: the sum or the product of two functions Riemann-integrable are
Riemann-integrable, and so on. 

Any function Riemann-integrable is bounded---by definition. A continuous
function on $P$ is Riemann-integrable. 

A function that is continuous by pieces is Riemann-integrable if the
``pieces'' are sufficiently simple. Such simple pieces may include
for instance, in dimension 2, polygons, disks, and more generally,
sets of the form
\[
\left\{ a\leqslant x\leqslant b,\,\,\,\,g\left(x\right)\leqslant y\leqslant h\left(x\right)\right\} 
\]
where $g$ and $h$ are continuous function on $\left[a,b\right],$
or sets of the form
\[
\left\{ c\leqslant y\leqslant d,\,\,\,\,k\left(y\right)\leqslant x\leqslant l\left(y\right)\right\} 
\]
where $k$ and $l$ are continuous function on $\left[c,d\right],$
and unions of such sets.

\subsubsection{Integral on a Jordan-measurable Bounded Part of $\mathbb{R}^{n}$ }

\begin{definition}{Jordan-measurable Bounded Subset of $\mathbb{R}^n$}{}

Let $A$ be a bounded subset of $\mathbb{R}^{n}.$ Consider a compact
rectangle $P$ containing $A.$ We say that $A$ is \textbf{Jordan-measurable\index{Jordan-measurable}}
if its characteristic function $\boldsymbol{1}_{A}$ is integrable
on $P.$ 

In this case, we define the \textbf{volume of $A$} \index{volume of a bounded subset}by\\
\boxeq{
\[
\text{vol}\left(A\right)=\intop_{P}\boldsymbol{1}_{A}\left(x\right)\textrm{d}x.
\]
}

\end{definition}

\begin{definition}{Riemann-Integrable Function on a Jordan-measurable Subset of $\mathbb{R}^n$}{}

Suppose that $A$ is Jordan-measurable, and let $f$ be a bounded
function defined on $A.$ 

We say that $f$ is \textbf{\index{Riemann-integrable}Riemann-integrable
on $A$} if the function $f^{\ast}$ obtained by extending $f$ by
0 outside of $A$ is integrable on $P.$ In that case, we define\\
\boxeq{
\[
\intop_{A}f\left(x\right)\textrm{d}x=\intop_{P}f^{\ast}\left(x\right)\textrm{d}x.
\]
}

\end{definition}

All these concepts do not depend on the particular choice of the rectangle
$P.$

For instance, in dimension 2, the sets of ``simple shape'' as defined
previously are Jordan-measurable.

Moreover, the interior $\mathring{A}$ and the closure $\overline{A}$
of a Jordan-measurable bounded set are also Jordan-measurable, and
\[
\text{vol}\left(\mathring{A}\right)=\text{vol}\left(A\right)=\text{vol}\left(\overline{A}\right).
\]

Similarly,\\
\boxeq{
\[
\intop_{\mathring{A}}f\left(x\right)\textrm{d}x=\intop_{A}f\left(x\right)\textrm{d}x.
\]
}

That is, in the computation of the integral, the contribution of the
boundary of $A$ can be neglected. This observation is important in
practice to apply the formula of change of variables.

\subsubsection{Integral on $\mathbb{R}^{n}$---Generalized Integral}

\begin{definition}{Riemann-Ingrable Nonnegative Function}{}

Let $f$ be a nonnegative function defined on $\mathbb{R}^{n},$ Riemann-integrable
on every rectangle. We say that $f$ is \textbf{Riemann-integrable\mindex{Riemann-integrable!nonnegative function}}
on $\mathbb{R}^{n}$ if the supremum of the integrals 
\[
\intop_{P}f\left(x\right)\textrm{d}x,
\]
 taken over all compact rectangles $P\subset\mathbb{R}^{n}$ is finite.

In this case, we define\\
\boxeq{
\[
\intop_{\mathbb{R}^{n}}f\left(x\right)\textrm{d}x=\underset{P\text{ compact rectangle}}{\sup}\intop_{P}f\left(x\right)\textrm{d}x.
\]
}

\end{definition}

\begin{proposition}{Sufficient and Necessary Condition of Riemann-Integrability}{}

Let $\left(A_{k}\right)$ be any nondecreasing sequence of Jordan-measurable
bounded sets such that
\[
\bigcup_{k}A_{k}=\mathbb{R}^{n}.
\]

For instance, one can take for $A_{k}$ the rectangle $\left[-k,k\right]\times\dots\times\left[-k,k\right],$
or the Euclidean ball of radius $k$ centered at the origin. 

Then, in order for $f$ to be Riemann-integrable on $\mathbb{R}^{n},$
it is necessary and sufficient\footnotemark that the sequence of
integrals 
\[
\intop_{A_{k}}f\left(x\right)\textrm{d}x
\]
 is bounded. In that case,\boxeq{
\[
\intop_{\mathbb{R}^{n}}f\left(x\right)\textrm{d}x=\lim_{k}\uparrow\intop_{A_{k}}f\left(x\right)\textrm{d}x.
\]
}

\end{proposition}

\footnotetext{The result is elementary to prove if we additionally
assume that $A_{k}\subset\mathring{A_{k+1}}$ for every $k=1,2,\dots.$}

\begin{definition}{Riemann-Integral of a Function of Arbitrary Sign}{}

Now, let $f$ be a function of arbitrary sign defined on $\mathbb{R}^{n},$
such that $\left|f\right|$ is Riemann-integrable. Then, for every
sequence $\left(A_{k}\right)$ of Jordan-measurable bounded sets with
\[
\bigcup_{k}A_{k}=\mathbb{R}^{n},
\]
the sequence of integrals 
\[
\left(\intop_{A_{k}}f\left(x\right)\textrm{d}x\right)
\]
is convergent\footnotemark. Moreover, its limit does not depend on
the choice of the sequence $\left(A_{k}\right).$

We then define\\
\boxeq{
\[
\intop_{\mathbb{R}^{n}}f\left(x\right)\textrm{d}x=\lim_{k\to+\infty}\intop_{A_{k}}f\left(x\right)\textrm{d}x.
\]
}

\end{definition}

\footnotetext{Here again, the proof is elementary only if we moreover
suppose that each of the $A_{k}$ is contained in the interior of
$A_{k+1}.$}

\subsubsection{Integral over a Non-Bounded Jordan-measurable Part of $\mathbb{R}^{n}$}

Finally, we can define the generalised integral of a function $f$
over a subset $A$ of $\mathbb{R}^{n},$ non-necessarily bounded,
provided that the intersection of $A$ with every rectangle is Jordan-measurable.
In this case, we simply say that $A$ is \textbf{\index{Jordan-measurable}Jordan-measurable}. 

\begin{definition}{Integral over a Non-Bounded Jordan-measurable Part of $\mathbb{R}^n$}{}

If $A$ is a non-bounded Jordan-measurable subset of $\mathbb{R}^{n},$
we say that $f$ is \textbf{Riemann-integrable on $A$\mindex{Riemann-integrable ! non-bounded subset}}
if the function $\overline{f}$ defined by extending $f$ by 0 outside
of $A,$ is Riemann-integrable on $\mathbb{R}^{n}.$ In that case,
we define

\boxeq{
\[
\intop_{A}f\left(x\right)\textrm{d}x=\intop_{\mathbb{R}^{n}}\overline{f}\left(x\right)\textrm{d}x.
\]
}

\end{definition}

\begin{remark}{}{}

For a function $f$ whose absolute value $\left|f\right|$ is \textbf{not}
Riemann-integrable, the so-called ``half-convergent'' or ``conditionally
convergent'' or ``improper'' integral is not defined when the dimension
$d>1.$ 

In such cases, the behaviour of the integrals
\[
\intop_{A_{k}}f\left(x\right)\textrm{d}x
\]
heavily depends on the choice of the sets $A_{k},$ whose shapes can
be arbitrarily chosen---and possibly in very weird and pathological
ways.

In particular, we can always choose the sequence $\left(A_{k}\right)$
in such a way that the limit of the integral becomes $+\infty,$ $-\infty$
or even any prescribed real number. 

In some cases, the limit may differ depending on wether we take $A_{k}$
to be the rectangle of side $2k$ centered at the origin, or the Euclidean
ball of radius $k$ centered at the origin.

\end{remark}

Tr.N. A classical example of such behavior where the behaviour of
the integrals depends of the sets $A_{k}$ is given in the following
example.

\begin{example}{}{}

Let $f:\mathbb{R}^{2}\backslash\left\{ \left(0,0\right)\right\} \to\mathbb{R}$
be the function defined for every $\left(x,y\right)\in\mathbb{R}^{2}$
such that $\left(x,y\right)\in\mathbb{R}^{2}\backslash\left\{ \left(0,0\right)\right\} $
by
\[
f\left(x,y\right)=\dfrac{xy}{x^{2}+y^{2}}.
\]

We consider the improper integral of $f$ over $\mathbb{R}^{2}\backslash\left\{ \left(0,0\right)\right\} $
and investigate the behaviour of
\[
\int_{A_{k}}f\left(x,y\right)\text{d}x\text{d}y
\]
where $A_{k}$ is a nondecreasing sequence of Jordan-measurable domains
such that
\[
\bigcup_{k\in\mathbb{N}}A_{k}=\mathbb{R}^{2}\backslash\left\{ \left(0,0\right)\right\} 
\]

Then:

1. Integrating over symmetric annuli centered at the origin gives
0.

2. Restricting to domain that are asymmetric regions, such a quadrant
makes the integral becomes $+\infty.$

\end{example}

\begin{solutionexample}{}{}

1. We consider the polar coordinates variable change: $x=r\cos\theta,$
$y=r\sin\theta,$ of Jacobian $r,$ we get
\[
f\left(x,y\right)=\dfrac{\cos\theta\sin\theta}{r^{2}}.
\]

And so, considering for $A_{k}=B\left(0,R_{k}\right)\backslash B\left(0,\epsilon\right)$
where $R_{0}<R_{1}<...<R_{k},$
\[
\intop_{B\left(0,R_{k}\right)\backslash B\left(0,\epsilon\right)}f\left(x,y\right)\text{d}x\text{d}y=\left(\intop^{2\pi}_{0}\cos\theta\sin\theta\text{d}\theta\right)\left(\intop^{R_{k}}_{\epsilon}\dfrac{1}{r}\text{d}r\right)=0\cdot\ln\left(\dfrac{R_{k}}{\epsilon}\right)=0.
\]
Hence, integrating over symmetric annuli centered at the origin gives
0.

2. Considering the following non-symmetric domain, such as
\[
D_{k}=\left\{ \left(x,y\right)\in\mathbb{R}^{2}\backslash\left\{ \left(0,0\right)\right\} :\,0<x<k,\,0<y<k\right\} ,
\]
then $f\left(x,y\right)$ becomes positive over the entire domain,
and the integral diverges to $+\infty:$ 
\[
\intop_{D_{k}}f\left(x,y\right)\text{d}x\text{d}y\underset{k\to+\infty}{\rightarrow}+\infty.
\]

\end{solutionexample}

\subsection{The Fubini Theorem}

The Fubini theorem allows us to compare integrals defined on a product
space with corresponding iterated integrals. For simplicity, we consider
only the case $\mathbb{R}^{2}=\mathbb{R}\times\mathbb{R}.$ The statements
for $\mathbb{R}^{n}=\mathbb{R}^{n_{1}}\times\mathbb{R}^{n_{2}},$
with $n=n_{1}+n_{2}$ are analogous.

\subsubsection{Case of a Continuous Function on a Closed Rectangle $\left[a,b\right]\times\left[c,d\right]$}

The most useful version of the \textbf{Fubini theorem\index{Fubini theorem}}
is the following.

\begin{theorem}{Fubini Theorem}{}

Let $f$ be a continuous function defined on the rectangle $\left[a,b\right]\times\left[c,d\right]$
in $\mathbb{R}^{2}.$

Then,\\
\boxeq{
\begin{align*}
\intop_{\left[a,b\right]\times\left[c,d\right]}f\left(x,y\right)\textrm{d}x\text{d}y & =\intop^{b}_{a}\left(\intop^{d}_{c}f\left(x,y\right)\text{d}y\right)\textrm{d}x\\
 & =\intop^{d}_{c}\left(\intop^{b}_{a}f\left(x,y\right)\textrm{d}x\right)\text{d}y
\end{align*}
}

\end{theorem}

\begin{remark}{}{}

For a continuous function $f$ defined on the rectangle $\left[a,b\right]\times\left[c,d\right]$
in $\mathbb{R}^{2},$ it is well known that the functions
\[
\begin{array}{ccccc}
x\mapsto\intop^{d}_{c}f\left(x,y\right)\text{d}y & \,\,\,\, & \text{and} & \,\,\,\, & y\mapsto\intop^{b}_{a}f\left(x,y\right)\textrm{d}x\end{array}
\]
are continuous on $\left[a,b\right]$ and $\left[c,d\right],$ respectively.

\end{remark}

\subsubsection{General Case}

There is no simple statement in the case where $f$ is only assumed
to be Riemann integrable. Indeed, even if $f$ is integrable on $\left[a,b\right]\times\left[c,d\right],$
we cannot guarantee that the function $y\mapsto f\left(x,y\right)$
is integrable on $\left[c,d\right]$ for every $x\in\left[a,b\right].$
And even when this is true, it does not ensure that the function 
\[
x\mapsto\intop^{d}_{c}f\left(x,y\right)\text{d}y,
\]
which is then defined for every $x\in\left[a,b\right],$ is itself
Riemann integrable on $\left[a,b\right].$

Similar difficulties arise when considering an extended integral over
the entire plane---even if $f$ is assumed to be continuous.

The Fubini theorem for the Riemann integral is stated as follows.

\begin{theorem}{Fubini Theorem for the Riemann Integral}{}

Let $f$ be a function defined and Riemann-integrable on $\mathbb{R}^{2}.$
Suppose that for every $x\in\mathbb{R},$ the function $y\mapsto f\left(x,y\right)$
is Riemann-integrable on $\mathbb{R}$---which, according to our
earlier convention, means that $y\mapsto f\left(x,y\right)$ is Riemann-integrable
on every bounded and closed interval, and that the integral $\intop^{+\infty}_{-\infty}\left|f\left(x,y\right)\right|\text{d}y$
converges.

If, in addition, the function
\[
x\mapsto\intop^{+\infty}_{-\infty}\left|f\left(x,y\right)\right|\text{d}y
\]
 is Riemann-integrable on $\mathbb{R},$ then\\
\boxeq{
\begin{equation}
\intop_{\mathbb{R}^{2}}f\left(x,y\right)\text{d}x\text{d}y=\intop^{+\infty}_{-\infty}\left(\intop^{+\infty}_{-\infty}f\left(x,y\right)\text{d}y\right)\text{d}x.\label{eq:Fubini_theorem_Riemann_integral}
\end{equation}
}

\end{theorem}

The corresponding statement for a rectangle follows immediately. Similar
results hold when the roles of $x$ and $y$ are interchanged, allowing---under
appropriate conditions---the order of integration to be reversed.

It is also worth noting that, if $f$ is continuous on the plane,
or even only piecewise continuous, it suffices to verify that $f$
is Riemann-integrable by checking that, for every nondecreasing sequence
of rectangles $\left(P_{k}\right)_{k\in\mathbb{N}}$ such that $\bigcup_{k\in\mathbb{N}}P_{k}=\mathbb{R}^{2},$
the integrals
\[
\intop_{P_{k}}\left|f\left(x,y\right)\right|\text{d}x\text{d}y
\]
form a bounded sequence.

\subsection{Change of Variable Formula}

Let $A$ and $B$ be two Jordan-measurable subsets of $\mathbb{R}^{2},$
and let $\phi:B\to A$ be a continuous function that defines a diffeomorphism
from the interior of $B$ onto the interior of $A.$

Then, for $f$ to be Riemann-integrable on $A,$ it is necessary and
sufficient that the function
\[
u\mapsto f\left(\phi\left(u\right)\right)\left|\text{det}\left(\phi^{\prime}\left(u\right)\right)\right|
\]
is Riemann-integrable on $B.$

In that case,\\
\boxeq{
\begin{equation}
\intop_{A}f\left(x\right)\text{d}x=\intop_{B}f\left(\phi\left(u\right)\right)\left|\text{det}\left(\phi^{\prime}\left(u\right)\right)\right|\text{d}u.\label{eq:formula_change}
\end{equation}
}

We say that the second part of the formula $\refpar{eq:formula_change}$
is obtained from the first by performing the change of variables $x=\phi\left(u\right)$---or,
in terms of coordinates
\[
\left\{ \begin{array}{c}
x_{1}=\phi_{1}\left(u_{1},u_{2}\right)\\
x_{2}=\phi_{2}\left(u_{1},u_{2}\right).
\end{array}\right.
\]

In the formula $\refpar{eq:formula_change},$ the notation $\phi^{\prime}\left(u\right)$
refers to the Jacobian matrix of $\phi$ at the point $u=\left(u_{1},u_{2}\right),$
\[
\phi^{\prime}\left(u\right)=\left(\begin{array}{cc}
\dfrac{\partial\phi_{1}}{\partial u_{1}}\left(u_{1},u_{2}\right) & \dfrac{\partial\phi_{1}}{\partial u_{2}}\left(u_{1},u_{2}\right)\\
\dfrac{\partial\phi_{2}}{\partial u_{1}}\left(u_{1},u_{2}\right) & \dfrac{\partial\phi_{2}}{\partial u_{2}}\left(u_{1},u_{2}\right)
\end{array}\right).
\]

This is often written more simply as
\[
\left(\begin{array}{cc}
\dfrac{\partial x_{1}}{\partial u_{1}} & \dfrac{\partial x_{1}}{\partial u_{2}}\\
\dfrac{\partial x_{2}}{\partial u_{1}} & \dfrac{\partial x_{2}}{\partial u_{2}}
\end{array}\right).
\]

Its determinant\\
\boxeq{
\[
\text{det}\left(\phi^{\prime}\left(u\right)\right)=\left|\begin{array}{cc}
\dfrac{\partial x_{1}}{\partial u_{1}} & \dfrac{\partial x_{1}}{\partial u_{2}}\\
\dfrac{\partial x_{2}}{\partial u_{1}} & \dfrac{\partial x_{2}}{\partial u_{2}}
\end{array}\right|
\]
}is called the \textbf{Jacobian\index{Jacobian}} of $\phi$ at the
point $u.$ In the formula $\refpar{eq:formula_change},$ it is the
absolute value of the Jacobian that appears.

\begin{example}{}{int_exp_min_xsqu}

Compute the integral 
\[
I=\intop^{+\infty}_{-\infty}\text{e}^{-\frac{x^{2}}{2}}\text{d}x.
\]

\end{example}

\begin{solutionexample}{}{}

One method is to compute in two different ways the double integral
\begin{equation}
\iintop_{\mathbb{R}^{2}}\text{e}^{-\frac{1}{2}\left(x^{2}+y^{2}\right)}\text{d}x\,\text{d}y.\label{eq:double_int_exp}
\end{equation}

\textbf{Using the Fubini theorem}

By applying the Fubini theorem previously stated, and by setting 
\[
f\left(x,y\right)=\text{e}^{-\frac{1}{2}\left(x^{2}+y^{2}\right)}=\text{e}^{-\frac{1}{2}x^{2}}\text{e}^{-\frac{1}{2}y^{2}}.
\]

Since, for instance,
\[
\text{e}^{-\frac{1}{2}y^{2}}\leqslant\dfrac{1}{1+\dfrac{1}{2}y^{2}},
\]
the function $y\mapsto f\left(x,y\right)$ is integrable over $\mathbb{\mathbb{R}},$
and
\[
\intop_{\mathbb{R}}f\left(x,y\right)\text{d}y=\text{e}^{-\frac{1}{2}x^{2}}\intop_{\mathbb{R}}\text{e}^{-\frac{1}{2}y^{2}}\text{d}y.
\]

The function $x\mapsto\intop_{\mathbb{R}}f\left(x,y\right)\text{d}y$
is integrable over $\mathbb{R}$ for the same reason. Noting that
$f$ is nonnegative, we may apply the Fubini theorem to obtain
\begin{align}
\iintop_{\mathbb{R}^{2}}\text{e}^{-\frac{1}{2}\left(x^{2}+y^{2}\right)}\text{d}x\,\text{d}y & =\intop_{\mathbb{R}}\text{e}^{-\frac{1}{2}x^{2}}\left[\intop_{\mathbb{R}}\text{e}^{-\frac{1}{2}y^{2}}\text{d}y\right]\text{d}x\nonumber \\
 & =\left(\intop_{\mathbb{R}}\text{e}^{-\frac{1}{2}x^{2}}\text{d}x\right)\left(\intop_{\mathbb{R}}\text{e}^{-\frac{1}{2}y^{2}}\text{d}y\right)\nonumber \\
 & =I^{2}.\label{eq:Fubini}
\end{align}

\textbf{Polar coordinates computation}

Let us now compute the double integral $\refpar{eq:double_int_exp}$
by switching to polar coordinates. More precisely, we perform in $\refpar{eq:double_int_exp}$
the change of variables for $r\in\left[0,+\infty\right[$ and $\theta\in\left[0,2\pi\right[$
\begin{equation}
\begin{array}{c}
\left\{ \begin{array}{c}
x=r\cos\theta\\
y=r\sin\theta
\end{array}\right.\end{array}\label{eq:polar_coord}
\end{equation}
 The function $\phi:\left[0,+\infty\right[\times\left[0,2\pi\right[\to\mathbb{R}^{2}$
which defines this change of variables, is a diffeomorphism from $\left[0,+\infty\right[\times\left[0,2\pi\right[$
onto the plane minus the nonnegative $x$-axis, that is $\left[0,+\infty\right[\times\left\{ 0\right\} .$ 

The Jacobian of $\phi$ can be computed by differentiating the system
$\refpar{eq:polar_coord}$
\[
\left\{ \begin{array}{c}
\text{d}x=\cos\theta\,\text{d}r-r\sin\theta\,\text{d\ensuremath{\theta}}\\
\text{d}y=\sin\theta\,\text{d}r+r\cos\theta\text{\,d}\theta.
\end{array}\right.
\]
Thus, the Jacobian is
\begin{align*}
J\left(\phi\right) & =\left|\begin{array}{cc}
\cos\theta & -r\sin\theta\\
\sin\theta & r\cos\theta
\end{array}\right|=r.
\end{align*}

Applying formula $\refpar{eq:formula_change},$ we obtain
\[
\iintop_{\mathbb{R}^{2}}\exp\left(-\dfrac{1}{2}\left(x^{2}+y^{2}\right)\right)\text{d}x\,\text{d}y=\iintop_{\left[0,+\infty\right[\times\left[0,2\pi\right[}\exp\left(-\dfrac{1}{2}r^{2}\right)r\text{d}r\,\text{d}\theta.
\]

Now applying the Fubini theorem to this expression, we get
\begin{align}
\iintop_{\mathbb{R}^{2}}\exp\left(-\dfrac{1}{2}\left(x^{2}+y^{2}\right)\right)\text{d}x\,\text{d}y & =\left(\intop_{\left[0,+\infty\right[}\exp\left(-\dfrac{1}{2}r^{2}\right)r\text{d}r\right)\left(\intop_{\left[0,2\pi\right[}\text{d}\theta\right)\nonumber \\
 & =2\pi\left[\exp\left(-\dfrac{1}{2}r^{2}\right)\right]^{+\infty}_{0}\nonumber \\
 & =2\pi.\label{eq:I_by_polar}
\end{align}

Hence, gathering $\refpar{eq:Fubini}$ and $\refpar{eq:I_by_polar}$,
we obtain $I^{2}=2\pi,$ and therefore
\begin{equation}
\intop^{+\infty}_{-\infty}\exp\left(-\dfrac{x^{2}}{2}\right)dx=\sqrt{2\pi}.\label{eq:int_normal}
\end{equation}

\end{solutionexample}

\begin{remark}{}{}

An alternate proof---often presented to avoid using the Fubini theorem
on an unbounded domain---consists in computing:
\begin{itemize}
\item On the one hand, the integral, denoted $I_{R},$ of the function $\left(x,y\right)\mapsto\exp\left(-\dfrac{1}{2}\left(x^{2}+y^{2}\right)\right)$
over the square $\left[-R,R\right]^{2},$ 
\item And on the other hand, by switching to polar coordinates to compute
the integral, denoted $J_{R}$ of the same function on the disk of
center $0$ with radius $R,$
\end{itemize}
and then letting $R$ tends to the infinite after noting the inequality
\[
J_{R}\leqslant I_{R}\leqslant I_{R\sqrt{2}}.
\]

\end{remark}

\section*{Exercises}

\addcontentsline{toc}{section}{Exercises}

\begin{exercise}{Triangular Law and Independence}{exercise6.1}

Let $a>0$ and $\alpha>0$ be two real numbers. Define the function
$f$ on $\mathbb{R}^{+}$ by
\[
\forall x\in\mathbb{R}^{+},\,\,\,\,f\left(x\right)=\alpha\left(x\boldsymbol{1}_{\left[0,\frac{a}{2}\right[}\left(x\right)+\left(a-x\right)\boldsymbol{1}_{\left[\frac{a}{2},a\right[}\left(x\right)\right).
\]

1. Compute the constant $\alpha_{0}$ such that $f$ is a probability
density function.

We fix $\alpha=\alpha_{0}$ for the remainder of this exercise.

2. Let $X$ be a random variable with density $f,$ and let $b\in\left]0,\dfrac{a}{2}\right[.$
Compute the probabilities 
\[
P\left(X>\dfrac{a}{2}\right)\,\,\,\,\text{and}\,\,\,\,\,P\left(\dfrac{a}{2}-b<X\leqslant\dfrac{a}{2}+b\right).
\]

3. Prove that for every $b\in\left]0,\dfrac{a}{2}\right[,$ the events
\[
A=\left(X>\dfrac{a}{2}\right)\,\,\,\,\text{and}\,\,\,\,B=\left(\dfrac{a}{2}-b<X\leqslant\dfrac{a}{2}+b\right)
\]
are independent.

\end{exercise}

\begin{exercise}{Exponential Law and Arrival Time}{exercise6.2}

The arrival time of the first customer after a shop opens is modelled
by a random variable $T$ defined on a probabilized space $\left(\Omega,\mathcal{A},P\right),$
and follows the exponential law with parameter $p>0.$

1. Compute the probability $P\left(T>\dfrac{1}{p}\right).$

2. Fix $\epsilon>0.$ 

For each $k\in\mathbb{N},$ consider the time interval $\left[k\epsilon,\left(k+1\right)\epsilon\right[$
and define the event 
\[
A_{k}=\left\{ T\in\left[k\epsilon,\left(k+1\right)\epsilon\right[\right\} ,
\]
corresponding to the statement ``the client arrives in the time interval
$\left[k\epsilon,\left(k+1\right)\epsilon\right[$''.

Compute the probability 
\[
P\left(T\in\left[k\epsilon,\left(k+1\right)\epsilon\right[\right).
\]

3. Define a random variable $X$ with values in $\mathbb{N},$ by
setting for every $\omega\in\Omega,$
\[
X\left(\omega\right)=k\Leftrightarrow\omega\in A_{k}.
\]

What is the law of the random variable $X?$

4. For every $t>0$ and for every $h>0,$ compute the probability
\[
P\left(t<T\right)
\]
 and the conditional probability 
\[
P\left(T>t+h\left|T>t\right.\right).
\]

\end{exercise}

\begin{exercise}{Uniform and Decimal Laws}{exercise6.3}

Let $X$ be a random variable defined on a probabilized space $\left(\Omega,\mathcal{A},P\right),$
following the uniform law on an interval $\left[0,1\right].$

We define two random variables $D_{1}$ and $D_{2}$---representing
the first and second decimal digits of $X$---by
\[
\begin{array}{ccc}
D_{1}=\left\lfloor 10X\right\rfloor  & \,\,\,\, & D_{2}=\left\lfloor 100X-10D_{1}\right\rfloor .\end{array}
\]

We recall that $\left\lfloor x\right\rfloor $ denotes the integer
part of the real number $x.$

1. What are the laws of the random variables $D_{1}$ and $D_{2}?$

2. Prove that the the random variables $D_{1}$ and $D_{2}$ are independent.

\begin{remark}{}{}This exercise can be generalized: one can show
that if $D_{n}$ denotes the $n-$th decimal digit of $X,$ then the
sequence of random variables $\left(D_{n}\right)_{n\in\mathbb{N}^{\ast}}$
forms a sequence of independent random variables, each uniformly distributed
over the set $\left\llbracket 0,9\right\rrbracket .$

\end{remark}

\end{exercise}

\begin{exercise}{Uniform and Triangular Law. Convolution.}{exercise6.4}

Let $X$ and $Y$ be two independent random variables defined on a
probabilized space $\left(\Omega,\mathcal{A},P\right),$ each following
the uniform law on the interval $\left[0,1\right].$

What is the law of the random variable $Z=X+Y?$

\end{exercise}

\begin{exercise}{Law of a Function of a Random Variable and Convolution}{exercise6.5}

Let $X$ be a random variable defined on a probabilized space $\left(\Omega,\mathcal{A},P\right),$
with a piecewise continuous density $f,$ and cumulative distribution
function $F.$

Define the random variable $Y=X^{2}.$

1. Prove that $Y$ admits a density $f_{Y},$ and express it in terms
of $f$ and $F.$ Give the explicit result in the particular case
where $X$ follows the Gauss law $\mathcal{N}\left(0,1\right).$

2. Let $X_{1}$ and $X_{2}$ be two independent real-valued random
variables defined on a probabilized space $\left(\Omega,\mathcal{A},P\right),$
both following the Gauss law $\mathcal{N}\left(0,1\right).$ Justify
the existence of a density for the random variable $Z=X^{2}_{1}+X^{2}_{2}$
and compute it.

3. Answer the same question for the random variable $T=\sqrt{Z}.$

\end{exercise}

\begin{exercise}{Expectation Computation}{exercise6.6}

Let $X$ be a real-valued random variable, defined on a probabilized
space $\left(\Omega,\mathcal{A},P\right)$ and following the Gauss
law $\mathcal{N}\left(0,1\right).$

Check that, for every $t\in\mathbb{R},$ the random variable $\exp\left(tX\right)$
admits an expectation, and compute it.

\end{exercise}

\begin{exercise}{Moments of Classical Continuous Laws}{exercise6.7}

Let $X$ be a real-valued random variable.

1. Compute the expectation and variance of $X$ in each of the following
cases, assuming that the law of $X$ is:

a. The uniform law on $\left[a,b\right];$

b. The exponential law with parameter $p>0;$

c. The chi-squared law with $n$ degrees of liberty.

2. Show that if $X$ follows the Cauchy law, then it does not admit
an expectation---and therefore has no variance either.

\end{exercise}

\begin{exercise}{Marginal Law and Independence}{exercise6.8}

Let $X=\left(X_{1},X_{2}\right)$ be a random variable defined on
a probabilized space $\left(\Omega,\mathcal{A},P\right),$ following
the uniform law on the square $\left[0,1\right]\times\left[0,1\right].$

Show that the random variables $X_{1}$ and $X_{2}$ are independent
and both follow the uniform law on $\left[0,1\right].$

\end{exercise}

\begin{exercise}{A Zero Covariance Does Not Imply Independence}{exercise6.9}

Let $X=\left(X_{1},X_{2}\right)$ be a random variable defined on
a probabilized space $\left(\Omega,\mathcal{A},P\right),$ following
the uniform law on the disk $\mathcal{D}\left(0,1\right).$

Show that the random variables $X_{1}$ and $X_{2}$ are not independent.

Compute the covariance of $X_{1}$ and $X_{2}.$

\end{exercise}

\section*{Solutions of Exercises}

\addcontentsline{toc}{section}{Solutions of Exercises}

\begin{solution}{}{solexercise6.1}

\textbf{1. Computation of $\alpha_{0}$ for $f$ to be a probability
density}

We compute
\begin{align*}
\intop^{+\infty}_{-\infty}f\left(x\right)\text{d}x & =\alpha\left[\intop^{\frac{a}{2}}_{0}x\text{d}x+\intop^{a}_{\frac{a}{2}}\left(a-x\right)\text{d}x\right]=2\alpha\left[\dfrac{x^{2}}{2}\right]^{\frac{a}{2}}_{0}=\alpha\dfrac{a^{2}}{4}.
\end{align*}

Hence $f$ is a probability density if and only if $\intop^{+\infty}_{-\infty}f\left(x\right)\text{d}x=1$
which occurs only for\\
\boxeq{
\[
\alpha_{0}=\dfrac{4}{a^{2}}.
\]
}

\textbf{2. Computation of the probabilities $P\left(X>\dfrac{a}{2}\right)$
and $P\left(\dfrac{a}{2}-b<X\leqslant\dfrac{a}{2}+b\right)$}

We have
\[
P\left(X>\dfrac{a}{2}\right)=\alpha_{0}\intop^{a}_{\frac{a}{2}}\left(a-x\right)\text{d}x=\dfrac{1}{2}.
\]

Similarly,
\[
P\left(\dfrac{a}{2}-b<X\leqslant\dfrac{a}{2}+b\right)=\alpha_{0}\left(\intop^{\frac{a}{2}}_{\frac{a}{2}-b}x\text{d}x+\intop^{\frac{a}{2}+b}_{\frac{a}{2}}\left(a-x\right)\text{d}x\right).
\]

Using the change of variable $y=a-x$ in the second integral, we obtain
\begin{align*}
P\left(\dfrac{a}{2}-b<X\leqslant\dfrac{a}{2}+b\right) & =2\alpha_{0}\intop^{\frac{a}{2}}_{\frac{a}{2}-b}x\text{d}x=2\alpha\left[\dfrac{x^{2}}{2}\right]^{\frac{a}{2}}_{\frac{a}{2}-b}=\alpha_{0}b\left(a-b\right).
\end{align*}

Thus,\\
\boxeq{
\[
P\left(\dfrac{a}{2}-b<X\leqslant\dfrac{a}{2}+b\right)=\dfrac{4b\left(a-b\right)}{a^{2}}.
\]
}

\textbf{3. Independence of $A=\left(X>\dfrac{a}{2}\right)$ and $B=\left(\dfrac{a}{2}-b<X\leqslant\dfrac{a}{2}+b\right)$}

We have
\[
P\left(A\cap B\right)=P\left(\dfrac{a}{2}<X\leqslant\dfrac{a}{2}+b\right)=\alpha_{0}\intop^{\frac{a}{2}+b}_{\frac{a}{2}}\left(a-x\right)\text{\,d}x.
\]

Making the change of variable $y=a-x,$ we obtain\\
\boxeq{
\[
P\left(A\cap B\right)=\alpha_{0}\intop^{\frac{a}{2}}_{\frac{a}{2}-b}y\text{\,d}y=\dfrac{1}{2}\alpha_{0}b\left(a-b\right)=P\left(A\right)P\left(B\right).
\]
}

This shows that the events $A$ and $B$ are independent.

\end{solution}

\begin{solution}{}{solexercise6.2}

\textbf{1. Computation of $P\left(T>\dfrac{1}{p}\right)$}

Since $T$ follows the exponential law with parameter $p,$ we have\\
\boxeq{
\[
P\left(T>\dfrac{1}{p}\right)=\intop^{+\infty}_{\frac{1}{p}}p\exp\left(-px\right)\text{d}x=\left[-\exp\left(-px\right)\right]^{+\infty}_{\frac{1}{p}}=\dfrac{1}{\text{e}}\approx0.367
\]
}

\textbf{2. Computation of $P\left(T\in\left[k\epsilon,\left(k+1\right)\epsilon\right[\right)$}

We compute
\begin{align*}
P\left(T\in\left[k\epsilon,\left(k+1\right)\epsilon\right[\right) & =\intop^{\left(k+1\right)\epsilon}_{k\epsilon}p\exp\left(-px\right)\text{\,d}x\\
 & =\left[-\exp\left(-px\right)\right]^{k\left(1+\epsilon\right)}_{k\epsilon}\\
 & =\text{e}^{-k\epsilon}-\text{e}^{-\left(k+1\right)\epsilon}
\end{align*}
Hence,\\
\boxeq{
\[
P\left(T\in\left[k\epsilon,\left(k+1\right)\epsilon\right[\right)=\text{e}^{-k\epsilon}\left(1-e^{-\epsilon}\right).
\]
}

\textbf{3. Law followed by the random variable $X$}

Since, for every $k\in\mathbb{N},$ we have
\[
\left(X=k\right)=\left(T\in\left[k\epsilon,\left(k+1\right)\epsilon\right[\right),
\]
it follows from the previous question that the law of $X$ is the
geometric law $\mathcal{G}_{\mathbb{N}}\left(1-\text{e}^{-\epsilon}\right).$
This illustrates the principle that a continuous-time phenomenon modelled
by an exponential law becomes a geometric law when time is discretized.

\textbf{4. Computation of $P\left(t<T\right)$ and $P\left(T>t+h\left|T>t\right.\right)$}

We compute
\[
P\left(t<T\right)=\intop^{+\infty}_{t}p\exp\left(-px\right)\text{\,d}x=\left[-\exp\left(-px\right)\right]^{+\infty}_{t}=\exp\left(-pt\right).
\]

Since for every $t>0$ and for every $h>0,$
\[
\left(T>t+h\right)\cap\left(T>t\right)=\left(T>t+h\right),
\]
we have
\[
P\left(T>t+h\right)=\exp\left(-p\left(t+h\right)\right).
\]

By definition of the conditional probability,\\
\boxeq{
\[
P\left(T>t+h\left|T>t\right.\right)=\dfrac{P\left(T>t+h,T>t\right)}{P\left(T>t\right)}=\exp\left(-ph\right)\equiv P\left(T>h\right).
\]
}

\begin{remark}{}{}

We say that the \textbf{exponential law\index{exponential law}\mindex{law ! exponential}}
is a \index{memoryless law}\textbf{\mindex{law ! memoryless}memoryless
law}. The probability of waiting an additional time $h,$ given that
one has already waited time $t,$ is independent of $t.$ That is,
this probability neither increases nor decrease with $t!$ 

Other physical phenomena are also modelled by an exponential law---for
instance the lifetime of a radioactive nucleus. Whether a nucleus
formed one minute ago or five billion years ago, its probability of
decaying in the next second remains the same.

It can be shown that the memoryless property characterizes the random
variables with an exponential law density---this is an interesting
exercise in itself. Similarly, the condition\\
\boxeq{
\[
P\left(T>n+h\left|T\geqslant n\right.\right)=P\left(T\geqslant h\right)
\]
}characterizes the random variables of geometric law among the integer-valued
random variables.

\end{remark}

\end{solution}

\begin{solution}{}{solexercise6.3}

\textbf{1. Laws of the random variables $D_{1}$ and $D_{2}$}

For every integer $k_{1}\in\left\llbracket 0,9\right\rrbracket ,$
\[
\left(D_{1}=k_{1}\right)=\left(k_{1}\leqslant10X<k_{1}+1\right).
\]

Since the random variable $X$ follows the uniform law on the interval
$\left[0,1\right],$
\[
P\left(D_{1}=k_{1}\right)=\dfrac{1}{10}.
\]

Now, for every integer $k_{2}\in\left\llbracket 0,9\right\rrbracket ,$
\begin{align*}
\left(D_{2}=k_{2}\right) & =\left(k_{2}\leqslant10^{2}X-10D_{1}<k_{2}+1\right)\\
 & =\biguplus^{9}_{k_{1}=0}\left[\left(D_{1}=k_{1}\right)\cap\left(k_{2}\leqslant10^{2}X-10D_{1}<k_{2}+1\right)\right]\\
 & =\biguplus^{9}_{k_{1}=0}\left[\left(D_{1}=k_{1}\right)\cap\left(\dfrac{k_{1}}{10}+\dfrac{k_{2}}{10^{2}}\leqslant X<\dfrac{k_{1}}{10}+\dfrac{k_{2}}{10^{2}}+\dfrac{1}{10^{2}}\right)\right]\\
 & =\biguplus^{9}_{k_{1}=0}\left(\dfrac{k_{1}}{10}+\dfrac{k_{2}}{10^{2}}\leqslant X<\dfrac{k_{1}}{10}+\dfrac{k_{2}}{10^{2}}+\dfrac{1}{10^{2}}\right).
\end{align*}

Each of this event is disjoint and of same probability $10^{-2},$
it follows that\\
\boxeq{
\[
P\left(D_{2}=k_{2}\right)=\dfrac{1}{10}.
\]
}

\textbf{2. Independence of $D_{1}$ and $D_{2}$} 

For every integers $k_{1},k_{2}\in\left\llbracket 0,9\right\rrbracket ,$
\begin{multline*}
\left(D_{1}=k_{1}\right)\cap\left(D_{2}=k_{2}\right)\\
=\left(\dfrac{k_{1}}{10}\leqslant X<\dfrac{k_{1}}{10}+\dfrac{1}{10}\right)\cap\left(\dfrac{k_{1}}{10}+\dfrac{k_{2}}{10^{2}}\leqslant X<\dfrac{k_{1}}{10}+\dfrac{k_{2}}{10^{2}}+\dfrac{1}{10^{2}}\right),
\end{multline*}
thus
\[
P\left(\left(D_{1}=k_{1}\right)\cap\left(D_{2}=k_{2}\right)\right)=\dfrac{1}{100}.
\]

We have shown that for $k_{1},k_{2}\in\left\llbracket 0,9\right\rrbracket ,$\\
\boxeq{
\[
P\left(\left(D_{1}=k_{1}\right)\cap\left(D_{2}=k_{2}\right)\right)=P\left(D_{1}=k_{1}\right)P\left(D_{2}=k_{2}\right).
\]
}

Hence, the random variables $D_{1}$ and $D_{2}$ are independent,
and each follows the uniform law on $\left\llbracket 0,9\right\rrbracket .$

\end{solution}

\begin{solution}{}{solexercise6.4}

\textbf{Law of the random variable $Z=X+Y$}

The random variable $Z$ admits a density $f_{Z},$ which is the convolution
of the densities of $X$ and $Y.$ That is, for every $z\in\mathbb{R},$
\[
f_{Z}\left(z\right)=\intop^{+\infty}_{-\infty}\boldsymbol{1}_{\left[0,1\right]}\left(x\right)\boldsymbol{1}_{\left[0,1\right]}\left(z-x\right)\text{\,d}x.
\]

Since
\[
\boldsymbol{1}_{\left[0,1\right]}\left(x\right)\boldsymbol{1}_{\left[0,1\right]}\left(z-x\right)=\boldsymbol{1}_{\left[0,1\right]}\left(z\right)\boldsymbol{1}_{\left[0,z\right]}\left(x\right)+\boldsymbol{1}_{\left[1,2\right]}\left(z\right)\boldsymbol{1}_{\left[z-1,1\right]}\left(x\right),
\]
we can write
\[
f_{Z}\left(z\right)=\boldsymbol{1}_{\left[0,1\right]}\left(z\right)\intop^{z}_{0}\text{d}x+\boldsymbol{1}_{\left[1,2\right]}\left(z\right)\intop^{1}_{z-1}\text{d}x.
\]

Thus the random variable $Z$ admits a density $f_{Z}$ given for
every $z\in\mathbb{R}$ by\\
\boxeq{
\[
f_{Z}\left(z\right)=z\boldsymbol{1}_{\left[0,1\right]}\left(z\right)+\left(2-z\right)\boldsymbol{1}_{\left[1,2\right]}\left(z\right).
\]
}

This means that $Z$ follows a \textbf{\index{triangular law}\mindex{law ! triangular}triangular
law}.

\end{solution}

\begin{solution}{}{solexercise6.5}

\textbf{1. Density existence of $Y$ and expression. Normal law particular
case}

We denote $F_{Y}$ the cumulative distribution function of the random
variable $Y.$ First, we remark that for every $y\leqslant0,$
\[
F_{Y}\left(y\right)\equiv P\left(Y\leqslant y\right)=0.
\]

Let $y$ be any positive real number. Then
\[
F_{Y}\left(y\right)=P\left(-\sqrt{y}\leqslant X\leqslant\sqrt{y}\right)=\intop^{\sqrt{y}}_{-\sqrt{y}}f\left(u\right)\text{d}u.
\]

Since the function $f$ is piecewise continuous, the cumulative distribution
function $F_{Y}$ is differentiable on every open interval where $f$
is continuous, and its derivative is the density $f_{Y}$ of $Y.$
It yields
\[
f_{Y}\left(y\right)=\dfrac{1}{2\sqrt{y}}\left(f\left(\sqrt{y}\right)+f\left(-\sqrt{y}\right)\right).
\]

In summary, $Y$ admits a piecewise continuous density given by
\[
f_{Y}\left(y\right)=\begin{cases}
\dfrac{1}{2\sqrt{y}}\left(f\left(\sqrt{y}\right)+f\left(-\sqrt{y}\right)\right), & \text{if\,}y>0,\\
0, & \text{otherwise.}
\end{cases}
\]

If $X$ follows the Gauss law $\mathcal{N}\left(0,1\right),$ the
function $f$ is even, and\\
\boxeq{
\[
f_{Y}\left(y\right)=\begin{cases}
\dfrac{1}{\sqrt{2\pi y}}\text{e}^{-\frac{y}{2}}, & \text{if\,}y>0,\\
0, & \text{otherwise.}
\end{cases}
\]
}

In other words, $Y$ follows the chi-squared law with one degree of
freedom.

\textbf{2. Density existence of $Z=X^{2}_{1}+X^{2}_{2}$ and computation}

The random variables $X^{2}_{1}$ and $X^{2}_{2}$ are independent,
and their sum $Z$ admits a density $f_{Z},$ which is the convolution
of the densities of $X^{2}_{1}$ and $X^{2}_{2}.$ It follows from
the previous question that, for every $z\leqslant0,$ we have $f_{Z}\left(z\right)=0,$
and for every $z>0,$
\begin{multline*}
f_{Z}\left(z\right)=\intop^{+\infty}_{-\infty}\boldsymbol{1}_{\mathbb{R}^{\ast+}}\left(y\right)\dfrac{1}{\sqrt{2\pi y}}\text{e}^{-\frac{y}{2}}\boldsymbol{1}_{\mathbb{R}^{\ast+}}\left(z-y\right)\dfrac{1}{\sqrt{2\pi\left(z-y\right)}}\text{e}^{-\frac{z-y}{2}}\text{d}y.
\end{multline*}

Thus, by simplifying and applying the change of variable $y=zu,$
it yields
\begin{align*}
f_{Z}\left(z\right) & =\dfrac{1}{2\pi}\text{e}^{-\frac{z}{2}}\intop^{z}_{0}\dfrac{1}{\sqrt{y\left(z-y\right)}}\text{d}y=\dfrac{1}{2\pi}\text{e}^{-\frac{z}{2}}\intop^{1}_{0}\dfrac{1}{\sqrt{u\left(1-u\right)}}\text{d}u.
\end{align*}

We recognize the analytical form of an exponential law density with
parameter $\dfrac{1}{2},$ which additionally yields, without computation
\[
\intop^{1}_{0}\dfrac{1}{\sqrt{u\left(1-u\right)}}\text{d}u=\pi.
\]

Therefore,\\
\boxeq{
\[
f_{Z}\left(z\right)=\begin{cases}
\dfrac{1}{2}\text{e}^{-\frac{z}{2}}, & \text{if\,}z>0,\\
0, & \text{otherwise.}
\end{cases}
\]
}

\textbf{3. Density existence of $T=\sqrt{Z}$ and computation}

Let $F_{T}$ be the cumulative distribution function of $T.$ For
every $t\leqslant0,$
\[
F_{T}\left(t\right)=0,
\]
and for every $t>0,$
\[
F_{T}\left(t\right)\equiv P\left(T\leqslant t\right)=P\left(Z\leqslant t^{2}\right)=\intop^{t^{2}}_{0}\dfrac{1}{2}\text{e}^{-\frac{z}{2}}\text{\,d}z.
\]

The function $F_{T}$ is differentiable everywhere except possibly
at 0. The random variable $T$ admits a density $f_{T},$ the derivative
of $F_{T},$ given by\\
\boxeq{
\[
f_{T}\left(t\right)=\begin{cases}
t\text{e}^{-\frac{t^{2}}{2}}, & \text{if\,}t>0,\\
0, & \text{otherwise.}
\end{cases}
\]
}

\end{solution}

\begin{solution}{}{solexercise6.6}

\textbf{Existence and Computation of the Expectation of $\exp\left(tX\right)$} 

Let $t$ be a real number. 

We denote $f_{X}$ the density of the Gaussian law $\mathscr{N}\left(0,1\right).$
For every $x\in\mathbb{R},$
\[
f_{X}\left(x\right)=\dfrac{1}{\sqrt{2\pi}}\text{e}^{-\frac{x^{2}}{2}}.
\]

For every $x\in\mathbb{R},$
\begin{align*}
\text{e}^{tx}f_{X}\left(x\right) & =\text{e}^{tx}\dfrac{1}{\sqrt{2\pi}}\text{e}^{-\frac{x^{2}}{2}}\\
 & =\dfrac{1}{\sqrt{2\pi}}\text{e}^{-\frac{\left(x-t\right)^{2}}{2}}\text{e}^{\frac{t^{2}}{2}},
\end{align*}
which proves that the function $x\mapsto\text{e}^{tx}\dfrac{1}{\sqrt{2\pi}}\text{e}^{-\frac{x^{2}}{2}}$
is Riemann-integrable. Therefore, the random variable $\exp\left(tX\right)$
admits an expectation given by
\begin{align*}
\mathbb{E}\left(\text{e}^{tX}\right) & =\intop^{+\infty}_{-\infty}\text{e}^{tx}f_{X}\left(x\right)\text{\,d}x\\
 & =\dfrac{1}{\sqrt{2\pi}}\text{e}^{\frac{t^{2}}{2}}\intop^{+\infty}_{-\infty}\text{e}^{-\frac{\left(x-t\right)^{2}}{2}}\text{\,d}x.
\end{align*}
Hence,\\
\boxeq{
\[
\mathbb{E}\left(\text{e}^{tX}\right)=\text{e}^{\frac{t^{2}}{2}}.
\]
}

\end{solution}

\begin{solution}{}{solexercise6.7}

\textbf{1. a. Computation of the expectation and ariance for the uniform
law on $\left[a,b\right]$}
\[
\mathbb{E}\left(X\right)=\intop^{+\infty}_{-\infty}\text{\ensuremath{\boldsymbol{1}_{\left[a,b\right]}\left(x\right)}}\dfrac{x}{b-a}\text{d}x=\dfrac{1}{b-a}\left[\dfrac{x^{2}}{2}\right]^{b}_{a}=\dfrac{a+b}{2}.
\]
\[
\sigma^{2}_{X}=\intop^{b}_{a}\dfrac{1}{b-a}\left[x-\dfrac{a+b}{2}\right]^{2}\text{d}x=\dfrac{1}{b-a}\intop^{\frac{b-a}{2}}_{\frac{a-b}{2}}y^{2}\text{dy}=\dfrac{1}{b-a}\left[\dfrac{y^{3}}{3}\right]^{\frac{b-a}{2}}_{\frac{a-b}{2}}.
\]
Hence,

\boxeq{
\[
\sigma^{2}_{X}=\dfrac{\left(b-a\right)^{2}}{12}.
\]
}

\textbf{b. Computation of the expectation and variance for the exponential
law with parameter $p>0$}
\begin{align*}
\mathbb{E}\left(X\right) & =\intop^{+\infty}_{-\infty}\boldsymbol{1}_{\mathbb{R}^{+}}\left(x\right)xp\text{e}^{-px}\textrm{d}x\\
 & =\intop^{+\infty}_{0}xp\text{e}^{-px}\textrm{d}x\\
 & =\left[-x\text{e}^{-px}\right]^{+\infty}_{0}+\intop^{+\infty}_{0}\text{e}^{-px}\textrm{d}x\\
 & =\dfrac{1}{p}.
\end{align*}

Using two integration by parts, in a similar way, we obtain
\begin{align*}
\sigma^{2}_{X} & =\intop^{+\infty}_{-\infty}\boldsymbol{1}_{\mathbb{R}^{+}}\left(x\right)\left(x-\dfrac{1}{p}\right)^{2}p\text{e}^{-px}\text{d}x=\cdots=\dfrac{1}{p^{2}}.
\end{align*}
Thus,\\
\boxeq{
\[
\sigma^{2}_{X}=\dfrac{1}{p^{2}}.
\]
}

\textbf{c. Computation of the expectation and variance for the Chi-Squared
law with $n$ degrees of freedom}

\begin{align*}
\mathbb{E}\left(X\right) & =\intop^{+\infty}_{-\infty}\boldsymbol{1}_{\mathbb{R}^{+}}\left(x\right)x\dfrac{1}{K_{n}}\text{e}^{-\frac{x}{2}}x^{\frac{n}{2}-1}\text{d}x\\
 & =\dfrac{1}{K_{n}}\intop^{+\infty}_{0}x^{\frac{n}{2}+1}\text{e}^{-\frac{x}{2}}\text{d}x\\
 & =...
\end{align*}

Hence,
\[
\mathbb{E}\left(X\right)=n.
\]

\begin{align*}
\sigma^{2}_{X} & =\intop^{+\infty}_{-\infty}\boldsymbol{1}_{\mathbb{R}^{+}}\left(x\right)\left(x-n\right)^{2}\dfrac{1}{K_{n}}\text{e}^{-\frac{x}{2}}x^{\frac{n}{2}-1}\text{d}x=\cdots
\end{align*}
Thus,\\
 \boxeq{
\[
\sigma^{2}_{X}=2n.
\]
}

\textbf{2. A random variable following the Cauchy law admits no expectation
and no variance}

The function 
\[
x\mapsto\dfrac{\left|x\right|}{\pi\left(1+x^{2}\right)}
\]
is not integrable on $\mathbb{R}.$ Therefore, the random variable
$X$ in this case does not have an expectation, and consequently no
variance either.

\end{solution}

\begin{solution}{}{solexercise6.8}

\textbf{Independence and law of $X_{1}$ and $X_{2}$}

The marginals $X_{1}$ and $X_{2}$ admit densities given by
\begin{itemize}
\item For $x_{1}\in\mathbb{R},$$\begin{array}{c}
f_{X_{1}}\left(x_{1}\right)=\intop^{+\infty}_{-\infty}\boldsymbol{1}_{\left[0,1\right]\times\left[0,1\right]}\left(x_{1},x_{2}\right)\text{\,d}x_{2}.\end{array}$
\item For $x_{2}\in\mathbb{R},$$\begin{array}{c}
f_{X_{2}}\left(x_{2}\right)=\intop^{+\infty}_{-\infty}\boldsymbol{1}_{\left[0,1\right]\times\left[0,1\right]}\left(x_{1},x_{2}\right)\text{\,d}x_{1}.\end{array}$
\end{itemize}
Since
\begin{equation}
\boldsymbol{1}_{\left[0,1\right]\times\left[0,1\right]}\left(x_{1},x_{2}\right)=\boldsymbol{1}_{\left[0,1\right]}\left(x_{1}\right)\boldsymbol{1}_{\left[0,1\right]}\left(x_{2}\right),\label{eq:caract_square}
\end{equation}
it follows that for every $x_{1}$ and $x_{2}$ in $\mathbb{R}$
\[
\begin{array}{ccccc}
f_{X_{1}}\left(x_{1}\right)=\boldsymbol{1}_{\left[0,1\right]}\left(x_{1}\right) & \,\,\,\, & \text{and} & \,\,\,\, & f_{X_{2}}\left(x_{2}\right)=\boldsymbol{1}_{\left[0,1\right]}\left(x_{2}\right).\end{array}
\]

The equality $\refpar{eq:caract_square}$ implies that the density
of $X$ is the direct product of the densities of the random variables
$X_{1}$ and $X_{2}$ and therefore, $X_{1}$ and $X_{2}$ are independent.

\end{solution}

\begin{solution}{}{solexercise6.9}

\textbf{Showing that the random variables $X_{1}$ and $X_{2}$ are
not independent}

The marginals $X_{1}$ and $X_{2}$ admit densities given by
\begin{itemize}
\item For $x_{1}\in\mathbb{R},$$\begin{array}{c}
\begin{array}{c}
f_{X_{1}}\left(x_{1}\right)=\intop^{+\infty}_{-\infty}\dfrac{1}{\pi}\boldsymbol{1}_{\mathcal{D}\left(0,1\right)}\left(x_{1},x_{2}\right)\text{\,d}x_{2}\end{array}.\end{array}$
\item For $x_{2}\in\mathbb{R},$$\begin{array}{c}
f_{X_{2}}\left(x_{2}\right)=\intop^{+\infty}_{-\infty}\dfrac{1}{\pi}\boldsymbol{1}_{\mathcal{D}\left(0,1\right)}\left(x_{1},x_{2}\right)\text{\,d}x_{1}.\end{array}$
\end{itemize}
Since $\left(x_{1},x_{2}\right)\in\mathcal{D}\left(0,1\right)$ is
equivalent to
\[
x^{2}_{1}+x^{2}_{2}\leqslant1
\]
 we have
\[
\begin{array}{c}
f_{X_{1}}\left(x_{1}\right)=\begin{cases}
\intop^{\sqrt{1-x^{2}_{1}}}_{-\sqrt{1-x^{2}_{1}}}\dfrac{1}{\pi}\text{\,d}x_{2}=\dfrac{2}{\pi}\sqrt{1-x^{2}_{1}}, & \text{if }-1<x_{1}<1,\\
0, & \text{otherwise.}
\end{cases}\end{array}
\]
\[
\begin{array}{c}
f_{X_{2}}\left(x_{2}\right)=\begin{cases}
\intop^{\sqrt{1-x^{2}_{2}}}_{-\sqrt{1-x^{2}_{2}}}\dfrac{1}{\pi}\text{\,d}x_{1}=\dfrac{2}{\pi}\sqrt{1-x^{2}_{2}}, & \text{if }-1<x_{2}<1,\\
0, & \text{otherwise}.
\end{cases}\end{array}
\]

The product $f_{X_{1}}\left(x_{1}\right)f_{X_{2}}\left(x_{2}\right)$
is positive on the open square $\left]-1,1\right[^{2}$ and equal
to zero outside of it. Therefore, it cannot equal the joint density
of $X,$ which is
\[
f_{X}\left(x_{1},x_{2}\right)=\dfrac{1}{\pi}\boldsymbol{1}_{D\left(0,1\right)}\left(x\right)
\]
which is in particular nonzero on $\mathcal{D}\left(0,1\right)\backslash\left]-1,1\right[^{2}.$ 

\boxeq{Thus, the random variables $X_{1}$ and $X_{2}$ are not independent.
}

\textbf{Computation of the covariance of $X_{1}$ and $X_{2}.$}

We have
\begin{align*}
\text{cov}\left(X_{1},X_{2}\right) & =\mathbb{E}\left(\left(X_{1}-\mathbb{E}\left(X_{1}\right)\right)\left(X_{2}-\mathbb{E}\left(X_{2}\right)\right)\right)\\
 & =\mathbb{E}\left(X_{1}X_{2}\right)-\mathbb{E}\left(X_{1}\right)\mathbb{E}\left(X_{2}\right).
\end{align*}

Since, the random variables $X_{1}$ and $X_{2}$ have even densities,
their expectations are zero.

By symmetry, we also have
\[
\mathbb{E}\left(X_{1}X_{2}\right)=\intop^{+\infty}_{-\infty}\intop^{+\infty}_{-\infty}\dfrac{1}{\pi}x_{1}x_{2}\boldsymbol{1}_{\mathcal{D}\left(0,1\right)}\left(x_{1},x_{2}\right)\text{\,d}x_{1}\text{d}x_{2}=0.
\]
This proves that\\
\boxeq{
\[
\text{cov}\left(X_{1},X_{2}\right)=0.
\]
}

\end{solution}

\chapter{Approximation of Laws. Weak Law of Large Numbers.}\label{chap:Approximation-of-laws.}

\begin{objective}{}{}

Chapter \ref{chap:Approximation-of-laws.} introduces the approximation
of probability laws and presents weak law of large numbers.
\begin{itemize}
\item Section \ref{sec:Approximation-of-laws} begins with several approximations
of classical laws: it includes the Poisson theorem, which approximates
a binomial law by a Poisson law; the approximation of a binomial law
by a Gauss law; and, the approximation of a hypergeometric law by
a binomial law. The central limit theorem is then presented, approximating
the sum of independent random variables of same law by the normal
law. A definition of the convergence in law is then given.
\item In Section \ref{sec:Weak-Law-of}, after having defined the convergence
in probability, the weak law of large numbers is stated. The chapter
concludes this part with the Bernoulli theorem, which shows that independent
events with the same probability allows to define a sequence of random
variables that converges in probability to that probability, thus
providing a compatible framework for the frequentist approach of probabilities.
\end{itemize}
\end{objective}

\section*{Introduction}

This part provides a first introduction to the concept of approximation
of laws, which will be depevelopped in greater depth in Part \ref{part:Deepening-Probability-Theory},
in terms of convergence of sequences of bounded measures. Different
notions of convergence for sequences of random variables will be studied
within the general framework. 

\section{Approximation of Laws}\label{sec:Approximation-of-laws}

The aim of this section is to introduce the problem of approximation
in laws\footnote{In the French edition, the Author mentions that it does not introduce
the concept of convergence in law on purpose, in order to comply to
the French curriculum for certain teaching qualification examinations.}.

The binomial law, the hypergeometric law and gambling game, when modelled
probabilistically, lead naturally to expressions involving factorials.
The following figures show how large these values can get, and illustrate
the potential benefit of replacing exact probability laws with approximations\footnote{This numerical complexity reflects the complexity of the formula involved.
Although modern software can compute with numbers of arbitrary size,
there are still other reasons to use approximations, as we will see
in this chapter, especially in the exercises.}:

\[
\begin{array}{ccc}
10!=3,628,800 & \,\,\,\, & 15!=1,307,674,368,000\end{array}
\]

\begin{multline*}
52!=80,658,175,170,943,878,571,660,636,856,403,766,975,289,505,\\
440,883,277,824,000,000,000,000.
\end{multline*}

\subsection{Poisson approximation}

The Poisson theorem gives an approximation of the binomial law $\mathcal{B}\left(n,p\right)$
when $n$ is ``large'' and $p$ is ``small''.

\begin{theorem}{Poisson Theorem}{Poisson-theorem}

Let $\left(p_{n}\right)_{n\in\mathbb{N}}$ be a sequence of real numbers
in the interval $\left]0;1\right[$ such that $\lim_{n\rightarrow+\infty}np_{n}=\lambda,$
with $\lambda>0.$

Let consider, for each integer $n,$ a random variable $S_{n}$ of
law $\mathcal{B}\left(n,p_{n}\right).$ 

Then
\[
P\left(S_{n}=k\right)=P_{n}\left(k\right)=\begin{cases}
\binom{n}{k}p^{k}_{n}\left(1-p_{n}\right)^{n-k}, & \text{if\,}0\leqslant k\leqslant n,\\
0, & \text{otherwise.}
\end{cases}
\]

Thus, for every integer $k\in\mathbb{N},$ the sequence with general
term $P\left(S_{n}=k\right)$ converges, and\\
\boxeq{
\[
\lim_{n\to+\infty}P\left(S_{n}=k\right)=\dfrac{\lambda^{k}}{k!}\text{e}^{-\lambda}.
\]
}

\end{theorem}

\begin{proof}{}{}

We have first to remark, as $n$ tends to infinity
\[
p_{n}=\dfrac{\lambda}{n}+o\left(\dfrac{1}{n}\right).
\]

Thus, for every fixed $k\in\mathbb{N},$ and for every $n\geqslant k,$
\[
P\left(S_{n}=k\right)=\dfrac{n\left(n-1\right)\dots\left(n-k+1\right)}{k!}\left(\dfrac{\lambda}{n}+o\left(\dfrac{1}{n}\right)\right)^{k}\left(1-\dfrac{\lambda}{n}+o\left(\dfrac{1}{n}\right)\right)^{n-k}.
\]

Now,
\begin{multline*}
n\left(n-1\right)\dots\left(n-k+1\right)\left(\dfrac{\lambda}{n}+o\left(\dfrac{1}{n}\right)\right)^{k}\\
=\dfrac{n\left(n-1\right)\dots\left(n-k+1\right)}{n^{k}}\left(\lambda+o\left(1\right)\right)^{k}\underset{n\to+\infty}{\longrightarrow}\lambda^{k}.
\end{multline*}

Moreover, we know that
\[
\left(1-\dfrac{\lambda}{n}+o\left(\dfrac{1}{n}\right)\right)^{n-k}\underset{n\to+\infty}{\longrightarrow}\text{e}^{-\lambda}.
\]

Hence, we can conclude that
\[
\lim_{n\to+\infty}P\left(S_{n}=k\right)=\text{e}^{-\lambda}\dfrac{\lambda^{k}}{k!}.
\]

\end{proof}

\begin{remark}{}{}

There exists a more precise result which gives a convergence rate,
uniform in $k.$ This result requires a technical and difficult proof.
The interested reader can refer for instance to the book \cite{shiryaev2016probability}\footnotemark.
We merely cite it.

If, for every $n\in\mathbb{N},$ we have $np_{n}=\lambda,$ where
$\lambda$ is a positive real number, then we have the estimate

\boxeq{
\[
\sum^{+\infty}_{k=0}\left|P_{n}\left(k\right)-\text{e}^{-\lambda}\dfrac{\lambda^{k}}{k!}\right|\leqslant\dfrac{2\lambda}{n}\min\left(2,\lambda\right).
\]
}

The Poisson theorem gives an approximation of the binomial law when
the parameter $p$ is ``small''. In practice, it is common to replace
the binomial law with the Poisson law when $n$ is sufficiently large---on
the order of 30---and $p$ is small enougth---on the order of $0.1.$
That is, if $X$ is a random variable following the binomial law with
first parameter $n$ and expectation $\lambda$ not too large, its
law is approximately the Poisson law $\mathcal{P}\left(\lambda\right).$

\end{remark}

\footnotetext{Tr. N: We cite here the most recent edition of this
book.}

\begin{example}{Probability of Birthday on January the First}{example7.1}

We want to determine the probability $P_{n}\left(k\right)$ that,
among $n$ people, exactly $k$ were born on the first of January.
We suppose that no one is born on February 29, and that all other
days of the year are equally probable. 

Compute the value of $P_{n}\left(k\right)$ for $n=500$ and $n=600,$
and in each case, for values of $k$ between 0 and 7, both directly
and using an approximation.

\end{example}

\begin{solutionexample}{}{}

The number of people born on January the first follows a binomial
law $\mathcal{B}\left(n,\dfrac{1}{365}\right).$ With such parameters
and values of $n,$ the binomial law can be approximated by a Poisson
law $\mathcal{P}_{\lambda},$ with parameter $\lambda,$ and expression
for every $k\in\mathbb{N},$ 
\[
\mathcal{P}_{\lambda}\left(k\right)=\dfrac{\lambda^{k}}{k!}\text{e}^{-\lambda}.
\]

We have:
\begin{itemize}
\item For $n=500,$ $\lambda_{1}=np\approx1.36986.$
\item For $n=600,$ $\lambda_{2}=np\approx1.64384.$
\end{itemize}
Moreover, we can compute inductively for every $k\in\mathbb{N},$
\begin{itemize}
\item $P_{n}\left(k+1\right)=P_{n}\left(k\right)\dfrac{n-k}{k+1}\dfrac{p}{1-p}.$
\item $\mathcal{P}_{\lambda}\left(k+1\right)=\mathcal{P}_{\lambda}\left(k\right)\dfrac{\lambda}{k+1}.$
\end{itemize}
We then deduce

\begin{table}[H]
\begin{center}%
\begin{tabular}{|c|c|c|c|c|}
\hline 
$k$ &
$P_{500}\left(k\right)$ &
$\mathcal{P}_{\lambda_{1}}\left(k\right)$ &
$P_{600}\left(k\right)$ &
$\mathcal{P}_{\lambda_{2}}\left(k\right)$\tabularnewline
\hline 
0 &
0.2536 &
0.2541 &
0.1928 &
0.1932\tabularnewline
\hline 
1 &
0.3465 &
0.3481 &
0.3178 &
0.3177\tabularnewline
\hline 
2 &
0.2375 &
0.2384 &
0.2615 &
0.2611\tabularnewline
\hline 
3 &
0.1083 &
0.1088 &
0.1432 &
0.1431\tabularnewline
\hline 
4 &
0.0369 &
0.0372 &
0.0587 &
0.0588\tabularnewline
\hline 
5 &
0.0100 &
0.0102 &
0.0192 &
0.0193\tabularnewline
\hline 
6 &
0.0020 &
0.0023 &
0.0052 &
0.0053\tabularnewline
\hline 
7 &
0.0004 &
0.0004 &
0.0012 &
0.0012\tabularnewline
\hline 
\end{tabular}\end{center}

\caption{Binomial law and its approximation for $n=500$ and $n=600$}
\end{table}

\begin{remark}{}{}

We have $\mathcal{P}_{\lambda}\left(k\right)$ small for $k\geqslant5.$
The Chebyshev inequality gives an upper estimate of the dispersion
of a random variable around its expectation. Applying this inequality
to a random variable $S$ following the law $\mathcal{P}\left(\lambda\right),$
with $\lambda=\lambda_{1},$ and using the fact that $\mathbb{E}\left(S\right)=\sigma^{2}_{S}=\lambda_{1},$
we get
\[
P\left(S\geqslant5\right)\leqslant P\left(S\leqslant\mathbb{E}\left(S\right)+2.6\sigma_{S}\right)\leqslant2.6^{-2}\leqslant0.15.
\]

This bound is approximately 10 times greater that the actual value
of $P\left(S\geqslant5\right),$ showing how loose the Chebyshev inequality
can be. 

\end{remark}

\end{solutionexample}

\subsection{Approximation of a Binomial Law by the Gauss Law}

We now focus on the case where $p$ is not sufficiently small and
$n$ is large, and how to approximate the probability $P_{n}\left(k\right).$
The answer is given by the \index{Moivre-Laplace theorems}\textbf{Moivre}{\bfseries\footnote{Abraham De Moivre (1667-1754), an English mathematician of French
origin, refined the principles of probability computation, and the
rules of compound probabilities. He also worked in the domain of finite
differences equations and introduced imaginary numbers in trigonometry.}}\textbf{-Laplace} theorems. We do not provide their proof here; the
interested reader may refer to the book \cite{shiryaev2016probability}.
The proof of the local theorem is based on the Stirling formula\footnote{$n!\sim_{+\infty}n^{n+\frac{1}{2}}\text{e}^{-n}\sqrt{2\pi}$},
and the global theorem relies on the approximation of integrals by
Riemann sums.

\begin{theorem}{Local Moivre-Laplace Theorem}{}

Let $p$ be a real number such that $0<p<1.$ Let $S_{n}$ be a random
variable following the law $\mathcal{B}\left(n,p\right).$ 

\boxeq{As $n$ tends to infinity
\[
P\left(S_{n}=k\right)\sim_{+\infty}\dfrac{1}{\sqrt{2npq}}\text{e}^{-\frac{\left(k-np\right)^{2}}{2npq}},
\]
 uniformly in $k$ such that
\[
\left|k-np\right|=o\left[\left(npq\right)^{\frac{2}{3}}\right].
\]
}

\end{theorem}

\begin{theorem}{Global Moivre-Laplace Theorem}{}

Let $p$ be a real number such that $0<p<1.$ Let $S_{n}$ be a random
variable following the law $\mathcal{B}\left(n,p\right),$ and let
$\widetilde{S_{n}}=\dfrac{S_{n}-\mathbb{E}\left(S_{n}\right)}{\sigma\left(S_{n}\right)}$
be the centered reduced variable associated with $S_{n}.$ Then

\boxeq{
\[
\lim_{n\to+\infty}{}_{-\infty\leqslant a<b\leqslant+\infty}\left|P\left(a<\widetilde{S_{n}}\leqslant b\right)-\dfrac{1}{\sqrt{2\pi}}\intop^{b}_{a}\text{e}^{-\frac{x^{2}}{2}}\text{\,d}x\right|=0.
\]
}

\end{theorem}

The law of $\widetilde{S_{n}}$ can then be approximated when $n$
is large, by the center reduced Gauss law $\mathcal{N}\left(0,1\right).$
A straightforward computation then yields an approximation for the
law of $S_{n}$ itself.

\begin{corollary}{}{}

If $S_{n}$ is a random variable following the binomial law $\mathcal{B}\left(n,p\right),$
then for every real numbers $a$ and $b$ such that $a<b,$ denoting
$q=1-p,$

\boxeq{
\[
\lim_{n\to+\infty}\left(P\left(a<S_{n}\leqslant b\right)-\dfrac{1}{\sqrt{2\pi}}\intop^{\frac{b-np}{\sqrt{npq}}}_{\frac{a-np}{\sqrt{npq}}}\text{e}^{-\frac{x^{2}}{2}}\text{\,d}x\right)=0.
\]
}

\end{corollary}

We give on Figure \ref{Fig Binomial and Gauss approx} an example
comparing the binomial law with a Gauss law having the same expectation
and the same standard deviation.

\begin{figure}
\begin{center}\begin{tikzpicture}
	\definecolor{orange}{HTML}{FFA200}
	\definecolor{purple}{HTML}{9500FF}
  \begin{axis}[
		axis x line=center,
		axis y line=center,
		xtick = {-2,0,2,4,6,8,10,12,14,16},
		ytick = {0,0.1,0.2},
		ymin=-0.05,
		ymax=0.25
	]
 \addplot[green, samples=80, smooth, thick, domain=-2.25:16.5 ]
      {1/(sqrt(pi*2*3.6))*exp(-(x-6)^2/(2*3.6)};
	\pgfplotsinvokeforeach {-1,0,1,2,3,4,5,6,7,8,9,10,11,12,13,14} {
    	\draw[orange, line width=3] (axis cs: #1,0)
      -- (axis cs: #1,{(1/(sqrt(2*3.6*pi)))*exp((-1/7.2)*(#1-6)^2)});
  }
  \end{axis}

\end{tikzpicture}\end{center}

\caption{Comparison of the binomial law $\mathcal{B}\left(15,0.4\right)$ and
of the Gauss law of same expectation and standard deviation.}
\label{Fig Binomial and Gauss approx}
\end{figure}

\begin{remark}{}{}

It can be shown---as a particular case of the Berry-Esseen theorem---that

\boxeq{
\[
_{x\in\mathbb{R}}\left|P\left(S_{n}\leqslant x\right)-\varphi\left(x\right)\right|\leqslant\dfrac{p^{2}+q^{2}}{\sqrt{npq}},
\]
}

where $\varphi$ denotes the cumulative distribution function of the
Gauss law $\mathcal{N}\left(0,1\right).$ However, this approximation
may be too imprecise for values in the neighborhood of 0 or 1. In
such cases, the Poisson theorem proves more useful.

Recall that for every $x\in\mathbb{R},$
\[
\Phi\left(x\right)=\dfrac{1}{\sqrt{2\pi}}\intop^{x}_{-\infty}\text{e}^{-\frac{t^{2}}{2}}\text{\,d}t,
\]
 and that if $X$ follows the law $\mathcal{N}\left(0,1\right),$
then
\[
P\left(a<X\leqslant b\right)=\dfrac{1}{\sqrt{2\pi}}\intop^{b}_{a}\text{e}^{-\frac{x^{2}}{2}}\text{\,d}x=\Phi\left(b\right)-\Phi\left(a\right).
\]

The function $\Phi$ is typically tabulated for $x>0,$ and for $x>0,$
we use the identity
\[
\Phi\left(-x\right)=1-\Phi\left(x\right).
\]

Below are three commonly nondecreasing values of the function $\Phi:$
\begin{align*}
\Phi\left(1.64\right)-\Phi\left(-1.64\right) & \approx0.9\\
\Phi\left(1.96\right)-\Phi\left(-1.96\right) & \approx0.95\\
\Phi\left(3.09\right)-\Phi\left(-3.09\right) & \approx0.99
\end{align*}

Figure \ref{Fig area under normal} illustrates the interpretation
of this integral when $X$ follows the standard normal law $\mathcal{N}\left(0,1\right).$

\end{remark}

\begin{figure}
\begin{center}    \begin{tikzpicture}
	\definecolor{purple}{HTML}{9500FF}
  \begin{axis}[
		axis x line=center,
		axis y line=center,
		xtick = \empty,
		ytick =\empty,
		xmin=-4.5,
		xmax=4.5,
		ymin=-0.05,
		ymax=0.45
	]
 \addplot[green, samples=80, smooth, thick, domain=-4.5:4.5, name path = dist ]
      {1/(sqrt(pi*2))*exp(-x^2/2)};

\pgfplotsinvokeforeach {-1.96,1.96} {
    	\draw[orange, thick, opacity=0.5] (axis cs: #1,0)
      -- (axis cs: #1,{(1/sqrt(2*pi))*exp((-1/2)*(#1)^2)});
  }

  % This is necessary for the filling later
\path [name path = base] (\pgfkeysvalueof{/pgfplots/xmin},0)
    -- (\pgfkeysvalueof{/pgfplots/xmax},0);
   
\node[orange] at (axis cs: -1.96,-0.03) {$-1.96$};
\node[orange] at (axis cs: 1.96,-0.03) {$1.96$};

\addplot [fill=orange, fill opacity=0.05] fill between [of = dist and base, soft clip = {domain=-2:2}];

  \end{axis}

\end{tikzpicture}\end{center}

\caption{$P\left(-1.96\leqslant X\leqslant1.96\right)=0.95$ when $X$ follows
the standard normal law $\mathcal{N}\left(0,1\right)$}
\label{Fig area under normal}
\end{figure}

\begin{example}{Die Rolling}{}

We roll a fair die 12,000 times. 

We are interested in the probability that the number of times a six
is obtained lies between 1,800 and 2,100.

\end{example}

\begin{solutionexample}{}{}

The number of sixes obtained is a random variable $X$ following the
binomial law $\mathcal{B}\left(12,000;\dfrac{1}{6}\right).$

The expectation of $X$ is
\[
\mathbb{E}\left(X\right)=12,000\times\dfrac{1}{6}=2,000
\]
and the standard deviation is $\sigma_{X}=\sqrt{12,000\times\dfrac{1}{6}\times\dfrac{5}{6}}=\sqrt{\dfrac{5,000}{3}}.$

Here, $a=1,800$ and $b=2,100,$ so the approximative value of the
desired probability is

\boxeq{
\[
\Phi\left(\dfrac{2,100-2,000}{\sqrt{\dfrac{10,000}{6}}}\right)-\Phi\left(\dfrac{1,800-2,000}{\sqrt{\dfrac{10,000}{6}}}\right)=\Phi\left(\sqrt{6}\right)-\Phi\left(-2\sqrt{6}\right)\approx0.992.
\]
}

\end{solutionexample}

\subsection{Approximation of a Hypergeometric Law by a Binomial Law}

We recall the definition of the hypergeometric law through an example.
A lake contains $r$ fishes\footnote{The reader may easily adapt this example to polling.},
among which $r_{1}$ belong to an interesting species $a.$ We catch
$n$ fishes, assuming that each fish has the same probability of being
caught. Let $X$ denote the number of caught fish of species $a.$
Then
\[
P\left(X=k\right)=\dfrac{\binom{r_{1}}{k}\binom{r-r_{1}}{n-k}}{\binom{r}{n}}
\]
when the right hand side is defined, and $P\left(X=k\right)=0$ otherwise.
We say that $X$ follows a hypergeometric law of parameters $n,$
$r$ and $r_{1}.$

If, instead, we consider a catch-and-release process---no-kill fishing---the
number $Y$ of caught fish of species $a$ follows the binomial law
$\mathcal{B}\left(n,\dfrac{r_{1}}{r}\right).$ For $0\leqslant k\leqslant n,$
\[
P\left(Y=k\right)=\binom{n}{k}\left(\dfrac{r_{1}}{r}\right)^{k}\left(1-\dfrac{r_{1}}{r}\right)^{n-k},
\]
 and $P\left(Y=k\right)$ otherwise.

It is natural to expect that if $r$ and $r_{1}$ are large compared
to $n,$ the difference between the two processes is slight. This
is confirmed by the following statement.

\begin{theorem}{The Binomial Law as an Approximation of the Hypergeometric Law}{}

Consider, on the one hand, for each $j\in\mathbb{N},$ a random variable
$X_{j}$ following the hypergeomtric law with parameters $n,r_{j},r^{1}_{j};$
and on the other hand, a random variable $Y$ following the binomial
law $\mathcal{B}\left(n,p\right).$

If
\[
\lim_{j\to+\infty}\dfrac{r^{1}_{j}}{r_{j}}=p,
\]
where $p\in\left]0,1\right[,$ then

\boxeq{
\[
\lim_{j\to+\infty}P\left(X_{j}=k\right)=P\left(Y=k\right)=\binom{n}{k}p^{k}\left(1-p\right)^{n-k}.
\]
}

\end{theorem}

\begin{proof}{}{}

As soon as we have $k\leqslant r^{1}_{j}$ and $n-k\leqslant r_{j}-r^{1}_{j},$
we obtain, after simplification of the binomial coefficients
\[
P\left(X_{j}=k\right)=\binom{n}{k}\prod^{k-1}_{l=0}\left(\dfrac{r^{1}_{j}-l}{r_{j}-l}\right)\prod^{\left(n-k\right)-1}_{l=0}\left(\dfrac{r_{j}-r^{1}_{j}-l}{r_{j}-l}\right),
\]
which readily ensures the anounced convergence and thus proves the
result.

\end{proof}

\subsection{Central Limit Theorem}

Results concerning the approximation of laws---or more precisely,
the behaviour of certain probability laws as a parameter (either an
integer or a real number) tends to infinity---are often designated
as \textbf{limit theorems}\index{limit theorems}. The most important
of these is known as the \index{central limit theorem}\textbf{central
limit theorem}. We give here an elementary version, without proof,
and show how it applies in practice. 

\begin{theorem}{Central Limit Theorem}{}

Let $\left(X_{n}\right)_{n\in\mathbb{N}}$ be a sequence of independent,
non-constant random variables, following the same law and admitting
a second-order moment, with expectation $\mathbb{E}\left(X_{1}\right)$
and standard deviation $\sigma_{X_{1}}.$ 

For each $n\in\mathbb{N},$ we define the random variable $S_{n},$
corresponding to the centered and reduced variable associated to $\sum^{n}_{j=1}X_{j},$
by
\[
S_{n}=\dfrac{\sum^{n}_{j=1}X_{j}-n\mathbb{E}\left(X_{1}\right)}{\sqrt{n}\sigma_{X_{1}}}.
\]

Then, for every pairs $\left(a,b\right)\in\overline{\mathbb{R}}^{2}$
such that $a<b,$

\boxeq{
\[
P\left(a<S_{n}\leqslant b\right)\underset{n\to+\infty}{\longrightarrow}\dfrac{1}{\sqrt{2\pi}}\intop^{b}_{a}\text{e}^{-\frac{x^{2}}{2}}\text{\,d}x.
\]
}

\end{theorem}

\begin{remark}{}{}

The global Moivre-Laplace theorem is a particular case of the central
limit theorem: indeed, a random variable $S_{n}$ following the binomial
law $\mathcal{B}\left(n,p\right)$ has the same law as the sum $X_{1}+\dots+X_{n},$
where $X_{1},\dots,X_{n}$ are independent random variables following
the Bernoulli law with parameter $p.$

\end{remark}

\begin{application}{}{}

If a random variable $X$ can be written $\sum^{n}_{j=1}X_{j},$ where
the random variables $X_{j}$ satisfy the central limit theorem hypotheses,
then the law of $X$ can be approximated by the Gauss law $\mathcal{N}\left(\mathbb{E}\left(X\right),\sigma^{2}_{X}\right)$
for $n$ sufficiently large.

Besides the case of random variables following the binomial law $\mathcal{B}\left(n,p\right),$
this also applies to random variables following the Poisson law $\mathcal{P}\left(\lambda\right)$---which
has the same law as the sum of $n$ random variables of law $\mathcal{P}\left(\dfrac{\lambda}{n}\right)$
as mentioned in Example $\ref{ex:sum_rv_poisson_law}$. This is illustrated
in Figure \ref{Fig Gauss-Poisson}. It also applies to random variables
following the Chi-Squared law $\chi^{2}_{n}$---which by definition
is the law of the sum of $n$ independent random variables following
the same law $\chi^{2}_{1}$---; this is illustrated in Figure \ref{Fig Gauss-Chi-Two}.

The approximations in Table \ref{Table: Approximation by Gauss law}
are valid when $n$ and $\lambda$ are sufficiently large.

\begin{table}[H]
\begin{center}%
\begin{tabular}{|c|c|c|c|}
\hline 
$P_{X}$ &
$P_{X_{j}}$ &
$\mathbb{E}\left(X\right)$ &
$\sigma^{2}_{X}$\tabularnewline
\hline 
$\mathcal{B}\left(n,p\right)$ &
$\mathcal{B}\left(1,p\right)$ &
$np$ &
$npq$\tabularnewline
\hline 
$\mathcal{P}\left(\lambda\right)$ &
$\mathcal{P}\left(\dfrac{\lambda}{n}\right)$ &
$\lambda$ &
$\lambda$\tabularnewline
\hline 
$\chi^{2}_{n}$ &
$\chi^{2}_{1}$ &
$n$ &
$2n$\tabularnewline
\hline 
\end{tabular}\end{center}

\caption{Approximation by a Gauss law for $n$ and $\lambda$ sufficiently
big.}
\label{Table: Approximation by Gauss law}
\end{table}

\end{application}

\begin{figure}
\begin{center}\begin{tikzpicture}[
    declare function={poissonlaw(\t,\n)=\t^\n/\n!*exp(-\t);}
]
	\definecolor{orange}{HTML}{FFA200}
	\definecolor{purple}{HTML}{9500FF}
  \begin{axis}[
		axis x line=center,
		axis y line=center,
		xtick = {-2,0,2,4,6,8,10,12,14},
		ytick = {0,0.05,0.10,0.15,0.20},
		ymin=-0.02,
		ymax=0.22
	]
\addplot[orange, samples=80, smooth, thick, domain=-4.5:14.5 ]{1/(sqrt(pi*2*3.7))*exp(-(x-3.7)^2/(2*3.7)};
\pgfplotsinvokeforeach {-4,-3,-2,-1,0,1,2,3,4,5,6,7,8,9,10} {
    	\draw[green!80!black, line width=3] (axis cs: #1,0)
      -- (axis cs: #1,{poissonlaw(3.7,#1)});
  }
  \end{axis}

\end{tikzpicture}\end{center}

\caption{Comparison of a Poisson law $\mathcal{P}\left(\lambda\right)$---with
$\lambda=3.7$---and a Gauss law of same expectation and variance.}
\label{Fig Gauss-Poisson}
\end{figure}

\begin{figure}
\begin{center}\begin{tikzpicture}[
    declare function={poissonlaw(\t,\n)=\t^\n/\n!*exp(-\t);}
]
	\definecolor{orange}{HTML}{FFA200}
	\definecolor{purple}{HTML}{9500FF}
  \begin{axis}[
		axis x line=center,
		axis y line=center,
		xtick = {-10,0,10,20,30,40,50},
		ytick = {0,0.01,0.02,0.03,0.04,0.05,0.06},
		ymin=-0.005,
		ymax=0.07,
		restrict y to domain = -1:1
	]
\addplot[green, samples=80, smooth, thick, domain=-10:50 ]{1/(sqrt(pi*2*40))*exp(-(x-20)^2/(2*40)};
\addlegendentry{$\mathcal{N}\left(20,40\right)$}
\addplot[purple, samples=80, smooth, thick, domain=0:50 ]{(x^9)/371226240*(exp(-(x/2))};
\addplot[purple, samples=80, smooth, thick, domain=-10:0 ]{0};
\addlegendentry{$\chi^2_{20}$}

  \end{axis}

\end{tikzpicture}\end{center}

\caption{Comparison of a Chi-Squared law with 20 degrees of freedom and a Gauss
law of same expectation and variance.}
\label{Fig Gauss-Chi-Two}
\end{figure}

\subsection{Convergence in Law}

\begin{definition}{Convergence in Law}{}

A sequence $\left(X_{n}\right)_{n\in\mathbb{N}}$ of real-valued random
variables is said to \index{converge in law}\textbf{converge in law}
to a random variable $X$ if, for every point where the cumulative
distribution function $F_{X}$ of the random variable $X$ is continuous,
\[
F_{X_{n}}\left(x\right)\underset{n\to+\infty}{\longrightarrow}F_{X}\left(x\right).
\]

\end{definition}

\begin{remarks}{}{}
\begin{itemize}
\item Convergence cannot be enforced at every point. Indeed, taking $X_{n}=\dfrac{1}{n},$
it follows naturally that the sequence $\left(X_{n}\right)_{n\in\mathbb{N}^{\ast}}$
converges in law to the random variable $X=0.$ Nonetheless,
\[
F_{X_{n}}\left(0\right)=P\left(X_{n}\leqslant0\right)=0
\]
and
\[
F_{X}\left(0\right)=P\left(X\leqslant0\right)=1.
\]
\item It may happen that the sequence of cumulative distribution function
$F_{X_{n}}$ converges pointwise to a function that is not a cumulative
distribution function. For instance, if $X_{n}=n,$ then
\[
F_{X_{n}}\left(x\right)=\begin{cases}
0, & \text{if }x<n,\\
1, & \text{if }x\geqslant n,
\end{cases}
\]
so for every fixed real number $x,$
\[
\lim_{n\to+\infty}F_{X_{n}}\left(x\right)=0.
\]
\item The random variables $X_{n}$ and $X$ need not to be defined on the
same probabilized space; it is a concept referring to the convergence
of the sequence of the laws followed by these random variables.
\item The conclusion of the central limit theorem can be reformulated as:
the sequence of random variables $S_{n}$ converges in law to the
Gauss law $\mathcal{N}\left(0,1\right).$
\end{itemize}
\end{remarks}

\section{Weak Law of Large Numbers}\label{sec:Weak-Law-of}

The study of random phenomena often leads to the analysis of the sequence
of arithmetic means of a sequence of independent random variables
sharing the same law. This arises in particular in statistics, when,
for example, when estimating a parameter of the law followed by a
random variable associated with a phenomenon, based on the sequence
of outcomes obtained from independent repetitions of the experiment. 

The study of the convergence of such sequences is addressed by results
known as laws of large numbers: weak laws, which concern convergence
in probability---a concept we will define shorly---and strong laws
which pertain to almost sure convergence---and that will be studied
later in Part \ref{part:Deepening-Probability-Theory}.

\begin{definition}{Congergence in Probability}{}

A sequence of real valued random variables $\left(X_{n}\right)_{n\in\mathbb{N}}$
is said to \index{converge in probability}\textbf{converge in probability}
to a random variable $X$ if, for every $\epsilon>0,$\\
\boxeq{
\[
\lim_{n\to+\infty}P\left(\left|X_{n}-X\right|>\epsilon\right)=0.
\]
}

\end{definition}

\begin{theorem}{Weak Law of Large Numbers}{weak_law}

Let $\left(X_{n}\right)_{n\in\mathbb{N}^{\ast}}$ be a sequence of
independent random variables defined on a probabilized space $\left(\Omega,\mathcal{A},P\right)$
each admitting a second-order moment. 

Suppose the convergence of the sequences
\[
\begin{array}{ccccc}
\dfrac{1}{n}\sum\limits^{n}_{j=1}\mathbb{E}\left(X_{j}\right)\underset{n\to+\infty}{\longrightarrow}m & \,\,\,\, & \text{and} & \,\,\,\, & \dfrac{1}{n^{2}}\sum\limits^{n}_{j=1}\sigma^{2}_{X_{j}}\underset{n\to+\infty}{\longrightarrow}0.\end{array}
\]

Then, the sequence of random variables $\overline{X_{n}}=\dfrac{1}{n}\sum\limits^{n}_{j=1}X_{j}$
\textbf{converges in probability} to $m.$

\end{theorem}

\begin{proof}{}{}

We have
\[
\mathbb{E}\left(\overline{X_{n}}\right)=\dfrac{1}{n}\sum\limits^{n}_{j=1}\mathbb{E}\left(X_{j}\right).
\]

Since the random variables $X_{n}$ are independent, it follows that
\[
\sigma^{2}_{\overline{X_{n}}}=\dfrac{1}{n^{2}}\sum^{n}_{j=1}\sigma^{2}_{X_{j}}.
\]

The triangular inequality gives
\[
\left|\overline{X_{n}}-m\right|\leqslant\left|\overline{X_{n}}-\dfrac{1}{n}\sum\limits^{n}_{j=1}\mathbb{E}\left(X_{j}\right)\right|+\left|\dfrac{1}{n}\sum\limits^{n}_{j=1}\mathbb{E}\left(X_{j}\right)-m\right|.
\]

By hypothesis, the sequence $\left(\dfrac{1}{n}\sum^{n}_{j=1}\mathbb{E}\left(X_{j}\right)\right)_{n\in\mathbb{N}}$
converges to $m.$ Therefore, for every $\epsilon>0,$ there exists
$N\left(\epsilon\right)\in\mathbb{N}^{\ast}$ such that for every
$n\geqslant N\left(\epsilon\right),$
\[
\left|\dfrac{1}{n}\sum^{n}_{j=1}\mathbb{E}\left(X_{j}\right)-m\right|\leqslant\dfrac{\epsilon}{2}.
\]

For every $n\geqslant N\left(\epsilon\right),$ this implies the inclusion
of sets
\[
\left(\left|\overline{X_{n}}-\dfrac{1}{n}\sum\limits^{n}_{j=1}\mathbb{E}\left(X_{j}\right)\right|\leqslant\dfrac{\epsilon}{2}\right)\subset\left(\left|\overline{X_{n}}-m\right|\leqslant\epsilon\right).
\]

Taking complement, we get
\[
\left(\left|\overline{X_{n}}-m\right|>\epsilon\right)\subset\left(\left|\overline{X_{n}}-\dfrac{1}{n}\sum\limits^{n}_{j=1}\mathbb{E}\left(X_{j}\right)\right|>\dfrac{\epsilon}{2}\right).
\]

Applying the Bienaymé-Chebishev inequality, we obtain
\[
P\left(\left|\overline{X_{n}}-\dfrac{1}{n}\sum\limits^{n}_{j=1}\mathbb{E}\left(X_{j}\right)\right|>\dfrac{\epsilon}{2}\right)\leqslant\dfrac{4}{\epsilon^{2}}\sigma^{2}_{\overline{X_{n}}}=\dfrac{4}{\epsilon^{2}}\dfrac{1}{n^{2}}\sum^{n}_{j=1}\sigma^{2}_{X_{j}}.
\]

Hence, for every $n\geqslant N\left(\epsilon\right),$
\[
P\left(\left|\overline{X_{n}}-m\right|>\epsilon\right)\leqslant\dfrac{4}{\epsilon^{2}}\dfrac{1}{n^{2}}\sum^{n}_{j=1}\sigma^{2}_{X_{j}},
\]

Using the second hypothesis, proves the result.

\end{proof}

\begin{remark}{}{}

In particular, the assumptions of the previous theorem are satisfied
when the random variables $X_{n}$ are independent and follow the
same law, and when $X_{1}$ admits a second-order moment.

We now study a particular case of the previous theorem, the \index{Bernoulli theorem}Bernoulli
theorem, even though it historically predates it.

\end{remark}

\begin{theorem}{Bernoulli Theorem}{bernoulli_th}

Let $\left(A_{n}\right)_{n\in\mathbb{N}}$ be a sequence of independent
events, each of the same probability $p.$ Then the sequence of random
variables 
\[
\dfrac{1}{n}\sum^{n}_{j=1}\boldsymbol{1}_{A_{n}}
\]
 converges in probability to $p.$

\end{theorem}

\begin{proof}{}{}

The random variables $\boldsymbol{1}_{A_{n}},\,\,n\in\mathbb{N}^{\ast}$
are independent and follow the same Bernoulli law. They admit a second-order
moment, and thus correspond to the situation of the previous remark. 

\end{proof}

\begin{remark}{}{}

This theorem ensures that, when studying a given random phenomenon,
if we perform a sequence of independent experiments, the sequence
of \textbf{relative frequencies\index{relative frequencies}} of a
given property associated with the phenomenon converges in probability---in
the sense of the probability $P$ induced by the chosen model---to
the probability of the event corresponding to that property.

The Bernoulli theorem is therefore a coherence theorem for the probabilistic
model within the frequentist interpretation of the probability of
an event, point of view that was not especially adopted in this book,
but hat is at the origin of the calculus of probabilities...!

\end{remark}

\section*{Exercises}

\addcontentsline{toc}{section}{Exercises}

\begin{exercise}{Central Limit Theorem and Poisson Law}{exercise7.1}

Prove, using the central limit theorem, that
\[
\lim_{n\to+\infty}\exp\left(-n\right)\sum^{n}_{k=0}\dfrac{n^{k}}{k!}=\dfrac{1}{2}.
\]

\end{exercise}

\begin{exercise}{Numerical Comparison: Central Limit Theorem and Chebyshev Inequality}{exercise7.2}

Let $\left(X_{n}\right)_{n\in\mathbb{N}}$ be a sequence of independent
random variables following the same law and admitting a second-order
moment.

For every $n\in\mathbb{N}^{\ast},$ define the random variable $S_{n}$
by
\[
S_{n}=\dfrac{\sum^{n}_{j=1}X_{j}-n\mathbb{E}\left(X_{1}\right)}{\sqrt{n}\sigma_{X_{1}}}.
\]

Compare, for $\epsilon=1,2,3$ the numerical information given by
the probability $P\left(\left|S_{n}\right|\geqslant\epsilon\right)$
estimated by the central limit theorem, and then by the Chebishev
inequality.

\end{exercise}

\begin{exercise}{Size of Sampling and Central Limit Theorem}{exercise7.3}

Let $X$ be a random variable with unknown expectation $m$ and standard
deviation $\sigma.$ Let $\left(X_{n}\right)_{n\in\mathbb{N}^{\ast}}$
be a sequence of independent random variables following the same law
as $X$---we say that $\left(X_{1},X_{2},\dots,X_{n}\right)$ is
a sample of size $n$ of $X$---; we estimate the expectation $m$
by the random variable $\overline{X_{n}}=\dfrac{1}{n}\sum^{n}_{j=1}X_{j}.$

1. Justify the approximation, for $\alpha>0$ and for large $n$
\begin{equation}
P\left(\left|\overline{X_{n}}-m\right|\geqslant\alpha\right)\approx1-\left[\varphi\left(\dfrac{\sqrt{n}}{\sigma}\alpha\right)-\varphi\left(-\dfrac{\sqrt{n}}{\sigma}\alpha\right)\right].\label{eq:approx_expectation}
\end{equation}

2. What is the minimum size of the sample such that, for fixed $\alpha>0$
and $\beta\in\left]0,1\right[,$ we have
\[
P\left(\left|\overline{X_{n}}-m\right|\geqslant\alpha\right)\leqslant\beta?
\]

Compute the numerical value in the case where $\sigma=3,$ $\alpha=0.05\sigma,$
and $\beta=0.05.$

\textit{Hint: Recall that
\[
\varphi^{-1}\left(0.975\right)\approx1.96.
\]
}

\end{exercise}

\begin{exercise}{Polls and Moivre-Laplace Theorem}{exercise7.4}

We want to estimate the percentage $p$ of positive answers in a referendum.
To achieve this goal, we conduct a poll on $n$ persons and estimate
$p$ by the relative frequency $F_{n}$ of ``yes'' answers among
the polled individuals. 

1. What is the smallest integer $n_{0}$ such that the probability
that $F_{n}$ does not differ from $p$ by more than $\alpha>0$ is
less than $\beta\in\left]0,1\right[?$

2. Numerical application:

a. We choose $\beta=0.05$ and $\alpha=0.01.$ 

Study the two following cases:

i. Where it is known that $0<p<0.3$

ii. When $p$ is completely unknown.

b. By how much is the sample size reduced when we choose $\alpha=0.05?$

\end{exercise}

\begin{exercise}{Roulette Game and Central Limit Theorem}{exercise7.5}

The probability of winning at a roulette game is $\dfrac{19}{37},$
and the gambler's stake is one token---which has a counter value
in the currency of the casino's location. We adopt the point of view
of the casino owner.

1. What is the minimum number $n_{0}$ of games to be played per day
so that the casino wins at least 1,000 tokens per day with probability
0.5?

2. What is the probability of a total loss for the casino during those
$n_{0}$ games?

\end{exercise}

\section*{Solutions of Exercises}

\addcontentsline{toc}{section}{Solutions of Exercises}

\begin{solution}{}{solexercise7.1}

\textbf{Proof of $\lim_{n\to+\infty}\text{e}^{-n}\sum^{n}_{k=0}\dfrac{n^{k}}{k!}=\dfrac{1}{2}$}

Let $\left(X_{n}\right)_{n\in\mathbb{N}^{\ast}}$ be a sequence of
independent random variables following the same Poisson law with parameter
1, and let $S_{n}$ be the random variable $S_{n}=\sum^{n}_{j=1}X_{j}.$

The law of $S_{n}$ is the Poisson law with parameter $n.$ Thus,
\[
P\left(S_{n}\leqslant n\right)=\text{e}^{-n}\sum^{n}_{k=0}\dfrac{n^{k}}{k!}.
\]

Applying the central limit theorem, it follows, with the notation
of this theorem, that\\
\boxeq{
\[
P\left(S_{n}\leqslant n\right)=P\left(\dfrac{S_{n}-n}{\sqrt{n}}\leqslant0\right)\underset{n\to+\infty}{\longrightarrow}\phi\left(0\right)=\dfrac{1}{2},
\]
} which proves the theorem.

\end{solution}

\begin{solution}{}{solexercise7.2}

The random variables, for $1\leqslant j\leqslant n,$
\[
\mathring{X_{j}}=\dfrac{X_{j}-\mathbb{E}\left(X_{1}\right)}{\sigma_{X_{1}}}
\]
are centered and reduced, and
\[
S_{n}=\dfrac{1}{\sqrt{n}}\sum^{n}_{j=1}\mathring{X_{j}}.
\]

Since the random variables $\mathring{X_{j}}$ are independent and
follow the same law,
\[
\sigma^{2}_{S_{n}}=\dfrac{1}{n}\left(n\sigma^{2}_{\mathring{X_{1}}}\right)=1.
\]

The Chebyshev inegality then gives, for every $\epsilon>0,$
\[
P\left(\left|S_{n}\right|\geqslant\epsilon\right)\leqslant\dfrac{1}{\epsilon^{2}}.
\]

We have:
\begin{itemize}
\item For $\epsilon=1,$ the Chebyshev inequality corresponds to
\[
P\left(\left|S_{n}\right|\geqslant1\right)\leqslant1
\]
which provides no information, while the central limit theorem yields
\[
P\left(\left|S_{n}\right|\geqslant1\right)\underset{n\to+\infty}{\longrightarrow}2\left(1-\phi\left(1\right)\right)\approx0.3173.
\]
\item For $\epsilon=2,$ the Chebyshev inequality corresponds to
\[
P\left(\left|S_{n}\right|\geqslant2\right)\leqslant0.25
\]
while the central limit theorem yields
\[
P\left(\left|S_{n}\right|\geqslant2\right)\underset{n\to+\infty}{\longrightarrow}2\left(1-\phi\left(2\right)\right)\approx0.0455.
\]
\item For $\epsilon=3,$ the Chebyshev inequality corresponds to
\[
P\left(\left|S_{n}\right|\geqslant3\right)\leqslant\dfrac{1}{9}\approx0.1111
\]
while the central limit theorem yields
\[
P\left(\left|S_{n}\right|\geqslant2\right)\underset{n\to+\infty}{\longrightarrow}2\left(1-\phi\left(3\right)\right)\approx0.0027.
\]
\end{itemize}
\end{solution}

\begin{solution}{}{solexercise7.3}

\textbf{1. Justification of the approximation $P\left(\left|\overline{X_{n}}-m\right|\geqslant\alpha\right)\approx1-\left[\varphi\left(\dfrac{\sqrt{n}}{\sigma}\alpha\right)-\varphi\left(-\dfrac{\sqrt{n}}{\sigma}\alpha\right)\right]$}

We set
\[
S_{n}=\dfrac{\sum^{n}_{j=1}X_{j}-nm}{\sqrt{n}\sigma}\equiv\sqrt{n}\dfrac{\overline{X_{n}}-m}{\sigma}.
\]

We then have the equality
\[
\left(\left|\overline{X_{n}}-m\right|\geqslant\alpha\right)=\left(\left|S_{n}\right|\geqslant\dfrac{\sqrt{n}}{\sigma}\alpha\right).
\]

By the central limit theorem, we have for $a>0,$
\[
P\left(\left|S_{n}\right|\geqslant a\right)\underset{n\to+\infty}{\longrightarrow}1-\left[\varphi\left(a\right)-\varphi\left(-a\right)\right].
\]

Thus, the approximation $\refpar{eq:approx_expectation}$ obtained
for $n$ sufficiently large.

\textbf{2. Minimum sample size to have $P\left(\left|\overline{X_{n}}-m\right|\geqslant\alpha\right)\leqslant\beta$}

Since
\[
1-\left[\Phi\left(a\right)-\Phi\left(-a\right)\right]=2\left[1-\Phi\left(a\right)\right],
\]
we choose $n$ as the smallest integer such that
\[
2\left[1-\Phi\left(\dfrac{\sqrt{n}}{\sigma}\alpha\right)\right]\leqslant\beta,
\]
which corresponds to the smallest integer such that\\
\boxeq{
\[
n\geqslant\dfrac{\sigma^{2}}{\alpha^{2}}\left[\Phi^{-1}\left(1-\dfrac{\beta}{2}\right)\right]^{2}
\]
}

\textbf{Numerical application with $\sigma=3,$ $\alpha=0.05\sigma$
and $\beta=0.05$}

The minimum size $n$ of the sample is the smallest ingeger such that
\[
n\geqslant\dfrac{10^{4}}{25}\left(\Phi^{-1}\left(0.975\right)\right)^{2},
\]
thus,\boxeq{
\[
n=1537.
\]
}

Tr.N.: This means that, to ensure to have less than a 5\% risk that
the sample mean deviates from the unknown expectation by more than
5\% of the standard deviation, ones need a sample of size at least
1537.

\end{solution}

\begin{solution}{}{solexercise7.4}

\textbf{1. Smallest integer $n_{0}$}

Let $\left(X_{n}\right)_{n\in\mathbb{N}^{\ast}}$ be a sequence of
independent random variables following the same Bernoulli law $\mathcal{B}\left(p\right).$
We estimate $p$ by the random variable 
\[
F_{n}=\dfrac{1}{n}\sum^{n}_{j=1}X_{j}.
\]
We are in the particular case of Exercise $\ref{exo:exercise7.3}$
to which we refer. 

Here,
\[
\begin{array}{ccccc}
\mathbb{E}\left(X_{1}\right)=p & \,\,\,\, & \text{and} & \,\,\,\, & \sigma_{X_{1}}=\sqrt{pq},\end{array}
\]
where we set $q=1-p.$

Setting
\[
Y_{n}=\dfrac{\sum^{n}_{j=1}X_{j}-np}{\sqrt{npq}}\equiv\sqrt{n}\dfrac{F_{n}-p}{\sqrt{pq}},
\]
it follows that
\[
\left(\left|F_{n}-p\right|\geqslant\alpha\right)=\left(\left|Y_{n}\right|\geqslant\dfrac{\sqrt{n}}{\sigma_{X_{1}}}\alpha\right).
\]

Applying the central limit theorem, $n_{0}$ is the smallest value
of $n$ such that
\begin{equation}
n\geqslant\dfrac{pq}{\alpha^{2}}\left(\varphi^{-1}\left(1-\dfrac{\beta}{2}\right)\right)^{2}.\label{eq:n_0 inequation}
\end{equation}

Tr.N.: Thus, writing the integer part with $\left\lfloor .\right\rfloor ,$

\boxeq{
\[
n_{0}=\left\lfloor \underset{p\in\left]p_{\text{min}},p_{\text{max}}\right[}{\max}\dfrac{p\left(1-p\right)}{\alpha^{2}}\left(\Phi^{-1}\left(1-\dfrac{\beta}{2}\right)\right)^{2}\right\rfloor +1
\]
} where $p_{\text{max}}$ is the maximum value $p$ can take, and
$p_{\text{min}}$ its least value.

2. Numerical application

Since we have chosen $\beta=0.05$ and $\alpha=0.01$, it yields
\[
n_{0}=\left\lfloor \underset{p\in\left]p_{\text{min}},p_{\text{max}}\right[}{\max}\dfrac{p\left(1-p\right)}{10^{-4}}\left(\varphi^{-1}\left(0.975\right)\right)^{2}\right\rfloor +1,
\]
 thus\boxeq{
\[
n_{0}=\left\lfloor \underset{p\in\left]p_{\text{min}},p_{\text{max}}\right[}{\max}\left(196\right)^{2}p\left(1-p\right)\right\rfloor +1.
\]
}
\begin{itemize}
\item When $p<0.3,$ 
\[
\max_{p\in\left]0,0.3\right[}p\left(1-p\right)=0.3\times0.7
\]
thus,
\[
n_{0}=8,068.
\]
\item If $p$ is unknown
\[
\max_{p\in\left]0,1\right[}p\left(1-p\right)=\dfrac{1}{4}
\]
thus,
\[
n_{0}=9,604.
\]
\end{itemize}
If we now choose $\beta=0.05$ and $\alpha=0.05,$ it follows from
the inequality $\refpar{eq:n_0 inequation}$ that $n_{0}$ is divided
by 25. Thus, the minimal sizes obtained are 323 and 384 respectively.

That is, in the case where $p$ is unknow (or close to 50\%), if we
want the probability that the estimation error exceeds 1\% to be less
than 1/20, we must poll approximately 10,000 people.

If we accept a lower precision and tolerate a 5\% margin of error
with a probability of 1/20, we can reduce the sample size by a factor
of 25, which leads to a sample of about 400 people.

\end{solution}

\begin{solution}{}{solexercise7.5}

1. Let $\left(X_{n}\right)_{n\in\mathbb{N}^{\ast}}$ be a sequence
of independent random variables following the same law determined
by
\[
\begin{array}{ccccc}
P\left(X_{1}=1\right)=\dfrac{18}{37} & \,\,\,\, & \text{and} & \,\,\,\, & P\left(X_{1}=-1\right)=\dfrac{19}{37}.\end{array}
\]

The random variable $X_{n}$ represents the algebraic gain of the
casino at the $n$-th game. The total algebraic gain during the first
$n$ games is
\[
G_{n}=\sum^{n}_{j=1}X_{j}.
\]
We have
\[
\mathbb{E}\left(X_{1}\right)=-\dfrac{1}{37},
\]
and since $\mathbb{E}\left(X^{2}_{1}\right)=1,$
\[
\sigma^{2}_{X_{1}}=1-\left(\dfrac{1}{37}\right)^{2}.
\]
For $n$ sufficiently large, the random variable $\dfrac{G_{n}-\dfrac{n}{37}}{\sqrt{n}\sigma_{X_{1}}}$
approximately follows a centered reduced Gauss law, and 
\[
\left(G_{n}\geqslant1,000\right)=\left(\dfrac{G_{n}-\dfrac{n}{37}}{\sqrt{n}\sigma_{X_{1}}}\geqslant\dfrac{1,000-\dfrac{n}{37}}{\sqrt{n}\sigma_{X_{1}}}\right).
\]

The number $n_{0}$ we look for is the smallest number of games that
have to be played daily such that
\[
P\left(G_{n}\geqslant1,000\right)=P\left(\dfrac{G_{n}-\dfrac{n}{37}}{\sqrt{n}\sigma_{X_{1}}}\geqslant\dfrac{1,000-\dfrac{n}{37}}{\sqrt{n}\sigma_{X_{1}}}\right)\geqslant\dfrac{1}{2}.
\]

This corresponds to finding the smallest $n$ such that
\[
1,000-\dfrac{n}{37}\leqslant0,
\]
which happens for\boxeq{
\[
n_{0}=37,000.
\]
}

2. The probability of a total loss for the casino during those $n_{0}$
games is then
\[
P\left(G_{n_{0}}<0\right)\equiv P\left(\dfrac{G_{n_{0}}-\dfrac{n_{0}}{37}}{\sqrt{n_{0}}\sigma_{X_{1}}}<\dfrac{-\dfrac{n_{0}}{37}}{\sqrt{n_{0}}\sigma_{X_{1}}}\right)
\]

By taking $\sigma_{X_{1}}\approx1,$
\[
P\left(G_{n_{0}}<0\right)\approx P\left(\dfrac{G_{n_{0}}-\dfrac{n_{0}}{37}}{\sqrt{n_{0}}\sigma_{X_{1}}}<-5.19\right).
\]

The random variable $\dfrac{G_{n_{0}}-\dfrac{n_{0}}{37}}{\sqrt{n_{0}}\sigma_{X_{1}}}$
approximately follows a centered reduced Gauss law. We deduce that\boxeq{
\[
P\left(G_{n_{0}}<0\right)\approx0.
\]
}

\end{solution}

\markright{Summary of Usual Probability Laws}
\markleft{Summary of Usual Probability Laws}

\chapter*{Summary of Usual Probability Laws}

We summarize below the different laws and the results attached to
them seen during the different chapters of this first part of the
book

\section*{Discrete Laws}

\begin{summary}{Bernoulli Law $\mathcal{B}\left(1,p\right)$ with $0<p<1$}{}  

\textit{A Bernoulli test is a random experiment with two possible
outcomes, often called ``success'' with probability $p$ and ``failure''
with probability $q=1-p.$ }

\textit{The Bernoulli law $\mathcal{B}\left(1,p\right)$ also abusively
denoted $\mathcal{B}\left(p\right)$ is the law followed by the random
variable associated to a Bernoulli test which takes 1 in case of success
and 0 otherwise.}

\begin{multicols}{2}
\begin{itemize}
\item $P\left(X=1\right)=p$ and $P\left(X=0\right)=q.$
\item $\mathbb{E}\left(X\right)=p$ and $\sigma^{2}_{X}=pq.$
\end{itemize}
\end{multicols}

\end{summary}

\begin{summary}{Binomial Law $\mathcal{B}\left(n,p\right)$ with $n\in\mathbb{N}^\ast$ and $0<p<1$}{}

\textit{The binomial law is the law of the number of success in a
sequence of $n$ independent Bernoulli tests of parameter $p.$ We
denote $q=1-p.$}

\begin{multicols}{2}
\begin{itemize}
\item $P\left(X=k\right)=\left(\begin{array}{c}
n\\
k
\end{array}\right)p^{k}q^{n-k}.$
\item $\mathbb{E}\left(X\right)=np$ and $\sigma^{2}_{X}=npq.$
\end{itemize}
\end{multicols}

\end{summary}

\begin{summary}{Hypergeometric Law with Parameters $n,r,r_1$}{}

\textit{A poll box contains $r$ balls, with $r$ tokens, with $r_{1}$
blue tokens and $r-r_{1}$ white tokens. We draw $n$ tokens without
reputting them back in the box\footnotemark. The law followed by
the random variable $X$ corresponding to the number of blue balls
obtained is the hypergeometric law of parameters $n,r,r_{1}.$ The
hypergeometric law is useful in poll modelling and quality control.}
\begin{itemize}
\item $P\left(X=k\right)=\dfrac{\left(\begin{array}{c}
r_{1}\\
k
\end{array}\right)\left(\begin{array}{c}
r-r_{1}\\
n-k
\end{array}\right)}{\left(\begin{array}{c}
r\\
n
\end{array}\right)}.$
\item $\mathbb{E}\left(X\right)=\dfrac{nr_{1}}{r}$ and $\sigma^{2}_{X}=n\dfrac{r-n}{r-1}\dfrac{r_{1}}{r}\left(1-\dfrac{r_{1}}{r}\right).$
\end{itemize}
\end{summary}

\footnotetext{In case the balls are put back within the box, the
law obtained to model the case is the binomial law $\mathcal{B}\left(n,\dfrac{r_{1}}{r}\right).$}

\begin{summary}{Poisson Law $\mathcal{P}\left(\lambda\right)$ with $\lambda>0$}{}

\textit{The Poisson law is useful in the modelling of counting phenomena.
The number of nucleus desintegration in a radioactive element during
a time interval, the number of arrivals at a client desk, the number
of drops of rain on a given surface, follow a Poisson law. The Poisson
law is obtained as a limit of binomial laws.}

\begin{multicols}{2}
\begin{itemize}
\item $P\left(X=n\right)=\dfrac{\lambda^{n}}{n!}\exp\left(-\lambda\right).$
\item $\mathbb{E}\left(X\right)=\lambda$ and $\sigma^{2}_{X}=\lambda.$
\end{itemize}
\end{multicols}

\end{summary}

\begin{summary}{Geometric Law on $\mathbb{N},$ $\mathcal{G}_\mathbb{N}\left(p\right)$ with $0<p<1$}{}

\textit{The geometric law on $\mathbb{N}$ is the law of the number
of failures met before obtaining a success in the repetition of independent
Bernoulli tests with parameter $p.$ We denote $q=1-p.$}

\begin{multicols}{2}
\begin{itemize}
\item $P\left(X=n\right)=pq^{n}.$
\item $\mathbb{E}\left(X\right)=\dfrac{q}{p}$ and $\sigma^{2}_{X}=\dfrac{q}{p^{2}}.$
\end{itemize}
\end{multicols}

\end{summary}

\begin{summary}{Geometric Law on $\mathbb{N}^\ast,$ $\mathcal{G}_\mathbb{N}^\ast\left(p\right)$ with $0<p<1$}{}

\textit{The geometric law on $\mathbb{N}^{\ast}$ is the law of the
number of Bernoulli tests before obtaining a success, including itself
in the repetition of independent Bernoulli tests with parameter $p.$
We denote $q=1-p.$}

\begin{multicols}{2}
\begin{itemize}
\item $P\left(X=n\right)=pq^{n-1}.$
\item $\mathbb{E}\left(X\right)=\dfrac{1}{p}$ and $\sigma^{2}_{X}=\dfrac{q}{p^{2}}.$
\end{itemize}
\end{multicols}

\end{summary}

\begin{summary}{Negative Binomial Law $\mathcal{B}^-\left(r,p\right)$ with $r\in\mathbb{N}^\ast$ and $0<p<1$}{}

\textit{The negative binomial law is the law of the number of failures
in a sequence of $n$ independent Bernoulli tests with parameter $p$
before obtaining $r$ success. We denote $q=1-p.$}

\begin{multicols}{2}
\begin{itemize}
\item $P\left(X=k\right)=\left(\begin{array}{c}
r+k-1\\
k
\end{array}\right)p^{r}q^{k}.$
\item $\mathbb{E}\left(X\right)=r\dfrac{q}{p}$ and $\sigma^{2}_{X}=r\dfrac{q}{p^{2}}.$
\end{itemize}
\end{multicols}

\end{summary}

\section*{Laws with Density}

\begin{summary}{Uniform Law $\mathcal{U}\left(\left[a,b\right]\right)$ on $\left[a,b\right]$ with $a<b$}{}

For $X$ following the uniform law on $\left[a,b\right]$, we have:
\begin{itemize}
\item Its density for every $x\in\mathbb{R},$ 
\[
f_{X}\left(x\right)=\dfrac{1}{b-a}\boldsymbol{1}_{\left[a,b\right]}\left(x\right)
\]
\item $\mathbb{E}\left(X\right)=\dfrac{a+b}{2}$ and $\sigma^{2}_{X}=\dfrac{\left(b-a\right)^{2}}{12}.$
\end{itemize}
\end{summary}

\begin{summary}{Gauss Law $\mathcal{N}\left(m,\sigma^2\right)$ with $m\in\mathbb{R},$,$\sigma^2\neq0$}{}

\textit{The Gauss law is also known as the normal law. }

\textit{Many measures follow approximately a normal law, the so call
error law. The central limit theorem is the main reason of the universal
intervention of the normal law in the natural phenomena and others.
The normal law is often used as an a priori model, or to approximate
known laws of same expectation and same variance.}

For $X$ following the normal law $\mathcal{N}\left(m,\sigma^{2}\right),$
with $m\in\mathbb{R}$ and $\sigma^{2}\neq0,$ we have:
\begin{itemize}
\item Its density for every $x\in\mathbb{R},$ 
\[
f_{X}\left(x\right)=\dfrac{1}{\sigma\sqrt{2\pi}}\exp\left(-\dfrac{\left(x-m\right)^{2}}{\sigma^{2}}\right)
\]
\item $\mathbb{E}\left(X\right)=m$ and $\sigma^{2}_{X}=\sigma^{2}.$
\end{itemize}
\end{summary}

\begin{summary}{Standard Cauchy Law}{}

\textit{The standard Cauchy law is often used in Probability theory
as the pathological case where the law does not have an expectation
and variance defined.}

For $X$ following the standard Cauchy law, we have:
\begin{itemize}
\item Its density for every $x\in\mathbb{R},$ $f_{X}\left(x\right)=\dfrac{1}{\pi}\dfrac{1}{1+x^{2}}.$
\item $X$ does not have neither an expectation defined nor a variance.
\end{itemize}
\end{summary}

\begin{summary}{Exponential Law $\exp\left(p\right)$ with Parameter $p>0$}{}

The exponential law intervenes in the modelling of waiting times and
in lifespans, for instance of radioactive nuclei. It is similar to
geometric laws.

For $X$ following the exponential law $\exp\left(p\right)$ with
parameter $p>0$, we have:
\begin{itemize}
\item Its density for every $x\in\mathbb{R},$ 
\[
f_{X}\left(x\right)=\boldsymbol{1}_{\mathbb{R}^{+}}\left(x\right)p\text{e}^{-px}.
\]
\item $\mathbb{E}\left(X\right)=\dfrac{1}{p}$ and $\sigma^{2}_{X}=\dfrac{1}{p^{2}}.$
\end{itemize}
\end{summary}

\begin{summary}{Chi-Squared Law $\chi^2_n$ with $n\in\mathbb{N}^\ast$}{}

\textit{The chi-squared law is the law followed by the random variable
$X^{2}_{1}+X^{2}_{2}+\dots+X^{2}_{n}$ when $X_{1},\dots,X_{n}$ are
independent random variables following the centered reduced normal
law $\mathcal{N}\left(0,1\right).$ The chi-squared law is used in
statistics in the adequation test of an empirical distribution to
a given probability law.}

For $X$ following the chi-squared law $\chi^{2}_{n}$ with $n\in\mathbb{N}$,
we have:
\begin{itemize}
\item Its density, for every $x\in\mathbb{R},$ 
\[
f_{X}\left(x\right)=\boldsymbol{1}_{\mathbb{R}^{+}}\left(x\right)\dfrac{1}{K_{n}}\exp\left(-\dfrac{x}{2}\right)x^{\frac{n}{2}-1}
\]
where, for $p\geqslant1,$
\[
\begin{array}{ccccc}
K_{2p}=2^{p}\left(p-1\right)! & \,\,\,\, & \text{and} & \,\,\,\, & K_{2p+1}=\dfrac{\left(2p-1\right)!}{2^{p-1}\left(p-1\right)!}\sqrt{2\pi}.\end{array}
\]
\item $\mathbb{E}\left(X\right)=n$ and $\sigma^{2}_{X}=2n.$
\end{itemize}
\end{summary}

\pagestyle{noheader}

\part{Deepening into Probability Theory}\label{part:Deepening-Probability-Theory}

\pagestyle{scrheadings}

\chapter{Introduction: A Measure Theory Summary}\label{chap:PartIIChap8}

\subsubsection*{Translator's Note}

\begin{leftbar}

In the French edition, this summary appeared as an Appendix. However,
we have chosen to include it here as the Introduction of Part \ref{part:Deepening-Probability-Theory},
since several concepts covered in this summary are essential for understanding
the remainder of the book. 

Readers with a solid background in measure theory may choose to skip
this chapter. The results presented are stated without proof.

\end{leftbar}

\begin{objective}{}{}

Chapter \ref{chap:PartIIChap8} aims to provide a summary of measure
theory that serves as a foundation for deeper study in probability
theory.
\begin{itemize}
\item Section \ref{sec:Measure-and-Probability} focuses on measures and
probabilities. It begins by introducing algebras and $\sigma-$algebras,
followed by the definition of $\sigma-$algebras generated by families
of subsets. The Borel $\sigma-$algebra is then introduced. Measurable
spaces are defined, followed by measurable functions---also referred
to as random variables. Basic properties of measurable functions are
presented, particularly those related to the case of Borel $\sigma-$algebra.
The concept of finite additivity is introduced, leading to the definition
of a measure, its mass, and a probability measure. Key properties
of measures are discussed, culminating in a description of how measures
can be generated.
\item Section \ref{sec:Integral} introduces integration, beginning with
nonnegative functions. It then explores the relationship between integration
and measures, before extending the concept to functions of arbitrary
sign.
\item Section \ref{sec:Three-Convergence-Theorems} presents three fundamental
convergence theorems, starting with Fatou lemma. It then covers the
Monotone Convergence Theorem and two versions of the Dominated Convergence
Theorem.
\item Section \ref{sec:Product-Measure-and} addresses product measures
and introduces the Fubini theorem. The section concludes with the
Kolomogorov Extension Theorem and the concept of product probability.
\end{itemize}
\end{objective}

We give in this chapter, the main results and statements of the measure
theory and of the integration, so we have available all the key theorems.
For a deeper study, we refer to the books of measure theory or of
probability such as those of \cite{durrett2019probability}---in
which exists a measure theory summary quite detailed. Tr.N. In the
French Edition three other references were given, all texts in French,
\cite{gramain1994,metivierneuveu}. We also point the books of \cite{tao2011introduction,ash2000probability,adams1996measure}. 

\section{Measure and Probability}\label{sec:Measure-and-Probability}

\begin{definition}{Algebra. Unitary Algebra. Semi-Algebra. $\sigma$-Algebra}{}

A family $\mathcal{A}$ of subsets of a set $\Omega$ is:
\begin{itemize}
\item An \textbf{\index{algebra}algebra} (or a \textbf{ring}\index{ring}),
if it is stable under (finite) union and set differences.
\item An \textbf{\index{unitary algebra}\mindex{algebra ! unitary}unitary
algebra} (or a \textbf{unitary ring}\index{unitary ring}), if it
is an algebra that contains the entire set $\Omega.$
\item A \textbf{semi-algebra\index{semi-algebra}\mindex{algebra ! semi}}
(or \textbf{unitary semi-ring}\index{unitary semi-ring}), if $\Omega\in\mathcal{A}$
and $\emptyset\in\mathcal{A},$ if it is stable under (finite) intersection
and if, for every $A\in\mathcal{A},$ the complement $A^{c}$ can
be written as a finite union of pairwise disjoint elements of $\mathcal{A}.$
\item A \textbf{$\sigma$-algebra\mindex{sigma-algebra @ $\sigma-$algebra}}
(or a \textbf{tribe}\index{tribe}), if it is a unitary algebra stable
under countable union. That is if for every sequence $\left(A_{n}\right)_{n\in\mathbb{N}}$
of elements of $\mathcal{A},$
\[
\bigcup_{n\in\mathbb{N}}A_{n}\in\mathcal{A}.
\]
\end{itemize}
\end{definition}

\begin{examples}{}{}
\begin{itemize}
\item The set of all finite unions of intervals of $\mathbb{R}$ is a unitary
algebra on $\mathbb{R}.$
\item The set of all rectangles in $\mathbb{R}^{n}$ of the form $\prod^{n}_{i=1}\left]a_{i},b_{i}\right]$
where $-\infty\leqslant a_{i}<b_{i}<+\infty$ forms a semi-algebra
on $\mathbb{R}^{n}.$
\end{itemize}
\end{examples}

\begin{remark}{}{}

If $\mathcal{A}$ is a $\sigma-$algebra, then
\begin{itemize}
\item $\Omega\in\mathcal{A}.$ 
\item $\mathcal{A}$ is stable under complementation;\\
If $A\in\mathcal{A},$ then $A^{c}\in\mathcal{A},$ 
\item $\mathcal{A}$ is stable under countable intersection: \\
For every sequence $\left(A_{n}\right)_{n\in\mathbb{N}}$ of elements
of $\mathcal{A},$ 
\[
\bigcap_{n\in\mathbb{N}}A_{n}\in\mathcal{A}.
\]
\end{itemize}
\end{remark}

\begin{proposition}{}{}

Let $\left\{ \mathcal{A}_{i},\,i\in I\right\} $ be a family of rings---respectively
of $\sigma-$algebra---on a set $\Omega.$ Then:
\begin{itemize}
\item Their intersection $\bigcap_{i\in I}\mathcal{A}_{i}$ is still a ring---respectively
a $\sigma-$algebra. 
\item In contrast, the union of $\sigma-$algebra $\bigcup_{i\in I}\mathcal{A}_{i}$
is not necessarily a $\sigma-$algebra.
\end{itemize}
The power set $\mathcal{P}\left(\Omega\right)$ of all subsets of
$\Omega$ is a $\sigma-$algebra---and, hence, a ring.

\end{proposition}

\begin{definition}{Generated Ring. Generated $\sigma-$algebra}{}

We now define the \textbf{ring} \textbf{generated}\index{generated ring}\mindex{ring ! generated}\mindex{generated!ring}---respectively
the \textbf{\mindex{generated!sigma-algebra @ $\sigma-$algebra}\mindex{sigma-algebra @ $\sigma-$algebra! generated}$\sigma-$algebra}
\textbf{generated}---by an arbitrary family $\mathcal{E}$ of subsets
of $\Omega$ as the intersection of all the rings---respectively
all the $\sigma-$algebra---containing $\mathcal{E}.$ The $\sigma-$algebra
generated by $\mathcal{E}$ is often denoted $\sigma\left(\mathcal{E}\right)$
and $\mathcal{E}$ is called the \textbf{generating system\index{generating system}}
of the $\sigma-$algebra $\sigma\left(\mathcal{E}\right).$

\end{definition}

\begin{examples}{}{}

The $\sigma-$algebra generated by a subset $A$ of $\Omega$ is the
family $\left\{ A,A^{c},\Omega,\emptyset\right\} .$

The $\sigma-$algebra $\left\{ \emptyset,\Omega\right\} $ is called
the \textbf{\mindex{trivial sigma-algebra @ trivial $\sigma-$algebra}\mindex{sigma-algebra @ $\sigma-$algebra ! trivial}trivial
$\sigma-$algebra}.

\end{examples}

We now give an important example of $\sigma-$algebra, in the case
where $\Omega$ is equipped as well of a topological structure. \footnote{Tr.N. We recall that a topology on a set $X$ may be defined, considering
a collection $\mathcal{C}$ of subsets of $X,$ called \textbf{open
sets\index{open sets}} and satisfying:
\begin{itemize}
\item $\emptyset\in\mathcal{C}$ and $X\in\mathcal{C}$
\item Any arbitrary (finite or infinite) union of elements of $\mathcal{C}$
is in $\mathcal{C}.$
\item Finite intersection of elements of $\mathcal{C}$ is in $\mathcal{C}.$
\end{itemize}
$\left(X,\mathcal{C}\right)$ is called a \textbf{\index{topological space}topological
space}.} 

\begin{definition}{Borel $\sigma-$algebra}{}

Let $\Omega$ be a topological space, and $\mathcal{O}$ its family
of open sets.

The \textbf{Borel $\sigma-$algebra}

\textbf{\mindex{topological space ! Borel sigma-algebra @ Borel $\sigma-$algebra}\mindex{Borel sigma-algebra @ Borel $\sigma-$algebra! topological space}\mindex{sigma-algebra @ $\sigma-$algebra ! Borel ! topological space}}
is the $\sigma-$algebra generated by all the open sets of $\Omega.$
It is denoted $\mathcal{B}_{\Omega,\mathcal{O}}$ or more simply $\mathcal{B}_{\Omega}$
if there is no ambiguity on the family of open sets used.

This is the smallest $\sigma-$algebra containing all the open subsets
of the topological space. 

\end{definition}

\begin{remark}{}{}

The Borel $\sigma-$algebra of a topological space are much ``bigger''
than only the open sets of the topological space, as the Borel $\sigma-$algebra
contains not only all the open sets, but also the closed sets, as
they are the complement of the open sets, the countable unions of
closed sets, the countable intersections of open sets, the countable
union of open sets, the countable intersection of union of closed
sets, etc...

\end{remark}

\begin{examples}{}{}
\begin{itemize}
\item The \textbf{\mindex{Borel sigma-algebra @ Borel $\sigma-$algebra! Rn @ $\mathbb{R}^n$}Borel
$\sigma-$algebra of $\mathbb{R}^{n},$} denoted $\mathcal{B}_{\mathbb{R}^{n}},$
is generated by the open subsets of $\mathbb{R}^{n}.$ \\
This $\sigma-$algebra is also generated, for instance, by the family
of rectangles of $\mathbb{R}^{n}$ of the form $\prod^{n}_{i=1}\left]a_{i},b_{i}\right]$
where $-\infty\leqslant a_{i}<b_{i}<+\infty.$
\item The Borel $\sigma-$algebra of the extended real line $\overline{\mathbb{R}},$
denoted $\mathcal{B}_{\overline{\mathbb{R}}},$ is generated by the
open subsets of $\overline{\mathbb{R}}.$ \\
This $\sigma-$algebra is also generated, for instance, by the family
of intervals of the shape $\left[a,b\right]$ where $-\infty\leqslant a<b\leqslant+\infty.$
\end{itemize}
\end{examples}

\begin{denotation}{Joins of $\sigma-$algebra}{}

Let $\left(\mathcal{A}_{i}\right)_{i\in I}$ be a family of $\sigma-$algebras
on $\Omega.$ 

The $\sigma-$algebra generated by the union of the $\mathcal{A}_{i},\,i\in I$
is denoted $\bigvee_{i\in I}\mathcal{A}_{i}.$

This is the smallest $\sigma-$algebra containing $\bigcup_{i\in I}\mathcal{A}_{i}.$

\end{denotation}

\begin{definition}{Measurable Space / Probabilizable Space}{}

If $\mathcal{A}$ is a $\sigma-$algebra, the pair $\left(\Omega,\mathcal{A}\right)$
is called a \textbf{measurable space\index{measurable space}} or
a \index{probabilizable space}\textbf{probabilizable space}.

\end{definition}

\begin{definition}{Measurable Function / Random Variable}{}

Let $f$ be a function from a set $E$ to another set $F,$ each equipped
with the $\sigma-$algebra $\mathcal{E}$ and $\mathcal{F},$ respectively. 

We say that $f$ is \textbf{measurable\index{measurable}}---or is
a \textbf{random variable\index{random variable}}---if, for every
$A\in\mathcal{F},$ the preimage 
\[
f^{-1}\left(A\right)\in\mathcal{E}.
\]

Recall that
\[
f^{-1}\left(A\right)=\left\{ x\in E\,:\,f\left(x\right)\in A\right\} .
\]

\end{definition}

\begin{proposition}{Basic Properties of Measurable Functions}{}

\textbf{(a) Composition}

The composition of two measurable functions is measurable.

\textbf{(b) Pullback of a Generated $\sigma-$Algebra}

Let $f$ be a function from $E$ to $F$ and let $\mathcal{C}$ be
a family of subsets of $F.$ 

Then, we have the equality of $\sigma-$algebra\boxeq{
\[
f^{-1}\left(\sigma\left(\mathcal{C}\right)\right)=\sigma\left(f^{-1}\left(\mathcal{C}\right)\right).
\]
}

We recall that for an arbitrary family $\mathcal{D}$ of subsets of
$F,$ $f^{-1}\left(\mathcal{D}\right)$ designates the family constituted
by $f^{-1}\left(A\right)$ for every $A$ in $\mathcal{D}.$

In particular, if $\mathcal{F}$ is a $\sigma-$algebra, then the
family $f^{-1}\left(\mathcal{F}\right)$ is also a $\sigma-$algebra:
it is called the\textbf{ $\sigma-$algebra generated by $f.$\mindex{sigma-algebra @ $\sigma-$algebra ! generated by a function}}

\textbf{(c) Criterion for Measurability via Generators}

Let $f$ be a function from $E$ to $F,$ sets respectively equipped
with the $\sigma-$algebras $\mathcal{E}$ and $\mathcal{F}.$ 

Suppose that $\mathcal{F}$ is generated by a family $\mathcal{C}$
of subsets of $F.$ Then
\[
f\text{\,is\,measurable}\,\,\,\,\Longleftrightarrow\,\,\,\,f^{-1}\left(\mathcal{C}\right)\subset\mathcal{E}.
\]

\end{proposition}

\begin{definition}{$\sigma$-Algebra Generated by a Family of Functions}{}

Let $\left(f_{i}\right)_{i\in I}$ be a family of functions, such
that for each $i\in I,$ $f_{i}$ is a function from $E$ to $F_{i},$
equipped with the $\sigma-$algebra $\mathcal{F}_{i}.$

The $\sigma-$algebra generated by the union of $\sigma-$algebra
$f^{-1}_{i}\left(\mathcal{F}_{i}\right),\,i\in I$ is called the\textbf{
$\sigma-$algebra generated by the family $\left(f_{i}\right)_{i\in I}$}\textbf{\mindex{sigma-algebra @ $\sigma-$algebra ! generated by a family of  functions}}
and it is denoted $\sigma\left(f_{i},i\in I\right).$ 

This is the smallest $\sigma-$algebra on $E$ that makes all the
functions $f_{i}$ measurable.

\end{definition}

\begin{proposition}{Properties of Measurable Function in a Borel $\sigma-$Algebra}{}

Let $\left(f_{n}\right)_{n\in\mathbb{N}}$ be a sequence of measurable
functions from a measurable space $\left(E,\mathcal{E}\right)$ to
$\mathbb{R}$---respectively $\overline{\mathbb{R}}$---equipped
with their Borel $\sigma-$algebras.

Then, whenever defined\footnotemark, the following functions $f_{1}+f_{2},$
$f_{1}f_{2},$ $f^{+}_{1},$ $f^{-}_{1},$ $\sup_{n\in\mathbb{N}}f_{n},$
$\inf_{n\in\mathbb{N}}f_{n},$ $\lim\sup_{n\in\mathbb{N}}f_{n}$ and
$\lim\inf_{n\in\mathbb{N}}f_{n}$ are measurable.

A continuous function from $\mathbb{R}^{n}$ to $\mathbb{R}^{p}$
is \textbf{\index{Borel measurable}Borel measurable,} when this function
is measurable with respect to the Borel $\sigma-$algebras on both
$\mathbb{R}^{n}$ and $\mathbb{R}^{p}.$

\end{proposition}

\footnotetext{By convention, for every $a\in\mathbb{R},$ we have
$+\infty+a=+\infty,$ $-\infty+a=-\infty,$ $+\infty+\left(+\infty\right)=+\infty,$
$-\infty+\left(-\infty\right)=-\infty.$ $0\times\left(\pm\infty\right)=0$
and for every $a\in\mathbb{R}^{\ast},$ $a\times\left(+\infty\right)=\text{sign\ensuremath{\left(a\right)}\ensuremath{\ensuremath{\infty}}},$
$a\times\left(-\infty\right)=-\text{sign\ensuremath{\left(a\right)}\ensuremath{\ensuremath{\infty}}}.$
Finally, summing $+\infty$ and $-\infty$ is not defined.}

\textbf{In what follows, it is understood that the spaces $\mathbb{R},$
$\overline{\mathbb{R}}$ and $\mathbb{R}^{n}$ are equipped with their
respective Borel $\sigma-$algebras.}

\begin{definition}{Step Function}{}

A function $f$ defined on a measurable space $\left(E,\mathcal{E}\right)$
taking values in $\mathbb{R}$---respectively $\overline{\mathbb{R}}$---is
a \textbf{step function\index{step function}} if:
\begin{itemize}
\item $f$ is measurable, and
\item $f$ takes only finitely many real finite values. 
\end{itemize}
Such a function can be written as\boxeq{
\[
f=\sum^{n}_{j=1}f_{j}\boldsymbol{1}_{A_{j}}
\]
} where:
\begin{itemize}
\item The $A_{j}$ belong to $\mathcal{E},$
\item The sets $A_{1},A_{2},\cdots,A_{n}$ are pairwise disjoint,
\item For every $j\in\left\llbracket 1,n\right\rrbracket ,$ $f_{j}\in\mathbb{R}.$ 
\end{itemize}
\end{definition}

\begin{lemma}{Approximation of Measurable Functions by Non-Decreasing Sequence of Step Functions}{}

Let $\left(E,\mathcal{E}\right)$ be a measurable space.
\begin{itemize}
\item Any measurable function defined on $E$ taking values in $\overline{\mathbb{R}}^{+}$
is the pointwise limit of a non-decreasing sequence of step functions
with values in $\mathbb{R}^{+}.$
\item Any measurable function defined on $E$ with values in $\mathbb{R}$
or $\overline{\mathbb{R}}$ is the pointwise limit of a sequence of
step functions.
\end{itemize}
\end{lemma}

\begin{definition}{Finite Additivity. $\sigma-$Additivity}{}

Let $\mathcal{F}$ be a family of subsets of a set $\Omega,$ and
let $\mu$ be a function from $\mathcal{F}$ to $\overline{\mathbb{R}}^{+}.$

We say that $\mu$ has the property of:
\begin{itemize}
\item \textbf{\index{finite additivity}Finite additivity} if for every
finite family $\left(A_{i}\right)_{i\in I}$ of pairwise disjoint
elements of $\mathcal{F}$ such that $\bigcup_{i\in I}A_{i}\in\mathcal{F},$\boxeq{
\[
\mu\left(\biguplus_{i\in I}A_{i}\right)=\sum_{i\in I}\mu\left(A_{i}\right).
\]
}
\item \textbf{Countable additivity}\index{countable additivity}, or that
$\mu$ is \textbf{\mindex{sigma-additive @ $\sigma-$additive}$\sigma-$additive}
if for every countable family $\left(A_{i}\right)_{i\in I}$ of pairwise
disjoint elements of $\mathcal{F}$ such that $\bigcup_{i\in I}A_{i}\in\mathcal{F},$\boxeq{
\[
\mu\left(\biguplus_{i\in I}A_{i}\right)=\sum_{i\in I}\mu\left(A_{i}\right).
\]
}
\end{itemize}
\end{definition}

\begin{examples}{}{}

1. The function $\mu$ defined by
\begin{align*}
\mu\left(A\right)=\begin{cases}
0, & \text{if }\left|A\right|<+\infty,\\
+\infty, & \text{otherwise,}
\end{cases}
\end{align*}
is finitely additive but not $\sigma-$additive.

2. The function $\mu$ defined on the family $\mathcal{I}$ of intervals
of $\mathbb{R}$ by $\mu\left(A\right)=\text{length}\left(A\right)$
for every $A\in\mathcal{I}$ is $\sigma-$additive.

\end{examples}

\begin{definition}{Measure. Mass. Probability. Dirac Measure. Discrete Measure}{}

Let $\left(\Omega,\mathcal{A}\right)$ be a measurable space.
\begin{itemize}
\item A \textbf{measure\index{measure}} $\mu$ on $\left(\Omega,\mathcal{A}\right)$
is a $\sigma-$additive function from $\mathcal{A}$ to $\overline{\mathbb{R}}^{+}$
such that $\mu\left(\emptyset\right)=0.$
\item A measure $\mu$ is \textbf{finite\index{finite measure}\mindex{measure ! finite}}
if it has its values in $\mathbb{R}^{+}.$
\item The \textbf{mass of a measure\index{mass of a measure}\mindex{measure ! mass}}
is the value of $\mu\left(\Omega\right).$
\item A measure is \textbf{$\sigma-$finite\mindex{sigma-finite@$\sigma-$finite}\mindex{measure ! sigma-finite@$\sigma-$finite}}
if there exists a countable covering of $\Omega$ by a family $\left(A_{n}\right)_{n\in\mathbb{N}}$
of finite measure elements of $\mathcal{A}.$
\item A \textbf{probability\index{probability}} $P$ on $\left(\Omega,\mathcal{A}\right)$
is a measure of mass 1.
\item The \textbf{Dirac\sindex[fam]{Dirac, Paul} measure\index{Dirac measure}\mindex{measure ! Dirac}}
in $\omega\in\Omega$ is a measure denoted $\delta_{\omega}$ defined
by $\delta_{\omega}\left(A\right)=1,$ if $\omega\in A$ and, $\delta_{\omega}\left(A\right)=0,$
otherwise.
\item A measure is a \textbf{\mindex{discrete!measure}discrete measure\mindex{measure ! discrete}},
if it is of the shape $\mu=\sum_{\omega\in D}\alpha_{\omega}\delta_{\omega}$
where $D$ is a countable subset of $\Omega$ and $\alpha_{\omega}\in\mathbb{R}^{+}.$
\end{itemize}
\end{definition}

\begin{proposition}{Measure Properties}{}

Let $\mu$ be a measure on a measurable space $\left(\Omega,\mathcal{A}\right)$
non identical to $+\infty.$ 

We have the following properties:

\textbf{(a) Additivity for Disjoint Sets}

For every $A,B\in\mathcal{A}$ with $A\cap B=\emptyset,$\boxeq{
\[
\mu\left(A\cup B\right)=\mu\left(A\right)+\mu\left(B\right).
\]
}

\textbf{(b) Monotonicity}

For every $A,B\in\mathcal{A}$ with $A\subset B,$\boxeq{
\[
\mu\left(A\right)\leqslant\mu\left(B\right).
\]
}

\textbf{(c) Sub-$\sigma-$additivity}

For every $A,B\in\mathcal{A},$\boxeq{
\[
\mu\left(A\cup B\right)\leqslant\mu\left(A\right)+\mu\left(B\right).
\]
}

\textbf{(d) Continuity from Below}

If $\left(A_{n}\right)_{n\in\mathbb{N}}$ is a non-decreasing sequence
of elements of $\mathcal{A}.$ Then\boxeq{
\[
\mu\left(\bigcup_{n\in\mathbb{N}}A_{n}\right)=\lim_{n\to+\infty}\mu\left(A_{n}\right).
\]
}

\textbf{(e) Continuity from Above}

If $\left(A_{n}\right)_{n\in\mathbb{N}}$ is a nonincreasing sequence
of elements of $\mathcal{A},$ such that there exists $n_{0}$ for
which $\mu\left(A_{n_{0}}\right)<+\infty.$ Then\boxeq{
\[
\mu\left(\bigcap_{n\in\mathbb{N}}A_{n}\right)=\lim_{n\to+\infty}\mu\left(A_{n}\right).
\]
}

\end{proposition}

\subsection*{Measure generation}

\begin{figure}[t]
\begin{center}\includegraphics[width=0.4\textwidth]{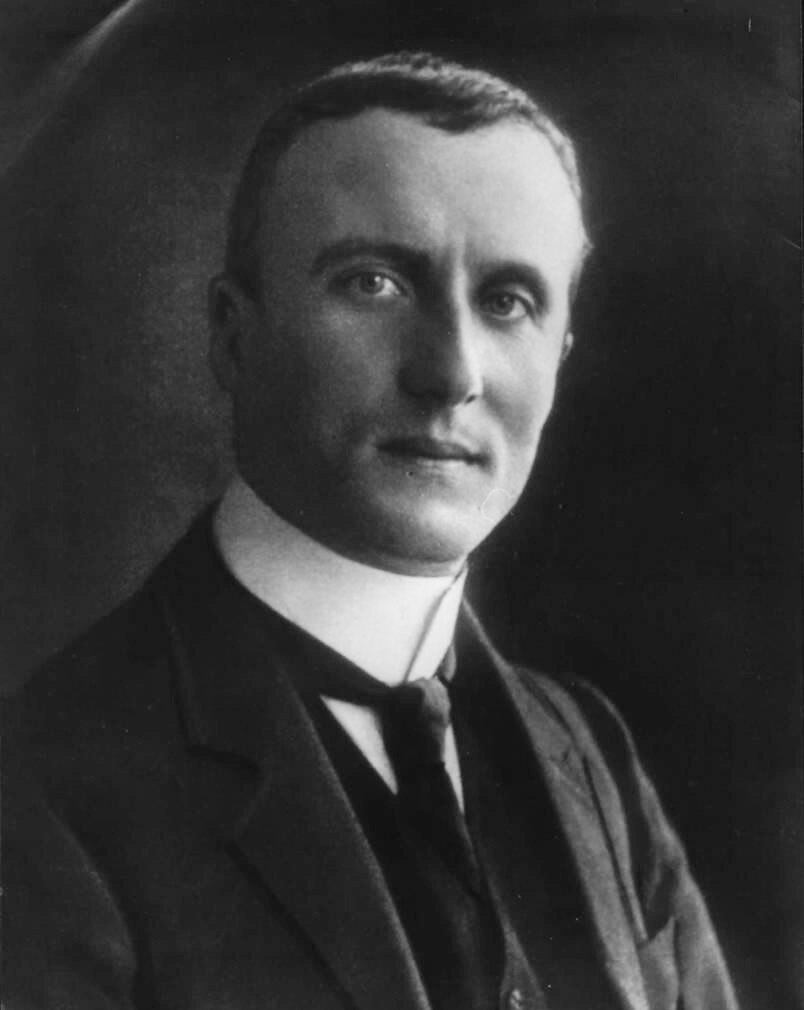}

{\tiny Public Domain}\end{center}

\caption{\textbf{\protect\href{https://en.wikipedia.org/wiki/Constantin_Carath\%25C3\%25A9odory}{Constantin Caratheodory}}
(1873-1950)}\sindex[fam]{Caratheodory, Constantin}
\end{figure}

We now state the Caratheodory\footnote{\textbf{\sindex[fam]{Caratheodory, Constantin}\href{https://en.wikipedia.org/wiki/Constantin_Carath\%25C3\%25A9odory}{Constantin Caratheodory}}
(1873-1950) is a greek mathematician with important works in theory
of functions with real-valued variables, in variation computation
and measure theory.} extension theorem.

\begin{theorem}{Caratheodory Extension Theorem}{}

Let $\mu$ be a $\sigma-$additive function on an unitary algebra
$\mathcal{A}$ and satisfying 
\[
\mu\left(\emptyset\right)=0.
\]
Then $\mu$ admits a unique extension to a measure $\mu$ on the $\sigma-$algebra
generated by $\mathcal{A}.$

\end{theorem}

\begin{theorem}{Extension from a Semi-Algebra}{}

Let $\mathcal{S}$ be a semi-algebra on $\Omega.$ The algebra $\overline{\mathcal{S}}$
generated by $\mathcal{S}$ is the family of finite unions of disjoint
elements of $\mathcal{S}.$ 

Let $\mu$ be a function, additive on the semi-algebra $\mathcal{S},$
such that $\mu\left(\emptyset\right)=0,$ and sub$-\sigma-$additive
on $\mathcal{S},$ that is such that, for every countable family $A_{i},$\boxeq{
\[
\mu\left(\biguplus_{i\in I}A_{i}\right)\leqslant\sum_{i\in I}\mu\left(A_{i}\right).
\]
}

Then $\mu$ can be uniquely extended to a $\sigma-$additive function
on $\overline{\mathcal{S}},$ and consequently, to a unique measure
$\mu$ on the $\sigma-$algebra generated by $\mathcal{S}.$

\end{theorem}

\begin{figure}[t]
\begin{center}\includegraphics[width=0.4\textwidth]{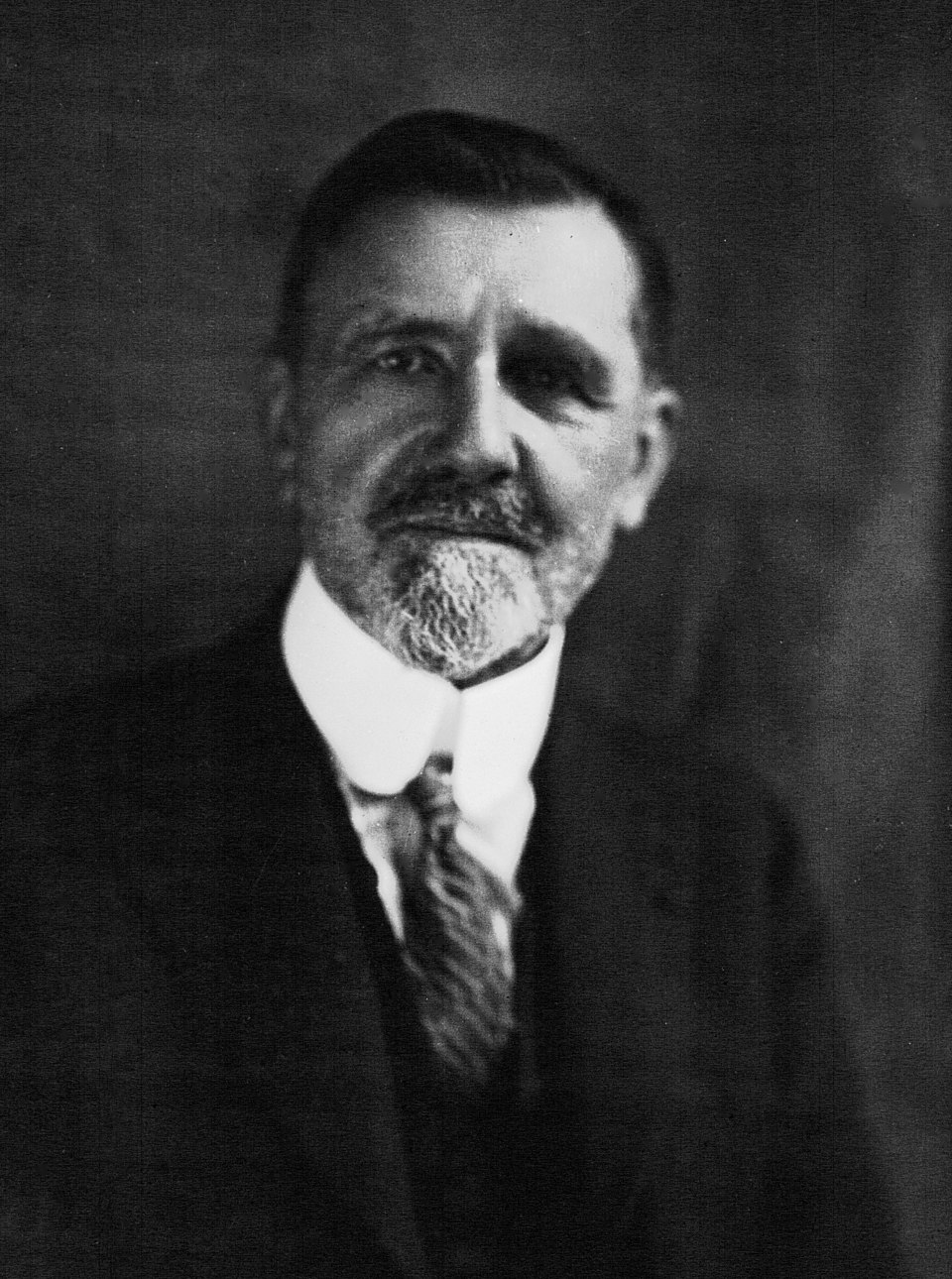}

{\tiny Public domain}\end{center}

\caption{\textbf{\protect\href{https://en.wikipedia.org/wiki/\%25C3\%2589mile_Borel}{Emile Borel}}
(1856 - 1994)}\sindex[fam]{Borel, Emile}
\end{figure}

\begin{figure}[t]
\begin{center}\includegraphics[width=0.4\textwidth]{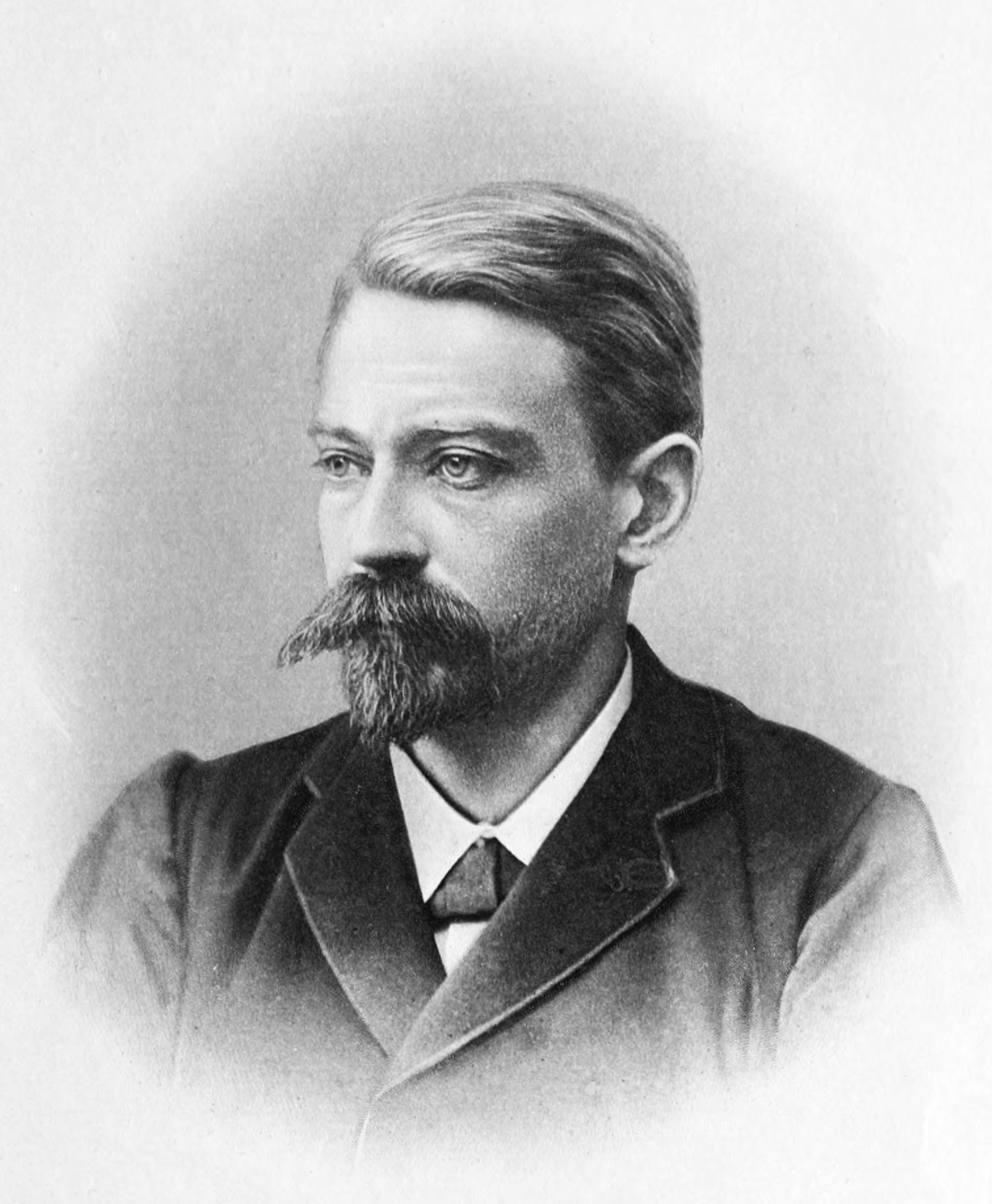}

{\tiny Public domain}\end{center}

\caption{\textbf{\protect\href{https://en.wikipedia.org/wiki/Thomas_Joannes_Stieltjes}{Thomas Joannes Stieltjes}}
(1856 - 1994)}\sindex[fam]{Stieltjes, Thomas Joannes}
\end{figure}

\begin{figure}[t]
\begin{center}\includegraphics[width=0.4\textwidth]{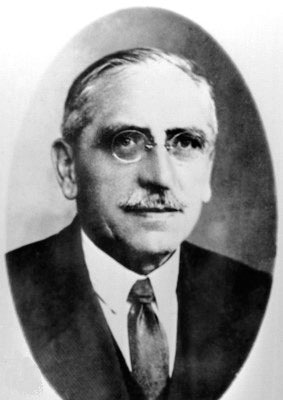}

{\tiny Public domain}\end{center}

\caption{\textbf{\protect\href{https://en.wikipedia.org/wiki/Henri_Lebesgue}{Henri-Léon Lebesgue}}
(1875 - 1941)}\sindex[fam]{Lebesgue, Henri-Léon}
\end{figure}

\begin{example}{}{}

Let $\mathcal{I}_{d}$ be the semi-algebra on $\mathbb{R}$ consisting
of intervals of the form $\left]a,b\right],$ and define $\mu$ the
function length defined on $\mathcal{I}_{d}.$ $\mu$ is \textbf{$\sigma-$additive\mindex{sigma-additive @ $\sigma-$additive}}
on $\mathcal{I}_{d}.$

By the Caratheodory extension theorem, $\mu$ admits a unique extension
to a measure on the Borel $\sigma-$algebra $\mathcal{B}_{\mathbb{R}}.$
This extension measure is called the \textbf{Borel}\footnotemark\textbf{
measure on $\mathbb{R}.$\index{Borel measure on ensuremath{mathbb{R}}@Borel measure on $\mathbb{R}$}}

More generally, let $F:\mathbb{R}\to\mathbb{R}$ be a right-continuous,
non-decreasing function. Then, there exists a unique measure $\mu$
on $\left(\mathbb{R},\mathcal{B}_{\mathbb{R}}\right)$ such that
\[
\forall a,b\in\mathbb{R},\,\,\,\,\mu\left(\left]a,b\right]\right)=F\left(b\right)-F\left(a\right).
\]
This measure is called the \textbf{Borel-Stieltjes\footnotemark measure\index{Borel-Stieltjes measure}\mindex{measure ! Borel-Stieljes}}
associated with $F,$ and it is $\sigma-$finite. 

When $F\left(x\right)=x,$ this measure is called the \textbf{\mindex{Lebesgue!measure}\mindex{measure ! Lebesgue}Lebesgue\footnotemark
measure} on $\mathbb{R}.$

Similarly, if $\mathcal{P}$ is the semi-algebra of rectangles in
$\mathbb{R}^{n}$ of the form 
\[
\prod^{n}_{i=1}\left]a_{i},b_{i}\right].
\]

Define the function volume $\mu$ on $\mathcal{P}$ by
\[
\mu\left(\prod^{n}_{i=1}\left]a_{i},b_{i}\right]\right)=\prod^{n}_{i=1}\left(b_{i}-a_{i}\right).
\]
Then $\mu$ is \textbf{\mindex{sigma-additive @ $\sigma-$additive}$\sigma-$additive}
on $\mathcal{P}.$ It extends uniquely to a measure on the Borel $\sigma-$algebra
and is called the \textbf{Borel measure\index{Borel measure on ensuremath{mathbb{R}}^{n}@Borel measure on $\mathbb{R}^{n}$}\mindex{measure ! Borel}
on $\mathbb{R}^{n}.$}

\end{example}

\addtocounter{footnote}{-2}

\footnotetext{\textbf{\href{https://en.wikipedia.org/wiki/\%25C3\%2589mile_Borel}{Emile Borel}}
(1871-1956) was a French mathematician, borned at Saint-Affrique.
He taught at École Normale Supérieure---a prestigious Research and
Teacher School in France---and then at the Sorbonne---also a prestigious
University in Paris. His research work are first focused on measure
theory---he introduced the concept of set of measure zero---on real
variable functions and on series summation. He then turn on working
on probability theory, on game theory and on mathematical physics.
In particular, he gave a probabilistic approach of gaze kinetic.}

\stepcounter{footnote}

\footnotetext{\textbf{\href{https://en.wikipedia.org/wiki/Thomas_Joannes_Stieltjes}{Thomas Joannes Stieltjes}}
(1856--1894) was a Dutch mathematician known for his work in analysis
and continued fractions. He is best known for the Stieltjes integral,
a generalization of the Riemann integral that laid the groundwork
for modern measure and probability theory. Stieltjes made significant
contributions to orthogonal polynomials and moment problems. Though
largely self-taught, he held a professorship at the University of
Toulouse. His work continues to influence real analysis, spectral
theory, and functional analysis.}

\stepcounter{footnote}

\footnotetext{\textbf{\href{https://en.wikipedia.org/wiki/Henri_Lebesgue}{Henri-Léon Lebesgue}}
(1875--1941) was a French mathematician who revolutionized integration
theory by introducing the Lebesgue integral, which extended the class
of functions that can be integrated and laid the groundwork for modern
probability theory.}

\begin{definition}{$\mu-$negligible Set. Complete Measured Space}{}

Let $\left(\Omega,\mathcal{A},\mu\right)$ be a measured space.
\begin{itemize}
\item A set $A$ is said to be\textbf{ \mindex{mu-negligible@$\mu-$negligible}$\mu-$negligible}
if there exists a set $B\in\mathcal{A}$ with $A\subset B$ and $\mu\left(B\right)=0.$
\item The measured space $\left(\Omega,\mathcal{A},\mu\right)$ is said
to be \textbf{\index{complete measured space}complete}, if every
$\mu-$neglectible set belongs to $\mathcal{A}.$ This way, any subset
of sets of measure zero is measurable.
\end{itemize}
\end{definition}

\begin{proposition}{Extension to a Complete Measure}{}

Let $\left(\Omega,\mathcal{A},\mu\right)$ be a measured space. 

Define the family of subsets of $\Omega,$
\[
\mathcal{A}^{\mu}=\left\{ X:\,\exists B_{1},B_{2}\in\mathcal{A}\text{ such that }B_{1}\subset X\subset B_{2}\text{ and }\mu\left(B_{2}-B_{1}\right)=0\right\} .
\]

Then
\begin{itemize}
\item $\mathcal{A}^{\mu}$ is a $\sigma-$algebra
\item There exists a unique extension of $\mu$ to a measure $\widehat{\mu}$
on $\mathcal{A}^{\mu},$ 
\item $\left(\Omega,\mathcal{A},\widehat{\mu}\right)$ is a complete measured
space.
\end{itemize}
$\widehat{\mu}$ is called the completed measure related to $\mu.$

\end{proposition}

\begin{example}{Lebesgue Measure}{}The completion of the Borel measure
on $\mathbb{R}^{n}$ is the \textbf{Lebesgue\sindex[fam]{Lebesgue, Henri-Léon}
measure\mindex{Lebesgue ! measure}} on $\mathbb{R}^{n}.$

\end{example}

\begin{definition}{$\mu-$Almost Everywhere True 
Property}{}

Let $\left(\Omega,\mathcal{A},\mu\right)$ be a measured space. 

A property $\mathcal{P}$ depending on $\omega\in\Omega$ is said
to be \textbf{true $\mu-$almost everywhere\mindex{mu-almost everywhere@$\mu-$almost everywhere}}---$\mu-$a.e.
for short---if the set $\omega$ where the property $\mathcal{P}\left(\omega\right)$
is false is $\mu-$negligible, i.e. contained in a measurable set
of measure zero.\footnotemark

\end{definition}

\footnotetext{We sometimes write for short $\mu-\text{a.e.}$ instead
of $\mu-$almost everywhere.}

\begin{example}{}{}Saying that a sequence of measurable functions
$\left(f_{n}\right)_{n\in\mathbb{N}^{\ast}}$ converges $\mu-$almost
everywhere is equivalent to say that the set of $\omega$ where the
sequence $\left(f_{n}\left(\omega\right)\right)_{n\in\mathbb{N}^{\ast}}$
does not converge is of measure zero.

\end{example}

\section{Integral}\label{sec:Integral}

\subsection{Integration of Nonnegative Functions}

Let $\left(\Omega,\mathcal{A}\right)$ be a measurable space, and
let $\mathcal{M}^{+}$ denote the set of functions with values in
$\overline{\mathbb{R}}^{+}$ and $\mathcal{A}-$measurable.

\begin{definition}{Integral on a Measurable Space}{}

An \textbf{integral on a measurable space $\left(\Omega,\mathcal{A}\right)$\index{integral on a measurable space}\mindex{measurable space!integral}}
is any function $\mathscr{I}$ from $\mathcal{M}^{+}$ to $\overline{\mathbb{R}}^{+}$
satisfying the following properties:
\begin{itemize}
\item $\mathscr{I}\left(0\right)=0.$ 
\item $\sigma-$additivity on $\mathcal{M}^{+}.$ That is, for every sequence
$\left(f_{n}\right)_{n\in\mathbb{N}}$ of elements of $\mathcal{M}^{+},$\boxeq{
\[
\mathscr{I}\left(\sum_{n\in\mathbb{N}}f_{n}\right)=\sum_{n\in\mathbb{N}}\mathscr{I}\left(f_{n}\right)
\]
}
\end{itemize}
\end{definition}

\begin{figure}[t]
\begin{center}\includegraphics[width=0.4\textwidth]{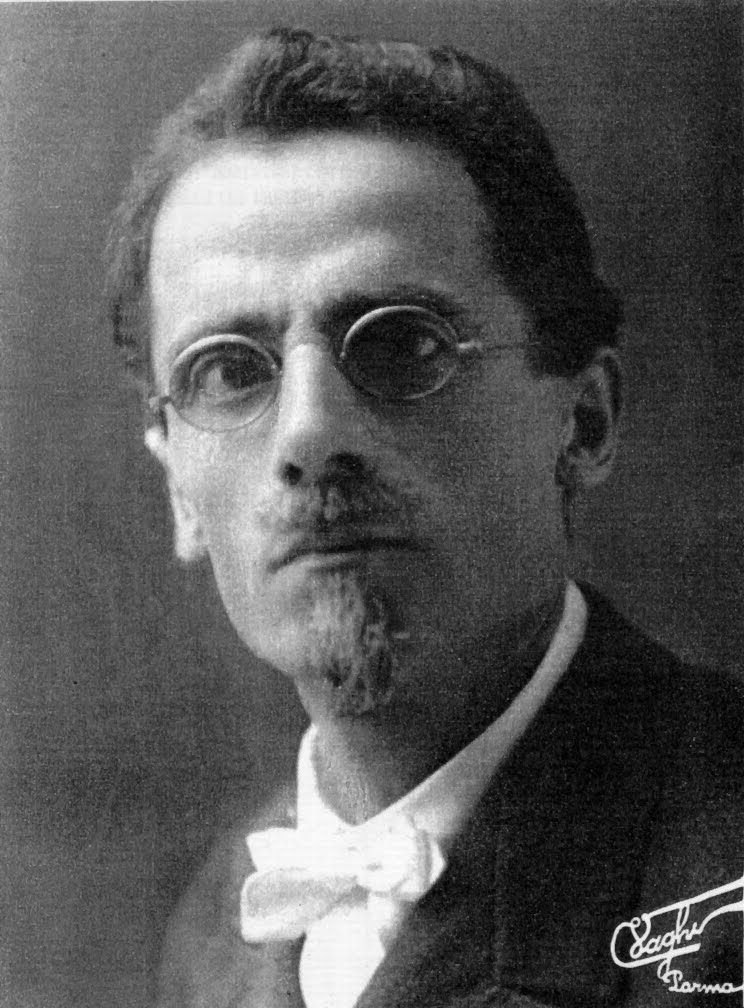}

{\tiny Public domain}\end{center}

\caption{\textbf{\protect\href{https://en.wikipedia.org/wiki/Beppo_Levi}{Beppo Levi}}
(1875 - 1961)}\sindex[fam]{Levi, Beppo}
\end{figure}

\begin{proposition}{Properties of the Integral}{integral_properties}

Let $\mathscr{I}$ be an integral on $\left(\Omega,\mathcal{A}\right)$,
and let $f,\,g$ and, $f_{n},\,n\in\mathbb{N}$ be elements of $\mathcal{M}^{+}.$ 

The following properties hold in $\overline{\mathbb{R}}^{+}:$

\textbf{(a) Additivity of the Integral}\boxeq{
\[
\mathscr{I}\left(f+g\right)=\mathscr{I}\left(f\right)+\mathscr{I}\left(g\right).
\]
}

\textbf{(b) Monotonicity (Order-Preserving)} \boxeq{
\[
f\leqslant g\Rightarrow\mathscr{I}\left(f\right)\leqslant\mathscr{I}\left(g\right).
\]
}

(c) \textbf{Monotone Convergence Theorem (Beppo Levi\footnotemark
Property)} \boxeq{
\[
\text{If }\lim_{n\to+\infty}\nearrow f_{n}=f,\text{ then }\lim_{n\to+\infty}\nearrow\mathscr{I}\left(f_{n}\right)=\mathscr{I}\left(f\right).
\]
}

\textbf{(d) Continuity from Above (under Finite Bound)}

If $\lim_{n\to+\infty}\searrow f_{n}=f$ pointwise and there exists
$n_{0}$ such that $\mathscr{I}\left(f_{n_{0}}\right)<+\infty$ then
\[
\lim_{n\to+\infty}\searrow\mathscr{I}\left(f_{n}\right)=\mathscr{I}\left(f\right).
\]

\textbf{(e) Linearity with scalars}\boxeq{
\[
\forall a\in\mathbb{R}^{+},\,\,\,\,\mathscr{I}\left(af\right)=a\mathscr{I}\left(f\right).
\]
}

\end{proposition}

\footnotetext{Beppo Levi (1875-1961), was an Italian mathematician,
borned in Turin, set up in Argentina in 1939, where he fled the fascist
regim in Italy. His work is related to the theory of integration and
also to early quantum mechanics.}

\subsection{Link Between Integral and Measure}

\begin{theorem}{Link Between Integral and Measure}{}

Let $\mathscr{I}$ be an integral on $\left(\Omega,\mathcal{A}\right).$ 

The function $A\mapsto\mathscr{I}\left(\boldsymbol{1}_{A}\right)$
is a measure on $\left(\Omega,\mathcal{A}\right).$

Conversely, let $\mu$ be a measure on $\left(\Omega,\mathcal{A}\right).$
There exists a unique integral $\mathscr{I}_{\mu}$ on $\left(\Omega,\mathcal{A}\right)$
such that for every $A\in\mathcal{A},$ $\mathscr{I}_{\mu}\left(\boldsymbol{1}_{A}\right)=\mu\left(A\right).$ 

Moreover, if $f\in\mathcal{M}^{+},$ then $\mathscr{I}_{\mu}\left(f\right)$
is given by
\[
\mathscr{I}_{\mu}\left(f\right)=\begin{cases}
\sum_{x\in f\left(\Omega\right)}x\mu\left(f=x\right), & \text{if \ensuremath{f} is a step function,}\\
\sup\left\{ \mathscr{I}_{\mu}\left(g\right):\,g\leqslant f,\,g\text{\,step function}\right\} , & \text{in the general case.}
\end{cases}
\]

\end{theorem}

\begin{denotation}{Integral of a Function in $\mathcal{M}^+$ with Respect to a Measure}{}

The integral of the function $f\in\mathcal{M}^{+}$ with respect to
a measure $\mu$ is denoted $\mathscr{I}_{\mu}\left(f\right),$ element
of $\overline{\mathbb{R}}^{+},$ is equivalently denoted by\boxeq{
\[
\intop_{\Omega}f\text{d}\mu,\,\,\,\,\intop_{\Omega}f\left(\omega\right)\text{d}\mu\left(\omega\right)\,\,\,\,\text{or}\,\,\,\,\intop_{\Omega}f\left(\omega\right)\mu\left(\text{d}\omega\right),
\]
}and is called the integral of $f$ with respect to $\mu.$

\end{denotation}

\begin{lemma}{Unicity Lemma}{unicity_lemma}

Two integrals $\mathscr{I}$ and $\mathscr{I}^{\prime}$ on $\left(\Omega,\mathcal{A}\right)$
such that, for every $A\in\mathcal{A},$ 
\[
\mathscr{I}\left(\boldsymbol{1}_{A}\right)=\mathscr{I}^{\prime}\left(\boldsymbol{1}_{A}\right)
\]
are equal.

\end{lemma}

\subsection{Integration of Functions of Arbitrary Sign}

\begin{definition}{$\mu$-integrability. Integral of a Function with Respect to a Measure}{}

Let $\left(\Omega,\mathcal{A},\mu\right)$ be a measured space. 

A function taking values in $\overline{\mathbb{R}}$ is said to be
\textbf{\mindex{mu-integrable@$\mu-$integrable}$\mu-$integrable}
if it is $\mathcal{A}-$measurable and satisfies \boxeq{
\[
\intop_{\Omega}\left|f\right|\text{d}\mu<+\infty,
\]
} or equivalently, \boxeq{
\[
\intop_{\Omega}f^{+}\text{d}\mu<+\infty\,\,\,\,\text{and}\,\,\,\,\intop_{\Omega}f^{-}\text{d}\mu<+\infty.
\]
}

If $f$ is $\mu-$integrable, the element of $\mathbb{R},$ 
\[
\intop_{\Omega}f\text{d}\mu=\intop_{\Omega}f^{+}\text{d}\mu-\intop_{\Omega}f^{-}\text{d}\mu
\]
is indifferently denoted by \boxeq{
\[
\intop_{\Omega}f\text{d}\mu,\,\,\,\,\intop_{\Omega}f\left(\omega\right)\text{d}\mu\left(\omega\right)\,\,\,\,\text{or}\,\,\,\,\intop_{\Omega}f\left(\omega\right)\mu\left(\text{d}\omega\right),
\]
}and is called the \textbf{integral of $f$ with respect to $\mu.$\index{integral of a function with respect to a measure}}

\end{definition}

\begin{definition}{$\mu$-semi-integrability. Integral of a Function with Respect to a Measure}{}

Let $\left(\Omega,\mathcal{A},\mu\right)$ be a measured space. 

A function with values in $\overline{\mathbb{R}}$ is said to be \textbf{\mindex{mu-semi-integrable@$\mu-$semi-integrable}$\mu-$semi-integrable}
if it is $\mathcal{A}-$measurable and satisfies 
\[
\intop_{\Omega}f^{+}\text{d}\mu<+\infty\,\,\,\,\text{or}\,\,\,\,\intop_{\Omega}f^{-}\text{d}\mu<+\infty.
\]

If $f$ is $\mu-$semi-integrable, the element of $\overline{\mathbb{R}},$
\[
\intop_{\Omega}f\text{d}\mu=\intop_{\Omega}f^{+}\text{d}\mu-\intop_{\Omega}f^{-}\text{d}\mu
\]
is indifferently denoted by \boxeq{
\[
\intop_{\Omega}f\text{d}\mu,\,\,\,\,\intop_{\Omega}f\left(\omega\right)\text{d}\mu\left(\omega\right)\,\,\,\,\text{or}\,\,\,\,\intop_{\Omega}f\left(\omega\right)\mu\left(\text{d}\omega\right),
\]
}and is called the \textbf{integral of $f$ with respect to $\mu.$\index{integral of a function with respect to a measure}}

\end{definition}

\begin{proposition}{}{integrable_vector_space}

(a) If $f$ and $g$ are semi-integrable and $f\leqslant g,$ then
\[
\intop_{\Omega}f\text{d}\mu\leqslant\intop_{\Omega}g\text{d}\mu.
\]

(b) If $f$ is measurable, if $g$ is $\mu-$integrable and if $\left|f\right|\leqslant g,$
then $f$ is $\mu-$integrable.

(c) If $f$ is $\mu-$integrable, then 
\[
\left|\intop_{\Omega}f\text{d}\mu\right|\leqslant\intop_{\Omega}\left|f\right|\text{d}\mu.
\]

(d) The set $\mathscr{L}^{1}\left(\Omega,\mathcal{A},\mu\right)$
of functions with values in $\mathbb{R}$ that are $\mu-$integrable
is a vector space and the mapping 
\[
f\mapsto\intop_{\Omega}f\text{d}\mu
\]
 is linear from $\mathscr{L}^{1}\left(\Omega,\mathcal{A},\mu\right)$
to $\mathbb{R}.$

\end{proposition}

\begin{figure}[t]
\begin{center}\includegraphics[width=0.4\textwidth]{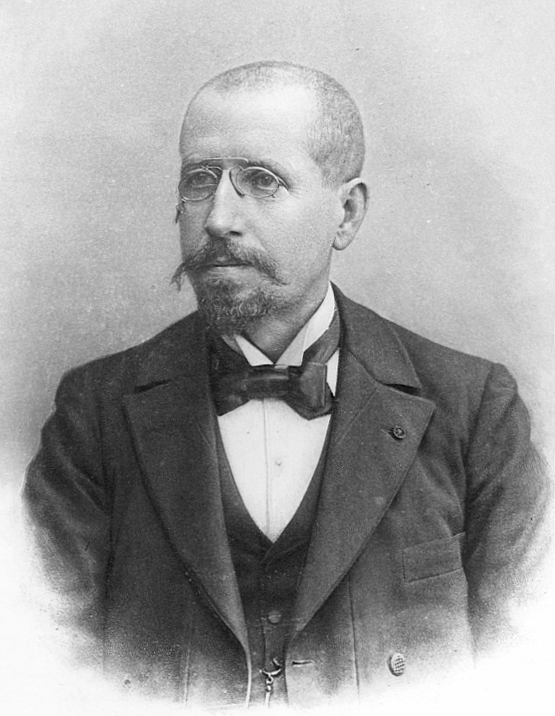}

{\tiny Public domain}\end{center}

\caption{\textbf{\protect\href{https://en.wikipedia.org/wiki/Jean_Gaston_Darboux}{Jean Gaston Darboux}}
(1887 - 1956)}\sindex[fam]{Darboux, Gaston}
\end{figure}

\begin{examples}{}{}
\begin{itemize}
\item Let $\delta_{\omega_{0}}$ be the Dirac measure at $\omega_{0}\in\Omega.$
The function $\mathscr{I}$ from $\mathcal{M}^{+}$ to $\overline{\mathbb{R}}^{+},$
defined by $f\mapsto f\left(\omega_{0}\right)$ is an integral. Since
for every $A\in\mathcal{A},$
\[
\mathscr{I}\left(\boldsymbol{1}_{A}\right)=\boldsymbol{1}_{A}\left(\omega_{0}\right)=\delta_{\omega_{0}}\left(A\right),
\]
and by using Lemma $\ref{lm:unicity_lemma}$, we have for every $f\in\mathcal{M}^{+},$
\[
\intop_{\Omega}f\text{d}\delta_{\omega_{0}}=f\left(\omega_{0}\right).
\]
Moreover, if $f$ is $\mathcal{A}-$measurable and of arbitrary sign,
then $f$ is $\mu-$integrable---respectively $\mu-$semi-integrable---if
and only if $f^{+}\left(\omega_{0}\right)<+\infty$ and $f^{-}\left(\omega_{0}\right)<+\infty$---respectively
$f^{+}\left(\omega_{0}\right)<+\infty$ or $f^{-}\left(\omega_{0}\right)<+\infty.$
In this case, we still have
\[
\intop_{\Omega}f\text{d}\delta_{\omega_{0}}=f\left(\omega_{0}\right).
\]
\item A similar argument shows that if $\mu$ is a discrete measure 
\[
\mu=\sum^{+\infty}_{n=1}\alpha_{n}\delta_{\omega_{n}}
\]
 where $\alpha_{n}\in\mathbb{R}^{+}$ and $\omega_{n}\in\Omega,$
then for every $f\in\mathcal{M}^{+},$
\[
\intop_{\Omega}f\text{d}\mu=\sum^{+\infty}_{n=1}\alpha_{n}f\left(\omega_{n}\right).
\]
If $f$ is $\mathcal{A}-$measurable and of arbitrary sign, then $f$
is $\mu-$integrable if and only if 
\[
\sum^{+\infty}_{n=1}\alpha_{n}\left|f\left(\omega_{n}\right)\right|<+\infty.
\]
 Again, in this case,
\[
\intop_{\Omega}f\text{d}\mu=\sum^{+\infty}_{n=1}\alpha_{n}f\left(\omega_{n}\right).
\]
\item \textbf{Lebesgue integral of a Riemann-integrable function on $\left[a,b\right].$}\\
Let $\mathcal{P}$ be a finite partition of $\left[a,b\right]$ into
intervals, and let $f$ be a bounded function of arbitrary sign defined
on $\left[a,b\right].$ \\
For each interval $P\in\mathcal{P},$ define
\[
\begin{array}{ccccc}
\underline{f_{P}}=\inf\left\{ f\left(x\right):\,x\in P\right\}  & \,\,\,\, & \text{and} & \,\,\,\, & \overline{f_{P}}=\sup\left\{ f\left(x\right):\,x\in P\right\} ,\end{array}
\]
and let $\left|P\right|$ denote the length of $P.$ \\
The \index{Darboux sums}\textbf{Darboux}\footnotemark\sindex[fam]{Darboux, Gaston}\textbf{
sums} are then defined as:
\[
\begin{array}{ccccc}
s_{P}=\sum\limits_{P\in\mathcal{P}}\underline{f_{P}}\left|P\right| & \,\,\,\, & \text{and} & \,\,\,\, & S_{P}=\sum\limits_{P\in\mathcal{P}}\overline{f_{P}}\left|P\right|.\end{array}
\]
By definition, the function $f$ is Riemann-integrable on $\left[a,b\right]$
if, for every sequence $\left(P_{n}\right)_{n\in\mathbb{N}}$ of nested
partitions with mesh tending to zero, the sequences $\left(s_{P_{n}}\right)_{n\in\mathbb{N}}$
and $\left(S_{P_{n}}\right)_{n\in\mathbb{N}}$ converge to the same
limit. This common limit is 
\[
\intop^{b}_{a}f\left(x\right)\text{d}x
\]
the Riemann integral of $f$ on $\left[a,b\right].$\\
If $f$ is Riemann-integrable and nonnegative on $\left[a,b\right],$
then the function $f\boldsymbol{1}_{\left[a,b\right]}$ is Legesgue-integrable---i.e.,
integrable with respect to the Lebesgue measure $\lambda$ on $\mathbb{R}$---and\boxeq{
\[
\intop_{\mathbb{R}}f\boldsymbol{1}_{\left[a,b\right]}\text{d}\lambda=\intop^{b}_{a}f\left(x\right)\text{d}x.
\]
}
\end{itemize}
\end{examples}

\footnotetext{\textbf{\sindex[fam]{Darboux, Gaston}\href{https://en.wikipedia.org/wiki/Jean_Gaston_Darboux}{Jean Gaston Darboux}}
(1842 - 1917) was a French mathematician known for his work in geometry
and analysis. He made significant contributions to differential geometry,
particularly through the concept of the Darboux frame. In real analysis,
he is recognized for Darboux sums, which form a foundational part
of the Riemann integral theory. He also contributed to the theory
of orthogonal systems and partial differential equations. Darboux
held professorships at prestigious institutions including the Collège
de France and the Sorbonne. Among his students were Emile Borel and
Elie Cartan.}

\section{Three Convergence Theorems}\label{sec:Three-Convergence-Theorems}

\begin{figure}[t]
\begin{center}\includegraphics[width=0.4\textwidth]{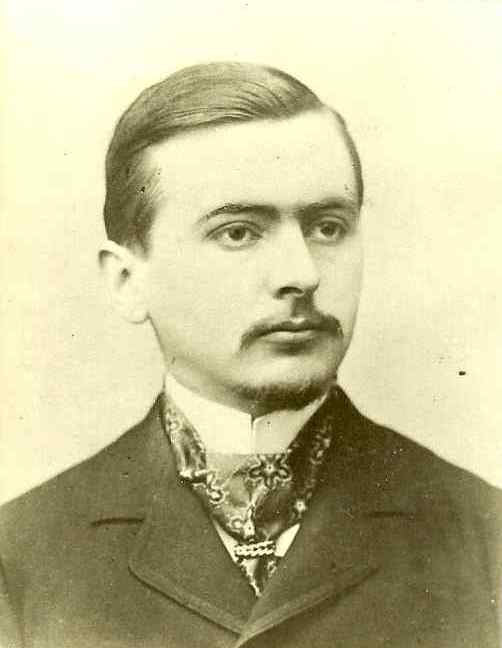}

{\tiny Public domain}\end{center}

\caption{\textbf{\protect\href{https://en.wikipedia.org/wiki/Pierre_Fatou}{Pierre Fatou}}
(1878 - 1929)}\sindex[fam]{Fatou, Pierre}
\end{figure}

We now introduce the fundamental \index{Fatou lemma}\textbf{Fatou}{\bfseries\footnote{\textbf{\href{https://en.wikipedia.org/wiki/Pierre_Fatou}{Pierre Fatou}}
(1878 -1829) is a French mathematician and astronomer. He contributed
in many branches of analysis. The Fatou lemma and the Fatou set are
named after him. After studying at Ecole Normale Supérieure in Paris,
he went graduated and appointed at Paris Observatory, from intern
to full tenure astronomer/}}\sindex[fam]{Fatou, Pierre}\textbf{ lemma}.

\begin{lemma}{Fatou Lemma}{}

Let $\left(\Omega,\mathcal{A},\mu\right)$ be a measured space, and
let $f_{n},n\in\mathbb{N},$ be elements of $\mathcal{M}^{+}.$ Then
the following inequality holds in $\overline{\mathbb{R}}^{+},$\boxeq{
\[
\intop_{\Omega}\liminf_{n\to+\infty}f_{n}\leqslant\liminf_{n\to+\infty}\intop_{\Omega}f_{n}\text{d}\mu.
\]
}

\end{lemma}

\begin{remark}{Utility of Fatou Lemma}{}

Fatou lemma is primarily used to establish the integrability of a
function that arises as the pointwise limit of a sequence of measurable
functions.

\end{remark}

\begin{theorem}{Monotone Convergence Theorem}{}

Let $\left(\Omega,\mathcal{A},\mu\right)$ be a measured space, and
let $\left(f_{n}\right)_{n\in\mathbb{N}}$ be a monotonic sequence
of measurable functions taking values in $\overline{\mathbb{R}},$
of limit $f.$

\textbf{1. Non-decreasing case}

If the sequence $\left(f_{n}\right)_{n\in\mathbb{N}}$ is non-decreasing
and if there exists $n_{0}$ such that
\[
\intop_{\Omega}f^{-}_{n_{0}}\text{d}\mu<+\infty,
\]
 then the following equality holds in $\overline{\mathbb{R}}$\boxeq{
\[
\lim_{n\to+\infty}\nearrow\intop_{\Omega}f_{n}\text{d}\mu=\intop_{\Omega}f\text{d}\mu.
\]
}

\textbf{2. nonincreasing case}

If the sequence $\left(f_{n}\right)_{n\in\mathbb{N}}$ is nonincreasing
and if there exists $n_{0}$ such that
\[
\intop_{\Omega}f^{+}_{n_{0}}\text{d}\mu<+\infty,
\]
 then the following equality holds in $\overline{\mathbb{R}}:$\boxeq{
\[
\lim_{n\to+\infty}\searrow\intop_{\Omega}f_{n}\text{d}\mu=\intop_{\Omega}f\text{d}\mu.
\]
}

\end{theorem}

\begin{theorem}{Dominated Convergence Theorem, First Version}{}

Let $\left(\Omega,\mathcal{A},\mu\right)$ be a measured space, and
let $\left(f_{n}\right)_{n\in\mathbb{N}}$ be a sequence of measurable
functions taking values in $\overline{\mathbb{R}},$ that converges
pointwise to a function $f.$

If $\sup\left|f_{n}\right|$ is $\mu-$integrable, then each function
$f_{n}$ and $f$ are $\mu-$integrable, and the sequence of integrals
$\left(\intop_{\Omega}f_{n}\text{d}\mu\right)_{n\in\mathbb{N}}$ converges
in $\mathbb{R}.$ 

Moreover, we have
\[
\begin{array}{ccccc}
\lim_{n\to+\infty}\intop_{\Omega}f_{n}\text{d}\mu=\intop_{\Omega}f\text{d}\mu & \,\,\,\, & \text{and} & \,\,\,\: & \lim_{n\to+\infty}\intop_{\Omega}\left|f_{n}-f\right|\text{d}\mu=0.\end{array}
\]

\end{theorem}

\begin{remark}{}{}

The hypothesis ``$\sup\left|f_{n}\right|$ is $\mu-$integrable''
is equivalent to the more commonly stated condition that gives the
theorem its name: ``There exists a $\mu-$integrable function $g$
such that for every $n\in\mathbb{N}$ et for every $\omega\in\Omega,$
$\left|f_{n}\left(\omega\right)\leqslant g\left(\omega\right)\right|$''.

\end{remark}

A direct application of this theorem yields the following essential
result: 

\begin{corollary}{Link Riemann-Integrability and Lebesgue-Integrability}{}

If $f$ is Riemann-integrable of arbitrary sign on $\left[a,b\right]$---and
thus bounded---the function $f.\boldsymbol{1}_{\left[a,b\right]}$
is Lebesgue-integrable, and\boxeq{
\[
\intop_{\mathbb{R}}f.\boldsymbol{1}_{\left[a,b\right]}\text{d}\lambda=\intop^{b}_{a}f\left(x\right)\text{d}x.
\]
}

\end{corollary}

\begin{definition}{Local Riemann-Integrability}{}

A function $f$ defined on an arbitrary open or semi-open interval
$I=\left]a,b\right[,$ where $-\infty\leqslant a<b\leqslant+\infty$
is said to be \textbf{locally Riemann-integrable\index{locally Riemann-integrable}}
on $I,$ if it is Riemann-integrable on every bounded closed subinterval
contained in I.

\end{definition}

\begin{proposition}{Necessary and Sufficient Condition for Lebesgue Integrability}{}

Let $f$ be a function defined on an arbitrary open or semi-open interval
$I=\left]a,b\right[,$ where $-\infty\leqslant a<b\leqslant+\infty,$
and assume that $f$ is locally Riemann-integrable on $I.$

Then the function $\boldsymbol{1}_{I}.f$ is Lebesgue-integrable if
and only if the generalized Riemann integral 
\[
\intop^{b}_{a}f\left(x\right)\text{d}x
\]
 is absolutely convergente. 

In this case,\boxeq{
\[
\intop_{\mathbb{R}}\boldsymbol{1}_{I}.f\text{d}\lambda=\intop^{b}_{a}f\left(x\right)\text{d}x.
\]
}

\end{proposition}

\subsubsection*{Integral on a Set. Integral of a $\mu-$almost Everywhere Defined
Function}

\begin{definition}{Integral of a Function on a Set}{}

Let $\left(\Omega,\mathcal{A},\mu\right)$ be a measured space, and
let $f$ a function defined on $\Omega,$ taking values in $\overline{\mathbb{R}},$
that is $\mu-$semi-integrable---respectively $\mu-$integrable.

For every $A\in\mathcal{A},$ the function $\boldsymbol{1}_{A}.f$
is also $\mu-$semi-integrable---respectively $\mu-$integrable.
The integral 
\[
\intop_{\Omega}\boldsymbol{1}_{A}.f\text{d}\mu
\]
is denoted 
\[
\intop_{A}f\text{d}\mu
\]
and called the\textbf{ integral of $f$ over the set $A.$}

\end{definition}

\begin{proposition}{Necessary and Sufficient Condition for an Integral of a Nonnegative Function on $\Omega$ to be Zero}{}

Let $\left(\Omega,\mathcal{A},\mu\right)$ be a measured space, and
let $f\in\mathcal{M}^{+}.$ 

Then 
\[
\intop_{\Omega}f\text{d}\mu=0
\]
if and only if $f=0$ $\mu-$almost everywhere.

\end{proposition}

\begin{proposition}{Sufficient Conditions To Have Equality of Integrals}{}

Let $\left(\Omega,\mathcal{A},\mu\right)$ be a measured space, and
let $f$ and $g$ be functions defined on $\Omega,$ taking values
in $\overline{\mathbb{R}},$ $\mathcal{A}-$measurable, such that
$f=g$ $\mu-$almost everywhere.

(a) If $f$ is nonnegative and $g$ is nonnegative $\mu-$almost everywhere,
then \boxeq{
\[
\intop_{\Omega}f\text{d}\mu=\intop_{\Omega}g\text{d}\mu.
\]
}

(b) If $f$ and $g$ are of arbitrary sign, and if $f$ is $\mu-$integrable,
then $g$ is also $\mu-$integrable, and\boxeq{
\[
\intop_{\Omega}f\text{d}\mu=\intop_{\Omega}g\text{d}\mu.
\]
}

\end{proposition}

\begin{definition}{Trace $\sigma-$algebra. Trace Measured Space}{}

Let $\left(\Omega,\mathcal{A},\mu\right)$ be a measured space, and
let $B\in\mathcal{A}.$

The family of subsets 
\[
B\cap\mathcal{A}=\left\{ B\cap A:\,A\in\mathcal{A}\right\} 
\]
is a $\sigma-$algebra called the \textbf{\mindex{trace sigma-algebra@trace $\sigma-$algebra}\mindex{sigma-algebra@$\sigma-$algebra ! trace}trace
$\sigma-$algebra}---or\textbf{ induced \salg}\mindex{induced \salg}\mindex{sigma-algebra@$\sigma-$algebra ! induced}---of
$\mathcal{A}$ on $B.$ 

We then define the measured space $\left(B,B\cap A,\mu_{\left|B\right.}\right),$
called the \textbf{trace measured space\index{trace measured space}}
on $B$ from $\left(\Omega,\mathcal{A},\mu\right),$ where $\mu_{\left|B\right.}$
is the \textbf{\index{trace measure}trace measure}, restriction of
$\mu$ to $B\cap\mathcal{A}.$ It is often also denoted by $\mu.$

\end{definition}

\begin{definition}{Measurable Extension. Function Defined $\mu-$almost Everywhere}{}

Let $\left(\Omega,\mathcal{A},\mu\right)$ be a measured space, and
let $f$ be a function defined on a subset $\Omega_{f}\subset\Omega.$

For every measurable function $g$ on $\left(\Omega,\mathcal{A},\mu\right),$
define an extension $f_{g}$ of $f$ to $\Omega$ for every $\omega\in\Omega,$
by:
\[
f_{g}\left(\omega\right)=\begin{cases}
f\left(\omega\right), & \text{if }\omega\in\Omega_{f},\\
g\left(\omega\right), & \text{otherwise.}
\end{cases}
\]

If $\Omega_{f}\in\mathcal{A},$ and if $f$ is measurable with respect
to the trace measured space on $\Omega_{f},$ then $f_{g}$ is measurable.
The function $f_{g}$ is called a \textbf{\index{measurable extension}measurable
extension} of $f.$ 

Moreover, if $\mu\left(\Omega_{f}\right)=0,$ then we say that $f$
is \textbf{\mindex{defined mu-almost everywhere@defined $\mu-$almost everywhere}defined
$\mu-$almost everywhere.} Two measurable extensions of $f$ are thus
equal $\mu-$almost everywhere.

\end{definition}

By the previous proposition, if $f$ is defined $\mu-$a.e. and admits
a measurable $\mu-$integrable extension, then any other measurable
extension of $f$ is also $\mu-$integrable and their integrals are
equal. This justifies defining the integral of $f$ as the integral
of any of its measurable extensions. We then say again that $f$ is
\textbf{$\mu-$integrable\mindex{mu-integrable@$\mu-$integrable}}
and we denote $\intop_{\Omega}f\text{d}\mu$ its integral. 

\begin{lemma}{Sufficient Condition to be a $\mu$-almost Everywhere Finite Function}{}

Let $\left(\Omega,\mathcal{A},\mu\right)$ be a measured space and
let $f$ be a function taking values in $\overline{\mathbb{R}},$
and $\mu-$integrable. Then $f$ is finite $\mu-$almost everywhere.

\end{lemma}

\begin{theorem}{Dominated Convergence Theorem, Second Version}{}

Let $\left(\Omega,\mathcal{A},\mu\right)$ be a measured space, and
let $\left(f_{n}\right)_{n\in\mathbb{N}}$ be a sequence of measurable
functions taking values in $\overline{\mathbb{R}},$ converging to
$f$ $\mu-$almost everywhere. 

Suppose there exists a $\mu-$integrable function $g$ such that,
for every $n\in\mathbb{N},$ 
\[
\left|f_{n}\right|\leqslant g,\,\,\,\,\mu-\text{a.e.}
\]

Then the functions $f_{n}$ and $f$ are $\mu-$integrable, and the
sequence of integrals $\left(\intop_{\Omega}f_{n}\text{d}\mu\right)_{n\in\mathbb{N}}$
converges in $\mathbb{R}.$ 

Moreover,
\[
\begin{array}{ccccc}
\lim\limits_{n\to+\infty}\intop_{\Omega}f_{n}\text{d\ensuremath{\mu}=\ensuremath{\intop_{\Omega}f\text{d}\mu}} & \,\,\,\, & \text{and} & \,\,\,\, & \lim\limits_{n\to+\infty}\intop_{\Omega}\left|f_{n}-f\right|\text{d\ensuremath{\mu=0.}}\end{array}
\]

\end{theorem}

\begin{remark}{}{}

The domination hypothesis is equivalent to the hypothesis ``there
exists a function $g$ $\mu-$integrable such that $\mu-$a.e. for
every $n\in\mathbb{N},$ $\left|f_{n}\right|\leqslant g.$ This equivalence
comes from the fact that, as all countable union of sets of measure
zero still has measure zero, allowing us to interchange the conditions
``$\mu-$a.e.'' and ``for every $n\in\mathbb{N}$''.

\end{remark}

\begin{corollary}{}{sum_measurable_functions_with_bounded_abs_value}

Let $\left(\Omega,\mathcal{A},\mu\right)$ be a measured space, and
let $\left(f_{n}\right)_{n\in\mathbb{N}}$ be a sequence of measurable
functions taking values in $\overline{\mathbb{R}},$ such that
\[
\sum^{+\infty}_{n=0}\intop_{\Omega}\left|f_{n}\right|\text{d}\mu<+\infty.
\]

Then the series of function $\sum f_{n}$ converges absolutely $\mu-$almost
everywhere, its sum is $\mu-$integrable, and\boxeq{
\[
\intop_{\Omega}\sum^{+\infty}_{n=0}f_{n}\text{d}\mu=\sum^{+\infty}_{n=0}\intop_{\Omega}f_{n}\text{d}\mu.
\]
}

\end{corollary}

\section{Product Measure and Fubini Theorem}\label{sec:Product-Measure-and}

Let $\left(\Omega_{1},\mathcal{A}_{1},\mu_{1}\right)$ and $\left(\Omega_{2},\mathcal{A}_{2},\mu_{2}\right)$
be two measured spaces. Denote by $\Pi_{1}$ and $\Pi_{2}$ the canonical
projections on $\Omega_{1}\times\Omega_{2}$ onto $\Omega_{1}$ and
$\Omega_{2},$ respectively.

\begin{definition}{Product $\sigma-$algebra}{}

On $\Omega_{1}\times\Omega_{2}$ the $\sigma-$algebra generated by
the semi-algebra of rectangles $A_{1}\times A_{2},$ where $A_{1}\in\mathcal{A}_{1}$
and $A_{2}\in\mathcal{A}_{2},$ is called the \textbf{product $\sigma-$algebra\mindex{product sigma-algebra@product $\sigma-$algebra}\mindex{sigma-algebra@$\sigma-$algebra ! product}}
of $\mathcal{A}_{1}$ and $\mathcal{A}_{2},$ and is denoted by $\mathcal{A}_{1}\otimes\mathcal{A}_{2}.$ 

It is the smallest $\sigma-$algebra that makes the canonical projections
measurable.

\end{definition}

The Caratheodory extension theorem guarantees the existence and uniqueness
of the product measure associated with $\mu_{1}$ and $\mu_{2}.$

\begin{proposition}{Product Measure}{}

If $\mu_{1}$ and $\mu_{2}$ are $\sigma-$finite measures, then there
exists a unique measure $\mu$ on the product measurable space $\left(\Omega_{1}\times\Omega_{2},\mathcal{A}_{1}\otimes\mathcal{A}_{2}\right)$
such that
\[
\forall A_{1}\in\mathcal{A}_{1},\forall A_{2}\in\mathcal{A}_{2},\,\,\,\,\mu\left(A_{1}\times A_{2}\right)=\mu\left(A_{1}\right)\mu\left(A_{2}\right).
\]

This measure $\mu$ is called the \textbf{product measure\index{product measure}}
and is denoted $\mu_{1}\otimes\mu_{2}.$

\end{proposition}

\subsection*{Sections of Sets}

Let $A$ be a subset of $\Omega_{1}\times\Omega_{2}.$ 

For each $\omega_{2}\in\Omega_{2},$ we define the section of $A$
in $\omega_{2},$ possibly empty, by
\[
A^{1}_{\omega_{2}}=\left\{ \omega_{1}\in\Omega_{1}:\,\left(\omega_{1},\omega_{2}\right)\in A\right\} 
\]
and similarly, for each $\omega_{1}\in\Omega_{1},$ we define the
section of $A$ in $\omega_{1},$ possibly empty, by
\[
A^{2}_{\omega_{1}}=\left\{ \omega_{2}\in\Omega_{2}:\,\left(\omega_{1},\omega_{2}\right)\in A\right\} .
\]

If $f$ is a function from $\Omega_{1}\times\Omega_{2}$ to $\overline{\mathbb{R}},$
for every $\omega_{2}\in\Omega_{2},$ we define the partial function
of $f$ in $\omega_{2},$ denoted $f^{1}_{\omega_{2}},$ from $\Omega_{1}$
to $\overline{\mathbb{R}},$ which for every $\omega_{1}\in\Omega_{1}$
maps it to $f\left(\omega_{1},\omega_{2}\right).$ 

Similarly, we define the partial function of $f$ in $\omega_{1},$
denoted $f^{2}_{\omega_{1}},$ from $\Omega_{2}$ to $\overline{\mathbb{R}},$
which for every $\omega_{2}\in\Omega_{2}$ maps it to $f\left(\omega_{1},\omega_{2}\right).$ 

\begin{lemma}{Measurability of Partial Functions on the Product Space}{}

(a) Let $A\in\mathcal{A}_{1}\otimes\mathcal{A}_{2}.$ 

Then:
\begin{itemize}
\item For every $\omega_{2}\in\Omega_{2},$ the section $A^{1}_{\omega_{2}}\in\mathcal{A}_{1}.$ 
\item For every $\omega_{1}\in\Omega_{1},$ the section $A^{2}_{\omega_{1}}\in\mathcal{A}_{2}.$
\end{itemize}
(b) Let $f$ be a $\mathcal{A}_{1}\otimes\mathcal{A}_{2}-$measurable
function from $\Omega_{1}\times\Omega_{2}$ to $\overline{\mathbb{R}}.$ 

Then:
\begin{itemize}
\item For every $\omega_{2}\in\Omega_{2},$ the partial function $f^{1}_{\omega_{2}}$
is $\mathcal{A}_{1}-$measurable.
\item For every $\omega_{1}\in\Omega_{1},$ the partial function $f^{2}_{\omega_{1}}$
is $\mathcal{A}_{2}-$measurable.
\end{itemize}
\end{lemma}

\begin{theorem}{Fubini Theorem}{}

Let $\left(\Omega_{1},\mathcal{A}_{1},\mu_{1}\right)$ and $\left(\Omega_{2},\mathcal{A}_{2},\mu_{2}\right)$
be two measured spaces, where $\mu_{1}$ and $\mu_{2}$ are $\sigma-$finite
measures.

Let $f$ be a $\mathcal{A}_{1}\otimes\mathcal{A}_{2}-$measurable
function from $\Omega_{1}\times\Omega_{2}$ to $\overline{\mathbb{R}}.$

1. \textbf{If $f$ is nonnegative}, then
\begin{itemize}
\item The function $\omega_{1}\mapsto\intop_{\Omega_{2}}f^{2}_{\omega_{1}}\left(\omega_{2}\right)\text{d}\mu_{2}\left(\omega_{2}\right)$
is measurable from $\left(\Omega_{1},\mathcal{A}_{1}\right)$ to $\left(\overline{\mathbb{R}}^{+},\mathcal{B}_{\overline{\mathbb{R}}^{+}}\right).$
\item The function $\omega_{2}\mapsto\intop_{\Omega_{1}}f^{1}_{\omega_{2}}\left(\omega_{1}\right)\text{d}\mu_{1}\left(\omega_{1}\right)$
is measurable from $\left(\Omega_{2},\mathcal{A}_{2}\right)$ to $\left(\overline{\mathbb{R}}^{+},\mathcal{B}_{\overline{\mathbb{R}}^{+}}\right).$
\end{itemize}
Moreover the integral of $f$ over the product space can be computed
by iterated integration through one of the following formula:
\begin{equation}
\intop_{\Omega_{1}\times\Omega_{2}}f\text{d}\mu_{1}\otimes\mu_{2}=\intop_{\Omega_{1}}\left[\intop_{\Omega_{2}}f^{2}_{\omega_{1}}\left(\omega_{2}\right)\text{d}\mu_{2}\left(\omega_{2}\right)\right]\text{d}\mu_{1}\left(\omega_{1}\right),\label{eq:fubini_1}
\end{equation}
or
\begin{equation}
\intop_{\Omega_{1}\times\Omega_{2}}f\text{d}\mu_{1}\otimes\mu_{2}=\intop_{\Omega_{2}}\left[\intop_{\Omega_{1}}f^{1}_{\omega_{2}}\left(\omega_{1}\right)\text{d}\mu_{1}\left(\omega_{1}\right)\right]\text{d}\mu_{2}\left(\omega_{2}\right).\label{eq:fubini_2}
\end{equation}

2. \textbf{If $f$ is of arbitrary sign, and if $f$ is $\mu_{1}\otimes\mu_{2}-$integrable},
then
\begin{itemize}
\item For $\mu_{1}-$a.e., for $\omega_{1}\in\Omega_{1},$ the partial function
$f^{2}_{\omega_{1}}$ is $\mu_{2}-$integrable. 
\item For $\mu_{2}-$a.e., for $\omega_{2}\in\Omega_{2},$ the partial function
$f^{1}_{\omega_{2}}$ is $\mu_{1}-$integrable.
\item The function 
\[
\omega_{1}\mapsto\intop_{\Omega_{2}}f^{2}_{\omega_{1}}\left(\omega_{2}\right)\text{d}\mu_{2}\left(\omega_{2}\right)
\]
 is defined $\mu_{1}-$a.e. and is $\mu_{1}-$integrable.
\item The function 
\[
\omega_{2}\mapsto\intop_{\Omega_{1}}f^{1}_{\omega_{2}}\left(\omega_{1}\right)\text{d}\mu_{1}\left(\omega_{1}\right)
\]
 is defined $\mu_{2}-$a.e. and is $\mu_{2}-$integrable.
\end{itemize}
And the computation by iterated integrals of the integral of $f$
by one of the formula $\refpar{eq:fubini_1}$ or $\refpar{eq:fubini_2}$
still holds.

\end{theorem}

\begin{remark}{}{}

To verify that a function $f$ is $\mu_{1}\otimes\mu_{2}-$integrable,
one often checks the integrability by computing the integral $\intop_{\Omega_{1}\times\Omega_{2}}\left|f\right|\text{d}\mu_{1}\otimes\mu_{2}$
with the help of the first part of Fubini theorem.

\end{remark}

It is straightforward to define the product of $n$ $\sigma-$finite
measures and to verify that this product is associative. The Kolmogorov
Extension Theorem enables the construction of probability measures
on the product space $\mathbb{R}^{\mathbb{N}}.$ This result is particularly
used to construct a probabilized space supporting a family of real-valued
independent random variables with arbitrary given laws.

\begin{definition}{Product $\sigma-$algebra}{}

Let $\left(\Omega_{n},\mathcal{A}_{n}\right)_{n\in\mathbb{N}}$ be
a sequence of measurable spaces. 

For $j\in\mathbb{N},$ denote by $\Pi_{j},$ the canonical projection
from the infinite Cartesian product $\Pi_{n\in\mathbb{N}}\Omega_{n}$
onto $\Omega_{j},$ i.e. the function $\omega\mapsto\omega_{j}$ where
$\omega_{j}$ is the $j-$th coordinate of the sequence $\omega=\left(\omega_{n}\right)_{n\in\mathbb{N}}.$

On $\Pi_{n\in\mathbb{N}}\Omega_{n},$ we define the \textbf{product}
\textbf{$\sigma-$algebra\mindex{product sigma-algebra @ product $\sigma-$algebra}\mindex{sigma-algebra @ $\sigma-$algebra ! product}}
of the $\mathcal{A}_{n},$ as the $\sigma-$algebra generated by the
semi-algebra of cylinders with finite basis, that is the infinite
Cartesian product of the form $\prod_{n\in\mathbb{N}}A_{n},$ where
$A_{n}\in\mathcal{A}_{n}$ and where $A_{n}=\Omega_{n},$ for every
but a finite many indices. 

The product $\sigma-$algebra is then denoted $\otimes_{n\in\mathbb{N}}\mathcal{A}_{n}.$
This is the smallest $\sigma-$algebra that makes all the canonical
projections measurable.

\end{definition}

Equip $\mathbb{R}^{\mathbb{N}}$ with the product $\sigma-$algebra
of the Borel $\sigma-$algebras on $\mathbb{R},$ denoted concisely
by $\mathcal{B}^{\otimes_{\mathbb{N}}}_{\mathbb{R}}.$ 

\begin{theorem}{Kolmogorov\footnotemark Extension Theorem}{kolmogorov_extension_th}

Consider a \textbf{consistent sequence of probabilities}\index{consistent sequence of probabilities},
that is a sequence $\left(P_{n}\right)_{n\in\mathbb{N}}$ such that
for every $n\in\mathbb{N},$ $P_{n}$ is a probability on $\left(\mathbb{R}^{n},\mathcal{B}_{\mathbb{R}^{n}}\right)$
and such that, for every rectangles $\prod^{n}_{j=1}\left]a_{j},b_{j}\right],$
\[
P_{n+1}\left(\prod^{n}_{j=1}\left]a_{j},b_{j}\right]\times\mathbb{R}\right)=P_{n}\left(\prod^{n}_{j=1}\left]a_{j},b_{j}\right]\right).
\]

Then, there exists a unique probability measure $P$ on the measurable
space $\left(\mathbb{R}^{\mathbb{N}},\mathcal{B}^{\otimes_{\mathbb{N}}}_{\mathbb{R}}\right)$
such that for every $n\in\mathbb{N},$
\[
P\left(\left\{ \omega\in\mathbb{R}^{\mathbb{N}}:\,\left(\omega_{1},\dots,\omega_{n}\right)\in\prod^{n}_{j=1}\left]a_{j},b_{j}\right]\right\} \right)=P_{n}\left(\prod^{n}_{j=1}\left]a_{j},b_{j}\right]\right).
\]
\end{theorem}

\footnotetext{See the footnote \ref{fn:Kolmogorov} of Chapter \ref{chap:Random-Phenomena}
for a short summary on Andrey Kolmogorov}

\begin{corollary}{Infinite Product Probability}{infinite_product_proba}

Let $\left(\mu_{n}\right)_{n\in\mathbb{N}^{\ast}}$ be a probability
sequence on $\left(\mathbb{R},\mathcal{B}_{\mathbb{R}}\right)$ and
consider for every $n\in\mathbb{N}^{\ast}$ the product probability
$P_{n}=\otimes^{n}_{j=1}\mu_{j}$ on $\left(\mathbb{R}^{n},\mathcal{B}_{\mathbb{R}^{n}}\right),$
unique probability such that for every rectangles $\prod^{n}_{j=1}\left]a_{j},b_{j}\right],$
\[
P_{n}\left(\prod^{n}_{j=1}\left]a_{j},b_{j}\right]\right)=\prod^{n}_{j=1}\mu_{j}\left(\left]a_{j},b_{j}\right]\right).
\]

The sequence of probabilities $\left(P_{n}\right)_{n\in\mathbb{N}}$
is consistent, and there exists a unique probability $P$ on the measurable
space $\left(\mathbb{R}^{n},\mathcal{B}_{\mathbb{R}^{n}}\right)$
such that for every $n\in\mathbb{N},$

\[
P\left(\left\{ \omega\in\mathbb{R}^{\mathbb{N}}:\,\left(\omega_{1},\dots,\omega_{n}\right)\in\prod^{n}_{j=1}\left]a_{j},b_{j}\right]\right\} \right)=\prod^{n}_{j=1}\mu_{j}\left(\left]a_{j},b_{j}\right]\right).
\]

This probability is called the infinite product probability of the
probabilities $\mu_{n},n\in\mathbb{N}^{\ast}.$

\end{corollary}

\chapter{Laws and Moments of Random Variables}\label{chap:PartIIChap9}

\begin{objective}{}{}

Chapter \ref{chap:PartIIChap9} aims to introduce the laws and moments
of random variables based on the measure theory.
\begin{itemize}
\item Section \ref{sec:Complements-in-Measure} provides additional material
to Chapter \ref{chap:PartIIChap8} on measure theory that is useful
for probability theory. It introduces the principle of extension by
measurability, following the definitions of $\pi-$ and $\lambda-$systems.
This principle is presented in both a set-theoretic and a functional
version, the latter requiring Radon measures. The Radon-Nikodym theorem
is then stated, followed by two key results: the transfer theorem
and the change of variables theorem.
\item Section \ref{sec:Law-of-a-random-variable} addresses the laws of
random variables. After defining the notion of a law, cumulative distribution
functions are introduced. The law obtained from another via a transormation
by a diffeomorphism is described. Marginals are defined, and the method
for computing their law is presented.
\item Section \ref{sec:Moments-of-Random-Variables} begins by extending
Hölder and Minkowski inequalities in order to derive the first properties
of the $\mathscr{L}^{p}$ spaces. It then introduces the moments of
random variables: first, the expectation, variance, and standard-deviation,
followed by higher order moments, as well as the covariance, together
with its relationships with expectation and variance. These notions
are generalized to the case of a finite-dimensional vector space by
introducing the covariance matrix. The Markov and Bienaymé-Chebyshev
inequalities are then presented. The correlation coefficient is defined,
before tackling the problem of linear regression. The chapter concludes
with a summary of the most common laws, either discrete or with density. 
\end{itemize}
\end{objective}

\section*{Introduction}

In this second part, we assume that the abstract theory of measure
and integration is known. A summary of this theory can be found in
Chapter \ref{chap:PartIIChap8}. To this summary, we add in this chapter
some complements that are generally omitted in integration courses,
but are needed in probabilities.

In this chapter, we give the final presentation, within the framework
of measure theory, of the concepts of the law and moments of a random
variable.

\section{Complements in Measure Theory}\label{sec:Complements-in-Measure}

We begin by studying the \textbf{principle of extension by measurability}\index{principle of the extension by measurability},
which is frequently used in probability theory.

\begin{definition}{$\pi-$System. $\lambda-$System}{}

A family $\mathcal{C}$ of subsets of a set $\Omega$ is called a
\textbf{$\pi-$system\mindex{pi-system@$\pi-$system}} if it is stable
under finite intersection.

A family $\mathscr{S}$ of subsets of a set $\Omega$ is called a
\textbf{$\lambda-$system\mindex{lambda-system@$\lambda-$system}}
if it satisfies the following two axioms:

\boxeq{$\left(\lambda_{1}\right)$ For every nondecreasing sequence
$\left(S_{n}\right)_{n\in\mathbb{N}}$ of elements of $\mathscr{S},$
\[
\bigcup_{n\in\mathbb{N}}S_{n}\in\mathscr{S}.
\]

$\left(\lambda_{2}\right)$ For every elements $A$ and $B$ of $\mathscr{S}$
with $A\subset B,$
\[
B\backslash A\in\mathscr{S}.
\]
}

\end{definition}

We will need the notion of the $\pi-$system, as well as $\lambda-$system,
generated by a family $\mathscr{L}$ of subsets of $\Omega.$ We begin
by observing that the intersection of an arbitrary family of $\pi-$system
is itself a $\pi-$system. Similarly the intersection of an arbitrary
family of $\lambda-$system is also a $\lambda-$system.

Moreover, the power set $\mathcal{P}\left(\Omega\right)$ is both
a $\pi-$system and a $\lambda-$system. Thus, there exists a $\pi-$system
containing $\mathscr{L}.$ We then define the \textbf{$\pi-$system}\textbf{generated
by $\mathcal{L}$\mindex{pi-system @ $\pi-$system ! generated by a family of subsets}}
as the intersection of all $\pi-$systems containing $\mathscr{L}.$
Similarly, the \textbf{$\lambda-$system generated by $\mathcal{L}$}\textbf{\mindex{lambda-system @$\lambda$-system ! generated by a family of subsets}}
is the intersection of all $\lambda-$systems containing $\mathscr{L}.$

We may also characterize the $\pi-$system generated by $\mathscr{L}$
as the smallest $\pi-$system with respect to inclusion, that contains
$\mathscr{L}.$ The same holds for a $\lambda-$system.

\begin{remark}{}{}

As in the context of $\sigma-$algebras, the process of defining a
$\lambda-$system by closure is not constructive: in general, we do
not have an explicit expression of a generic element of the $\lambda-$system
generated by $\mathscr{L}.$ However, the $\pi-$system generated
by $\mathcal{L}$ is simply the family of all finite intersections
of elements of $\mathscr{L}.$ This family is the smallest $\pi-$system
containing $\mathscr{L}.$

\end{remark}

\begin{examples}{}{}

On $\mathbb{R},$ the following families are all $\pi-$systems:
\begin{itemize}
\item Intervals $\left]a,b\right[,$ with $a\leqslant b;$
\item Intervals $\left]a,b\right],$ with $a\leqslant b;$
\item Intervals $\left[a,b\right],$ with $a\leqslant b,$ together with
the empty set;
\item Half-lines $\left[a,+\infty\right[,$ with $a\in\mathbb{R};$
\item Half-lines $\left]-\infty,a\right],$ with $a\in\mathbb{R}.$
\end{itemize}
Similarly, on $\mathbb{R}^{d},$ the following families are all $\pi-$systems:
\begin{itemize}
\item Open sets;
\item Bounded open sets;
\item Closed sets;
\item Rectangles of the form $\prod^{d}_{i=1}\left[a_{i},b_{i}\right],\,a_{i}\leqslant b_{i}$
together with the empty set.
\end{itemize}
\end{examples}

In what follows, we will frequently encounter $\lambda-$systems.
To get an idea, here is an example of a $\lambda-$system that is
not a $\sigma-$algebra: if $\Omega$ is an uncountable set, the family
of its countable subsets is both a $\pi-$system and a $\lambda-$system.
Nonetheless, this family does not contain $\Omega$ and, in general,
is not stable under complementation. Hence, it is not a $\sigma-$algebra.

The relations between these different structures are stated more precisely
in the following lemma.

\begin{lemma}{}{extension_by_measurability}

A $\lambda-$system $\mathscr{S}$ on $\Omega$ is a $\sigma-$algebra
on $\Omega$ if and only if $\mathscr{S}$ is a $\pi-$system and
$\Omega\in\mathscr{S}.$

\end{lemma}

\begin{proof}{}{}

The necessary condition is straightforward.

Let us prove the sufficient condition. 

Suppose that $\mathscr{S}$ is both a $\lambda-$system and a $\pi-$system
on $\Omega,$ and that $\Omega\in\mathscr{S}.$ 
\begin{itemize}
\item Then $\mathscr{S}$ is stable under complementation---since $\Omega\in\mathscr{S}$
implies that for every $A\in\mathscr{S},$ $A\subset\Omega$ and by
the second axiom of the definition of a $\lambda-$system, $A^{c}=\Omega\backslash A\in\mathscr{S}.$
\item $\mathscr{S}$ is also stable by finite unions. To prove it, it is
enough to note that if $A$ and $B$ are some elements of $\mathscr{S},$
$A^{c}$ and $B^{c}$ are also elements of $\mathscr{S}$ using the
first point.\\
Since $\mathscr{S}$ is a $\pi-$system, $A^{c}\cap B^{c}\in\mathscr{S}.$
Moreover $A^{c}\cap B^{c}=\left(A\cup B\right)^{c},$ and thus,
\[
\left(A\cup B\right)^{c}\in\mathscr{S}.
\]
Consequently,
\[
A\cup B\in\mathscr{S}.
\]
\item It remains to show that $\mathscr{S}$ is stable by countable union.
\\
Let $\left(A_{n}\right)_{n\in\mathbb{N}^{\ast}}$ be a sequence of
elements of $\mathscr{S}.$ We construct a nondecreasing sequence
$\left(B_{n}\right)_{n\in\mathbb{N}^{\ast}}$ of elements of $\mathscr{S}$
with the same union as the sequence $\left(A_{n}\right)_{n\in\mathbb{N}^{\ast}}.$
It is sufficient to define $B_{n}$ by induction, setting 
\[
B_{1}=A_{1}
\]
 and for every $n\geqslant2,$ 
\[
B_{n}=\bigcup^{n}_{j=1}A_{j}.
\]
Since $\mathscr{S}$ is a $\lambda-$system, $B_{n}\in\mathscr{S}.$
\end{itemize}
\end{proof}

\textbf{The following lemma is frequently used in probability theory.
It allows to extend a property satisfied by a family of events forming
a structure of $\lambda-$system to the $\sigma-$algebra generated
from this family.}

\begin{lemma}{Principle of Extension by Measurability---Set Version}{meas_ext_set}

Let $\mathscr{S}$ be a $\lambda-$system on $\Omega$ which contains
a $\pi-$system $\mathscr{C}$ and such that $\Omega\in\mathscr{S}.$
Then $\mathscr{S}$ contains the $\sigma-$algebra $\sigma\left(\mathscr{C}\right)$
generated by $\mathscr{C}.$

\end{lemma}

\begin{proof}{}{}

It is enough to prove that the $\lambda-$system $\Lambda,$ generated
by the $\pi-$system $\mathscr{C}$ and $\Omega,$ is equal to the
\salg $\sigma\left(\mathscr{C}\right).$

We will show that \textbf{$\Lambda$ is a $\pi-$system}. To do this,
define, for each $A\in\mathcal{P}\left(\Omega\right),$ the set family
\[
\Lambda_{A}=\left\{ B\in\Lambda:\,B\cap A\in\Lambda\right\} .
\]

\begin{itemize}
\item \textbf{$\Lambda$ being a $\lambda-$system, it is the same for $\Lambda_{A}.$}
Indeed, let $\left(B_{i}\right)_{1\leqslant i\leqslant n}$ be a nondecreasing
sequence of elements of $\Lambda_{A},$ then, since $B_{i}\in\Lambda$
and $\Lambda$ is a $\lambda-$system, $\bigcup^{n}_{i=1}B_{i}\in\Lambda.$
Moreover, as $B_{i}\in\Lambda_{A},$ then $B_{i}\in\Lambda$ and $B_{i}\cap A\in\Lambda.$
Hence, $\left(B_{i}\cap A\right)_{1\leqslant i\leqslant n}$ is also
a nondecreasing sequence, and thus $\bigcup^{n}_{i=1}\left(B_{i}\cap A\right)$
is in $\Lambda$ since $\Lambda$ is a $\lambda-$system, and thus
as
\[
\bigcup^{n}_{i=1}\left(B_{i}\cap A\right)=\left(\bigcup^{n}_{i=1}B_{i}\right)\cap A\in\Lambda
\]
which shows that $\bigcup^{n}_{i=1}B_{i}\in\Lambda_{A}.$ So $\Lambda_{A}$
is a $\lambda-$system.
\item Clearly, $\Lambda_{A}\subset\Lambda.$ Furthermore, for every $A\in\mathscr{C},$
by definition of $\Lambda,$ $\Omega\in\Lambda_{A}$ and $\mathscr{C}\subset\Lambda_{A},$
hence $\Lambda_{A}=\Lambda,$ the family $\Lambda$---defined as
the $\lambda-$system generated by $\mathscr{C}$ and $\Omega$---being
the smallest $\lambda-$system containing $\mathscr{C}$ and $\Omega.$\\
Consequently,
\[
\forall A\in\mathscr{C},\forall B\in\Lambda,\,\,\,\,B\cap A\in\Lambda.
\]
In other words,
\[
\forall B\in\Lambda,\mathscr{C}\subset\Lambda_{B},
\]
which implies that for every $B\in\Lambda,$ the family $\Lambda_{B}$
is a $\lambda-$system containing $\mathscr{C}$ and $\Omega,$ thus
\[
\Lambda_{B}=\Lambda.
\]
Hence, we proved that
\[
\forall B\in\Lambda,\forall C\in\Lambda,\,\,\,\,B\cap C\in\Lambda
\]
showing that $\Lambda$ is a $\pi-$system.
\end{itemize}
Since $\Lambda$ is both a $\pi-$system and a $\lambda-$system containing
$\Omega,$ the previous lemma implies that $\Lambda$ is a $\sigma-$algebra.
Therefore
\[
\Lambda\supset\sigma\left(\mathscr{C}\right).
\]

Since $\sigma\left(\mathscr{C}\right)$ is a $\lambda-$system containing
$\mathscr{C}$ and $\Omega,$
\[
\Lambda\subset\sigma\left(\mathscr{C}\right).
\]

We conclude on the equality
\[
\Lambda=\sigma\left(\mathscr{C}\right).
\]

\end{proof}

Next, we state an important application of the previous lemma.

\begin{theorem}{Uniqueness of Measures}{unicity_measures}

Let $\mu_{1}$ and $\mu_{2}$ be two nonnegative measures on the probabilizable
space $\left(\Omega,\mathcal{A}\right)$ such that
\[
\forall A\in\mathscr{C},\,\,\,\,\mu_{1}\left(A\right)=\mu_{2}\left(A\right),
\]
where $\mathscr{C}$ is a $\pi-$system that generates the $\sigma-$algebra
$\mathcal{A}.$

1. If $\mu_{1}$ and $\mu_{2}$ are bounded and have the same mass,
then $\mu_{1}=\mu_{2}.$

2. If one of the measures $\mu_{1}$ or $\mu_{2}$ is unbounded, and
there exists a sequence $\left(E_{n}\right)_{n\in\mathbb{N}}$ of
elements of $\mathscr{C}$ such that 
\[
\Omega=\cup_{n\in\mathbb{N}}E_{n}
\]
and
\[
\forall n\in\mathbb{N},\,\,\,\,\mu_{1}\left(E_{1}\right)=\mu_{2}\left(E_{2}\right)<+\infty,
\]
then $\mu_{1}=\mu_{2}.$

\end{theorem}

\begin{proof}{}{}

1. If $\mu_{1}$ and $\mu_{2}$ are bounded and with same mass, consider
the family
\[
\mathscr{S}=\left\{ A\in\mathcal{A}:\,\mu_{1}\left(A\right)=\mu_{2}\left(A\right)\right\} .
\]
This is a $\lambda-$system containing $\mathscr{C}$ and $\Omega.$
It thus contains the $\sigma-$algebra $\mathcal{A}$ generated by
$\mathscr{C}.$ This shows that on $\mathcal{A}$ 
\[
\mu_{1}=\mu_{2}.
\]

2. If one of the measures $\mu_{1}$ or $\mu_{2}$ is unbounded, consider
their restrictions to the sets $E_{n}.$ From Part \ref{part:Introduction-to-Probability},
they are equal for every $n\in\mathbb{N}.$ The Poincaré formula---Part
\ref{part:Introduction-to-Probability}, Proposition $\ref{pr:poincareform}$---is
still valid for finite measures. The restriction of $\mu_{1}$ and
$\mu_{2}$ to the sets $E_{n}$ are finite. Hence, by the Poincaré
formula, the restrictions of $\mu_{1}$ and $\mu_{2}$ to the sets
$F_{n}=\bigcup_{0\leqslant j\leqslant n}E_{j}$ are equal. Therefore,
the measures $\mu_{1}$ and $\mu_{2}$ coincide on $\Omega,$ since
the sequence $F_{n}$ is nondecreasing and that $\bigcup_{n}F_{n}=\Omega.$

\end{proof}

\begin{remark}{}{}

The theorem hypotheses imply that $\mu_{1}$ and $\mu_{2}$ are $\sigma-$finite.
Moreover, the theorem implies that if two probability measures agree
on a $\pi-$system generating $\mathcal{A},$ then they are equal.

\end{remark}

\begin{example}{Application}{}

If two measures on $\mathbb{R}$ agree on all half-lines $\left]-\infty,x\right],$
where $x\in\mathbb{R},$ then they are equal.

\end{example}

The following corollary is frequently used in probability computations.

\begin{denotation}{}{}

We denote $\mathscr{C}_{\mathscr{K}}\left(\mathbb{R}^{d}\right)$
the set of continuous functions from $\mathbb{R}^{d}$ to $\mathbb{R}$
with compact support and $\mathscr{C}^{+}_{\mathscr{K}}\left(\mathbb{R}^{d}\right)$
the subset of $\mathscr{C}_{\mathscr{K}}\left(\mathbb{R}^{d}\right)$
of nonnegative functions.

\end{denotation}

\begin{figure}[t]
\begin{center}\includegraphics[width=0.4\textwidth]{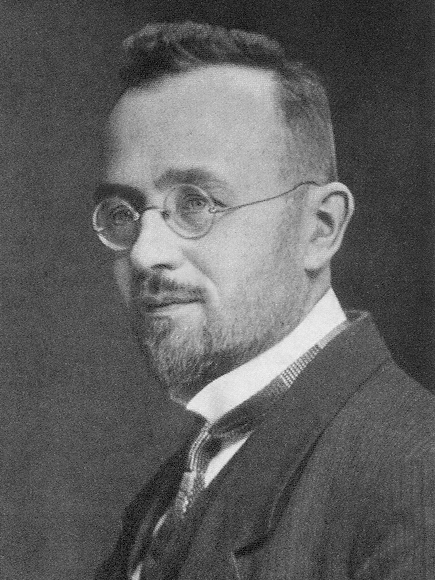}

{\tiny Public domain}\end{center}

\caption{\textbf{\protect\href{https://en.wikipedia.org/wiki/Johann_Radon}{Johann Radon}}
(1887 - 1956)}\sindex[fam]{Radon, Johann}
\end{figure}

\begin{corollary}{Equality of Radon Measures}{equality_radon_measures}

Let $\mu_{1}$ and $\mu_{2}$ be two nonnegative measures on $\left(\mathbb{R}^{d},\mathcal{B}_{\mathbb{R}^{d}}\right)$
finite on every compact set---we say that they are \textbf{\index{Radon measure}\sindex[fam]{Radon, Johann}Radon\footnotemark
measures}. If
\[
\forall f\in\mathscr{C}^{+}_{\mathscr{K}}\left(\mathbb{R}^{d}\right),\,\,\,\,\intop_{\mathbb{R}^{d}}f\text{d}\mu_{1}=\intop_{\mathbb{R}^{d}}f\text{d}\mu_{2},
\]
then the measures $\mu_{1}$ and $\mu_{2}$ are equal.

\end{corollary}

\footnotetext{\href{https://en.wikipedia.org/wiki/Johann_Radon}{Johann Radon}
(1887 - 1956) is an Austrian mathematician. He worked on the calculus
of variations during his PhD, and gave many lasting contributions
such as the Radon-Nikodym theorem, the Radon measure, the Radon transform
in integral geometry, the Radon's theorem and other results.}

\begin{figure}[t]
\begin{center}\includegraphics[width=0.4\textwidth]{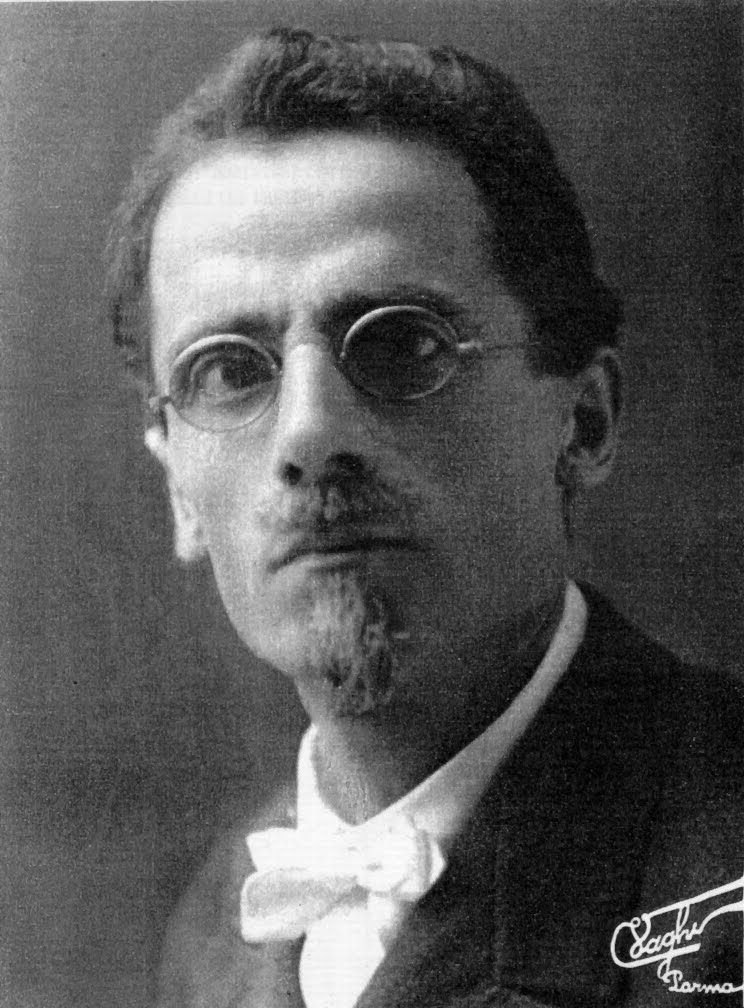}

{\tiny Public domain}\end{center}

\caption{\textbf{\protect\href{https://en.wikipedia.org/wiki/Beppo_Levi}{Beppo Levi}}
(1875 - 1961)}\sindex[fam]{Beppo, Levi}
\end{figure}

\begin{proof}{}{}

The class $\mathscr{C}$ of bounded open subsets of $\mathbb{R}^{d}$
is a $\pi-$system. The measures $\mu_{1}$ and $\mu_{2}$ coincide
on $\mathscr{C}.$ Indeed, for every $O\in\mathscr{C},$ there exists
a sequence $\left(f_{n}\right)_{n\in\mathbb{N}}$ of functions in
$\mathscr{C}^{+}_{\mathscr{K}}\left(\mathbb{R}^{d}\right)$ that converges
pointwise to $\boldsymbol{1}_{O}.$ 

By the \textbf{\sindex[fam]{Beppo, Levi}\href{https://en.wikipedia.org/wiki/Beppo_Levi}{Beppo Levi}}
property\footnotemark and the assumption that the integrals of any
function in $\mathscr{C}^{+}_{\mathscr{K}}\left(\mathbb{R}^{d}\right)$---and
thus the $f_{n}$---agree under $\mu_{1}$ and $\mu_{2},$
\[
\mu_{1}\left(O\right)=\lim_{n\to+\infty}\intop_{\mathbb{R}^{d}}f_{n}\text{d}\mu_{1}=\lim_{n\to+\infty}\intop_{\mathbb{R}^{d}}f_{n}\text{d}\mu_{2}=\mu_{2}\left(O\right)<+\infty.
\]

We conclude by applying the uniqueness theorem for measures.

\end{proof}

%\addtocounter{footnote}{-1}

\footnotetext{\textbf{\href{https://en.wikipedia.org/wiki/Beppo_Levi}{Beppo Levi}}
(1875-1971) was an Italian mathematician. He left Italy, being expelled
from his position at the University of Bologna because he was Jewish,
and migrated to Argentina for the rest of his life. He contributed
in algebraic and diophantine geometry as well as in Lebesgue integration.}

\stepcounter{footnote}

\footnotetext{The Beppo Levi property also known as the monotone
convergence theorem.}

We now give a functional version of Lemma $\ref{lm:meas_ext_set}$.

\begin{theorem}{Measurability Extension Principle. Functional Version}{meas_ext_pple}

Let $\mathscr{C}$ be a $\pi-$system on $\Omega,$ and let $\mathscr{H}$
be a vector space of real-valued functions on $\Omega$ such that:

(i) For every nondecreasing sequence $\left(h_{n}\right)_{n\in\mathbb{N}}$
of nonnegative elements of $\mathscr{H}$ such that $h\equiv\sup_{n\in\mathbb{N}}h_{n}$
is finite---respectively bounded---$h\in\mathscr{H}.$

(ii) $\boldsymbol{1}_{\Omega}\in\mathscr{H},$ and for every $C\in\mathscr{C},\,\boldsymbol{1}_{C}\in\mathscr{H}.$

Then $\mathscr{H}$ contains all the $\sigma\left(\mathscr{C}\right)-$measurable---respectively
$\sigma\left(\mathscr{C}\right)-$measurable and bounded---real-valued
functions.

\end{theorem}

\begin{proof}{}{}
\begin{itemize}
\item It is enough to show that for every $A\in\sigma\left(\mathscr{C}\right),$
the functions $\boldsymbol{1}_{A}\in\mathscr{H}.$ Indeed, under this
assumption, the vector space $\mathscr{H}$ contains all $\sigma\left(\mathscr{C}\right)-$measurable
step functions, since any such function can be written 
\[
\sum_{i\in I}a_{i}\boldsymbol{1}_{A_{i}},\,\,\,\,\text{with}\,I\,\text{finite,\,\ensuremath{a_{i}\in}\ensuremath{\mathbb{R}}}\,\,\,\,\text{and}\,\,\,\,A_{i}\in\sigma\left(\mathscr{C}\right).
\]
Then, by hypothesis (i), any nonnegative and finite---respectively,
bounded---$\sigma\left(\mathscr{C}\right)-$measurable function belongs
to $\mathscr{H}.$ Finally, $\mathscr{H}$ is going to contain any
finite---respectively bounded---$\sigma\left(\mathscr{C}\right)-$measurable
function $h,$ since such a function can be decomposed under the form:
$h=h^{+}-h^{-},$ where $h^{+}$ and $h^{-}$ are nonnegative finite---respectively
bounded---,and $\sigma\left(\mathscr{C}\right)-$measurable.
\item It remains to show that for every $A\in\sigma\left(\mathscr{C}\right),$
the functions $\boldsymbol{1}_{A}\in\mathscr{H},$ i.e. that 
\[
\mathscr{S}=\left\{ A\in\mathcal{P}\left(\Omega\right):\,\boldsymbol{1}_{A}\in\mathscr{H}\right\} 
\]
contains $\sigma\left(\mathscr{C}\right).$ \\
By assumption,
\[
\mathscr{S}\supset\mathscr{C}\,\,\,\,\text{and\,\,\,\,\ensuremath{\Omega\in}\ensuremath{\mathscr{S}}}.
\]
We now show that $\mathscr{S}$ is a $\lambda-$system:
\begin{itemize}
\item Since, from the one hand, $\mathscr{H}$ being a vector space, for
every $S_{1}$ and $S_{2}$ such that $S_{1}S_{2},$ 
\[
\boldsymbol{1}_{S_{1}\backslash S_{2}}=\boldsymbol{1}_{S_{1}}-\boldsymbol{1}_{S_{2}}\in\mathscr{H}.
\]
\item And, on the other hand, by the first hypothesis, for every nondecreasing
sequence $\left(S_{n}\right)_{n\in\mathbb{N}}\subset\mathscr{S},$
\[
\boldsymbol{1}_{\bigcup_{n\in\mathbb{N}}S_{n}}=\sup_{n\in\mathbb{N}}\boldsymbol{1}_{S_{n}}\in\mathscr{H}.
\]
\end{itemize}
\item By the principle of extension by measurability---Lemma $\ref{lm:meas_ext_set}$---it
follows that 
\[
\sigma\left(\mathscr{C}\right)\subset\mathscr{S}.
\]
 This proves the theorem.
\end{itemize}
\end{proof}

We end this section by recalling, without proof, the statements of
some commonly used theorems.

\begin{definition}{Density of a Measure}{}

Let $\mu$ be a nonnegative measure on the probabilizable space $\left(\Omega,\mathcal{A}\right).$

Let $f$ be a nonnegative, measurable, real-valued function defined
on this space.

The measure defined by \boxeq{
\[
A\mapsto\intop_{A}f\text{d}\mu
\]
}is called \textbf{measure with density $f$ with respect to $\mu,$\index{measure with density a function with respect to a measure}}
and is denoted by\footnotemark \boxeq{
\[
f\cdot\mu.
\]
}

\end{definition}

\footnotetext{A notation justified by formula $\refpar{eq:integration_vs_density_measure}$}

\begin{definition}{Absolutely Continuous Measure. Foreign Measures}{abs_cont_foreign_meas}

A measure $\nu$ on $\left(\Omega,\mathcal{A}\right)$ is said \textbf{absolutely
continuous\mindex{measure ! absolutely continuous}} with respect
to $\mu$ if\boxeq{
\[
\forall A\in\mathcal{A},\,\,\,\,\mu\left(A\right)=0\Longrightarrow\nu\left(A\right)=0.
\]
}We write $\nu\ll\mu.$

The \textbf{measures} $\mu$ and $\nu$ on $\left(\Omega,\mathcal{A}\right)$
are said to be \textbf{\mindex{measures ! foreign}foreign} or in
a more modern terminology are \textbf{\mindex{measures ! mutually singular}mutually
singular} if there exists $N\in\mathcal{A}$ such that 
\[
\mu\left(N\right)=0\,\,\,\,\text{and}\,\,\,\,\nu\left(N^{c}\right)=0.
\]
 We write $\nu\perp\mu.$

\end{definition}

\begin{figure}[t]
\begin{center}\includegraphics[width=0.4\textwidth]{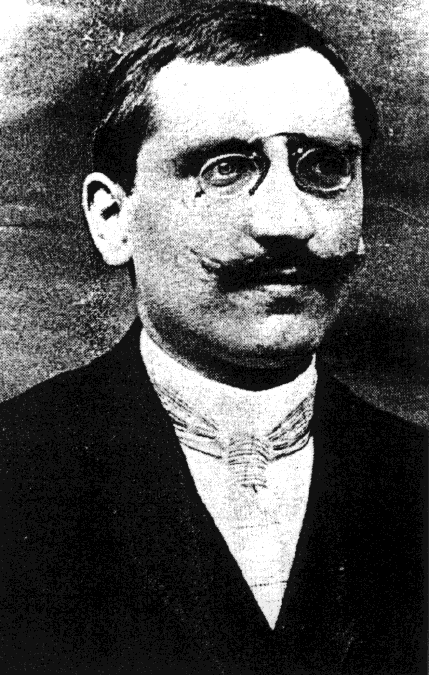}

{\tiny Public domain}\end{center}

\caption{\textbf{\protect\href{https://fr.wikipedia.org/wiki/Henri-L\%25C3\%25A9on_Lebesgue}{Henri Lebesgue}}
(1875 - 1941)}\sindex[fam]{Lebesgue, Henri-Léon}
\end{figure}

\begin{figure}[t]
\begin{center}\includegraphics[width=0.4\textwidth]{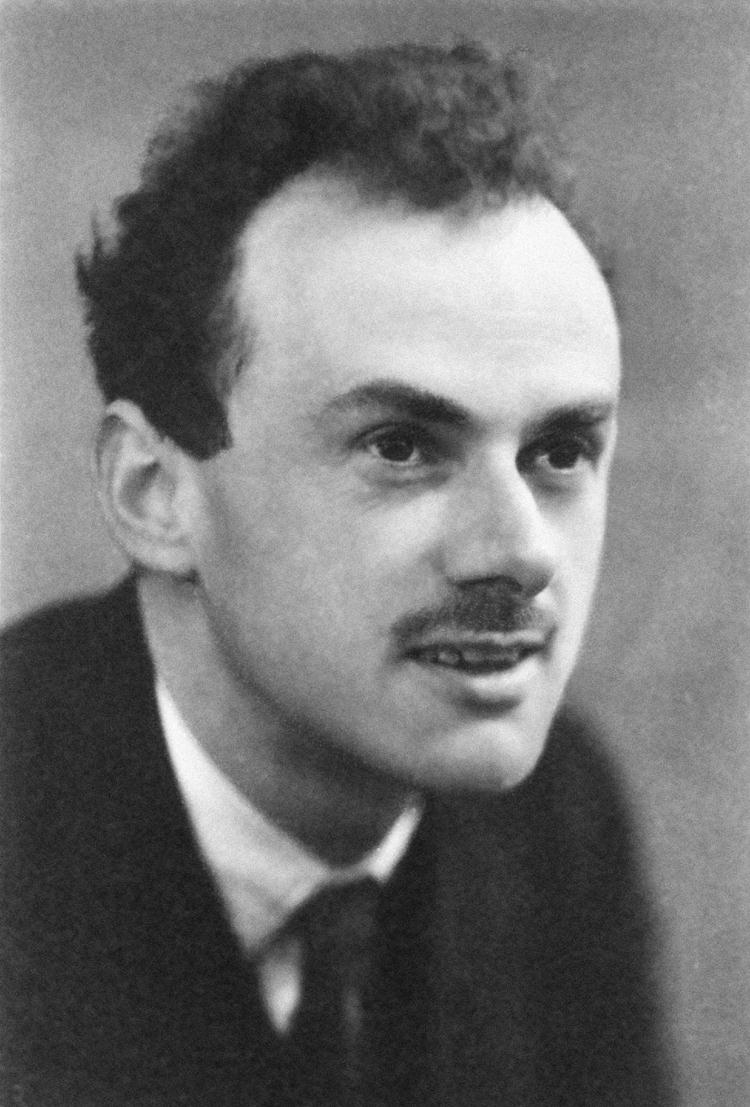}

{\tiny\href{https://en.wikipedia.org/wiki/Nobel_Foundation}{Nobel Foundation}
Public domain}\end{center}

\caption{\textbf{\protect\href{https://fr.wikipedia.org/wiki/Paul_Dirac}{Paul Dirac}}
(1902 - 1984)}\sindex[fam]{Dirac, Paul}
\end{figure}

\begin{example}{}{}

The \textbf{\href{https://fr.wikipedia.org/wiki/Henri-L\%25C3\%25A9on_Lebesgue}{Lebesgue}\footnotemark\sindex[fam]{Lebesgue, Henri-Léon}}
measure $\lambda$ on $\mathbb{R}$ and the \textbf{Dirac\footnotemark\sindex[fam]{Dirac, Paul}}
measure $\delta_{0}$ at 0 are mutually singular, since 
\[
\lambda\left(\left\{ 0\right\} \right)=\delta\left(\left\{ 0\right\} ^{c}\right)=0.
\]

\end{example}

\addtocounter{footnote}{-1}

\footnotetext{\textbf{\href{https://fr.wikipedia.org/wiki/Henri-L\%25C3\%25A9on_Lebesgue}{Henri-Léon Lebesgue}\sindex[fam]{Lebesgue, Henri-Léon}}
(1875-1941) is one of the most famous French mathematician from the
first half of the XX-th century. He is famous for his integration
theory published in his Thesis in 1902, which extends the Borel work.
He is also at the origin of the Fourier transform.}

\stepcounter{footnote}

\footnotetext{\textbf{\href{https://fr.wikipedia.org/wiki/Paul_Dirac}{Paul Dirac}}\sindex[fam]{Dirac, Paul}
(1902-1984) is an English mathematician and physicit. He is one of
the ``fathers'' of the quantum mecanic and has foreseen the existence
of the anti-matter. He won with Erwin Schrödinger the Nobel price
in Physics in 1933.}

\begin{figure}[t]
\begin{center}\includegraphics[width=0.4\textwidth]{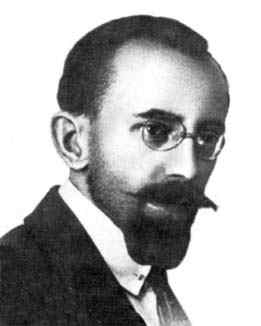}

{\tiny\href{https://mathshistory.st-andrews.ac.uk/Biographies/Nikodym/}{MacTutor}
CC-BY-SA 4.0}\end{center}

\caption{\textbf{\protect\href{https://en.wikipedia.org/wiki/Otto_M._Nikodym}{Otto M. Nikodym}}
(1887 - 1974)}\sindex[fam]{Nikodym, Otto}
\end{figure}

If $\nu=f\cdot\mu,$ then clearly $\nu\ll\mu.$ The converse---whether
every measure $\nu\ll\mu$ can be written in this form---is the object
of the \textbf{Radon-Nikodym}{\bfseries\footnote{\textbf{\href{https://en.wikipedia.org/wiki/Otto_M._Nikodym}{Otto M. Nikodym}}
(1887-1974) was a Polish mathematician. He worked in many areas of
mathematics, but his most well-know contribution relates to the development
of the Lebesgue-Radon-Nikodym integral. He worked also on the theory
of operators in Hilbert space, based on Boolean lattices.}}\textbf{\sindex[fam]{Nikodym, Otto}} theorem. For a proof of this
theorem, the interested reader may refer to \cite{neveu1972martingales}
or \cite{metivier1968notions}.

\begin{theorem}{Radon-Nikodym Theorem}{}

Let $\left(\Omega,\mathcal{A}\right)$ be a probabilizable space,
and let $\mu$ be a $\sigma-$finite measure on this space. Let $\nu$
be a measure such that $\nu\ll\mu.$ Then, there exists a nonnegative
measurable function $f$---unique up to $\mu-$almost everywhere
equivalence---such that \boxeq{
\[
\nu=f\cdot\mu.
\]
}

\end{theorem}

\begin{proposition}{Integration with Respect to a Measure with Density}{}

Let $\mu$ be a nonnegative measure on the probabilizable space $\left(\Omega,\mathcal{A}\right),$
and let $f$ be a nonnegative measurable real-valued function defined
on this space. 

Let $\nu=f\cdot\mu$ be the measure with density $f$ with respect
to $\mu.$ 

Let $h$ be a measurable function on $\left(\Omega,\mathcal{A}\right).$
\begin{itemize}
\item If $h$ is nonnegative, then\boxeq{
\begin{equation}
\intop_{\Omega}h\text{d}\nu=\intop_{\Omega}h\cdot f\text{d}\mu.\label{eq:integration_vs_density_measure}
\end{equation}
}
\item If $h$ is of arbitrary sign, for $h$ to be $\nu-$integrable, it
must and it suffices that $h\cdot f$ is $\mu-$integrable and, in
this case, the equality $\refpar{eq:integration_vs_density_measure}$
still holds.
\end{itemize}
\end{proposition}

\begin{definition}{Measure Image}{}

Let $T$ be a measurable function from the probabilizable space $\left(E_{1},\mathcal{E}_{1}\right)$
to the probabilizable space $\left(E_{2},\mathcal{E}_{2}\right).$ 

Let $\mu_{1}$ be a measure on $\left(E_{1},\mathcal{E}_{1}\right).$ 

The measure $\mu_{2}$ on $\left(E_{2},\mathcal{E}_{2}\right)$ defined
by\boxeq{
\[
\forall B\in\mathcal{E}_{2},\,\,\,\,\mu_{2}\left(B\right)=\mu_{1}\left(T^{-1}\left(B\right)\right)
\]
}is called the \textbf{image measure\index{image measure}} of $\mu_{1}$
by $T$ and denoted $T\left(\mu_{1}\right).$

\end{definition}

\begin{theorem}{Measure Image Theorem---Transfer Theorem}{}

Let $T$ be a measurable function from the probabilizable space $\left(E_{1},\mathcal{E}_{1}\right)$
to the probabilizable space $\left(E_{2},\mathcal{E}_{2}\right).$ 

Let $T\left(\mu_{1}\right)$ be the image measure of $\mu_{1}$ by
$T.$ 

Let $h$ be a measurable function on $\left(E_{2},\mathcal{E}_{2}\right).$
\begin{itemize}
\item If $h$ is nonnegative, then\boxeq{
\begin{equation}
\intop_{E_{2}}h\text{d}T\left(\mu_{1}\right)=\intop_{E_{1}}h\circ T\text{d}\mu_{1}.\label{eq:measure_image_formula}
\end{equation}
}
\item If $h$ is of arbitrary sign, then for $h$ to be $T\left(\mu_{1}\right)-$integrable
it must and it suffices that $h\circ T$ is $\mu_{1}-$integrable;
in that case, the equality $\refpar{eq:measure_image_formula}$ still
holds.
\end{itemize}
\end{theorem}

\begin{theorem}{Change of Variables}{}

Let $T$ be a $C^{1}$-diffeomorphism from an open set $U$ of $\mathbb{R}^{d}$
onto an open set $V$ of $\mathbb{R}^{d}.$

Let $f$ be a measurable real-valued function defined on $U.$

Then $f$ is Lebesgue-integrable on $U$ if and only if the function
$v\mapsto\left|\text{det}\left(T^{-1}\right)^{\prime}\left(v\right)\right|f\left(T^{-1}\left(v\right)\right)$
is Lebesgue-integrable on $V.$

In this case,\boxeq{
\begin{equation}
\intop_{U}f\left(x\right)\text{d}\lambda_{d}\left(x\right)=\intop_{V}\left|\text{det}\left(T^{-1}\right)^{\prime}\left(v\right)\right|f\left(T^{-1}\left(v\right)\right)\text{d}\lambda_{d}\left(v\right).\label{eq:variables_change_formula}
\end{equation}
}

\end{theorem}

\begin{figure}[t]
\begin{center}\includegraphics[width=0.4\textwidth]{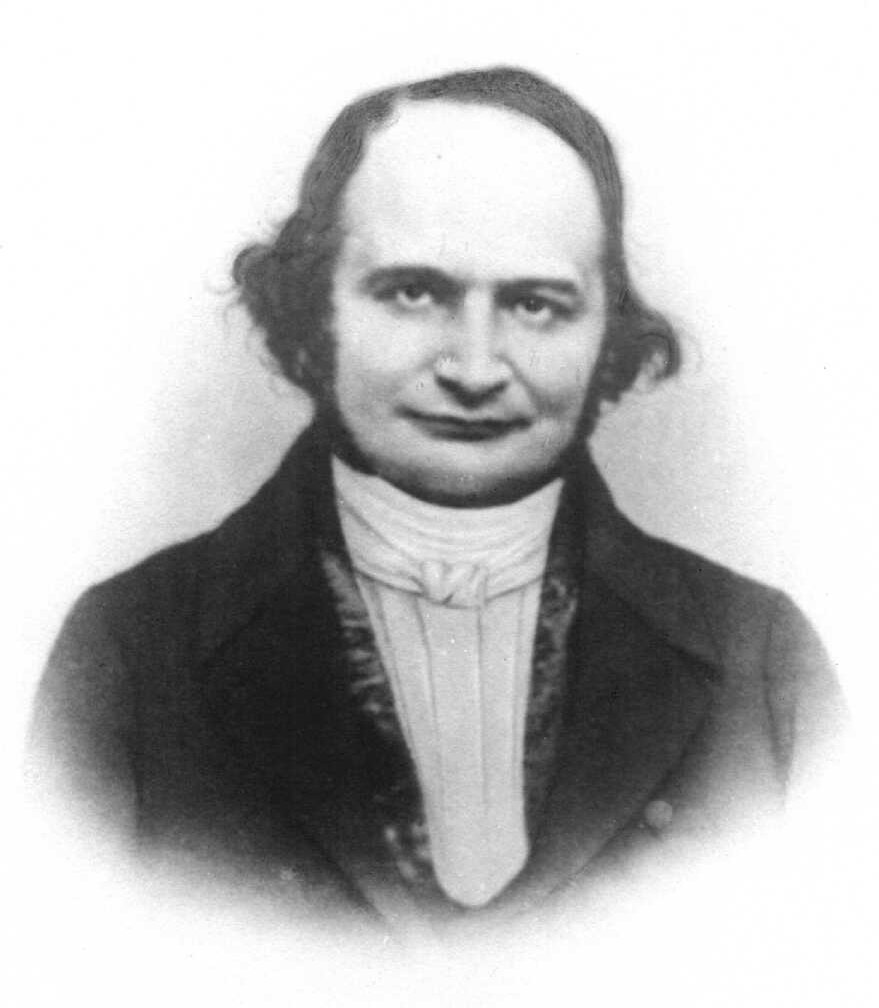}

{\tiny\href{http://www.sil.si.edu/digitalcollections/hst/scientific-identity/explore.htm}{Smithsonian Libraries}
Public Domain}\end{center}

\caption{\textbf{\protect\href{https://en.wikipedia.org/wiki/Carl_Gustav_Jacob_Jacobi}{Carl Jacobi}}
(1804 - 1851)}\sindex[fam]{Jacobi, Carl}
\end{figure}

\begin{remark}{}{}

We often say that the right-hand side of the equality $\refpar{eq:variables_change_formula}$
is obtained from the left-hand side by doing the change of variable
$v=T\left(x\right),$ or equivalently $x=T^{-1}\left(v\right),$ where
$v$ is considered the new variable and $x$ the old one. 

Moreover, $\text{det}\left(T^{-1}\right)^{\prime}\left(v\right)$
is often denoted by $\dfrac{\text{D}v}{\text{D}x}$ and is called
the \textbf{\index{Jacobian}\sindex[fam]{Jacobi, Carl}Jacobian\footnotemark}
of the variable change.

\end{remark}

\footnotetext{The Jacobian roots his name from \href{https://en.wikipedia.org/wiki/Carl_Gustav_Jacob_Jacobi}{Carl Jacobi}
(1804 - 1851) a Prussian mathematician, famous for his work on elliptic
integrals, on determinants, on equations with partial derivatives,
and their application to analytical mechanic.}

\section{Law of a Random Variable}\label{sec:Law-of-a-random-variable}

All random variables in this section are defined on the same probabilized
space $\left(\Omega,\mathcal{A},P\right).$ 

\begin{definition}{Random Variable (extended version)}{}

A \textbf{random variable\index{random variable}} $X$ taking values
in the probabilizable space $\left(E,\mathcal{E}\right)$ is, by definition,
a measurable function from $\left(\Omega,\mathcal{A}\right)$ to $\left(E,\mathcal{E}\right);$
that is, a function such that\boxeq{
\[
\forall B\in\mathcal{E},\,\,\,\,X^{-1}\left(B\right)\in\mathcal{A}.
\]
}

\end{definition}

\begin{definition}{Probability Law}{}

The \textbf{law}\index{law}---or also \textbf{probability law}\mindex{probability!law}---of
a random variable $X$ taking values in the probabilizable space $\left(E,\mathcal{E}\right)$
is the image measure $P_{X}$ of $P$ under $X.$

\end{definition}

To extend the concept of a cumulative distribution function to random
variables taking values in $\mathbb{R}^{d}$ for $d>1,$ we introduce
a partial order on $\mathbb{R}^{d}$ defined by
\[
x\leqslant y\,\,\,\,\Leftrightarrow\,\,\,\,\forall i\in\left\llbracket 1,d\right\rrbracket ,\,\,\,\,x_{i}\leqslant y_{i}.
\]

For $d=1,$ this coincides with the usual---total---order on $\mathbb{R}.$

\begin{definition}{Cumulative Distribution Function. Density}{}

Let $X$ be a random variable taking values in $\left(\mathbb{R}^{d},\mathcal{B}_{\mathbb{R}^{d}}\right).$

(a) The \textbf{cumulative distribution function\index{cumulative distribution function}}
of $X$ is the function $F_{X}$ from $\mathbb{R}^{d}$ to $\mathbb{R}^{+}$
defined by\boxeq{
\[
\forall x\in\mathbb{R}^{d},\,\,\,\,F_{X}\left(x\right)=P\left(X\leqslant x\right),
\]
} where $\leqslant$ denotes the usual partial order on $\mathbb{R}^{d}.$

(b) $X$ is said to admit the function $f$ as \textbf{a \index{density}density},
if its law $P_{X}$ admits $f$ as a density with respect to the Lebesgue
measure $\lambda_{d}$ on $\mathbb{R}^{d}.$

\end{definition}

If a random variable $X$ has a density $f,$ then any function $\lambda_{d}-$almost
everywhere equal to $f$ is also a density of $X.$ Conversely, any
density of $X$ is $\lambda_{d}-$almost everywhere equal to $f.$

Thus, the density of $X$ is defined up to $\lambda_{d}-$almost everywhere
equality, and is often identified with its equivalence class under
this relation, denoted by $f_{X}.$ The density of $X$ then satisfies
\[
\forall A\in\mathcal{B}_{\mathbb{R}^{d}},\,\,\,\,P_{X}\left(A\right)=\intop_{A}f_{X}\left(x\right)\text{d}\lambda_{d}\left(x\right).
\]

The uniqueness theorem for measures ensures that $X$ admits a density
if and only if there exists a nonnegative function $f_{X}$ of $\mathcal{L}^{1}\left(\mathbb{R}^{d},\text{\ensuremath{\mathcal{B}_{\mathbb{R}^{d}}}},\lambda_{d}\right)$
such that
\begin{equation}
\forall x\in\mathbb{R}^{d},\,\,\,\,F_{X}\left(x\right)=\intop_{\left\{ u\leqslant x\right\} }f_{X}\left(u\right)\text{d}\lambda_{d}\left(u\right).\label{eq:cns_density_existence}
\end{equation}

In particular, when $d=1,$ if there exists a nonnegative Riemann-integrable
function $f_{X}$ such that
\[
\forall x\in\mathbb{R},\,\,\,\,F_{X}\left(x\right)=\intop_{\left]-\infty,x\right]}f_{X}\left(u\right)\text{d}\lambda_{d}\left(u\right),
\]
by noticing that, in this case, this integral coincides with the Riemann
integral, we recover the basic definition of a density given previously
in Part \ref{part:Introduction-to-Probability} Chapter \ref{chap:Random-Variables-with}.

If the law $P_{X}$ is a measure with density with respect to the
counting measure\footnote{Tr.N.: Let $S$ be a set. For each set $E\subseteq S,$ define the
\textbf{\index{counting measure}counting measure} on $X$ as the
function $\mu:S\rightarrow\overline{\mathbb{N}}$ such that
\[
\mu\left(E\right)=\begin{cases}
\#E, & \text{if }E\text{ is finite},\\
\infty, & \text{if }E\text{ is infinite.}
\end{cases}
\]
} of $\mathbb{R}^{d},$ then the \textbf{random variable} is said to
be \textbf{discrete}\mindex{random variable ! discrete}---this definition
is a bit more general than the one given in Part \ref{part:Introduction-to-Probability}. 

In this case, the set of values
\[
\text{val}\left(X\right)=\left\{ x:\,P\left(X=x\right)\neq0\right\} 
\]
 is countable, and if $\delta_{x}$ denotes the Dirac measure in $x,$
then
\[
\forall A\in\mathcal{B}_{\mathbb{R}^{d}},\,\,\,\,P_{X}\left(A\right)=\sum_{x\in\text{val}\left(X\right)}P\left(X=x\right)\delta_{x}\left(A\right),
\]
which, in probabilistic notation, is written as\boxeq{
\[
P_{X}=\sum_{x\in\text{val}\left(X\right)}P\left(X=x\right)\delta_{x}.
\]
}

We recall that for such a discrete random variable $X,$ a function
$f$ belongs to $\mathscr{L}^{1}\left(\mathbb{R}^{d},\mathcal{B}_{\mathbb{R}^{d}},P_{X}\right)$
if and only if
\[
\sum_{x\in\text{val}\left(X\right)}\left|f\left(x\right)\right|P\left(X=x\right)<+\infty
\]
and in this case,
\[
\intop_{\mathbb{R}^{d}}f\text{d}P_{X}=\sum_{x\in\text{val}\left(X\right)}f\left(x\right)P\left(X=x\right).
\]

\begin{remark}{Important!}{}

If $X$ takes values in $\left(\mathbb{R}^{d},\mathcal{B}_{\mathbb{R}^{d}}\right),$
it follows from Corollary $\ref{co:equality_radon_measures}$ that
its law is fully determined by the family of integrals $\intop_{\Omega}f\left(X\right)\text{d}P\equiv\intop_{\mathbb{R}^{d}}f\text{d}P_{X}$
where $f$ is in $\mathbb{\mathscr{C}}^{+}_{\mathscr{K}}\left(\mathbb{R}^{d}\right).$
This provides an effective method for studying the law of a random
variable taking values in $\left(\mathbb{R}^{d},\mathcal{B}_{\mathbb{R}^{d}}\right),$
since all standard integration theorems can be applied without difficulty.

\end{remark}

\begin{example}{A Gamma Law Example}{}

Let $X$ be a real-valued random variable following the Gaussian law
$\mathcal{N}_{\mathbb{R}}\left(0,1\right),$ with density $f_{X}$
defined by
\[
\forall x\in\mathbb{R},\,\,\,\,f_{X}\left(x\right)=\dfrac{1}{\sqrt{2\pi}}\text{e}^{-\frac{x^{2}}{2}}.
\]

We aim to study the law of the random variable $X^{2}.$

\end{example}

\begin{solutionexample}{}{}

For every $f\in\mathscr{C}^{+}_{\mathscr{K}}\left(\mathbb{R}^{d}\right),$
we compute, using the transfer theorem\footnotemark and the one of
integration relatively to a measure with density
\begin{align*}
\intop_{\Omega}f\left(X^{2}\right)\text{d}P & =\intop_{\mathbb{R}}f\left(x^{2}\right)\text{d}P_{X}\left(x\right)\\
 & =\intop_{\mathbb{R}}f\left(x^{2}\right)\dfrac{1}{\sqrt{2\pi}}\text{e}^{-\frac{x^{2}}{2}}\text{d}\lambda\left(x\right)\\
 & =\dfrac{2}{\sqrt{2\pi}}\intop_{\mathbb{R}^{+\ast}}f\left(x^{2}\right)\text{e}^{-\frac{x^{2}}{2}}\text{d}\lambda\left(x\right).
\end{align*}

The change of variables using the $C^{1}-$diffeomorphism $T$ from
$\left]0,+\infty\right[$ onto itself, defined by $T\left(x\right)=x^{2}$
yields
\[
\intop_{\Omega}f\left(X^{2}\right)\text{d}P=\dfrac{1}{\sqrt{2\pi}}\intop_{\mathbb{R}^{+\ast}}f\left(y\right)\text{e}^{-\frac{y}{2}}y^{-\frac{1}{2}}\text{d}\lambda\left(y\right).
\]
Since, for every $y\in\left]0,+\infty\right[,$ $T^{-1}\left(y\right)=\sqrt{y}$
and $\left(T^{-1}\right)^{\prime}\left(y\right)=\dfrac{1}{2}y^{-\frac{1}{2}},$
it follows that
\[
\intop_{\Omega}f\left(X^{2}\right)\text{d}P=\intop_{\mathbb{R}}fg\text{d}\lambda,
\]
where the function $g$ is defined on the whole $\mathbb{R}$ by
\[
\forall y\in\mathbb{R},\,\,\,\,g\left(y\right)=\dfrac{1}{\sqrt{2\pi}}\boldsymbol{1}_{\mathbb{R}^{+\ast}}\text{e}^{-\frac{y}{2}}y^{-\frac{1}{2}}.
\]

This proves that the random variable admits $g$ as a density. The
probability law with density $g$ is called the \textbf{Gamma law\index{Gamma law}}
with parameters $\left(\dfrac{1}{2},\dfrac{1}{2}\right),$ and denoted
$\gamma\left(\dfrac{1}{2},\dfrac{1}{2}\right).$ We further talk again
about the Gamma laws, as well as in the classical laws summary at
the end of this chapter.

\end{solutionexample}

\footnotetext{Note that formally, it suffices to substitute the uppercase
denoting the random variable by the lowercase corresponding to the
random variable taken values.}

\begin{example}{Cauchy Law}{}

Let $U=\left(U_{1},U_{2}\right)$ be a random variable taking values
in $\mathbb{R}^{2}$ following the standard bivariate normal law $\mathcal{N}_{\mathbb{R}^{2}}\left(0,\boldsymbol{1}_{\mathbb{R}^{2}}\right),$
that is, with density $f_{U},$ the function defined on $\mathbb{R}^{2}$
by
\[
\forall u\in\mathbb{R}^{2},\,\,\,\,f_{U}\left(u\right)=\dfrac{1}{2\pi}\text{e}^{-\frac{\left\Vert u\right\Vert ^{2}}{2}},
\]
where $\left\Vert \cdot\right\Vert $ denotes the standard Euclidean
norm. 

Let $g$ be the function from $\mathbb{R}^{2}$ to $\mathbb{R}$ defined
by
\[
\forall u\in\mathbb{R}^{2},\,\,\,\,g\left(u_{1},u_{2}\right)=\begin{cases}
\dfrac{u_{1}}{u_{2}}, & \text{if }u_{2}\neq0,\\
0, & \text{if }u_{2}=0.
\end{cases}
\]

Let $X=g\left(U\right).$

We aim to study the law of the random variable $X.$

\end{example}

\begin{solutionexample}{}{}

For every $f\in\mathscr{C}^{+}_{\mathscr{K}}\left(\mathbb{R}^{2}\right),$
it follows from the transfer theorem and the one of integration with
respect to a measure with density that---by recalling that the Lebesgue
measure on a $\mathbb{R}^{2}$ line is zero---
\begin{align*}
\intop_{\Omega}f\left(X\right)\text{d}P & =\intop_{\mathbb{R}^{2}}\left(f\circ g\right)\left(x\right)\text{d}P_{X}\left(x\right)\\
 & =\intop_{\mathbb{R}^{2}\backslash\left\{ u_{2}=0\right\} }f\left(\dfrac{u_{1}}{u_{2}}\right)\dfrac{1}{2\pi}\text{e}^{-\frac{u^{2}_{1}+u^{2}_{2}}{2}}\text{d}\lambda_{2}\left(u_{1},u_{2}\right).
\end{align*}

Consider $T$ the diffeomorphism from $\mathbb{R}^{2}\backslash\left\{ u_{2}=0\right\} $
onto itself defined by
\[
T\left(u_{1},u_{2}\right)=\left(\dfrac{u_{1}}{u_{2}},u_{2}\right).
\]

Its invert $T^{-1}$ is determined by solving the system with variables
$u_{1}$ and $u_{2},$
\[
\left\{ \begin{array}{c}
x=\dfrac{u_{1}}{u_{2}}\\
y=u_{2}
\end{array}\right.\Leftrightarrow\left\{ \begin{array}{c}
u_{1}=xy\\
u_{2}=y
\end{array}\right.,
\]
which implies
\[
T^{-1}\left(x,y\right)=\left(xy,y\right).
\]

By doing the change of variables associated to $T,$ of Jacobian determinant
(or simply abusively called Jacobian), $\det\left(T^{-1}\right)^{\prime}\left(x,y\right),$
also denoted $\dfrac{D\left(u_{1},u_{2}\right)}{D\left(x,y\right)}.$
\[
\dfrac{D\left(u_{1},u_{2}\right)}{D\left(x,y\right)}=\left|\begin{array}{cc}
y & x\\
0 & 1
\end{array}\right|=y,
\]
it follows---without forgetting the absolute value of the Jacobian
determinant
\[
\intop_{\Omega}f\left(X\right)\text{d}P=\intop_{\mathbb{R}^{2}\backslash\left\{ y=0\right\} }f\left(x\right)\dfrac{1}{2\pi}\text{e}^{-\frac{\left(1+x^{2}\right)y^{2}}{2}}\left|y\right|\text{d}\lambda_{2}\left(x,y\right).
\]

Since $\left\{ y=0\right\} $ is a line and has Lebesgue measure zero,
we can extend the domain

\[
\intop_{\Omega}f\left(X\right)\text{d}P=\intop_{\mathbb{R}^{2}}f\left(x\right)\dfrac{1}{2\pi}\text{e}^{-\frac{\left(1+x^{2}\right)y^{2}}{2}}\left|y\right|\text{d}\lambda_{2}\left(x,y\right).
\]

By the Fubini theorem, which straightforwardly applies to nonnegative
measurable functions,
\begin{align*}
\intop_{\Omega}f\left(X\right)\text{d}P & =\intop_{\mathbb{R}}f\left(x\right)\left(\intop_{\mathbb{R}}\dfrac{1}{2\pi}\text{e}^{-\frac{\left(1+x^{2}\right)y^{2}}{2}}\left|y\right|\text{d}\lambda\left(y\right)\right)\text{d}\lambda\left(x\right)\\
 & =\intop_{\mathbb{R}}f\left(x\right)\dfrac{1}{\pi}\dfrac{1}{1+x^{2}}\text{d}\lambda\left(x\right).
\end{align*}

Thus, the random variable $X$ admits a density $f_{X}$ defined by\boxeq{
\[
\forall x\in\mathbb{R},\,\,\,\,f_{X}\left(x\right)=\dfrac{1}{\pi}\dfrac{1}{1+x^{2}}.
\]
}It means that $X$ follows the Cauchy law. 

\end{solutionexample}

\begin{remark}{}{}

The family of subsets of $\mathbb{R}^{d}$ of the form 
\[
\left\{ y\in\mathbb{R}^{d}:\,y\leqslant x\right\} 
\]
where $x\in\mathbb{R}^{d}$ is a $\pi-$system. 

Thus, if we know the cumulative distribution function of a random
variable $X$ taking values in $\mathbb{R}^{d},$ then we know its
law $P_{X}$ on this $\pi-$system, and hence, by the uniqueness theorem
for measures, we know its law entirely.

The cumulative distribution function is thus a useful tool for identifying
the law of a random variable. Its definition relies on an order structure:
so it is particularly suited to studying random variables through
operations that preserve or relate to this order, such as $\text{sup,}\,\inf,\,\max,\,\min.$

\end{remark}

The following proposition allows us to determine the law of a random
variable transformed by another among a diffeomorphism.

\begin{proposition}{Law of a Random Variable Transformed by Another via a Diffeomorphism}{diffeomorphism_trans_law}

Let $X$ be a random variable taking values in $\mathbb{R}^{d},$
and let $T$ be a diffeomorphism from $\mathbb{R}^{d}$ onto itself.

If $X$ admits a density $f_{X},$ then the random variable $Y=T\circ X,$
also denoted $T\left(X\right),$ admits a density $f_{Y}$ defined
by\boxeq{
\[
\forall y\in\mathbb{R}^{d},\,\,\,\,f_{Y}\left(y\right)=\left|\det\left(T^{-1}\right)^{\prime}\left(y\right)\right|f_{X}\left(T^{-1}\left(y\right)\right).
\]
}

\end{proposition}

\begin{proof}{}{}

For every fixed $f\in\mathscr{C}^{+}_{\mathscr{K}}\left(\mathbb{R}^{d}\right),$
by the transfer theorem and the theorem of integration related to
measures with density,
\begin{align*}
\intop_{\Omega}f\left(Y\right)\text{d}P & =\intop_{\mathbb{R}^{d}}\left(f\circ T\right)\left(x\right)\text{d}P_{X}\left(x\right)\\
 & =\intop_{\mathbb{R}^{d}}\left(f\circ T\right)\left(x\right)f_{X}\left(x\right)\text{d}\lambda_{d}\left(x\right).
\end{align*}

By doing the change of variables $y=T\left(x\right)$ from $\mathbb{R}^{d}$
onto itself, defined by the diffeomorphism $T,$
\[
\intop_{\Omega}f\left(Y\right)\text{d}P=\intop_{\mathbb{R}^{d}}f\left(y\right)\left|\det\left(T^{-1}\right)^{\prime}\left(y\right)\right|f_{X}\left(T^{-1}\left(y\right)\right)\text{d}\lambda_{d}\left(y\right).
\]

The result follows.

\end{proof}

\begin{definition}{Marginals of a Random Variable}{}

Let $X$ be a random variable taking values in 
\[
\mathbb{R}^{d}=\prod^{k}_{i=1}\mathbb{R}^{d_{i}}.
\]
 Let $\Pi_{i}$ be the canonical projection from $\mathbb{R}^{d}$
to $\mathbb{R}^{d_{i}}.$

The random variable $X_{i}=\Pi_{i}\circ X$ taking values in $\mathbb{R}^{d_{i}}$
is called the $i$-th \textbf{\mindex{random variable ! marginal}marginal\index{marginal}}
of $X.$

\end{definition}

The following propositions describe how to compute marginale laws.
To simplify the presentation, we state them for the special case $k=2;$
the general case is entirely analogous.

\begin{proposition}{Computation of Marginals Law}{comput_marginale_laws}

Let $X$ be a random variable taking values in $\mathbb{R}^{d}=\mathbb{R}^{d_{1}}\times\mathbb{R}^{d_{2}}.$

(a) If $X$ admits a density $f_{X},$ then the marginals $X_{1}$
and $X_{2}$ also admit densities $f_{X_{1}}$ and $f_{X_{2}},$ respectively,
given by\boxeq{
\begin{align*}
\forall x_{1}\in\mathbb{R}^{d_{1}},\,\,\,\, & f_{X_{1}}\left(x_{1}\right)=\intop_{\mathbb{R}^{d_{2}}}f_{X}\left(x_{1},x_{2}\right)\text{d}\lambda_{d_{2}}\left(x_{2}\right),\\
\forall x_{2}\in\mathbb{R}^{d_{2}},\,\,\,\, & f_{X_{2}}\left(x_{2}\right)=\intop_{\mathbb{R}^{d_{1}}}f_{X}\left(x_{1},x_{2}\right)\text{d}\lambda_{d_{1}}\left(x_{1}\right).
\end{align*}
}

(b) If $X$ is a discrete random variable, then $X_{1}$ and $X_{2}$
are also discrete, and\boxeq{
\begin{align*}
\forall x_{1}\in\text{val}\left(X_{1}\right),\,\,\,\, & P\left(X_{1}=x_{1}\right)=\sum_{x_{2}\in\text{val}\left(X_{2}\right)}P\left(X_{1}=x_{1},X_{2}=x_{2}\right),\\
\forall x_{2}\in\text{val}\left(X_{2}\right),\,\,\,\, & P\left(X_{2}=x_{2}\right)=\sum_{x_{1}\in\text{val}\left(X_{1}\right)}P\left(X_{1}=x_{1},X_{2}=x_{2}\right).
\end{align*}
}

\end{proposition}

\begin{proof}{}{}

(a) For every fixed $f\in\mathscr{C}^{+}_{\mathscr{K}}\left(\mathbb{R}^{d_{1}}\right),$
by the transfer theorem and by the theorem of integration with respect
to a measure with density,
\begin{align*}
\intop_{\Omega}f\left(X_{1}\right)\text{d}P & =\intop_{\mathbb{R}^{d}}\left(f\circ\Pi_{1}\right)\left(x\right)\text{d}P_{X}\left(x\right)\\
 & =\intop_{\mathbb{R}^{d}}\left(f\circ\Pi_{1}\right)\left(x_{1},x_{2}\right)f_{X}\left(x_{1},x_{2}\right)\text{d}\lambda_{d}\left(x_{1},x_{2}\right).
\end{align*}

By the Fubini theorem---applicable, since $f$ is nonnegative and
measurable---it follows that
\begin{align*}
\intop_{\Omega}f\left(X_{1}\right)\text{d}P & =\intop_{\mathbb{R}^{d_{1}}}f\left(x_{1}\right)\left(\intop_{\mathbb{R}^{d_{2}}}f_{X}\left(x_{1},x_{2}\right)\text{d}\lambda_{d_{2}}\left(x_{2}\right)\right)\text{d}\lambda_{d_{1}}\left(x_{1}\right).
\end{align*}

Thus, the enounced result.

(b) Recall that if $X$ is discrete, then $\text{val}\left(X\right)$
is countable, and that
\[
\forall x_{1}\in\text{val}\left(X_{1}\right),\,\,\,\,\left(X_{1}=x_{1}\right)\overset{P-a.e.}{=}\biguplus_{x_{2}\in\text{val}\left(X_{2}\right)}\left[\left(X_{1}=x_{1}\right)\cap\left(X_{2}=x_{2}\right)\right].
\]

Hence, the result by taking probabilities on both sides.

\end{proof}

\begin{remark}{}{}

Propositions $\ref{pr:diffeomorphism_trans_law}$ and $\ref{pr:comput_marginale_laws}$
are often used consecutively. The next example will illustrates this
combination perfectly.

\end{remark}

\begin{example}{}{}

Let $p>0$ be a real number.

Let $X=\left(X_{1},X_{2}\right)$ be a random variable taking values
in $\mathbb{R}^{2}$ with density $f_{X}$ defined by
\[
\forall\left(x_{1},x_{2}\right)\in\mathbb{R}^{2},\,\,\,\,f_{X}\left(x_{1},x_{2}\right)=p^{2}\boldsymbol{1}_{\mathbb{R}^{+}}\left(x_{1}\right)\boldsymbol{1}_{\mathbb{R}^{+}}\left(x_{2}\right)\text{e}^{-p\left(x_{1}+x_{2}\right)}.
\]

We seek the law of the random variable $Y=2X_{1}-X_{2}.$

\end{example}

\begin{solutionexample}{}{}

To find the law of the random variable $Y=2X_{1}-X_{2},$ we introduce
the random variable $\left(Y,X_{2}\right),$ obtained as a transformation
of $X$ by the diffeomorphism $T$ from $\mathbb{R}^{2}$ onto itself,
and defined by
\[
T\left(x_{1},x_{2}\right)=\left(2x_{1}-x_{2},x_{2}\right).
\]

Then $Y$ is its first marginal of $\left(Y,X_{2}\right)$. 

The random variable $\left(Y,X_{2}\right)$ admits the density $f_{\left(Y,X_{2}\right)}$
defined by
\[
\forall\left(y,x_{2}\right)\in\mathbb{R}^{2},\,\,\,\,f_{\left(Y,X_{2}\right)}\left(y,x_{2}\right)=\dfrac{1}{2}f_{X}\left(\dfrac{y+x_{2}}{2},x_{2}\right).
\]

Thus, the marginal $Y$ admits the density $f_{Y}$ defined, for every
$y\in\mathbb{R}$ by
\begin{align*}
f_{Y}\left(y\right) & =\intop_{\mathbb{R}}f_{\left(Y,X_{2}\right)}\left(y,x_{2}\right)\text{d}\lambda\left(x_{2}\right)\\
 & =\dfrac{p}{2}\text{e}^{-\frac{py}{2}}\intop_{\mathbb{R}}p\boldsymbol{1}_{\mathbb{R}^{+}}\left(y+x_{2}\right)\boldsymbol{1}_{\mathbb{R}^{+}}\left(x_{2}\right)\text{e}^{-\frac{3px_{2}}{2}}\text{d}\lambda\left(x_{2}\right)\\
 & =\dfrac{p}{2}\text{e}^{-\frac{py}{2}}\intop^{+\infty}_{\max\left(-y,0\right)}p\text{e}^{-\frac{3px_{2}}{2}}\text{d}x_{2}\\
 & =\dfrac{p}{3}\text{e}^{-\frac{py}{2}}\text{e}^{-\frac{3p\max\left(-y,0\right)}{2}}.
\end{align*}
Thus,
\[
\forall y\in\mathbb{R},\,\,\,\,f_{Y}\left(y\right)=\dfrac{p}{3}\left[\boldsymbol{1}_{\mathbb{R}^{+}}\left(y\right)\text{e}^{-\frac{py}{2}}+\boldsymbol{1}_{\mathbb{R}^{-\ast}}\left(y\right)\text{e}^{py}\right].
\]

\end{solutionexample}

\section{Moments of Random Variables}\label{sec:Moments-of-Random-Variables}

The moments of a random variable, when they exist, are parameters
that provide information about the law of this random variable. In
some cases, they can even determine the law entirely. Before defining
moments and studying their properties, we establish the \sindex[fam]{Hölder, Otto}\textbf{Hölder\index{Hölder inequality}}
and \textbf{\sindex[fam]{Minkowski, Hermann}Minkowski\index{Minkowski inequality}
inequalities}, from which we deduce fundamental properties of the
spaces $\mathscr{L}^{p}.$

\begin{definition}{Essentially Bounded Random Variable}{}

Let $p\geqslant1$ be a real number. Denote $\mathscr{L}^{p}\left(\Omega,\mathcal{A},P\right)$
the set of random variables $X$ defined $P-$almost everywhere, and
taking values in $\mathbb{R}$ or $\overline{\mathbb{R}},$ such that
\[
\intop_{\Omega}\left|X\right|^{p}\text{d}P<+\infty.
\]

For such a random variable $X,$ define the $p-$norm
\[
\left\Vert X\right\Vert _{p}=\left(\intop_{\Omega}\left|X\right|^{p}\text{d}P\right)^{1/p}.
\]

Denote $\mathscr{L}^{\infty}\left(\Omega,\mathcal{A},P\right)$ the
set of random variables $X,$ defined $P-$almost everywhere and taking
values in $\mathbb{R}$ or $\overline{\mathbb{R}},$ such that
\[
\left\{ x:\,P\left(\left|X\right|>x\right)>0\right\} <+\infty.
\]

For such a random variable $X,$ define the essential supremum
\[
\left\Vert X\right\Vert _{\infty}=\sup\left\{ x:\,P\left(\left|X\right|>x\right)>0\right\} =\inf\left\{ x:\,P\left(\left|X\right|>x\right)=0\right\} .
\]

In this case, we then say that $X$ is \textbf{essentially bounded\index{essentially bounded}---or
bounded $P-$almost everywhere}\mindex{P-almost everywhere bounded @ $P-$almost everywhere bounded}.

\end{definition}

\begin{remark}{}{}

If $X\in\mathscr{L}^{p}\left(\Omega,\mathcal{A},P\right),$ then $X$
is finite $P-$almost everywhere.

If $X\in\mathscr{L}^{\infty}\left(\Omega,\mathcal{A},P\right),$ then
$\left|X\right|\leqslant\left\Vert X\right\Vert _{\infty}$ $P-$almost
everywhere.

\end{remark}

We refer to Part \ref{part:Introduction-to-Probability} Chapter \ref{chap:Moments-of-a}
for the notion of conjugate real numbers, as introduced in Definition
$\ref{df:conjugate_real_numbers}$ and for Lemma $\ref{lm:lemma_pre_holder}$
that is instrumental in proving the following generalization of Proposition
$\ref{pr:holder_inequality}$

\begin{proposition}{Hölder Inequality}{}

Let $p$ and $q$ be two conjugate real numbers, finite or not.

(a) For any random variables $X$ and $Y$ taking values in $\overline{\mathbb{R}}^{+}$
and defined $P-$almost surely, the following inequality holds,
\begin{equation}
\intop_{\Omega}XY\text{d}P\leqslant\left(\intop_{\Omega}X^{p}\text{d}P\right)^{\frac{1}{p}}\left(\intop_{\Omega}Y^{q}\text{d}P\right)^{\frac{1}{q}}.\label{eq:holder_ineq_rv_r_bar}
\end{equation}

(b) If $X\in\mathscr{L}^{p}\left(\Omega,\mathcal{A},P\right)$ and
$Y\in\mathscr{L}^{q}\left(\Omega,\mathcal{A},P\right),$ then the
product $XY$ is integrable and the following inequality, called \textbf{Hölder}\footnotemark\textbf{\sindex[fam]{Hölder, Otto}
inequality}\index{Hölder inequality} holds,
\begin{equation}
\left\Vert XY\right\Vert _{1}\leqslant\left\Vert X\right\Vert _{p}\left\Vert Y\right\Vert _{q}.\label{eq:holder_ineq_for_lp_lq}
\end{equation}

In the special case where $p=q=2,$ this inequality is called the
Schwarz\footnotemark inequality, and becomes
\[
\left|\intop_{\Omega}XY\text{d}P\right|\leqslant\left(\intop_{\Omega}X^{2}\text{d}P\right)^{\frac{1}{2}}\left(\intop_{\Omega}Y^{2}\text{d}P\right)^{\frac{1}{2}}.
\]

\end{proposition}

\addtocounter{footnote}{-1}

\footnotetext{See Chapter \ref{chap:Moments-of-a} Footnote \ref{fn:Holder}}

\stepcounter{footnote}

\footnotetext{See Chapter \ref{chap:Moments-of-a} Footnote \ref{fn:schwartz}}

\begin{proof}{}{}

\textbf{(a) Case of nonnegative random variables}
\begin{itemize}
\item If $p$ and $q$ are finite:
\begin{itemize}
\item If one of the two terms on the right-hand side of inequality $\refpar{eq:holder_ineq_rv_r_bar}$
is zero, say the first, then $X=0$ $P-$almost everywhere. \\
Therefore,
\[
\intop_{\Omega}XY\text{d}P=0.
\]
\item If both terms on the right-hand side are nonzero, then it suffices
to show the inequality $\refpar{eq:holder_ineq_rv_r_bar}$ when the
two factors of the right-hand side member are finite. From Lemma $\ref{lm:lemma_pre_holder},$
\[
\dfrac{X}{\left\Vert X\right\Vert _{p}}\,\dfrac{Y}{\left\Vert Y\right\Vert _{q}}\leqslant\dfrac{1}{p}\dfrac{X^{p}}{\left\Vert X\right\Vert ^{p}_{p}}+\dfrac{1}{q}\dfrac{Y^{q}}{\left\Vert Y\right\Vert ^{q}_{q}},
\]
 which, by integrating both sides, shows the inequality $\refpar{eq:holder_ineq_rv_r_bar}.$
\end{itemize}
\item If $p=1$ and $q=+\infty,$ then we have $P-$almost surely
\[
0\leqslant Y\leqslant\left\Vert Y\right\Vert _{\infty}
\]
which, after integration yields $\refpar{eq:holder_ineq_rv_r_bar}.$
\end{itemize}
\textbf{(b) Case of random variables with arbitrary sign}

We apply inequality $\refpar{eq:holder_ineq_rv_r_bar}$ to the absolute
values.

\end{proof}

We deduce the \textbf{\sindex[fam]{Minkowski, Hermann}Minkowski inequality.}\index{Minkowski inequality}

\begin{proposition}{Minkowski Inequality}{minkowski_ineq_9}

Let $p>1$ be a real number, finite or not. 

(a) For any random variables $X$ and $Y$ taking values in $\overline{\mathbb{R}}^{+}$
and defined $P-$almost surely, the following inequality holds in
$\overline{\mathbb{R}}^{+},$
\begin{equation}
\left(\intop_{\Omega}\left(X+Y\right)^{p}\text{d}P\right)^{\frac{1}{p}}\leqslant\left(\intop_{\Omega}X^{p}\text{d}P\right)^{\frac{1}{p}}+\left(\intop_{\Omega}Y^{p}\text{d}P\right)^{\frac{1}{p}}.\label{eq:pre_minkovski_r_bar}
\end{equation}

(b) If $X\in\mathscr{L}^{p}\left(\Omega,\mathcal{A},P\right)$ and
$Y\in\mathscr{L}^{p}\left(\Omega,\mathcal{A},P\right),$ the sum $X+Y$
is in $\mathscr{L}^{p}\left(\Omega,\mathcal{A},P\right)$ and the
inequality called the \textbf{\sindex[fam]{Minkowski, Hermann}Minkowski\footnotemark
inequality}, holds\index{Minkowski inequality}
\begin{equation}
\left\Vert X+Y\right\Vert _{p}\leqslant\left\Vert X\right\Vert _{p}+\left\Vert Y\right\Vert _{p}.\label{eq:minkowski_inequality_lp}
\end{equation}

\end{proposition}

\footnotetext{See Chapter \ref{chap:Moments-of-a} Footnote \ref{fn:Minkowski}.}

\begin{proof}{}{}

\textbf{(a) Case of nonnegative random variables}
\begin{itemize}
\item If $p>1$ is finite, then by linearity,
\[
\intop_{\Omega}\left(X+Y\right)^{p}\text{d}P=\intop_{\Omega}\left(X+Y\right)^{p-1}X\text{d}P+\intop_{\Omega}\left(X+Y\right)^{p-1}Y\text{d}P.
\]
\\
By applying the Hölder inequality $\refpar{eq:holder_ineq_for_lp_lq}$
to each factor on the right-hand side yields
\begin{multline*}
\intop_{\Omega}\left(X+Y\right)^{p}\text{d}P\leqslant\left(\intop_{\Omega}\left(X+Y\right)^{q\left(p-1\right)}\text{d}P\right)^{1/q}\left(\intop_{\Omega}X^{p}\text{d}P\right)^{1/p}\\
+\left(\intop_{\Omega}\left(X+Y\right)^{q\left(p-1\right)}\text{d}P\right)^{1/q}\left(\intop_{\Omega}Y^{p}\text{d}P\right)^{1/p}.
\end{multline*}
Since $q\left(p-1\right)=p,$
\begin{multline*}
\intop_{\Omega}\left(X+Y\right)^{p}\text{d}P\leqslant\left(\intop_{\Omega}\left(X+Y\right)^{p}\text{d}P\right)^{1/q}\left[\left(\intop_{\Omega}X^{p}\text{d}P\right)^{1/p}+\left(\intop_{\Omega}Y^{p}\text{d}P\right)^{1/p}\right]
\end{multline*}
and by dividing both sides of the inequality by $\left(\intop_{\Omega}\left(X+Y\right)^{p}\text{d}P\right)^{1/q}=\left(\intop_{\Omega}\left(X+Y\right)^{p}\text{d}P\right)^{1-\frac{1}{p}}$
when this is nonzero (otherwise, this is immediate), we obtain
\[
\left(\intop_{\Omega}\left(X+Y\right)^{p}\text{d}P\right)^{\frac{1}{p}}\leqslant\left(\intop_{\Omega}X^{p}\text{d}P\right)^{\frac{1}{p}}+\left(\intop_{\Omega}Y^{p}\text{d}P\right)^{\frac{1}{p}}.
\]
\item If $p=+\infty,$ then
\[
\left|X+Y\right|\leqslant\left|X\right|+\left|Y\right|\leqslant\left\Vert X\right\Vert _{\infty}+\left\Vert Y\right\Vert _{\infty}\,\,\,\,P-\text{a.s.},
\]
and thus
\[
\left\Vert X+Y\right\Vert _{\infty}\leqslant\left\Vert X\right\Vert _{\infty}+\left\Vert Y\right\Vert _{\infty}.
\]
\end{itemize}
\textbf{(b) Case of random variables with arbitrary sign}

The random variables $X$ and $Y$ being defined and finite $P-$almost
surely, it is also the same for $X+Y.$ Then it suffices to apply
the inequality $\refpar{eq:minkowski_inequality_lp}$ to the absolute
values to obtain that 
\[
\left(\intop_{\Omega}\left|X+Y\right|^{p}\right)^{1/p}<+\infty,
\]
and the inequality follows directly.

\end{proof}

From these two inequalities, we deduce the properties of the sets
$\mathscr{L}^{p}\left(\Omega,\mathcal{A},P\right).$

\begin{proposition}{Properties of the Sets $\mathcal{L}^p\left(\Omega,\mathcal{A},P\right)$}{}

(a) If $p\geqslant1,$ $\mathscr{L}^{p}\left(\Omega,\mathcal{A},P\right)$
is a vector space equipped with a seminorm.

(b) If $p$ and $q$ are two integers such that
\[
1\leqslant p\leqslant q\leqslant+\infty.
\]

Then the inclusion of sets holds 
\[
\mathscr{L}^{q}\left(\Omega,\mathcal{A},P\right)\subset\mathscr{L}^{p}\left(\Omega,\mathcal{A},P\right)
\]
 as well as the seminorm inequality
\begin{equation}
\left\Vert X\right\Vert _{p}\leqslant\left\Vert X\right\Vert _{q}.\label{eq:inequality_seminorms}
\end{equation}

\end{proposition}

\begin{proof}{}{}

(a) This first point follows from the Minkowski inequality and the
homogeneity of the seminorm: for any real number\footnotemark $c,$
\[
\left\Vert cX\right\Vert _{p}=\left|c\right|\left\Vert X\right\Vert _{p}.
\]

(b) It suffices to consider the case where $p$ and $q$ are distinct.
Consider $X\in\mathscr{L}^{q}\left(\Omega,\mathcal{A},P\right),$
and let $r$ be the conjugate of $\dfrac{q}{p},$ that is
\[
r=\dfrac{q}{q-p}.
\]

By the Hölder inequality $\refpar{eq:holder_ineq_rv_r_bar},$ we obtain
\[
\intop_{\Omega}\left|X\right|^{p}\cdot1\text{d}P\leqslant\left(\intop_{\Omega}\left(\left|X\right|^{p}\right)^{\frac{q}{p}}\text{d}P\right)^{\frac{p}{q}}\left(\intop_{\Omega}1^{r}\text{d}P\right)^{\frac{1}{r}}=\left\Vert X\right\Vert ^{p}_{q}<+\infty.
\]

This shows that $X\in\mathscr{L}^{p}\left(\Omega,\mathcal{A},P\right),$
and also yields the inequality $\refpar{eq:inequality_seminorms}.$

\end{proof}

\footnotetext{We recall the convention $0\times\pm\infty=0.$}

\begin{remark}{}{}

We have thus proved that if $1\leqslant p\leqslant q\leqslant+\infty,$
then\footnotemark\boxeq{
\begin{equation}
\mathscr{L}^{\infty}\left(\Omega,\mathcal{A},P\right)\subset\mathscr{L}^{q}\left(\Omega,\mathcal{A},P\right)\subset\mathscr{L}^{p}\left(\Omega,\mathcal{A},P\right)\subset\mathscr{L}^{1}\left(\Omega,\mathcal{A},P\right)\label{eq:set_order_of_linf_lp_l1}
\end{equation}
}

If $p\geqslant1,$ the function $X\mapsto\left\Vert X\right\Vert _{p}$
defines a seminorm on $\mathscr{L}^{p}\left(\Omega,\mathcal{A},P\right),$
and $\left\Vert X\right\Vert _{p}=0$ if and only if $X=0,$ $P-$almost
surely. 

The quotient vector space of $\mathscr{L}^{p}\left(\Omega,\mathcal{A},P\right)$
by the equivalence relation of equality $P-$almost surely is denoted
$L^{P}\left(\Omega,\mathcal{A},P\right).$ It is a normed vector space,
with the norm obtained by taking the quotient of the seminorm $X\mapsto\left\Vert X\right\Vert _{p}$---commonly
referred to as the \textbf{\mindex{p-norm @ $p-$norm}$p-$norm}
of $X.$ 

It is standard practice to refer to a random variable and its equivalence
class using the same notation. We shall adopt this convention, and
likewise, we will use the same notation for the seminorm and its induced
quotient norm.

\end{remark}

\footnotetext{One have to remember that these inclusion relationships
fail when the measure is unbounded.}

\begin{definition}{Mean or Mathematical Expectation. Moment of Order $\alpha$. Variance. Standard Deviation}{}

Let $X$ be a real-valued random variable defined on a probabilized
space $\left(\Omega,\mathcal{A},P\right).$

(a) If $X\in\mathscr{L}^{1}\left(\Omega,\mathcal{A},P\right),$ the
quantity denoted $\mathbb{E}\left(X\right)$ or $\mathbb{E}X$ \boxeq{
\[
\mathbb{E}\left(X\right)=\intop_{\Omega}X\text{d}P
\]
}is called the \textbf{mean},\textbf{\index{mean}} \textbf{mathematical
expectation},\textbf{\index{mathematical expectation}} or simply
\textbf{expectation\index{expectation}} of $X.$

(b) If $X\in\mathscr{L}^{\alpha}\left(\Omega,\mathcal{A},P\right),$
for some $\alpha>0,$ the quantity, corresponding to the expectation
$\mathbb{E}\left(X^{\alpha}\right)$ of the random variable $X^{\alpha},$\boxeq{
\[
\mathbb{E}\left(X^{\alpha}\right)=\intop_{\Omega}X^{\alpha}\text{d}P
\]
}is called the \textbf{moment of order $\alpha$\mindex{moment ! order alpha@ order $\alpha$}}
of $X.$ In particular, if $\alpha>1$ and if $X\in\mathscr{L}^{\alpha}\left(\Omega,\mathcal{A},P\right),$
the quantity $\mathbb{E}\left(\left[X-\mathbb{E}\left(X\right)\right]^{\alpha}\right)$
is called \textbf{centered moment of order $\alpha$\mindex{centered moment of order alpha@centered moment of order $\alpha$}~\mindex{moment ! order alpha@ order $\alpha$ ! centered}}
of $X.$

(c) When $\alpha=2,$ the second-order centered moment is called the
\textbf{variance\index{variance}} of $X$ and is denoted $\sigma^{2}_{X}.$
Its nonnegative square root $\sigma_{X}$ is called the \textbf{\index{standard-deviation}standard-deviation}
of $X.$ 

\end{definition}

\begin{proposition}{Expectation is a Linear Form}{}

The mapping $\mathbb{E}$ is a linear---and continuous---form on
the vector space $\mathscr{L}^{1}\left(\Omega,\mathcal{A},P\right).$

\end{proposition}

\begin{proof}{}{}

This is a result from the integration theory\,---\,Proposition
$\ref{pr:integrable_vector_space}.$

\end{proof}

\begin{definition}{Associated Centered Random Variable. Reduced and Centered Random Variable}{}

If $X\in\mathscr{L}^{1}\left(\Omega,\mathcal{A},P\right),$ the random
variable \boxeq{
\[
\mathring{X}=X-\mathbb{E}\left(X\right)
\]
}is called the \textbf{centered random variable\mindex{random variable ! centered}}
associated with $X.$

If $X\in\mathscr{L}^{2}\left(\Omega,\mathcal{A},P\right),$ the random
variable \boxeq{
\[
\dfrac{X}{\sigma_{X}}\,\,\,\,\,\text{—}\text{respectively\,}\dfrac{\mathring{X}}{\sigma_{X}}\text{—}
\]
}is called the \textbf{reduced\mindex{random variable ! reduced}}---respectively,
\textbf{centered reduced\mindex{random variable ! centered reduced}}---\textbf{random
variable} associated with $X.$

\end{definition}

\begin{remarks}{}{}

1. These definitions do not provide any specific information about
the form of the law of the random variable $X.$ In particular, they
encompass the definitions introduced in previous chapters of Part
\ref{part:Introduction-to-Probability} for discrete and density-type
random variables. The interested reader may refer back to those chapters
for classical results and exercises related to such random variables.

2. As a consequence of the inclusion relations in $\refpar{eq:set_order_of_linf_lp_l1},$
if a random variable admits a moment of order $p\geqslant1,$ then
it also admits moments of every order between 1 and $p,$ inclusive.

\end{remarks}

\subsection{Moments Computation}

If $X\in\mathscr{L}^{\alpha}\left(\Omega,\mathcal{A},P\right),$ then
the transfer theorem gives\boxeq{
\[
\mathbb{E}\left(X^{\alpha}\right)=\intop_{\mathbb{R}}x^{\alpha}\text{d}P_{X}\left(x\right).
\]
}

In particular:
\begin{itemize}
\item If $X$ is a \textbf{discrete random variable}\mindex{random variable!discrete},
then since the law of $X$ satisfies
\[
P_{X}=\sum_{x\in\text{val}\left(X\right)}P\left(X=x\right)\delta_{x},
\]
it follows that\boxeq{
\[
\mathbb{E}\left(X^{\alpha}\right)=\sum_{x\in\text{val}\left(X\right)}x^{\alpha}P\left(X=x\right).
\]
}
\item If $X$ is a \textbf{random variable with density}\mindex{random variable!with density}
$f_{X},$ then since the law of $X$ satisfies $P_{X}=f_{X}\cdot\lambda,$
the integration theorem for a measure with density yields\boxeq{
\[
\mathbb{E}\left(X^{\alpha}\right)=\intop_{\mathbb{R}}x^{\alpha}f_{X}\left(x\right)\text{d}\lambda\left(x\right).
\]
}
\end{itemize}
\begin{proposition}{Properties of the Variance of a Random Variable}{}

If $X\in\mathscr{L}^{2}\left(\Omega,\mathcal{A},P\right),$ then its
variance satisfies\boxeq{
\[
\begin{array}{ccccccc}
\sigma^{2}_{X}=\mathbb{E}\left(X^{2}\right)-\left(\mathbb{E}\left(X\right)\right)^{2} & \,\,\,\, & \text{and} & \,\,\,\, & \forall\left(a,b\right)\in\mathbb{R}, & \,\,\,\, & \sigma^{2}_{aX+b}=a^{2}\sigma^{2}_{X}.\end{array}
\]
}

\end{proposition}

\begin{proof}{}{}

The results follows by expanding the square and applying the linearity
of the expectation. The proof is similar to that given in Part \ref{part:Introduction-to-Probability}.

\end{proof}

The standard deviation is thus invariant and nonnegatively homogeneous.

\begin{definition}{Covariance of Two Random Variables}{}

If $X$ and $Y$ belong to $\mathscr{L}^{2}\left(\Omega,\mathcal{A},P\right),$
then, by the Schwarz inequality, the random variable $\left(X-\mathbb{E}\left(X\right)\right)\left(Y-\mathbb{E}\left(Y\right)\right)$
belongs to $\mathscr{L}^{1}\left(\Omega,\mathcal{A},P\right).$ The
quantity, denoted $\text{cov}\left(X,Y\right),$ defined by\boxeq{
\[
\text{cov}\left(X,Y\right)=\mathbb{E}\left(\left(X-\mathbb{E}\left(X\right)\right)\left(Y-\mathbb{E}\left(Y\right)\right)\right)
\]
}is called \textbf{covariance\index{covariance}} of $X$ and $Y.$

\end{definition}

\begin{proposition}{Properties of the Covariance}{covariance_prop}

If $X$ and $Y$ belong to $\mathscr{L}^{2}\left(\Omega,\mathcal{A},P\right),$
then\boxeq{
\[
\begin{array}{ccccc}
\text{cov}\left(X,Y\right)=\mathbb{E}\left(XY\right)-\mathbb{E}\left(X\right)\mathbb{E}\left(Y\right) & \,\,\,\, & \text{and} & \,\,\,\, & \sigma^{2}_{X+Y}=\end{array}\sigma^{2}_{X}+\sigma^{2}_{Y}+2\text{cov}\left(X,Y\right).
\]
}

\end{proposition}

\begin{proof}{}{}

For the first equality, expand the product and use the linearity of
the expectation.

For the second equality, note that
\[
\sigma^{2}_{X+Y}=\mathbb{E}\left(\left(\mathring{X}+\mathring{Y}\right)^{2}\right).
\]

Then expand the square and use again the expectation linearity.

\end{proof}

We now extend these concepts to the case where the random variable
takes values in a finite-dimensional vector space. For such a space
$F,$ we denote by $F^{\ast}$ its \textbf{algebraic dual space\index{algebraic dual space}}
and $\left\langle \cdot,\cdot\right\rangle $ the associated \textbf{dual
bilinear form}\index{dual bilinear form}. In finite dimension, we
recall that all norms on $F$ are equivalent. We denote generically
by $\left\Vert \cdot\right\Vert $ one such norm. 

The vector space $F$ is equipped with its Borel $\sigma-$algebra,
generated by the open subsets of $F.$ In most applications, $F$
will be a Euclidean space, and the bilinear form $\left\langle \cdot,\cdot\right\rangle $
will coincide with the scalar product; in this case, $F$ is canonically
identified with its dual.

The presentation we adopt here is intended to define the moments in
an intrinsic manner. The reader may always replace $F$ by $\mathbb{R}^{d}$
equipped with the canonical scalar product if desired.

\begin{proposition}{Equivalence of Random Variable Norm and Bilinear Form Belonging to $\mathscr{L}^p\left(\Omega,\mathscr{A},P\right)$}{norm_and_bilinear_form_Lp}

Let $X$ be a random variable taking values in a finite-dimensional
vector space $F,$ and let $p\in\left[1,+\infty\right].$ The following
statements are equivalent:

(i) $\left\Vert X\right\Vert \in\mathscr{L}^{p}\left(\Omega,\mathcal{A},P\right)$

(ii) $\forall x^{\ast}\in F^{\ast},\,\,\,\,\left\langle X,x^{\ast}\right\rangle \in\mathscr{L}^{p}\left(\Omega,\mathcal{A},P\right).$

\end{proposition}

\begin{proof}{}{}

Equip $F^{\ast}$ with the norm defined by
\[
\left\Vert x^{\ast}\right\Vert =\sup_{\left\Vert x\right\Vert \leqslant1}\left\langle x,x^{\ast}\right\rangle .
\]

To prove that (i) implies (ii), observe that for every $x^{\ast}\in F^{\ast},$
\[
\left|\left\langle X,x^{\ast}\right\rangle \right|\leqslant\left\Vert X\right\Vert \left\Vert x^{\ast}\right\Vert 
\]
so if $\left\Vert X\right\Vert \in\mathscr{L}^{p}\left(\Omega,\mathcal{A},P\right),$
then $\left\langle X,x^{\ast}\right\rangle \in\mathscr{L}^{p}\left(\Omega,\mathcal{A},P\right)$
for every $x^{\ast}\in F^{\ast}.$

Conversely, suppose (ii) holds. Let $\left(e_{i}\right)_{1\leqslant i\leqslant d}$
be a basis of $F$ and let $\left(e^{\ast}_{i}\right)_{1\leqslant i\leqslant d}$
be its dual basis in $F^{\ast}.$ Since $X=\sum^{d}_{i=1}\left\langle X,e^{\ast}_{i}\right\rangle e_{i},$
\[
\left\Vert X\right\Vert \leqslant\sum^{d}_{i=1}\left|\left\langle X,e^{\ast}_{i}\right\rangle \right|\left\Vert e_{i}\right\Vert .
\]

Because the function $x\mapsto x^{p}$ is a nondecreasing function
on $\mathbb{R}^{+},$ applying the Minkowski inequality yields $\left\Vert X\right\Vert \in\mathscr{L}^{p}\left(\Omega,\mathcal{A},P\right).$

\end{proof}

\begin{denotation}{}{}

$\mathscr{L}^{p}_{F}\left(\Omega,\mathcal{A},P\right)$ is the vector
space of random variables taking values in $F$ such that $\left\Vert X\right\Vert \in\mathscr{L}^{p}\left(\Omega,\mathcal{A},P\right).$ 

\end{denotation}

\begin{definition}{Expectation of a Random Variable in $\mathcal{L}^1_F\left(\Omega,\mathcal{A},P\right)$}{}

Let $X\in\mathscr{L}^{1}_{F}\left(\Omega,\mathcal{A},P\right).$ The
linear form on $F^{\ast},$ denoted $\mathbb{E}\left(X\right),$ and
defined by \boxeq{
\[
\mathbb{E}\left(X\right):x^{\ast}\mapsto\intop_{\Omega}\left\langle X,x^{\ast}\right\rangle \text{d}P
\]
}is called the \textbf{expectation of $X.$\mindex{expectation!random variable}}

Identifying $F$ with its bidual---that is the algebraic dual\footnotemark
of $F^{\ast}$---the expectation $\mathbb{E}\left(X\right)$ is the
unique element of $F$ which satisfies\boxeq{
\begin{equation}
\forall x^{\ast}\in F^{\ast},\,\,\,\,\left\langle \mathbb{E}\left(X\right),x^{\ast}\right\rangle =\mathbb{E}\left(\left\langle X,x^{\ast}\right\rangle \right).\label{eq:expectation_and_duality}
\end{equation}
} 

\end{definition}

\footnotetext{The fact that $F$ is finite-dimensional ensures that
$F$ and its bidual $\left(F^{\ast}\right)^{\ast}$ are naturally
isomorph. This identification does not generally hold in infinite-dimensional
spaces, where the definition of the expectation as an element of $F$
may no longer be well-defined.}

\begin{remark}{}{}

1. We denote by the same symbol the operator expectation on both $\mathscr{L}^{1}\left(\Omega,\mathcal{A},P\right)$
and $\mathscr{L}^{1}_{F}\left(\Omega,\mathcal{A},P\right).$ When
$F=\mathbb{R},$ these two operators coincide. The operator $\mathbb{E}$
remains linear on $\mathscr{L}^{1}_{F}\left(\Omega,\mathcal{A},P\right).$

2. If $F=\mathbb{R}^{d},$ equipped with its canonical basis, then
by the relation $\refpar{eq:expectation_and_duality},$ the expectation
$\mathbb{E}\left(X\right)$ is the vector in $\mathbb{R}^{d}$ whose
$i-$th component is $\mathbb{E}\left(X_{i}\right).$

\end{remark}

\begin{proposition}{}{expectation_linearity}

Let $X\in\mathscr{L}^{1}_{F}\left(\Omega,\mathcal{A},P\right).$ Let
$A\in\mathscr{L}\left(F,G\right)$ and $b\in G,$ where $G$ is another
finite-dimensional vector space.

Then the random variable $AX+b$ belongs to $\mathscr{L}^{1}_{G}\left(\Omega,\mathcal{A},P\right),$
and\boxeq{
\[
\mathbb{E}\left(AX+b\right)=A\mathbb{E}\left(X\right)+b.
\]
}

\end{proposition}

\begin{proof}{}{}

We have
\[
\left\Vert AX+b\right\Vert \leqslant\left\Vert A\right\Vert \left\Vert X\right\Vert +\left\Vert b\right\Vert ,
\]
which shows that $AX+b\in\mathscr{L}^{1}_{G}\left(\Omega,\mathcal{A},P\right).$

Moreover, using the transpose $A^{\ast}$ of $A,$ we compute
\begin{align*}
\forall y^{\ast}\in G^{\ast},\,\,\,\,\left\langle \mathbb{E}\left(AX+b\right),y^{\ast}\right\rangle  & =\mathbb{E}\left(\left\langle AX+b,y^{\ast}\right\rangle \right)\\
 & =\mathbb{E}\left(\left\langle X,A^{\ast}y^{\ast}\right\rangle +\left\langle b,y^{\ast}\right\rangle \right)\\
 & =\left\langle \mathbb{E}\left(X\right),A^{\ast}y^{\ast}\right\rangle +\left\langle b,y^{\ast}\right\rangle \\
 & =\left\langle A\mathbb{E}\left(X\right)+b,y^{\ast}\right\rangle ,
\end{align*}
which proves the result.

\end{proof}

\begin{definition}{Variance of a Random Variable in $\mathcal{L}^2_F\left(\Omega,\mathcal{A},P\right)$. Auto-Covariance Operator. Matrix of Covariances}{}

Let $X\in\mathscr{L}^{2}_{F}\left(\Omega,\mathcal{A},P\right).$ The
nonnegative quadratic form on $F^{\ast},$ denoted $\sigma^{2}_{X}\left(\cdot\right),$\boxeq{
\[
\sigma^{2}_{X}:x^{\ast}\mapsto\intop_{\Omega}\left\langle X-\mathbb{E}\left(X\right),x^{\ast}\right\rangle ^{2}\text{d}P
\]
}is called \textbf{variance\index{variance}} of $X.$

The variance of $X$ is bijectively mapped to the nonnegative self-adjoint
linear operator $\Lambda_{X}\in\mathscr{L}\left(F^{\ast},F\right)$
by the relationship\boxeq{
\[
\forall x^{\ast}\in F^{\ast},\,\,\,\,\left\langle \Lambda_{X}x^{\ast},x^{\ast}\right\rangle =\sigma^{2}_{X}\left(x^{\ast}\right).
\]
}

This operator is called \textbf{\index{auto-covariance operator}auto-covariance
operator} of $X.$ 

Let $\left(e_{i}\right)_{1\leqslant i\leqslant d}$ be a basis of
$F$ and let $\left(e^{\ast}_{i}\right)_{1\leqslant i\leqslant d}$
be its dual basis in $F^{\ast}.$ The matrix representation $C_{X}$
of the operator $\Lambda_{X}$ in this basis is called the covariance
matrix of $X.$

This matrix satisfies\boxeq{
\[
\forall i,j\in\left\llbracket 1,d\right\rrbracket ,\,\,\,\,\left(C_{X}\right)_{ij}=\text{cov}\left(\left\langle X,e^{\ast}_{i}\right\rangle ,\left\langle X,e^{\ast}_{j}\right\rangle \right).
\]
}

\end{definition}

\begin{remark}{}{}

When $F=\mathbb{R}^{d},$ equipped with its canonical basis, the matrix
$C_{X}$ is the nonnegative symmetric matrix of size $d\times d:$
\[
C_{X}=\left(\begin{array}{ccccccc}
\sigma^{2}_{X_{1}}\\
 & \ddots\\
 &  & \sigma^{2}_{X_{i}} & \dots & \text{cov}\left(X_{j},X_{i}\right)\\
 &  & \vdots & \ddots & \vdots\\
 &  & \text{cov}\left(X_{i},X_{j}\right) &  & \sigma^{2}_{X_{j}}\\
 &  &  &  &  & \ddots\\
 &  &  &  &  &  & \sigma^{2}_{X_{d}}
\end{array}\right)
\]

\end{remark}

\begin{proposition}{}{}

Let $F$ and $G$ be two finite-dimensional vector spaces. Let $X\in\mathscr{L}^{2}_{F}\left(\Omega,\mathcal{A},P\right).$
Let $A\in\mathscr{L}\left(F,G\right)$ and $b\in G.$

Then, the random variable $AX\in\mathscr{L}^{2}_{G}\left(\Omega,\mathcal{A},P\right)$
and\boxeq{
\[
\Lambda_{AX+b}=A\Lambda_{X}A^{\ast},
\]
}which, in terms of covariance matrices, reads\boxeq{
\[
C_{AX+b}=AC_{X}A^{\ast}.
\]
}

\end{proposition}

\begin{proof}{}{}

We have
\[
\left\Vert AX+b\right\Vert ^{2}\leqslant\dfrac{1}{2}\left(\left\Vert A\right\Vert ^{2}\left\Vert X\right\Vert ^{2}+\left\Vert b\right\Vert ^{2}\right).
\]
which proves that $AX+b\in\mathscr{L}^{2}_{G}\left(\Omega,\mathcal{A},P\right).$

Moreover, by the definition of the transpose of $A$ and Proposition
$\ref{pr:expectation_linearity},$ it holds
\begin{align*}
\forall y^{\ast}\in G^{\ast},\,\,\,\,\left\langle \Lambda_{AX+b}y^{\ast},y^{\ast}\right\rangle  & =\mathbb{E}\left(\left\langle A\mathring{X},y^{\ast}\right\rangle ^{2}\right)\\
 & =\mathbb{E}\left(\left\langle \mathring{X},A^{\ast}y^{\ast}\right\rangle ^{2}\right)\\
 & =\left\langle \Lambda_{X}A^{\ast}y^{\ast},A^{\ast}y^{\ast}\right\rangle \\
 & =\left\langle A\Lambda_{X}A^{\ast}y^{\ast},y^{\ast}\right\rangle .
\end{align*}

By bilinearization of the quadratic form, it follows that
\[
\forall y^{\ast}\in G^{\ast},\,\,\,\,\left\langle \Lambda_{AX+b}x^{\ast},y^{\ast}\right\rangle =\left\langle A\Lambda_{X}A^{\ast}x^{\ast},y^{\ast}\right\rangle ,
\]
which completes the proof.

\end{proof}

\begin{figure}[t]
\begin{center}\includegraphics[width=0.4\textwidth]{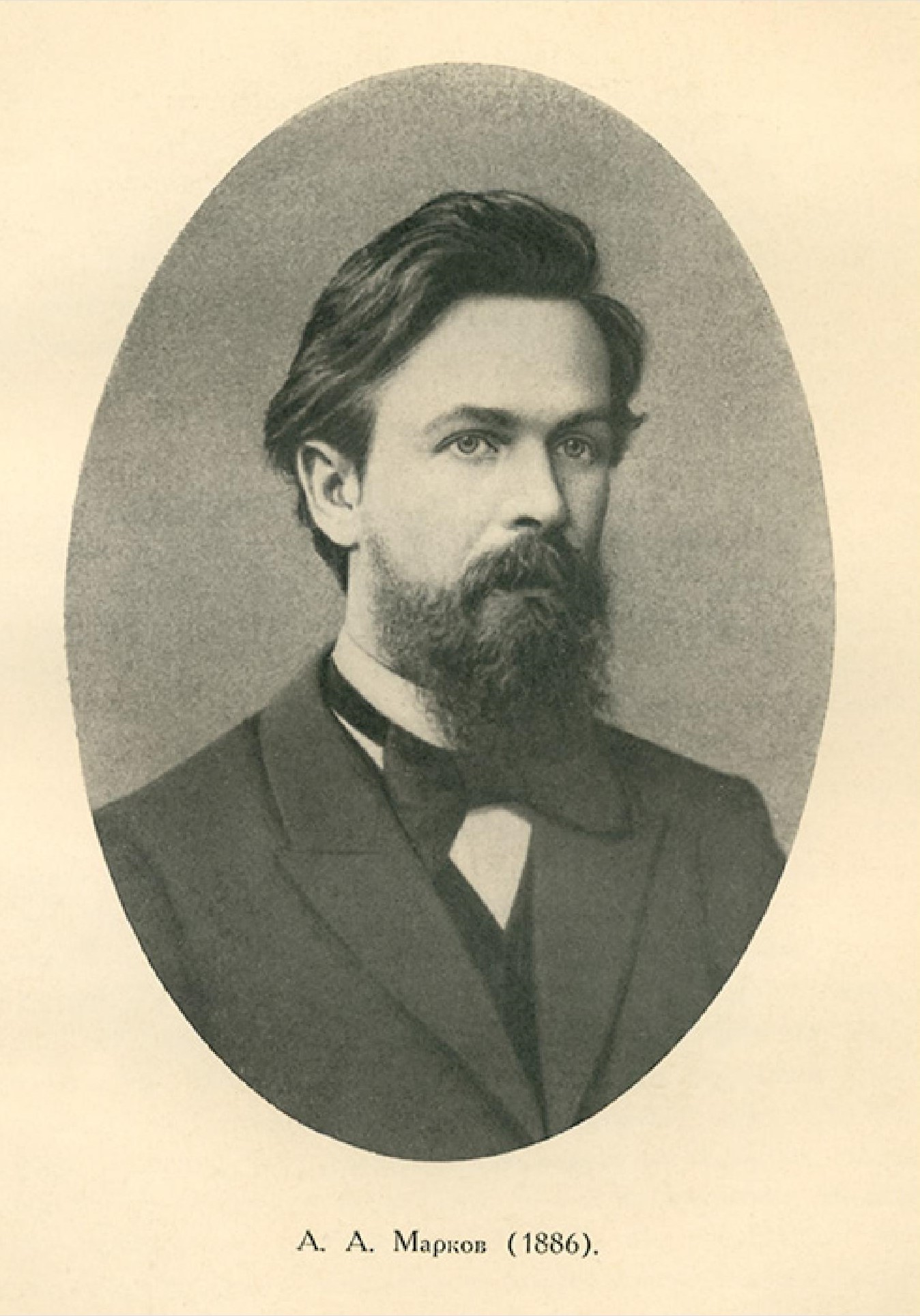}

{\tiny Credits: Public Domain}\end{center}

\caption{\textbf{\protect\href{https://en.wikipedia.org/wiki/Andrey_Markov}{Andreï Markov}}
(1856-1922)}\sindex[fam]{Markov, Andreï}
\end{figure}

We now present two classical inequalities---the \textbf{Markov}\footnote{Andreï Markov (1856 - 1922) was a Russian mathematician, famous for
his work on stochastic processes and what became known as Markov chain.
His brother Vladimir and his son, also named Andreï, were two mathematicians.}\sindex[fam]{Markov, Andreï}\textbf{ inequality}\index{Markov inequality}
and the \textbf{Bienaymé\sindex[fam]{Bienaymé, Irénée-Jules}-Chebyshev}\sindex[fam]{Chebyshev, Pafnuty}\textbf{
inequality}\index{Bienaymé-Chebyshev inequality}. Though rather
coarse, they nevertheless provide basic information on the concentration
of the values taken by a random variable, particularly around its
expectation. Their numerical looseness is not surprising upon examining
their proofs. These inequalities are mainly useful for proving convergence
in probability---see Chapter \ref{chap:PartIIChap11}.

\begin{proposition}{Markov Inequality}{markov_ineq}

Let $X\in\mathscr{L}^{1}\left(\Omega,\mathcal{A},P\right)$ be a nonnegative
random variable. Then, for every $\epsilon>0,$\boxeq{
\[
P\left(X\geqslant\epsilon\right)\leqslant\dfrac{\mathbb{E}\left(X\right)}{\epsilon},
\]
}and in particular\boxeq{
\[
P\left(X>\epsilon\right)\leqslant\dfrac{\mathbb{E}\left(X\right)}{\epsilon}.
\]
}

As a consequence, if $X\in\mathscr{L}^{1}_{F}\left(\Omega,\mathcal{A},P\right),$
then for every $\epsilon>0,$\boxeq{
\[
P\left(\left\Vert X\right\Vert >\epsilon\right)\leqslant\dfrac{\mathbb{E}\left(\left\Vert X\right\Vert \right)}{\epsilon}.
\]
}

\end{proposition}

\begin{proof}{}{}

Let $X\in\mathscr{L}^{1}\left(\Omega,\mathcal{A},P\right)$ be a nonnegative
random variable. Let $D=\left\{ X\geqslant\epsilon\right\} $ for
some fixed $\epsilon>0.$

Since $X\geqslant\epsilon$ on $D,$ the following chain of inequalities
holds:
\[
\mathbb{E}\left(X\right)=\intop_{\Omega}X\text{d}P\geqslant\intop_{D}X\text{d}P\geqslant\epsilon P\left(D\right),
\]
hence, the first inequality.

The second inequality follows from the inclusion $\left(X>\epsilon\right)\subset\left(X\geqslant\epsilon\right).$

When $X\in\mathscr{L}^{1}_{F}\left(\Omega,\mathcal{A},P\right),$
the second inequality is obtained by applying the first inequality
to the nonnegative random variable $\left\Vert X\right\Vert .$ 

\end{proof}

\begin{proposition}{Bienaymé-Tchebitchev Inequality}{}

Let $F$ be an Euclidean space, and consider $X\in\mathscr{L}^{2}_{F}\left(\Omega,\mathcal{A},P\right).$

Then, for every $\epsilon>0,$
\[
P\left(\left\Vert X-\mathbb{E}\left(X\right)\right\Vert >\epsilon\right)\leqslant\dfrac{\text{tr}\left(\Lambda_{X}\right)}{\epsilon^{2}}.
\]

In particular, if $F=\mathbb{R},$ this inequality becomes
\[
P\left(\left\Vert X-\mathbb{E}\left(X\right)\right\Vert >\epsilon\right)\leqslant\dfrac{\sigma^{2}_{X}}{\epsilon^{2}}.
\]

\end{proposition}

\begin{proof}{}{}

It suffices to apply the Markov inequality to the nonnegative random
variable $\left\Vert X-\mathbb{E}\left(X\right)\right\Vert ^{2},$
with the real number $\epsilon^{2},$ noting that
\[
\left(\left\Vert X-\mathbb{E}\left(X\right)\right\Vert >\epsilon\right)=\left(\left\Vert X-\mathbb{E}\left(X\right)\right\Vert ^{2}>\epsilon^{2}\right)
\]
and that
\[
\mathbb{E}\left(\left\Vert \mathring{X}\right\Vert ^{2}\right)=\text{tr}\left(\Lambda_{X}\right).
\]

In particular, if $F=\mathbb{R},$ then
\[
\mathbb{E}\left(\mathring{X}^{2}\right)=\sigma^{2}_{X}.
\]

\end{proof}

\begin{figure}[t]
\begin{center}\includegraphics[width=0.4\textwidth]{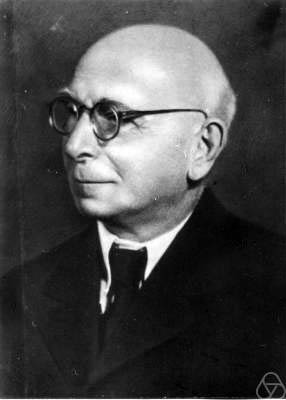}

{\tiny Credits: \href{https://commons.wikimedia.org/wiki/File:Sergei_Natanowitsch_Bernstein.jpg}{Konrad Jacobs, Erlangen}.
Free to use}\end{center}

\caption{\textbf{\protect\href{https://en.wikipedia.org/wiki/Sergei_Bernstein}{Sergei Bernstein}}
(1880-1968)}\sindex[fam]{Bernstein, Sergei}
\end{figure}

The \textbf{Markov}\sindex[fam]{Markov, Andreï} inequality has for
consequence the following inequalities, from which the \textbf{Bernstein}{\bfseries\footnote{Sergei Berstein (1880 - 1968) was a Ukrainian and Soviet mathematician,
known for his contributions in partial differential equations, differential
geometry, probability theory and approximation theory.}}\textbf{\sindex[fam]{Bernstein, Sergei} inequalities}\index{Bernstein inequalities}
can be derived. These improve upon the Chebyshev inequalities, and
form the starting point of the large deviations theory.

\begin{proposition}{}{}

Let $f$ be a nondecreasing and positive function from $\mathbb{R}$
to $\mathbb{R}.$ Let $X$ be a real-valued random variable such that
$f\circ X\in\mathscr{L}^{1}\left(\Omega,\mathcal{A},P\right).$

Then, for any real number $\epsilon,$\boxeq{
\[
P\left(X>\epsilon\right)\leqslant\dfrac{\mathbb{E}\left(f\left(X\right)\right)}{f\left(\epsilon\right)}.
\]
}

In particular, if $X$ is a real-valued random variable such that
for every $a>0,$ we have $\exp\left(aX\right)\in\mathscr{L}^{1}\left(\Omega,\mathcal{A},P\right),$
then\boxeq{
\[
\forall\epsilon>0,\,\,\,\,P\left(X>\epsilon\right)\leqslant\text{e}^{-a\epsilon}\mathbb{E}\left(\text{e}^{aX}\right).
\]
}

\end{proposition}

\begin{proof}{}{}

Since $f$ is nondecreasing,
\[
\left(X>\epsilon\right)\subset\left(f\left(X\right)\geqslant f\left(\epsilon\right)\right).
\]

And thus
\[
P\left(X>\epsilon\right)\leqslant P\left(f\left(X\right)\geqslant f\left(\epsilon\right)\right).
\]

By the Markov inequality applied to the random variable $f\left(X\right)$
and to the positive real number $f\left(\epsilon\right),$ we conclude.

The second inequality is obtained from the first by taking $f\left(x\right)=\text{e}^{ax}.$

\end{proof}

We introduce the concept of \textbf{\index{correlation coefficient}correlation
coefficient}, which, as we will see later, serves to quantify a certain
degree of link between random variables. 

\begin{definition}{Correlation Coefficient of Two Random Variables}{}

Let $X$ and $Y\in\mathscr{L}^{2}\left(\Omega,\mathcal{A},P\right)$
be two random variables with nonzero variance.

The \textbf{correlation coefficient\index{correlation coefficient}}
of $X$ and $Y$ is the real number denoted $\rho_{X,Y}$ and defined
by\boxeq{
\[
\rho_{X,Y}=\dfrac{\text{cov}\left(X,Y\right)}{\sigma_{X}\sigma_{Y}}.
\]
}

\end{definition}

\begin{proposition}{Characterization of Perfect Correlation}{}

Let $X$ and $Y\in\mathscr{L}^{2}\left(\Omega,\mathcal{A},P\right)$
be two random variables with nonzero variance. Then their correlation
coefficient $\rho_{X,Y}$ satisfies the following properties:

(a) The inequality\boxeq{
\[
\left|\rho_{X,Y}\right|\leqslant1.
\]
}

(b) The equality $\left|\rho_{X,Y}\right|=1$ occurs if and only if
there exists three real numbers $a,b,c,$ not all zero, such that\boxeq{
\[
P\left(aX+bY+c=0\right)=1.
\]
}

\end{proposition}

\begin{proof}{}{}

(a) The Schwarz inequality gives
\[
\left|\mathbb{E}\left(\mathring{X}\mathring{Y}\right)\right|\leqslant\mathbb{E}\left(\left|\mathring{X}\mathring{Y}\right|\right)\leqslant\left[\mathbb{E}\left(\mathring{X}^{2}\right)\right]^{\frac{1}{2}}\left[\mathbb{E}\left(Y^{2}\right)\right]^{\frac{1}{2}}=\sigma_{X}\sigma_{Y},
\]
which proves the desired inequality.

(b) Suppose $\left|\rho_{X,Y}\right|=1.$ Then the second-degree polynomial
in $\lambda,$ $\mathbb{E}\left(\left(\mathring{X}+\lambda\mathring{Y}\right)^{2}\right),$
has its reduced discriminant equals to zero, and admits a double root
$\lambda_{0}.$ It follows
\[
\mathbb{E}\left(\left(\mathring{X}+\lambda\mathring{Y}\right)^{2}\right)=0.
\]
and hence
\[
P\left(\mathring{X}+\lambda\mathring{Y}=0\right)=1.
\]

Conversely, suppose there exist three real numbers $a,b$ and $c,$
not all zero, such that
\begin{equation}
P\left(aX+bY+c=0\right)=1.\label{eq:P_linear}
\end{equation}

If $c\neq0,$ then necessarily $a$ and $b$ are different from zero.
Indeed, if for instance $a=0,$ then we would have $P\left(bY+c=0\right)=1,$
so that $\sigma^{2}_{bY+c}=\sigma^{2}_{0}=0,$ which implies $b^{2}\sigma^{2}_{Y}=0,$
and thus $b=0,$ contradicting $\refpar{eq:P_linear}$. 

In that case, we can write
\begin{equation}
P\left(X=\alpha Y+\beta\right)=1,\label{eq:P_linear_rel_2}
\end{equation}
for some $\alpha\neq0.$

If $c=0,$ then either $a\neq0$ or $b\neq0.$ If for instance $a\neq0,$
the equality $\refpar{eq:P_linear_rel_2}$ still holds with $\beta=0$---otherwise,
if it was $b$ that was different from 0, a similar computation would
be done.

Anyway, in both cases, we obtain
\[
\text{cov}\left(X,Y\right)=\mathbb{E}\left(\left(\alpha\mathring{Y}\right)\mathring{Y}\right)=\alpha\sigma^{2}_{Y}\,\,\,\,\text{and}\,\,\,\,\sigma^{2}_{X}=\sigma^{2}_{\alpha Y+\beta}=\alpha^{2}\sigma^{2}_{Y}.
\]

This yields
\[
\rho_{X,Y}=\dfrac{\alpha\sigma^{2}_{Y}}{\left|\alpha\right|\sigma^{2}_{Y}}
\]
and thus
\[
\left|\rho_{X,Y}\right|=1.
\]

\end{proof}

\subsection{Linear Regression Problem}

Let $X$ and $Y\in\mathscr{L}^{2}\left(\Omega,\mathcal{A},P\right)$
be two random variables. We seek the ``best'' approximation of $Y$
as a linear function of $X,$ in the least squares sense. That is,
an element $\left(\widehat{a},\widehat{b}_{\widehat{a}}\right)\in\mathbb{R}^{2}$
of the following set corresponding to the minimization problem
\[
\inf\left\{ \phi\left(a,b\right):\,\left(a,b\right)\in\mathbb{R}^{2}\right\} ,
\]
where
\[
\phi\left(a,b\right)=\mathbb{E}\left(\left[Y-\left(aX+b\right)\right]^{2}\right).
\]

This problem is known as the \textbf{linear regression problem}\index{linear regression problem}.

The analysis carried out in Part \ref{part:Introduction-to-Probability}
Subsection \ref{subsec:The-problem-of-linear-regression} applies
identically here in a more general setting. We recall the essential
steps.

We first write
\begin{align*}
\phi\left(a,b\right) & =\mathbb{E}\left(\left[\mathring{Y}-a\mathring{X}+\left(\mathbb{E}\left(Y\right)-a\mathbb{E}\left(X\right)-b\right)\right]^{2}\right)\\
 & =\mathbb{E}\left(\left[\mathring{Y}-a\mathring{X}\right]^{2}\right)+\left[\mathbb{E}\left(Y\right)-a\mathbb{E}\left(X\right)-b\right]^{2}.
\end{align*}

For a fixed $a,$ the function $\phi\left(a,b\right)$ is minimized
in $b$ when 
\[
\hat{b_{a}}=\mathbb{E}\left(Y\right)-a\mathbb{E}\left(X\right),
\]
 i.e., when $b$ is solution of
\[
\mathbb{E}\left(Y\right)-a\mathbb{E}\left(X\right)-b=0
\]

It remains to minimize in $a,$ the polynomial
\begin{align*}
f(a) & =\phi\left(a,\hat{b}_{a}\right)=\mathbb{E}\left(\left[\mathring{Y}-a\mathring{X}\right]^{2}\right)\\
 & =\sigma^{2}_{Y}-2a\text{cov}\left(X,Y\right)+a^{2}\sigma^{2}_{X}.
\end{align*}

Taking the derivative, we obtain
\[
f^{\prime}\left(a\right)=2a\sigma^{2}_{X}-2\text{cov}\left(X,Y\right).
\]
and
\[
f'\left(a\right)\geqslant0\Leftrightarrow a\geqslant\dfrac{\text{cov\ensuremath{\left(X,Y\right)}}}{\sigma^{2}_{X}}.
\]

Define
\[
\hat{a}=\dfrac{\text{cov\ensuremath{\left(X,Y\right)}}}{\sigma^{2}_{X}}
\]

Then $f$ reaches a minimum at $\hat{a}.$ The solution to the linear
regression problem is thus the pair $\left(\hat{a},\hat{b}_{\hat{a}}\right),$
given by:
\[
\begin{cases}
\hat{a}=\rho_{X,Y}\dfrac{\sigma_{Y}}{\sigma_{X}}\\
\hat{b}_{\hat{a}}=\mathbb{E}\left(Y\right)-\hat{a}\mathbb{E}\left(X\right)=\mathbb{E}\left(Y\right)-\mathbb{E}\left(X\right)\cdot\rho_{X,Y}\dfrac{\sigma_{Y}}{\sigma_{X}}.
\end{cases}
\]

The line $D$ with equation\boxeq{
\[
\left(y-\mathbb{E}\left(Y\right)\right)-\rho_{X,Y}\dfrac{\sigma_{Y}}{\sigma_{X}}\left(x-\mathbb{E}\left(X\right)\right)=0
\]
}is called the \textbf{linear regression line of $Y$ in $X.$\index{linear regression line of $Y$ in $X$}}
The best approximation $\widetilde{Y}$ of $Y$ as a linear function
of $X,$ in the least squares sense, is\boxeq{
\[
\widetilde{Y}=\mathbb{E}\left(Y\right)+\rho_{X,Y}\dfrac{\sigma_{Y}}{\sigma_{X}}\left(X-\mathbb{E}\left(X\right)\right).
\]
}

Moreover,
\[
P\left(\left(X,Y\right)\in D\right)=1
\]
if and only if
\[
\phi\left(\hat{a},\hat{b}_{\hat{a}}\right)=0.
\]

\textbf{Special Case}

If the random variable follows the uniform law on a finite set of
$n$ points in the plane $\left\{ \left(x_{i},y_{i}\right):\,i\in\left\llbracket 1,n\right\rrbracket \right\} ,$
then the quantity $\phi\left(a,b\right)$ becomes\boxeq{
\[
\phi\left(a,b\right)=\dfrac{1}{n}\sum^{n}_{i=1}\left(y_{i}-\left(ax_{i}+b\right)\right)^{2}.
\]
}

We thus recover the \textbf{least squares approximation line} commonly
used in physics and data fitting. The reader is invited to determine
the equation of this line as an exercise.

\subsection{Most Usual Laws}

We present here the most commonly encountered probability laws $\mu,$
along with their expectation $m$ and variance $\sigma^{2}$ when
these moments exist. We also provide their Fourier transform $\hat{\mu}$---whose
precise definition will be introduced in Chapter \ref{chap:PartIIChap13}.

\subsubsection{Discrete Laws}

\begin{tabular}{|c|c|c|c|c|}
\hline 
\begin{cellvarwidth}[t]
\centering
Name of the law

(parameters)
\end{cellvarwidth} &
\begin{cellvarwidth}[t]
\centering
Probability measure

$\mu$
\end{cellvarwidth} &
\begin{cellvarwidth}[t]
\centering
$\hat{\mu}\left(t\right)$

$\left(t\in\mathbb{R}\right)$
\end{cellvarwidth} &
$m$ &
$\sigma^{2}$\tabularnewline
\hline 
\hline 
\begin{cellvarwidth}[t]
\centering
Bernoulli law $\mathcal{B}\left(1,p\right)$

$\left(0<p<1,\,q=1-p\right)$
\end{cellvarwidth} &
$p\delta_{1}+q\delta_{0}$ &
$p\text{e}^{\text{i}t}+q$ &
$p$ &
$pq$\tabularnewline
\hline 
\begin{cellvarwidth}[t]
\centering
Binomial law $\mathcal{B}\left(n,p\right)$

$\left(0<p<1,\,n\in\mathbb{N}^{\ast}\right)$
\end{cellvarwidth} &
$\sum\limits^{n}_{k=0}\binom{n}{k}p^{k}q^{n-k}\delta_{k}$ &
$\left(p\text{e}^{\text{i}t}+q\right)^{n}$ &
$np$ &
$npq$\tabularnewline
\hline 
\begin{cellvarwidth}[t]
\centering
Poisson law $\mathcal{P}\left(\lambda\right)$

$\left(\lambda>0\right)$
\end{cellvarwidth} &
$\sum\limits^{+\infty}_{n=0}\dfrac{\lambda^{n}\text{e}^{-\lambda}}{n!}\delta_{n}$ &
$\text{e}^{\lambda\left(\text{e}^{\text{i}t}-1\right)}$ &
$\lambda$ &
$\lambda$\tabularnewline
\hline 
\begin{cellvarwidth}[t]
\centering
Geometric law on $\mathbb{N},\,$$\mathcal{G}_{\mathbb{N}}\left(p\right)$

$\left(0<p<1\right)$
\end{cellvarwidth} &
$\sum\limits^{+\infty}_{n=0}pq^{n}\delta_{n}$ &
$\dfrac{p}{1-q\text{e}^{\text{i}t}}$ &
$\dfrac{q}{p}$ &
$\dfrac{q}{p^{2}}$\tabularnewline
\hline 
\begin{cellvarwidth}[t]
\centering
Geometric law on $\mathbb{N}^{\ast},\,$$\mathcal{G}_{\mathbb{N}^{\ast}}\left(p\right)$

$\left(0<p<1\right)$
\end{cellvarwidth} &
$\sum\limits^{+\infty}_{n=1}pq^{n-1}\delta_{n}$ &
$\dfrac{p\text{e}^{\text{i}t}}{1-q\text{e}^{\text{i}t}}$ &
$\dfrac{1}{p}$ &
$\dfrac{q}{p^{2}}$\tabularnewline
\hline 
\end{tabular}

\subsubsection{Laws with Density $\mu=f.\lambda$}

\begin{tabular}{|c|c|c|c|c|}
\hline 
\begin{cellvarwidth}[t]
\centering
Name of the law

(parameters)
\end{cellvarwidth} &
\begin{cellvarwidth}[t]
\centering
Density $f\left(x\right)$

$\left(x\in\mathbb{R}\right)$
\end{cellvarwidth} &
\begin{cellvarwidth}[t]
\centering
$\hat{\mu}\left(t\right)$

$\left(t\in\mathbb{R}\right)$
\end{cellvarwidth} &
$m$ &
$\sigma^{2}$\tabularnewline
\hline 
\hline 
\begin{cellvarwidth}[t]
\centering
Uniform law on $\left[a,b\right]$

$\left(a<b\right)$
\end{cellvarwidth} &
$\dfrac{1}{b-a}\boldsymbol{1}_{\left[a,b\right]}\left(x\right)$ &
$\dfrac{\text{e}^{\text{i}tb}-\text{e}^{\text{i}ta}}{\text{i}t\left(b-a\right)}$ &
$\dfrac{a+b}{2}$ &
$\dfrac{\left(b-a\right)^{2}}{12}$\tabularnewline
\hline 
Cauchy law &
$\dfrac{1}{\pi}\dfrac{1}{1+x^{2}}$ &
$\text{e}^{-\left|t\right|}$ &
\multicolumn{2}{c|}{Does not exist}\tabularnewline
\hline 
\begin{cellvarwidth}[t]
\centering
Gauss law $\mathscr{N}\left(m,\sigma^{2}\right)$

$\left(m\in\mathbb{R},\,\sigma^{2}>0\right)$
\end{cellvarwidth} &
$\dfrac{1}{\sigma\sqrt{2\pi}}\text{e}^{-\frac{\left(x-m\right)^{2}}{2\sigma^{2}}}$ &
$\text{e}^{\text{i}mt-\frac{t^{2}-\sigma^{2}}{2}}$ &
$m$ &
$\sigma^{2}$\tabularnewline
\hline 
First Laplace law &
$\dfrac{1}{2}\text{e}^{-\left|x\right|}$ &
$\dfrac{1}{1+t^{2}}$ &
0 &
2\tabularnewline
\hline 
\begin{cellvarwidth}[t]
\centering
Exponential law $\exp\left(p\right)$

$\left(p>0\right)$
\end{cellvarwidth} &
$\boldsymbol{1}_{\mathbb{R}^{+}}\left(x\right)p\text{e}^{-px}$ &
$\dfrac{1}{1-\dfrac{\text{i}t}{p}}$ &
$\dfrac{1}{p}$ &
$\dfrac{1}{p^{2}}$\tabularnewline
\hline 
\begin{cellvarwidth}[t]
\centering
Gamma law $\gamma\left(a,p\right)$

$\left(a>0,\,p>0\right)$
\end{cellvarwidth} &
$\boldsymbol{1}_{\mathbb{R}^{+}}\left(x\right)\dfrac{p^{\alpha}}{\Gamma\left(a\right)}\text{e}^{-px}x^{a-1}$ &
$\left(1-\dfrac{\text{i}t}{p}\right)^{-a}$ &
$\dfrac{a}{p}$ &
$\dfrac{a}{p^{2}}$\tabularnewline
\hline 
\begin{cellvarwidth}[t]
\centering
Chi-squared law

with $n$ degrees 

of freedom $\chi^{2}_{n}$
\end{cellvarwidth} &
$\boldsymbol{1}_{\mathbb{R}^{+}}\left(x\right)\dfrac{\text{e}^{-\frac{x}{2}}x^{\frac{n}{2}-1}}{2^{\frac{n}{2}}\Gamma\left(\dfrac{n}{2}\right)}$ &
$\left(1-2\text{i}t\right)^{-\frac{n}{2}}$ &
$n$ &
$2n$\tabularnewline
\hline 
\begin{cellvarwidth}[t]
\centering
Beta law 

of first kind

$\left(a>0,\,b>0\right)$
\end{cellvarwidth} &
\begin{cellvarwidth}[t]
\centering
$\dfrac{1}{B\left(a,b\right)}x^{a-1}\left(1-x\right)^{b-1}$

$\left(x\in\left[0,1\right]\right)$
\end{cellvarwidth} &
 &
$\dfrac{a}{a+b}$ &
$\dfrac{ab}{\left(a+b\right)^{2}\left(a+b+1\right)}$\tabularnewline
\hline 
\begin{cellvarwidth}[t]
\centering
Beta law 

of second kind

$\left(a>0,\,b>0\right)$
\end{cellvarwidth} &
\begin{cellvarwidth}[t]
\centering
$\dfrac{1}{B\left(a,b\right)}\dfrac{x^{a-1}}{\left(1+x\right)^{a+b}}$

$\left(x\in\mathbb{R}^{+}\right)$
\end{cellvarwidth} &
 &
\begin{cellvarwidth}[t]
\centering
If $b>1,$

$\dfrac{a}{b-1}$

If $b\leqslant1,$

does not 

exist
\end{cellvarwidth} &
\begin{cellvarwidth}[t]
\centering
If $b>2,$

$\dfrac{a\left(a+b-1\right)}{\left(b-1\right)^{2}\left(b-2\right)}$

If $b\leqslant2,$

does not exist
\end{cellvarwidth}\tabularnewline
\hline 
\begin{cellvarwidth}[t]
\centering
Student law

with $n$ degrees 

of freedom
\end{cellvarwidth} &
$\dfrac{\Gamma\left(\dfrac{n+1}{2}\right)}{\sqrt{n\pi}\Gamma\left(\dfrac{n}{2}\right)}\left(1+\dfrac{t^{2}}{n}\right)^{-\frac{n+1}{2}}$ &
 &
$0$ &
\begin{cellvarwidth}[t]
\centering
If $n>2,$$\dfrac{n}{n-2}$

If $n=2,$

does not exist
\end{cellvarwidth}\tabularnewline
\hline 
\end{tabular}

\begin{remark}{}{}
\begin{itemize}
\item An exponential law $\text{e}^{p}$ is a gamma law with parameters
$\left(1,p\right).$
\item A chi-squared law with $n$ degrees of freedom is a gamma law with
parameters $\left(\dfrac{n}{2},\dfrac{1}{2}\right).$
\item A uniform law on $\left[0,1\right]$ is a beta law of first kind with
parameters $\left(1,1\right).$
\end{itemize}
\end{remark}

\begin{reminders}{Eulerian Functions B (Bêta) and $\Gamma$ (Gamma)}{}
\begin{itemize}
\item Gamma function $\Gamma$\\
For $a>0,$
\[
\Gamma\left(a\right)=\intop^{+\infty}_{0}\text{e}^{-x}x^{a-1}\text{d}x.
\]
For $a>1,$
\[
\Gamma\left(a\right)=\left(a-1\right)\Gamma\left(a-1\right).
\]
\[
\Gamma\left(\dfrac{1}{2}\right)=\sqrt{\pi}.
\]
\item Beta function $\text{B}$\\
For $a,b>0,$
\[
\text{B}\left(a,b\right)=\intop^{+\infty}_{0}\dfrac{x^{a-1}}{\left(1+x\right)^{a+b}}\text{d}x=\intop^{1}_{0}x^{a-1}\left(1-x\right)^{b-1}\text{d}x.
\]
\[
\text{B}\left(a,b\right)=\dfrac{\Gamma\left(a\right)\Gamma\left(b\right)}{\Gamma\left(a+b\right)}.
\]
\end{itemize}
\end{reminders}

\section*{Exercises}

\addcontentsline{toc}{section}{Exercises}

All the random variables introduced below are defined on a probabilized
space $\left(\Omega,\mathcal{A},P\right).$

\begin{exercise}{Fundamental Result For The Simulation of Probability Laws}{exercise9.1}

Let $X$ be a random variable with cumulative distribution function
$F.$

Define the function $G$ with real-valued variable by
\[
\forall t\in\mathbb{R},\,\,\,\,G\left(t\right)=\text{inf}\left\{ x\in\mathbb{R}:\,F\left(x\right)\geqslant t\right\} .
\]

This function $G$ is called the \textbf{\index{pseudo-inverse}pseudo-inverse\mindex{cumulative distribution function ! pseudo-inverse}}
of $F.$

1. Prove successively that:

(a) If $F$ is continuous, then
\[
\forall t\in\left]0,1\right[,\,\,\,\,F\left(G\left(t\right)\right)=t.
\]

(b) If $F$ is increasing, then
\[
\forall x\in\mathbb{R},\,\,\,\,\,G\left(F\left(x\right)\right)=x.
\]

(c) If $F$ is continuous and increasing, then $F$ is a bijection
from $\mathbb{R}$ onto $\left]0,1\right[$ and $G=F^{-1}.$

2. Prove that if $F$ is continuous and increasing, then $F\left(X\right)$
follows the uniform law on $\left[0,1\right].$

3. Prove that if $Y$ is a uniform law on $\left[0,1\right],$ then
the random variable $G\left(Y\right)$ admits $F$ as cumulative distribution
function.

\end{exercise}

\begin{exercise}{Simulation of Law of Discrete Random Variables}{exercise9.2}

Let $X$ be a discrete real-valued random variable taking values in
an increasing sequence $\left(x_{n}\right)_{n\in\mathbb{N}},$ such
that for every $n\in\mathbb{N},$ $P\left(X=x_{n}\right)=p_{n}\geqslant0,$
with $\sum^{+\infty}_{n=0}p_{n}=1.$

Let $U$ be a random variable following the uniform law on $\left[0,1\right],$
and define the random variable $Y$ by
\[
Y=x_{0}\boldsymbol{1}_{\left(U<p_{0}\right)}+\sum^{+\infty}_{n=1}x_{n}\boldsymbol{1}_{\left(p_{0}+p_{1}+\dots+p_{n-1}<U<p_{0}+p_{1}+\dots+p_{n}\right)}.
\]

Verify that $X$ and $Y$ have the same law.

\end{exercise}

\begin{exercise}{Simulation of the Exponential Law}{exercise9.3}

Let $U$ be a random variable following the uniform law on $\left[0,1\right].$ 

Let $X$ be the random variable defined by
\[
X=-\dfrac{1}{p}\ln\left(U\right),
\]
 where $p>0.$

Determine the law of $X.$

\end{exercise}

\begin{figure}[t]
\begin{center}\includegraphics[width=0.4\textwidth]{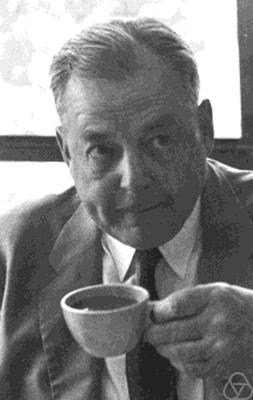}

{\tiny Source: \href{https://opc.mfo.de/person_detail?id=10946}{Institute of Mathematical Statistics}
CC-BY-SA 2.0}\end{center}

\caption{\textbf{\protect\href{https://en.wikipedia.org/wiki/Harold_Hotelling}{Harold Hotelling}}
(1895-1973)}\sindex[fam]{Hotelling, Harold}
\end{figure}

\begin{exercise}{Normal Laws in $\mathbb{R}^2$. Exponential and Hotelling Laws}{exercise9.4}

Let $X=\left(X_{1},X_{2}\right)$ be a random variable taking values
in $\mathbb{R}^{2}$ following the standard normal law $\mathcal{N}_{\mathbb{R}^{2}}\left(0,1\right),$
that is, with density function $f_{X}$ given by
\[
\forall x\in\mathbb{R}^{2},\,\,\,\,f_{X}\left(x\right)=\dfrac{1}{2\pi}\text{e}^{-\frac{\left\Vert x\right\Vert ^{2}}{2}},
\]
where $\left\Vert \cdot\right\Vert $ denotes the usual Euclidean
norm.

1. Determine the law of the random variable $\left\Vert X\right\Vert ^{2}.$

2. Let $D=\left\{ \left(x_{1},x_{2}\right)\in\mathbb{R}^{2}:\,x_{1}=x_{2}\right\} .$
Prove that the random variable $T$ defined by
\[
T=\begin{cases}
\left(\dfrac{X_{1}+X_{2}}{X_{1}-X_{2}}\right)^{2}, & \text{if\,}X\notin D,\\
0, & \text{otherwise,}
\end{cases}
\]
admits a density. Compute it. 

The law of $T$ is called the \textbf{Hotelling\footnotemark\footnotemark\sindex[fam]{Hotelling, Harold}
law}\index{Hotelling law}\mindex{law ! Hotelling}.

\end{exercise}

\addtocounter{footnote}{-1}

\footnotetext{\textbf{\href{https://en.wikipedia.org/wiki/Harold_Hotelling}{Harold Hotelling}\sindex[fam]{Hotelling, Harold}}
(1895 - 1973) was an American mathematical statistician and an economic
theorist. He is known for Hotelling law and lemma, and Hotelling rule
in economics, as well as Hotelling T-squared distribution. He developed
the principal component analysis method, widely used in many area,
including machine learning.}

\stepcounter{footnote}

\footnotetext{This law appears in the study of the comparison tests
between an empirical mean and a theoretical mean---the interested
reader may consult \cite{fourgeaud1967statistique}.}

\begin{exercise}{Moments of Order $\alpha$ and Fubini Theorem}{exercise9.5}

Let $X$ be a nonnegative random variable.

Define the function $G$ by
\[
\forall x\in\mathbb{R},\,\,\,\,G\left(x\right)=P\left(X>x\right).
\]

Prove that for $X$ to admit a moment of order $\alpha\geqslant1,$
it is necessary and sufficient that the function $x\mapsto x^{\alpha-1}G\left(x\right)$
is Lebesgue-integrable on $\mathbb{R}^{+};$ in that case,\boxeq{
\[
\mathbb{E}\left(X^{\alpha}\right)=\alpha\intop_{\mathbb{R}^{+}}x^{\alpha-1}G\left(x\right)\text{d}\lambda\left(x\right).
\]
}

\end{exercise}

\begin{exercise}{Hölder Inequality Equivalent Formulation}{exercise9.6}

Let $p,\,q$ and $r$ be three nonnegative real numbers such that
$\dfrac{1}{p}+\dfrac{1}{q}=\dfrac{1}{r}.$

Prove that:

1. For every random variables $X,Y$ taking values in $\overline{\mathbb{R}}^{+}$
and defined $P-$almost surely, the following inequality on $\overline{\mathbb{R}}^{+}$
holds
\begin{equation}
\left(\intop_{\Omega}\left(XY\right)^{r}\text{d}P\right)^{\frac{1}{r}}\leqslant\left(\intop_{\Omega}X^{p}\text{d}P\right)^{\frac{1}{p}}\left(\intop_{\Omega}Y^{q}\text{d}P\right)^{\frac{1}{q}}.\label{eq:holder_equiv}
\end{equation}

2. If $X\in\mathscr{L}^{p}\left(\Omega,\mathcal{A},P\right)$ and
$Y\in\mathscr{L}^{q}\left(\Omega,\mathcal{A},P\right),$ then $XY\in\mathscr{L}^{r}\left(\Omega,\mathcal{A},P\right),$
and\boxeq{
\begin{equation}
\left\Vert XY\right\Vert _{r}\leqslant\left\Vert X\right\Vert _{p}\left\Vert Y\right\Vert _{q}.\label{eq:inequality_norm_Lr}
\end{equation}
}

Deduce that if $p,q$ and $r$ are three nonnegative real numbers
such that $\dfrac{1}{p}+\dfrac{1}{q}+\dfrac{1}{r}=1,$ if $X\in\mathscr{L}^{p}\left(\Omega,\mathcal{A},P\right),$
$Y\in\mathscr{L}^{q}\left(\Omega,\mathcal{A},P\right),$ and $Z\in\mathscr{L}^{r}\left(\Omega,\mathcal{A},P\right),$
then $XYZ\in\mathscr{L}^{1}\left(\Omega,\mathcal{A},P\right)$ and\boxeq{
\begin{equation}
\left\Vert XYZ\right\Vert _{1}\leqslant\left\Vert X\right\Vert _{p}\left\Vert Y\right\Vert _{q}\left\Vert Z\right\Vert _{r}.\label{eq:inequality_holder_equiv}
\end{equation}
}

\end{exercise}

\begin{exercise}{Variance, Covariance Operator and Support of the Law}{exercise9.7}

Let $X\in\mathscr{L}^{2}\left(\Omega,\mathcal{A},P\right)$ be a random
variable taking values in an Euclidean space $F,$ and let $\Lambda_{X}$
denote its covariance operator.

Prove first that if $F=\mathbb{R},$ then\boxeq{
\[
X=\mathbb{E}\left(X\right)\,\,\,\,P-\text{almost surely}\,\,\,\,\Leftrightarrow\,\,\,\,\,\sigma_{X}=0.
\]
}

Deduce that, in the general case,\boxeq{
\[
P-\text{almost surely}\,\,\,\,\left(X-\mathbb{E}\left(X\right)\right)\in\left(\text{ker}\Lambda_{X}\right)^{\perp}.
\]
}

\end{exercise}

\begin{exercise}{Generalization of the Linear Regression Problem to the Case of Random Variables with Values in an Euclidean Space}{exercise9.8}

Let $F$ and $G$ be two Euclidean space. 

Let $X\in\mathscr{L}^{2}_{F}\left(\Omega,\mathcal{A},P\right)$ and
$Y\in\mathscr{L}^{2}_{G}\left(\Omega,\mathcal{A},P\right)$ be two
random variables.

Assume that the auto-covariance operator $\Lambda_{X}$ is invertible.

Seek the best approximation of $Y$ by a linear function of $X$ in
the sense of the least squares; that is, an element $\left(A,b\right)\in\mathscr{L}\left(F,G\right)\times G$
of the set corresponding to the minimization problem
\[
\underset{\left(A,b\right)\in\mathscr{L}\left(F,G\right)\times G}{\text{argmin}}\Phi\left(A,b\right)
\]
 where
\[
\Phi\left(A,b\right)=\mathbb{E}\left(\left\Vert Y-\left(AX+b\right)\right\Vert ^{2}\right).
\]

\textit{Hint: Introduce the intercovariance operator of $X$ and $Y,$
unique operator such that
\[
\forall\left(x,y\right)\in F\times G,\,\,\,\,\left\langle \Lambda_{X,Y}x,y\right\rangle =\mathbb{E}\left(\left\langle \mathring{X},x\right\rangle \left\langle \mathring{Y},y\right\rangle \right).
\]
}

\textit{Note the symmetry relation
\[
\Lambda_{X,Y}=\left(\Lambda_{Y,X}\right)^{\ast}.
\]
}

\end{exercise}

\section*{Solutions of Exercises}

\addcontentsline{toc}{section}{Solutions of Exercises}

\begin{solution}{}{solexercise9.1}

We distinguish carefully between assumptions of monotonocity and continuity.

\textbf{1. (a) $F$ is continuous implies that $\forall t\in\left]0,1\right[,\,\,\,\,F\left(G\left(t\right)\right)=t$}
\begin{itemize}
\item Let us consider, for every $t\in\left]0,1\right[,$ the set 
\[
A_{t}=\left\{ x:\,F\left(x\right)\geqslant t\right\} .
\]
Then $A_{t}$ is a half-line: since $F$ is nondecreasing, for every
$x_{0}\in A_{t}$ and for every $y\geqslant x_{0},$
\[
F\left(y\right)\geqslant F\left(x_{0}\right)\geqslant t
\]
 and thus $y\in A_{t}.$
\item We begin by proving that\boxeq{
\[
F\left[G\left(t\right)\right]\geqslant t.
\]
}Since $G\left(t\right)=\inf_{t\in\left]0,1\right[}A_{t}$ and that
$A_{t}$ is a half-line, for every $y>G\left(t\right),$ $y\in A_{t}$
and hence $F\left(y\right)\geqslant t.$ \\
As $F$ is continuous from the right, we consider a sequence of real
numbers $y_{n}\in A_{t}$ that converges nonincreasingly to $G\left(t\right),$
and obtain
\[
F\left(G\left(t\right)\right)=\lim_{n\to+\infty}\searrow F\left(y_{n}\right)\geqslant t.
\]
\item Next, since $F$ is continuous,\boxeq{
\[
F\left[G\left(t\right)\right]\leqslant t.
\]
}Indeed, by the definition of $G,$ for every $y\in\mathbb{R}$ such
that $y<G\left(t\right),$ $F\left(y\right)<t.$ As the function $F$
is continuous from the left, we now take a sequence of real numbers
$y_{n}<G\left(t\right)$ that converges nondecreasingly to $G\left(t\right),$
and obtain
\[
F\left(G\left(t\right)\right)=\lim_{n\to+\infty}\nearrow F\left(y_{n}\right)\leqslant t.
\]
\item Hence, if $F$ is continuous, then for every $t\in\left]0,1\right[,$
$F\left(G\left(t\right)\right)=t.$
\end{itemize}
\textbf{(b) $F$ increasing implies $\forall x\in\mathbb{R},\,\,\,\,\,G\left(F\left(x\right)\right)=x$}
\begin{itemize}
\item For every $x\in\mathbb{R},$ by definition of $G,$
\[
G\left(F\left(x\right)\right)=\inf\left\{ y\in\mathbb{R}:\,F\left(y\right)\geqslant F\left(x\right)\right\} 
\]
which implies
\[
G\left(F\left(x\right)\right)\leqslant x.
\]
\item Moreover, if $F$ is increasing, for every $y\in\mathbb{R}$ such
that $F\left(y\right)\geqslant F\left(x\right),$ $y\geqslant x$---otherwise,
we would have $y<x$ and thus $F\left(y\right)<F\left(x\right),$
which would contradict the assumption. This shows that
\[
G\left(F\left(x\right)\right)\geqslant x.
\]
Hence, we conclude that
\[
G\left(F\left(x\right)\right)=x.
\]
\end{itemize}
\textbf{(c) $F$ continuous and increasing implies $F$ is bijective
and $G=F^{-1}$}
\[
\forall t\in\left]0,1\right[,\,\,\,\,F\left(G\left(t\right)\right)=t\,\,\,\,\,\text{and}\,\,\,\,\forall x\in\mathbb{R},\,\,\,\,G\left(F\left(x\right)\right)=x.
\]

Hence, $F$ is a bijection from $\mathbb{R}$ onto $\left]0,1\right[$
and $G=F^{-1}.$

\textbf{2. $F$ continuous and increasing implies that $F\left(X\right)$
follows the uniform law on $\left[0,1\right]$}

If $F$ is continuous, then for every $y\in\left]0,1\right[,$ $F\left(G\left(y\right)\right)=y.$
Hence, using the strict growth of $F,$ we obtain
\[
P\left(F\left(X\right)\leqslant y\right)=P\left(F\left(X\right)\leqslant F\left(G\left(y\right)\right)\right)=P\left(X\leqslant G\left(y\right)\right).
\]

Since $F$ is the cumulative distribution function of $X,$
\[
P\left(X\leqslant G\left(y\right)\right)=F\left(G\left(y\right)\right)=y.
\]

Moreover,
\[
P\left(F\left(X\right)\leqslant y\right)=\begin{cases}
0, & \text{if}\,y<0,\\
1, & \text{if\,}y>1.
\end{cases}
\]

Therefore, $F\left(X\right)$ follows the uniform law on $\left[0,1\right].$

\textbf{3. $Y$ is a uniform law on $\left[0,1\right]$ implies that
the random variable $G\left(Y\right)$ admits $F$ as cumulative distribution
function}

We have the equivalence
\[
F\left(x\right)\geqslant t\Leftrightarrow x\geqslant G\left(t\right).
\]

Hence, if $Y$ follows the uniform law on $\left[0,1\right],$ then
\[
\forall x\in\mathbb{R},\,\,\,\,P\left(G\left(Y\right)\leqslant x\right)=P\left(Y\leqslant F\left(x\right)\right)=F\left(x\right).
\]

This shows that $G\left(Y\right)$ admits $F$ as cumulative distribution
function.

\end{solution}

\begin{remark}{}{}

The results of this exercise theoretically allow to simulate any law
on $\mathbb{R}$ from a random variable following the uniform law.
Indeed, a call to the ``random'' function on a computer is assumed
to return a realization $y$ of a random variable $Y$ following the
uniform law on $\left[0,1\right],$ this random number being generated
by a uniform generator.

If we want to simulate a real-valued random variable of cumulative
distribution function $F,$ we compute, when feasible, its pseudo-inverse
$G.$ Then $G\left(y\right)$ is a realization of the random variable
$G\left(Y\right)$ which admits $F$ as cumulative distribution function.

This method can be numerically intensive, or even intractable. Specific
methods exist to simulate several classical laws, as the following
exercise will show.

\end{remark}

\begin{solution}{}{solexercise9.2}

We have
\[
P\left(Y=x_{0}\right)=P\left(U<p_{0}\right)=p_{0}.
\]

And, for every $n\geqslant1,$
\[
P\left(Y=x_{n}\right)=P\left(p_{0}+p_{1}+\dots+p_{n-1}<U<p_{0}+p_{1}+\dots+p_{n}\right)=p_{n}.
\]

\end{solution}

\begin{remark}{}{}

The result of this exercise allows to simulate any discrete law on
$\mathbb{R}$ from a random variable following the uniform law.

\end{remark}

\begin{solution}{}{solexercise9.3}

By applying the transfer theorem and then performing the change of
variable associated with the diffeomorphism from $\left]0,1\right[$
to $\mathbb{R}^{+\ast}$ defined by $v=-\dfrac{1}{p}\ln\left(u\right),$
we obtain
\begin{align*}
\forall f\in\mathscr{C}^{+}_{\mathscr{K}}\left(\mathbb{R}\right),\,\,\,\,\mathbb{E}\left(f\left(X\right)\right) & =\mathbb{E}\left(f\left(-\dfrac{1}{p}\ln U\right)\right)\\
 & =\intop_{\mathbb{R}}f\left(-\dfrac{1}{p}\ln U\right)\boldsymbol{1}_{\left]0,1\right[}\left(u\right)\text{d}\lambda\left(u\right)\\
 & =\intop_{\mathbb{R}^{+\ast}}f\left(v\right)\text{e}^{-pv}p\text{d}\lambda\left(v\right).
\end{align*}

This shows that $P_{X}$ is the law $\exp\left(p\right).$ 

\end{solution}

\begin{remark}{}{}

The result of this exercise is often used to simulate an exponential
law from a random variable following the uniform law, without resorting
to the general but heavier method involving the pseudo-inverse of
the cumulative distribution function.

\end{remark}

\begin{solution}{}{solexercise9.4}

\textbf{1. Law of $\left\Vert X\right\Vert ^{2}.$}

By applying the transfer theorem and the integration theorem related
to a density measure, for every $f\in\mathscr{C}^{+}_{\mathscr{K}}\left(\mathbb{R}\right),$
\begin{align*}
\mathbb{E}\left(f\left(\left\Vert X\right\Vert ^{2}\right)\right) & =\intop_{\mathbb{R}^{2}}f\left(\left\Vert x\right\Vert ^{2}\right)\dfrac{1}{2\pi}\text{e}^{-\frac{\left\Vert x\right\Vert ^{2}}{2}}\text{d}\lambda_{2}\left(x\right)\\
 & =\intop_{\mathbb{R}^{2}\backslash\mathbb{R}^{+}\times\left\{ 0\right\} }f\left(\left\Vert x\right\Vert ^{2}\right)\dfrac{1}{2\pi}\text{e}^{-\frac{\left\Vert x\right\Vert ^{2}}{2}}\text{d}\lambda_{2}\left(x\right).
\end{align*}
Performing the change of variables in polar coordinates associated
with the diffeomorphism from $\mathbb{R}^{+\ast}\times\left]0,2\pi\right[$
onto $\mathbb{R}^{2}\backslash\mathbb{R}^{+}\times\left\{ 0\right\} $
defined by
\[
x=\rho\cos\theta\,\,\,\,y=\rho\sin\theta,
\]
and whose Jacobian is $\rho,$ and applying the Fubini theorem, we
obtain
\begin{align*}
\mathbb{E}\left(f\left(\left\Vert X\right\Vert ^{2}\right)\right) & =\intop_{\mathbb{R}^{+\ast}\times\left]0,2\pi\right[}f\left(\rho^{2}\right)\dfrac{1}{2\pi}\text{e}^{-\frac{\rho^{2}}{2}}\rho\text{d}\left(\lambda\otimes\lambda\right)\left(\rho,\theta\right)\\
 & =\intop_{\mathbb{R}^{+\ast}}f\left(\rho^{2}\right)\text{e}^{-\frac{\rho^{2}}{2}}\rho\text{d}\lambda\left(\rho\right).
\end{align*}

Using a final change of variable associated with the diffeomorphism
from $\mathbb{R}^{+\ast}$ onto itself, defined by 
\[
u=\rho^{2},
\]
we obtain
\[
\forall f\in\mathscr{C}^{+}_{\mathscr{K}}\left(\mathbb{R}\right),\,\,\,\,\mathbb{E}\left(f\left(\left\Vert X\right\Vert ^{2}\right)\right)=\intop_{\mathbb{R}}f\left(u\right)\boldsymbol{1}_{\mathbb{R}^{+}}\left(u\right)\dfrac{1}{2}\text{e}^{-\frac{u}{2}}\text{d}\lambda\left(u\right),
\]
which shows that $\left\Vert X\right\Vert ^{2}$ follows an exponential
law $\exp\left(\dfrac{1}{2}\right).$

\textbf{2. Computation of Hotteling law}

By applying the transfer theorem and the integration theorem related
to a density measure, we have, for every $f\in\mathscr{C}^{+}_{\mathscr{K}}\left(\mathbb{R}\right),$
and since $\lambda_{2}\left(D\right)=0,$
\[
\mathbb{E}\left(f\left(T\right)\right)=\intop_{\mathbb{R}^{2}\backslash D}f\left(\left(\dfrac{x_{1}+x_{2}}{x_{1}-x_{2}}\right)^{2}\right)\dfrac{1}{2\pi}\text{e}^{-\frac{x^{2}_{1}+x^{2}_{2}}{2}}\text{d}\lambda_{2}\left(x_{1},x_{2}\right).
\]

We perform the change of variables associated with the diffeomorphism
from $\mathbb{R}^{2}\backslash D$ onto $\mathbb{R}^{2}\backslash\left(\left\{ 0\right\} \times\mathbb{R}\right)$
defined by
\[
\begin{cases}
u= & \dfrac{x_{1}+x_{2}}{x_{1}-x_{2}}\\
v= & x_{1}+x_{2}
\end{cases}\Leftrightarrow\begin{cases}
x_{1}= & \dfrac{1}{2}\left(v+\dfrac{v}{u}\right)\\
x_{2}= & \dfrac{1}{2}\left(v-\dfrac{v}{u}\right)
\end{cases}
\]
whose Jacobian is
\begin{align*}
\dfrac{D\left(x_{1},x_{2}\right)}{D\left(u,v\right)} & =\det\left(\begin{array}{cc}
-\dfrac{v}{2u^{2}} & \dfrac{1}{2}\left(1+\dfrac{1}{u}\right)\\
\dfrac{v}{2u^{2}} & \dfrac{1}{2}\left(1-\dfrac{1}{u}\right)
\end{array}\right)\\
 & =-\dfrac{v}{2u^{2}}.
\end{align*}

By applying the Fubini theorem, we obtain, for every $f\in\mathscr{C}^{+}_{\mathscr{K}}\left(\mathbb{R}\right),$
\begin{align*}
\mathbb{E}\left(f\left(T\right)\right) & =\dfrac{1}{2\pi}\intop_{\mathbb{R}^{2}\backslash\left\{ 0\right\} \times\mathbb{R}}f\left(u^{2}\right)\text{e}^{-\frac{1}{4}v^{2}\left(1+\frac{1}{u^{2}}\right)}\left|\dfrac{v}{2u^{2}}\right|\text{d}\lambda\otimes\lambda\left(u,v\right)\\
 & =\dfrac{1}{4\pi}\intop_{\mathbb{R}^{\ast}}\dfrac{f\left(u^{2}\right)}{u^{2}}\left[\intop_{\mathbb{R}}\left|v\right|\text{e}^{-\frac{1}{4}v^{2}\left(1+\frac{1}{u^{2}}\right)}\text{d}\lambda\left(v\right)\right]\text{d}\lambda\left(u\right).
\end{align*}

Now, by comparison of the Lebesgue and generalized Riemann integrals,
we have, and using the variable change on the second line 
\[
w=\dfrac{1}{4}v^{2}\left(1+\dfrac{1}{u^{2}}\right),
\]
we compute
\begin{align*}
\intop_{\mathbb{R}}\left|v\right|\text{e}^{-\frac{1}{4}v^{2}\left(1+\frac{1}{u^{2}}\right)}\text{d}\lambda\left(v\right) & =2\intop^{+\infty}_{0}v\text{e}^{-\frac{1}{4}v^{2}\left(1+\frac{1}{u^{2}}\right)}\text{d}v\\
 & =\intop^{+\infty}_{0}\dfrac{8}{1+\dfrac{1}{u^{2}}}\text{e}^{-w}\text{d}w\\
 & =\dfrac{8}{1+\dfrac{1}{u^{2}}}.
\end{align*}

Hence, for every $f\in\mathscr{C}^{+}_{\mathscr{K}}\left(\mathbb{R}\right),$
performing a second change of variable $t=u^{2},$ we obtain
\begin{align*}
\mathbb{E}\left(f\left(T\right)\right) & =\dfrac{1}{\pi}\intop_{\mathbb{R}^{\ast}}f\left(u^{2}\right)\dfrac{2}{u^{2}+1}\text{d}\lambda\left(u\right)\\
 & =\dfrac{1}{\pi}\intop_{\mathbb{R}}\boldsymbol{1}_{\mathbb{R}^{+\ast}}\left(t\right)f\left(t\right)\dfrac{1}{\left(t+1\right)\sqrt{t}}\text{d}\lambda\left(t\right).
\end{align*}

This shows that $T$ admits a density $f_{T}$ given by\boxeq{
\[
\forall t\in\mathbb{R},\,\,\,\,f_{T}\left(t\right)=\dfrac{1}{\pi}\boldsymbol{1}_{\mathbb{R}^{+\ast}}\left(t\right)\dfrac{1}{\left(t+1\right)\sqrt{t}}.
\]
}

This is the Beta law of second kind, $\text{B}\left(\dfrac{1}{2},\dfrac{1}{2}\right).$

\end{solution}

\begin{solution}{}{solexercise9.5}

By the Fubini theorem, for every nonnegative measurable functions,
\begin{align*}
\intop_{\mathbb{R}^{+}}x^{\alpha-1}G\left(x\right)\text{d}\lambda\left(x\right) & =\intop_{\mathbb{R}^{+}}x^{\alpha-1}\left(\intop_{\Omega}\boldsymbol{1}_{\left(X>x\right)}\text{d}P\right)\text{d}\lambda\left(x\right)\\
 & =\intop_{\Omega}\left(\intop_{\mathbb{R}^{+}}x^{\alpha-1}\boldsymbol{1}_{\left(X>x\right)}\text{d}\lambda\left(x\right)\right)\text{d}P\\
 & =\intop_{\Omega}\dfrac{X^{\alpha}}{\alpha}\text{d}P.
\end{align*}

Hence, the result.

\end{solution}

\begin{solution}{}{solexercise9.6}

1. This follows from the inequality $\refpar{eq:Holder_inequality}$
applied to the random variables $X^{r}$ and $Y^{r},$ with the conjugate
real numbers $\dfrac{p}{r}$ and $\dfrac{q}{r},$ giving
\[
\intop_{\Omega}X^{r}Y^{r}\text{d}P\leqslant\left(\intop_{\Omega}\left(X^{r}\right)^{p/r}\text{d}P\right)^{r/p}\left(\intop_{\Omega}\left(Y^{r}\right)^{q/r}\text{d}P\right)^{r/q}.
\]

Raising both sides to the power $\dfrac{1}{r}$ yields $\refpar{eq:holder_equiv}.$

2. If $X\in\mathscr{L}^{p}\left(\Omega,\mathcal{A},P\right)$ and
$Y\in\mathscr{L}^{q}\left(\Omega,\mathcal{A},P\right),$ then it follows
that
\[
\intop_{\Omega}\left|XY\right|^{r}\text{d}P\leqslant\left(\intop_{\Omega}\left|X^{r}\right|^{p/r}\text{d}P\right)^{r/p}\left(\intop_{\Omega}\left|Y^{r}\right|^{q/r}\text{d}P\right)^{r/q}<+\infty,
\]
and thus 
\[
XY\in\mathscr{L}^{r}\left(\Omega,\mathcal{A},P\right).
\]
The inequality $\refpar{eq:inequality_norm_Lr}$ follows.

\begin{remark}{}{}

This formulation is therefore equivalent to the one yielding the Hölder
inequality---for the converse, take $r=1.$

\end{remark}

Lastly, let $p,q$ and $r$ be three nonnegative real numbers such
that 
\[
\dfrac{1}{p}+\dfrac{1}{q}+\dfrac{1}{r}=1.
\]
 Define $\alpha$ by the relationship 
\[
\dfrac{1}{\alpha}=\dfrac{1}{q}+\dfrac{1}{r}.
\]

Let $X\in\mathscr{L}^{p}\left(\Omega,\mathcal{A},P\right),$ $Y\in\mathscr{L}^{q}\left(\Omega,\mathcal{A},P\right),$
and $Z\in\mathscr{L}^{r}\left(\Omega,\mathcal{A},P\right).$ From
what precedes, we have $YZ\in\mathscr{L}^{\alpha}\left(\Omega,\mathcal{A},P\right)$
and
\[
\left\Vert YZ\right\Vert _{\alpha}\leqslant\left\Vert Y\right\Vert _{q}\left\Vert Z\right\Vert _{r}.
\]

Since the real numbers $\alpha$ and $p$ are conjugate, the product
$X\left(YZ\right)$ is integrable, the Hölder inequality then yields
\[
\left\Vert X\left(YZ\right)\right\Vert _{1}\leqslant\left\Vert X\right\Vert _{p}\left\Vert YZ\right\Vert _{\alpha},
\]
which, combined with the previous inequality proves the inequality
$\refpar{eq:inequality_holder_equiv}.$

\end{solution}

\begin{solution}{}{solexercise9.7}

If $F=\mathbb{R},$ then
\[
\sigma^{2}_{X}=0\Leftrightarrow\mathbb{E}\left(\left(X-\mathbb{E}\left(X\right)\right)^{2}\right)=0\Leftrightarrow X=\mathbb{E}\left(X\right)\,\,\,\,P-\text{almost surely}.
\]

In the general case,
\[
x\in\text{ker}\Lambda_{X}\Leftrightarrow\mathbb{E}\left(\left\langle \mathring{X},x\right\rangle ^{2}\right)=0\Leftrightarrow\left\langle \mathring{X},x\right\rangle =0\,\,\,\,P-\text{almost surely}.
\]

Let $\left(e_{i}\right)_{i=1,\cdots,l}$ be a basis of $\text{ker}\Lambda_{X}.$

For every $i=1,\cdots,l,$ there exists a set of zero probability
$N_{i}$ such that, for every $\omega\notin N_{i},$
\[
\left\langle \mathring{X}\left(\omega\right),e_{i}\right\rangle =0.
\]

Let $N=\bigcup^{l}_{i=1}N_{i}.$ We have $P\left(N\right)=0$ and
\[
\forall\omega\notin N,\,\,\,\,\forall i=1,\cdots,l,\,\,\,\,\left\langle \mathring{X}\left(\omega\right),e_{i}\right\rangle =0.
\]

Hence,
\[
\forall\omega\notin N,\,\,\,\,\mathring{X}\left(\omega\right)\in\left(\text{ker}\Lambda_{X}\right)^{\perp}.
\]

\end{solution}

\begin{solution}{}{solexercise9.8}

Since $\mathring{Y}-A\mathring{X}$ is centered,
\begin{align*}
\Phi\left(A,b\right) & =\mathbb{E}\left(\left\Vert \mathring{Y}-A\mathring{X}+\left(\mathbb{E}\left(Y\right)-A\mathbb{E}\left(X\right)-b\right)\right\Vert ^{2}\right)\\
 & =\mathbb{E}\left(\left\Vert \mathring{Y}-A\mathring{X}\right\Vert ^{2}\right)+\left\Vert \mathbb{E}\left(Y\right)-A\mathbb{E}\left(X\right)-b\right\Vert ^{2}+2\mathbb{E}\left[\left\langle \mathring{Y}-A\mathring{X},\mathbb{E}\left(Y\right)-A\mathbb{E}\left(X\right)-b\right\rangle \right]\\
 & =\mathbb{E}\left(\left\Vert \mathring{Y}-A\mathring{X}\right\Vert ^{2}\right)+\left\Vert \mathbb{E}\left(Y\right)-A\mathbb{E}\left(X\right)-b\right\Vert ^{2}.
\end{align*}
For each fixed $A,$ this quantity reaches its minimum when 
\[
\hat{b}_{A}=\mathbb{E}\left(Y\right)-A\mathbb{E}\left(X\right).
\]
It remains to minimize with respect to $A,$
\[
\Phi\left(A,\hat{b}_{A}\right)\equiv\mathbb{E}\left(\left\Vert \mathring{Y}-A\mathring{X}\right\Vert ^{2}\right).
\]

Now, we compute
\begin{align*}
\mathbb{E}\left(\left\Vert \mathring{Y}-A\mathring{X}\right\Vert ^{2}\right) & =\text{tr}\left(\mathbb{E}\left(\left[\mathring{Y}-A\mathring{X}\right]\left[\mathring{Y}-A\mathring{X}\right]^{\ast}\right)\right)\\
 & =\text{tr}\left(\mathbb{E}\left(\mathring{Y}\mathring{Y}^{\ast}+A\left(\mathring{X}\mathring{X}^{\ast}\right)A^{\ast}-\left(\mathring{Y}\mathring{X}^{\ast}\right)A^{\ast}-A\mathring{X}\mathring{Y}^{\ast}\right)\right)\\
 & =\text{tr}\left(\Lambda_{Y}+A\Lambda_{X}A^{\ast}-\Lambda_{X,Y}A^{\ast}-A\Lambda_{Y,X}\right).
\end{align*}

Define $\Psi$ as the differentiable function from $\mathscr{L}\left(F,G\right)$
to $\mathbb{R},$ given by
\[
\Psi\left(A\right)=\text{tr}\left(\Lambda_{Y}+A\Lambda_{X}A^{\ast}-\Lambda_{X,Y}A^{\ast}-A\Lambda_{Y,X}\right).
\]

Its differential at $A$ is given by
\begin{align*}
\forall H\in\mathscr{L}\left(F,G\right),\,\,\,\,\Psi^{\prime}\left(A\right)\left(H\right) & =\text{tr}\left(H\Lambda_{X}A^{\ast}+A\Lambda_{X}H^{\ast}-\Lambda_{X,Y}H^{\ast}-H\Lambda_{Y,X}\right)\\
 & =2\text{tr}\left(\left(A\Lambda_{X}-\Lambda_{X,Y}\right)H^{\ast}\right).
\end{align*}
A stationary point $\widehat{A}$ is given by
\[
\widehat{A}=\Lambda_{X,Y}\Lambda^{-1}_{X}.
\]

This stationary point corresponds to a minimum. Indeed,
\[
\forall H\in\mathscr{L}\left(F,G\right),\,\,\,\,\Psi^{\prime\prime}\left(A\right)\left(H\right)\left(H\right)=2\text{tr}\left(H\Lambda_{X}H^{\ast}\right)\geqslant0.
\]

Moreover, if $\left(e_{i}\right)_{1\leqslant i\leqslant d}$ is an
orthonormal basis of $G,$ then
\begin{align*}
\text{tr}\left(H\Lambda_{X}H^{\ast}\right) & =\sum^{d}_{j=1}\left\langle H\Lambda_{X}H^{\ast}e_{i},e_{i}\right\rangle =\sum^{d}_{j=1}\left\langle \Lambda_{X}H^{\ast}e_{i},H^{\ast}e_{i}\right\rangle =\sum^{d}_{j=1}\left\Vert \Lambda^{1/2}_{X}H^{\ast}e_{i}\right\Vert ^{2}.
\end{align*}

Since $\Lambda_{X}$ is invertible, it results that $\Psi^{\prime\prime}\left(A\right)\left(H\right)\left(H\right)=0$
if and only if $H^{\ast}e_{i}=0$ for every $i=1,\cdots,d.$ Otherwisely
said, if and only if $H^{\ast},$ thus $H$ is zero. This ensures
that $\widehat{A}$ corresponds to a minimum.

The solution of the linear regression problem is the pair
\[
\left(\widehat{A},\widehat{b}_{\widehat{A}}\right)=\left(\Lambda_{X,Y}\Lambda^{-1}_{X},\mathbb{E}\left(Y\right)-\Lambda_{X,Y}\Lambda^{-1}_{X}\left(\mathbb{E}\left(X\right)\right)\right).
\]

The surface $D$---an affine subspace---of equation\boxeq{
\[
\left(y-\mathbb{E}\left(Y\right)\right)-\Lambda_{X,Y}\Lambda^{-1}\left(x-\mathbb{E}\left(X\right)\right)=0
\]
}is called the surface of linear regression of $Y$ in $X.$ The
best approximation of $Y$ as a linear function of $X$ in the sense
of the least square, is\boxeq{
\[
\mathbb{E}\left(Y\right)+\Lambda_{X,Y}\Lambda^{-1}\left(X-\mathbb{E}\left(X\right)\right)=0.
\]
}

We have $P\left(\left(X,Y\right)\in D\right)=1$ if and only if 
\[
\Phi\left(\widehat{A},\widehat{b}_{\widehat{A}}\right)=0.
\]

\end{solution}

\chapter{Independence of $\sigma-$Algebra, of Random Variables}\label{chap:PartIIChap10}

\begin{objective}{}{}

Chapter \ref{chap:PartIIChap10} aims to introduce the concept of
independence of $\sigma-$algebra and random variables.
\begin{itemize}
\item Section \ref{sec:Independence-of-Family} introduces the concept of
independence for families of events and random variables. After defining
the concept of independence of events, a necessary and sufficient
condition of independence is given in the case of \salg generated
by $\pi-$systems. A similar approach is taken for random variables.
A general criterion for independence is given for two random variables
based on their laws, which is then stated in a functional form. Various
criteria for independence are provided based on cumulative distribution
functions, densities and the discrete case. The impact of independence
on expectation, covariance, and variance is then analyzed, as well
as on the covariance operator and the covariance matrix. The section
concludes by extending this concept to arbitrary families indexed
by a set $I.$
\item Section \ref{sec:Independence-and-Asymptotic} addresses the concept
of asymptotic $\sigma-$algebras. After providing a definition, the
zero-one law is stated. The section ends with the Borel-Cantelli lemma.
\item Section \ref{sec:Some-Results-Linked} focuses on the heads and tails
model, with the aim of generalizing it to infinite sequences of the
game. To this end, it begins by recalling the dyadic expansion of
a real number, before presenting results on this model.
\item Section \ref{sec:Convolution-and-Law} begins by defining the convolution
product of two bounded measures, and then addresses the law of the
sum of two independent random variables.
\end{itemize}
\end{objective}

Independence is a fundamental concept in probability theory and statistics.
It is always defined with respect to a given probabilized space $\left(\Omega,\mathcal{A},P\right).$
All random variables are assumed to be defined on this space.

The elementary concepts on independence studied in Part \ref{part:Introduction-to-Probability}
Chapter \ref{chap:Independence} are assumed to be known.

\section{Independence of Family of Events and of Random Variables}\label{sec:Independence-of-Family}

\begin{definition}{Independent Event. Independent Families}{}

Let $\left(\Omega,\mathcal{A},P\right)$ be a probabilized space.

(a) The events $A\in\mathcal{A}$ and $B\in\mathcal{A}$ are \textbf{independent\index{independent}\mindex{independent ! events}}
if\boxeq{
\[
P\left(A\cap B\right)=P\left(A\right)P\left(B\right).
\]
}

(b) Two \textbf{families of events} $\mathcal{A}_{1}$ and $\mathcal{A}_{2}$
are \textbf{independent\mindex{independent ! families}} if every
element of $\mathcal{A}_{1}$ is independent of every element of $\mathcal{A}_{2}.$

\end{definition}

\begin{remark}{}{}

As established in Part \ref{part:Introduction-to-Probability} Chapter
\ref{chap:Independence}, we verify that, for $A\in\mathcal{A}$ and
$B\in\mathcal{A}$ to be independent, it is necessary and sufficient
that the $\sigma-$algebras $\sigma\left(\left\{ A\right\} \right)$
and $\sigma\left(\left\{ B\right\} \right)$ are independent. The
next theorem generalizes this property and is commonly used.

\end{remark}

\begin{theorem}{Necessary and Sufficient Condition for the Independence of Family of Events}{nsc_independence_family_events}

Let $\mathscr{C}_{1}$ and $\mathscr{C}_{2}$ be two $\pi-$systems
contained in $\mathcal{A}.$ Let $\mathscr{F}_{1}$---respectively
$\mathscr{F}_{2}$---denote the $\sigma-$algebra generated by $\mathscr{C}_{1}$---respectively
$\mathscr{C}_{2}.$ For the family of events $\mathscr{F}_{1}$ and
$\mathscr{F}_{2}$ to be independent, it is necessary and sufficient
that the families of events $\mathscr{C}_{1}$ and $\mathscr{C}_{2}$
are independent.

\end{theorem}

\begin{proof}{}{}

The necessary condition follows directly from the definition of the
independence of two families of events.

For the sufficient condition, suppose that $\mathscr{C}_{1}$ and
$\mathscr{C}_{2}$ are independent. We use the principle of the extension
by measurability, in its set-theoretic version---Lemma $\ref{lm:extension_by_measurability}.$

Let $\mathscr{D}$ be the set of independent events of every element
of $\mathscr{C}_{2},$ defined by
\[
\mathscr{D}=\left\{ D\in\mathscr{A}:\,\forall A_{2}\in\mathscr{C}_{2},\,\,\,\,P\left(D\cap A_{2}\right)=P\left(D\right)P\left(A_{2}\right)\right\} .
\]

By assumption, $\mathscr{C}_{1}\subset\mathscr{D}$ and $\Omega\in\mathscr{D}.$
We now prove that $\mathscr{D}$ is a $\lambda-$system.
\begin{itemize}
\item If $D_{1}$ and $D_{2}\in\mathscr{D}$ with $D_{1}\supset D_{2},$
then for every $A_{2}\in\mathscr{C}_{2},$
\begin{align*}
P\left(\left(D_{1}\backslash D_{2}\right)\cap A_{2}\right) & =P\left(\left(D_{1}\cap A_{2}\right)\backslash\left(D_{2}\cap A_{2}\right)\right)\\
 & =P\left(D_{1}\cap A_{2}\right)-P\left(D_{2}\cap A_{2}\right).
\end{align*}
Since $D_{1}$ and $D_{2}\in\mathscr{D},$
\begin{align*}
P\left(\left(D_{1}\backslash D_{2}\right)\cap A_{2}\right) & =P\left(D_{1}\right)P\left(A_{2}\right)-P\left(D_{2}\right)P\left(A_{2}\right)\\
 & =P\left(D_{1}\backslash D_{2}\right)P\left(A_{2}\right).
\end{align*}
Moreover, as $D_{1}\backslash D_{2}\in\mathcal{A},$ it follows that
$D_{1}\backslash D_{2}\in\mathscr{D}.$
\item Additionally, for every nondecreasing sequence $\left(D_{n}\right)_{n\in\mathbb{N}}\subset\mathscr{D},$
\[
\forall n\in\mathbb{N},\,\forall A_{2}\in\mathscr{C}_{2},\,\,\,\,P\left(D_{n}\cap A_{2}\right)=P\left(D_{n}\right)P\left(A_{2}\right).
\]
Taking the limit by monotonicity
\[
\forall A_{2}\in\mathscr{C}_{2},\,\,\,\,P\left(\left[\bigcup_{n\in\mathbb{N}}D_{n}\right]\cap A_{2}\right)=P\left(\bigcup_{n\in\mathbb{N}}D_{n}\right)P\left(A_{2}\right).
\]
As $\bigcup_{n\in\mathbb{N}}D_{n}\in\mathscr{A},$ we conclude that
$\bigcup_{n\in\mathbb{N}}D_{n}\in\mathscr{D}.$
\end{itemize}
By the principle of extension by measurability, it follows that $\mathscr{D\supset\mathscr{F}}_{1},$
i.e.,
\begin{equation}
\forall F_{1}\in\mathscr{F}_{1},\,\forall A_{2}\in\mathscr{C}_{2},\,\,\,\,P\left(F_{1}\cap A_{2}\right)=P\left(F_{1}\right)P\left(A_{2}\right).\label{eq:independence_F_1_C_2}
\end{equation}

Now, define $\mathscr{E}$ as the set of events that are independent
of every element of $\mathscr{F}_{1}$ by
\[
\mathscr{E}=\left\{ E\in\mathscr{A}:\,\forall F_{1}\in\mathscr{F}_{1},\,\,\,\,P\left(F_{1}\cap E\right)=P\left(F_{1}\right)\cap P\left(E\right)\right\} .
\]

From equation $\refpar{eq:independence_F_1_C_2},$ we have $\mathscr{E}\supset\mathscr{C}_{2}.$
We have: $\Omega\in\mathscr{E}$ and $\mathscr{E}$ is still a $\lambda-$system.
Therefore, by the principle of extension by measurability, we obtain
$\mathscr{E}\supset\mathscr{F}_{2},$ which completes the proof of
the theorem.

\end{proof}

\begin{denotation}{}{}

Let $X$ be a function from $\Omega$ to $E.$ 

If $\mathscr{G}$ is a family of subsets of $E,$ then we denote by
$X^{-1}\left(\mathscr{G}\right)$ the family of subsets of $\Omega$
defined by\boxeq{
\[
X^{-1}\left(\mathscr{G}\right)=\left\{ X^{-1}\left(G\right)\in\mathcal{P}\left(\Omega\right):\,G\in\mathscr{G}\right\} .
\]
}

In particular, if $\mathscr{G}$ is a $\sigma-$algebra on $E,$ the
family $X^{-1}\left(\mathscr{G}\right)$ is a $\sigma-$algebra on
$\Omega,$ called the \textbf{$\sigma-$algebra generated by the function
$X.$}

\end{denotation}

\begin{definition}{Independent Random Variables}{}

Let $X_{i},$ $i=1,2,$ be two random variables defined on a probabilized
space $\left(\Omega,\mathcal{A},P\right),$ taking values in the respective
probabilizable spaces $\left(E_{i},\mathcal{E}_{i}\right).$

The random variables $X_{1}$ and $X_{2}$ are \textbf{independent\index{independent}}
if the $\sigma-$algebra $X^{-1}_{1}\left(\mathcal{E}_{1}\right)$
and $X^{-1}_{2}\left(\mathcal{E}_{2}\right)$ generated respectively
by the functions $X_{1}$ and $X_{2}$ are independent.

\end{definition}

\begin{remark}{Preservation of Independence by Measurable Functions}{}

Let $X_{i},$ for $i=1,2,$ be two random variables defined on a probabilized
space $\left(\Omega,\mathcal{A},P\right),$ taking values in the respective
probabilizable spaces $\left(E_{i},\mathscr{E}_{i}\right).$ Let $f_{i},\,i=1,2,$
be two measurable functions from $\left(E_{i},\mathscr{E}_{i}\right)$
to the probabilizable space $\left(F_{i},\mathcal{F}_{i}\right).$ 

If the random variables $X_{i},\,i=1,2$ are independent, then the
random variables $f_{i}\circ X_{i},\,i=1,2$---often denoted $f_{i}\left(X_{i}\right)$---are
also independent. 

\end{remark}

\begin{example}{}{}

If the random variables $X_{i},\,i=1,2$ take values in $\mathbb{R}^{d_{i}}$
and are independent, then any marginal of $X_{1}$ is independent
of any marginal of $X_{2}.$ 

\end{example}

\begin{proposition}{Necessary and Sufficient Condition for Independence and $\pi$-systems}{}

With the previous notations, let $\mathscr{C}_{i},$ for each $i=1,2,$
be $\pi-$systems generating the respective $\mathcal{E}_{i}.$

For the random variables $X_{1}$ and $X_{2}$ to be independent,
it is necessary and sufficient for the $\pi-$systems $X^{-1}_{i}\left(\mathscr{C}_{i}\right),\,i=1,2$
to be independent.

\end{proposition}

\begin{proof}{}{}

Recall that the $\sigma-$algebras generated by the functions $X_{i},\,i=1,2$
satisfy the relations
\[
\sigma\left(X^{-1}_{i}\left(\mathscr{C}_{i}\right)\right)=X^{-1}_{i}\left(\sigma\left(\mathscr{C}_{i}\right)\right)\equiv X^{-1}_{i}\left(\mathcal{E}_{i}\right).
\]

Since the family of events $X^{-1}_{i}\left(\mathscr{C}_{i}\right)$
are $\pi-$systems, the result follows directly from Theorem $\ref{th:nsc_independence_family_events}.$

\end{proof}

We now present the \textbf{general criterion for the independence}
of two random variables in terms of their \textbf{laws}, in the following
corollary.

\begin{corollary}{General Criterion of Independence for Two Random Variables}{general_criterion_independence}

Let $X_{i},$ for $i=1,2,$ be two random variables defined on a probabilized
space $\left(\Omega,\mathcal{A},P\right)$ taking values in the respective
probabilizable spaces $\left(E_{i},\mathscr{E}_{i}\right).$

For the random variables $X_{1}$ and $X_{2}$ to be independent it
is necessary and sufficient that\boxeq{
\begin{equation}
P_{X_{1},X_{2}}=P_{X_{1}}\otimes P_{X_{2}}\label{eq:general_independence_criteria}
\end{equation}
} where $P_{X_{1}}\otimes P_{X_{2}}$ denotes the product probability
of the probabilities $P_{X_{1}}$ and $P_{X_{2}},$ laws of $X_{1}$
and $X_{2}.$

\end{corollary}

\begin{proof}{}{}

By definition, for $X_{1}$ and $X_{2}$ to be independent, it is
necessary and sufficient that
\[
\forall A_{1}\in\mathcal{E}_{1},\,\,\,\,\forall A_{2}\in\mathcal{E}_{2},\,\,\,\,P\left(X^{-1}_{1}\left(A_{1}\right)\cap X^{-1}_{2}\left(A_{2}\right)\right)=\prod^{2}_{i=1}P\left(X^{-1}_{i}\left(A_{i}\right)\right).
\]

But, for every $A_{1}\in\mathcal{E}_{1}$ and $A_{2}\in\mathcal{E}_{2},$
\begin{align*}
P\left(X^{-1}_{1}\left(A_{1}\right)\cap X^{-1}_{2}\left(A_{2}\right)\right) & =P\left(\left(X_{1}\cdot X_{2}\right)^{-1}\left(A_{1}\times A_{2}\right)\right)\\
 & =P_{X_{1},X_{2}}\left(A_{1}\times A_{2}\right)
\end{align*}
and
\[
\prod^{2}_{i=1}P\left(X^{-1}_{i}\left(A_{i}\right)\right)=P_{X_{1}}\otimes P_{X_{2}}\left(A_{1}\times A_{2}\right).
\]

It follows that $X_{1}$ and $X_{2}$ are independent if and only
if
\[
\forall A_{1}\in\mathcal{E}_{1},\,\,\,\,\forall A_{2}\in\mathcal{E}_{2},\,\,\,\,P_{X_{1},X_{2}}\left(A_{1}\times A_{2}\right)=P_{X_{1}}\otimes P_{X_{2}}\left(A_{1}\times A_{2}\right).
\]

By the uniqueness theorem for measures---Theorem $\ref{th:unicity_measures}$---this
is equivalent to
\[
P_{X_{1},X_{2}}=P_{X_{1}}\otimes P_{X_{2}},
\]
since the set of rectangles $A_{1}\times A_{2},$ with $A_{1}\in\mathcal{E}_{1}$
and $A_{2}\in\mathcal{E}_{2},$ forms a $\pi-$system generating the
product $\sigma-$algebra $\mathcal{E}_{1}\times\mathcal{E}_{2}.$

\end{proof}

This criterion takes the following functional form.

\begin{proposition}{Functional Form of the General Criterion of Independence}{func_form_gal_criterion_indep}

Let $X_{i},$ $i=1,2$ be two random variables defined on the probabilized
space $\left(\Omega,\mathcal{A},P\right)$ taking values in the respective
probabilizable spaces $\left(E_{i},\mathcal{E}_{i}\right).$

The following three assertions are equivalent:

(i) The random variables $X_{1}$ and $X_{2}$ are independent.

(ii) For any nonnegative real-valued $\mathcal{E}_{i}-$measurable
function $f_{i},$ with $i=1,2,$\boxeq{
\[
\intop_{\Omega}\left(\prod^{2}_{i=1}f_{i}\circ X_{i}\right)\text{d}P=\prod^{2}_{i=1}\intop_{\Omega}f_{i}\circ X_{i}\text{d}P.
\]
}

(iii) For any real bounded $\mathcal{E}_{i}-$measurable function
$f_{i},$ with $i=1,2,$\boxeq{
\[
\intop_{\Omega}\left(\prod^{2}_{i=1}f_{i}\circ X_{i}\right)\text{d}P=\prod^{2}_{i=1}\intop_{\Omega}f_{i}\circ X_{i}\text{d}P,
\]
}which can also be written as\boxeq{
\[
\mathbb{E}\left(f_{1}\left(X_{1}\right)f_{2}\left(X_{2}\right)\right)=\mathbb{E}\left(f_{1}\left(X_{1}\right)\right)\mathbb{E}\left(f_{2}\left(X_{2}\right)\right).
\]
}

\end{proposition}

\begin{proof}{}{}

First, note that for the functions considered in (ii) and (iii), we
have, by the transfer theorem
\[
\intop_{\Omega}\left(\prod^{2}_{i=1}f_{i}\circ X_{i}\right)\text{d}P=\intop_{E_{1}\times E_{2}}\left[\prod^{2}_{i=1}f_{i}\left(x_{i}\right)\right]\text{d}P_{X_{1},X_{2}}\left(x_{1},x_{2}\right)
\]
and
\[
\prod^{2}_{i=1}\intop_{\Omega}f_{i}\circ X_{i}\text{d}P=\prod^{2}_{i=1}\intop_{E_{i}}f_{i}\left(x_{i}\right)\text{d}P_{X_{i}}\left(x_{i}\right).
\]

We successively show the implications:
\begin{itemize}
\item $\text{(i)}\Rightarrow\text{(ii)}\,\text{and}\,\text{(iii)}:$ \\
It suffices to apply the relation $\refpar{eq:general_independence_criteria}$
together with the Fubini theorem.
\item $\text{(iii)}\Rightarrow\text{(ii)}:$ \\
The relation holds in particular for the bounded, nonnegative, $\mathcal{E}_{i}-$measurable
functions $f_{i}.$ \\
For arbitrary nonnegative $\mathcal{E}_{i}-$measurable functions
$f_{i},$ consider the sequence of bounded nonnegative $\mathcal{E}_{i}-$measurable
$\min\left(f_{i},n\right)$ and apply the Beppo Levi property---Proposition
$\ref{pr:integral_properties}.$
\item $\text{(iii)}\Rightarrow\text{(i)}:$ \\
Taking $f_{i}=\boldsymbol{1}_{A_{i}}$ where $A_{i}\in\mathcal{E}_{i},$
we immediately recover relation $\refpar{eq:general_independence_criteria},$
which is equivalent to the independence of $X_{1}$ and $X_{2}.$
\end{itemize}
\end{proof}

From the independence general criterion, we deduce the following particular
criterion.

\begin{corollary}{Independence Criterion in Terms of Cumulative Distribution Function and Densities, and for Discrete Random Variables}{}

Let $X_{i},$ $i=1,2$ be two random variables defined on the probabilized
space $\left(\Omega,\mathcal{A},P\right)$ taking values in the respective
probabilizable spaces $\left(\text{\ensuremath{\mathbb{R}}}^{d_{i}},\mathscr{B}_{\mathbb{R}^{d_{i}}}\right).$

\textbf{1. Independence criterion in terms of cumulative distribution
functions}

\boxeq{For the random variables $X_{1}$ and $X_{2}$ to be independent,
it is necessary and sufficient that
\begin{equation}
\forall x_{1}\in\mathbb{R}^{d_{1}},\,\,\,\,\forall x_{2}\in\mathbb{R}^{d_{2}},\,\,\,\,F_{X_{1},X_{2}}\left(x_{1},x_{2}\right)=F_{X_{1}}\left(x_{1}\right)F_{X_{2}}\left(x_{2}\right).\label{eq:cum_dist_function_indep_crit}
\end{equation}
}

\textbf{2. Independence criterion in terms of densities}

(a) If the random variables $X_{1}$ and $X_{2}$ admit respective
densities $f_{X_{1}}$ and $f_{X_{2}},$ and are independent, then
the random variable $\left(X_{1},X_{2}\right)$ admits a density $f_{X_{1},X_{2}}$
\textbf{\index{direct product}direct product} of $f_{X_{1}}$ and
$f_{X_{2}},$ that is such that\boxeq{
\begin{equation}
\forall x_{1}\in\mathbb{R}^{d_{1}},\,\,\,\,\forall x_{2}\in\mathbb{R}^{d_{2}},\,\,\,\,f_{X_{1},X_{2}}\left(x_{1},x_{2}\right)=f_{X_{1}}\left(x_{1}\right)f_{X_{2}}\left(x_{2}\right).\label{eq:density_direct_prod_indep_crit}
\end{equation}
}

(b) Conversely, if the random variable $\left(X_{1},X_{2}\right)$
admits a density $f_{X_{1},X_{2}}$ direct product of two nonnegative
integrable functions $f_{1}$ and $f_{2},$ that is such that it verifies\boxeq{
\begin{equation}
\forall x_{1}\in\mathbb{R}^{d_{1}},\,\,\,\,\forall x_{2}\in\mathbb{R}^{d_{2}},\,\,\,\,f_{X_{1},X_{2}}\left(x_{1},x_{2}\right)=f_{1}\left(x_{1}\right)f_{2}\left(x_{2}\right)\label{eq:density_prod_function_indep_crit}
\end{equation}
}then $f_{1}$ and $f_{2}$ are, up to a nonnegative factor, the
respective densities of $X_{1}$ and $X_{2}$ and the random variables
$X_{1}$ and $X_{2}$ are independent.

\textbf{3. Independence criterion for the discrete random variables}

If $X_{1}$ and $X_{2}$ are discrete random variables, then so is
for the random variable $\left(X_{1},X_{2}\right).$

\boxeq{The random variables $X_{1}$ and $X_{2}$ are independent
if and only if
\begin{align}
\forall x_{1} & \in X_{1}\left(\Omega\right),\,\,\,\,\forall x_{2}\in X_{2}\left(\Omega\right),\label{eq:discrete_var_independence_criterion}\\
 & P\left(X_{1}=x_{1},X_{2}=x_{2}\right)=P\left(X_{1}=x_{1}\right)P\left(X_{2}=x_{2}\right).\nonumber 
\end{align}
}

\end{corollary}

\begin{proof}{}{}

1. If the random variables $X_{1}$ and $X_{2}$ are independent,
then the relation $\refpar{eq:cum_dist_function_indep_crit}$ follows
directly from the general criterion. 

Conversely, if the relation $\refpar{eq:cum_dist_function_indep_crit}$
is satisfied, then the probabilities $P_{X_{1},X_{2}}$ and $P_{X_{1}}\otimes P_{X_{2}}$
coincide on the $\pi-$system
\[
\left\{ \left\{ u_{1}\leqslant x_{1}\right\} \times\left\{ u_{2}\leqslant x_{2}\right\} :\,x_{1}\in\mathbb{R}^{d_{1}}_{1},\,x_{2}\in\mathbb{R}^{d_{2}}_{2}\right\} 
\]
which generates the product $\sigma-$algebra $\mathscr{B}_{\mathbb{R}^{d_{1}}_{1}}\otimes\mathscr{B}_{\mathbb{R}^{d_{2}}_{2}}.$ 

By the uniqueness theorem for measures, we conclude that
\[
P_{X_{1},X_{2}}=P_{X_{1}}\otimes P_{X_{2}}
\]
 and thus the random variables $X_{1}$ and $X_{2}$ are independent.

2. Let $\lambda_{i}$ be the Lebesgue measure on $\mathbb{R}^{d_{i}},$
for $i=1,2.$ In the two considered cases, it holds
\[
P_{X_{i}}=f_{X_{i}}\cdot\lambda_{i},
\]
and thus, by the Fubini theorem,
\begin{equation}
P_{X_{1}}\otimes P_{X_{2}}=f_{X_{1}}f_{X_{2}}\cdot\lambda_{1}\otimes\lambda_{2}.\label{eq:direct_product_P_X_i}
\end{equation}

(a) Therefore, if $X_{1}$ and $X_{2}$ are independent, Corollary
$\ref{co:general_criterion_independence},$ together with the above
relation, ensures that $\left(X_{1},X_{2}\right)$ admits a density
$f_{X_{1},X_{2}},$ which is the direct product of $f_{X_{1}}$ and
$f_{X_{2}}.$

(b) Conversely, if the relation $\refpar{eq:density_prod_function_indep_crit}$
holds, then the random variable $X_{i},$ for $i=1,2,$ admit the
respective density $f_{X_{i}}$ given by
\[
\forall x_{i}\in\mathbb{R}^{d_{i}},\,\,\,\,f_{X_{i}}\left(x_{i}\right)=f_{i}\left(x_{i}\right)\intop_{\mathbb{R}^{d_{j}}}f_{j}\left(x_{j}\right)\text{d}\lambda_{j}\left(x_{j}\right),
\]
where $j=1,2$ and $i\neq j.$

Integrating this identity---for instance with $i=1$ and $j=2$---yields
\[
1=\left(\intop_{\mathbb{R}^{d_{1}}}f_{1}\left(x_{1}\right)\text{d}\lambda_{1}\left(x_{1}\right)\right)\left(\intop_{\mathbb{R}^{d_{2}}}f_{2}\left(x_{2}\right)\text{d}\lambda_{2}\left(x_{2}\right)\right).
\]

It follows that
\[
\forall x_{1}\in\mathbb{R}^{d_{1}},\,\,\,\,\forall x_{2}\in\mathbb{R}^{d_{2}},\,\,\,\,f_{X_{1}}\left(x_{1}\right)f_{X_{2}}\left(x_{2}\right)=f_{1}\left(x_{1}\right)f_{2}\left(x_{2}\right),
\]
which, by the hypothesis, implies that
\begin{equation}
P_{X_{1},X_{2}}=f_{X_{1}}f_{X_{2}}\cdot\lambda_{1}\otimes\lambda_{2}.\label{eq:P_X_1_X_2_direct_prod_lambda}
\end{equation}

Combining $\refpar{eq:direct_product_P_X_i}$ and $\refpar{eq:P_X_1_X_2_direct_prod_lambda},$
it follows that
\[
P_{X_{1},X_{2}}=P_{X_{1}}\otimes P_{X_{2}},
\]
which establishes the independence of the random variables $X_{1}$
and $X_{2}.$ 

\textbf{It is also worth noting that the relation $\refpar{eq:density_direct_prod_indep_crit}$
remains valid.}

3. If the random variables $X_{1}$ and $X_{2}$ are discrete, the
necessary condition is straightforward. 

Let us prove that this condition is sufficient. 

Assume that the relation $\refpar{eq:discrete_var_independence_criterion}$
holds. 

Then, for every $A_{i}\in\mathscr{B}_{\mathbb{R}^{d_{i}}},i=1,2,$
\[
P_{X_{1},X_{2}}\left(A_{1}\times A_{2}\right)=\sum_{\substack{x_{1}\in X_{1}\left(\Omega\right)\\
x_{2}\in X_{2}\left(\Omega\right)
}
}P\left(\left(X_{1}=x_{1}\right),\left(X_{2}=x_{2}\right)\right)\delta_{x_{1}}\otimes\delta_{x_{2}}\left(A_{1}\times A_{2}\right),
\]
and thus, by the Fubini theorem
\[
P_{X_{1},X_{2}}\left(A_{1}\times A_{2}\right)=\left(\sum_{\substack{x_{1}\in X_{1}\left(\Omega\right)}
}P\left(X_{1}=x_{1}\right)\delta_{x_{1}}\left(A_{1}\right)\right)\left(\sum_{\substack{x_{2}\in X_{2}\left(\Omega\right)}
}P\left(X_{2}=x_{2}\right)\delta_{x_{2}}\left(A_{2}\right)\right)
\]
This gives
\[
P_{X_{1},X_{2}}\left(A_{1}\times A_{2}\right)=P_{X_{1}}\left(A_{1}\right)P_{X_{2}}\left(A_{2}\right).
\]

By the theorem of uniqueness of measures---Theorem $\ref{th:unicity_measures}$---it
follows that
\[
P_{X_{1},X_{2}}=P_{X_{1}}\otimes P_{X_{2}}.
\]

This proves the independence of the random variables $X_{1}$ and
$X_{2}.$

\end{proof}

The following properties, linking independence and moments of random
variables are commonly used.

\begin{proposition}{Independence and Product. Independence and Covariance}{}

Let $X_{i},$ $i=1,2$ be two independent real-valued random variables
defined on a probabilized space $\left(\Omega,\mathcal{A},P\right).$

(a) If $X_{1}$ and $X_{2}$ admit an expectation, then the product
random variable $X_{1}X_{2}$ also admits an expectation, and\boxeq{
\[
\mathbb{E}\left(X_{1}X_{2}\right)=\mathbb{E}\left(X_{1}\right)\mathbb{E}\left(X_{2}\right).
\]
}

(b) If $X_{1}$ and $X_{2}$ admit second-order moments, then\boxeq{
\[
\text{cov}\left(X_{1},X_{2}\right)=0\,\,\,\,\text{and}\,\,\,\,\sigma^{2}_{X_{1}+X_{2}}=\sigma^{2}_{X_{1}}+\sigma^{2}_{X_{2}}.
\]
}

\end{proposition}

\begin{proof}{}{}

(a) By Proposition $\ref{pr:func_form_gal_criterion_indep},$
\[
\intop_{\Omega}\left|X_{1}X_{2}\right|\text{d}P=\left(\intop_{\Omega}\left|X_{1}\right|\text{d}P\right)\left(\intop_{\Omega}\left|X_{2}\right|\text{d}P\right)<+\infty.
\]

Thus, the product $X_{1}X_{2}$ admits an expectation.

By the transfer theorem,
\[
\intop_{\Omega}X_{1}X_{2}\text{d}P=\intop_{\mathbb{R}}x_{1}x_{2}\text{d}P_{X_{1}X_{2}}\left(x_{1},x_{2}\right).
\]

The independence of the random variables $X_{1}$ and $X_{2}$ is
equivalent to the equality $P_{X_{1}X_{2}}=P_{X_{1}}\otimes P_{X_{2}}.$
Hence, by the Fubini theorem,
\begin{align*}
\mathbb{E}\left(X_{1}X_{2}\right) & =\intop_{\Omega}X_{1}X_{2}\text{d}P\\
 & =\intop_{\mathbb{R}}x_{1}x_{2}\text{d}\left(P_{X_{1}}\otimes P_{X_{2}}\right)\left(x_{1},x_{2}\right)\\
 & =\prod^{2}_{i=1}\intop_{\mathbb{R}}x\text{d}P_{X_{i}}\left(x\right)\\
 & =\mathbb{E}\left(X_{1}\right)\mathbb{E}\left(X_{2}\right).
\end{align*}

(b) Applying this last result and using the general formula for the
variance of a sum of random variables---Proposition $\ref{pr:covariance_prop}$---we
obtain
\[
\text{cov}\left(X_{1},X_{2}\right)=\mathbb{E}\left(X_{1}X_{2}\right)-\mathbb{E}\left(X_{1}\right)\mathbb{E}\left(X_{2}\right)=0,
\]
and thus
\[
\sigma^{2}_{X_{1}+X_{2}}=\sigma^{2}_{X_{1}}+\sigma^{2}_{X_{2}}+2\text{cov}\left(X_{1},X_{2}\right)=\sigma^{2}_{X_{1}}+\sigma^{2}_{X_{2}}.
\]

\end{proof}

\begin{corollary}{Independence and Covariance Operators and Covariance Matrix}{}

Let $X_{i}\in\mathscr{L}^{2}_{E}\left(\Omega,\mathcal{A},P\right),\,i=1,2,$
where $E$ is a Euclidean space.

If the random variables $X_{1}$ and $X_{2}$ are \textbf{independent},
then the covariance operators satisfy\boxeq{
\[
\Lambda_{X_{1}+X_{2}}=\Lambda_{X_{1}}+\Lambda_{X_{2}}.
\]
}

Moreover, if $E=\mathbb{R}^{d},$ this relation expresses in term
of the covariance matrices as\boxeq{
\[
C_{X_{1}+X_{2}}=C_{X_{1}}+C_{X_{2}}.
\]
}

\end{corollary}

\begin{proof}{}{}

For every $x\in E,$ the random variables $\left\langle X_{1},x\right\rangle $
and $\left\langle X_{2},x\right\rangle $ are independent.

Thus,

$\left\langle \Lambda_{X_{1}+X_{2}}x,x\right\rangle \equiv\sigma^{2}_{\left\langle X_{1}+X_{2},x\right\rangle }=\sigma^{2}_{\left\langle X_{1},x\right\rangle }+\sigma^{2}_{\left\langle X_{2},x\right\rangle }\equiv\left\langle \Lambda_{X_{1}}x,x\right\rangle +\left\langle \Lambda_{X_{2}}x,x\right\rangle ,$
hence, the result.

\end{proof}

\textbf{We now generalize the concept of independence of events, family
of events, random variables of arbitrary families indexed on a set
$I.$}

\begin{definition}{Independence and Events of an Arbitrary Family}{}

Let $\left(A_{i}\right)_{i\in I}$ be a family of events. The events
$A_{i},i\in I$ are \textbf{independent\mindex{independent ! family of events}}
if\boxeq{
\[
\forall j\in\mathscr{P}_{f}\left(I\right),\,\,\,\,P\left(\bigcap_{j\in J}A_{j}\right)=\prod_{j\in J}P\left(A_{j}\right),
\]
}where $\mathscr{P}_{f}\left(I\right)$ denotes the set of all finite
subsets of $I.$

\end{definition}

\begin{remark}{}{}

In this context, we also speak of a ``family of independent events.''
This expresses a notion of global independence. The reader may refer
to Part \ref{part:Introduction-to-Probability} Chapter \ref{chap:Independence}
for the concept of $n$ to $n$ independance and its relations with
this last.

\end{remark}

\begin{definition}{Independence of an Arbitrary Family of Families of Events}{}

Let $\left(\mathcal{A}_{i}\right)_{i\in I}$ be a family of families
of events. 

We say that the family of events $\mathcal{A}_{i},\,i\in I$ are independent
if, for any choice of events $A_{i}$ in $\mathcal{A}_{i},\,i\in I,$
the events $A_{i},\,i\in I$ are independent.

\end{definition}

\begin{remark}{}{}

In this definition, the families $\mathcal{A}_{i}$ are not assumed
to have any specific structure. Nonetheless, this notion is particularly
useful when these families of events are $\pi-$systems or are $\sigma-$algebras. 

For instance, one may consider a sequence of sub-$\sigma-$algebras
of $\mathcal{A}$ that are independent. \textbf{Theorem $\ref{th:nsc_independence_family_events}$
can then be extended to arbitrary families of independent $\pi-$systems}.

\end{remark}

\begin{theorem}{Independence of $\sigma$-algebras Generated by Union of $\pi-$systems}{indep_salg_gen_union_pis}

Let $\left(\mathscr{C}_{i}\right)_{i\in I}$ be a family of $\pi-$systems
contained in $\mathcal{A},$ assumed to be independent. 

Denote $\mathscr{F}_{i}$ the $\sigma-$algebra generated by $\mathscr{C}_{i},$
for each $i\in I.$

Let $\left\{ I_{j}\right\} _{j\in J}$ be an arbitrary partition of
$I.$

For every $j\in J,$ let $\mathcal{A}_{j}$ be the $\sigma-$algebra
generated by the family of events $\bigcup_{i\in I_{j}}\mathscr{C}_{i},$
that is the smallest $\sigma-$algebra containing $\bigcup_{i\in I_{j}}\mathscr{C}_{i}.$

Then the $\sigma-$algebras $\mathcal{A}_{j},j\in J$ are independent.

In particular, the $\sigma-$algebras generated by the $\mathscr{C}_{i},\,i\in I$
are independent.

\end{theorem}

\begin{proof}{}{}

We provide only an outline of the proof. 

It is sufficient to prove the result when $I$ is finite, say $I=\left\llbracket 1,n\right\rrbracket .$
It can be done by induction.

To avoid complications with indexing, we only show the case where
$\mathscr{C}_{i},i=1,2,3$ are independent $\pi-$systems, and we
prove that the $\sigma-$algebra generated by $\mathscr{C}_{1}$ and
$\mathscr{C}_{2}\cup\mathscr{C}_{3}$ are independent.

Note that the family of events $\mathscr{C}_{2}\cup\mathscr{C}_{3}$
is no more a $\pi-$system. 

Let $\mathscr{C}_{4}$ be the $\pi-$system generated by $\mathscr{C}_{2}\cup\mathscr{C}_{3}$
and $\Omega;$ that is
\[
\left\{ C_{2}\cap C_{3}:\,C_{2}\in\mathscr{C}_{2}\cup\left\{ \Omega\right\} \,\text{and\,}C_{3}\in\mathscr{C}_{3}\cup\left\{ \Omega\right\} \right\} .
\]

We clearly have
\[
\mathscr{C}_{2}\cup\mathscr{C}_{3}\subset\mathscr{C}_{4}\subset\sigma\left(\mathscr{C}_{2}\cup\mathscr{C}_{3}\right).
\]

And therefore, 
\[
\sigma\left(\mathscr{C}_{4}\right)=\sigma\left(\mathscr{C}_{2}\cup\mathscr{C}_{3}\right).
\]

It is straightforward that the $\pi-$systems $\mathscr{C}_{1}$ and
$\mathscr{C}_{4}$ are independent. By Theorem $\ref{th:nsc_independence_family_events},$
the $\sigma-$algebras are also independent.

\end{proof}

\begin{definition}{Independence of a Family of Random Variables}{}

A family $\left(X_{i}\right)_{i\in I}$ of random variables, taking
values respectively in the probabilizable spaces $\left(E_{i},\mathcal{E}_{i}\right),i\in I,$
is said to be \textbf{a family of independent random variables} if
the $\sigma-$algebras $X^{-1}_{i}\left(\mathcal{E}_{i}\right),i\in I,$
generated by these random variables, are independent.

\end{definition}

We then often say more concisely that ``the random variables $X_{i},\,i\in I$
are independent.''

All the criteria previously studied for families where $\left|I\right|=2$
extend easily to the case where $I$ is finite. For an arbitrary index
set $I,$ however, one must define the concept of a product measure
on $\prod_{i\in I}E_{i}.$ This is possible, for instance, when $E_{i}=\mathbb{R}$
for every $i\in I,$ by applying the Kolmogorov extension theorem---see
Chapter \ref{chap:PartIIChap8} Theorem $\ref{th:kolmogorov_extension_th}$
and Corollary $\ref{co:infinite_product_proba}$ for a statement.

\section{Independence and Asymptotic Events}\label{sec:Independence-and-Asymptotic}

In this Section, we study two well-known theorems that are commonly
used in the analysis of almost sure convergence for sequences or series
of random variables---a topic we will examine later.

\begin{definition}{Asymptotic $\sigma-$algebra. Asymptotic Events. Asymptotic Random Variable}{}

Let $\left(\mathcal{A}_{n}\right)_{n\in\mathbb{N}}$ be a sequence
of sub-$\sigma-$algebras of the $\sigma-\text{algebra}$ $\mathcal{A}$
defined on a set $\Omega.$

For each $n\in\mathbb{N},$ denote $\bigvee_{p\geqslant n}\mathcal{A}_{p}$
the $\sigma-$algebra generated by $\bigcup_{p\geqslant n}\mathcal{A}_{p},$
i.e., the smallest $\sigma-\text{algebra}$---with respect to families
of subsets inclusion---that contains $\bigcup_{p\geqslant n}\mathcal{A}_{p}.$

The \mindex{sigma-algebra @ $\sigma-$algebra ! asymptotic}\mindex{asymptotic!sigma-algebra @ asymtotic!$\sigma-$algebra}\textbf{asymptotic
$\sigma-\text{algebra},$} denoted $\mathcal{A}_{\infty},$ is defined
by
\[
\mathcal{A}_{\infty}=\bigcap_{n\in\mathbb{N}}\left(\bigvee_{p\geqslant n}\mathcal{A}_{p}\right).
\]

The elements of $\mathcal{A}_{\infty}$ are called the \textbf{asymptotic
events}\mindex{asymptotic!events}\mindex{events ! asymptotic}. And
any \textbf{random variable} $\mathcal{A}_{\infty}-$measurable is
called an\mindex{asymptotic!random variable}\mindex{random variable ! asymptotic}\textbf{
asymptotic} \textbf{random variable}.

\end{definition}

\begin{example}{Asympotic $\sigma-$algebra Associated to a Sequence of Random Variables}{asymptotic_seq_rv}

Let $\left(X_{n}\right)_{n\in\mathbb{N}}$ be a sequence of random
variables defined on $\left(\Omega,\mathcal{A},P\right)$ taking values
in the respective probabilizable space $\left(E_{n},\mathcal{E}_{n}\right).$

For each $p\in\mathbb{N},$ we define $\mathcal{A}_{p}$ as the $\sigma-\text{algebra}$
$X^{-1}_{p}\left(\mathcal{E}_{p}\right)$ generated by the random
variable $X_{p},$ also denoted $\sigma\left(X_{p}\right).$ This
is the $\sigma-\text{algebra}$ of events of the form $\left(X_{p}\in A_{p}\right),\,\left(A_{p}\in\mathcal{E}_{p}\right);$
i.e., for a given $\omega,$ as soon as we know the value $X_{p}\left(\omega\right)$
we can determine whether the event occurs or not.

The $\sigma-$algebra $\bigvee_{p\geqslant n}X^{-1}_{p}\left(\mathcal{E}_{p}\right),$
the smallest \salg on $\Omega$ making simultaneously measurable
all the functions $X_{p},\,p\geqslant n,$ and still denoted $\sigma\left(X_{p}:\,p\geqslant n\right),$
is, by definition, the \salg generated by the family of random variables
$\left(X_{p}\right)_{p\geqslant n}.$ It is constituted of events
whose realization, for a given $\omega,$ depends on the sequence
$\left(X_{p}\left(\omega\right),X_{p+1}\left(\omega\right),\cdots\right).$
Nonetheless, this dependence is not necessarily explicitly described. 

Then, the asymptotic \salg is
\[
\mathcal{A}_{\infty}=\bigcap_{n\in\mathbb{N}}\sigma\left(X_{p}:\,p\geqslant n\right).
\]

An event belongs to $\mathcal{A}_{\infty}$ if its outcome, for a
given $\omega$---even if it depends on the sequence $\left(X_{i}\left(\omega\right)\right)_{i\in\mathbb{N}}$---does
not depend on the first $n$ values, regardless the value of the integer
$n$ is. This is the case, for instance, when the $X_{n}$ take real
values, of the event
\[
\left\{ \text{the sequence }\left(X_{n}\right)_{n\in\mathbb{N}}\text{ converges in }\mathbb{R}\right\} .
\]

\end{example}

\begin{solutionexample}{}{}

Let us show that this event is precisely asymptotic. 

First, recall that a numerical sequence $\left(a_{n}\right)_{n\in\mathbb{N}}$
converges in $\overline{\mathbb{R}}$ if and only if
\[
\liminf_{n\to+\infty}a_{n}=\limsup_{n\to+\infty}a_{n}.
\]

Thus, such a sequence converges in $\mathbb{R}$ if and only if
\[
-\infty<\liminf_{n\to+\infty}a_{n}=\limsup_{n\to+\infty}a_{n}<+\infty.
\]

Hence, we have the equality of events
\[
\left\{ \text{the sequence }\left(X_{n}\right)_{n\in\mathbb{N}}\text{ converges in }\mathbb{R}\right\} =\left\{ \liminf_{n\to+\infty}X_{n}=\limsup_{n\to+\infty}X_{n}\in\mathbb{R}\right\} .
\]

Thus, it suffices to prove that the random variables $\limsup_{n\to+\infty}X_{n}$
and $\liminf_{n\to+\infty}X_{n}$ are asymptotic. It suffices to cover
the case of $\limsup_{n\to+\infty}X_{n}.$ 

To this aim, let $\mathscr{B}_{n}$ be the \salg $\sigma\left(X_{p}:\,p\geqslant n\right).$
We note that, for every $p$ and for every $n\geqslant p,$ $\sup_{k\geqslant n}X_{k}$
is $\mathscr{B}_{p}-\text{measurable.}$ It follows that, for each
$p,$ $\lim_{n\to+\infty}\left(\sup_{k\geqslant n}X_{k}\right)$ is
$\mathscr{B}_{p}-$measurable, that is $\limsup_{n\to+\infty}X_{n}$
is $\mathcal{A}_{\infty}-$measurable.

Then, we deduce immediately that the event $\left(\text{the series }\sum X_{n}\text{ converges in }\mathbb{R}\right)$
is also asymptotic. Another argument is also to say that, by the Cauchy
criterion, we have
\[
\left(\sum_{n\in\mathbb{N}}X_{n}\text{ converges}\right)=\left(\lim_{n,m\to+\infty}\sum^{m}_{k=n}X_{k}=0\right)\in\bigcap_{n\in\mathbb{N}}\mathscr{B}_{n}=\mathcal{A}_{\infty}.
\]

\end{solutionexample}

\begin{theorem}{Law of the All or Nothing or Law of 0,1}{}

On the probabilized space \preds, let $\left(\mathcal{A}_{n}\right)_{n\in\mathbb{N}}$
be a sequence of independent sub-$\sigma-$algebras of the $\sigma-$algebra
$\mathcal{A}$, and let $\mathcal{A}_{\infty}$ be the associated
asymptotic $\sigma-$algebra.

Then\boxeq{
\[
\forall A\in\mathcal{A}_{\infty},\,\,\,\,P\left(A\right)=0\text{ or }1.
\]
}

\end{theorem}

\begin{proof}{}{}

By Theorem $\ref{th:indep_salg_gen_union_pis},$ for any nonnegative
integers $M$ and $N$ such that $M<N,$ the \salg $\mathcal{A}_{n},\,n\leqslant M,$
and $\bigvee_{p\geqslant N}\mathcal{A}_{p}$ are independent. 

But for every $N,$ 
\[
\mathcal{A}_{\infty}\subset\bigvee_{p\geqslant N}\mathcal{A}_{p},
\]
Thus, for every $M,$ the \salg $\mathcal{A}_{n},\,n\leqslant M,$
and $\mathcal{A}_{\infty}$ are independent. Similarly, the \salg
$\mathcal{A}_{n},\,n\in\mathbb{N}$ and $\mathcal{A}_{\infty}$ are
independent. But, then the \salg $\bigvee_{p\geqslant1}\mathcal{A}_{p}$
and $\mathcal{A}_{\infty}$ are independent, thus, since
\[
\mathcal{A}_{\infty}\subset\bigvee_{p\geqslant1}\mathcal{A}_{p},
\]
the \salg is independent of itself. 

In particular,
\[
\forall A\in\mathcal{A}_{\infty},\,\,\,\,P\left(A\cap A\right)=P\left(A\right)P\left(A\right),
\]
hence, the result.

\end{proof}

\begin{corollary}{}{}

With the same notations as in the previous theorem, any $\mathcal{A}_{\infty}-$measurable
random variable is almost surely constant.

\end{corollary}

\begin{example}{}{}

Let $\left(X_{n}\right)_{n\in\mathbb{N}}$ be a sequence of independent
real-valued random variables defined on \preds.

Then the series $\sum_{n\in\mathbb{N}}X_{n}$ converges or diverges
almost surely, that is
\[
P\left(\sum_{n\in\mathbb{N}}X_{n}\text{ converges}\right)=0\text{ or }1.
\]

Indeed, we saw in Example $\ref{ex:asymptotic_seq_rv}$ that the event
$\left(\sum_{n\in\mathbb{N}}X_{n}\text{ converges}\right)$ is asymptotic.
The previous theorem yields the result.

\end{example}

\begin{remark}{}{}

Consequentely, if $\left(X_{n}\right)_{n\in\mathbb{N}}$ is a sequence
of independent real-valued random variables defined on \preds, then
to prove that the series $\sum_{n\in\mathbb{N}}X_{n}$ converges almost
surely, it is sufficient to show that
\[
P\left(\sum_{n\in\mathbb{N}}X_{n}\text{ converges}\right)>0.
\]

\end{remark}

\begin{remark}{}{}

A special case of the previous example is when 
\[
X_{n}=\dfrac{\epsilon_{n}}{n},
\]
where the random variables $\epsilon_{n},\,n\in\mathbb{N}$ are independent,
take values -1 or 1, following the same law given by
\[
P\left(\epsilon_{n}=1\right)=P\left(\epsilon_{n}=-1\right)=\dfrac{1}{2}.
\]

By the law of the all or nothing, exactly one of the following two
assertions is true---without saying which one:

(i) The series with general term $\dfrac{\epsilon_{n}}{n}$ converges
$P-$almost surely.

(ii) The series with general term $\dfrac{\epsilon_{n}}{n}$ digerges
$P-$almost surely.

\end{remark}

\begin{figure}[t]
\begin{center}\includegraphics[width=0.4\textwidth]{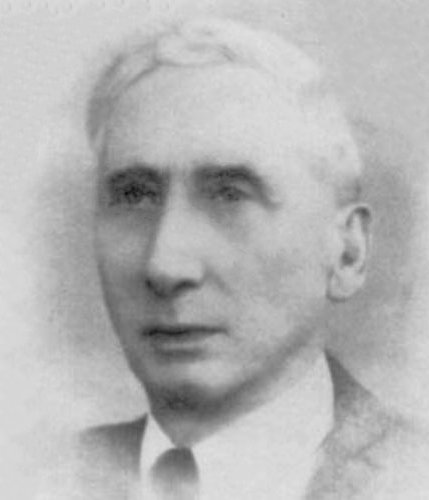}

{\tiny Public domain}\end{center}

\caption{\textbf{\protect\href{https://en.wikipedia.org/wiki/Francesco_Paolo_Cantelli}{Francesco Paolo Cantelli}}
(1875 - 1966)}\sindex[fam]{Cantelli, Francesco}
\end{figure}

\begin{reminder}{}{}

Let $\left(A_{n}\right)_{n\in\mathbb{N}}$ be a sequence of subsets
of $\Omega.$ We define the two subsets of $\Omega,$ superior and
inferior limit\footnotemark of the sequence of sets $\left(A_{n}\right)_{n\in\mathbb{N}}$
by
\[
\limsup_{n\to+\infty}A_{n}=\bigcap_{n\in\mathbb{N}}\bigcup_{p\geqslant n}A_{p}\,\,\,\,\,\,\,\,\liminf_{n\to+\infty}A_{n}=\bigcup_{n\in\mathbb{N}}\bigcap_{p\geqslant n}A_{p}.
\]

The set $\limsup_{n\to+\infty}A_{n}$ consists of the $\omega\in\Omega$
that belong to infinitely many of the $A_{n}.$ 

The set $\liminf_{n\to+\infty}A_{n}$ consists of the $\omega\in\Omega$
that belong from some rank onward---depending on $\omega$---belong
to all the $A_{n}.$

We thus always have the inclusion\boxeq{
\[
\liminf_{n\to+\infty}A_{n}\subset\limsup_{n\to+\infty}A_{n}.
\]
}

Moreover, the following equalities hold\boxeq{
\[
\left(\limsup_{n\to+\infty}A_{n}\right)^{c}=\liminf_{n\to+\infty}\left(A^{c}_{n}\right)\,\,\,\,\text{and}\,\,\,\,\left(\liminf_{n\to+\infty}A_{n}\right)^{c}=\limsup_{n\to+\infty}A^{c}_{n}.
\]
}

This sets play a key role in the study of convergence of sequences
of random variables. The fundamental tool is the \textbf{\sindex[fam]{Borel, Emile}Borel-Cantelli\sindex[fam]{Cantelli, Francesco}\footnotemark
lemma}\index{Borel-Cantelli lemma}, stated next.

\end{reminder}

\footnotetext{\textbf{\href{https://en.wikipedia.org/wiki/Francesco_Paolo_Cantelli}{Francesco Paolo Cantelli}}
(1875 - 1966) was an Italian mathematician. He contributed to celestial
mechanics, probability theory and actuarial science. He started his
carreer at the Pension Fund for the Government Deposits and Loans
Bank in Italy, before becomin in 1923, professor of actuarial mathematics,
first at Catania University then at Naples and Rome Sapienza Universities.}

\footnotetext{Other notations: $\limsup_{n\to+\infty}A_{n}=\overline{\lim}_{n\to+\infty}A_{n}$
and $\liminf_{n\to+\infty}A_{n}=\underline{\lim}_{n\to+\infty}A_{n}.$}

\begin{lemma}{Borel-Cantelli Lemma}{borel-cantelli}

Let $\left(A_{n}\right)_{n\in\mathbb{N}}$ be a sequence of events.

(a) We have the implication\boxeq{
\[
\sum^{+\infty}_{n=0}P\left(A_{n}\right)<+\infty\Rightarrow P\left(\limsup_{n\to+\infty}A_{n}\right)=0.
\]
}

(b) If the events $A_{n}$ are independent, then\boxeq{
\[
\sum^{+\infty}_{n=0}P\left(A_{n}\right)=+\infty\Rightarrow P\left(\limsup_{n\to+\infty}A_{n}\right)=1.
\]
}

\end{lemma}

\begin{proof}{}{}

(a) For every $n\in\mathbb{N},$
\[
P\left(\limsup_{n\to+\infty}A_{n}\right)\leqslant P\left(\bigcup_{p\geqslant n}A_{p}\right)\leqslant\sum_{p\geqslant n}P\left(A_{p}\right),
\]
which yields the result, the second member being the remainder of
order $n$ of a convergent series.

(b) We start by writing
\[
P\left(\limsup_{n\to+\infty}A_{n}\right)=1-P\left(\liminf_{n\to+\infty}\left(A^{c}_{n}\right)\right).
\]

Since the events $A^{c}_{n}$ are also independent, we have, from
the properties of a probability for monotonic sequences---properties
said of sequential monotonicity---
\begin{align*}
P\left(\liminf_{n\to+\infty}\left(A^{c}_{n}\right)\right) & =\lim_{n\to+\infty}\lim_{q\to+\infty}P\left(\bigcap^{q}_{p=n}A^{c}_{p}\right)\\
 & =\lim_{n\to+\infty}\lim_{q\to+\infty}\prod^{q}_{p=n}P\left(A^{c}_{p}\right)\\
 & =\lim_{n\to+\infty}\lim_{q\to+\infty}\prod^{q}_{p=n}\left(1-P\left(A_{p}\right)\right).
\end{align*}

Now, since $\text{e}^{-x}\geqslant1-x,$ we obtain
\[
0\leqslant\prod^{q}_{p=n}\left(1-P\left(A_{p}\right)\right)\leqslant\text{e}^{-\sum^{q}_{p=n}P\left(A_{p}\right)}.
\]

As the sum $\sum^{q}_{p=n}P\left(A_{p}\right)$ diverges, the right-hand
side converges to 0 when $q$ tends to $+\infty,$ so
\[
\lim_{q\to+\infty}\prod^{q}_{p=n}\left(1-P\left(A_{p}\right)\right)=0,
\]

Hence,
\[
P\left(\liminf_{n\to+\infty}\left(A^{c}_{n}\right)\right)=0,
\]
so
\[
P\left(\limsup_{n\to+\infty}A_{n}\right)=1.
\]

\end{proof}

\begin{remark}{}{}

The first implication of the Borel-Cantelli lemma always holds. However,
its converse is not generally true. As a counterexample, consider
the probabilized space $\left(\left[0,1\right],\mathscr{B}_{\left[0,1\right]},\lambda\right),$
where $\lambda$ denotes the restriction of the Lebesgue measure to
$\left[0,1\right].$ Define, for every $n\in\mathbb{N}^{\ast},$ $A_{n}=\left]0,\dfrac{1}{n}\right].$

We have 
\[
\limsup_{n\to+\infty}A_{n}=\emptyset,
\]
 and thus 
\[
P\left(\limsup_{n\to+\infty}A_{n}\right)=0.
\]
However, 
\[
\sum^{+\infty}_{n=1}P\left(A_{n}\right)=+\infty.
\]

This example clearly shows that the second implication of the Borell-Cantelli
lemma does not hold in general without the additional assumption of
independence. 

\end{remark}

\section{Some Results Related to Independence and the Heads and Tails Model}\label{sec:Some-Results-Linked}

We have seen how to construct a probabilistic model for a game of
heads and tails over $n$ tosses with a fair coin. Mathematically,
this is similar to construct a probabilized space \preds and $n$
independent random variables $X_{i},\,i\in\left\llbracket 1,n\right\rrbracket $
following the same uniform law on $\left[0,1\right[,$ defined on
this space. One can take $\Omega=\left\{ 0,1\right\} ^{n},$ equipped
with the uniform probability, and define $X_{i}$ as the projection
on the $i-$th coordinate,
\[
\left(x_{1},\cdots,x_{i},\cdots,x_{n}\right)\mapsto x_{i}.
\]

The analogous problem in the case of an infinite sequence of tosses
would naturally lead to choose as outcome space $\Omega=\left\{ 0,1\right\} ^{\mathbb{N}^{\ast}}.$
Nonetheless, the existence on this space of a probability under which,
for each $i,$ the $i-$th marginal follows the uniform law on $\left\{ 0,1\right\} $
is not straightforward. If we want to obtain such a probability by
invoking a general result, we must appeal to the Kolmogorov theorem---see
Corollary $\ref{co:infinite_product_proba}.$

Nevertheless, the problem in question admits another solution, that,
at first glance, seems more elementary: we take for $\Omega$ the
interval $\left[0,1\right[,$ equipped with the Lebesgue measure.
This solution, which we study in what follows, is in fact not so far
from the previous one: to each $x\in\left[0,1\right[,$ we can indeed
associate a sequence belonging to $\left\{ 0,1\right\} ^{\mathbb{N}^{\ast}}$
by writing the binary expansion of $x$ in base 2 
\[
x=0,x_{1}x_{2}\cdots x_{n}\cdots.
\]

The $x_{i},i\in\mathbb{N}^{\ast}$ are called the dyadic digits of
$x.$

\begin{remark}{}{}

There is an ambiguity for the rational numbers of the form $\dfrac{p}{2^{q}},$
since such numbers have two binary expansions. 

For instance,
\[
\dfrac{1}{2}=0.1000...=0.01111....
\]

Similarly, rational numbers of the form $\dfrac{p}{10^{q}}$ have
two decimal expansions. 

For instance,
\[
\dfrac{7}{10}=0.7000...=0.6999...
\]

Here, we agree to select the expansion that ends with an infinite
sequence of zeros. 

\end{remark}

Hence, on $\Omega=\left[0,1\right[,$ equipped with its Borel \salg
and the Lebesgue measure, we define a sequence of random variables
$D_{n},\,n\in\mathbb{N}^{\ast}$ by taking for $D_{n}\left(x\right)$
the $n-$th dyadic digit of $x.$ We then show---see Proposition
$\ref{pr:infinite_sequence_bernoulli_rv}$---that these random variables
are independent and follow the uniform law on $\left\{ 0,1\right\} ,$
which provides an answer to the stated problem.

\subsection{Dyadic Expansion of a Real Number in $\left[0,1\right[$}

Define for $x\in\left[0,1\right[,$ the sequences of general term
$D_{n}\left(x\right)$ and $R_{n}\left(x\right)$ by
\[
R_{0}\left(x\right)=x
\]
 and, for $n\in\mathbb{N}^{\ast},$
\[
D_{n}=\left[2R_{n-1}\left(x\right)\right]\,\,\,\,\,\,\,\,R_{n}\left(x\right)=2R_{n-1}\left(x\right)-D_{n}\left(x\right).
\]

By construction,
\[
D_{n}\left(x\right)\in\left\{ 0,1\right\} \,\,\,\,\text{and}\,\,\,\,R_{n}\left(x\right)\in\left[0,1\right[.
\]

An immediate induction shows that
\[
\forall n\in\mathbb{N}^{\ast},\,\,\,\,x=\sum^{n}_{j=1}\dfrac{D_{j}\left(x\right)}{2^{j}}+\dfrac{1}{2^{n}}R_{n}\left(x\right).
\]

When $n$ tends to $+\infty,$ we obtain
\begin{equation}
x=\sum^{+\infty}_{j=1}\dfrac{D_{j}\left(x\right)}{2^{j}}.\label{eq:dyadic_dev_x}
\end{equation}

In general, when we have
\begin{equation}
x=\sum^{n}_{j=1}\dfrac{x_{j}}{2^{j}},\label{eq:general_dyadic_dev_x}
\end{equation}
we write symbolically
\begin{equation}
x=0,x_{1}x_{2}\cdots x_{j}\cdots\label{eq:symbolic_dyadic_dev_x}
\end{equation}
and say that the right-hand side of $\refpar{eq:general_dyadic_dev_x}$
or $\refpar{eq:symbolic_dyadic_dev_x}$ is a \textbf{\index{dyadic expansion}dyadic
expansion} of $x.$

We have just shown that every $x\in\left[0,1\right[$ admits a dyadic
expansion.

The dyadic expansion of a real number is not unique; indeed,
\begin{equation}
\forall n\in\mathbb{N}^{*},\,\,\,\,\dfrac{1}{2^{n-1}}=\sum^{+\infty}_{j=n}\dfrac{1}{2^{j}}.\label{eq:1/2^n-1}
\end{equation}

Hence, for every $k\in\mathbb{N}^{\ast}$ and for every finite sequences
$\left(x_{j}\right)_{1\leqslant j\leqslant k}\in\left\{ 0,1\right\} ^{k},$
\begin{equation}
\sum^{k}_{j=1}\dfrac{x_{j}}{2^{j}}+\dfrac{1}{2^{k+1}}=\sum^{k}_{j=1}\dfrac{x_{j}}{2^{j}}+\sum^{+\infty}_{j=k+2}\dfrac{1}{2^{j}}.\label{eq:two_dyadic_dev}
\end{equation}

Hence,
\begin{equation}
0.x_{1}x_{2}\cdots x_{k}10000\cdots=0.x_{1}x_{2}\cdots x_{k}01111\cdots\label{eq:two_dyadic_dev_symbolic}
\end{equation}

Now, since any odd integer can be written under the form 
\[
x_{1}2^{k}+x_{2}2^{k-1}+\cdots+x_{k}2+1
\]
---which corresponds to writing an integer in base 2---we can easily
verify that the numbers of the form $\refpar{eq:two_dyadic_dev}$
are exactly the rational numbers of the form $\dfrac{p}{2^{q}},$
called \textbf{dyadic rational numbers}\index{dyadic rational numbers},
from the interval $\left]0,1\right[.$ Thus, we have now two expansions
for the dyadic rational numbers: 
\begin{itemize}
\item The first, is a finite expansion, and ends with an infinite sequence
of $0,$
\item The second finishes with an infinite sequence of $1.$
\end{itemize}
There is no other case of non-uniqueness. Indeed, suppose there were.
Then
\[
0.x_{1}x_{2}\cdots x_{j}\cdots=0.y_{1}y_{2}\cdots y_{j}\cdots
\]

Let $k$ be the first index such that $x_{j}\neq y_{j}.$ Without
loss of generality, suppose $x_{k}=1$ and $y_{k}=0.$ Then
\[
\dfrac{1}{2^{k}}+\sum^{+\infty}_{j=k+1}\dfrac{x_{j}}{2^{j}}=\sum^{+\infty}_{j=k+1}\dfrac{y_{j}}{2^{j}}.
\]

From the relation $\refpar{eq:1/2^n-1},$ the only possibility is
that for every $j\geqslant k+1:$ $x_{j}=0$ and $y_{j}=1.$ That
is, we are in the dyadic rational case $\refpar{eq:two_dyadic_dev}.$

Finally, we can retain that for $x\in\left[0,1\right[,$ the sequence
with general term $D_{n}\left(x\right)$ gives the digit of the dyadic
expansion of $x$ when this one is unique. When $x\in\left]0,1\right[$
is a dyadic rational number, this sequence yields the digit of the
finite expansion. Indeed, a simple computation, setting $d_{n}=D_{n}\left(x\right)$
for brevity, shows that
\[
R_{n}\left(x\right)=0.d_{1}d_{2}\cdots d_{n}\cdots
\]
Since $0.111\cdots=1$ and $R_{n}\left(x\right)<1,$ the expansion
$x=0.d_{1}d_{2}\cdots d_{n}\cdots$ cannot contain an infinite sequence
of 1s from some rank onward.

\begin{proposition}{}{infinite_sequence_bernoulli_rv}

Let $\left(\left[0,1\right[,\mathscr{B}_{\left[0,1\right[},P\right)$
be a probabilized space where $P$ is the restriction of the Lebesgue
measure to $\left[0,1\right[.$

Define for $x\in\left[0,1\right[$ the sequences of general term $D_{n}\left(x\right)$
and $R_{n}\left(x\right)$ by
\[
R_{0}\left(x\right)=x
\]
 and, for $n\in\mathbb{N}^{\ast},$
\[
D_{n}=\left[2R_{n-1}\left(x\right)\right]\,\,\,\,\,\,\,\,R_{n}\left(x\right)=2R_{n-1}\left(x\right)-D_{n}\left(x\right).
\]

Then, on this probabilized space, the sequence $\left(D_{n}\right)_{n\in\mathbb{N}^{\ast}}$
is a sequence of independent random variables, each following the
same Bernoulli law $\mathcal{B}\left(1,\dfrac{1}{2}\right)\equiv\dfrac{1}{2}\left(\delta_{0}+\delta_{1}\right).$

Moreover, for every $n\in\mathbb{N}^{\ast},$ the random variable
$R_{n}$ follows the uniform law on $\left[0,1\right[,$ and the random
variables $R_{n}$ and $\left(D_{1},D_{2},\cdots,D_{n}\right)$ are
independent.

\end{proposition}

\begin{figure}

\begin{centering}
    \begin{tikzpicture}[every edge/.style={shorten <=1pt, shorten >=1pt}]
    \draw (0,0)  node [below=3mm] {0} -- (2,0) node [below=3mm] {$\frac{1}{4}$} -- (4,0) node [below=3mm] {$\frac{1}{2}$} -- (6,0) node [below=3mm] {$\frac{3}{4}$} -- (8,0) node [below=3mm] {1};
	\draw[>={Bracket[width=6mm,line width=0.5pt,length=1.5mm]},<-] (0, 0) -- (2, 0);
\draw[>={Bracket[width=6mm,line width=0.5pt,length=1.5mm]},<-] (2, 0) -- (4, 0);
\draw[>={Bracket[width=6mm,line width=0.5pt,length=1.5mm]},<-] (4, 0) -- (6, 0);
\draw[>={Bracket[width=6mm,line width=0.5pt,length=1.5mm]},<-] (6, 0) -- (8, 0);
\draw[>={Bracket[width=6mm,line width=0.5pt,length=1.5mm]},<-] (8, 0) -- (8.001, 0);
    % draw the tick marks
    \coordinate (p) at (0,2pt);
    \foreach \myprop/\mytext [count=\n] in {2/$I_{00}^2$,2/$I_{01}^2$,2/$I_{10}^2$,2/$I_{11}^2$}
    \draw [decorate,decoration={brace,amplitude=4,raise=10pt}] (p)   +(0,0pt) -- ++(\myprop,0) coordinate (p) node [midway, above=5mm, anchor=south] {\mytext} ;
    %\path (10,2pt) edge [draw]  ++(0,-4pt);
    \end{tikzpicture}
    
\end{centering}

\caption{Dyadic intervals for $n=2$}
\label{fig:dyadic_interval_n2}
\end{figure}

\begin{proof}{}{}
\begin{itemize}
\item For every $n\in\mathbb{N}^{\ast}$ and for every $n-$uples $\uuline{\epsilon_{n}}=\left(\epsilon_{1},\epsilon_{2},\cdots,\epsilon_{n}\right)\in\left\{ 0,1\right\} ^{n},$
denote by $I^{n}_{\uuline{\epsilon_{n}}}$ the dyadic interval, defined
by
\[
I^{n}_{\uuline{\epsilon_{n}}}=\left[\sum^{n}_{j=1}\dfrac{\epsilon_{j}}{2^{j}},\sum^{n}_{j=1}\dfrac{\epsilon_{j}}{2^{j}}+\dfrac{1}{2^{n}}\right[.
\]
This interval consists of the real numbers in $\left[0,1\right[$
whose dyadic expansion begins with $0.\epsilon_{1}\epsilon_{2}\cdots\epsilon_{n}.$\\
Figure \ref{fig:dyadic_interval_n2} shows the dyadic intervals for
$n=2.$\\
We have
\[
\bigcap^{n}_{j=1}\left(D_{j}=\epsilon_{j}\right)=I^{n}_{\uuline{\epsilon_{n}}}.
\]
Hence,
\[
P\left(\bigcap^{n}_{j=1}\left(D_{j}=\epsilon_{j}\right)\right)=\dfrac{1}{2^{n}}.
\]
But then, for every non-empty subset $J$ of $\left\llbracket 1,n\right\rrbracket ,$
we obtain, by summing over all the $\epsilon_{j}$ for $j\in J^{c},$
\[
P\left(\bigcap_{j\in J}\left(D_{j}=\epsilon_{j}\right)\right)=\dfrac{1}{2^{\left|J\right|}}.
\]
In particular, for every $j\in\left\llbracket 1,n\right\rrbracket ,$
\[
P\left(D_{j}=\epsilon_{j}\right)=\dfrac{1}{2}.
\]
We then obtain
\[
P\left(\bigcap_{j\in J}\left(D_{j}=\epsilon_{j}\right)\right)=\prod_{j\in J}P\left(D_{j}=\epsilon_{j}\right),
\]
which shows---since $n$ and $J$ are arbitrary---that the random
variables $D_{j}$ constitute a sequence of independent random variables
following the Bernoulli law $\mathscr{B}\left(1,\dfrac{1}{2}\right).$
\item Let $I$ denote the identity function $x\mapsto x$ from $\left[0,1\right[$
to $\mathbb{R}.$ We have
\[
R_{n}=2^{n}I-\sum^{n}_{j=1}2^{n-j}D_{j}.
\]
Then, for every $f\in\mathscr{C}^{+}_{\mathscr{K}}\left(\mathbb{R}\right)$
and for $\uuline{\epsilon_{n}}=\left(\epsilon_{1},\epsilon_{2},\cdots,\epsilon_{n}\right)\in\left\{ 0,1\right\} ^{n},$
\begin{align*}
\mathbb{E}\left(f\left(R_{n}\right)\prod^{n}_{j=1}\boldsymbol{1}_{\left(D_{j}=\epsilon_{j}\right)}\right) & =\mathbb{E}\left(f\left(2^{n}I-\sum^{n}_{j=1}2^{n-j}\epsilon_{j}\right)\prod^{n}_{j=1}\boldsymbol{1}_{\left(D_{j}=\epsilon_{j}\right)}\right)\\
 & =\intop_{\mathbb{R}}\boldsymbol{1}_{I^{n}_{\uuline{\epsilon_{n}}}}\left(x\right)f\left(2^{n}x-\sum^{n}_{j=1}2^{n-j}\epsilon_{j}\right)\text{d}\lambda\left(x\right).
\end{align*}
\\
Thus, by making the change of variables in the Lebesgue integral defined
by $y=2^{n}x-\sum^{n}_{j=1}2^{n-j}\epsilon_{j},$
\begin{align*}
\mathbb{E}\left(f\left(R_{n}\right)\prod^{n}_{j=1}\boldsymbol{1}_{\left(D_{j}=\epsilon_{j}\right)}\right) & =\intop_{\mathbb{R}}\boldsymbol{1}_{I^{n}_{\uuline{\epsilon_{n}}}}\left(x\right)\left(\dfrac{1}{2^{n}}y+\sum^{n}_{j=1}\dfrac{\epsilon_{j}}{2^{j}}\right)f\left(y\right)\text{d}\lambda\left(y\right).\\
 & =\dfrac{1}{2^{n}}\intop_{\mathbb{R}}\boldsymbol{1}_{\left[0,1\right[}\left(y\right)f\left(y\right)\text{d}\lambda\left(y\right).
\end{align*}
This yields
\begin{align}
\mathbb{E}\left(f\left(R_{n}\right)\prod^{n}_{j=1}\boldsymbol{1}_{\left(D_{j}=\epsilon_{j}\right)}\right) & =\left[P\left(\bigcap^{n}_{j=1}\left(D_{j}=\epsilon_{j}\right)\right)\right]\left[\intop_{\mathbb{R}}\boldsymbol{1}_{\left[0,1\right[}\left(y\right)f\left(y\right)\text{d}\lambda\left(y\right)\right].\label{eq:expect_f_r_n}
\end{align}
And thus, by summing on $\uuline{\epsilon_{n}}\in\left\{ 0,1\right\} ^{n}$
on each side of the previous equality
\[
\mathbb{E}\left(f\left(R_{n}\right)\right)=\intop_{\mathbb{R}}\boldsymbol{1}_{\left[0,1\right[}\left(y\right)f\left(y\right)\text{d}\lambda\left(y\right).
\]
This shows that $R_{n}$ follows the uniform law on $\left[0,1\right[.$
Moreover, for every $f\in\mathscr{C}^{+}_{\mathscr{K}}\left(\mathbb{R}\right),$
for every subset $J$ of $\left\llbracket 1,n\right\rrbracket $ and
for every $\left(\epsilon_{j}\right)_{j\in J}\in\left\{ 0,1\right\} ^{J},$
it follows, by summing in each side of the equation $\refpar{eq:expect_f_r_n}$
over all values $\epsilon_{j}$ for $j\in J^{c}$
\begin{align}
\mathbb{E}\left(f\left(R_{n}\right)\prod_{j\in J}\boldsymbol{1}_{\left(D_{j}=\epsilon_{j}\right)}\right) & =\left[\prod_{j\in J}P\left(D_{j}=\epsilon_{j}\right)\right]\mathbb{E}\left(f\left(R_{n}\right)\right)\label{eq:expect_f_r_n-1}
\end{align}
which proves that $R_{n}$ and $\left(D_{1},D_{2},\cdots,D_{n}\right)$
are independent. 
\end{itemize}
\end{proof}

\begin{remark}{}{}

The sequence of random variables $\left(R_{n}\right)$ does not form
a family of independent random variables. Indeed, since:
\[
D_{n}=-R_{n}+2R_{n-1}
\]
and since both $R_{n-1}$ and $R_{n}$ admit a density, if they were
independent, then $D_{n}$ would also admit a density---see the next
section, Proposition $\ref{pr:sum_two_rv_density}$---which is false,
since $D_{n}$ takes only the values $0$ and $1$ and therefore does
not have a density.

\end{remark}

\begin{figure}[t]
\begin{center}\includegraphics[width=0.4\textwidth]{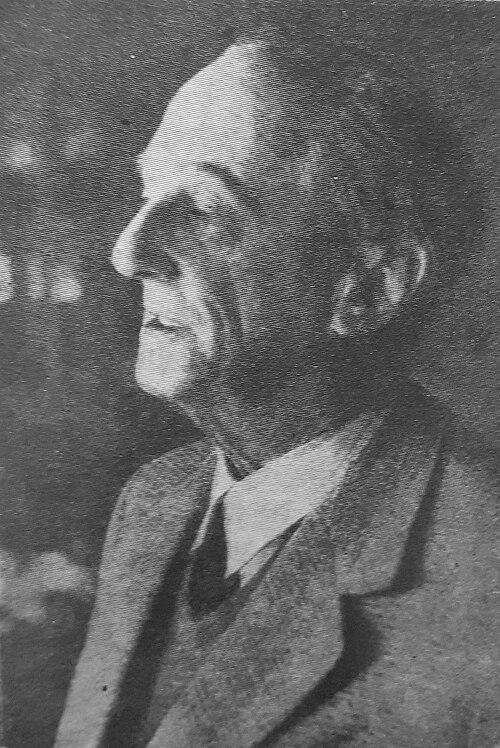}

{\tiny\href{https://en.wikipedia.org/wiki/Nobel_Foundation}{Nobel Foundation}
Public domain}\end{center}

\caption{\textbf{\protect\href{https://en.wikipedia.org/wiki/Hugo_Steinhaus}{Hugo Steinhaus}}
(1887 - 1972)}\sindex[fam]{Steinhaus, Hugo}
\end{figure}

As a corollary, we give a constructive proof of the existence of a
sequence of independent real-valued random variables with arbitrary
prescribed laws\footnote{The first rigourous mathematical presentation of the sequences of
independent random variables---and in particular of the game of heads
and tails---has been done by \textbf{\sindex[fam]{Steinhaus, Hugo}\href{https://en.wikipedia.org/wiki/Hugo_Steinhaus}{Hugo Steinhaus}}
(1887 - 1972) in 1923 and 1930. Hugo Steinhaus is a Polish mathematician
and educator. He proposed considering random variables as measurable
functions defined on $\left[0,1\right].$ Steinhaus's work precedes
by a few years the publication by \textbf{Kolmogorov} of its axiomatic
construction of probability theory, which is based on measure theory
and uses arbitrary probabilized spaces \preds, formulated in 1929
and 1933. He cpntributed in many branches of mathematics, such as
functional analysis, geometry, mathematical logic and trigonometry,
and is considered as one of the early founders of game theory and
probability theory.}.

\begin{corollary}{Existence of a Sequence of Independent Real-Valued Random Variables with Given Arbitrary Laws}{exist_seq_indep_rv_given_laws}

Let $\left(\mu_{j}\right)_{j\in\mathbb{N}^{\ast}}$ be a sequence
of probability on $\left(\mathbb{R},\mathscr{B}_{\mathbb{R}}\right).$

There exists a sequence of independent real-valued random variables
$\left(X_{j}\right)_{j\in\mathbb{N}^{\ast}}$ defined on the probabilized
space $\left(\left[0,1\right[,\mathscr{B}_{\left[0,1\right[},P\right),$
where $P$ is the probability, restriction of the Lebesgue measure
to $\left[0,1\right[,$ such that, for each $j\in\mathbb{N}^{\ast},$
$X_{j}$ has law $\mu_{j}.$

\end{corollary}

\begin{proof}{}{}

We begin by proving the existence of a sequence of independent random
variables that all follow the uniform law on the interval $\left[0,1\right[.$
The general case will then follow easily.

We retain the same notations as in Proposition $\ref{pr:infinite_sequence_bernoulli_rv}.$
The functions $D_{n},n\in\mathbb{N}^{\ast}$ are considered as random
variables defined on $\left[0,1\right[,$ equipped with its Borel
\salg and the Lebesgue measure.

Let $\left(N_{j}\right)_{j\in\mathbb{N}^{\ast}}$ be a sequence of
infinite subsets partitioning $\mathbb{N}^{\ast},$ that is
\[
\mathbb{N}^{\ast}=\biguplus_{j\in\mathbb{N}^{\ast}}N_{j}
\]
and let $\varphi_{j}$ be a sequence obtained by enumerating the elements
of $N_{j}$ in nondecreasing order. To construct such a partition,
we may begin with a bijection 
\[
\phi:\mathbb{N}^{\ast}\times\mathbb{N}^{\ast}\to\mathbb{N}^{\ast},
\]
 for instance, the one defined by
\[
\phi\left(j,k\right)=\dfrac{\left(j+k-2\right)\left(j+k-1\right)}{2}+k.
\]

We denote $\varphi_{j}$ the nondecreasing sequence $k\mapsto\phi\left(j,k\right)$
and we denote $N_{j}$ the image of $\mathbb{N}^{\ast}$ by this sequence.

For each $j\in\mathbb{N}^{\ast},$ we set
\[
Y_{j}=\sum^{+\infty}_{k=1}\dfrac{1}{2^{k}}D_{\varphi_{j}\left(k\right)}.
\]

In other words, we redistribute the dyadic digits of $x$ into infinitely
disjoint subsequences, each of which is then used to build a new real
number, denoted $Y_{j}\left(x\right).$

The random variables $Y_{j},j\in\mathbb{N},$ are independent. Indeed,
each $Y_{j}$ is mesurable with respect to the \salg $\sigma\left(D_{n},n\in N_{j}\right)$
and these \salg are independent, since the sets $N_{j}$ form a partition
of $\mathbb{N}^{\ast},$ and the $D_{n},n\in\mathbb{N}^{\ast}$ are
independent---by Theorem $\ref{th:indep_salg_gen_union_pis}.$

Moreover, for each $j\in\mathbb{N}^{\ast},$ the law of $Y_{j}$ is
uniform on $\left[0,1\right[.$ Indeed, for every $n\in\mathbb{N}^{\ast},$
define the truncated sum
\[
Y_{j,n}=\sum^{n}_{k=1}\dfrac{1}{2^{k}}D_{\varphi_{j}\left(k\right)}.
\]

The law of a sum of independent random variables does not depend on
the laws of these random variables---see further, Proposition $\ref{pr:sum_two_rv_density}$.
Since the random variables $D_{1},D_{2},\cdots,D_{k}$ are independent
and with same law, the law of $Y_{j,n}$ is thus the same as that
of
\[
Z_{n}=\sum^{n}_{k=1}\dfrac{1}{2^{k}}D_{k}.
\]

Note that
\[
Z=\lim_{n\to+\infty}\nearrow Z_{n}
\]
is simply the identity function on $\left[0,1\right[.$

Since $Y_{j}=\lim_{n\to+\infty}\nearrow Y_{j,n},$
\[
\left(Y_{j}\leqslant y\right)=\lim_{n\to+\infty}\searrow\left(Y_{j,n}\leqslant y\right).
\]

Similarly,
\[
\left(Z\leqslant y\right)=\lim_{n\to+\infty}\searrow\left(Z_{n}\leqslant y\right).
\]

Hence, we obtain for $y\in\left[0,1\right],$
\[
P\left(Y_{j}\leqslant y\right)=\lim_{n\to+\infty}P\left(Y_{j,n}\leqslant y\right)=\lim_{n\to+\infty}\left(Z_{n}\leqslant y\right)=P\left(Z\leqslant y\right)=y.
\]

Finally, if $F_{j}$ is the cumulative distribution function associated
with the probability $\mu_{j},$ defined by
\[
\forall x\in\mathbb{R},\,\,\,\,F_{j}\left(x\right)=\mu_{j}\left(\left]-\infty,x\right]\right),
\]
and if $G_{j}$ denotes its ``pseudo-inverse'' defined by
\[
\forall t\in\mathbb{R},\,\,\,\,G_{j}\left(t\right)=\inf\left\{ x\in\mathbb{R}:\,F_{j}\left(x\right)\geqslant t\right\} ,
\]
and if $X_{j}=G_{j}\left(Y_{j}\right),$ then by Exercise $\ref{exo:exercise9.1}$
in Chapter \ref{chap:PartIIChap9}---a fundamental result for the
simulation of probability laws---we obtain that the law of $X_{j}$
is $\mu_{j},$ which completes the proof.

\end{proof}

\subsection{Complement. Probabilities Product on $\left\{ 0,1\right\} ^{\mathbb{N}^{\ast}}$}

\subsubsection*{Modelling the Heads and Tails Tossing Game with A Space of Sequences}

In practice, it is possible, from the constructed model on the interval
$\left[0,1\right[$ to build a model where the fundamental space is
the space of sequences $\left\{ 0,1\right\} ^{\mathbb{N}^{\ast}},$
which is the ``natural `` model alluded to at the beginning of this
section.

Denote $D$ the function from $\left[0,1\right]$ to $\left\{ 0,1\right\} ^{\mathbb{N}^{\ast}}$
defined by
\[
\forall x\in\left[0,1\right[,\,\,\,\,D\left(x\right)=\left(D_{n}\left(x\right)\right)_{n\in\mathbb{N}^{\ast}}.
\]
It follows from the relation $\refpar{eq:dyadic_dev_x}$ that $D$
is injective. Based on our analysis of the dyadic expansion, the image
of $D$ is $\left\{ 0,1\right\} ^{\mathbb{N}^{\ast}},$ minus the
sub-set $\Omega_{1},$ consisting of the sequences that are equal
to 1 from a certain rank.

We equip $\left\{ 0,1\right\} ^{\mathbb{N}^{\ast}}$ with the product
\salg $\mathcal{A},$ product of the \salg on the subsets of each
of its components. More precisely, $\mathcal{A}$ is the \salg generated
by the family of subsets of the form
\[
\prod_{n\in\mathbb{N}^{\ast}}A_{n}\equiv\left\{ \omega\in\left\{ 0,1\right\} ^{\mathbb{N}^{\ast}}:\,\forall n\in\mathbb{N}^{\ast},\,\,\,\,\omega_{n}\in A_{n}\right\} ,
\]
where $A_{n}$ is a subset of $\left\{ 0,1\right\} ^{\mathbb{N}^{\ast}}$
equal to $\left\{ 0,1\right\} $ except a finite number of indices
$n.$

The function $D$ from $\left(\left[0,1\right[,\mathscr{B}_{\left[0,1\right[}\right)$
to $\left(\left\{ 0,1\right\} ^{\mathbb{N}^{\ast}},\mathcal{A}\right)$
is measurable. 

Indeed, it suffices to note that for every $I\in\mathscr{P}_{f}\left(\mathbb{N}^{\ast}\right),$
for every $n\in I,$ and for every $A_{n}=\left\{ \epsilon_{n}\right\} ,$
where $\epsilon_{n}\in\left\{ 0,1\right\} ,$
\[
D^{-1}\left(\prod_{n\in\mathbb{N}^{\ast}}A_{n}\right)=\bigcap_{n\in I}\left(D_{n}=\epsilon_{n}\right)\in\mathscr{B}_{\left[0,1\right[}.
\]

Let $Q$ be the probability on $\left(\Omega,\mathcal{A}\cap\Omega\right),$
image of $P$ under the measurable function $D.$ 

By taking for $Z_{j}$ the projection from $\left\{ 0,1\right\} ^{\mathbb{N}^{\ast}}$
on the $j-$th coordinate, we obtain a sequence of random variables
defined on $\left(\left\{ 0,1\right\} ^{\mathbb{N}^{\ast}},\mathcal{A},Q\right),$
independent and with same Bernoulli law $\mathcal{B}\left(1,\dfrac{1}{2}\right).$ 

Indeed, we have $D_{j}=Z_{j}\circ D$ for $j\in\mathbb{N}^{\ast}.$
It follows, by the definition itself of the probability $Q,$ that
\[
Q\left(Z_{j}=\epsilon\right)=P\left(D_{j}=\epsilon\right),
\]
and more generally that
\[
Q\left(Z_{j_{1}}=\epsilon_{j_{1}},\cdots,Z_{j_{n}}=\epsilon_{j_{n}}\right)=P\left(D_{j_{1}}=\epsilon_{j_{1}},\cdots,D_{j_{n}}=\epsilon_{j_{n}}\right).
\]
Hence, the random variables $Z_{j}$ have the same law as the $D_{j},$
and since the $D_{j}$ are independent, the $Z_{j}$ are independent.

We can observe that $Q$ is precisely the product measure of the uniform
probabilities on the factors $\left\{ 0,1\right\} $---see Corollary
$\ref{co:infinite_product_proba}$ for the formal definition. The
properties required to establish this identification translates exactly
the fact that the random variables $Z_{j},j\in\mathbb{N}^{\ast}$
are independent and with same law $\mathcal{B}\left(1,\dfrac{1}{2}\right).$

\begin{example}{Usage Example}{}

In the context of an infinite heads and tails tossing game, the probability
that the finite sequence $\left(\epsilon_{1},\dots,\epsilon_{n}\right)$
appears infintely often is equal to 1. 

Indeed, the considered event is denoted $\limsup_{n\to+\infty}\left(A_{n}\right),$
where
\[
A_{j}=\left\{ \left(Z_{j+1},Z_{j+2},\cdots,Z_{j+n}\right)=\left(\epsilon_{1},\cdots,\epsilon_{n}\right)\right\} .
\]

Now define
\[
B_{j}=\left\{ \left(Z_{jn+1},Z_{jn+2},\cdots,Z_{\left(j+1\right)n}\right)=\left(\epsilon_{1},\cdots,\epsilon_{n}\right)\right\} ,
\]

The events $B_{j}$ are $Q-$independent, because they depend on disjoint
blocks of coordinates, and
\[
\limsup_{j\to+\infty}B_{j}\subset\limsup_{j\to+\infty}A_{j}.
\]

Moreover, since $Q\left(B_{j}\right)=\dfrac{1}{2^{n}},$
\[
\sum^{+\infty}_{j=1}Q\left(B_{j}\right)=+\infty
\]
and by the Borel-Cantelli lemma, $Q\left(\limsup_{j\to+\infty}B_{j}\right)=1.$

A fortiori,
\[
Q\left(\limsup_{j\to+\infty}A_{j}\right)=1.
\]

\end{example}

\begin{remark}{}{}

To answer the question, it is enough to construct a model allowing
to handle an infinite sequence of independent random variables, following
the same Bernoulli law $\mathcal{B}\left(1,\dfrac{1}{2}\right).$

\end{remark}

\subsubsection*{Canonical Model for a Sequence of Independent Random Variables of
Bernoulli}

If now we consider a sequence of random variables of given laws $\mu_{j},$
we can transport to $\left\{ 0,1\right\} ^{\mathbb{N}^{\ast}}$ the
solution to the problem ``construct a sequence of independent random
variables with law $\mu_{j}$'' given in Corollary $\ref{co:exist_seq_indep_rv_given_laws}.$
Since the random variables $X_{j}$ constructed there are defined
on the space $\left[0,1\right[,$ and, after having observed that
the function $D^{-1}$ defined on $\left\{ 0,1\right\} ^{\mathbb{N}^{\ast}}\backslash\Omega_{1}$
is measurable, it suffices to consider the variables $X_{j}\circ D^{-1}$---we
arbitrarily extend them to the set $\Omega_{1},$ which is of zero
$Q-$probability.

This yields a model where the fundamental space is a sequence space,
but not fully adapted to the problem. The random variables $X_{i}$
are not like in the previous model of heads and tails in relation
with the coordinate projections---also called the coordinate functions---
on the factor spaces.

We now describe a better-adapted model, assuming for simplicity that
the laws $\mu_{j}$ are all supported on $\left\{ 0,1\right\} $---i.e.,
they are Bernoulli laws with parameter $p_{j}$. This model is suited,
for instance, to the description of a heads and tails game with an
unfair coin, in which the $p_{j}$ are all equal to some fixed parameter
$p\in\left[0,1\right].$

Let $X:\left[0,1\right[\to\left\{ 0,1\right\} ^{\mathbb{N}^{\ast}}$
be the function which maps each $x$ to the sequence $\left(X_{j}\left(x\right)\right)_{j\in\mathbb{N}^{\ast}}.$
Note that this function depends on the choice of the sequence $\mu=\left(\mu_{j}\right)_{j\in\mathbb{N}^{\ast}}.$
The function $X$ is measurable: by the definition of the \salg $\mathcal{A}$
on $\left\{ 0,1\right\} ^{\mathbb{N}^{\ast}},$ it suffices to verify
that, for every $n,$ the function $x\mapsto\left(X_{1}\left(x\right),\cdots,X_{n}\left(x\right)\right)$
is measurable, which it is.

Let $P_{\mu}$ be the image of $X$ by the probability $P.$ Then,
if $\left\{ 0,1\right\} ^{\mathbb{N}^{\ast}}$ is equipped with the
\salg $\mathcal{A}$ and the probability $P_{\mu},$ the projections
$Z_{j}=\text{pr}_{j}$ are independent and of law $\mu_{j},j\in\mathbb{N}^{\ast}.$
As previously, by definition of the image probability, we have 
\[
P_{\mu}\left(Z_{j_{1}}=\epsilon_{j_{1}},\cdots,Z_{j_{n}}=\epsilon_{j_{n}}\right)=P\left(X_{j_{1}}=\epsilon_{j_{1}},\cdots,X_{j_{n}}=\epsilon_{j_{n}}\right),
\]
from which it follows that $Z_{j}$ has the same law than $X_{j}$
and that the $Z_{j}$ are independent.

We have shown in passing that $P_{\mu}$ is the product probability
of the probabilities $\mu_{j},j\in\mathbb{N}^{\ast}.$

In particular, the law $P_{\mu}$ depends only on the $\mu_{j},$
and not on the $X_{j}$ which may involve arbitrary choices in their
construction. This is why the model we have just built---namely the
fundamental space $\left(\left\{ 0,1\right\} ^{\mathbb{N}^{\ast}},\mathcal{A},P_{\mu}\right)$
and the random variables $Z_{j},j\in\mathbb{N}^{\ast}$---can be
called the \textbf{canonical model} for the realization of a sequence
of independent Bernoulli random variables with prescribed laws. We
also observe that this provides a proof, in a particular case, of
the Kolmogorov theorem.

Finally, let us emphasize once again that we have obtained two mathematically
equivalent but structurally different solutions to the same modelling
problem: the one described here, and the one presented in Corollary
$\ref{co:exist_seq_indep_rv_given_laws}.$

\begin{remark}{}{}

Consider the case where all the $\mu_{j}$ are equal to $\mathcal{B}\left(1,p\right),$
with $0<p<1.$ We then write $P_{p}$ instead of $P_{\mu}$---for
instance, $P_{1/2}=Q.$

The measure $P_{p}$ is diffuse\footnotemark, and mutually singular\index{mutually singular}
(or foreign\index{foreign})\footnotemark to $Q$ if $p\neq\dfrac{1}{2}.$

Indeed, for $\omega\in\left\{ 0,1\right\} $ and for every $n\in\mathbb{N}^{\ast},$
\begin{align*}
P_{p}\left(\left\{ \omega\right\} \right) & =P_{p}\left(Z_{1}=\omega_{1},\cdots,Z_{n}=\omega_{n},\cdots\right)\\
 & \leqslant P_{p}\left(Z_{1}=\omega_{1},\cdots,Z_{n}=\omega_{n}\right)\leqslant\rho^{n}
\end{align*}
where $\rho=\sup\left(p,1-p\right).$ Hence, 
\[
P_{p}\left(\left\{ \omega\right\} \right)\leqslant\lim_{n\to+\infty}\rho^{n}=0.
\]

To demonstrate that $P_{p}$ and $Q$ are mutually singular, the simplest
approach is to exhibit an event that has at the same time probability
1 for $P_{p}$ and probability 0 for $Q.$ Such is the case for the
event
\[
\lim_{n\to+\infty}\dfrac{1}{n}\left(Z_{1}+\cdots+Z_{n}\right)=p
\]
by the strong law of large numbers---Theorem $\ref{th:strong_law_large_nb}.$

As a consequence, we obtain a construction of a measure on $\left[0,1\right]$
that is diffuse and mutually singular to the Lebesgue measure. Both
the probabilities $P_{p}$ and $Q$ assign a measure 0 to the countable
subset $\Omega_{1}.$ Therefore, we may restrict these probabilities
to $\left\{ 0,1\right\} ^{\mathbb{N}^{\ast}}\backslash\Omega_{1}.$
Additionally, $D$ defines a bijection from $\left[0,1\right[$ onto
this last set and the inverse of this bijection, denoted $D^{-1},$
is measurable. 

Indeed, the \salg $\mathscr{B}_{\left[0,1\right[}$ is generated
by dyadic intervals of the form $\left[\dfrac{k}{2^{n}},\dfrac{k+1}{2^{n}}\right[.$
The inverse image of such an interval by $D^{-1}$ is the set of the
$\omega\in\Omega$ which $n$ first coordinates are equal to the $n$
first digits of the diadyc expansion of $\dfrac{k}{2^{n}}.$ 

By now considering the images by $D^{-1}$ of the probabilities $Q$
and $P_{p},$ we obtain from the one hand the Lebesgue measure on
$\left[0,1\right[$, and on the other hand the diffuse probability
on $\left[0,1\right[,$ necessarily mutually singular with respect
to the Lebesgue measure.

\end{remark}

\addtocounter{footnote}{-1}

\footnotetext{Tr.N. In measure theory, given a measurable space $\left(X,\Sigma\right)$
and a measure $\mu$ on that space, a set $A\subset X$ is called
an \textbf{atom\index{atom}} if $\mu\left(A\right)>0$ and for any
measurable subset $B$ of $A,$ either $\mu\left(B\right)=0$ or $\mu\left(B\right)=\mu\left(A\right).$
A measure with no atoms is called a \textbf{non-atomic measure\index{non-atomic measure}\mindex{measure ! non-atomic}}
or a \textbf{diffuse measure\mindex{measure ! diffuse}}\index{diffuse measure}.}

\stepcounter{footnote}

\footnotetext{The notion of \textbf{foreign measures\index{foreign measures}\mindex{measures ! foreign}
}(\textbf{or mutually singular\mindex{measures ! mutally singular}})
is defined in Definition $\ref{df:abs_cont_foreign_meas}.$}

\section{Convolution and Law of the Sum of Independent Random Variables}\label{sec:Convolution-and-Law}

\begin{definition}{Convolution Product of Two Measures}{}

Let $\mu_{1}$ and $\mu_{2}$ be two bounded measures---or, respectively,
two probabilities---on $\left(\mathbb{R}^{d},\mathscr{B}_{\mathbb{R}^{d}}\right).$
Let $S$ be the sum function on $\mathbb{R}^{d}.$ The measure image
by $S$ of $\mu_{1}\otimes\mu_{2}$ is called the \textbf{convolution
product\index{convolution product}} of $\mu_{1}$ and $\mu_{2}$
and denoted $\mu_{1}\ast\mu_{2}.$ This is a bounded measure---respectively
a probability---on $\left(\mathbb{R}^{d},\mathscr{B}_{\mathbb{R}^{d}}\right).$

\end{definition}

\begin{proposition}{}{}

For any measurable nonnegative function $f$ on $\mathbb{R}^{d},$\boxeq{
\begin{equation}
\intop_{\mathbb{R}^{d}}f\text{d}\left(\mu_{1}\ast\mu_{2}\right)=\intop_{\mathbb{R}^{d}\times\mathbb{R}^{d}}f\left(x_{1}+x_{2}\right)\text{d}\left(\mu_{1}\otimes\mu_{2}\right)\left(x_{1},x_{2}\right).\label{eq:int_convolution_product_two_measures}
\end{equation}
}Moreover, $f\in\mathscr{L}^{1}\left(\mathbb{R}^{d},\mathscr{B}_{\mathbb{R}^{d}},\mu_{1}\ast\mu_{2}\right)$
if and only if the function $\left(x_{1},x_{2}\right)\mapsto f\left(x_{1}+x_{2}\right)$
is $\left(\mu_{1}\otimes\mu_{2}\right)-$integrable. In this case,
the equality above still holds.

\end{proposition}

\begin{proof}{}{}

The proof is classical. 

If $f=\boldsymbol{1}_{A},$ where $A\in\mathscr{B}_{\mathbb{R}^{d}},$
this follows directly from the definition of $\mu_{1}\ast\mu_{2}.$
The equality $\refpar{eq:int_convolution_product_two_measures}$ then
holds by linearity for any nonnegative step function. By the Beppo
Levi property, it follows for any nonnegative measurable function
by taking a nondecreasing sequence of nonnegative step functions converging
to $f.$

The remaining of the proof is also standard: for integrability, we
take the absolute values and apply $\refpar{eq:int_convolution_product_two_measures};$
finally, we decompose $f$ into its nonnegative and nonpositive parts.

\end{proof}

\begin{proposition}{Law of the Sum of Two Random Variables}{sum_two_rv_density}

Let $X_{1}$ and $X_{2}$ be two independent random variables taking
values in $\mathbb{R}^{d}$ on the probabilized space $\left(\Omega,\mathscr{A},P\right).$

The law of $X_{1}+X_{2}$ is the convolution product of the laws $X_{1}$
and $X_{2},$\boxeq{
\[
P_{X_{1}+X_{2}}=P_{X_{1}}\ast P_{X_{2}}.
\]
}

\end{proposition}

\begin{proof}{Density of the Sum of Two Random Variables with Densities}{}

Since $X_{1}$ and $X_{2}$ are independent, $P_{X_{1},X_{2}}=P_{X_{1}}\otimes P_{X_{2}}.$

It follows that
\[
P_{X_{1}+X_{2}}=P_{S\circ\left(X_{1},X_{2}\right)}=S\left(P_{X_{1}}\otimes P_{X_{2}}\right)=P_{X_{1}}\ast P_{X_{2}}.
\]

\end{proof}

\begin{proposition}{Convolution Product of Two Densities}{}If, additionally
to the assumptions of Proposition $\ref{pr:sum_two_rv_density},$
we suppose that $X_{1}$ and $X_{2}$ admit the respective densities
$f_{X_{1}}$ and $f_{X_{2}},$ then $X_{1}+X_{2}$ admits a density
$f_{X_{1}+X_{2}}$ defined for $y\in\mathbb{R}^{d}$ by\boxeq{
\begin{align*}
f_{X_{1}+X_{2}}\left(y\right) & =\intop_{\mathbb{R}^{d}}f_{X_{1}}\left(x_{1}\right)f_{X_{2}}\left(y-x_{1}\right)\text{d}\lambda_{d}\left(x_{1}\right)\\
 & =\intop_{\mathbb{R}^{d}}f_{X_{1}}\left(y-x_{2}\right)f_{X_{2}}\left(x_{2}\right)\text{d}\lambda_{d}\left(x_{2}\right).
\end{align*}
}

We say that $f_{X_{1}+X_{2}}$ is the \textbf{\index{convolution product}convolution
product} of the functions $f_{X_{1}}$ and $f_{X_{2}}.$

\end{proposition}

\begin{proof}{}{}

For every $A\in\mathscr{B}_{\mathbb{R}^{d}},$ by Proposition $\ref{pr:sum_two_rv_density}$
and the fact that $X_{1}$ and $X_{2}$ admit densities and are independent,
\[
P_{X_{1}+X_{2}}\left(A\right)=\intop_{\mathbb{R}^{d}}\boldsymbol{1}_{A}\left(x_{1}+x_{2}\right)f_{X_{1}}\left(x_{1}\right)f_{X_{2}}\left(x_{2}\right)\text{d}\lambda_{d}\otimes\lambda_{d}\left(x_{1},x_{2}\right).
\]
We now make the change of variables
\[
\left\{ \begin{array}{l}
y_{1}=x_{1}+x_{2}\\
y_{2}=x_{2}.
\end{array}\right.
\]

The associated diffeomorphism has a Jacobian determinant equal to
1, it follows that
\[
P_{X_{1}+X_{2}}\left(A\right)=\intop_{\mathbb{R}^{d}}\boldsymbol{1}_{A}\left(y_{1}\right)f_{X_{1}}\left(y_{1}-y_{2}\right)f_{X_{2}}\left(y_{2}\right)\text{d}\lambda_{d}\otimes\lambda_{d}\left(y_{1},y_{2}\right).
\]
Then, by the Fubini theorem,
\[
P_{X_{1}+X_{2}}\left(A\right)=\intop_{A}\left[\intop_{\mathbb{R}^{d}}f_{X_{1}}\left(y_{1}-y_{2}\right)f_{X_{2}}\left(y_{2}\right)\text{d}\lambda_{d}\left(y_{2}\right)\right]\text{d}\lambda_{d}\left(y_{1}\right),
\]
which proves the result.

\end{proof}

\begin{example}{}{}

Let $X_{1}$ and $X_{2}$ be two independent real-valued random variables
defined on a probabilized space \preds, with respective laws $\gamma\left(a_{1},p\right)$
and $\gamma\left(a_{2},p\right),$ where $a_{1},a_{2}$ and $p$ are
positive real numbers. The law of $X_{1}+X_{2}$ is the law $\gamma\left(a_{1}+a_{2},p\right).$

\end{example}

\begin{solutionexample}{}{}

Indeed, for $i=1,2,$
\[
\forall x\in\mathbb{R},\,\,\,\,f_{X_{i}}\left(x\right)=\boldsymbol{1}_{\mathbb{R}^{+}}\left(x\right)\dfrac{p^{a_{i}}}{\Gamma\left(a_{i}\right)}\text{e}^{-px}x^{a_{i}-1},
\]
and $X_{1}+X_{2}$ admits as density the convolution product of the
functions $f_{X_{1}}$ and $f_{X_{2}}.$

Since
\begin{multline*}
f_{X_{1}}\left(y-x_{2}\right)f_{X_{2}}\left(x_{2}\right)=\boldsymbol{1}_{\mathbb{R}^{+}}\left(y-x_{2}\right)\boldsymbol{1}_{\mathbb{R}^{+}}\left(x_{2}\right)\dfrac{p^{a_{1}}p^{a_{2}}}{\Gamma\left(a_{1}\right)\Gamma\left(a_{2}\right)}\text{e}^{-py}\left(y-x_{2}\right)^{a_{1}-1}x^{a_{2}-1}_{2},
\end{multline*}
and, since\footnotemark
\[
\boldsymbol{1}_{\mathbb{R}^{+}}\left(y-x_{2}\right)\boldsymbol{1}_{\mathbb{R}^{+}}\left(x_{2}\right)=\boldsymbol{1}_{\mathbb{R}^{+}}\left(y\right)\boldsymbol{1}_{\left[0,y\right]}\left(x_{2}\right),
\]
it follows
\[
f_{X_{1}+X_{2}}\left(y\right)=\boldsymbol{1}_{\mathbb{R}^{+}}\left(y\right)\dfrac{p^{a_{1}+a_{2}}}{\Gamma\left(a_{1}\right)\Gamma\left(a_{2}\right)}\text{e}^{-py}\intop^{y}_{0}\left(y-x_{2}\right)^{a_{1}-1}x^{a_{2}-1}_{2}\text{d}x_{2}.
\]
Thus, after making the change of variable $u=1+x_{2}-y,$
\[
f_{X_{1}+X_{2}}\left(y\right)=\boldsymbol{1}_{\mathbb{R}^{+}}\left(y\right)\dfrac{p^{a_{1}+a_{2}}}{\Gamma\left(a_{1}\right)\Gamma\left(a_{2}\right)}\text{e}^{-py}y^{a_{1}+a_{2}-1}\intop^{1}_{0}\left(1-u\right)^{a_{1}-1}u^{a_{2}-1}\text{d}u.
\]
But, since the function 
\[
y\mapsto\boldsymbol{1}_{\mathbb{R}^{+}}\left(y\right)\dfrac{p^{a_{1}+a_{2}}}{\Gamma\left(a_{1}+a_{2}\right)}\text{e}^{-py}y^{a_{1}+a_{2}-1}
\]
is a density of probability, the function $f_{X_{1}+X_{2}}$ is equal
to this density---the normalizing coefficient is unique---and $X_{1}+X_{2}$
follows the law $\gamma\left(a_{1}+a_{2},p\right).$ 

Moreover, we obtain the relation
\[
\intop^{1}_{0}\left(1-u\right)^{a_{1}-1}u^{a_{2}-1}\text{d}u=\dfrac{\Gamma\left(a_{1}\right)\Gamma\left(a_{2}\right)}{\Gamma\left(a_{1}+a_{2}\right)}.
\]

\end{solutionexample}

\footnotetext{The meticulous usage of the indicator functions avoids
many errors in the determination of the support of the laws.}

\section*{Exercises}

\addcontentsline{toc}{section}{Exercises}

Unless otherwise specified, all random variables are defined on the
same probabilized space \preds.

\begin{exercise}{Independent Gaussian Random Variables. Independence of the Radius and the Polar Angle}{exercise10.1}

Let $X$ and $Y$ be two independent real-valued random variables
following the same Gauss law $\mathcal{N}\left(0,1\right).$ We denote
$\phi$ the function from $\mathbb{R}^{+\ast}\times\left]0,2\pi\right[$
to $\mathbb{R}^{2}\backslash\mathbb{R}^{+}\times\left\{ 0\right\} $---polar
coordinate transformation---defined by
\[
\forall\left(\rho,\theta\right)\in\mathbb{R}^{+\ast}\times\left]0,2\pi\right[,\,\,\,\,\phi\left(\rho,\theta\right)=\left(\rho\cos\theta,\rho\sin\theta\right),
\]
 and $\psi$ its inverse.

Now, define the function $\widehat{\psi}$ from $\mathbb{R}^{2}$
to $\mathbb{R}\times\left[0,2\pi\right[$ by
\[
\widehat{\psi}\left(x,y\right)=\begin{cases}
\psi\left(x,y\right), & \text{if }\left(x,y\right)\in\mathbb{R}^{2}\backslash\mathbb{R}^{+}\times\left\{ 0\right\} ,\\
\left(0,0\right), & \text{otherwise.}
\end{cases}
\]

Define the random variable $\left(R,\Theta\right)$ by
\[
\left(R,\Theta\right)=\widehat{\psi}\left(X,Y\right).
\]

1. Determine the law of the random variable $\left(R,\Theta\right)$
expressed through its density. 

2. Justify, on the one hand, the independence of $R$ and $\Theta,$
and on the other hand, the independence of $S\equiv R^{2}$ and $\Theta.$ 

3. Specify the law of $S.$

\end{exercise}

\begin{exercise}{On the Simulation Path  of Two Independent Gaussian Random Variables}{exercise10.2}

Let $S$ and $\Theta$ be two independent real-valued random variables
of respective law the exponential law $\text{e}^{\frac{1}{2}}$ and
the uniform law on $\left[0,2\pi\right].$ 

Define two real-valued random variables $X$ and $Y$ by
\[
X=\sqrt{S}\cos\Theta\,\,\,\,\,\,\,\,Y=\sqrt{S}\sin\Theta.
\]

What is the law of the random variable $\left(X,Y\right)?$

\end{exercise}

\begin{exercise}{Law and Moments of the Maximum of Exponential Law Random Variables}{exercise10.3}

Let $\left(X_{n}\right)_{n\in\mathbb{N}^{\ast}}$ be a sequence of
independent random variables, following the same exponential law $\text{e}^{1}.$

Define for each $n\in\mathbb{N}^{\ast},$ the random variables $Y_{n}$
and $Z_{n}$ by
\[
Y_{n}=\max_{1\leqslant i\leqslant n}X_{i}\,\,\,\,\,\,\,\,Z_{n}=\sum^{n}_{i=1}\dfrac{X_{i}}{i}.
\]

Prove by induction that $Y_{n}$ and $Z_{n}$ follow the same law.
Deduce the expectation and variance of $Y_{n}.$

\end{exercise}

\begin{exercise}{Chi-Squared Law and Student Law}{exercise10.4}

Let $\left(X_{n}\right)_{n\in\mathbb{N}^{\ast}}$ be a sequence of
independent random variables following the same Gauss law $\mathcal{N}_{\mathbb{R}}\left(0,1\right).$

Define for every $n\in\mathbb{N}^{\ast},$ the random variables $Y_{n},$
$T^{\prime}_{n}$ and $T_{n}$ by
\[
Y_{n}=\sum^{n}_{j=1}X^{2}_{j}\,\,\,\,\,\,\,\,T^{\prime}_{n}=\dfrac{X_{n+1}}{\sqrt{Y_{n}}}\,\,\,\,\,\,\,\,T_{n}=\sqrt{n}T^{\prime}_{n}.
\]

1. a. Prove by induction that $Y_{n}$ admits a density given by
\begin{equation}
\forall t\in\mathbb{R},\,\,\,\,f_{Y_{n}}\left(y\right)=\boldsymbol{1}_{\mathbb{R}^{+}}\left(y\right)\dfrac{1}{2^{\frac{n}{2}}\Gamma\left(\dfrac{n}{2}\right)}\text{e}^{-\frac{y}{2}}y^{\frac{n}{2}-1}.\label{eq:density_chi_two}
\end{equation}

b. Compute the expectation and variance of $Y_{n}.$

The law of $Y_{n}$ is the \textbf{Chi-squared law with $n$ degrees
of freedom}\index{Chi-squared law}\mindex{law ! Chi-Squared}.

2. a. Prove that $T^{\prime}_{n}$ admits a density, and compute it.

b. Deduce that $T_{n}$ admits a density $f_{T_{n}}$ given by
\begin{equation}
\forall t\in\mathbb{R},\,\,\,\,f_{T_{n}}\left(t\right)=\dfrac{1}{\sqrt{n\pi}}\dfrac{\Gamma\left(\dfrac{n+1}{2}\right)}{\Gamma\left(\dfrac{n}{2}\right)}\dfrac{1}{\left(1+\dfrac{t^{2}}{n}\right)^{\frac{n+1}{2}}}.\label{eq:Student_density}
\end{equation}

The law of $T_{n}$ is the \textbf{\index{Student law}\mindex{law! Student}law
of Student with parameter $n.$}

c. Without using the density of $T_{n},$ determine for which values
of $n$ the expectation and variance of $T_{n}$ exist, and compute
when they do exist. 

\end{exercise}

\begin{exercise}{Law of the Sum of Two Independent Random Variables---One with Density, One Discrete}{exercise10.5}

Let $X$ and $Y$ be two independent random variables of respective
laws the geometric law $\mathscr{G}_{\mathbb{N}}\left(p\right),$
where $0<p<1$ and the exponential law $\text{e}^{1}.$ 

Study the law of the random variable $Z=X+Y.$

\end{exercise}

\begin{exercise}{The Independence of Random Variables is not Always Intuitive! Gamma and Beta Laws}{exercise10.6}

Let $X$ and $Y$ be two independent random variables with respective
laws $\gamma\left(a,p\right)$ and $\gamma\left(b,p\right)$ where
$a,b,p$ are positive real numbers. 

Define the random variables $U,T,Z$ by
\[
U=X+Y,\,\,\,\,\,\,\,\,T=\dfrac{X}{X+Y},\,\,\,\,\,\,\,\,\text{and}\,\,\,\,\,\,\,\,Z=\dfrac{X}{Y}.
\]

Study the law of the random variable $\left(U,T\right)$ and deduce
from it the independence of the random variables $U$ and $T.$ Be
more precise on the laws of the random variables $U,T$ and $Z.$

\end{exercise}

\begin{exercise}{Characterization of the Exponential Laws}{exercise10.7}

Let $X$ and $Y$ be two independent real-valued random variables
with the same law $\mu.$ Suppose that $\mu$ admits a density $f$---that
is $\mu=f\cdot\lambda,$ where $\lambda$ is the measure of Lebesgue
on $\mathbb{R}$---such that $f\left(x\right)>0$ if and only if
$x\in\mathbb{R}^{+}.$ 

Define the random variables $U$ and $W$ by
\[
U=\min\left(X,Y\right)\,\,\,\,\,\,\,\,W=\left|X-Y\right|.
\]

1. Prove that the random variable $\left(U,W\right)$ admits a density,
and express it in terms of $f.$

\textit{Hint: To do this, compute, for any nonnegative measurable
function $g$ defined on $\mathbb{R}^{2},$ possibly bounded, the
integral 
\[
\intop_{\Omega}g\left(U,W\right)\text{d}P.
\]
}

Deduce that the random variables $U$ and $W$ each admit densities,
denoted respectively $f_{U}$ and $f_{W}.$

2. To simplify the analysis, assume further that $f$ is bounded and
that its restriction to $\mathbb{R}^{+}$ is continuous. 

Prove that $U$ and $W$ admit densities denoted $f_{U}$ and $f_{W},$
and that these densities have restrictions to $\mathbb{R}^{+}$ that
are continuous.

Deduce from the previous results that \textbf{$\mu$ is an exponential
law if and only if $U$ and $W$ are independent.}

3. Finally, prove that even if we do not assume anymore that $f$
is bounded and that its restriction to $\mathbb{R}^{+}$ is continuous,
\textbf{$\mu$ is an exponential law if and only if $U$ and $W$
are independent.}

\end{exercise}

\begin{figure}[t]
\begin{center}\includegraphics[width=0.4\textwidth]{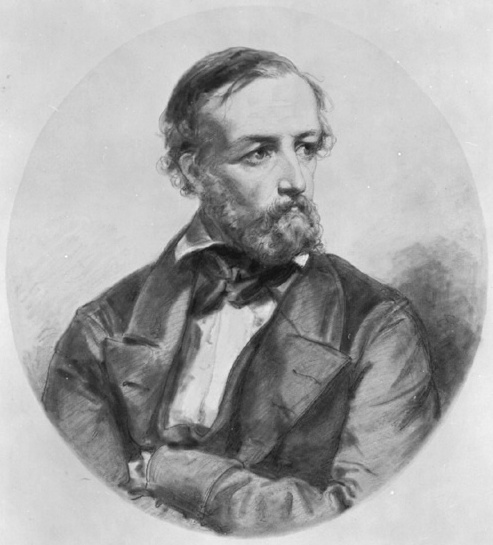}

{\tiny Unknown author---Unknown source---Public Domain}\end{center}

\caption{\textbf{\protect\href{https://en.wikipedia.org/wiki/Peter_Gustav_Lejeune_Dirichlet}{Peter Gustav Lejeune Dirichlet}}
(1805 - 1859)}\sindex[fam]{Dirichlet, Gustav}
\end{figure}

\begin{exercise}{Dirichlet Law and Order Statistics}{exercise10.8}

Let $\left(X_{i}\right)_{1\leqslant i\leqslant n}$ be a finite sequence
of independent real-valued random variables, each following the same
law $\mu$ and admitting a piecewise continuous density $f.$ Denote
by $F$ the cumulative distribution function. 

Define the function $r$ from $\mathbb{R}^{n}$ onto itself by
\[
\forall\left(x_{1},x_{2},\cdots,x_{n}\right)\in\mathbb{R}^{n},\,\,\,\,r\left(x_{1},x_{2},\cdots,x_{n}\right)=\left(x_{\left(1\right)},x_{\left(2\right)},\cdots,x_{\left(n\right)}\right),
\]
where the $x_{\left(i\right)}$ are the real numbers $x_{i}$ sorted
by nondecreasing values, that is the real numbers defined by
\[
\left\{ x_{1},x_{2},\cdots,x_{n}\right\} =\left\{ x_{\left(1\right)},x_{\left(2\right)},\cdots,x_{\left(n\right)}\right\} \,\,\,\,\text{and}\,\,\,\,x_{\left(1\right)}\leqslant x_{\left(2\right)}\leqslant\cdots\leqslant x_{\left(n\right)}.
\]
We define the random variable 
\[
\left(X_{\left(1\right)},X_{\left(2\right)},\cdots,X_{\left(n\right)}\right)=r\left(X_{1},X_{2},\cdots,X_{n}\right).
\]
If $1\leqslant k\leqslant n,$ $X_{\left(k\right)}$ is called the
\textbf{$k^{\text{th}}$ order statistic}\mindex{kth order statistic@$k^\textrm{th}$ order statistic}.

1. a. Compute the cumulative distribution function of $X_{\left(k\right)}$
as a function of $F$ and $f.$ 

b. Justify the existence of a density $f_{X_{\left(k\right)}}$ for
the random variable $X_{\left(k\right)},$ and give its explicit expression.

c. Identify the law of $X_{\left(k\right)}$ in the special case where
$\mu$ is the uniform law on $\left[0,t\right],$ with $t>0.$

2. Determine the law of the random variable $\left(X_{\left(1\right)},X_{\left(2\right)},\cdots,X_{\left(n\right)}\right)$
in the particular case where $\mu$ is the uniform law on $\left[0,t\right],$
with $t>0.$ 

This law is called the \textbf{\index{Dirichlet law}\mindex{law ! Dirichlet}Dirichlet\footnotemark\sindex[fam]{Dirichlet, Gustav}
law}.

\end{exercise}

\footnotetext{\textbf{\href{https://en.wikipedia.org/wiki/Peter_Gustav_Lejeune_Dirichlet}{Peter Gustav Lejeune Dirichlet}\sindex[fam]{Dirichlet, Gustav}}
(1805 - 1859) was a Prussian mathematician, who gave deep contributions
to number theory, by founding the number analytical theory and Fourier
seies theory. He also contributed to the mathematical analysis, and
gave the modern formal definition of a function. He taught in different
universities, from Berlin to Breslau and Götingen.}

\begin{exercise}{Multinomial Law}{exercise10.9}

Let $k\in\mathbb{N}^{\ast}$ be fixed. For each $n\in\mathbb{N}^{\ast},$
consider a partition $\left(A^{n}_{i}\right)_{1\leqslant i\leqslant k}$
of $\Omega,$ where $A^{n}_{i}\in\mathcal{A}.$ Assume that the families,
indexed on $n,$ constituted by the elements of these partitions are
independent. 

Also, suppose that
\[
\forall n\in\mathbb{N}^{\ast},\,\,\,\,\forall i\in\left\llbracket 1,k\right\rrbracket ,\,\,\,\,P\left(A^{n}_{i}\right)=p_{i},
\]
where $p_{i}>0$ and $\sum^{k}_{i=1}p_{i}=1.$ 

Define the random variables $X^{n}$ and $Y^{n}$ taking values in
$\mathbb{R}^{k}$ by
\[
X^{n}=\left(\begin{array}{c}
1_{A^{n}_{1}}\\
\vdots\\
1_{A^{n}_{k}}
\end{array}\right)\,\,\,\,\,\,\,\text{and}\,\,\,\,\,\,\,Y^{n}=\sum^{n}_{j=1}X^{j}.
\]

Denote $S=\left\{ e_{1},e_{2},\cdots,e_{k}\right\} $ the canonical
basis of $\mathbb{R}^{k}$ and
\[
D_{n}=\left\{ y\in\left\llbracket 0,n\right\rrbracket ^{k}:\,\sum^{k}_{j=1}y_{j}=n\right\} .
\]

1. Determine $P\left(X^{n}=e_{j}\right)$ for every $j$ such that
$1\leqslant j\leqslant k.$ 

Deduce, that for $y\in D_{n},$ the probability $P\left(Y^{n}=y\right).$ 

The law of $Y^{n}$ is called the \textbf{multinomial law with parameters
$n,$ $p_{1},p_{2},\cdots,p_{k-1}$\index{multinomial law}\mindex{law!multinomial}}
and is denoted 
\[
M\left(n,p_{1},p_{2},\cdots,p_{k-1}\right).
\]
 It generalizes the binomial law.

2. Determine the expectation and covariance matrix of $Y^{n}.$

\end{exercise}

\begin{exercise}{A Result Related To The Central Limit Theorem}{exercise10.10}

Let $\left(X_{n}\right)_{n\in\mathbb{N}^{\ast}}$ be a sequence of
independent real-valued random variables of same law, admitting a
second-order moment, centered.

For each $n\in\mathbb{N}^{\ast},$ define the random variable $Z_{n}$
by
\[
Z_{n}=\dfrac{1}{\sqrt{n}}\sum^{n}_{j=1}X_{j}.
\]

Assume that if the sequence $\left(Z_{n}\right)_{n\in\mathbb{N}^{\ast}}$
converges almost surely to a random variable $Z,$ it follows from
the central limit theorem---see Chapter \ref{chap:PartIIChap11}
on the convergences in law---that the law of $Z$ is the law $\mathcal{N}_{\mathbb{R}}\left(0,1\right).$ 

Deduce from this that the sequence $\left(Z_{n}\right)_{n\in\mathbb{N}^{\ast}}$
cannot converge almost surely.

\end{exercise}

\begin{exercise}{A Function of the Borel-Cantelli Lemma}{exercise10.11}

Let $\left(X_{n}\right)_{n\in\mathbb{N}^{\ast}}$ be a sequence of
independent, real-valued random variables. Prove the equivalence
\[
P\left(\sup_{n\in\mathbb{N}^{\ast}}X_{n}<+\infty\right)=1\Leftrightarrow\exists A>0:\sum^{+\infty}_{j=1}P\left(X_{n}>A\right)<+\infty.
\]

\end{exercise}

\begin{exercise}{Borel-Cantelli Lemma and Almost Sure Convergence}{exercise10.12}

Let $\left(X_{n}\right)_{n\in\mathbb{N}^{\ast}}$ be a sequence of
independent real-valued random variables, with same law. Prove the
equivalence
\[
P\left(\dfrac{X_{n}}{n}\to0\right)=1\Leftrightarrow\intop_{\Omega}\left|X_{1}\right|\text{d}P<+\infty.
\]

\begin{remark}{}{}

This property is used when studying the strong laws of large numbers.

\end{remark}

\end{exercise}

\begin{exercise}{Duplication and Symmetrisation}{exercise10.13}

To a real-valued random variable $X$ defined in a probabilized space
\preds, we associate the functions $\widehat{X}$ and $X^{s}$ defined
on the Cartesian product $\Omega\times\Omega$ and with values respectively
in $\mathbb{R}^{2}$ and $\mathbb{R},$ for every $\left(\omega,\omega^{\prime}\right)\in\Omega\times\Omega,$
by
\[
\widehat{X}\left(\omega,\omega^{\prime}\right)=\left(X\left(\omega\right),X\left(\omega^{\prime}\right)\right)\,\,\,\,\text{and}\,\,\,\,X^{s}\left(\omega,\omega^{\prime}\right)=X\left(\omega\right)-X\left(\omega^{\prime}\right).
\]
1. Verify that $\widehat{X}$ and $X^{s}$ are random variables defined
on the probabilized space $\left(\Omega\times\Omega,\mathcal{A}\otimes\mathcal{A},P\otimes P\right)$
and that the marginals $\widehat{X_{1}}$ and $\widehat{X_{2}}$ of
$\widehat{X}$ are independent and each has the same law than $X.$

2. a. Let $p\geqslant1.$ Prove that if $X\in\mathscr{L}^{p}\left(\Omega,\mathcal{A},P\right),$
then the symmetrized $X^{s}$ belongs to $\mathscr{L}^{p}\left(\Omega\times\Omega,\mathcal{A}\otimes\mathcal{A},P\otimes P\right).$ 

b. In particular, if $X\in\mathscr{L}^{2},$ compute the expectation
and variance of $X^{s}.$

3. Let $\left(X_{i}\right)_{i\in I}$ be a family of independent real-valued
random variables defined on $\left(\Omega,\mathcal{A},P\right).$
Prove that the associated random variables $\widehat{X_{i}}$---respectively
the symmetrized $X^{s}_{i}$---are $P\otimes P-$independent.

\begin{remark}{}{}

This property is used in particular when studying the convergence
$P-$almost surely of series of independent random variables.

\end{remark}

\end{exercise}

\section*{Solutions of Exercises}

\addcontentsline{toc}{section}{Solutions of Exercises}

\begin{solution}{}{solexercise10.1}

\textbf{1. Law of the random variable $\left(R,\Theta\right)$}

For every $f\in\mathscr{C}^{+}_{\mathscr{K}}\left(\mathbb{R}^{2}\right),$
by the transfer theorem,
\[
\mathbb{E}\left(f\left(R,\Theta\right)\right)=\mathbb{E}\left(f\circ\widehat{\psi}\left(X,Y\right)\right)=\intop_{\mathbb{R}^{2}}f\circ\widehat{\psi}\left(x,y\right)\text{d}P_{\left(X,Y\right)}\left(x,y\right).
\]
The random variables $X$ and $Y$ are independent and admit densities:
thus,
\[
\mathbb{E}\left(f\left(R,\Theta\right)\right)=\intop_{\mathbb{R}^{2}}f\circ\widehat{\psi}\left(x,y\right)f_{X}\left(x\right)f_{Y}\left(y\right)\text{d}\lambda_{2}\left(x,y\right).
\]

Since, in the plane, a half-line has Lebesgue measure zero,
\begin{align*}
\mathbb{E}\left(f\left(R,\Theta\right)\right) & =\intop_{\mathbb{R}^{2}\backslash\mathbb{R}^{+}\times\left\{ 0\right\} }f\circ\psi\left(x,y\right)f_{X}\left(x\right)f_{Y}\left(y\right)\text{d}\lambda_{2}\left(x,y\right)\\
 & =\intop_{\mathbb{R}^{2}\backslash\mathbb{R}^{+}\times\left\{ 0\right\} }f\circ\psi\left(x,y\right)\dfrac{1}{2\pi}\text{e}^{-\frac{x^{2}+y^{2}}{2}}\text{d}\lambda_{2}\left(x,y\right)
\end{align*}
Let make the change of variables to polar coordinates defined by the
diffeomorphism $\psi,$ of Jacobian $r;$ it follows that
\begin{align*}
\mathbb{E}\left(f\left(R,\Theta\right)\right) & =\intop_{\mathbb{R}^{+\ast}\times\left]0,2\pi\right[}f\left(r,\theta\right)\dfrac{1}{2\pi}\text{e}^{-\frac{r^{2}}{2}}r\text{d}\lambda_{2}\left(r,\theta\right),
\end{align*}
This is equivalent to
\begin{align*}
\mathbb{E}\left(f\left(R,\Theta\right)\right) & =\intop_{\mathbb{R}^{2}}f\left(r,\theta\right)\boldsymbol{1}_{\mathbb{R}^{+\ast}\times\left[0,2\pi\right]}\left(r,\theta\right)\dfrac{1}{2\pi}\text{e}^{-\frac{r^{2}}{2}}r\text{d}\lambda_{2}\left(r,\theta\right).
\end{align*}
Then, we conclude that $\left(R,\Theta\right)$ admits the density
$f_{\left(R,\Theta\right)}$ given by\boxeq{
\[
\forall\left(r,\theta\right)\in\mathbb{R}^{2},\,\,\,\,f_{\left(R,\Theta\right)}\left(r,\theta\right)=\boldsymbol{1}_{\left[0,2\pi\right]}\left(\theta\right)\dfrac{1}{2\pi}\boldsymbol{1}_{\mathbb{R}^{+}}\text{e}^{-\frac{r^{2}}{2}}r.
\]
}

\textbf{2. Independence of $R$ and $\Theta$} 

The density $f_{\left(R,\Theta\right)}$ is a direct product of two
nonnegative measurable functions. The random variables $R$ and $\Theta$
are thus independent and of respective density $f_{R}$ and $f_{\Theta}$
given by\boxeq{
\[
\forall r\in\mathbb{R},\,\,\,\,f_{R}\left(r\right)=\boldsymbol{1}_{\mathbb{R}^{+}}\left(r\right)\text{e}^{-\frac{r^{2}}{2}}r.
\]
\[
\forall\theta\in\mathbb{R},\,\,\,\,f_{\Theta}\left(\theta\right)=\boldsymbol{1}_{\left[0,2\pi\right]}\left(\theta\right)\dfrac{1}{2\pi}.
\]
}

\textbf{The law of $\Theta$ is the uniform law on $\left[0,2\pi\right].$}

\textbf{Independence of $S=R^{2}$ and $\Theta$ }

\boxeq{The random variables $S$ and $\Theta$ are also independent
as measurable functions of such random variables.}

\textbf{3. Law of $S$}

Finally, for every $f\in\mathscr{C}^{+}_{\mathscr{K}}\left(\mathbb{R}\right),$
by following the same steps as before,
\[
\mathbb{E}\left(f\left(S\right)\right)=\mathbb{E}\left(f\left(\mathbb{R}^{2}\right)\right)=\intop_{\mathbb{R}}f\left(r^{2}\right)\text{d}P_{R}\left(r\right).
\]

Since the Lebesgue measure of a singleton is zero,
\[
\mathbb{E}\left(f\left(S\right)\right)=\intop_{\mathbb{R}^{+\ast}}f\left(r^{2}\right)\text{e}^{-\frac{r^{2}}{2}}r\text{d}\lambda\left(r\right).
\]

By making the coordinate change defined by $s=r^{2},$ which is a
diffeomorphism from $\mathbb{R}^{+\ast}$ onto itself, of Jacobian
$\dfrac{1}{2\sqrt{s}},$\boxeq{
\[
\mathbb{E}\left(f\left(S\right)\right)=\intop_{\mathbb{R}}f\left(s\right)\boldsymbol{1}_{\mathbb{R}^{+\ast}}\left(s\right)\dfrac{1}{2}\text{e}^{-\frac{s}{2}}\text{d}\lambda\left(s\right),
\]
}Therefore, $S$ follows the exponential law $\exp\left(-\dfrac{1}{2}\right).$

\end{solution}

\begin{solution}{}{solexercise10.2}

Note that $P\left(S\geqslant0\right)=1.$ $\left(X,Y\right)$ is thus
defined almost surely.

For every $f\in\mathscr{C}^{+}_{\mathscr{K}}\left(\mathbb{R}^{2}\right),$
\[
\mathbb{E}\left(f\left(X,Y\right)\right)=\mathbb{E}\left(f\left(\sqrt{S}\cos\Theta,\sqrt{S}\sin\Theta\right)\right).
\]

By the transfer theorem,
\[
\mathbb{E}\left(f\left(X,Y\right)\right)=\intop_{\mathbb{R}^{2}}f\left(\sqrt{s}\cos\Theta,\sqrt{s}\sin\Theta\right)\text{d}P_{\left(S,\Theta\right)}\left(s,\theta\right).
\]

The random variables $S$ and $\theta$ are independent and both admit
densities; thus,
\[
\mathbb{E}\left(f\left(X,Y\right)\right)=\intop_{\mathbb{R}^{2}}f\left(\sqrt{s}\cos\Theta,\sqrt{s}\sin\Theta\right)f_{S}\left(s\right)f_{\Theta}\left(\theta\right)\text{d}\lambda_{2}\left(s,\theta\right).
\]

Hence,
\[
\mathbb{E}\left(f\left(X,Y\right)\right)=\intop_{\mathbb{R}^{+\ast}\times\left]0,2\pi\right[}f\left(\sqrt{s}\cos\Theta,\sqrt{s}\sin\Theta\right)\dfrac{1}{2\pi}\dfrac{1}{2}\text{e}^{-\frac{s}{2}}\text{d}\lambda_{2}\left(s,\theta\right).
\]

We now make the change of variables associated to the diffeomorphism
$H$ from $\mathbb{R}^{+\ast}\times\left]0,2\pi\right[$ onto $\mathbb{R}^{2}\backslash\mathbb{R}^{+}\times\left\{ 0\right\} $
defined by
\[
\forall\left(s,\theta\right)\in\mathbb{R}^{+\ast}\times\left]0,2\pi\right[,\,\,\,\,x=\sqrt{s}\cos\theta\,\,\,\,y=\sqrt{s}\sin\theta.
\]

Since
\[
H^{\prime}\left(s,\theta\right)=\left[\begin{array}{cc}
\dfrac{1}{2\sqrt{s}}\cos\theta & \,\,-\sqrt{s}\sin\theta\\
\dfrac{1}{2\sqrt{s}}\sin\theta & \,\,\sqrt{s}\cos\theta
\end{array}\right]
\]
and
\[
\left(H^{-1}\right)^{\prime}\left(x,y\right)=\left[H^{\prime}\left(H^{-1}\left(x,y\right)\right)\right]^{-1},
\]
the Jacobian of the diffeomorphism is equal to
\[
\det\left[\left(H^{-1}\right)^{\prime}\left(x,y\right)\right]=\dfrac{1}{\det H^{\prime}\left(H^{-1}\left(x,y\right)\right)}=2.
\]

Then, it follows that
\[
\mathbb{E}\left(f\left(X,Y\right)\right)=\intop_{\mathbb{R}^{2}\backslash\mathbb{R}^{+}\times\left\{ 0\right\} }f\left(x,y\right)\dfrac{1}{2\pi}\text{e}^{-\frac{x^{2}+y^{2}}{2}}\text{d}\lambda_{2}\left(x,y\right).
\]
Thus, since any half-line has a Lebesgue measure zero
\[
\mathbb{E}\left(f\left(X,Y\right)\right)=\intop_{\mathbb{R}^{2}}f\left(x,y\right)\dfrac{1}{2\pi}\text{e}^{-\frac{x^{2}+y^{2}}{2}}\text{d}\lambda_{2}\left(x,y\right).
\]
This proves that $\left(X,Y\right)$ admits the density $f_{\left(X,Y\right)}$
given by\boxeq{
\[
\forall\left(x,y\right)\in\mathbb{R}^{2},\,\,\,\,f_{\left(X,Y\right)}\left(x,y\right)=\dfrac{1}{\sqrt{2\pi}}\text{e}^{-\frac{x^{2}}{2}}\dfrac{1}{\sqrt{2\pi}}\text{e}^{-\frac{y^{2}}{2}}.
\]
}It follows that $X$ and $Y$ are two independent real-valued random
variables of same Gaussian law $\mathcal{N}_{\mathbb{R}}\left(0,1\right).$

\begin{remark}{}{}

We saw in exercise that if $U$ is a random variable following the
uniform law on $\left]0,1\right[,$ the random variable $-2\ln U$
follows the exponential law$\exp\left(1/2\right).$ 

This remark and the exercise then give a simulation method of two
independent random variables following a Gaussian law $\mathcal{N}_{\mathbb{R}}\left(0,1\right).$ 

Tr.N: This is the basis of the \index{Box-Muller transform}\href{https://en.wikipedia.org/wiki/Box\%25E2\%2580\%2593Muller_transform}{Box-Muller transform}
used to generate pair of random numbers distributed in a centered
reduced normal law, from a source of uniform law random numbers. This
method is used in the GCC compiler in the standard library of C++.

\end{remark}

\end{solution}

\begin{solution}{}{solexercise10.3}
\begin{itemize}
\item The cumulative distribution function of $Y_{n}$ is given, for every
$y\in\mathbb{R},$ by
\[
F_{Y_{n}}\left(y\right)=P\left(\bigcap^{n}_{i=1}\left(X_{i}\leqslant y\right)\right)=\begin{cases}
\left(\intop^{y}_{0}\text{e}^{-u}\text{d}u\right)^{n}, & \text{if }y>0,\\
0, & \text{otherwise.}
\end{cases}
\]
Thus, the random variable $Y_{n}$ admits a density given, except
in 0, by the derivative of $F_{Y_{n}},$ that is\boxeq{
\[
\forall y\in\mathbb{R},\,\,\,\,f_{Y_{n}}\left(y\right)=\boldsymbol{1}_{\mathbb{R}^{+}}\left(y\right)n\text{e}^{-y}\left(1-\text{e}^{-y}\right)^{n-1}.
\]
}
\item We now prove by induction on $n$ that $Y_{n}$ and $Z_{n}$ have
the same law. 
\begin{itemize}
\item $Y_{1}=Z_{1},$ thus these variables have the same law. 
\item Now, assume that at the rank $n,$ $Y_{n}$ and $Z_{n}$ have the
same law, and let us show that it still holds at the rank $n+1.$\\
We have
\[
Z_{n+1}=Z_{n}+\dfrac{X_{n+1}}{n+1}.
\]
\\
The random variables $Z_{n}$ and $\dfrac{X_{n+1}}{n+1}$ are independent---Tr.N:
$Z_{n}$ is the sum of variables that are independent of $X_{n+1}$.
Then $Z_{n+1}$ has a density which is the convolution product of
$Z_{n}$ and $\dfrac{X_{n+1}}{n+1}.$ As, by a classical computation,
the law of $\dfrac{X_{n+1}}{n+1}$ is the law $\exp\left(n+1\right),$
for every $z\in\mathbb{R},$
\begin{multline*}
f_{Z_{n+1}}\left(z\right)=\intop_{\mathbb{R}}\boldsymbol{1}_{\mathbb{R}^{+}}\left(y\right)n\text{e}^{-y}\left(1-\text{e}^{-y}\right)^{n-1}\\
\left[\left(n+1\right)\boldsymbol{1}_{\mathbb{R}^{+}}\left(z-y\right)\text{e}^{-\left(n+1\right)\left(z-y\right)}\right]\text{d}\lambda\left(y\right).
\end{multline*}
Thus, since
\[
\boldsymbol{1}_{\mathbb{R}^{+}}\left(y\right)\boldsymbol{1}_{\mathbb{R}^{+}}\left(z-y\right)=\boldsymbol{1}_{\mathbb{R}^{+}}\left(z\right)\boldsymbol{1}_{\left[0,z\right]}\left(y\right),
\]
it holds that
\[
f_{Z_{n+1}}\left(z\right)=\boldsymbol{1}_{\mathbb{R}^{+}}\left(z\right)n\left(n+1\right)\text{e}^{-\left(n+1\right)z}\intop^{z}_{0}\text{e}^{ny}\left(1-\text{e}^{-y}\right)^{n-1}\text{d}y.
\]
However,
\begin{align*}
\intop^{z}_{0}\text{e}^{ny}\left(1-\text{e}^{-y}\right)^{n-1}\text{d}y & =\intop^{z}_{0}\text{e}^{y}\left(\text{e}^{y}-1\right)^{n-1}\text{d}y\\
 & =\left[\dfrac{\left(\text{e}^{y}-1\right)^{n}}{n}\right]^{z}_{0}\\
 & =\dfrac{\left(\text{e}^{z}-1\right)^{n}}{n},
\end{align*}
which shows that
\[
f_{Z_{n+1}}\left(z\right)=\boldsymbol{1}_{\mathbb{R}^{+}}\left(z\right)\left(n+1\right)\text{e}^{-z}\left(1-\text{e}^{-z}\right)^{n}=f_{Y_{n+1}}\left(z\right),
\]
that is
\[
f_{Z_{n+1}}=f_{Y_{n+1}}.
\]
\end{itemize}
\item Then\boxeq{
\[
\mathbb{E}\left(Y_{n}\right)=\mathbb{E}\left(Z_{n}\right)=\sum^{n}_{i=1}\dfrac{1}{i}\mathbb{E}\left(X_{i}\right)=\sum^{n}_{i=1}\dfrac{1}{i}.
\]
}The random variables $Y_{n}$ and $Z_{n}$ have also same variance
and the random variables $X_{i}$ being independent,
\[
\sigma^{2}_{Y_{n}}=\sigma^{2}_{Z_{n}}=\sum^{n}_{i=1}\dfrac{1}{i^{2}}\sigma^{2}_{X_{i}},
\]
that is\boxeq{
\[
\sigma^{2}_{Y_{n}}=\sum^{n}_{i=1}\dfrac{1}{i^{2}}.
\]
}
\end{itemize}
\end{solution}

\begin{solution}{}{solexercise10.4}

\textbf{1. a. $Y_{n}$ admits a density}

We prove this by induction on $n.$
\begin{itemize}
\item \textbf{Initialization step} $\left(n=1\right)$\\
By the same classical method than in the previous exercices, we prove
that, for every $f\in\mathscr{C}^{+}_{\mathscr{K}}\left(\mathbb{R}\right),$
\[
\mathbb{E}\left(f\left(X^{2}_{1}\right)\right)=\intop_{\mathbb{R}}f\left(x\right)\boldsymbol{1}_{\mathbb{R}^{+}}\left(x\right)\dfrac{1}{\sqrt{2\pi}}\text{e}^{-\frac{x}{2}}x^{-\frac{1}{2}}\text{d}\lambda\left(x\right),
\]
which proves that $Y_{1}$ admits a density given by the relation
$\refpar{eq:density_chi_two}$ at the order 1. This is the chi-squared
with 1 degree of freedom.
\item \textbf{Inductive step}\\
Suppose that $Y_{n}$ admits a density $f_{Y_{n}}$ given by the relation
$\refpar{eq:density_chi_two}$ at the order n.\\
The random variables $Y_{n}$ and $X^{2}_{n+1}$ are independent and
$Y_{n+1}=Y_{n}+X^{2}_{n+1}.$ Then $Y_{n+1}$ has a density obtained
by convolution of the ones of $Y_{n}$ and $X^{2}_{n+1},$ and, moreover
as $X^{2}_{n+1}$ and $X^{2}_{1}$ have the same law, it follows that
\begin{multline*}
\forall y\in\mathbb{R},\,\,\,\,f_{Y_{n+1}}\left(y\right)=\intop_{\mathbb{R}}\boldsymbol{1}_{\mathbb{R}^{+}}\left(x\right)\dfrac{1}{2^{\frac{n}{2}}\Gamma\left(\frac{n}{2}\right)}\text{e}^{-\frac{x}{2}}x^{\frac{n}{2}-1}\\
\left[\boldsymbol{1}_{\mathbb{R}^{+}}\left(y-x\right)\dfrac{1}{\sqrt{2\pi}}\text{e}^{-\frac{y-x}{2}}\left(y-x\right)^{-\frac{1}{2}}\right]\text{d}\lambda\left(x\right).
\end{multline*}
It follows that, for every $y\in\mathbb{R},$
\[
f_{Y_{n+1}}\left(y\right)=\boldsymbol{1}_{\mathbb{R}^{+}}\left(y\right)\text{e}^{-\frac{y}{2}}\dfrac{1}{\sqrt{2\pi}}\dfrac{1}{2^{\frac{n}{2}}\Gamma\left(\frac{n}{2}\right)}\intop_{\left[0,y\right]}x^{\frac{n}{2}-1}\left(y-x\right)^{-\frac{1}{2}}\text{d}\lambda\left(x\right).
\]
Making the change of variables $x=uy$ for $y>0,$
\[
\intop_{\left[0,y\right]}x^{\frac{n}{2}-1}\left(y-x\right)^{-\frac{1}{2}}\text{d}\lambda\left(x\right)=y^{\frac{n+1}{2}-1}\intop_{\left[0,1\right]}u^{\frac{n}{2}-1}\left(1-u\right)^{-\frac{1}{2}}\text{d}\lambda\left(u\right).
\]
Thus,
\begin{equation}
\forall y\in\mathbb{R},\,\,\,\,f_{Y_{n+1}}\left(y\right)=K\boldsymbol{1}_{\mathbb{R}^{+}}\left(y\right)\text{e}^{-\frac{y}{2}}y^{\frac{n+1}{2}-1}\label{eq:f_Y_n+1}
\end{equation}
where
\[
K=\dfrac{1}{\sqrt{2\pi}}\dfrac{1}{2^{\frac{n}{2}}\Gamma\left(\frac{n}{2}\right)}\intop_{\left[0,1\right]}u^{\frac{n}{2}-1}\left(1-u\right)^{-\frac{1}{2}}\text{d}\lambda\left(u\right).
\]
To identify $K,$ it suffices to write that $f_{Y_{n+1}}$ is a density,
to integrate on $\mathbb{R}$ and to make the variable change $t=\dfrac{y}{2}.$
\\
This confirms that $Y_{n+1}$ follows the chi-squared distribution
with $n+1$ degrees of freedom, which completes the inductive step.
Hence, the asked result for every $n.$
\end{itemize}
\textbf{b. Expectation and variance of $Y_{n}.$}

We have
\[
\mathbb{E}\left(Y_{1}\right)=\mathbb{E}\left(X^{2}\right)=\sigma^{2}_{X_{1}}+\left(\mathbb{E}\left(X_{1}\right)\right)^{2}=1.
\]
Thus, by linearity of expectation,\boxeq{
\[
\mathbb{E}\left(Y_{n}\right)=\sum^{n}_{j=1}\mathbb{E}\left(X^{2}_{1}\right)=n.
\]
}

The random variables $X^{2}_{j}$ being independent,
\[
\sigma^{2}_{Y_{n}}=\sum^{n}_{j=1}\sigma^{2}_{X^{2}_{j}}.
\]
Now,
\[
\sigma^{2}_{X^{2}_{j}}=\mathbb{E}\left(X^{4}_{1}\right)-\left(\mathbb{E}\left(X^{2}_{1}\right)\right)^{2},
\]
 and a simple computation implies $\mathbb{E}\left(X^{4}_{1}\right)=3,$
it then follows that\boxeq{
\[
\sigma^{2}_{Y_{n}}=2n.
\]
}

\textbf{2. a. $T^{\prime}_{n}$ admits a density}

Note first that, since $Y_{n}$ admits a density, $T^{\prime}_{n}$
is defined $P-$almost surely. Moreover, $Y_{n}$ and $X_{n+1}$ are
independent. 

Then, for every $f\in\mathscr{C}^{+}_{\mathscr{K}}\left(\mathbb{R}\right),$
\[
\mathbb{E}\left(f\left(T^{\prime}_{n}\right)\right)=\intop_{\left(Y_{n}\neq0\right)}f\left(\dfrac{X_{n+1}}{\sqrt{Y_{n}}}\right)\text{d}P.
\]

By the transfer theorem,
\[
\mathbb{E}\left(f\left(T^{\prime}_{n}\right)\right)=\intop_{\mathbb{R}\times\mathbb{R}^{\ast}}f\left(\dfrac{x}{\sqrt{y}}\right)f_{X_{n+1}}\left(x\right)f_{Y_{n}}\left(y\right)\text{d}\lambda_{2}\left(x,y\right).
\]

Hence,
\begin{multline*}
\mathbb{E}\left(f\left(T^{\prime}_{n}\right)\right)=\intop_{\mathbb{R}\times\mathbb{R}^{+\ast}}f\left(\dfrac{x}{\sqrt{y}}\right)\dfrac{1}{\sqrt{2\pi}}\text{e}^{-\frac{x^{2}}{2}}\dfrac{1}{2^{\frac{n}{2}}\Gamma\left(\dfrac{n}{2}\right)}\text{e}^{-\frac{y}{2}}y^{\frac{n}{2}-1}\text{d}\lambda_{2}\left(x,y\right).
\end{multline*}
Making the change of variables associated to the diffeomorphism from
the open $\mathbb{R}\times\mathbb{R}^{+\ast}$ onto itself, defined
by
\[
\left\{ \begin{array}{l}
t=\dfrac{x}{\sqrt{y}}\\
z=y
\end{array}\right.\Leftrightarrow\left\{ \begin{array}{l}
x=t\sqrt{z}\\
y=z.
\end{array}\right.
\]
and with Jacobian
\[
\dfrac{D\left(x,y\right)}{D\left(t,z\right)}=\left|\begin{array}{cc}
\sqrt{z} & \dfrac{t}{2\sqrt{z}}\\
0 & 1
\end{array}\right|=\sqrt{z}.
\]

It follows that
\begin{align*}
\mathbb{E}\left(f\left(T^{\prime}_{n}\right)\right) & =\intop_{\mathbb{R}\times\mathbb{R}^{+\ast}}f\left(t\right)\dfrac{1}{\sqrt{2\pi}}\dfrac{1}{2^{\frac{n}{2}}\Gamma\left(\dfrac{n}{2}\right)}\text{e}^{-\frac{z}{2}\left(1+t^{2}\right)}z^{\frac{n+1}{2}-1}\text{d}\lambda_{2}\left(t,z\right).
\end{align*}
By the Fubini theorem
\begin{align*}
\mathbb{E}\left(f\left(T^{\prime}_{n}\right)\right) & =\intop_{\mathbb{R}}f\left(t\right)\left[\intop_{\mathbb{R}^{+\ast}}\dfrac{1}{\sqrt{2\pi}}\dfrac{1}{2^{\frac{n}{2}}\Gamma\left(\dfrac{n}{2}\right)}\text{e}^{-\frac{z}{2}\left(1+t^{2}\right)}z^{\frac{n+1}{2}-1}\text{d}\lambda\left(z\right)\right]\text{d}\lambda\left(t\right).
\end{align*}
Hence, after performing a new change of variables defined by $u=\dfrac{z}{2}\left(1+t^{2}\right)$
and after reduction,
\begin{align*}
\mathbb{E}\left(f\left(T^{\prime}_{n}\right)\right) & =\intop_{\mathbb{R}}f\left(t\right)\dfrac{1}{\sqrt{\pi}}\dfrac{1}{\Gamma\left(\dfrac{n}{2}\right)}\dfrac{1}{\left(1+t^{2}\right)^{\frac{n+1}{2}}}\left[\intop_{\mathbb{R}^{+\ast}}\text{e}^{-u}u^{\frac{n+1}{2}-1}\text{d}\lambda\left(z\right)\right]\text{d}\lambda\left(t\right).
\end{align*}
This proves that $T^{\prime}_{n}$ admits a density given by\boxeq{
\[
\forall t\in\mathbb{R},\,\,\,\,f_{T^{\prime}_{n}}\left(t\right)=\dfrac{1}{\sqrt{\pi}}\dfrac{\Gamma\left(\dfrac{n+1}{2}\right)}{\Gamma\left(\dfrac{n}{2}\right)}\dfrac{1}{\left(1+t^{2}\right)^{\frac{n+1}{2}}}.
\]
}

\textbf{b. $T_{n}$ admits a density}

Lastly, since $T_{n}=\sqrt{n}T^{\prime}_{n},$
\[
\forall t\in\mathbb{R},\,\,\,\,f_{T_{n}}\left(t\right)=f_{T^{\prime}_{n}}\left(\dfrac{t}{\sqrt{n}}\right)\dfrac{1}{\sqrt{n}}.
\]

This proves the relation $\refpar{eq:Student_density}.$

\textbf{c. Values of $n$ for which the expectation and variance of
$T_{n}$ exist and computation}

The random variables $\left|X_{n+1}\right|$ and $\dfrac{1}{\sqrt{Y_{n}}}$
are independent and nonnegative. Then, in $\overline{\mathbb{R}}^{+},$
\[
\intop_{\Omega}\left|T_{n}\right|\text{d}P=\sqrt{n}\left(\intop_{\Omega}\left|X_{n+1}\right|\text{d}P\right)\left(\intop_{\Omega}\dfrac{1}{\sqrt{Y_{n}}}\text{d}P\right).
\]
The first factor of the right-hand side is finite; the second, by
the transfer theorem, becomes
\[
\intop_{\Omega}\dfrac{1}{\sqrt{Y_{n}}}\text{d}P=\intop_{\mathbb{R}^{+\ast}}\dfrac{1}{\sqrt{y}}f_{Y_{n}}\left(y\right)\text{d}\lambda\left(y\right).
\]

It follows that
\[
\intop_{\Omega}\left|T_{n}\right|\text{d}P<+\infty
\]
if and only if $n>1.$

In this case, still by independence,\boxeq{
\[
\mathbb{E}\left(T_{n}\right)=\sqrt{n}\mathbb{E}\left(X_{n+1}\right)\mathbb{E}\left(\dfrac{1}{\sqrt{Y_{n}}}\right)=0.
\]
}

In a similar way, the random variables $X^{2}_{n+1}$ and $\dfrac{1}{Y_{n}}$
are independent and nonnegative. Then, in $\overline{\mathbb{R}}^{+},$
\[
\intop_{\Omega}T^{2}_{n}\text{d}P=n\left(\intop_{\Omega}X^{2}_{n+1}\text{d}P\right)\left(\intop_{\Omega}\dfrac{1}{Y_{n}}\text{d}P\right).
\]

The first factor of the right-hand side is finite; while the second,
by the transfer theorem can be written
\[
\intop_{\Omega}\dfrac{1}{Y_{n}}\text{d}P=\intop_{\mathbb{R}^{+\ast}}\dfrac{1}{y}f_{Y_{n}}\left(y\right)\text{d}\lambda\left(y\right).
\]

It follows that
\[
\intop_{\Omega}T^{2}_{n}\text{d}P<+\infty
\]
if and only if $n>2.$ 

If $n>2,$ still by independence
\[
\mathbb{E}\left(T^{2}_{n}\right)=n\mathbb{E}\left(X^{2}_{n+1}\right)\mathbb{E}\left(\dfrac{1}{Y_{n}}\right).
\]

Since
\[
\mathbb{E}\left(\dfrac{1}{Y_{n}}\right)=\intop^{+\infty}_{0}\dfrac{1}{y}\dfrac{1}{2^{\frac{n}{2}}\Gamma\left(\dfrac{n}{2}\right)}\text{e}^{-\frac{y}{2}}y^{\frac{n}{2}-1}\text{d}y,
\]
it follows that, after a change of variables,
\[
\mathbb{E}\left(\dfrac{1}{Y_{n}}\right)=\dfrac{1}{n-2},
\]
and thus\boxeq{
\[
\sigma^{2}_{T_{n}}=\dfrac{n}{n-2}.
\]
}

\end{solution}

\begin{solution}{}{solexercise10.5}

For every $f\in\mathscr{C}^{+}_{\mathscr{K}}\left(\mathbb{R}\right),$
by taking into account successively the independence of the variables
$X$ and $Y,$ and thus the independence of the random variables $\boldsymbol{1}_{\left(X=n\right)}$
and $f\left(n+Y\right),$ and by the transfer theorem,
\begin{align*}
\mathbb{E}\left(f\left(Z\right)\right) & =\sum_{n\in\mathbb{N}}\intop_{\left(X=n\right)}f\left(n+Y\right)\text{d}P\\
 & =\sum_{n\in\mathbb{N}}\left[P\left(X=n\right)\intop_{\Omega}f\left(n+Y\right)\text{d}P\right]\\
 & =\sum_{n\in\mathbb{N}}\left[P\left(X=n\right)\intop_{\mathbb{R}}f\left(n+y\right)\boldsymbol{1}_{\mathbb{R}^{+}}\left(y\right)\text{e}^{-y}\text{d}\lambda\left(y\right)\right].
\end{align*}
Let $q=1-p.$ Thus, by making a change of variables in each integral,
\begin{align*}
\mathbb{E}\left(f\left(Z\right)\right) & =\sum_{n\in\mathbb{N}}\left[P\left(X=n\right)\intop_{\mathbb{R}}f\left(z-n\right)\boldsymbol{1}_{\mathbb{R}^{+}}\left(z-n\right)\text{e}^{-\left(z-n\right)}\text{d}\lambda\left(z\right)\right]\\
 & =\intop_{\mathbb{R}}f\left(z\right)\boldsymbol{1}_{\mathbb{R}^{+}}\left(z\right)\left[\sum^{\left\lfloor z\right\rfloor }_{n=0}p\left(q\text{e}\right)^{n}\right]\text{e}^{-z}\text{d}\lambda\left(z\right).
\end{align*}
Thus, the random variable $Z$ admits a density $f_{Z}$ given by\boxeq{
\[
\forall z\in\mathbb{R},\,\,\,\,f_{Z}\left(z\right)=p\boldsymbol{1}_{\mathbb{R}^{+}}\left(z\right)\dfrac{1-\left(q\text{e}\right)^{\left\lfloor z\right\rfloor +1}}{1-q\text{e}}\text{e}^{-z}.
\]
}

\end{solution}

\begin{solution}{}{solexercise10.6}

The random variables $X$ and $Y$ are independent and admit densities.
Then $X+Y$ also admits a density, and consequently
\[
P\left(X+Y=0\right)=1
\]
so the function $\left(U,T\right)$ is defined almost surely.

For every $f\in\mathscr{C}^{+}_{\mathscr{K}}\left(\mathbb{R}^{2}\right),$
the transfer theorem can then be written as
\[
\mathbb{E}\left(f\left(U,T\right)\right)=\intop_{\mathbb{R}^{2}\backslash\Delta}f\left(x+y,\dfrac{x}{x+y}\right)f_{X}\left(x\right)f_{Y}\left(y\right)\text{d}\lambda_{2}\left(x,y\right),
\]
where $\Delta$ is the line 
\[
\left\{ \left(x,y\right)\in\mathbb{R}^{2}:\,x+y=0\right\} .
\]
This gives
\[
\mathbb{E}\left(f\left(U,T\right)\right)=K\intop_{\left(\mathbb{R}^{+\ast}\right)^{2}}f\left(x+y,\dfrac{x}{x+y}\right)\text{e}^{-p\left(x+y\right)}x^{a-1}y^{b-1}\text{d}\lambda_{2}\left(x,y\right)
\]
where $K$ is a positive constant depending on $a,b,p.$

Let us make the change of variables associated to the diffeomorphism
from $\left(\mathbb{R}^{+\ast}\right)^{2}$ onto $\mathbb{R}^{+\ast}\times\left]0,1\right[$
defined by
\[
\left\{ \begin{array}{c}
u=x+y\\
t=\dfrac{x}{x+y}
\end{array}\right.\Leftrightarrow\left\{ \begin{array}{c}
x=ut\\
y=u\left(1-t\right),
\end{array}\right.
\]
with Jacobian
\[
\dfrac{D\left(x,y\right)}{D\left(u,t\right)}=\left|\begin{array}{ccc}
t &  & u\\
1-t &  & -u
\end{array}\right|=-u.
\]

It follows that
\[
\mathbb{E}\left(f\left(U,T\right)\right)=K\intop_{\mathbb{R}^{+\ast}\times\left]0,1\right[}f\left(u,t\right)\text{e}^{-u}\left(ut\right)^{a-1}\left(u\left(1-t\right)\right)^{b-1}\left|-u\right|\text{d}\lambda_{2}\left(u,t\right).
\]

Hence,
\[
\mathbb{E}\left(f\left(U,T\right)\right)=K\intop_{\mathbb{R}^{2}}f\left(u,t\right)\boldsymbol{1}_{\mathbb{R}^{+\ast}}\left(u\right)\text{e}^{-u}ut^{a+b-1}\boldsymbol{1}_{\left]0,1\right[}\left(t\right)t^{a-1}\left(1-t\right)^{b-1}\text{d}\lambda_{2}\left(u,t\right).
\]

It follows that $\left(U,T\right)$ admits a density $f_{\left(U,T\right)}$
defined in every $\left(u,t\right)\in\mathbb{R}^{2}$ by\boxeq{
\[
f_{\left(U,T\right)}\left(u,t\right)=K\left(\boldsymbol{1}_{\mathbb{R}^{+\ast}}\left(u\right)\text{e}^{-u}u^{a+b-1}\right)\left(\boldsymbol{1}_{\left]0,1\right[}\left(t\right)t^{a-1}\left(1-t\right)^{b-1}\right).
\]
}

The function $f_{\left(U,T\right)}$ is a direct product of two nonnegative
measurable functions. The random variables $U$ and $T$ are thus
independent. Moreover, the law of $U$ is the Gamma law $\gamma\left(a+b,p\right)$
and the one of $T$ is the Beta law $B\left(a,b\right)$ of first
kind---supported by $\left[0,1\right].$

Concerning the law of $Z,$ one may note that
\[
Z=\dfrac{T}{1-T},
\]
and apply the now usual law computation method. We find that $Z$
admits a density $f_{Z}$ given by\boxeq{
\[
\forall z\in\mathbb{R},\,\,\,\,\,\,\,\,f_{Z}\left(z\right)=\boldsymbol{1}_{\mathbb{R}^{+}}\left(z\right)\dfrac{1}{B\left(a,b\right)}\dfrac{z^{a-1}}{\left(1+z\right)^{a+b}}.
\]
}

The law of $Z$ is the Beta law $\text{B}\left(a,b\right)$ of second
kind---supported by $\mathbb{R}^{+}.$ Of course, the random variables
$U$ and $Z$ are also independent!

\end{solution}

\begin{solution}{}{solexercise10.7}

\textbf{1. $\left(U,W\right)$ admits a density}

For any $g\in\mathscr{C}^{+}_{\mathscr{K}}\left(\mathbb{R}^{2}\right),$
by the independence of $X$ and $Y$ and the transfer theorem\footnotemark,
\begin{align*}
\mathbb{E}\left(g\left(U,W\right)\right) & =\intop_{\mathbb{R}^{2}}g\left(x\wedge y,\left|x-y\right|\right)f_{X}\left(x\right)f_{Y}\left(y\right)\text{d}\lambda_{2}\left(x,y\right)\\
 & =\intop_{\left(x\leqslant y\right)}g\left(x,y-x\right)f_{X}\left(x\right)f_{Y}\left(y\right)\text{d}\lambda_{2}\left(x,y\right)\\
 & \qquad\qquad+\intop_{\left(x>y\right)}g\left(y,x-y\right)f_{X}\left(x\right)f_{Y}\left(y\right)\text{d}\lambda_{2}\left(x,y\right).
\end{align*}

But, since $f_{X}=f_{Y}=f$ and that $\lambda_{2}\left(\left\{ x=y\right\} \right)=0,$
\[
\intop_{\left(x>y\right)}g\left(y,x-y\right)f_{X}\left(x\right)f_{Y}\left(y\right)\text{d}\lambda_{2}\left(x,y\right)=\intop_{\left(x\geqslant y\right)}g\left(y,x-y\right)f_{X}\left(y\right)f_{Y}\left(x\right)\text{d}\lambda_{2}\left(x,y\right),
\]
and thus,
\[
\mathbb{E}\left(g\left(U,W\right)\right)=2\intop_{\left(x\leqslant y\right)}g\left(y,x-y\right)f\left(x\right)f\left(y\right)\text{d}\lambda_{2}\left(x,y\right).
\]

Making the change of variables associated to the diffeomorphism from
the open $\mathbb{R}^{2}$ onto itself defined by
\[
\left\{ \begin{array}{l}
u=x\\
w=y-x
\end{array}\right.\Leftrightarrow\left\{ \begin{array}{l}
x=u\\
y=u+w
\end{array}\right.
\]
with Jacobian 1. It follows that
\[
\mathbb{E}\left(g\left(U,W\right)\right)=2\intop_{\left(0\leqslant w\right)}g\left(u,w\right)f\left(u\right)f\left(u+w\right)\text{d}\lambda_{2}\left(u,w\right),
\]
which shows that the random variable $\left(U,W\right)$ admits a
density $f_{\left(U,W\right)}$ given by\boxeq{
\[
\forall\left(u,w\right)\in\mathbb{R}^{2},\,\,\,\,f_{\left(U,W\right)}\left(u,w\right)=2\cdot\boldsymbol{1}_{\mathbb{R}^{+}}\left(w\right)f\left(u\right)f\left(u+w\right).
\]
}Then the random variables $U$ and $W$ admit densities, denoted
$f_{U}$ and $f_{W}$ given by\boxeq{
\[
\forall u\in\mathbb{R},\,\,\,\,f_{U}\left(u\right)=2f\left(u\right)\intop_{\mathbb{R}^{+}}f\left(u+w\right)\text{d}\lambda\left(w\right),
\]
}and\boxeq{
\[
\forall w\in\mathbb{R},\,\,\,\,f_{W}\left(w\right)=2\boldsymbol{1}_{\mathbb{R}^{+}}\left(w\right)\intop_{\mathbb{R}^{+}}f\left(u\right)f\left(u+w\right)\text{d}\lambda\left(u\right).
\]
}

\textbf{2. Proof that $U$ and $W$ admit densities when $f$ is bounded
and its restriction to $\mathbb{R}^{+}$ is continuous}

It follows that for the random variables $U$ and $W$ to be independent,
it is necessary and sufficient that for $\lambda_{2}-$almost any
$\left(u,w\right)\in\mathbb{R}\times\mathbb{R}^{+},$
\begin{equation}
f\left(u\right)f\left(u+w\right)=2\left(f\left(u\right)\intop_{\mathbb{R}^{+}}f\left(u+\alpha\right)\text{d}\lambda\left(\alpha\right)\right)\left(\intop_{\mathbb{R}}f\left(\alpha\right)f\left(\alpha+w\right)\text{d}\lambda\left(\alpha\right)\right).\label{eq:fufu+w}
\end{equation}
 We can easily check that this holds when $\mu$ is an exponential
law.

Conversely, suppose that the relation $\refpar{eq:fufu+w}$ is satisfied
for $\lambda_{2}-$almost every $\left(u,w\right)\in\mathbb{R}\times\mathbb{R}^{+}.$ 

Since $f\left(u\right)>0$ when $u\geqslant0,$ we have, after a change
of variables, for $\lambda_{2}-$almost every $\left(u,w\right)\in\mathbb{R}^{+}\times\mathbb{R}^{+},$
\begin{equation}
f\left(u+w\right)=2\left(\intop_{\left[u,+\infty\right[}f\left(\alpha\right)\text{d}\lambda\left(\alpha\right)\right)\left(\intop_{\mathbb{R}}f\left(\alpha\right)f\left(\alpha+w\right)\text{d}\lambda\left(\alpha\right)\right)\label{eq:fu+w}
\end{equation}

We further suppose that $f$ is bounded and that its restriction to
$\mathbb{R}^{+}$ is continuous. 

We are going to show, that this equality actually holds for every
$\left(u,w\right)\in\mathbb{R}^{+}\times\mathbb{R}^{+}.$ 

Since $f$ is bounded and its restriction to $\mathbb{R}^{+}$ is
continuous, the function 
\[
w\mapsto\intop_{\mathbb{R}}f\left(\alpha\right)f\left(\alpha+w\right)\text{d}\lambda\left(\alpha\right)
\]
is continuous on $\mathbb{R}^{+}$---by the dominated convergence
theorem. Since the functions $f$ and $u\mapsto\intop_{\left[u,+\infty\right[}f\left(\alpha\right)\text{d}\lambda\left(\alpha\right)$
are continuous on $\mathbb{R}^{+},$ the function
\[
\left(u,w\right)\mapsto f\left(u+w\right)-2\left(\intop_{\left[u,+\infty\right[}f\left(\alpha\right)\text{d}\lambda\left(\alpha\right)\right)\left(\intop_{\mathbb{R}}f\left(\alpha\right)f\left(\alpha+w\right)\text{d}\lambda\left(\alpha\right)\right)
\]
is continuous on $\mathbb{R}^{+}\times\mathbb{R}^{+}.$ 

It follows\footnotemark that the equality $\refpar{eq:fu+w}$ is
valid for every $\left(u,w\right)\in\mathbb{R}^{+}\times\mathbb{R}^{+}.$ 

In particular, we can take $w=0,$ which yields, for every $u\in\mathbb{R}^{+},$
\begin{equation}
f\left(u\right)=2C\intop_{\left[u,+\infty\right[}f\left(\alpha\right)\text{d}\lambda\left(\alpha\right),\label{eq:fu}
\end{equation}
where
\[
C=\intop_{\mathbb{R}}\left[f\left(\alpha\right)\right]^{2}\text{d}\lambda\left(\alpha\right)>0
\]
---the equality $\refpar{eq:fu}$ ensures that $C>0,$ since for
$u\geqslant0,$ $f\left(u\right)>0.$ 

Since $f$ is continuous on $\mathbb{R}^{+},$ the function 
\[
u\mapsto\intop_{\left[u,+\infty\right[}f\left(\alpha\right)\text{d}\lambda\left(\alpha\right)
\]
 is differentiable on $\mathbb{R}^{+}$ of derivative $f.$ 

We then have, by the equality $\refpar{eq:fu},$
\[
\forall u\in\mathbb{R}^{+},\,\,\,\,f^{\prime}\left(u\right)=-2Cf\left(u\right).
\]

The general nonzero solution of this differential equation is given
for every $u\in\mathbb{R}^{+},$ by
\[
f\left(u\right)=p\text{e}^{-2Cu},\,\,\,\,\text{where\,}p>0.
\]
The solution, that makes $f$ a probability density is then given,
for every $u\in\mathbb{R}^{+},$ by
\[
f\left(u\right)=2C\text{e}^{-2Cu},
\]
that is $\mu$ is an exponential law.

\begin{lemma}{}{}

A function $g$ continuous on $\mathbb{R}^{+}\times\mathbb{R}^{+}$
taking as value zero $\lambda_{d}-$almost everywhere, is in fact
equal to zero everywhere on $\mathbb{R}^{+}\times\mathbb{R}^{+}.$

\end{lemma}

\begin{proof}{}{}

Indeed, since $g$ is continuous in particular on $\mathbb{R}^{+\ast}\times\mathbb{R}^{+\ast},$
the set $O=\left\{ t\in\mathbb{R}^{+\ast}\times\mathbb{R}^{+\ast}:\,g\left(t\right)\neq0\right\} $
is an open subset of $\mathbb{R}^{+\ast}\times\mathbb{R}^{+\ast}.$ 

If $O$ is empty, $g$ is equal to 0 on $\mathbb{R}^{+\ast}\times\mathbb{R}^{+\ast}.$
Otherwise, there exists an open ball---thus, of positive measure---contained
in $O,$ which contradicts the fact that $g$ is either equal to zero
$\lambda_{d}-$almost everywhere. Hence, $g$ is zero on $\mathbb{R}^{+\ast}\times\mathbb{R}^{+\ast}.$
It follows from the continuity of $g$ on $\mathbb{R}^{+}\times\mathbb{R}^{+},$
that $g$ is also equal to zero on the axis.

\end{proof}

\textbf{3. $\mu$ is still an exponential if and only if $U$ and
$W$ are independent when $f$ is bounded and that its restriction
to $\mathbb{R}^{+}$ is continuous.}

If the random variables $U$ and $V$ are independent, for $\lambda_{2}-$almost
any $\left(u,v\right)\in\mathbb{R}^{+}\times\mathbb{R}^{+},$ we still
have the equality $\refpar{eq:fu+w}$ written as
\begin{equation}
f\left(u+w\right)=G\left(u\right)f_{W}\left(w\right),\label{eq:f_u+w_Gf_w}
\end{equation}
where $G$ is the continuous function defined on $\mathbb{R}^{+}$
by, for every $u\geqslant0,$
\begin{equation}
G\left(u\right)=\intop_{\left[u,+\infty\right[}f\left(\alpha\right)\text{d}\lambda\left(\alpha\right).\label{eq:G_u}
\end{equation}

Then, for $\lambda-$almost any $w\in\mathbb{R}^{+},$
\[
G\left(w\right)=\intop_{\left[w,+\infty\right[}f\left(\alpha\right)\text{d}\lambda\left(\alpha\right)=\intop_{\mathbb{R}^{+}}f\left(u+w\right)\text{d}\lambda\left(u\right)=\intop_{\mathbb{R}^{+}}G\left(u\right)f_{W}\left(w\right)\text{d}\lambda\left(u\right).
\]

Hence,
\begin{equation}
G\left(w\right)=mf_{W}\left(w\right),\label{eq:G_w}
\end{equation}
 where $m=\intop_{\mathbb{R}^{+}}G\left(u\right)\text{d}\lambda\left(u\right).$

Since the support of the density $f_{W}$ is included in $\mathbb{R}^{+},$
it follows from the equality $\refpar{eq:G_w}$ that the function
$G$ is integrable over $\mathbb{R}^{+},$ with a nonzero integral.
Hence, from the equality $\refpar{eq:f_u+w_Gf_w},$ for $\lambda_{2}-$almost
any $\left(u,w\right)\in\mathbb{R}^{+}\times\mathbb{R}^{+},$
\begin{equation}
f\left(u+w\right)=\dfrac{G\left(u\right)G\left(w\right)}{m}.\label{eq:f_u+wGuGw}
\end{equation}
Then taking into account the definition of $G,$ for every $u\in\mathbb{R}^{+}$
and for every $v\in\mathbb{R}^{+},$
\[
G\left(u+v\right)=\intop_{\left[u+v,+\infty\right[}f\left(w\right)\text{d}\lambda\left(w\right)=\intop_{\left[v,+\infty\right[}f\left(u+w\right)\text{d}\lambda\left(w\right),
\]
 and thus, by the equality $\refpar{eq:f_u+wGuGw},$ for $\lambda-$almost
any $u\in\mathbb{R}^{+}$ and for every $v\in\mathbb{R}^{+}.$
\[
G\left(u+v\right)=\intop_{\left[v,+\infty\right[}\dfrac{G\left(u\right)G\left(w\right)}{m}\text{d}\lambda\left(w\right)=\dfrac{G\left(u\right)}{m}\intop_{\left[v,+\infty\right[}G\left(w\right)\text{d}\lambda\left(w\right).
\]
By continuity of $G,$ we then have for every $u\in\mathbb{R}^{+}$
and every $v\in\mathbb{R}^{+}.$
\[
G\left(u+v\right)=\dfrac{G\left(u\right)}{m}\intop_{\left[v,+\infty\right[}G\left(w\right)\text{d}\lambda\left(w\right).
\]

Since $G\left(0\right)=1,$ it follows that for every $v\in\mathbb{R}^{+},$
\[
G\left(v\right)=\dfrac{1}{m}\intop_{\left[v,+\infty\right[}G\left(w\right)\text{d}\lambda\left(w\right),
\]
which implies that $G$ is differentiable---G is continous---and
that
\[
G^{\prime}\left(v\right)=-\dfrac{G\left(v\right)}{m}.
\]

Taking into account that $G\left(0\right)=1,$ it follows that for
every $v\in\mathbb{R}^{+},$
\[
G\left(v\right)=\text{e}^{-\frac{v}{m}}.
\]

Hence, by equality $\refpar{eq:G_u},$ for $\lambda-$almost every
$u\in\mathbb{R}^{+},$
\[
f\left(u\right)=-G^{\prime}\left(u\right)=\dfrac{1}{m}\text{e}^{-\frac{u}{m}},
\]
which proves that $\mu$ is still the exponential law.

\end{solution}

\addtocounter{footnote}{-1}

\footnotetext{Classically, we denote: $x\wedge y=\min\left(x,y\right)$
and $x\vee y=\max\left(x,y\right)$}

\stepcounter{footnote}

\footnotetext{See the lemma.}

\begin{solution}{}{solexercise10.8}

\textbf{1. a. Computation of the cumulative distribution function
of $X_{\left(k\right)}$ as a function of $F$ and $f$ }

We first note that, since the random variable $\left(X_{1},X_{2},\cdots,X_{n}\right)$
admits a density,
\[
P\left(X_{\left(1\right)}<X_{\left(2\right)}<\cdots<X_{\left(n\right)}\right)=1.
\]

Then, for every $y\in\mathbb{R},$
\[
P\left(X_{\left(k\right)}\leqslant y\right)=P\left[\bigcup_{\substack{J\in\mathcal{P}_{f}\left(\left\{ 1,2,\cdots,n\right\} \right)\\
\left|J\right|\geqslant k
}
}B_{J}\right]=\sum^{n}_{i=k}P\left[\biguplus_{\substack{J\in\mathcal{P}_{f}\left(\left\{ 1,2,\cdots,n\right\} \right)\\
\left|J\right|=i
}
}B_{J}\right]
\]
where
\[
B_{J}=\left[\bigcap_{j\in J}\left(X_{j}\leqslant y\right)\right]\bigcap\left[\bigcap_{j\in J^{c}}\left(X_{j}>y\right)\right],
\]
and note that if $\left|J\right|=\left|J^{\prime}\right|,$ then $B_{J}\bigcap B_{J^{\prime}}=\emptyset.$

Since, for fixed $i,$ all subsets $J$ of size $i$ yield the same
probability $P\left(B_{J}\right),$ it follows that\boxeq{
\[
\forall y\in\mathbb{R},\,\,\,\,F_{X_{\left(k\right)}}\left(y\right)=\sum^{n}_{i=k}\binom{n}{i}\left[F\left(y\right)\right]^{i}\left[1-F\left(y\right)\right]^{n-i}.
\]
}

\textbf{b. Existence of a density $f_{X_{\left(k\right)}}$ and explicit
expression}

The function $F$ is differentiable with derivative $f;$ the same
holds for $F_{X_{\left(k\right)}}.$ Thus, the random variable $X_{\left(k\right)}$
admits a density $f_{X_{\left(k\right)}},$ the derivative of $F_{X_{\left(k\right)}},$
which is given, for every $y\in\mathbb{R},$ by
\begin{align*}
f_{X_{\left(k\right)}}\left(y\right) & =\sum^{n}_{i=k}i\binom{n}{i}f\left(y\right)\left[F\left(y\right)\right]^{i-1}\left[1-F\left(y\right)\right]^{n-i}\\
 & \qquad\qquad-f\left(y\right)\sum^{n-1}_{i=k}\left(n-i\right)\binom{n}{i}\left[F\left(y\right)\right]^{i}\left[1-F\left(y\right)\right]^{n-i-1},
\end{align*}
Now, using the identities
\[
i\binom{n}{i}=n\binom{n-1}{i-1}\,\,\,\,\,\,\,\,\text{and}\,\,\,\,\,\,\,\,\left(n-i\right)\binom{n}{i}=n\binom{n-1}{i}
\]
and changing the index in the first sum via $j=i-1,$ we get
\begin{multline*}
f_{X_{\left(k\right)}}\left(y\right)=nf\left(y\right)\sum^{n-1}_{j=k-1}\binom{n-1}{j}\left[F\left(y\right)\right]^{j}\left[1-F\left(y\right)\right]^{n-j-1}\\
-nf\left(y\right)\sum^{n-1}_{i=k}\binom{n-1}{i}\left[F\left(y\right)\right]^{i}\left[1-F\left(y\right)\right]^{n-i-1}.
\end{multline*}
After simplification, it follows that\boxeq{
\[
\forall y\in\mathbb{R},\,\,\,\,f_{X_{\left(k\right)}}\left(y\right)=nf\left(y\right)\left\{ \binom{n-1}{k-1}\left[F\left(y\right)\right]^{k-1}\left[1-F\left(y\right)\right]^{n-k}\right\} .
\]
}

\textbf{c. Case where $\mu$ is the uniform law on $\left[0,t\right]$
where $t>0$}

In particular, if $\mu$ is the uniform law on $\left[0,t\right],$
then\boxeq{
\[
\forall y\in\mathbb{R},\,\,\,\,f_{X_{\left(k\right)}}\left(y\right)=\dfrac{n!}{\left(k-1\right)!\left(n-k\right)!}\boldsymbol{1}_{\left]0,t\right]}\left(y\right)\dfrac{1}{t}\left(\dfrac{y}{t}\right)^{k-1}\left(1-\dfrac{y}{t}\right)^{n-k}
\]
}which shows that the random variable $X_{\left(k\right)}$ follows
a Beta law of the first kind on the interval $\left[0,t\right],$
with parameters $k$ and $n-k+1.$

\textbf{2. Law of the random variable $\left(X_{\left(1\right)},X_{\left(2\right)},\cdots,X_{\left(n\right)}\right)$
in the special case, where $\mu$ is the uniform law on $\left[0,t\right]$
where $t>0$}

Now let us determine the law of the random variable $\left(X_{\left(1\right)},X_{\left(2\right)},\cdots,X_{\left(n\right)}\right)$
in the special case, where $\mu$ is the uniform law on $\left[0,t\right]$
where $t>0.$

Denote $\mathscr{S}_{n}$ the simplex of $\mathbb{R}^{n},$ 
\[
\left\{ x_{1}\leqslant x_{2}\leqslant\cdots\leqslant x_{n}\right\} .
\]

For any permutation $\sigma\in\Sigma_{n},$ define the isometry $\mathscr{P}_{\sigma}$
of $\mathbb{R}^{n}$ onto itself by
\[
\forall\left(x_{1},x_{2},\cdots,x_{n}\right)\in\mathbb{R}^{n},\,\,\,\,\mathscr{P_{\sigma}}\left(x_{1},x_{2},\cdots,x_{n}\right)=\left(x_{\sigma\left(1\right)},x_{\sigma\left(2\right)},\cdots,x_{\sigma\left(n\right)}\right).
\]

We have the equality
\begin{equation}
\mathbb{R}^{n}=\biguplus_{\sigma\in\Sigma_{n}}\mathscr{P}^{-1}_{\sigma}\left(\mathscr{S}_{n}\right).\label{eq:part_r_n}
\end{equation}

The random variable $\left(X_{1},X_{2},\cdots,X_{n}\right)$ admits
the density $f_{C_{n}}=\dfrac{1}{t^{n}}\boldsymbol{1}_{\left[0,t\right]^{n}}.$ 

For every $g\in\mathscr{C}^{+}_{\mathscr{K}}\left(\mathbb{R}^{n}\right),$
\[
\mathbb{E}\left(g\left(X_{\left(1\right)},X_{\left(2\right)},\cdots,X_{\left(n\right)}\right)\right)=\mathbb{E}\left(\left(g\circ r\right)\left(X_{1},X_{2},\cdots,X_{n}\right)\right).
\]
Thus, by the transfer theorem,
\begin{align*}
 & \mathbb{E}\left(g\left(X_{\left(1\right)},X_{\left(2\right)},\cdots,X_{\left(n\right)}\right)\right)\\
 & \qquad=\intop_{\mathbb{R}^{n}}\left(g\circ r\right)\left(x_{1},x_{2},\cdots,x_{n}\right)f_{C_{n}}\left(x_{1},x_{2},\cdots,x_{n}\right)\text{d}\lambda_{n}\left(x_{1},x_{2},\cdots,x_{n}\right).
\end{align*}

By denoting $\uuline{x_{n}}=\left(x_{1},x_{2},\cdots,x_{n}\right),$
and by taking into account the equality $\refpar{eq:part_r_n},$ we
obtain

\begin{align*}
 & \mathbb{E}\left(g\left(X_{\left(1\right)},X_{\left(2\right)},\cdots,X_{\left(n\right)}\right)\right)\\
 & \qquad\qquad=\sum_{\sigma\in\Sigma_{n}}\intop_{\mathscr{P}^{-1}_{\sigma}\left(\mathscr{S}_{n}\right)}\left(g\circ r\right)\left(\uuline{x_{n}}\right)f_{C_{n}}\left(\uuline{x_{n}}\right)\text{d}\lambda_{n}\left(\uuline{x_{n}}\right)\\
 & \qquad\qquad=\sum_{\sigma\in\Sigma_{n}}\intop_{\mathbb{R}^{n}}\boldsymbol{1}_{\mathscr{S}_{n}}\mathscr{P}_{\sigma}\left(\uuline{x_{n}}\right)\left(g\circ r\right)\left(\uuline{x_{n}}\right)f_{C_{n}}\left(\uuline{x_{n}}\right)\text{d}\lambda_{n}\left(\uuline{x_{n}}\right).
\end{align*}

Now, make the change of variables $\uuline{y_{n}}=\mathscr{P}_{\sigma}\left(\uuline{x_{n}}\right)$
defined by the isometry $\mathscr{P}_{\sigma}$---of Jacobian $\pm1$---then
yields

\begin{multline*}
\mathbb{E}\left(g\left(X_{\left(1\right)},X_{\left(2\right)},\cdots,X_{\left(n\right)}\right)\right)\\
=\sum_{\sigma\in\Sigma_{n}}\intop_{\mathbb{R}^{n}}\boldsymbol{1}_{\mathscr{S}_{n}}\left(\uuline{y_{n}}\right)\left(g\circ r\right)\left(\mathscr{P}^{-1}_{\sigma}\left(\uuline{y_{n}}\right)\right)f_{C_{n}}\left(\mathscr{P}^{-1}_{\sigma}\left(\uuline{y_{n}}\right)\right)\text{d}\lambda_{n}\left(\uuline{y_{n}}\right).
\end{multline*}
Since the function $\left(g\circ r\right)\cdot f_{C_{n}}$ is invariant
under any permutation $\mathscr{P}^{-1}_{\sigma},$ the integrals
in the sum are independent of $\sigma.$ Thus,
\[
\mathbb{E}\left(g\left(X_{\left(1\right)},X_{\left(2\right)},\cdots,X_{\left(n\right)}\right)\right)=n!\intop_{\mathscr{S}_{n}}\left(g\circ r\right)\left(\uuline{y_{n}}\right)f_{C_{n}}\left(\uuline{y_{n}}\right)\text{d}\lambda_{n}\left(\uuline{y_{n}}\right).
\]
Taking into account the value of the density $f_{C_{n}},$ we conclude
that the random variable $\left(X_{\left(1\right)},X_{\left(2\right)},\cdots,X_{\left(n\right)}\right)$
admits the density $f_{\left(X_{\left(1\right)},X_{\left(2\right)},\cdots,X_{\left(n\right)}\right)}$
given by\boxeq{
\[
f_{\left(X_{\left(1\right)},X_{\left(2\right)},\cdots,X_{\left(n\right)}\right)}=\dfrac{n!}{t^{n}}\boldsymbol{1}_{\left(0\leqslant x_{1}\leqslant x_{2}\leqslant\cdots\leqslant x_{n}\leqslant t\right)}.
\]
}

\end{solution}

\begin{solution}{}{solexercise10.9}

\textbf{1. Determination of $P\left(X^{n}=e_{j}\right)$ and $P\left(Y^{n}=y\right)$}

We have
\[
P\left(X^{n}=e_{j}\right)=P\left(A^{n}_{j}\right)=p_{j}.
\]

Now, since
\[
\sum^{k}_{i=1}Y^{n}_{i}=\sum^{k}_{i=1}\sum^{n}_{j=1}X^{j}_{i}=\sum^{n}_{j=1}\left(\sum^{k}_{i=1}\boldsymbol{1}_{A^{j}_{i}}\right)=\sum^{n}_{j=1}\boldsymbol{1}_{\Omega}=n,
\]
it follows that $Y^{n}\left(\Omega\right)\subset D_{n}.$

Let $y\in D_{n}.$ Denoting $\uuline{J}^{y_{1},y_{2},\cdots,y_{k}}$
the set of partitions of $\left\{ 0,1,2,\cdots,n\right\} $ with $k$
elements $J_{1},J_{2},\cdots,J_{k}$ such that $\left|J_{1}\right|=y_{1},\left|J_{2}\right|=y_{2},\cdots,\left|J_{k}\right|=y_{k},$
\begin{multline*}
\left(Y^{n}=y\right)=\\
\biguplus_{\left(J_{1},J_{2},\cdots,J_{k}\right)\in\uuline{J}^{y_{1},y_{2},\cdots,y_{k}}}\left[\bigcap_{j_{1}\in J_{1}}\left(X^{j_{1}}=e_{1}\right)\cap\bigcap_{j_{2}\in J_{2}}\left(X^{j_{2}}=e_{2}\right)\cap\cdots\cap\bigcap_{j_{k}\in J_{k}}\left(X^{j_{k}}=e_{k}\right)\right],
\end{multline*}

Using the independence of the random variables $X^{j},$ it yields
\[
P\left(Y^{n}=y\right)=\sum_{\left\{ J_{1},J_{2},\cdots,J_{k}\right\} \in\uuline{J}^{y_{1},y_{2},\cdots,y_{k}}}p^{y_{1}}_{1}p^{y_{2}}_{2}\cdots p^{y_{k}}_{k}.
\]

The number of such partitions is
\[
\#\uuline{J}^{y_{1},y_{2},\cdots,y_{k}}=\binom{n}{y_{1}}\binom{n-y_{1}}{y_{2}}\binom{n-\left(y_{1}+y_{2}\right)}{y_{3}}\cdots\binom{n-\left(y_{1}+y_{2}+\cdots+y_{k-2}\right)}{y_{k-1}}
\]
As, $y_{1}+y_{2}+\cdots+y_{k-1}=n-y_{k},$
\[
\#\uuline{J}^{y_{1},y_{2},\cdots,y_{k}}=\dfrac{n!}{y_{1}!y_{2}!\cdots y_{k}!}.
\]

---we recognize the multinomial coefficient. Thus\boxeq{
\[
\forall y\in D_{n},\,\,\,\,P\left(Y^{n}=y\right)=\dfrac{n!}{y_{1}!y_{2}!\cdots y_{k}!}p^{y_{1}}_{1}p^{y_{2}}_{2}\cdots p^{y_{k}}_{k}.
\]
}

\textbf{2. Expectation and covariance matrix of $Y^{n}$}

The law of $Y^{n}$ is
\[
P_{Y^{n}}=\sum_{y\in D_{n}}P\left(Y^{n}=y\right)\delta_{y}.
\]

We compute the expectation of $X^{n}$
\[
\mathbb{E}\left(X^{n}\right)=\left(\begin{array}{c}
\mathbb{E}\left(1_{A^{n}_{1}}\right)\\
\vdots\\
\mathbb{E}\left(1_{A^{n}_{k}}\right)
\end{array}\right)=\left(\begin{array}{c}
p_{1}\\
\vdots\\
p_{k}
\end{array}\right),
\]
so\boxeq{
\[
\mathbb{E}\left(Y^{n}\right)=\left(\begin{array}{c}
np_{1}\\
\vdots\\
np_{k}
\end{array}\right).
\]
}

The covariance matric $C_{X^{n}}$ has for entries
\begin{itemize}
\item Diagonal:
\[
\left(C_{X^{n}}\right)_{i,i}=\sigma^{2}_{X^{n}_{i}}=\mathbb{E}\left(\boldsymbol{1}^{2}_{A^{n}_{i}}\right)-\left[\mathbb{E}\left(\boldsymbol{1}_{A^{n}_{i}}\right)\right]^{2}=p_{i}\left(1-p_{i}\right),
\]
\item Off-diagonal for $i\neq j:$
\[
\left(C_{X^{n}}\right)_{i,j}=\text{cov}\left(X^{n}_{i},X^{n}_{j}\right)=\left[\mathbb{E}\left(\boldsymbol{1}_{A^{n}_{i}}\boldsymbol{1}_{A^{n}_{j}}\right)\right]-\left[\mathbb{E}\left(\boldsymbol{1}_{A^{n}_{i}}\right)\right]\left[\mathbb{E}\left(\boldsymbol{1}_{A^{n}_{j}}\right)\right]=-p_{i}p_{j}.
\]
\end{itemize}
Since $Y^{n}$ is the sum of independent random variables $X^{j},$
the covariance matrix of $Y^{n}$ is the sum of $X^{j}$ covariance
matrices, thus\boxeq{
\[
C_{Y^{n}}=nC_{X^{n}}.
\]
}

\end{solution}

\begin{solution}{}{solexercise10.10}

Define, for every $n\in\mathbb{N}^{\ast},$ the \salg 
\[
\mathcal{A}_{n}=\sigma\left(X_{p}:\,p\geqslant n\right)
\]
 and the asymptotic \salg 
\[
\mathcal{A}_{\infty}=\bigcap_{n\in\mathbb{N}^{\ast}}\mathcal{A}_{n}.
\]

For any fixed $n_{0},$ and for $p\geqslant n_{0},$ define the random
variable 
\[
Y_{n_{0},p}=\dfrac{1}{\sqrt{p}}\sum^{p}_{j=n_{0}}X_{j}.
\]

Let $\left(U_{n}\longrightarrow\right)$ denote the set where the
sequence of random variables $\left(U_{n}\right)_{n\in\mathbb{N}^{\ast}}$
converges, we have $\left(Y_{n_{0},p}\longrightarrow\right)\in\mathcal{A}_{n_{0}}.$ 

However
\[
\left(Z_{p}\longrightarrow\right)=\left(Y_{n_{0},p}\longrightarrow\right),
\]
so for any fixed $n_{0},$ $\left(Z_{p}\longrightarrow\right)\in\mathcal{A}_{n_{0}}.$

Since this holds for any $n_{0},$ it follows that $\left(Z_{p}\longrightarrow\right)\in\mathcal{A}_{\infty}.$

Therefore the sequence $\left(Z_{n}\right)_{n\in\mathbb{N}^{\ast}}$
converges almost surely to a random variable $Z,$ then $Z$ is almost
surely equal to a random variable $\mathcal{A}_{\infty}-$measurable.
By the law of the all or nothing, $Z$ equals almost surely to a constant,
which is in contradiction to the fact that the law of $Z$ is the
law $\mathcal{N}_{\mathbb{R}}\left(0,1\right).$

\end{solution}

\begin{solution}{}{solexercise10.11}

Suppose there exists $A>0$ such that
\[
\sum^{+\infty}_{j=1}P\left(X_{n}>A\right)<+\infty,
\]

Then, by the Borel-Cantelli lemma---Lemma $\ref{lm:borel-cantelli}$---we
have:
\[
P\left(\limsup_{n\to+\infty}\left(X_{n}>A\right)\right)=0,
\]
 which implies
\[
P\left(\liminf_{n\to+\infty}\left(X_{n}\leqslant A\right)\right)=1.
\]

Now, since
\[
\liminf_{n\to+\infty}\left(X_{n}\leqslant A\right)\subset\left(\sup_{n\in\mathbb{N}}X_{n}<+\infty\right),
\]

it follows that
\[
P\left(\sup X_{n}<+\infty\right)=1.
\]

Conversely, suppose now, that for every $A>0,$ then
\[
\sum^{+\infty}_{j=1}P\left(X_{n}>A\right)<+\infty.
\]

Since the random variables $X_{n}$ are independent, the ``converse''
of the Borel-Cantelli lemma applies. Thus
\[
\forall A>0,\,\,\,\,P\left(\limsup_{n\to+\infty}\left(X_{n}>A\right)\right)=1.
\]

Because $\mathbb{N}^{\ast}$ is countable, we can write that
\[
P\left(\bigcap_{A\in\mathbb{N}^{\ast}}\left(\limsup_{n\to+\infty}\left(X_{n}>A\right)\right)\right)=1.
\]

Now, observe the sequence of implications
\begin{align*}
\omega & \in\bigcap_{A\in\mathbb{N}^{\ast}}\left(\limsup_{n\to+\infty}\left(X_{n}>A\right)\right)\\
 & \Leftrightarrow\forall A\in\mathbb{N}^{\ast},\,\forall n\in\mathbb{N}^{\ast},\,\exists p\geqslant n:\,X_{p}\left(\omega\right)\geqslant A\\
 & \Rightarrow\forall A\in\mathbb{N}^{\ast},\,\,\,\,\sup X_{n}\left(\omega\right)>A\Rightarrow X_{n}\left(\omega\right)=+\infty.
\end{align*}

It follows that
\[
P\left(\sup X_{n}=+\infty\right)=1,
\]
which gives
\[
P\left(\sup X_{n}<+\infty\right)=0.
\]

Hence, if $P\left(\sup X_{n}=+\infty\right)=1,$ there exists $A>0,$
such that $\sum^{+\infty}_{j=1}P\left(X_{n}>A\right)<+\infty.$

\end{solution}

\begin{solution}{}{solexercise10.12}

Recall first---see Exercise $\ref{exo:exercise9.5}$---that using
the Fubini theorem, we can establish the equality, in $\mathbb{R}^{+},$
\[
\intop_{\Omega}\left|X_{1}\right|\text{d}P=\intop_{\mathbb{R}^{+}}P\left(\left|X_{1}\right|>x\right)\text{d}\lambda\left(x\right).
\]

Since the function $x\mapsto P\left(\left|X_{1}\right|>x\right)$
is nonincreasing, we obtain the following double inequality
\begin{equation}
\epsilon\sum^{+\infty}_{n=0}P\left(\left|X_{1}\right|>\left(n+1\right)\epsilon\right)\leqslant\intop_{\Omega}\left|X_{1}\right|\text{d}P\leqslant\epsilon\sum^{+\infty}_{n=0}P\left(\left|X_{1}\right|>n\epsilon\right).\label{eq:double_ineq_intopom_X_1}
\end{equation}
In particular, taking $\epsilon=1,$ it yields
\[
\intop_{\Omega}\left|X_{1}\right|\text{d}P<+\infty\,\,\,\,\Leftrightarrow\,\,\,\,\sum^{+\infty}_{n=1}P\left(\left|X_{1}\right|>n\right)<+\infty.
\]
Because the random variables $X_{n}$ follow the same law, this equivalence
becomes
\[
\intop_{\Omega}\left|X_{1}\right|\text{d}P<+\infty\,\,\,\,\Leftrightarrow\,\,\,\,\sum^{+\infty}_{n=1}P\left(\left|X_{n}\right|>n\right)<+\infty.
\]

Furthermore, since the random variables $X_{n}$ are independent,
we apply the Borel-Cantelli lemma to obtain the equivalences
\[
\intop_{\Omega}\left|X_{1}\right|\text{d}P<+\infty\,\,\,\,\Leftrightarrow\,\,\,\,\sum^{+\infty}_{n=1}P\left(\left|X_{n}\right|>n\right)<+\infty\,\,\,\,\Leftrightarrow\,\,\,\,P\left(\limsup_{n\to+\infty}\left(\left|X_{n}\right|>n\right)\right)=0.
\]
and
\[
\intop_{\Omega}\left|X_{1}\right|\text{d}P=+\infty\,\,\,\,\Leftrightarrow\,\,\,\,\sum^{+\infty}_{n=1}P\left(\left|X_{n}\right|>n\right)=+\infty\,\,\,\,\Leftrightarrow\,\,\,\,P\left(\limsup_{n\to+\infty}\left(\left|X_{n}\right|>n\right)\right)=1.
\]

Hence, we obtain the implications
\[
P\left(\dfrac{X_{n}}{n}\to0\right)=1\,\,\,\,\Rightarrow\,\,\,\,P\left(\limsup_{n\to+\infty}\left(\left|X_{n}\right|>n\right)\right)=0\,\,\,\,\Leftrightarrow\,\,\,\,\intop_{\Omega}\left|X_{1}\right|\text{d}P<+\infty.
\]

Conversely, suppose $\intop_{\Omega}\left|X_{1}\right|\text{d}P<+\infty.$
Then, from $\refpar{eq:double_ineq_intopom_X_1}$, we know that for
every $\epsilon>0,$
\[
\sum^{+\infty}_{n=0}P\left(\left|X_{1}\right|>\left(n+1\right)\epsilon\right)<+\infty.
\]

Since the random variables $X_{n}$ follow the same law,
\[
\sum^{+\infty}_{n=0}P\left(\left|X_{n}\right|>\left(n+1\right)\epsilon\right)<+\infty.
\]

By the Borel-Cantelli lemma again,
\[
\forall\epsilon>0,\,\,\,\,P\left(\limsup_{n\to+\infty}\left(\left|X_{n}\right|>\left(n+1\right)\epsilon\right)\right)=0.
\]

Since $\mathbb{Q}^{+\ast}$ is countable,
\[
P\left(\bigcup_{\epsilon\in\mathbb{Q}^{+\ast}}\limsup_{n\to+\infty}\left(\left|X_{n}\right|>\left(n+1\right)\epsilon\right)\right)=0.
\]

Therefore,
\[
P\left(\dfrac{X_{n}}{n}\to0\right)\geqslant P\left(\bigcap_{\epsilon\in\mathbb{Q}^{+\ast}}\liminf_{n\to+\infty}\left(\dfrac{\left|X_{n}\right|}{n+1}\leqslant\epsilon\right)\right)=1,
\]
 and consequentely,\boxeq{
\[
P\left(\dfrac{X_{n}}{n}\to0\right)=1.
\]
}

The equivalence is then proved.

\end{solution}

\begin{solution}{}{solexercise10.13}

\textbf{1. $\widehat{X_{1}}$ and $\widehat{X_{2}}$ are independent
and follow the same law than $X$}

For every Borel subsets $A$ and $B$ of $\mathbb{R},$
\[
\widehat{X}^{-1}\left(A\times B\right)=X^{-1}\left(A\right)\times X^{-1}\left(B\right)\in\mathcal{A}\otimes\mathcal{A}
\]
which shows that $\widehat{X}$ is a random variable. Therefore, the
same holds for $X^{s},$ which is a measurable function of $\widehat{X}.$ 

Moreover, for every Borel subsets $A$ and $B$ of $\mathbb{R},$
\[
\widehat{X_{1}}^{-1}\left(A\right)\cap\widehat{X_{2}}^{-1}\left(B\right)=X^{-1}\left(A\right)\times X^{-1}\left(B\right),
\]
which implies, using successively the definition of the product measure,
and the equality $P\left(\Omega\right)=1,$
\begin{align*}
P\otimes P\left(\widehat{X_{1}}^{-1}\left(A\right)\cap\widehat{X_{2}}^{-1}\left(B\right)\right) & =P\left(X^{-1}\left(A\right)\right)\cdot P\left(X^{-1}\left(B\right)\right)\\
 & =P\otimes P\left(X^{-1}\left(A\right)\times\Omega\right)\cdot P\otimes P\left(\Omega\times X^{-1}\left(B\right)\right)\\
 & =P\otimes P\left(\widehat{X}^{-1}_{1}\left(A\right)\right)\cdot P\otimes P\left(\widehat{X}^{-1}_{2}\left(B\right)\right).
\end{align*}
This proves the $P\otimes P-$independence of $\widehat{X}_{1}$ and
$\widehat{X}_{2}.$

Now, taking for $B$ the set $\mathbb{R},$
\[
P\otimes P\left(\widehat{X}^{-1}_{1}\left(A\right)\right)=P\otimes P\left(\widehat{X}^{-1}_{1}\left(A\right)\cap\widehat{X}^{-1}_{2}\left(\mathbb{R}\right)\right)=P\left(X^{-1}\left(A\right)\right)\cdot P\left(X^{-1}\left(\mathbb{R}\right)\right),
\]
it follows that
\[
P\otimes P\left(\widehat{X}^{-1}_{1}\left(A\right)\right)=P\left(X^{-1}\left(A\right)\right),
\]
which shows that the random variables $\widehat{X_{1}}$ and $X$
have the same law---and thus also $\widehat{X}_{2}.$

\textbf{2. a. $X^{s}$ is in $\mathscr{L}^{p},$ if $X\in\mathscr{L}^{p}$}

By the Fubini theorem, we compute
\begin{align*}
\intop_{\Omega\times\Omega}\left|\widehat{X}_{1}\right|^{p}\text{d}P\otimes P & =\intop_{\Omega}\left[\intop_{\Omega}\left|\widehat{X}_{1}\left(\omega,\omega^{\prime}\right)\right|^{p}\text{d}P\left(\omega\right)\right]\text{d}P\left(\omega^{\prime}\right)\\
 & =\intop_{\Omega}\left[\intop_{\Omega}\left|X\left(\omega\right)\right|^{p}\text{d}P\left(\omega\right)\right]\text{d}P\left(\omega^{\prime}\right)\\
 & =\intop_{\Omega}\left|X\left(\omega\right)\right|^{p}\text{d}P\left(\omega\right).
\end{align*}
Thus if $X\in\mathscr{L}^{p},$ then the random variables $\widehat{X}_{1}$
and $\widehat{X}_{2}$ are in $\mathscr{L}^{p},$ and thus also $X^{s}.$

\textbf{b. Expectation and variance of $X^{s}$ when $X\in\mathscr{L}^{2}$}

If $X\in\mathscr{L}^{2},$ then, since $X^{s}=\widehat{X}_{1}-\widehat{X}_{2},$
and since $\widehat{X}_{1}$ and $\widehat{X}_{2}$ follows the same
law, and thus have same moments than $X,$\boxeq{
\[
\mathbb{E}\left(X^{s}\right)=\mathbb{E}\left(\widehat{X}_{1}\right)-\mathbb{E}\left(\widehat{X}_{2}\right)=0.
\]
}Moreover $\widehat{X}_{1}$ and $\widehat{X}_{2}$ are independent,
we get\boxeq{
\[
\sigma^{2}_{X^{s}}=\sigma^{2}_{\widehat{X}_{1}}+\sigma^{2}_{\widehat{X}_{2}}=2\sigma^{2}_{X}.
\]
}

\textbf{3. Independence of $\widehat{X}_{i},$ and $X^{s}_{i}$ for
independent family $\left(X_{i}\right)_{i\in I}$}

It suffices to show it for a \textbf{finite} index set $I.$ For every
Borel subsets $A_{i}$ and $B_{i},$ $i\in I,$ from $\mathbb{R},$
\begin{align*}
P\otimes P\left(\bigcap_{i\in I}\widehat{X}^{-1}_{i}\left(A_{i}\times B_{i}\right)\right) & =P\otimes P\left[\left(\bigcap_{i\in I}X^{-1}_{i}\left(A_{i}\right)\right)\times\left(\bigcap_{j\in I}X^{-1}_{j}\left(B_{j}\right)\right)\right]\\
 & =P\left(\bigcap_{i\in I}X^{-1}_{i}\left(A_{i}\right)\right)\cdot P\left(\bigcap_{j\in I}X^{-1}_{j}\left(B_{j}\right)\right).
\end{align*}
Hence, by independence of the $X_{j},$
\begin{align*}
P\otimes P\left(\bigcap_{i\in I}\widehat{X}^{-1}_{i}\left(A_{i}\times B_{i}\right)\right) & =\prod_{i\in I}P\left(X^{-1}_{i}\left(A_{i}\right)\right)\cdot\prod_{j\in I}P\left(X^{-1}_{j}\left(B_{j}\right)\right)\\
 & =\prod_{i\in I}P\otimes P\left(X^{-1}_{i}\left(A_{i}\right)\times X^{-1}_{i}\left(B_{i}\right)\right)\\
 & =\prod_{i\in I}P\otimes P\left(\widehat{X}^{-1}_{i}\left(A_{i}\times B_{i}\right)\right).
\end{align*}
This proves the independence of the random variables $\widehat{X}_{i}.$
Since the symmetrised variables $X^{s}_{i}$ are measurable functions
of the random variables $\widehat{X}_{i}$ are also $P\otimes P-$independent.

\end{solution}

\chapter{Convergences and Laws of Large Numbers}\label{chap:PartIIChap11}

\begin{objective}{}{}

Chapter \ref{chap:PartIIChap11} is devoted to understanding how random
variables behave when we observe them repeatedly. It introduces several
important types of convergences and explains how they lead to the
laws of large numbers.
\begin{itemize}
\item Section \ref{sec:Convergence-in-Probability} introduces two fundamental
notions: \textbf{convergence in probability} and \textbf{almost sure
convergence}. It explains how these concepts differ, provides simple
criteria that guarantee $P-$almost sure convergence, and compares
the various modes of convergence used in probability theory.
\item Section \ref{sec:Convergence--and} defines the idea of \textbf{equi-integrability},
a key tool for controlling families of random variables in $\mathbb{R}^{d},$
that extends the dominated convergence theorem. After presenting a
basic sufficient condition for equi-integrability, the section introduced
the concept of \textbf{equi-continuity} and shows how it completely
characterizes equi-integrability. The notion of \textbf{$\mathscr{L}^{p}-$convergence}
is then introduced, and the section ends with an important theorem
describing how convergence in probability relates to $\mathscr{L}^{p}-$convergence. 
\item Section \ref{sec:Series-of-Independent} extends Chebychev inequality
through \textbf{Kolmogorov inequality}, a powerful estimate for partial
sums of independent random variables. Using this inequality, the section
presents improved conditions ensuring that a series of independent
random variables converges $P-$almost surely.
\item Section \ref{sec:Laws-of-Large} addresses the \textbf{laws of large
numbers}, which describe how averages of random variables behaves
when the number of obervations grows. The section begins with the
\textbf{Cesaro and Kronecker lemmas}, which serve as foundational
tools. It then introduces the \textbf{weak law of large numbers},
followed by Bernouilli and Khintchine theorems--- two classical form
of the weak law. Next comes the \textbf{strong law of large numbers},
including the Kolmogorov-Khintchine theorem, which applies when the
variables are only integrable. The chapter concludes with the notions
of \textbf{samples} and empirical distribution function, leading to
the fundamental theorem of statistics: the \textbf{Glivenko-Cantelli}
\textbf{theorem} which states that the statistical mean of the empirical
distribution converges $P-$almost surely to the true distribution
expectation. 
\end{itemize}
\end{objective}

In the first part of this Chapter, we study the concepts of convergence
in probability, almost sure convergence and convergence in $\text{L}^{p},$
as well as the relationships between the different modes of convergence.
To this end, we introduce the concept of equi-integrability. The second
part deals with the weak and strong laws of large numbers.

\section{Convergence in Probability and Almost Sure}\label{sec:Convergence-in-Probability}

In this Section, all random variables are defined on the same probabilized
space $\left(\Omega,\mathcal{A},P\right)$ and take values in $\mathbb{R}^{d},\,d\geqslant1$
or in $\overline{\mathbb{R}}.$

We use the notation $\left|.\right|$ interchangeably to denote the
absolute value in $\mathbb{R}$---possible extended to $\overline{\mathbb{R}}$---or
a norm on $\mathbb{R}^{d}.$

\begin{definition}{Almost Sure Convergence. Convergence in Probability. Convergence Almost Certain}{}

a. (i) A sequence $\left(X_{n}\right)_{n\in\mathbb{N}}$ of random
variables \textbf{converges almost surely\mindex{convergence almost sure ! to a random variable}
or $P-$almost surely\mindex{P-almost sure convergence@$P-$almost sure convergence ! to a random variable}}
to a random variable $X$ if there exists a set $C\in\mathcal{A}$
of probability 1 on which the sequence converges pointwise (or simply)
to $X.$ We denote this by
\[
X_{n}\stackrel[n\to+\infty]{P-\text{a.s.}}{\longrightarrow}X.
\]

(ii) A sequence $\left(X_{n}\right)_{n\in\mathbb{N}}$ of random variables
is said to \textbf{\mindex{P-almost sure convergence@$P-$almost sure convergence ! sequence}\mindex{convergence almost sure ! sequence}converge
$\left(P-\right)$almost surely} if there exists a random variable
$X$ such that this\textbf{ sequence converges $\left(P-\right)$almost
surely} to $X.$

b. (i) A sequence $\left(X_{n}\right)_{n\in\mathbb{N}}$ of random
variables \textbf{converges in probability\mindex{convergence in probability ! to a random variable}}
to a random variable $X$ if, for every $\epsilon>0,$ the sequence
with general term $P\left(\left|X_{n}-X\right|>\epsilon\right)$ converges
to 0 when $n$ tends to infinity. We denote this by
\[
X_{n}\stackrel[n\to+\infty]{P}{\longrightarrow}X.
\]

(ii) A sequence $\left(X_{n}\right)_{n\in\mathbb{N}}$ of random variables
is said to \textbf{converge} \textbf{in probability\index{convergence in probability}}
if there exists a random variable $X$ such that the\textbf{ sequence
converges} \textbf{in probability} to $X.$

\end{definition}

\begin{denotation}{}{}

We denote $\left(X_{n}\longrightarrow\right)$---respectively $\left(X_{n}\longrightarrow X\right)$---the
set of $\omega$ for which the sequence $\left(X_{n}\left(\omega\right)\right)_{n\in\mathbb{N}}$
converges---respectively converges to $X\left(\omega\right).$

\end{denotation}

\begin{remark}{}{}

1. If a sequence $\left(X_{n}\right)_{n\in\mathbb{N}}$ of random
variables converges $\left(P-\right)$almost surely---respectively
in probability---the limit $X$ is $P-$almost surely unique. 

This is clear for almost sure convergence. 

For convergence in probability, let $X$ and $X'$ be two limits in
probability. For every $\mathbb{\epsilon}>0$ and every $n\in\mathbb{N},$
by the triangle inequality,
\[
\left(\left|X-X^{\prime}\right|>\epsilon\right)\subset\left(\left|X-X_{n}\right|>\dfrac{\epsilon}{2}\right)\cup\left(\left|X_{n}-X^{\prime}\right|>\dfrac{\epsilon}{2}\right),
\]
and therefore,
\[
P\left(\left|X-X^{\prime}\right|>\epsilon\right)\leqslant P\left(\left|X-X_{n}\right|>\dfrac{\epsilon}{2}\right)+P\left(\left|X_{n}-X^{\prime}\right|>\dfrac{\epsilon}{2}\right).
\]

Taking the limit,
\[
\forall\epsilon>0,\,\,\,\,P\left(\left|X-X^{\prime}\right|>\epsilon\right)=0.
\]

The result follows by noting that
\[
\left(X\neq X^{\prime}\right)=\bigcup_{n\in\mathbb{N}^{\ast}}\left(\left|X-X^{\prime}\right|>\dfrac{1}{n}\right).
\]

2. The convergence in probability $X_{n}\stackrel[n\to+\infty]{P}{\longrightarrow}X$
can be written in the quantified form
\[
\forall\epsilon>0,\,\forall\delta>0,\,\exists N\left(\epsilon,\delta\right):\,\,\,\,n\geqslant N\left(\epsilon,\delta\right)\,\,\,\,\Rightarrow\,\,\,\,P\left(\left|X_{n}-X\right|>\epsilon\right)\leqslant\delta.
\]
This is equivalent to the assertion
\[
\forall\epsilon>0,\,\exists N\left(\epsilon,\delta\right):\,\,\,\,n\geqslant N\left(\epsilon\right)\,\,\,\,\Rightarrow\,\,\,\,P\left(\left|X_{n}-X\right|>\epsilon\right)\leqslant\epsilon.
\]

It is clear that the first assertion implies the second. 

Conversely, assume the second assertion holds and let $\epsilon>0$
and $\delta>0.$ 
\begin{itemize}
\item If $\delta\geqslant\epsilon,$ we take for $N\left(\epsilon,\delta\right)=N\left(\epsilon\right).$ 
\item If $\delta<\epsilon,$ we take for $N\left(\epsilon,\delta\right)=N\left(\delta\right).$ 
\end{itemize}
Then, for every $n\geqslant N\left(\delta\right),$
\[
P\left(\left|X_{n}-X\right|>\delta\right)\leqslant\delta,
\]
 and the result follows from the inclusion
\[
\left(\left|X_{n}-X\right|>\epsilon\right)\subset\left(\left|X_{n}-X\right|>\delta\right).
\]

3. If the random variables take values in $\mathbb{R}^{d}$ with $d\geqslant2,$
the choice of the norm is irrelevant, since all norms are equivalent
in finite dimension. 

Moreover, $X_{n}\stackrel[n\to+\infty]{P}{\longrightarrow}X$ holds
if and only if, for every $j\in\left\llbracket 1,d\right\rrbracket ,$
\[
X^{j}_{n}\stackrel[n\to+\infty]{P}{\longrightarrow}X^{j},
\]
where $X^{j}_{n}$ denotes the $j-$th component of $X_{n}.$

The necessity is straightforward. The sufficiency follows from the
inequalities---using the max norm---
\[
P\left(\max_{1\leqslant j\leqslant d}\left|X^{j}_{n}-X^{j}\right|>\epsilon\right)\leqslant P\left(\bigcup^{d}_{i=1}\left(\left|X^{j}_{n}-X^{j}\right|>\epsilon\right)\right)\leqslant\sum^{d}_{j=1}P\left(\left|X^{j}_{n}-X^{j}\right|>\epsilon\right).
\]
\end{remark}

The sufficient conditions for $P-$almost sure convergence presented
below are commonly used.

\begin{theorem}{A First Sufficient Condition for Almost Sure Convergence of A Sequence of Random Variables}{sc_as_cv_seq_rv}

If there exists a series with nonnegative terms, of general term $\epsilon_{n}$
that converges and satisfies
\[
\sum^{+\infty}_{n=0}P\left(\left|X_{n+1}-X_{n}\right|>\epsilon_{n}\right)<+\infty,
\]
then the sequence of random variables $\left(X_{n}\right)_{n\in\mathbb{N}}$
converges almost surely.

\end{theorem}

\begin{proof}{}{}

By the Borel-Cantelli lemma,
\[
P\left(\limsup_{n\to+\infty}\left(\left|X_{n+1}-X_{n}\right|>\epsilon_{n}\right)\right)=0.
\]

The set
\[
C=\liminf_{n\to+\infty}\left(\left|X_{n+1}-X_{n}\right|\leqslant\epsilon_{n}\right)
\]
therefore has probability 1. For every $\omega\in C,$ the series
with general term $\left|X_{n+1}\left(\omega\right)-X_{n}\left(\omega\right)\right|$
converges. Consequently, the sequence with general term $X_{n}\left(\omega\right)$
also converges.

\end{proof}

\begin{remark}{}{}

This theorem will be used in particular to compare convergence in
probability and almost sure convergence---Theorem $\ref{th:comp_var_mode_cv}.$

\end{remark}

\begin{theorem}{A Second Sufficient Condition of Almost Certain Convergence of A Sequence of Random Variables}{second_th_sc_cv_seq}

Let $X$ be a random variable such that, for every $\epsilon>0,$
\[
\sum^{+\infty}_{n=0}P\left(\left|X_{n}-X\right|>\epsilon\right)<+\infty.
\]

Then the sequence of random variables $\left(X_{n}\right)_{n\in\mathbb{N}}$
converges almost surely to $X.$

\end{theorem}

\begin{proof}{}{}

By the Borel-Cantelli lemma, for every $\epsilon\in\mathbb{Q}^{\ast+},$
\[
P\left(\limsup_{n\to+\infty}\left(\left|X_{n}-X\right|>\epsilon\right)\right)=0.
\]

Since $\mathbb{Q}^{\ast+}$ is countable,
\[
P\left(\bigcup_{\epsilon\in\mathbb{Q}^{*+}}\limsup_{n\to+\infty}\left(\left|X_{n}-X\right|>\epsilon\right)\right)=0.
\]

Thus, the set $C=\bigcap_{\epsilon\in\mathbb{Q}^{\ast+}}\left(\liminf_{n\to+\infty}\left(\left|X_{n}-X\right|\leqslant\epsilon\right)\right)$
has probability 1.

The set $C$ corresponds to the $\omega$ for which the sequence with
general term $X_{n}\left(\omega\right)$ converges to $X\left(\omega\right).$

\end{proof}

\begin{remark}{}{}

The previous theorem provides a sufficient but not a necessarry condition
for $P-$almost sure convergence. Indeed, consider the probabilized
space $\left(\left[0,1\right],\mathscr{B}_{\left[0,1\right]},P\right),$
where $P$ is the Lebesgue measure on $\left[0,1\right],$ and the
random variables $X_{n}=\boldsymbol{1}_{\left[0,1/n\right[}.$ For
every $\epsilon>0,$
\[
P\left(\left|X_{n}\right|>\epsilon\right)=\dfrac{1}{n},
\]
and thus
\[
\sum^{+\infty}_{j=1}P\left(\left|X_{n}\right|>\epsilon\right)=+\infty,
\]
whereas the sequence with general term $X_{n}$ converges $P-$almost
surely to $0.$ A partial converse will be studied in Exercise $\ref{exo:exercise11.2}.$

We now compare the different modes of convergence.

\end{remark}

\begin{theorem}{Comparison of the Various Modes of Convergence}{comp_var_mode_cv}

Let $\left(X_{n}\right)_{n\in\mathbb{N}}$ be a sequence of random
variables.

a. If the sequence $\left(X_{n}\right)_{n\in\mathbb{N}}$ converges
almost surely, then it also converges in probability, and the limits
are $P-$almost surely equal.

b. If the sequence $\left(X_{n}\right)_{n\in\mathbb{N}}$ converges
in probability to a random variable $X,$ then there exists a subsequence
$\left(X_{n_{i}}\right)_{i\in\mathbb{N}}$ that converges almost surely
to $X.$

c. The sequence $\left(X_{n}\right)_{n\in\mathbb{N}}$ converges in
probability to $X,$ if and only if the sequence is Cauchy with respect
to convergence in probability, that is, for $\epsilon>0,$ the double
sequence with general term $P\left(\left|X_{n}-X_{m}\right|>\epsilon\right)$
converges to 0.

\end{theorem}

\begin{proof}{}{}

a. Suppose that $\left(X_{n}\right)_{n\in\mathbb{N}}$ converges almost
surely. Let $X$ be the almost sure limit of the sequence $\left(X_{n}\right)_{n\in\mathbb{N}}.$
For every $\epsilon>0,$
\[
\left(X_{n}\rightarrow X\right)\subset\liminf_{n\to+\infty}\left(\left|X_{n}-X\right|\leqslant\epsilon\right).
\]

Taking complements,
\begin{align*}
0\leqslant\limsup_{n\to+\infty}P\left(\left|X_{n}-X\right|>\epsilon\right) & \leqslant P\left(\limsup_{n\to+\infty}\left(\left|X_{n}-X\right|>\epsilon\right)\right)\\
 & \leqslant P\left(\left(X_{n}\rightarrow X\right)^{c}\right)=0
\end{align*}

Thus, $\left(X_{n}\right)_{n\in\mathbb{N}}$ converges in probability
to $X,$ and the limits are $P-$almost surely identical.

b. \textbf{We prove this statement in two steps}

\textbf{Step 1: Suppose that $\left(X_{n}\right)_{n\in\mathbb{N}}$
converges in probability to $X.$ We show that this sequence is Cauchy
for the convergence in probability.} 

Indeed, for every $\epsilon>0,$ and every $n,m\in\mathbb{N},$
\[
\left(\left|X_{n}-X_{m}\right|>\epsilon\right)\subset\left(\left|X-X_{n}\right|>\dfrac{\epsilon}{2}\right)\cup\left(\left|X_{m}-X\right|>\dfrac{\epsilon}{2}\right).
\]

This also shows the necessary part of the third assertion stated in
(c).

\textbf{Step 2: Now suppose that $\left(X_{n}\right)_{n\in\mathbb{N}}$
is Cauchy for the convergence in probability. We show there exists
a subsequence $\left(X_{n_{j}}\right)_{j\in\mathbb{N}}$ that converges
almost surely. }

Construct the sequence of integers $\left(n_{j}\right)_{j\in\mathbb{N}}$
by setting $n_{0}=1,$ and for every $j\in\mathbb{N}^{\ast},$
\[
n_{j}=\inf\left\{ n>n_{j-1}:\,\forall p,q\geqslant n,\,\,\,\,P\left(\left|X_{p}-X_{q}\right|>\dfrac{1}{2^{j}}\right)<\dfrac{1}{3^{j}}\right\} .
\]

Since the sequence $\left(X_{n}\right)_{n\in\mathbb{N}}$ is Cauchy
for the convergence in probability, the constructed sequence tends
to $+\infty$ by increasing. Moreover,
\[
\sum^{+\infty}_{n=0}P\left(\left|X_{n_{j+1}}-X_{n_{j}}\right|>\dfrac{1}{2^{j}}\right)<+\infty.
\]

By Theorem $\ref{th:sc_as_cv_seq_rv},$ the subsequence $\left(X_{n_{j}}\right)_{j\in\mathbb{N}}$
converges almost surely.

\textbf{Summary}: if the sequence $\left(X_{n}\right)_{n\in\mathbb{N}}$
converges in probability, then it is Cauchy, and one can extract a
$P-$almost surely convergent subsequence.

(c) It remains to prove the sufficient condition: \textbf{if $\left(X_{n}\right)_{n\in\mathbb{N}}$
is Cauchy for the convergence in probability, then it converges in
probability}. Let $X$ be the almost sure limit of the subsequence
$\left(X_{n_{j}}\right)_{j\in\mathbb{N}}.$ By the assertion (a),
this subsequence converges in probability to $X.$ 

Moreover, for every integers $n$ and $j,$
\[
P\left(\left|X_{n}-X\right|>\epsilon\right)\leqslant P\left(\left|X_{n}-X_{n_{j}}\right|>\dfrac{\epsilon}{2}\right)+P\left(\left|X_{n_{j}}-X\right|>\dfrac{\epsilon}{2}\right).
\]

Then the fact that the sequence $\left(X_{n}\right)_{n\in\mathbb{N}}$
converges to $X$ in probability results from the fact that this sequence
is Cauchy and that the sequence $\left(X_{n_{j}}\right)_{j\in\mathbb{N}}$
converges in probability to $X.$

\end{proof}

\begin{counterexample}{}{}

We now give an example of sequence that converges in probability but
not almost surely.

Let $\left(X_{n}\right)_{n\in\mathbb{N}^{\ast}}$ be a sequence of
independent random variables with values $0$ or $1$ such that
\[
\forall n\in\mathbb{N}^{\ast},\,\,\,\,P\left(X_{n}=1\right)=\dfrac{1}{n}\,\,\,\,\text{and}\,\,\,\,P\left(X_{n}=0\right)=1-\dfrac{1}{n}.
\]

The sequence $\left(X_{n}\right)_{n\in\mathbb{N}^{\ast}}$ converges
in probability to 0, since for every $\epsilon>0,$
\[
P\left(\left|X_{n}>\epsilon\right|\right)=P\left(X_{n}=1\right)=\dfrac{1}{n}.
\]

Nonetheless, we will show that this sequence does not converge $P-$almost
surely. Indeed, if it did converge $P-$almost surely, it would necessarily
converge to 0. However, this is not the case, as the following argument
shows. The events $\left(X_{n}=1\right)$ are independent and satisfy
\[
\sum^{+\infty}_{n=1}P\left(X_{n}=1\right)=+\infty.
\]

By the Borel-Cantelli lemma,
\[
P\left(\limsup_{n\to+\infty}\left(X_{n}=1\right)\right)=1,
\]
which means that the sequence takes the value 1 infinitely often,
$P-$almost surely. Therefore, $\left(X_{n}\right)_{n\in\mathbb{N}^{\ast}}$
cannot converge $P-$almost surely to 0.

\end{counterexample}

\begin{remark}{}{}

When $f$ is a continous function from $\mathbb{R}^{d}$ to $\mathbb{R}^{k},$
and when the sequence $\left(X_{n}\right)_{n\in\mathbb{N}}$ converges
almost surely to $X,$ it is straightforward that the sequence $\left(f\left(X_{n}\right)\right)_{n\in\mathbb{N}}$
converges almost surely to $f\left(X\right).$ We will prove a similar
result for the convergence in probability.

\end{remark}

\begin{proposition}{Convergence in Probability of the Image of Random Variables by A Continuous Function}{}

Let $f$ be a continuous function from $\mathbb{R}^{d}$ to $\mathbb{R}^{k}.$
If the sequence $\left(X_{n}\right)_{n\in\mathbb{N}}$ converges in
probability to $X,$ then the sequence $\left(f\left(X_{n}\right)\right)_{n\in\mathbb{N}}$
converges in probability to $f\left(X\right).$

\end{proposition}

\begin{proof}{}{}

Fix an arbitrary $\delta>0.$ Choose $a>0$ such that
\[
P_{X}\left(B\left(0,a\right)^{c}\right)\leqslant\dfrac{\delta}{2},
\]
where $B\left(0,a\right)$ denotes the open ball centered at the origin
with radius $a.$ This choice is possible since
\[
\lim_{n\to+\infty}P_{X}\left(B\left(0,n\right)^{c}\right)=P_{X}\left(\emptyset\right)=0.
\]

Since $f$ is countinuous on $\mathbb{R}^{d},$ it is uniformly continuous
on the closed ball $B_{f}\left(0,2a\right),$ thus,
\[
\forall\epsilon>0,\,\exists\eta\left(\epsilon\right)>0:\,\forall x,y\in B_{f}\left(0,2a\right),\,\,\,\,\left|x-y\right|\leqslant\eta\left(\epsilon\right)\,\,\,\,\Rightarrow\,\,\,\,\left|f\left(x\right)-f\left(y\right)\right|\leqslant\epsilon.
\]

In particular, by the triangle inequality, for any fixed $\epsilon>0,$
\[
\left|x\right|\leqslant a\,\,\,\,\text{and}\,\,\,\,\left|x-y\right|\ensuremath{\leqslant}\ensuremath{\eta}\left(\epsilon\right)\ensuremath{\wedge}a\,\,\,\,\Rightarrow\,\,\,\,\left|f\left(x\right)-f\left(y\right)\right|\leqslant\epsilon.
\]

Thus, taking the contrapositive of this implication, we obtain
\[
\left|f\left(x\right)-f\left(y\right)\right|>\epsilon\,\,\,\,\Rightarrow\,\,\,\,\left|x\right|>a\,\,\,\,\text{or\,\,\,\,}\left|x-y\right|>\ensuremath{\eta}\left(\epsilon\right)\ensuremath{\wedge}a.
\]

Therefore, for every $n\in\mathbb{N},$ we have the inclusion
\[
\left(\left|f\left(X_{n}\right)-f\left(X\right)\right|>\epsilon\right)\subset\left(\left|X\right|>a\right)\cup\left(\left|X_{n}-X\right|>\ensuremath{\eta}\left(\epsilon\right)\ensuremath{\wedge}a\right).
\]

Since the sequence $\left(X_{n}\right)_{n\in\mathbb{N}}$ converges
in probability to $X,$ there exists an integer $N$ such that, for
every $n\geqslant N,$
\[
P\left(\left|X_{n}-X\right|>\ensuremath{\eta}\left(\epsilon\right)\ensuremath{\wedge}a\right)\leqslant\dfrac{\delta}{2}.
\]

Hence, for every $n\geqslant N,$
\[
P\left(\left|f\left(X_{n}\right)-f\left(X\right)\right|>\epsilon\right)\subset P\left(\left|X\right|>a\right)+P\left(\left|X_{n}-X\right|>\ensuremath{\eta}\left(\epsilon\right)\ensuremath{\wedge}a\right)\leqslant\delta.
\]

Thus, for every $\epsilon>0,$ 
\[
\lim_{n\to+\infty}P\left(\left|f\left(X_{n}\right)-f\left(X\right)\right|>\epsilon\right)=0,
\]
which shows that $\left(f\left(X_{n}\right)\right)_{n\in\mathbb{N}}$
converges in probability to $f\left(X\right).$

\end{proof}

\begin{example}{}{}

Let $\left(X_{n}\right)_{n\in\mathbb{N}}$ and $\left(Y_{n}\right)_{n\in\mathbb{N}}$
be two sequences of random variables taking values in $\mathbb{R}^{d}$
and converging in probability to $X$ and $Y,$ respectively.

Then
\[
\left\langle X_{n},Y_{n}\right\rangle \stackrel[n\to+\infty]{P}{\longrightarrow}\left\langle X,Y\right\rangle .
\]

Indeed, the sequence $\left(\left(X_{n},Y_{n}\right)\right)_{n\in\mathbb{N}},$
taking values in $\mathbb{R}^{2d},$ converges in probability to $\left(X,Y\right),$
and the scalar product is a continuous function. The conclusion then
follows from the previous proposition. 

By the same argument, if $d=1,$ then
\[
\max\left(X_{n},Y_{n}\right)\stackrel[n\to+\infty]{P}{\longrightarrow}\max\left(X,Y\right),
\]
since the function $\left(x,y\right)\mapsto\max\left(x,y\right)$
is continuous on $\mathbb{R}^{2}.$

\end{example}

\begin{remark}{}{}

From what we have presented, it appears that we may modify random
variables on a set of probability zero without affecting either the
definitions, or the results. This suggests the construction of a ``theory''
of convergence on \textbf{classes of random variables\mindex{random variable ! class}}\index{classes of random variables}.

More precisely, let $X$ be a function from $D_{X}\in\mathcal{A}$
taking values in $\mathbb{R}^{d}$ or $\overline{\mathbb{R}}.$ We
say that $X$ is a \textbf{random variable defined $P-$almost surely\mindex{random variable ! P-almost surely @ $P-$almost surely ! defined}}
if:
\begin{itemize}
\item $X$ is measurable with respect to the trace \salg $D_{X}\cap\mathcal{A}$
from $D_{X}$ to $\mathcal{A}$ 
\item And, $P\left(D_{X}\right)=1.$
\end{itemize}
We sat that the random variables $X$ and $Y$ defined $P-$almost
surely are \textbf{$P-$almost equal\mindex{random variables ! P-almost equal @ $P-$almost equal}}
if
\[
P\left(\left\{ \omega\in D_{X}\cap D_{Y}:X\left(\omega\right)=Y\left(\omega\right)\right\} \right)=1.
\]
Depending on whether the space of values taken by the random variables
is $\mathbb{R}^{d}$ or $\overline{\mathbb{R}},$ we define the set
$G$ of random variables defined $P-$almost surely---defined $P-$almost
surely and defined $P-$almost surely finitely, respectively. This
is a vector space. The subset $K$ of random variables that are $P-$almost
surely equal to 0 is a vector subspace. The relation ``equal $P-$almost
surely'' is an equivalence relation. 

The quotient set $L^{0}\left(\Omega,\mathcal{A},P\right)$ of $G$
by this relation is the quotient vector space from $G$ by $K.$ It
is called the \textbf{space of classes of random variables defined
$P-$almost surely\mindex{random variable ! P-almost surely @ $P-$almost surely ! space of classes}}\index{space of classes of random variables defined $P-$almost surely}---and
in the case of $\overline{\mathbb{R}},$ \textbf{$P-$almost surely
finite}\index{set of classes of random variables defined $P-$almost surely finite}.
All statements previously made about convergence then extend naturally
to $L^{0}\left(\Omega,\mathcal{A},P\right).$

On $L^{0}\left(\Omega,\mathcal{A},P\right),$ one can define a metric
that makes it complete and that has the property that the convergence
of a sequence in the sense of this metric is equivalent to the convergence
in probability---see Exercise $\ref{exo:exercise11.1}.$

\end{remark}

\section{Convergence $\text{L}^{p}$ and Equi-Integrability}\label{sec:Convergence--and}

\textbf{In this section, all random variables are defined on the same
probabilized space $\left(\Omega,\mathcal{A},P\right)$ and take values
in $\mathbb{R}^{d},\,d\geqslant1,$ or in $\overline{\mathbb{R}}.$}

If $X$ is an integrable random variable, the dominated convergence
theorem implies that
\[
\lim_{a\to+\infty}\intop_{\left(\left|X\right|>a\right)}\left|X\right|\text{d}P=0.
\]

The concept of \textbf{equi-integrability} extends this property to
an arbitrary family of random variables in an uniform way.

\begin{definition}{Equi-integrability of a Family of Random Variables}{}

A family of random variables $\left(X_{i}\right)_{i\in I}$ where
$I$ is an arbitrary index set, is said to be \textbf{\index{equi-integrable}equi-integrable}
if
\[
\lim_{a\to+\infty}\sup_{i\in I}\intop_{\left(\left|X_{i}\right|>a\right)}\left|X_{i}\right|\text{d}P=0.
\]

\end{definition}

\textbf{We now give a sufficient condition of equi-integrability.}

\begin{proposition}{Sufficient Condition of Equi-Integrability}{}

If the family $\left(X_{i}\right)_{i\in I}$ is uniformly dominated
by a nonnegative and integrable random variable $X,$ that is
\[
\forall i\in I,\,\,\,\,\left|X_{i}\right|\leqslant X\,\,\,\,P-\text{almost surely}
\]
then the family is equi-integrable. 

In particular, any finite family of integrable random variable is
equi-integrable.

\end{proposition}

\begin{proof}{}{}

We have
\[
\forall i\in I,\,\forall a>0,\,\,\,\,\,\left(\left|X_{i}\right|>a\right)\subset\left(\left|X\right|>a\right).
\]

Hence,
\[
\forall a>0,\,\forall i\in I,\,\,\,\,\intop_{\left(\left|X_{i}\right|>a\right)}\left|X_{i}\right|\text{d}P\leqslant\intop_{\left(\left|X\right|>a\right)}\left|X\right|\text{d}P.
\]

Then
\[
\forall a>0,\,\,\,\,\sup_{i\in I}\intop_{\left(\left|X_{i}\right|>a\right)}\left|X_{i}\right|\text{d}P\leqslant\intop_{\left(\left|X\right|>a\right)}\left|X\right|\text{d}P.
\]

Since the right-hand side tends to 0 when $a$ tends to infinity,
the family $\left(X_{i}\right)_{i\in I}$ is equi-integrable.

If $I$ is finite, the nonnegative random variable $X=\max_{i\in I}\left|X_{i}\right|$
is integrable: the conclusion follows from the first part.

\end{proof}

We now provide a \textbf{necessary and sufficient condition of equi-integrability}.
First, we need to define the notion of \textbf{equi-continuity}\index{equi-continuity}.

\begin{definition}{Equi-Continuity  of a Family of Random Variables}{}

Let $I$ be an arbitrary set.

A family of random variables $\left(X_{i}\right)_{i\in I}$ is said
to be\textbf{ \index{equi-continuous}equi-continuous} if
\[
\forall\epsilon>0,\,\,\,\,\exists\eta>0:\,\,\,\,P\left(A\right)\leqslant\eta\Rightarrow\sup_{i\in I}\intop_{A}\left|X_{i}\right|\text{d}P\leqslant\epsilon.
\]

\end{definition}

\begin{remark}{}{}

This concept is precisely the usual equi-continuity at a point for
a family of functions. Indeed, it is classical to define the \textbf{metric
algebra} $\mathcal{A},$ that is, the set $\mathcal{A}$ equipped
with the gap---or pseudo-distance---defined by the function
\[
\phi:\left(A,B\right)\mapsto P\left(A\triangle B\right),
\]
where $\Delta$ denotes the symmetric difference\footnotemark.

In this setting, $P\left(A\right)=\phi\left(A,\emptyset\right)$ represents
the ``distance'' from $A$ to the empty set. Thus, equi-continuity
in $\emptyset$ of the family of functions $A\mapsto\intop_{A}\left|X_{i}\right|\text{d}P$
corresponds exactly to the definition given above.

\end{remark}

\footnotetext{$A\Delta B=\left(A\backslash B\right)\cup\left(B\backslash A\right).$

}

\begin{proposition}{Necessary and Sufficient Condition of Equi-Integrability}{}

A family of random variables $\left(X_{i}\right)_{i\in I}$ is equi-integrable
if and only if it is equi-continuous and bounded in $\text{L}^{1},$
that is, such that 
\[
\sup_{i\in I}\intop_{\Omega}\left|X_{i}\right|\text{d}P<+\infty.
\]

\end{proposition}

\begin{proof}{}{}
\begin{itemize}
\item \textbf{Necessary condition}\\
Suppose that the family $\left(X_{i}\right)_{i\in I}$ is equi-integrable.
\\
For every $A\in\mathcal{A},$ every $a>0,$ and every $i\in I,$
\begin{align*}
\intop_{A}\left|X_{i}\right|\text{d}P & \leqslant\intop_{A\cap\left(\left|X_{i}\right|\leqslant a\right)}\left|X_{i}\right|\text{d}P+\intop_{A\cap\left(\left|X_{i}\right|>a\right)}\left|X_{i}\right|\text{d}P\\
 & \leqslant aP\left(A\right)+\sup_{i\in I}\intop_{\left(\left|X_{i}\right|>a\right)}\left|X_{i}\right|\text{d}P.
\end{align*}
\\
Hence, for every $A\in\mathcal{A}$ and every $a>0,$
\[
\sup_{i\in I}\intop_{A}\left|X_{i}\right|\text{d}P\leqslant aP\left(A\right)+\sup_{i\in I}\intop_{\left(\left|X_{i}\right|>a\right)}\left|X_{i}\right|\text{d}P.
\]
By taking $A=\Omega$ in the previous inequaliy, we obtain that the
family is bounded in $\text{L}^{1}.$ Moreover, let $\epsilon>0.$
Choose $a>0$ such that 
\[
\sup_{i\in I}\intop_{\left(\left|X_{i}\right|>a\right)}\left|X_{i}\right|\text{d}P\leqslant\dfrac{\epsilon}{2}
\]
and then $\eta=\dfrac{\epsilon}{2a}.$ \\
Whenever $P\left(A\right)\leqslant\eta,$ we get 
\[
\sup_{i\in I}\intop_{A}\left|X_{i}\right|\text{d}P\leqslant\epsilon
\]
and therefore the family is equi-continuous.
\item \textbf{Sufficient condition}\\
Suppose that the family $\left(X_{i}\right)_{i\in I}$ is equi-continuous
and bounded in $\text{L}^{1}.$ By Markov inequality---Proposition
$\ref{pr:markov_ineq}$---we have, for every $a>0$ and every $i\in I,$
\[
P\left(\left|X_{i}\right|>a\right)\leqslant\dfrac{1}{a}\intop\left|X_{i}\right|\text{d}P\leqslant\dfrac{1}{a}\sup_{i\in I}\intop\left|X_{i}\right|\text{d}P.
\]
Since the family $\left(X_{i}\right)_{i\in I}$ is bounded in $\text{L}^{1},$
it follows that
\begin{equation}
\lim_{a\to+\infty}\sup_{i\in I}P\left(\left|X_{i}\right|>a\right)=0.\label{eq:limsup_p_absXi}
\end{equation}
\\
Consider an arbitrary $\epsilon>0.$ Since the family is equi-continuous,
choose $\eta>0$ such that
\begin{equation}
P\left(A\right)\leqslant\eta\Rightarrow\sup_{i\in I}\intop_{A}\left|X_{i}\right|\text{d}P\leqslant\epsilon.\label{eq:equi-cont}
\end{equation}
\\
Choose $M>0$ such that
\[
\forall a\geqslant M,\,\,\,\,\sup_{i\in I}P\left(\left|X_{i}\right|>a\right)\leqslant\eta,
\]
which is possible by $\refpar{eq:limsup_p_absXi}.$ \\
Then, by $\refpar{eq:equi-cont}$
\[
\forall a\geqslant M,\,\forall i\in I,\,\,\,\,\intop_{\left(\left|X_{i}\right|>a\right)}\left|X_{i}\right|\text{d}P\leqslant\epsilon.
\]
Thus
\[
\forall a\geqslant M,\,\,\,\,\sup_{i\in I}\intop_{\left(\left|X_{i}\right|>a\right)}\left|X_{i}\right|\text{d}P\leqslant\epsilon.
\]
This proves the equi-integrability of the family.
\end{itemize}
\end{proof}

\begin{definition}{Convergence in $\mathscr{L}^p$ of a Sequence of Random Variables}{}

Let $p\geqslant1.$

A sequence $\left(X_{n}\right)_{n\in\mathbb{N}}$ of random variables,
admitting a moment of order $p$ \textbf{converges in $\mathscr{L}^{p}$\mindex{convergence in Lp@convergence in $\mathcal{L}^p$!to a random variable}
to a random variable $X$} if $X\in\mathscr{L}^{p}\left(\Omega,\mathcal{A},P\right)$
and if
\[
\lim_{n\to+\infty}\mathbb{E}\left(\left|X_{n}-X\right|^{p}\right)=0.
\]

We write
\[
X_{n}\stackrel[n\to+\infty]{\mathscr{L}^{p}}{\longrightarrow}X.
\]
 The sequence $\left(X_{n}\right)_{n\in\mathbb{N}}$ of random variables
\textbf{converges in $\mathscr{L}^{p}$\mindex{convergence in Lp@convergence in $\mathcal{L}^p$}}
if there exists a random variable $X\in\mathscr{L}^{p}\left(\Omega,\mathcal{A},P\right)$
such that this sequence converges in $\mathscr{L}^{p}$ to $X.$

\end{definition}

\begin{remark}{}{}
\begin{itemize}
\item If $p=1,$ we say that the sequence \textbf{\index{convergence in mean}converges
in mean}.
\item If $p=2,$ we say that the sequence \textbf{converges\index{convergence in quadratic mean}
in quadratic mean}.
\item If $p\geqslant1,$ by the Minkowski inequality---Proposition $\ref{pr:minkowski_ineq_9}$---the
map 
\[
X\mapsto\left[\mathbb{E}\left(\left|X\right|^{p}\right)\right]^{1/p}
\]
is a \textbf{semi-norm} on $\mathscr{L}^{p}\left(\Omega,\mathcal{A},P\right):$
the notions of $\mathscr{L}^{p}-$convergence are therefore notions
of convergence with respect to this semi-norm.\\
In particular, if a sequence $\left(X_{n}\right)_{n\in\mathbb{N}}$
of random variables \textbf{converges in $\mathscr{L}^{p},$ its limit
is $P-$almost surely unique}. \\
The quotient set of $\mathscr{L}^{p}\left(\Omega,\mathcal{A},P\right)$
by the equivalence relation of equality $P-$almost surely is denoted
$\text{L}^{p}\left(\Omega,\mathcal{A},P\right).$ This is, by the
Minkowski inequallity, a normed vector space, whose norm is obtained
by taking the quotient of the semi-norm 
\[
X\mapsto\left[\mathbb{E}\left(\left|X\right|^{p}\right)\right]^{1/p},
\]
called the $p-$norm of $X,$ usually denoted $\left\Vert X\right\Vert _{p}=\left[\mathbb{E}\left(\left|X\right|^{p}\right)\right]^{1/p}.$
In this space of classes of random variables, the limit is then unique.
It is a common usage to denote the same way the random variable and
its equivalence class. \\
We do the same for the semi-norm and the associated norm and we talk
indifferently of \textbf{$\mathscr{L}^{p}-$convergence\mindex{Lp-convergence@$\mathscr{L}^p-$convergence}}
or \textbf{$\text{L}^{p}-$convergence\mindex{Lp-convergence@$\mathrm{L}^p-$convergence}}.
\end{itemize}
\end{remark}

The next theorem states the relations between convergence in probability
and $L^{p}-$convergence, and shows that if $p\geqslant1,$ the set
$\mathscr{L}^{p}\left(\Omega,\mathcal{A},P\right)$ is complete---and
non-separated. The space $L^{p}\left(\Omega,\mathcal{A},P\right)$
is then a Banach space. To prove this theorem, we will use the following
convexity inequality.

\begin{lemma}{Convexity Inequality}{}

Let $p\geqslant1.$ For every real numbers $a,b,c,$
\begin{equation}
\left|a-b\right|^{p}\leqslant2^{p-1}\left[\left|a-c\right|^{p}+\left|c-b\right|^{p}\right].\label{eq:convex_xp}
\end{equation}

\end{lemma}

\begin{proof}{}{}

The function $x\mapsto x^{p}$ is convex on $\mathbb{R}^{+}.$ Therefore,
for every nonnegative $u,v,$
\[
\left[\dfrac{1}{2}\left(u+v\right)\right]^{p}\leqslant\dfrac{1}{2}\left(u^{p}+v^{p}\right).
\]

Hence,
\[
\left(u+v\right)^{p}\leqslant2^{p-1}\left(u^{p}+v^{p}\right).
\]

By the triangle inequality and the growth of the function $x\mapsto x^{p}$
on $\mathbb{R}^{+},$ it follows that
\[
\left|a-b\right|^{p}\leqslant\left(\left|a-c\right|+\left|c-b\right|\right)^{p}\leqslant2^{p-1}\left[\left|a-c\right|^{p}+\left|c-b\right|^{p}\right].
\]

\end{proof}

\begin{theorem}{Relation between $\textrm{L}^p-$convergence and Convergence in Probability. Completeness of $\mathcal{L}^p\left(\Omega,\mathcal{A},P\right)$}{}

Let $p\geqslant1.$ Let $\left(X_{n}\right)_{n\in\mathbb{N}}$ be
a sequence of random variables admitting a moment of order $p.$ The
following assertions are equivalent:

(i) The sequence $\left(X_{n}\right)_{n\in\mathbb{N}}$ converges
in $\text{L}^{p}.$

(ii) The sequence $\left(X_{n}\right)_{n\in\mathbb{N}}$ is Cauchy
in $\text{L}^{p},$ that is
\[
\lim_{m,n\to+\infty}\mathbb{E}\left(\left|X_{n}-X_{m}\right|^{p}\right)=0.
\]

(iii) The sequence $\left(\left|X_{n}\right|^{p}\right)_{n\in\mathbb{N}}$
is equi-integrable and there exists $X\in\mathscr{L}^{p}\left(\Omega,\mathcal{A},P\right)$
such that $X_{n}\stackrel[n\to+\infty]{P}{\longrightarrow}X.$

\end{theorem}

\begin{proof}{}{}
\begin{itemize}
\item \textbf{(i)$\Rightarrow$(ii)}\\
If the sequence $\left(X_{n}\right)_{n\in\mathbb{N}}$ converges in
$L^{p},$ there exists $X\in\mathscr{L}^{p}\left(\Omega,\mathcal{A},P\right)$
such that
\[
\lim_{n\to+\infty}\mathbb{E}\left(\left|X_{n}-X\right|^{p}\right)=0.
\]
By the Minkowski inequality, for every $m,n\in\mathbb{N},$
\[
\left\Vert X_{n}-X_{m}\right\Vert _{p}\leqslant\left\Vert X_{n}-X\right\Vert _{p}+\left\Vert X-X_{m}\right\Vert _{p},
\]
which shows that the sequence is Cauchy---we just proved, in this
special case, the general fact that any convergent sequence for a
semi-norm is Cauchy with respect to this half-norm.
\item \textbf{(ii)$\Rightarrow$(iii)}\\
Let $\epsilon>0.$ Choose an integer $N$ such that for every $n,m\geqslant N,$
\[
\mathbb{E}\left(\left|X_{n}-X_{m}\right|^{p}\right)\leqslant\dfrac{\epsilon}{2^{p}}.
\]
It follows from the inequality $\refpar{eq:convex_xp}$ that, for
every $A\in\mathcal{A}$ and every $n\geqslant N,$
\begin{align*}
\intop_{A}\left|X_{n}\right|^{p}\text{d}P & \leqslant2^{p-1}\left[\intop_{A}\left|X_{N}\right|^{p}\text{d}P+\intop_{A}\left|X_{n}-X_{N}\right|^{p}\text{d}P\right]\\
 & \leqslant2^{p-1}\intop_{A}\left|X_{N}\right|^{p}\text{d}P+\dfrac{\epsilon}{2}.
\end{align*}
Hence, for every $A\in\mathcal{A},$
\[
\sup_{n\geqslant N}\intop_{A}\left|X_{n}\right|^{p}\text{d}P\leqslant2^{p-1}\intop_{A}\left|X_{N}\right|^{p}\text{d}P+\dfrac{\epsilon}{2}.
\]
Thus,
\begin{equation}
\sup_{n\geqslant N}\intop_{A}\left|X_{n}\right|^{p}\text{d}P\leqslant\sup_{n\leqslant N}\intop_{A}\left|X_{n}\right|^{p}\text{d}P+2^{p-1}\intop_{A}\left|X_{n}\right|^{p}\text{d}P+\dfrac{\epsilon}{2}.\label{eq:sup_intAX_np}
\end{equation}
It follows that the family $\left\{ \left|X_{n}\right|^{p}:\,n\in\mathbb{N}\right\} $
is bounded in $\text{L}^{1}.$ Moreover, the finite family $\left\{ \left|X_{n}\right|^{p}:\,n\leqslant N\right\} $
being equi-integrable, is in particular equi-continuous. The upper-bound
$\refpar{eq:sup_intAX_np}$ then implies that the family $\left\{ \left|X_{n}\right|^{p}:\,n\in\mathbb{N}\right\} $
is also equi-continuous, and hence equi-integrable since bounded in
$\text{L}^{1}.$ \\
Finally, by the growth of the function $x\mapsto x^{p}$ on $\mathbb{R}^{+}$
and by the Markov inequality, for every $\epsilon>0,$ for every $n$
and $m,$
\[
P\left(\left|X_{n}-X_{m}\right|>\epsilon\right)\leqslant\epsilon^{-p}\mathbb{E}\left(\left|X_{n}-X_{m}\right|^{p}\right).
\]
Thus, the sequence $\left(X_{n}\right)_{n\in\mathbb{N}}$ is Cauchy
in probability, and therefore converges in probability to a random
variable $X.$ Since the sequence $\left(\left|X_{n}\right|^{p}\right)_{n\in\mathbb{N}}$
is bounded in $\text{L}^{1},$ by Fatou lemma,
\[
\intop_{\Omega}\left|X\right|^{p}\text{d}P\leqslant\liminf_{n\to+\infty}\mathbb{E}\left(\left|X_{n}\right|^{p}\right)\leqslant\sup_{n\in\mathbb{N}}\mathbb{E}\left(\left|X_{n}\right|^{p}\right)<+\infty,
\]
which shows that $X\in\mathscr{L}^{p}\left(\Omega,\mathcal{A},P\right).$
\item \textbf{(iii)$\Rightarrow$(i)}\\
For every $\epsilon>0,$ by the inequality $\refpar{eq:sup_intAX_np},$
\[
\mathbb{E}\left(\left|X_{n}-X\right|^{p}\right)\leqslant\intop_{\left(\left|X_{n}-X\right|\leqslant\epsilon^{1/p}\right)}\left|X_{n}-X\right|^{p}\text{d}P+2^{p-1}\left[\intop_{\left(\left|X_{n}-X\right|>\epsilon^{1/p}\right)}\left(\left|X_{n}\right|^{p}+\left|X\right|^{p}\right)\text{d}P\right].
\]
Hence,
\begin{equation}
\mathbb{E}\left(\left|X_{n}-X\right|^{p}\right)\leqslant\epsilon+2^{p-1}\left[\intop_{\left(\left|X_{n}-X\right|>\epsilon^{1/p}\right)}\left|X_{n}\right|^{p}\text{d}P+\intop_{\left(\left|X_{n}-X\right|>\epsilon^{1/p}\right)}\left|X\right|^{p}\text{d}P\right].\label{eq:11_abs_exp_x_n-x_p}
\end{equation}
The equi-continuity of the family $\left\{ \left|X_{n}\right|^{p},n\in\mathbb{N};\left|X\right|^{p}\right\} $
allows us to find $\eta>0$ such that
\[
\sup{}_{n\in\mathbb{N}}\left(\intop_{A}\left|X_{n}\right|^{p}\text{d}P\right)+\intop_{A}\left|X\right|^{p}\text{d}P\leqslant\dfrac{\epsilon}{2^{p-1}}
\]
whenever $P\left(A\right)\leqslant\eta.$\\
Moreover, the convergence in probability of the sequence $\left(X_{n}\right)_{n\in\mathbb{N}}$
towards $X$ implies that there exists $N$ such that, for every $n\geqslant N,$
\[
P\left(\left|X_{n}-X\right|>\epsilon^{\frac{1}{p}}\right)\leqslant\eta.
\]
By the inequality $\refpar{eq:11_abs_exp_x_n-x_p},$
\[
\forall\epsilon>0,\,\,\,\,\limsup_{n\to+\infty}\mathbb{E}\left(\left|X_{n}-X\right|^{p}\right)\leqslant2\epsilon.
\]
This proves that the sequence $\left(X_{n}\right)_{n\in\mathbb{N}}$
converges in $L^{p}$ to $X.$
\end{itemize}
\end{proof}

\begin{counterexample}{}{}

If, for every $n\in\mathbb{N}^{\ast},$ $X_{n}$ has the law $\dfrac{1}{n}\delta_{n^{2}}+\left(1-\dfrac{1}{n}\right)\delta_{0},$
then for every $\epsilon>0,$
\[
P\left(\left|X_{n}\right|>\epsilon\right)=\dfrac{1}{n}\,\,\,\,\text{and}\,\,\,\,\mathbb{E}\left(X_{n}\right)=n.
\]

The sequence $\left(X_{n}\right)_{n\in\mathbb{N}}$ converges in probability
to $0,$ but does not converge in $\mathscr{L}^{1}.$

\end{counterexample}

\section{Series of Independent Random Variables}\label{sec:Series-of-Independent}

We now study a sufficient condition $P-$almost sure convergence and
$\text{L}^{2}-$convergence of series of independent random variables
admitting a second-order moment. We first give the Kolmogorov inequality
which extends the Chebyshev inequality.

\begin{theorem}{Kolmogorov Inequality}{kolmogorov_ineq}

Let $X_{1},X_{2},\cdots,X_{n}$ be $n$ independent centered random
variables, admitting a second-order moment.

Then, for every $\epsilon>0,$
\[
P\left(\max_{1\leqslant k\leqslant n}\left|\sum^{k}_{i=1}X_{i}\right|\geqslant\epsilon\right)\leqslant\dfrac{1}{\epsilon^{2}}\left(\sum^{n}_{i=1}\sigma^{2}_{X_{i}}\right).
\]

\end{theorem}

\begin{proof}{}{}

Let us denote, for every nonnegative integer $j$ and $k$ such that
$1\leqslant j,\,k\leqslant n,$
\[
S_{j}=\sum^{j}_{i=1}X_{i}\,\,\,\,\text{and}\,\,\,\,M_{k}=\underset{1\leqslant j\leqslant k}{\max}\left|S_{j}\right|.
\]

We want to find an upper-bound for the probability of the set $E=\left(M_{n}\geqslant\epsilon\right).$
If this set is empty, the inequality is obvious; thus we consider
the case where the set is non-empty. Let us make appear the index
$i$ for which $\left|S_{j}\right|$ exceeds $\epsilon$ for the first
time. For this purpose, we consider the sets
\[
E_{1}=\left(\left|S_{1}\right|\geqslant\epsilon\right)
\]
and, if $2\leqslant i\leqslant n,$
\[
E_{i}=\left(\left|S_{i}\right|\geqslant\epsilon\right)\cap\left[\bigcap^{i-1}_{j=1}\left(\left|S_{j}\right|<\epsilon\right)\right]
\]

These sets $\left(E_{i}\right)_{1\leqslant i\leqslant n}$ form a
partition of $E,$ and therefore,
\[
P\left(E\right)=\sum^{n}_{i=1}P\left(E_{i}\right).
\]

By definition of the $E_{i},$ $\left|S_{i}\right|\geqslant\epsilon$
which implies that $S^{2}_{i}\geqslant\epsilon^{2},$ which is equivalent
to $\boldsymbol{1}_{E_{i}}S^{2}_{i}\geqslant\boldsymbol{1}_{E_{i}}\epsilon^{2}.$
Hence, $\mathbb{E}\left(\boldsymbol{1}_{E_{i}}S^{2}_{i}\right)\geqslant\mathbb{E}\left(\boldsymbol{1}_{E_{i}}\epsilon^{2}\right).$
Since $\mathbb{E}\left(\boldsymbol{1}_{E_{i}}\epsilon^{2}\right)=\epsilon^{2}\mathbb{E}\left(\boldsymbol{1}_{E_{i}}\right)=\epsilon^{2}P\left(E_{i}\right),$
it follows by rearranging that
\[
P\left(E_{i}\right)\leqslant\dfrac{1}{\epsilon^{2}}\mathbb{E}\left(\boldsymbol{1}_{E_{i}}S^{2}_{i}\right).
\]
Hence,
\begin{equation}
P\left(E\right)\leqslant\dfrac{1}{\epsilon^{2}}\sum^{n}_{i=1}\mathbb{E}\left(\boldsymbol{1}_{E_{i}}S^{2}_{i}\right).\label{eq:PE_bound}
\end{equation}

Now let us prove that, if $1\leqslant i\leqslant n,$ then
\begin{equation}
\mathbb{E}\left(\boldsymbol{1}_{E_{i}}S^{2}_{i}\right)\leqslant\mathbb{E}\left(\boldsymbol{1}_{E_{i}}S^{2}_{n}\right).\label{eq:exp_1_E_iS_i^2}
\end{equation}

Indeed,
\begin{align*}
\mathbb{E}\left(\boldsymbol{1}_{E_{i}}S^{2}_{n}\right) & =\mathbb{E}\left(\boldsymbol{1}_{E_{i}}\left(S_{i}+\sum^{n}_{j=i+1}X_{j}\right)^{2}\right)\\
 & =\mathbb{E}\left(\boldsymbol{1}_{E_{i}}S^{2}_{i}\right)+2\mathbb{E}\left(\boldsymbol{1}_{E_{i}}S_{i}\sum^{n}_{j=i+1}X_{j}\right)+\mathbb{E}\left(\boldsymbol{1}_{E_{i}}\left(\sum^{n}_{j=i+1}X_{j}\right)^{2}\right).
\end{align*}
The random variables $\boldsymbol{1}_{E_{i}}S_{i}$ and $\sum^{n}_{j=i+1}X_{j}$
are independent and the random variable $\sum^{n}_{j=i+1}X_{j}$ is
centered. The middle term on the right-hand side is therefore zero.
Since the third term is nonnegative, we obtain the inequality $\refpar{eq:exp_1_E_iS_i^2}.$

By substituting this bound into the equation $\refpar{eq:exp_1_E_iS_i^2},$
and using the fact that the sets $E_{i}$ form a partition of $E,$
we obtain
\[
P\left(E\right)\leqslant\dfrac{1}{\epsilon^{2}}\sum^{n}_{i=1}\mathbb{E}\left(\boldsymbol{1}_{E_{i}}S^{2}_{i}\right)=\dfrac{1}{\epsilon^{2}}\mathbb{E}\left(\boldsymbol{1}_{E}S^{2}_{n}\right)\leqslant\dfrac{1}{\epsilon^{2}}\mathbb{E}\left(S^{2}_{n}\right).
\]

The random variables $X_{i}$ being centered and independent,
\[
\mathbb{E}\left(S^{2}_{n}\right)=\sum^{n}_{i=1}\sigma^{2}_{X_{1}}
\]
which concludes the proof.

\end{proof}

We deduce from the Kolmogorov inequality a sufficient condition of
$P-$almost sure convergence of a series of independent random variables.

\begin{proposition}{Sufficient Condition P-Almost Sure for a Series of Independent Random Variables}{sc_p-as_series_indep_rv}

Let $\left(X_{n}\right)_{n\in\mathbb{N}}$ be a sequence of independent
centered real-valued random variables admitting a second-order moment.

If $\sum^{+\infty}_{n=1}\sigma^{2}_{X_{n}}<+\infty,$ then the series
$\sum X_{n}$ with general term $X_{n}$ converges $P-$almost surely
and in $\text{L}^{2}.$

\end{proposition}

\begin{proof}{}{}

We first prove the $P-$almost sure convergence. For $m\in\mathbb{N}^{\ast},$
let us denote
\[
S_{m}=\sum^{m}_{i=1}X_{i}\,\,\,\,\,\,\,A_{m}=\sup_{k\in\mathbb{N}^{\ast}}\left|S_{m+k}-S_{m}\right|\,\,\,\,\text{and}\,\,\,\,A=\inf_{m\in\mathbb{N}^{\ast}}A_{m}.
\]

By the Cauchy criterion for numerical series,
\[
\left\{ \sum X_{n}\,\text{converges}\right\} =\left\{ A=0\right\} .
\]

But,
\[
\left\{ A\neq0\right\} \subset\bigcup_{n\in\mathbb{N}^{\ast}}\left\{ A>\dfrac{1}{n}\right\} 
\]
and, for every $n\in\mathbb{N}^{\ast},$
\[
\left\{ A>\dfrac{1}{n}\right\} \subset\bigcap_{m\in\mathbb{N}^{\ast}}\left\{ A_{m}>\dfrac{1}{n}\right\} ,
\]
which yields the inclusion
\begin{equation}
\left\{ A\neq0\right\} \subset\bigcup_{n\in\mathbb{N}^{\ast}}\bigcap_{m\in\mathbb{N}^{\ast}}\left\{ A_{m}>\dfrac{1}{n}\right\} .\label{eq:Aneq0subset}
\end{equation}

Since
\[
\sup_{k\in\mathbb{N}^{\ast}}\left|S_{m+k}-S_{m}\right|=\lim_{r\to+\infty}\nearrow\sup_{1\leqslant k\leqslant r}\left|S_{m+k}-S_{m}\right|,
\]
the sequence of sets
\[
\left(\left\{ \sup_{1\leqslant k\leqslant r}\left|S_{m+k}-S_{m}\right|>\dfrac{1}{n}\right\} \right)_{r\in\mathbb{N}^{\ast}}
\]
is nondecreasing, and
\begin{equation}
\left\{ A_{m}>\dfrac{1}{n}\right\} =\bigcup_{r\in\mathbb{N}^{\ast}}\left\{ \sup_{1\leqslant k\leqslant r}\left|S_{m+k}-S_{m}\right|>\dfrac{1}{n}\right\} .\label{eq:A_m_as_union_sup}
\end{equation}

By the Kolmogorov inequality, it follows that
\[
P\left(\sup_{1\leqslant k\leqslant r}\left|S_{m+k}-S_{m}\right|>\dfrac{1}{n}\right)\leqslant n^{2}\sum^{r}_{j=m+1}\sigma^{2}_{X_{j}}.
\]

Since the equality $\refpar{eq:A_m_as_union_sup}$ holds for a nondecreasing
sequence of sets,
\[
P\left(A_{m}>\dfrac{1}{n}\right)=\lim_{r\to+\infty}P\left(\sup_{1\leqslant k\leqslant r}\left|S_{m+k}-S_{m}\right|>\dfrac{1}{n}\right)\leqslant n^{2}\sum^{+\infty}_{j=m+1}\sigma^{2}_{X_{j}}.
\]

Hence, for every $m\in\mathbb{N}^{\ast},$
\[
0\leqslant P\left(\bigcap_{p\in\mathbb{N}^{\ast}}\left(A_{p}>\dfrac{1}{n}\right)\right)\leqslant P\left(A_{m}>\dfrac{1}{n}\right)\leqslant n^{2}\sum^{+\infty}_{j=m+1}\sigma^{2}_{X_{j}}.
\]

The right-hand term converges to 0 as $m$ tends to infinity, since
it is the remainder of a convergent series. It follows that, for every
$n\in\mathbb{N}^{\ast},$
\[
P\left(\bigcap_{p\in\mathbb{N}^{\ast}}\left(A_{p}>\dfrac{1}{n}\right)\right)=0.
\]

By the inclusion $\refpar{eq:Aneq0subset},$
\[
P\left(A\neq0\right)=0,
\]
that is, the series with general term $X_{n}$ converges $P-$almost
surely.

Since the sequence of partial sums is Cauchy for the $\text{L}^{2}-$norm,
there is also convergence in $\text{L}^{2}.$ Indeed, the random variables
$X_{n}$ being centered and independent, we have, if $m<n,$
\[
\mathbb{E}\left(\left(S_{n}-S_{m}\right)^{2}\right)=\sum^{n}_{j=m+1}\sigma^{2}_{X_{j}}.
\]

This proves the result, since the series of variances converges. 

\end{proof}

\section{Laws of Large Numbers}\label{sec:Laws-of-Large}

When studying random phenomena, one often needs to analyze the convergence
of sequence of \textbf{arithmetic means} of a sequence of \textbf{independent}
random variables having the \textbf{same law}. This occurs, for instance,
in statistics and estimation theory: if $X$ is a random variable
modelling a characteristic linked to a random phenomenon, one often
needs to estimate its law, or some parameters of this law, from a
sequence of outcomes of this phenomenon obtained during \textbf{independent}
experiences. This leads to introduce a sequence $\left(X_{n}\right)_{n\in\mathbb{N}^{\ast}}$
of independent random variables with the same law than $X,$ and to
studying, for a given function $f,$ the sequence whose general term
is of the form
\[
\dfrac{1}{n}\sum^{n}_{j=1}f\left(X_{j}\right).
\]

If $\left(X_{n}\right)_{n\in\mathbb{N}^{\ast}}$ is a sequence of
real-valued random variables, we denote for every $n\in\mathbb{N}^{\ast},$
\[
\overline{X}_{n}=\dfrac{1}{n}\sum^{n}_{j=1}X_{j}.
\]
In statistical terminology, $\overline{X}_{n}$ is called the \textbf{\mindex{empirical!mean}empirical
mean} of the sample $\left(X_{1},X_{2},\cdots,X_{n}\right).$ We give
the name of \textbf{\mindex{law ! large numbers}\index{law of large numbers}law
of large numbers} to refer to two main theorems asserting the convergence
of the sequence with general term $\overline{X}_{n}$ under certain
assumptions. For the \textbf{weak law}\index{weak law}\textbf{\mindex{law ! weak}},
it is a \textbf{convergence in probability}. For the \textbf{strong
law}\index{strong law}\textbf{\mindex{law ! strong}}, it is the
\textbf{almost sure convergence}. 

The name ``law of large numbers''---weak or strong, depending on
the case---is also given to many variants of these two results, obtained
under stronger or weaker assumptions.

It is worth noting that for the weak laws, the hypothesis of independence
of the random variables $X_{n}$ is not necessary\footnote{We assumed mutual independence in the statement of Theorem $\ref{th:weak_law}$
only for simplicity.}: uncorrelation or pairwise independence of the random variables is
sufficient.

By contrast, the \textbf{mutal independence}\footnote{The laws of large numbers appear in many contexts and generate an
extensive litterature. In particular, in the framework of our study,
the independence hypothesis can be removed using martingal theory.} is required in the classical versions of strong laws for the random
variables. We leave it to the interested reader to formulate a version
of the strong law for random variables that are pairwise independent.

To begin, let us recall two elementary lemmas from analysis that will
be used several times in establishing the laws of large numbers.

\begin{figure}[t]
\begin{center}\includegraphics[width=0.4\textwidth]{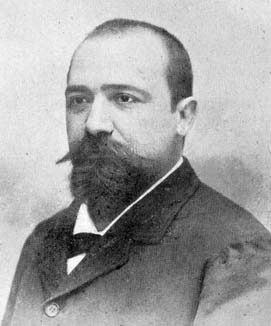}

{\tiny Public domain}\end{center}

\caption{\textbf{\protect\href{https://it.wikipedia.org/wiki/Ernesto_Cesaro}{Ernesto Cesàro}}
(1805-1906)}\sindex[fam]{Cesàro, Ernesto}
\end{figure}

\begin{lemma}{Cesàro\footnotemark Lemma}{}

Let $\left(x_{n}\right)_{n\in\mathbb{N}^{\ast}}$ be a sequence of
real numbers converging to $x$ when $n$ tends to infinity. 

Then the sequence whose general term is $\dfrac{1}{n}\sum^{n}_{j=1}x_{j}$
also converges, and its limit is $x.$

\end{lemma}

\footnotetext{\textbf{\sindex[fam]{Cesàro, Ernesto}\href{https://it.wikipedia.org/wiki/Ernesto_Cesaro}{Ernesto Cesàro}}
(1805-1906), born in Naples. He became professor at the University
of Naples in 1883. He worked in several areas of mathematics: in particular,
he studied the relationships between arithmetic and integral calculus,
and the behaviour of entire series on the circle of convergence.}

\begin{proof}{}{}

Let $\epsilon>0$ be an arbitrary real number. Let $N$ be an integer
such that, for every $n\geqslant N,$ $\left|x_{n}-x\right|\leqslant\epsilon.$

Since
\[
\left|\dfrac{1}{n}\sum^{n}_{j=1}x_{j}-x\right|\leqslant\dfrac{1}{n}\sum^{N}_{j=1}\left|x_{j}-x\right|+\dfrac{1}{n}\sum^{n}_{j=N+1}\left|x_{j}-x\right|\leqslant\dfrac{1}{n}\sum^{N}_{j=1}\left|x_{j}-x\right|+\epsilon,
\]
we obtain
\[
\limsup_{n\to+\infty}\left|\dfrac{1}{n}\sum^{n}_{j=1}x_{j}-x\right|\leqslant\epsilon
\]
Since $\epsilon$ is arbitrary, this proves the result.

\end{proof}

\begin{figure}[t]
\begin{center}\includegraphics[width=0.4\textwidth]{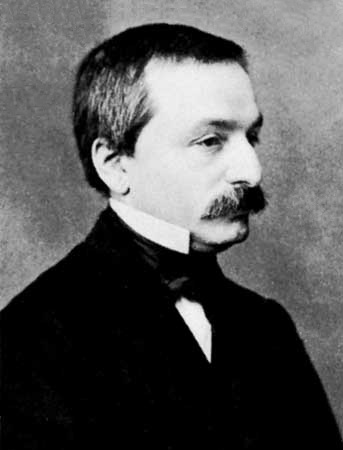}

{\tiny Public domain}\end{center}

\caption{\textbf{\protect\href{https://fr.wikipedia.org/wiki/Leopold_Kronecker}{Leopold Kronecker}}
(1823-1891)}\sindex[fam]{Kronecker, Leopold}
\end{figure}

\begin{lemma}{Kronecker\footnotemark Lemma}{}

Consider a series with general real-valued term $x_{n}$ that converges,
and a nondecreasing sequence $\left(b_{n}\right)_{n\in\mathbb{N}^{\ast}}$
of real numbers tending to infinity with $n.$

Then
\[
\lim_{n\to+\infty}\dfrac{1}{b_{n}}\sum^{n}_{j=1}b_{j}x_{j}=0.
\]

\end{lemma}

\footnotetext{\textbf{\sindex[fam]{Kronecker, Leopold}\href{https://fr.wikipedia.org/wiki/Leopold_Kronecker}{Leopold Kronecker}}
(1823-1891) was born in Liegnitz, in Poland. After studying in Berlin
and Bonn, he became wealthy in finance, which allowed him to devote
himself entirely to mathematics. From 1861 onward, he taught in Berlin.
His work concerned the theory of equations, elliptic functions and
algebraic number theory. He strongly opposed Cantor set theory, and
Weierstrass construction of the real numbers.}

\begin{proof}{}{}

We denote $S=\sum^{+\infty}_{j=1}x_{j}$ and, for every $n\in\mathbb{N}^{\ast},$
\[
S_{n}=-S+\sum^{n}_{j=1}x_{j}.
\]

The sequence $\left(S_{n}\right)_{n\in\mathbb{N^{\ast}}}$ tends to
0 when $n$ tends to infinity.

We have
\[
x_{n}=S_{n}-S_{n-1}
\]
and thus, by Abel transformation---or summation---for every integer
$n$ and $N$ such that $n>N\geqslant2,$
\[
\sum^{n}_{j=N}b_{j}x_{j}=\sum^{n}_{j=N}b_{j}\left(S_{j}-S_{j-1}\right)=S_{n}b_{n}-b_{N}S_{N-1}-\sum^{n-1}_{j=N}S_{j}\left(b_{j+1}-b_{j}\right).
\]
Hence, when $b_{n}\neq0,$
\[
\dfrac{1}{b_{n}}\sum^{n}_{j=1}b_{j}x_{j}=\dfrac{1}{b_{n}}\sum^{N}_{j=1}b_{j}x_{j}-\dfrac{b_{N}S_{N-1}}{b_{n}}+S_{n}-\dfrac{1}{b_{n}}\sum^{n-1}_{j=N}S_{j}\left(b_{j+1}-b_{j}\right).
\]

The sequence whose general term is $\dfrac{1}{b_{n}}\sum^{N}_{j=1}b_{j}x_{j}-\dfrac{b_{N}S_{N-1}}{b_{n}}+S_{n}$
tends to 0 when $n$ tends to infinity. 

Hence, fix $\epsilon>0$ and choose $N$ such that, for every $n\geqslant N,$
\[
\left|\dfrac{1}{b_{n}}\sum^{N}_{j=1}b_{j}x_{j}-\dfrac{b_{N}S_{N-1}}{b_{n}}+S_{n}\right|\leqslant\dfrac{\epsilon}{2}\,\,\,\,\,\,\,\,\text{and}\,\,\,\,\,\,\,\,\left|S_{n}\right|\leqslant\dfrac{\epsilon}{2}.
\]

Since the sequence $\left(b_{n}\right)_{n\in\mathbb{N}^{\ast}}$ is
nondecreasing,
\[
\left|\sum^{n-1}_{j=N}S_{j}\left(b_{j+1}-b_{j}\right)\right|\leqslant\sum^{n-1}_{j=N}\left|S_{j}\right|\left(b_{j+1}-b_{j}\right)\leqslant\dfrac{\epsilon}{2}\sum^{n-1}_{j=N}\left(b_{j+1}-b_{j}\right)=\dfrac{\epsilon}{2}\left(b_{n}-b_{N}\right).
\]

It follows that
\[
\limsup_{n\to+\infty}\dfrac{1}{b_{n}}\left|\sum^{n-1}_{j=N}S_{j}\left(b_{j+1}-b_{j}\right)\right|\leqslant\dfrac{\epsilon}{2}
\]
and thus,
\[
\limsup_{n\to+\infty}\left|\dfrac{1}{b_{n}}\sum^{n}_{j=1}b_{j}x_{j}\right|\leqslant\epsilon,
\]
which shows the result, since $\epsilon$ is arbitrary.

\end{proof}

\begin{theorem}{Weak Law of Large Numbers}{weak_law}

Let $\left(X_{n}\right)_{n\in\mathbb{N}^{\ast}}$ be a sequence of
random variables defined on the probabilized space $\left(\Omega,\mathscr{A},P\right),$
admitting a second-order moment, and pairwise uncorrelated\footnotemark.
We suppose the convergence of the sequences
\[
\dfrac{1}{n}\sum^{n}_{j=1}\mathbb{E}\left(X_{j}\right)\underset{n\to+\infty}{\longrightarrow}m\,\,\,\,\,\,\,\text{and}\,\,\,\,\,\,\,\,\dfrac{1}{n^{2}}\sum^{n}_{j=1}\sigma^{2}_{X_{j}}\underset{n\to+\infty}{\longrightarrow}0.
\]

Then, the sequence of random variables $\overline{X}_{n}=\dfrac{1}{n}\sum^{n}_{j=1}X_{j}$
\textbf{converges in probability to $m.$}

\end{theorem}

\footnotetext{Two real-valued random variables admitting a second-order
moment are said \textbf{uncorrelated\index{uncorrelated}} if their
correlation coefficient is zero---which is equivalent to saying that
their covariance is zero.}

\begin{proof}{}{}

We have
\[
\mathbb{E}\left(\overline{X}_{n}\right)=\dfrac{1}{n}\sum^{n}_{j=1}\mathbb{E}\left(X_{j}\right).
\]

Since the random variables $X_{n}$ are pairwise uncorrelated, we
also have
\[
\sigma^{2}_{\overline{X}_{n}}=\dfrac{1}{n^{2}}\sum^{n}_{j=1}\sigma^{2}_{X_{j}}.
\]

By the triangle inequality,
\[
\left|\overline{X}_{n}-m\right|\leqslant\left|\overline{X}_{n}-\dfrac{1}{n}\sum^{n}_{j=1}\mathbb{E}\left(X_{j}\right)\right|+\left|\dfrac{1}{n}\sum^{n}_{j=1}\mathbb{E}\left(X_{j}\right)-m\right|.
\]

But, for every $\epsilon>0,$ there exists $N\left(\epsilon\right)\in\mathbb{N}^{\ast}$
such that, for every $n\geqslant N\left(\epsilon\right),$
\[
\left|\dfrac{1}{n}\sum^{n}_{j=1}\mathbb{E}\left(X_{j}\right)-m\right|\leqslant\dfrac{\epsilon}{2}.
\]

For every $n\geqslant N\left(\epsilon\right),$ we therefore have
the inclusion of sets
\[
\left(\left|\overline{X}_{n}-\dfrac{1}{n}\sum^{n}_{j=1}\mathbb{E}\left(X_{j}\right)\right|\leqslant\dfrac{\epsilon}{2}\right)\subset\left(\left|\overline{X}_{n}-m\right|\leqslant\epsilon\right).
\]

Taking complements gives
\[
\left(\left|\overline{X}_{n}-m\right|>\epsilon\right)\subset\left(\left|\overline{X}_{n}-\dfrac{1}{n}\sum^{n}_{j=1}\mathbb{E}\left(X_{j}\right)\right|>\dfrac{\epsilon}{2}\right).
\]

By the Bienaymé-Chebyshev inequality,
\[
P\left(\left|\overline{X}_{n}-\dfrac{1}{n}\sum^{n}_{j=1}\mathbb{E}\left(X_{j}\right)\right|>\dfrac{\epsilon}{2}\right)\leqslant\dfrac{4}{\epsilon^{2}}\sigma^{2}_{\overline{X}_{n}}=\dfrac{4}{\epsilon^{2}}\dfrac{1}{n^{2}}\sum^{n}_{j=1}\sigma^{2}_{X_{j}}.
\]

It follows that, for every $n\geqslant N\left(\epsilon\right),$
\[
P\left(\left|\overline{X}_{n}-m\right|>\epsilon\right)\leqslant\dfrac{4}{\epsilon^{2}}\dfrac{1}{n^{2}}\sum^{n}_{j=1}\sigma^{2}_{X_{j}},
\]
which, by using the second assumption, shows the result.

\end{proof}

\begin{figure}[t]
\begin{center}\includegraphics[width=0.4\textwidth]{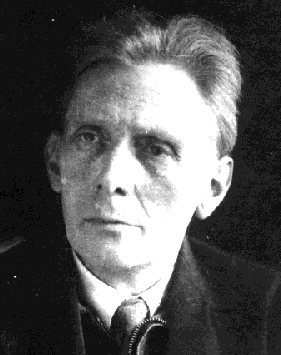}

{\tiny Public domain}\end{center}

\caption{\textbf{\protect\href{https://en.wikipedia.org/wiki/Aleksandr_Khinchin}{Aleksandr Yakovlevich Khinchin}}
(1894 - 1959)}\sindex[fam]{Khinchin, Alexander}
\end{figure}

\begin{remark}{}{}

The assumption on the variances restricts the random variables to
be not too dispersed around their average. 

All the assumptions of the previous theorem are satisfied, in particular,
if the random variables $X_{n}$ are independent and with same law,
and if $X_{1}$ admits a second-order moment. 

In fact, if the random variables are independent and with same law,
the mere existence of a first-order moment is sufficient, as shown
by the \textbf{Khinchin}\footnotemark theorem stated below.

\end{remark}

\footnotetext{\textbf{\href{https://en.wikipedia.org/wiki/Aleksandr_Khinchin}{Alexander Khinchin}\sindex[fam]{Khinchin, Alexander}}
(1894-1959) was a professor at Moscow University from 1922 onward.
His work concerned real analysis, number theory and probability theory.
He introduced, at the same time as P. Lévy the notion of a random
variable. He also introduced the definition of a stationary random
process.}

Before studying this theorem, the reader who wishes may recall a special
case of the previous theorem---indeed, historically earlier---the
Bernoulli theorem stated in Theorem $\ref{th:bernoulli_th},$ along
with the remark following it. 

\begin{theorem}{Khinchin Theorem. Weak Law of Large Numbers}{}

Let $\left(X_{n}\right)_{n\in\mathbb{N}^{\ast}}$ be a sequence of
random variables defined on a probabilized space $\left(\Omega,\mathcal{A},P\right),$
pairwise independent, following the same law $\mu,$ and admitting
an expectation.

Then the sequence of random variables
\[
\overline{X}_{n}=\dfrac{1}{n}\sum^{n}_{j=1}X_{j}
\]
converges in probability to their common expectation $\mathbb{E}\left(X_{1}\right).$

\end{theorem}

\begin{proof}{}{}

We return to Theorem $\ref{th:weak_law}$ by means of a truncation
argument.

Let
\[
Y_{n}=\boldsymbol{1}_{\left(\left|X_{n}\right|\leqslant n\right)}X_{n}\,\,\,\,\,\text{and}\,\,\,\,\overline{Y}_{n}=\dfrac{1}{n}\sum^{n}_{i=1}Y_{i}.
\]

The random variables $Y_{n}$ are pairwise independent and bounded.
We now prove that they satisfy the assumptions of Theorem $\ref{th:weak_law}.$ 

We have
\begin{align*}
\sum^{n}_{j=1}\mathbb{E}\left(Y_{j}\right) & =\sum^{n}_{j=1}\intop_{\left(\left|x\right|\leqslant j\right)}x\text{d}\mu\left(x\right)\\
 & =\sum^{n}_{j=1}\sum^{j-1}_{k=0}\intop_{\left(k<\left|x\right|\leqslant k+1\right)}x\text{d}\mu\left(x\right)\\
 & =\sum^{n-1}_{k=0}\sum^{n}_{j=k+1}\intop_{\left(k<\left|x\right|\leqslant k+1\right)}x\text{d}\mu\left(x\right).
\end{align*}
Hence
\[
\sum^{n}_{j=1}\mathbb{E}\left(Y_{j}\right)=\sum^{n-1}_{k=0}\left(n-k\right)\intop_{\left(k<\left|x\right|\leqslant k+1\right)}x\text{d}\mu\left(x\right),
\]
which yields,
\[
\dfrac{1}{n}\sum^{n}_{j=1}\mathbb{E}\left(Y_{j}\right)=\intop_{\left(\left|x\right|\leqslant n\right)}x\text{d}\mu\left(x\right)-\dfrac{1}{n}\sum^{n-1}_{k=0}k\intop_{\left(k<\left|x\right|\leqslant k+1\right)}x\text{d}\mu\left(x\right).
\]

The existence of the mean of $X_{1}$ implies that
\[
\lim_{n\to+\infty}\int_{\left(\left|x\right|\leqslant n\right)}x\text{d}\mu\left(x\right)=\mathbb{E}\left(X_{1}\right),
\]
and that the series with general term
\[
\intop_{\left(k<\left|x\right|\leqslant k+1\right)}x\text{d}\mu\left(x\right)
\]
converges. 

The Kronecker lemma then yields
\[
\lim_{n\to+\infty}\dfrac{1}{n}\sum^{n-1}_{k=0}k\intop_{\left(k<\left|x\right|\leqslant k+1\right)}x\text{d}\mu\left(x\right)=0,
\]
which gives
\[
\lim_{n\to+\infty}\dfrac{1}{n}\sum^{n}_{j=1}\mathbb{E}\left(Y_{j}\right)=\mathbb{E}\left(X_{1}\right).
\]

Furthermore,
\[
\sigma^{2}_{Y_{j}}\leqslant\mathbb{E}\left(Y^{2}_{j}\right)=\int_{\left(\left|x\right|\leqslant j\right)}x^{2}\text{d}\mu\left(x\right)\leqslant\int_{\left(\left|x\right|\leqslant n\right)}x^{2}\text{d}\mu\left(x\right),
\]
and thus
\begin{align*}
0 & \leqslant\dfrac{1}{n^{2}}\sum^{n}_{j=1}\sigma^{2}_{Y_{j}}\\
 & \leqslant\dfrac{1}{n}\intop_{\left(\left|x\right|\leqslant n\right)}x^{2}\text{d}\mu\left(x\right)\\
 & \leqslant\dfrac{1}{n}\left[\intop_{\left(\left|x\right|\leqslant\sqrt{n}\right)}x^{2}\text{d}\mu\left(x\right)+\intop_{\left(\sqrt{n}<\left|x\right|\leqslant n\right)}x^{2}\text{d}\mu\left(x\right)\right].
\end{align*}

It follows that
\[
0\leqslant\dfrac{1}{n^{2}}\sum^{n}_{j=1}\sigma^{2}_{Y_{j}}\leqslant\dfrac{1}{\sqrt{n}}\intop_{\left(\left|x\right|\leqslant\sqrt{n}\right)}\left|x\right|\text{d}\mu\left(x\right)+\intop_{\left(\sqrt{n}<\left|x\right|\right)}\left|x\right|\text{d}\mu\left(x\right).
\]

Since $\intop\left|x\right|\text{d}\mu\left(x\right)<+\infty,$ the
right-hand side tends to 0, and
\[
\lim_{n\to+\infty}\dfrac{1}{n^{2}}\sum^{n}_{j=1}\sigma^{2}_{Y_{j}}=0.
\]

By Theorem $\ref{th:weak_law},$
\begin{equation}
\overline{Y}_{n}\stackrel[n\to+\infty]{P}{\longrightarrow}\mathbb{E}\left(X_{1}\right).\label{eq:overlineY_P_conv}
\end{equation}

If $n>r,$ denote
\[
\overline{Y}_{n,r}=\dfrac{1}{n}\left(\sum^{r}_{j=1}X_{j}+\sum^{n}_{j=r+1}Y_{j}\right).
\]

We have
\begin{align*}
P\left(\overline{Y}_{n,r}\neq\overline{X}_{n}\right) & \leqslant P\left(\bigcup^{n}_{j=r+1}\left(Y_{j}\neq X_{j}\right)\right)\\
 & \leqslant\sum^{n}_{j=r+1}P\left(Y_{j}\neq X_{j}\right)=\sum^{n}_{j=r+1}\intop_{\mathbb{R}}\boldsymbol{1}_{\left[0,\left|x\right|\right[}\left(j\right)\text{d}\mu\left(x\right).
\end{align*}

Since
\[
\sum^{n}_{j=r+1}\boldsymbol{1}_{\left[0,\left|x\right|\right[}\left(j\right)\leqslant\sum^{+\infty}_{j=r+1}\boldsymbol{1}_{\left[0,\left|x\right|\right[}\left(j\right)\leqslant\boldsymbol{1}_{\left(\left|x\right|>r\right)}\left|x\right|,
\]
it follows that
\[
P\left(\overline{Y}_{n,r}\neq\overline{X}_{n}\right)\leqslant\intop_{\left(\left|x\right|>r\right)}\left|x\right|\text{d}\mu\left(x\right).
\]

Thus, for every $\epsilon>0,$ there exists $r>0$ such that, as soon
as $n>r,$
\[
P\left(\overline{Y}_{n,r}\neq\overline{X}_{n}\right)\leqslant\dfrac{\epsilon}{2}.
\]
Then, for every $\delta>0,$
\begin{multline*}
P\left(\left|\overline{X}_{n}-\mathbb{E}\left(X_{1}\right)\right|>\delta\right)=P\left(\left(\left|\overline{X}_{n}-\mathbb{E}\left(X_{1}\right)\right|>\delta\right)\cap\left(\overline{Y}_{n,r}\neq\overline{X}_{n}\right)\right)\\
+P\left(\left(\left|\overline{X}_{n}-\mathbb{E}\left(X_{1}\right)\right|>\delta\right)\cap\left(\overline{Y}_{n,r}=\overline{X}_{n}\right)\right).
\end{multline*}
Hence
\begin{align*}
P\left(\left|\overline{X}_{n}-\mathbb{E}\left(X_{1}\right)\right|>\delta\right) & \leqslant P\left(\overline{Y}_{n,r}\neq\overline{X}_{n}\right)+P\left(\left|\overline{Y}_{n,r}-\mathbb{E}\left(X_{1}\right)\right|>\delta\right)\\
 & \leqslant\dfrac{\epsilon}{2}+P\left(\left|\overline{Y}_{n,r}-\mathbb{E}\left(X_{1}\right)\right|>\delta\right).
\end{align*}

By $\refpar{eq:overlineY_P_conv},$ 
\[
\lim_{n\to+\infty}P\left(\left|\overline{Y}_{n,r}-\mathbb{E}\left(X_{1}\right)\right|>\delta\right)=0.
\]
Therefore there exists $N>r$ such that, for every $n\geqslant N,$
\[
P\left(\left|\overline{Y}_{n,r}-\mathbb{E}\left(X_{1}\right)\right|>\delta\right)\leqslant\dfrac{\epsilon}{2}.
\]

In summary, we have shown that, for every $\epsilon>0$ and every
$\delta>0,$ there exists $N$ such that, for every $n\geqslant N,$
\[
P\left(\left|\overline{X}_{n}-\mathbb{E}\left(X_{1}\right)\right|>\delta\right)\leqslant\epsilon,
\]
which means that the sequence $\left(\overline{X_{n}}\right)_{n\in\mathbb{N}^{\ast}}$
converges in probability to $\mathbb{E}\left(X_{1}\right).$

\end{proof}

\begin{theorem}{Strong Law of Large Numbers}{strong_law_large_nb}

Let $\left(X_{n}\right)_{n\in\mathbb{N}^{\ast}}$ be a sequence of
random variables defined on the probabilized space $\left(\Omega,\mathscr{A},P\right),$
independent and admitting a second-order moment. We suppose that
\[
\mathbb{E}\left(X_{n}\right)\underset{n\to+\infty}{\longrightarrow}m\,\,\,\,\,\,\,\text{and}\,\,\,\,\,\,\,\,\sum^{n}_{j=1}\dfrac{1}{j^{2}}\sigma^{2}_{X_{j}}<+\infty.
\]

Then, the sequence of random variables 
\[
\overline{X}_{n}=\dfrac{1}{n}\sum^{n}_{j=1}X_{j}
\]
 \textbf{converges $P-$almost surely and in $\text{L}^{2}$ to $m.$}

\end{theorem}

\begin{proof}{}{}

By Cesaro lemma, it follows that
\begin{equation}
\mathbb{E}\left(\overline{X}_{n}\right)=\dfrac{1}{n}\sum^{n}_{j=1}\mathbb{E}\left(X_{j}\right)\underset{n\to+\infty}{\longrightarrow}m.\label{eq:overlineX_n_to_m}
\end{equation}

The random variables $Y_{n}=\dfrac{X_{n}-\mathbb{E}\left(X_{n}\right)}{n}$
are independent, centered and have variance $\dfrac{1}{n^{2}}\sigma^{2}_{X_{n}}.$
Thus,
\[
\sum^{+\infty}_{j=1}\sigma^{2}_{Y_{j}}=\sum^{+\infty}_{j=1}\dfrac{1}{j^{2}}\sigma^{2}_{X_{j}}<+\infty.
\]

By Proposition $\ref{pr:sc_p-as_series_indep_rv},$ this proves the
$P-$almost surely convergence of the series with general term $Y_{n}.$
Then Kronecker lemma ensures that the sequence of arithmetic means
of the random variables $nY_{n}$ converges $P-$almost surely to
$0,$ and thus that \textbf{the sequence of random variables $\overline{X}_{n}$
converges $P-$almost surely to $m.$}

For the $\text{L}^{2}$ convergence, we note that, by independence
of the random variables $X_{n}-\mathbb{E}\left(X_{n}\right),$
\[
\left\Vert \dfrac{1}{n}\sum^{n}_{j=1}\left(X_{j}-\mathbb{E}\left(X_{j}\right)\right)\right\Vert ^{2}_{2}=\dfrac{1}{n^{2}}\sum^{n}_{j=1}\sigma^{2}_{X_{j}}.
\]
By the Kronecker lemma and the assumption, it holds that
\[
\lim_{n\to+\infty}\dfrac{1}{n^{2}}\sum^{n}_{j=1}\sigma^{2}_{X_{j}}=0.
\]

Since for every $n\in\mathbb{N}^{\ast},$
\[
\overline{X}_{n}-m=\dfrac{1}{n}\sum^{n}_{j=1}\left(X_{j}-\mathbb{E}\left(X_{j}\right)\right)+\left[\dfrac{1}{n}\sum^{n}_{j=1}\mathbb{E}\left(X_{j}\right)-m\right],
\]
the relation $\refpar{eq:overlineX_n_to_m}$ and the triangle inequality
yield the $\text{L}^{2}-$convergence of the sequence with general
term $\overline{X}_{n}$ to $m.$ 

\end{proof}

As the following example shows, a sequence of random variables may
satisfy the weak law of large numbers without sattisfying the strong
law.

\begin{counterexample}{Following the Weak Law Does Not Imply Following the Strong Law}{}

Let $\left(X_{n}\right)_{n\geqslant2}$ be a sequence of independent
random variables with laws given by
\[
P_{X_{n}}=\dfrac{1}{2n\ln n}\left(\delta_{n}+\delta_{-n}\right)+\left(1-\dfrac{1}{2n\ln n}\right)\delta_{0}.
\]

The random variables $X_{n}$ are centered. If $\dfrac{S_{n}}{n}$
converges $P-$almost surely, then it must converge to $0.$ However,
\[
\sum^{+\infty}_{n=1}P\left(\left|X_{n}\right|\geqslant n\right)=\sum^{+\infty}_{n=1}\dfrac{1}{n\ln n}=+\infty.
\]

The $X_{n}$ are independent, and therefore the Borel-Cantelli lemma
ensures that
\[
P\left(\limsup_{n\to+\infty}\left(\left|X_{n}\right|\geqslant n\right)\right)=1.
\]

Since
\[
\dfrac{X_{n}}{n}=\dfrac{S_{n}}{n}-\dfrac{n-1}{n}\dfrac{S_{n-1}}{n-1},
\]
it follows that, \textbf{$P-$almost surely, the sequence with general
term $\dfrac{S_{n}}{n}$ does not converge to 0.}

Nevertheless, the weak law holds. Indeed,
\[
\sigma^{2}_{X_{n}}=\dfrac{n}{\ln n}.
\]

Since the function $x\mapsto\dfrac{x}{\ln x}$ is nondecreasing on
$\left[\text{e},+\infty\right[,$ we have
\[
\sum^{n}_{k=3}\dfrac{k}{\ln k}\leqslant\sum^{n}_{k=3}\dfrac{n}{\ln n}=\dfrac{n\left(n-2\right)}{\ln n}.
\]

Hence, 
\[
0\leqslant\dfrac{1}{n^{2}}\sum^{n}_{k=2}\sigma^{2}_{X_{n}}\leqslant\dfrac{2}{n^{2}\ln2}+\dfrac{n-2}{n\ln n},
\]
which proves that
\[
\lim_{n\to+\infty}\dfrac{1}{n^{2}}\sum^{n}_{k=2}\sigma^{2}_{X_{n}}=0.
\]

By Theorem $\ref{th:weak_law},$ \textbf{the sequence with general
term $\dfrac{S_{n}}{n}$ converges to $0$ in probability.}

\end{counterexample}

\textbf{If the random variables $X_{n}$ are only integrable}, we
still have a strong law of large numbers provided we add an additional
assumption, namely that the $X_{n}$ are identically distributed.
This is the content of the next theorem.

\begin{theorem}{Kolmogorov-Khinchin Theorem}{kolm-khin-th}

Let $\left(X_{n}\right)_{n\in\mathbb{N}^{\ast}}$ be a sequence of
random variables defined on the probabilized space $\left(\Omega,\mathscr{A},P\right),$
independent and with same law.

Then the two following assertions are equivalent:

(i) There exists a real number $c$ such that the sequence of random
variables 
\[
\overline{X}_{n}=\dfrac{1}{n}\sum^{n}_{j=1}X_{j}
\]
 converges $P-$almost surely to $c.$

(ii) $X_{1}\in\mathscr{L}^{1}\left(\Omega,\mathcal{A},P\right).$ 

If the assertion (i) holds, then $c=\mathbb{E}\left(X_{1}\right).$

\end{theorem}

\begin{proof}{}{}

We have shown in Exercise $\ref{exo:exercise10.12}$ that, if the
$X_{n}$ are independent and with same law, then the sequence with
general term $\dfrac{X_{n}}{n}$ converges $P-$almost surely to 0
if and only if $X_{1}$ is integrable.
\begin{itemize}
\item Suppose that the sequence $\left(\overline{X}_{n}\right)_{n\in\mathbb{N}^{\ast}}$
converges $P-$almost surely to $c.$ Since for $n\in\mathbb{N}^{\ast},$
\[
\dfrac{X_{n}}{n}=\overline{X}_{n}-\dfrac{n-1}{n}\overline{X}_{n-1},
\]
the sequence with general term $\dfrac{X_{n}}{n}$ converges $P-$almost
surely to 0, and $X_{1}$ is therefore integrable.
\item Conversely, suppose that $X_{1}$ is integrable. Then the sequence
with general term $\dfrac{X_{n}}{n}$ converges $P-$almost surely
to 0 and 
\[
P\left(\limsup_{n\to+\infty}\left(\left|X_{n}\right|>n\right)\right)=0.
\]
By introducing for every $n$ the random variable $\widetilde{X}_{n}=\boldsymbol{1}_{\left(\left|X_{n}\right|\leqslant n\right)}X_{n},$
it follows that
\[
P\left(\liminf_{n\to+\infty}\left(X_{n}=\widetilde{X}_{n}\right)\right)=1.
\]
Let
\[
S_{n}=\sum^{n}_{j=1}X_{j}\,\,\,\,\text{and}\,\,\,\,\widetilde{S}_{n}=\sum^{n}_{j=1}\widetilde{X}_{j}.
\]
Then
\begin{align*}
P\left(\left(\dfrac{S_{n}}{n}\right)_{n\in\mathbb{N}^{\ast}}\,\text{converges}\right) & =P\left(\left[\left(\dfrac{S_{n}}{n}\right)_{n\in\mathbb{N}^{\ast}}\,\text{converges}\right]\cap\liminf_{n\to+\infty}\left(X_{n}=\widetilde{X}_{n}\right)\right)\\
 & =P\left(\left[\left(\dfrac{\widetilde{S}_{n}}{n}\right)_{n\in\mathbb{N}^{\ast}}\,\text{converges}\right]\cap\liminf_{n\to+\infty}\left(X_{n}=\widetilde{X}_{n}\right)\right).
\end{align*}
Hence
\[
P\left(\left(\dfrac{S_{n}}{n}\right)_{n\in\mathbb{N}^{\ast}}\,\text{converges}\right)=P\left(\left(\dfrac{\widetilde{S}_{n}}{n}\right)_{n\in\mathbb{N}^{\ast}}\,\text{converges}\right).
\]
Thus, it suffices to prove the $P-$almost sure convergence of the
sequence $\left(\dfrac{\widetilde{S}_{n}}{n}\right)_{n\in\mathbb{N}^{\ast}},$
which follows from Theorem $\ref{th:strong_law_large_nb}.$ We now
verify the two conditions on the moments:
\begin{itemize}
\item The $X_{n}$ have the same law, so it is for the $\widetilde{X}_{n};$
these last have thus same expectation and $X_{1}$ being integrable,
\[
\lim_{n\to+\infty}\mathbb{E}\left(\widetilde{X}_{n}\right)=\lim_{n\to+\infty}\mathbb{E}\left(X_{1}\boldsymbol{1}_{\left(\left|X_{1}\right|\leqslant n\right)}\right)=\mathbb{E}\left(X_{1}\right).
\]
\item We have
\[
0\leqslant\sigma^{2}_{\widetilde{X}_{n}}\leqslant\mathbb{E}\left(\widetilde{X}^{2}_{n}\right)=\mathbb{E}\left(X^{2}_{1}\boldsymbol{1}_{\left(\left|X_{1}\right|\leqslant n\right)}\right),
\]
and, by the Beppo Levi property,
\begin{equation}
0\leqslant\sum^{+\infty}_{n=1}\dfrac{1}{n^{2}}\sigma^{2}_{\widetilde{X}_{n}}\leqslant\mathbb{E}\left(\sum^{+\infty}_{n=1}\dfrac{1}{n^{2}}X^{2}_{1}\boldsymbol{1}_{\left(\left|X_{1}\right|\leqslant n\right)}\right).\label{eq:sig2_X_til_over_n2}
\end{equation}
But, since 
\[
X^{2}_{1}\boldsymbol{1}_{\left(\left|X_{1}\right|\leqslant n\right)}=X^{2}_{1}\sum^{n}_{m=1}\boldsymbol{1}_{\left(m-1\leqslant\left|X_{1}\right|\leqslant m\right)},
\]
we have in $\overline{\mathbb{R}^{+}},$
\begin{align*}
\sum^{+\infty}_{n=1}\dfrac{1}{n^{2}}X^{2}_{1}\boldsymbol{1}_{\left(\left|X_{1}\right|\leqslant n\right)} & =\sum^{+\infty}_{n=1}\dfrac{1}{n^{2}}X^{2}_{1}\left[\sum^{n}_{m=1}\boldsymbol{1}_{\left(m-1\leqslant\left|X_{1}\right|\leqslant m\right)}\right]\\
 & =\sum^{+\infty}_{m=1}\left|X_{1}\right|\left[\sum^{+\infty}_{n=m}\left|X_{1}\right|\dfrac{1}{n^{2}}\boldsymbol{1}_{\left(m-1\leqslant\left|X_{1}\right|\leqslant m\right)}\right],
\end{align*}
which gives the inequality
\[
\sum^{+\infty}_{n=1}\dfrac{1}{n^{2}}X^{2}_{1}\boldsymbol{1}_{\left(\left|X_{1}\right|\leqslant n\right)}\leqslant\sum^{+\infty}_{m=1}\left|X_{1}\right|\boldsymbol{1}_{\left(m-1\leqslant\left|X_{1}\right|\leqslant m\right)}\left[m\sum^{+\infty}_{n=m}\dfrac{1}{n^{2}}\right].
\]
Moreover, we have the upper-bound
\begin{align*}
\sum^{+\infty}_{n=m}\dfrac{1}{n^{2}} & =\dfrac{1}{m^{2}}+\sum^{+\infty}_{n=m}\dfrac{1}{\left(n+1\right)^{2}}\\
 & \leqslant\dfrac{1}{m^{2}}+\sum^{+\infty}_{n=m}\intop^{n+1}_{n}\dfrac{1}{x^{2}}\text{d}x=\dfrac{1}{m^{2}}+\intop^{+\infty}_{m}\dfrac{1}{x^{2}}\text{d}x.
\end{align*}
Hence
\[
m\sum^{+\infty}_{n=m}\dfrac{1}{n^{2}}\leqslant\dfrac{1}{m}+1\leqslant2.
\]
It follows the inequality
\[
\sum^{+\infty}_{n=1}\dfrac{1}{n^{2}}X^{2}_{1}\boldsymbol{1}_{\left(\left|X_{1}\right|\leqslant n\right)}\leqslant2\sum^{+\infty}_{m=1}\left|X_{1}\right|\boldsymbol{1}_{\left(m-1\leqslant\left|X_{1}\right|\leqslant m\right)}=2\left|X_{1}\right|.
\]
We then have, by the inequality $\refpar{eq:sig2_X_til_over_n2}$,
\[
\sum^{+\infty}_{n=1}\dfrac{1}{n^{2}}\sigma^{2}_{\widetilde{X}_{n}}\leqslant\mathbb{E}\left(\sum^{+\infty}_{n=1}\dfrac{1}{n^{2}}X^{2}_{1}\boldsymbol{1}_{\left(\left|X_{1}\right|\leqslant n\right)}\right)\leqslant2\mathbb{E}\left(\left|X_{1}\right|\right)<+\infty.
\]
We have then shown the convergence $P-$almost sure of the sequence
$\left(\dfrac{\widetilde{S}_{n}}{n}\right)_{n\in\mathbb{N}^{\ast}}$
and, thus also of the sequence $\left(\overline{X}_{n}\right)_{n\in\mathbb{N}^{\ast}}$
to $\mathbb{E}\left(X_{1}\right).$
\end{itemize}
\end{itemize}
\end{proof}

An important application of the previous theorem is the \textbf{fundamental
theorem in statistics} concerning the convergence of \textbf{empirical
cumulative distribution functions}.

Let $X$ be a real-valued random variable of law $\mu.$ Let $\left(X_{n}\right)_{n\in\mathbb{N}^{\ast}}$
be a sequence of real-valued random variables, independent and with
same law $\mu$ and with cumulative distribution function $F.$

\begin{definition}{Sample. Empirical Cumulative Distribution Function}{}

The random vector $\left(X_{1},X_{2},\cdots,X_{n}\right)$ is called
\textbf{\mindex{sample ! size}sample of size $n$} of $X.$

The function $F_{n}$ from $\mathbb{R}\times\Omega$ to $\left[0,1\right]$
defined by
\[
\forall\left(x,\omega\right)\in\mathbb{R}\times\Omega,\,\,\,\,\,\,\,F_{n}\left(x,\omega\right)=\dfrac{1}{n}\sum^{n}_{j=1}\boldsymbol{1}_{\left(X_{j}\leqslant x\right)}\left(\omega\right)
\]
is called \textbf{empirical cumulative distribution function}\mindex{cumulative distribution function!empirical}---associated
with $X$---based on the sample $\left(X_{1},X_{2},\cdots,X_{n}\right).$

\end{definition}

\begin{remark}{}{}

For a realization $\omega,$ the vector $\left(X_{1}\left(\omega\right),X_{2}\left(\omega\right),\cdots,X_{n}\left(\omega\right)\right)$
is called the \textbf{\mindex{sample!empirical}empirical sample}. 

For every real number $x,$ the quantity $nF_{n}\left(x,\omega\right)$
is the number of indices such that $X_{k}\left(\omega\right)\leqslant x.$

\end{remark}

\begin{theorem}{Fundamental Theorem of  Statistics, or Glivenko-Cantelli Theorem}{gli-canth}

With the previous notation, for $P-$almost every $\omega,$ the sequence
of cumulative distribution function $F_{n}\left(\cdot,\omega\right)$
converges uniformly to $F,$ that is, we have 
\[
P-\text{almost surely}\,\,\,\,\lim_{n\to+\infty}\sup_{x\in\mathbb{R}}\left|F_{n}\left(x,\cdot\right)-F\left(x\right)\right|=0.
\]

\end{theorem}

\begin{proof}{}{}

Since for every $\omega,$ $F_{n}\left(\cdot,\omega\right)$ and $F$
are right-continuous, we have
\[
\sup_{x\in\mathbb{R}}\left|F_{n}\left(x,\cdot\right)-F\left(x\right)\right|=\sup_{x\in\mathbb{Q}}\left|F_{n}\left(x,\cdot\right)-F\left(x\right)\right|,
\]
which shows that $\sup_{x\in\mathbb{R}}\left|F_{n}\left(x,\cdot\right)-F\left(x\right)\right|$
is indeed a random variable.

For every real number $x,$ the sequence $\left(\boldsymbol{1}_{\left(X_{j}\leqslant x\right)}\right)_{j\in\mathbb{N}^{\ast}}$
---the sequence $\left(\boldsymbol{1}_{\left(X_{j}<x\right)}\right)_{j\in\mathbb{N}^{\ast}}$---is
a sequence of independent random variables with same law and integrable.

Moreover\footnotemark,
\[
\mathbb{E}\left(\boldsymbol{1}_{\left(X_{j}\leqslant x\right)}\right)=P\left(X_{j}\leqslant x\right)=F\left(x\right)\,\,\,\,\text{and}\,\,\,\,\mathbb{E}\left(\boldsymbol{1}_{\left(X_{j}<x\right)}\right)=P\left(X_{j}<x\right)=F\left(x-\right).
\]
By Theorem $\ref{th:kolm-khin-th},$ it follows that
\begin{equation}
\lim_{n\to+\infty}F_{n}\left(x,\cdot\right)\overset{P-\text{a.s.}}{=}F\left(x\right)\,\,\,\,\text{and}\,\,\,\,\lim_{n\to+\infty}F_{n}\left(x-,\cdot\right)\overset{P-\text{a.s.}}{=}F\left(x-\right)\label{eq:P-as-lim_F_nxdot}
\end{equation}
Let $D$ be the union of the set of rationals $\mathbb{Q}$ and of
the set---countable and possibly empty---of discontituity points
of $F.$ This set $D$ is countable and dense. By assertion $\refpar{eq:P-as-lim_F_nxdot},$
there exist, for every $x\in D,$ two sets of zero probability, $N^{1}_{x}$
and $N^{2}_{x},$ such that
\[
\forall\omega\notin N^{1}_{x},\,\,\,\,\lim_{n\to+\infty}F_{n}\left(x,\omega\right)=F\left(x\right)\,\,\,\,\text{and}\,\,\,\,\forall\omega\notin N^{2}_{x},\,\,\,\,\lim_{n\to+\infty}F_{n}\left(x-,\omega\right)=F\left(x-\right)
\]
Let 
\[
N=\left(\bigcup_{x\in D}N^{1}_{x}\right)\cup\left(\bigcup_{x\in D}N^{2}_{x}\right)
\]
Since $D$ is countable, $N$ also has probability zero, and
\[
\forall\omega\notin N,\,\forall x\in D,\,\,\,\,\lim_{n\to+\infty}F_{n}\left(x,\omega\right)=F\left(x\right)\,\,\,\,\text{and}\,\,\,\,\lim_{n\to+\infty}F_{n}\left(x-,\omega\right)=F\left(x-\right).
\]
By Lemma $\ref{lm:cum_dist_unif_cv},$ we obtain
\[
\forall\omega\notin N,\,\,\,\,\lim_{n\to+\infty}\sup_{x\in\mathbb{R}}\left|F_{n}\left(x,\omega\right)-F\left(x\right)\right|=0,
\]
which is exactly the result stated in the theorem.

\end{proof}

\footnotetext{Tr.N: $F(x-):=\lim_{\begin{array}{c}
t\to x\\
t<x
\end{array}}F\left(t\right)$}

For completeness, we just now need to enounce and prove the lemma.

\begin{lemma}{}{cum_dist_unif_cv}

Let $f$ and $f_{n},\,n\in\mathbb{N}^{\ast}$ be nonnegative, nondecreasing
functions defined on $\mathbb{R}$ bounded by $1.$

(a) If the sequence $\left(f_{n}\right)_{n\in\mathbb{N}^{\ast}}$
pointwise converges on a countable and dense set $D$ of $\mathbb{R},$
this sequence converges simply on the set of points of continuity
of $f.$

(b) Suppose moreover that the functions $f$ and $f_{n}$ are cumulative
distribution functions. Let $D$ be the union of the set of rationals
$\mathbb{Q}$ and of the set of discontituity points of $f.$ If
\[
\forall x\in D,\,\,\,\,\lim_{n\to+\infty}f_{n}\left(x\right)=f\left(x\right)\,\,\,\,\text{and}\,\,\,\,\lim_{n\to+\infty}f_{n}\left(x-\right)=f\left(x-\right),
\]
then the sequence $\left(f_{n}\right)_{n\in\mathbb{N}^{\ast}}$ converges
to $f$ uniformly on $\mathbb{R}.$

\end{lemma}

\begin{proof}{}{}

(a) Let $x$ be a point of continuity of $f.$ 

Let $\epsilon>0$ be arbitrary and choose $\eta>0$ such that
\[
x^{\prime}\in\left[x-\eta,x+\eta\right]\Longrightarrow\left|f\left(x\right)-f\left(x^{\prime}\right)\right|\leqslant\epsilon.
\]

Let $y,y^{\prime}\in D$ satisfy $x-\eta<y<x<y^{\prime}<x+\eta.$ 

By the growth of the functions $f$ and $f_{n}$ and the assumption
of sequence convergence, we have
\[
f\left(y\right)=\lim_{n\to+\infty}f_{n}\left(y\right)\leqslant\liminf_{n\to+\infty}f_{n}\left(x\right)\leqslant\limsup_{n\to+\infty}f_{n}\left(x\right)\leqslant\lim_{n\to+\infty}f_{n}\left(y^{\prime}\right)=f\left(y^{\prime}\right).
\]
Thus, for every $\epsilon>0,$
\[
\left|\liminf_{n\to+\infty}f_{n}\left(x\right)-\limsup_{n\to+\infty}f_{n}\left(x\right)\right|\leqslant\epsilon\,\,\,\,\text{and}\,\,\,\,\left|\limsup_{n\to+\infty}f_{n}\left(x\right)-f\left(x\right)\right|\leqslant\epsilon,
\]
Since $\epsilon$ is arbitrary, we conclude that
\[
\liminf_{n\to+\infty}f_{n}\left(x\right)=\limsup_{n\to+\infty}f_{n}\left(x\right)=f\left(x\right).
\]
That is, the sequence $\left(f_{n}\left(x\right)\right)_{n\in\mathbb{N}^{\ast}}$
converges to $f\left(x\right).$

(b) For every integers $j,k$ such that $1\leqslant j\leqslant k,$
define
\[
x_{j,k}=\sup\left\{ x\in\mathbb{R}:\,f\left(x-\right)\leqslant\dfrac{j}{k}\leqslant f\left(x\right)\right\} 
\]
---$\sup\emptyset=+\infty$ by convention---and $x_{0,k}=-\infty.$ 

Since $f$ is a cumulative distribution function, $\lim_{x\to-\infty}f\left(x\right)=0$
and $\lim_{x\to+\infty}f\left(x\right)=1.$ 

Hence $x_{j,k}<x_{j+1,k}$ for $k\in\mathbb{N}^{\ast}$ and $0\leqslant j\leqslant k-1.$ 

Thus, the intervals $\left[x_{j,k},x_{j+1,k}\right]$ form a partition
of $\mathbb{R}.$ 

Define
\[
\Delta^{1}_{n}\left(k\right)=\max_{0\leqslant j\leqslant k-1}\left|f_{n}\left(x_{j,k}\right)-f\left(x_{j,k}\right)\right|,
\]
\[
\Delta^{2}_{n}\left(k\right)=\max_{0\leqslant j\leqslant k-1}\left|f_{n}\left(x_{j,k}-\right)-f\left(x_{j,k}-\right)\right|\,\,\,\,\text{and}\,\,\,\,\Delta_{n}=\sup_{x\in\mathbb{R}}\left|f_{n}\left(x\right)-f\left(x\right)\right|.
\]

For every $k\in\mathbb{N}^{\ast},$
\begin{equation}
\Delta_{n}\leqslant\max\left(\Delta^{1}_{n}\left(k\right),\Delta^{2}_{n}\left(k\right)\right)+\dfrac{1}{k}.\label{eq:upper_bound_delta_n}
\end{equation}

Indeed, if $x\in\left]x_{j,k},x_{j+1,k}\right[,$ then
\[
f\left(x_{j,k}\right)\leqslant f\left(x\right)\leqslant f\left(x_{j+1,k}-\right)\,\,\,\,\,\,\,\,f_{n}\left(x_{j,k}\right)\leqslant f_{n}\left(x\right)\leqslant f_{n}\left(x_{j+1,k}-\right)
\]
\[
\text{and}\,\,\,\,0\leqslant f\left(x_{j+1,k}-\right)-f\left(x_{j,k}\right)\leqslant\dfrac{1}{k}.
\]

Hence,
\[
f_{n}\left(x\right)-f\left(x\right)\leqslant f_{n}\left(x_{j+1,k}-\right)-f\left(x_{j,k}\right)\leqslant f_{n}\left(x_{j+1,k}-\right)-f\left(x_{j+1,k}-\right)+\dfrac{1}{k}
\]
and
\[
f_{n}\left(x\right)-f\left(x\right)\geqslant f_{n}\left(x_{j,k}\right)-f\left(x_{j+1,k}-\right)\geqslant f_{n}\left(x_{j,k}\right)-f\left(x_{j,k}\right)-\dfrac{1}{k},
\]
which shows the inequality $\refpar{eq:upper_bound_delta_n},$ since
for every $k\in\mathbb{N}^{\ast},$
\[
\Delta_{n}=\max_{1\leqslant j\leqslant k}\sup_{x\in\left]x_{j,k},x_{j+1,k}\right]}\left|f_{n}\left(x\right)-f\left(x\right)\right|.
\]

Next, for every $k\in\mathbb{N}^{\ast},$
\[
\lim_{n\to+\infty}\Delta^{1}_{n}\left(k\right)=\lim_{n\to+\infty}\Delta^{2}_{n}\left(k\right)=0.
\]

Indeed,
\[
\lim_{n\to+\infty}f_{n}\left(x_{j,k}\right)=f\left(x_{j,k}\right)\,\,\,\,\text{and}\,\,\,\,\lim_{n\to+\infty}f_{n}\left(x_{j,k}-\right)=f\left(x_{j,k}-\right)=f\left(x_{j,k}\right),
\]
as either $x_{j,k}$ is a continuity point of $f,$ and it follows
from the first part of the lemma, or $x_{j,k}$ is a discontinuity
point of $f,$ and this is the assumption. It is then enough to note
that in $\Delta^{1}_{n}\left(k\right)$ and $\Delta^{2}_{n}\left(k\right)$
there is only a finite number of quantities $\left|f_{n}\left(x_{j,k}\right)-f\left(x_{j,k}\right)\right|$
or $\left|f_{n}\left(x_{j,k}-\right)-f\left(x_{j,k}-\right)\right|.$
It follows that, for all $k\in\mathbb{N}^{\ast},$
\[
0\leqslant\limsup_{n\to+\infty}\Delta_{n}=0,
\]
which shows that the sequence with general term $f_{n}$ converges
to 0.

\end{proof}

The Glivenko-Cantelli theorem---Theorem $\ref{th:gli-canth}$---suggests
the idea behind the \textbf{\mindex{test!Kolmogorov-Smirnov}\index{Kolmogorov-Smirnov test}Kolmogorov-Smirnov
test}. 

With the same notations than in this theorem, the goal is the following:
given a sample of size $n,$ we want to test the hypothesis that the
random variable $X$ has as cumulative distribution function the continuous
function $F.$

The testing method constists of defining an \textbf{acceptance region\index{acceptance region}}
for the hypothesis with a probability of error $\alpha,$ where $\alpha$
is called the \textbf{significance level\index{significance level}}
(or threshold) of the test. 

This test is \textbf{non-parametric}\index{non-parametric test},
in the sense that the hypothesis only requires that $F$ belongs to
a class of functions---here, the continuous functions. 

This is in contrast with a \textbf{parametric test},\textbf{\index{parametric test}\mindex{test! parametric}}
where we assume that $F$ belongs to a family of functions determined
by parameters---for instance, the family of Gaussian laws with parameters
$m$ and $\sigma^{2}$---and where the hypothesis concerns the values
of these parameters.

The test is based on the observation that the random variable $D_{n},$
called \textbf{\index{Kolmogorov-Smirnov statistic}\mindex{statistic!Kolmogorov-Smirnov}Kolmogorov-Smirnov
statistic}, and defined by
\[
D_{n}=\sup_{x\in\mathbb{R}}\left|F_{n}\left(x,\cdot\right)-F\left(x\right)\right|,
\]
has an independent law from $F.$

Let us prove this. 

Let $G$ be the pseudo-inverse of $F,$ that is the function defined
by
\[
G\left(x\right)=\inf\left\{ x\in\mathbb{R}:\,F\left(x\right)\geqslant y\right\} .
\]

Recall---see Exercise $\ref{exo:exercise9.1}$---that $G$ is defined
on $\left[0,1\right],$ and that when $F$ is continuous, we have,
for all $y\in\left[0,1\right],$ $F\left(G\left(y\right)\right)=y.$ 

Moreover, the law of $F\left(X\right)$ is the uniform law on the
interval $\left[0,1\right].$ Finally, since $F$ is continuous, at
each point $x$ of strict increase\footnote{A point $x$ is a \textbf{point of strict growth \mindex{point!strict growth}}
for a function $f,$ if there exists an open interval containing $x$
on which $f$ is strictly increasing.} of $F,$ there is the equivalence 
\[
F\left(x\right)\leqslant F\left(y\right)\,\,\,\,\Leftrightarrow\,\,\,\,x\leqslant y.
\]
Let $C$ be the set of points where $F$ is increasing. Its complement
$C^{c}$ is a countable union of intervals $\left]a_{i},b_{i}\right[,i\in I.$
These invervals correspond to the plateaus of $F,$ which in turn
correspond to the jumps of the nondecreasing function $G.$ 

For $1\leqslant j\leqslant n,$
\[
P\left(X_{j}\in C^{c}\right)\leqslant\sum_{i\in I}P\left(X_{j}\in\left(a_{i},b_{i}\right)\right)=\sum_{i\in I}\left(F\left(b_{i}\right)-F\left(a_{i}\right)\right)=0.
\]

Since the random variables $X_{j}$ are independent, it follows $P-$almost
surely that $\left(X_{1},X_{2},\cdots,X_{n}\right)\in C^{n}.$

For real numbers $x_{1},x_{2},\cdots,x_{n},x,$ define
\[
\nu_{n}\left(x_{1},x_{2},\cdots,x_{n},x\right)=\sum^{n}_{j=1}\boldsymbol{1}_{\left(x_{j}\leqslant x\right)}
\]
the number of $x_{j}$ less than or equal to $x.$ 

We have $P-$almost surely
\[
D_{n}=\sup_{x\in\mathbb{R}}\left|F_{n}\left(x,\cdot\right)-F\left(x\right)\right|=\sup_{x\leqslant G\left(1\right)}\left|F_{n}\left(x,\cdot\right)-F\left(x\right)\right|,
\]
since for $x>G\left(1\right),$ we have $P-$almost surely $F_{n}\left(x,\cdot\right)=F\left(x\right)=1$.

Thus $P-$almost surely
\begin{align*}
D_{n} & =\sup_{x\leqslant G\left(1\right)}\left|\dfrac{1}{n}\nu_{n}\left(X_{1},X_{2},\cdots,X_{n},x\right)-F\left(x\right)\right|\\
 & =\sup_{y\in\left[0,1\right]}\left|\dfrac{1}{n}\nu_{n}\left(X_{1},X_{2},\cdots,X_{n},G\left(y\right)\right)-F\left(G\left(y\right)\right)\right|.
\end{align*}
Therefore, $P-$almost surely,
\[
D_{n}=\sup_{y\in\left[0,1\right]}\left|\dfrac{1}{n}\nu_{n}\left(X_{1},X_{2},\cdots,X_{n},G\left(y\right)\right)-y\right|.
\]

But, since $P-$almost surely, $\left(X_{1},X_{2},\cdots,X_{n}\right)\in C^{n},$
$P-$almost surely
\begin{align*}
\nu_{n}\left(X_{1},X_{2},\cdots,X_{n},G\left(y\right)\right) & =\nu_{n}\left(F\left(X_{1}\right),F\left(X_{2}\right),\cdots,F\left(X_{n}\right),F\left(G\left(y\right)\right)\right)\\
 & =\nu_{n}\left(F\left(X_{1}\right),F\left(X_{2}\right),\cdots,F\left(X_{n}\right),y\right).
\end{align*}
Hence, \boxeq{
\[
D_{n}\overset{P-\text{a.s.}}{=}\sup_{y\in\left[0,1\right]}\left|\dfrac{1}{n}\nu_{n}\left(F\left(X_{1}\right),F\left(X_{2}\right),\cdots,F\left(X_{n}\right),y\right)-y\right|.
\]
}

Since the random variables $F\left(X_{1}\right),F\left(X_{2}\right),\cdots,F\left(X_{n}\right)$
are independent, of same uniform law on $\left[0,1\right],$ we have
shown that $D_{n}$ has a law independent of $F.$ 

This law is tabulated\footnote{For instance, one can find the table of the law of $D_{n}$ in \cite{trivedi2001probability},
---Tr.N: p.599 in the edition of 1982.}. 

Then the test constists, given a significance level $\alpha,$ to
find in the table the value $d_{\alpha}$ such that 
\[
P\left(D_{n}\leqslant d_{\alpha}\right)=1-\alpha.
\]
We \textbf{accept the hypothesis\index{accept the hypothesis}} that
$X$ has for cumulative distribution function $F,$ if 
\[
\forall x\in\mathbb{R},\,\,\,\,\left|F_{n}\left(x,\cdot\right)-F\left(x\right)\right|\leqslant d_{\alpha},
\]
that is, if the graph of $F$ lies entirely within the band bounded
by the graphs of the empirical cumulative distribution function shifted
by $\pm d_{\alpha},$ using the sample $\left(x_{1},x_{2},\cdots,x_{n}\right).$

\begin{example}{Numerical Example}{}

We want to determine whether, at the 0.05 significance level, if we
can accept the hypothesis that the following sample of size 15, obtained
from a random generator, comes from a random variable of uniform law
on $\left[0,1\right]:$
\[
\begin{array}{ccccccccccccccc}
0.8 & 0.4 & 0.25 & 0.7 & 0.6 & 0.2 & 0.5 & 0.3 & 0.15 & 0.1 & 0.65 & 0.9 & 0.45 & 0.85 & 0.55.\end{array}
\]

\end{example}

\begin{solutionexample}{}{}

The empirical cumulative distribution function increases in steps
of $\dfrac{1}{15}$ at the abscissas $x_{i}.$ 

From the table, the critical value at the 0.05 significance level
is $d_{0.05}=0.34.$ 

A graphical representation shows easily that the diagonal lies entirely
within the band delimited by the functions $F_{n}\pm0.34$---a direct
computation leads to the same conclusion.

Hence, with a probability 0.05 of error, we accept the assumption
that the sample comes from a random variable following the uniform
law on $\left[0,1\right].$

\end{solutionexample}

\section*{Exercises}

\addcontentsline{toc}{section}{Exercises}

\textbf{Unless otherwise specified, all random variables are defined
on the same probabilized space $\left(\Omega,\mathcal{A},P\right).$}

\begin{exercise}{Metrics and Convergence in Probability}{exercise11.1}

Let $\mathscr{L}^{0}$---respectively $\text{L}^{0}$---denote the
set of random variables---respectively the set of equivalence classes
of random variables---defined $P-$almost surely, and, in the case
of $\overline{\mathbb{R}},$ $P-$almost surely finite.

Define, for every $X$ and $Y$ of $\mathscr{L}^{0}$---respectively
of $\text{L}^{0}$---the functions
\[
d\left(X,Y\right)=\mathbb{E}\left(\dfrac{\left|X-Y\right|}{1+\left|X-Y\right|}\right)\,\,\,\,\text{and}\,\,\,\,\delta\left(X,Y\right)=\mathbb{E}\left(\min\left(1,\left|X-Y\right|\right)\right).
\]

Prove that $d$ and $\delta$ defined on the vector space $\mathscr{L}^{0}$---respectively
on $\text{L}^{0}$---are two pseudo-metrics\footnotemark---respectively
metrics---equivalent and such that the convergence of the sequences
with respect to these metrics is equivalent to the convergence in
probability. Conclude that the ``metric'' spaces are complete.

\end{exercise}

\footnotetext{A pseudo-metric is also called \textbf{gap}.}

\begin{exercise}{A Partial Converse of Theorem $\ref{th:second_th_sc_cv_seq}$}{exercise11.2}

Prove that if $\left(X_{n}\right)_{n\in\mathbb{N}}$ is a sequence
of independent random variables that converges $P-$almost surely
to 0, then for every $\epsilon>0,$
\[
\sum^{+\infty}_{n=0}P\left(\left(\left|X_{n}\right|>\epsilon\right)\right)<+\infty.
\]

\end{exercise}

\begin{exercise}{Equi-integrability and Convergence $\mathcal{L}^p$ of Sequences of Gaussian Random Variables}{exercise11.3}

Let $X$ be a real-valued random variable with Gaussian law $\mathcal{N}_{\mathbb{R}}\left(m,\sigma^{2}\right).$

1. Prove the inequality
\begin{equation}
\mathbb{E}\left(\text{e}^{\left|X\right|}\right)\leqslant2\text{ch}\left(m\right)\text{e}^{\frac{\sigma^{2}}{2}}.\label{eq:exp_normal_bound}
\end{equation}

2. Let $I$ be an arbitrary index set and, for each $i\in I,$ let
$X_{i}$ be a random variable with Gaussian law $\mathcal{N}_{\mathbb{R}}\left(m_{i},\sigma^{2}_{i}\right).$
Prove that if the families of real numbers $\left(m_{i}\right)_{i\in I}$
and $\left(\sigma^{2}_{i}\right)_{i\in I}$ are bounded, then the
family of random variables $\left(\left|X_{i}\right|^{p}\right)_{i\in I}$
is equi-integrable for every $p\geqslant1.$

3. If moreover $I=\mathbb{N}$ and the sequence $\left(X_{n}\right)_{n\in\mathbb{N}}$
converges in probability to a random variable $X,$ then prove that
$X$ is Gaussian and that the convergence holds in every $\mathscr{L}^{p},\,p\geqslant1.$

\end{exercise}

\begin{exercise}{Necessary Condition to the P-almost Sure Convergence of Series of Random Variables Uniformly Bounded}{exercise11.4}

Let $\left(X_{n}\right)_{n\in\mathbb{N}}$ be a sequence of centered
and independent random variables, uniformly bounded $P-$almost surely
by $c>0.$ Set $S_{n}=\sum^{n}_{k=0}X_{k}.$

For a fixed integer $l\geqslant1,$ define the sets
\[
A=\left\{ \sup_{n\in\mathbb{N}}\left|S_{n}\right|\leqslant l\right\} \,\,\,\,\text{and}\,\,\,\,A_{p}=\left\{ \sup_{0\leqslant n\leqslant p}\left|S_{n}\right|\leqslant l\right\} .
\]

1. Prove the inequality
\begin{equation}
\mathbb{E}\left(\boldsymbol{1}_{A_{p}}S^{2}_{p+1}\right)\geqslant\mathbb{E}\left(\boldsymbol{1}_{A_{p}}S^{2}_{p}\right)+P\left(A\right)\sigma^{2}_{X_{p+1}}.\label{eq:low_bound_exp_ind_Ap_S_p+1^2}
\end{equation}

2. Deduce that the condition $P\left(\sup_{n\in\mathbb{N}}\left|S_{n}\right|<+\infty\right)>0$
implies $\sum^{+\infty}_{n=0}\sigma^{2}_{X_{n}}<+\infty:$ in particular,
this holds whenever the series with general term $X_{n}$ converges
$P-$almost surely.

\begin{remark}{}{}

This is a partial converse of Proposition $\ref{pr:sc_p-as_series_indep_rv}$ 

\end{remark}

\end{exercise}

\begin{exercise}{Kolmogorov Three Series Theorem}{exercise11.5}

Let $\left(X_{n}\right)_{n\in\mathbb{N}}$ be a sequence of independent
real-valued random variables. Define $S_{n}=\sum^{n}_{k=0}X_{k}$
and $Y_{n}=\boldsymbol{1}_{\left(\left|X_{n}\right|\leqslant c\right)}X_{n},$
where $c$ is an arbitrary positive real number.

Prove that the series $\sum X_{n}$ converges $P-$almost surely if
and only if the three series $\sum\mathbb{E}\left(Y_{n}\right),$
$\sum\sigma^{2}_{Y_{n}}$ and $\sum P\left(\left|X_{n}\right|>c\right)$
all converge.

\textit{Hint: For the necessary condition, use symmetrisation to reduce
to centered random variables---see Exercise $\ref{exo:exercise10.13}$
and apply the preceding Exercise $\ref{exo:exercise11.4}.$}

\end{exercise}

\begin{exercise}{There is No Cesàro Lemma for the Convergence in Probability}{exercise11.6}

Let $\left(X_{n}\right)_{n\in\mathbb{N}^{\ast}}$ be a sequence of
independent real-valued random variables, where each $X_{n}$ has
cumulative distribution function $F_{n}$ defined by
\[
F_{n}\left(x\right)=\begin{cases}
0, & \text{if}\,x\leqslant0,\\
1-\dfrac{1}{x+n}, & \text{otherwise.}
\end{cases}
\]
Let $S_{n}=\sum^{n}_{k=1}X_{k}$ and $Y_{n}=\dfrac{S_{n}}{n}.$

Prove that the sequence $\left(X_{n}\right)_{n\in\mathbb{N}^{\ast}}$
converges to 0 in probability, but that the sequence ~$\left(Y_{n}\right)_{n\in\mathbb{N}^{\ast}}$
does not converge to 0 in probability.

\end{exercise}

\begin{exercise}{Theorem $\ref{th:strong_law_large_nb}$ Gives a Sufficient but not Necessary Condition to the Strong Law of Large Numbers}{exercise11.7}

Let $\left(X_{n}\right)_{n\in\mathbb{N}^{\ast}}$ be a sequence of
independent real-valued random variables, each $X_{n}$ being of law
\[
P_{X_{n}}=\dfrac{1}{2}\left(1-\dfrac{1}{2^{n}}\right)\left(\delta_{1}+\delta_{-1}\right)+\dfrac{1}{2^{n+1}}\left(\delta_{2^{n}}+\delta_{-2^{n}}\right).
\]

Prove that 
\[
\sum^{+\infty}_{n=1}\dfrac{1}{n^{2}}\sigma^{2}_{X_{n}}=+\infty
\]
and nevertheless, the sequence $\left(X_{n}\right)_{n\in\mathbb{N}^{\ast}}$
satisfies the strong law of large numbers.

\end{exercise}

\begin{exercise}{An Application of the Strong Law of Large Numbers: The Monte-Carlo Method for Numerical Integration}{exercise11.8}

Let $D$ be a domain of $\mathbb{R}^{d}$ and $f$ a measureable real-valued
function defined on $\mathbb{R}^{d}$ such that $\boldsymbol{1}_{D}\cdot f$
is Lebesgue-integrable.

Let $\left(U_{n}\right)_{n\in\mathbb{N}^{\ast}}$ be a sequence of
independent real-valued random variables, of uniform law on $\left[0,1\right].$
Define for each $n\in\mathbb{N}^{\ast},$ the random variable $\uuline{U_{n}}$
taking values in $\mathbb{R}^{d}$ by 
\[
\uuline{U_{n}}=\left(U_{nd+1},U_{nd+2},\cdots,U_{\left(n+1\right)d}\right)
\]
and the real-valued random variable 
\[
X_{n}=\left(\boldsymbol{1}_{D}\cdot f\right)\circ\uuline{U_{n}}.
\]
 Prove that the sequence with general term
\[
S_{n}=\dfrac{1}{n}\sum^{n}_{i=1}X_{i}
\]
converges $P-$almost surely to the integral $I=\intop_{D\cap\left[0,1\right]^{d}}f\left(x\right)\text{d}x$
and that, if $f$ is bounded by $c>0,$ then for every $\epsilon>0,$
\begin{equation}
P\left(\left|S_{n}-1\right|>\epsilon\right)\leqslant\dfrac{c^{2}}{n\epsilon^{2}}.\label{eq:PabsS_n-1overeps}
\end{equation}

\end{exercise}

\begin{exercise}{Ottaviani Inequality}{exercise11.9}

Let $X_{1},X_{2},\cdots,X_{n}$ be $n$ independent random variables.
For each $k$ such that $1\leqslant k\leqslant n,$ let
\[
S_{k}=\sum^{k}_{i=1}X_{i}\,\,\,\,\text{and}\,\,\,\,M_{k}=\max_{1\leqslant i\leqslant k}\left|S_{i}\right|
\]
and for $0\leqslant k\leqslant n-1,$
\[
S_{k,n}=\sum^{n}_{i=k+1}X_{i}.
\]

For every $\epsilon>0,$ we introduce the sets $E=\left(M_{n}>2\epsilon\right),$
$E_{I}=\left(\left|S_{1}\right|>2\epsilon\right)$ and, if $2\leqslant k\leqslant n,$
\[
E_{k}=\left(\left|S_{k}\right|>2\epsilon\right)\cap\left[\bigcap^{k-1}_{i=1}\left(\left|S_{i}\right|\leqslant2\epsilon\right)\right].
\]

Prove the inequality
\[
P\left(\left|S_{n}\right|>\epsilon\right)\geqslant\sum^{n}_{k=1}P\left(\left(\left|S_{k,n}\right|\leqslant\epsilon\right)\cap E_{k}\right).
\]

Deduce the \textbf{\sindex[fam]{Ottaviani, Giorgio}Ottaviani\footnotemark
inequality
\[
\min_{1\leqslant k\leqslant n}P\left(\left|S_{k,n}\right|\leqslant\epsilon\right)P\left(\left|M_{n}\right|>2\epsilon\right)\leqslant P\left(\left|S_{n}\right|>\epsilon\right).
\]
}

\end{exercise}

\footnotetext{Giorgio Ottaviani is an Italian mathematician, Professor
at University of Florence. He obtained his PhD in 1987. He made several
contribution in particular in algebraic geometry, tensors and their
applications, as well as algebraic statistics.}

\begin{exercise}{Equivalence of the Convergences in Probability and P-Almost Sure of Series of Independent Random Variables---Levy's Theorem---}{exercise11.10}

Let $\left(X_{n}\right)_{n\in\mathbb{N}^{\ast}}$ be a sequence of
independent real-valued random variables.

Prove that if the series $\sum X_{n}$ of general term $X_{n}$ converges
in probability, then, in fact, it converges $P-$almost surely.

\textit{Hint: Use the Ottaviani inequality proved in Exercise $\ref{exo:exercise11.9}.$}

\end{exercise}

\begin{figure}[t]
\begin{center}\includegraphics[width=0.4\textwidth]{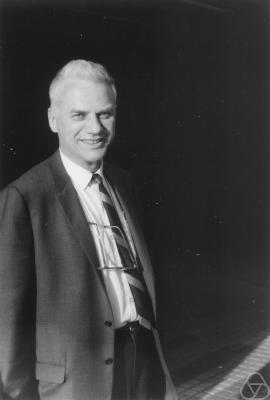}

{\tiny Unknown author---\href{https://mathshistory.st-andrews.ac.uk/Biographies/Hoeffding/}{https://mathshistory.st-andrews.ac.uk/Biographies/Hoeffding/}---Public
domain}\end{center}

\caption{\textbf{\protect\href{https://en.wikipedia.org/wiki/Wassily_Hoeffding}{Wassily Hoeffding}}
(1914 - 1991)}\sindex[fam]{Hoeffding, Wassily}
\end{figure}

\begin{exercise}{Large-Deviation-Type Inequality: Hoeffding Inequality}{exercise11.11}

1. Let $X$ be a real-valued random variable that is $P-$almost surely
bounded by 1, i.e. such that $\left|X\right|\leqslant1$ $P-$almost
surely. Assume moreover that $X$ is centered.

(a) Let $t\in\mathbb{R}.$ Justify the following convexity inequality
\[
\forall x\in\left[-1,1\right],\,\,\,\,\text{e}^{tx}\leqslant\dfrac{1}{2}\left(1-x\right)\text{e}^{-t}+\dfrac{1}{2}\left(1+x\right)\text{e}^{t}.
\]

(b) After justifying the existence of the expectation of the random
variable $\text{e}^{tX},$ deduce the inequality\footnotemark
\[
\mathbb{E}\left(\text{e}^{tX}\right)\leqslant\dfrac{1}{2}\left(\text{e}^{-t}+\text{e}^{t}\right).
\]

Then, prove the inequality
\begin{equation}
\mathbb{E}\left(\text{e}^{tX}\right)\leqslant\text{e}^{\frac{t^{2}}{2}}.\label{eq:bound_Laplace_transform_X_by_exp_t^2over2}
\end{equation}

\textit{Hint: Compare the general term in the series expansion of
both sides.}

2. Let $\left(X_{n}\right)_{n\in\mathbb{N}^{\ast}}$ be a sequence
of independent, centered, real-valued random variables, $P-$almost
surely bounded. Assume that $\left|X_{n}\right|\leqslant c_{n}$ $P-$almost
surely, with $c_{n}>0.$ For each $n\in\mathbb{N}^{\ast},$ set $S_{n}=\sum^{n}_{j=1}X_{j}.$

(a) Prove that, for every $t\in\mathbb{R},$
\[
\mathbb{E}\left(\text{e}^{tS_{n}}\right)\leqslant\text{e}^{\frac{t^{2}}{2}\sum^{n}_{j=1}c^{2}_{j}}
\]

(b) Deduce from the Markov inequality, that for every $t>0$ and $\epsilon>0,$
\begin{equation}
P\left(S_{n}>\epsilon\right)\leqslant\text{e}^{-t\epsilon+\frac{t^{2}}{2}\sum^{n}_{j=1}c^{2}_{j}}.\label{eq:P_S_noverepsilon_first_bound}
\end{equation}

(c) By minimising in $t$ the right-hand side of the inequality $\refpar{eq:P_S_noverepsilon_first_bound},$
deduce that for every $\epsilon>0,$
\begin{equation}
P\left(S_{n}>\epsilon\right)\leqslant\text{e}^{-\frac{\epsilon^{2}}{2\sum^{n}_{j=1}c^{2}_{j}}}.\label{eq:P_S_nnoverepsilon_min_bound}
\end{equation}

(d) Prove that for every $\epsilon>0,$ the \textbf{\sindex[fam]{Hoeffding, Wassily}Hoeffding\footnotemark}
inequality holds:
\begin{equation}
P\left(\left|S_{n}\right|>\epsilon\right)\leqslant2\text{e}^{-\frac{\epsilon^{2}}{2\sum^{n}_{j=1}c^{2}_{j}}}.\label{eq:Hoeffing inequality}
\end{equation}

(e) Let $\alpha,\beta>0.$ Suppose that the sequence with general
term $c_{n}$ satisfies
\[
\sum^{n}_{j=1}c^{2}_{j}\leqslant n^{2\alpha-\beta}.
\]

Prove that, for every $\epsilon>0,$ the series with general term
$P\left(\left|S_{n}\right|>n^{\alpha}\epsilon\right)$ converges.
Deduce that
\[
P\left(\bigcup_{\epsilon\in\mathbb{Q}^{+\ast}}\limsup_{n\to+\infty}\left(\left|S_{n}\right|>n^{\alpha}\epsilon\right)=0\right).
\]

What can then be concluded about the $P-$almost sure convergence
of the sequence with general term $n^{-\alpha}S_{n}?$

3. Now assume additionally that the $X_{n}$ follow the same triangle
law, or more precisely they admit the density $g$ defined by
\[
g\left(x\right)=\boldsymbol{1}_{\left[0,1\right]}\left(\left|x\right|\right)\left(1-\left|x\right|\right).
\]

a. Compute, for every real number $t,$ $\mathbb{E}\left(\text{e}^{tX_{1}}\right).$

b. Prove that the function $\phi:t\mapsto\mathbb{E}\left(\text{e}^{tX_{1}}\right)$
is infinitely differentiable and that for every $k\in\mathbb{N}^{\ast},$
\[
\phi^{\left(k\right)}\left(0\right)=\mathbb{E}\left(X^{k}_{1}\right).
\]

Deduce the variance of $S_{n}.$

c. By using the results of the second question, prove that for every
$\alpha>\dfrac{1}{2},$
\[
\lim_{n\to+\infty}n^{-\alpha}S_{n}=0\,\,\,\,P-\text{almost surely.}
\]

\end{exercise}

\addtocounter{footnote}{-1}

\footnotetext{The function $t\mapsto\mathbb{E}\left(\text{e}^{tX}\right)$
is called the Laplace transform or generating function of the random
variable $X.$

}

\stepcounter{footnote}

\footnotetext{\textbf{\href{https://en.wikipedia.org/wiki/Wassily_Hoeffding}{Wassily Hoeffding}\sindex[fam]{Hoeffding, Wassily}}
(1914 - 1991) was a Finnish and American statistician and probabilist.
He is one of the founders of the nonparametric statistics. He let
his name to the Hoeffding inequality, a concentration inequality,
and to the independence test of Hoeffding. He worked at North Carolina
University, first as research assistant, and obtain the title of Professor
in 1973.}

\section*{Solutions of Exercises}

\addcontentsline{toc}{section}{Solutions of Exercises}

\begin{solution}{}{solexercise11.1}

The function $x\mapsto\dfrac{x}{1+x}$ is increasing on $\mathbb{R}^{+}.$
By the triangle inequality, for all nonnegative real numbers $x,y$
and $z,$
\[
\dfrac{\left|x-y\right|}{1+\left|x-y\right|}\leqslant\dfrac{\left|x-z\right|+\left|z-y\right|}{1+\left|x-z\right|+\left|z-y\right|}\leqslant\dfrac{\left|x-z\right|}{1+\left|x-z\right|}+\dfrac{\left|z-y\right|}{1+\left|z-y\right|}.
\]
The triangle inequality for $d$ follows from the nondecreasing monotonicity
of the integral. Moreover, $d\left(X,Y\right)=0$ if and only if $X=Y$
$P-$almost surely, and $d$ is symmetric, hence $d$ is a pseudo-metric
on $\mathscr{L}^{0}.$

By the triangle inequality, for all nonnegative real numbers $x,y$
and $z,$
\[
\min\left(1,\left|x-y\right|\right)\leqslant\min\left(1,\left|x-z\right|+\left|z-y\right|\right)\leqslant\min\left(1,\left|x-z\right|\right)+\min\left(1,\left|z-y\right|\right).
\]
The triangle inequality for $\delta$ again follows from the nondecreasing
monotonicity of the integral. Moreover, $\delta\left(X,Y\right)=0$
if and only if $X=Y$ $P-$almost surely, and $\delta$ is symmetric,
so $\delta$ is a pseudo-metric on $\mathscr{L}^{0}.$

We easily check that for all $x\geqslant0,$
\[
\dfrac{1}{2}\min\left(1,x\right)\leqslant\dfrac{x}{1+x}\leqslant\min\left(1,x\right),
\]
 Hence,
\[
\dfrac{1}{2}\delta\left(X,Y\right)\leqslant d\left(X,Y\right)\leqslant\delta\left(X,Y\right).
\]

Thus, the metrics $d$ and $\delta$ are \textbf{equivalent}.

Since the function $x\mapsto\dfrac{x}{1+x}$ is increasing on $\mathbb{R}^{+}$
and bounded by 1, we obtain for every $\epsilon>0,$
\[
\dfrac{\epsilon}{1+\epsilon}P\left(\left|X-Y\right|>\epsilon\right)\leqslant d\left(X,Y\right)\leqslant\dfrac{\epsilon}{1+\epsilon}\mathbb{E}\left(\boldsymbol{1}_{\left(\left|X-Y\right|\leqslant\epsilon\right)}\right)+\mathbb{E}\left(\boldsymbol{1}_{\left(\left|X-Y\right|>\epsilon\right)}\right),
\]
which yields
\[
\dfrac{\epsilon}{1+\epsilon}P\left(\left|X-Y\right|>\epsilon\right)\leqslant d\left(X,Y\right)\leqslant\epsilon+P\left(\left|X-Y\right|>\epsilon\right).
\]
Hence, the convergence of a sequence with respect to the metric $d$
is equivalent to the convergence in probability. 

Since a sequence converges in probability if and only if it is Cauchy
with respect to convergence in probability, it follows that the \textbf{pseudo-metric
space} $\left(\mathscr{L}^{0},d\right)$ and therefore also $\left(\mathscr{L}^{0},\delta\right),$
\textbf{is} \textbf{complete}.

\begin{remark}{}{}

We could have proceeded in the reverse order: first prove that $\left(\mathscr{L}^{0},d\right)$
is complete, and then deduce that a sequence converges in probability
if and only if it is Cauchy for the convergence in probability.

Let us prove directly that $\left(\mathscr{L}^{0},d\right)$ is complete.

If $\left(X_{n}\right)_{n\in\mathbb{N}}$ is a Cauchy sequence for
$d,$ we may extract a subsequence $\left(X_{n_{k}}\right)_{k\in\mathbb{N}}$
such that
\[
\sum^{+\infty}_{k=0}d\left(X_{n_{k}},X_{n_{k+1}}\right)<+\infty.
\]

By the Beppo-Levi property,
\[
\mathbb{E}\left(\sum^{+\infty}_{k=0}\dfrac{\left|X_{n_{k}}-X_{n_{k+1}}\right|}{1+\left|X_{n_{k}}-X_{n_{k+1}}\right|}\right)<+\infty.
\]
Hence, $P-$almost surely,
\[
\sum^{+\infty}_{k=0}\dfrac{\left|X_{n_{k}}-X_{n_{k+1}}\right|}{1+\left|X_{n_{k}}-X_{n_{k+1}}\right|}<+\infty.
\]
However, since $P-$almost surely 
\[
\lim_{k\to+\infty}\left|X_{n_{k}}-X_{n_{k+1}}\right|=0,
\]
we also have $P-$almost surely
\[
\sum^{+\infty}_{k=0}\left|X_{n_{k}}-X_{n_{k+1}}\right|<+\infty.
\]
Consequently, the subsequence $\left(X_{n_{k}}\right)_{k\in\mathbb{N}}$
converges $P-$almost surely and, by the dominated convergence theorem,
it also converges with respect to the metric $d.$ We conclude since
any Cauchy sequence that admits a convergent subsequence is itself
convergent. 

\end{remark}

\end{solution}

\begin{solution}{}{solexercise11.2}

The set
\[
C=\bigcap_{\epsilon>0}\liminf_{n\to+\infty}\left(\left|X_{n}\right|\leqslant\epsilon\right)
\]
is the set of all $\omega$ such that the sequence with general term
$X_{n}\left(\omega\right)$ converges to 0. It has probability 1.

Hence, for every $\epsilon>0,$
\[
P\left(\liminf_{n\to+\infty}\left(\left|X_{n}\right|\leqslant\epsilon\right)\right)=1,
\]
so
\[
P\left(\limsup_{n\to+\infty}\left(\left|X_{n}\right|>\epsilon\right)\right)=0.
\]
Since the $X_{n}$ are independent, the Borel-Cantelli lemma ensures
that for every $\epsilon>0,$
\[
\sum^{+\infty}_{n=0}P\left(\left(\left|X_{n}\right|>\epsilon\right)\right)<+\infty.
\]

\end{solution}

\begin{solution}{}{solexercise11.3}

\textbf{1. Proof of $\mathbb{E}\left(\text{e}^{\left|X\right|}\right)\leqslant2\text{ch}\left(m\right)\text{e}^{\frac{\sigma^{2}}{2}}$}

We have
\begin{align*}
\mathbb{E}\left(\text{e}^{X}\right) & =\intop_{\mathbb{R}}\text{e}^{x}\dfrac{1}{\sigma\sqrt{2\pi}}\text{e}^{-\frac{\left(x-m\right)^{2}}{2\sigma^{2}}}\text{d}x\\
 & =\text{e}^{m+\frac{\sigma^{2}}{2}}\intop_{\mathbb{R}}\dfrac{1}{\sigma\sqrt{2\pi}}\text{e}^{-\frac{1}{2\sigma^{2}}\left(x-\left(m+\sigma^{2}\right)\right)^{2}}\text{d}x.
\end{align*}

Thus
\[
\mathbb{E}\left(\text{e}^{X}\right)=\text{e}^{m+\frac{\sigma^{2}}{2}}.
\]

Since $-X$ is of law $\mathcal{N}_{\mathbb{R}}\left(-m,\sigma^{2}\right),$
it follows that
\[
\mathbb{E}\left(\text{e}^{\left|X\right|}\right)\leqslant\mathbb{E}\left(\text{e}^{X}\right)+\mathbb{E}\left(\text{e}^{-X}\right)\leqslant\text{e}^{m+\frac{\sigma^{2}}{2}}+\text{e}^{-m+\frac{\sigma^{2}}{2}},
\]
 which yields the formula $\refpar{eq:exp_normal_bound}.$

\textbf{2. $\left(\left|X_{i}\right|^{p}\right)_{i\in I}$ is equi-integrable
for every $p\geqslant1$ }

Let $p\geqslant1.$ There exists $M>0$ such that, for every $x\geqslant M,$
$\left|x\right|^{p}\leqslant\text{e}^{\frac{x}{2}}.$ Hence, for every
$i\in I,$
\[
\left|X_{i}\right|^{p}\leqslant M^{p}\boldsymbol{1}_{\left(\left|X_{i}\right|\leqslant M\right)}+\boldsymbol{1}_{\left(\left|X_{i}\right|>M\right)}\text{e}^{\frac{\left|X_{i}\right|}{2}}.
\]

Thus, for every $A\in\mathcal{A},$ by the Schwarz inequality,
\[
\intop_{A}\left|X_{i}\right|^{p}\text{d}P\leqslant M^{p}P\left(A\right)+\left(P\left(A\right)\right)^{\frac{1}{2}}\mathbb{E}\left(\text{e}^{\left|X_{i}\right|}\right).
\]

By $\refpar{eq:exp_normal_bound},$ by the boundedness of the real
number families $\left(m_{i}\right)_{i\in I}$ and $\left(\sigma^{2}_{i}\right)_{i\in I}$
and the fact that $\sup_{i\in I}\mathbb{E}\left(\text{e}^{\left|X_{i}\right|}\right)<+\infty,$
the family \textbf{$\left(\left|X_{i}\right|^{p}\right)_{i\in I}$
}is equi-integrable.

\textbf{3. $X$ is Gaussian, convergence in all $\mathscr{L}^{p},\,p\geqslant1.$}

The sequence $\left(X^{p}_{n}\right)_{n\in\mathbb{N}}$ is then equi-integrable
for every $p\geqslant1.$ Hence, since the sequence $\left(X_{n}\right)_{n\in\mathbb{N}}$
converges in probability to $X,$ it also converges in every $\mathscr{L}^{p}.$ 

In particular, by taking $p=1,2,$ it follows that the sequences of
general term $m_{n}$ and $\sigma_{n}$ converge and that 
\[
\lim_{n\to+\infty}m_{n}=\mathbb{E}\left(X\right)\,\,\,\,\,\,\,\text{and}\,\,\,\,\,\,\,\,\lim_{n\to+\infty}\sigma^{2}_{n}=\sigma^{2}_{X}.
\]

Now let $f\in\mathscr{C}^{+}_{\mathscr{K}}\left(\mathbb{R}\right)$
arbitrary. Then
\[
\mathbb{E}\left(f\left(X_{n}\right)\right)=\intop_{\mathbb{R}}f\left(x\right)\dfrac{1}{\sigma_{n}\sqrt{2\pi}}\text{e}^{-\frac{\left(x-m_{n}\right)^{2}}{2\sigma^{2}_{n}}}\text{d}x.
\]

By the change of variables $y=\dfrac{x-m_{n}}{\sigma_{n}},$
\[
\mathbb{E}\left(f\left(X_{n}\right)\right)=\intop_{\mathbb{R}}f\left(\sigma_{n}y+m_{n}\right)\dfrac{1}{\sqrt{2\pi}}\text{e}^{-\frac{y^{2}}{2}}\text{d}y.
\]
Since $f$ is bounded, the dominated convergence theorem ensures that
\[
\lim_{n\to+\infty}\mathbb{E}\left(f\left(X_{n}\right)\right)=\intop_{\mathbb{R}}f\left(\sigma y+m\right)\dfrac{1}{\sqrt{2\pi}}\text{e}^{-\frac{y^{2}}{2}}\text{d}y.
\]
By the change of variables $y=\dfrac{x-m}{\sigma},$
\[
\lim_{n\to+\infty}\mathbb{E}\left(f\left(X_{n}\right)\right)=\intop_{\mathbb{R}}f\left(x\right)\dfrac{1}{\sigma\sqrt{2\pi}}\text{e}^{-\frac{\left(x-m\right)^{2}}{2\sigma^{2}}}\text{d}x.
\]

Since $f$ is continuous and bounded, the sequence with general term
$f\left(X_{n}\right),$ on the one hand converges in probability to
$f\left(X\right),$ and on the other hand, is equi-integrable. It
therefore also converges in $\mathscr{L}^{1},$ which shows that $\lim_{n\to+\infty}\mathbb{E}\left(f\left(X_{n}\right)\right)=\mathbb{E}\left(f\left(X\right)\right),$
and thus
\[
\mathbb{E}\left(f\left(X\right)\right)=\intop_{\mathbb{R}}f\left(x\right)\dfrac{1}{\sigma\sqrt{2\pi}}\text{e}^{-\frac{\left(x-m\right)^{2}}{2\sigma^{2}}}\text{d}x.
\]

As $f$ was arbitrary, the law of $X$ is $\mathcal{N}_{\mathbb{R}}\left(m,\sigma^{2}\right).$

\end{solution}

\begin{solution}{}{solexercise11.4}

\textbf{1. Proof of $\mathbb{E}\left(\boldsymbol{1}_{A_{p}}S^{2}_{p+1}\right)\geqslant\mathbb{E}\left(\boldsymbol{1}_{A_{p}}S^{2}_{p}\right)+P\left(A\right)\sigma^{2}_{X_{p+1}}$}

We have
\begin{align*}
\mathbb{E}\left(\boldsymbol{1}_{A_{p}}S^{2}_{p+1}\right) & =\mathbb{E}\left(\boldsymbol{1}_{A_{p}}\left(S_{p}+X_{p+1}\right)^{2}\right)\\
 & =\mathbb{E}\left(\boldsymbol{1}_{A_{p}}S^{2}_{p}\right)+2\mathbb{E}\left(\boldsymbol{1}_{A_{p}}S_{p}X_{p+1}\right)+\mathbb{E}\left(\boldsymbol{1}_{A_{p}}X^{2}_{p+1}\right).
\end{align*}

Since the random variables $\boldsymbol{1}_{A_{p}}S_{p}$ and $X_{p+1}$
are independent,
\[
\mathbb{E}\left(\boldsymbol{1}_{A_{p}}S^{2}_{p+1}\right)=\mathbb{E}\left(\boldsymbol{1}_{A_{p}}S^{2}_{p}\right)+2\mathbb{E}\left(\boldsymbol{1}_{A_{p}}S_{p}\right)\mathbb{E}\left(X_{p+1}\right)+\mathbb{E}\left(\boldsymbol{1}_{A_{p}}X^{2}_{p+1}\right).
\]
Since $X_{p+1}$ is centered and that $\boldsymbol{1}_{A_{p}}$ and
$X_{p+1}$ are independent, we obtain
\[
\mathbb{E}\left(\boldsymbol{1}_{A_{p}}S^{2}_{p+1}\right)=\mathbb{E}\left(\boldsymbol{1}_{A_{p}}S^{2}_{p}\right)+\mathbb{E}\left(\boldsymbol{1}_{A_{p}}\right)\mathbb{E}\left(X^{2}_{p+1}\right).
\]

Thus, it remains only to note that $A\subset A_{p}$ to obtain the
inequality $\refpar{eq:low_bound_exp_ind_Ap_S_p+1^2}.$ 

\textbf{2. $P\left(\sup_{n\in\mathbb{N}}\left|S_{n}\right|<+\infty\right)>0$
implies $\sum^{+\infty}_{n=0}\sigma^{2}_{X_{n}}<+\infty$ }

Since the set $\left\{ \sup_{n\in\mathbb{N}}\left|S_{n}\right|<+\infty\right\} $
is the nondecreasing union of the sequence of sets $\left\{ \sup_{n\in\mathbb{N}}\left|S_{n}\right|\leqslant l\right\} ,l\in\mathbb{N}^{\ast},$
\[
P\left(\sup_{n\in\mathbb{N}}\left|S_{n}\right|<+\infty\right)=\lim_{l\to+\infty}P\left(\sup_{n\in\mathbb{N}}\left|S_{n}\right|\leqslant l\right).
\]

Hence, we can choose an integer $l$ such that 
\[
P\left(\sup_{n\in\mathbb{N}}\left|S_{n}\right|\leqslant l\right)>0,
\]
that is, with the previous notations, such that $P\left(A\right)>0.$
It then follows from the inequality $\refpar{eq:low_bound_exp_ind_Ap_S_p+1^2}$
and from the inclusion $A_{p}\supset A_{p+1}$ that
\[
P\left(A\right)\sigma^{2}_{X_{p+1}}\leqslant\mathbb{E}\left(\boldsymbol{1}_{A_{p}\backslash A_{p+1}}S^{2}_{p+1}\right)+\mathbb{E}\left(\boldsymbol{1}_{A_{p+1}}S^{2}_{p+1}\right)-\mathbb{E}\left(\boldsymbol{1}_{A_{p}}S^{2}_{p}\right).
\]
Since the sequence $\left(X_{n}\right)_{n\in\mathbb{N}}$ is $P-$almost
surely uniformly bounded by $c>0,$ we have on $A_{p}\backslash A_{p+1},$
\[
\left|S_{p+1}\right|\leqslant\left|S_{p}\right|+\left|X_{p+1}\right|\leqslant l+c,
\]
which leads to the inequality
\[
P\left(A\right)\sigma^{2}_{X_{p+1}}\leqslant\left(l+c\right)^{2}P\left(A_{p}\backslash A_{p+1}\right)+\mathbb{E}\left(\boldsymbol{1}_{A_{p+1}}S^{2}_{p+1}\right)-\mathbb{E}\left(\boldsymbol{1}_{A_{p}}S^{2}_{p}\right).
\]
By summing term by term, we obtain for every $n\geqslant2,$
\[
P\left(A\right)\sum^{n}_{p=1}\sigma^{2}_{X_{p+1}}\leqslant\left(l+c\right)^{2}+\mathbb{E}\left(\boldsymbol{1}_{A_{n}}S^{2}_{n}\right),
\]
and thus, by the definition of $A_{n},$
\[
P\left(A\right)\sum^{n}_{p=1}\sigma^{2}_{X_{p+1}}\leqslant\left(l+c\right)^{2}+\text{L}^{2}.
\]

Since $P\left(A\right)>0,$ it follows that 
\[
\sum^{+\infty}_{n=0}\sigma^{2}_{X_{n}}<+\infty.
\]

In particular, if the series with general term $X_{n}$ converges
$P-$almost surely, then $P\left(\sup_{n\in\mathbb{N}}\left|S_{n}\right|<+\infty\right)>0$
and thus $\sum^{+\infty}_{n=0}\sigma^{2}_{X_{n}}<+\infty.$

\begin{remark}{}{}

As the following example shows, the boundedness assumption is necessary.
If the $X_{n},\,n\in\mathbb{N}^{\ast}$ are independent and of law
\[
\dfrac{1}{2n^{3}}\left(\delta_{n}+\delta_{-n}\right)+\left(1-\dfrac{1}{n^{3}}\right)\delta_{0},
\]
then they are centered. We have $\mathbb{E}\left(\left|X_{n}\right|\right)=\dfrac{1}{n^{2}}.$
Hence, we have $\mathbb{E}\left(\sum^{+\infty}_{n=1}\left|X_{n}\right|\right)<+\infty$
and thus $P\left(\sum^{+\infty}_{n=0}\left|X_{n}\right|<+\infty\right)=1.$
Moreover, $\sigma^{2}_{X_{n}}=\mathbb{E}\left(X^{2}_{n}\right)=\dfrac{1}{n};$
hence, $\sum^{+\infty}_{n=0}\sigma^{2}_{X_{n}}=+\infty.$ But the
sequence is not $P-$almost surely uniformly bounded for a constant
$c>0!$

\end{remark}

\end{solution}

\begin{solution}{}{solexercise11.5}

\textbf{Necessary condition}

If the series $\sum X_{n}$ converges $P-$almost surely, the sequence
with general term $X_{n}$ converges $P-$almost surely to 0. Thus
\begin{equation}
P\left(\liminf_{n\to+\infty}\left(\left|X_{n}\right|\leqslant c\right)\right)=1,\label{eq:Pliminf_absX_nlessthanc}
\end{equation}
or also 
\[
P\left(\limsup_{n\to+\infty}\left(\left|X_{n}\right|>c\right)\right)=0.
\]

Since the events $\left(\left|X_{n}\right|>c\right)$ are independent,
it follows from the Borel-Cantelli lemma that
\[
\sum^{+\infty}_{n=0}P\left(\left|X_{n}\right|>c\right)<+\infty.
\]

Moreover, the equality $\refpar{eq:Pliminf_absX_nlessthanc}$ may
also be written
\[
P\left(\liminf_{n\to+\infty}\left(X_{n}=Y_{n}\right)\right)=1.
\]
Hence, the series $\sum Y_{n}$ converges $P-$almost surely. 

On the product probabilized space, let $Y^{s}_{n}$ be the symmetrized
of $Y_{n}$---see Exercise $\ref{exo:exercise10.13}.$ The $Y^{s}_{n}$
are independent and centered. Since the series $\sum Y_{n}$ converges
$P-$almost surely, it follows from the Fubini theorem that the series
$\sum Y^{s}_{n}$ converges $P\otimes P-$almost surely. As the random
variables $Y^{s}_{n}$ are bounded by $2c,$ it follows from Exercise
$\ref{exo:exercise11.4}$ that $\sum^{+\infty}_{n=0}\sigma^{2}_{Y^{s}_{n}}<+\infty.$
As $\sigma^{2}_{Y^{s}_{n}}=2\sigma^{2}_{Y_{n}},$ we also have $\sum^{+\infty}_{n=0}\sigma^{2}_{Y_{n}}<+\infty.$

Finally, the independent and centered random variables $\mathring{Y_{n}}=Y_{n}-\mathbb{E}\left(Y_{n}\right)$
satisfy $\sum^{+\infty}_{n=0}\sigma^{2}_{0}<+\infty.$ It follows
that the series $\sum\mathring{Y_{n}}$ converges $P-$almost surely.
Since $\sum Y_{n}$ converges $P-$almost surely, then the series
$\sum\mathbb{E}\left(Y_{n}\right)$ also converges.

\textbf{Sufficient condition}

Suppose that the three series $\sum\mathbb{E}\left(Y_{n}\right),$
$\sum\sigma^{2}_{Y_{n}}$ and $\sum P\left(\left|X_{n}\right|>c\right)$
converge. Then $\sum^{+\infty}_{n=0}\sigma^{2}_{0}<+\infty$ and the
series $\sum\mathring{Y_{n}}$ converges $P-$almost surely. The same
holds for the series $\sum Y_{n}.$ 

Moreover, since
\[
\sum^{+\infty}_{n=0}P\left(\left(X_{n}\neq Y_{n}\right)\right)=\sum^{+\infty}_{n=0}P\left(\left|X_{n}\right|>c\right)<+\infty,
\]
it follows from the Borel-Cantelli lemma that 
\[
P\left(\limsup_{n\to+\infty}\left(X_{n}\neq Y_{n}\right)\right)=0,
\]
which can also be written 
\[
P\left(\liminf_{n\to+\infty}\left(X_{n}=Y_{n}\right)\right)=1.
\]
 The $P-$almost sure convergence of the series $\sum X_{n}$ follows.

\end{solution}

\begin{solution}{}{solexercise11.6}

The random variables $X$ are $P-$almost surely nonnegative. For
every $\epsilon>0,$
\[
P\left(\left|X_{n}\right|>\epsilon\right)=1-F_{n}\left(\epsilon\right)=\dfrac{1}{\epsilon+n},
\]
and thus
\[
\lim_{n\to+\infty}P\left(\left|X_{n}\right|>\epsilon\right)=0,
\]
which means \textbf{the sequence $\left(X_{n}\right)_{n\in\mathbb{N}^{\ast}}$
converges to 0 in probability.}

Moreover, if $M_{n}=\max_{1\leqslant k\leqslant n}X_{k},$ then, since
the $X_{n}$ are $P-$almost surely nonnegative, we have $\dfrac{M_{n}}{n}\leqslant Y_{n},$
$P-$almost surely. Thus, for every $\epsilon>0,$
\[
P\left(\epsilon<\dfrac{M_{n}}{n}\right)\leqslant P\left(\epsilon<Y_{n}\right).
\]

As the $X_{n}$ are independent, we have, for every $x>0,$
\begin{align*}
P\left(M_{n}\leqslant x\right) & =P\left(\bigcap^{n}_{k=1}\left(X_{k}\leqslant x\right)\right)\\
 & =\prod^{n}_{k=1}P\left(\left(X_{k}\leqslant x\right)\right)\\
 & =\prod^{n}_{k=1}\left(1-\dfrac{1}{x+k}\right)\\
 & \leqslant\left(1-\dfrac{1}{x+n}\right)^{n}.
\end{align*}

It follows that
\[
1-\left(1-\dfrac{1}{n\left(\epsilon+1\right)}\right)^{n}\leqslant P\left(\epsilon<\dfrac{M_{n}}{n}\right)\leqslant P\left(\epsilon<Y_{n}\right),
\]
and thus, by taking the inferior limit,
\[
0<1-\exp\left(-\dfrac{1}{\epsilon+1}\right)\leqslant\liminf_{n\to+\infty}P\left(\epsilon<Y_{n}\right),
\]
which proves that the sequence $\left(Y_{n}\right)_{n\in\mathbb{N}^{\ast}}$
does not converge in probability to 0.

\end{solution}

\begin{solution}{}{solexercise11.7}

The law of $X_{n}$ is symmetric, we have $\mathbb{E}\left(X_{n}\right)=0.$
Thus,
\[
\sigma^{2}_{X_{n}}=\mathbb{E}\left(X^{2}_{n}\right)=1-\dfrac{1}{2^{n}}+\left(2^{n}\right)^{2}\dfrac{1}{2^{n}}=1-\dfrac{1}{2^{n}}+2^{n},
\]
which proves that
\[
\sum^{+\infty}_{n=1}\dfrac{1}{n^{2}}\sigma^{2}_{X_{n}}=+\infty.
\]
Moreover, 
\[
Y_{n}=\boldsymbol{1}_{\left(\left|X_{n}\right|\leqslant1\right)}X_{n}.
\]

The random variable $Y_{n}$ takes the values 0 or $\pm1,$ and
\[
P\left(Y_{n}=+1\right)=P\left(X_{n}=\pm1\right)=1-\dfrac{1}{2^{n}},
\]
while
\[
P\left(Y_{n}=0\right)=P\left(\left|X_{n}\right|=2^{n}\right)=\dfrac{1}{2^{n}}.
\]
Thus, 
\[
P\left(X_{n}\neq Y_{n}\right)=P\left(\left|X_{n}\right|>1\right)=\dfrac{1}{2^{n}}
\]
which implies
\[
\sum^{+\infty}_{n=1}P\left(X_{n}\neq Y_{n}\right)<+\infty.
\]

By the Borel-Cantelli lemma,
\[
P\left(\limsup_{n\to+\infty}\left(X_{n}\neq Y_{n}\right)\right)=0
\]
which is equivalent to
\begin{equation}
P\left(\liminf_{n\to+\infty}\left(X_{n}=Y_{n}\right)\right)=1.\label{eq:liminfX_n_is_Y_n}
\end{equation}

Moreover,
\[
\sigma^{2}_{Y_{n}}=\mathbb{E}\left(Y^{2}_{n}\right)=1-\dfrac{1}{2^{n}},
\]
 and therefore
\[
\sum^{+\infty}_{n=1}\dfrac{1}{n^{2}}\sigma^{2}_{Y_{n}}<+\infty.
\]
By Theorem $\ref{th:strong_law_large_nb},$ it follows that the sequence
with general term $\dfrac{1}{n}\sum^{n}_{j=1}Y_{j}$ converges $P-$almost
surely to 0. Using the probability obtained in $\refpar{eq:liminfX_n_is_Y_n},$
the sequence with general term $\dfrac{1}{n}\sum^{n}_{j=1}X_{j}$
also converges $P-$almost surely to 0.

\end{solution}

\begin{solution}{}{solexercise11.8}

The random variables $\uuline{U_{n}}$ are independent; hence, the
random variables $X_{n}$ are also independent. Moreover, the $X_{n}$
follow the same law and, by the transfer theorem together with the
hypothesis that $\boldsymbol{1}_{D}\cdot f$ is Lebesgue integrable,
they admit an expectation. Therefore, the strong law theorem---Theorem
$\ref{th:kolm-khin-th}$---applies. 

It remains to compute the expectation of $X_{1}.$ By the transfer
theorem
\[
\mathbb{E}\left(X_{1}\right)=\intop_{\mathbb{R}^{d}}\left(\boldsymbol{1}_{D}\cdot f\right)\left(x\right)\text{d}P_{\uuline{U_{n}}}\left(x\right),
\]
and, since $\uuline{U_{n}}$ is of uniform law on $\left[0,1\right]^{d},$
it follows
\[
\mathbb{E}\left(X_{1}\right)=\intop_{\mathbb{R}^{d}}\left(\boldsymbol{1}_{D}\cdot f\right)\left(x\right)\boldsymbol{1}_{\left[0,1\right]^{d}}\left(x\right)\text{d}x=I.
\]
Hence,
\[
\dfrac{1}{n}\sum^{n}_{j=1}X_{j}\stackrel[n\to+\infty]{P-\text{a.s.}}{\longrightarrow}\intop_{D\cap\left[0,1\right]^{d}}f\left(x\right)\text{d}x.
\]

Now, assume that $f$ is bounded by $c>0.$ Since $D\cap\left[0,1\right]^{d}$
is a bounded domain, the random variables $X_{n}$ are in $\mathscr{L}^{2}.$
As they have the same law and are independent, 
\[
\mathbb{E}\left(S_{n}\right)=\mathbb{E}\left(X_{1}\right)=1\,\,\,\,\text{and}\,\,\,\,\sigma^{2}_{S_{n}}=\dfrac{\sigma^{2}_{X_{1}}}{n}.
\]
The Chebyshev inequality applied to $S_{n}$ together with the upper-bound
\[
\sigma^{2}_{X_{1}}=\mathbb{E}\left(X^{2}_{1}\right)-\left(\mathbb{E}\left(X_{1}\right)\right)^{2}\leqslant\mathbb{E}\left(X^{2}_{1}\right)\leqslant c^{2}
\]
yields the inequality $\refpar{eq:PabsS_n-1overeps}.$

\begin{remark}{}{}

In dimension 1 and for sufficiently regular functions, this Monte-Carlo
method cannot compete with classical numerical analysis methods. However,
it is particularly useful when the function is very irregular---here
we only require measurability---or when $d\geqslant2.$ The upper-bound
$\refpar{eq:PabsS_n-1overeps}$ can also be improved using a Bernstein-type
inequality.

\end{remark}

\end{solution}

\begin{solution}{}{solexercise11.9}

Let $\epsilon>0.$

Since the sets $E_{k}$ constitute a partition of $E,$
\begin{equation}
\left(\left|S_{n}\right|>\epsilon\right)\supset\left(\left|S_{n}\right|>\epsilon\right)\cap E=\biguplus^{n}_{k=1}\left(\left(\left|S_{n}\right|>\epsilon\right)\cap E_{k}\right)\label{eq:abs_S_n_eps_contain}
\end{equation}

Since $S_{n}=S_{k}+S_{k,n},$
\[
\left|S_{k}\right|>2\epsilon\,\,\,\,\text{and}\,\,\,\,\left|S_{k,n}\right|\leqslant\epsilon\Longrightarrow\left|S_{n}\right|>\epsilon,
\]
because otherwise we would have $\left|S_{k}\right|\leqslant\left|S_{n}\right|+\left|S_{k,n}\right|\leqslant2\epsilon$
which would contradicts $\left|S_{k}\right|>2\epsilon.$

Thus,
\[
\left(\left|S_{n}\right|>\epsilon\right)\cap E_{k}\supset\left(\left|S_{k,n}\right|\leqslant\epsilon\right)\cap E_{k},
\]
and consequently, by $\refpar{eq:abs_S_n_eps_contain}$,
\[
P\left(\left(\left|S_{n}\right|>\epsilon\right)\right)\geqslant\sum^{n}_{k=1}P\left(\left(\left|S_{k,n}\right|\leqslant\epsilon\right)\cap E_{k}\right).
\]

Since the events $\left(\left|S_{k,n}\right|\leqslant\epsilon\right)$
and $E_{k}$ are independent, we have
\[
P\left(\left|S_{n}\right|>\epsilon\right)\geqslant\sum^{n}_{k=1}P\left(\left|S_{k,n}\right|\leqslant\epsilon\right)P\left(E_{k}\right)\geqslant\min_{1\leqslant k\leqslant n}P\left(\left|S_{k,n}\right|\leqslant\epsilon\right)\sum^{n}_{k=1}P\left(E_{k}\right).
\]
Taking into account the equality $\sum^{n}_{k=1}P\left(E_{k}\right)=P\left(E\right),$
we obtain the Ottaviani inequality.

\begin{remark}{}{}

In contrast to Kolmogorov inequality, Ottaviani inequality does not
require the existence of any moment of the random variables.

\end{remark}

\end{solution}

\begin{solution}{}{solexercise11.10}

For $m\in\mathbb{N}^{\ast},$ set
\[
S_{m}=\sum^{m}_{i=1}X_{i}\,\,\,\,A_{m}=\sup_{k\in\mathbb{N}^{\ast}}\left|S_{m+k}-S_{m}\right|\,\,\,\,\text{and}\,\,\,\,A=\inf_{m\in\mathbb{N}^{\ast}}A_{m}.
\]

By the Cauchy criterion for numerical series,
\[
\left\{ \sum X_{n}\,\text{converges}\right\} =\left\{ A=0\right\} .
\]

But,
\[
\left\{ A\neq0\right\} =\bigcup_{\epsilon\in\mathbb{Q}^{+\ast}}\left\{ A>\epsilon\right\} ,
\]
and, for every $n\in\mathbb{N}^{\ast},$
\[
\left\{ A>\epsilon\right\} \subset\bigcap_{m\in\mathbb{N}^{\ast}}\left\{ A_{m}>\epsilon\right\} ,
\]
which yields
\begin{equation}
\left\{ A\neq0\right\} \subset\bigcup_{\epsilon\in\mathbb{Q}^{+\ast}}\bigcap_{m\in\mathbb{N}^{\ast}}\left\{ A_{m}>\epsilon\right\} .\label{eq:A_neq_0_contained_in}
\end{equation}

Since
\[
\sup{}_{k\in\mathbb{N}^{\ast}}\left|S_{m+k}-S_{m}\right|=\lim_{r\to+\infty}\nearrow\sup_{1\leqslant k\leqslant r}\left|S_{m+k}-S_{m}\right|,
\]
the sequence of sets $\left\{ \sup_{1\leqslant k\leqslant r}\left|S_{m+k}-S_{m}\right|>\epsilon\right\} $
increases with $r,$ and we have
\begin{equation}
\left\{ A_{m}>\epsilon\right\} =\bigcup_{r\in\mathbb{N}^{\ast}}\left\{ \sup_{1\leqslant k\leqslant r}\left|S_{m+k}-S_{m}\right|>\epsilon\right\} .\label{eq:A_m>eps}
\end{equation}

Noticing that $S_{m+k}-S_{m}=\sum^{k}_{j=1}X_{j+m},$ and applying
the Ottaviani inequality to the sequence $\left(X_{j+m}\right)_{j\in\mathbb{N}^{\ast}},$
it yields
\begin{multline}
\min_{1\leqslant k\leqslant r}P\left(\left|S_{m+r}-S_{m+k}\right|\leqslant\epsilon\right)P\left(\max_{1\leqslant k\leqslant r}\left|S_{m+k}-S_{m}\right|>2\epsilon\right)\\
\leqslant P\left(\left|S_{r+m}-S_{m}\right|>\epsilon\right).\label{eq:ottaviani_X_j+m}
\end{multline}

Since the series $\sum X_{n}$ converges in probability, the sequence
with general term $S_{n}$ is Cauchy in probability. Thus, for a given
$\eta>0,$ we can choose an integer $N_{\epsilon,\eta}$ such that
for every $m\geqslant N_{\epsilon,\eta},$ 
\begin{equation}
P\left(\left|S_{m+r}-S_{m+k}\right|>\epsilon\right)\leqslant\eta,\,\,\,\,\text{whenever }0\leqslant k\leqslant r,\label{eq:P_abs_dif_Sm+r-Sm+k>eps}
\end{equation}
which implies
\[
1-\eta\leqslant P\left(\left|S_{m+r}-S_{m+k}\right|\leqslant\epsilon\right),\,\,\,\,\text{whenever }0\leqslant k\leqslant r.
\]

Hence, for such a choice of $m,$
\[
1-\eta\leqslant\min_{1\leqslant k\leqslant r}P\left(\left|S_{m+r}-S_{m+k}\right|\leqslant\epsilon\right).
\]

By the inequalities $\refpar{eq:ottaviani_X_j+m}$ and $\refpar{eq:P_abs_dif_Sm+r-Sm+k>eps},$
\[
P\left(\max_{1\leqslant k\leqslant r}\left|S_{m+k}-S_{m+n}\right|>2\epsilon\right)\leqslant\dfrac{1}{1-\eta}P\left(\left|S_{m+r}-S_{m}\right|>\epsilon\right)\leqslant\dfrac{\eta}{1-\eta}.
\]

Since in the equality $\refpar{eq:A_m>eps}$ there is a nondecreasing
sequence of sets, it follows that
\[
P\left(A_{m}>\epsilon\right)=\lim_{r\to+\infty}P\left(\sup_{1\leqslant k\leqslant r}\left|S_{m+k}-S_{m}\right|>\epsilon\right)\leqslant\dfrac{\eta}{1-\eta}.
\]
It then follows, that for every $m\geqslant N_{\epsilon,\eta},$
\[
0\leqslant P\left(\bigcap_{p\in\mathbb{N}^{\ast}}\left(A_{p}>\epsilon\right)\right)\leqslant P\left(A_{m}>\epsilon\right)\leqslant\dfrac{\eta}{1-\eta}.
\]
 To conclude, we proved that, for every $\eta>0,$
\[
0\leqslant P\left(\bigcap_{p\in\mathbb{N}^{\ast}}\left(A_{p}>\epsilon\right)\right)\leqslant\dfrac{\eta}{1-\eta}
\]
which shows that, for every $\epsilon>0,$ $P\left(\bigcap_{p\in\mathbb{N}^{\ast}}\left(A_{p}>\epsilon\right)\right)=0.$
It then follows from the inclusion $\refpar{eq:A_neq_0_contained_in}$
that $P\left(A\neq0\right)=0,$ i.e. \textbf{the series with general
term $X_{n}$ converges $P-$almost surely}.

\begin{remark}{}{}

As the converse is always true, we indeed obtain the \textbf{equivalence}
of convergence in probability and $P-$almost sure convergence for
\textbf{series} of \textbf{independent} random variables.

\end{remark}

\end{solution}

\begin{solution}{}{solexercise11.11}

1. (a) Let $t$ be an arbitrary real number. Suppose that, for every
$x$ such that $\left|x\right|\leqslant1,$
\[
0\leqslant\dfrac{1}{2}\left(1-x\right)\leqslant1,\,\,\,\,0\leqslant\dfrac{1}{2}\left(1+x\right)\leqslant1\,\,\,\,\text{and,}\,\,\,\,\text{\ensuremath{\dfrac{1}{2}\left(1-x\right)+\dfrac{1}{2}\left(1+x\right)=1.}}
\]

Moreover, since
\[
tx=\dfrac{1}{2}\left(1-x\right)\left(-t\right)+\dfrac{1}{2}\left(1+x\right)t,
\]
and as the function $x\mapsto\text{e}^{tx}$ is convex, since its
second derivative is positive, we obtain, for every $x\in\left[-1,1\right],$\boxeq{
\[
\text{e}^{tx}\leqslant\dfrac{1}{2}\left(1-x\right)\text{e}^{-t}+\dfrac{1}{2}\left(1+x\right)\text{e}^{t}.
\]
}

(b) The random variable $X$ is $P-$almost surely bounded by 1, hence
the random variable $\text{e}^{tX}$ is $P-$almost surely bounded
and admits an expectation. By the result of the previous question,
we have $P-$almost surely
\[
\text{e}^{tX}\leqslant\dfrac{1}{2}\left(1-X\right)\text{e}^{-t}+\dfrac{1}{2}\left(1+X\right)\text{e}^{t}.
\]

Therefore,
\[
\mathbb{E}\left(\text{e}^{tX}\right)\leqslant\dfrac{1}{2}\mathbb{E}\left(1-X\right)\text{e}^{-t}+\dfrac{1}{2}\mathbb{E}\left(1+X\right)\text{e}^{t}.
\]

Since the random variable $X$ is centered, this simplifies to
\[
\mathbb{E}\left(\text{e}^{tX}\right)\leqslant\dfrac{1}{2}\left(\text{e}^{-t}+\text{e}^{t}\right)=\text{ch}\left(t\right).
\]

Now 
\[
\text{ch}\left(t\right)=\sum^{+\infty}_{n=0}\dfrac{t^{2n}}{\left(2n\right)!}\,\,\,\,\,\,\,\text{and}\,\,\,\,\,\,\,\,\text{e}^{\frac{t^{2}}{2}}=\sum^{+\infty}_{n=0}\dfrac{t^{2n}}{n!2^{n}},
\]
and for every $n\in\mathbb{N},$
\[
n!2^{n}=\prod^{n}_{k=1}2k\leqslant\left(2n\right)!.
\]
Hence
\[
\text{ch}\left(t\right)\leqslant\text{e}^{\frac{t^{2}}{2}}.
\]

and it follows that\boxeq{
\begin{equation}
\mathbb{E}\left(\text{e}^{tX}\right)\leqslant\text{e}^{\frac{t^{2}}{2}}.\label{eq:upper_bound_expect_lapl_transf_X}
\end{equation}
}

2. (a) Let $t$ be an arbitrary real number. By the inequality $\refpar{eq:upper_bound_expect_lapl_transf_X}$
applied to the random variable $\dfrac{X_{n}}{c_{n}},$ we get, for
every $t^{\prime}\in\mathbb{R},$
\[
\mathbb{E}\left(\text{e}^{t^{\prime}\frac{X_{n}}{c_{n}}}\right)\leqslant\text{e}^{\frac{t^{\prime2}}{2}}.
\]

Taking $t^{\prime}=tc_{n},$ we obtain
\[
\mathbb{E}\left(\text{e}^{tX_{n}}\right)\leqslant\text{e}^{\frac{t^{2}}{2}c^{2}_{n}}.
\]

Since the random variables $\exp\left(tX_{n}\right)$ are independent,
\[
\mathbb{E}\left(\text{e}^{tS_{n}}\right)=\prod^{n}_{j=1}\mathbb{E}\left(\text{e}^{tX_{j}}\right).
\]

Thus, we have proved that for every $t\in\mathbb{R},$\boxeq{
\begin{equation}
\mathbb{E}\left(\text{e}^{tS_{n}}\right)\leqslant\exp\left(\frac{t^{2}}{2}\sum^{n}_{j=1}c^{2}_{j}\right).\label{eq:expect_lapl_transf_S_n}
\end{equation}
}

(b) Let $t>0$ and $\epsilon>0$ be arbitrary. Since the function
$x\mapsto\text{e}^{tx}$ is increasing with $x$,
\[
\left(S_{n}>\epsilon\right)\subset\left(\text{e}^{tS_{n}}>\text{e}^{t\epsilon}\right).
\]
By the Markov inequality,
\[
P\left(S_{n}>\epsilon\right)\leqslant P\left(\text{e}^{tS_{n}}>\text{e}^{t\epsilon}\right)\leqslant\dfrac{\mathbb{E}\left(\text{e}^{tS_{n}}\right)}{\text{e}^{t\epsilon}}.
\]
Using the inequality $\refpar{eq:expect_lapl_transf_S_n},$ we obtain\boxeq{
\begin{equation}
P\left(S_{n}>\epsilon\right)\leqslant\exp\left(-t\epsilon+\frac{t^{2}}{2}\sum^{n}_{j=1}c^{2}_{j}\right).\label{eq:P_S_noverepsilon}
\end{equation}
}

(c) Let $\epsilon>0$ be arbitrary and set $a=\sum^{n}_{j=1}c^{2}_{j}.$
The function $t\mapsto a\dfrac{t^{2}}{2}-t\epsilon$ reaches its minimum
at $t=\dfrac{\epsilon}{a}>0.$ The value of this minimum is $-\dfrac{\epsilon^{2}}{2a}.$
Since the exponential is strictly increasing, it follows that\boxeq{
\begin{equation}
P\left(S_{n}>\epsilon\right)\leqslant\exp\left(\min_{t>0}\left(-t\epsilon+\frac{t^{2}}{2}\sum^{n}_{j=1}c^{2}_{j}\right)\right)=\exp\left(-\frac{\epsilon^{2}}{2\sum^{n}_{j=1}c^{2}_{j}}\right).\label{eq:P_S_nover_eps_min_bound}
\end{equation}
}

(d) Let $\epsilon>0$ be arbitrary. We have the equalities
\[
\left(\left|S_{n}\right|>\epsilon\right)=\left(S_{n}>\epsilon\right)\cup\left(S_{n}<-\epsilon\right)=\left(S_{n}>\epsilon\right)\cup\left(-S_{n}>\epsilon\right),
\]
and thus the inequality
\[
P\left(\left|S_{n}\right|>\epsilon\right)\leqslant P\left(S_{n}>\epsilon\right)+P\left(-S_{n}>\epsilon\right).
\]

Applying the inequality $\refpar{eq:P_S_nover_eps_min_bound}$ to
the random variables $-S_{n},$ it yields
\[
P\left(-S_{n}>\epsilon\right)\leqslant\exp\left(-\frac{\epsilon^{2}}{2\sum^{n}_{j=1}c^{2}_{j}}\right).
\]

Hence,\boxeq{
\begin{equation}
P\left(\left|S_{n}\right|>\epsilon\right)\leqslant2\exp\left(-\frac{\epsilon^{2}}{2\sum^{n}_{j=1}c^{2}_{j}}\right).\label{eq:P_abs_S_novereps_min_upper_bound}
\end{equation}
}

(e) Let $\epsilon^{\prime}>0$ be arbitrary. Taking $\epsilon=n^{\alpha}\epsilon^{\prime}$
in the inequality $\refpar{eq:P_abs_S_novereps_min_upper_bound},$
we get
\[
P\left(\left|S_{n}\right|>n^{\alpha}\epsilon^{\prime}\right)\leqslant2\exp\left(-\frac{n^{2\alpha}\epsilon^{\prime2}}{2\sum^{n}_{j=1}c^{2}_{j}}\right).
\]
Since $\sum^{n}_{j=1}c^{2}_{j}\leqslant n^{2\alpha-\beta}$ for some
$\beta>0,$ it follows that
\begin{equation}
P\left(\left|S_{n}\right|>n^{\alpha}\epsilon^{\prime}\right)\leqslant2\text{e}^{-n^{\beta}\epsilon^{\prime2}}.\label{eq:P_S_n_n_alpha_eps_prime_upper_bound}
\end{equation}
The series with general term $\text{e}^{-n^{\beta}\epsilon^{\prime2}}$
converges: indeed, from some rank onward,
\[
n^{\beta}\epsilon^{\prime2}\geqslant2\ln n,
\]
and thus,
\[
0\leqslant\text{e}^{-n^{\beta}\epsilon^{\prime2}}\leqslant n^{-2}.
\]
By the inequality $\refpar{eq:P_S_n_n_alpha_eps_prime_upper_bound}$
the series with general term $P\left(\left|S_{n}\right|>n^{\alpha}\epsilon^{\prime}\right)$
also converges, and the Borel Cantelli lemma ensures that
\[
P\left(\limsup_{n\to+\infty}\left(\left|S_{n}\right|>n^{\alpha}\epsilon^{\prime}\right)\right)=0.
\]
Since $\mathbb{Q}^{+\ast}$ is countable, it follows that
\[
P\left(\bigcup_{\epsilon\in\mathbb{Q}^{+\ast}}\limsup_{n\to+\infty}\left(\left|S_{n}\right|>n^{\alpha}\epsilon\right)\right)=0.
\]
Taking complements,
\[
P\left(\bigcap_{\epsilon\in\mathbb{Q}^{+\ast}}\liminf_{n\to+\infty}\left(\left|S_{n}\right|\leqslant n^{\alpha}\epsilon\right)\right)=1,
\]
which means that the series with general term $n^{-\alpha}S_{n}$
converges $P-$almost surely to 0.

3. (a) Let $t$ be an arbitrary real number. By the transfer theorem
and the theorem of integration with respect to a density measure,
\[
\mathbb{E}\left(\text{e}^{tX_{1}}\right)=\intop_{\mathbb{R}}\text{e}^{tx}\text{d}P_{X}\left(x\right)=\intop_{\mathbb{R}}\text{e}^{tx}\boldsymbol{1}_{\left[0,1\right]}\left(\left|x\right|\right)\left(1-\left|x\right|\right)\text{d}x,
\]
Thus
\[
\mathbb{E}\left(\text{e}^{tX_{1}}\right)=\intop^{0}_{-1}\text{e}^{tx}\left(1+x\right)\text{d}x+\intop^{1}_{0}\text{e}^{tx}\left(1-x\right)\text{d}x.
\]
By making the change of variables $y=-x$ in the first integral, and
combining terms, we get
\[
\mathbb{E}\left(\text{e}^{tX_{1}}\right)=\intop^{1}_{0}\left(\text{e}^{-tx}+\text{e}^{tx}\right)\left(1-x\right)\text{d}x,
\]
hence
\[
\mathbb{E}\left(\text{e}^{tX_{1}}\right)=2\intop^{1}_{0}\text{ch}\left(tx\right)\left(1-x\right)\text{d}x.
\]

\begin{itemize}
\item If $t=0,$ then $\mathbb{E}\left(\text{e}^{tX_{1}}\right)=1.$
\item If $t\neq0,$ then an integration by parts yields
\[
\mathbb{E}\left(\text{e}^{tX_{1}}\right)=2\left[\left|\left(1-x\right)\dfrac{\text{sh}\left(tx\right)}{t}\right|^{1}_{0}+\intop^{1}_{0}\dfrac{\text{sh}\left(tx\right)}{t}\text{d}x\right].
\]
Hence\boxeq{
\[
\mathbb{E}\left(\text{e}^{tX_{1}}\right)=\begin{cases}
2\dfrac{\text{ch}\left(t\right)-1}{t^{2}}, & \text{if}\,t\neq0,\\
1, & \text{if}\,t=0.
\end{cases}
\]
}
\end{itemize}
(b) The function $t\mapsto\text{e}^{tX_{1}}$ is indefinitely differentiable,
and for every $k\in\mathbb{N}$ and every $t\in\mathbb{R},$
\[
\left|X^{k}_{1}\text{e}^{tX_{1}}\right|\leqslant\text{e}^{\left|t\right|}\,\,\,\,P-\text{almost surely.}
\]

The theorem of differentiation under the integral depending of a parameter---corollary
of the dominated convergence theorem---ensures that $\phi$ is indefinitely
differentiable, and that for every $k\in\mathbb{N}^{\ast},$
\[
\phi^{\left(k\right)}\left(t\right)=\mathbb{E}\left(X^{k}_{1}\text{e}^{tX_{1}}\right).
\]
In particular,\boxeq{
\[
\phi^{\left(k\right)}\left(0\right)=\mathbb{E}\left(X^{k}_{1}\right).
\]
}

By the Taylor-Young formula, for every $n\in\mathbb{N}^{\ast},$
\[
\phi\left(t\right)=1+\sum^{n}_{k=1}\dfrac{\phi^{\left(k\right)}\left(0\right)}{k!}t^{k}+o\left(t^{n}\right)=1+\sum^{n}_{k=1}\dfrac{\mathbb{E}\left(X^{k}_{1}\right)}{k!}t^{k}+o\left(t^{n}\right).
\]
In particular, for $n=2,$
\[
\phi\left(t\right)=1+t\mathbb{E}\left(X_{1}\right)+\dfrac{t^{2}}{2}\mathbb{E}\left(X_{1}\right)+o\left(t^{2}\right).
\]

Since by the expression of $\phi,$
\[
\phi\left(t\right)=1+2\dfrac{t^{2}}{4!}+o\left(t^{2}\right),
\]
so by uniqueness of the Taylor expansion, $\mathbb{E}\left(X_{1}\right)=0,$
and we obtain that
\[
\sigma^{2}_{X_{1}}=\mathbb{E}\left(X^{2}_{1}\right)=2\cdot2\dfrac{1}{4!}=\dfrac{1}{6}.
\]
The random variables $X_{n}$ being independent with same law, then\boxeq{
\[
\sigma^{2}_{S_{n}}=n\mathbb{E}\left(X^{2}_{1}\right)=\dfrac{n}{6}.
\]
}

(c) In this case, we can take $c_{n}=1$ for every $n\in\mathbb{N}^{\ast}.$
Then $\sum^{n}_{j=1}c^{2}_{j}=n$ and, to satisfy the condition in
Question 2.e., it is enough to find $\beta>0$ such that $1=2\alpha-\beta>0.$ 

If $\alpha>\dfrac{1}{2},$ then $\beta=2\alpha-1>0$ works. It follows
that, for every $\alpha>\dfrac{1}{2},$
\[
\lim_{n\to+\infty}n^{-\alpha}S_{n}=0\,\,\,\,P-\text{almost surely.}
\]

\begin{remark}{}{}As will be seen in Chapter \ref{chap:PartIIChapConvMeasuresandInLaw15},
the central limit theorem implies that, for every $x\in\mathbb{R},$
\[
\lim_{n\to+\infty}P\left(n^{-\frac{1}{2}}S_{n}\leqslant x\right)=\intop^{x}_{-\infty}\dfrac{1}{\sqrt{2\pi}}\text{e}^{-\frac{u^{2}}{2}}\text{d}u,
\]
so the sequence with general term $n^{-\frac{1}{2}}S_{n}$ does not
converge $P-$almost surely to 0.

\end{remark}

\end{solution}

\chapter{Probabilities and Conditional Expectations}\label{chap:PartIIChap12}

\begin{objective}{}{}

Chapter \ref{chap:PartIIChap12} aims at extending the concept of
conditional law.
\begin{itemize}
\item Section \ref{sec:Kernel-and-Conditional} begins by introducing probability
kernel as mapping that make it possible to define conditional laws
in a more general framework. Once a kernel is defined, a probability
measure can be defined on the product space. A generalized Fubini
theorem is then stated, followed by the definition of a conditional
law using a kernel. The Jirina theorem ensures the existence of such
conditional laws, and a version of the Fubini theorem formulated in
terms of conditional laws is presented. The section concludes with
the conditional transfer theorem, which is frequently used in the
computation of conditional laws.
\item Section \ref{sec:Conditional-Moments} focuses on the definition of
conditional moments of a real-valued random variable admitting a conditional
law, whenever these moments exist, in particular the conditional mean
and the conditional variance. Using the conditional transfer theorem,
the computation of the conditional mean is carried out. The section
ends with an application of the conditional mean to the problem of
regression.
\item Section \ref{sec:Conditional-Expectation} tackles the definition
of conditional expectation. It starts with its definition as an orthogonal
projector in the Hilbert space $\text{L}^{2}\left(\Omega,\mathscr{A},P\right),$
which is defined on a class of random variables, where an arbitrary
representative is chosen. The operator is shown to be linear, continuous,
and positive, and several of its properties are established. \\
The definition of conditional expectation is then extended to $\text{L}^{1}\left(\Omega,\mathscr{A},P\right),$
together with the corresponding properties. The case of conditional
expectation of a random variable given an independent \salg is then
examined. The definition is further extended to random variables taking
values in $\overline{\mathbb{R}}^{+}$ and that are $\mathscr{A}-$measurable,
leading to the statement of the conditional Beppo Lévi property. \\
A conditional Fatou lemma and a conditional Fatou-Lebesgue dominated
convergence are then presented. The case of continuous functions on
closed convex sets is then studied, yielding Jensen inequality. The
chapter concludes with the computation of conditional expectations
and highlights the connection between conditional mean and conditional
expectation.
\end{itemize}
\end{objective}

\section{Kernel and Conditional Laws}\label{sec:Kernel-and-Conditional}

In Part \ref{part:Introduction-to-Probability} Chapter \ref{chap:Probabilities-and-Conditional},
we introduced the conditional law of a random variable $Y$ with respect
to a \textbf{discrete} random variable $X$ taking values in a probabilizable
space $\left(E,\mathcal{E}\right).$ In that setting, conditioning
was defined, for $x\in E,$ using expressions of the form 
\[
P\left(Y\in B\mid X=x\right)=\dfrac{P\left(Y\in B,X=x\right)}{P\left(X=x\right)}
\]
However, this approach breaks down when $X$ is no more discrete---most
notably when $P\left(X=x\right)=0$ for every $x\in E.$ Since division
by zero is not defined, a more general notion of conditional law is
required. 

The purpose of this chapter is to introduce the concept of a probability
kernel---also called a transition probability---which provides the
appropriate framework for defining conditional laws in full generality.

Throughout this section, $\left(E,\mathscr{E}\right)$ and $\left(F,\mathscr{F}\right)$
denote two arbitrary probabilizable spaces.

\begin{definition}{Probability Kernel or Transition Probability}{}

A function $\nu$ from $E\times\mathscr{F}$ to $\left[0,1\right]$
is called a\textbf{ probability kernel\mindex{probability!kernel}}
on $E\times\mathscr{F}$---also called a \textbf{transition probability\index{transition probability}}
or more simply a \textbf{kernel\index{kernel}---}from $\left(E,\mathscr{E}\right)$
to $\left(F,\mathscr{F}\right),$ if it satisfies the following two
properties:

(i) For every $x\in E,$ the function $\nu\left(x,\cdot\right):\mathscr{F}\longrightarrow\left[0,1\right]$
is a probability measure on $\left(F,\mathscr{F}\right).$

(ii) For every $B\in\mathscr{F},$ the function $\nu\left(\cdot,B\right):E\longrightarrow\left[0,1\right]$
is $\mathscr{E}-$measurable.

\end{definition}

\begin{example}{Kernel from a Fixed Probability Measure}{kernel_from_proba}

Let $P$ be a probability measure on $\left(F,\mathscr{F}\right).$
The function $\nu$ from $E\times\mathscr{F}$ into $\left[0,1\right],$
defined by
\[
\forall B\in\mathscr{F},\,\,\,\:\nu\left(\cdot,B\right)=P\left(B\right)
\]
is a probability kernel on $E\times\mathscr{F}.$

\end{example}

\begin{example}{Kernel from a Dirac Measure}{kernel_from_dirac}

Let $p$ be a measurable function from $\left(E,\mathscr{E}\right)$
to $\left(F,\mathscr{F}\right).$ Denote by $\delta_{y}$ the Dirac
measure at $y.$ The function $\nu$ from $E\times\mathscr{F}$ into
$\left[0,1\right],$ defined by
\[
\forall x\in E,\,\,\,\,\,\nu\left(x,\cdot\right)=\delta_{p\left(x\right)},
\]
is a probability kernel on $E\times\mathscr{F}.$ For every $B\in\mathscr{F},$
the measurability of $\nu\left(\cdot,B\right)$ follows from the identity
\[
\nu\left(\cdot,B\right)=\boldsymbol{1}_{B}\circ p.
\]

\end{example}

\begin{example}{}{kernel_from_app_meas}

Let $f$ be an $\mathscr{E}\otimes\mathscr{F}-$measurable function
from $E\times F$ into $\overline{\mathbb{R}}^{+},$ and let $\mu$
be a nonnegative measure on $\left(F,\mathscr{F}\right).$ Moreover,
let $\rho$ be a probability measure on $\left(F,\mathscr{F}\right).$ 

Show that the function $\nu$ from $E\times\mathscr{F}$ to $\left[0,1\right]$
defined for every $\left(x,B\right)\in E\times\mathscr{F}$ by
\[
\nu\left(x,B\right)=\begin{cases}
\dfrac{\intop_{B}f\left(x,y\right)\text{d}\mu\left(y\right)}{\intop_{F}f\left(x,y\right)\text{d}\mu\left(y\right)}, & \text{if }\intop_{F}f\left(x,y\right)\text{d}\mu\left(y\right)\in\left]0,+\infty\right[,\\
\rho\left(B\right), & \text{otherwise,}
\end{cases}
\]
is a probability kernel on $E\times\mathscr{F}.$

\end{example}

\begin{solutionexample}{}{}

The first property follows directly from the definition. For each
fixed $x\in E,$ if 
\[
\intop_{F}f\left(x,y\right)\text{d}\mu\left(y\right)\in\left]0,+\infty\right[,
\]
then $\nu\left(x,\cdot\right)$ is a probability measure on $\left(F,\mathscr{F}\right),$
since it is nonnegative, countably additive, and normalized to one.
In the remaining case, $\nu\left(x,\cdot\right)=\rho\left(\cdot\right),$
which is a probability measure by assumption.

For the second property, the measurability of $\nu\left(\cdot,B\right)$
for every $B\in\mathscr{F}$ follows from an intermediate result of
Fubini theorem which ensures the measurability of functions $x\mapsto\intop_{B}f\left(x,y\right)\text{d}\mu\left(y\right)$
and $x\mapsto\intop_{F}f\left(x,y\right)\text{d}\mu\left(y\right).$

Combining these observations yields that $\nu$ is a probability kernel
on $E\times\mathscr{F}.$

\end{solutionexample}

\begin{remark}{}{}

If 
\[
\intop_{F}f\left(x,y\right)\text{d}\mu\left(y\right)\in\left]0,+\infty\right[,
\]
 then the function $B\mapsto\nu\left(x,B\right)$ is the measure with
density
\[
\dfrac{f\left(x,\cdot\right)}{\intop_{F}f\left(x,y\right)\text{d}\mu\left(y\right)}
\]
with respect to $\mu.$ 

A commonly encountered case is when $E=\mathbb{R}^{n}$ and $F=\mathbb{R}^{m},$
both equipped with their Borel $\sigma-$algebras, and when the reference
measure $\mu$ is the Lebesgue measure on $\mathbb{R}^{m}.$

\end{remark}

\textbf{A probability kernel, together with a probability measure,
allows one to define a probability measure on the product space.}

\begin{theorem}{Probability of the Product Space from a Kernel and a Probability}{}

Let $\lambda$ be a probability on the probabilizable space $\left(E,\mathscr{E}\right),$
and let $\nu$ be a probability kernel on $E\times\mathscr{F}.$ Define
the function $\lambda\cdot\nu$ on the semi-algebra of rectangles,
abusively denoted $\mathscr{E}\times\mathscr{F},$ by
\begin{equation}
\forall A\times B\in\mathscr{E}\times\mathscr{F},\,\,\,\,\lambda\cdot\nu\left(A\times B\right)=\intop_{A}\nu\left(x,B\right)\text{d}\lambda\left(x\right).\label{eq:probability_product_space}
\end{equation}

Then $\lambda\cdot\nu$ is $\sigma-$additive on $\mathscr{E}\times\mathscr{F}.$
Morevoer, there exists a unique extension of $\lambda\cdot\nu$ to
a probability measure on the product probabilizable space $\left(E\times F,\,\mathscr{E}\otimes\mathscr{F}\right)$
which is also denoted by $\lambda\cdot\nu.$

\end{theorem}

\begin{proof}{}{}

The proof relies on the theorem ensuring the existence and uniqueness
of the extension of a $\sigma-$additive function defined on a semi-algebra
to a probability measure on the generated $\sigma-$algebra.

We first prove the $\sigma-$additivity of $\lambda\cdot\nu$ on the
semi-algebra $\mathscr{E}\times\mathscr{F}.$ Let $A\times B\in\mathscr{E}\times\mathscr{F},$
and let $\left(A_{n}\times B_{n}\right)_{n\in\mathbb{N}}\subset\mathscr{E}\times\mathscr{F}$
be a sequence of disjoint sets such that
\[
A\times B=\biguplus_{n\in\mathbb{N}}\left(A_{n}\times B_{n}\right).
\]
This is equivalent to stating that
\[
\forall\left(x,y\right)\in E\times F,\,\,\,\,\boldsymbol{1}_{\left(A\times B\right)}\left(x,y\right)=\sum^{+\infty}_{n=0}\boldsymbol{1}_{A_{n}}\left(x\right)\boldsymbol{1}_{B_{n}}\left(y\right).
\]
Fix $x\in E$ and integrate with respect to $y$ under the probability
measure $\nu\left(x,\cdot\right).$ It follows that
\[
\forall x\in E,\,\,\,\,\boldsymbol{1}_{A}\left(x\right)\nu\left(x,B\right)=\sum^{+\infty}_{n=0}\boldsymbol{1}_{A_{n}}\left(x\right)\nu\left(x,B_{n}\right).
\]

It remains to integrate both sides with respect to the probability
measure $\lambda$ and to use the $\sigma-$additivity of the integral.
\[
\intop_{E}\boldsymbol{1}_{A}\left(x\right)\nu\left(x,B\right)\text{d}\lambda\left(x\right)=\sum^{+\infty}_{n=0}\intop_{E}\boldsymbol{1}_{A_{n}}\left(x\right)\nu\left(x,B_{n}\right)\text{d}\lambda\left(x\right).
\]
This yields
\[
\lambda\cdot\nu\left(A\times B\right)=\sum^{+\infty}_{n=0}\lambda\cdot\nu\left(A_{n}\times B_{n}\right),
\]
which proves the $\sigma-$additivity of $\lambda\cdot\nu$ on $\mathscr{E}\times\mathscr{F}.$

\end{proof}

\begin{remarks}{}{}

1. In Example $\ref{ex:kernel_from_proba},$ the measure $\lambda\cdot\nu$
coincides with the product probability measure $\lambda\otimes P.$

2. In Example $\ref{ex:kernel_from_dirac},$ the measure $\lambda\cdot\nu$
is defined by
\[
\forall A\times B\in\mathscr{E}\times\mathscr{F},\,\,\,\,\lambda\cdot\nu\left(A\times B\right)=\intop_{A}\boldsymbol{1}_{B}\circ p\text{d}\lambda.
\]

3. In Example $\ref{ex:kernel_from_app_meas},$ if $\mu$ is a measure
$\sigma-$finite, the measure $\lambda\cdot\mu$ is defined for every
$A\times B\in\mathscr{E}\times\mathscr{F}$ by
\[
\lambda\cdot\nu\left(A\times B\right)=\intop_{\left(A\cap C\right)\times B}\dfrac{f\left(x,y\right)}{g\left(x\right)}\text{d}\left(\lambda\otimes\mu\right)\left(x,y\right)+\lambda\left(A\cap C^{c}\right)p\left(B\right),
\]
where $g$ is defined by
\[
\forall x\in E,\,\,\,\,g\left(x\right)=\intop_{F}f\left(x,y\right)\text{d}\mu\left(y\right),
\]
and
\[
C=\left\{ x\in E:\,g\left(x\right)\in\left]0,+\infty\right[\right\} .
\]

If $\lambda\left(C\right)=1,$ then 
\[
\lambda\otimes\mu\left(\left(A\cap C\right)\times B\right)=\lambda\otimes\mu\left(A\times B\right),
\]
 and the measure $\lambda\cdot\mu$ is defined by
\[
\forall A\times B\in\mathscr{E}\times\mathscr{F},\,\,\,\,\lambda\cdot\nu\left(A\times B\right)=\intop_{A\times B}\dfrac{f\left(x,y\right)}{g\left(x\right)}\text{d}\lambda\otimes\mu\left(x,y\right).
\]
It is therefore a \textbf{measure with density} with respect to the
product measure $\lambda\otimes\mu.$

\end{remarks}

\begin{denotation}{}{}

Depending on the context, if $f$ is a function on $F$ that is integrable
with respect to the probability measure $\nu\left(x,\cdot\right),$
then its integral may be denoted interchangeably by $\nu\left(x,f\right),$
$\intop_{F}f\left(y\right)\nu\left(x,\text{d}y\right)$ or $\nu f\left(x\right).$

\end{denotation}

We now present an integration theorem with respect to the measure
$\lambda\cdot\mu,$ which extends Fubini theorem---which corresponds
to the situation of Example $\ref{ex:kernel_from_proba}.$

\begin{theorem}{Extended Fubini Theorem}{extended_fubini}

Let $f$ be a measurable function from the product probabilizable
space $\left(E\times F,\mathscr{E}\otimes\mathscr{F}\right)$ to $\left(\overline{\mathbb{R}},\mathscr{B}_{\overline{\mathbb{R}}}\right).$

(a) If $f$ is nonnegative, then the function 
\[
x\mapsto\intop_{F}f\left(x,y\right)\nu\left(x,\text{d}y\right)
\]
is $\mathscr{E}-$measurable, and
\begin{equation}
\intop_{E\times F}f\text{d}\lambda\cdot\nu=\intop_{E}\left(\intop_{F}f\left(x,y\right)\nu\left(x,\text{d}y\right)\right)\text{d}\lambda\left(x\right).\label{eq:extended_fubini}
\end{equation}

(b) If $f$ is $\lambda\cdot\nu-$integrable, then for $\lambda-$almost
every $x\in E,$ the partial function $f\left(x,\cdot\right)$ is
$\nu\left(x,\cdot\right)-$integrable. Moreover, the function defined
for $\lambda-$almost every $x\in E$ by 
\[
x\mapsto\intop_{F}f\left(x,y\right)\nu\left(x,\text{d}y\right)
\]
is $\lambda-$integrable and the equality $\refpar{eq:extended_fubini}$
still holds.

\end{theorem}

\begin{proof}{}{}

The proof follows standard arguments.

(a) Consider the family
\[
\mathscr{S}=\left\{ A\in\mathscr{E}\otimes\mathscr{F}:\,x\mapsto\nu\left(x,A^{2}_{x}\right)\text{\,is\,\,}\mathscr{E-}\text{measurable}\right\} ,
\]
where $A^{2}_{x}$ is the section\footnotemark ~of $A$ at $x.$
We now show that $\mathscr{S}$ is a $\lambda-$system. 
\begin{itemize}
\item Let $A,B\in\mathscr{S}$ with $A\subset B.$ Then $A^{2}_{x}\subset B^{2}_{x}$
and 
\[
\left(B\backslash A\right)^{2}_{x}=B^{2}_{x}\backslash A^{2}_{x}.
\]
Since $\nu\left(x,\cdot\right)$ is a probability measure,
\[
\nu\left(x,\left(B\backslash A\right)^{2}_{x}\right)=\nu\left(x,B^{2}_{x}\right)-\nu\left(x,A^{2}_{x}\right),
\]
and therefore the function $x\mapsto\nu\left(x,\left(B\backslash A\right)^{2}_{x}\right)$
is $\mathscr{E}-$measurable. Hence, $B\backslash A\in\mathscr{S}.$
\item Let $\left(A_{n}\right)_{n\in\mathbb{N}}$ be an increasing sequence
of elements of $\mathscr{S}.$ For every $n\in\mathbb{N},$ 
\[
\left(A_{n}\right)^{2}_{x}\subset\left(A_{n+1}\right)^{2}_{x}\,\,\,\,\text{and}\,\,\,\,\left[\bigcup_{n\in\mathbb{N}}A_{n}\right]^{2}_{x}=\bigcup_{n\in\mathbb{N}}\left(A_{n}\right)^{2}_{x}.
\]
By the continuity from below of the probability measure $\nu\left(x,\cdot\right),$
\[
\nu\left(x,\left[\bigcup_{n\in\mathbb{N}}A_{n}\right]^{2}_{x}\right)=\lim_{n\to+\infty}\nearrow\nu\left(x,\left(A_{n}\right)^{2}_{x}\right).
\]
Hence, the function $x\mapsto\nu\left(x,\left[\bigcup_{n\in\mathbb{N}}A_{n}\right]^{2}_{x}\right)$
is $\mathscr{E}-$measurable, and thus $\bigcup_{n\in\mathbb{N}}A_{n}\in\mathscr{S}.$
\end{itemize}
It is immediate that $\mathscr{S}$ contains the $\pi-$system of
rectangles $\mathscr{E}\times\mathscr{F}.$ Therefore, $\mathscr{S}$
contains the $\sigma-$algebra $\mathscr{E}\otimes\mathscr{F}$ generated
by $\mathscr{E}\times\mathscr{F}.$ We conclude that $\mathscr{S}=\mathscr{E}\otimes\mathscr{F}.$

Thus, the function $x\mapsto\intop_{F}f\left(x,y\right)\nu\left(x,\text{d}y\right)$
is $\mathscr{E}-$measurable for every indicator function $f$ on
a set $A\in\mathscr{E}\otimes\mathscr{F}.$ By linearity, the same
holds for all measurable step functions. For a general nonnegative
measurable function $f$, we take an increasing sequence of nonnegative
step functions converging pointwise to $f.$ We can then define for
any nonnegative measurable function, the element of $\overline{\mathbb{R}}^{+},$
\[
\phi\left(f\right)=\intop_{E}\left(\intop_{F}f\left(x,y\right)\nu\left(x,\text{d}y\right)\right)\text{d}\lambda\left(x\right).
\]

The functional $\phi$ is an integral. By definition of $\lambda\cdot\nu,$
it coincides with the integral 
\[
f\mapsto\intop_{E\times F}f\text{d}\left(\lambda\cdot\nu\right)
\]
on the set of indicator functions of measurable rectangles. Hence,
the two integrals are equal for all nonnegative measurable function
$f.$

(b) If $f$ is $\lambda\cdot\nu-$integrable, then the equality $\refpar{eq:extended_fubini}$
holds for $\left|f\right|.$ 

This implies that, for $\lambda-$almost every $x\in E,$ the partial
function $f\left(x,\cdot\right)$ is $\nu\left(x,\cdot\right)-$integrable,
and that the function, defined $\lambda-$almost everywhere, 
\[
x\mapsto\intop_{F}f\left(x,y\right)\nu\left(x,\text{d}y\right)
\]
 is $\lambda-$integrable. 

The equality $\refpar{eq:extended_fubini}$ also holds for $f^{+}$
and $f^{-}.$ Thus, by the definition of the integral on $f$ with
respect to the measure $\lambda\cdot\nu,$
\begin{align*}
\intop_{E\times F}f\text{d}\left(\lambda\cdot\nu\right) & =\intop_{E\times F}f^{+}\text{d}\left(\lambda\cdot\nu\right)-\intop_{E\times F}f^{-}\text{d}\left(\lambda\cdot\nu\right)\\
 & =\intop_{E}\left(\intop_{F}f^{+}\left(x,y\right)\nu\left(x,\text{d}y\right)\right)\text{d}\lambda\left(x\right)-\intop_{E}\left(\intop_{F}f^{-}\left(x,y\right)\nu\left(x,\text{d}y\right)\right)\text{d}\lambda\left(x\right)\\
 & =\intop_{E}\left(\intop_{F}f\left(x,y\right)\nu\left(x,\text{d}y\right)\right)\text{d}\lambda\left(x\right).
\end{align*}

\end{proof}

\footnotetext{If $A\in\mathcal{P}\left(E\times F\right),$ we define
its---possibly empty---\textbf{sections\index{sections}} as follows: 
\begin{itemize}
\item For $y\in F,$ $A^{1}_{y}=\left\{ x\in E:\,\left(x,y\right)\in A\right\} ,$
\item For $x\in E,$ $A^{2}_{x}=\left\{ y\in F:\,\left(x,y\right)\in A\right\} .$
\end{itemize}
}

\begin{corollary}{}{}With the notations of Theorem $\ref{th:extended_fubini},$
the function defined on $\mathscr{F}$ by
\[
\forall B\in\mathscr{F},\,\,\,\,\mu\left(B\right)=\intop_{E}\nu\left(x,B\right)\text{d}\lambda\left(x\right)
\]
is a probability measure on the probabilizable space $\left(F,\mathscr{F}\right).$ 

Let $g$ be a measurable function from the probabilizable space $\left(F,\mathscr{F}\right)$
into $\left(\overline{\mathbb{R}},\mathscr{B}_{\overline{\mathbb{R}}}\right).$

(a) If $g$ is nonnegative, then
\begin{equation}
\intop_{F}g\text{d}\mu=\intop_{E}\nu\left(x,g\right)\text{d}\lambda\left(x\right).\label{eq:integral_of_g_related_to_mu}
\end{equation}

(b) If $g$ is of arbitrary sign and $\mu-$integrable, then the function
$\nu\left(\cdot,g\right)$ is defined $\lambda-$almost everywhere,
coincides $\lambda-$almost everywhere with an $\mathscr{E}-$measurable
function, and the equality $\refpar{eq:integral_of_g_related_to_mu}$
still holds.

\end{corollary}

\begin{remark}{}{}

Without any difficulty, one can obtain versions of Theorem $\ref{th:extended_fubini}$
and of its corollary in the framework of functions taking values in
$\mathbb{R}^{d},$ or more generally in an Euclidean space.

\end{remark}

\textbf{In what follows, $X$ and $Y$ denote two random variables
taking values respectively, in the arbitrary probabilizable spaces
$\left(E,\mathscr{E}\right)$ and $\left(F,\mathscr{F}\right).$}

Assume that $X$ is discrete and that $\mathscr{E}$ contains all
singletons, that is, $\left\{ x\right\} \in\mathscr{E},$ for every
$x\in E.$ We then define a kernel $\nu$ on $E\times\mathscr{F}$
by
\[
\forall\left(x,B\right)\in E\times\mathscr{F},\,\,\,\,\nu\left(x,B\right)=\begin{cases}
P^{\left(X=x\right)}\left(Y\in B\right), & \text{if }x\in\text{val}\left(X\right),\\
\rho\left(B\right), & \text{otherwise,}
\end{cases}
\]
where $\rho$ is an arbitrary probability on $\left(F,\mathscr{F}\right).$ 

Recall that 
\[
\text{val}\left(X\right)=\left\{ x\in E:\,P\left(X=x\right)\neq0\right\} ,
\]
which is a countable union of singletons, and therefore belongs to
$\mathscr{E},$ and that $P_{X}\left(\text{val}\left(X\right)\right)=1.$

For every $A\in\mathscr{E}$ and $B\in\mathscr{F},$
\begin{align*}
P\left[\left(X\in A\right)\cap\left(Y\in B\right)\right] & =\sum_{x\in\text{val}\left(X\right)\cap A}P\left(\left(X=x\right)\cap\left(Y\in B\right)\right)\\
 & =\sum_{x\in\text{val}\left(X\right)\cap A}\nu\left(x,B\right)P_{X}\left(\left\{ x\right\} \right).
\end{align*}
 Since the probability $P_{X}$ can be written as
\[
P_{X}=\sum_{x\in\text{val}\left(X\right)}P\left(X=x\right)\delta_{x},
\]
it follows that
\[
P_{\left(X,Y\right)}\left(A\times B\right)=\intop_{A}\nu\left(x,B\right)\text{d}P_{X}\left(x\right).
\]
This is equivalent to 
\[
P_{\left(X,Y\right)}=P_{X}\cdot\nu.
\]
This relationship serves as the starting point for defining conditional
laws in the general case.

\begin{definition}{Conditional Law}{}

A \textbf{conditional law of $Y$ given $X$\index{conditional law of $Y$ given $X$}}
is a probability kernel $\nu$ on $E\times\mathscr{F}$ such that\boxeq{
\[
P_{\left(X,Y\right)}=P_{X}\cdot\nu.
\]
}This conditional law is often denoted $P^{X=\cdot}_{Y},$ and the
preceding identity is then written\boxeq{
\begin{equation}
P_{\left(X,Y\right)}=P_{X}\cdot P^{X=\cdot}_{Y}.\label{eq:conditional_law_def}
\end{equation}
}

\end{definition}

\begin{example}{}{cond_dens_from_kernel_app_meas}

Let $\lambda$ and $\mu$ be $\sigma-$finite measures on $\left(E,\mathscr{E}\right)$
and $\left(F,\mathscr{F}\right),$ respectively. Assume that 
\[
P_{\left(X,Y\right)}=f\cdot\left(\lambda\otimes\mu\right),
\]
where $f$ is a nonnegative measurable function on $\left(E\times F,\mathscr{E}\otimes\mathscr{F}\right),$
whose $\lambda\otimes\mu-$integral equals 1. Then the kernel $\nu$
defined in Example $\ref{ex:kernel_from_app_meas}$ is a conditional
law on $Y$ given $X.$ 

The function 
\[
x\mapsto\intop_{F}f\left(x,y\right)\text{d}\mu\left(y\right)
\]
is simply the density of $P_{X}$ with respect to $\lambda.$

\end{example}

\begin{definition}{Conditional Density}{}

In the case where $E=\mathbb{R}^{n}$ and $F=\mathbb{R}^{m},$ both
equipped with their Borel $\sigma-$algebras, let $\nu$ be a conditional
law of $Y$ given $X$ such that, for $P_{X}-$almost every $x\in\mathbb{R}^{n},$
the measure $\nu\left(x,\cdot\right)$ admits a density $f^{X=x}_{Y}$
with respect to the Lebesgue measure on $\mathbb{R}^{m},$ then $f^{X=x}_{Y}$
is called a \textbf{\index{conditional density of $Y$ given $X=x$}conditional
density of $Y$ given $X=x.$} 

\end{definition}

\begin{example}{Common Special Case}{common_sp_case_example}

Let $E=\mathbb{R}^{n}$ and $F=\mathbb{R}^{m},$ both equipped with
their Borel $\sigma-$algebras, and take the Lebesgue measures on
these spaces as reference measures.

\textbf{Assume} that the random variable $\left(X,Y\right)$ admits
a density $f_{\left(X,Y\right)}.$ Then $X$ admits a density $f_{X},$
and for any probability $\rho$ on $\left(\mathbb{R}^{m},\mathscr{B}_{\mathbb{R}^{m}}\right),$
the kernel $\nu$ defined for $\left(x,B\right)\in\mathbb{R}^{n}\times\mathscr{B}_{\mathbb{R}^{m}}$
by
\[
\nu\left(x,B\right)=\begin{cases}
\intop_{B}\dfrac{f_{\left(X,Y\right)}\left(x,y\right)}{f_{X}\left(x\right)}\text{d}\lambda_{m}\left(y\right), & \text{if }f_{X}\left(x\right)>0,\\
\rho\left(B\right), & \text{if }f_{X}\left(x\right)=0,
\end{cases}
\]
is a conditional law of $Y$ given $X,$ which is a particular case
of Example $\ref{ex:cond_dens_from_kernel_app_meas}.$ 

If $f_{X}\left(x\right)>0,$ then the measure $\nu\left(x,\cdot\right)$
admits the density 
\[
\dfrac{f_{\left(X,Y\right)}\left(x,\cdot\right)}{f_{X}\left(x\right)}
\]
with respect to the Lebesgue measure $\lambda_{m}.$ 

Thus, the marginal random variable $Y$ admits a conditional density
with respect to $X$---or given $X$---denoted $f^{X=\cdot}_{Y},$
which satisfies, for every $\left(x,y\right)\in\mathbb{R}^{n}\times\mathbb{R}^{m}$
such that $f_{X}\left(x\right)>0,$\boxeq{
\begin{equation}
f^{X=x}_{Y}\left(y\right)=\dfrac{f_{\left(X,Y\right)}\left(x,y\right)}{f_{X}\left(x\right)}.\label{eq:conditional_density_marginal}
\end{equation}
}

\textbf{Conversely}, if, for $P_{X}-$almost every $x\in\mathbb{R}^{n},$
there exists a conditional density $f^{X=x}_{Y}$ of $Y$ given $X=x,$
and if $X$ admits a density $f_{X},$ then the random variable $\left(X,Y\right)$
admits a density $f_{\left(X,Y\right)}$ which verifies, for $P_{\left(X,Y\right)}-$almost
every $\left(x,y\right)\in\mathbb{R}^{n}\times\mathbb{R}^{m},$\boxeq{
\begin{equation}
f_{\left(X,Y\right)}\left(x,y\right)=f_{X}\left(x\right)f^{X=x}_{Y}\left(y\right).\label{eq:density_from_cond_density}
\end{equation}
}

Indeed, by the definition of the conditional law, for every $A\in\mathscr{B}_{\mathbb{R}^{n}}$
and $B\in\mathscr{B}_{\mathbb{R}^{m}},$
\[
P_{\left(X,Y\right)}\left(A\times B\right)=\intop_{A}\left(\intop_{B}f^{X=x}_{Y}\left(y\right)\text{d}\lambda_{m}\left(y\right)\right)f_{X}\left(x\right)\text{d}\lambda_{n}\left(x\right),
\]
which yields the result by Fubini theorem.

\end{example}

Let us now show how these two situations may arise simultaneously.
Consider two independent real-valued random variables following the
same exponential law $\exp\left(\lambda\right).$ Let $S=X+Y,$ and
we seek the conditional law of $X$ given $S.$ The random variables
$X$ and $Y$ may represent, for instance, the waiting times of two
customers arriving as a service counter.

By independence of $X$ and $Y,$ the random variable $\left(X,Y\right)$
admits a density, equals to the direct product of the marginal densities.
Using the change of variables, defined on $\mathbb{R}^{2}$ by
\[
\left\{ \begin{array}{l}
x=t\\
y=s-t,
\end{array}\right.
\]
whose Jacobian is 1, the random variable $\left(X,S\right)$ admits
a density $f_{\left(X,S\right)}$ given by
\[
\forall\left(t,s\right)\in\mathbb{R}^{2},\,\,\,\,f_{\left(X,S\right)}\left(t,s\right)=f_{X}\left(t\right)f_{Y}\left(t,s-t\right).
\]
It follows---as seen previously---that $S$ admits a density $f_{S}$
given by
\[
\forall s\in\mathbb{R},\,\,\,\,f_{S}\left(s\right)=\boldsymbol{1}_{\mathbb{R}^{+}}\left(s\right)\lambda^{2}s\text{e}^{-\lambda s},
\]
and that, for every $s>0,$ $X$ admits a conditional density given
$S=s,$ which, after simplification, is\boxeq{
\[
f^{S=s}_{X}\left(x\right)=\dfrac{1}{s}\boldsymbol{1}_{\left[0,s\right]}\left(x\right).
\]
}Hence, for every $s$ in the interior of the support of the law
of $S,$ the conditional law of $X$ given $S$ is the uniform law
on the interval $\left[0,s\right].$

\textbf{A concrete example will now illustrate that the notion of
conditional law introduced above agrees with our intuitive understanding.}

\begin{example}{}{}

The random variable $\left(X,Y\right)$ represents a point drawn at
random in the square $\left[0,1\right]^{2}.$ That is, $\left(X,Y\right)$
is of uniform law on $\left[0,1\right]^{2}$ and admits the density
$\boldsymbol{1}_{\left[0,1\right]^{2}}.$ Let $S=X+Y.$ What is the
conditional law of $X$ given $S?$

\textit{Hint: Note that $X$ and $Y$ are independent in this case
and both have the same uniform law on $\left[0,1\right].$}

\end{example}

\begin{solutionexample}{}{}

We are in the same setting as the previous example, but with a different
law. We now carry out the computation in more detail. Let $T$ be
the diffeomorphism on $\mathbb{R}^{2}$ onto itself defined by
\[
\forall\left(x,y\right)\in\mathbb{R}^{2},\,\,\,\,T\left(x,y\right)=\left(x,x+y\right).
\]

Its inverse is given by
\[
\forall\left(u,s\right)\in\mathbb{R}^{2},\,\,\,\,T^{-1}\left(u,s\right)=\left(u,s-u\right).
\]
The Jacobian of this diffeomorphism has absolute value 1 and we have
$\left(X,S\right)=T\circ\left(X,Y\right).$ Therefore, the random
variable $\left(X,S\right)$ has a density $f_{\left(X,S\right)}$
given by
\[
\forall\left(u,s\right)\in\mathbb{R}^{2},\,\,\,\,f_{\left(X,S\right)}\left(u,s\right)=f_{\left(X,Y\right)}\left(u,s-u\right).
\]
Hence,
\[
\forall\left(u,s\right)\in\mathbb{R}^{2},\,\,\,\,f_{\left(X,S\right)}\left(u,s\right)=\boldsymbol{1}_{\left[0,1\right]}\left(u\right)\boldsymbol{1}_{\left[0,1\right]}\left(s-u\right).
\]

Thus, the marginal random variable $S$ admits a density given by
\[
\forall s\in\mathbb{R},\,\,\,\,f_{S}\left(s\right)=\intop_{\mathbb{R}}f_{\left(X,S\right)}\left(u,s\right)\text{d}u=\intop_{\mathbb{R}}\boldsymbol{1}_{\left[0,1\right]}\left(u\right)\boldsymbol{1}_{\left[0,1\right]}\left(s-u\right)\text{d}u.
\]
Decomposing the product of indicator as
\[
\boldsymbol{1}_{\left[0,1\right]}\left(u\right)\boldsymbol{1}_{\left[0,1\right]}\left(s-u\right)=\boldsymbol{1}_{\left[0,1\right]}\left(s\right)\boldsymbol{1}_{\left[0,s\right]}\left(u\right)+\boldsymbol{1}_{\left]1,2\right]}\left(s\right)\boldsymbol{1}_{\left[s-1,1\right]}\left(u\right),
\]
we obtain
\[
\forall s\in\mathbb{R},\,\,\,\,f_{S}\left(s\right)=s\boldsymbol{1}_{\left[0,1\right]}\left(s\right)+\left(2-s\right)\boldsymbol{1}_{\left]1,2\right]}\left(s\right).
\]
Thus, the law of $S$ is the \textbf{triangular law}. 

For every $s$ in the interior of the support of $f_{S},$ the random
variable $X$ admits a conditional density given $S=s,$ $f^{S=s}_{X},$
given by\boxeq{
\[
\forall x\in\mathbb{R},\,\,\,\,f^{S=s}_{X}\left(x\right)=\begin{cases}
\dfrac{1}{s}\boldsymbol{1}_{\left[0,s\right]}\left(x\right), & \text{if\,\,\,\,}0<s\leqslant1,\\
\dfrac{1}{2-s}\boldsymbol{1}_{\left[s-1,1\right]}\left(x\right), & \text{if\,\,\,\,}1<s<2.
\end{cases}
\]
}

Hence, the conditional law of the random variable $X$ given $S=s$
is the uniform law on $\left[0,s\right],$ if $0<s\leqslant1,$ and
the uniform law on $\left[s-1,1\right]$ if $1<s<2$---which is a
perfectly natural result. 

\end{solutionexample}

\begin{example}{}{}

If $X$ and $Y$ are two independent random variables taking values
respectively in the probabilizable spaces $\left(E,\mathscr{E}\right)$
and $\left(F,\mathscr{F}\right)$ respectively, show that the ``constant''
kernel $\nu$ defined by
\[
\forall x\in E,\,\,\,\,\nu\left(x,\cdot\right)=P_{Y}
\]
 is a conditional law of $Y$ given $X.$

\end{example}

\begin{solutionexample}{}{}

Since $X$ and $Y$ are independent, $P_{\left(X,Y\right)}=P_{X}\otimes P_{Y}.$
Therefore, for every $A\in\mathscr{E}$ and $B\in\mathscr{F},$
\[
P_{\left(X,Y\right)}\left(A\times B\right)=\intop_{A}\nu\left(x,B\right)\text{d}P_{X}\left(x\right),
\]
which shows that $P_{\left(X,Y\right)}=P_{X}\cdot\nu.$ Hence, $\nu$
is a conditional law of $Y$ given $X.$

\end{solutionexample}

\begin{remark}{}{}

It is immediate that any other kernel $\nu'$ on $E\times F$ such
that, for every $B\in\mathscr{F},$ 
\[
\nu'\left(\cdot,B\right)=\nu\left(\cdot,B\right),\,\,\,\,P_{X}-\text{almost surely},
\]
is also a conditional law of $Y$ given $X.$ Consequently, the conditional
law is not unique.

What remains is the question of existence. This issue has already
been partially addressed in the preceding examples. We now state,
for reference, a rather general existence theorem. Its proof is omitted,
as it lies beyond the scope of this book.

\end{remark}

\begin{theorem}{Jirina Theorem}{}

Let $E$ and $F$ be two complete separable metric spaces---in particular,
Euclidean spaces---both equipped with their Borel $\sigma-$algebras.
Let $X$ and $Y$ be two random variables taking values in $E$ and
$F,$ respectively. Then, there exists a conditional law of $Y$ given
$X.$

\end{theorem}

We now state an extended version of Fubini theorem in terms of conditional
laws; no proof is required. 

\begin{theorem}{Extended Fubini Theorem in Term of Conditional Laws}{gen_fubini_cond}

Let $\left(X,Y\right)$ be a random variable taking values in an arbitrary
probabilizable space $\left(E\times F,\,\mathscr{E}\otimes\mathscr{F}\right),$
and assume that there exists a conditional law $P^{X=\cdot}_{Y}$
of $Y$ given $X.$ Let $f$ be a measurable function from the probabilizable
space $\left(E\times F,\,\mathscr{E}\otimes\mathscr{F}\right)$ into
$\left(\overline{\mathbb{R}},\mathscr{B}_{\overline{\mathbb{R}}}\right).$

(a) If $f$ is nonnegative, then the function $x\mapsto\intop_{F}f\left(x,y\right)\text{d}P^{X=x}_{Y}\left(y\right)$
is $\mathscr{E}-$measurable, and
\begin{equation}
\intop_{E\times F}f\text{d}P_{\left(X,Y\right)}=\intop_{E}\left(\intop_{F}f\left(x,y\right)\text{d}P^{X=x}_{Y}\left(y\right)\right)\text{d}P_{X}\left(x\right).\label{eq:extended_fubini_th_non_neg}
\end{equation}

(b) If $f$ is of arbitrary sign and $P_{\left(X,Y\right)}-$integrable,
then for $P_{X}-$almost every $x\in E,$ the partial function $f\left(x,\cdot\right)$
is $P^{X=\cdot}_{Y}-$integrable. Moreover the function $x\mapsto\intop_{F}f\left(x,y\right)\text{d}P^{X=x}_{Y}\left(y\right)$
defined $P_{X}-$almost surely, is $P_{X}-$integrable, and the equality
$\refpar{eq:extended_fubini_th_non_neg}$ still holds.

\end{theorem}

This result yields a ``conditional'' transfer theorem, which is
frequently used in the computation of conditional laws.

\begin{theorem}{Conditional Transfer Theorem}{cond_transfer}

Let $\left(X,Y\right)$ be a random variable taking values in a probabilizable
space $\left(E\times F,\mathscr{E}\otimes\mathscr{F}\right).$ Assume
that there exists a conditional law $\nu=P^{X=\cdot}_{Y}$ of $Y$
given $X.$ Let $f$ be a measurable function from $\left(E\times F,\mathscr{E}\otimes\mathscr{F}\right)$
to another probabilizable space $\left(G,\mathscr{G}\right).$

A conditional law of $f\left(X,Y\right)$ given $X$ is given by the
kernel $\mu,$ defined in terms of image measures by
\[
\forall x\in E,\,\,\,\,\mu\left(x,\cdot\right)=f\left(x,\cdot\right)\left[\nu\left(x,\cdot\right)\right].
\]
This can be written more suggestively as\boxeq{
\begin{equation}
\forall x\in E,\,\,\,\,P^{X=x}_{f\left(X,Y\right)}=P^{X=x}_{f\left(x,Y\right)}.\label{eq:cond_law_from_kernel}
\end{equation}
}

In particular, if $X$ and $Y$ are independent,\boxeq{
\begin{equation}
\forall x\in E,\,\,\,\,P^{X=x}_{f\left(X,Y\right)}=P_{f\left(x,Y\right)}.\label{eq:cond_law_from_kernel_independence}
\end{equation}
}

\end{theorem}

\begin{proof}{}{}

Let $A\in\mathscr{E}$ and $B\in\mathscr{G}.$ Using standard notations
and applying Theorem $\ref{th:gen_fubini_cond},$ we obtain
\begin{align*}
P_{\left(X,f\left(X,Y\right)\right)}\left(A\times B\right) & =P_{\left(X,Y\right)}\left(\left(A\times F\right)\cap f^{-1}\left(B\right)\right)\\
 & =\intop_{E}\left(\intop_{F}\boldsymbol{1}_{A\times F}\left(x,y\right)\boldsymbol{1}_{f^{-1}\left(B\right)}\left(x,y\right)\nu\left(x,\text{d}y\right)\right)\text{d}P_{X}\left(x\right).
\end{align*}
Since
\[
\boldsymbol{1}_{A\times F}\left(x,y\right)\boldsymbol{1}_{f^{-1}\left(B\right)}\left(x,y\right)=\boldsymbol{1}_{A}\left(x\right)\boldsymbol{1}_{\left[f\left(x,\cdot\right)\right]^{-1}\left(B\right)}\left(y\right),
\]
it follows, by the definition of the image measure $\mu\left(x,\cdot\right)$
of $\nu\left(x,\cdot\right),$ that
\[
P_{\left(X,f\left(X,Y\right)\right)}\left(A\times B\right)=\intop_{A}\mu\left(x,B\right)\text{d}P_{X}\left(x\right).
\]
This proves that $P_{\left(X,f\left(X,Y\right)\right)}=P_{X}\cdot\mu,$
and therefore that $\mu$ is a conditional law of $f\left(X,Y\right)$
given $X.$

\end{proof}

The following example illustrates the different concepts and theorems
that have been introduced so far in two non-standard situations, where
the laws and conditional laws combine laws with density and punctual
measures---weighted Dirac measures.

\begin{example}{}{2_cases_mixed_laws}

Let $\left(X,Y\right)$ be a random variable taking values in $\left(\mathbb{R}^{2},\mathscr{B}_{\mathbb{R}^{2}}\right).$
Let $\mathbb{\epsilon\in}\left]0,1\right[.$ Assume that $X$ follows
the uniform law on the interval $\left[0,1\right].$ Denote by $\lambda$
the Lebesgue measure on $\mathbb{R}.$ Consider the following two
cases:
\begin{itemize}
\item \textbf{Case 1:} A conditional law $P^{X=\cdot}_{Y}$ of $Y$ given
$X$ is defined by
\[
\forall x\in\left[0,1\right],\,\,\,\,P^{X=x}_{Y}=\boldsymbol{1}_{\left[\epsilon,1\right]}\cdot\lambda+\epsilon\delta_{x},
\]
that is, the probability $P^{X=x}_{Y}$ is a mixture of the uniform
probability measure on the interval $\left[\epsilon,1\right]$ and
a Dirac mass at $x.$
\item \textbf{Case 2:} A conditional law $P^{X=\cdot}_{Y}$ of $Y$ given
$X$ is defined by
\[
P^{X=x}_{Y}=\begin{cases}
\boldsymbol{1}_{\left[\epsilon,1\right]}\cdot\lambda+\epsilon\delta_{x}, & \forall x\in\left[0,\epsilon\right[,\\
\mathscr{U}\left(\left[0,1\right]\right), & \forall x\in\left[\epsilon,1\right[,
\end{cases}
\]
that is, if $0\leqslant x<\epsilon,$ the probability $P^{X=x}_{Y}$
is again a mixture of the uniform probability measure on the interval
$\left[\epsilon,1\right]$ and a Dirac mass at $x;$ whereas if $\epsilon\leqslant x\leqslant1,$
the law $P^{X=x}_{Y}$ is the uniform law on $\left[0,1\right].$
\end{itemize}
1. Study, in both cases, the law of the random variable $Y.$

2. Compute, in both cases, the covariance of the random variables
$X$ and $Y.$

\end{example}

\begin{solutionexample}{}{}

\textbf{1. Study of the law of $Y$}

By the definition of a conditional law, for every Borel sets $A$
and $B$ of $\mathbb{R},$
\[
P_{\left(X,Y\right)}\left(A\times B\right)=\intop_{A}P^{X=x}_{Y}\left(B\right)\text{d}P_{X}\left(x\right).
\]

\begin{itemize}
\item \textbf{Case 1}\\
For $A,B\in\mathscr{B}_{\mathbb{R}},$
\begin{align*}
P_{\left(X,Y\right)}\left(A\times B\right) & =\intop_{A}\boldsymbol{1}_{\left[0,1\right]}\left(x\right)\left[\lambda\left(B\cap\left[\epsilon,1\right]\right)+\epsilon\boldsymbol{1}_{B}\left(x\right)\right]\text{d}\lambda\left(x\right)\\
 & =\lambda\left(A\cap\left[0,1\right]\right)\lambda\left(B\cap\left[\epsilon,1\right]\right)+\epsilon\lambda\left(A\cap B\cap\left[0,1\right]\right).
\end{align*}
Taking $A=\mathbb{R}.$ we obtain, for every $B\in\mathscr{B}_{\mathbb{R}},$
the law of $Y,$
\begin{align*}
P_{Y}\left(B\right) & =\lambda\left(B\cap\left[\epsilon,1\right]\right)+\epsilon\lambda\left(B\cap\left[0,1\right]\right)\\
 & =\intop_{B}\left[\boldsymbol{1}_{\left[\epsilon,1\right]}+\epsilon\boldsymbol{1}_{\left[0,1\right]}\right]\text{d}\lambda.
\end{align*}
Hence, the random variable $Y$ admits a density $f_{Y}$ given by
\[
f_{Y}=\boldsymbol{1}_{\left[\epsilon,1\right]}+\epsilon\boldsymbol{1}_{\left[0,1\right]},
\]
or equivalently,\boxeq{
\[
f_{Y}=\epsilon\boldsymbol{1}_{\left[0,\epsilon\right[}+\left(1+\epsilon\right)\boldsymbol{1}_{\left[\epsilon,1\right]}.
\]
}
\item \textbf{Case 2}\\
For $A,B\in\mathscr{B}_{\mathbb{R}},$
\begin{multline*}
P_{\left(X,Y\right)}\left(A\times B\right)=\intop_{A\cap\left[0,\epsilon\right[}\boldsymbol{1}_{[0,1[}\left(x\right)\left[\lambda\left(B\cap\left[\epsilon,1\right]\right)+\epsilon\boldsymbol{1}_{B}\left(x\right)\right]\text{d}\lambda\left(x\right)\\
+\intop_{A\cap\left[\epsilon,1\right[}\boldsymbol{1}_{\left[0,1\right]}\left(x\right)\lambda\left(B\cap\left[0,1\right]\right)\text{d}\lambda\left(x\right).
\end{multline*}
This becomes
\begin{multline*}
P_{\left(X,Y\right)}\left(A\times B\right)=\lambda\left(A\cap\left[0,\epsilon\right[\right)\lambda\left(B\cap\left[\epsilon,1\right]\right)+\epsilon\lambda\left(A\cap B\cap\left[0,\epsilon\right[\right)\\
+\lambda\left(A\cap\left[\epsilon,1\right]\right)\lambda\left(B\cap\left[0,1\right]\right).
\end{multline*}
Taking $A=\mathbb{R},$ yields, after simplification,
\[
\forall B\in\mathscr{B}_{\mathbb{R}},\,\,\,\,P_{Y}\left(B\right)=\lambda\left(B\cap\left[0,1\right]\right).
\]
Thus, in this case, the random variable $Y$ follows the uniform law
on the inverval $\left[0,1\right].$\\
Thus this is \textbf{an example of random variables $X$ and $Y$
each following the uniform lawn in the interval $\left[0,1\right],$
and such that the law followed by the pair $\left(X,Y\right)$ is
not the uniform law on the square $\left[0,1\right]^{2}.$}
\end{itemize}
\textbf{We compute in these two cases the covariance of the random
variables $X$ and Y.}

\textbf{2. Covariance of $X$ and $Y$}

The random variables $X$ and $Y$ are bounded by 1. By the transfer
theorem, the function $\left(x,y\right)\mapsto xy$ is $P_{\left(X,Y\right)}-$integrable,
and by applying Theorem $\ref{th:gen_fubini_cond},$
\[
\mathbb{E}\left(XY\right)=\intop_{\mathbb{R}}\left[\intop_{\mathbb{R}}xy\text{d}P^{X=x}_{Y}\left(y\right)\right]\text{d}P_{X}\left(x\right).
\]

\begin{itemize}
\item \textbf{Case 1}
\[
\mathbb{E}\left(XY\right)=\intop_{\mathbb{R}}x\left[\intop_{\left[\epsilon,1\right]}y\text{d}\lambda\left(y\right)+\epsilon x\right]\boldsymbol{1}_{\left[0,1\right]}\left(x\right)\text{d}\lambda\left(x\right).
\]
Identifying Lebesgue and Riemann integrals,
\[
\mathbb{E}\left(XY\right)=\intop^{1}_{0}x\left[\dfrac{1}{2}\left(1-\epsilon^{2}\right)+\epsilon x\right]\text{d}x=\dfrac{1}{4}\left(1-\epsilon^{2}\right)+\dfrac{\epsilon}{3}.
\]
Moreover,
\[
\mathbb{E}\left(X\right)=\dfrac{1}{2}\,\,\,\,\text{and}\,\,\,\,\mathbb{E}\left(Y\right)=\dfrac{1}{2}\left(1-\epsilon^{2}\right)+\dfrac{\epsilon}{2}.
\]
Since $\text{cov}\left(X,Y\right)=\mathbb{E}\left(XY\right)-\mathbb{E}\left(X\right)\mathbb{E}\left(Y\right),$
we obtain\\
\boxeq{
\[
\text{cov}\left(X,Y\right)=\dfrac{\epsilon}{12}.
\]
}
\item \textbf{Case 2}\\
Similarly,
\[
\mathbb{E}\left(XY\right)=\intop_{\left[0,\epsilon\right[}x\left[\intop_{\left[\epsilon,1\right]}y\text{d}\lambda\left(y\right)+\epsilon x\right]\text{d}\lambda\left(x\right)+\intop_{\left[\epsilon,1\right]}x\left[\intop_{\left[0,1\right]}y\text{d}\lambda\left(y\right)\right]\text{d\ensuremath{\lambda}\ensuremath{\left(x\right)}},
\]
Thus, identifying the Lebesgue and Riemann integrals,
\[
\mathbb{E}\left(XY\right)=\intop^{\epsilon}_{0}x\left(\dfrac{1}{2}\left(1-\epsilon^{2}\right)+\epsilon x\right)\text{d}x+\dfrac{1}{2}\intop^{1}_{\epsilon}x\text{d}x.
\]
Hence, a direct computation yields 
\[
\mathbb{E}\left(XY\right)=\dfrac{1}{4}+\dfrac{\epsilon^{4}}{12}.
\]
Since the random variables $X$ and $Y$ are both of uniform law on
$\left[0,1\right],$ 
\[
\mathbb{E}\left(X\right)=\mathbb{E}\left(Y\right)=\dfrac{1}{2}.
\]
Hence,\boxeq{
\[
\text{cov}\left(X,Y\right)=\dfrac{\epsilon^{4}}{12}.
\]
}
\end{itemize}
\end{solutionexample}

\section{Conditional Moments}\label{sec:Conditional-Moments}

We now define\textbf{ \index{conditional moments}conditional moments},
when they exist.

\begin{proposition}{Conditional Moments}{}

Let $X$ be a random variable taking values in an arbitrary probabilizable
space $\left(E,\mathscr{E}\right).$ Let $Y$ be a real-valued random
variable such that there exists a conditional law $P^{X=\cdot}_{Y}$
of $Y$ given $X.$ If, for some $p\in\mathbb{N}^{\ast},$ the random
variable $Y$ admits a moment of order $p,$ then $P_{X}-$almost
surely
\[
\intop_{\mathbb{R}}\left|y\right|^{p}\text{d}P^{X=\cdot}_{Y}\left(y\right)<+\infty.
\]

\begin{itemize}
\item If $p=1,$ a \textbf{conditional mean\index{conditional mean}}\footnotemark
~is any $\mathscr{E}-$measurable function that is $P_{X}-$almost
surely equal to 
\[
m^{X=\cdot}_{Y}\equiv\intop_{\mathbb{R}}y\text{d}P^{X=\cdot}_{Y}\left(y\right).
\]
\item If $p=2,$ a \textbf{conditional variance\index{conditional variance}}
is any $\mathcal{E-}$measurable function, $P_{X}-$almost surely
equal to 
\[
\intop_{\mathbb{R}}\left(y-\intop_{\mathbb{R}}y\text{d}P^{X=\cdot}_{Y}\left(y\right)\right)^{2}\text{d}P^{X=\cdot}_{Y}\left(y\right).
\]
\end{itemize}
\end{proposition}

\footnotetext{Tr.N. Here, the term mean is used deliberately to distinguish
the conditional mean of a real-valued random variable from the conditional
expectation, which is an operator and will be defined in the next
section.}

\begin{proof}{}{}

By assumption,
\[
\mathbb{E}\left(\left|Y\right|^{p}\right)=\intop_{E\times\mathbb{R}}\left|y\right|^{p}\text{d}P^{X=\cdot}_{Y}\left(y\right)<+\infty.
\]
Applying the extended Fubini theorem in terms of conditional laws
yields
\[
\intop_{E}\left(\intop_{\mathbb{R}}\left|y\right|^{p}\text{d}P^{X=\cdot}_{Y}\left(y\right)\right)\text{d}P_{X}\left(x\right)<+\infty.
\]
The result follows.

\end{proof}

\begin{remark}{}{}

By the extended Fubini theorem and the transfer theorem, we obtain,
for every $C\in\mathscr{E},$
\begin{align*}
\intop_{C}m^{X=\cdot}_{Y}\text{d}P_{X}\left(x\right) & =\intop_{E}\boldsymbol{1}_{C}\left(x\right)\left(\intop_{\mathbb{R}}y\text{d}P^{X=\cdot}_{Y}\left(y\right)\right)\text{d}P_{X}\left(x\right)\\
 & =\intop_{E\times\mathbb{R}}\boldsymbol{1}_{C}\left(x\right)y\text{d}P_{\left(X,Y\right)}\left(x,y\right)\\
 & =\intop_{\Omega}\boldsymbol{1}_{C}\left(X\right)Y\text{d}P.
\end{align*}
Hence,\boxeq{
\begin{equation}
\forall C\in\mathscr{E},\,\,\,\,\intop_{C}m^{X=x}_{Y}\text{d}P_{X}\left(x\right)=\intop_{X^{-1}\left(C\right)}Y\text{d}P.\label{eq:cond_exp_from_law}
\end{equation}
}

That is, the conditional mean $m^{X=\cdot}_{Y},$ together with the
law of $X,$ is sufficient to compute the expectation of $Y$ over
any element $X^{-1}\left(C\right),$ of the $\sigma-$algebra generated
by $X.$

\end{remark}

\textbf{We now give an example illustrating the computation of a conditional
mean.}

\begin{example}{}{}

Using Example $\ref{ex:2_cases_mixed_laws},$ compute the conditional
mean of $Y$ given $X$ in both cases.

\end{example}

\begin{solutionexample}{}{}
\begin{itemize}
\item \textbf{Case 1}\\
For every $x\in\left[0,1\right],$
\[
m^{X=x}_{Y}=\intop_{\left[\epsilon,1\right]}y\text{d}\lambda\left(y\right)+\epsilon x.
\]
Hence,\boxeq{
\[
m^{X=x}_{Y}=\dfrac{1}{2}\left(1-\epsilon^{2}\right)+\epsilon x.
\]
}Thus the conditional mean of $Y$ given $X$ is a linear function
on $\left[0,1\right].$
\item \textbf{Case 2}\\
Using the previous computation, for every $x\in\left[0,\epsilon\right[,$
\[
m^{X=x}_{Y}=\dfrac{1}{2}\left(1-\epsilon^{2}\right)+\epsilon x.
\]
For every $x\in\left[\epsilon,1\right],$ $m^{X=x}_{Y}=\dfrac{1}{2}.$\\
Therefore, the conditional mean of $Y$ given $X$ is a piecewise
linear function on $\left[0,1\right],$ with a discontinuity at $\epsilon.$
It can be written\boxeq{
\[
m^{X=x}_{Y}=\boldsymbol{1}_{\left[0,\epsilon\right[}\left(x\right)\left[\dfrac{1}{2}\left(1-\epsilon^{2}\right)+\epsilon x\right]+\boldsymbol{1}_{\left[\epsilon,1\right]}\left(x\right)\dfrac{1}{2}.
\]
}
\end{itemize}
\end{solutionexample}

Taking the situation of conditional transfer theorem---Theorem $\ref{th:cond_transfer}$---we
now present a useful method for computing conditional means.

\begin{lemma}{Formula to Compute the Conditional Mean}{}

Let $\left(X,Y\right)$ be a random variable taking values in the
probabilizable space $\left(E\times F,\mathscr{E}\otimes\mathscr{F}\right),$
and assume that there exists a conditional law $\nu=P^{X=\cdot}_{Y}$
of $Y$ given $X.$ Let $f$ be a measurable function from $\left(E\times F,\mathscr{E}\otimes\mathscr{F}\right)$
into $\left(\mathbb{R},\mathscr{B}_{\mathbb{R}}\right).$ Assume that
$f\left(X,Y\right)\in\mathscr{L}^{1}\left(\Omega,\mathscr{A},P\right).$ 

Then, for $P_{X}-$almost every $x\in E,$\boxeq{
\[
m^{X=x}_{f\left(X,Y\right)}=m^{X=x}_{f\left(x,Y\right)}.
\]
}

In particular, if $X$ and $Y$ are \textbf{independent}, then for
$P_{X}-$almost every $x\in E,$\boxeq{
\[
m^{X=x}_{f\left(X,Y\right)}=\mathbb{E}\left(f\left(x,Y\right)\right).
\]
}

\end{lemma}

\begin{proof}{}{}

By definition of the conditional mean $m^{X=\cdot}_{f\left(X,Y\right)}$
and by the conditional transfer theorem---Theorem $\ref{th:cond_transfer}$---for
$P_{X}-$almost every $x\in E,$
\[
m^{X=x}_{f\left(X,Y\right)}=\intop_{\mathbb{R}}z\text{d}P^{X=x}_{f\left(x,Y\right)}\left(z\right),
\]
which yields the stated identity. In the case where $X$ and $Y$
are independent, it suffices to use the relation $\refpar{eq:cond_law_from_kernel_independence}.$

\end{proof}

\subsubsection*{Application of the Concept of Conditional Mean to the Regression
Problem}

We conclude this section by presenting an application of the concept
of conditional mean to the regression problem. This is a least-square
problem that generalizes the linear regression problem studied in
Part \ref{part:Introduction-to-Probability}.

\subsubsection*{The general problem}

\begin{leftbar}

Let $\left(X,Y\right)$ be a random variable taking values in the
probabilizable space $\left(E\times F,\mathscr{E}\otimes\mathscr{F}\right),$
and assume that there exists a conditional law $\nu=P^{X=\cdot}_{Y}$
of $Y$ given $X.$ We wish to estimate to what extent $Y$ is in
the neighborhood of a functional $X.$ This heuristic idea leads to
the following precise minimization problem, to which we restrict our
attention.

Assume that $F=\mathbb{R}$ equipped with its Borel $\sigma-$algebra---a
natural extension would be to take $F$ as a Euclidean space---and
that the random variable $Y$ admits a second-order moment. We seek
to solve the minimization problem\boxeq{
\begin{equation}
\min\left\{ \mathbb{E}\left(\left[Y-f\circ X\right]^{2}\right):\,f\in\mathscr{L}^{2}\left(E,\mathscr{E},P_{X}\right)\right\} .\label{eq:cond_linear_reg_problem}
\end{equation}
}

\end{leftbar}

\begin{remark}{}{}

To obtain a geometric interpretation of the problem, we may reformulate
it as a projection problem in the Hilbert space $\text{L}^{2}\left(\Omega,\mathscr{A},P\right).$
Temporarily assume that the subspace
\[
\Pi_{X}\equiv\left\{ \widetilde{f\circ X}:\,f\in\mathscr{L}^{2}\left(E,\mathscr{E},P_{X}\right)\right\} 
\]
is a closed subspace of $\text{L}^{2}\left(\Omega,\mathscr{A},P\right),$
where $\widetilde{f\circ X}$ denotes the equivalence class of $f\circ X.$
The solutions to the problem $\refpar{eq:cond_linear_reg_problem}$
are then precisely the representatives of the orthogonal projection
of the class of $Y$ onto $\Pi_{X}.$

\end{remark}

\begin{lemma}{}{}

The subspace $\Pi_{X}$ is closed in $\text{L}^{2}\left(\Omega,\mathcal{A},P\right).$

\end{lemma}

\begin{proof}{}{}

Let $\left(f_{n}\right)_{n\in\mathbb{N}}$ be a sequence of elements
in $\mathscr{L}^{2}\left(E,\mathscr{E},P_{X}\right)$ such that the
sequence $\left(\widetilde{f_{n}\circ X}\right)_{n\in\mathbb{N}}$
converges to $Z\in\text{L}^{2}\left(\Omega,\mathcal{A},P\right).$
The sequence $\left(f_{n}\circ X\right)_{n\in\mathbb{N}}$ is therefore
bounded in $\mathscr{L}^{2}\left(\Omega,\mathscr{A},P\right)$ by
a real number $c>0.$ There exists a subsequence $\left(f_{n_{k}}\right)_{k\in\mathbb{N}}$
such that the sequence $\left(f_{n_{k}}\circ X\right)_{k\in\mathbb{N}}$
converges $P-$almost surely---toward a representative of $Z.$ 

In particular, define $f=\limsup_{k\in\mathbb{N}}f_{n_{k}}.$ Then
the sequence $\left(f_{n_{k}}\circ X\right)_{k\in\mathbb{N}}$ converges
$P-$almost surely to $f\circ X.$ 

By the transfer theorem and the Fatou lemma,
\begin{align*}
\intop_{E}\left(f\left(x\right)\right)^{2}\text{d}P_{X}\left(x\right) & =\intop_{\Omega}\left(f\circ X\right)^{2}\text{d}P\\
 & =\intop_{\Omega}\lim_{k\to+\infty}\left(\left(f_{n_{k}}\circ X\right)^{2}\right)\text{d}P\\
 & \leqslant\lim_{k\to+\infty}\intop_{\Omega}\left(f_{n_{k}}\circ X\right)^{2}\text{d}P\leqslant c,
\end{align*}
which shows that $f\in\mathscr{L}^{2}\left(E,\mathscr{E},P_{X}\right).$ 

Moreover $\widetilde{f_{n}\circ X}=Z,$ hence, the lemma is proved.

\end{proof}

\begin{proposition}{Solution to the Regression Problem}{}

The conditional mean $m^{X=\cdot}_{Y}$ is a solution of the regression
problem $\refpar{eq:cond_linear_reg_problem}.$

\end{proposition}

\begin{proof}{}{}

By the extended Fubini theorem, for every $f\in\mathscr{L}^{2}\left(E,\mathscr{E},P_{X}\right),$
\[
\mathbb{E}\left(\left(Y-f\circ X\right)^{2}\right)=\intop_{E}\left[\intop_{F}\left(y-f\left(x\right)\right)^{2}\text{d}P^{X=x}_{Y}\left(y\right)\right]\text{d}P_{X}\left(x\right).
\]

Any solution $f_{0}$ to the problem $\refpar{eq:cond_linear_reg_problem}$
must therefore satisfy, for $P_{X}-$almost every $x,$
\[
\intop_{F}\left(y-f_{0}\left(x\right)\right)^{2}\text{d}P^{X=x}_{Y}\left(y\right)=\min\left\{ \intop_{F}\left(y-f\left(x\right)\right)^{2}\text{d}P^{X=x}_{Y}\left(y\right):\,f\in\mathscr{L}^{2}\left(E,\mathscr{E},P_{X}\right)\right\} .
\]
Thus, for $P_{X}-$almost every $x,$ the value $f_{0}\left(x\right)$
must be a stationnary point of the quadratic polynomial $Q,$ in the
real variable $z,$
\[
Q\left(z\right)=z^{2}-2z\intop_{F}y\text{d}P^{X=x}_{Y}\left(y\right).
\]
Hence,
\[
f_{0}\left(x\right)=\intop_{F}y\text{d}P^{X=x}_{Y}\left(y\right),
\]
and this stationary point clearly corresponds to a minimum.

\end{proof}

\section{Conditional Expectation}\label{sec:Conditional-Expectation}

The study of a random phenomenon leads us, given a certain level of
information, to adopt as a basic model, a probabilized space $\left(\Omega,\mathcal{A},P\right).$
If the available information is more limited, we may instead need
to work with a probabilized space $\left(\Omega,\mathscr{B},P\right),$
where $\mathscr{B}$ is a sub-$\sigma-$algebra of $\mathcal{A},$
that is, a $\sigma-$algebra satisfying $\mathscr{B}\subset\mathcal{A}.$
This situation arises in particular when the random phenomenon evolves
over time, since the information available typically increases with
time. 

If $Y$ is a random variable defined on the probabilized space $\left(\Omega,\mathcal{A},P\right),$
how can we compute its expectation on the elements of $\mathscr{B},$
using only $\mathscr{B}-$measurable random variables? The \textbf{conditional
expectation}, a fundamental tool in probability theory, provides precisely
this answer. One may say that its use enables a ``progressive''
computation of expectations, a feature that will appear repeatedly
in the study of martingales---see Chapter \ref{chap:Processus-and-Discrete}---and
Markov chains---see Chapter \ref{chap:PartIIChap17}.

\textbf{In this section, we consider a probabilized space $\left(\Omega,\mathcal{A},P\right)$
together with a sub-$\sigma-$algebra $\mathscr{B}$ of $\mathcal{A}.$
Unless explicitly stated otherwise, we use the same symbol to denote
a random variable $X$ and its equivalence class $\widetilde{X}.$}

\subsection{Conditional Expectation as an Orthogonal Projector in $\text{L}^{2}\left(\Omega,\mathcal{A},P\right)$}

\begin{lemma}{Characterization of the Orthogonal Projection}{}

The subspace $\text{L}^{2}\left(\Omega,\mathscr{B},P\right)$ is closed
in the Hilbert space $\text{L}^{2}\left(\Omega,\mathcal{A},P\right).$
The orthogonal projector onto $\text{L}^{2}\left(\Omega,\mathscr{B},P\right)$
is denoted by $\mathbb{E}^{\mathscr{B}}$---Tr.N. Some authors denote
it $\mathbb{E}\left(\cdot\,|\,\mathscr{B}\right).$ 

For $Y\in\text{L}^{2}\left(\Omega,\mathcal{A},P\right)$, its orthogonal
projection $\mathbb{E}^{\mathscr{B}}\left(Y\right)$ is characterized
by the orthogonality relation\boxeq{
\begin{equation}
\forall Z\in\text{L}^{2}\left(\Omega,\mathscr{B},P\right),\,\,\,\,\text{\ensuremath{\mathbb{E}\left(ZY\right)}}=\mathbb{E}\left(Z\mathbb{E}^{\mathscr{B}}\left(Y\right)\right).\label{eq:exp_of_prod_and__orthog_project}
\end{equation}
}

\end{lemma}

\begin{proof}{}{}

The subspace $\text{L}^{2}\left(\Omega,\mathscr{B},P\right)$ is complete
and therefore closed in $\text{L}^{2}\left(\Omega,\mathcal{A},P\right).$
The relation $\refpar{eq:exp_of_prod_and__orthog_project}$ is simply
the expression of the fact that $Y-\mathbb{E}^{\mathscr{B}}\left(Y\right)$
is orthogonal to the subspace $\text{L}^{2}\left(\Omega,\mathscr{B},P\right).$

\end{proof}

\begin{remark}{}{}

The uniqueness of the orthogonal projection onto a closed subspace
of a Hilbert space implies that $\mathbb{E}^{\mathscr{B}}\left(Y\right)$
is the unique equivalence class $U$ of $\mathscr{B}-$measurable
random variables such that
\[
\forall Z\in\text{L}^{2}\left(\Omega,\mathscr{B},P\right),\,\,\,\,\mathbb{E}\left(ZY\right)=\mathbb{E}\left(ZU\right).
\]

\textbf{This uniqueness property is often what allows one to identify
the conditional expectation explicitly.}

\end{remark}

\begin{definition}{Conditional Expectation as an Orthogonal Projector in $\textrm{L}^2\left(\Omega,\mathcal{A},P\right)$}{}

Let $Y\in\text{L}^{2}\left(\Omega,\mathcal{A},P\right)$ be an equivalence
class of random variables. The class of random variables $\mathbb{E}^{\mathscr{B}}\left(Y\right)$
is called the \textbf{conditional expectation of $Y$ given $\mathscr{B}.$\index{conditional expectation of ensuremath{Y} given $ensuremath{mathscr{B}}@conditional expectation of $Y$ given $\mathscr{B}$}}

Let $Y\in\mathscr{L}^{2}\left(\Omega,\mathcal{A},P\right)$ be an
actual random variable. The class of random variables $\mathbb{E}^{\mathscr{B}}\left(\widetilde{Y}\right)$
is also called the \textbf{conditional expectation of $Y$ given $\mathscr{B}$\index{conditional expectation of ensuremath{Y} given $ensuremath{mathscr{B}}@conditional expectation of $Y$ given $\mathscr{B}$}}
and is denoted by $\mathbb{E}^{\mathscr{B}}\left(Y\right).$

\end{definition}

\begin{denotation}{}{}

We therefore speak interchangeably of the conditional expectation
of a random variable or of an equivalence class of random variables:
in both cases, what is meant is an equivalence class. By abuse of
language, when there is no risk of confusion, $\mathbb{E}^{\mathscr{B}}\left(Y\right)$
may denote an arbitrary representative of this equivalence class---often
called \textbf{a version of the conditional expectation}\index{version of a conditional expectation}\mindex{conditional expectation!version}. 

To indicate that a random variable $U$ is a version of $\mathbb{E}^{\mathscr{B}}\left(Y\right)$,
we write
\[
U=\mathbb{E}^{\mathscr{B}}\left(Y\right)\,\,\,\,P-\text{almost surely.}
\]

\end{denotation}

\begin{proposition}{}{}

Let $Y\in\text{L}^{2}\left(\Omega,\mathcal{A},P\right).$ The relation
$\refpar{eq:exp_of_prod_and__orthog_project}$ is equivalent to the
relation\boxeq{
\begin{equation}
\forall B\in\mathscr{B},\,\,\,\,\mathbb{E}\left(\boldsymbol{1}_{B}Y\right)=\mathbb{E}\left(\boldsymbol{1}_{B}\mathbb{E}^{\mathscr{B}}\left(Y\right)\right),\label{eq:exp_of_prod_with_indicator_orth_proj}
\end{equation}
}which can also be written as\boxeq{
\begin{equation}
\forall B\in\mathscr{B},\,\,\,\,\intop_{B}Y\text{d}P=\intop_{B}\mathbb{E}^{\mathscr{B}}\left(Y\right)\text{d}P.\label{eq:exp_of_prod_int_vers_orthog_project}
\end{equation}
}

\end{proposition}

\begin{proof}{}{}

Let $Y\in\text{L}^{2}\left(\Omega,\mathcal{A},P\right).$ The implication
$\refpar{eq:exp_of_prod_and__orthog_project}\Rightarrow\refpar{eq:exp_of_prod_with_indicator_orth_proj}$
is immediate. 

Conversely, assume that $\refpar{eq:exp_of_prod_with_indicator_orth_proj}$
holds. Then by linearity, the relation $\refpar{eq:exp_of_prod_and__orthog_project}$
holds for every $\mathscr{B-}$measurable step random variable $Z.$
Since such simple random variables form a dense subset of $\text{L}^{2}\left(\Omega,\mathscr{B},P\right),$
the result extends by continuity to all $Z\in\text{L}^{2}\left(\Omega,\mathscr{B},P\right).$
Indeed the functions $Z\mapsto\mathbb{E}\left(ZY\right)\,\,\,\,\text{and}\,\,\,\,Z\mapsto\mathbb{E}\left(Z\mathbb{E}^{\mathscr{B}}\left(Y\right)\right)$
are continuous linear functionals on $\text{L}^{2}\left(\Omega,\mathscr{B},P\right)$
by the Cauchy-Schwarz inequality. 

\end{proof}

\begin{remark}{}{}

The conditional expectation is therefore characterized as the unique
equivalence class $U\in\text{L}^{2}\left(\Omega,\mathscr{B},P\right)$
of $\mathscr{B}-$measurable random variables such that
\[
\forall B\in\mathscr{B},\,\,\,\,\intop_{B}Y\text{d}P=\intop_{B}U\text{d}P.
\]

\end{remark}

\subsubsection*{Preliminary Note on the Order Defined on the Set of Equivalence Classes
of Random Variables}

The equivalence relation of ``$P-$almost sure equality'' is compatible
with the usual partial order on real-valued random variables---or
random variables taking values in $\overline{\mathbb{R}}.$ It therefore
induces a partial order on the corresponding equivalence classes,
which we continue to denote by $\leqslant.$ 

In particular, let $\left(X_{n}\right)_{n\in\mathbb{N}}$ and $\left(Y_{n}\right)_{n\in\mathbb{N}}$
be sequences of random variables such that
\[
\forall n\in\mathbb{N},\,\,\,\,\,P-\text{almost surely,}\,\,\,\,X_{n}=Y_{n}.
\]
Since a countable union of sets of probability zero still has probability
zero, it follows that
\[
P-\text{almost surely,}\,\,\,\,\forall n\in\mathbb{N},\,\,\,\,X_{n}=Y_{n}.
\]
Consequently, for random variables taking values in $\overline{\mathbb{R}},$
we obtain the $P-$almost sure equalities 
\[
P-\text{almost surely,}\,\,\,\,\sup_{n\in\mathbb{N}}X_{n}=\sup_{n\in\mathbb{N}}Y_{n}\,\,\,\,\text{and}\,\,\,\,\inf_{n\in\mathbb{N}}X_{n}=\inf_{n\in\mathbb{N}}Y_{n}.
\]

\begin{remark}{}{}

It is worth noting that this is no longer true when dealing with uncountable
families of random variables, since the supremum and infinimum may
then may fail to be random variables at all!

\end{remark}

\begin{proposition}{Non-Negativity of $\mathbb{E}^\mathscr{B}$}{}

The operator $\mathbb{E}^{\mathscr{B}}$ is a continuous linear operator
on $\text{L}^{2}\left(\Omega,\mathscr{A},P\right)$ with norm 1. Moreover,
it is \textbf{nonnegative}, in the sense that it satisfies the implication\boxeq{
\[
Y\geqslant0\Rightarrow\mathbb{E}^{\mathscr{B}}\left(Y\right)\geqslant0.
\]
}

In particular, if $Y_{1},Y_{2}\in\text{L}^{2}\left(\Omega,\mathscr{A},P\right)$
are such that $Y_{1}\leqslant Y_{2},$ then\boxeq{
\[
\mathbb{E}^{\mathscr{B}}\left(Y_{1}\right)\leqslant\mathbb{E}^{\mathscr{B}}\left(Y_{2}\right).
\]
}

\end{proposition}

\begin{proof}{}{}

This is a standard property of the orthogonal projectors. The non-negativity
follows from the fact that if $Y\geqslant0,$ then for every $B\in\mathscr{B},$
\[
\intop_{B}\mathbb{E}^{\mathscr{B}}\left(Y\right)\text{d}P\geqslant0,
\]
which is equivalent to saying that $\mathbb{E}^{\mathscr{B}}\left(Y\right)\geqslant0.$

\end{proof}

\begin{proposition}{Properties of the Orthogonal Projector $\mathbb{E}^\mathscr{B}$}{orthog_proj_prop}

Let $Y\in\text{L}^{2}\left(\Omega,\mathscr{A},P\right).$ The following
properties hold:

(a) $\mathbb{E}\left(\mathbb{E}^{\mathscr{B}}\left(Y\right)\right)=\mathbb{E}\left(Y\right).$

(b) If $Y$ is $\mathscr{B}-$measurable, then\boxeq{
\[
\mathbb{E}^{\mathscr{B}}\left(Y\right)=Y.
\]
}

(c) If $Z$ is $\mathscr{B}-$measurable and bounded, then\boxeq{
\[
\mathbb{E}^{\mathscr{B}}\left(ZY\right)=Z\mathbb{E}^{\mathscr{B}}\left(Y\right)\,\,\,\,\,P-\text{almost surely.}
\]
}

(d) \textbf{Theorem of the three perpendiculars}

If $\mathscr{B}_{1}$ and $\mathscr{B}_{2}$ are two sub-$\sigma$-algebra
with $\mathscr{B}_{1}\subset\mathscr{B}_{2},$ then\boxeq{
\begin{equation}
E^{\mathscr{B}_{1}}\left(Y\right)=E^{\mathscr{B}_{1}}\left(E^{\mathscr{B}_{2}}\left(Y\right)\right).\label{eq:three_perp_cond_oper}
\end{equation}
}

(e) $\left|\mathbb{E}^{\mathscr{B}}\left(Y\right)\right|\leqslant\mathbb{E}^{\mathscr{B}}\left(\left|Y\right|\right).$

(f) The operator $\mathbb{E}^{\mathscr{B}}$ from $\text{L}^{2}\left(\Omega,\mathscr{A},P\right)$
to $\text{L}^{2}\left(\Omega,\mathscr{B},P\right)$ has norm 1 with
respect to the $\text{L}^{1}-$norm. That is, for every $Y\in\text{L}^{2}\left(\Omega,\mathscr{A},P\right),$
\begin{equation}
\left\Vert \mathbb{E}^{\mathscr{B}}\left(Y\right)\right\Vert _{1}\leqslant\left\Vert Y\right\Vert _{1}.\label{eq:norm_1_cond_op_bound}
\end{equation}

\end{proposition}

\begin{proof}{}{}

(a) It suffices to take $Z=1,$ which is indeed $\mathscr{B}-$measurable,
in the relation $\refpar{eq:exp_of_prod_and__orthog_project}.$

(b) Since $Y$ is $\mathscr{B}-$measurable, then $Y$ is in the subspace
$\text{L}^{2}\left(\Omega,\mathscr{B},P\right),$ so its orthogonal
projection onto this subspace is itself.

(c) If $Z$ is bounded, then $YZ\in\text{L}^{2}\left(\Omega,\mathscr{A},P\right).$
For every $T\in\text{L}^{2}\left(\Omega,\mathscr{B},P\right),$ the
definition of $\mathbb{E}^{\mathscr{B}}\left(ZY\right)$ gives 
\[
\mathbb{E}\left(T\mathbb{E}^{\mathscr{B}}\left(ZY\right)\right)=\mathbb{E}\left(TZY\right).
\]
Since $ZT$ is $\mathscr{B}-$measurable, the definition of $\mathbb{E}^{\mathscr{B}}\left(Y\right)$
yields 
\[
\mathbb{E}\left(T\mathbb{E}^{\mathscr{B}}\left(ZY\right)\right)=\mathbb{E}\left(TZ\mathbb{E}^{\mathscr{B}}\left(Y\right)\right),
\]
which can be also read as
\[
\mathbb{E}\left(T\mathbb{E}^{\mathscr{B}}\left(ZY\right)\right)=\mathbb{E}\left(T\left[Z\mathbb{E}^{\mathscr{B}}\left(Y\right)\right]\right).
\]

Since $Z\mathbb{E}^{\mathscr{B}}\left(Y\right)\in\text{L}^{2}\left(\Omega,\mathscr{B},P\right)$
and combined with the first remark with $T=1,$ it allows to conclude.

(d) This is a general Hilbert space property---known as the theorem
of the three perpendiculars.

Since $\mathscr{B}_{1}\subset\mathscr{B}_{2},$ the space $\text{L}^{2}\left(\Omega,\mathscr{B}_{1},P\right)$
is a closed subspace of $\text{L}^{2}\left(\Omega,\mathscr{B}_{2},P\right).$
Let us show this property in this context: for every $Z\in\text{L}^{2}\left(\Omega,\mathscr{B}_{1},P\right),$
the variable $Z$ is $\mathscr{B}_{2}-$measurable and thus
\[
\mathbb{E}\left(ZY\right)=\mathbb{E}\left(Z\mathbb{E}^{\mathscr{B}}\left(Y\right)\right).
\]
By the definition of the projection $E^{\mathscr{B}_{1}}\left(E^{\mathscr{B}_{2}}\left(Y\right)\right),$
we obtain
\[
\mathbb{E}\left(ZY\right)=\mathbb{E}\left(ZE^{\mathscr{B}_{1}}\left(E^{\mathscr{B}_{2}}\left(Y\right)\right)\right),
\]
which establishes the identity.

(e) Since the conditional expectation is linear, we use the convexity
of the absolute value function by writing that it is the superior
enveloppe of its linear lower-bounding. More precisely, by taking
only the extremes, if $A=\left\{ -1,1\right\} ,$ then
\[
\forall x\in\mathbb{R},\,\,\,\,\left|x\right|=\sup_{a\in A}\left(ax\right).
\]
This argument of convexity will be used later to prove the \textbf{\index{Jensen inequality}Jensen
inequality}. Hence,
\[
\forall a\in A,\,\,\,\,P-\text{almost surely,}\,\,\,\,a\mathbb{E}^{\mathscr{B}}\left(Y\right)=\mathbb{E}^{\mathscr{B}}\left(aY\right)\leqslant\mathbb{E}^{\mathscr{B}}\left(\left|Y\right|\right),
\]
and thus,
\[
P-\text{almost surely,}\,\,\,\,\forall a\in A,\,\,\,\,a\mathbb{E}^{\mathscr{B}}\left(Y\right)\leqslant\mathbb{E}^{\mathscr{B}}\left(\left|Y\right|\right),
\]
which implies that
\[
P-\text{almost surely,}\,\,\,\,\left|\mathbb{E}^{\mathscr{B}}\left(Y\right)\right|=\sup_{a\in A}\left(a\mathbb{E}^{\mathscr{B}}\left(Y\right)\right)\leqslant\mathbb{E}^{\mathscr{B}}\left(\left|Y\right|\right).
\]

(f) This follows by integrating the inequality obtained in (e).

\end{proof}

\subsection{Extension of the Definition of the Conditional Expectation to $\text{L}^{1}\left(\Omega,\mathscr{A},P\right)$}

\begin{proposition}{Conditional Expectation of an $\textrm{L}^1\left(\Omega,\mathcal{A},P\right)$ Random Variable}{}

Let $Y\in\text{L}^{1}\left(\Omega,\mathscr{A},P\right)$---or $\mathscr{L}^{1}\left(\Omega,\mathscr{A},P\right).$
Then, there exists a unique equivalence class of $\mathscr{B}-$measurable
random variables $U\in\text{L}^{1}\left(\Omega,\mathscr{B},P\right)$
satisfying
\begin{equation}
\forall B\in\mathscr{B},\,\,\,\,\intop_{B}Y\text{d}P=\intop_{B}U\text{d}P.\label{eq:cond_exp_Y_knwoing_B_L1}
\end{equation}

This random variable $U$ is denoted by $\mathbb{E}^{\mathscr{B}}\left(Y\right)$
and is called the \textbf{\index{conditional expectation of ensuremath{Y} given $ensuremath{mathscr{B}}@conditional expectation of $Y$ given $\mathscr{B}$}conditional
expectation of $Y$ given $\mathscr{B}.$} Moreover, it satisfies
the inequality
\begin{equation}
\left\Vert \mathbb{E}^{\mathscr{B}}\left(Y\right)\right\Vert _{1}\leqslant\left\Vert Y\right\Vert _{1}.\label{eq:cond_exp_norm1_bound}
\end{equation}

\end{proposition}

\begin{proof}{}{}

We reduce the problem to the $\text{L}^{2}\left(\Omega,\mathscr{A},P\right)$
setting as follows. Let $\left(Y_{n}\right)_{n\in\mathbb{N}}$ be
the sequence defined by
\[
\forall n\in\mathbb{N},\,\,\,\,Y_{n}=\boldsymbol{1}_{\left(\left|Y\right|\leqslant n\right)}Y.
\]
For every $n\in\mathbb{N},$ $Y_{n}\in\text{L}^{2}\left(\Omega,\mathscr{A},P\right),$
and
\[
\left|Y_{n}-Y\right|\leqslant\left|Y\right|.
\]
Since the sequence $\left(Y_{n}\right)_{n\in\mathbb{N}}$ converges
$P-$almost surely to $Y,$ the dominated convergence theorem implies
that it also converges to $Y$ in $\text{L}^{1}\left(\Omega,\mathscr{A},P\right).$ 

For each $n\in\mathbb{N},$ define 
\[
Z_{n}=\mathbb{E}^{\mathscr{B}}\left(Y_{n}\right)\in\text{L}^{2}\left(\Omega,\mathscr{B},P\right).
\]
By inequality $\refpar{eq:norm_1_cond_op_bound},$ for every $n,m\in\mathbb{N},$
\[
\left\Vert Z_{n}-Z_{m}\right\Vert _{1}\leqslant\left\Vert Y_{n}-Y_{m}\right\Vert _{1}.
\]

Since the sequence $\left(Y_{n}\right)_{n\in\mathbb{N}}$ converges
in $\text{L}^{1}\left(\Omega,\mathscr{A},P\right),$ it is a Cauchy
sequence, and therefore so is the sequence $\left(Z_{n}\right)_{n\in\mathbb{N}}.$
Because the space $\text{L}^{1}\left(\Omega,\mathscr{A},P\right)$
is complete, the sequence $\left(Z_{n}\right)_{n\in\mathbb{N}}$ converges
in $\text{L}^{1}\left(\Omega,\mathscr{A},P\right)$ to some limit
$Z.$ More precisely, since for each $n\in\mathbb{N},$ $Z_{n}$ is
$\mathscr{B}-$measurable, $Z\in\text{L}^{1}\left(\Omega,\mathscr{B},P\right).$
Moreover, for every $B\in\mathscr{B}$ and every $n\in\mathbb{N},$
\[
\intop_{B}Y_{n}\text{d}P=\intop_{B}Z_{n}\text{d}P.
\]

The convergence of sequences $\left(Y_{n}\right)_{n\in\mathbb{N}}$
and $\left(Z_{n}\right)_{n\in\mathbb{N}}$ in $\text{L}^{1}\left(\Omega,\mathscr{A},P\right)$
allows to pass to the limit, which yields the relation
\[
\forall B\in\mathscr{B},\,\,\,\,\intop_{B}Y\text{d}P=\intop_{B}Z\text{d}P.
\]
This proves existence. Uniqueness is immediate. It remains to show
the inequality $\refpar{eq:cond_exp_norm1_bound}.$ For each $n,$
we apply the inequality $\refpar{eq:norm_1_cond_op_bound}$ to $Y_{n}$
and then pass to the limit. We obtain
\[
\mathbb{E}\left(\left|Z\right|\right)=\lim_{n\to+\infty}\mathbb{E}\left(\left|Z_{n}\right|\right)\leqslant\lim_{n\to+\infty}\mathbb{E}\left(\left|Y_{n}\right|\right)=\mathbb{E}\left(\lim_{n\to+\infty}\left|Y_{n}\right|\right)=\mathbb{E}\left(\left|Y\right|\right).
\]

\end{proof}

\begin{proposition}{Properties of $\mathbb{E}^\mathscr{B}$}{ext_cond_exp_prop}

The operator $\mathbb{E}^{\mathscr{B}}$ is linear and continuous
on $\text{L}^{1}\left(\Omega,\mathscr{A},P\right)$ with norm 1. Moreover,
it is nonnegative.

Let $Y\in\text{L}^{1}\left(\Omega,\mathscr{A},P\right).$ Then the
following properties hold:

(a) $\mathbb{E}\left(\mathbb{E}^{\mathscr{B}}\left(Y\right)\right)=\mathbb{E}\left(Y\right).$

(b) If $Y$ is $\mathscr{B}-$measurable, then\boxeq{
\[
\mathbb{E}^{\mathscr{B}}\left(Y\right)=Y.
\]
}

(c) If $Z$ is $\mathscr{B}-$measurable and bounded, then\boxeq{
\[
\mathbb{E}^{\mathscr{B}}\left(ZY\right)=Z\mathbb{E}^{\mathscr{B}}\left(Y\right)\,\,\,\,\,P-\text{almost surely.}
\]
}

(d) Theorem of the three perpendiculars

If $\mathscr{B}_{1}$ and $\mathscr{B}_{2}$ are two sub$-\sigma-$algebra
such that $\mathscr{B}_{1}\subset\mathscr{B}_{2},$ then\boxeq{
\begin{equation}
E^{\mathscr{B}_{1}}\left(Y\right)=E^{\mathscr{B}_{1}}\left(E^{\mathscr{B}_{2}}\left(Y\right)\right).\label{eq:three_perp_cond_oper-L1}
\end{equation}
}

(e) $\left|\mathbb{E}^{\mathscr{B}}\left(Y\right)\right|\leqslant\mathbb{E}^{\mathscr{B}}\left(\left|Y\right|\right).$

\end{proposition}

\begin{proof}{}{}

Linearity follows from the characterization $\refpar{eq:cond_exp_Y_knwoing_B_L1}$
of $\mathbb{E}^{\mathscr{B}}\left(Y\right).$ The inequality $\refpar{eq:norm_1_cond_op_bound}$
shows that $\mathbb{E}^{\mathscr{B}}$ is continuous with operator
norm less than or equal to 1. This norm is in fact equal to 1, since,
if $Y\in\text{L}^{1}\left(\Omega,\mathscr{B},P\right),$ the characterization
$\refpar{eq:cond_exp_Y_knwoing_B_L1}$ implies that $\mathbb{E}^{\mathscr{B}}\left(Y\right)=Y.$ 

All remaining properties follow either by continuity from the corresponding
properties in $\text{L}^{2}\left(\Omega,\mathscr{A},P\right),$ or
directly by applying the characterization $\refpar{eq:cond_exp_Y_knwoing_B_L1}$
together with the same arguments as in Proposition $\ref{pr:orthog_proj_prop}.$

\end{proof}

\begin{remark}{}{}

Since $\mathbb{E}^{\mathscr{B}}$ is continuous on $\text{L}^{1}\left(\Omega,\mathscr{A},P\right),$
if a sequence $\left(X_{n}\right)_{n\in\mathbb{N}}$ converges to
$X$ in $\text{L}^{1}\left(\Omega,\mathscr{A},P\right),$ then the
sequence $\left(\mathbb{E}^{\mathscr{B}}\left(X_{n}\right)\right)_{n\in\mathbb{N}}$
converges to $\mathbb{E}^{\mathscr{B}}\left(X\right)$ in $\text{L}^{1}\left(\Omega,\mathscr{A},P\right).$

\end{remark}

\begin{proposition}{Conditional Expectation and Independence}{}

Let $Y\in\text{L}^{1}\left(\Omega,\mathscr{A},P\right).$ If $Y$
and $\mathscr{B}$ are independent---that is, if the $\sigma-$algebras
$\sigma\left(Y\right)$ and $\mathscr{B}$ are independent---then
\[
\mathbb{E}^{\mathscr{B}}\left(Y\right)=E\left(Y\right)\,\,\,\,P-\text{almost surely}.
\]

\end{proposition}

\begin{proof}{}{}

For every $B\in\mathscr{B},$ the random variables $\boldsymbol{1}_{B}$
and $Y$ are independent, 
\[
\mathbb{E}\left(\boldsymbol{1}_{B}Y\right)=\mathbb{E}\left(\boldsymbol{1}_{B}\right)\mathbb{E}\left(Y\right).
\]
Therefore,
\[
\mathbb{E}\left(\boldsymbol{1}_{B}Y\right)=\mathbb{E}\left(\boldsymbol{1}_{B}\mathbb{E}\left(Y\right)\right).
\]
Invoking the characterization $\refpar{eq:cond_exp_Y_knwoing_B_L1}$
of the conditional expectation $\mathbb{E}^{\mathscr{B}}\left(Y\right)$
concludes the proof.

\end{proof}

\begin{remark}{}{}

Note that the equality obtained above is an equality of equivalence
classes, and that $\mathbb{E}\left(Y\right)$ denotes the class of
random variables $P-$almost surely equal to $\mathbb{E}\left(Y\right).$

\end{remark}

The following extension of this property is frequently used in computations.

\begin{proposition}{Extension of  Conditional Expectation under Independence}{cond_exp_exist_indep}

Let $\left(X,Y\right)$ be a random variable taking values in a probabilizable
space $\left(E\times F,\mathscr{E}\otimes\mathscr{F}\right),$ and
let $f\in\mathscr{L}^{1}\left(E\times F,\mathscr{E}\otimes\mathscr{F},P_{\left(X,Y\right)}\right).$
Suppose that $X$ is $\mathscr{B}-$measurable and that $Y$ and $\mathscr{B}$
are independent. The function $\widehat{f}$ defined by
\[
\forall x\in E,\,\,\,\,\widehat{f}\left(x\right)=\mathbb{E}\left(f\left(x,Y\right)\right)
\]
is $\mathscr{E}-$measurable, and
\begin{equation}
\mathbb{E}^{\mathscr{B}}\left(f\left(x,Y\right)\right)=\widehat{f}\circ X\,\,\,\,P-\text{almost surely}.\label{eq:conditional_exp_and_L1_app}
\end{equation}

\end{proposition}

\begin{proof}{}{}

Since the random variables $X$ and $Y$ are independent,
\[
P_{\left(X,Y\right)}=P_{X}\otimes P_{Y}.
\]

Because
\[
\forall x\in E,\,\,\,\,\widehat{f}\left(x\right)=\intop_{F}f\left(x,y\right)\text{d}P_{Y}\left(y\right),
\]
the measurability property of $\widehat{f}$ follows from the Fubini
theorem. Moreover, the vector space generated by the functions $\left(x,y\right)\mapsto g\left(x\right)h\left(y\right)$
where $g\in\mathscr{L}^{1}\left(E,\mathscr{E},P_{X}\right)$ and $h\in\mathscr{L}^{1}\left(F,\mathscr{F},P_{Y}\right)$
is dense in $\mathscr{L}^{1}\left(E\times F,\mathscr{E}\otimes\mathscr{F},P_{\left(X,Y\right)}\right).$
Additionally, the functions
\[
f\mapsto\widehat{f}\circ X\,\,\,\,\text{and}\,\,\,\,f\mapsto\mathbb{E}^{\mathscr{B}}\left(f\left(X,Y\right)\right)
\]
are continuous from $\mathscr{L}^{1}\left(E\times F,\mathscr{E}\otimes\mathscr{F},P_{\left(X,Y\right)}\right)$
into $\text{L}^{1}\left(\Omega,A,P\right).$ Indeed, by the transfer
theorem and the Fubini theorem,
\begin{align*}
\left\Vert \widehat{f}\circ X\right\Vert _{1} & =\intop_{\Omega}\left|\intop_{F}f\left(X,y\right)\text{d}P_{Y}\left(y\right)\right|\text{d}P\\
 & \leqslant\intop_{E}\left[\left|f\left(x,y\right)\right|\text{d}P_{Y}\left(y\right)\right]\text{d}P_{X}\left(x\right)=\intop_{E\times F}\left|f\right|\text{d}\left(P_{X}\otimes P_{Y}\right).
\end{align*}

Since $P_{\left(X,Y\right)}=P_{X}\otimes P_{Y},$
\[
\left\Vert \widehat{f}\circ X\right\Vert _{1}\leqslant\left\Vert f\right\Vert _{1}.
\]

The continuity of the function $f\mapsto\mathbb{E}^{\mathscr{B}}\left(f\left(X,Y\right)\right)$
follows from the transfer theorem and the inequality
\[
\left\Vert \mathbb{E}^{\mathscr{B}}\left(f\left(X,Y\right)\right)\right\Vert _{1}\leqslant\left\Vert f\left(X,Y\right)\right\Vert _{1}=\left\Vert f\right\Vert _{1}.
\]

By linearity, it therefore suffices to establish the relation $\refpar{eq:conditional_exp_and_L1_app}$
for functions $f$ direct product of $g\in\mathscr{L}^{1}\left(E,\mathscr{E},P_{X}\right)$
and $h\in\mathscr{L}^{1}\left(F,\mathscr{F},P_{Y}\right).$ 

Since $g\circ X$ is $\mathscr{B}-$measurable and $h\circ Y$ is
independent of $\mathscr{B},$
\[
\mathbb{E}^{\mathscr{B}}\left(f\left(X,Y\right)\right)=\mathbb{E}^{\mathscr{B}}\left(\left(g\circ X\right)\left(h\circ Y\right)\right)=\left(g\circ X\right)\mathbb{E}^{\mathscr{B}}\left(h\circ Y\right)=\left(g\circ X\right)\mathbb{E}\left(h\circ Y\right).
\]
Hence,\\
\[
\mathbb{E}^{\mathscr{B}}\left(f\left(X,Y\right)\right)=\left(g\circ X\right)\intop_{F}\left[h\left(y\right)\right]\text{d}P_{Y}\left(y\right)=\intop_{F}\left(g\circ X\right)\left[h\left(y\right)\right]\text{d}P_{Y}\left(y\right)=\widehat{f}\circ X.
\]

\end{proof}

\begin{example}{}{Poisson_and_B_meas_indep}

Let $X$ and $Y$ be two real-valued random variables. Assume that
$Y$ follows a Poisson law with parameter $\lambda>0,$ and that $X$
is $\mathscr{B}-$measurable, and that $Y$ and $\mathscr{B}$ are
independent.

Compute $\mathbb{E}^{\mathscr{B}}\left(\cos\left(XY\right)\right).$

\end{example}

\begin{solutionexample}{}{}

We have
\[
\widehat{f}\left(x\right)=\mathbb{E}\left(\cos\left(xY\right)\right)=\text{e}^{-\lambda}\sum^{+\infty}_{k=0}\dfrac{\lambda^{k}}{k!}\cos\left(kx\right).
\]

As
\[
\sum^{+\infty}_{k=0}\dfrac{\lambda^{k}}{k!}\cos\left(kx\right)=\text{Re}\left(\sum^{+\infty}_{k=0}\dfrac{\lambda^{k}}{k!}\text{e}^{\text{i}kx}\right)=\text{Re}\left(\lambda\text{e}^{\text{i}x}\right).
\]
This yields\boxeq{
\[
\mathbb{E}^{\mathscr{B}}\left[\cos\left(XY\right)\right]=\text{e}^{-\lambda\left[1-\cos\left(X\right)\right]}\cos\left(\lambda\sin X\right).
\]
}

\end{solutionexample}

\subsection{Extension of the Definition of the Conditional Expectation to $\mathscr{M}^{+}\left(\mathscr{A}\right)$}

Denote $\mathscr{M}^{+}\left(\mathscr{B}\right)$ the set of $\mathscr{B}-$measurable
random variables taking values in $\overline{\mathbb{R}}^{+}.$

\begin{proposition}{Conditional Expectation of a Random Variable in $\mathcal{M}^+\left(\mathcal{A}\right)$}{}

Let $Y\in\mathscr{M}^{+}\left(\mathscr{A}\right).$ There exists a
unique equivalence class $U$ of elements of $\mathscr{M}^{+}\left(\mathscr{B}\right)$
satisfying
\begin{equation}
\forall B\in\mathscr{B},\,\,\,\,\intop_{B}Y\text{d}P=\intop_{B}U\text{d}P.\label{eq:unique_class_exp_cond}
\end{equation}
This class is still denoted $\mathbb{E}^{\mathscr{B}}\left(Y\right)$
and is called the \textbf{conditional expectation of Y given $\mathscr{B}.$\index{conditional expectation of ensuremath{Y} given $ensuremath{mathscr{B}}@conditional expectation of $Y$ given $\mathscr{B}$}}

\end{proposition}

\begin{proof}{}{}

Uniqueness is proved exactly as in the $\text{L}^{1}$ case. For existence,
define, for each $n\in\mathbb{N},$ the bounded random variable $Y_{n}=\inf\left(Y,n\right)$
and let $U_{n}$ be a version of $\mathbb{E}^{\mathscr{B}}\left(Y_{n}\right).$ 

The sequence $\left(Y_{n}\right)_{n\in\mathbb{N}}$ is nondecreasing
and converges $P-$almost surely to $Y.$ The sequence $\left(U_{n}\right)_{n\in\mathbb{N}}$
non-decreases, thus converges in $\overline{\mathbb{R}}^{+}$ to a
$\mathscr{B}-$measurable limit $U.$ 

By the Beppo Levi property, for every $B\in\mathscr{B},$
\[
\intop_{B}Y\text{d}P=\lim_{n\to+\infty}\intop_{B}Y_{n}\text{d}P=\lim_{n\to+\infty}\intop_{B}U_{n}\text{d}P=\intop_{B}U\text{d}P.
\]

\end{proof}

\begin{proposition}{Conditional Beppo Levi Property}{}

Let $Y,Z\in\mathscr{M}^{+}\left(\mathscr{A}\right)$ be such that
$Y\leqslant Z.$ Then 
\[
\mathbb{E}^{\mathscr{B}}\left(Y\right)\leqslant\mathbb{E}^{\mathscr{B}}\left(Z\right).
\]
Moreover, the conditional \textbf{Beppo Levi} property holds: if a
sequence $\left(Y_{n}\right)_{n\in\mathbb{N}}$ of elements of $\mathscr{M}^{+}\left(\mathscr{A}\right)$
converges by nondecreasing to $Y,$ the sequence $\left(\mathbb{E}^{\mathscr{B}}\left(Y_{n}\right)\right)_{n\in\mathbb{N}}$
converges by nondecreasing to $\mathbb{E}^{\mathscr{B}}\left(Y\right).$

\end{proposition}

\begin{proof}{}{}

For the first property, note that for every $n\in\mathbb{N},$ $\inf\left(Y,n\right)\leqslant\inf\left(Z,n\right)$
and therefore
\[
\mathbb{E}^{\mathscr{B}}\left(\inf\left(Y,n\right)\right)\leqslant\mathbb{E}^{\mathscr{B}}\left(\inf\left(Z,n\right)\right).
\]
Then, it is sufficient to pass to the limit in $\overline{\mathbb{R}}^{+}$
and use the definitions of $\mathbb{E}^{\mathscr{B}}\left(Y\right)$
and $\mathbb{E}^{\mathscr{B}}\left(Z\right).$

For the conditional Beppo Levi property, the nondecreasing monotonicity
of the sequence $\left(\mathbb{E}^{\mathscr{B}}\left(Y_{n}\right)\right)_{n\in\mathbb{N}}$
follows from the first property. Hence, this sequence converges in
$\overline{\mathbb{R}}^{+}.$ Moreover, for every $B\in\mathscr{B},$
the usual Beppo Levi property yields
\[
\intop_{B}\mathbb{E}^{\mathscr{B}}\left(Y\right)\text{d}P=\intop_{B}Y\text{d}P=\lim_{n\to+\infty}\intop_{B}Y_{n}\text{d}P=\lim_{n\to+\infty}\intop_{B}\mathbb{E}^{\mathscr{B}}\left(Y_{n}\right)\text{d}P=\intop_{B}\lim_{n\to+\infty}\mathbb{E}^{\mathscr{B}}\left(Y_{n}\right)\text{d}P.
\]
This implies that
\[
\mathbb{E}^{\mathscr{B}}\left(Y\right)=\lim_{n\to+\infty}\mathbb{E}^{\mathscr{B}}\left(Y_{n}\right).
\]

\end{proof}

\begin{remark}{}{}As an immediate consequence, the first four properties
of $\mathbb{E}^{\mathscr{B}}$ stated in Proposition $\ref{pr:ext_cond_exp_prop}$
remain valid on $\mathscr{M}^{+}\left(\mathscr{A}\right).$

\end{remark}

\subsection{Convergence Theorems}

As we obtained a conditional Beppo Levi property, one can derive---exactly
as in classical integration theory---a conditional version both of
Fatou lemma and the dominated convergence theorem.

\begin{lemma}{Conditional Fatou Lemma}{}

Let $\left(X_{n}\right)_{n\in\mathbb{N}}$ be a sequence of elements
of $\mathscr{M}^{+}\left(\mathscr{A}\right).$ Then,
\[
\mathbb{E}^{\mathscr{B}}\left(\liminf_{n\to+\infty}X_{n}\right)\leqslant\liminf_{n\to+\infty}\mathbb{E}^{\mathscr{B}}\left(X_{n}\right).
\]

\end{lemma}

\begin{proof}{}{}

For every $n\in\mathbb{N},$
\[
\forall p\geqslant n,\,\,\,\,\inf_{k\geqslant n}X_{k}\leqslant X_{p}
\]
and thus, by conditional expectation growth,
\[
\forall p\geqslant n,\,\,\,\,\mathbb{E}^{\mathscr{B}}\left(\inf_{k\geqslant n}X_{k}\right)\leqslant\mathbb{E}^{\mathscr{B}}\left(X_{p}\right).
\]

It follows that
\[
\mathbb{E}^{\mathscr{B}}\left(\inf_{k\geqslant n}X_{k}\right)\leqslant\inf_{p\geqslant n}\mathbb{E}^{\mathscr{B}}\left(X_{p}\right).
\]
By the conditional Beppo Levi property, we conclude.

\end{proof}

\begin{theorem}{Conditional Fatou-Lebesgue Dominated Convergence}{cond_fatou_lebes_dom_cv}

Let $\left(X_{n}\right)_{n\in\mathbb{N}}$ be a sequence of $P-$almost
surely finite random variables, and let $Y\in\mathscr{L}^{1}_{\overline{\mathbb{R}}}\left(\Omega,\mathscr{A},P\right)$
be such that
\[
\forall n\in\mathbb{N},\,\,\,\,\left|X_{n}\right|\leqslant Y\,\,\,\,P-\text{almost surely.}
\]

(a) The following chain of inequalities holds:\boxeq{
\[
\mathbb{E}^{\mathscr{B}}\left(\liminf_{n\to+\infty}X_{n}\right)\leqslant\liminf_{n\to+\infty}\mathbb{E}^{\mathscr{B}}\left(X_{n}\right)\leqslant\limsup_{n\to+\infty}\mathbb{E}^{\mathscr{B}}\left(X_{n}\right)\leqslant\mathbb{E}^{\mathscr{B}}\left(\limsup_{n\to+\infty}X_{n}\right).
\]

}

(b) Moreover, if the sequence $\left(X_{n}\right)_{n\in\mathbb{N}}$
converges $P-$almost surely, the sequence $\left(\mathbb{E}^{\mathscr{B}}\left(X_{n}\right)\right)_{n\in\mathbb{N}}$
also converges $P-$almost surely, and\boxeq{
\[
\mathbb{E}^{\mathscr{B}}\left(\liminf_{n\to+\infty}X_{n}\right)=\lim_{n\to+\infty}\mathbb{E}^{\mathscr{B}}\left(X_{n}\right)\,\,\,\,P-\text{almost surely}
\]
}

\end{theorem}

\begin{proof}{}{}

Since $Y$ and $X_{n}$ are $P-$almost surely finite, the random
variables $Y+X_{n}$ and $Y-X_{n}$ are well defined and nonnegative
$P-$almost surely. Note similarly their measurable extension by 0.
Then applying the conditional Fatou lemma to the sequence $\left(Y+X_{n}\right)_{n\in\mathbb{N}}$
yields
\[
\mathbb{E}^{\mathscr{B}}\left(\liminf_{n\to+\infty}\left(Y+X_{n}\right)\right)\leqslant\liminf_{n\to+\infty}\mathbb{E}^{\mathscr{B}}\left(Y+X_{n}\right).
\]
Hence,
\[
\mathbb{E}^{\mathscr{B}}\left(Y\right)+\mathbb{E}^{\mathscr{B}}\left(\liminf_{n\to+\infty}X_{n}\right)\leqslant\mathbb{E}^{\mathscr{B}}\left(Y\right)+\liminf_{n\to+\infty}\mathbb{E}^{\mathscr{B}}\left(X_{n}\right).
\]
Since $\mathbb{E}^{\mathscr{B}}\left(Y\right)$ is integrable, thus
finite $P-$almost surely, this gives the first inequality. 

The second inequality is obtained analogously by applying the conditional
Fatou with the sequence $\left(Y-X_{n}\right)_{n\in\mathbb{N}}.$

Lastly, assume that the sequence $\left(X_{n}\right)_{n\in\mathbb{N}}$
is $P-$almost surely convergent. This is equivalent to
\[
\liminf_{n\to+\infty}X_{n}=\limsup_{n\to+\infty}X_{n}=\lim_{n\to+\infty}X_{n}.
\]
By the dominated convergence theorem, $\lim_{n\to+\infty}X_{n}\in\mathscr{L}^{1}_{\overline{\mathbb{R}}}\left(\Omega,\mathscr{A},P\right)$
and
\[
\mathbb{E}^{\mathscr{B}}\left(\lim_{n\to+\infty}X_{n}\right)\leqslant\liminf_{n\to+\infty}\mathbb{E}^{\mathscr{B}}\left(X_{n}\right)\leqslant\limsup_{n\to+\infty}\mathbb{E}^{\mathscr{B}}\left(X_{n}\right)\leqslant\mathbb{E}^{\mathscr{B}}\left(\lim_{n\to+\infty}X_{n}\right),
\]
which yields the announced result.

\end{proof}

\begin{corollary}{}{}

Let $\left(X_{n}\right)_{n\in\mathbb{N}}$ be a sequence of $P-$almost
surely finite random variables such that
\[
\sum^{+\infty}_{n=0}E\left(\left|X_{n}\right|\right)<+\infty.
\]
Then, $P-$almost surely the series $\sum X_{n}$ converges absolutely.
Its sum belongs to $\mathscr{L}^{1}_{\overline{\mathbb{R}}}\left(\Omega,\mathscr{A},P\right),$
and
\[
\mathbb{E}^{\mathscr{B}}\left(\sum^{+\infty}_{n=0}X_{n}\right)=\sum^{+\infty}_{n=0}\mathbb{E}^{\mathscr{B}}\left(X_{n}\right).
\]

\end{corollary}

\begin{proof}{}{}

The conditional Lebesgue theorem is applied to the sequence of partial
sums.

\end{proof}

\begin{example}{}{}

Returning to Example $\ref{ex:Poisson_and_B_meas_indep},$ assume
first that $X$ is bounded by a constant $M,$ and compute $\mathbb{E}^{\mathscr{B}}\left(\cos\left(XY\right)\right)$
by expanding the cosine into its power series. The general case is
then obtained by passing to the limit.

\end{example}

\begin{remark}{}{}

This method is considerably longer and illustrates the simplifications
achieved by establishing Proposition $\refpar{pr:cond_exp_exist_indep}.$

\end{remark}

\begin{solutionexample}{}{}

We start from the power series expansion
\[
\cos\left(XY\right)=\sum^{+\infty}_{r=0}\left(-1\right)^{r}\dfrac{\left(XY\right)^{2r}}{\left(2r\right)!}.
\]
We are going to check, that under the hypothesis $\left|X\right|\leqslant M,$
\[
\sum^{+\infty}_{r=0}\mathbb{E}\left(\dfrac{\left|XY\right|^{2r}}{\left(2r\right)!}\right)<+\infty.
\]
Since $Y$ follows a Poisson law with parameter $\lambda$, we have,
in $\overline{\mathbb{R}}^{+},$
\begin{align*}
\sum^{+\infty}_{r=0}\mathbb{E}\left(\dfrac{\left|XY\right|^{2r}}{\left(2r\right)!}\right) & \leqslant\sum^{+\infty}_{r=0}\dfrac{M^{2r}}{\left(2r\right)!}\left(\sum^{+\infty}_{k=0}k^{2r}\text{e}^{-\lambda}\dfrac{\lambda^{k}}{k!}\right)\\
 & \leqslant\text{e}^{-\lambda}\sum^{+\infty}_{k=0}\dfrac{\lambda^{k}}{k!}\left(\sum^{+\infty}_{r=0}\dfrac{\left(Mk\right)^{2r}}{\left(2r\right)!}\right)\\
 & =\text{e}^{-\lambda}\sum^{+\infty}_{k=0}\dfrac{\lambda^{k}}{k!}\text{ch}\left(Mk\right)\\
 & =\dfrac{\text{e}^{-\lambda}}{2}\left(\text{e}^{\lambda\text{e}^{M}}+\text{e}^{\lambda\text{e}^{-M}}\right)<+\infty,
\end{align*}
Hence,
\[
\mathbb{E}^{\mathscr{B}}\left(\cos\left(XY\right)\right)=\sum^{+\infty}_{r=0}\left(-1\right)^{r}\mathbb{E}^{\mathscr{B}}\left(\dfrac{\left|XY\right|^{2r}}{\left(2r\right)!}\right).
\]

Since $X^{2r}$ is $\mathscr{B}-$measurable and $Y^{2r}$ and $\mathscr{B}$
are independent,
\[
\mathbb{E}^{\mathscr{B}}\left(X^{2r}Y^{2r}\right)=X^{2r}\mathbb{E}^{\mathscr{B}}\left(Y^{2r}\right)=X^{2r}\mathbb{E}\left(Y^{2r}\right).
\]

A computation identical to the previous, the interversions of the
sum signs being justified by the absolute convergence of the double
series, yields
\begin{align*}
\mathbb{E}^{\mathscr{B}}\left(\cos\left(XY\right)\right) & =\sum^{+\infty}_{r=0}\dfrac{\left(-1\right)^{r}}{\left(2r\right)!}X^{2r}\left(\sum^{+\infty}_{k=0}k^{2r}\text{e}^{-\lambda}\dfrac{\lambda^{k}}{k!}\right)\\
 & =\text{e}^{-\lambda}\sum^{+\infty}_{k=0}\dfrac{\lambda^{k}}{k!}\left(\sum^{+\infty}_{r=0}\dfrac{\left(-1\right)^{r}\left(kX\right)^{2r}}{\left(2r\right)!}\right)\\
 & =\text{e}^{-\lambda}\sum^{+\infty}_{k=0}\dfrac{\lambda^{k}}{k!}\cos\left(kX\right).
\end{align*}

Since
\[
\sum^{+\infty}_{k=0}\dfrac{\lambda^{k}}{k!}\cos\left(kX\right)=\text{Re}\left(\sum^{+\infty}_{k=0}\dfrac{\lambda^{k}}{k!}\text{e}^{\text{i}kx}\right)=\text{Re}\left(\text{e}^{\lambda\text{e}^{\text{i}x}}\right),
\]
it yields\boxeq{
\begin{equation}
\mathbb{E}^{\mathscr{B}}\left(\cos\left(XY\right)\right)=\text{e}^{-\lambda\left(1-\cos X\right)}\times\cos\left(\lambda\sin X\right).\label{eq:cond_exp_cosXY}
\end{equation}
}

Now, let $X$ be arbitrary. For each $n\in\mathbb{N},$ define $X_{n}=\boldsymbol{1}_{\left\{ \left|X\right|\leqslant n\right\} }X.$
Then, the sequence $\left(\cos\left(X_{n}Y\right)\right)_{n\in\mathbb{N}}$
is $P-$almost surely convergent, and, for every $n\in\mathbb{N},$
\[
\left|\cos\left(X_{n}Y\right)\right|\leqslant1.
\]
By the conditional dominated convergence theorem, the sequence $\left(\mathbb{E}^{\mathscr{B}}\left(\cos\left(X_{n}Y\right)\right)\right)_{n\in\mathbb{N}}$
is $P-$almost surely convergent and
\[
\mathbb{E}^{\mathscr{B}}\left(\cos\left(XY\right)\right)=\lim_{n\to+\infty}\mathbb{E}^{\mathscr{B}}\left(\cos\left(X_{n}Y\right)\right)\,\,\,\,P-\text{almost surely.}
\]
Therefore, the formula $\refpar{eq:cond_exp_cosXY}$ holds for an
arbitrary random variable $X.$

\end{solutionexample}

\subsection{Jensen Inequality}

Jensen inequality is a convexity inequality often used. We begin with
an elementary version, before presenting a more refined formulation.

\begin{proposition}{Jensen Inequality}{jensen_ineq}

Let $g$ be a convex\footnotemark function on $\mathbb{R}$ and let
$Y\in\mathscr{L}^{1}\left(\Omega,\mathscr{A},P\right)$ be such that
$g\circ Y\in\mathscr{L}^{1}\left(\Omega,\mathscr{A},P\right).$ Then,\boxeq{
\begin{equation}
g\left(\mathbb{E}^{\mathscr{B}}\left(Y\right)\right)\leqslant\mathbb{E}^{\mathscr{B}}\left(g\circ Y\right).\label{eq:jensen_ineq}
\end{equation}
}

\end{proposition}

\footnotetext{Recall that any convex real-valued function defined
on an open interval of $\mathbb{R}$ is continuous. This property
fail when the interval is not open. For instance, consider the function
$g$ defined on $\left[0,+\infty\right)$ by $g\left(0\right)=1$
and $g\left(x\right)=0$ if $x>0.$

}

\begin{proof}{}{}

Since the function $g$ is convex, there exists two sequences of real
numbers such that
\[
\forall x\in\mathbb{R},\,\,\,\,g\left(x\right)=\sup_{n\in\mathbb{N}}\left(a_{n}x+b_{n}\right).
\]

Then,
\[
\forall n\in\mathbb{N},\,\,\,\,P-\text{almost surely,}\,\,\,\,a_{n}\mathbb{E}^{\mathscr{B}}\left(Y\right)+b_{n}=\mathbb{E}^{\mathscr{B}}\left(a_{n}Y+b_{n}\right)\leqslant\mathbb{E}^{\mathscr{B}}\left(g\circ Y\right).
\]
Since a countable union of sets of null probability is of null probability,
\[
P-\text{almost surely,}\,\,\,\,\forall n\in\mathbb{N},\,\,\,\,a_{n}\mathbb{E}^{\mathscr{B}}\left(Y\right)+b_{n}\leqslant\mathbb{E}^{\mathscr{B}}\left(g\circ Y\right),
\]
which implies
\[
P-\text{almost surely,}\,\,\,\,\sup_{n\in\mathbb{N}}\left(a_{n}\mathbb{E}^{\mathscr{B}}\left(Y\right)+b_{n}\right)\leqslant\mathbb{E}^{\mathscr{B}}\left(g\circ Y\right),
\]
which proves the desired inequality.

\end{proof}

\begin{corollary}{$\mathbb{E}^{\mathscr{B}}$ as Contraction}{}

Let $p\in\mathbb{N}^{\ast}.$ If $Y\in\text{L}^{p}\left(\Omega,\mathscr{A},P\right),$
then $\mathbb{E}^{\mathscr{B}}\left(Y\right)\in\text{L}^{p}\left(\Omega,\mathscr{B},P\right)$
and
\begin{equation}
\left\Vert \mathbb{E}^{\mathscr{B}}\left(Y\right)\right\Vert _{p}\leqslant\left\Vert Y\right\Vert _{p}.\label{eq:contraction}
\end{equation}
That is, $\mathbb{E}^{\mathscr{B}}$ is a \textbf{contraction\index{contraction}}
of $\text{L}^{p}\left(\Omega,\mathscr{A},P\right)$ into $\text{L}^{p}\left(\Omega,\mathscr{B},P\right).$

\end{corollary}

\begin{proof}{}{}

We apply the Jensen inequality to the convex function $x\mapsto\left\Vert x\right\Vert _{p}.$

\end{proof}

\begin{proposition}{Continuous Function on Closed Convex of $\mathbb{R}$ and Jensen Inequality}{}

(a) Let $Y\in L^{p}\left(\Omega,\mathscr{A},P\right)$ take values
in a closed convex $K$ of $\mathbb{R}$---that is, in a closed interval.
The conditional expectation $\mathbb{E}^{\mathscr{B}}\left(Y\right)$
takes values in $K,$ $P-$almost surely. 

For every continuous convex function $g$ on $K,$ taking nonnegative
values in $\mathbb{R}\cup\left\{ +\infty\right\} ,$ or such that
$g\circ Y\in\mathscr{L}^{1}_{\overline{\mathbb{R}}}\left(\Omega,\mathscr{A},P\right),$
the Jensen inequality $\refpar{eq:jensen_ineq}$ is satisfied.

(b) Let $Y$ be a random variable taking values in $\overline{\mathbb{R}}^{+}.$
For every continuous convex function $g$ on $\overline{\mathbb{R}}^{+},$
such that $g\left(+\infty\right)=+\infty,$ and which is non negative
or satisfies $g\circ Y\in\mathscr{L}^{1}_{\overline{\mathbb{R}}}\left(\Omega,\mathscr{A},P\right),$
the Jensen inequality $\refpar{eq:jensen_ineq}$ holds.

\end{proposition}

\begin{proof}{}{}

The proof is exactly the same as that of the previous proposition. 

\end{proof}

\begin{remark}{}{}

It is not difficult to generalize the concept of conditional expectation
to the case where the random variable $Y$ takes values in a Euclidean
space. In this setting, the previous proposition remains valid, since
every closed convex can be written as a countable intersection of
closed half-spaces.

\end{remark}

\subsection{Computation of the Conditional Expectation}

We have already given an example of such a computation. \textbf{A
particularly frequent case is when the sub-$\sigma-$algebra $\mathscr{B}$
is generated by a random variable $X,$ and when there exists a conditional
law of $Y$ given $X.$}

\begin{proposition}{}{}

Let $X$ be a random variable taking values in an arbitrary probabilizable
space $\left(E,\mathscr{E}\right),$ and let $Y\in\mathscr{L}^{1}\left(\Omega,\mathscr{A},P\right).$
Assume that there exists a conditional law $P^{X=\cdot}_{Y}$ of $Y$
given $X.$ Then $m^{X=\cdot}_{Y}\circ X$ is a version of the conditional
expectation $\mathbb{E}^{\sigma\left(X\right)}\left(Y\right),$ where
$m^{X=\cdot}_{Y}$ denotes the conditional mean of $Y$ given $X,$
which is written\footnotemark
\[
\mathbb{E}^{\sigma\left(X\right)}\left(Y\right)=m^{X=\cdot}_{Y}\circ X\,\,\,\,P-\text{almost\,surely.}
\]

\end{proposition}

\footnotetext{Some authors write the conditional expectation $\mathbb{E}^{\sigma\left(X\right)}\left(Y\right)$
in the form $\mathbb{E}\left(Y\mid X\right).$ This notation will
be used in this book when the typographical context makes it appropriate.}

\begin{proof}{}{}

Under the present assumptions, the conditional expectation exists.
Recall that 
\[
\sigma\left(X\right)=\left\{ X^{-1}\left(C\right):\,C\in\mathscr{E}\right\} ,
\]
and that $\boldsymbol{1}_{X^{-1}}\left(C\right)=\boldsymbol{1}_{C}\circ X.$ 

For every $C\in\mathscr{E},$ by the transfer theorem and the expanded
Fubini theorem,
\begin{align*}
\intop_{X^{-1}\left(C\right)}m^{X=\cdot}_{Y}\circ X\text{d}P & =\intop_{E}\boldsymbol{1}_{C}\left(x\right)m^{X=x}_{Y}\text{d}P_{X}\left(x\right)=\intop_{E}\boldsymbol{1}_{C}\left(x\right)\left[\intop_{\mathbb{R}}y\text{d}P^{X=x}_{Y}\left(y\right)\right]\text{d}P_{X}\left(x\right)\\
 & =\intop_{E\times\mathbb{R}}\boldsymbol{1}_{C}\left(x\right)y\text{d}P_{\left(X,Y\right)}\left(x,y\right)=\intop_{X^{-1}\left(C\right)}Y\text{d}P.
\end{align*}
This proves the result.

\end{proof}

\begin{example}{}{}

We return to Example $\ref{ex:common_sp_case_example},$ where $X$
and $Y$ are independent of same law $\exp\left(\lambda\right).$
We have seen that, for every $s>0,$ $P^{S=s}_{X}$ is the uniform
law on $\left[0,s\right].$ Consequently, 
\[
m^{S=s}_{X}=\dfrac{s}{2}.
\]
It follows that $\dfrac{S}{2}$ is a version of the conditional expectation
$\mathbb{E}^{\sigma\left(S\right)}\left(X\right).$

\end{example}

As the following example illustrates, this result does not depend
on the specific nature of the laws of the random variables.

\begin{example}{}{cond_exp_know_sigma-alg_gen_sum}

Let $X_{1}$ and $X_{2}$ be two independent real-valued random variables
with the same law $\mu.$ Let $S=X_{1}+X_{2}.$ Prove that $\mathbb{E}^{\sigma\left(S\right)}\left(X_{1}\right)=\mathbb{E}^{\sigma\left(S\right)}\left(X_{2}\right)$
and deduce $\mathbb{E}^{\sigma\left(S\right)}\left(X_{1}\right).$

\end{example}

\begin{remark}{}{}

This example will be developed further as an exercise.

\end{remark}

\begin{solutionexample}{}{}

Taking into account the independence of $X_{1}$ and $X_{2},$ we
have, for every Borel subset $C$ of $\mathbb{R},$
\[
\intop_{S^{-1}\left(C\right)}X_{1}\text{d}P=\intop_{\mathbb{R}^{2}}\boldsymbol{1}_{C}\left(x_{1}+x_{2}\right)x_{1}\text{d}\left(P_{X_{1}}\otimes P_{X_{2}}\right)\left(x_{1},x_{2}\right).
\]
Since $X_{1}$ and $X_{2}$ have the same law,
\[
\intop_{S^{-1}\left(C\right)}X_{1}\text{d}P=\intop_{\mathbb{R}^{2}}\boldsymbol{1}_{C}\left(x_{1}+x_{2}\right)x_{1}\text{d}\left(P_{X_{2}}\otimes P_{X_{1}}\right)\left(x_{1},x_{2}\right).
\]
It follows that
\[
\forall C\in\mathscr{B}_{\mathbb{R}},\,\,\,\,\intop_{S^{-1}\left(C\right)}X_{1}\text{d}P=\intop_{S^{-1}\left(C\right)}X_{2}\text{d}P,
\]
which shows the equality\boxeq{
\[
\mathbb{E}^{\sigma\left(S\right)}\left(X_{1}\right)=\mathbb{E}^{\sigma\left(S\right)}\left(X_{2}\right).
\]
}

Since $S$ is $\sigma\left(S\right)-$measurable,
\[
\mathbb{E}^{\sigma\left(S\right)}\left(X_{1}+X_{2}\right)=S=\mathbb{E}^{\sigma\left(S\right)}\left(X_{1}\right)+\mathbb{E}^{\sigma\left(S\right)}\left(X_{2}\right)\,\,\,\,P-\text{almost surely.}
\]
Therefore,\boxeq{
\[
\mathbb{E}^{\sigma\left(S\right)}\left(X_{1}\right)=\dfrac{S}{2}\,\,\,\,P-\text{almost surely},
\]
}a result that is entirely natural.

\end{solutionexample}

\section*{Exercises}

\textbf{Unless otherwise specified, all the random variables are defined
on the same probabilized space $\left(\Omega,\mathscr{A},P\right).$}

\addcontentsline{toc}{section}{Exercises}

\begin{exercise}{Poisson and Multinomial Laws}{exercise12.1}

Let $X_{1},X_{2},\cdots,X_{n}$ be $n$ independent random variables,
where $X_{i}$ follows a Poisson law $\mathscr{P}\left(\lambda_{i}\right)$
for $i=1,2,\cdots,n.$ Let $X=\left(X_{1},X_{2},\cdots,X_{n}\right)$
be a random variable taking values in $\mathbb{N}^{n},$ and let $S_{n}=\sum^{n}_{i=1}X_{i}.$ 

Determine a conditional law $P^{S_{n}=\cdot}_{X}$ of $X$ given $S_{n}.$

\end{exercise}

\begin{exercise}{Bernoulli and Uniform Laws}{exercise12.2}

Let $X_{1},X_{2},\cdots,X_{n}$ be $n$ independent random variables,
all following the same Bernoulli law $\mathcal{B}\left(1,p\right)$
where $0<p<1.$ Let $X=\left(X_{1},X_{2},\cdots,X_{n}\right)$ be
the random variable taking values in $\mathbb{N}^{n},$ and let $S_{n}=\sum^{n}_{i=1}X_{i}.$ 

Determine a conditional law $P^{S_{n}=\cdot}_{X}$ of $X$ given $S_{n}.$

\end{exercise}

\begin{exercise}{Poisson Process}{exercise12.3}

Let $\left(W_{n}\right)_{n\in\mathbb{N}^{\ast}}$ be a nondecreasing
sequence of nonnegative random variables such that $W_{0}=0.$ For
$n\in\mathbb{N}^{\ast},$ define the random variable $T_{n}=W_{n}-W_{n-1}.$
Assume that the random variables $T_{n},\,n\in\mathbb{N}^{\ast},$
constitute a family of independent random variables, of same \textbf{exponential
law} $\exp\left(\lambda\right),$ where $\lambda>0.$ Set $X_{0}=0,$
and, for every $t>0,$
\[
X_{t}=\sum_{n\in\mathbb{N}^{\ast}}\boldsymbol{1}_{\left(W_{n}\leqslant t\right)}.
\]

The family of random variables $\left(X_{t}\right)_{t\in\mathbb{R}^{+}}$
is called a \textbf{\index{Poisson process of intensity λ@Poisson process of intensity $\lambda$}Poisson
process of intensity $\lambda.$}

\textbf{1.} Let $s,t$ be two real numbers such that $0\leqslant s<t.$
Compute by induction the integral defined, for every $n\in\mathbb{N}^{\ast},$
by
\[
I_{n}\left(s,t\right)=\intop_{\mathbb{R}^{n}}\boldsymbol{1}_{\left(s\leqslant x_{1}\leqslant x_{2}\leqslant\cdots\leqslant x_{n}\leqslant t\right)}\text{d}\lambda_{n}\left(x_{1},x_{2},\cdots,x_{n}\right).
\]

\textbf{2.} Compute, for every $n\in\mathbb{N^{\ast}}$ and for every
family $\left(f_{j}\right)_{1\leqslant j\leqslant n}$ of nonnegative
measurable functions bounded on $\mathbb{R},$ the quantity
\[
\mathbb{E}\left(\boldsymbol{1}_{\left(X_{t}=n\right)}\prod^{n}_{j=1}f_{j}\left(W_{j}\right)\right).
\]
Deduce the law of $X_{t}$ and a conditional law of $\left(W_{1},W_{2},\cdots,W_{n}\right)$
given $\left(X_{t}=n\right).$

\textbf{3.} Consider $t>0,$ and an arbitrary integer $k\geqslant1,$
and an arbitrary finte sequence of real numbers such that $0=t_{0}\leqslant t_{1}\leqslant\cdots\leqslant t_{k}=t.$ 

\textbf{a.} Determine the law of the random variable $\left(X_{t_{1}},X_{t_{2}}-X_{t_{1}},\cdots,X_{t_{k}}-X_{t_{k-1}}\right)$
and justify the independence of the random variables $X_{t_{1}},X_{t_{2}}-X_{t_{1}},\cdots,X_{t_{k}}-X_{t_{k-1}}.$ 

We say that the \textbf{process} $\left(X_{t}\right)_{t\in\mathbb{R}^{+}}$
is \textbf{\index{process with independent increments}with independent
increments}.

\textbf{b.} For every $s,t$ such that $0\leqslant s<t,$ determine
the law of the random variable $X_{t}-X_{s}.$ Deduce from it its
expectation $\mathbb{E}\left(X_{t}-X_{s}\right).$

\textbf{4.} Let $k\in\mathbb{N}^{\ast}$ such that $1\leqslant k\leqslant n.$
Determine a conditional law $P^{X_{t}=n}_{W_{k}}$ of $W_{k}$ given
$X_{t}=n$ and identify it.

\begin{remark}{}{}The Poisson process is a particular case of a \textbf{\index{counting process}counting
process}. A property occurs randomly over time. The random variable
$W_{n}$ corresponds to the time of the $n-$th outcome and $T_{n}$
to the time elapsed between the $(n-1)-$th and the $n-$th realization.
The random variable $X_{t}$ represents the number of occurrences
of the property in the time interval $\left[0,t\right].$ This process
arises in particular in the queuing theory.

\end{remark}

\end{exercise}

\begin{exercise}{Uniform Drawing and Interval of Random Length}{exercise12.4}

Let $\left\{ L,\left(X_{n}\right)_{n\in\mathbb{N}^{\ast}}\right\} $
be a family of independent real-valued random variables, all following
the same uniform law on $\left[0,1\right].$ For every $n\in\mathbb{N}^{\ast},$
define the function $S_{n}$ by
\[
\forall\omega\in\Omega,\,\,\,\,S_{n}\left(\omega\right)=\sum^{n}_{j=1}\boldsymbol{1}_{\left[0,L\left(\omega\right)\right]}\left(X_{j}\left(\omega\right)\right).
\]

1. Verify that $S_{n}$ is a random variable and determine a conditional
law $P^{L=\cdot}_{S_{n}}$ of $S_{n}$ given $L.$ 

2. Deduce the law of $S_{n}.$

3. Determine a conditional law $P^{S_{n}=\cdot}_{L}$ of $L$ given
$S_{n}.$ Compute the conditional mean $m^{S_{n}=\cdot}_{L}$ and
recover the expectation $\mathbb{E}\left(L\right)$ of $L.$

\end{exercise}

\begin{exercise}{Whole and Decimal Parts}{exercise12.5}

Let $X$ be a nonnegative random variable of density $f_{X}.$ Define
$Y=X-\left\lfloor X\right\rfloor $ where $\left\lfloor \cdot\right\rfloor $
denotes the integer part.

1. Determine the law of the random variable $\left(\left\lfloor X\right\rfloor ,Y\right)$
as a function of $f_{X}.$ Deduce the law of $\left\lfloor X\right\rfloor $
and of $Y.$ Retrieve directly the law of $\left\lfloor X\right\rfloor .$

2. Determine the conditional laws of $Y$ given $\left\lfloor X\right\rfloor $
and of $\left\lfloor X\right\rfloor $ given $Y,$ respectively $P^{\left\lfloor X\right\rfloor =\cdot}_{Y}$
and $P^{Y=\cdot}_{\left\lfloor X\right\rfloor }.$

3. Suppose that $X$ follows a gamma law $\gamma\left(a,p\right)$
where $a>0$ and $p>0.$ For which values of the couple $\left(a,p\right)$
are the random variables $\left\lfloor X\right\rfloor $ and $Y$
independent? Determine the laws of $\left\lfloor X\right\rfloor $
and $Y$ in the case where $X$ follows the exponential law $\exp\left(p\right),$
$p>0.$

4. Suppose the density $f_{X}$ of $X$ to be given by
\[
\forall x\in\mathbb{R},\,\,\,\,f_{X}\left(x\right)=\sum^{+\infty}_{n=0}\boldsymbol{1}_{\left[n,n+1\right[}\left(x\right)\text{e}^{-\lambda}\dfrac{\lambda^{n}}{n!}.
\]
Prove that the random variables $\left\lfloor X\right\rfloor $ and
$Y$ are independent and determine their law.

5. Suppose again that $X$ follows the gamma law $\gamma\left(a,p\right).$
Compute the conditional expectations $m^{\left\lfloor X\right\rfloor =\cdot}_{Y}$
and $m^{Y=\cdot}_{\left\lfloor X\right\rfloor }$---no explicitation
of the integrals and series sum is required. 

\end{exercise}

\begin{exercise}{Conditional Expectation and Gaussian Random Variable. Different Computation Methods}{exercise12.6}

Let $\mathscr{B}$ be a sub-$\sigma-$algebra of $\mathscr{A},$ and
let $X$ and $Y$ be two real-valued random variables such that $X$
is $\mathscr{B}-$measurable and the $\sigma-$algebras $\mathscr{B}$
and $\sigma\left(Y\right)$ are independent. Moreover, we suppose
that $Y$ follows the normal law $\mathscr{N}_{\mathbb{R}}\left(0,1\right).$

1. Prove the equivalence of the three following properties:

(i) $\text{e}^{\frac{X^{2}}{2}}$ is $P-$integrable

(ii) $\text{e}^{XY}$ is $P-$integrable

(iii) $\text{e}^{\left|XY\right|}$ is $P-$integrable

2. Suppose that $\text{e}^{\frac{X^{2}}{2}}$ is $P-$integrable.

(a) Without computing the conditional expectation, prove that
\[
\mathbb{E}^{\mathscr{B}}\left(\text{e}^{XY}\right)\geqslant1\,\,\,\,P-\text{almost surely}
\]

(b) In the case where $\mathscr{B}=\sigma\left(X\right),$ compute
$\mathbb{E}^{\mathscr{B}}\left(\text{e}^{XY}\right).$

(c) In the general case---a priori $\sigma\left(X\right)\underset{\neq}{\subset}\mathscr{B}$---compute,
using two different methods, the conditional expectation $\mathbb{E}^{\mathscr{B}}\left(\text{e}^{XY}\right):$
first by expanding the exponential into a series, and then by applying
Proposition $\ref{pr:cond_exp_exist_indep}.$

3. We no longer assume that $\text{e}^{\frac{X^{2}}{2}}$ is $P-$integrable. 

Compute $\mathbb{E}^{\mathscr{B}}\left(\text{e}^{XY}\right).$

\end{exercise}

\begin{exercise}{Conditional Expectation and Independence}{exercise12.7}

Let $X\in\mathscr{L}^{1}\left(\Omega,\mathscr{A},P\right)$ and let
$\mathscr{A}_{1}$ and $\mathscr{A}_{2}$ be two sub$-$\salg of $\mathscr{A}$
such that the $\sigma-$algebra $\mathscr{A}_{1}\vee\sigma\left(X\right)$
and $\mathscr{A}_{2}$ are independent. Here, $\mathscr{A}_{1}\vee\sigma\left(X\right)$
denotes the $\sigma-$algebra generated by $\mathscr{A}_{1}$ and
$\sigma\left(X\right),$ that is, the smallest $\sigma-$algebra containing
$\mathscr{A}_{1}$ and $\sigma\left(X\right).$

Prove the equality
\[
\mathbb{E}^{\mathscr{A}_{1}\vee\mathscr{A}_{2}}\left(X\right)=\mathbb{E}^{\mathscr{A}_{1}}\left(X\right).
\]

\end{exercise}

\begin{exercise}{On the Path of a Strong Law of Large Numbers}{exercise12.8}

Let $\left(X_{n}\right)_{n\in\mathbb{N}^{\ast}}$ be a sequence of
independent real-valued random variables, all following the same law
$\mu.$

Let $S_{n}=\sum^{n}_{i=1}X_{i}.$ Prove that, for every $i\in\left\llbracket 1,n\right\rrbracket $
\[
\mathbb{E}^{\sigma\left(S_{n}\right)}\left(X_{i}\right)=\mathbb{E}^{\sigma\left(S_{n}\right)}\left(X_{1}\right)
\]
and deduce $\mathbb{E}^{\sigma\left(S_{n}\right)}\left(X_{1}\right),$
and then $\mathbb{E}^{\mathscr{A}_{n}}\left(X_{1}\right),$ where
$\mathscr{A}_{n}=\sigma\left(S_{n+j}:\,j\in\mathbb{N}\right).$ 

\textit{Hint: Use the previous exercise.}

This exercise expands Example $\ref{ex:cond_exp_know_sigma-alg_gen_sum}.$ 

\end{exercise}

\begin{exercise}{Method of Simulation by Rejection and Introduction to Markov Chains Methods}{exercise12.9}

Let $f$ and $g$ be two probability densities on $\mathbb{R}$ with
respect to the Lebesgue measure. The objective is to simulate a random
variable $Y$ with density $f,$ whose analytical form is at first
sight ``complicated'', by using a density $g$ that is chosen sufficiently
``close to'' $f$ and has a ``simpler'' analytical form. With
the convention $\dfrac{0}{0}=0$, set
\[
t\left(x\right)=\dfrac{f\left(x\right)}{g\left(x\right)}
\]
and assume that
\[
1<\overline{t}=\sup_{x\in\mathbb{R}}t\left(x\right)<+\infty.
\]

Consider a family of independent random variables 
\[
\left\{ X_{n},Y_{p}:\,n\in\mathbb{N}^{\ast},p\in\mathbb{N}^{\ast}\right\} 
\]
such that, for every $n\in\mathbb{N}^{\ast},$ $X_{n}$ has density
$g$ and $Y_{n}$ follows the uniform law on the interval $\left[0,\overline{t}\right].$ 

Define the random variables, taking values respectively in $\mathbb{R}^{2}$
and $\mathbb{R}\cup\left\{ +\infty\right\} ,$
\[
M_{n}=\left(X_{n},Y_{n}\right)\,\,\,\,\text{and}\,\,\,\,X_{\infty}=\limsup_{n\to+\infty}X_{n}.
\]

Consider the set 
\[
G=\left\{ \left(x,y\right)\in\mathbb{R}^{2}:\,t\left(x\right)\geqslant y\right\} 
\]
and define the functions $T$ and $X_{T},$ taking values respectively
in $\overline{\mathbb{N}}$ and $\mathbb{R}\cup\left\{ +\infty\right\} ,$
for every $\omega\in\Omega,$ by
\[
T\left(\omega\right)=\inf\left\{ n\in\mathbb{N}^{\ast}:\,M_{n}\left(\omega\right)\in G\right\} \,\,\,\,\text{and}\,\,\,\,X_{T}\left(\omega\right)=X_{T\left(\omega\right)}\left(\omega\right),
\]
with the convention $\inf\emptyset=+\infty.$

1. Prove that, for every bounded and measurable function $\varphi$
on $\mathbb{R},$ the quantity
\[
I\left(\varphi\right)\equiv\mathbb{E}\left(\boldsymbol{1}_{\left(M_{n}\in G\right)}\varphi\left(X_{n}\right)\right)
\]
is equal, up to a constant that will be determined to the integral
$\intop_{\mathbb{\mathbb{R}}}\varphi f\text{d}\lambda.$ 

Deduce from this the probability $P\left(M_{n}\in G\right).$

2. Verify that $T$ and $X_{T}$ are random variables. Determine the
law of $T$ and deduce that $X_{T}$ is $P-$almost surely finite.
Compute the expectation of $T.$

3. Determine the law of $X_{T}.$

4. Still with the convention $\inf\emptyset=+\infty,$ define the
sequence of functions with values in $\overline{\mathbb{N}}$ by $T_{1}=T$
and, for every $k\in\mathbb{N}^{\ast},$
\[
\forall\omega\in\Omega,\,\,\,\,T_{k+1}\left(\omega\right)=\inf\left\{ n>T_{k}\left(\omega\right):\,M_{n}\in G\right\} .
\]

Denote by $\mathscr{A}_{n}$ the $\sigma-$algebra 
\[
\mathscr{A}_{n}=\sigma\left(M_{j}:\,1\leqslant j\leqslant n\right).
\]

(a) Prove, by induction on $k,$ that for every $k\in\mathbb{N}^{\ast}$
and every $n\in\mathbb{N^{\ast}},$ the event $\left(T_{k}=n\right)$
belongs to $\mathscr{A}_{n},$ and that $T_{k}$ is $P-$almost surely
finite. 

For every $k\in\mathbb{N}^{\ast},$ define the family of events
\[
\mathscr{A}_{T_{k}}=\left\{ A\in\mathscr{A}:\,\forall n\in\mathbb{N}^{\ast},\,\,\,\,A\cap\left(T_{k}=n\right)\in\mathscr{A}_{n}\right\} .
\]

Verify that, for every $k\in\mathbb{N}^{\ast},$
\begin{equation}
\mathscr{A}_{T_{k}}\subset\mathscr{A}_{T_{k+1}}.\label{eq:A_T_k_subset_A_t_k_plus_1}
\end{equation}

(b) Let $f_{k}$, for each $k\in\mathbb{N}^{\ast},$ be an arbitrary
bounded nonnegative measurable function. Compute, for every $n\in\mathbb{N}^{\ast},$
the conditional expectation
\[
\mathbb{E}^{\mathscr{A}_{n}}\left[\boldsymbol{1}_{\left(T_{k}=n\right)}f_{k+1}\left(X_{T_{k+1}}\right)\right].
\]
Deduce from this the conditional expectation
\[
\mathbb{E}^{\mathscr{A}_{T_{k}}}\left[f_{k+1}\left(X_{T_{k+1}}\right)\right],
\]
and then the law of $X_{T_{k+1}}.$

(c) Prove that for every $k\in\mathbb{N}^{\ast},$ the random variable
$f_{k}\left(X_{T_{k}}\right)$ is $\mathscr{A}_{T_{k}}-$measurable
and deduce that $\left(X_{T_{k}}\right)_{k\in\mathbb{N}^{\ast}}$
is a sequence of independent random variables.

(d) Numerical application

Let $a>2$ and $b>0$ such that $ab>1,$ and define $f$ and $g$
by
\[
f\left(x\right)=\boldsymbol{1}_{\mathbb{R}^{+}}\left(x\right)\dfrac{b^{a}}{\Gamma\left(a\right)}x^{a-1}\text{e}^{-bx}\,\,\,\,\text{and}\,\,\,\,g\left(x\right)=\boldsymbol{1}_{\mathbb{R}^{+}}\left(x\right)\dfrac{1}{a}\text{e}^{-\frac{x}{a}}.
\]
Verify that $t$ is bounded, and determine $\widehat{x}$ such that
$t\left(\widehat{x}\right)=\widehat{t}.$ 

In the two following case, $b=1$ and $a=\dfrac{5}{2},$ then $b=1$
and $a=\dfrac{9}{2},$ verify that $\overline{t}>1$ and determine
the numerical values of $P\left(M_{n}\in G\right)$ and $\mathbb{E}\left(T\right)$
with an accuracy of $10^{-2}.$ 

\textit{Hint: Recall that $\Gamma\left(\dfrac{1}{2}\right)=\sqrt{\pi}.$}

\end{exercise}

\section*{Solutions of Exercises}

\addcontentsline{toc}{section}{Solutions of Exercises}

\begin{solution}{}{solexercise12.1}

Recall that $S_{n}$ follows a Poisson law $\mathscr{P}\left(\sum^{n}_{i=1}\lambda_{i}\right).$
Moreover, for every $\left(k_{1},k_{2},\cdots,k_{n},x\right)\in\mathbb{N}^{n+1},$
\[
P\left[\left(\bigcap^{n}_{i=1}\left(X_{i}=k_{i}\right)\right)\cap\left(S_{n}=x\right)\right]=\boldsymbol{1}_{\left(\sum^{n}_{i=1}l_{i}=x\right)}\left(k_{1},\cdots,k_{n}\right)P\left[\bigcap^{n}_{i=1}\left(X_{i}=k_{i}\right)\right].
\]
Hence, by independence of the random variables $X_{i},$

\[
P\left[\left(\bigcap^{n}_{i=1}\left(X_{i}=k_{i}\right)\right)\cap\left(S_{n}=x\right)\right]=\boldsymbol{1}_{\left(\sum^{n}_{i=1}l_{i}=x\right)}\left(k_{1},\cdots,k_{n}\right)\text{e}^{-\sum^{n}_{i=1}\lambda_{i}}\prod^{n}_{i=1}\dfrac{\lambda^{k_{i}}_{i}}{k_{i}!}.
\]
Thus, for every $\left(k_{1},k_{2},\cdots,k_{n},x\right)\in\mathbb{N}^{n+1},$
\[
P^{S_{n}=x}_{X}\left[\left(k_{1},k_{2},\cdots,k_{n}\right)\right]=\boldsymbol{1}_{\left(\sum^{n}_{i=1}l_{i}=x\right)}\left(k_{1},\cdots,k_{n}\right)\dfrac{\text{e}^{-\sum^{n}_{i=1}\lambda_{i}}\prod^{n}_{i=1}\dfrac{\lambda^{k_{i}}_{i}}{k_{i}!}}{\dfrac{\left(\sum^{n}_{i=1}\lambda_{i}\right)^{x}}{x!}\text{e}^{-\sum^{n}_{i=1}\lambda_{i}}}.
\]
Hence,\boxeq{
\[
P^{S_{n}=x}_{X}\left[\left(k_{1},k_{2},\cdots,k_{n}\right)\right]=\boldsymbol{1}_{\left(\sum^{n}_{i=1}l_{i}=x\right)}\left(k_{1},\cdots,k_{n}\right)\dfrac{x!}{k_{1}!k_{2}!\cdots k_{n}!}\prod^{n}_{i=1}\left(\dfrac{\lambda^{k_{i}}_{i}}{\sum^{n}_{j=1}\lambda_{j}}\right)^{k_{i}}.
\]
}That is, for every $x\in\mathbb{N}^{\ast},$ $P^{S_{n}=x}_{X}$
is the multinomial law
\[
M\left(x;\dfrac{\lambda_{1}}{\sum^{n}_{j=1}\lambda_{j}},\dfrac{\lambda_{2}}{\sum^{n}_{j=1}\lambda_{j}},\cdots,\dfrac{\lambda_{n}}{\sum^{n}_{j=1}\lambda_{j}}\right).
\]
If $x=0,$ $P^{S_{n}=x}_{X}$ is the Dirac measure in 0.

\begin{remark}{}{}

This property of the Poisson law will be recovered later, in Exercise
$\ref{exo:exercise12.3}$ devoted to the Poisson process: it will
then receive a natural interpretation.

\end{remark}

\end{solution}

\begin{solution}{}{solexercise12.2}

Recall that $S_{n}$ follows a binomial law $\mathcal{B}\left(n,p\right).$
Moreover, for every $\left(k_{1},k_{2},\cdots,k_{n},x\right)\in\left\{ 0,1\right\} ^{n}\times\mathbb{N},$
\[
P\left[\bigcap^{n}_{i=1}\left(X_{i}=k_{i}\right)\cap\left(S_{n}=x\right)\right]=\boldsymbol{1}_{\left(\sum^{n}_{i=1}l_{i}=x\right)}\left(k_{1},k_{2},\cdots,k_{n}\right)P\left[\bigcap^{n}_{i=1}\left(X_{i}=k_{i}\right)\right].
\]
Hence, by independence of the random variables $X_{i},$
\begin{align*}
P\left[\bigcap^{n}_{i=1}\left(X_{i}=k_{i}\right)\cap\left(S_{n}=x\right)\right] & =\boldsymbol{1}_{\left(\sum^{n}_{i=1}l_{i}=x\right)}\left(k_{1},k_{2},\cdots,k_{n}\right)\prod^{n}_{i=1}\left[p^{k_{i}}\left(1-p\right)^{1-k_{i}}\right]\\
 & =\boldsymbol{1}_{\left(\sum^{n}_{i=1}l_{i}=x\right)}\left(k_{1},k_{2},\cdots,k_{n}\right)p^{x}\left(1-p\right)^{1-x}.
\end{align*}
Thus, for every $\left(k_{1},k_{2},\cdots,k_{n},x\right)\in\left\{ 0,1\right\} ^{n}\times\mathbb{N},$
\begin{align*}
P^{S_{n}=x}_{X}\left[\left(k_{1},k_{2},\cdots,k_{n}\right)\right] & =\boldsymbol{1}_{\left(\sum^{n}_{i=1}l_{i}=x\right)}\left(k_{1},k_{2},\cdots,k_{n}\right)\dfrac{p^{x}\left(1-p\right)^{1-x}}{\binom{n}{x}p^{x}\left(1-p\right)^{x}}.
\end{align*}
Hence,\boxeq{
\[
P^{S_{n}=x}_{X}\left[\left(k_{1},k_{2},\cdots,k_{n}\right)\right]=\boldsymbol{1}_{\left(\sum^{n}_{i=1}l_{i}=x\right)}\left(k_{1},k_{2},\cdots,k_{n}\right)\dfrac{1}{\binom{n}{x}}.
\]
}

That is, for every $x\in\mathbb{N}^{\ast},$ $P^{S_{n}=x}_{X}$ is
the \textbf{uniform law} on the set $\left\{ \left(k_{1},k_{2},\cdots,k_{n}\right)\in\left\{ 0,1\right\} ^{n}:\,\sum^{n}_{i=1}k_{i}=x\right\} .$ 

If $x=0,$ $P^{S_{n}=x}_{X}$ is the \textbf{Dirac measure} at 0.

\end{solution}

\begin{solution}{}{solexercise12.3}

\textbf{1. Computation of $I_{n}\left(s,t\right)=\intop_{\mathbb{R}^{n}}\boldsymbol{1}_{\left(s\leqslant x_{1}\leqslant x_{2}\leqslant\cdots\leqslant x_{n}\leqslant t\right)}\text{d}\lambda_{n}\left(x_{1},x_{2},\cdots,x_{n}\right)$}
\begin{itemize}
\item \textbf{Initialization step}\\
Immediately, we have
\[
I_{1}\left(s,t\right)=\intop_{\mathbb{R}}\boldsymbol{1}_{\left(s\leqslant x_{1}\leqslant t\right)}\text{d\ensuremath{\lambda}}\left(x_{1}\right)=t-s.
\]
By the Fubini theorem,
\[
I_{2}\left(s,t\right)=\intop_{\left[s,t\right]}\left(\intop_{\left[x_{1},t\right]}\text{d}\lambda\left(x_{2}\right)\right)\text{d}\lambda\left(x_{1}\right)=\intop^{t}_{s}\left(t-x_{1}\right)\text{d}x_{1}=\dfrac{\left(t-s\right)^{2}}{2}.
\]
\item \textbf{Induction step}\\
Assume that, for every $s,t$ such that $0\leqslant s\leqslant t,$
\begin{equation}
I_{n}\left(s,t\right)=\dfrac{\left(t-s\right)^{n}}{n!}.\label{eq:I_n_s_t}
\end{equation}
\\
By the Fubini theorem,
\begin{align*}
I_{n+1}\left(s,t\right) & =\intop_{\mathbb{R}}\boldsymbol{1}_{\left(s\leqslant x_{1}\leqslant t\right)}\left(\intop_{\mathbb{R}^{n}}\boldsymbol{1}_{\left(s\leqslant x_{1}\leqslant x_{2}\leqslant\cdots\leqslant x_{n+1}\leqslant t\right)}\text{d}\lambda\left(x_{2},\cdots,x_{n+1}\right)\right)\text{d}\lambda\left(x_{1}\right)\\
 & =\intop_{\mathbb{R}}\boldsymbol{1}_{\left(s\leqslant x_{1}\leqslant t\right)}I_{n}\left(x_{1},t\right)\text{d}\lambda\left(x_{1}\right)=\intop^{t}_{s}\dfrac{\left(t-x_{1}\right)^{n}}{n!}\text{d}x_{1}=\dfrac{\left(t-s\right)^{n+1}}{\left(n+1\right)!}.
\end{align*}
It follows that, for every $n\in\mathbb{N}^{\ast},$\boxeq{
\[
I_{n}\left(s,t\right)=\dfrac{\left(t-s\right)^{n}}{n!}.
\]
}
\end{itemize}
\textbf{2. a. Computation of $\mathbb{E}\left(\boldsymbol{1}_{\left(X_{t}=n\right)}\prod^{n}_{j=1}f_{j}\left(W_{j}\right)\right)$}
\[
\mathbb{E}\left(\boldsymbol{1}_{\left(X_{t}=n\right)}\prod^{n}_{j=1}f_{j}\left(W_{j}\right)\right)=\mathbb{E}\left(\boldsymbol{1}_{\left(W_{n}\leqslant t\right)\cap\left(W_{n+1}>t\right)}\prod^{n}_{j=1}f_{j}\left(W_{j}\right)\right).
\]
Hence, using the random variables $T_{n}$ which carry the probabilistic
information,
\[
\mathbb{E}\left(\boldsymbol{1}_{\left(X_{t}=n\right)}\prod^{n}_{j=1}f_{j}\left(W_{j}\right)\right)=\mathbb{E}\left(\boldsymbol{1}_{\left(\sum^{N}_{i=1}T_{i}\leqslant t\right)\cap\left(\sum^{n+1}_{i=1}T_{i}>t\right)}\prod^{n}_{j=1}f_{j}\left(\sum^{j}_{i=1}T_{i}\right)\right).
\]
By the transfer theorem and since the random variables $T_{n}$ are
independent following an exponential law with parameter $\lambda,$
we obtain
\begin{multline*}
\mathbb{E}\left(\boldsymbol{1}_{\left(X_{t}=n\right)}\prod^{n}_{j=1}f_{j}\left(W_{j}\right)\right)=\intop_{\mathbb{R}^{n+1}}\left[\boldsymbol{1}_{\left(\sum^{N}_{i=1}t_{i}\leqslant t\right)\cap\left(\sum^{n+1}_{i=1}t_{i}>t\right)}\prod^{n}_{j=1}f_{j}\left(\sum^{j}_{i=1}t_{i}\right)\right]\\
\times\left[\prod^{n+1}_{j=1}\boldsymbol{1}_{\mathbb{R}^{+}}\left(t_{j}\right)\right]\lambda^{n+1}\text{e}^{-\lambda\sum^{n+1}_{i=1}t_{i}}\text{d}\lambda_{n+1}\left(t_{1},t_{2},\cdots,t_{n+1}\right).
\end{multline*}
By the change of variables on $\mathbb{R}^{n+1},$ with Jacobian 1,
defined by
\[
\left\{ \begin{array}{l}
w_{1}=t_{1}\\
w_{2}=t_{1}+t_{2}\\
\cdots\cdots\\
w_{n+1}=t_{1}+t_{2}+\cdots+t_{n+1}
\end{array}\right.\Leftrightarrow\left\{ \begin{array}{l}
t_{1}=w_{1}\\
t_{2}=w_{2}-w_{1}\\
\cdots\cdots\\
t_{n+1}=w_{n+1}-w_{n},
\end{array}\right.
\]
it yields
\begin{multline*}
\mathbb{E}\left(\boldsymbol{1}_{\left(X_{t}=n\right)}\prod^{n}_{j=1}f_{j}\left(W_{j}\right)\right)=\intop_{\mathbb{R}^{n+1}}\left(\boldsymbol{1}_{\left(w_{n}\leqslant t\right)\cap\left(w_{n+1}>t\right)}\prod^{n}_{j=1}f_{j}\left(w_{j}\right)\right)\\
\times\left[\boldsymbol{1}_{\mathbb{R}^{+}}\left(w_{1}\right)\prod^{n+1}_{j=2}\boldsymbol{1}_{\mathbb{R}^{+}}\left(w_{j}-w_{j-1}\right)\right]\lambda^{n+1}\text{e}^{-\lambda w_{n+1}}\text{d}\lambda_{n+1}\left(w_{1},w_{2},\cdots,w_{n+1}\right).
\end{multline*}
By the Fubini theorem,
\[
\mathbb{E}\left(\boldsymbol{1}_{\left(X_{t}=n\right)}\prod^{n}_{j=1}f_{j}\left(W_{j}\right)\right)=\intop_{\left[t,+\infty\right[}\lambda\text{e}^{-\lambda w_{n+1}}\times\Phi\left(t\right)\text{d}\lambda\left(w_{n+1}\right),
\]
where
\[
\Phi\left(t\right)=\intop_{\mathbb{R}^{n}}\lambda^{n}\boldsymbol{1}_{\left(w_{n}\leqslant t\right)}\boldsymbol{1}_{\left(0<w_{1}\leqslant w_{2}\leqslant\cdots\leqslant w_{n}\leqslant t\right)}\prod^{n}_{j=1}f_{j}\left(w_{j}\right)\text{d}\lambda_{n}\left(w_{1},w_{2},\cdots,w_{n}\right).
\]
Hence,\boxeq{
\begin{multline}
\mathbb{E}\left(\boldsymbol{1}_{\left(X_{t}=n\right)}\prod^{n}_{j=1}f_{j}\left(W_{j}\right)\right)=\lambda^{n}\text{e}^{-\lambda t}\\
\times\intop_{\mathbb{R}^{n}}\prod^{n}_{j=1}f_{j}\left(w_{j}\right)\boldsymbol{1}_{\left(0<w_{1}\leqslant w_{2}\leqslant\cdots\leqslant w_{n}\leqslant t\right)}\text{d}\lambda_{n}\left(w_{1},w_{2},\cdots,w_{n}\right).\label{eq:exp_of_prod_f_j}
\end{multline}
}

\textbf{b. Law of $X_{t}$ }

In particular, if $f_{j}=1$ for every $j,$ it follows
\[
\mathbb{E}\left(\boldsymbol{1}_{\left(X_{t}=n\right)}\right)=\lambda^{n}\text{e}^{-\lambda t}\times I_{n}\left(0,t\right),
\]
which yields\boxeq{
\[
P\left(X_{t}=n\right)=\text{e}^{-\lambda t}\dfrac{\left(\lambda t\right)^{n}}{n!}.
\]
}That is $X_{t}$ follows a Poisson law with parameter $\lambda t.$ 

Then, \boxeq{
\[
\mathbb{E}\left(X_{t}\right)=\lambda t.
\]
}

\textbf{c. Conditional law of $\left(W_{1},W_{2},\cdots,W_{n}\right)$
given $\left(X_{t}=n\right)$}

Moreover, by taking for instance $f_{j}=\boldsymbol{1}_{A_{j}}$ where
$A_{j}\in\mathscr{B}_{\mathbb{R}},$ it follows from the equality
$\refpar{eq:exp_of_prod_f_j}$ that the random variable $\left(W_{1},W_{2},\cdots,W_{n}\right)$
admits a conditional density given $\left(X_{t}=n\right),$ $f^{\left(X_{t}=n\right)}_{\left(W_{1},W_{2},\cdots,W_{n}\right)},$
given for every $\left(w_{1},w_{2},\cdots,w_{n}\right)\in\mathbb{R}^{n},$
by
\[
f^{\left(X_{t}=n\right)}_{\left(W_{1},W_{2},\cdots,W_{n}\right)}\left(w_{1},w_{2},\cdots,w_{n}\right)=\dfrac{n!}{t^{n}}\boldsymbol{1}_{\left(0\leqslant w_{1}\leqslant w_{2}\leqslant\cdots\leqslant w_{n}\leqslant t\right)},
\]
that is the conditional law of $\left(W_{1},W_{2},\cdots,W_{n}\right)$
given $\left(X_{t}=n\right)$ is the \textbf{\mindex{law!Dirichlet}\index{Dirichlet law}Dirichlet
law}.

\textbf{3.} \textbf{a. Law of $\left(X_{t_{1}},X_{t_{2}}-X_{t_{1}},\cdots,X_{t_{k}}-X_{t_{k-1}}\right).$
Independence of the random variables $X_{t_{1}},X_{t_{2}}-X_{t_{1}},\cdots,X_{t_{k}}-X_{t_{k-1}}$ }

Let $\alpha_{1},\alpha_{2},\cdots,\alpha_{k}$ be arbitrary nonnegative
integers. Denote by $n$ their sum. For $j$ such that $1\leqslant j\leqslant k,$
define
\[
l_{j}=\alpha_{1}+\alpha_{2}+\cdots+\alpha_{j}.
\]
Note that, since $l_{k}=n,$ and that $\sum^{k}_{j=1}\left(t_{j}-t_{j-1}\right)=t,$
\[
\bigcap^{k}_{j=1}\left(X_{t_{j}}-X_{t_{j-1}}=\alpha_{j}\right)\subset\left(X_{t}=n\right).
\]
Thus,
\begin{align*}
\bigcap^{k}_{j=1}\left(X_{t_{j}}-X_{t_{j-1}}=\alpha_{j}\right) & =\left(X_{t}=n\right)\cap\left[\bigcap^{k-1}_{j=1}\left(X_{t_{j}}=l_{j}\right)\right]\\
 & =\left(X_{t}=n\right)\cap\left[\bigcap^{k-1}_{j=1}\left(W_{l_{j}}\leqslant t_{j}\right)\cap\left(W_{l_{j}+1}>t_{j}\right)\right].
\end{align*}
By the equality $\refpar{eq:exp_of_prod_f_j},$ it follows that
\[
P\left(\bigcap^{k}_{j=1}\left(X_{t_{j}}-X_{t_{j-1}}=\alpha_{j}\right)\right)=\lambda^{n}\text{e}^{-\lambda t}\times\psi\left(t_{1},t_{2},\cdots,t_{k}\right),
\]
where
\[
\psi\left(t_{1},t_{2},\cdots,t_{k}\right)=\intop_{\mathbb{R}^{n}}\prod^{k-1}_{j=1}\boldsymbol{1}_{\left(w_{l_{j}}\leqslant t_{j}\right)\cap\left(w_{l_{j}+1}>t_{j}\right)}\boldsymbol{1}_{\left(0\leqslant w_{1}\leqslant w_{2}\leqslant\cdots\leqslant w_{n}\leqslant t\right)}\text{d}\lambda_{n}\left(w_{1},w_{2},\cdots,w_{n}\right).
\]
Setting $l_{0}=0,$
\[
\prod^{k-1}_{j=1}\boldsymbol{1}_{\left(w_{l_{j}}\leqslant t_{j}\right)\cap\left(w_{l_{j}+1}>t_{j}\right)}\boldsymbol{1}_{\left(0\leqslant w_{1}\leqslant w_{2}\leqslant\cdots\leqslant w_{n}\leqslant t\right)}=\prod^{k}_{j=1}\boldsymbol{1}_{\left(t_{j-1}\leqslant w_{l_{j-1}+1}\leqslant\cdots\leqslant w_{l_{j}}\leqslant t_{j}\right).}
\]
Taking into account that $l_{j}-l_{j-1}=\alpha_{j},$ the Fubini theorem
yields
\[
P\left(\bigcap^{k}_{j=1}\left(X_{t_{j}}-X_{t_{j-1}}=\alpha_{j}\right)\right)=\lambda^{n}\text{e}^{-\lambda t}\times\prod^{k}_{j=1}I_{\alpha_{j}}\left(t_{j-1},t_{j}\right).
\]
Thus,\boxeq{
\begin{equation}
P\left(\bigcap^{k}_{j=1}\left(X_{t_{j}}-X_{t_{j-1}}=\alpha_{j}\right)\right)=\prod^{k}_{j=1}\text{e}^{-\lambda\left(t_{j}-t_{j-1}\right)}\dfrac{\left[\lambda\left(t_{j}-t_{j-1}\right)\right]^{\alpha_{j}}}{\alpha_{j}!}.\label{eq:proba_inter_diff_X_t_j}
\end{equation}
}This shows that the random variables $X_{t_{j}}-X_{t_{j-1}}$ are
independent and follow Poisson laws with respective parameter $\lambda\left(t_{j}-t_{j-1}\right).$

\textbf{b. Law of $X_{t}-X_{s}.$ Computation of its expectation $\mathbb{E}\left(X_{t}-X_{s}\right)$ }

Since the $t_{j}$ are arbitrary, it follows that the law of $X_{t}-X_{s}$
is a Poisson law with parameter $\lambda\left(t-s\right)$ and that\boxeq{
\[
\lambda=\dfrac{\mathbb{E}\left(X_{t}-X_{s}\right)}{t-s}.
\]
}

Hence the name, for the parameter $\lambda,$ of \textbf{process intensity}. 

The random variable $X_{t}$ follows a Poisson law with parameter
$\lambda t,$ and by the equality $\refpar{eq:proba_inter_diff_X_t_j},$
after simplifications, we obtain\boxeq{
\[
P^{X_{t}=n}\left(\bigcap^{k}_{j=1}\left(X_{t_{j}}-X_{t_{j-1}}=\alpha_{j}\right)\right)=\dfrac{n!}{\alpha_{1}!\alpha_{2}!\cdots\alpha_{k}!}\prod^{k}_{j=1}\left(\dfrac{t_{j}-t_{j-1}}{t}\right)^{\alpha_{j}}.
\]
}

Hence, the conditional law of $\left(X_{t_{1}},X_{t_{2}}-X_{t_{1}},\cdots,X_{t_{k}}-X_{t_{k-1}}\right)$
given $\left(X_{t}=n\right)$ is the multinomial law $M\left(n,\dfrac{t_{1}}{t},\dfrac{t_{2}-t_{1}}{t},\cdots,\dfrac{t_{k}-t_{k-1}}{t}\right).$

\textbf{Intuitive interpretation}

Let $X_{1},X_{2},\cdots,X_{n}$ be $n$ independent random variables
with uniform law on $\left]0,t\right],$ and let 
\[
\left(Y_{1},Y_{2},\cdots,Y_{k}\right)=\left(\sum^{n}_{j=1}\boldsymbol{1}_{\left]0,t_{1}\right]}\left(X_{j}\right),\sum^{n}_{j=1}\boldsymbol{1}_{\left]t_{1},t_{2}\right]}\left(X_{j}\right),\cdots,\sum^{n}_{j=1}\boldsymbol{1}_{\left]t_{k-1},t_{k}\right]}\left(X_{j}\right)\right)
\]
be the random variable giving the \textbf{number of ``points''}
in every interval $\left]t_{j-1},t_{j}\right].$ Its law is a multinomial
law
\[
M\left(n,\dfrac{t_{1}}{t},\dfrac{t_{2}-t_{1}}{t},\cdots,\dfrac{t_{k}-t_{k-1}}{t}\right).
\]

We have just proved that, for the Poisson process, knowing that the
property has been realized exactly $n$ times in the interval of time
$\left]0,t\right],$ the random variable giving the number of realizations
of the property in each sub-intervals $\left]t_{j-1},t_{j}\right]$---which
form a partition of $\left]0,t\right]$---has the same law as $\left(Y_{1},Y_{2},\cdots,Y_{k}\right),$
and this holds for any chosen partition! This illustrates a form of
\textbf{uniformity} in the times at which the property occurs.

\textbf{4. Determination of a conditional law $P^{X_{t}=n}_{W_{k}}$
of $W_{k}$ given $X_{t}=n$}

Let $f\in\mathscr{C}^{+}_{\mathscr{K}}\left(\mathbb{R}\right).$ Taking
in the equality $\refpar{eq:exp_of_prod_f_j}$ all the functions $f_{j}$
equal to 1, except $f_{k}$ which we take equal to $f,$ we obtain
\[
\mathbb{E}\left(\boldsymbol{1}_{\left(X_{t}=n\right)}f\left(W_{k}\right)\right)=\lambda^{n}\text{e}^{-\lambda t}\times\left(\intop_{\mathbb{R}^{n}}f\left(w_{k}\right)\boldsymbol{1}_{\left(0\leqslant w_{1}\leqslant w_{2}\leqslant\cdots\leqslant w_{n}\leqslant t\right)}\text{d}\lambda_{n}\left(w_{1},w_{2},\cdots,w_{n}\right)\right).
\]
Note that
\[
f\left(w_{k}\right)\boldsymbol{1}_{\left(0\leqslant w_{1}\leqslant w_{2}\leqslant\cdots\leqslant w_{n}\leqslant t\right)}=\boldsymbol{1}_{\left(0\leqslant w_{1}\leqslant w_{2}\leqslant\cdots\leqslant w_{k}\right)}f\left(w_{k}\right)\boldsymbol{1}_{\left(w_{k}\leqslant w_{k+1}\leqslant\cdots\leqslant w_{n}\leqslant t\right)}.
\]
Hence, by the Fubini theorem, integrating first with respect to the
first $k-1$ variables, and using the equality $\refpar{eq:I_n_s_t},$
\begin{equation}
\mathbb{E}\left(\boldsymbol{1}_{\left(X_{t}=n\right)}f\left(W_{k}\right)\right)=\lambda^{n}\text{e}^{-\lambda t}\times R\left(t\right),\label{eq:exp_ind_x_t_is_n_f_w_k}
\end{equation}
where
\[
R\left(t\right)=\intop_{\mathbb{R}^{n-k+1}}\boldsymbol{1}_{\left(0\leqslant w_{k}\leqslant t\right)}\dfrac{\left(w_{k}\right)^{k-1}}{\left(k-1\right)!}f\left(w_{k}\right)\boldsymbol{1}_{\left(w_{k}\leqslant\cdots\leqslant w_{n}\leqslant t\right)}\text{d}\lambda_{n-k+1}\left(w_{k},w_{k+1},\cdots,w_{n}\right).
\]
Again by the Fubini theorem,
\[
R\left(t\right)=\intop_{\mathbb{R}}f\left(w_{k}\right)\boldsymbol{1}_{\left(0\leqslant w_{k}\leqslant t\right)}\dfrac{\left(w_{k}\right)^{k-1}}{\left(k-1\right)!}\times\left[\intop_{\mathbb{R}^{n-k}}\boldsymbol{1}_{\left(w_{k}\leqslant\cdots\leqslant w_{n}\leqslant t\right)}\text{d}\lambda_{n-k+1}\left(w_{k+1},\cdots,w_{n}\right)\right]\text{d}\lambda\left(w_{k}\right).
\]
Hence, by definition of the integrals $I_{n}\left(s,t\right),$
\[
R\left(t\right)=\intop_{\mathbb{R}}f\left(w_{k}\right)\boldsymbol{1}_{\left(0\leqslant w_{k}\leqslant t\right)}\dfrac{\left(w_{k}\right)^{k-1}}{\left(k-1\right)!}\times I_{n-k}\left(w_{k},t\right)\text{d}\lambda\left(w_{k}\right).
\]
Equivalently, using the equality $\refpar{eq:exp_of_prod_f_j},$ 
\[
R\left(t\right)=\intop_{\mathbb{R}}f\left(w_{k}\right)\boldsymbol{1}_{\left(0\leqslant w_{k}\leqslant t\right)}\dfrac{\left(w_{k}\right)^{k-1}}{\left(k-1\right)!}\times\dfrac{\left(t-w_{k}\right)^{n-k}}{\left(n-k\right)!}\text{d}\lambda\left(w_{k}\right).
\]
Substituting into the equality $\refpar{eq:exp_ind_x_t_is_n_f_w_k},$
\[
\mathbb{E}\left(\boldsymbol{1}_{\left(X_{t}=n\right)}f\left(W_{k}\right)\right)=\lambda^{n}\text{e}^{-\lambda t}\times\left[\intop_{\mathbb{R}}f\left(w_{k}\right)\boldsymbol{1}_{\left(0\leqslant w_{k}\leqslant t\right)}\dfrac{\left(w_{k}\right)^{k-1}}{\left(k-1\right)!}\times\dfrac{\left(t-w_{k}\right)^{n-k}}{\left(n-k\right)!}\text{d}\lambda\left(w_{k}\right)\right].
\]
This can also be written
\[
\mathbb{E}\left(\boldsymbol{1}_{\left(X_{t}=n\right)}f\left(W_{k}\right)\right)=\dfrac{n!}{t^{n}}P\left(X_{t}=n\right)\times\left[\intop_{\mathbb{R}}f\left(w_{k}\right)\boldsymbol{1}_{\left(0\leqslant w_{k}\leqslant t\right)}\dfrac{\left(w_{k}\right)^{k-1}}{\left(k-1\right)!}\times\dfrac{\left(t-w_{k}\right)^{n-k}}{\left(n-k\right)!}\text{d}\lambda\left(w_{k}\right)\right].
\]
Thus, there exists a \textbf{conditional density} $f^{X_{t}=n}_{W_{k}}$
of $W_{k}$ given $\left(X_{t}=n\right)$ given by\boxeq{
\[
\forall w\in\mathbb{R},\,\,\,\,f^{X_{t}=n}_{W_{k}}\left(w\right)=\boldsymbol{1}_{\left[0,t\right]}\left(w\right)\dfrac{n!}{\left(k-1\right)!\left(n-k\right)!}\dfrac{1}{t}\left(\dfrac{w}{t}\right)^{k-1}\left(1-\dfrac{w}{t}\right)^{n-k}.
\]
}

\textbf{The conditional law of $W_{k}$ given $\left(X_{t}=n\right)$
is thus a beta law of the first kind $B\left(k,n-k+1\right)$ on $\left[0,t\right].$}

\textbf{Intuitive interpretation}

Consider $n$ independent random variables $X_{1},X_{2},\cdots,X_{n}$
following a uniform law on $\left]0,t\right].$ Let $X_{\left(k\right)}$
be the $k-$th ``\textbf{order statistic}''---see Exercise $\ref{exo:exercise10.8}$
on \textbf{Dirichlet} law and order statistics. We have just proved
that, for the Poisson process, knowing that the property has been
realized exactely $n$ times in the time interval $\left]0,t\right],$
the time of the $k-$th realization of the property is a random variable
following the same law as the one of $X_{\left(k\right)}.$ This again
highlights a \textbf{uniformity} in the times of occurrence of the
property.

\end{solution}

\begin{solution}{}{solexercise12.4}

\textbf{1. $S_{n}$ is a random variable. Conditional law $P^{L=\cdot}_{S_{n}}$}

The random variables $Y=\left(X_{1},X_{2},\cdots,X_{n}\right)$ and
$L$ are independent. Let $f$ be the function defined on $\mathbb{R}^{n}\times\left[0,1\right]$
by
\[
\forall\left(y,l\right)\in\mathbb{R}^{n}\times\left[0,1\right],\,\,\,\,f\left(y,l\right)=\sum^{n}_{j=1}\boldsymbol{1}_{\left[0,l\right]}\left(y_{j}\right).
\]
Then, $P-$almost surely, $S_{n}=f\left(Y,L\right),$ and therefore,
for $P_{L}-$almost every $l,$
\[
P^{L=l}_{S_{n}}=P^{L=l}_{f\left(Y,L\right)}=P_{f\left(Y,L\right)},
\]
where the last equality follows from the independence of $Y$ and
$L.$ Nonetheless, for every $l\in\left]0,1\right],$ $P_{f\left(Y,l\right)}$
is the binomial law $\mathscr{B}\left(n,l\right).$ Hence, for $P_{L}-$almost
every $l$,
\[
P^{L=l}_{S_{n}}=\mathscr{B}\left(n,l\right).
\]

\textbf{2. Law of $S_{n}$}

Then, for every $A\in\mathscr{B}_{\mathbb{R}},$
\begin{align*}
P_{S_{n}}\left(A\right) & =\intop_{\mathbb{R}}P^{L=l}_{S_{n}}\left(A\right)\text{d}P_{L}\left(l\right)\\
 & =\intop_{\mathbb{R}}\boldsymbol{1}_{\left[0,1\right]}\left(l\right)\left[\sum^{n}_{k=0}\binom{n}{k}l^{k}\left(1-l\right)^{n-k}\delta_{k}\left(A\right)\right]\text{d}\lambda\left(l\right)\\
 & =\sum^{n}_{k=0}\binom{n}{k}\delta_{k}\left(A\right)\left[\intop_{\left[0,1\right]}l^{k}\left(1-l\right)^{n-k}\text{d}\lambda\left(l\right)\right]\\
 & =\sum^{n}_{k=0}\binom{n}{k}\delta_{k}\left(A\right)B\left(k+1,n-k+1\right).
\end{align*}
Since,
\[
\binom{n}{k}B\left(k+1,n-k+1\right)=\dfrac{n!}{k!\left(n-k\right)!}\dfrac{\Gamma\left(k+1\right)\Gamma\left(n-k+1\right)}{\Gamma\left(n+2\right)}=\dfrac{1}{n+1},
\]
it follows that\boxeq{
\[
P_{S_{n}}\left(A\right)=\sum^{n}_{k=0}\dfrac{1}{n+1}\delta_{k}\left(A\right),
\]
}that is the law of $S_{n}$ is the uniform law on $\left\llbracket 0,n\right\rrbracket .$

\textbf{3. Conditional law $P^{S_{n}=\cdot}_{L}.$ Conditional expectation
$m^{S_{n}=\cdot}_{L}.$ Computation of $\mathbb{E}\left(L\right)$}

The law of the couple $\left(S_{n},L\right)$ is determined by the
values, for every $A,B\in\mathscr{B}_{\mathbb{R}},$ of $P_{\left(S_{n},L\right)}\left(A\times B\right).$
Hence,
\begin{align*}
P_{\left(S_{n},L\right)}\left(A\times B\right) & =\intop_{B}P^{L=l}_{S_{n}}\left(A\right)\text{d}P_{L}\left(l\right)\\
 & =\intop_{B}\left[\sum^{n}_{k=0}\binom{n}{k}l^{k}\left(1-l\right)^{n-k}\delta_{k}\left(A\right)\right]\boldsymbol{1}_{\left[0,1\right]}\left(l\right)\text{d}\lambda\left(l\right)\\
 & =\sum^{n}_{k=0}\dfrac{1}{n+1}\delta_{k}\left(A\right)\left[\intop_{B}\boldsymbol{1}_{\left[0,1\right]}\left(l\right)\binom{n}{k}\left(n+1\right)l^{k}\left(1-l\right)^{n-k}\text{d}\lambda\left(l\right)\right]\\
 & =\intop_{A}\beta_{1}\left(k+1,n-k+1\right)\left(B\right)\text{d}P_{S_{n}}\left(k\right),
\end{align*}
where $\beta_{1}\left(k+1,n-k+1\right)\left(\cdot\right)$ denotes
the beta law of first kind on $\left[0,1\right]$ with parameters
$k+1$ and $n-k+1.$

For every $k\in\left\llbracket 0,n\right\rrbracket ,$ the conditional
law $P^{S_{n}=k}_{L}$ of $L$ given $S_{n}=k$ is the law $\beta_{1}\left(k+1,n-k+1\right).$

The conditional expectation of $L$ given $S_{n}$ is thus given---refer
to the standard tables of law---by\boxeq{
\[
\forall k\in\left\llbracket 0,n\right\rrbracket ,\,\,\,\,m^{S_{n}=k}_{L}=\dfrac{k+1}{n+2}.
\]
}

We recover the expectation of $L,$ since
\[
\mathbb{E}\left(L\right)=\intop_{\mathbb{R}}m^{S_{n}=k}_{L}\text{d}P_{S_{n}}\left(k\right)=\dfrac{1}{n+1}\sum^{n}_{k=0}\dfrac{k+1}{n+2}.
\]

Hence, since $\sum^{n}_{k=0}\left(k+1\right)=\dfrac{\left(n+1\right)\left(n+2\right)}{2},$\boxeq{
\[
\mathbb{E}\left(L\right)=\dfrac{1}{2}.
\]
}

\end{solution}

\begin{solution}{}{solexercise12.5}

\textbf{1. Law of $\left(\left\lfloor X\right\rfloor ,Y\right)$ as
a function of $f_{X}.$ Laws of $\left\lfloor X\right\rfloor $ and
$Y.$ Direct differentiation of the law of $\left\lfloor X\right\rfloor $}

By the transfer theorem, for every $A,B\in\mathscr{B}_{\mathbb{R}},$
\begin{align*}
P_{\left(\left\lfloor X\right\rfloor ,Y\right)}\left(A\times B\right) & =\intop_{\mathbb{R}}\boldsymbol{1}_{A}\left(\left\lfloor x\right\rfloor \right)\boldsymbol{1}_{B}\left(x-\left\lfloor x\right\rfloor \right)f_{X}\left(x\right)\text{d}\lambda\left(x\right)\\
 & =\sum^{+\infty}_{n=0}\intop_{\left[n,n+1\right[}\boldsymbol{1}_{A}\left(\left\lfloor x\right\rfloor \right)\boldsymbol{1}_{B}\left(x-\left\lfloor x\right\rfloor \right)f_{X}\left(x\right)\text{d}\lambda\left(x\right)\\
 & =\sum^{+\infty}_{n=0}\delta_{n}\left(A\right)\intop_{\left[n,n+1\right[}\boldsymbol{1}_{B}\left(x-n\right)f_{X}\left(x\right)\text{d}\lambda\left(x\right).
\end{align*}
Hence, for every $A,B\in\mathscr{B}_{\mathbb{R}},$\boxeq{
\begin{equation}
P_{\left(\left\lfloor X\right\rfloor ,Y\right)}\left(A\times B\right)=\sum^{+\infty}_{n=0}\delta_{n}\left(A\right)\intop_{\left[0,1\right[}\boldsymbol{1}_{B}\left(x\right)f_{X}\left(x+n\right)\text{d}\lambda\left(x\right).\label{eq:P_couple_intX_and_Y}
\end{equation}
}

The law of $\left\lfloor X\right\rfloor $ is obtained by taking $B=\mathbb{R}.$
For every $A\in\mathscr{B}_{\mathbb{R}},$
\[
P_{\left\lfloor X\right\rfloor }\left(A\right)=P_{\left(\left\lfloor X\right\rfloor ,Y\right)}\left(A\times\mathbb{R}\right)=\sum^{+\infty}_{n=0}\delta_{n}\left(A\right)\intop_{\left[0,1\right[}f_{X}\left(x+n\right)\text{d}\lambda\left(x\right).
\]
This shows that $\left\lfloor X\right\rfloor $ is a discrete random
variable taking values in $\mathbb{N}$ such that\boxeq{
\[
\forall n\in\mathbb{N},\,\,\,\,P\left(\left\lfloor X\right\rfloor =n\right)=\intop_{\left[0,1\right[}f_{X}\left(x+n\right)\text{d}\lambda\left(x\right).
\]
}

This result can be obtained directly by observing that
\[
\left(\left\lfloor X\right\rfloor =n\right)=\left(n\leqslant X<n+1\right),
\]
which yields
\[
P\left(\left\lfloor X\right\rfloor =n\right)=P_{X}\left(\left[n,n+1\right[\right)=\intop_{\left[n,n+1\right[}f_{X}\left(x\right)\text{d}\lambda\left(x\right).
\]

A change of variables then gives,
\[
P\left(\left\lfloor X\right\rfloor =n\right)=\intop_{\left[0,1\right[}f_{X}\left(x+n\right)\text{d}\lambda\left(x\right).
\]

The law of $Y$ is obtained by taking $A=\mathbb{R}.$ For every $B\in\mathscr{B}_{\mathbb{R}},$
\[
P_{Y}\left(B\right)=P_{\left(\left\lfloor X\right\rfloor ,Y\right)}\left(\mathbb{R}\times B\right)=\intop_{B}\boldsymbol{1}_{\left[0,1\right[}\left(x\right)\left(\sum^{+\infty}_{n=0}f_{X}\left(x+n\right)\right)\text{d}\lambda\left(x\right).
\]
Thus, the marginal $Y$ admits a density $f_{Y}$ given by\boxeq{
\[
\forall y\in\mathbb{R},\,\,\,\,f_{Y}\left(y\right)=\boldsymbol{1}_{\left[0,1\right[}\left(y\right)\sum^{+\infty}_{n=0}f_{X}\left(y+n\right).
\]
}

\textbf{2. Determination of $P^{\left\lfloor X\right\rfloor =\cdot}_{Y}$
and $P^{Y=\cdot}_{\left\lfloor X\right\rfloor }$}

For every $A,B\in\mathscr{B}_{\mathbb{R}},$ the equality $\refpar{eq:P_couple_intX_and_Y}$
can also be written as
\[
P_{\left(\left\lfloor X\right\rfloor ,Y\right)}\left(A\times B\right)=\sum_{n\in\text{val}\left(X\right)}\delta_{n}\left(A\right)P_{\left\lfloor X\right\rfloor }\left(\left\{ n\right\} \right)\intop_{B}\dfrac{\boldsymbol{1}_{\left[0,1\right[}\left(x\right)f_{X}\left(x+n\right)}{P_{\left\lfloor X\right\rfloor }\left(\left\{ n\right\} \right)}\text{d}\lambda\left(x\right).
\]

Thus, for every $n\in\text{val}\left(\left\lfloor X\right\rfloor \right),$
$Y$ admits a conditional density given $\left\lfloor X\right\rfloor =n$
given by\boxeq{
\[
\forall y\in\mathbb{R},\,\,\,\,f^{\left\lfloor X\right\rfloor =n}_{Y}\left(y\right)=\boldsymbol{1}_{\left[0,1\right[}\left(y\right)\dfrac{f_{X}\left(y+n\right)}{\intop_{\left[0,1\right[}f_{X}\left(x+n\right)\text{d}\lambda\left(x\right)}.
\]
}

Similarly, for every $A,B\in\mathscr{B}_{\mathbb{R}},$ the equality
$\refpar{eq:P_couple_intX_and_Y}$ can be written, with the convention
$\dfrac{0}{0}=0,$ as
\[
P_{\left(\left\lfloor X\right\rfloor ,Y\right)}\left(A\times B\right)=\intop_{B}\boldsymbol{1}_{\left[0,1\right[}\left(y\right)\dfrac{\sum^{+\infty}_{n\in0}\delta_{n}\left(A\right)f_{X}\left(y+n\right)}{\sum^{+\infty}_{n\in0}f_{X}\left(y+n\right)}f_{Y}\left(y\right)\text{d}\lambda\left(y\right).
\]

Thus, for every $y$ such that $f_{Y}\left(y\right)\neq0,$ $\left\lfloor X\right\rfloor $
admits a conditional density given $Y=y$ such that\boxeq{
\[
\forall A\in\mathscr{B}_{\mathbb{R}},\,\,\,\,P^{Y=y}_{\left\lfloor X\right\rfloor }\left(A\right)=\sum^{+\infty}_{n\in0}\delta_{n}\left(A\right)\dfrac{f_{X}\left(y+n\right)}{\sum^{+\infty}_{n\in0}f_{X}\left(y+n\right)}.
\]
}

\textbf{3. Values of the couple $\left(a,p\right)$ for $\left\lfloor X\right\rfloor $
and $Y$ to be independent. Laws of $\left\lfloor X\right\rfloor $
and $Y$ when $X$ follows the law $\exp\left(p\right),$ $p>0$}

If $P_{X}=\gamma\left(a,p\right),$ then for every $n\in\mathbb{N},$
\[
\intop_{\left[0,1\right[}f_{X}\left(x+n\right)\text{d}\lambda\left(x\right)>0,
\]
which implies that $\text{val}\left(\left\lfloor X\right\rfloor \right)=\mathbb{N}.$ 

Thus, for every $n\in\mathbb{N}$ and for every $y\in\mathbb{R},$
\[
f^{\left\lfloor X\right\rfloor =n}_{Y}\left(y\right)=\boldsymbol{1}_{\left[0,1\right[}\left(y\right)\dfrac{\text{e}^{-py}\left(y+n\right)^{a-1}}{\intop_{\left[0,1\right[}\text{e}^{-px}\left(x+n\right)^{a-1}\text{d}\lambda\left(x\right)}.
\]
This expression is independent of $n$ only if $a=1.$ \textbf{In
this case, that is, if $P_{X}=\text{e}^{p},$ the random variables
$Y$ and $\left\lfloor X\right\rfloor $ are independent.}

If $P_{X}=\text{e}^{p},$ then
\[
P\left(\left\lfloor X\right\rfloor =n\right)=\text{e}^{-pn}\intop_{\left[0,1\right[}p\text{e}^{-px}\text{d}\lambda\left(x\right).
\]
Hence,\boxeq{
\[
P\left(\left\lfloor X\right\rfloor =n\right)=\text{e}^{-pn}\left(1-\textrm{e}^{-p}\right).
\]
} 

Therefore $\left\lfloor X\right\rfloor $ follows the geometric law
on $\mathbb{N}$ with parameter $1-\text{e}^{-p},$ and
\[
\forall y\in\mathbb{R},\,\,\,\,f_{Y}\left(y\right)=\boldsymbol{1}_{\left[0,1\right[}\left(y\right)\sum^{+\infty}_{n=0}p\text{e}^{-p\left(y+n\right)}.
\]
Hence,\boxeq{
\[
\forall y\in\mathbb{R},\,\,\,\,f_{Y}\left(y\right)=\boldsymbol{1}_{\left[0,1\right[}\left(y\right)\dfrac{p}{1-\text{e}^{-p}}\text{e}^{-py}.
\]
}

\textbf{4. Indepencence of $\left\lfloor X\right\rfloor $ and $Y$
for a given $f_{X}$}

In this case, by the equality $\refpar{eq:P_couple_intX_and_Y},$
for every $A,B\in\mathscr{B}_{\mathbb{R}},$
\[
P_{\left(\left\lfloor X\right\rfloor ,Y\right)}\left(A\times B\right)=\sum^{+\infty}_{n=0}\delta_{n}\left(A\right)\intop_{B}\boldsymbol{1}_{\left[0,1\right[}\left(x\right)\times\left[\sum^{+\infty}_{k=0}\boldsymbol{1}_{\left[k,k+1\right[}\left(x+n\right)\text{e}^{-\lambda}\dfrac{\lambda^{k}}{k!}\right]\text{d}\lambda\left(x\right).
\]
Hence,\boxeq{
\[
P_{\left(\left\lfloor X\right\rfloor ,Y\right)}\left(A\times B\right)=\left[\sum^{+\infty}_{n=0}\delta_{n}\left(A\right)\text{e}^{-\lambda}\dfrac{\lambda^{n}}{n!}\right]\left[\intop_{B}\boldsymbol{1}_{\left[0,1\right[}\left(x\right)\text{d}\lambda\left(x\right)\right],
\]
}which shows that \textbf{$Y$ and $\left\lfloor X\right\rfloor $
are independent, with respective laws the uniform law on $\left[0,1\right]$
and the Poisson law $\mathscr{P}\left(\lambda\right).$}

\textbf{5. Computation of $m^{\left\lfloor X\right\rfloor =\cdot}_{Y}$
and $m^{Y=\cdot}_{\left\lfloor X\right\rfloor }$ when $X$ follows
the gamma law $\gamma(a,p)$}

\[
\forall n\in\text{val}\left(\left\lfloor X\right\rfloor \right),\,\,\,\,m^{\left\lfloor X\right\rfloor =n}_{Y}=\intop_{\mathbb{R}}yf^{\left\lfloor X\right\rfloor =n}_{Y}\left(y\right)\text{d}\lambda\left(y\right).
\]
Hence,\boxeq{
\[
\forall n\in\mathbb{N},\,\,\,\,m^{\left\lfloor X\right\rfloor =n}_{Y}=\dfrac{\intop_{\left[0,1\right[}y\text{e}^{-py}\left(y+n\right)^{a-1}\text{d}\lambda\left(y\right)}{\intop_{\left[0,1\right[}\text{e}^{-px}\left(x+n\right)^{a-1}\text{d}\lambda\left(x\right)}.
\]
}

Moreover, for every $y$ such that $f_{Y}\left(y\right)\neq0,$ 
\[
m^{Y=y}_{\left\lfloor X\right\rfloor }=\intop_{\mathbb{R}}x\text{d}P^{Y=y}_{\left\lfloor X\right\rfloor }\left(x\right).
\]
Hence,
\[
\forall y\in\left[0,1\right[,\,\,\,\,m^{Y=y}_{\left\lfloor X\right\rfloor }=\dfrac{\sum^{+\infty}_{n=0}n\text{e}^{-pn}\left(y+n\right)^{a-1}}{\sum^{+\infty}_{n=0}\text{e}^{-pn}\left(y+n\right)^{a-1}}.
\]

\begin{remark}{}{}

In the case $a=2,$ $m^{Y=\cdot}_{\left\lfloor X\right\rfloor }$
is the restriction to the interval $\left[0,1\right[$ of a homographic
function.

\end{remark}

\end{solution}

\begin{solution}{}{solexercise12.6}

\textbf{1. Equivalence between (i) $\text{e}^{\frac{X^{2}}{2}}$ is
$P-$integrable, (ii) $\text{e}^{XY}$ is $P-$integrable, and (iii)
$\text{e}^{\left|XY\right|}$ is $P-$integrable}
\begin{itemize}
\item \textbf{Equivalence of (i) and (ii)}\\
The random vairables $X$ and $Y$ are independent. It follows from
the transfer theorem and the Fubini theorem, in $\overline{\mathbb{R}}^{+},$
\begin{align*}
\intop_{\Omega}\text{e}^{XY}\text{d}P & =\intop_{\mathbb{R}^{2}}\text{e}^{xy}\text{d}P_{X}\otimes P_{Y}\left(x,y\right)\\
 & =\intop_{\mathbb{R}}\left[\intop_{\mathbb{R}}\text{e}^{xy}\text{d}P_{Y}\left(y\right)\right]\text{d}P_{X}\left(x\right).
\end{align*}
But, since $Y$ is Gaussian, using the identity
\[
\text{e}^{xy}\text{e}^{-\frac{y^{2}}{2}}=\text{e}^{\frac{x^{2}}{2}}\text{e}^{-\frac{\left(y-x\right)^{2}}{2}},
\]
we have, for every $x\in\mathbb{R},$\\
\boxeq{
\begin{equation}
\intop_{\mathbb{R}}\text{e}^{xy}\dfrac{1}{\sqrt{2\pi}}\text{e}^{-\frac{y^{2}}{2}}\text{d}\lambda\left(y\right)=\text{e}^{\frac{x^{2}}{2}}.\label{eq:exp_xy_exp-y2over2}
\end{equation}
}Hence,
\begin{equation}
\intop_{\Omega}\text{e}^{XY}\text{d}P=\intop_{\Omega}\text{e}^{\frac{X^{2}}{2}}\text{d}P.\label{eq:int_omega_expXY}
\end{equation}
Thus, the equivalence of (i) and (ii) follows immediately.
\item \textbf{(i) implies (iii)}\\
Moreover, using the following identity, for every $x\in\mathbb{R},$
\[
\text{e}^{x\left|y\right|}\text{e}^{-\frac{y^{2}}{2}}=\boldsymbol{1}_{\mathbb{R}^{+}}\left(y\right)\text{e}^{\frac{x^{2}}{2}}\text{e}^{-\frac{\left(y-x\right)^{2}}{2}}+\boldsymbol{1}_{\mathbb{R}^{-\ast}}\left(y\right)\text{e}^{\frac{x^{2}}{2}}\text{e}^{-\frac{\left(y+x\right)^{2}}{2}},
\]
it yields by integrating over $\mathbb{R}^{+}$ and $\mathbb{R}^{-\ast},$
\[
\intop_{\mathbb{R}}\text{e}^{\left|xy\right|}\dfrac{1}{\sqrt{2\pi}}\text{e}^{-\frac{y^{2}}{2}}\text{d}\lambda\left(y\right)\leqslant2\text{e}^{\frac{x^{2}}{2}}.
\]
It follows that
\[
\intop_{\Omega}\text{e}^{\left|XY\right|}\text{d}P\leqslant2\intop_{\Omega}\text{e}^{\frac{X^{2}}{2}}\text{d}P,
\]
which shows that (i) implies (iii).
\item \textbf{(iii) implies (i)} since $\text{e}^{XY}\leqslant\text{e}^{\left|XY\right|}.$
\end{itemize}
\textbf{2. (a) Proof of $\mathbb{E}^{\mathscr{B}}\left(\text{e}^{XY}\right)\geqslant1\,\,\,\,P-\text{almost surely}$}

By the Jensen inequality,
\[
\text{e}^{\mathbb{E}^{\mathscr{B}}\left(XY\right)}\leqslant\mathbb{E}^{\mathscr{B}}\left(\text{e}^{XY}\right).
\]

But, since $X$ is $\mathscr{B}-$measurable,
\[
\mathbb{E}^{\mathscr{B}}\left(XY\right)=X\mathbb{E}^{\mathscr{B}}\left(Y\right).
\]
The $\sigma-$algebra $\mathscr{B}$ and $\sigma\left(Y\right)$ are
independent and, as $Y$ is centered, it follows that
\[
\mathbb{E}^{\mathscr{B}}\left(XY\right)=X\mathbb{E}\left(Y\right)=0.
\]
Thus,\boxeq{
\[
1\leqslant\mathbb{E}^{\mathscr{B}}\left(\text{e}^{XY}\right)\,\,\,\,P-\text{almost surely.}
\]
}

\textbf{(b) Computation of $\mathbb{E}^{\mathscr{B}}\left(\text{e}^{XY}\right),$
when $\mathscr{B}=\sigma\left(X\right)$}

If $\mathscr{B}=\sigma\left(X\right),$ then a version of $\mathbb{E}^{\mathscr{B}}\left(\text{e}^{XY}\right)$
is obtained by composing with $X$ the conditional expectation of
$\text{e}^{XY}$ given $X.$ By the conditional transfer theorem,
for $P_{X}-$almost every $x,$
\[
m^{X=x}_{\text{e}^{XY}}=m^{X=x}_{\text{e}^{xY}}=\mathbb{E}\left(\text{e}^{xY}\right),
\]
where the last equality follows from the independence of $X$ and
$Y.$ Hence, by the equality $\refpar{eq:exp_xy_exp-y2over2},$ it
yields
\[
m^{X=s}_{\text{e}^{XY}}=\text{e}^{\frac{x^{2}}{2}}.
\]
It follows that\boxeq{
\[
\mathbb{E}^{\mathscr{B}}\left(\text{e}^{XY}\right)=\text{e}^{\frac{X^{2}}{2}}\,\,\,\,P-\text{almost surely}
\]
}

\textbf{(c) Computation of $\mathbb{E}^{\mathscr{B}}\left(\text{e}^{XY}\right),$
general case (2 methods)}
\begin{itemize}
\item \textbf{First method}\\
Since, for every $p\in\mathbb{N},$
\[
\left|\sum^{p}_{n=0}\dfrac{\left(XY\right)^{n}}{n!}\right|\leqslant\sum^{+\infty}_{n=0}\dfrac{\left|XY\right|^{n}}{n!}=\text{e}^{\left|XY\right|},
\]
and $\text{e}^{\left|XY\right|}$ is $P-$integrable, the dominated
convergence theorem for conditional expectations implies that
\[
\mathbb{E}^{\mathscr{B}}\left(\text{e}^{XY}\right)=\sum^{+\infty}_{n=0}\mathbb{E}^{\mathscr{B}}\left(\dfrac{\left(XY\right)^{n}}{n!}\right).
\]
But, since $X^{n}$ is $\mathscr{B}-$measurable,
\[
\mathbb{E}^{\mathscr{B}}\left(\left(XY\right)^{n}\right)=X^{n}\mathbb{E}^{\mathscr{B}}\left(Y^{n}\right).
\]
As the $\sigma-$algebra $\mathscr{B}$ and $\sigma\left(Y\right)$
are independent, it follows that
\[
\mathbb{E}^{\mathscr{B}}\left(\left(XY\right)^{n}\right)=X^{n}\mathbb{E}\left(Y^{n}\right).
\]
Since the random variable $Y$ is Gaussian, centered and reduced,
a classical computation of the moments---via integration by parts---yields,
for every $p\in\mathbb{N}^{\ast},$\\
\boxeq{
\[
\mathbb{E}\left(Y^{2p+1}\right)=0\,\,\,\,\text{and}\,\,\,\,\mathbb{E}\left(Y^{2p}\right)=\dfrac{\left(2p-1\right)!}{2^{p-1}\left(p-1\right)!}.
\]
}Consequently,
\[
\mathbb{E}^{\mathscr{B}}\left(\text{e}^{XY}\right)=\sum^{+\infty}_{p=0}\dfrac{\left(2p-1\right)!}{2^{p-1}\left(p-1\right)!\left(2p\right)!}X^{2p}=\sum^{+\infty}_{p=0}\dfrac{1}{p!}\left(\dfrac{X^{2}}{2}\right)^{p}.
\]
Thus,\boxeq{
\[
\mathbb{E}^{\mathscr{B}}\left(\text{e}^{XY}\right)=\text{e}^{\frac{X^{2}}{2}}\,\,\,\,P-\text{almost surely}.
\]
}
\item \textbf{Second method}\\
Let $f$ be the function defined by
\[
\forall x\in\mathbb{R},\,\,\,\,f\left(x\right)=\mathbb{E}\left(\text{e}^{xY}\right),
\]
Then equality $\refpar{eq:exp_xy_exp-y2over2}$ can be written as
\[
\forall x\in\mathbb{R},\,\,\,\,\,f\left(x\right)=\text{e}^{\frac{x^{2}}{2}},
\]
and Proposition $\ref{pr:cond_exp_exist_indep}$ states that $f\circ X$
is a version of $\mathbb{E}^{\mathscr{B}}\left(\text{e}^{XY}\right).$\\
\begin{remark}{}{}We recover the previous result. Indeed, by the
inclusion of $\sigma-$algebras $\mathscr{B}\supset\sigma\left(X\right)$
and the $\sigma\left(X\right)-$measurability of the random variable
$\text{e}^{\frac{X^{2}}{2}},$
\[
\mathbb{E}^{\sigma\left(X\right)}\left(\text{e}^{XY}\right)=\mathbb{E}^{\sigma\left(X\right)}\left(\mathbb{E}^{\mathscr{B}}\left(\text{e}^{XY}\right)\right)=\text{e}^{\frac{X^{2}}{2}}\,\,\,\,P-\text{almost surely.}
\]
\end{remark}
\end{itemize}
\textbf{3. Computation $\mathbb{E}^{\mathscr{B}}\left(\text{e}^{XY}\right),$
when $\text{e}^{\frac{X^{2}}{2}}$ is not necessarily $P-$integrable}

For every $n\in\mathbb{N},$ consider the random variable $X_{n}=\boldsymbol{1}_{\left(\left|X\right|\leqslant n\right)}X.$
Then $X_{n}$ is bounded by $n,$ and therefore $\text{e}^{\frac{X^{2}_{n}}{2}}$
is $P-$integrable. Hence, for every $n\in\mathbb{N},$
\[
\mathbb{E}^{\mathscr{B}}\left(\text{e}^{X_{n}Y}\right)=\text{e}^{\frac{X^{2}_{n}}{2}}\,\,\,\,P-\text{almost surely.}
\]
Moreover, the sequence of nonnegative general term $\text{e}^{\frac{X^{2}_{n}}{2}}$
converges $P-$almost surely by nondecreasing to $\text{e}^{\frac{X^{2}}{2}}.$
Applying the conditional Fatou-Lebesgue theorem $\ref{th:cond_fatou_lebes_dom_cv}$
yields
\begin{multline*}
\mathbb{E}^{\mathscr{B}}\left(\text{e}^{XY}\right)=\mathbb{E}^{\mathscr{B}}\left(\liminf_{n\to+\infty}\text{e}^{X_{n}Y}\right)\leqslant\liminf_{n\to+\infty}\left[\mathbb{E}^{\mathscr{B}}\left(\text{e}^{X_{n}Y}\right)\right]\\
\leqslant\limsup_{n\to+\infty}\left[\mathbb{E}^{\mathscr{B}}\left(\text{e}^{X_{n}Y}\right)\right]\leqslant\mathbb{E}^{\mathscr{B}}\left(\limsup_{n\to+\infty}\text{e}^{X_{n}Y}\right)=\mathbb{E}^{\mathscr{B}}\left(\text{e}^{XY}\right).
\end{multline*}
It follows that the sequence with general term $\mathbb{E}^{\mathscr{B}}\left(\text{e}^{X_{n}Y}\right)$
converges $P-$almost surely and that
\[
\mathbb{E}^{\mathscr{B}}\left(\text{e}^{XY}\right)=\lim_{n\to+\infty}\mathbb{E}^{\mathscr{B}}\left(\text{e}^{X_{n}Y}\right).
\]
Therefore,\boxeq{
\[
\mathbb{E}^{\mathscr{B}}\left(\text{e}^{XY}\right)=\text{e}^{\frac{X^{2}}{2}}\,\,\,\,P-\text{almost surely.}
\]
}

\end{solution}

\begin{solution}{}{solexercise12.7}

For every $A_{1}\in\mathscr{A}_{1}$ and $A_{2}\in\mathscr{A}_{2},$
since $\boldsymbol{1}_{A_{1}}X$ and $\boldsymbol{1}_{A_{2}}$ are
independent,
\[
\mathbb{E}\left(\boldsymbol{1}_{A_{1}\cap A_{2}}X\right)=\mathbb{E}\left(\boldsymbol{1}_{A_{1}}\boldsymbol{1}_{A_{2}}X\right)=\mathbb{E}\left(\boldsymbol{1}_{A_{1}}X\right)\mathbb{E}\left(\boldsymbol{1}_{A_{2}}\right),
\]
 and thus, by the definition of $\mathbb{E}^{\mathscr{A}_{1}}\left(X\right),$
\[
\mathbb{E}\left(\boldsymbol{1}_{A_{1}\cap A_{2}}X\right)=\mathbb{E}\left(\boldsymbol{1}_{A_{1}}\mathbb{E}^{\mathscr{A}_{1}}\left(X\right)\right)\mathbb{E}\left(\boldsymbol{1}_{A_{2}}\right).
\]

But $\boldsymbol{1}_{A_{1}}\mathbb{E}^{\mathscr{A}_{1}}\left(X\right)$
and $\boldsymbol{1}_{A_{2}}$ are independent, hence
\[
\mathbb{E}\left(\boldsymbol{1}_{A_{1}\cap A_{2}}X\right)=\mathbb{E}\left(\left(\boldsymbol{1}_{A_{1}}\mathbb{E}^{\mathscr{A}_{1}}\left(X\right)\right)\boldsymbol{1}_{A_{2}}\right)=\mathbb{E}\left(\boldsymbol{1}_{A_{1}\cap A_{2}}\mathbb{E}^{\mathscr{A}_{1}}\left(X\right)\right).
\]

Since $\left\{ A_{1}\cap A_{2}:\,A_{1}\in\mathscr{A}_{1}\,\text{and}\,A_{2}\in\mathscr{A}_{2}\right\} $
is a $\pi-$system that generates the $\sigma-$algebra $\mathscr{A}_{1}\vee\mathscr{A}_{2},$
it follows from the measurability expansion theorem that
\[
\forall A\in\mathscr{A}_{1}\vee\mathscr{A}_{2},\,\,\,\,\mathbb{E}\left(\boldsymbol{1}_{A}X\right)=\mathbb{E}\left(\boldsymbol{1}_{A}\mathbb{E}^{\mathscr{A}_{1}}\left(X\right)\right).
\]

To conclude, it suffices to note that $\mathbb{E}^{\mathscr{A}_{1}}\left(X\right)$
is $\mathscr{A}_{1}\vee\mathscr{A}_{2}-$measurable.

\end{solution}

\begin{solution}{}{solexercise12.8}

Taking into account the independence of the $X_{i},$ for every Borel
subset $C$ of $\mathbb{R},$
\[
\intop_{S^{-1}_{n}\left(C\right)}X_{i}\text{d}P=\intop_{\mathbb{R}^{n}}\boldsymbol{1}_{C}\left(x_{1}+x_{2}+\cdots+x_{n}\right)x_{i}\text{d}P_{X_{1}}\otimes\text{d}P_{X_{2}}\otimes\cdots\otimes\text{d}P_{X_{n}}\left(x_{1},x_{2},\cdots,x_{n}\right).
\]
Since the $X_{i}$ have the same law $\mu,$
\begin{align*}
\intop_{S^{-1}_{n}\left(C\right)}X_{i}\text{d}P & =\intop_{\mathbb{R}^{n}}\boldsymbol{1}_{C}\left(x_{1}+x_{2}+\cdots+x_{n}\right)x_{i}\text{d}\mu^{\otimes n}\left(x_{1},x_{2},\cdots,x_{n}\right)\\
 & =\intop_{\mathbb{R}^{n}}\boldsymbol{1}_{C}\left(x_{1}+x_{2}+\cdots+x_{n}\right)x_{1}\text{d}\mu^{\otimes n}\left(x_{1},x_{2},\cdots,x_{n}\right).
\end{align*}
It follows that
\[
\forall C\in\mathscr{B}_{\mathbb{R}},\,\,\,\,\intop_{S^{-1}_{n}\left(C\right)}X_{i}\text{d}P=\intop_{S^{-1}_{n}\left(C\right)}X_{1}\text{d}P,
\]
which established the required equality.

Next, since $S_{n}$ is $\sigma\left(S_{n}\right)-$measurable,
\[
\mathbb{E}^{\sigma\left(S_{n}\right)}\left(S_{n}\right)=S_{n}=n\mathbb{E}^{\sigma\left(S_{n}\right)}\left(X_{1}\right),
\]
 and therefore\boxeq{
\[
\mathbb{E}^{\sigma\left(S_{n}\right)}\left(X_{1}\right)=\dfrac{S_{n}}{n}.
\]
}

Finally, since for every $k\in\mathbb{N}^{\ast},$
\[
S_{n+k}=S_{n}+\sum^{n+k}_{j=n+1}X_{j},
\]
we have the equality of $\sigma-$algebras
\[
\mathscr{A}_{n}=\sigma\left(S_{n}\right)\vee\sigma\left(\sum^{n+k}_{j=n+1}X_{j}:\,k\in\mathbb{N}^{\ast}\right).
\]
It follows from the previous exercise that
\[
\mathbb{E}^{\mathscr{A}_{n}}\left(X_{1}\right)=\mathbb{E}^{\sigma\left(S_{n}\right)}\left(X_{1}\right),
\]
which yields\boxeq{
\[
\mathbb{E}^{\mathscr{A}_{n}}\left(X_{1}\right)=\dfrac{S_{n}}{n}.
\]
}

\end{solution}

\begin{solution}{}{solexercise12.9}

\textbf{1. Computation of $I\left(\varphi\right)\equiv\mathbb{E}\left(\boldsymbol{1}_{\left(M_{n}\in G\right)}\varphi\left(X_{n}\right)\right).$
Computation of $P\left(M_{n}\in G\right)$}

Since $\varphi$ is bounded, we may apply the transfer theorem. Moreover,
$X_{n}$ and $Y_{n}$ are independent and have respective densities
$g$ and $\boldsymbol{1}_{\left[0,\overline{t}\right]}.$ Hence,
\[
I\left(\varphi\right)\equiv\mathbb{E}\left(\boldsymbol{1}_{\left(M_{n}\in G\right)}\varphi\left(X_{n}\right)\right)=\intop_{\mathbb{R}^{2}}\boldsymbol{1}_{G}\left(x,y\right)\varphi\left(x\right)\dfrac{1}{\overline{t}}\boldsymbol{1}_{\left[0,\overline{t}\right]}\left(y\right)g\left(x\right)\text{d}\lambda_{2}\left(x,y\right).
\]
By the Fubini theorem, since $\varphi$ is bounded,
\[
I\left(\varphi\right)=\intop_{\mathbb{R}}\varphi\left(x\right)g\left(x\right)\left[\intop_{\left[0,t\left(x\right)\right]}\dfrac{1}{\overline{t}}\text{d}\lambda\left(y\right)\right]\text{d}\lambda\left(x\right)=\dfrac{1}{\overline{t}}\intop_{\left\{ g\neq0\right\} }\varphi\left(x\right)f\left(x\right)\text{d}\lambda\left(x\right).
\]
The inclusion $\left\{ g=0\right\} \subset\left\{ f=0\right\} $ follows
from inequalities: for every $x\in\mathbb{R},$ $0\leqslant f\left(x\right)\leqslant\overline{t}g\left(x\right),$
and $\overline{t}>0.$ It follows that
\[
\intop_{\left\{ g=0\right\} }\left|\varphi\left(x\right)\right|f\left(x\right)\text{d}\lambda\left(x\right)\leqslant\intop_{\left\{ f=0\right\} }\left|\varphi\left(x\right)\right|f\left(x\right)\text{d}\lambda\left(x\right)=0,
\]
 and thus,\boxeq{
\[
I\left(\varphi\right)=\dfrac{1}{\overline{t}}\intop_{\mathbb{R}}\varphi\left(x\right)f\left(x\right)\text{d}\lambda\left(x\right).
\]
}In particular, if $\varphi=1,$ then\boxeq{
\[
P\left(M_{n}\in G\right)=\dfrac{1}{\overline{t}}.
\]
}

\textbf{2. $T$ and $X_{T}$ are random variables. Law of $T.$ $X_{T}$
is $P-$almost surely finite. Expectation of $T$}

The random variable $T$ takes values in $\overline{\mathbb{N}},$
since $\left(T=1\right)=\left(M_{1}\in G\right)\in\mathscr{A}$ and,
for every $n\geqslant2,$
\[
\left(T=n\right)=\left[\bigcap^{n-1}_{j=1}\left(M_{j}\notin G\right)\cap\left(M_{n}\in G\right)\right]\in\mathscr{A}.
\]

Hence $\left(T=+\infty\right)=\left(T\in\mathbb{N}^{\ast}\right)^{c}\in\mathscr{A}.$
Moreover, since the random variables $M_{n}$ are independent and
following the same law,
\[
P\left(T=1\right)=P\left(M_{1}\in G\right)=\dfrac{1}{\overline{t}},
\]
and, for every $n\geqslant2,$
\[
P\left(T=n\right)=\prod^{n-1}_{j=1}P\left(M_{j}\notin G\right)\times P\left(M_{n}\in G\right)=\dfrac{1}{\overline{t}}\left(1-\dfrac{1}{\overline{t}}\right)^{n-1},
\]
Thus $T$ follows a geometric law on $\mathbb{N}^{\ast},$ $\mathscr{G}_{\mathbb{N}^{\ast}}\left(\dfrac{1}{\overline{t}}\right).$
In particular $T$ is $P-$almost surely finite, and since
\[
\left(X_{T}=+\infty\right)\subset\left(T=+\infty\right),
\]
it follows that $X_{T}$ is $P-$almost surely finite.

Moreover,\boxeq{
\[
\mathbb{E}\left(T\right)=\overline{t}.
\]
}

\textbf{3. Law of $X_{T}$}

For every $\varphi\in\mathscr{C}^{+}_{\mathscr{K}}\left(\mathbb{R}\right),$
\begin{align*}
\mathbb{E}\left(\varphi\left(X_{T}\right)\right) & =\sum_{n\in\mathbb{N}^{\ast}}\intop_{\left(T=n\right)}\varphi\left(X_{n}\right)\text{d}P\\
 & =\intop_{\left(M_{1}\in G\right)}\varphi\left(X_{1}\right)\text{d}P+\sum_{n\geqslant2}\intop_{\Omega}\left[\prod^{n-1}_{j=1}\boldsymbol{1}_{\left(M_{j}\notin G\right)}\right]\boldsymbol{1}_{\left(M_{n}\in G\right)}\varphi\left(X_{n}\right)\text{d}P.
\end{align*}
Since the random variables $\prod^{n-1}_{j=1}\boldsymbol{1}_{\left(M_{j}\notin G\right)}$
and $\boldsymbol{1}_{\left(M_{n}\in G\right)}\varphi\left(X_{n}\right)$
are independent, it follows from the Question 1 that
\begin{align*}
\mathbb{E}\left(\varphi\left(X_{T}\right)\right) & =I\left(\varphi\right)\left[1+\sum_{n\geqslant2}\prod^{n-1}_{j=1}P\left(M_{j}\notin G\right)\right]\\
 & =I\left(\varphi\right)\left[1+\sum_{n\geqslant2}\left(1-\dfrac{1}{\overline{t}}\right)^{n-1}\right]\\
 & =I\left(\varphi\right)\dfrac{1}{1-\left(1-\dfrac{1}{\overline{t}}\right)}\\
 & =\overline{t}I\left(\varphi\right).
\end{align*}
Hence, we proved that\boxeq{
\[
\forall\varphi\in\mathscr{C}^{+}_{\mathscr{K}}\left(\mathbb{R}\right),\,\,\,\,\mathbb{E}_{\varphi}\left(X_{T}\right)=\intop_{\mathbb{R}}\varphi\left(x\right)f\left(x\right)\text{d}\lambda\left(x\right).
\]
}That is, \textbf{$X_{T}$ admits the density $f.$}

\textbf{4. (a) $\left(T_{k}=n\right)\in\mathscr{A}_{n}$ and $T_{k}$
is $P-$almost surely finite. $\mathscr{A}_{T_{k}}\subset\mathscr{A}_{T_{k+1}}$}

Let us show by induction on $k$ the property $\left(T_{k}=n\right)\in\mathscr{A}_{n}$
and $T_{k}$ is $P-$almost surely finite. 
\begin{itemize}
\item \textbf{Initialization step}\\
The property has already been established for $T_{1}.$
\item \textbf{Induction step}\\
Assume the property hold up to $k.$ If $n<k,$ then $\left(T_{k+1}=n\right)=\emptyset\in\mathscr{A}_{n}.$\\
If $n\geqslant k+1,$ then
\begin{equation}
\left(T_{k+1}=n\right)=\biguplus^{n-1}_{j=1}\left[\left(T_{k}=j\right)\cap\bigcap^{n-1}_{i=j+1}\left(M_{i}\notin G\right)\cap\left(M_{n}\in G\right)\right],\label{eq:T_k_plus_1_eq_n}
\end{equation}
with the convention $\bigcap^{n-1}_{i=n}\left(M_{i}\notin G\right)=\Omega.$\\
By the induction hypothesis and by the fact that for $1\leqslant j\leqslant n-1,$
we have $\left(T_{k}=j\right)\in\mathscr{A}_{j}\subset\mathscr{A}_{n},$
it follows that $\left(T_{k+1}=n\right)\in\mathscr{A}_{n}.$\\
We already know that $T_{1}$ is $P-$almost surely finite. Assume
that $T_{k}$ is $P-$almost surely finite. From the previous equalities,
for every $n\in\mathbb{N}^{\ast},$
\[
P\left(T_{k+1}=n\right)=\dfrac{1}{\overline{t}}P\left(T_{k}=n-1\right)+\sum^{n-2}_{j=1}\left[P\left(T_{k}=j\right)\dfrac{1}{\overline{t}}\left(1-\dfrac{1}{\overline{t}}\right)^{n-j-1}\right].
\]
Summing over $\mathbb{N}^{\ast},$ we obtain
\begin{align*}
\sum_{n\in\mathbb{N}^{\ast}}P\left(T_{k+1}=n\right) & =\dfrac{1}{\overline{t}}\sum_{n\in\mathbb{N}^{\ast}}P\left(T_{k}=n-1\right)+\sum_{n\in\mathbb{N}^{\ast}}\sum^{n-2}_{j=1}\left[P\left(T_{k}=j\right)\dfrac{1}{\overline{t}}\left(1-\dfrac{1}{\overline{t}}\right)^{n-j-1}\right]\\
 & =\dfrac{1}{\overline{t}}P\left(T_{k}<+\infty\right)+\dfrac{1}{\overline{t}}\sum_{j\in\mathbb{N}^{\ast}}\left[P\left(T_{k}=j\right)\sum_{n\geqslant j+2}\left(1-\dfrac{1}{\overline{t}}\right)^{n-j-1}\right]\\
 & =\dfrac{1}{\overline{t}}P\left(T_{k}<+\infty\right)+\dfrac{1}{\overline{t}}\sum_{j\in\mathbb{N}^{\ast}}\left[P\left(T_{k}=j\right)\left(1-\dfrac{1}{\overline{t}}\right)\overline{t}\right]\\
 & =P\left(T_{k}<+\infty\right).
\end{align*}
Hence $P\left(T_{k}<+\infty\right)=1,$ and the claim holds for every
$k.$
\end{itemize}
Let $A\in\mathscr{A}_{T_{k}},$ the equality $\refpar{eq:T_k_plus_1_eq_n}$
shows that for every $n\in\mathbb{N}^{\ast},$
\[
A\cap\left(T_{k+1}=n\right)\in\mathscr{A}_{n},
\]
that is, $A\in\mathscr{A}_{T_{k+1}}.$ Consequently, 
\[
\mathscr{A}_{T_{k}}\subset\mathscr{A}_{T_{k+1}}.
\]

\textbf{(b) Computation of $\mathbb{E}^{\mathscr{A}_{n}}\left(\boldsymbol{1}_{\left(T_{k}=n\right)}f_{k+1}\left(X_{T_{k+1}}\right)\right).$
Computation of $\mathbb{E}^{\mathscr{A}_{T_{k}}}\left(f_{k+1}\left(X_{T_{k+1}}\right)\right).$
Law of $X_{T_{k+1}}$}

We write
\begin{align*}
\boldsymbol{1}_{\left(T_{k}=n\right)}f_{k+1}\left(X_{T_{k+1}}\right) & =\sum_{l\in\mathbb{N}^{\ast}}\boldsymbol{1}_{\left(T_{k}=n\right)}f_{k+1}\left(X_{n+l}\right)\boldsymbol{1}_{\left(T_{k+1}=n+l\right)}\\
 & =\sum_{l\in\mathbb{N}^{\ast}}\left[\boldsymbol{1}_{\left(T_{k}=n\right)}f_{k+1}\left(X_{n+l}\right)\prod^{l-1}_{j=1}\boldsymbol{1}_{\left(M_{n+j}\notin G\right)}\boldsymbol{1}_{\left(M_{n+l}\in G\right)}\right],
\end{align*}
with the convention $\prod^{0}_{j=1}\boldsymbol{1}_{\left(M_{n+j}\notin G\right)}=1.$

Taking the conditional expectation with respect to $\mathscr{A}_{n},$
and using that $\left(T_{k}=n\right)\in\mathscr{A}_{n},$ we obtain
\begin{align*}
 & \mathbb{E}^{\mathscr{A}_{n}}\left(\boldsymbol{1}_{\left(T_{k}=n\right)}f_{k+1}\left(X_{T_{k+1}}\right)\right)\\
 & \qquad\qquad=\sum_{l\in\mathbb{N}^{\ast}}\boldsymbol{1}_{\left(T_{k}=n\right)}\mathbb{E}^{\mathscr{A}_{n}}\left(f_{k+1}\left(X_{n+l}\right)\prod^{l-1}_{j=1}\boldsymbol{1}_{\left(M_{n+j}\notin G\right)}\boldsymbol{1}_{\left(M_{n+l}\in G\right)}\right)\\
 & \qquad\qquad=\sum_{l\in\mathbb{N}^{\ast}}\boldsymbol{1}_{\left(T_{k}=n\right)}\mathbb{E}\left(f_{k+1}\left(X_{n+l}\right)\prod^{l-1}_{j=1}\boldsymbol{1}_{\left(M_{n+j}\notin G\right)}\boldsymbol{1}_{\left(M_{n+l}\in G\right)}\right).
\end{align*}
The last equality above comes from the independence of the $\sigma-$algebras
$\sigma\left(M_{n+j}:\,j\in\mathbb{N}^{\ast}\right)$ and $\mathscr{A}_{n}.$
It then yields
\[
\boldsymbol{1}_{\left(T_{k}=n\right)}\mathbb{E}^{\mathscr{A}_{n}}\left(f_{k+1}\left(X_{T_{k+1}}\right)\right)=\boldsymbol{1}_{\left(T_{k}=n\right)}\sum_{l\in\mathbb{N}^{\ast}}\left(1-\dfrac{1}{\overline{t}}\right)^{l-1}I\left(f_{k+1}\right).
\]

Hence,\boxeq{
\[
\boldsymbol{1}_{\left(T_{k}=n\right)}\mathbb{E}^{\mathscr{A}_{n}}\left(f_{k+1}\left(X_{T_{k+1}}\right)\right)=\boldsymbol{1}_{\left(T_{k}=n\right)}\overline{t}I\left(f_{k+1}\right).
\]
}

Let $A\in\mathscr{A}_{T_{k}}.$ Since $A\cap\left(T_{k}=n\right)\in\mathscr{A}_{n}$
and that $T_{k}$ is finite $P-$almost surely,
\begin{align*}
\mathbb{E}\left(\boldsymbol{1}_{A}f_{k+1}\left(X_{T_{k+1}}\right)\right) & =\sum_{n\in\mathbb{N}^{\ast}}\mathbb{E}\left(\boldsymbol{1}_{A\cap\left(T_{k}=n\right)}f_{k+1}\left(X_{T_{k+1}}\right)\right)\\
 & =\sum_{n\in\mathbb{N}^{\ast}}\mathbb{E}\left(\boldsymbol{1}_{A}\boldsymbol{1}_{\left(T_{k}=n\right)}\overline{t}I\left(f_{k+1}\right)\right).
\end{align*}
Thus, for every $A\in\mathscr{A}_{T_{k}},$
\[
\mathbb{E}\left(\boldsymbol{1}_{A}f_{k+1}\left(X_{T_{k+1}}\right)\right)=\mathbb{E}\left(\boldsymbol{1}_{A}\overline{t}I\left(f_{k+1}\right)\right),
\]
which proves---since $\overline{t}I\left(f_{k+1}\right)$ is constant
and hence $\mathscr{A}_{T_{k}}-$measurable---that\boxeq{
\[
\mathbb{E}^{\mathscr{A}_{T_{k}}}\left(f_{k+1}\left(X_{T_{k+1}}\right)\right)=\overline{t}I\left(f_{k+1}\right),
\]
}Taking expectation,
\[
\mathbb{E}\left(f_{k+1}\left(X_{T_{k+1}}\right)\right)=\mathbb{E}\left(\mathbb{E}^{\mathscr{A}_{T_{k}}}\left(f_{k+1}\left(X_{T_{k+1}}\right)\right)\right)=\overline{t}I\left(f_{k+1}\right).
\]
That is,
\[
\mathbb{E}\left(f_{k+1}\left(X_{T_{k+1}}\right)\right)=\intop_{\mathbb{R}}f_{k+1}\left(x\right)f\left(x\right)\text{d}\lambda\left(x\right).
\]
Thus, \textbf{$X_{T_{k+1}}$ still admits a density $f.$}

\textbf{(c) $f_{k}\left(X_{T_{k}}\right)$ is $\mathscr{A}_{T_{k}}-$measurable.
$\left(X_{T_{k}}\right)_{k\in\mathbb{N}^{\ast}}$ is a sequence of
independent random variables}

Note that for every $k\in\mathbb{N}^{\ast},$ $f_{k}\left(X_{T_{k}}\right)$
is $\mathscr{A}_{T_{k}}-$measurable. Indeed, for every Borel subset
$B$ of $\mathbb{R},$ and every $n\in\mathbb{N}^{\ast},$
\[
\left(T_{k}=n\right)\cap\left[f_{k}\left(X_{T_{k}}\right)\right]^{-1}\left(B\right)=\left(T_{k}=n\right)\cap\left[f_{k}\left(X_{n}\right)\right]^{-1}\left(B\right)\in\mathscr{A}_{n},
\]
since $\left(T_{k}=n\right)\in\mathscr{A}_{n}$ and $\left[f_{k}\left(X_{n}\right)\right]^{-1}\left(B\right)\in\mathscr{A}_{n}.$
Taking into account the inclusion $\refpar{eq:A_T_k_subset_A_t_k_plus_1},$
the product $\prod^{k-1}_{j=1}f_{j}\left(X_{T_{j}}\right)$ is then
$\mathscr{A}_{T_{k-1}}-$measurable. Then,
\begin{align*}
\mathbb{E}\left(\prod^{k}_{j=1}f_{j}\left(X_{T_{j}}\right)\right) & =\mathbb{E}\left(\mathbb{E}^{\mathscr{A}_{T_{k-1}}}\left(\prod^{k}_{j=1}f_{j}\left(X_{T_{j}}\right)\right)\right)\\
 & =\mathbb{E}\left(\prod^{k-1}_{j=1}f_{j}\left(X_{T_{j}}\right)\mathbb{E}^{\mathscr{A}_{T_{k-1}}}\left(f_{k}\left(X_{T_{k}}\right)\right)\right)\\
 & =\mathbb{E}\left(\prod^{k-1}_{j=1}f_{j}\left(X_{T_{j}}\right)\right)\overline{t}I\left(f_{k}\right),
\end{align*}
Hence,
\[
\mathbb{E}\left(\prod^{k}_{j=1}f_{j}\left(X_{T_{j}}\right)\right)=\mathbb{E}\left(\prod^{k-1}_{j=1}f_{j}\left(X_{T_{j}}\right)\right)\mathbb{E}\left(f_{k}\left(X_{T_{k}}\right)\right),
\]
and, by backward iteration,\boxeq{
\[
\mathbb{E}\left(\prod^{k}_{j=1}f_{j}\left(X_{T_{j}}\right)\right)=\prod^{k}_{j=1}\mathbb{E}\left(f_{j}\left(X_{T_{j}}\right)\right).
\]
}Since this holds for every bounded nonnegative measurable function
$f_{k},$ \textbf{the sequence $\left(X_{T_{j}}\right)$ is a sequence
of independent random variables with density $f.$ }

\textbf{(d) Numerical application}

We have
\[
\forall x\in\mathbb{R}^{+},\,\,\,\,t\left(x\right)=\dfrac{ab^{a}}{\Gamma\left(a\right)}x^{a-1}\text{e}^{-x\left(b-\frac{1}{a}\right)}.
\]
Let $h=\ln\circ t.$
\[
\forall x\in\mathbb{R}^{+\ast},\,\,\,\,h^{\prime}\left(x\right)=\dfrac{a-1}{x}-\left(b-\dfrac{1}{a}\right)\,\,\,\,\text{and}\,\,\,\,h^{\prime\prime}\left(x\right)=\dfrac{1-a}{x^{2}}<0.
\]
Moreover,
\[
h^{\prime}\left(0^{+}\right)=+\infty\,\,\,\,\text{and}\,\,\,\,\lim_{x\to+\infty}h'\left(x\right)=-\dfrac{ab-1}{a}<0.
\]
It follows that $h,$ and thus $t,$ admits a unique maximum on $\mathbb{R}^{+}$
in $\widehat{x},$ solution of the equation $h^{\prime}\left(x\right)=0,$
namely
\[
\widehat{x}=\dfrac{a\left(a-1\right)}{ab-1},
\]
 and therefore\boxeq{
\[
\overline{t}=t\left(\widehat{x}\right)=\dfrac{\left(ab\right)^{a}}{\Gamma\left(a\right)}\left(\dfrac{a-1}{ab-1}\right)^{a-1}\text{e}^{-a+1}.
\]
}
\begin{itemize}
\item If $b=1$ and $a=\dfrac{5}{2},$ then
\[
\overline{t}=\dfrac{\left(\dfrac{5}{2}\right)^{\frac{5}{2}}}{\dfrac{3}{2}\dfrac{1}{2}\sqrt{\pi}}\text{e}^{-\frac{3}{2}}\approx1.66,
\]
which yields\boxeq{
\[
P\left(M_{n}\in G\right)\approx0.6\,\,\,\,\text{and}\,\,\,\,\mathbb{E}\left(T\right)\approx1.66.
\]
}
\item If $b=1$ and $a=\dfrac{9}{2},$ a similar computation yields\boxeq{
\[
P\left(M_{n}\in G\right)\approx0.44\,\,\,\,\text{and}\,\,\,\,\mathbb{E}\left(T\right)\approx2.26.
\]
}
\end{itemize}
\end{solution}

\chapter{Fourier Transforms and Characteristic Functions}\label{chap:PartIIChap13}

\begin{objective}{}{}

Chapter \ref{chap:PartIIChap13} deals with Fourier transforms and
characteristic functions.
\begin{itemize}
\item Section \ref{sec:Definition-and-Immediate} introduces the Fourier
transform of a bounded measure on $\mathbb{R}^{d}.$ First properties
are then stated.
\item Section \ref{sec:Injectivity-Theorem} starts by introducing the convolution
of a Borel function and a bounded measure on $\mathbb{R}^{d}.$ The
Gaussian kernel is then defined and properties are given, as well
as its convolution with a bounded measure. This then allows to state
the Fourier transform injectivity theorem, a fundamental result first
proved by Paul Lévy.
\item Section \ref{sec:Properties-Relative-to} concerns properties related
to independence. The Fourier transform of a product measure is given,
and an independence criterion via the characteristic functions is
deduced, together with a second independence criterion for characteristic
functions. The Fourier transform of a convolution product of two bounded
measures is also given, and application to characteristic functions
is deduced.
\item Section \ref{sec:Caracteristic-Function-and} focuses on the differentiability
of characteristic functions and how it relates to the existence of
moments of random variables. The case of real-valued random variables
is first examined before extending to $\mathbb{R}^{d}.$ The section
ends with the Taylor expansion at zero and the power-series expansion
of a characteristic function.
\end{itemize}
\end{objective}

The Fourier transform associates a function to every bounded measure
defined on $\mathbb{R}^{d}.$ Since it operates on the set of bounded
measures defined on $\mathbb{R}^{d},$ this transformation is injective.
Consequently, it allows, without loss of information, to replace the
study of a family of measures by the study of the family of associated
functions. More precisely, the strength of the Fourier transform comes
from the fact that it turns a convolution product of measures into
a pointwise product of functions, and that properties of convergence
of measures translate into corresponding convergence properties of
their Fourier transforms.

\section{Definition and Immediate Properties}\label{sec:Definition-and-Immediate}

Unless otherwise stated, in this section, \textbf{$\mu$ denotes a
bounded measure on $\mathbb{R}^{d},$ equipped with its Borel $\sigma-$algebra,
and $X$ denotes a random variable defined on the probabilized space
$\left(\Omega,\mathscr{A},P\right)$ with values in $\mathbb{R}^{d}.$}

We denote by $\left\langle \cdot,\cdot\right\rangle $ the usual Euclidean
inner product on $\mathbb{R}^{d}.$

Since, for every $t\in\mathbb{R}^{d},$ $\left|\text{e}^{\text{i}\left\langle x,t\right\rangle }\right|=1$
and since $\mu$ is a bounded measure, the function $x\mapsto\text{e}^{\text{i}\left\langle x,t\right\rangle }$
is $\mu-$integrable.

\begin{definition}{Fourier Transform of a Bounded Measure}{}

The function $\widehat{\mu}$ from $\mathbb{\mathbb{R}}^{d}$ into
$\mathbb{C},$ defined by\\
\boxeq{
\begin{equation}
\forall t\in\mathbb{R}^{d},\,\,\,\,\text{\ensuremath{\widehat{\mu}\left(t\right)=}\ensuremath{\intop_{\mathbb{R}^{d}}}}\text{e}^{\text{i}\left\langle x,t\right\rangle }\text{d}\mu\left(x\right)\label{eq:Fourier_transform_definition}
\end{equation}
}is called the \textbf{Fourier transform\index{Fourier transform}}
of the bounded measure $\mu.$

The Fourier transform of the law $P_{X}$ of a random variable $X$
is called the \textbf{\index{characteristic function}characteristic
function}. It is denoted by $\varphi_{X}.$

\end{definition}

\begin{remark}{}{}

It is worth noting that the notion of characteristic function depends
only on the law of the random variable $X,$ and not on the function
$X$ itself. 

\end{remark}

The following fundamental formula follows directly from the transfer
theorem,\boxeq{
\begin{equation}
\forall t\in\mathbb{R}^{d},\,\,\,\,\varphi_{X}\left(t\right)=\mathbb{E}\left(\text{e}^{\text{i}\left\langle X,t\right\rangle }\right).\label{eq:charact_function_and_exp}
\end{equation}
}

The notions of Fourier transform and characteristic function extend
immediately, without any change in the formulae, to the case where
$\mathbb{R}^{d}$ is replaced by a finite dimensional vector space
$E.$ The bracket corresponds then to the duality bilinear form\footnote{Recall that if $E$ is equipped with an Euclidean structure, $E^{\ast}$
can be naturally identified to $E.$} between $E$ and its dual $E^{\ast}.$ In this setting, $\widehat{\mu}$
is defined as a function on $E^{\ast}.$ If $\phi$ is an isorphism
from $E$ to $\mathbb{R}^{d},$ and if $\phi\left(\mu\right)$ denotes
the measure image of $\mu$ under $\phi,$ an immediate computation
shows that 
\[
\widehat{\phi\left(\mu\right)}=\widehat{\mu}\circ\phi^{\ast},
\]
where $\phi^{\ast}$ denotes the adjoint of $\phi,$ defined on $\mathbb{R}^{d}$
and taking values on $E^{\ast}.$ 

If $X$ is a random variable taking values in $E,$ its characteristic
function is then defined, by the same formula $\refpar{eq:charact_function_and_exp}$
as a function on the dual space $E^{\ast}.$ All properties established
in the case of $\mathbb{R}^{d}$ extend to the general case. This
extension may be treated as an exercise.

We now present, in parallel, the basic properties of $\widehat{\mu}$
and $\varphi_{X}.$

\begin{proposition}{Basic Properties of the Fourier Transform and the Characteristic Function}{}

With the above notations, the following properties hold:

1. $\widehat{\mu}\left(0\right)=\mu\left(\mathbb{R}^{d}\right)$ and
$\varphi_{X}\left(0\right)=1.$

2. $\forall t\in\mathbb{\mathbb{R}}^{d},\,\,\,\,\left|\widehat{\mu}\left(t\right)\right|\leqslant\mu\left(\mathbb{R}^{d}\right)\,\,\,\,\text{and}\,\,\,\,\left|\varphi_{X}\left(t\right)\right|\leqslant1,$

3. $\forall t\in\mathbb{\mathbb{R}}^{d},\,\,\,\,\widehat{\mu}\left(-t\right)=\overline{\widehat{\mu}\left(t\right)}\,\,\,\,\text{and}\,\,\,\,\varphi_{X}\left(-t\right)=\overline{\varphi_{X}\left(t\right)}.$

4. Let $A\in\mathscr{L}\left(\mathbb{R}^{d},\mathbb{R}^{k}\right)$
and $b\in\mathbb{R}^{k}.$ It holds\boxeq{
\begin{equation}
\forall t\in\mathbb{R}^{k},\,\,\,\,\varphi_{AX+b}\left(t\right)=\varphi_{X}\left(A^{\ast}t\right)\text{e}^{\text{i}\left\langle b,t\right\rangle }\label{eq:charact_affine_translation_rv}
\end{equation}
}where $A^{\ast}$ denotes the adjoint of $A.$

5. The functions $\widehat{\mu}$ and $\varphi_{X}$ are uniformly
continuous on $\mathbb{R}^{d}.$

\end{proposition}

\begin{proof}{}{}

1.2.3.The first three properties are straightforward and left to the
reader.

4. For the fourth property, for every $t\in\mathbb{R}^{k},$
\[
\varphi_{AX+b}\left(t\right)=\mathbb{E}\left(\text{e}^{\text{i}\left\langle AX+b,t\right\rangle }\right).
\]
By definition of the adjoint of $A,$
\[
\left\langle AX+b,t\right\rangle =\left\langle X,A^{\ast}t\right\rangle +\left\langle b,t\right\rangle .
\]
Therefore,
\[
\varphi_{AX+b}\left(t\right)=\text{e}^{\text{i}\left\langle b,t\right\rangle }\mathbb{E}\left(\text{e}^{\text{i}\left\langle X,A^{\ast}t\right\rangle }\right)=\varphi_{X}\left(A^{\ast}t\right)\text{e}^{\text{i}\left\langle b,t\right\rangle }.
\]

5. Now, we prove that $\widehat{\mu}$ is uniformly continuous. Let
$\epsilon>0$ be fixed. Since $\mu$ is a bounded measure, there exists
an integer $n$ such that 
\[
\mu\left(B\left(0,n\right)^{c}\right)\leqslant\dfrac{\epsilon}{4},
\]
where $B\left(0,n\right)$ denotes the open ball centered at $0$
and radius $n.$

For every $u,t\in\mathbb{R}^{d},$
\[
\left|\widehat{\mu}\left(u\right)-\widehat{\mu}\left(t\right)\right|\leqslant\intop_{B\left(0,n\right)}\left|\text{e}^{\text{i}\left\langle x,u\right\rangle }-\text{e}^{\text{i}\left\langle x,t\right\rangle }\right|\text{d}\mu\left(x\right)+2\mu\left(B\left(0,n\right)^{c}\right).
\]
By the mean value inequality applied to the complex exponential,
\[
\left|\text{e}^{\text{i}\left\langle x,u\right\rangle }-\text{e}^{\text{i}\left\langle x,t\right\rangle }\right|\leqslant\left\Vert u-t\right\Vert \left\Vert x\right\Vert .
\]
Hence,
\[
\left|\widehat{\mu}\left(u\right)-\widehat{\mu}\left(t\right)\right|\leqslant n\mu\left(\mathbb{R}^{d}\right)\left\Vert u-t\right\Vert +2\mu\left(B\left(0,n\right)^{c}\right).
\]

Let $\eta=\dfrac{\epsilon}{2n\mu\left(\mathbb{R}^{d}\right)},$ then
\[
\forall u,t\in\mathbb{R}^{d},\,\,\,\,\left\Vert u-t\right\Vert \leqslant\eta\Longrightarrow\left|\widehat{\mu}\left(u\right)-\widehat{\mu}\left(t\right)\right|\leqslant\epsilon.
\]
This proves the result, since $\epsilon$ is arbitrary.

\end{proof}

\section{Injectivity Theorem}\label{sec:Injectivity-Theorem}

\begin{denotation}{}{}

Throughout this chapter, we denote $\intop_{\mathbb{R}^{d}}g\left(x\right)\text{d}x$
the integral of a function $g$ over $\mathbb{R}^{d}.$

\end{denotation}

\begin{definition}{Convolution of a Borel Function with a Bounded Measure on $\mathbb{R}^d$. Fourier Transform of a Lebesgue-Integrable Function}{}

Let $\mu$ be a bounded measure on $\mathbb{R}^{d},$ and let $f$
be a Borel function such that, for every $x,$ the function $y\mapsto f\left(x-y\right)$
is $\mu-$integrable. The \textbf{convolution\index{convolution}}
of $f$ and $\mu$ is the function denoted $f\ast\mu$ defined by\\
\boxeq{
\[
\forall x\in\mathbb{R}^{d},\,\,\,\,\left(f\ast\mu\right)\left(x\right)=\intop_{\mathbb{R}^{d}}f\left(x-y\right)\text{d}\mu\left(y\right).
\]
}

If $g$ is Lebesgue-integrable, we denote $\widehat{g}$ its \textbf{Fourier
transform}\index{Fourier transform}, that is the function defined
on $\mathbb{R}^{d}$ by\\
\boxeq{
\[
\forall t\in\mathbb{R}^{d},\,\,\,\,\widehat{g}\left(t\right)=\intop_{\mathbb{R}^{d}}g\left(x\right)\text{e}^{\text{i}\left\langle x,t\right\rangle }\text{d}x.
\]
}

\end{definition}

For our purposes, \textbf{the fundamental property of the Fourier
transform is its injectivity: a bounded measure $\mu$ on $\mathbb{R}^{d}$
is uniquely determined by its Fourier transform $\widehat{\mu}.$}
This property is established through a sequence of lemmas.

\begin{lemma}{Gaussian Kernel Property}{}

For every $\sigma>0,$ consider the function $g_{\sigma},$ called
the \textbf{\mindex{Gaussian!kernel}\mindex{kernel!Gaussian}Gaussian
kernel}, defined on $\mathbb{R}^{d}$ by
\[
\forall x\in\mathbb{R}^{d},\,\,\,\,g_{\sigma}\left(x\right)=\dfrac{1}{\left(\sigma\sqrt{2\pi}\right)^{d}}\text{e}^{-\frac{\left\Vert x\right\Vert ^{2}}{2\sigma^{2}}},
\]
where $\left\Vert \cdot\right\Vert $ denotes the Euclidean norm on
$\mathbb{R}^{d}.$

(a) The function $g_{\sigma}$ is a probability density on $\mathbb{R}^{d}.$

(b) For every $\epsilon>0,$
\[
\lim_{\sigma\to0}\intop_{\left\{ \left\Vert x\right\Vert >\epsilon\right\} }g_{\sigma}\left(x\right)\text{d}x=0.
\]

(c) For every $f\in\mathscr{C}_{b}\left(\mathbb{R}^{d}\right)$ and
for every $x\in\mathbb{R}^{d},$
\begin{equation}
\left(f\ast g_{\sigma}\right)\left(x\right)\underset{\sigma\to0}{\longrightarrow}f\left(x\right).\label{eq:lim_sig_0_conv_f_gaussian_kernel}
\end{equation}

(d) The Fourier transform of $g_{1}$ is given by
\begin{equation}
\forall t\in\mathbb{R}^{d},\,\,\,\,\widehat{g}_{1}\left(t\right)=\text{e}^{-\frac{\left\Vert t\right\Vert ^{2}}{2}}\equiv\left(\sqrt{2\pi}\right)^{d}g_{1}\left(t\right).\label{eq:fourier_transform_g_1}
\end{equation}

\end{lemma}

\begin{proof}{}{}

(a) We begin with the case $d=1.$ The change of variables $x=\sigma y$
yields
\[
\intop_{\mathbb{R}}\dfrac{1}{\sigma\sqrt{2\pi}}\text{e}^{-\frac{x^{2}}{2\sigma^{2}}}\text{d}x=\intop_{\mathbb{R}}\dfrac{1}{\sqrt{2\pi}}\text{e}^{-\frac{y^{2}}{2}}\text{d}y=1,
\]
where the last equality was established in Example $\ref{ex:int_exp_min_xsqu}.$ 

In the general case, the Fubini theorem gives
\[
\intop_{\mathbb{R}^{d}}g_{\sigma}\left(x\right)\text{d}x=\prod^{d}_{j=1}\left(\intop_{\mathbb{R}}\dfrac{1}{\sigma\sqrt{2\pi}}\text{e}^{-\frac{x^{2}_{j}}{2\sigma^{2}}}\text{d}x_{j}\right)=1.
\]

(b) Temporarily, denote by $\left\Vert \,\cdot\,\right\Vert _{2}$
the Euclidean norm of $\mathbb{R}^{d}$ and by $\left\Vert \,\cdot\,\right\Vert _{\infty}$
the max norm. There exists a constant $c>0,$ such that, for every
$x\in\mathbb{R}^{d},$ $c\left\Vert x\right\Vert _{2}\leqslant\left\Vert x\right\Vert _{\infty}.$
Consequently, for every $\epsilon>0,$
\[
\left\{ \left\Vert x\right\Vert _{\infty}\leqslant c\epsilon\right\} \subset\left\{ \left\Vert x\right\Vert _{2}\leqslant\epsilon\right\} .
\]
Therefore,
\[
\intop_{\left\{ \left\Vert x\right\Vert _{\infty}\leqslant c\epsilon\right\} }g_{\sigma}\left(x\right)\text{d}x\leqslant\intop_{\left\{ \left\Vert x\right\Vert _{2}\leqslant\epsilon\right\} }g_{\sigma}\left(x\right)\text{d}x\leqslant1.
\]

By the Fubini theorem
\begin{align*}
\intop_{\left\{ \left\Vert x\right\Vert _{\infty}\leqslant c\epsilon\right\} }g_{\sigma}\left(x\right)\text{d}x & =\prod^{d}_{j=1}\left(\intop_{\left\{ \left|x_{j}\right|\leqslant c\epsilon\right\} }\dfrac{1}{\sigma\sqrt{2\pi}}\text{e}^{-\frac{x^{2}_{j}}{2\sigma^{2}}}\text{d}x_{j}\right)\\
 & =\left(\intop_{\left\{ \left|y\right|\leqslant\frac{c\epsilon}{\sigma}\right\} }\dfrac{1}{\sqrt{2\pi}}\text{e}^{-\frac{y^{2}}{2}}\text{d}y\right)^{d},
\end{align*}
Letting $\sigma\to0$ shows that
\[
\intop_{\left\{ \left\Vert x\right\Vert _{\infty}\leqslant c\epsilon\right\} }g_{\sigma}\left(x\right)\text{d}x\underset{\sigma\to0}{\longrightarrow}1,
\]
which proves the claim.

\begin{remark}{}{}

This result may also be obtained by switching to spherical coordinates.
Indeed, a direct computation shows that there exists a constant $c_{d}>0$
such that
\[
\intop_{\left\{ \left\Vert x\right\Vert _{2}>\epsilon\right\} }g_{\sigma}\left(x\right)\text{d}x=c_{d}\intop^{+\infty}_{\frac{\epsilon}{\sigma}}\text{e}^{-\frac{r^{2}}{2}}r^{d-1}\text{d}r.
\]

\end{remark}

(c) Let $x\in\mathbb{R}^{d}.$ The change of variables $z=\dfrac{x-y}{\sigma},$
whose Jacobian is $\sigma^{d},$ yields
\begin{align}
\left(f\ast g_{\sigma}\right)\left(x\right) & =\intop_{\mathbb{R}^{d}}f\left(y\right)\dfrac{1}{\left(\sigma\sqrt{2\pi}\right)^{d}}\text{e}^{-\frac{\left\Vert x-y\right\Vert ^{2}}{2\sigma^{2}}}\text{d}y\label{eq:conv_prod_f_g_sigma}\\
 & =\intop_{\mathbb{R}^{d}}f\left(x-\sigma z\right)g_{1}\left(z\right)\text{d}z.\nonumber 
\end{align}

Since $f$ is continuous and bounded, 
\[
\lim_{\sigma\to0}f\left(x-\sigma z\right)=f\left(x\right),
\]
and 
\[
\left|f\left(x-\sigma z\right)\right|\leqslant\left\Vert f\right\Vert _{\infty},
\]
which is an integrable constant with respect to the density probability
$g_{1}.$ The dominated convergence theorem, applied to an arbitrary
nonnegative sequence that converges to 0, yields the result.

(d) Let $t\in\mathbb{R}^{d}.$ Since the function $x\mapsto\text{e}^{\text{i}\left\langle x,t\right\rangle }$
is integrable with respect to the probability measure with density
$g_{1}$, the Fubini theorem applies and
\[
\widehat{g_{1}}\left(t\right)=\prod^{d}_{j=1}\left(\intop_{\mathbb{R}}\dfrac{1}{\sqrt{2\pi}}\text{e}^{\text{i}x_{j}t_{j}-\frac{x^{2}_{j}}{2}}\text{d}x_{j}\right).
\]
Thus,
\[
\widehat{g_{1}}\left(t\right)=\text{e}^{-\frac{\left\Vert t\right\Vert ^{2}}{2}},
\]
provided that\\
\boxeq{
\[
\forall u\in\mathbb{R},\,\,\,\,\intop_{\mathbb{R}}\dfrac{1}{\sqrt{2\pi}}\text{e}^{\text{i}xu-\frac{x^{2}}{2}}\text{d}x=\text{e}^{-\frac{u^{2}}{2}}.
\]
}

We now establish this identity, which will also be useful in other
applications.

It holds that
\[
\forall z\in\mathbb{R},\,\,\,\,\intop_{\mathbb{R}}\dfrac{1}{\sqrt{2\pi}}\text{e}^{-\frac{\left(x-z\right)^{2}}{2}}\text{d}x=1,
\]
which, after expanding the square, implies
\begin{equation}
\forall z\in\mathbb{R},\,\,\,\,\intop_{\mathbb{R}}\dfrac{1}{\sqrt{2\pi}}\text{e}^{zx-\frac{x^{2}}{2}}\text{d}x=\text{e}^{\frac{z^{2}}{2}}.\label{eq:int_on_r_exp_zxminus_xsquare2}
\end{equation}

Moreover, for every $z\in\mathbb{C},$
\[
\intop_{\mathbb{R}}\dfrac{1}{\sqrt{2\pi}}\text{e}^{\left|zx\right|-\frac{x^{2}}{2}}\text{d}x<+\infty.
\]

This shows both the integrability of the function $x\mapsto\text{e}^{zx-\frac{x^{2}}{2}},$
since
\[
\left|\text{e}^{zx-\frac{x^{2}}{2}}\right|\leqslant\text{e}^{\left|zx\right|-\frac{x^{2}}{2}},
\]
and the following inequality, by monotonic convergence,
\[
\sum^{+\infty}_{n=0}\intop_{\mathbb{R}}\dfrac{1}{\sqrt{2\pi}}\dfrac{\left|zx\right|^{n}}{n!}\text{e}^{-\frac{x^{2}}{2}}\text{d}x<+\infty.
\]

By the corollary of the dominated convergence related to series theorem---see
Corollary $\ref{co:sum_measurable_functions_with_bounded_abs_value}$,---it
follows that, for every $z\in\mathbb{C},$ 
\[
\intop_{\mathbb{R}}\dfrac{1}{\sqrt{2\pi}}\text{e}^{zx-\frac{x^{2}}{2}}\text{d}x=\sum^{+\infty}_{n=0}\dfrac{z^{n}}{n!}\left(\intop_{\mathbb{R}}\dfrac{1}{\sqrt{2\pi}}x^{n}\text{e}^{-\frac{x^{2}}{2}}\text{d}x\right).
\]

Hence, the function 
\[
z\mapsto\intop_{\mathbb{R}}\dfrac{1}{\sqrt{2\pi}}\text{e}^{zx-\frac{x^{2}}{2}}\text{d}x
\]
is thus entire and coincides on $\mathbb{R}$ with the entire function
$z\mapsto\text{e}^{\frac{z^{2}}{2}}.$ 

By the analytical expansion principle, these functions coincide on
$\mathbb{C}.$ Substituting $z=\text{i}u$ in the equality $\refpar{eq:int_on_r_exp_zxminus_xsquare2}.$
yields the desired result.

\end{proof}

The following lemma is the key ingredient in the proof of the injectivity
theorem. It shows that knowledge of the Fourier transform $\widehat{\mu}$
uniquely determines the convolution product $g_{\sigma}\ast\mu$ $\left(\sigma>0\right).$
The proof of the injectivity theorem---Theorem $\ref{th:fourier_transf_inject}$---then
consists in showing that, once the family $\left(g_{\sigma}\ast\mu\right)_{\sigma>0}$
is given, the measure $\mu$ is uniquely determined.

\begin{lemma}{Gaussian Kernel and Bounded Measure: Convolution and Fourier Transform}{}

Let $\mu$ be a bounded measure on $\mathbb{R}^{d}.$ For every $\sigma>0$
and every $y\in\mathbb{R}^{d},$ the function $g_{\sigma}\left(\cdot-y\right)$
is $\mu-$integrable, and
\begin{equation}
\left(g_{\sigma}\ast\mu\right)\left(y\right)=\left(\sqrt{2\pi}\right)^{-d}\intop_{\mathbb{R}^{d}}\widehat{\mu}\left(v\right)g_{1}\left(\sigma_{v}\right)\text{e}^{-\text{i}\left\langle y,v\right\rangle }\text{d}v.\label{eq:conv_gaussian_kernel_bounded_meas}
\end{equation}

\end{lemma}

\begin{proof}{}{}

Let $y\in\mathbb{R}^{d}.$ The function $g_{\sigma}\left(\cdot-y\right)$
is $\mu-$integrable since it is bounded and $\mu$ is a bounded measure.
Using the relation $\refpar{eq:fourier_transform_g_1}$ between $g_{1}$
and its Fourier transform $\widehat{g}_{1},$ we may write
\[
g_{\sigma}\left(y-x\right)=g_{\sigma}\left(x-y\right)=\dfrac{1}{\left(\sigma\sqrt{2\pi}\right)^{d}}\widehat{g}_{1}\left(\dfrac{x-y}{\sigma}\right).
\]
By the change of variable $v=\dfrac{z}{\sigma}$ whose Jacobian is
$\dfrac{1}{\sigma^{d}},$ this becomes
\[
g_{\sigma}\left(y-x\right)=\dfrac{1}{\left(\sqrt{2\pi}\right)^{d}}\intop_{\mathbb{R}^{d}}g_{1}\left(\sigma v\right)\text{e}^{\text{i}\left\langle x-y,v\right\rangle }\text{d}v.
\]
Consequently,
\[
\left(g_{\sigma}\ast\mu\right)\left(y\right)=\intop_{\mathbb{R}^{d}}\left[\dfrac{1}{\left(\sqrt{2\pi}\right)^{d}}\intop_{\mathbb{R}^{d}}g_{1}\left(\sigma v\right)\text{e}^{\text{i}\left\langle x-y,v\right\rangle }\text{d}v\right]\text{d}\mu\left(x\right).
\]

Since
\[
\left|g_{1}\left(\sigma v\right)\text{e}^{\text{i}\left\langle x-y,v\right\rangle }\right|\leqslant\text{e}^{-\sigma^{2}\frac{\left\Vert v\right\Vert ^{2}}{2}},
\]
and since the function determined on the right-hand side is $\mu\otimes\lambda_{d}-$integrable---once
again, the boundedness of $\mu$ is essential---the Fubini theorem
applies. We therefore obain
\[
\left(g_{\sigma}\ast\mu\right)\left(y\right)=\left(\sqrt{2\pi}\right)^{-d}\intop_{\mathbb{R}^{d}}\left[\intop_{\mathbb{R}^{d}}\text{e}^{\text{i}\left\langle x,v\right\rangle }\text{d}\mu\left(x\right)\right]g_{1}\left(\sigma v\right)\text{e}^{-\text{i}\left\langle y,v\right\rangle }\text{d}v,
\]
 which yields the stated result.

\end{proof}

\begin{theorem}{Fourier Transform Injectivity Theorem}{fourier_transf_inject}

Two bounded measures on $\mathbb{R}^{d}$ that have the same Fourier
transform are equal.

\end{theorem}

\begin{proof}{}{}

Recall that a bounded measure on $\mathbb{R}^{d}$ is uniquely determined,
for every $f\in\mathscr{C}_{b}\left(\mathbb{R}^{d}\right),$ when
the integrals $\intop_{\mathbb{R}^{d}}f\text{d}\mu$---refer to Chapter
\ref{chap:PartIIChap9}, Corollary $\ref{co:equality_radon_measures}.$ 

Let $f\in\mathscr{C}_{b}\left(\mathbb{R}^{d}\right)$ be arbitrary.
We show that its integral with respect to $\mu$ is a function of
$\widehat{\mu},$ which ensures the injectivity of the Fourier transform.

Let $\left(\sigma_{n}\right)_{n\in\mathbb{N}}$ be an arbitrary nonnegative
sequence converging to 0. By the relation $\refpar{eq:lim_sig_0_conv_f_gaussian_kernel},$
for every $x\in\mathbb{R}^{d},$
\begin{equation}
\lim_{n\to+\infty}\left(f\ast g_{\sigma_{n}}\right)\left(x\right)=f\left(x\right).\label{eq:lim_f_conv_g_sigma_n}
\end{equation}

Moreover, by $\refpar{eq:conv_prod_f_g_sigma},$ for every $n\in\mathbb{N},$
\begin{equation}
\left|\left(f\ast g_{\sigma_{n}}\right)\left(x\right)\right|\leqslant\left\Vert f\right\Vert _{\infty}.\label{eq:bound_f_conv_g_sigma_n}
\end{equation}

Since every constant function is $\mu-$integrable, the dominated
convergence theorem yields
\[
\intop_{\mathbb{R}^{d}}f\text{d}\mu=\intop_{\mathbb{R}^{d}}\lim_{n\to+\infty}\left(f\ast g_{\sigma_{n}}\right)\left(x\right)\text{d}\mu\left(x\right)=\lim_{n\to+\infty}\intop_{\mathbb{R}^{d}}\left(f\ast g_{\sigma_{n}}\right)\left(x\right)\text{d}\mu\left(x\right).
\]
Hence,
\begin{equation}
\intop_{\mathbb{R}^{d}}f\text{d}\mu=\lim_{n\to+\infty}\intop_{\mathbb{R}^{d}}\left[\intop_{\mathbb{R}^{d}}f\left(y\right)g_{\sigma_{n}}\left(x-y\right)\text{d}y\right]\text{d}\mu\left(x\right).\label{eq:int_on_r_d_of_f_dmu_as_lim_int_prod_fandkernel_sigma_n}
\end{equation}
In this expression, the order of integration may be exchanged since
the function $\left(x,y\right)\mapsto f\left(y\right)g_{\sigma_{n}}\left(x-y\right)$
is $\lambda_{d}\otimes\mu-$integrable. Indeed, applying the Fubini
theorem to the nonnegative function
\[
\intop_{\mathbb{R}^{d}\times\mathbb{R}^{d}}\left|f\left(y\right)g_{\sigma_{n}}\left(x-y\right)\right|\text{d}\lambda_{d}\otimes\mu\left(x,y\right)=\intop_{\mathbb{R}^{d}}\left|f\left(y\right)\right|\left[\intop_{\mathbb{R}^{d}}g_{\sigma_{n}}\left(x-y\right)\text{d}\lambda_{d}\left(x\right)\right]\text{d}\mu\left(y\right).
\]
It follows that
\[
\intop_{\mathbb{R}^{d}\times\mathbb{R}^{d}}\left|f\left(y\right)g_{\sigma_{n}}\left(x-y\right)\right|\text{d}\lambda_{d}\otimes\mu\left(x,y\right)\leqslant\intop_{\mathbb{R}^{d}}\left\Vert f\right\Vert _{\infty}\left[\intop_{\mathbb{R}^{d}}g_{\sigma_{n}}\left(x-y\right)\text{d}\lambda_{d}\left(x\right)\right]\text{d}\mu\left(y\right).
\]
Since, for every $y\in\mathbb{R}^{d},$
\[
\intop_{\mathbb{R}^{d}}g_{\sigma_{n}}\left(x-y\right)\text{d}\lambda_{d}\left(x\right)=\intop_{\mathbb{R}^{d}}g_{\sigma_{n}}\left(x\right)\text{d}\lambda_{d}\left(x\right)=1,
\]
we obtain
\[
\intop_{\mathbb{R}^{d}\times\mathbb{R}^{d}}\left|f\left(y\right)g_{\sigma_{n}}\left(x-y\right)\right|\text{d}\lambda_{d}\otimes\mu\left(x,y\right)\leqslant\left\Vert f\right\Vert _{\infty}\mu\left(\mathbb{R}^{d}\right)<+\infty.
\]

Hence, the Fubini theorem applies to the integrals appearing in the
equality $\refpar{eq:int_on_r_d_of_f_dmu_as_lim_int_prod_fandkernel_sigma_n}.$
Taking into account the parity of $g_{\sigma_{n}},$ we obtain
\[
\intop_{\mathbb{R}^{d}}f\text{d}\mu=\lim_{n\to+\infty}\intop_{\mathbb{R}^{d}}f\left(y\right)\left[\intop_{\mathbb{R}^{d}}g_{\sigma_{n}}\left(x-y\right)\text{d}\mu\left(x\right)\right]\text{d}y=\lim_{n\to+\infty}\intop_{\mathbb{R}^{d}}f\left(y\right)\left(g_{\sigma_{n}}\ast\mu\right)\left(y\right)\text{d}y.
\]
Finally, using equality $\refpar{eq:conv_gaussian_kernel_bounded_meas},$
we conclude that
\[
\intop_{\mathbb{R}^{d}}f\text{d}\mu=\left(\sqrt{2\pi}\right)^{-d}\lim_{n\to+\infty}\intop_{\mathbb{R}^{d}}f\left(y\right)\left[\intop_{\mathbb{R}^{d}}\widehat{\mu}\left(v\right)g_{1}\left(\sigma_{n}v\right)\text{e}^{-\text{i}\left\langle y,v\right\rangle }\text{d}v\right]\text{d}y.
\]

This proves, as stated, that $\intop_{\mathbb{R}^{d}}f\text{d}\mu$
is function of $\widehat{\mu},$ and the theorem follows.

\end{proof}

\begin{figure}[t]
\begin{center}\includegraphics[width=0.4\textwidth]{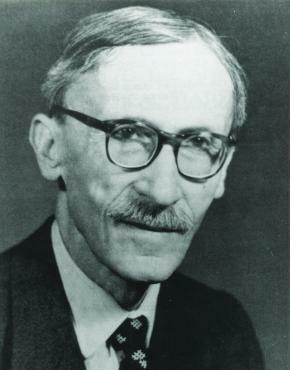}

{\tiny Credits: Konrad Jacobs (Erlangen), CC-BY-SA 2.0}\end{center}

\caption{\textbf{\protect\href{https://en.wikipedia.org/wiki/Paul_L\%25C3\%25A9vy_(mathematician)}{Paul Lévy}}
(1886 - 1971)}
\end{figure}

\begin{remark}{}{}

It follows from the injectivity property that the characteristic function
of a random variable taking values in $\mathbb{R}^{d}$ completely
characterizes the law of this random variable---hence its name. In
particular, the Table in Chapter \ref{chap:PartIIChap9} which presents
probability laws together with their Fourier transforms, can be read
in both directions. This property was first proved by \textbf{\sindex[fam]{Lévy, Paul}Paul
Lévy}\footnotemark, in 1922, for real-valued random variables, in
terms of their cumulative distribution functions.

\end{remark}

\footnotetext{\textbf{\href{https://en.wikipedia.org/wiki/Paul_L\%25C3\%25A9vy_(mathematician)}{Paul Lévy}}
(1886-1971), was born in Paris. He taught at École des mines in Saint-Etienne,
then at the École des Mines in Paris, and concurrently at the École
Polytechnique. His early work focused on functional analysis, but
he rapidly turned to probability theory. His contributions to probability
theory are fundamental, in particular in the study of stochastic processes
and Brownian motion. His works were collected in a series of three
volumes published by Gauthier-Villars in 1976 and 1980.}

To refine the injectivity theorem, there exists an explicit formula
expressing the cumulative function of a bounded measure on $\mathbb{R}^{d}$
in terms of its Fourier transform---see Exercise $\ref{exo:exercise13.10}$
of this chapter. Here, we restrict ourselves to the case where the
measure is \textbf{absolutely continuous}.

\begin{proposition}{Fourier Transform in the Case The Measure is Absolutely Continuous}{abs_cont_meas_four_tr}

Let $\mu$ be a bounded measure on $\mathbb{R}^{d}$ such that its
Fourier transform $\widehat{\mu}$ is Lebesgue-integrable. Then $\mu$
is absolutely continuous with respect to the Lebesgue measure, and
its density is given $\lambda_{d}-$almost surely by the continuous
function $h$ defined by\boxeq{
\[
\forall x\in\mathbb{R}^{d},\,\,\,\,h\left(x\right)=\dfrac{1}{\left(2\pi\right)^{d}}\intop_{\mathbb{R}^{d}}\widehat{\mu}\left(t\right)\text{e}^{-\text{i}\left\langle x,t\right\rangle }\text{d}t.
\]
}

\end{proposition}

\begin{proof}{}{}

To identify the measure $\mu,$ it is sufficient to compute the integrals
$\intop_{\mathbb{R}^{d}}f\text{d}\mu$ for every $f\in\mathscr{C}_{\mathscr{K}}\left(\mathbb{R}^{d}\right).$
Returning to the end of the proof of the injectivity theorem, we recall
that for every $f\in\mathscr{C}_{\mathscr{K}}\left(\mathbb{R}^{d}\right)$
and for every nonnegative sequence $\left(\sigma_{n}\right)_{n\in\mathbb{N}}$
converging to 0,
\[
\intop_{\mathbb{R}^{d}}f\text{d}\mu=\lim_{n\to+\infty}\intop_{\mathbb{R}^{d}}f\left(y\right)\left(g_{\sigma_{n}}\ast\mu\right)\left(y\right)\text{d}y.
\]

Moreover, by $\refpar{eq:conv_gaussian_kernel_bounded_meas},$
\[
\left|\left(g_{\sigma_{n}}\ast\mu\right)\left(y\right)\right|\leqslant\intop_{\mathbb{R}^{d}}\left|\widehat{\mu}\left(v\right)\right|g_{1}\left(\sigma_{n}v\right)\text{d}v\leqslant\intop_{\mathbb{R}^{d}}\left|\widehat{\mu}\left(v\right)\right|\text{d}v<+\infty,
\]
and therefore
\[
\left|f\left(y\right)\left(g_{\sigma_{n}}\ast\mu\right)\left(y\right)\right|\leqslant\left|f\left(y\right)\right|\left\Vert \widehat{\mu}\right\Vert _{L_{1}}.
\]

The right-hand side is Lebesgue-integrable, since $f$ is continuous
with compact support. Hence, by the dominated convergence theorem
\[
\intop_{\mathbb{R}^{d}}f\text{d}\mu=\intop_{\mathbb{R}^{d}}f\left(y\right)\left[\lim_{n\to+\infty}\left(g_{\sigma_{n}}\ast\mu\right)\left(y\right)\right]\text{d}y.
\]
This shows that $\mu$ is absolutely continuous with respect to the
Lebesgue measure and that its density is given $\lambda_{d}-$almost
surely by the continuous function $h$ defined as 
\[
h\left(y\right)=\lim_{n\to+\infty}\left(g_{\sigma_{n}}\ast\mu\right)\left(y\right).
\]
 The result then follows by noting that
\[
\lim_{n\to+\infty}\widehat{\mu}\left(v\right)g_{1}\left(\sigma_{n}v\right)\text{e}^{-i\left\langle y,v\right\rangle }=\left(\sqrt{2\pi}\right)^{-d}\widehat{\mu}\left(v\right)\text{e}^{-i\left\langle y,v\right\rangle }
\]
and by applying once again the dominated convergence theorem. This
is justified since $\widehat{\mu}$ is integrable and, by the relation
$\refpar{eq:conv_gaussian_kernel_bounded_meas},$
\[
\left|\widehat{\mu}\left(v\right)g_{1}\left(\sigma_{n}v\right)\text{e}^{-i\left\langle y,v\right\rangle }\right|\leqslant\left|\widehat{\mu}\left(v\right)\right|.
\]

Finally, the continuity of $h$ follows from the theorem on continuity
of integrals depending on a parameter, which is a corollary of the
dominated convergence theorem.

\end{proof}

\section{Properties Relative to Independence}\label{sec:Properties-Relative-to}

On $\mathbb{R}^{d_{1}}\times\mathbb{R}^{d_{2}}$ the usual Euclidean
scalar product satisfies\footnote{Tr.N: $\mathbb{R}^{d_{1}}\times\mathbb{R}^{d_{2}}$ is identified
to $\mathbb{R}^{d_{1}+d_{2}},$ where the usual Euclidean scalar product
for $x=\sum^{d_{1}+d_{2}}_{i=1}x_{i}e_{i}$ and $t=\sum^{d_{1}+d_{2}}_{i=1}t_{i}e_{i}$
is $\left\langle x,t\right\rangle =\sum^{d_{1}+d_{2}}_{i=1}x_{i}t_{i}=\sum^{d_{1}}_{i=1}x_{i}t_{i}+\sum^{d_{2}}_{j=1}x_{d_{1}+j}t_{d_{1}+j}.$}, for every $\left(x_{1},x_{2}\right)$ and $\left(t_{1},t_{2}\right)$
of $\mathbb{R}^{d_{1}}\times\mathbb{R}^{d_{2}},$
\[
\left\langle \left(x_{1},x_{2}\right),\left(t_{1},t_{2}\right)\right\rangle =\left\langle x_{1},t_{1}\right\rangle +\left\langle x_{2},t_{2}\right\rangle .
\]

Then, the next proposition is immediate.

\begin{proposition}{Fourier Transform of the Product Measure}{}

Let $\mu_{1}$ and $\mu_{2}$ be two bounded measures respectively
on $\mathbb{R}^{d_{1}}$ and $\mathbb{R}^{d_{2}}.$ The Fourier transform
of the product measure is the product of the Fourier transforms of
$\mu_{1}$ and $\mu_{2};$ that is\boxeq{
\[
\forall\left(t_{1},t_{2}\right)\in\mathbb{R}^{d_{1}}\times\mathbb{R}^{d_{2}},\,\,\,\,\widehat{\mu_{1}\otimes\mu_{2}}\left(t_{1},t_{2}\right)=\widehat{\mu_{1}}\left(t_{1}\right)\widehat{\mu_{2}}\left(t_{2}\right).
\]
}

\end{proposition}

\begin{proof}{}{}

The function $\left(x_{1},x_{2}\right)\mapsto\text{e}^{\text{i}\left\langle \left(x_{1},x_{2}\right),\left(t_{1},t_{2}\right)\right\rangle }$
is bounded and therefore $\mu_{1}\otimes\mu_{2}-$integrable. Applying
the Fubini theorem and using the identity
\[
\text{e}^{\text{i}\left\langle \left(x_{1},x_{2}\right),\left(t_{1},t_{2}\right)\right\rangle }=\text{e}^{\text{i}\left\langle x_{1},t_{1}\right\rangle }\text{e}^{\text{i}\left\langle x_{2},t_{2}\right\rangle },
\]
yields the desired result.

\end{proof}

We then obtain a\textbf{ criterion for independence} of random variables
expressed in terms of characteristic functions.

\begin{corollary}{Independence Criterion via Characteristic Functions}{indep_crit_char_function}

Let $X=\left(X_{1},X_{2}\right)$ be a random variable taking values
in $\mathbb{R}^{d_{1}}\times\mathbb{R}^{d_{2}}.$ The random variables
$X_{1}$ and $X_{2}$ are independent if and only if\boxeq{
\begin{equation}
\forall\left(t_{1},t_{2}\right)\in\mathbb{R}^{d_{1}}\times\mathbb{R}^{d_{2}},\,\,\,\,\varphi_{\left(X_{1},X_{2}\right)}\left(t_{1},t_{2}\right)=\varphi_{X_{1}}\left(t_{1}\right)\varphi_{X_{2}}\left(t_{2}\right).\label{eq:independence_via_charact_function}
\end{equation}
}

\end{corollary}

\begin{proof}{}{}

The random variables $X_{1}$ and $X_{2}$ are independent if and
only if $P_{\left(X_{1},X_{2}\right)}=P_{X_{1}}\otimes P_{X_{2}}.$
By the injectivity theorem, this is equivalent to $\widehat{P_{\left(X_{1},X_{2}\right)}}=\widehat{P_{X_{1}}\otimes P_{X_{2}}}.$
The result then follows directly from the previous proposition.

\end{proof}

\begin{remark}{}{}

The characteristic function of a marginal is immediate to obtain.
With the notations of Corollary $\ref{co:indep_crit_char_function},$\boxeq{
\begin{align}
\forall t_{1}\in\mathbb{R}^{d_{1}},\,\,\,\, & \varphi_{X_{1}}\left(t_{1}\right)=\varphi_{\left(X_{1},X_{2}\right)}\left(t_{1},0\right),\label{eq:marginal_charact_function}\\
\forall t_{2}\in\mathbb{R}^{d_{2}},\,\,\,\, & \varphi_{X_{2}}\left(t_{2}\right)=\varphi_{\left(X_{1},X_{2}\right)}\left(0,t_{2}\right).\label{eq:marginal_charact_function_2}
\end{align}
}

\end{remark}

Corollary $\ref{co:indep_crit_char_function}$ can also be stated
in the following form.

\begin{corollary}{Second Independence Criterion for Characteristic Functions}{}

Let $X=\left(X_{1},X_{2}\right)$ be a random variable taking values
in $\mathbb{R}^{d_{1}}\times\mathbb{R}^{d_{2}}.$ The random variables
$X_{1}$ and $X_{2}$ are independent, if and only if\boxeq{
\begin{equation}
\forall\left(t_{1},t_{2}\right)\in\mathbb{R}^{d_{1}}\times\mathbb{R}^{d_{2}},\,\,\,\,\varphi_{\left(X_{1},X_{2}\right)}\left(t_{1},t_{2}\right)=\varphi_{\left(X_{1},X_{2}\right)}\left(t_{1},0\right)\varphi_{\left(X_{1},X_{2}\right)}\left(0,t_{2}\right).\label{eq:independence_via_charact_function_marginals}
\end{equation}
}

\end{corollary}

\begin{example}{}{}

Let $X_{1}$ and $X_{2}$ be two inpendent real-valued random variables
with the same Laplace law, with characteristic function $\varphi$
given by
\[
\forall t\in\mathbb{R},\,\,\,\,\varphi\left(t\right)=\dfrac{1}{1+t^{2}}.
\]

Define the random variables $Y_{1}$ and $Y_{2}$ by
\[
Y_{1}=X_{1}-X_{2}\,\,\,\,\,\,\,\,\text{and}\,\,\,\,\,\,\,\,Y_{2}=X_{1}+X_{2}.
\]

Are the random variables $Y_{1}$ and $Y_{2}$ independent? Are they
uncorrelated?

\end{example}

\begin{solutionexample}{}{}

We have
\[
\left(\begin{array}{l}
Y_{1}\\
Y_{2}
\end{array}\right)=\left(\begin{array}{cc}
1 & -1\\
1 & 1
\end{array}\right)\left(\begin{array}{l}
X_{1}\\
X_{2}
\end{array}\right).
\]

By the equality $\refpar{eq:independence_via_charact_function},$
the characteristic function of $\left(X_{1},X_{2}\right)$ is defined
by
\[
\forall\left(t_{1},t_{2}\right)\in\mathbb{R}^{2},\,\,\,\,\varphi_{\left(X_{1},X_{2}\right)}\left(t_{1},t_{2}\right)=\dfrac{1}{1+t^{2}_{1}}\dfrac{1}{1+t^{2}_{2}}.
\]
Moreover, by the equality $\refpar{eq:charact_affine_translation_rv},$
the characteristic function of $\left(Y_{1},Y_{2}\right)$ is
\begin{align*}
\forall\left(t_{1},t_{2}\right)\in\mathbb{R}^{2},\,\,\,\,\varphi_{\left(Y_{1},Y_{2}\right)}\left(t_{1},t_{2}\right) & =\varphi_{\left(X_{1},X_{2}\right)}\left(t_{1}+t_{2},-t_{1}+t_{2}\right)\\
 & =\dfrac{1}{1+\left(t_{1}+t_{2}\right)^{2}}\dfrac{1}{1+\left(-t_{1}+t_{2}\right)^{2}}.
\end{align*}
We obtain the characteristic functions of the marginals $Y_{1}$ and
$Y_{2}$ by the relations $\refpar{eq:marginal_charact_function}$
and $\refpar{eq:marginal_charact_function_2}.$
\[
\forall t_{1}\in\mathbb{R},\,\,\,\,\varphi_{Y_{1}}\left(t_{1}\right)=\dfrac{1}{\left(1+t^{2}_{1}\right)^{2}}\,\,\,\,\text{and}\,\,\,\,\forall t_{2}\in\mathbb{R},\,\,\,\,\varphi_{Y_{2}}\left(t_{2}\right)=\dfrac{1}{\left(1+t^{2}_{2}\right)^{2}}.
\]
Thus, the random variables $Y_{1}$ and $Y_{2}$ have the same characteristic
function, and therefore the same law. They are not independent, since
$\varphi_{\left(Y_{1},Y_{2}\right)}\left(1,1\right)=\dfrac{1}{5}$
whereas $\varphi_{Y_{1}}\left(1\right)\varphi_{Y_{2}}\left(1\right)=\dfrac{1}{16},$
and these two quantities are not equal.

Nonetheless, we may note that these random variables are uncorrelated.
Indeed, since, $X_{1}$ and $X_{2}$ have the same law, they have
the same moments, and therefore
\[
\mathbb{E}\left(Y_{1}Y_{2}\right)=\mathbb{E}\left(X^{2}_{1}\right)-\mathbb{E}\left(X^{2}_{2}\right)=0\,\,\,\,\text{and}\,\,\,\,\text{\ensuremath{\mathbb{E}}\ensuremath{\left(Y_{1}\right)}}=\mathbb{E}\left(X_{1}\right)-\mathbb{E}\left(X_{2}\right)=0.
\]

\end{solutionexample}

One of the key point of this theory is that the Fourier transform
of the \textbf{convolution} of bounded measures, as well as the characteristic
function of \textbf{independent random variables sum,} can be computed
in a particularly simple way, as we now show.

\begin{proposition}{Fourier Transform of the Convolution Product of Two Bounded Measures}{}

Let $\mu_{1}$ and $\mu_{2}$ be two bounded measures on $\mathbb{R}^{d}.$
The Fourier transform of the convolution product of $\mu_{1}$ and
$\mu_{2}$ is the product of their Fourier transforms\boxeq{
\[
\widehat{\mu_{1}\ast\mu_{2}}=\widehat{\mu_{1}}\widehat{\mu_{2}}.
\]
}

\end{proposition}

\begin{proof}{}{}

Since the convolution $\mu_{1}\ast\mu_{2}$ is the measure image of
$\mu_{1}\otimes\mu_{2}$ by the sum function, and since the complex
exponential function is bounded, the transfer theorem yields
\[
\forall t\in\mathbb{R}^{d},\,\,\,\,\widehat{\mu_{1}\ast\mu_{2}}\left(t\right)=\intop_{\mathbb{R}^{d}}\text{e}^{\text{i}\left\langle x_{1}+x_{2},t\right\rangle }\text{d}\mu_{1}\otimes\mu_{2}\left(x_{1},x_{2}\right).
\]
The result then follows by factoring the exponential term and applying
the Fubini theorem.

\end{proof}

\begin{corollary}{}{}

Let $X_{1}$ and $X_{2}$ be two independent random variables taking
values in $\mathbb{R}^{d}.$ The characteristic function of their
sum satisfies‌\boxeq{
\[
\forall t\in\mathbb{R}^{d},\,\,\,\,\varphi_{X_{1}+X_{2}}\left(t\right)=\varphi_{X_{1}}\left(t\right)\varphi_{X_{2}}\left(t\right).
\]
}

\end{corollary}

\begin{proof}{}{}

It suffices to recall that, by independence, the law of $X_{1}+X_{2}$
is the convolution of the laws of $X_{1}$ and $X_{2}.$

\end{proof}

Taking into account the injectivity of the Fourier transform, we thus
obtain a \textbf{method to determine the law of a finite sum of independent
random variables}.

\begin{example}{Multinomial Law and Characteristic Function}{mult_law_char_fct}

We recall the modelling framework of the \textbf{\mindex{law!multinomial}\index{multinomial law}multinomial
law}.

Fix $k\in\mathbb{N^{\ast}}.$ For every $n\in\mathbb{N}^{\ast},$
consider a partition $\left(A^{n}_{j}\right)_{1\leqslant j\leqslant k}$
of $\Omega,$ where $A^{n}_{j}\in\mathscr{A}.$ Assume that the families
of events, indexed on $n,$ formed by the elements of these partitions
are independent. Moreover, suppose that
\[
\forall n\in\mathbb{N}^{\ast},\,\forall j\in\left\llbracket 1,k\right\rrbracket ,\,\,\,\,P\left(A^{n}_{j}\right)=p_{j},
\]
where $p_{j}>0$ and $\sum^{k}_{j=1}p_{j}=1.$ Define the random variables
$X^{n}$ and $Y^{n}$ taking values in $\mathbb{R}^{k}$ by
\[
X^{n}=\left(\begin{array}{c}
\boldsymbol{1}_{A^{n}_{1}}\\
\vdots\\
\boldsymbol{1}_{A^{n}_{k}}
\end{array}\right)\,\,\,\,\text{and}\,\,\,\,Y^{n}=\sum^{n}_{j=1}X^{j}.
\]

Determine the characteristic function of $Y^{n}.$

\end{example}

\begin{solutionexample}{}{}

The random variables $X^{j}$ are independent and of same law. Therefore,
the characteristic function of $Y^{n}$ is
\[
\forall t\in\mathbb{R}^{k},\,\,\,\,\varphi_{Y^{n}}\left(t\right)=\left[\varphi_{X^{1}}\left(t\right)\right]^{n}.
\]

Moreover, since $\left(A^{1}_{j}\right)_{1\leqslant j\leqslant k}$
is a partition of $\Omega,$
\[
\forall t\in\mathbb{R}^{k},\,\,\,\,\varphi_{X^{1}}\left(t\right)=\sum^{k}_{j=1}\intop_{A^{1}_{j}}\text{e}^{\text{i}\left\langle X^{1},t\right\rangle }\text{d}P=\sum^{k}_{j=1}p_{j}\text{e}^{\text{i}t_{j}},
\]
and hence,\boxeq{
\begin{equation}
\forall t\in\mathbb{R}^{k},\,\,\,\,\varphi_{X^{n}}\left(t\right)=\left[\sum^{k}_{j=1}p_{j}\text{e}^{\text{i}t_{j}}\right]^{n}.\label{eq:char_funct_Yn}
\end{equation}
}

\end{solutionexample}

\begin{remark}{}{}

It follows that if $Z_{1}$ and $Z_{2}$ are two independent $\mathbb{R}^{k}-$valued
random variables with respective multinomial laws $M\left(n_{1},p_{1},p_{2},\cdots,p_{k-1}\right)$
and $M\left(n_{2},p_{1},p_{2},\cdots,p_{k-1}\right),$ then the random
variable $Z_{1}+Z_{2}$ follows the multinomial law $M\left(n_{1}+n_{2},p_{1},p_{2},\cdots,p_{k-1}\right).$
In other words, \textbf{the family of multinomial laws $M\left(n,p_{1},p_{2},\cdots,p_{k-1}\right)$
is stable under convolution.}

\end{remark}

\section{Characteristic Functions and Moments}\label{sec:Caracteristic-Function-and}

\textbf{We study the relationship between the differentiabilty properties
of characteristic functions and the existence of random variable moments.}
We begin by recalling, without proof and in a vector-valued setting,
the theorem on differentiation of an integral with a parameter.

\begin{theorem}{Differentiation of An Integral With a  Parameter}{der_int_param}

Let $\mu$ be a $\sigma-$finite measure on a probabilizable space
$\left(\Omega,\mathscr{A}\right).$ Let $E$ and $F$ be two finite
dimensional normed vector spaces, and let $O$ be an open set of $E.$
Let $f$ be a function from $O\times\Omega$ into $F$ satisfying:

1. For every $\omega\in\Omega,$ the partial function $f\left(\cdot,\omega\right)$
is of class $C^{1}$ on $O,$ and there exists a function $g\in\mathscr{L}^{1}\left(\Omega,\mathscr{A},\mu\right)$
such that
\[
\forall x\in E,\,\,\,\,\left\Vert \dfrac{\partial}{\partial x}f\left(x,\cdot\right)\right\Vert _{\mathscr{L}\left(E,F\right)}\leqslant g.
\]

2. For every $x\in O,$ the partial function $f\left(x,\cdot\right)$
is $\mu-$integrable.

Then the function from $O$ to $F$ defined by 
\[
x\mapsto\intop_{\Omega}f\left(x,\omega\right)\text{d}\mu\left(\omega\right)
\]
is differentiable, and for every $x\in O,$\boxeq{
\[
\dfrac{\partial}{\partial x}\intop_{\Omega}f\left(x,\omega\right)\text{d}\mu\left(\omega\right)=\intop_{\Omega}\dfrac{\partial}{\partial x}f\left(x,\omega\right)\text{d}\mu\left(\omega\right).
\]
}

\end{theorem}

\begin{proof}{}{}

The proof follows directly from the finite increment inequality and
the dominated convergence theorem.

\end{proof}

To avoid the technicalities of multivariate differential calculus,
we first restrict our attention to the case of the real-valued random
variables.

\begin{proposition}{Moments and Differentiability of a Real-Valued Random Variable}{mom_der_real_rv}

Let $X$ be a real-valued random variable and let $\varphi_{X}$ denote
its characteristic function.

(a) If $X$ admits a moment of order $n\in\mathbb{N}^{\ast},$ then
$\varphi_{X}$ is of class $C^{n}$ and, for every $k$ such that
$1\leqslant k\leqslant n,$
\begin{equation}
\forall t\in\mathbb{R},\,\,\,\,\varphi^{\left(k\right)}_{X}\left(t\right)=\text{i}^{k}\intop_{\Omega}X^{k}\text{e}^{\text{i}tX}\text{d}P,\label{eq:k-th_deriv_charact_funct}
\end{equation}
and in particular, \boxeq{
\begin{equation}
\varphi^{\left(k\right)}_{X}\left(0\right)=\text{i}^{k}\mathbb{E}\left(X^{k}\right).\label{eq:k-th_der_char_funct_in_zero}
\end{equation}
}

(b) Conversely, if $\varphi_{X}$ is $k$ times differentiable at
0 $\left(k\geqslant2\right),$ then $X$ admits moments up to the
order $2\left\lfloor \dfrac{k}{2}\right\rfloor .$ These moments are
given by the formula $\refpar{eq:k-th_deriv_charact_funct}$

\end{proposition}

\begin{proof}{}{}

(a) Since
\[
\dfrac{\text{d}^{k}}{\text{d}t^{k}}\text{e}^{\text{i}tX}=\left(\text{i}X\right)^{k}\text{e}^{\text{i}tX},
\]
we have
\[
\left|\dfrac{\text{d}^{k}}{\text{d}t^{k}}\text{e}^{\text{i}tX}\right|\leqslant\left|X\right|^{k}.
\]
Hence, Theorem $\ref{th:der_int_param}$ can be applied $n$ times.

(b) Let us prove the result by induction on $k.$
\begin{itemize}
\item \textbf{Initialization step.}\\
Let us prove the result for $k=2.$ In this case, $\varphi_{X}$ admits
a second-order Taylor-Young expansion, and therefore
\[
\lim_{t\to0}\dfrac{\varphi_{X}\left(t\right)+\varphi_{X}\left(-t\right)-2}{t^{2}}=\varphi^{\prime\prime}_{X}\left(0\right).
\]
Moreover, since
\[
\varphi_{X}\left(t\right)+\varphi_{X}\left(-t\right)=2\text{Re}\left(\varphi_{X}\left(t\right)\right)=2\mathbb{E}\left(\cos\left(tX\right)\right),
\]
it follows that
\[
\lim_{t\to0}\mathbb{E}\left(\dfrac{1-\cos\left(tX\right)}{t^{2}}\right)=-\dfrac{1}{2}\varphi^{\prime\prime}_{X}\left(0\right).
\]
Since the integrand is nonnegative, the Fatou lemma implies that if
$\left(t_{n}\right)_{n\in\mathbb{N}}$ is a sequence converging to
0, then
\begin{align*}
\intop_{\Omega}X^{2}\text{d}P & =\mathbb{E}\left(2\liminf_{n\to+\infty}\dfrac{1-\cos\left(t_{n}X\right)}{t^{2}_{n}}\right)\\
 & \leqslant2\liminf_{n\to+\infty}\mathbb{E}\left(\dfrac{1-\cos\left(t_{n}X\right)}{t^{2}_{n}}\right)<+\infty.
\end{align*}
\item \textbf{Induction step}.\\
Suppose that the existence of all moments up to the order $2\left(n-1\right)=2\left\lfloor \dfrac{k}{2}\right\rfloor -2$
has been established. We must show that the moment of order $2n=2\left\lfloor \dfrac{k}{2}\right\rfloor $
exists. By the direct proposition
\[
\varphi^{\left(2\left(n-1\right)\right)}_{X}\left(t\right)+\varphi^{\left(2\left(n-1\right)\right)}_{X}\left(-t\right)=\left(-1\right)^{n-1}2\mathbb{E}\left(X^{2\left(n-1\right)}\cos\left(tX\right)\right)
\]
and
\[
\varphi^{\left(2\left(n-1\right)\right)}_{X}\left(0\right)=\left(-1\right)^{n-1}\mathbb{E}\left(X^{2\left(n-1\right)}\right).
\]
Moereover, $\varphi^{2\left(n-1\right)}_{X}$ is, by hypothesis, twice
differentiable at 0. Hence it admits a second-order Taylor-Young expansion,
and thus
\[
\lim_{t\to0}\dfrac{\varphi^{\left(2\left(n-1\right)\right)}_{X}\left(t\right)+\varphi^{\left(2\left(n-1\right)\right)}_{X}\left(-t\right)-2\varphi^{\left(2\left(n-1\right)\right)}_{X}\left(0\right)}{t^{2}}=\varphi^{\left(2n\right)}_{X}\left(0\right).
\]
Combining these three relations yields
\[
\lim_{t\to0}\mathbb{E}\left(X^{\left(2\left(n-1\right)\right)}\dfrac{1-\cos\left(tX\right)}{t^{2}}\right)=\dfrac{\left(-1\right)^{n}}{2}\varphi^{\left(2n\right)}_{X}\left(0\right).
\]
We conclude by applying the Fatou lemma, exactly as in the initialization
step.
\end{itemize}
\end{proof}

\begin{remark}{}{}

As the next example shows it---see \cite{stoyanov2013counterexamples}
and \cite{jeulin2006semi} Chapter 2, p.20---\textbf{the characteristic
function may be differentiable at the origin---and even differentiable
everywhere---without the random variable admitting an expectation.}

\end{remark}

\begin{example}{}{}

Let $X$ be a real-valued random variable with law $P_{X}=\sum_{k\in\mathbb{Z}}a_{k}\delta_{k},$
assumed to be symmetric, that is such that $a_{k}=a_{-k},$ and such
that $\sum_{k\in\mathbb{N}}ka_{k}=+\infty.$

Prove that $\varphi_{X}$ is differentiable everywhere, but that $X$
has no expectation, for a suitable choice of the sequence $\left(a_{k}\right)_{k\in\mathbb{Z}}.$

\end{example}

\begin{solutionexample}{}{}

We have
\[
\intop_{\Omega}\left|X\right|\text{d}P=2\sum_{k\in\mathbb{N}^{\ast}}ka_{k}=+\infty\,\,\,\,\text{and}\,\,\,\,\varphi_{X}\left(t\right)=a_{0}+2\sum^{+\infty}_{k=1}a_{k}\cos\left(kt\right).
\]

We choose the sequence $\left(a_{k}\right)_{k\in\mathbb{Z}}$ so that
the sequence $\left(ka_{k}\right)_{k\in\mathbb{N}^{\ast}}$ is nonincreasing
and converges to 0. Recall the upper bound
\[
\left|\sum^{n}_{k=0}\sin\left(kx\right)\right|\leqslant\left|\sum^{n}_{k=0}\text{e}^{\text{i}kx}\right|=\left|\dfrac{1-\text{e}^{\text{i}\left(n+1\right)x}}{1-\text{e}^{\text{i}x}}\right|\leqslant\dfrac{2}{\left|\sin\left(\dfrac{x}{2}\right)\right|}.
\]
Thus, for $\alpha\in\left]0,2\pi\right[,$
\[
\forall n\in\mathbb{N}^{\ast},\,\forall x\in\left[\alpha,2\pi-\alpha\right],\,\,\,\,\left|\sum^{n}_{k=0}\sin\left(kx\right)\right|\leqslant\dfrac{2}{\left|\sin\left(\dfrac{\alpha}{2}\right)\right|}.
\]

The Abel criteria ensures the uniform convergence on the interval
$\left[\alpha,2\pi-\alpha\right]$ of the series of functions with
general term $ka_{k}\sin kt.$ Therefore the function $\varphi_{X}$
is differentiable on this interval, and hence also on $\mathbb{R}\backslash2\pi\mathbb{Z},$
since it is $2\pi-$periodic. 

It remains to choose the sequence $\left(a_{k}\right)_{k\in\mathbb{Z}}$
so as to obtain the differentiability at 0. We take the sequence defined
by
\[
a_{0}=a_{1}=a_{-1}=0\,\,\,\,\text{and}\,\,\,\,\forall k\geqslant2,\,\,\,\,a_{k}=a_{-k}=\dfrac{c}{k^{2}\ln k},
\]
where 
\[
c=\dfrac{1}{2}\left(\sum^{+\infty}_{k=2}\dfrac{1}{k^{2}\ln k}\right)^{-1}
\]
---the series with general term $a_{k}$ is a convergent Bertrand
series. All the previous required conditions are then satisfied. 

Moreover, for every $t\neq0,$ by the transfer theorem
\[
0\leqslant\dfrac{1-\varphi_{X}\left(t\right)}{t}=\dfrac{1}{t}\mathbb{E}\left(1-\cos\left(tX\right)\right)=\dfrac{2c}{t}\sum^{+\infty}_{k=2}\dfrac{1}{k^{2}\ln k}\left(1-\cos\left(tk\right)\right).
\]
For every $t$ such that $0<t<\dfrac{1}{2},$ we split the sum in
two parts, depending on whether $k$ is smaller or larger than $t^{-1}.$
Since the functions $x\mapsto\left(\ln x\right)^{-1}$ and $x\mapsto x^{-2}$
are decreasing, comparison with series and integrals yields
\begin{multline*}
\dfrac{1}{t}\sum_{k\geqslant\frac{1}{t}}\dfrac{1}{k^{2}\ln k}\left(1-\cos\left(tk\right)\right)\leqslant-\dfrac{2}{t\ln t}\sum_{k\geqslant\left\lfloor \frac{1}{t}\right\rfloor }\dfrac{1}{k^{2}}\leqslant-\dfrac{2}{t\ln t}\intop^{+\infty}_{\left\lfloor \frac{1}{t}\right\rfloor -1}\dfrac{1}{x^{2}}\text{d}x\\
=-\dfrac{2}{t\left(\left\lfloor \frac{1}{t}\right\rfloor -1\right)\ln t}\leqslant-2\dfrac{\left\lfloor \frac{1}{t}\right\rfloor +1}{\left\lfloor \frac{1}{t}\right\rfloor -1}\cdot\dfrac{1}{\ln t}\underset{t\to0}{\longrightarrow}0.
\end{multline*}
Morevoer, using the inequality
\[
\forall x\in\mathbb{R},\,\,\,\,1-\cos x\leqslant\dfrac{x^{2}}{2},
\]
we similarly obtain, 
\begin{multline*}
\dfrac{1}{t}\sum_{2\leqslant k<\frac{1}{t}}\dfrac{1}{k^{2}\ln k}\left(1-\cos\left(tk\right)\right)\leqslant t\sum_{2\leqslant k<\frac{1}{t}}\dfrac{1}{\ln k}\leqslant\dfrac{t}{\ln2}+t\sum_{3\leqslant k\leqslant\left\lfloor \frac{1}{t}\right\rfloor }\dfrac{1}{\ln k}\\
\leqslant\dfrac{t}{\ln2}+t\sum_{3\leqslant k\leqslant\left\lfloor \frac{1}{t}\right\rfloor }\intop^{k}_{k-1}\dfrac{1}{\ln x}\text{d}x\leqslant\dfrac{t}{\ln2}+t\intop^{\frac{1}{t}}_{2}\dfrac{1}{\ln x}\text{d}x.
\end{multline*}

Finally, by integration by parts, 
\[
\intop^{y}_{2}\dfrac{1}{\ln x}\text{d}x=\left[\dfrac{x}{\ln x}\right]^{y}_{2}+\intop^{y}_{2}\dfrac{1}{\left(\ln x\right)^{2}}\text{d}x.
\]

Since as $x$ tends to infinity, $\frac{1}{\left(\ln x\right)^{2}}=o\left(\frac{1}{\ln x}\right),$
when $y$ tends to infinity, we also have
\[
\intop^{y}_{2}\dfrac{1}{\left(\ln x\right)^{2}}\text{d}x=o\left(\intop^{y}_{2}\dfrac{1}{\ln x}\text{d}x\right).
\]
It follows that
\[
\lim_{t\to0}t\intop^{\frac{1}{t}}_{2}\dfrac{1}{\ln x}\text{d}x=0.
\]
This proves that
\[
\lim_{t\to0}\dfrac{1-\varphi_{X}\left(t\right)}{t}=0,
\]
that is, $\varphi_{X}$ is differentiable at 0 with derivative 0.\textbf{
In summary, for this choice of law, $X$ does not admit an expectation,
while $\varphi_{X}$ is differentiable everywhere.}

\end{solutionexample}

\textbf{Proposition $\ref{pr:mom_der_real_rv}$ expands to the case
of random variables taking values in $\mathbb{R}^{d}.$} To simplify
the presentation, we restrict ourselves to moments of order at most
two.

\begin{proposition}{}{}

Let $X$ be a random variable taking values in $\mathbb{R}^{d},$
and let $\varphi_{X}$ denote its characteristic function.

(a) If $X$ admits an expectation---that is, $\left\Vert X\right\Vert $
is integrable---then $\varphi_{X}$ is differentiable. Its differential
at $t,$ viewed as a linear function from $\mathbb{R}^{d}$ to $\mathbb{C},$
is given by
\[
\forall t\in\mathbb{R}^{d},\,\forall x\in\mathbb{R}^{d},\,\,\,\,\varphi^{\prime}_{X}\left(t\right)\left(x\right)=\text{i}\intop_{\Omega}\left\langle X,x\right\rangle \text{e}^{\text{i}\left\langle X,t\right\rangle }\text{d}P.
\]
In particular,\boxeq{
\[
\forall x\in\mathbb{R}^{d},\,\,\,\,\varphi^{\prime}_{X}\left(0\right)\left(x\right)=\text{i}\left\langle \mathbb{E}\left(X\right),x\right\rangle .
\]
}

(b) If the norm of $X$ is of integrable square, $\varphi_{X}$ is
twice differentiable. Its second-order differential at $t,$ viewed
as a bilinear function from $\mathbb{R}^{d}\times\mathbb{R}^{d}$
to $\mathbb{C},$ is given by
\[
\forall t\in\mathbb{R}^{d},\,\forall x,y\in\mathbb{R}^{d},\,\,\,\,\varphi^{\prime\prime}_{X}\left(t\right)\left(x,y\right)=-\intop_{\Omega}\left\langle X,x\right\rangle \left\langle X,y\right\rangle \text{e}^{\text{i}\left\langle X,t\right\rangle }\text{d}P.
\]
In particular, \boxeq{
\[
\forall x\in\mathbb{R}^{d},\,\,\,\,\varphi^{\prime\prime}_{X}\left(0\right)\left(x,x\right)=-\mathbb{E}\left(\left\langle X,x\right\rangle ^{2}\right).
\]
}The variance of $X$ is then given by the relation\boxeq{
\[
\forall x\in\mathbb{R}^{d},\,\,\,\,\sigma^{2}_{X}\left(x\right)=-\varphi^{\prime\prime}_{X}\left(0\right)\left(x,x\right)+\left[\varphi^{\prime}_{X}\left(0\right)\left(x\right)\right]^{2}.
\]
}

Moreover, the covariance matrix $C_{X}$ of $X$ is given by
\[
C_{X}=\left(-\dfrac{\partial^{2}}{\partial t_{i}\partial t_{j}}\varphi_{X}\left(0\right)\right)+\left(\dfrac{\partial}{\partial t_{i}}\varphi_{X}\left(0\right)\dfrac{\partial}{\partial t_{j}}\varphi_{X}\left(0\right)\right).
\]

\end{proposition}

\begin{proof}{}{}

(a) Since
\[
\dfrac{\partial}{\partial t}\text{e}^{\text{i}\left\langle X,t\right\rangle }=\text{i}\left\langle X,\cdot\right\rangle \text{e}^{\text{i}\left\langle X,t\right\rangle },
\]
we have
\[
\left\Vert \dfrac{\partial}{\partial t}\text{e}^{\text{i}\left\langle X,t\right\rangle }\right\Vert _{\mathscr{L}\left(\mathbb{R}^{d},\mathbb{C}\right)}\leqslant\left\Vert X\right\Vert ,
\]
and we may apply Theorem $\ref{th:der_int_param}.$

(b) Similarly, for $x,y\in\mathbb{R}^{d},$
\[
\left[\dfrac{\partial^{2}}{\partial t^{2}}\text{e}^{\text{i}\left\langle X,t\right\rangle }\right]\left(x,y\right)=-\left\langle X,x\right\rangle \left\langle X,y\right\rangle \text{e}^{\text{i}\left\langle X,t\right\rangle },
\]
and thus
\[
\left\Vert \dfrac{\partial^{2}}{\partial t^{2}}\text{e}^{\text{i}\left\langle X,t\right\rangle }\right\Vert _{\mathscr{L}\left(\mathbb{R}^{d}\times\mathbb{R}^{d},\mathbb{C}\right)}\leqslant\left\Vert X\right\Vert ^{2},
\]
so we may again apply Theorem $\ref{th:der_int_param}.$

Finally, to compute the variance, it suffices to recall that
\[
\forall x\in\mathbb{R}^{d},\,\,\,\,\sigma^{2}_{X}\left(x\right)=\mathbb{E}\left(\left\langle X,x\right\rangle ^{2}\right)-\left(\mathbb{E}\left(\left\langle X,x\right\rangle \right)\right)^{2}.
\]
The auto-covariance operator of $X,$ denoted $\Lambda_{X},$ obtained
by bilinearization of the variance, satisfies
\[
\forall x,y\in\mathbb{R}^{d},\,\,\,\,\left\langle \Lambda_{X}x,y\right\rangle =-\varphi^{\prime\prime}_{X}\left(0\right)\left(x,y\right)+\left[\varphi^{\prime}_{X}\left(0\right)\left(x\right)\right]\left[\varphi^{\prime}_{X}\left(0\right)\left(y\right)\right],
\]
which yields the covariance matrix $C_{X},$ namely the matrix representation
of $\Lambda_{X}$ in the canonical basis---The partial derivatives
at 0 are precisely the values of the differential at 0 evaluated on
the vectors of this basis.

\end{proof}

\begin{example}{Follow-up of the Example $\ref{ex:mult_law_char_fct}$}{}

Compute the expectation and the covariance matrix of $Y_{n}$ for
$n\geqslant2.$

\end{example}

\begin{solutionexample}{}{}

By the equality $\refpar{eq:char_funct_Yn},$
\[
\dfrac{\partial}{\partial t_{j}}\varphi_{Y_{n}}\left(t\right)=\text{i}np_{j}\text{e}^{\text{i}t_{j}}\left[\sum^{k}_{m=1}p_{m}\text{e}^{\text{i}t_{m}}\right]^{n-1}.
\]
Hence
\[
\dfrac{\partial}{\partial t_{j}}\varphi_{Y_{n}}\left(0\right)=\text{i}np_{j},
\]
and therefore \boxeq{
\[
\mathbb{E}\left(Y_{n}\right)=np_{j}.
\]
}

If $j\neq l,$ then
\[
\dfrac{\partial^{2}}{\partial t_{l}\partial t_{j}}\varphi_{Y_{n}}\left(t\right)=-n\left(n-1\right)p_{j}p_{l}\text{e}^{\text{i}t_{j}}\text{e}^{\text{i}t_{l}}\left[\sum^{k}_{m=1}p_{m}\text{e}^{\text{i}t_{m}}\right]^{n-2},
\]
and thus
\[
\dfrac{\partial^{2}}{\partial t_{l}\partial t_{j}}\varphi_{Y_{n}}\left(t\right)=-n\left(n-1\right)p_{j}p_{l},
\]
which yields, after computation,\boxeq{
\[
\left(C_{Y_{n}}\right)_{jl}=-np_{j}p_{l},\,\,\,\,\text{if}\,j\neq l.
\]
}

Finally, 
\[
\dfrac{\partial^{2}}{\partial t^{2}_{j}}\varphi_{Y_{n}}\left(t\right)=\text{i}np_{j}\text{e}^{\text{i}t_{j}}\left[\text{i}\left(\sum^{k}_{m=1}p_{m}\text{e}^{\text{i}t_{m}}\right)^{n-1}+\text{i}\left(n-1\right)p_{j}\text{e}^{\text{i}t_{j}}\left(\sum^{k}_{m=1}p_{m}\text{e}^{\text{i}t_{m}}\right)^{n-2}\right],
\]
and thus,
\[
\dfrac{\partial^{2}}{\partial t^{2}_{j}}\varphi_{Y_{n}}\left(0\right)=-np_{j}\left[1+\left(n-1\right)p_{j}\right],
\]
which yields, after computation,\boxeq{
\[
\left(C_{Y_{n}}\right)_{jj}=np_{j}\left(1-p_{j}\right).
\]
}

\end{solutionexample}

The existence of moments makes it possible to obtain a Taylor expansion
of the characteristic function at 0, which is particularly useful
for the study of convergences in law---see Chapter \ref{chap:PartIIChapConvMeasuresandInLaw15}
for the proof of the central limit theorem---but also for computing
moments, by identifying the coefficients in the expansion. We also
provide a sufficient condition ensuring that a characteristic function
admits a power series expansion.

\begin{proposition}{Taylor Expansion in Zero and Power Series of A Characteristic Function}{tayl_exp_pow_ser_ch_f}

Let $X$ be a real-valued random variable and let $\varphi_{X}$ be
its characteristic function.

(a) If $X$ admits a moment of order $n\in\mathbb{N}^{\ast},$ then
$\varphi_{X}$ admits a Taylor expansion at 0 with an integral remainder,
which can be written, for every $t\in\mathbb{R},$ as\boxeq{
\begin{equation}
\varphi_{X}\left(t\right)=\sum^{n-1}_{k=0}\dfrac{\left(\text{i}t\right)^{k}}{k!}\mathbb{E}\text{\ensuremath{\left(X^{k}\right)}}+\dfrac{\left(\text{i}t\right)^{n}}{\left(n-1\right)!}\mathbb{E}\left[X^{n}\intop^{1}_{0}\left(1-u\right)^{n-1}\exp\left(\text{i}tuX\right)\text{d}u\right].\label{eq:Taylor_exp_zero_integ_rem_ch_f}
\end{equation}
}

It follows that\boxeq{
\begin{equation}
\forall t\in\mathbb{R},\,\,\,\,\varphi_{X}\left(t\right)=\sum^{n}_{k=0}\dfrac{\left(\text{i}t\right)^{k}}{k!}\mathbb{E}\left(X^{k}\right)+\dfrac{\left(\text{i}t\right)^{n}}{n!}\epsilon_{n}\left(t\right),\label{eq:Taylor_exp_order_n_in_zero}
\end{equation}
} where
\[
\left|\epsilon_{n}\left(t\right)\right|\leqslant2\mathbb{E}\left(\left|X^{n}\right|\right)\,\,\,\,\text{and}\,\,\,\,\lim_{t\to0}\epsilon_{n}\left(t\right)=0.
\]

Hence, in particular, we obtain a Taylor expansion of order $n$ of
$\varphi_{X}$ at the neighborhood of $0.$

(b) If $X$ admits moments of all orders, and if 
\begin{equation}
\limsup_{n\to+\infty}\dfrac{\left\Vert X\right\Vert _{n}}{n}=\dfrac{1}{R}<+\infty\label{eq:hyp_norm_Xovern_finite}
\end{equation}
where $\left\Vert X\right\Vert _{n}$ is the norm\footnotemark $n$
of $X$---this condition is satisfied in particular when the random
variable $X$ is bounded---then $\varphi_{X}$ admits a power series
expansion in a neighborhood of any real number, with radius of convergence
at least $\dfrac{R}{\text{e}}.$ In other words, the function $\varphi_{X}$
is analytic. In particular, on the interval $\left]-\dfrac{R}{\text{e}},\dfrac{R}{\text{e}}\right[,$
$\varphi_{X}$ admits the expansion\boxeq{
\begin{equation}
\forall t\in\left]-\dfrac{R}{\text{e}},\dfrac{R}{\text{e}}\right[,\,\,\,\,\varphi_{X}\left(t\right)=\sum^{+\infty}_{k=0}\dfrac{\left(\text{i}t\right)^{k}}{k!}\text{\ensuremath{\mathbb{E}}}\left(X^{k}\right).\label{eq:charact_fct_power_ser_dev}
\end{equation}
}

\end{proposition}

\footnotetext{Recall that the $n-$norm $\left\Vert X\right\Vert _{n}$
of a random variable $X$ is defined by 
\[
\left\Vert X\right\Vert _{n}=\left(\intop_{\Omega}\left|X\right|^{n}\text{d}P\right)^{\frac{1}{n}},
\]
whenever this quantity is finite.}

\begin{proof}{}{}

(a) The Taylor formula with integral remainder applied to the complex
exponential yields, for every $y\in\mathbb{R},$
\[
\text{e}^{\text{i}y}=\sum^{n-1}_{k=0}\dfrac{\left(\text{i}y\right)^{k}}{k!}+\dfrac{\left(\text{i}y\right)^{n}}{\left(n-1\right)!}\intop^{1}_{0}\left(1-u\right)^{n-1}\text{e}^{\text{i}uy}\text{d}u.
\]
Therefore, for every $t\in\mathbb{R},$
\begin{equation}
\text{e}^{\text{i}tX}=\sum^{n-1}_{k=0}\dfrac{\left(\text{i}t\right)^{k}}{k!}X^{k}+\dfrac{\left(\text{i}t\right)^{n}}{\left(n-1\right)!}X^{n}\intop^{1}_{0}\left(1-u\right)^{n-1}\text{e}^{\text{i}utX}\text{d}u.\label{eq:exp_itX_taylor_exp_int_rem}
\end{equation}

Integrating with respect to $P$ gives then the formula $\refpar{eq:Taylor_exp_zero_integ_rem_ch_f}.$
Note, moreover, that
\[
\intop^{1}_{0}\left(1-u\right)^{n-1}\text{d}u=\dfrac{1}{n}.
\]

For every $y\in\mathbb{R},$
\[
\text{e}^{\text{i}y}=\sum^{n}_{k=0}\dfrac{\left(\text{i}y\right)^{k}}{k!}+\dfrac{\left(\text{i}y\right)^{n}}{\left(n-1\right)!}\intop^{1}_{0}\left(1-u\right)^{n-1}\left(\text{e}^{\text{i}uy}-1\right)\text{d}u,
\]
which yields, for every $t\in\mathbb{R},$

\[
\text{e}^{\text{i}tX}=\sum^{n}_{k=0}\dfrac{\left(\text{i}t\right)^{k}}{k!}X^{k}+\dfrac{\left(\text{i}t\right)^{n}}{\left(n-1\right)!}X^{n}\intop^{1}_{0}\left(1-u\right)^{n-1}\left(\text{e}^{\text{i}utX}-1\right)\text{d}u.
\]

By integrating with respect to $P$, it follows, for every $t\in\mathbb{R},$
\[
\varphi_{X}\left(t\right)=\sum^{n}_{k=0}\dfrac{\left(\text{i}t\right)^{k}}{k!}\mathbb{E}\left(X^{k}\right)+\dfrac{\left(\text{i}t\right)^{n}}{n!}\epsilon_{n}\left(t\right),
\]
where
\[
\epsilon_{n}\left(t\right)=n\intop_{\Omega}\left(X^{n}\intop^{1}_{0}\left(1-u\right)^{n-1}\left(\text{e}^{\text{i}utX}-1\right)\text{d}u\right)\text{d}P.
\]

By the Fubini theorem, we obtain the upper bound
\[
\left|\epsilon_{n}\left(t\right)\right|\leqslant2n\mathbb{E}\left(\left|X^{n}\right|\right)\intop^{1}_{0}\left(1-u\right)^{n-1}\text{d}u=2\mathbb{E}\left(\left|X^{n}\right|\right).
\]

Moreover,
\[
\left|X^{n}\left(1-u\right)^{n-1}\left(\text{e}^{\text{i}utX}-1\right)\right|\leqslant2\left|X^{n}\right|\left|\left(1-u\right)^{n-1}\right|,
\]
which provides an upper bound independent of $t$ by a function that
is $\lambda_{\left[0,1\right]}\otimes P-$integrable. By the Fubini
theorem, and the dominated convergence theorem---with an arbitrary
sequence that converges to 0---it follows that
\[
\lim_{t\to0}\epsilon_{n}\left(t\right)=0.
\]

(b) Let $t_{0}$ be an arbitrary real number. Since the random variable
$X$ admits moments of all orders, its characteristic function $\varphi_{X}$---with
real valued variable---is of class $\mathscr{C}^{\infty}.$ It therefore
admits a Taylor expansion of order $n$ at $t_{0},$ namely, for every
real number $t,$
\[
\varphi_{X}\left(t\right)=\varphi_{X}\left(t_{0}\right)+\sum^{n}_{k=1}\dfrac{\left(t-t_{0}\right)^{k}}{k!}\varphi^{\left(k\right)}_{X}\left(t_{0}\right)+R_{n}\left(t_{0},t\right),
\]
 with remainder
\[
R_{n}\left(t_{0},t\right)=\intop^{t}_{t_{0}}\dfrac{\left(t-u\right)^{n}}{n!}\varphi^{\left(n+1\right)}_{X}\left(u\right)\text{d}u.
\]

Now, we want to prove that this remainder converges to 0 as $n$ tends
to $+\infty.$ From relation $\refpar{eq:k-th_deriv_charact_funct},$
\[
\left|R_{n}\left(t_{0},t\right)\right|\leqslant\dfrac{\left(\left|t-t_{0}\right|\left\Vert X\right\Vert _{n+1}\right)^{n+1}}{\left(n+1\right)!}.
\]

Let consider an arbitrary $\epsilon>0.$ Assumption $\refpar{eq:hyp_norm_Xovern_finite}$
implies that there exists $N$ such that, for every $n\geqslant N,$
\[
\dfrac{\left\Vert X\right\Vert _{n}}{n}\leqslant\dfrac{1}{R}+\epsilon.
\]

Using the Stirling formula, we obtain
\[
\left(\dfrac{\left|t-t_{0}\right|^{n}\left\Vert X\right\Vert ^{n}_{n}}{n!}\right)^{\frac{1}{n}}\leqslant\left|t-t_{0}\right|\left(\dfrac{1}{R}+\epsilon\right)n\left[\left(\dfrac{n}{\text{e}}\right)\sqrt{2\pi n}\left(1+\dfrac{1}{12n}+o\left(\dfrac{1}{n}\right)\right)\right]^{-\frac{1}{n}}.
\]
The right member converges to $\left|t-t_{0}\right|\left(\dfrac{1}{R}+\epsilon\right)\text{e},$
and it follows from the arbitrary of $\epsilon$ that
\[
\limsup_{n\to+\infty}\left(\dfrac{\left|t-t_{0}\right|^{n}\left\Vert X\right\Vert ^{n}_{n}}{n!}\right)^{\frac{1}{n}}\leqslant\left|t-t_{0}\right|\dfrac{\text{e}}{R}.
\]

Then, for every $t\in\mathbb{R},$ such that 
\[
\left|t-t_{0}\right|<\dfrac{R}{\text{e}},
\]
the Cauchy test shows that the series with general term $\dfrac{\left|t-t_{0}\right|^{n}\left\Vert X\right\Vert ^{n}_{n}}{n!}$
converges. Thus, the Taylor remainder $R_{n}\left(t_{0},t\right)$
tends to $0$ when $n$ tends to infinity, which proves the existence
of the Taylor expansion of $\varphi_{X}$ at $t_{0}$ and also the
analyticity of $\varphi_{X}.$ 

Taking $t_{0}=0$ and using the values of the derivatives of $\varphi_{X}$
at 0 given by $\refpar{eq:k-th_der_char_funct_in_zero},$ we obtain
the expansion $\refpar{eq:charact_fct_power_ser_dev}.$

\end{proof}

\section*{Exercises}

\addcontentsline{toc}{section}{Exercises}

All the introduced random variables introduced below are defined on
the same probabilized space $\left(\Omega,\mathscr{A},P\right).$

\begin{exercise}{Characteristic Function and Injectivity. Triangular Law}{exercise13.1}

Characteristics functions may coincide on an interval without being
equal---see \cite{stoyanov2013counterexamples}.

Let $\phi$ be a function defined for every $t\in\mathbb{R},$ by
\[
\phi\left(t\right)=\begin{cases}
1-\left|t\right|, & \text{if}\,\,\,\,0\leqslant\left|t\right|\leqslant1,\\
0, & \text{otherwise.}
\end{cases}
\]

Let $X$ be a random variable of the triangular law on the interval
$\left[-1,1\right],$ that is, such that $X$ has density $f_{X}=\phi.$

1. Compute the characteristic function $\varphi_{X}$ of $X.$ Denote
by $\rho$ the uniform law on $\left[-\dfrac{1}{2},\dfrac{1}{2}\right].$
Justify the fact that the law $P_{X}$ of $X$ is the convolution
product $\rho\ast\rho.$

2. Prove that $\phi$ is the Fourier transform of a probability $\mu=f\cdot\lambda$
on $\mathbb{R}$ where $f$ is a density of probability to be determined.

3. Let $Y$ and $Z$ be two random variables such that $Y$ has density
$f$ and $Z$ takes values in the set of integers $\mathbb{Z},$ with
law
\[
P_{Z}=\dfrac{1}{2}\delta_{0}+\sum_{k\in\mathbb{Z}}\dfrac{2}{\left(2k-1\right)^{2}\pi^{2}}\delta_{\left(2k-1\right)\pi}.
\]

The aim is to prove that the characteristic functions of $Y$ and
$Z$ coincide on the interval $\left[-1,1\right]$ but are not equal.
To this end, expand into a Fourier series the function $\psi,$ periodic
with period 2, equal to $1-\phi$ on the set $\left[-1,1\right],$
and conclude.

\end{exercise}

\begin{exercise}{Characteristic Function of a Product of Independent Random Variables}{exercise13.2}

Let $X$ and $Y$ be two independent real-valued random variables.

1. Prove that the characteristic function of the product $XY$ is
given by the relation: for every $t\in\mathbb{R},$
\[
\varphi_{XY}\left(t\right)=\intop_{\mathbb{R}}\varphi_{X}\left(ty\right)\text{d}P_{Y}\left(y\right).
\]

Moreover, if $X$ and $Y$ have same normal law $\mathscr{N}_{\mathbb{R}}\left(0,1\right),$
determine the characteristic function of $XY.$

2. Let $X_{1},X_{2},X_{3},$ and $X_{4}$ be four independent real-valued
random variables, with the same normal law $\mathscr{N}_{\mathbb{R}}\left(0,1\right).$
Determine the characteristic function and the law of the random variable
$X_{1}X_{2}+X_{3}X_{4}.$

3. What is the law of the random variable $\left|X_{1}X_{2}+X_{3}X_{4}\right|?$

\end{exercise}

\begin{exercise}{Characteristic Function, Convolution and Moments}{exercise13.3}

Let $U$ be a real-valued random variable with density $f_{U}$ defined,
for every $u\in\mathbb{R},$ by
\[
f_{U}\left(u\right)=\sum^{+\infty}_{n=0}\boldsymbol{1}_{\left[n,n+1\right[}\left(u\right)\text{e}^{-\lambda}\dfrac{\lambda^{n}}{n!}.
\]

1. Compute the characteristic function $\varphi_{U}$ of $U,$ and
deduce from it that the law of $U$ is the convolution of two laws
to be specified.

2. Determine, without computation, the expectation and the variance
of $U.$

3. Let $T$ be a random variable independent of $U,$ with the uniform
law on the interval $\left[0,1\right].$ Determine the characteristic
function of the random variable $W=T+U.$\\
\textit{Hint: refer to the first question. }

Justify its differentiability and give, without computation its derivative
in 0.

\end{exercise}

\begin{exercise}{}{exercise13.4}

Let $X,\,Y$ and $Z$ be three independent random variables with the
same normal law $\mathscr{N}_{\mathbb{R}}\left(0,1\right).$

1. Determine the moments of every order of $X$ from its characteristic
function.

2. Deduce from this the moments of every order of the random variable
$U=XY.$ Using this approach, determine the characteristic function
of $U.$ The reader may refer to Exercise $\ref{exo:exercise13.1}$
for another method, which is indeed faster.

3. Let $V$ be the random variable $YZ.$ Determine the characteristic
function of the pair $\left(U,V\right).$ Deduce from it the characteristic
function of the random variable $\dfrac{U+V}{\sqrt{2}}.$ Compare
the law of this last random variable with that of $U.$

4. Are the random variables $U$ and $V$ independent? Are they correlated?

\end{exercise}

\begin{exercise}{Independence Criterion of Bounded Random Variables---M. Kac}{exercise13.5}

Let $X$ and $Y$ be two bounded real-valued random variables. Prove
that for $X$ and $Y$ are independent if and only if
\begin{equation}
\forall\left(k,l\right)\in\mathbb{N}^{2},\,\,\,\,\mathbb{E}\left(X^{k}Y^{l}\right)=\mathbb{E}\left(X^{k}\right)\mathbb{E}\left(Y^{l}\right).\label{eq:bounded_real_val_var_indep_criteria}
\end{equation}

\end{exercise}

\begin{exercise}{Case Where $\varphi_{X+Y}=\varphi_X\varphi_Y$ Without The Random Variables $X$ and $Y$ Being Independent}{exercise13.6}

Let $f$ be the function defined on $\mathbb{R}^{2}$ by
\[
\forall\left(x,y\right)\in\mathbb{R}^{2},\,\,\,\,f\left(x,y\right)=\dfrac{1}{4}\boldsymbol{1}_{C}\left(x,y\right)\left[1+xy\left(x^{2}-y^{2}\right)\right],
\]
where $C=\left[-1,1\right]^{2}.$

1. Check that $f$ is a probability density on $\mathbb{R}^{2}$ with
respect to the Lebesgue measure.

2. Let $\left(X,Y\right)$ be a random variable taking values in $\mathbb{R}^{2},$
with density $f.$ Define the random variable $Z=X+Y.$ Compute the
densities of the random variables $X,\,Y,\,Z.$ Specify the laws of
$X$ and $Y.$

3. Prove that if $U$ is a real-valued random variable admitting an
even density, then its characteristic function satisfies
\[
\forall t\in\mathbb{R},\,\,\,\,\varphi_{U}\left(t\right)=2\intop^{+\infty}_{0}\cos\left(tu\right)f_{U}\left(u\right)\text{d}u.
\]

Then compute the characteristic functions of $X$ and $Y.$ Express
the characteristic function of $Z$ as a function of those of $X$
and $Y.$

4. Compute the correlation coefficient of $X$ and $Y.$

5. Remarks on this exercise.

\end{exercise}

\begin{exercise}{Another Example Where $\varphi_{X+Y}=\varphi_X\varphi_Y$ Without The Random Variables $X$ and $Y$ Being Independent. Cauchy Laws}{exercise13.7}

The Cauchy law $\mu_{a}$ with parameter $a>0$ is the probability
on $\mathbb{R}$ with density $f_{a}$ defined by, for every real
number $x,$
\[
f_{a}\left(x\right)=\dfrac{a}{\pi\left(a^{2}+x^{2}\right)}.
\]

1. Prove the following relation for the Fourier transform:
\[
\forall t\in\mathbb{R},\,\,\,\,\widehat{\mu_{a}}\left(t\right)=\widehat{\mu_{1}}\left(at\right).
\]

Deduce from this relation that if a random variable $Z$ follows a
Cauchy law with parameter 1, then the random variable $aZ$ follows
a Cauchy law with parameter $a.$

2. Let $U$ and $V$ be two random variables of Cauchy law with parameter
1. Recall that 
\[
\widehat{\mu_{1}}\left(t\right)=\text{e}^{-\left|t\right|},
\]
as can be shown by a simple computation using residues. Let $a,\,b,\,c$
and $d$ be four positive real numbers. Let $X$ and $Y$ be the random
variables defined by
\[
X=aU+bV\,\,\,\,\,\,\,\,Y=cU+dV.
\]

Compute the characteristic function of the random variable $\left(X,Y\right)$
and deduce that $X$ and $Y$ are not independent.

3. Compute the characteristic function of $X+Y$ and deduce from it
the equality of laws
\[
P_{X+Y}=P_{X}\ast P_{Y}.
\]

\end{exercise}

\begin{exercise}{Characteristic Function and Support of Law}{exercise13.8}

Let $X$ be a real-valued random variable with characteristic function
$\varphi_{X}.$

1. Prove that if there exists a real number $t_{0}\neq0$ such that
$\left|\varphi_{X}\left(t_{0}\right)\right|=1,$ then there exists
a real number $a$ such that 
\[
P_{X}\left(a+\mathbb{Z}\dfrac{2\pi}{t_{0}}\right)=1.
\]

2. Prove that if there exists two nonzero real numbers $t_{1}$ and
$t_{2}$ such that $\dfrac{t_{1}}{t_{2}}$ is irrational and 
\[
\left|\varphi_{X}\left(t_{1}\right)\right|=\left|\varphi_{X}\left(t_{2}\right)\right|=1,
\]
then the random variable $X$ is degenerate---that is, $X$ is $P-$almost
surely equal to a constant.

3. Prove that for the random variable $X$ to be degenerate, it is
necessary and sufficient that $\left|\varphi_{X_{0}}\right|=1.$

\end{exercise}

\begin{exercise}{Characteristic Function and Conditional Expectation. Game of Heads and Tails---Variant}{exercise13.9}

Let $\left(U_{n}\right)_{n\in\mathbb{N}}$ be a sequence of independent
real-valued random variables of same law $\dfrac{\delta_{-1}+\delta_{1}}{2}.$
Set $U_{-1}=0.$ Define the sequence of random variables $\left(Y_{n}\right)_{n\in\mathbb{N}}$
by
\[
Y_{n}=\sum^{n}_{j=0}U_{j-1}U_{j}.
\]

For every $n\in\mathbb{N},$ denote by $\mathscr{F}_{n}$ the $\sigma-$algebra
generated by the random variables $U_{j},\,0\leqslant j\leqslant n.$ 

The random variable $Y_{n}$ may represent the algebraic gain, after
the $n-$th toss, of a player playing heads and tails with a fair
coin, with the following gain rule: the play wins one unit after the
$n-$th toss if the outcome is the same as the previous toss. Otherwise,
the player loses one unit. This corresponds to the usual strategy
of a player always betting on the face that just appeared.

1. Compute, for every real number $t,$ the conditional expectation
\[
\mathbb{E}^{\mathscr{F}_{n-1}}\left(\text{e}^{\text{i}tU_{n-1}U_{n}}\right).
\]

2. Compute, for every real number $t$ and for every integer $l$
such that $1\leqslant l\leqslant n,$ the conditional expectation
\[
\mathbb{E}^{\mathscr{F}_{n-l}}\left(\text{e}^{\text{i}tY_{n}}\right).
\]

\textit{Hint: One may proceed by induction on $l.$}

3. Deduce the characteristic function $\varphi_{Y_{n}}$ of $Y_{n}.$

4. Determine the law of $Y_{n}.$

5. Compute, for every real number $t,$ the sequence with general
term $\varphi_{\frac{Y_{n}}{n}}\left(t\right).$

\end{exercise}

\begin{exercise}{Inversion Formula of a Probability Fourier Transform}{exercise13.10}

Let $\mu$ be a probability on $\left(\mathbb{R},\mathscr{B}_{\mathbb{R}}\right)$
with Fourier transform $\varphi.$ We begin by proving the following
inversion formula: for every real numbers $a$ and $b$ such that
$a<b,$
\begin{equation}
\lim_{T\to+\infty}\dfrac{1}{2\pi}\intop^{T}_{-T}\dfrac{\text{e}^{-\text{i}ta}-\text{e}^{-\text{i}tb}}{\text{i}t}\varphi\left(t\right)\text{d}t=\dfrac{1}{2}\mu\left(\left\{ a,b\right\} \right)+\mu\left(\left]a,b\right[\right).\label{eq:first_inversion_formula}
\end{equation}

The results of questions 2 and 3 are particularly interesting and
useful.

1. For every nonnegative real number $T,$ define the complex-valued
function of the real variable $x$ by
\[
\forall x\in\mathbb{R},\,\,\,\,I_{T}\left(x\right)=\dfrac{1}{2\pi}\intop^{+T}_{-T}\dfrac{\text{e}^{\text{i}t\left(x-a\right)}-\text{e}^{\text{i}t\left(x-b\right)}}{\text{i}t}\text{d}t.
\]

Prove that, for every real number $x,$
\[
\lim_{T\to+\infty}I_{T}\left(x\right)=\dfrac{1}{2}\boldsymbol{1}_{\left\{ a,b\right\} }\left(x\right)+\boldsymbol{1}_{\left]a,b\right[}\left(x\right).
\]

2. Deduce from the first question the inversion formula $\refpar{eq:first_inversion_formula}.$

3. Follow the same steps, prove that for every real number $b,$
\[
\mu\left(\left\{ b\right\} \right)=\lim_{T\to+\infty}\dfrac{1}{2T}\intop^{+T}_{-T}\text{e}^{-\text{i}tb}\varphi\left(t\right)\text{d}t.
\]

4. Prove the equality in $\overline{\mathbb{R}}^{+},$
\[
\lim_{T\to+\infty}\dfrac{1}{2T}\intop^{+T}_{-T}\left|\varphi\left(t\right)\right|^{2}\text{d}t=\sum_{x\in\mathbb{R}}\mu\left(\left\{ x\right\} \right)^{2}.
\]
 To this end, introduce two independent random variables $X$ and
$Y$ with law $\mu,$ and apply the previous results to the characteristic
function of the random variable $X-Y.$

\end{exercise}

\begin{exercise}{Characteristic Function of Vectorial Random Variables. Computation of Conditional Expectation and Injectivity of the Fourier Transform}{exercise13.11}

Let $X$ and $Y$ be two real valued random variables such that $Y$
follows a Bernoulli law $\mathscr{B}\left(1,1-\rho\right),$ and such
that a conditional law of $X$ given $Y$ is defined by
\[
P^{Y=0}_{X}=\delta_{0}\,\,\,\,\text{and}\,\,\,\,P^{Y=1}_{X}=\text{e}^{\lambda}
\]
where $\delta_{0}$ is the Dirac measure at 0, and $\exp\left(\lambda\right)$
the exponential law with parameter $\lambda>0.$

1. Compute the characteristic function $\varphi_{X}$ of $X.$ Denote
by $\mu$ the law of $X$---it is not required to compute it explicitly.

2. Deduce from the first question the expectation and the variance
of $X.$

Now consider a family of independent random variables $X_{0},\epsilon_{n},n\in\mathbb{N}^{\ast}.$
Assume that the $\epsilon_{n}$ have common law $\mu,$ and that $X_{0}$
has law $\exp\left(\lambda\right).$ Define inductively the random
variables $X_{n}$ by
\[
\forall n\in\mathbb{N}^{\ast},\,\,\,\,X_{n}=\rho X_{n-1}+\epsilon_{n}.
\]

3. Justify the independence, for every $n\in\mathbb{N}^{\ast},$ of
the random variables $\epsilon_{n}$ and $\left(X_{0},X_{1},\cdots,X_{n-1}\right).$
Prove by induction that the random variables $X_{n}$ have the same
characteristic function. Identify their law. 

4. Compute the conditional mean $m^{X_{n-1}=x_{n-1}}_{X_{n}}$ of
$X_{n}$ given $X_{n-1},$ for every real number $x_{n-1}.$

5. Express the characteristic function $\varphi_{\left(X_{n-1},X_{n}\right)}$
of the random variable $\left(X_{n-1},X_{n}\right),$ for every $\left(u,v\right)\in\mathbb{R}^{2},$
as a function of $u,\,v,\,\lambda$ and $\rho.$

6. Show that
\[
\forall v\in\mathbb{R},\,\,\,\,\dfrac{\partial}{\partial u}\varphi_{\left(X_{n-1},X_{n}\right)}\left(0,v\right)=\text{i}\intop_{\mathbb{R}}\text{e}^{\text{i}vx}m^{X_{n}=x}_{X_{n-1}}\text{d}P_{X_{n}}\left(x\right).
\]

Using the injectivity theorem of the Fourier transform, show that

\[
\forall x\in\mathbb{R},\,\,\,\,m^{X_{n}=x}_{X_{n-1}}=\boldsymbol{1}_{\mathbb{R}^{+}}\left(x\right)\dfrac{1}{\lambda\left(1-\rho\right)}\left[1-\text{e}^{-\frac{\lambda\left(1-\rho\right)}{\rho}x}\right].
\]

\end{exercise}

\section*{Solutions of Exercises}

\addcontentsline{toc}{section}{Solutions of Exercises}

\begin{solution}{}{solexercise13.1}

\textbf{1. Computation of $\varphi_{X}.$ Expression of $P_{X}$ as
a convolution product $\rho\ast\rho$}

Since the density of $X$ is even, then
\[
\varphi_{X}\left(t\right)=\intop^{+1}_{-1}\left(1-\left|t\right|\right)\text{e}^{\text{i}tx}\text{d}x=2\intop^{+1}_{0}\left(1-t\right)\cos\left(tx\right)\text{d}x,
\]
An integration by parts, valid for every $t\neq0,$ yields\boxeq{
\[
\forall t\neq0,\,\,\,\,\varphi_{X}\left(t\right)=2\dfrac{1-\cos t}{t^{2}}\,\,\,\,\text{and}\,\,\,\,\varphi_{X}\left(0\right)=1.
\]
}

This expression can also be written
\[
\forall t\neq0,\,\,\,\,\varphi_{X}\left(t\right)=4\dfrac{\sin^{2}\dfrac{t}{2}}{t^{2}}=\left(\widehat{\rho}\left(t\right)\right)^{2}\,\,\,\,\text{and}\,\,\,\,\varphi_{X}\left(0\right)=1,
\]
 since the Fourier transform of $\rho$ is equal, for every $\rho\neq0,$
\[
\widehat{\rho}\left(t\right)=\intop^{\frac{1}{2}}_{-\frac{1}{2}}\text{e}^{\text{i}tx}\text{d}x=2\dfrac{\sin\dfrac{t}{2}}{t}.
\]

By the injectivity property of the Fourier transform, it follows that
the law of $X$ is the convolution of the uniform law on $\left[-\dfrac{1}{2},\dfrac{1}{2}\right]$
with itself.

\textbf{2. Expression of $\phi$ as the Fourier transform of $\mu=f\cdot\lambda$}

Since $\phi$ is integrable, Proposition $\ref{pr:abs_cont_meas_four_tr}$
implies that, if $\mu$ is a probability such that $\widehat{\mu}=\phi,$
then it is with density $f$ given by
\[
\forall x\in\mathbb{R},\,\,\,\,f\left(x\right)=\dfrac{1}{\sqrt{2\pi}}\intop_{\mathbb{R}}\widehat{\mu}\left(t\right)\text{e}^{-\text{i}xt}\text{d}t\equiv\dfrac{1}{\sqrt{2\pi}}\widehat{\phi}\left(-x\right).
\]

Since the random variable $X$ has density $f_{X}=\phi,$ its characteristic
function $\varphi_{X}$ is equal to $\widehat{\phi}.$ Therefore,
\[
f\left(x\right)=\dfrac{1}{\sqrt{2\pi}}\varphi_{X}\left(-x\right),
\]
and since $\varphi_{X}$ is even,
\[
f\left(x\right)=\dfrac{1}{\sqrt{2\pi}}\varphi_{X}\left(x\right).
\]
It remains to verify that the probability $\mu=f\cdot\lambda$ thus
defined satisfies $\widehat{\mu}=\phi.$ For every real number $t,$
\[
\widehat{\mu}\left(t\right)=\intop_{\mathbb{R}}f\left(x\right)\text{e}^{\text{i}xt}\text{d}x=\dfrac{1}{\sqrt{2\pi}}\intop_{\mathbb{R}}\varphi_{X}\left(x\right)\text{e}^{\text{i}xt}\text{d}x.
\]
Since $\varphi_{X}$ is integrable, Proposition $\ref{pr:abs_cont_meas_four_tr}$
yields 
\[
\widehat{\mu}\left(t\right)=f_{X}\left(-t\right),
\]
and since $f_{X}$ is even, we conclude that 
\[
\widehat{\mu}=f_{X}=\phi.
\]

\textbf{3. Coincidence without equality of the characteristic functions
of $Y$ and $Z$ on $\left[-1,1\right]$}

For every real number $t,$ the characteristic function of $Z$ is
given by
\[
\varphi_{Z}\left(t\right)=\dfrac{1}{2}+\dfrac{2}{\pi^{2}}\sum_{k\in\mathbb{Z}}\dfrac{2}{\left(2k-1\right)^{2}}\text{e}^{\text{i}t\left(2k-1\right)\pi}.
\]
Equivalently,
\[
\varphi_{Z}\left(t\right)=\dfrac{1}{2}+\dfrac{4}{\pi^{2}}\sum^{+\infty}_{k=1}\dfrac{1}{\left(2k-1\right)^{2}}\cos\left(\left(2k-1\right)\pi t\right).
\]

Moreover, the function $\psi$ is even, periodic with period 2, continuous,
and piecewise $\mathscr{C}^{1}.$ By the Dirichlet theorem, for every
real number $t,$
\[
\psi\left(t\right)=\dfrac{a_{0}}{2}+\sum^{+\infty}_{n=1}a_{n}\cos\left(2\pi n\dfrac{t}{2}\right),
\]
where
\[
a_{0}=\intop^{1}_{-1}\left|t\right|\text{d}t=1
\]
and, for $n\in\mathbb{N}^{\ast},$
\begin{gather*}
a_{n}=\intop^{1}_{-1}\left|t\right|\cos\left(2\pi n\dfrac{t}{2}\right)\text{d}t=2\intop^{1}_{0}t\cos\left(\pi nt\right)\text{d}t.
\end{gather*}
 An integration by parts yields, for $n\in\mathbb{N}^{\ast},$
\[
a_{n}=-\dfrac{2}{\pi^{2}n^{2}}\left(1-\cos\left(\pi n\right)\right)=\begin{cases}
-\dfrac{4}{\pi^{2}n^{2}}, & \text{if }n\text{ is even,}\\
0, & \text{otherwise.}
\end{cases}
\]

Therefore, for every real number $t,$
\[
\psi\left(t\right)=\dfrac{1}{2}-\dfrac{4}{\pi^{2}}\sum^{+\infty}_{k=1}\dfrac{1}{\left(2k-1\right)^{2}}\cos\left[\left(2k-1\right)\pi t\right].
\]
It follows that, on the interval $\left[-1,1\right],$
\[
\varphi_{Z}\left(t\right)=1-\psi\left(t\right)=\phi\left(t\right)=\varphi_{Y}\left(t\right).
\]

Finally, note that the random variable $Z$ is discrete, whereas the
random variable $Y$ admits a density, which proves that their characteristic
functions are not equal everywhere.

\end{solution}

\begin{solution}{}{solexercise13.2}

\textbf{1. Proof that for every $t\in\mathbb{R},$ $\varphi_{XY}\left(t\right)=\intop_{\mathbb{R}}\varphi_{X}\left(ty\right)\text{d}P_{Y}\left(y\right).$
Characteristic function of $XY$ when $X$ and $Y$ follow the centered
reduced normal law}

Since the function $\left(x,y\right)\mapsto\text{e}^{\text{i}txy}$
is bounded, the transfer theorem ensures that
\[
\varphi_{XY}\left(t\right)=\intop_{\mathbb{R}^{2}}\text{e}^{\text{i}txy}\text{d}P_{\left(X,Y\right)}\left(x,y\right).
\]

As the random variables $X$ and $Y$ are independent, the law of
the pair $\left(X,Y\right)$ is the product of the laws of $X$ and
$Y.$ The function $\left(x,y\right)\mapsto\text{e}^{\text{i}txy}$
being bounded and therefore $P_{X}\otimes P_{Y}-$integrable, the
Fubini theorem yields,
\[
\varphi_{XY}\left(t\right)=\intop_{\mathbb{R}}\left[\intop_{\mathbb{R}}\text{e}^{\text{i}txy}\text{d}P_{X}\left(x\right)\right]\text{d}P_{Y}\left(y\right),
\]
which proves the requested formula. 

In the case where the random variables have the same centered reduced
normal law, the characteristic function of $XY$ is then given, for
every $t\in\mathbb{R},$ by
\begin{align*}
\varphi_{XY}\left(t\right) & =\intop_{\mathbb{R}}\text{e}^{-\frac{t^{2}y^{2}}{2}}\dfrac{1}{\sqrt{2\pi}}\text{e}^{-\frac{y^{2}}{2}}\text{d}y\\
 & =\intop_{\mathbb{R}}\dfrac{1}{\sqrt{2\pi}}\text{e}^{-\frac{\left(1+t^{2}\right)y^{2}}{2}}\text{d}y.
\end{align*}
Thus,\boxeq{
\[
\varphi_{XY}\left(t\right)=\dfrac{1}{\sqrt{1+t^{2}}}.
\]
}

\textbf{2. Characteristic function and law of $X_{1}X_{2}+X_{3}X_{4}$}

The random variables $X_{1}X_{2}$ and $X_{3}X_{4}$ are independent
and follow the same law as $XY.$ Consequently, for every $t\in\mathbb{R},$
\[
\varphi_{X_{1}X_{2}+X_{3}X_{4}}\left(t\right)=\varphi_{X_{1}X_{2}}\left(t\right)\varphi_{X_{3}X_{4}}\left(t\right)=\dfrac{1}{1+t^{2}}.
\]
By the injectivity theorem, it follows that the Fourier transform
of $X_{1}X_{2}+X_{3}X_{4}$ corresponds to the \textbf{\mindex{law!Laplace}\index{Laplace law}Laplace\sindex[fam]{Laplace, Pierre-Simon}
law} with density the function $x\mapsto\dfrac{1}{2}\text{e}^{-\left|x\right|},$
as listed in the table of laws.

\textbf{3. Law of the random variable $\left|X_{1}X_{2}+X_{3}X_{4}\right|$}

For every $f\in\mathscr{C}^{+}_{\mathscr{K}}\left(\mathbb{R}\right),$
\[
\intop_{\mathbb{R}}f\left(\left|u\right|\right)\text{d}P_{\left|X_{1}X_{2}+X_{3}X_{4}\right|}\left(u\right)=\intop_{\mathbb{R}}f\left(\left|u\right|\right)\dfrac{1}{2}\text{e}^{-\left|u\right|}\text{d}u.
\]
Using the parity of the integrand, this becomes
\[
\intop_{\mathbb{R}}f\left(\left|u\right|\right)\text{d}P_{\left|X_{1}X_{2}+X_{3}X_{4}\right|}\left(u\right)=\intop_{\mathbb{R}}f\left(\left|u\right|\right)\boldsymbol{1}_{\mathbb{R}^{+}}\text{e}^{-\left|u\right|}\text{d}u.
\]
Therefore, the random variable $\left|X_{1}X_{2}+X_{3}X_{4}\right|$
follows the exponential law with parameter 1, denoted $\exp\left(1\right).$

\end{solution}

\begin{solution}{}{solexercise13.3}

\textbf{1. Computation of $\varphi_{U}.$ Law of $U$ as the convolution
of two laws to be specified.}

Let $t$ be an arbitrary real number. We have
\[
\varphi_{U}\left(t\right)=\intop_{\mathbb{R}}\text{e}^{\text{i}tu}f_{U}\left(u\right)\text{d}u.
\]

Since
\[
\sum^{+\infty}_{n=0}\intop_{\mathbb{R}}\left|\text{e}^{\text{i}tu}\right|\boldsymbol{1}_{\left[n,n+1\right[}\left(u\right)\text{e}^{-\lambda}\dfrac{\lambda^{n}}{n!}\text{d}u=\sum^{+\infty}_{n=0}\text{e}^{-\lambda}\dfrac{\lambda^{n}}{n!}=1,
\]
it follows, by the corollary of the dominated convergence theorem
for series of functions---see Chapter \ref{chap:PartIIChap8}, Corollary
$\ref{co:sum_measurable_functions_with_bounded_abs_value},$---that,
for every $t\neq0,$
\begin{align*}
\varphi_{U}\left(t\right) & =\sum^{+\infty}_{n=0}\left[\text{e}^{-\lambda}\dfrac{\lambda^{n}}{n!}\intop_{\mathbb{R}}\text{e}^{\text{i}tu}\boldsymbol{1}_{\left[n,n+1\right[}\left(u\right)\text{d}u\right]\\
 & =\sum^{+\infty}_{n=0}\text{e}^{-\lambda}\dfrac{\lambda^{n}}{n!}\dfrac{\text{e}^{\text{i}t\left(n+1\right)}-\text{e}^{\text{i}tn}}{\text{i}t}\\
 & =\dfrac{\text{e}^{-\lambda}\left(\text{e}^{\text{i}t}-1\right)}{\text{i}t}\sum^{+\infty}_{n=0}\dfrac{\left[\lambda\text{e}^{\text{i}t}\right]^{n}}{n!}.
\end{align*}
Hence,\boxeq{
\[
\varphi_{U}\left(t\right)=\dfrac{\text{e}^{\text{i}t}-1}{\text{i}t}\text{e}^{\lambda\left(\text{e}^{\text{i}t}-1\right)}.
\]
}

By the Fourier transform injectivity property, it follows that the
law of $U$ is the convolution of the uniform law on $\left[0,1\right]$
and the Poisson law $\mathscr{P}\left(\lambda\right).$

\textbf{2. Determination of the expectation and variance of $U$}

Let $X$ and $N$ be two independent random variables, where $X$
follows the uniform law on the interval $\left[0,1\right]$ and $N$
the Poisson law $\mathscr{P}\left(\lambda\right).$ Then the law of
$X+N$ coincides with the law of $U.$

Thus,

\boxeq{
\[
\mathbb{E}\left(U\right)=\mathbb{E}\left(X\right)+\mathbb{E}\left(N\right)=\dfrac{1}{2}+\lambda,
\]
}and, since $X$ and $N$ are independent,

\boxeq{
\[
\sigma^{2}_{U}=\sigma^{2}_{X}+\sigma^{2}_{N}=\dfrac{1}{12}+\lambda.
\]
}

\textbf{3. Characteristic function of $W=T+U.$ Derivative at 0}

Assume moreover that the random variables $X$ and $N$ are also independent
of $T.$ Then
\[
P_{W}=P_{T}\ast P_{U}=P_{T}\ast\left(P_{X}\ast P_{N}\right).
\]
Taking the Fourier transforms, we obtain
\[
\varphi_{W}=\varphi_{T}\varphi_{X}\varphi_{N}.
\]
Note that this argument also shows, again via the injectivity of the
Fourier transform, that the \textbf{convolution product is associative}.
Hence, for every $t\neq0,$\boxeq{
\[
\varphi_{W}\left(t\right)=-\dfrac{\left(\text{e}^{\text{i}t}-1\right)^{2}}{t^{2}}\text{e}^{\lambda\left(\text{e}^{\text{i}t}-1\right)}.
\]
}

The random variables $T,X$ and $N$ all admit expectations; hence,
the same is true for $W.$ The characteristic function $\varphi_{W}$
of $W$ is thus differentiable, and\boxeq{
\[
\varphi^{\prime}_{W}\left(0\right)=\text{i}\mathbb{E}\left(W\right)=\text{i}\mathbb{E}\left(T+X+N\right)=\text{i}\left(1+\lambda\right).
\]
}

\end{solution}

\begin{solution}{}{solexercise13.4}

\textbf{1. Moments of every order of $X$}

The characteristic function $\varphi_{X}$ of $X$ admits a Taylor
expansion of arbitrary order. For every $n\in\mathbb{N}$ and every
$t\in\mathbb{R},$ 
\[
\varphi_{X}\left(t\right)=\text{e}^{-\frac{t^{2}}{2}}=\sum^{n}_{k=0}\dfrac{\left(-1\right)^{k}t^{2k}}{2^{k}k!}+o\left(t^{2n}\right).
\]

Thus, the random variable $X$ admits moments of every order---this
could also be seen directly. They are given, for every $k\in\mathbb{N},$
by\boxeq{
\[
\mathbb{E}\left(X^{2k+1}\right)=0\,\,\,\,\text{and}\,\,\,\,\mathbb{E}\left(X^{2k}\right)=\dfrac{\left(2k\right)!}{2^{k}k!}
\]
}

Moreover, we note that $\varphi_{X}$ is analytic on $\mathbb{R}.$

\textbf{2. Moments of every order of $U=XY.$ Characteristic function
of $U$}

The characteristic function of $U$ admits a Taylor expansion of every
order $n,$ given for every real number $t,$ by
\[
\varphi_{U}\left(t\right)=1+\sum^{n}_{k=1}\dfrac{\left(\text{i}t\right)^{k}}{k!}\mathbb{E}\left(U^{k}\right)+R_{n}\left(t\right),
\]
 where the remainder is defined by
\[
R_{n}\left(t\right)=\intop^{t}_{0}\dfrac{\left(t-u\right)^{n}}{n!}\varphi^{\left(n+1\right)}_{U}\left(u\right)\text{d}u.
\]
Since the random variables $X$ and $Y$ are independent, $U$ admits
moments of every order given by
\[
\mathbb{E}\left(U^{n}\right)=\mathbb{E}\left(X^{n}\right)\mathbb{E}\left(Y^{n}\right).
\]
Hence, for every $k\in\mathbb{N},$\boxeq{
\[
\mathbb{E}\left(U^{2k+1}\right)=0\,\,\,\,\text{and}\,\,\,\,\mathbb{E}\left(U^{2k}\right)=\left(\dfrac{\left(2k\right)!}{2^{k}k!}\right)^{2}.
\]
}

It remains to show that the Taylor remainder $R_{n}\left(t\right)$
tends to 0. By the relation $\refpar{eq:k-th_deriv_charact_funct},$
\[
\left|R_{n}\left(t\right)\right|\leqslant\dfrac{\left|t\right|^{n+1}\mathbb{E}\left(\left|U\right|^{n+1}\right)}{\left(n+1\right)!}.
\]

By the Schwarz inequality, 
\[
\mathbb{E}\left(\left|U\right|^{2k+1}\right)\leqslant\left[\mathbb{E}\left(\left|U\right|^{2k}\right)\right]^{\frac{1}{2}}\left[\mathbb{E}\left(\left|U\right|^{2k+2}\right)\right]^{\frac{1}{2}},
\]
Taking into account the values of these moments, we obtain 
\[
\mathbb{E}\left(\left|U\right|^{2k+1}\right)\leqslant\left(2k+1\right)\mathbb{E}\left(U^{2k}\right),
\]
which yields the upper bound
\[
\left|R_{2k}\left(t\right)\right|\leqslant\dfrac{\left|t\right|^{2k+1}\mathbb{E}\left(U^{2k}\right)}{\left(2k\right)!}=\left|t\right|^{2k+1}\dfrac{\left(2k\right)!}{2^{2k}\left(k!\right)^{2}}.
\]
Moreover, we also have the upper bound
\[
\left|R_{2k-1}\left(t\right)\right|\leqslant\left|t\right|^{2k}\dfrac{\left(2k\right)!}{2^{2k}\left(k!\right)^{2}}.
\]
Therefore, the Taylor remainder $R_{n}\left(t\right)$ tends to $0$
as soon as $\left|t\right|<1$ and
\[
\varphi_{U}\left(t\right)=\sum^{+\infty}_{k=0}\left(-1\right)^{k}\dfrac{\left(2k\right)!}{2^{2k}\left(k!\right)^{2}}t^{2k}.
\]
Hence\boxeq{
\[
\varphi_{U}\left(t\right)=\dfrac{1}{\sqrt{1+t^{2}}}.
\]
}

By the principle of analytical continuation, this formula holds for
every real number.

\textbf{3. Characteristic function of $\left(U,V\right).$ Characteristic
function of $\dfrac{U+V}{\sqrt{2}}.$ Comparison of the laws of $\dfrac{U+V}{\sqrt{2}}$
and $U$}

Let $\left(a,b\right)\in\mathbb{R}^{2}$ be arbitrary. Since the random
variables $X,\,Y$ and $Z$ are independent, the transfer and Fubini
theorems---the integrand being bounded---yield the following expression
for the characteristic function of $\left(U,V\right),$
\begin{align*}
\varphi_{\left(U,V\right)}\left(a,b\right) & =\intop_{\mathbb{R}^{2}}\left[\intop_{\mathbb{R}}\text{e}^{\text{i}\left(ax+bz\right)y}\text{d}P_{Y}\left(y\right)\right]\text{d}P_{X}\otimes\text{d}P_{Z}\left(x,z\right)\\
 & =\intop_{\mathbb{R}^{2}}\varphi_{Y}\left(ax+bz\right)\text{d}P_{X}\otimes\text{d}P_{Z}\left(x,z\right).
\end{align*}
That is
\begin{align*}
\varphi_{\left(U,V\right)}\left(a,b\right) & =\dfrac{1}{2\pi}\intop_{\mathbb{R}^{2}}\text{e}^{-\frac{\left(ax+bz\right)^{2}}{2}}\text{e}^{-\frac{x^{2}+z^{2}}{2}}\text{d}x\text{d}z\\
 & =\dfrac{1}{2\pi}\intop_{\mathbb{R}^{2}}\text{e}^{-\frac{\left(1+a^{2}\right)x^{2}+2abxz+\left(1+b^{2}\right)z^{2}}{2}}\text{d}x\text{d}z.
\end{align*}
Applying the Fubini theorem
\[
\varphi_{\left(U,V\right)}\left(a,b\right)=\dfrac{1}{2\pi}\intop_{\mathbb{R}}\text{e}^{-\frac{\left(1+a^{2}\right)x^{2}}{2}}\left[\intop_{\mathbb{R}}\text{e}^{-\frac{\left(1+b^{2}\right)z^{2}+2abxz}{2}}\text{d}z\right]\text{d}x.
\]

A straightforward computation shows that, for every $\left(a,b\right)\in\mathbb{R}^{2},$
\[
\intop_{\mathbb{R}}\text{e}^{-\frac{a^{2}u^{2}+2bu}{2}}\text{d}u=\text{e}^{\frac{b^{2}}{2a}}\intop_{\mathbb{R}}\text{e}^{-\frac{v^{2}}{2}}\dfrac{1}{\sqrt{a}}\text{d}v=\sqrt{\dfrac{2\pi}{a}}\text{e}^{\frac{b^{2}}{2a}}.
\]
It follows that
\[
\varphi_{\left(U,V\right)}\left(a,b\right)=\dfrac{1}{2\pi}\intop_{\mathbb{R}}\text{e}^{-\frac{\left(1+a^{2}\right)x^{2}}{2}}\sqrt{\dfrac{2\pi}{1+b^{2}}}\text{e}^{\frac{a^{2}b^{2}x^{2}}{2\left(1+b^{2}\right)}}\text{d}x,
\]
which yields, after simplification,\boxeq{
\[
\varphi_{\left(U,V\right)}\left(a,b\right)=\dfrac{1}{\sqrt{1+a^{2}+b^{2}}}.
\]
}

The characteristic function of $\dfrac{U+V}{\sqrt{2}}$ is therefore
given, for every $t\in\mathbb{R},$ by\boxeq{
\[
\varphi_{\frac{U+V}{\sqrt{2}}\left(t\right)}=\varphi_{\left(U,V\right)}\left(\dfrac{t}{\sqrt{2}},\dfrac{t}{\sqrt{2}}\right)=\dfrac{1}{\sqrt{1+t^{2}}}.
\]
}

The Fourier transform injectivity theorem ensures that \textbf{the
random variables $\dfrac{U+V}{\sqrt{2}}$ and $U$ have the same law}.

\textbf{4. Independence of $U$ and $V.$ Correlation of $U$ and
$V$}

We have
\[
\varphi_{\left(U,V\right)}\left(a,b\right)\neq\varphi_{\left(U,V\right)}\left(a,0\right)\varphi_{\left(U,V\right)}\left(0,b\right)=\varphi_{U}\left(a\right)\varphi_{V}\left(b\right),
\]
which proves that \textbf{the random variables $U$ and $V$ are not
independent}. 

Nevertheless, since the random variables $X,$ $Y$ and $Z$ are independent,
\[
\mathbb{E}\left(UV\right)=\mathbb{E}\left(XY^{2}Z\right)=\mathbb{E}\left(X\right)\mathbb{E}\left(Y^{2}\right)\mathbb{E}\left(Z\right)=0.
\]

As $U$ and $V$ are centered, it follows that \textbf{the random
variables $U$ and $V$ are uncorrelated}.

\end{solution}

\begin{solution}{}{solexercise13.5}

The condition is necessary, since when $X$ and $Y$ are independent,
the random variables $X^{k}$ and $Y^{l}$ are also independent, which
yields the relation $\refpar{eq:bounded_real_val_var_indep_criteria}.$

Conversely, assume the relation $\refpar{eq:bounded_real_val_var_indep_criteria}$
holds. The characteristic function of $\left(X,Y\right)$ is given,
for every $\left(u,v\right)\in\mathbb{R}^{2},$ by
\[
\varphi_{\left(X,Y\right)}\left(u,v\right)=\mathbb{E}\left(\text{e}^{\text{i}uX}\text{e}^{\text{i}vY}\right)=\mathbb{E}\left[\left(\sum^{+\infty}_{k=0}\dfrac{\left(\text{i}u\right)^{k}X^{k}}{k!}\right)\left(\sum^{+\infty}_{l=0}\dfrac{\left(\text{i}v\right)^{l}Y^{l}}{l!}\right)\right].
\]
Let $C$ be an upper bound for $\left|X\right|$ and $\left|Y\right|.$
Then
\[
\sum_{\left(k,l\right)\in\mathbb{N}^{2}}\dfrac{\left|u\right|^{k}\left|v\right|^{l}C^{k+l}}{k!l!}=\text{e}^{\left|u\right|C}\text{e}^{\left|v\right|C}<+\infty.
\]
It follows that the family 
\[
\left\{ \dfrac{\left(\text{i}u\right)^{k}X^{k}}{k!}\cdot\dfrac{\left(\text{i}v\right)^{l}Y^{l}}{l!}:\,\left(k,l\right)\in\mathbb{N}^{2}\right\} 
\]
is summable. Since this family is countable, an application of the
dominated convergence theorem yields
\[
\varphi_{\left(X,Y\right)}\left(u,v\right)=\sum_{\left(k,l\right)\in\mathbb{N}^{2}}\text{i}^{k+l}\dfrac{u^{k}v^{l}}{k!l!}\mathbb{E}\left(X^{k}Y^{l}\right).
\]
Taking the hypothesis into account, it follows that
\[
\varphi_{\left(X,Y\right)}\left(u,v\right)=\sum_{\left(k,l\right)\in\mathbb{N}^{2}}\text{i}^{k+l}\dfrac{u^{k}v^{l}}{k!l!}\mathbb{E}\left(X^{k}\right)\mathbb{E}\left(Y^{l}\right).
\]
This last family is still summable, and therefore, by application
of the Fubini property, followed again by the dominated convergence
theorem,
\begin{align*}
\varphi_{\left(X,Y\right)}\left(u,v\right) & =\left[\sum_{k\in\mathbb{N}}\text{i}^{k}\dfrac{u^{k}}{k!}\mathbb{E}\left(X^{k}\right)\right]\left[\sum_{l\in\mathbb{N}}\text{i}^{l}\dfrac{v^{l}}{l!}\mathbb{E}\left(Y^{l}\right)\right]\\
 & =\mathbb{E}\left(\sum_{k\in\mathbb{N}}\text{i}^{k}\dfrac{u^{k}}{k!}X^{k}\right)\mathbb{E}\left(\sum_{l\in\mathbb{N}}\text{i}^{l}\dfrac{v^{l}}{l!}Y^{l}\right).
\end{align*}
That is,
\[
\varphi_{\left(X,Y\right)}\left(u,v\right)=\varphi_{X}\left(u\right)\varphi_{Y}\left(v\right).
\]

This proves the independence of $X$ and $Y.$

\end{solution}

\begin{solution}{}{solexercise13.6}

\textbf{1. Proof that $f$ is a probability density on $\mathbb{R}^{2}$ }

The function $f$ is nonnegative. Indeed, for every $\left(x,y\right)\in C,$
\[
-1\leqslant-y^{2}\leqslant x^{2}-y^{2}\leqslant x^{2}\leqslant1
\]
and thus
\[
\left|xy\left(x^{2}-y^{2}\right)\right|\leqslant1.
\]
It follows that
\[
1+xy\left(x^{2}-y^{2}\right)\geqslant0.
\]
The function $f$ is measurable. Moreover, it is a probability density,
using the symmetry arguments, we obtain
\[
\intop_{\mathbb{R}^{2}}f\left(x,y\right)\text{d}x\text{d}y=\dfrac{1}{4}\intop_{C}\left[1+xy\left(x^{2}-y^{2}\right)\right]\text{d}x\text{d}y=1.
\]

\textbf{2. Densities of $X,\,Y,\,Z.$ Laws of $X$ and $Y$}

The marginal random variable $X$ admits the density $f_{X}$ given,
for every real number $x,$ by
\[
f_{X}\left(x\right)=\intop_{\mathbb{R}}f\left(x,y\right)\text{d}y=\dfrac{1}{4}\boldsymbol{1}_{\left[-1,1\right]}\left(x\right)\intop^{+1}_{-1}\left[1+xy\left(x^{2}-y^{2}\right)\right]\text{d}y.
\]

Hence,\boxeq{
\[
f_{X}=\dfrac{1}{2}\boldsymbol{1}_{\left[-1,1\right]}.
\]
}

\textbf{The random variable $X$ follows the uniform law on the interval
$\left[-1,1\right].$} By symmetry, the same holds for the random
variable $Y.$ Note that $f$ is not the direct product of the densities
$X$ and $Y.$ Hence, \textbf{the random variables $X$ and $Y$ are
not independent}.

Now, let study the law of $Z.$ The random variable $\left(X+Y,X\right)$
is obtained from $\left(X,Y\right)$ via a linear diffeomorphism whose
Jacobian has absolute value 1. Consequently, the random variable $\left(Z,X\right)$
admits the density $f_{\left(Z,X\right)}$ given, for every $\left(z,t\right)\in\mathbb{R}^{2},$
by
\[
f_{\left(Z,X\right)}\left(z,t\right)=f_{\left(X,Y\right)}\left(t,z-t\right).
\]
The random variable $Z$ therefore admits the density $f_{Z}$ given
for every $z\in\mathbb{R}$ by
\begin{align*}
f_{Z}\left(z\right) & =\intop^{+\infty}_{-\infty}f_{\left(X,Y\right)}\left(t,z-t\right)\text{d}t\\
 & =\dfrac{1}{4}\intop^{+\infty}_{-\infty}\left[1+t\left(z-t\right)\left(t^{2}-\left(z-t\right)^{2}\right)\boldsymbol{1}_{C}\left(t,z-t\right)\right]\text{d}t.
\end{align*}
Thus, 
\[
f_{Z}\left(z\right)=\dfrac{1}{4}\boldsymbol{1}_{\left[-2,2\right]}\left(z\right)\intop^{\min\left(z+1,1\right)}_{\max\left(z-1,-1\right)}\left[1+t\left(z-t\right)\left(t^{2}-\left(z-t\right)^{2}\right)\right]\text{d}t.
\]

If $0<z<2,$ then
\[
f_{Z}\left(z\right)=\dfrac{1}{4}\intop^{1}_{z-1}\left[1+zt\left(z-t\right)\left(2t-z\right)\right]\text{d}t.
\]

Making the change of variables $u=2t-z,$ we obtain
\begin{align*}
f_{Z}\left(z\right) & =\dfrac{1}{8}\intop^{2-z}_{-\left(2-z\right)}\left[1+zu\dfrac{u+z}{2}\dfrac{z-u}{2}\right]\text{d}u\\
 & =\dfrac{1}{8}\intop^{2-z}_{-\left(2-z\right)}\dfrac{z^{2}-u^{2}}{4}zu\text{d}u+\dfrac{1}{8}\intop^{2-z}_{-\left(2-z\right)}\text{d}u.
\end{align*}
Since the first integral is zero, we obtain
\[
f_{Z}\left(z\right)=\dfrac{2-z}{4}.
\]
The random variables $\left(X,Y\right)$ and $\left(-X,-Y\right)$
have the same law, and therefore so do the random variables $Z$ and
$-Z.$ It follows that $f_{Z}$ is even. Thus,\boxeq{
\[
\forall z\in\mathbb{R},\,\,\,\,f_{Z}\left(z\right)=\boldsymbol{1}_{\left[-2,2\right]}\left(z\right)\dfrac{2-\left|z\right|}{4}.
\]
}

\textbf{3. Case of an even density. Characteristic functions of $X,\,Y$
and $Z$}

If a random variable $U$ has an even density, then its characteristic
function $\varphi_{U}$ verifies, for every $t\in\mathbb{R},$
\[
\varphi_{U}\left(t\right)=\intop_{\mathbb{R}}\left[\cos\left(tu\right)+\text{i}\sin\left(tu\right)\right]f_{U}\left(u\right)\text{d}u.
\]
Hence,
\[
\varphi_{U}\left(t\right)=2\intop^{+\infty}_{0}\cos\left(tu\right)f_{U}\left(u\right)\text{d}u.
\]
The characteristic function of $X$ is therefore given by
\[
\forall t\in\mathbb{R},\,\,\,\,\varphi_{X}\left(t\right)=2\intop^{1}_{0}\cos\left(tu\right)\dfrac{1}{2}\text{d}u.
\]
Hence, \boxeq{
\[
\forall y\neq0,\,\,\,\,\varphi_{X}\left(t\right)=\varphi_{Y}\left(t\right)=\dfrac{\sin t}{t}\,\,\,\,\text{and}\,\,\,\,\varphi_{X}\left(0\right)=\varphi_{Y}\left(0\right)=1.
\]
}

Similarly, the characteristic function of $Z$ is given by
\[
\forall t\neq0,\,\,\,\,\varphi_{Z}\left(t\right)=2\intop^{2}_{0}\cos\left(tz\right)\dfrac{2-z}{4}\text{d}z=\dfrac{\sin2t}{t}-\dfrac{1}{2}\intop^{2}_{0}z\cos\left(tz\right)\text{d}z.
\]
After an integration by parts, we obtain
\[
\forall t\neq0,\,\,\,\,\varphi_{Z}\left(t\right)=\dfrac{1-\cos2t}{2t^{2}}=\dfrac{\sin^{2}t}{t^{2}}.
\]
Thus,\boxeq{
\[
\varphi_{Z}=\varphi_{X}\varphi_{Y}.
\]
}

By the Fourier transform injectivity property, this implies that
\[
P_{X+Y}=P_{X}\ast P_{Y}.
\]

\textbf{4. Correlation coefficient of $X$ and Y}

We have 
\[
\mathbb{E}\left(X\right)=\mathbb{E}\left(Y\right)=0.
\]
Moreover,
\[
\mathbb{E}\left(XY\right)=\dfrac{1}{4}\intop_{C}xy\left[1+xy\left(x^{2}-y^{2}\right)\right]\text{d}x\text{d}y,
\]
thus,
\[
\mathbb{E}\left(XY\right)=\dfrac{1}{4}\left[\intop_{C}x^{4}y^{2}\text{d}x\text{d}y-\intop_{C}x^{2}y^{4}\text{d}x\text{d}y\right]=0.
\]
\textbf{The correlation coefficient of $X$ and $Y$ is therefore
equal to zero.}

\textbf{5. Remarks}

In summary, we obtain an example of a random variable $\left(X,Y\right)$
with a non-uniform law on $C,$ whose marginals are of uniform law
but not independent, while still being uncorrelated. Nonetheless,
these marginals satisfy 
\[
\varphi_{X+Y}=\varphi_{X}\varphi_{Y},
\]
and therefore 
\[
P_{X+Y}=P_{X}\ast P_{Y}.
\]

\end{solution}

\begin{solution}{}{solexercise13.7}

\textbf{1. Proof that $\forall t\in\mathbb{R},\,\,\,\,\widehat{\mu_{a}}\left(t\right)=\widehat{\mu_{1}}\left(at\right).$
If $Z$ follows a Cauchy law with parameter 1, then $aZ$ follows
a Cauchy law with parameter $a$}

For every real number $t,$ making a change of variables
\[
\widehat{\mu_{a}}\left(t\right)=\intop_{\mathbb{R}}\text{e}^{\text{i}tx}\dfrac{a}{\pi\left(a^{2}+x^{2}\right)}\text{d}x=\intop_{\mathbb{R}}\text{e}^{\text{i}at\frac{x}{a}}\dfrac{1}{\pi\left(1+\left(\dfrac{x}{a}\right)^{2}\right)}\dfrac{\text{d}x}{a}.
\]
This shows that\boxeq{
\[
\widehat{\mu_{a}}\left(t\right)=\widehat{\mu_{1}}\left(at\right).
\]
}

The characteristic function $\varphi_{aZ}$ of $aZ$ is then given,
for every real number $t,$ by 
\[
\varphi_{aZ}\left(t\right)=\varphi_{Z}\left(at\right)=\widehat{\mu_{1}}\left(at\right)=\widehat{\mu_{a}}\left(t\right).
\]
By the injectivity of the Fourier transform, the law of $aZ$ is thus
the law $\mu_{a}.$ 

\textbf{2. Characteristic function $\left(X,Y\right).$ $X$ and $Y$
are not independent}

We have
\[
\left(\begin{array}{c}
X\\
Y
\end{array}\right)=A\left(\begin{array}{c}
U\\
V
\end{array}\right),
\]
where $A$ is the matrix
\[
A=\left(\begin{array}{cc}
a & b\\
c & d
\end{array}\right).
\]
For every $\left(\alpha,\beta\right)\in\mathbb{R}^{2},$
\[
A^{\ast}\left(\begin{array}{c}
\alpha\\
\beta
\end{array}\right)=\left(\begin{array}{c}
a\alpha+c\beta\\
b\alpha+d\beta
\end{array}\right).
\]
It follows that, for every $\left(\alpha,\beta\right)\in\mathbb{R}^{2},$
\[
\varphi_{\left(X,Y\right)}\left(\alpha,\beta\right)=\varphi_{\left(U,V\right)}\left(a\alpha+c\beta,b\alpha+d\beta\right).
\]

Since the random variables $U$ and $V$ are independent, we obtain
\[
\varphi_{\left(X,Y\right)}\left(\alpha,\beta\right)=\varphi_{U}\left(a\alpha+c\beta\right)\varphi_{V}\left(b\alpha+d\beta\right)=\text{e}^{-\left[\left|a\alpha+c\beta\right|+\left|b\alpha+d\beta\right|\right]}.
\]

Hence,
\[
\varphi_{X}\left(\alpha\right)=\varphi_{\left(X,Y\right)}\left(\alpha,0\right)=\text{e}^{-\left(a+b\right)\left|\alpha\right|}.
\]
Thus,\boxeq{
\[
\varphi_{\left(X,Y\right)}\left(\alpha,\beta\right)\neq\varphi_{X}\left(\alpha\right)\varphi_{Y}\left(\beta\right),
\]
}which proves that \textbf{the random variables $X$ and $Y$ are
not independent}.

\textbf{3. Characteristic function of $X+Y.$ Proof of $P_{X+Y}=P_{X}\ast P_{Y}$}

Nevertheless, for every real number $\alpha,$
\[
\varphi_{X+Y}\left(\alpha\right)=\varphi_{\left(X,Y\right)}\left(\alpha,\alpha\right)=\text{e}^{-\left(a+b+c+d\right)\left|\alpha\right|}=\varphi_{X}\left(\alpha\right)\varphi_{Y}\left(\alpha\right).
\]
By the injectivity of the Fourier transform, this implies\boxeq{
\[
P_{X+Y}=P_{X}\ast P_{Y}.
\]
}

\end{solution}

\begin{solution}{}{solexercise13.8}

\textbf{1. Proof that $\exists t_{0}\neq0:$ $\left|\varphi_{X}\left(t_{0}\right)\right|=1,$
implies $\exists a\in\mathbb{R},$$P_{X}\left(a+\mathbb{Z}\dfrac{2\pi}{t_{0}}\right)=1$}

Let $t_{0}$ be a nonzero real number such that $\varphi_{X}\left(t_{0}\right)=\text{e}^{\text{i}t_{0}a},$
that is, such that $\mathbb{E}\left(\text{e}^{\text{i}t_{0}\left(X-a\right)}\right)=1.$
Then,
\[
\mathbb{E}\left(1-\text{e}^{\text{i}t_{0}\left(X-a\right)}\right)=0
\]
 and, taking the real part,
\[
\mathbb{E}\left(1-\cos\left(t_{0}\left(X-a\right)\right)\right)=0.
\]
Since the integrand is nonnegative, it follows that $P-$almost surely
$\cos\left(t_{0}\left(X-a\right)\right)=1,$ and 
\[
P_{X}\left(a+\mathbb{Z}\dfrac{2\pi}{t_{0}}\right)=1.
\]

\textbf{2. Proof that $\exists t_{1},t_{2}\in\mathbb{R}^{\ast}:\,\dfrac{t_{1}}{t_{2}}\in\mathbb{R}\backslash\mathbb{Q}$
and $\left|\varphi_{X}\left(t_{1}\right)\right|=\left|\varphi_{X}\left(t_{2}\right)\right|=1,$
implies that $X$ is degenerate}

Assume that there exist two nonzero real numbers $t_{1}$ and $t_{2}$
such that $\dfrac{t_{1}}{t_{2}}$ is irrationnal and 
\[
\left|\varphi_{X}\left(t_{1}\right)\right|=\left|\varphi_{X}\left(t_{2}\right)\right|=1.
\]
By the first question, there exist two real numbers $a$ and $b$
such that
\[
P_{X}\left(a+\mathbb{Z}\dfrac{2\pi}{t_{1}}\right)=P_{X}\left(a+\mathbb{Z}\dfrac{2\pi}{t_{2}}\right)=1.
\]

If the random variable $X$ were non-degenerated, the sets $a+\mathbb{Z}\dfrac{2\pi}{t_{1}}$
and $a+\mathbb{Z}\dfrac{2\pi}{t_{2}}$ would have at least two distinct
common points, such that there would exist integers $k\neq k^{\prime}$
and $l\neq l^{\prime}$ such that $a+k\dfrac{2\pi}{t_{1}}=b+l\dfrac{2\pi}{t_{2}}$
and $a+k^{\prime}\dfrac{2\pi}{t_{1}}=b+l^{\prime}\dfrac{2\pi}{t_{2}}.$

Subtracting there equalities would yield
\[
a-b=l\dfrac{2\pi}{t_{2}}-k\dfrac{2\pi}{t_{1}}=l^{\prime}\dfrac{2\pi}{t_{2}}-k^{\prime}\dfrac{2\pi}{t_{1}},
\]
and thus
\[
\dfrac{l-l^{\prime}}{k-k^{\prime}}=\dfrac{t_{2}}{t_{1}},
\]
which is impossible if $\dfrac{t_{1}}{t_{2}}$ is irrational.\textbf{
Therefore, the random variable $X$ must be degenerated.}

\textbf{3. Proof $X$ is degenerated if and only if $\left|\varphi_{X_{0}}\right|=1.$}

If the random variable $X$ is degenerated and is equal $P-$almost
surely to $a,$ then, for every real number $t,$ $\varphi_{X}\left(t\right)=\text{e}^{\text{i}ta}$
and hence $\left|\varphi_{X}\right|=1.$ The converse follows from
the previous question.

\end{solution}

\begin{solution}{}{solexercise13.9}

\textbf{1. Computation of $\mathbb{E}^{\mathscr{F}_{n-1}}\left(\text{e}^{\text{i}tU_{n-1}U_{n}}\right)$}

As the random variables $\left(U_{0},U_{1},\cdots,U_{n-1}\right)$
and $U_{n}$ are independent, we have 
\[
\mathbb{E}^{\mathscr{F}_{n-1}}\left(\text{e}^{\text{i}tU_{n-1}U_{n}}\right)=f\left(U_{0},U_{1},\cdots,U_{n-1}\right),
\]
where the function $f$ is defined on $\mathbb{R}^{n}$ by
\[
\forall\left(u_{0},u_{1},\cdots,u_{n-1}\right)\in\mathbb{R}^{n},\,\,\,\,f\left(u_{0},u_{1},\cdots,u_{n-1}\right)=\mathbb{E}\left(\text{e}^{\text{i}tu_{n-1}U_{n}}\right).
\]
Thus,
\[
f\left(u_{0},u_{1},\cdots,u_{n-1}\right)=\dfrac{1}{2}\left[\text{e}^{\text{i}tu_{n-1}}+\text{e}^{-\text{i}tu_{n-1}}\right]=\cos\left(tu_{n-1}\right).
\]
It follows that \boxeq{
\[
\mathbb{E}^{\mathscr{F}_{n-1}}\left(\text{e}^{\text{i}tU_{n-1}U_{n}}\right)=\cos\left(tU_{n-1}\right).
\]
}

\textbf{2. Computation of $\mathbb{E}^{\mathscr{F}_{n-l}}\left(\text{e}^{\text{i}tY_{n}}\right)$}

Since $Y_{n-1}$ is $\mathscr{F}_{n-1}-$measurable, 
\begin{align*}
\mathbb{E}^{\mathscr{F}_{n-1}}\left(\text{e}^{\text{i}tY_{n}}\right) & =\text{e}^{\text{i}tY_{n}-1}\mathbb{E}^{\mathscr{F}_{n-1}}\left(\text{e}^{\text{i}tU_{n-1}U_{n}}\right)\\
 & =\text{e}^{\text{i}tY_{n}-1}\cos\left(tU_{n-1}\right).
\end{align*}

We now compute by induction on $l,$ $\mathbb{E}^{\mathscr{F}_{n-l}}\left(\text{e}^{\text{i}tY_{n}}\right)$.
Computations for $l=2$ and $l=3$ suggest the following induction
hypothesis at the order $l,$
\[
\left(\text{IH}_{l}\right)\,\,\,\,\mathbb{E}^{\mathscr{F}_{n-l}}\left(\text{e}^{\text{i}tY_{n}}\right)=\text{e}^{\text{i}tY_{n-l}}\cos\left(tU_{n-l}\right)\cos^{l-1}\left(t\right).
\]

Assume that $\left(\text{IH}_{l}\right)$ holds at rank $l,$ we now
prove that it holds at rank $l+1.$

Since $\mathscr{F}_{n-l-1}\subset\mathscr{F}_{n-l},$
\[
\mathbb{E}^{\mathscr{F}_{n-\left(l+1\right)}}\left(\text{e}^{\text{i}tY_{n}}\right)=\mathbb{E}^{\mathscr{F}_{n-l-1}}\left(\mathbb{E}^{\mathscr{F}_{n-l}}\left(\text{e}^{\text{i}tY_{n}}\right)\right).
\]
By the induction hypothesis $\left(\text{IH}_{l}\right)$,
\[
\mathbb{E}^{\mathscr{F}_{n-\left(l+1\right)}}\left(\text{e}^{\text{i}tY_{n}}\right)=\mathbb{E}^{\mathscr{F}_{n-l-1}}\left(\text{e}^{\text{i}tY_{n-l}}\cos\left(tU_{n-l}\right)\cos^{l-1}\left(t\right)\right).
\]

Also, since $Y_{n-\left(l+1\right)}$ is $\mathscr{F}_{n-l-1}-$measurable,
\[
\mathbb{E}^{\mathscr{F}_{n-\left(l+1\right)}}\left(\text{e}^{\text{i}tY_{n}}\right)=\text{e}^{\text{i}tY_{n-\left(l+1\right)}}\mathbb{E}^{\mathscr{F}_{n-l-1}}\left(\text{e}^{\text{i}tU_{n-\left(l+1\right)}U_{n-l}}\cos\left(tU_{n-l}\right)\cos^{l-1}\left(t\right)\right).
\]
A computation similar to that of the previous question yields
\begin{align*}
\mathbb{E}^{\mathscr{F}_{n-\left(l+1\right)}}\left(\text{e}^{\text{i}tY_{n}}\right) & =\text{e}^{\text{i}tY_{n-\left(l+1\right)}}\\
 & \times\dfrac{1}{2}\left[\text{e}^{\text{i}tU_{n-\left(l+1\right)}}\cos\left(t\right)+\text{e}^{-\text{i}tU_{n-\left(l+1\right)}}\cos\left(-t\right)\right]\cos^{l-1}\left(t\right).
\end{align*}
Hence,
\[
\mathbb{E}^{\mathscr{F}_{n-\left(l+1\right)}}\left(\text{e}^{\text{i}tY_{n}}\right)=\text{e}^{\text{i}tY_{n-\left(l+1\right)}}\cos\left(tU_{n-\left(l+1\right)}\right)\cos^{l}\left(t\right).
\]

Thus, the formula holds at rank $l+1.$

This proves the formula for every $l$ such that $1\leqslant l\leqslant n.$

\textbf{3. Characteristic function $\varphi_{Y_{n}}$ of $Y_{n}$}

The characteristic function of $Y_{n}$ is given, for every real number
$t,$ by
\[
\varphi_{Y_{n}}\left(t\right)=\mathbb{E}\left(\mathbb{E}^{\mathscr{F}_{0}}\left(\text{e}^{\text{i}tY_{n}}\right)\right)=\mathbb{E}\left(\text{e}^{\text{i}tY_{0}}\cos\left(tU_{0}\right)\right)\cos^{n-1}\left(t\right).
\]
Since $Y_{0}=0,$ this yields
\[
\varphi_{Y_{n}}\left(t\right)=\dfrac{1}{2}\left[\cos\left(t\right)+\cos\left(-t\right)\right]\cos^{n-1}\left(t\right).
\]
Thus,\boxeq{
\[
\varphi_{Y_{n}}\left(t\right)=\cos^{n}\left(t\right).
\]
}

\textbf{4. Law of $Y_{n}$}

The Fourier transform in $t$ of the probability $\dfrac{\delta_{-1}+\delta_{1}}{2}$
is $\cos\left(t\right).$ By the injectivity of the Fourier transform,
it follows that $Y_{n}$ has the same law as the sum of $n$ independent
random variables with law $\dfrac{\delta_{-1}+\delta_{1}}{2},$ that
is, the same law as $\sum^{n}_{j=1}U_{j}.$ 

Thus, for every integer $k$ such that $-n\leqslant k\leqslant n,$
using the independence of the $U_{j},$ we obtain
\[
P\left(Y_{n}=k\right)=P\left[\biguplus_{\substack{J\subset\left\{ 1,\cdots,n\right\} \\
\left|J\right|=\frac{k+n}{2}
}
}\left\{ \bigcap_{j\in J}\left(U_{j}=1\right)\cap\bigcap_{j\in J^{c}}\left(U_{j}=-1\right)\right\} \right].
\]
Hence, 
\[
P\left(Y_{n}=k\right)=\sum_{\substack{J\subset\left\{ 1,\cdots,n\right\} \\
\left|J\right|=\frac{k+n}{2}
}
}\dfrac{1}{2^{n}}.
\]
Thus, for every integer $k$ with $-n\leqslant k\leqslant n,$\boxeq{
\[
P\left(Y_{n}=k\right)=\begin{cases}
\binom{n}{\frac{k+n}{2}}, & \text{if }k+n\text{ is even,}\\
0, & \text{if }k+n\text{ is odd.}
\end{cases}
\]
}

\textbf{5. Computation of $\varphi_{\frac{Y_{n}}{n}}\left(t\right)$}

We have
\[
\varphi_{\frac{Y_{n}}{n}\left(t\right)}=\varphi_{Y_{n}}\left(\dfrac{t}{n}\right)=\cos^{n}\left(\dfrac{t}{n}\right).
\]
Thus,
\[
\ln\left(\varphi_{\frac{Y_{n}}{n}\left(t\right)}\right)=n\ln\left[1-\dfrac{t^{2}}{2n^{2}}+o\left(\dfrac{1}{n^{2}}\right)\right],
\]
which yields
\[
\varphi_{\frac{Y_{n}}{n}\left(t\right)}=\text{e}^{-\frac{t^{2}}{2n}+o\left(\frac{1}{n}\right)}.
\]
Therefore\boxeq{
\[
\lim_{n\to+\infty}\varphi_{\frac{Y_{n}}{n}\left(t\right)}=1.
\]
}

\begin{remark}{}{}

We have thus shown that the sequence with general term $\frac{Y_{n}}{n}$
converges in law to $0$---the notion of convergence in law is studied
in the next chapter. 

\end{remark}

\end{solution}

\begin{solution}{}{solexercise13.10}

\textbf{1. Proof that $\forall x\in\mathbb{R},\,\,\,\,\lim_{T\to+\infty}I_{T}\left(x\right)=\dfrac{1}{2}\boldsymbol{1}_{\left\{ a,b\right\} }\left(x\right)+\boldsymbol{1}_{\left]a,b\right[}\left(x\right).$}

For every $x\notin\left\{ a,b\right\} ,$ a change of variables gives
\[
I_{T}\left(x\right)=\dfrac{1}{2\pi}\intop^{+T\left(x-a\right)}_{-T\left(x-a\right)}\dfrac{\sin\left(u\right)}{u}\text{d}u-\dfrac{1}{2\pi}\intop^{+T\left(x-b\right)}_{-T\left(x-b\right)}\dfrac{\sin\left(u\right)}{u}\text{d}u.
\]
Using the equality
\[
\lim_{\alpha\to+\infty}\intop^{\alpha}_{0}\dfrac{\sin\left(u\right)}{u}\text{d}u=\dfrac{\pi}{2}
\]
we obtain:
\begin{itemize}
\item If $x<a$ or $x>b,$ then $\lim_{T\to+\infty}I_{T}\left(x\right)=0.$
\item If $a<x<b,$ then $\lim_{T\to+\infty}I_{T}\left(x\right)=\dfrac{1}{2\pi}\left(\pi+\pi\right)=1.$
\item If $x=a,$ then
\begin{align*}
I_{T}\left(a\right) & =\dfrac{1}{2\pi}\intop^{+T}_{-T}\dfrac{1-\text{e}^{\text{i}t\left(a-b\right)}}{\text{i}t}\text{d}t\\
 & =\dfrac{1}{2\pi}\left[\intop^{+T}_{-T}\dfrac{1-\cos\left(t\left(b-a\right)\right)}{\text{i}t}\text{d}t+\intop^{+T}_{-T}\dfrac{\sin\left(t\left(b-a\right)\right)}{\text{i}}\text{d}t\right].
\end{align*}
By a change of variables and a parity argument,
\[
I_{T}\left(a\right)=\dfrac{1}{2\pi}\intop^{T\left(b-a\right)}_{-T\left(b-a\right)}\dfrac{\sin\left(u\right)}{u}\text{d}u,
\]
so that
\[
\lim_{T\to+\infty}I_{T}\left(a\right)=\dfrac{1}{2}.
\]
Similarly, 
\[
\lim_{T\to+\infty}I_{T}\left(b\right)=\dfrac{1}{2}.
\]
\end{itemize}
\textbf{2. Inversion formula $\refpar{eq:first_inversion_formula}.$}

By the finite increments inequality, the measurable function 
\[
\left(t,x\right)\mapsto\dfrac{\text{e}^{-\text{i}ta}-\text{e}^{-\text{i}tb}}{\text{i}t}\text{e}^{\text{i}tx}
\]
is bounded on $\left[-T,T\right]\times\mathbb{R}.$ Hence, by the
Fubini theorem
\[
\dfrac{1}{2\pi}\intop^{T}_{-T}\dfrac{\text{e}^{-\text{i}ta}-\text{e}^{-\text{i}tb}}{\text{i}t}\varphi\left(t\right)\text{d}t=\dfrac{1}{2\pi}\intop_{\mathbb{R}}I_{T}\left(x\right)\text{d}\mu\left(x\right).
\]
Moreover, the function $\alpha\mapsto\intop^{\alpha}_{0}\dfrac{\sin\left(u\right)}{u}\text{d}u$
is uniformly continuous on $\mathbb{R}$ and converges to $\dfrac{\pi}{2}$
as $\alpha$ tends to $+\infty.$ Therefore, there exists a real number
$M$ such that, for every $\left(x,T\right),$ 
\[
\left|I_{T}\left(x\right)\right|\leqslant M.
\]
The dominated convergence theorem then yields
\[
\lim_{T\to+\infty}\dfrac{1}{2\pi}\intop^{+T}_{-T}\dfrac{\text{e}^{-\text{i}ta}-\text{e}^{-\text{i}tb}}{\text{i}t}\varphi\left(t\right)\text{d}t=\dfrac{1}{2\pi}\intop_{\mathbb{R}}\lim_{T\to+\infty}I_{T}\left(x\right)\text{d}\mu\left(x\right),
\]
and, using the result of the previous question, we obtain the inversion
formula $\refpar{eq:first_inversion_formula}.$

\textbf{3. Proof that $\forall b\in\mathbb{R},\,\,\,\,\mu\left(\left\{ b\right\} \right)=\lim_{T\to+\infty}\dfrac{1}{2T}\intop^{+T}_{-T}\text{e}^{-\text{i}tb}\varphi\left(t\right)\text{d}t.$}

By the Fubini theorem,
\[
\intop^{+T}_{-T}\text{e}^{-\text{i}tb}\varphi\left(t\right)\text{d}t=\intop_{\mathbb{R}}\text{\ensuremath{\left[\intop^{+T}_{-T}\text{e}^{\text{i}t\left(x-b\right)}\text{d}t\right]}d}\mu\left(x\right).
\]
If $x\neq b,$ then
\[
\intop^{+T}_{-T}\text{e}^{\text{i}t\left(x-b\right)}\text{d}t=2\dfrac{\sin\left(T\left(x-b\right)\right)}{x-b},
\]
which yields
\[
\dfrac{1}{2T}\intop^{+T}_{-T}\text{e}^{-\text{i}tb}\varphi\left(t\right)\text{d}t=\mu\left(\left\{ b\right\} \right)+\intop_{\mathbb{R}\backslash\left\{ b\right\} }\dfrac{\sin\left(T\left(x-b\right)\right)}{T\left(x-b\right)}\text{d}\mu\left(x\right).
\]

Since the function $u\mapsto\dfrac{\sin u}{u}$---extended by 1 at
0---is bounded and tends to $0$ as $u$ tends to $\pm\infty,$ the
dominated convergence theorem implies that the integral on the right-hand
side tends to $0$ as $T$ tends to $+\infty.$ Hence,\boxeq{
\[
\lim_{T\to+\infty}\dfrac{1}{2T}\intop^{+T}_{-T}\text{e}^{-\text{i}tb}\varphi\left(t\right)\text{d}t=\mu\left(\left\{ b\right\} \right).
\]
}

\textbf{4. Proof that $\lim_{T\to+\infty}\dfrac{1}{2T}\intop^{+T}_{-T}\left|\varphi\left(t\right)\right|^{2}\text{d}t=\sum_{x\in\mathbb{R}}\mu\left(\left\{ x\right\} \right)^{2}.$ }

Let $X$ and $Y$ be two independent random variables with law $\mu.$
The characteristic function $\varphi_{X-Y}$ of the random variable
$X-Y$ is given by
\[
\forall t\in\mathbb{R},\,\,\,\,\varphi_{X-Y}\left(t\right)=\varphi_{X}\left(t\right)\varphi_{Y}\left(-t\right)=\left|\varphi\left(t\right)\right|^{2}.
\]
The result of the previous question ensures that
\[
P\left(X-Y=0\right)=\lim_{T\to+\infty}\dfrac{1}{2T}\intop^{+T}_{-T}\left|\varphi\left(t\right)\right|^{2}\text{d}t.
\]

On the other hand, since the random variables $X$ and $Y$ are independent
with law $\mu,$ the Fubini theorem gives
\[
P\left(X-Y=0\right)=\mu\otimes\mu\left(\left\{ x=y\right\} \right)=\intop_{\mathbb{R}}\left[\intop_{\mathbb{R}}\boldsymbol{1}_{\left(x=y\right)}\text{d}\mu\left(y\right)\right]\text{d}\mu\left(x\right)=\intop_{\mathbb{R}}\mu\left(\left\{ x\right\} \right)\text{d}\mu\left(x\right).
\]
As the set $S=\left\{ x:\,\mu\left(\left\{ x\right\} \right)\neq0\right\} $
is countable,
\[
\intop_{\mathbb{R}}\mu\left(\left\{ x\right\} \right)\text{d}\mu\left(x\right)=\intop_{S}\mu\left(\left\{ x\right\} \right)\text{d}\mu\left(x\right)=\sum_{x\in S}\mu\left(\left\{ x\right\} \right)^{2}=\sum_{x\in\mathbb{R}}\mu\left(\left\{ x\right\} \right)^{2}.
\]
Thus,\boxeq{
\[
\lim_{T\to+\infty}\dfrac{1}{2T}\intop^{+T}_{-T}\left|\varphi\left(t\right)\right|^{2}\text{d}t=\sum_{x\in\mathbb{R}}\mu\left(\left\{ x\right\} \right)^{2}.
\]
}

\end{solution}

\begin{solution}{}{solexercise13.11}

\textbf{1. Computation of $\varphi_{X}$}

The transfer theorem yields
\[
\forall t\in\mathbb{R},\,\,\,\,\varphi_{X}\left(t\right)=\intop_{\mathbb{R}^{2}}\text{e}^{\text{i}tx}\text{d}P_{\left(X,Y\right)}\left(x,y\right).
\]

Since the function $\left(x,y\right)\mapsto\text{e}^{\text{i}tx}$
is bounded, we may apply the generalized Fubini theorem---Theorem
$\ref{th:extended_fubini}.$ As the law of $Y$ is $P_{Y}=\rho\delta_{0}+\left(1-\rho\right)\delta_{1},$
we obtain
\begin{align*}
\varphi_{X} & =\intop_{\mathbb{R}}\left[\intop_{\mathbb{R}}\text{e}^{\text{i}tx}\text{d}P^{Y=y}_{X}\left(x\right)\right]\text{d}P_{Y}\left(y\right)\\
 & =\rho\intop_{\mathbb{R}}\text{e}^{\text{i}tx}\text{d}P^{Y=0}_{X}\left(x\right)+\left(1-\rho\right)\intop_{\mathbb{R}}\text{e}^{\text{i}tx}\text{d}P^{Y=1}_{X}\left(x\right)\\
 & =\rho\text{e}^{\text{i}t\cdot0}+\left(1-\rho\right)\intop_{\mathbb{R}}\text{e}^{\text{i}tx}\boldsymbol{1}_{\mathbb{R}^{+}}\left(x\right)\lambda\text{e}^{-\lambda x}\text{d}x.
\end{align*}
Thus,\boxeq{
\[
\forall t\in\mathbb{R},\,\,\,\,\varphi_{X}\left(t\right)=\rho+\left(1-\rho\right)\dfrac{\lambda}{\lambda-\text{i}t}.
\]
}

\textbf{2. Expectation and variance of $X.$}

The characteristic function of $X$ is twice differentiable. Hence,
the random variable $X$ admits a second-order moment, and
\[
\mathbb{E}\left(X\right)=-\text{i}\varphi^{\prime}_{X}\left(0\right)\,\,\,\,\text{and}\,\,\,\,\mathbb{E}\left(X^{2}\right)=-\varphi^{\prime\prime}_{X}\left(0\right).
\]
Since
\[
\varphi^{\prime}_{X}\left(t\right)=\left(1-\rho\right)\dfrac{\lambda\text{i}}{\left(\lambda-\text{i}t\right)^{2}}\,\,\,\,\text{and}\,\,\,\,\varphi^{\prime\prime}_{X}\left(t\right)=-2\left(1-\rho\right)\dfrac{\lambda\text{i}}{\left(\lambda-\text{i}t\right)^{3}},
\]
we obtain\boxeq{
\[
\mathbb{E}\left(X\right)=\dfrac{1-\rho}{\lambda}\,\,\,\,\,\,\,\,\mathbb{E}\left(X^{2}\right)=2\dfrac{1-\rho}{\lambda^{2}}\,\,\,\,\,\,\,\,\sigma^{2}_{X}=\dfrac{1-\rho^{2}}{\lambda^{2}}.
\]
}

\textbf{3. Independence, for every $n\in\mathbb{N}^{\ast},$ of $\epsilon_{n}$
and $\left(X_{0},X_{1},\cdots,X_{n-1}\right).$ Induction proof that
the $X_{n}$ have same characteristic function. Law of $X_{n}$} 

The random variable $\left(X_{0},X_{1},\cdots,X_{n-1}\right)$ is
a linear function of $\left(X_{0},\epsilon_{1},\cdots,\epsilon_{n-1}\right).$
Hence, it is independent of $\epsilon_{n},$ by the assumed independence.

We now prove by induction that the $X_{n}$ have same characteristic
function.
\begin{itemize}
\item \textbf{Initialization step}\\
The characteristic function of $X_{0}$ is for every real number $t,$
\[
\varphi_{X_{0}}\left(t\right)=\dfrac{\lambda}{\lambda-\text{i}t}.
\]
\item \textbf{Induction step}\\
Assume that $\varphi_{X_{n-1}}=\varphi_{X_{0}}.$ The random variables
$X_{n-1}$ and $\epsilon_{n}$ are independent; hence, by applying
the induction hypothesis, for every real number $t,$
\[
\varphi_{X_{n}}\left(t\right)=\varphi_{X_{n-1}}\left(t\right)\varphi_{\epsilon_{n}}\left(t\right)=\varphi_{X_{n-1}}\left(\rho t\right)\varphi_{\epsilon_{n}}\left(t\right).
\]
Thus,
\[
\varphi_{X_{n}}\left(t\right)=\dfrac{\lambda}{\lambda-\text{i}\rho t}\left[\rho+\left(1-\rho\right)\dfrac{\lambda}{\lambda-\text{i}t}\right].
\]
\end{itemize}
Hence,\boxeq{
\[
\varphi_{X_{n}}\left(t\right)=\dfrac{\lambda}{\lambda-\text{i}t}.
\]
}

It follows that, for every $n\in\mathbb{N},$ $\varphi_{X_{n}}=\varphi_{X_{0}},$
and thus the random variables $X_{n}$ follow the same exponential
law $\exp\left(\lambda\right).$

\textbf{4. Expression of $m^{X_{n-1}=\cdot}_{X_{n}}$}

Since $X_{n-1}$ and $\epsilon_{n}$ are independent, the conditional
mean $m^{X_{n-1}=\cdot}_{X_{n}}$ of $X_{n}$ knowing $X_{n-1}$ is
given, for every real number $x_{n-1},$ by
\[
m^{X_{n-1}=x_{n-1}}_{X_{n}}=\rho x_{n-1}+\mathbb{E}\left(\epsilon_{n}\right),
\]
thus,\boxeq{
\[
m^{X_{n-1}=x_{n-1}}_{X_{n}}=\rho x_{n-1}+\dfrac{1-\rho}{\lambda}.
\]
}

\textbf{5. Expression of $\varphi_{\left(X_{n-1},X_{n}\right)}$}

The characteristic function $\varphi_{\left(X_{n-1},X_{n}\right)}$
of the random variable $\left(X_{n-1},X_{n}\right)$ is, for every
$\left(u,v\right)\in\mathbb{R}^{2},$
\[
\varphi_{\left(X_{n-1},X_{n}\right)}\left(u,v\right)=\mathbb{E}\left(\text{e}^{\text{i}\left(uX_{n-1}+vX_{n}\right)}\right)=\mathbb{E}\left(\text{e}^{\text{i}\left(u+\rho v\right)X_{n-1}}\text{e}^{\text{i}v\epsilon_{n}}\right).
\]
By independence of the random variables $X_{n-1}$ and $\epsilon_{n},$
\[
\varphi_{\left(X_{n-1},X_{n}\right)}\left(u,v\right)=\varphi_{X_{n-1}}\left(u+\rho v\right)\varphi_{\epsilon_{n}}\left(v\right),
\]
which yields\boxeq{
\[
\forall\left(u,v\right)\in\mathbb{R}^{2},\,\,\,\,\varphi_{\left(X_{n-1},X_{n}\right)}\left(u,v\right)=\dfrac{\lambda}{\lambda-\text{i}\left(u+\rho v\right)}\left[\rho+\left(1-\rho\right)\dfrac{\lambda}{\lambda-\text{i}v}\right].
\]
}

\textbf{6. Determination of $m^{X_{n}=x}_{X_{n-1}}$}

Since the pair $\left(X_{n-1},X_{n}\right)$ admits an expectation,
its characteristic function is differentiable and
\[
\dfrac{\partial}{\partial u}\varphi_{\left(X_{n-1},X_{n}\right)}\left(u,v\right)=\text{i}\mathbb{E}\left(X_{n-1}\text{e}^{\text{i}\left(uX_{n-1}+vX_{n}\right)}\right).
\]
Hence,
\[
\dfrac{\partial}{\partial u}\varphi_{\left(X_{n-1},X_{n}\right)}\left(0,v\right)=\text{i}\mathbb{E}\left(X_{n-1}\text{e}^{\text{i}vX_{n}}\right).
\]
By the generalized Fubini theorem\boxeq{
\[
\forall v\in\mathbb{R},\,\,\,\,\dfrac{\partial}{\partial u}\varphi_{\left(X_{n-1},X_{n}\right)}\left(0,v\right)=\text{i}\intop_{\mathbb{R}}\text{e}^{\text{i}vx}m^{X_{n}=x}_{X_{n-1}}\text{d}P_{X_{n}}\left(x\right),
\]
}or, since $X_{n}$ follows the exponential law $\exp\left(\lambda\right),$
\[
\forall v\in\mathbb{R},\,\,\,\,\dfrac{\partial}{\partial u}\varphi_{\left(X_{n-1},X_{n}\right)}\left(0,v\right)=\text{i}\intop_{\mathbb{R}}\text{e}^{\text{i}vx}f\left(x\right)\text{d}x=\text{i}\widehat{f}\left(v\right)
\]
where the nonnegative integrable function $f,$ with Fourier transform
$\widehat{f},$ is defined, for every real number $x,$ by 
\[
f\left(x\right)=\boldsymbol{1}_{\mathbb{R}^{+}}\left(x\right)\lambda\text{e}^{-\lambda x}m^{X_{n}=x}_{X_{n-1}}.
\]
Moreover, a direct computation of the partial derivative $\varphi_{\left(X_{n-1},X_{n}\right)}$
gives
\[
\forall v\in\mathbb{R},\,\,\,\,\dfrac{\partial}{\partial u}\varphi_{\left(X_{n-1},X_{n}\right)}\left(0,v\right)=\dfrac{\lambda\text{i}}{\left(\lambda-\text{i}\rho v\right)\left(\lambda-\text{i}v\right)}.
\]
It follows that the Fourier transform $\widehat{f}$ of $f$ is equal
for every real number $v$ to
\[
\widehat{f}\left(v\right)=\dfrac{\lambda}{\left(\lambda-\text{i}\rho v\right)\left(\lambda-\text{i}v\right)},
\]
which is a Lebesgue-integrable function. The inversion formula then
yields that
\[
f\left(x\right)=\dfrac{1}{2\pi}\intop_{\mathbb{R}}\dfrac{\lambda}{\left(\lambda-\text{i}\rho v\right)\left(\lambda-\text{i}v\right)}\text{e}^{-\text{i}xv}\text{d}v.
\]

It remains to compute this integral. We have
\[
\dfrac{\lambda}{\left(\lambda-\text{i}\rho v\right)\left(\lambda-\text{i}v\right)}=-\dfrac{1}{1-\rho}\cdot\dfrac{1}{\dfrac{\lambda}{\rho}-\text{i}v}+\dfrac{1}{1-\rho}\cdot\dfrac{1}{\lambda-\text{i}v}
\]
and, by the Fourier transform injectivity theorem, applied to the
Fourier transform of the exponential law $\exp\left(\lambda\right),$
\[
\dfrac{1}{2\pi}\intop_{\mathbb{R}}\dfrac{\lambda}{\lambda-\text{i}v}\text{e}^{-\text{i}xv}\text{d}v=\boldsymbol{1}_{\mathbb{R}^{+}}\left(x\right)\text{e}^{-\lambda x}.
\]
It follows that
\[
f\left(x\right)=-\dfrac{1}{1-\rho}\boldsymbol{1}_{\mathbb{R}^{+}}\left(\dfrac{x}{\rho}\right)\text{e}^{-\frac{\lambda}{\rho}x}+\dfrac{1}{1-\rho}\boldsymbol{1}_{\mathbb{R}^{+}}\left(x\right)\text{e}^{-\lambda x}.
\]

The definition of $f$ and a straightforward computation shows that\boxeq{
\[
\forall x\in\mathbb{R},\,\,\,\,m^{X_{n}=x}_{X_{n-1}}=\boldsymbol{1}_{\mathbb{R}^{+}}\left(x\right)\dfrac{1}{\lambda\left(1-\rho\right)}\left[1-\text{e}^{-\frac{\lambda\left(1-\rho\right)}{\rho}x}\right].
\]
}

\end{solution}

\chapter{Gaussian Random Variables}\label{chap:PartIIChap14}

\begin{objective}{}{}

Chapter \ref{chap:PartIIChap14} is devoted to the study of Gaussian
random variables taking values in a vector space. It begins with a
brief review of the Gaussian law on $\mathbb{R}.$
\begin{itemize}
\item Section \ref{sec:Definition-and-Properties} presents the definition
and main properties of a Gaussian law on a vector space. The expectation
and variance of a random variable taking values in a vector space
are introduced, together with its characteristic function and Fourier
transform. The effect of affine transformation on Gaussian random
variables is also discussed.
\item Section \ref{sec:Gaussian-Measure-Existence.} focuses on the existence
of Gaussian measures and conditions for absolute continuity. After
examining existence issues, the notion of a degenerate Gaussian random
variable is introduced. The section concludes with sufficient condition
for a Gaussian random variable to admit a density.
\item Section \ref{sec:Marginals} focuses on the independence properties
of marginals of a Gaussian random variable. It defines the cross-covariance
operator of two random variables and introduces the notion of uncorrelated
random variables. Sufficient and necessary conditions for a couple
of random variables to be Gaussian are stated, before addressing the
characterization of Gaussian random variables.
\item Section \ref{sec:Regression:-the-linear} concludes the chapter by
studying the linear regression model on a vector space. The estimation
of regression parameters is then carried out, starting with the estimation
of the regression line of two random variables. The definition of
a linear estimator is given before enunciating the Gauss-Markov theorem.
The section ends with the Gaussian linear model and its statistical
formulation, stating a theorem on estimators of laws and their expectation
and variance. Applications to hypothesis testing, confidence intervals
and prediction conclude the chapter.
\end{itemize}
\end{objective}

\textbf{This chapter is dedicated to the study of Gaussian random
variables taking values in a finite-dimensional vector space $E.$
The study is done intrinsically, that is, independently of the choice
of a basis of $E.$ Although it plays no role in the definition of
Gaussian random variables taking values in $E,$ it is useful to assume
that $E$ is equipped with a scalar product: this avoids the explicit
use of the dual space and, in particular, allows the variance to be
viewed as a quadratic form on $E$---see Chapter \ref{chap:PartIIChap9}.
The Euclidean structure may also arise naturally, for instance, in
the study of estimation problems, or statistical tests, which leads
to the study of random variables with matrix values.}

\textbf{The reader may, if desired, assume that $E=\mathbb{R}^{d}.$
This case can always be reduced to choosing an orthonormal basis of
$E.$}

\textbf{At the end of the chapter, we study a regression problem in
the Gaussian setting. In particular, within the framework of the Gaussian
linear model, we solve the problem of parameters estimation as well
as problems of testing and of determining confidence intervals for
these parameters.}

In this chapter, unless explicitly mentioned, all random variables
are defined on the same probabilized space $\left(\Omega,\mathscr{A},P\right).$
We denote by $E$ a fixed real Euclidean space of dimension $d,$
whose scalar product is denoted $\left\langle \,\cdot\,,\,\cdot\,\right\rangle .$
The space $E$ is identified with its dual. Hence, a linear form $u$
on $E$ is denoted $\left\langle \,\cdot\,,u\right\rangle .$ We denote
by $\mathscr{E}$ the Borel $\sigma-$algebra of $E,$ that is, the
$\sigma-$algebra generated by the family of open subsets of $E.$
The space $E$ is always assumed to be equipped with its Borel $\sigma-$algebra
$\mathscr{E}.$

\begin{reminders}{}{}

The Gaussian law on $\mathbb{R},$ also called the Laplace-Gauss law
or normal law, with parameters $m\in\mathbb{R}$ and $\sigma^{2}\neq0,$
is the probability with density $f$ with respect to the Lebesgue
measure, where $f$ is defined by\boxeq{
\[
\forall x\in\mathbb{R},\,\,\,\,f\left(x\right)=\dfrac{1}{\sigma\sqrt{2\pi}}\text{e}^{-\frac{\left(x-m\right)^{2}}{2\sigma^{2}}}.
\]
}This law is denoted indifferently by $\mathcal{N}\left(m,\sigma^{2}\right)$
or $\mathcal{N}_{\mathbb{R}}\left(m,\sigma^{2}\right).$ 

Its Fourier transform is given by the relation\boxeq{
\begin{equation}
\forall t\in\mathbb{R},\,\,\,\,\widehat{\mathcal{N}\left(m,\sigma^{2}\right)}\left(t\right)=\text{e}^{\text{i}tm}\text{e}^{-\frac{t^{2}\sigma^{2}}{2}}.\label{eq:Fourier_transform_normal_law}
\end{equation}
}

A real-valued random variable $X$ is called Gaussian---or normal---if
its law is Gaussian. If the random variable $X$ has law $\mathcal{N}\left(m,\sigma^{2}\right),$
its characteristic function is defined by\boxeq{
\begin{equation}
\forall t\in\mathbb{R},\,\,\,\,\varphi_{X}\left(t\right)=\text{e}^{\text{i}tm}\text{e}^{-\frac{t^{2}\sigma^{2}}{2}}.\label{eq:charact_fct_normal_law}
\end{equation}
}

Its expectation is $m$ and its variance $\sigma^{2}.$ This random
variable admits moments of every order, which can be obtained, for
instance, by Taylor expansion of $\varphi_{X}.$ In particular---see
Chapter \ref{chap:PartIIChap13}, Exercise $\ref{exo:exercise13.4}$---if
$X$ has law $\mathcal{N}\left(0,1\right),$ then, for every $n\in\mathbb{N},$\boxeq{
\[
\mathbb{E}\left(X^{2n+1}\right)=0,\,\,\,\,\text{and if }n\geqslant1,\,\,\,\,\text{\ensuremath{\mathbb{E}\left(X^{2n}\right)}}=\dfrac{\left(2n\right)!}{2^{n}n!}
\]
}

Expanding the concept of Gaussian law---or measure---to a Euclidean
space leads to considering a Dirac measure as a degenerate Gaussian
law. Gaussian random variables supported on an affine space arrise
naturally in this context. In particular, a random variable that is
$P-$almost surely constant---hence with a zero variance---is still
Gaussian.

\end{reminders}

\section{Definition and Properties}\label{sec:Definition-and-Properties}

\begin{definition}{Gaussian Law. Gaussian Random Variable}{}

We call \textbf{\mindex{Gaussian!law}\mindex{law!Gaussian}Gaussian
law}---or \textbf{normal law}\index{normal law}---on $E$ a probability
$\mu$ on $\left(E,\mathscr{E}\right)$ such that the image measure
of $\mu$ under any linear form on $E$ is a Gaussian law on $\mathbb{R}.$

A random variable $X$ taking values in $E$ is said to be \textbf{\mindex{Gaussian!random variable}\mindex{random variable!Gaussian}Gaussian}---or
\textbf{\mindex{normal!random variable}\mindex{random variable!normal}normal}---if
its law $P_{X}$ is Gaussian on $E.$

\end{definition}

\begin{remark}{}{}A random variable $X$ taking values in $E$ is
Gaussian if and only if, for every $u\in E,$ the real-valued random
variable $\left\langle X,u\right\rangle $ is Gaussian. This follows
from the equality of probabilities $\left\langle \,\cdot\,,u\right\rangle \left(P_{X}\right)$
and $P_{\left\langle X,u\right\rangle }.$

\end{remark}

We now study some immediate and important consequences of these definitions,
that will help us to prove the existence of Gaussian laws and Gaussian
random variables with given expectation and variance.

\begin{denotation}{}{}

We denote $\mathscr{L}^{+}\left(E\right)$ the set of \textbf{operators}---or
\textbf{\mindex{endomorphism!self-adjoint and positive}endomorphisms}---\textbf{self-adjoint
and positive\mindex{operator!self-adjoint and positive}} on $E.$
An endomorphism $\Lambda$ of $E$ belongs to $\mathscr{L}^{+}\left(E\right)$
if and only if it satisfies $\Lambda=\Lambda^{\ast}$ and if $\left\langle \Lambda x,x\right\rangle \geqslant0,$
for every $x\in E.$

\end{denotation}

\begin{proposition}{Expectation and Variance, Chracteristic Function, Fourier Transform, Affine Transformation of a Gaussian Random Variable}{prop_after_def_gaussian_rv}Let
$X$ be a random variable taking values in $E.$

\textbf{(a) Expectation and variance}

If $X$ is Gaussian, then the random variable $\left\Vert X\right\Vert $
is in $\mathscr{L}^{2},$ and thus $X$ admits an expectation $m$---an
element of $E$---and a variance $\sigma^{2}_{X},$ which is a quadratic
form on $E.$ We denote $\Lambda_{X}$ the autocovariance operator
of $X,$ that is, the unique positive self-adjoint operator such that
\[
\forall x\in E,\,\,\,\,\left\langle \Lambda_{X}x,x\right\rangle =\sigma^{2}_{X}\left(x\right).
\]
We recall that we have\boxeq{
\begin{equation}
\forall x\in E,\,\,\,\,\mathbb{E}\left(\left\langle X,x\right\rangle \right)=\left\langle m,x\right\rangle \,\,\,\,\text{and}\,\,\,\,\sigma^{2}_{X}\left(x\right)=\sigma^{2}_{\left\langle X,x\right\rangle }=\left\langle \Lambda_{X}x,x\right\rangle .\label{eq:expectation_variance_Gaussian_rv}
\end{equation}

}

\textbf{(b) Gaussian random variable and characteristic function}

The random variable $X$ is Gaussian if and only if the characteristic
function $\varphi_{X}$ is given by\boxeq{
\begin{equation}
\forall t\in E,\,\,\,\,\varphi_{X}\left(t\right)=\text{e}^{\text{i}\left\langle m,t\right\rangle }\text{e}^{-\frac{1}{2}\left\langle Ct,t\right\rangle },\label{eq:charact_function_gaussian_rv}
\end{equation}
}where $m\in E$ and $C\in\mathscr{L}^{+}\left(E\right).$ In this
case, we have $m=\mathbb{E}\left(X\right)$ and $C=\Lambda_{X}.$

\textbf{Consequently, $\mathbb{E}\left(X\right)$ and $\Lambda_{X}$
characterize fully the law of the Gaussian random variable $X.$}

\textbf{(c) Necessary and sufficient condition for a measure to be
Gaussian }

A measure $\mu$ on $\left(E,\mathscr{E}\right)$ is Gaussian if and
only if its Fourier transform $\widehat{\mu}$ is given by the relation\boxeq{
\begin{equation}
\forall t\in E,\,\,\,\,\widehat{\mu}\left(t\right)=\text{e}^{\text{i}\left\langle m,t\right\rangle }\text{e}^{-\frac{1}{2}\left\langle Ct,t\right\rangle }.\label{eq:fourier_transform_gaussian_measure}
\end{equation}
}where $m\in E$ and $C\in\mathscr{L}^{+}\left(E\right).$ The parameters
$m$ and $C$ are then uniquely determined. This defines a probability
on $E.$ This probability law is denoted $\mathcal{N}_{E}\left(m,C\right),$
and is called\footnotemark \textbf{Gaussian law}\mindex{Gaussian!law}\mindex{law!Gaussian}---or
\textbf{Gaussian measure}\mindex{Gaussian!measure}---\textbf{with
parameters $m$ and $C.$}

\textbf{(d) Affine transformation}

Let $f$ be an Euclidean space. If $X$ is Gaussian with law $\mathscr{N}_{E}\left(m,\Lambda_{X}\right),$
then for every $A\in\mathscr{L}\left(E,F\right)$ and for every $b\in F,$
the random variable $AX+b,$ taking values in $F,$ is Gaussian with
law $\mathscr{N}_{F}\left(Am+b,A\Lambda_{X}A^{\ast}\right).$

\end{proposition}

\footnotetext{For the moment, nothing ensures that this measure exists.
Its existence will be proved in Theorem $\ref{th:existence_theorem}.$

}

\begin{proof}{}{}

\textbf{(a) Expectation and variance}

For every $x\in E,$ the real-valued random variable $\left\langle X,x\right\rangle $
is Gaussian, hence, in $\mathscr{L}^{2}.$ This is equivalent to saying
that $\left\Vert X\right\Vert $ belongs to $\mathscr{L}^{2}$---see
Chapter \ref{chap:PartIIChap9}, Proposition $\ref{pr:norm_and_bilinear_form_Lp}.$
The remainder of the statement follows directly from the definitions.

\textbf{(b) Gaussian random variable and characteristic function}

For every $t\in E$ and every $\alpha\in\mathbb{R},$
\begin{equation}
\varphi_{X}\left(\alpha t\right)=\mathbb{E}\left(\text{e}^{\text{i}\left\langle X,\alpha t\right\rangle }\right)=\varphi_{\left\langle X,t\right\rangle }\left(\alpha\right).\label{eq:char_function_alpha_t_of_X}
\end{equation}
If $X$ is Gaussian, then by $\refpar{eq:expectation_variance_Gaussian_rv},$
$\left\langle X,t\right\rangle $ is of Gaussian law
\[
\mathcal{N}\left(\left\langle m,t\right\rangle ,\left\langle \Lambda_{X}t,t\right\rangle \right),
\]
which, by taking $\alpha=1$ and using $\refpar{eq:charact_fct_normal_law}$
yields the result.

Conversely, if $\varphi_{X}$ is given by $\refpar{eq:charact_function_gaussian_rv},$
then it follows from $\refpar{eq:char_function_alpha_t_of_X}$ that,
for every $t\in E$ and for every $\alpha\in\mathbb{R},$
\[
\varphi_{\left\langle X,t\right\rangle }\left(\alpha\right)=\text{e}^{\text{i}\alpha\left\langle m,t\right\rangle }\text{e}^{-\frac{\alpha^{2}}{2}\left\langle Ct,t\right\rangle }.
\]

This proves that $\left\langle X,t\right\rangle $ is Gaussian, and
consequently that $X$ is Gaussian. More precisely, the law of $\left\langle X,t\right\rangle $
is the law $\mathcal{N}\left(\left\langle m,t\right\rangle ,\left\langle Ct,t\right\rangle \right).$
Hence,
\[
\forall t\in E,\,\,\,\,\mathbb{E}\left(X,t\right)=\left\langle m,t\right\rangle \,\,\,\,\text{and}\,\,\,\,\sigma^{2}_{\left\langle X,t\right\rangle }=\left\langle Ct,t\right\rangle ,
\]
which proves, taking into account $\refpar{eq:expectation_variance_Gaussian_rv},$
that
\[
\mathbb{E}\left(X\right)=m\,\,\,\,\text{and\,\,\,\,}\forall t\in E,\,\,\,\,\sigma^{2}_{\left\langle X,t\right\rangle }=\left\langle Ct,t\right\rangle =\left\langle \Lambda_{X}t,t\right\rangle .
\]
That is, $\mathbb{E}\left(X\right)=m$ and $C=\Lambda_{X}.$

\textbf{(c) Necessary and sufficient condition for a measure to be
Gaussian }

The identity function $I$ on $E,$ viewed as a random variable defined
on the probabilized space $\left(E,\mathscr{E},\mu\right)$ with values
in $\left(E,\mathscr{E}\right)$ is of law $\mu$ and of characteristic
function $\widehat{\mu}.$ Hence, $I$ is a Gaussian random variable
if and only if the measure $\mu$ is Gaussian. It therefore suffices
to apply the previous characterization.

\textbf{(d) Affine translation}

By using the definition of the transpose, for every $t\in F,$
\[
\varphi_{AX+b}\left(t\right)=\text{e}^{\text{i}\left\langle b,t\right\rangle }\varphi_{X}\left(A^{\ast}t\right)=\text{e}^{\text{i}\left\langle Am+b,t\right\rangle }\text{e}^{-\frac{1}{2}\left\langle A\Lambda_{X}A^{\ast}t,t\right\rangle },
\]
which yields the result, as shown above.

\end{proof}

\begin{remark}{}{}

We recall---see Chapter \ref{chap:PartIIChap9}---that the matrix
$C_{X},$ which represents the operator $\Lambda_{X}$ in the orthonormal
basis $\left(e_{j}\right)_{1\leqslant j\leqslant d},$ is the covariance
matrix of $X$ in this basis and that
\[
\left(C_{X}\right)_{ij}=\text{cov}\left(\left\langle X,e_{i}\right\rangle ,\left\langle X,e_{j}\right\rangle \right).
\]

In particular, if $E=\mathbb{R}^{d},$ the usual choice of basis is
the canonical basis. In this case, the operator $\Lambda_{X}$ is
represented by the matrix $C_{X}$ of covariances $\text{cov}\left(X_{i},X_{j}\right)$
of the marginals $X_{i},\,1\leqslant i\leqslant d.$ The Gaussian
law of $X$ is then simply denoted $\mathscr{N}_{\mathbb{R}^{d}}\left(m,C_{X}\right).$

\end{remark}

\section{Gaussian Measure Existence. Condition of Absolute Continuity}\label{sec:Gaussian-Measure-Existence.}

\begin{lemma}{Product Measure of Centered Reduced Gaussian Laws on $\left(\mathbb{R}^d,\mathscr{B}_{\mathbb{R}^d}\right)$}{meas_prod_c_r_normal_law}

The product measure $\left[\mathcal{N}_{\mathbb{R}}\left(0,1\right)\right]^{\otimes d}$
on $\left(\mathbb{R}^{d},\mathscr{B}_{\mathbb{R}^{d}}\right)$ is
the Gaussian measure $\mathcal{N}_{\mathbb{R}^{d}}\left(0,\text{I}_{\mathbb{R}^{d}}\right),$
where $\text{I}_{\mathbb{R}^{d}}$ is the identity matrix on $\mathbb{R}^{d}.$

\end{lemma}

\begin{proof}{}{}

Since the Fourier transform of a product measure is the product of
the Fourier transforms of its factors, for every $t\in\mathbb{R}^{d},$
\begin{align*}
\widehat{\left[\mathcal{N}_{\mathbb{R}}\left(0,1\right)\right]^{\otimes d}}\left(t\right) & =\prod^{d}_{j=1}\left[\mathcal{N}_{\mathbb{R}}\left(0,1\right)\right]\left(t_{j}\right)\\
 & =\prod^{d}_{j=1}\text{e}^{-\frac{t^{2}_{j}}{2}}=\text{e}^{-\frac{\left\Vert t\right\Vert ^{2}}{2}}=\widehat{\mathcal{N}_{\mathbb{R}^{d}}\left(0,\text{I}_{\mathbb{R}^{d}}\right)}\left(t\right),
\end{align*}
which proves the result by injectivity of the Fourier transform.

\end{proof}

\begin{remark}{}{}

This ensures the existence of the Gaussian measure $\mathcal{N}_{\mathbb{R}^{d}}\left(0,\text{I}_{\mathbb{R}^{d}}\right).$
Moreover, the product measure $\left[\mathcal{N}_{\mathbb{R}}\left(0,1\right)\right]^{\otimes d}$
admits a density $f,$ given by the product of densities of the marginal
measures. It is defined by
\[
\forall x\in\mathbb{R}^{d},\,\,\,\,f\left(x\right)=\dfrac{1}{\left(2\pi\right)^{\frac{d}{2}}}\text{e}^{-\frac{\left\Vert x\right\Vert ^{2}}{2}}.
\]
Hence, the Gaussian measure $\mathcal{N}_{\mathbb{R}^{d}}\left(0,\text{I}_{\mathbb{R}^{d}}\right)$
admits the density $f.$

\end{remark}

From this lemma, we deduce the next theorem which ensures the existence
of a Gaussian measure of given expectation $m$ and autocovariance
operator $\Lambda.$ This theorem has a purely algebraic content:
it is about showing the existence of an operator $B$ such that $BB^{\ast}=\Lambda.$
We give two proofs of this fact: the first relies on the spectral
theorem for self-adjoint operators, and the second uses the decomposition
of quadratic forms into squares.

\begin{theorem}{Existence Theorem}{existence_theorem}

For every vector $m\in E$ and for every \textbf{self-adjoint} and
\textbf{positive}\footnotemark operator, there exists a unique Gaussian
measure $\mathcal{N}_{E}\left(m,\Lambda\right).$

\end{theorem}

\footnotetext{But not necessarily positive \textit{definite}!

}

\begin{proof}{}{}

It is sufficient to exhibit a Gaussian random variable $X$ taking
values in $E,$ with expectation $m$ and covariance operator $\Lambda.$

By Proposition $\ref{pr:prop_after_def_gaussian_rv},$ every random
variable of the form
\[
X=m+BX_{0},
\]
where $X_{0}$ is a random variable taking values in $\mathbb{R}^{k},$
with law $\mathscr{N}_{\mathbb{R}^{k}}\left(0,\text{I}_{\mathbb{R}^{k}}\right),$
and where $B\in\mathscr{L}\left(\mathbb{R}^{k},E\right)$ satisfies
\begin{equation}
BB^{\ast}=\Lambda,\label{eq:bb_ast_equal_lambda}
\end{equation}
answers the question.

It is always possible to take for $X_{0}$ as the identity function
from $\mathbb{R}^{k}$ onto itself, considered as a random variable
defined on the probabilizable space $\left(\mathbb{R}^{k},\mathscr{B}_{\mathbb{R}^{k}}\right)$
equipped with the probability $\mathcal{N}_{\mathbb{R}^{k}}\left(0,\text{I}_{\mathbb{R}^{k}}\right).$

To prove the existence of the operator $B$ satisfying the equality
$\refpar{eq:bb_ast_equal_lambda},$ there are two methods: the first
is based on the fact that self-adjoint operators can be diagonalized
in an orthonormal basis, and the second relies on the decomposition
of quadratic forms into squares, known as Gauss decomposition.
\begin{itemize}
\item \textbf{Method 1: Via diagonalization of self-adjoint operators}\\
We first define a positive self-adjoint operator $\Delta$ such that
\[
\Delta^{2}=\Lambda.
\]
The uniqueness of such an operator $\Delta$ can be shown. $\Delta$
is called the nonnegative square root of $\Lambda.$ \\
Indeed, since $\Lambda$ is self-adjoint, there exists an orthonormal
basis $\left(e_{i}\right)_{1\leqslant i\leqslant d}$ of $E$ consisting
of eigenvectors of $\Lambda.$ The operator $\Lambda$ can thus be
written as
\[
\Lambda=\sum^{d}_{i=1}\lambda_{i}\left\langle \,\cdot\,,e_{i}\right\rangle e_{i},
\]
where the $\lambda_{i}$ are the eigenvalues of $\Lambda,$ repeated
according to their multiplicity. They are nonnegative. \\
Define the operator
\[
\Delta=\sum^{d}_{i=1}\sqrt{\lambda_{i}}\left\langle \,\cdot\,,e_{i}\right\rangle e_{i}.
\]
It is self-adjoint, positive, and satisfies $\Delta^{2}=\Lambda.$
\\
Let $\Phi$ be the isomorphism from $\mathbb{R}^{d}$ on $E$ associated
to the basis $\left(e_{i}\right)_{1\leqslant i\leqslant d},$ defined
by
\[
\forall\left(a_{1},\cdots,a_{d}\right)\in\mathbb{R}^{d},\,\,\,\,\Phi\left(a_{1},\cdots,a_{d}\right)=\sum^{d}_{i=1}a_{i}e_{i}.
\]
We can take $B=\Delta\Phi.$\\
Indeed, the adjoint of $\Phi,$ isomorphism from $E$ to $\mathbb{R}^{d},$
is defined by
\[
\forall y\in E,\,\,\,\,\Phi^{\ast}\left(y\right)=\left(\left\langle y,e_{1}\right\rangle ,\cdots,\left\langle y,e_{d}\right\rangle \right),
\]
which follows immediately from the definition. Moreover, $\Phi\Phi^{\ast}=\text{Id}_{E},$
where $\text{Id}_{E}$ denotes the identity on $E.$ Therefore,
\[
BB^{\ast}=\Delta\Phi\Phi^{\ast}\Delta^{\ast}=\Delta^{2}=\Delta.
\]
\item \textbf{Method 2: Via Gauss decomposition of quadratic forms}\\
The square decomposition theorem for quadratic forms---Gauss decomposition---applied
to the quadratic form $x\mapsto\left\langle \Lambda x,x\right\rangle ,$
asserts that, for every $x\in E,$
\[
\left\langle \Lambda x,x\right\rangle =\sum^{r}_{i=1}\left\langle u_{i},x\right\rangle ^{2},
\]
where $r$ is the rank of $\Lambda$ and where the $u_{i}$ are independent
linear forms on $E,$ identified with elements of $E$.\\
Define $A\in\mathscr{L}\left(E,\mathbb{R}^{r}\right)$ by setting,
for every $x\in E,$
\[
Ax=\left(\left\langle u_{1},x\right\rangle ,\cdots,\left\langle u_{r},x\right\rangle \right).
\]
Then, for every $x\in E,$
\[
\left\langle \Lambda x,x\right\rangle =\left\langle Ax,Ax\right\rangle =\left\langle A^{\ast}Ax,x\right\rangle .
\]
Observe that $A$ admits as adjoint the linear function $B\in\mathscr{L}\left(\mathbb{R}^{r},E\right)$
defined, for every $\alpha\in\mathbb{R}^{r},$ by
\[
B\alpha=\sum^{r}_{i=1}\alpha_{i}u_{i}.
\]
We then obtain, for every $x\in E,$
\[
\left\langle \Lambda x,x\right\rangle =\left\langle BB^{\ast}x,x\right\rangle ,
\]
and hence $\Lambda=BB^{\ast}.$
\end{itemize}
\end{proof}

\begin{remarks}{Important!}{}

1. By the two methods, for the first with $k=d$ and for the second
with $k=r,$ we observe that $X$ is supported by the affine subspace
$m+\text{Im}\left(B\right).$
\begin{itemize}
\item With Method 1, it is straightforward that $\text{Im}\left(B\right)=\text{Im}\left(\Lambda\right).$
\item With Method 2, we have $\text{Im}\left(B\right)=\left(\text{Ker}\left(B^{\ast}\right)\right)^{\perp}=\left(\text{Ker}\left(A\right)\right)^{\perp}.$
But, since
\[
\text{Ker}\left(A\right)=\left\{ x:\,\left\langle u_{1},x\right\rangle =\cdots=\left\langle u_{r},x\right\rangle =0\right\} =\left\{ x:\,\left\langle \Lambda x,x\right\rangle =0\right\} ,
\]
hence\footnotemark,
\[
\text{\text{Ker}\ensuremath{\left(A\right)}}=\text{Ker}\left(\Lambda\right)\,\,\,\,\text{and}\,\,\,\,\text{Im}\left(B\right)=\text{Im}\left(\Lambda\right).
\]
\end{itemize}
2. With the second method, moreover, we additionally obtain immediately:
\begin{itemize}
\item The Gaussian measure $\mathcal{N}_{E}\left(m,\Lambda\right)$ is the
image of the standard measure $\mathcal{N}_{\mathbb{R}^{r}}\left(0,\text{Id}_{\mathbb{R}^{r}}\right)$
by the function
\begin{equation}
\alpha\mapsto m+\sum^{r}_{i=1}\alpha_{i}u_{i}.\label{eq:app_stand_measure_for_gaussian}
\end{equation}
\\
If $r=d,$ that is, if the operator $\Lambda$ is positive definite,
then this function is a diffeomorphism, and the measure $\mathcal{N}_{E}\left(m,\Lambda\right)$
admits a density with respect to the Lebesgue\footnotemark measure
on $E$---see Proposition $\ref{pr:suff_cond_gaus_rv_density}$ for
the computation of this density.
\item If $r<d,$ then the function $\refpar{eq:app_stand_measure_for_gaussian}$
is a diffeomorphism from $\mathbb{R}^{r}$ onto the affine subspace
$m+\text{Im}\Lambda.$ In this case, the measure $\mathcal{N}_{E}\left(m,\Lambda\right)$
admits a density with respect to the $r-$dimensional Lebesgue measure
on $m+\text{Im}\left(\Lambda\right)$---this Lebesgue measure is
well defined thanks to the Euclidean structure.
\end{itemize}
3. The two methods proposed to prove this theorem are constructive
and make it straightforward to derive two simulation algorithms for
a Gaussian random variable with law $\mathcal{N}_{\mathbb{R}^{d}}\left(m,\Lambda\right).$
When $r<d,$ the algorithm induced by the second method requires,
at least at first glance, fewer calls to the random number generator
than the one induced by the first method; however, this does not necessarily
mean it is faster. 

\end{remarks}

\addtocounter{footnote}{-1}

\footnotetext{Recall that the \textbf{isotropic cone\index{isotropic cone}}
of a quadratic form $q,$ ---that is, the set of vectors $x$ such
that $q\left(x\right)=0$---and the \textbf{kernel\index{kernel}}
of $q$---that is, the set of vectors $x$ such that $\varphi\left(x,y\right)=0,$
for every $y,$ where $\varphi$ is the bilinear form associated with
$q$---must not be confused. Nevertheless, these two sets coincide
in the case of a positive form, by the Schwarz inequality $\left|\varphi\left(x,y\right)\right|\leqslant q\left(x\right)q\left(y\right).$
}

\stepcounter{footnote}

\footnotetext{For a definition of the \textbf{Lebesgue measure on
an Euclidean space $E,$} the interested reader may refer to the add-on
at the end of this section.}

\begin{definition}{Degenerated Gaussian Random Variable }{}

With the previous notation, a Gaussian random variable with law $\mathcal{N}_{E}\left(m,\Lambda\right)$
is said \textbf{degenerate\mindex{random variable!Gaussian!degenerate}}
if the affine subspace $m+\text{Im}\Lambda$---which is still equal
to $m+\left(\text{Ker}\left(\Lambda\right)\right)^{\perp}$---is
a strict subspace of $E.$

\end{definition}

In Chapter \ref{chap:PartIIChap9}, Exercise $\ref{exo:exercise9.7},$
we saw that any random variable $X$ with square-integrable norm takes
its values $P-$almost surely in the affine subspace $\mathbb{E}\left(X\right)+\left(\text{Ker}\Lambda_{X}\right)^{\perp},$
wich coincides with $\mathbb{E}\left(X\right)+\text{Im}\left(\Lambda_{X}\right).$
Hence, for a Gaussian random variable with values in $\mathbb{R}^{d}$
to admit a density, it is necessary that the kernel reduces to $\left\{ 0\right\} .$
We now show in the next proposition that this condition is also sufficient.

\begin{proposition}{Sufficient Condition for a Gaussian Random Variable to Admit a Density}{suff_cond_gaus_rv_density}

Let $m\in\mathbb{R}^{d}.$ Let $C$ be a symmetric positive matrix
of size $d\times d.$ Let $X$ be an $\mathbb{R}^{d}-$valued random
variable with Gaussian law $\mathscr{N}_{\mathbb{R}^{d}}\left(m,C\right).$ 

(a) If $C$ is positive definite, then $X$ admits a density $f_{X}$
given, for every $x\in\mathbb{R}^{d},$ by\boxeq{
\begin{equation}
f_{X}\left(x\right)=\dfrac{1}{\left(\sqrt{2\pi}\right)^{d}}\left(\text{det}C\right)^{-\frac{1}{2}}\text{e}^{-\frac{1}{2}\left\langle C^{-1}\left(x-m\right),\left(x-m\right)\right\rangle }.\label{eq:density_of_Gaussian_rv_in_rd}
\end{equation}
}

(b) If $C$ is not positive definite, then $X$ takes its values $P-$almost
surely in the affine subspace $m+\text{Im}\left(C\right)$ and, consequently,
does not admit any density---that is, its law is not absolutely continuous.
In fact, it is even singular with respect to the Lebesgue measure
on $\mathbb{R}^{d}.$

\end{proposition}

\begin{proof}{}{}

(a) Assume that $C$ is positive definite. Let $Y$ be an $\mathbb{R}^{d}-$valued
random variable with Gaussian law $\mathscr{N}_{\mathbb{R}^{d}}\left(0,\text{Id}_{\mathbb{R}^{d}}\right).$
By Lemma $\ref{lm:meas_prod_c_r_normal_law},$ the marginals of $Y$
are independent and all have the same Gaussian law $\mathscr{N}_{\mathbb{R}}\left(0,1\right).$
Moreover, the random variable $Y$ admits a density $f_{Y},$ given
by the direct product of its marginal densities. For every $y\in\mathbb{R}^{d},$
\[
f_{Y}\left(y\right)=\dfrac{1}{\left(\sqrt{2\pi}\right)^{d}}\text{e}^{-\frac{\left\Vert y\right\Vert ^{2}}{2}},
\]
where $\left\Vert \,\cdot\,\right\Vert $ denotes the Euclidean norm
on $\mathbb{R}^{d}.$

Let $B$ be the nonnegative root square of $C.$ The random variable
$Z=m+BY$ follows the same Gaussian law $\mathcal{N}_{\mathbb{R}^{d}}\left(m,C\right)$
as $X.$ The function $y\mapsto m+By$ is a diffeomorphism, since
the matrices $C,$ and therefore $B,$ are invertible. Consequently,
the random variable $Z$ admits a density $f_{Z}$ given, for every
$z\in\mathbb{R}^{d},$ by
\begin{align*}
f_{Z}\left(z\right) & =f_{Y}\left(B^{-1}\left(z-m\right)\right)\left|\text{det}\left(B^{-1}\right)\right|\\
 & =\dfrac{1}{\left(\sqrt{2\pi}\right)^{d}}\text{e}^{-\frac{\left\Vert B^{-1}\left(z-m\right)\right\Vert ^{2}}{2}}\left|\text{det}\left(B^{-1}\right)\right|.
\end{align*}
This proves the result, after observing that
\[
\left\Vert B^{-1}\left(z-m\right)\right\Vert ^{2}=\left\langle C^{-1}\left(x-m\right),x-m\right\rangle \,\,\,\,\text{and}\,\,\,\,\text{det}\left(B^{-1}\right)=\left[\text{det}\left(C\right)\right]^{-\frac{1}{2}},
\]
using the fact that $B$ is self-adjoint and that $B^{2}=C.$

(b) If $C$ is not positive definite, then 
\[
P_{X}\left(m+\text{Im}\left(C\right)\right)=1,
\]
since $X$ takes $P-$almost surely its values in the affine subspace
$m+\text{Im}\left(C\right).$ On the other hand, 
\[
\lambda_{d}\left(m+\text{Im}\left(C\right)\right)=0,
\]
because the affine subspace $m+\text{Im}\left(C\right)$ is strict.
Hence, the law of $X$ is foreign to the Lebesgue measure $\lambda_{d}.$

\end{proof}

\begin{complem}{Lebesgue Measure on a Euclidean Space}{}

We now define the Lebesgue measure on a Euclidean space $E.$ After
having chosen an orthonormal basis on $E,$ we identify $\mathbb{R}^{d}$
with $E$ via the isomorphism $\Phi$ introduced in the proof of Theorem
$\ref{th:existence_theorem}.$ The Lebesgue measure on $E$ is then
defined as the measure image $\mu$ of the Lebesgue measure $\lambda_{d}$
on $\mathbb{R}^{d}$ by $\Phi.$ This definition is, in fact, independent
of the choice of the orthonormal basis.

Indeed, let $\Psi$ be another isomorphism corresponding to a---possibly
different---choice of the orthonormal basis on $E,$ and let $\nu$
denote the image measure of the Lebesgue measure $\lambda_{d}$ on
$\mathbb{R}^{d}$ by $\Psi.$ Then, for every $B\in\mathscr{E},$
\[
\nu\left(B\right)=\lambda_{d}\left(\Psi^{-1}\left(B\right)\right)=\intop_{\mathbb{R}^{d}}\boldsymbol{1}_{B}\circ\Psi\text{d}\lambda_{d}.
\]
Since $\Phi\Phi^{\ast}=\text{Id}_{E},$ 
\[
\nu\left(B\right)=\intop_{\mathbb{R}^{d}}\left(\boldsymbol{1}_{B}\circ\Phi\right)\circ\left(\Phi^{\ast}\circ\Psi\right)\text{d}\lambda_{d}.
\]
Making the change of variables defined by the diffeomorphism $\Phi^{\ast}\circ\Psi,$
whose Jacobian is $\pm1$---since $\Phi^{\ast}\circ\Psi$ is an isometry---we
obtain
\[
\nu\left(B\right)=\intop_{\mathbb{R}^{d}}\boldsymbol{1}_{B}\circ\Phi\text{d}\lambda_{d}=\mu\left(B\right).
\]

\end{complem}

\section{Marginals}\label{sec:Marginals}

We now focus on the \textbf{independence} properties of marginals.
We begin with the simple case where $E=\mathbb{R}^{d}$ and where
the marginals under consideration are all one-dimensional.

\begin{proposition}{Marginals of a Gaussian Law of a Random Variable in $\mathbb{R}^d$. Independence of Marginals}{marginals_gaussian_law}

Let $X=\left(X_{1},\cdots,X_{d}\right)$ be a random variable taking
values in $\mathbb{R}^{d}$ with Gaussian law $\mathcal{N}_{\mathbb{R}^{d}}\left(m,C\right),$
where $m\in\mathbb{R}^{d}$ and $C$ is a positive symmetric $d\times d$
matrix. Then the random variables $X_{j},\,1\leqslant j\leqslant d,$
are Gaussian.

Moreover, the random variables $X_{j},\,1\leqslant j\leqslant d,$
are independent if and only if they are pairwise uncorrelated, which
is equivalent to saying that the covariance matrix $C$ of $X$ is
diagonal.

\end{proposition}

\begin{proof}{}{}

The random variables $X_{j},\,1\leqslant j\leqslant d,$ are Gaussian,
since they are linear transformations of the Gaussian random variable
$X.$

If the random variables $X_{j},\,1\leqslant j\leqslant d,$ are independent,
then their covariances vanish pairwise, and the covariance matrix
$C$ of $X$ is therefore diagonal.

Conversely, suppose that the covariance matrix $C$ of $X$ is diagonal.
The characteristic function of $X$ then satisfies, for every $u\in\mathbb{R}^{d},$
\begin{align*}
\varphi_{X}\left(u\right) & =\text{e}^{\text{i}\left\langle m,u\right\rangle }\text{e}^{-\frac{\left\langle Cu,u\right\rangle }{2}}\\
 & =\prod^{d}_{j=1}\text{e}^{\text{i}\sum^{d}_{j=1}m_{j}u_{i}}\text{e}^{-\frac{1}{2}\sum^{d}_{j=1}C_{jj}u^{2}_{j}}.
\end{align*}
This shows that, for every $u$ in $\mathbb{R}^{d},$
\[
\varphi_{X}\left(u\right)=\prod^{d}_{j=1}\varphi_{X_{j}}\left(u_{j}\right).
\]

The marginal characteristic function is indeed obtained by the following
computation
\[
\varphi_{X_{j}}\left(u_{j}\right)=\varphi_{X}\left(0,\cdots,0,u_{j},0,\cdots,0\right)=\text{e}^{\text{i}\sum^{d}_{j=1}m_{j}u_{i}}\text{e}^{-\frac{1}{2}\sum^{d}_{j=1}C_{jj}u^{2}_{j}}.
\]
Hence the random variables $X_{j}$ are independent.

\end{proof}

We now turn to the marginal independence properties of a random variable
taking values in a Euclidean space $E.$

To this end, we recall the definition of the cross-covariance operator
between two random variables taking values in Euclidean spaces---see
Chapter \ref{chap:PartIIChap9}, Exercise $\ref{exo:exercise9.8}.$

\begin{definition}{Cross-covariance Operator. Uncorrelated Random Variables}{}

Let $F$ and $G$ be two Euclidean spaces and let $X\in\mathscr{L}^{2}_{F}\left(\Omega,\mathscr{A},P\right)$
and $Y\in\mathscr{L}^{2}_{G}\left(\Omega,\mathscr{A},P\right)$ be
two random variables. The \textbf{cross-covariance operator\index{cross-covariance operator}}
of $X$ and $Y$ is the unique operator $\Lambda_{X,Y}\in\mathscr{L}\left(F,G\right)$
satisfying
\[
\forall\left(x,y\right)\in F\times G,\,\,\,\,\left\langle \Lambda_{X,Y}x,y\right\rangle =\mathbb{E}\left(\left\langle \mathring{X},x\right\rangle \left\langle \mathring{Y},y\right\rangle \right)=\text{cov}\left(\left\langle X,x\right\rangle ,\left\langle Y,y\right\rangle \right).
\]

The random variables $X$ and $Y$ are said to be uncorrelated if
$\Lambda_{X,Y}=0.$

\end{definition}

\begin{remark}{}{}

The operator $\Lambda_{X,X}$ is nothing other than the autocovariance
operator of $X.$ Moreover, this concept of uncorrelated random variables
coincides, in the case where $E=F=\mathbb{R},$ with the notion of
uncorrelated real-valued random variables, as defined in Chapter \ref{chap:PartIIChap9}.

\end{remark}

\begin{proposition}{Necessary and Sufficient Condition for Direct Sum Components of a Gaussian Random Variable to be Independent}{cns_dir_sum_comp_gaus_rv_indep}

Let $X$ be a random variable taking values in a Euclidean space $E$
with Gaussian law $\mathcal{N}_{E}\left(m,\Lambda\right),$ where
$m\in E$ and $\Lambda\in\mathscr{L}^{+}\left(E\right).$

Let 
\[
E=\bigoplus^{n}_{j=1}E_{j},\,n\leqslant d,
\]
be a direct sum decomposition of $E,$ that is, each $x\in E$ can
be written uniquely as $x=\sum^{n}_{j=1}x_{j},$ where $x_{j}\in E_{j},$
for every $j\in\left\llbracket 1,n\right\rrbracket .$ Accordingly,
the random variable can be written as 
\[
X=\sum^{n}_{j=1}X_{j},
\]
where, for each $j\in\left\llbracket 1,n\right\rrbracket ,$ $X_{j}$
is a random variable taking values in the subspace $E_{j}.$

Then, the random variables $X_{j}$ are Gaussian.

Moreover, for $1\leqslant j\leqslant n,$ the random variables $X_{j}$
are pairwise uncorrelated.

\end{proposition}

\begin{proof}{}{}

The random variables $X_{j}$ are Gaussian, since they are linear
transforms---via the projections onto the subspaces $E_{j}$---of
the Gaussian random variable $X.$

Let $j$ and $k$ be two distinct integers in $\left\llbracket 1,n\right\rrbracket .$
If the random variables $X_{j},\,j\in\left\llbracket 1,n\right\rrbracket $
are independent, then for every $x_{j}\in E_{j}$ and for every $y_{k}\in E_{k},$
the random variables $\left\langle X_{j},x_{j}\right\rangle $ and
$\left\langle X_{k},y_{k}\right\rangle $ are independent. Consequently,
their covariance vanishes, and we obtain $\Lambda_{X_{j},X_{k}}=0.$

Conversely, suppose that the random variables $X_{j},\,j\in\left\llbracket 1,n\right\rrbracket $
are pairwise uncorrelated. For every choice of $u_{j}$ in $E_{j},\,j\in\left\llbracket 1,n\right\rrbracket ,$
the random variable $\left(\left\langle X_{1},u_{1}\right\rangle ,\left\langle X_{2},u_{2}\right\rangle ,\cdots,\left\langle X_{n},u_{n}\right\rangle \right)$
is Gaussian and takes values in $\mathbb{R}^{n}.$ By hypothesis,
its covariance matrix is diagonal. It then follows from Proposition
$\ref{pr:marginals_gaussian_law}$ that the random variables $\left\langle X_{j},u_{j}\right\rangle ,\,j\in\left\llbracket 1,n\right\rrbracket ,$
are independent. 

We therefore have
\begin{align*}
\varphi_{\left(X_{1},X_{2},\cdots,X_{n}\right)}\left(u_{1},u_{2},\cdots,u_{n}\right) & =\mathbb{E}\left(\text{e}^{\text{i}\sum^{n}_{j=1}\left\langle X_{j},u_{j}\right\rangle }\right)\\
 & =\mathbb{E}\left(\prod^{n}_{j=1}\text{e}^{\text{i}\left\langle X_{j},u_{j}\right\rangle }\right).
\end{align*}
By independence, this yields
\[
\varphi_{\left(X_{1},X_{2},\cdots,X_{n}\right)}\left(u_{1},u_{2},\cdots,u_{n}\right)=\prod^{n}_{j=1}\mathbb{E}\left(\text{e}^{\text{i}\left\langle X_{j},u_{j}\right\rangle }\right).
\]
Since the characteristic function of $X_{j}$ satisfies, for every
$u_{j}\in E_{j},$
\[
\varphi_{X_{j}}\left(u_{j}\right)=\varphi_{\left(X_{1},X_{2},\cdots,X_{n}\right)}\left(0,\cdots,0,u_{j},0,\cdots,0\right)=\mathbb{E}\left(\text{e}^{\text{i}\left\langle X_{j},u_{j}\right\rangle }\right),
\]
we conclude that, for every $\left(u_{1},u_{2},\cdots,u_{n}\right)\in\prod^{n}_{j=1}E_{j},$
\[
\varphi_{\left(X_{1},X_{2},\cdots,X_{n}\right)}\left(u_{1},u_{2},\cdots,u_{n}\right)=\prod^{n}_{j=1}\varphi_{X_{j}}\left(u_{j}\right),
\]
which is equivalent to the independence of the random variables $X_{j}.$

\end{proof}

\begin{remark}{}{}

A similar proposition to Proposition $\ref{pr:cns_dir_sum_comp_gaus_rv_indep}$
is obtained by replacing, in that proposition, the direct sum decomposition
$\bigoplus^{n}_{j=1}E_{j}$ by the Cartesian product $\prod^{n}_{j=1}E_{j},$
since these spaces are isomorphic as Euclidean spaces.

In particular, we obtain the following corollary concerning the marginals
of a Gaussian random variable taking values in the space $\mathbb{R}^{d}.$
We state it---of course, without a new proof---since it is of important
practical use.

\end{remark}

\begin{corollary}{Marginals of a Gaussian Random Variable Taking Values in $\mathbb{R}^d$}{marginals_gaussian_rv}

Let $X_{j},\,j\in\left\llbracket 1,n\right\rrbracket ,$ be random
variables defined on a probabilized space $\left(\Omega,\mathscr{A},P\right)$
taking values in $\mathbb{R}^{d_{j}}.$ If the random variable $X=\left(X_{1},X_{2},\cdots,X_{n}\right),$
taking values in $\mathbb{R}^{d_{1}+d_{2}+\cdots+d_{n}},$ is Gaussian,
and if the random variables $X_{j}$ are uncorrelated, then the random
variables $X_{j}$ are Gaussian and independent. 

\end{corollary}

\begin{remark}{}{}

As shown by the following counter-example, in Proposition $\ref{pr:cns_dir_sum_comp_gaus_rv_indep}$---and
also in its Corollary $\ref{co:marginals_gaussian_rv}$---it is essential
not to omit the hypothesis that the \textbf{global} random variable
$X$ is Gaussian. Moreover, saying that the $X_{i}$ are uncorrelated
is equivalent to saying that the covariance matrix $X$ is block-diagonal.

\end{remark}

\begin{counterexample}{Marginals are Gaussian without Being a Gaussian Random Variable}{marg_gauss_wo_gauss_rv}

Let $X$ be a real-valued random variable with a symmetric law, with
density $f_{X},$ and admitting a second-order moment. For every nonnegative
real number $a,$ define the random variable $Y_{a}$ by
\[
Y_{a}=-X\boldsymbol{1}_{\left(\left|X\right|\leqslant a\right)}+X\boldsymbol{1}_{\left(\left|X\right|>a\right)}\equiv X\left(2\boldsymbol{1}_{\left(\left|X\right|>a\right)}-1\right).
\]

1. Prove that $Y_{a}$ and $X$ have the same law.

2. Prove that, even if $X$ has law $\mathscr{N}_{\mathbb{R}}\left(0,1\right),$
the pair $\left(X,Y_{a}\right)$ is not Gaussian.

3. Prove that there exists $a,$ such that $\text{cov}\left(X,Y_{a}\right)=0.$
Conclude.

\end{counterexample}

\begin{solutioncounterexample}{}{}

\textbf{1. $Y_{a}$ and $X$ have the same law}

Indeed, for every $f\in\mathscr{C}^{+}_{\mathscr{K}}\left(\mathbb{R}\right),$
by the transfer theorem
\[
\mathbb{E}\left(f\left(Y_{a}\right)\right)=\intop_{\left(\left|x\right|\leqslant a\right)}f\left(-x\right)f_{X}\left(x\right)\text{d}x+\intop_{\left(\left|x\right|>a\right)}f\left(x\right)f_{X}\left(x\right)\text{d}x.
\]
Making the change of variables $x\mapsto-x$ in the first integral,
and using the parity of the density $f_{X},$ we obtain 
\begin{align*}
\mathbb{E}\left(f\left(Y_{a}\right)\right) & =\intop_{\left(\left|y\right|\leqslant a\right)}f\left(y\right)f_{X}\left(y\right)\text{d}y+\intop_{\left(\left|x\right|>a\right)}f\left(x\right)f_{X}\left(x\right)\text{d}x\\
 & =\intop_{\mathbb{R}}f\left(x\right)f_{X}\left(x\right)\text{d}x=\mathbb{E}\left(f\left(X\right)\right).
\end{align*}
This shows that \textbf{$Y_{a}$ and $X$ have the same law.} 

\textbf{2. Even if $X$ has law $\mathscr{N}_{\mathbb{R}}\left(0,1\right),$
$\left(X,Y_{a}\right)$ is not Gaussian}

Since the law of $X$ is symmetric, both $Y_{a}$ and $X$ are centered.
Their covariance is given by
\[
\text{cov}\left(X,Y_{a}\right)=\mathbb{E}\left(XY_{a}\right)=\mathbb{E}\left(X^{2}\left(2\boldsymbol{1}_{\left(\left|X\right|>a\right)}-1\right)\right)=4\intop_{\left(x>a\right)}x^{2}f_{X}\left(x\right)\text{d}x-\mathbb{E}\left(X^{2}\right).
\]
In particular, if $X$ has law $\mathscr{N}_{\mathbb{R}}\left(0,1\right),$
then the same holds for $Y_{a}.$ Nevertheless, the random variable
$X+Y_{a}=2X\boldsymbol{1}_{\left(\left|X\right|>a\right)}$ is not
Gaussian, since 
\[
P\left(X+Y_{a}=0\right)=P\left(\left|X\right|\leqslant a\right)>0.
\]
 Consequently, \textbf{the random variable $\left(X,Y_{a}\right)$
is not Gaussian}.

\textbf{3. There exists $a,$ such that $\text{cov}\left(X,Y_{a}\right)=0.$
Conclusion}

Finally, we can choose $a$ nonnegative such that $\text{cov}\left(X,Y_{a}\right)=0.$
Indeed, under the previous assumptions,
\[
\text{cov}\left(X,Y_{a}\right)=0\Leftrightarrow4\intop_{\left(\left|x\right|>a\right)}x^{2}f_{X}\left(x\right)\text{d}x=1.
\]

Since 
\[
\intop_{\mathbb{R}^{+}}x^{2}f\left(x\right)\text{d}x=\dfrac{1}{2},
\]
and since the function 
\[
a\mapsto\intop_{\left(x>a\right)}x^{2}f_{X}\left(x\right)\text{d}x
\]
is strictly decreasing on $\mathbb{R}^{+}$ and tends to 0 when $a$
tends to $+\infty,$ there exists a unique nonnegative $a$ such that
$\text{cov}\left(X,Y_{a}\right)=0.$ 

For this value of $a,$ the random variables $Y_{a}$ and $X$ have
the same Gaussian law and zero covariance, yet the couple $\left(X,Y_{a}\right)$
is not Gaussian.

\end{solutioncounterexample}

Nevertheless, we have the following important proposition

\begin{proposition}{Necessary and Sufficient Condition for a Couple of Random Variables to be Gaussian}{cns_couple_gaus}

Let $Z=\left(X,Y\right)$ be a random variable taking values in $E\times F,$
where $E$ and $F$ are Euclidean spaces. Suppose that $X$ and $Y$
are independent. For $Z$ to be Gaussian, it is necessary and sufficient
that $X$ and $Y$ are Gaussian. 

\end{proposition}

\begin{proof}{}{}

If $Z$ is Gaussian, then $X$ and $Y$ are Gaussian, since they are
linear transforms of $Z.$

Conversely, assume that $X$ and $Y$ are Gaussian. Their characteristic
functions are then given by
\[
\forall u\in E,\,\,\,\,\varphi_{X}\left(u\right)=\text{e}^{\text{i}\left\langle \mathbb{E}\left(X\right),u\right\rangle _{E}}\text{e}^{-\frac{1}{2}\left\langle \Lambda_{X}u,u\right\rangle _{E}}
\]
and
\[
\forall v\in F,\,\,\,\,\varphi_{Y}\left(v\right)=\text{e}^{\text{i}\left\langle \mathbb{E}\left(Y\right),v\right\rangle _{F}}\text{e}^{-\frac{1}{2}\left\langle \Lambda_{Y}v,v\right\rangle _{F}}.
\]
The independence of $X$ and $Y$ implies that the characteristic
function $\varphi_{Z}$ of $Z$ is the direct product of $\varphi_{X}$
and $\varphi_{Y}.$ Hence, for every $\left(u,v\right)\in E\times F,$
\[
\varphi_{Z}\left(u,v\right)=\text{e}^{\text{i}\left\langle \mathbb{E}\left(X\right),u\right\rangle _{E}+\text{i}\left\langle \mathbb{E}\left(Y\right),v\right\rangle _{F}}\text{e}^{-\frac{1}{2}\left[\left\langle \Lambda_{X}u,u\right\rangle _{E}+\left\langle \Lambda_{Y}v,v\right\rangle _{F}\right]}.
\]
The scalar product on $E\times F$ is defined by
\[
\forall\left(u,v\right)\in E\times F,\,\forall\left(u^{\prime},v^{\prime}\right)\in E\times F,\,\,\,\,\left\langle \left(u,v\right),\left(u^{\prime},v^{\prime}\right)\right\rangle _{E\times F}=\left\langle u,u^{\prime}\right\rangle _{E}+\left\langle v,v^{\prime}\right\rangle _{F},
\]

Define the operator $\Lambda\in\mathscr{L}^{+}\left(E\times F\right)$
by
\[
\left\langle \Lambda\left(u,v\right),\left(u,v\right)\right\rangle _{E\times F}=\left\langle \Lambda_{X}u,u\right\rangle _{E}+\left\langle \Lambda_{Y}v,v\right\rangle _{F}.
\]
Then
\[
\varphi_{Z}\left(u,v\right)=\text{e}^{\text{i}\left\langle \left(\mathbb{E}\left(X\right),\mathbb{E}\left(Y\right)\right),\left(u,v\right)\right\rangle _{E\times F}}\text{e}^{-\frac{1}{2}\left\langle \Lambda\left(u,v\right),\left(u,v\right)\right\rangle _{E\times F}},
\]
which proves that $Z$ is Gaussian.

\end{proof}

We now state two corollaries of Proposition $\ref{pr:cns_dir_sum_comp_gaus_rv_indep}$

\begin{corollary}{Necessary and Sufficient Condition of Independence of the Projections onto an Orthogonal Basis}{}

Let $X$ be a random variable taking values in the Euclidean space
$E$ with Gaussian law $\mathscr{N}_{E}\left(m,\Lambda\right),$ where
$m\in E$ and $\Lambda\in\mathscr{L}^{+}\left(E\right).$ Let $\left(e_{1},\cdots,e_{d}\right)$
be an orthogonal basis of $E.$

For the random variables $\left\langle X,e_{i}\right\rangle ,\,1\leqslant i\leqslant d$
to be independent, it is necessary and sufficent that, for every $i\in\left\llbracket 1,d\right\rrbracket ,$
the vector $e_{i}$ is an eigenvector of $\Lambda,$ which is equivalent
to say that the covariance matrix $C_{X},$ that is, the matrix representation
of $\Lambda$ in the basis $\left(e_{1},\cdots,e_{d}\right),$ is
diagonal.

\end{corollary}

\begin{proof}{}{}

Note that
\begin{equation}
\left\langle \Lambda e_{i},e_{j}\right\rangle =\text{cov}\left(\left\langle X,e_{i}\right\rangle ,\left\langle X,e_{j}\right\rangle \right).\label{eq:scal_prod_lambda_e___and_ei}
\end{equation}

If the random variables $\left\langle X,e_{i}\right\rangle ,\,1\leqslant i\leqslant d,$
are independent, then they are pairwise uncorrelated. Hence, for $i\neq j,$
\[
\left\langle \Lambda e_{i},e_{j}\right\rangle =0.
\]
It implies that, for every $i\in\left\llbracket 1,d\right\rrbracket ,$
\[
\Lambda e_{i}=\sum^{d}_{j=1}\left\langle \Lambda e_{i},e_{j}\right\rangle e_{j}=\left\langle \Lambda e_{i},e_{i}\right\rangle e_{i}.
\]
Therefore, $e_{i}$ is an eigenvector of $\Lambda$ associated with
the eigenvalue $\left\langle \Lambda e_{i},e_{i}\right\rangle =\sigma^{2}_{\left\langle X,e_{i}\right\rangle }.$

Conversely, suppose that, for every $i\in\left\llbracket 1,d\right\rrbracket ,$
$e_{i}$ is an eigenvector of $\Lambda$ associated with the nonnegative
eigenvalue $\lambda_{i},$ that is 
\[
\Lambda e_{i}=\lambda_{i}e_{i}.
\]
Since the basis $\left(e_{1},\cdots,e_{d}\right)$ is orthogonal,
it follows that, for $i\neq j,$ 
\[
\left\langle \Lambda e_{i},e_{j}\right\rangle =0.
\]
Thus, by the relation $\refpar{eq:scal_prod_lambda_e___and_ei},$
the random variables $\left\langle X,e_{i}\right\rangle ,\,1\leqslant i\leqslant d,$
are pairwise uncorrelated. Their independence then follows from Proposition
$\ref{pr:cns_dir_sum_comp_gaus_rv_indep}.$

\end{proof}

\begin{corollary}{Marginals of a Gaussian Random Variable on the Eigenspaces of its  Autocovariance Operator}{}

Let $X$ be a random variable taking values in a Euclidean space $E,$
with Gaussian law $\mathscr{N}_{E}\left(m,\Lambda\right),$ where
$m\in E$ and $\Lambda\in\mathscr{L}^{+}\left(E\right).$ Let $E_{j},j\in\left\llbracket 1,n\right\rrbracket ,$
denote the eigenspaces of $\Lambda.$ These subspaces are orthogonal,
invariant under $\Lambda,$ and form a direct sum decomposition of
$E.$ 

Denote by $\Pi_{j}$ the orthogonal projector onto $E_{j}.$ Then
the marginals $X_{j}=\Pi_{j}\circ X$ of $X$ on the subspaces $E_{j},\,j\in\left\llbracket 1,n\right\rrbracket $
are independent, and have respective laws $\mathscr{N}_{E_{j}}\left(\Pi_{j}m,\Lambda_{j}\right),$
where $\Lambda_{j}\in\mathscr{L}^{+}\left(E_{j}\right)$ is the restriction
of $\Lambda$ to $E_{j}.$

\end{corollary}

\begin{proof}{}{}

The independence follows from the fact that the eigenspaces $E_{j},\,j\in\left\llbracket 1,n\right\rrbracket $
form a direct sum decomposition of $E$---see Proposition $\refpar{eq:scal_prod_lambda_e___and_ei}.$
The random variables $X_{j}$ are Gaussian, with expectation $\Pi_{j}m$
and autocovariance operator $\Pi_{j}\Lambda\Pi^{\ast}_{j}=\Lambda_{j},$
since the subspaces $E_{j}$ are stable under $\Lambda.$

\end{proof}

In the previous counterexample $\ref{ex:marg_gauss_wo_gauss_rv}$,
we saw that a random variable may have Gaussian marginals without
itself being Gaussian. The next proposition provides a \textbf{characterization
of Gaussian random variables}.

\begin{proposition}{Characterization of Gaussian Random Variables}{}

Let $Z=\left(X,Y\right)$ be a random variable taking values in $E\times F,$
where $E$ and $F$ are Euclidean spaces. Suppose that $Z$ has a
square-integrable norm. Denote by $\Lambda_{X}$---respectively $\Lambda_{Y}$---the
covariance operator of $X$---respectively $Y$---and by $\Lambda_{X,Y}\in\mathscr{L}\left(E,F\right)$
the cross-covariance operator of $X$ and $Y.$ Assume that $\Lambda_{X}$
is invertible.

Then, the random variable $Z$ is Gaussian if and only if the following
two conditions are satisfied:

(i) The marginal $X$ is Gaussian.

(ii) There exist $A\in\mathscr{L}\left(E,F\right),$ $b\in F$ and
$\Lambda\in\mathscr{L}^{+}\left(F\right)$ such that, for $P_{X}-$almost
every $x\in E,$ the conditional law $P^{X=x}_{Y}$ of $Y$ given
$X=x$ is the Gaussian law $\mathscr{N}_{F}\left(Ax+b,\Lambda\right).$

In this case,\boxeq{
\[
m^{X=x}_{Y}=\mathbb{E}\left(Y\right)+\Lambda_{X,Y}\Lambda^{-1}_{X}\left(x-\mathbb{E}\left(X\right)\right)\,\,\,\,\text{and}\,\,\,\,\Lambda=\Lambda_{Y}-\Lambda_{X,Y}\Lambda^{-1}_{X}\Lambda^{\ast}_{X,Y},
\]
}where $m^{X=x}_{Y}$ denotes the conditional mean of $Y$ given
$X=x.$ 

\end{proposition}

\begin{proof}{}{}

\textbf{Suppose that $Z$ is Gaussian.} Then its marginal $X$ is
Gaussian. For $C\in\mathscr{L}\left(E,F\right),$ define $Y^{\prime}=Y-CX.$
The random variable $\left(X,Y^{\prime}\right),$ being a linear transform
of $Z,$ is then Gaussian. A straightforward computation shows that
\[
\Lambda_{X,Y^{\prime}}=\Lambda_{X,Y}-C\Lambda_{X}.
\]

Consequently, $\Lambda_{X,Y^{\prime}}=0$ if and only if $C=\Lambda_{X,Y}\Lambda^{-1}_{X}.$
We choose $C$ to be this value. It then follows from the first remark
following Proposition $\ref{co:marginals_gaussian_rv}$ that the random
variables $X$ and $Y^{\prime}$ are independent. Thus, for $P_{X}-$almost
every $x\in E,$ the following conditional law equalities hold
\[
P^{X=x}_{Y}=P^{X=x}_{Y^{\prime}+CX}=P^{X=x}_{Y^{\prime}+Cx}.
\]
Since $X$ and $Y^{\prime}$ are independent, we obtain 
\[
P^{X=x}_{Y}=P_{Y^{\prime}+Cx}=\mathscr{N}_{F}\left(\mathbb{E}\left(Y^{\prime}\right)+Cx,\Lambda_{Y^{\prime}}\right),
\]
which establishes the necessary condition. In fact, for $P_{X}-$almost
every $x\in E,$\boxeq{
\[
P^{X=x}_{Y}=\mathscr{N}_{F}\left(\mathbb{E}\left(Y^{\prime}\right)+\Lambda_{X,Y}\Lambda^{-1}_{X}\left(x-\mathbb{E}\left(X\right)\right),\Lambda_{Y}-\Lambda_{X,Y}\Lambda^{-1}_{X}\Lambda^{\ast}_{X,Y}\right).
\]
}

Indeed, we have $\mathbb{E}\left(Y^{\prime}\right)=\mathbb{E}\left(Y\right)-C\mathbb{E}\left(X\right)$
and, by independence of $Y^{\prime}$ and $CX,$
\[
\Lambda_{Y}=\Lambda_{Y^{\prime}}+\Lambda_{CX}=\Lambda_{Y^{\prime}}+C\Lambda_{X}C^{\ast}=\Lambda_{Y^{\prime}}+\left(\Lambda_{X,Y}\Lambda^{-1}_{X}\right)\Lambda_{X}\left(\Lambda_{X,Y}\Lambda^{-1}_{X}\right)^{\ast}.
\]
Taking into account the fact that $\Lambda_{X}$ is self-transposed,
this yields
\[
\Lambda_{Y^{\prime}}=\Lambda_{Y}-\Lambda_{X,Y}\Lambda^{-1}_{X}\Lambda^{\ast}_{X,Y}.
\]

\textbf{Conversely}, suppose that $X$ is Gaussian, and that there
exist $A\in\mathscr{L}\left(E,F\right),$ $b\in F$ and $\Lambda\in\mathscr{L}^{+}\left(F\right)$
such that, for every $P_{X}-$almost every $x\in E,$ the conditional
law $P^{X=x}_{Y}$ of $Y$ given $X=x$ is the Gaussian law $\mathscr{N}_{F}\left(Ax+b,\Lambda\right).$
Define $Y^{\prime\prime}=Y-AX-b.$ By the conditional transfer theorem,
for $P_{X}-$almost every $x\in E,$
\[
P^{X=x}_{Y^{\prime\prime}}=P^{X=x}_{Y-Ax-b}=\mathscr{N}_{F}\left(0,\Lambda\right),
\]
independent law of $x.$ The random variables $Y^{\prime\prime}$
and $X$ are therefore independent, and the law of $Y^{\prime\prime}$
is the Gaussian law $\mathscr{N}_{F}\left(0,\Lambda\right).$ 

It then follows from Proposition $\ref{pr:cns_couple_gaus}$ that
the random variable $\left(X,Y^{\prime\prime}\right)$ is Gaussian,
and hence so is the random variable $\left(X,Y\right),$ which is
a linear transform of $\left(X,Y^{\prime\prime}\right).$

\end{proof}

\section{Regression: the Linear Model}\label{sec:Regression:-the-linear}

In Chapter \ref{chap:PartIIChap9}, we studied the linear regression
problem for real-valued random vairables, as well as its generalization
to random variables taking values in a Euclidean space---see \ref{chap:PartIIChap9},
Exercise $\ref{exo:exercise9.8}.$ We now recall the statement and
the solution of these problems.

\subsubsection*{Case of Real-Valued Random Variables}

Let $X$ and $Y\in\mathscr{L}^{2}\left(\Omega,\mathscr{A},P\right)$
be real-valued random variables. The ``best'' approximation $\widetilde{Y}$
of $Y$ by an affine function of $X,$ in the least squares sense,
is obtained as the solution $\left(a,b\right)\in\mathbb{R}^{2}$ of
the minimization problem
\[
\text{inf}\left\{ \Phi\left(a,b\right):\,\left(a,b\right)\in\mathbb{R}^{2}\right\} ,
\]
where $\Phi\left(a,b\right)=\mathbb{E}\left(\left[Y-\left(aX+b\right)\right]^{2}\right),$
is\boxeq{
\[
\widetilde{Y}=\mathbb{E}\left(Y\right)+\rho_{X,Y}\dfrac{\sigma_{Y}}{\sigma_{X}}\left(X-\mathbb{E}\left(X\right)\right).
\]

}

The optimal couple $\left(\widehat{a},\widehat{b}\right)$ is given
by
\begin{equation}
\left\{ \begin{array}{l}
\widehat{a}=\rho_{X,Y}\dfrac{\sigma_{Y}}{\sigma_{X}}\\
\widehat{b}=\mathbb{E}\left(Y\right)-\mathbb{E}\left(X\right)\cdot\rho_{X,Y}\dfrac{\sigma_{Y}}{\sigma_{X}}.
\end{array}\right.\label{eq:reg_lin_real_opt_sol}
\end{equation}
The linear regression line of $Y$ on $X$ is therefore the line with
equation\boxeq{
\[
\left(y-\mathbb{E}\left(Y\right)\right)-\rho_{X,Y}\dfrac{\sigma_{Y}}{\sigma_{X}}\left(x-\mathbb{E}\left(X\right)\right)=0.
\]
}

The corresponding prediction error is
\[
\Phi\left(\widehat{a},\widehat{b}\right)=\mathbb{E}\left(\left[\mathring{Y}-\widehat{a}\mathring{X}\right]^{2}\right)=\sigma^{2}_{Y}-2\widehat{a}\text{cov}\left(X,Y\right)+\widehat{a}^{2}\sigma^{2}_{X}=\sigma_{Y^{2}}\left(1-\rho^{2}_{X,Y}\right).
\]

In particular, if the random variable follows a uniform law on the
set of $n$ points $\left\{ \left(x_{i},y_{i}\right)\right\} _{1\leqslant i\leqslant n},$
then 
\[
\Phi\left(a,b\right)=\dfrac{1}{n}\sum^{n}_{i=1}\left[y_{i}-\left(ax_{i}+b\right)\right]^{2}
\]
and we recover the classical least-squares approximation line used
in physics.

\subsubsection*{Case of Random Variables With Values in a Euclidean Space}

Let $X\in\mathscr{L}^{2}_{F}\left(\Omega,\mathscr{A},P\right)$ and
$Y\in\mathscr{L}^{2}_{G}\left(\Omega,\mathscr{A},P\right)$ be two
random vairables taking their values in the Euclidean spaces $F$
and $G,$ respectively. Suppose that the autocovariance operator $\Lambda_{X}$
is invertible. The ``best'' approximation $\widetilde{Y}$ of $Y$
by an affine function of $X,$ in the least-squares sense, is identified
as a solution in the couple $\left(A,b\right)\in\mathscr{L}\left(F,G\right)\times G$
of the minimization problem
\[
\inf\left\{ \Phi\left(A,b\right):\,\left(A,b\right)\in\mathscr{L}\left(F,G\right)\times G\right\} ,
\]
 where $\Phi\left(A,b\right)=\mathbb{E}\left(\left\Vert Y-\left(AX+b\right)\right\Vert ^{2}\right).$ 

This approximation is given by\boxeq{
\[
\widetilde{Y}=\mathbb{E}\left(Y\right)+\Lambda_{X,Y}\Lambda^{-1}_{X}\left(X-\mathbb{E}\left(X\right)\right),
\]
}where $\Lambda_{X,Y}$ denotes the \textbf{cross-covariance operator
of $X$ and $Y.$}

The optimal couple $\left(\widehat{A},\widehat{b}\right)$ is given
by\boxeq{
\[
\left(\widehat{A},\widehat{b}\right)=\left(\Lambda_{X,Y}\Lambda^{-1}_{X},\mathbb{E}\left(Y\right)-\Lambda_{X,Y}\Lambda^{-1}_{X}\mathbb{E}\left(X\right)\right).
\]
}The linear regression surface of $Y$ with respect to $X$ is the
surface---an affine subspace---defined by the equation\boxeq{
\[
\left(y-\mathbb{E}\left(Y\right)\right)-\Lambda_{X,Y}\Lambda^{-1}_{X}\left(x-\mathbb{E}\left(X\right)\right)=0.
\]
 }

The \textbf{prediction error} is
\[
\Phi\left(\widehat{A},\widehat{b}\right)=\mathbb{E}\left(\left\Vert \mathring{Y}-\widehat{A}\mathring{X}\right\Vert ^{2}\right)=\text{tr}\left(\Lambda_{Y}+\widehat{A}\Lambda_{X}\widehat{A}^{\ast}-\Lambda_{X,Y}\widehat{A}^{\ast}-\widehat{A}\Lambda_{Y,X}\right).
\]
 Taking into account that $\widehat{A}=\Lambda_{X,Y}\Lambda^{-1}_{X},$
we obtain\boxeq{
\[
\Phi\left(\widehat{A},\widehat{b}\right)=\text{tr}\left(\Lambda_{Y}-\Lambda_{X,Y}\Lambda^{-1}_{X}\Lambda^{\ast}_{X,Y}\right).
\]
}

These results are theoretical and require knowledge of the second-order
``moments'' of the couple $\left(X,Y\right).$ This leads to several
statistical problems. Giving experimental observations 
\[
\left(x_{1},y_{1}\right),\left(x_{2},y_{2}\right),\cdots,\left(x_{n},y_{n}\right),
\]
corresponding to realizations of the couple $\left(X,Y\right)$ through
independent experiences, one can obtain information about the \textbf{linear
regression line\index{linear regression line} of $Y$ on $X.$} Different
approaches can then be adopted: giving an \textbf{\mindex{linear regression line!estimate}estimate
of the linear regression line}, that is, providing suitable estimators
of the parameters $a$ and $b;$ testing admissible values of $a$
and $b;$ or constructing a \textbf{confidence interval\index{confidence interval}}
for these parameters. One may also address the problem of \textbf{predicting}
values of $Y$ given an observation of $X.$

We now address these different problems. The concepts introduced to
formulate them precisely and to solve them have a broad scope in statistics;
however, here we provide definitions only within the limited framework
of the regression problem. The estimation problem will be presented
more systematically in Chapter \ref{chap:PartIIChapConvMeasuresandInLaw15},
Section \ref{sec:Estimation}. The concept of hypothesis testing\footnote{Concerning the practical aspects of these statistical problems, the
interested reader may refer to the French book by Gilbert Saporta
(1990). Tr.N. We give here the 2006 reference, in French \cite{saporta2006probabilites}} is addressed at various points in this book---see, in particular,
the index entries related to the chi-squared test, the Student test,
and the Kolmogorov test. In what follows, for the sake of simplicity,
we consider only real-valued random variables.

\subsection{Estimate of regression parameters}

\subsubsection*{Problem 1}

Following the probabilistic modelling of a random phenomenon, we focus
on the couple of real-valued random variables $Z=\left(X,Y\right)$
which are assumed to represent two ``real'' measurements linked
to this phenomenon. The law of $Z$ is unknown by the experimenter;
nevertheless, after the computations and reasoning, the experimenter
is led to formulate \textbf{hypotheses} about this law---in particular,
that the random variables $X$ and $Y$ are of second order. 

Our objective is to \textbf{estimate the regression line of $Y$ on
$X$} from a \textbf{sample} of size $n$ drawn from $Z,$ namely
the vector $\uuline{z_{n}}=\left[\left(x_{1},y_{1}\right),\left(x_{2},y_{2}\right),\cdots,\left(x_{n},y_{n}\right)\right]$
in $\mathbb{R}^{2n},$ obtained by observing $n$ independent realizations
of this phenomenon.

The sample is assumed to represent a realization---that is the value
for a given outcome $\omega$---of $n$ \textbf{independent} random
variables $Z_{1},Z_{2},\cdots,Z_{n},$ all having the \textbf{same
law}{\bfseries\footnote{\textbf{The random variable $\uuline{Z_{n}}=\left(Z_{1},Z_{2},\cdots,Z_{n}\right)$
is called the ``empirical sample'' of size $n$ of the random variable
$Z.$}}}\textbf{ as $Z.$} The method used is the one of least-squares. For
a given $\omega,$ it consists in choosing the line with equation
\[
y=\widehat{a}_{n}\left(\omega\right)x+\widehat{b}_{n}\left(\omega\right)
\]
where the couple $\left(\widehat{a}_{n}\left(\omega\right),\widehat{b}_{n}\left(\omega\right)\right)$
is the solution of the minimization problem
\[
\inf\left\{ \Theta\left(a,b\right):\,\left(a,b\right)\in\mathbb{R}^{2}\right\} ,
\]
where
\[
\Theta\left(a,b\right)=\sum^{n}_{j=1}\left[Y_{j}\left(\omega\right)-\left(aX_{j}\left(\omega\right)+b\right)\right]^{2}.
\]
The line with equation $y=\widehat{a}_{n}\left(\omega\right)x+\widehat{b}_{n}\left(\omega\right)$
is called the least-square \textbf{estimate}\index{least-square estimate}
of the \textbf{regression line\index{regression line}} in $Y$ on
$X.$ Another justification for the use of this estimator is provided
by the Gauss-Markov theorem, which is presented below. This estimator
is determined---for every $\omega$---by the following proposition.

\begin{proposition}{Estimate of The Regression Line of Two Random Variables}{est_reg_line}

Let $X$ and $Y\in\mathscr{L}^{2}\left(\Omega,\mathscr{A},P\right)$
be two real-valued random variables. Let $\uuline{Z_{n}}=\left[\left(X_{1},Y_{1}\right),\left(X_{2},Y_{2}\right),\cdots,\left(X_{n},Y_{n}\right)\right]$
be an empirical sample of size $n$ of the random variable $Z=\left(X,Y\right).$
The \textbf{coefficients of the estimate}\index{coefficients of the estimate}---in
the least-squares sense---of the regression line of $Y$ on $X$
are given by\boxeq{
\begin{equation}
\left\{ \begin{array}{l}
\widehat{a}_{n}=r_{n}\dfrac{s_{Y,n}}{s_{X,n}}\\
\widehat{b}_{n}=\overline{Y}_{n}-\overline{X}_{n}\cdot r_{n}\dfrac{s_{Y,n}}{s_{X,n}},
\end{array}\right.\label{eq:coeff_estimate}
\end{equation}
}where the empirical moments associated to this sample are defined
as follows: 
\begin{itemize}
\item $\overline{X}_{n}$ and $\overline{Y}_{n}$ denote the \textbf{empirical
means\mindex{empirical!mean}} of $X$ and $Y,$ 
\item $s_{X,n}$ and $s_{Y,n}$ denote the \textbf{empirical variances\mindex{empirical!variance}}
of $X$ and $Y,$ 
\item and $r_{n}$ is the \textbf{empirical correlation coefficient\mindex{empirical!correlation coefficient}}
of $X$ and $Y.$
\end{itemize}
These quantities are defined by\boxeq{
\begin{align*}
\overline{X}_{n}=\dfrac{1}{n}\sum^{n}_{j=1}X_{j},\,\,\,\,\,\,\,\,\,\,\,\,\,\,\,\,\, & \qquad\qquad\overline{Y}_{n}=\dfrac{1}{n}\sum^{n}_{j=1}Y_{j}\\
s^{2}_{X,n}=\dfrac{1}{n}\sum^{n}_{j=1}X^{2}_{j}-\left(\overline{X}_{n}\right)^{2}, & \qquad\qquad s^{2}_{Y,n}=\dfrac{1}{n}\sum^{n}_{j=1}Y^{2}_{j}-\left(\overline{Y}_{n}\right)^{2},\\
r_{n}= & \dfrac{\dfrac{1}{n}\sum^{n}_{j=1}X_{j}Y_{j}-\overline{X}_{n}\overline{Y}_{n}}{s_{X,n}s_{Y,n}}.
\end{align*}
}

\end{proposition}

\begin{proof}{}{}

It is enough to apply, in the following manner, the results on linear
regression recalled earlier. For each fixed $\omega,$ consider the
probabilized space $\left(\mathbb{R}^{2},\mathscr{B}_{\mathbb{R}^{2}},\mu_{\omega}\right),$
where 
\[
\mu_{\omega}=\dfrac{1}{n}\sum^{n}_{j=1}\delta_{\left(X_{j}\left(\omega\right),Y_{j}\left(\omega\right)\right)}
\]
is the empirical measure associated with the sample. On this new probabilized
space, consider the random variables $U$ and $V,$ defined as the
canonical projections from $\mathbb{R}^{2}$ onto $\mathbb{R}.$

Since
\[
\Theta\left(a,b\right)=\sum^{n}_{j=1}\left[Y_{j}\left(\omega\right)-\left(aX_{j}\left(\omega\right)+b\right)\right]^{2}=n\intop_{\mathbb{R}^{2}}\left[V-\left(aU+b\right)\right]^{2}\text{d}\mu_{\omega},
\]
the formula $\refpar{eq:reg_lin_real_opt_sol}$ yields the result.

\end{proof}

\begin{remark}{}{}

The coefficients $\widehat{a}_{n}$ and $\widehat{b}_{n}$ introduced
in this proposition are in fact random variables whose values, for
each $\omega,$ determine an estimate of the regression line. They
are estimators---that is, measurable functions of the sample $\uuline{Z_{n}}$---of
the true coefficients $a$ and $b$ of the linear regression. 

\end{remark}

\subsubsection*{Problem 2}

In most situations, the random variable $X$ is deterministic. For
instance, during a chemical reaction, $X$ represents the catalyst
dose and $Y$ the quantity of a product formed by this reaction. The
linear model then consists in assuming that $Y$ can be written in
the form\boxeq{
\begin{equation}
Y=ax+b+\epsilon,\label{eq:linear_model_reaction}
\end{equation}
}where $\epsilon$ is a centered second-order random variable, representing
an approximation error or a measurement error. The problem is then,
given the observed results $y_{1},y_{2},\cdots,y_{n}$ of independent
experiments performed at the respective levels $x_{1},x_{2},\cdots,x_{n}$
of the values $x,$ to obtain an estimate of the coefficients $a$
and $b.$

The associated statistical model is as follows. We define $n$ independent
observations $Y_{1},Y_{2},\cdots,Y_{n}$ of $Y$ done at the levels
$x_{1},x_{2},\cdots,x_{n}$ of values of $x.$ That is, the $Y_{i}$
are random variables given by
\begin{equation}
Y_{i}=ax_{i}+b+\epsilon_{i},\label{eq:system_of_linear_model}
\end{equation}
where the random variables $\epsilon_{i},1\leqslant i\leqslant n,$
are independent, centered, second order random variables with the
same unknown variance $\sigma^{2}.$

We seek estimates of $a$ and $b$ in terms of the observations $Y_{i}.$
A least-squares technique can be used by transforming this problem
into one of the type considered in Problem 1: we regard the random
variable $X_{i}$ as constants equal to $x_{i}$ and minimize the
sum of the squared errors. The Gauss-Markov theorem provides a justification
for the use of the least-squares estimators obtained in this method. 

In what follows, unless explicitly stated otherwise, $\left\langle \cdot,\cdot\right\rangle $
denotes the usual scalar product of $\mathbb{R}^{n},$ and $\left\Vert \,\cdot\,\right\Vert $
the associated norm. 

\begin{definition}{Linear Estimator}{}

A \textbf{linear estimator\index{linear estimator}} of the unknown
parameter $\left(a,b\right)\in\mathbb{R}^{2}$ is a linear transformation
of the vector $\uuline{Y}=\left(Y_{1},Y_{2},\cdots,Y_{n}\right)$
of the form 
\[
T_{u,v}=\left(\left\langle \uuline{Y},u\right\rangle ,\left\langle \uuline{Y},v\right\rangle \right),
\]
where $u,v\in\mathbb{R}^{n}.$ 

An estimator of the unknown parameter $\left(a,b\right)\in\mathbb{R}^{2}$
is said to be\textbf{ unbiased\index{unbiased}} if its expectation
is equal to $\left(a,b\right),$ for every $\left(a,b\right)\in\mathbb{R}^{2}.$
A \textbf{linear estimator} $T_{u,v}$ of the unknown parameter $\left(a,b\right)$
is said to be of \textbf{minimum variance\index{minimum variance}}
among all \textbf{unbiased} linear estimator of $\left(a,b\right)$
if it is solution of the minimization problem
\begin{equation}
\min\left\{ \sigma^{2}_{\left\langle \uuline{Y},u\right\rangle }+\sigma^{2}_{\left\langle \uuline{Y},v\right\rangle }:\,u,v\in\mathbb{R}^{n}\right\} .\label{eq:minimisation_problem}
\end{equation}

\end{definition}

\begin{theorem}{Gauss-Markov Theorem}{}

Let $Y_{1},Y_{2},\cdots,Y_{n}$ be $n$ independent observations of
$Y$ taken at the levels $x_{1},x_{2},\cdots,x_{n}$ of the values
of $x.$ That is, suppose that the random variables $Y_{i}$ are given
by
\begin{equation}
Y_{i}=ax_{i}+b+\epsilon_{i},\label{eq:Y_i_as_linear_fct_of_x_i}
\end{equation}
where the random variables $\epsilon_{i},\,1\leqslant i\leqslant n$
are independent, centered, second-order random variables with the
same unknown variance $\sigma^{2}.$

The linear estimator of minimum variance among all unbiased linear
estimators of $\left(a,b\right)$ is the estimator $\left(\widehat{a}_{n},\widehat{b}_{n}\right),$
where $\widehat{a}_{n}$ and $\widehat{b}_{n}$ are estimators of
the least-squares of $a$ and $b$ given by Proposition $\ref{pr:est_reg_line}.$
They are given by\boxeq{
\begin{equation}
\left\{ \begin{array}{l}
\widehat{a}_{n}=r_{n}\dfrac{s_{Y,n}}{s_{x,n}}\\
\widehat{b}_{n}=\overline{Y}_{n}-\overline{x}_{n}\cdot r_{n}\dfrac{s_{Y,n}}{s_{x,n}},
\end{array}\right.\label{eq:estimates_least_squ}
\end{equation}
}where we denote\boxeq{
\begin{align*}
\overline{x}_{n}=\dfrac{1}{n}\sum^{n}_{j=1}x_{j},\,\,\,\,\,\,\,\,\,\,\,\,\,\,\,\,\, & \qquad\qquad\overline{Y}_{n}=\dfrac{1}{n}\sum^{n}_{j=1}Y_{j}\\
s^{2}_{X,n}=\dfrac{1}{n}\sum^{n}_{j=1}x^{2}_{j}-\left(\overline{x}_{n}\right)^{2}, & \qquad\qquad s^{2}_{Y,n}=\dfrac{1}{n}\sum^{n}_{j=1}Y^{2}_{j}-\left(\overline{Y}_{n}\right)^{2},\\
r_{n}= & \dfrac{\dfrac{1}{n}\sum^{n}_{j=1}x_{j}Y_{j}-\overline{x}_{n}\overline{Y}_{n}}{s_{x,n}s_{Y,n}}.
\end{align*}
}

\end{theorem}

\begin{proof}{}{}

Let $e=\left(1,1,\cdots,1\right),\,x=\left(x_{1},\cdots,x_{n}\right)$
and $\uuline{\epsilon}=\left(\epsilon_{1},\cdots,\epsilon_{n}\right).$ 

Then
\begin{equation}
\uuline{Y}=ax+be+\uuline{\epsilon}\label{eq:Y_underline_estimate}
\end{equation}
and since the random variable $\uuline{\epsilon}$ is centered,
\begin{equation}
\mathbb{E}\left(\uuline{Y}\right)=ax+be.\label{eq:expectation_Y_underline}
\end{equation}
A linear estimator $T_{u,v}=\left(\left\langle \uuline{Y},u\right\rangle ,\left\langle \uuline{Y},v\right\rangle \right)$
of $\left(a,b\right)$ is therefore unbiased if and only if, for every
$\left(a,b\right)\in\mathbb{R}^{2},$
\[
\mathbb{E}\left(\left\langle \uuline{Y},u\right\rangle \right)=a\,\,\,\,\text{and}\,\,\,\,\mathbb{E}\left(\left\langle \uuline{Y},v\right\rangle \right)=b.
\]
Hence, by the relation $\refpar{eq:expectation_Y_underline}$, this
holds if and only if $u\in A$ and $v\in B,$ where we define
\[
\left\{ \begin{array}{l}
A=\left\{ \left\langle x,u\right\rangle =1\,\,\,\,\text{and}\,\,\,\,\left\langle e,u\right\rangle =0\right\} \\
B=\left\{ \left\langle x,v\right\rangle =0\,\,\,\,\text{and}\,\,\,\,\left\langle e,v\right\rangle =1\right\} .
\end{array}\right.
\]
Noting that
\[
\sigma^{2}_{\left\langle \uuline{Y},u\right\rangle }=\left\langle \Lambda_{\uuline{Y}}u,u\right\rangle =\sigma^{2}\left\Vert u\right\Vert ^{2},
\]
the unbiased linear estimator $T_{\widehat{u},\widehat{v}}$ is therefore
of minimum variance---whatever the value of $\sigma$---among all
unbiased linear estimators of $\left(a,b\right)$ when $\widehat{u}$
and $\widehat{v}$ are solutions of the two constrained minimization
problems
\begin{equation}
\min\left\{ \left\Vert u\right\Vert ^{2}:\,u\in A\right\} \label{eq:min_norm_u_square_in_A}
\end{equation}
and
\begin{equation}
\min\left\{ \left\Vert v\right\Vert ^{2}:\,v\in B\right\} .\label{eq:min_norm_v_sq_in_B}
\end{equation}

We study the first extremum problem associated to $\refpar{eq:min_norm_u_square_in_A}.$
To the Lagrange multipliers $\lambda$ and $\mu,$ we associate the
function $\Phi_{\lambda,\mu}$ defined for every $u\in\mathbb{R}^{n}$
by
\begin{equation}
\Phi_{\lambda,\mu}\left(u\right)=\left\Vert u\right\Vert ^{2}-\lambda\left(\left\langle x,u\right\rangle -1\right)-\mu\left\langle e,u\right\rangle .\label{eq:phi_lambdau_mu_u}
\end{equation}
A point $\widehat{u}$ is a solution of the relative extremum problem
linked to $\refpar{eq:min_norm_u_square_in_A}$ if there exist $\lambda$
and $\mu$ such that
\begin{equation}
\Phi^{\prime}_{\lambda,\mu}\left(\widehat{u}\right)=0,\,\,\,\,\left\langle x,\widehat{u}\right\rangle =1\,\,\,\,\text{and}\,\,\,\,\left\langle e,\widehat{u}\right\rangle =0.\label{eq:deriv_phi_lamba_mu}
\end{equation}

Since
\[
\Phi^{\prime}_{\lambda,\mu}\left(u\right)=2\left\langle u,\cdot\right\rangle -\lambda\left\langle x,\cdot\right\rangle -\mu\left\langle e,\cdot\right\rangle =2\left\langle u,\cdot\right\rangle -\left\langle \lambda x+\mu e,\cdot\right\rangle ,
\]
we have $\Phi^{\prime}_{\lambda,\mu}\left(\widehat{u}\right)=0$ if
and only if 
\[
\widehat{u}=\dfrac{1}{2}\left(\lambda x+\mu e\right).
\]
Therefore, $\widehat{u}$ is a solution of the extremum problem if
there exist $\lambda$ and $\mu$ solving the system
\[
\left\{ \begin{array}{l}
\left\langle x,\dfrac{1}{2}\left(\lambda x+\mu e\right)\right\rangle =1\\
\left\langle e,\dfrac{1}{2}\left(\lambda x+\mu e\right)\right\rangle =0,
\end{array}\right.
\]
which is equivalent, since $\left\langle e,e\right\rangle =n,$ to
\[
\left\{ \begin{array}{l}
\lambda\left\Vert x\right\Vert ^{2}+\mu\left\langle x,e\right\rangle =2\\
\lambda\left\langle x,e\right\rangle +n\mu=0.
\end{array}\right.
\]
This system admits the solutions
\[
\lambda=\dfrac{2n}{n\left\Vert x\right\Vert ^{2}-\left\langle x,e\right\rangle ^{2}},\,\,\,\,\mu=-\dfrac{2\left\langle x,e\right\rangle }{n\left\Vert x\right\Vert ^{2}-\left\langle x,e\right\rangle ^{2}}.
\]
The unique solution $\widehat{u}$ is thus
\[
\widehat{u}=\dfrac{1}{2}\left(\lambda x+\mu e\right)=\dfrac{1}{2}\left(\dfrac{2n}{n\left\Vert x\right\Vert ^{2}-\left\langle x,e\right\rangle ^{2}}x-\dfrac{2\left\langle x,e\right\rangle }{n\left\Vert x\right\Vert ^{2}-\left\langle x,e\right\rangle ^{2}}e\right).
\]
Hence, 
\[
\widehat{u}=\dfrac{1}{n\left\Vert x\right\Vert ^{2}-\left\langle x,e\right\rangle ^{2}}\left(nx-\left\langle x,e\right\rangle e\right).
\]
The point $\widehat{u}$ corresponds to a \textbf{global extremum}.
It remains to show that this extremum is a \textbf{minimum}. 

To this end, decompose any $u\in A$ as $u=\widehat{u}+\delta.$ Since
$\widehat{u}\in A,$ 
\[
\left\langle x,\delta\right\rangle =0\,\,\,\,\text{and}\,\,\,\,\left\langle e,\delta\right\rangle =0,
\]
which implies
\[
\left\langle \widehat{u},\delta\right\rangle =\left\langle \dfrac{1}{2}\left(\lambda x+\mu e\right),\delta\right\rangle =\dfrac{1}{2}\left[\lambda\left\langle x,\delta\right\rangle +\mu\left\langle e,\delta\right\rangle \right]=0.
\]

Hence, for every $u\in A,$
\[
\left\Vert u\right\Vert ^{2}=\left\Vert \widehat{u}\right\Vert ^{2}+\left\Vert \delta\right\Vert ^{2}\geqslant\left\Vert \widehat{u}\right\Vert ^{2},
\]
which shows that $\widehat{u}$ is the unique solution of the extremum
problem linked to $\refpar{eq:min_norm_u_square_in_A}.$ 

Finally,
\begin{equation}
\left\langle \uuline{Y},\widehat{u}\right\rangle =\dfrac{n\left\langle \uuline{Y},x\right\rangle -\left\langle x,e\right\rangle \left\langle \uuline{Y},e\right\rangle }{n\left\Vert x\right\Vert ^{2}-\left\langle x,e\right\rangle ^{2}},\label{eq:scal_Y_u_widehat}
\end{equation}
which is precisely the estimator $\widehat{a}_{n},$ as can be verified
by a straightforward computation.

We now study the second extremum problem associated with $\refpar{eq:min_norm_v_sq_in_B}.$ 

To the Lagrange multipliers $\lambda$ and $\mu,$ we associate the
function $\psi_{\lambda,\mu}$ define for every $v\in\mathbb{R}^{n}$
by
\begin{equation}
\psi_{\lambda,\mu}\left(v\right)=\left\Vert v\right\Vert ^{2}-\lambda\left\langle x,v\right\rangle -\mu\left(\left\langle e,v\right\rangle -1\right).\label{eq:psi_lambda_mu_v}
\end{equation}

A point $\widehat{v}$ is a solution of the relative extremum problem
linked to $\refpar{eq:min_norm_v_sq_in_B}$ if there exist $\lambda$
and $\mu$ such that
\begin{equation}
\psi^{\prime}_{\lambda,\mu}\left(\widehat{v}\right)=0,\,\,\,\,\left\langle x,\widehat{v}\right\rangle =0\,\,\,\,\text{and}\,\,\,\,\left\langle e,\widehat{v}\right\rangle =1.\label{eq:derivative_of_psi_lambda_mu}
\end{equation}
Since
\[
\psi^{\prime}_{\lambda,\mu}\left(v\right)=2\left\langle v,\cdot\right\rangle -\lambda\left\langle v,\cdot\right\rangle -\mu\left\langle e,\cdot\right\rangle =2\left\langle v,\cdot\right\rangle -\left\langle \lambda x+\mu e,\cdot\right\rangle ,
\]
we have $\psi^{\prime}_{\lambda,\mu}\left(\widehat{v}\right)=0$ if
and only if 
\[
\widehat{v}=\dfrac{1}{2}\left(\lambda x+\mu e\right).
\]
Thus, $\widehat{v}$ is a solution of the extremum problem if there
exist $\lambda$ and $\mu$ solving the system
\[
\left\{ \begin{array}{l}
\left\langle x,\dfrac{1}{2}\left(\lambda x+\mu e\right)\right\rangle =0\\
\left\langle e,\dfrac{1}{2}\left(\lambda x+\mu e\right)\right\rangle =1,
\end{array}\right.
\]
which is equivalent, since $\left\langle e,e\right\rangle =n,$ to
\[
\left\{ \begin{array}{l}
\lambda\left\Vert x\right\Vert ^{2}+\mu\left\langle x,e\right\rangle =0\\
\lambda\left\langle x,e\right\rangle +n\mu=2.
\end{array}\right.
\]
This system admits the solutions
\[
\lambda=\dfrac{-2\left\langle x,e\right\rangle }{n\left\Vert x\right\Vert ^{2}-\left\langle x,e\right\rangle ^{2}}\,\,\,\,\,\,\,\,\mu=\dfrac{2\left\Vert x\right\Vert ^{2}}{n\left\Vert x\right\Vert ^{2}-\left\langle x,e\right\rangle ^{2}}.
\]
The unique solution $\widehat{v}$ found is therefore
\[
\widehat{v}=\dfrac{1}{2}\left(\lambda x+\mu e\right)=\dfrac{1}{2}\left(\dfrac{-2\left\langle x,e\right\rangle }{n\left\Vert x\right\Vert ^{2}-\left\langle x,e\right\rangle ^{2}}x-\dfrac{2\left\Vert x\right\Vert ^{2}}{n\left\Vert x\right\Vert ^{2}-\left\langle x,e\right\rangle ^{2}}e\right),
\]
hence, 
\[
\widehat{v}=\dfrac{1}{n\left\Vert x\right\Vert ^{2}-\left\langle x,e\right\rangle ^{2}}\left(\left\Vert x\right\Vert ^{2}e-\left\langle x,e\right\rangle x\right).
\]
The point $\widehat{v}$ corresponds to a global extremum. It remains
to show that this extremum is a minimum. 

To this end, decompose any $v\in B,$ as $v=\widehat{v}+\delta.$
Since $\widehat{v}\in B,$ 
\[
\left\langle x,\delta\right\rangle =0\,\,\,\,\text{and}\,\,\,\,\left\langle e,\delta\right\rangle =0,
\]
which implies
\[
\left\langle \widehat{v},\delta\right\rangle =\left\langle \dfrac{1}{2}\left(\lambda x+\mu e\right),\delta\right\rangle =\dfrac{1}{2}\left[\lambda\left\langle x,\delta\right\rangle +\mu\left\langle e,\delta\right\rangle \right]=0.
\]

Hence, for every $v\in B,$
\[
\left\Vert v\right\Vert ^{2}=\left\Vert \widehat{v}\right\Vert ^{2}+\left\Vert \delta\right\Vert ^{2}\geqslant\left\Vert \widehat{v}\right\Vert ^{2},
\]
which shows that $\widehat{v}$ is the unique solution of the extremum
problem linked to $\refpar{eq:min_norm_v_sq_in_B}.$ 

Finally,
\begin{equation}
\left\langle \uuline{Y},\widehat{v}\right\rangle =\dfrac{\left\Vert x\right\Vert ^{2}\left\langle \uuline{Y},e\right\rangle -\left\langle x,e\right\rangle \left\langle \uuline{Y},x\right\rangle }{n\left\Vert x\right\Vert ^{2}-\left\langle x,e\right\rangle ^{2}},\label{eq:scal_Y_u_widehat-1}
\end{equation}
which is nothing but the estimator $\widehat{b}_{n},$ as can be verified
by a direct, albeit somewhat lengthy, computation.

\end{proof}

\subsection{The Gaussian Linear Model}

We now study in greater detail the linear model $\refpar{eq:linear_model_reaction}$
and its statistical formulation $\refpar{eq:Y_i_as_linear_fct_of_x_i}.$ 

To obtain quantitative information for the estimators previously derived,
we must strengthen our assumptions by specifying the law of the \textbf{independent}
errors $\epsilon_{i}.$ We assume that the $\epsilon_{i}$ are following
the same Gaussian law $\mathscr{N_{\mathbb{R}}}\left(0,\sigma^{2}\right),$
where the variance $\sigma^{2}$ is unknown. In this case, the model
is referred to as the \textbf{\mindex{Gaussian!linear model}Gaussian
linear model}. 

With the notation of the previous section, the random variable $\uuline{Y}$
defined in $\refpar{eq:Y_underline_estimate}$ then follows the law
$\mathscr{N}_{\mathbb{R}^{n}}\left(ax+be,\sigma^{2}I_{n}\right)$
where $I_{n}$ denotes the identy matrix of $\mathbb{R}^{n}.$ The
density of $\uuline{Y}$---called the \textbf{\index{likelihood}likelihood}
by statisticians---is therefore given, for every $y\in\mathbb{R}^{n},$
by
\[
f_{\uuline{Y}}\left(y\right)=\dfrac{1}{\left(2\pi c\right)^{\frac{n}{2}}}\text{e}^{-\frac{\left\Vert y-m\left(a,b\right)\right\Vert ^{2}}{2c}},
\]
where we write $m\left(a,b\right)=ax+be$ and\footnote{The parameter that needs to be estimated is the variance, not the
standard deviation.} $c=\sigma^{2}.$ 

We now define and give \textbf{estimators}\mindex{estimator!likelihood maximum}\footnote{The notion of \textbf{estimator of maximum likelihood} is introduced
more systematically in Chapter \ref{chap:PartIIChapConvMeasuresandInLaw15}.} of the \textbf{likelihood maximum\mindex{likelihood!maximum}} of
$a,\,b$ and $c.$ These estimators are obtained as follows: for each
fixed $y,$ we determine the parameters that maximize the likelihood
at $y.$ In this setting, such maximizers exist and are unique; they
are denoted respectively $\widehat{a}\left(y\right),\,\widehat{b}\left(y\right)$
and $\widehat{c}\left(y\right).$ 

The maximum likelihood estimators of $a,\,b$ and $c$ are then the
random variables $\widehat{a}\left(\uuline{Y}\right),\,\widehat{b}\left(\uuline{Y}\right)$
and $\widehat{c}\left(\uuline{Y}\right).$

In this setting, it is more convenient to maximize, with respect to
$a,\,b$ and $c,$ the quantity known as the \textbf{log-likelihood\index{log-likelihood}}
at $y,$ namely
\[
\ln\left[f_{\uuline{Y}}\left(y\right)\right]=-\dfrac{n}{2}\ln\left(2\pi c\right)-\dfrac{\left\Vert y-m\left(a,b\right)\right\Vert ^{2}}{2c}.
\]
We now determine the stationary points. We have:
\begin{itemize}
\item $\dfrac{\partial}{\partial a}\ln\left[f_{\uuline{Y}}\left(y\right)\right]=0$
if and only if $\dfrac{\partial}{\partial a}\left\Vert y-m\left(a,b\right)\right\Vert ^{2}=0.$
\item $\dfrac{\partial}{\partial b}\ln\left[f_{\uuline{Y}}\left(y\right)\right]=0$
if and only if $\dfrac{\partial}{\partial b}\left\Vert y-m\left(a,b\right)\right\Vert ^{2}=0.$
\end{itemize}
Since 
\[
\left\Vert y-m\left(a,b\right)\right\Vert ^{2}=\left\Vert y\right\Vert ^{2}-2\left\langle y,m\left(a,b\right)\right\rangle +\left\Vert m\left(a,b\right)\right\Vert ^{2},
\]
we obtain:
\begin{itemize}
\item $\dfrac{\partial}{\partial a}\left\Vert y-m\left(a,b\right)\right\Vert ^{2}=-2\left\langle y,\cdot\right\rangle x+2\left\langle m\left(a,b\right),\cdot\right\rangle x=2\left\langle m\left(a,b\right)-y,\cdot\right\rangle x,$
\item $\dfrac{\partial}{\partial b}\left\Vert y-m\left(a,b\right)\right\Vert ^{2}=-2\left\langle y,\cdot\right\rangle e+2\left\langle m\left(a,b\right),\cdot\right\rangle e=2\left\langle m\left(a,b\right)-y,\cdot\right\rangle e.$
\end{itemize}
Moreover,
\[
\dfrac{\partial}{\partial c}\ln\left[f_{\uuline{Y}}\left(y\right)\right]=-\dfrac{n}{2c}+\dfrac{\left\Vert y-m\left(a,b\right)\right\Vert ^{2}}{2c^{2}}.
\]
A stationary point $\left(\widehat{a}\left(y\right),\widehat{b}\left(y\right),\widehat{c}\left(y\right)\right)$
must therefore satisfy
\begin{equation}
\widehat{a}\left(y\right)x+\widehat{b}\left(y\right)e=y,\label{eq:stat_point_cond1}
\end{equation}
and
\begin{equation}
\widehat{c}\left(y\right)=\dfrac{\left\Vert y-m\left(\widehat{a}\left(y\right),\widehat{b}\left(y\right)\right)\right\Vert ^{2}}{n}.\label{eq:stat_point_cond2}
\end{equation}

Taking the scalar product of both sides of $\refpar{eq:stat_point_cond1}$
successively with $x$ and with $e,$ we find that $\widehat{a}\left(y\right)$
and $\widehat{b}\left(y\right)$ must solve the system
\begin{equation}
\left\{ \begin{array}{l}
\widehat{a}\left(y\right)\left\Vert x\right\Vert ^{2}+\widehat{b}\left(y\right)\left\langle x,e\right\rangle =\left\langle y,x\right\rangle \\
\widehat{a}\left(y\right)\left\langle x,e\right\rangle +\widehat{b}\left(y\right)\left\Vert e\right\Vert ^{2}=\left\langle y,e\right\rangle ,
\end{array}\right.\label{eq:stat_point_cond1_scal}
\end{equation}
which admits the unique solution---using the identity $\left\Vert e\right\Vert ^{2}=n$---
\begin{equation}
\widehat{a}\left(y\right)=\dfrac{n\left\langle y,x\right\rangle -\left\langle x,e\right\rangle \left\langle y,e\right\rangle }{n\left\Vert x\right\Vert ^{2}-\left\langle x,e\right\rangle ^{2}}\qquad\qquad\widehat{b}\left(y\right)=\dfrac{\left\Vert x\right\Vert ^{2}\left\langle y,e\right\rangle -\left\langle x,e\right\rangle \left\langle y,x\right\rangle }{n\left\Vert x\right\Vert ^{2}-\left\langle x,e\right\rangle ^{2}}.\label{eq:sol_a_b_stationary_point}
\end{equation}
It follows, by comparing the equalities $\refpar{eq:scal_Y_u_widehat}$
and $\refpar{eq:scal_Y_u_widehat-1}$ that $\widehat{a}\left(\uuline{Y}\right)=\widehat{a}_{n}$
and $\widehat{b}\left(\uuline{Y}\right)=\widehat{b}_{n}.$

Thus, \textbf{in the Gaussian linear model, the maximum likelihood
estimators of $a$ and $b$ coincide with the least squares estimators
and, moreover, with the minimum-variance linear unbiased estimators.
The maximum likelihood estimator $\widehat{c}_{n}$ of the variance
is therefore given by
\begin{equation}
\widehat{c}_{n}=\dfrac{\left\Vert \uuline{Y}-\left(\widehat{a}_{n}x+\widehat{b}_{n}e\right)\right\Vert ^{2}}{n}.\label{eq:estimator_likelihood_max_variance}
\end{equation}
}

\begin{theorem}{Estimator Laws, Expectation and Variance of these Estimators}{est_exp_var_laws}

The random variable $\left(\widehat{a}_{n},\widehat{b}_{n},\uuline{Y}-\left(\widehat{a}_{n}x+\widehat{b}_{n}e\right)\right)$
taking values in $\mathbb{R}^{n+2},$ is Gaussian. The estimator $\left(\widehat{a}_{n},\widehat{b}_{n}\right)$
is a Gaussian random variable independent of $\widehat{c}_{n}.$ 

The expectations and variances of $\widehat{a}_{n}$ and $\widehat{b}_{n}$
are given by\boxeq{
\begin{equation}
\left\{ \begin{array}{lcc}
\mathbb{E}\left(\widehat{a}_{n}\right)=a, & \,\,\,\, & \sigma^{2}_{\widehat{a}_{n}}=\dfrac{n}{n\left\Vert x\right\Vert ^{2}-\left\langle x,e\right\rangle ^{2}}\sigma^{2}=\dfrac{\sigma^{2}}{ns^{2}_{x,n}}\\
\mathbb{E}\left(\widehat{b}_{n}\right)=b, & \,\,\,\, & \sigma^{2}_{\widehat{b}_{n}}=\dfrac{\left\Vert x\right\Vert ^{2}}{n\left\Vert x\right\Vert ^{2}-\left\langle x,e\right\rangle ^{2}}\sigma^{2}=\dfrac{\left\Vert x\right\Vert ^{2}}{n^{2}s^{2}_{x,n}}\sigma^{2}.
\end{array}\right.\label{eq:exp_var_estimators_an_bn}
\end{equation}
}The random variable $\dfrac{n\widehat{c}_{n}}{\sigma^{2}}$ follows
the chi-squared law $\chi^{2}_{n-2}.$ Consequently,\boxeq{
\begin{equation}
\mathbb{E}\left(\widehat{c}_{n}\right)=\dfrac{n-2}{n}\sigma^{2}\,\,\,\,\,\,\,\,\sigma^{2}_{\widehat{c}_{n}}=\dfrac{2\left(n-2\right)}{n^{2}}\sigma^{4}.\label{eq:exp_var_estimators_c_n}
\end{equation}
}

\end{theorem}

\begin{proof}{}{}

The random variable $\left(\widehat{a}_{n},\widehat{b}_{n},\uuline{Y}-\left(\widehat{a}_{n}x+\widehat{b}_{n}e\right)\right)$
is a linear transform of the Gaussian random variable $\uuline{Y},$
as can be seen by examining the equalities $\refpar{eq:scal_Y_u_widehat}$
and $\refpar{eq:scal_Y_u_widehat-1}.$ It is therefore itself Gaussian,
and the same holds, of course, for the random variable $\left(\widehat{a}_{n},\widehat{b}_{n}\right).$
The expectations and variances of $\widehat{a}_{n}$ and $\widehat{b}_{n}$
are easily computed from the equalities $\refpar{eq:scal_Y_u_widehat}$
and $\refpar{eq:scal_Y_u_widehat-1}.$ 

Let $V$ be the vector subspace of $\mathbb{R}^{n}$ generated by
the vectors $x$ and $e,$ and let us determine the orthogonal projection
$y_{V}=\widehat{\alpha}x+\widehat{\beta}e$ of an arbitrary vector
$y$ of $\mathbb{R}^{n}.$ This orthogonal projection is characterized
by the orthogonality relation
\[
\forall\left(u,v\right)\in\mathbb{R}^{2},\,\,\,\,\left\langle y-\left(\widehat{\alpha}x+\widehat{\beta}e\right),ux+ve\right\rangle =0,
\]
which is equivalent to
\[
\forall\left(u,v\right)\in\mathbb{R}^{2},\,\,\,\,u\left[\left\langle y,x\right\rangle -\widehat{\alpha}\left\Vert x\right\Vert ^{2}-\widehat{\beta}\left\langle x,e\right\rangle \right]+v\left[\left\langle y,e\right\rangle -\widehat{\alpha}\left\langle x,e\right\rangle -\widehat{\beta}\left\Vert e\right\Vert ^{2}\right]=0,
\]
and this, in turn, is equivalent to the system
\[
\left\{ \begin{array}{l}
\widehat{\alpha}\left\Vert x\right\Vert ^{2}+\widehat{\beta}\left\langle x,e\right\rangle =\left\langle y,x\right\rangle \\
\widehat{\alpha}\left\langle x,e\right\rangle +\widehat{\beta}\left\Vert e\right\Vert ^{2}=\left\langle y,e\right\rangle .
\end{array}\right.
\]
Hence, $\widehat{\alpha}$ and $\widehat{\beta}$ solve the system
$\refpar{eq:stat_point_cond1_scal},$ which proves that the random
variable $\uuline{Y}_{V}=\Pi_{V}\uuline{Y}$---where $\Pi_{V}$ denotes
the orthogonal projector onto $V$---satisfies
\[
\uuline{Y}_{V}=\widehat{a}_{n}x+\widehat{b}_{n}e.
\]

It follows that, pointwise,
\[
\uuline{Y}-\left(\widehat{a}_{n}x+\widehat{b}_{n}e\right)=\uuline{Y}_{V^{\perp}},
\]
that is, the orthogonal projection of $\uuline{Y}$ onto $V^{\perp}.$
Since $\uuline{Y}$ follows a Gaussian law $\mathscr{N}_{\mathbb{R}^{n}}\left(ax+be,\sigma^{2}I_{n}\right),$
the random variables $\uuline{Y}_{V}$ and $\uuline{Y}_{V^{\perp}}$
are independent---by Proposition $\ref{pr:cns_dir_sum_comp_gaus_rv_indep}$---and
thus the random variables $\uuline{Y}_{V^{\perp}}$ and $\left(\widehat{a}_{n},\widehat{b}_{n}\right)$
are also independent, since the latter is a measurable function of
$\uuline{Y}_{V}.$ This proves the independence of $\left(\widehat{a}_{n},\widehat{b}_{n}\right)$
and $\widehat{c}_{n}.$

Moreover, the dimension of $V^{\perp}$ is $n-2,$ and
\[
\Pi_{V^{\perp}}\left(ax+be\right)=0\,\,\,\,\text{and}\,\,\,\,\Lambda_{\uuline{Y}_{V^{\perp}}}=\Pi_{V^{\perp}}\Lambda_{\uuline{Y}}\Pi_{V^{\perp}}=\sigma^{2}\Pi_{V^{\perp}}.
\]
Then the law of $\uuline{Y}_{V^{\perp}}$ is the law $\mathscr{N}_{\mathbb{R}^{n}}\left(0,\sigma^{2}\Pi_{V^{\perp}}\right),$
which implies that the law of $\left\Vert \dfrac{\uuline{Y}_{V^{\perp}}}{\sigma}\right\Vert ^{2}$
is the chi-squared law $\chi^{2}_{n-2}.$ This is also the law of
$n\dfrac{\widehat{c}_{n}}{\sigma^{2}},$ since these two random variables
are equal. 

Finally,
\[
\mathbb{E}\left(n\dfrac{\widehat{c}_{n}}{\sigma^{2}}\right)=n-2\,\,\,\,\,\,\,\,\sigma^{2}\left(n\dfrac{\widehat{c}_{n}}{\sigma^{2}}\right)=2\left(n-2\right),
\]
 which immediately yields the equalities $\refpar{eq:exp_var_estimators_an_bn}.$

\end{proof}

Theorem $\ref{th:est_exp_var_laws}$ makes it possible to construct
\textbf{hypothesis} \textbf{tests} and \textbf{confidence intervals}
for the different regression parameters and, once the model has been
estimated, to perform predictions.

\subsubsection{Hypothesis Tests}

Suppose, for instance, that we wish to test the null hypothesis $H_{0}$
stating that the true variance value of the error $\epsilon$ is $\sigma^{2}.$
Fix a significant \textbf{threshold} $\alpha.$ Using a statistical
table or software, we determine the value $c_{\alpha}$ such that
$\chi^{2}_{n-2}\left(\left[c_{\alpha},+\infty\right[\right)=\alpha.$
We \textbf{reject the hypothesis\index{hypothesis rejection}} $H_{0}$
if $n\dfrac{\widehat{c}_{n}}{\sigma^{2}}>c_{\alpha},$ that is, if
$\widehat{c}_{n}>\dfrac{\sigma^{2}}{n}c_{\alpha}.$ 

Using the expressions of $\refpar{eq:scal_Y_u_widehat},$ $\refpar{eq:scal_Y_u_widehat-1}$
and $\refpar{eq:estimator_likelihood_max_variance},$ the \textbf{rejection
region\index{rejection region}} of the hypothesis $H_{0}$ can be
written explicitly as the subset of $\mathbb{R}^{n},$
\[
\left\{ \uuline{y}\in\mathbb{R}^{n}:\,\left\Vert \uuline{y}-\left(\dfrac{n\left\langle y,x\right\rangle -\left\langle x,e\right\rangle \left\langle y,e\right\rangle }{n\left\Vert x\right\Vert ^{2}-\left\langle x,e\right\rangle ^{2}}x+\dfrac{\left\Vert x\right\Vert ^{2}\left\langle y,e\right\rangle -\left\langle x,e\right\rangle \left\langle y,e\right\rangle }{n\left\Vert x\right\Vert ^{2}-\left\langle x,e\right\rangle ^{2}}e\right)\right\Vert >\sigma^{2}c_{\alpha}\right\} .
\]

We may proceed similarly to test an hypothesis on the parameter $b.$
To this end, we introduce the random variable $\widehat{B}_{n}\left(b\right),$
defined as the centered and reduced random variable associated to
the estimator $\widehat{b}_{n},$ where the unknown variance $\sigma^{2}$
is replaced by an unbiased estimator $\dfrac{n}{n-2}\widehat{c}_{n}.$
Using the expression of $\refpar{eq:exp_var_estimators_an_bn},$ this
variable is defined by
\[
\widehat{B}_{n}\left(b\right)=\left[\dfrac{\left\Vert x\right\Vert ^{2}}{n\left\Vert x\right\Vert ^{2}-\left\langle x,e\right\rangle ^{2}}\dfrac{n}{n-2}\widehat{c}_{n}\right]^{-\frac{1}{2}}\left(\widehat{b}_{n}-b\right),
\]
or equivalently
\begin{equation}
\widehat{B}_{n}\left(b\right)=\left[\dfrac{\left\Vert x\right\Vert ^{2}}{n^{2}s^{2}_{x,n}}\dfrac{n}{n-2}\widehat{c}_{n}\right]^{-\frac{1}{2}}\left(\widehat{b}_{n}-b\right).\label{eq:rv_cent_red_assoc_estimator_widehat_b_n}
\end{equation}
Recall that if $X$ and $Y$ are two independent random variables
with respective laws the Gaussian law $\mathscr{N}_{\mathbb{R}}\left(0,1\right)$
and the chi-squared law $\chi^{2}_{n},$ then the random variable
$\dfrac{\sqrt{n}X}{\sqrt{Y}}$ follows the Student law with $n$ degrees
of freedom---cf. Exercise $\ref{exo:exercise10.4}.$ Applying Theorem
$\ref{th:est_exp_var_laws},$ the random variable $\widehat{B}_{n}\left(b\right)$
follows the \textbf{Student law} $t_{n-2}$---the detailed verification
is left as an exercise for the interested reader.

To test the null hypothesis $H_{0}:$ ``the true value of the parameter
$b$ is $b_{0}$'' against the one-sided alternative $H_{1}:$ ``$b>b_{0}$'',
we proceed as follows. 

We fix a \textbf{significance threshold\index{significance threshold}}
$\alpha.$ Using a statiscal table or software, we determine the value
$b_{1-\alpha}$ such that $t_{n-2}\left(\left]-\infty,b_{1-\alpha}\right]\right)=1-\alpha.$
We \textbf{reject the hypothesis} $H_{0}$ if $\widehat{B}_{n}\left(b_{0}\right)>b_{1-\alpha}.$ 

Using the explicit expressions $\refpar{eq:scal_Y_u_widehat},$ $\refpar{eq:scal_Y_u_widehat-1}$
and $\refpar{eq:estimator_likelihood_max_variance},$ the \textbf{rejection
region} of the hypothesis $H_{0}$ against $H_{1}$ is there the subset
of $\mathbb{R}^{n},$
\[
\left\{ \uuline{y}\in\mathbb{R}^{n}:\,\left(\dfrac{\left\Vert x\right\Vert ^{2}}{n\left\Vert x\right\Vert ^{2}-\left\langle x,e\right\rangle ^{2}}\dfrac{\left\Vert \uuline{y}-\left[\widehat{a}\left(y\right)x+\widehat{b}\left(\uuline{y}\right)e\right]\right\Vert }{n-2}\right)^{-\frac{1}{2}}\left(\widehat{b}\left(\uuline{y}\right)-b_{0}\right)>b_{1-\alpha}\right\} .
\]

If we now test the same hypothesis $H_{0}$ against the two-sided
hypothesis $H_{2}:$ ``$b\neq b_{0}$'', we proceed significantly
differently. We choose $p,\,0<p<1,$ and we determine, with the help
of a table or with a statistical software, the values:
\begin{itemize}
\item $b_{1-\alpha p}$ such that $t_{n-2}\left(\left]-\infty,b_{1-\alpha p}\right]\right)=1-\alpha p$ 
\item and $b_{1-\alpha\left(1-p\right)}$ such that $t_{n-2}\left(\left]-\infty,b_{1-\alpha\left(1-p\right)}\right]\right)=1-\alpha\left(1-p\right).$ 
\end{itemize}
Since the law $t_{n-2}$ is symmetric, then
\begin{align*}
t_{n-2}\left(\left]-\infty,-b_{1-\alpha\left(1-p\right)}\right]\right) & =t_{n-2}\left(\left[b_{1-\alpha\left(1-p\right)},+\infty\right[\right)\\
 & =1-t_{n-2}\left(\left]-\infty,b_{1-\alpha\left(1-p\right)}\right]\right)\\
 & =\alpha\left(1-p\right).
\end{align*}
Then,
\[
t_{n-2}\left(\left]-\infty,-b_{1-\alpha\left(1-p\right)}\right]\cup\left[b_{1-\alpha p},+\infty\right[\right)=\alpha\left(1-p\right)+\alpha p=\alpha.
\]
We \textbf{reject the hypothesis} $H_{0}$ if $\widehat{B}_{n}\left(b_{0}\right)>b_{1-\alpha p}$
or if $\widehat{B}_{n}\left(b_{0}\right)<-b_{1-\alpha\left(1-p\right)}.$
The explicit form of the \textbf{rejection region} for this hypothesis
test of $H_{0}$ against $H_{2}$ can be written as before, but it
brings no additional insight.

Finally, we can test a hypothesis on the true value of the parameter
$a$ in exactly the same manner.

\subsubsection{Confidence Intervals}

Let us present in details the example of construction of a confidence
interval for $b$ at the level $\beta.$ Using a table or statistical
software, we determine the value $b_{1-\frac{\beta}{2}}$ such that
\[
t_{n-2}\left(\left]-\infty,b_{1-\frac{\beta}{2}}\right]\right)=1-\dfrac{\beta}{2}.
\]
Then
\[
t_{n-2}\left(\left]-\infty,-b_{1-\frac{\beta}{2}}\right]\cup\left[b_{1-\frac{\beta}{2}},+\infty\right[\right)=\beta.
\]
By Theorem $\ref{th:est_exp_var_laws},$ it follows once again that
the random variable $\widehat{B}_{n}\left(b\right)$ follows the Student
law $t_{n-2}.$ Hence,
\begin{equation}
P\left(\left|\widehat{B}_{n}\left(b\right)\right|\leqslant b_{1-\frac{\beta}{2}}\right)=1-\beta.\label{eq:P_widehatB_nlessthat b_1-alphaover2}
\end{equation}
This equality yields, at level $\beta,$ the confidence interval $\left[I,S\right],$
where
\[
I=\widehat{b}_{n}-b_{1-\frac{\beta}{2}}\left[\dfrac{\left\Vert x\right\Vert ^{2}}{n\left\Vert x\right\Vert ^{2}-\left\langle x,e\right\rangle ^{2}}\dfrac{n}{n-2}\widehat{c}_{n}\right]^{\frac{1}{2}},
\]
and
\[
S=\widehat{b}_{n}+b_{1-\frac{\beta}{2}}\left[\dfrac{\left\Vert x\right\Vert ^{2}}{n\left\Vert x\right\Vert ^{2}-\left\langle x,e\right\rangle ^{2}}\dfrac{n}{n-2}\widehat{c}_{n}\right]^{\frac{1}{2}}.
\]
Still by Theorem $\ref{th:est_exp_var_laws},$ we can construct confidence
intervals for the parameters $a$ and $\sigma^{2}.$ For instance,
for $a,$ we introduce the random variable
\begin{equation}
\widehat{A}_{n}\left(a\right)=\left[\dfrac{1}{ns^{2}_{x,n}}\dfrac{n}{n-2}\widehat{c}_{n}\right]^{-\frac{1}{2}}\left(\widehat{a}_{n}-a\right),\label{eq:rv_Ana}
\end{equation}
and we determine the value $a_{1-\frac{\alpha}{2}}$ such that 
\[
t_{n-2}\left(\left]-\infty,a_{1-\frac{\alpha}{2}}\right]\right)=1-\frac{\alpha}{2}.
\]
 Then,
\begin{equation}
P\left(\left|\widehat{A}_{n}\left(a\right)\right|\leqslant a_{1-\frac{\alpha}{2}}\right)=1-\alpha,\label{eq:Pabs_A_na_less_a1-alphaover2}
\end{equation}
and the corresponding confidence interval for $a$ is obtained in
a similar way.

\subsubsection{Prediction}

Since the theoretical model is still described by the equality $\refpar{eq:linear_model_reaction},$
we now consider the prediction of the outcome of an experiment performed
at the level $\widetilde{x}$ of the value of $x.$ We work with the
model estimated from a sample of size $n,$ and defined by
\begin{equation}
Y_{n+1}=\widehat{a}_{n}\widetilde{x}+\widehat{b}_{n}+\epsilon_{n+1},\label{eq:next_value_at_level_x_tilde}
\end{equation}
where the random variables $\widehat{a}_{n},\,\widehat{b}_{n}$ and
$\epsilon_{n}$ are those previously introduced. Our goal is to obtain
a confidence interval for $Y_{n+1}.$

First, observe that the random variables $\widehat{a}_{n}\widetilde{x}+\widehat{b}_{n}$
and $\epsilon_{n+1}$ are Gaussian and independent. Consequently,
the random variable $Y_{n+1}$ is also Gaussian. We now compute its
expectation and variance. 

Since the estimators $\widehat{a}_{n}$ and $\widehat{b}_{n}$ are
unbiased and $\epsilon_{n}$ is centered, 
\begin{equation}
\mathbb{E}\left(Y_{n+1}\right)=a\widetilde{x}+b.\label{eq:expectation_X_n_1}
\end{equation}

Moreover, by the equalities $\refpar{eq:estimates_least_squ},$
\[
\widehat{b}_{n}=\overline{Y}_{n}-\overline{x}_{n}\widehat{a}_{n}.
\]
Hence,
\[
\widehat{a}_{n}\widetilde{x}+\widehat{b}_{n}=\widehat{a}_{n}\left(\widetilde{x}-\overline{x}_{n}\right)+\overline{Y}_{n}.
\]
The random variables $\widehat{a}_{n}$ and $\overline{Y}_{n}$ are
not independent. Now, we can obtain an upper-bound of $\widehat{a}_{n}\widetilde{x}+\widehat{b}_{n}$
as follows: 
\[
\sigma^{2}_{\widehat{a}_{n}\widetilde{x}+\widehat{b}_{n}}\leqslant2\left(\sigma^{2}_{\widehat{a}_{n}\left(\widetilde{x}-\overline{x}_{n}\right)}+\sigma^{2}_{\overline{Y_{n}}}\right)=2\left(\left(\widetilde{x}-\overline{x}_{n}\right)^{2}\sigma^{2}_{\widehat{a}_{n}}+\sigma^{2}_{\overline{Y}_{n}}\right).
\]

Using the expression of the variance of $\widehat{a}_{n}$---see
the equality $\refpar{eq:exp_var_estimators_an_bn}$---and the identity
$\sigma^{2}_{\overline{Y}_{n}}=\dfrac{\sigma^{2}}{n},$ we obtain
the upper bound
\[
\sigma^{2}_{\widehat{a}_{n}\widetilde{x}+\widehat{b}_{n}}\leqslant2\sigma^{2}\left(\dfrac{\left(\widetilde{x}-\overline{x}_{n}\right)^{2}}{ns^{2}_{x,n}}+\dfrac{1}{n}\right).
\]

Since the random variables $\widehat{a}_{n}\widetilde{x}+\widehat{b}_{n}$
and $\epsilon_{n+1}$ are independent,
\begin{equation}
\sigma^{2}_{Y_{n+1}}\leqslant\sigma^{2}\left(1+\dfrac{2\left(\widetilde{x}-\overline{x}_{n}\right)^{2}}{ns^{2}_{x,n}}+\dfrac{2}{n}\right).\label{eq:upper_boud_sigma_2_x_n_and_1}
\end{equation}

Denote $\mathring{Y}_{n+1}$ the Gaussian centered random variable
$Y_{n+1}-\left(a\widetilde{x}+b\right).$ By Theorem $\ref{th:est_exp_var_laws},$
the random variables $\left(\widehat{a}_{n}\widetilde{x}+\widehat{b}_{n},\epsilon_{n+1}\right)$
and $\dfrac{n\widehat{c}_{n}}{\sigma^{2}},$ and therefore, also the
random variables $\mathring{Y}_{n+1}$ and $\dfrac{n\widehat{c}_{n}}{\sigma^{2}},$
are independent. Since $\dfrac{n\widehat{c}_{n}}{\sigma^{2}}$ follows
the chi-squared law $\chi^{2}_{n-2},$ the random variable
\[
Z_{n}=\sqrt{n-2}\dfrac{\mathring{Y}_{n+1}}{\sigma_{Y_{n+1}}}\left[\dfrac{n\widehat{c}_{n}}{\sigma^{2}}\right]^{-\frac{1}{2}}
\]
follows the Student law $t_{n-2}.$

We can then construct, in the following manner, a confidence interval
for $Y_{n+1}$ where the overall level is bounded above by $\alpha+\beta+\gamma$
and where $0<\alpha+\beta+\gamma<1.$ 

As before, using a statistical table or software, we determine the
value $z_{1-\frac{\gamma}{2}}$ such that
\[
t_{n-2}\left(\left]-\infty,-z_{1-\frac{\gamma}{2}}\right]\cup\left[z_{1-\frac{\gamma}{2}},+\infty\right[\right)=\gamma.
\]
Then,
\[
P\left(\left|Z_{n}\right|\leqslant z_{1-\frac{\gamma}{2}}\right)=1-\gamma.
\]
By the definition of $Z_{n},$ we have the equivalence
\[
\left|Z_{n}\right|\leqslant z_{1-\frac{\gamma}{2}}\,\,\,\,\Leftrightarrow\,\,\,\,\left|\mathring{Y}_{n+1}\right|\leqslant\dfrac{\sigma_{Y_{n+1}}}{\sigma}\left[\dfrac{n}{n-2}\widehat{c}_{n}\right]^{\frac{1}{2}}z_{1-\frac{\gamma}{2}}.
\]
Taking into account the upper-bound $\refpar{eq:upper_boud_sigma_2_x_n_and_1},$
we obtain the implication
\[
\left|Z_{n}\right|\leqslant z_{1-\frac{\gamma}{2}}\,\,\,\,\Longrightarrow\,\,\,\,\left|\mathring{Y}_{n+1}\right|\leqslant z_{1-\frac{\gamma}{2}}\left[\left(1+\dfrac{2\left(\widetilde{x}-\overline{x}_{n}\right)^{2}}{ns^{2}_{x,n}}+\dfrac{2}{n}\right)\left(\dfrac{n}{n-2}\widehat{c}_{n}\right)\right]^{\frac{1}{2}}.
\]
Thus,
\begin{equation}
P\left(\left|\mathring{Y}_{n+1}\right|\leqslant z_{1-\frac{\gamma}{2}}\left[\left(1+\dfrac{2\left(\widetilde{x}-\overline{x}_{n}\right)^{2}}{ns^{2}_{x,n}}+\dfrac{2}{n}\right)\left(\dfrac{n}{n-2}\widehat{c}_{n}\right)\right]^{\frac{1}{2}}\right)\geqslant1-\gamma.\label{eq:lower_bound_proba_abs_centered_Y_nand1}
\end{equation}
Hence, setting
\[
I_{n}=a\widetilde{x}+b-z_{1-\frac{\gamma}{2}}\left[\left(1+\dfrac{2\left(\widetilde{x}-\overline{x}_{n}\right)^{2}}{ns^{2}_{x,n}}+\dfrac{2}{n}\right)\left(\dfrac{n}{n-2}\widehat{c}_{n}\right)\right]^{\frac{1}{2}},
\]
and

\[
S_{n}=a\widetilde{x}+b+z_{1-\frac{\gamma}{2}}\left[\left(1+\dfrac{2\left(\widetilde{x}-\overline{x}_{n}\right)^{2}}{ns^{2}_{x,n}}+\dfrac{2}{n}\right)\left(\dfrac{n}{n-2}\widehat{c}_{n}\right)\right]^{\frac{1}{2}},
\]
we obtain
\[
P\left(Y_{n+1}\in\left[I_{n},S_{n}\right]\right)\geqslant1-\gamma,
\]
However, this does not yet provide a confidence interval for $Y_{n+1},$
since the interval $\left[I_{n},S_{n}\right]$ still depends on the
unknown parameters $a$ and $b.$ To obtain an explicit interval,
we must get the estimated values of these parameters. To keep the
presentation concise, we only state here the underlying principle
of the construction method.

After determining, using the previous method, the confidence intervals
for $a$ and $b$ at the respective levels $\alpha$ and $\beta$---by
means of the equalities $\refpar{eq:Pabs_A_na_less_a1-alphaover2}$
and $\refpar{eq:P_widehatB_nlessthat b_1-alphaover2}$---and taking
into account $\refpar{eq:lower_bound_proba_abs_centered_Y_nand1},$
we arrive at the following situation.

We have constructed the random variables $u_{n}\left(\alpha\right),\,v_{n}\left(\beta\right)$
and $w_{n}\left(\gamma\right)$ such that the following inequalities
hold simultaneously
\begin{gather*}
P\left(\left|a-\widehat{a}_{n}\right|\leqslant u_{n}\left(\alpha\right)\right)\geqslant1-\alpha,\\
P\left(\left|b-\widehat{b}_{n}\right|\leqslant v_{n}\left(\beta\right)\right)\geqslant1-\beta,\\
P\left(\left|Y_{n+1}-\left(a\widetilde{x}+b\right)\right|\leqslant w_{n}\left(\gamma\right)\right)\geqslant1-\gamma,
\end{gather*}
where
\[
\left\{ \begin{array}{l}
u_{n}\left(\alpha\right)=a_{1-\frac{\alpha}{2}}\left[\dfrac{1}{ns^{2}_{x,n}}\dfrac{n}{n-2}\widehat{c}_{n}\right]^{\frac{1}{2}}\\
v_{n}\left(\beta\right)=b_{1-\frac{\beta}{2}}\left[\dfrac{\left\Vert x\right\Vert ^{2}}{b^{2}s^{2}_{x,n}}\dfrac{n}{n-2}\widehat{c}_{n}\right]^{\frac{1}{2}}\\
w_{n}\left(\gamma\right)=z_{1-\frac{\gamma}{2}}\left[\left(1+\dfrac{2\left(\widetilde{x}-\overline{x}_{n}\right)^{2}}{ns^{2}_{x,n}}+\dfrac{2}{n}\right)\left(\dfrac{n}{n-2}\widehat{c}_{n}\right)\right]^{\frac{1}{2}}.
\end{array}\right.
\]
Now, let $A,B$ and $C$ be events satisfying
\[
P\left(A\right)\geqslant1-\alpha,\,\,\,\,P\left(B\right)\geqslant1-\beta,\,\,\,\,P\left(C\right)\geqslant1-\gamma.
\]
Then,
\[
P\left(A^{c}\cup B^{c}\cup C^{c}\right)\leqslant P\left(A^{c}\right)+P\left(B^{c}\right)+P\left(C^{c}\right)\leqslant\alpha+\beta+\gamma,
\]
and consequently
\[
P\left(A\cap B\cap C\right)\geqslant1-\left(\alpha+\beta+\gamma\right).
\]
Hence, by the triangle inequality, we obtain, with probability greater
than or equal to $1-\left(\alpha+\beta+\gamma\right),$
\begin{align*}
\left|Y_{n+1}-\left(\widehat{a}_{n}\widetilde{x}+\widehat{b}_{n}\right)\right| & \leqslant\left|Y_{n+1}-\left(a\widetilde{x}+b\right)\right|+\left|a-\widehat{a}_{n}\right|\left|\widetilde{x}\right|+\left|b-\widehat{b}_{n}\right|\\
 & \leqslant w_{n}\left(\gamma\right)+u_{n}\left(\alpha\right)\left|\widetilde{x}\right|+v_{n}\left(\beta\right),
\end{align*}
It follows that the interval $\left[\widehat{I}_{n},\widehat{S}_{n}\right]$
is a confidence interval for $Y_{n+1}$ with confidence level less
than or equal to $\alpha+\beta+\gamma.$ Here we define
\[
I_{n}=a\widetilde{x}+b-\left(u_{n}\left(\alpha\right)\left|\widetilde{x}\right|+v_{n}\left(\beta\right)+w_{n}\left(\gamma\right)\right),
\]
and
\[
S_{n}=a\widetilde{x}+b+\left(u_{n}\left(\alpha\right)\left|\widetilde{x}\right|+v_{n}\left(\beta\right)+w_{n}\left(\gamma\right)\right).
\]

We now provide several values of the inverse cumulative distribution
function of a random variable $X$ following a Student law $t_{n}$
for different values of $n.$ For fixed $n$ and $y,$ Table \ref{tab:student_val_pX}
gives the value $x$ such that $P\left(X\leqslant x\right)=y.$

\begin{table}
\begin{center}%
\begin{tabular}{|c|c|c|c|c|c|}
\hline 
n\textbackslash y &
0.75 &
0.90 &
0.95 &
0.990 &
0.995\tabularnewline
\hline 
5 &
0.727 &
1.476 &
2.015 &
3.365 &
4.032\tabularnewline
\hline 
10 &
0.700 &
1.372 &
1.812 &
2.764 &
3.169\tabularnewline
\hline 
15 &
0.691 &
1.341 &
1.753 &
2.602 &
2.947\tabularnewline
\hline 
20 &
0.687 &
1.325 &
1.725 &
2.528 &
2.845\tabularnewline
\hline 
\end{tabular}\end{center}

\caption{Values $x$ such that $P\left(X\leqslant x\right)=y$ for $X$ random
variable of Student law $t_{n}$ for given $n$ and $y.$}\label{tab:student_val_pX}
\end{table}

To conclude, observe that, for the sake of simplicity, we have restricted
our analysis to linear models with a single explanatory variable $x.$
A linear model with $k$ explanatory variables $x_{j},1\leqslant j\leqslant k,$
is a theoretical model of the form
\begin{equation}
Y=\sum^{k}_{j=1}a_{j}x_{j}+b+\epsilon,\label{eq:linear_model_k_factors}
\end{equation}
where $\epsilon$ is a centered second-order random variable, representing
approximation or measurement error. A comprehensive statistical treatment
of linear models---together with an extensive bibliography on the
subject---can be found, for example, in the first chapter of the
book by A. Antoniadis \cite{antoniadis1992regression}, which is primarily
devoted to nonlinear models.

\section*{Exercises}

\addcontentsline{toc}{section}{Exercises}

\textbf{Unless explicitly stated otherwise, all random variables are
defined on the same probabilized space $\left(\Omega,\mathscr{A},P\right).$}

\begin{exercise}{A Non Gaussian Measure with Gaussian Marginals}{exercise14.1}

Let $X$ be a real-valued random variable with Gaussian law $\mathscr{N}_{\mathbb{R}}\left(0,1\right).$
Consider the random variables taking values in $\mathbb{R}^{2},$
$Y=\left(X,-X\right)$ and $Z=\left(X,X\right).$ We study the probability
measure $\mu$ on $\mathbb{R}^{2}$ defined by 
\[
\mu=\dfrac{P_{Y}+P_{Z}}{2}.
\]
Denote $\Pi_{1}$ and $\Pi_{2}$ the coordinate functions defined
by $\Pi_{1}\left(x,y\right)=x$ and $\Pi_{2}\left(x,y\right)=y,$
for every $\left(x,y\right)\in\mathbb{R}^{2}.$ Finally, denote by
$\mu_{1}=\Pi_{1}\left(\mu\right)$ and $\mu_{2}=\Pi_{2}\left(\mu\right)$
the marginals of $\mu,$ that is, the image measures of $\mu$ by
$\Pi_{1}$ and $\Pi_{2}.$

Prove that $\mu_{1}$ and $\mu_{2}$ are equal to the Gaussian measure
$\mathscr{N}_{\mathbb{R}}\left(0,1\right).$ Compute the Fourier transform
of $\mu$ and deduce that $\mu$ is not Gaussian.

\end{exercise}

\begin{exercise}{A Non Linear Transform of a Gaussian Random Variable can be Gaussian}{exercise14.2}

Let $X,Y$ and $Z$ be three independent real-valued random variables,
with Gaussian law $\mathscr{N}_{\mathbb{R}}\left(0,1\right).$ Define
the random variable $U$ by
\[
U=\dfrac{X+YZ}{\sqrt{1+Z^{2}}}.
\]

Determine the conditional law $P^{Z=\cdot}_{U}$ of $U$ given $Z.$
Deduce from this that $U$ and $Z$ are independent and determine
the law of $U.$ Conclude.

\end{exercise}

\begin{exercise}{Characterization of Gaussian Laws on $\mathbb{R}$}{exercise14.3}

Let $X$ and $Y$ be two independent real-valued random variables,
with second order moments and common law $\mu.$ Assume that
\[
\intop_{\mathbb{R}}x\text{d}\mu\left(x\right)=0\,\,\,\,\text{and}\,\,\,\,\intop_{\mathbb{R}}x^{2}\text{d}\mu\left(x\right)=\sigma^{2}.
\]

Prove that if $\mu$ is the law $\mathscr{N}_{\mathbb{R}}\left(0,\sigma^{2}\right),$
then the random variable $\dfrac{X+Y}{\sqrt{2}}$ also has for law
$\mathscr{N}_{\mathbb{R}}\left(0,\sigma^{2}\right).$ Conversely,
prove that if the random variable $\dfrac{X+Y}{\sqrt{2}}$ has for
law $\mu,$ then $\mu$ must be the law $\mathscr{N}_{\mathbb{R}}\left(0,\sigma^{2}\right).$
For the converse, assume that $\sigma=1.$

Prove that, for every real number $t$ and every integer $n,$
\[
\widehat{\mu}\left(t\right)=\left[\widehat{\mu}\left(\dfrac{t}{2^{n}}\right)\right]^{4^{n}},
\]
and then deduce that $\widehat{\mu}\left(t\right)\neq0.$ Next, for
every $t\neq0,$
\[
h\left(t\right)=\dfrac{\ln\left|\widehat{\mu}\left(t\right)\right|}{t^{2}},
\]
and show that the function $h$ is constant. Deduce the expression
of $\left|\widehat{\mu}\left(t\right)\right|,$ and then determine
$\widehat{\mu}\left(t\right).$

\end{exercise}

We now study, in the following exercise, another characterization
of Gaussian random variables. It is a version of the Bernstein theorem
that is slightly more general than the one usually stated.

\begin{exercise}{Characterization of the Gaussian Random Variables: Bernstein Theorem}{exercise14.4}

Let $X$ and $Y$ be two independent real-valued random variables,
such that the random variables $X+Y$ and $X-Y$ are independent.
The aim of this exercise is to prove that $X$ and $Y$ are two Gaussian
random variables. To this end, we denote $\mu=P_{X},$ $\nu=P_{Y}$
and $\gamma=\mu\ast\nu.$

1. Prove that the Fourier transform $\widehat{\gamma}$ of $\gamma$
satisfies the relation
\begin{equation}
\forall\left(u,v\right)\in\mathbb{R}^{2},\,\,\,\,\widehat{\gamma}\left(u+v\right)\widehat{\gamma}\left(u-v\right)=\left[\widehat{\gamma}\left(u\right)\right]^{2}\left[\widehat{\gamma}\left(v\right)\right]^{2}\label{eq:fourier_transform_gamma_relation}
\end{equation}

2. Let $\overline{\gamma}$ be the probability defined by, for every
Borel subset $A,$ $\overline{\gamma}\left(A\right)=\gamma\left(-A\right)$
and $\delta=\gamma\ast\overline{\gamma}.$

Prove that the Fourier transform $\widehat{\delta}$ of $\delta$
verifies the relation
\begin{equation}
\forall\left(u,v\right)\in\mathbb{R}^{2},\,\,\,\,\widehat{\delta}\left(u+v\right)\widehat{\delta}\left(u-v\right)=\left[\widehat{\delta}\left(u\right)\right]^{2}\left[\widehat{\delta}\left(v\right)\right]^{2}\label{eq:fourier_transform_delta_relation}
\end{equation}
and that the set $G=\left\{ t\in\mathbb{R}:\,\widehat{\delta}\left(t\right)\neq0\right\} $
is a group. Deduce from this that the Fourier transform $\widehat{\delta}$
never vanishes. Determine $\widehat{\delta}$ and then $\left|\widehat{\gamma}\right|.$

3. For every real number $t,$ define 
\[
g\left(t\right)=\dfrac{\widehat{\gamma}\left(t\right)}{\text{\ensuremath{\left|\widehat{\gamma}\left(t\right)\right|}}}.
\]
 Prove that $g$ satisfies the relation
\begin{equation}
\forall\left(u,v\right)\in\mathbb{R}^{2},\,\,\,\,g^{2}\left(u+v\right)=g^{2}\left(u\right)g^{2}\left(v\right).\label{eq:relation_on_g}
\end{equation}

4. Let $\Phi$ be a Borel function from $\mathbb{R}$ to $\mathbb{C}$
such that:
\begin{itemize}
\item For every real number $t,$ $\left|\Phi\left(t\right)\right|=1,$
\item For every real numbers $s$ and $t,$
\begin{equation}
\Phi\left(s+t\right)=\Phi\left(s\right)\Phi\left(t\right).\label{eq:phi_s_plus_t}
\end{equation}
\end{itemize}
Prove that there exists a real number $c$ such that, for every real
number $t,$ $\Phi\left(t\right)=\text{e}^{\text{i}ct}.$ To this
end, show that there exists a real number $a$ such that $f\left(a\right)\neq0,$
where $f$ is the function defined, for every real number $x,$ by
\[
f\left(x\right)=\intop^{x}_{0}\Phi\left(t\right)\text{d}t,
\]
and then note that, for every real number $x,$
\[
\Phi\left(x\right)=\dfrac{f\left(x+a\right)-f\left(x\right)}{f\left(a\right)}.
\]
Deduce that $\Phi$ is differentiable, which allows one to conclude.

5. Deduce that there exists a real number $m$ and a real number $a>0$
such that, for every real number $t,$
\begin{equation}
\widehat{\gamma}\left(t\right)=\text{e}^{\text{i}mt-a\frac{t^{2}}{2}}.\label{eq:expression_transform_gamma}
\end{equation}
Then prove that the random variables $X$ and $Y$ are Gaussian.

6. Generalize this result to random variables taking values in $\mathbb{R}^{d}.$

\end{exercise}

\begin{exercise}{A Characterization of the Gaussian Law in Terms of Expectation and Empirical Variance}{exercise14.5}

Let $X_{1},X_{2},\cdots,X_{n}$ be independent real-valued random
variables with common law $\mu$ such that
\[
\intop_{\mathbb{R}}x^{2}\text{d}\mu\left(x\right)<+\infty.
\]

Define the random variables, respectively called empirical mean and
variance, by
\[
M_{n}=\dfrac{1}{n}\sum^{n}_{i=1}X_{i}\,\,\,\,\text{and}\,\,\,\,\Sigma_{n}=\dfrac{1}{n}\sum^{n}_{i=1}X^{2}_{i}-M^{2}_{n}.
\]
Denote by $X$ the random variable in $\mathbb{R}^{n}$ given by 
\[
X=\left(X_{1},X_{2},\cdots,X_{n}\right).
\]

1. Suppose that $\mu$ is a Gaussian law $\mathscr{N}_{\mathbb{R}}\left(m,\sigma^{2}\right).$
What is the law of $X?$ Let $C$ be an orthogonal $n\times n$ matrix
such that, for every $j\in\left\llbracket 1,n\right\rrbracket ,$
\[
C_{1j}=\dfrac{1}{\sqrt{n}}.
\]
Express $M_{n}$ and $\Sigma_{n}$ in terms of the components of $CX,$
and deduce that $M_{n}$ and $\Sigma_{n}$ are independent random
variables. 

In the case, where $m=0$ and $\sigma=1,$ specify the laws of the
random variables $M_{n}$ and $n\Sigma_{n}.$

2. We now consider the converse. Suppose that the random variables
$M_{n}$ and $\Sigma_{n}$ are independent. For simplicity, assume
that the random variables $X_{i}$ are centered. Denote 
\[
\sigma^{2}=\intop_{\mathbb{R}}x^{2}\text{d}\mu\left(x\right).
\]
Let $\varphi$ be the Fourier transform of $\mu,$ then $S_{n}=nM_{n}$
and $V_{n}=n\Sigma_{n}.$

(a) Compute the expectation $\mathbb{E}\left(V_{n}\right)$ as a function
of $\sigma.$

(b) Prove that the function from $\mathbb{R}^{2}$ to $\mathbb{C},$
defined by
\[
\left(u,v\right)\mapsto\mathbb{E}\left(\text{e}^{\text{i}\left(uS_{n}+vV_{n}\right)}\right)
\]
is differentiable. Justify the relation
\begin{equation}
\forall u\in\mathbb{R},\,\,\,\,\mathbb{E}\left(V_{n}\text{e}^{\text{i}\left(uS_{n}\right)}\right)=\left[\varphi\left(u\right)\right]^{n}\mathbb{E}\left(V_{n}\right).\label{eq:expectation_V_n_exp_iuS_n}
\end{equation}

(c) Compute $\mathbb{E}\left(V_{n}\text{e}^{\text{i}\left(uS_{n}\right)}\right)$
using the first and second derivatives of $\varphi.$

(d) Deduce, from relation $\refpar{eq:expectation_V_n_exp_iuS_n},$
that $\varphi$ satisfies the differential equation
\[
\dfrac{\varphi^{\prime\prime}}{\varphi}-\left(\dfrac{\varphi^{\prime}}{\varphi}\right)^{2}=-\sigma^{2},
\]
and conclude that $\mu$ is the Gaussian law $\mathscr{N}_{\mathbb{R}}\left(0,\sigma^{2}\right).$

(e) Prove that the same conclusion holds when the random variables
$X_{i}$ are no longer assumed to be centered.

\end{exercise}

\begin{exercise}{Signal Theory Initiation}{exercise14.6}

Let $S$---the signal---and $V$---the noise---be two independent
real-valued random variables, with respective Gaussian laws $\mathscr{N}_{\mathbb{R}}\left(m,\sigma^{2}\right)$
and $\mathscr{N}_{\mathbb{R}}\left(0,t\right),$ where $m$ is an
arbitrary real number, and $\sigma$ and $t$ are positive real numbers.
The observation is the random variable
\[
R=tS+V.
\]

Compute the best least-squares approximation of the signal given the
observation, that is, the conditional expectation $\mathbb{E}^{\sigma\left(R\right)}\left(S\right),$
where $\sigma\left(R\right)$ denotes the $\sigma-$algebra generated
by $R.$ To this end, determine constants $a$ and $b$ such that
the random variables $aR+bS$ and $R$ are independent, and deduce
from this the conditional expectation $\mathbb{E}^{\sigma\left(R\right)}\left(S\right).$

\end{exercise}

\begin{exercise}{Signal Theory Initiation (following)}{exercise14.7}

Let $t_{1},t_{2},\cdots,t_{n}$ be $n$ nonnegative real numbers.
Let $S$---the signal---and, for $i\in\left\llbracket 1,n\right\rrbracket ,\,W_{t_{i}}$---the
noise at time $t_{i}$---be Gaussian, independent, real-valued random
variables with respective laws $\mathscr{N}_{\mathbb{R}}\left(m,\sigma^{2}\right)$
and $\mathscr{N}_{\mathbb{R}}\left(0,t_{i}\right),$ where $m$ is
an arbitrary real number and $\sigma$ is a positive real number.
The observation at time $t_{i}$ is the random variable
\[
R_{i}=t_{i}S+W_{t_{i}}.
\]
Equip $\mathbb{R}^{n}$ with its canonical basis and denote by $t$
the vector with components $t_{1},t_{2},\cdots,t_{n}.$ Denote by
$R$ and $W$ the $\mathbb{R}^{n}-$valued random variables with respective
components $R_{1},R_{2},\cdots,R_{n}$ and $W_{t_{1}},W_{t_{2}},\cdots,W_{t_{n}},$
so that
\[
R=St+W.
\]
Compute the best least-squares approximation of the signal given the
observations $R_{1},R_{2},\cdots,R_{n},$ that is, the conditional
expectation $\mathbb{E}^{\sigma\left(R\right)}\left(S\right),$ where
$\sigma\left(R\right)$ denotes the $\sigma-$algebra generated by
$R.$ To this end, choose a vector $u$ of $\mathbb{R}^{n}$ and a
constant $b$ such that the random variables $\left\langle u,R\right\rangle +bS$
and $R$ are independent, and deduce from this the conditional expectation
$\mathbb{E}^{\sigma\left(R\right)}\left(S\right).$

\end{exercise}

\begin{exercise}{Quadratic Form of a Gaussian Random Variable: Cochran Theorem}{exercise14.8}

Let $X$ be a random variable taking values in an Euclidean space
$E$ of dimension $d,$ with Gaussian law $\mathscr{N}_{E}\left(0,\text{I}\right),$
where $\text{I}$ denotes the identity operator on $E.$ Let $a$
be a unit vector of $E,$ and let $U$ and $V$ be real-valued random
variables defined by
\[
U=\left\langle X,a\right\rangle \,\,\,\,\text{and}\,\,\,\,V=\left\Vert X\right\Vert ^{2}-\left\langle X,a\right\rangle ^{2}.
\]

1. Prove that the random variables $U$ and $V$ are independent and
determine their law.

2. Let $Y$ be a random variable taking values in $E$ with Gaussian
law $\mathscr{N}_{E}\left(m,\text{I}\right),$ where $m\in E.$ Deduce
from the previous question that the law of $\left\Vert Y\right\Vert ^{2}$
is the convolution of the chi-squared law with $d-1$ degrees of freedom
and the law of the square of a Gaussian real-valued random variable
with law $\mathscr{N}_{\mathbb{R}}\left(\left\Vert m\right\Vert ,\text{I}\right).$

\end{exercise}

\begin{exercise}{Empirical Mean and Variance}{exercise14.9}

Let $Y_{1},Y_{2},\cdots,Y_{n}$ be $n$ independent random variables
taking values in a Euclidean space $E$ of dimension $d,$ each following
a Gaussian law $\mathscr{N}_{E}\left(0,\text{I}\right)$ where $\text{I}$
denotes the identity operator on $E.$ Denote by $Y=\left(Y_{1},Y_{2},\cdots,Y_{n}\right)$
the random variable taking values in $E^{n},$ equipped with the Euclidean
product structure.

1. Prove that the random variables $\left\Vert Y\right\Vert ^{2}$
and $\dfrac{Y}{\left\Vert Y\right\Vert }$ are independent.

\textit{Hint: First consider the case where $E=\mathbb{R},$ by performing
a change of variables to spherical coordinates, and then deduce the
general case.}

Specify their laws in the case where $E=\mathbb{R}.$

2. Let $X_{1},X_{2},\cdots,X_{n}$ be $n$ independent real-valued
random variables, with Gaussian law $\mathscr{N}_{\mathbb{R}}\left(m,\sigma^{2}\right),$
where $m\in\mathbb{R}$ and $\sigma>0.$ Define the real-valued random
variables $M$ and $V,$ and the random variables $X$ and $X^{\prime}$
taking values in $\mathbb{R}^{n},$ by
\begin{gather*}
M=\dfrac{1}{n}\sum^{n}_{j=1}X_{j}\,\,\,\,\text{and}\,\,\,\,V=\sum^{n}_{j=1}\left(X_{j}-M\right)^{2},\\
X=\left(X_{1},X_{2},\cdots,X_{n}\right)\,\,\,\,\text{and}\,\,\,\,X^{\prime}=\left(X_{1}-M,X_{2}-M,\cdots,X_{n}-M\right).
\end{gather*}
Set $Z=\dfrac{X^{\prime}}{\sqrt{V}}.$ Note that $M=\dfrac{1}{n}\left\langle X,e\right\rangle ,$
where $e$ is the vector of $\mathbb{R}^{n}$ of components $\left(1,1,\cdots,1\right).$

(a) Prove that the random variables $M$ and $X^{\prime}$ are independent.

(b) Compute the auto-covariance operator of $X^{\prime}.$

(c) Let $H$ be the hyperplane orthogonal to $e.$ Prove that there
exists an isometry $B$ from $\mathbb{R}^{n-1}$ onto $H$ and a random
variable $U$ taking values in $\mathbb{R}^{n-1}$ with law $\mathscr{N}_{\mathbb{R}^{n-1}}\left(0,\sigma^{2}\boldsymbol{1}_{\mathbb{R}^{n-1}}\right)$
such that $X^{\prime}=BU$ $P-$almost surely.

(d) Deduce that the random variables $M,V,$ and $Z$ are independent.

(e) Determine the laws of $M$ and $\dfrac{1}{\sigma^{2}}V.$

\end{exercise}

\section*{Solutions of Exercises}

\addcontentsline{toc}{section}{Solutions of Exercises}

\begin{solution}{}{solexercise14.1}

For every $f\in\mathscr{C}^{+}_{\mathscr{K}}\left(\mathbb{R}\right),$
it follows from the transfer theorem, from the definition of $\mu$
and from the fact that $\Pi_{1}\left(Y\right)=\Pi_{1}\left(Z\right)=X,$
that
\[
\intop_{\mathbb{R}}f\text{d}\mu_{1}=\intop_{\mathbb{R}^{2}}f\circ\Pi_{1}\text{d}\mu=\dfrac{1}{2}\intop_{\Omega}\left[f\circ\Pi_{1}\left(Y\right)+f\circ\Pi_{1}\left(Z\right)\right]\text{d}P=\intop_{\Omega}f\left(X\right)\text{d}P,
\]
Thus,
\[
\intop_{\mathbb{R}}f\text{d}\mu_{1}=\intop_{\mathbb{R}}f\text{d}P_{X},
\]
which proves that $\mu_{1}=P_{X}=\mathscr{N}_{\mathbb{R}}\left(0,1\right).$

Similarly,
\[
\intop_{\mathbb{R}}f\text{d}\mu_{2}=\intop_{\mathbb{R}^{2}}f\circ\Pi_{2}\text{d}\mu=\dfrac{1}{2}\intop_{\Omega}\left[f\circ\Pi_{2}\left(Y\right)+f\circ\Pi_{2}\left(Z\right)\right]\text{d}P,
\]
hence
\[
\intop_{\mathbb{R}}f\text{d}\mu_{2}=\dfrac{1}{2}\intop_{\Omega}\left[f\left(-X\right)+f\left(X\right)\right]\text{d}P.
\]
Since $P_{X}$ is symmetric,
\[
\intop_{\mathbb{R}}f\text{d}\mu_{2}=\intop_{\mathbb{R}}f\text{d}P_{X},
\]
which proves that $\mu_{2}=P_{X}=\mathscr{N}_{\mathbb{R}}\left(0,1\right).$

Moreover, the Fourier transform $\widehat{\mu}$ of $\mu$ is given
by
\[
\widehat{\mu}=\dfrac{1}{2}\left(\widehat{P_{Y}}+\widehat{P_{Z}}\right)=\dfrac{1}{2}\left(\varphi_{Y}+\varphi_{Z}\right).
\]
Thus, for every $u\in\mathbb{R}^{2},$
\[
\widehat{\mu}\left(u\right)=\dfrac{1}{2}\left[\text{e}^{-\frac{\left(u_{1}-u_{2}\right)^{2}}{2}}+\text{e}^{-\frac{\left(u_{1}+u_{2}\right)^{2}}{2}}\right],
\]
After simplification, we obtain
\[
\widehat{\mu}\left(u\right)=\dfrac{1}{2}\text{e}^{-\frac{\left(u_{1}+u_{2}\right)^{2}}{2}}\left[\text{e}^{u_{1}u_{2}}+\text{e}^{-u_{1}u_{2}}\right],
\]
or equivalentely,\boxeq{
\[
\widehat{\mu}\left(u\right)=\text{e}^{-\frac{\left(u_{1}+u_{2}\right)^{2}}{2}}\text{ch}\left(u_{1}u_{2}\right).
\]
}

Hence, \textbf{the measure $\mu$ is not Gaussian.}

\begin{remark}{}{}

By injectivity of the Fourier transform, we recover that $\mu_{1}$
and $\mu_{2}$ are Gaussian, since
\[
\widehat{\mu_{1}}\left(u_{1}\right)=\widehat{\mu}\left(u_{1},0\right)=\text{e}^{-\frac{u^{2}_{1}}{2}}\,\,\,\,\text{and}\,\,\,\,\widehat{\mu_{2}}\left(u_{2}\right)=\widehat{\mu}\left(0,u_{2}\right)=\text{e}^{-\frac{u^{2}_{2}}{2}}.
\]

\end{remark}

\end{solution}

\begin{solution}{}{solexercise14.2}

A conditional law $P^{Z=\cdot}_{U}$ of $U$ given $Z$ is provided
by the kernel defined, for $P_{Z}-$almost every real number $z,$
by
\[
P^{Z=z}_{U}=P^{Z=z}_{\frac{X+Yz}{\sqrt{1+z^{2}}}}.
\]
Hence, by the independence of the random variables $X+Yz$ and $Z,$
\[
P^{Z=z}_{U}=P_{\frac{X+Yz}{\sqrt{1+z^{2}}}}=\mathscr{N}_{\mathbb{R}}\left(0,\sigma^{2}_{\frac{X+Yz}{\sqrt{1+z^{2}}}}\right).
\]

Since $X$ and $Y$ are independent,
\[
\sigma^{2}_{\frac{X+Yz}{\sqrt{1+z^{2}}}}=\dfrac{1}{1+z^{2}}\left(\sigma^{2}_{X}+z^{2}\sigma^{2}_{Y}\right)=1.
\]
It follows that, for $P_{Z}-$almost every real number $z,$
\[
P^{Z=z}_{U}=\mathscr{N}_{\mathbb{R}}\left(0,1\right),
\]
law that does not depend on $z.$ This shows that \textbf{the random
variables $U$ and $Z$ are independent} and that
\[
P^{Z=z}_{U}=P_{U}=\mathscr{N}_{\mathbb{R}}\left(0,1\right),
\]
that is, the law of $U$ is a Gaussian law $\mathscr{N}_{\mathbb{R}}\left(0,1\right).$

To conclude, \textbf{the random variable $U,$ which is a non linear
transformation of the Gaussian random variable $\left(X,Y,Z\right)$
taking values in $\mathbb{R}^{3}$ is Gaussian.} Moreover, \textbf{the
random variable $\left(U,Z\right)$ is Gaussian}, taking values in
$\mathbb{R}^{2},$ since $U$ and $Z$ are independent and both follows
the Gaussian law $\mathscr{N}_{\mathbb{R}}\left(0,1\right).$

\end{solution}

\begin{solution}{}{solexercise14.3}

If $\mu$ is the law $\mathscr{N}_{\mathbb{R}}\left(0,\sigma^{2}\right),$
and since $X$ and $Y$ are independent random variables with law
$\mathscr{N}_{\mathbb{R}}\left(0,\sigma^{2}\right),$ the random variable
$\left(X,Y\right)$ is Gaussian. Consequentely, the random variable
$\dfrac{X+Y}{\sqrt{2}}$ is also Gaussian. 

Moreover
\[
\mathbb{E}\left(\dfrac{X+Y}{\sqrt{2}}\right)=0\,\,\,\,\text{and}\,\,\,\,\sigma^{2}_{\frac{X+Y}{\sqrt{2}}}=\dfrac{1}{2}\left(\sigma^{2}_{X}+\sigma^{2}_{Y}\right)=\sigma^{2}.
\]
Therefore, 
\[
P_{\frac{X+Y}{\sqrt{2}}}=\mu.
\]

Conversely, we may assume without loss of generality that $\sigma=1.$
By independence of the random variables $X$ and $Y,$ the characteristic
function $\varphi_{\frac{X+Y}{\sqrt{2}}}$ satisfies, for every real
number $t,$
\[
\varphi_{\frac{X+Y}{\sqrt{2}}}\left(t\right)=\varphi_{X}\left(\dfrac{t}{\sqrt{2}}\right)\varphi_{Y}\left(\dfrac{t}{\sqrt{2}}\right)=\left[\widehat{\mu}\left(\dfrac{t}{\sqrt{2}}\right)\right]^{2}.
\]
The assumption that the random variable $\dfrac{X+Y}{\sqrt{2}}$ has
law $\mu$ is therefore equivalent to
\[
\forall t\in\mathbb{R},\,\,\,\,\widehat{\mu}\left(t\right)=\left[\widehat{\mu}\left(\dfrac{t}{\sqrt{2}}\right)\right]^{4}.
\]
which implies that, for every $t\in\mathbb{R},$
\[
\widehat{\mu}\left(t\right)=\left[\widehat{\mu}\left(\dfrac{t}{2}\right)\right]^{4}.
\]
By iteration, it follows that for every real number $t$ and every
integer $n,$
\begin{equation}
\widehat{\mu}\left(t\right)=\left[\widehat{\mu}\left(\dfrac{t}{2^{n}}\right)\right]^{4^{n}}.\label{eq:fourier_transform_mu_n}
\end{equation}

Suppose now that there exists $t_{0}$ such that $\widehat{\mu}\left(t_{0}\right)=0.$
Then, by $\refpar{eq:fourier_transform_mu_n},$ for every nonnegative
integer $n,$ we would have 
\[
\widehat{\mu}\left(\dfrac{t_{0}}{2^{n}}\right)=0.
\]
By continuity of $\widehat{\mu},$ this would imply 
\[
\widehat{\mu}\left(0\right)=0,
\]
which is impossible since $\widehat{\mu}\left(0\right)=1.$ Hence,
the Fourier transform $\widehat{\mu}$ never vanishes. 

Then, by $\refpar{eq:fourier_transform_mu_n},$ it follows that, for
every nonnegative integer $n,$
\[
h\left(t\right)=\dfrac{\ln\left|\widehat{\mu}\left(t\right)\right|}{t^{2}}=\dfrac{4^{n}\ln\left|\widehat{\mu}\left(\dfrac{t}{2^{n}}\right)\right|}{t^{2}}=\dfrac{\ln\left|\widehat{\mu}\left(\dfrac{t}{2^{n}}\right)\right|}{\left(\dfrac{t}{2^{n}}\right)^{2}}=h\left(\dfrac{t}{2^{n}}\right).
\]
Since the random variable $X$ has expectation 0 and variance 1, the
Fourier transform $\widehat{\mu}$ admits the second-order Taylor
expansion near 0 given by
\[
\widehat{\mu}\left(t\right)=1-\dfrac{t^{2}}{2}+o\left(t^{2}\right).
\]
Hence,
\[
\left|\widehat{\mu}\left(t\right)\right|=1-\dfrac{t^{2}}{2}+o\left(t^{2}\right),
\]
which implies
\[
\ln\left|\widehat{\mu}\left(t\right)\right|=-\dfrac{t^{2}}{2}+o\left(t^{2}\right)\,\,\,\,\text{and}\,\,\,\,h\left(t\right)=-\dfrac{1}{2}+o\left(1\right).
\]

Using again $\refpar{eq:fourier_transform_mu_n},$ we obtain for every
real number $t\neq0,$ the asymptotic expansion in $n,$
\[
h\left(t\right)=h\left(\dfrac{t}{2^{n}}\right)=-\dfrac{1}{2}+o\left(1\right).
\]

Thus, for every real number $t\neq0,$ $h\left(t\right)=-\dfrac{1}{2},$
and consequently\boxeq{
\[
\left|\widehat{\mu}\left(t\right)\right|=\text{e}^{-\frac{t^{2}}{2}}.
\]
}

Finally, for every real number $t\neq0,$ write $\widehat{\mu}\left(t\right)$
in polar form
\[
\widehat{\mu}\left(t\right)=\left|\widehat{\mu}\left(t\right)\right|\text{e}^{\text{i}g\left(t\right)}.
\]
Using $\refpar{eq:fourier_transform_mu_n},$ we have for every integer
$n,$
\[
\text{e}^{\text{i}g\left(t\right)}=\text{e}^{\frac{t^{2}}{2}}\widehat{\mu}\left(t\right)=\text{e}^{\frac{t^{2}}{2}}\left[\widehat{\mu}\left(\dfrac{t}{2^{n}}\right)\right]^{4^{n}}.
\]
which yields the asymptotic expansion as $n$ tends to infinity
\[
\text{e}^{\text{i}g\left(t\right)}=\text{e}^{\frac{t^{2}}{2}}\left[1-\dfrac{t^{2}}{2\cdot4^{n}}+o\left(\dfrac{t^{2}}{4^{n}}\right)\right]^{4^{n}}.
\]
Since the right-hand side converges to 1, we conclude that for every
real number $t\neq0,$ $\text{e}^{\text{i}g\left(t\right)}=1.$ Therefore,
for every real number $t,$\boxeq{
\[
\widehat{\mu}\left(t\right)=\text{e}^{-\frac{t^{2}}{2}},
\]
}which shows that $\mu$ is the Gaussian law $\mathscr{N}_{\mathbb{R}}\left(0,1\right).$

\end{solution}

\begin{solution}{}{solexercise14.4}

\textbf{1. Proof that $\widehat{\gamma}$ satisfies $\forall\left(u,v\right)\in\mathbb{R}^{2},\,\,\,\,\widehat{\gamma}\left(u+v\right)\widehat{\gamma}\left(u-v\right)=\left[\widehat{\gamma}\left(u\right)\right]^{2}\left[\widehat{\gamma}\left(v\right)\right]^{2}$}

By independence of the random variables $X+Y$ and $X-Y$ from the
one hand, and $X$ and $Y$ on the other hand, their characteristic
functions satisfies for every real numbers $u$ and $v,$
\[
\varphi_{\left(X+Y,X-Y\right)}\left(u,v\right)=\varphi_{X+Y}\left(u\right)\varphi_{X-Y}\left(v\right)=\varphi_{X}\left(u\right)\varphi_{Y}\left(u\right)\varphi_{X}\left(v\right)\overline{\varphi_{Y}\left(v\right)}.
\]
Since for every real numbers $u$ and $v,$ we have the relation
\[
\varphi_{\left(X+Y,X-Y\right)}\left(u,v\right)=\varphi_{\left(X,Y\right)}\left(u+v,u-v\right)=\varphi_{X}\left(u+v\right)\varphi_{Y}\left(u-v\right),
\]
we deduce for every real numbers $u$ and $v,$
\[
\varphi_{X}\left(u+v\right)\varphi_{Y}\left(u-v\right)=\varphi_{X}\left(u\right)\varphi_{Y}\left(u\right)\varphi_{X}\left(v\right)\overline{\varphi_{Y}\left(v\right)}.
\]
Replacing $v$ by $-v$ in this identity gives, for every $u,v\in\mathbb{R},$
\[
\varphi_{X}\left(u-v\right)\varphi_{Y}\left(u+v\right)=\varphi_{X}\left(u\right)\varphi_{Y}\left(u\right)\overline{\varphi_{X}\left(v\right)}\varphi_{Y}\left(v\right).
\]
Multiplying these two last equalities term by term, and using
\[
\widehat{\gamma}=\widehat{\mu}\widehat{v}=\varphi_{X+Y}=\varphi_{X}\varphi_{Y},
\]
we obtain the relation\boxeq{
\begin{equation}
\forall\left(u,v\right)\in\mathbb{R}^{2},\,\,\,\,\widehat{\gamma}\left(u+v\right)\widehat{\gamma}\left(u-v\right)=\left[\widehat{\gamma}\left(u\right)\right]^{2}\left[\widehat{\gamma}\left(v\right)\right]^{2}.\label{eq:tr_gamma_u+v_tr_gamma_u-v}
\end{equation}
}

\textbf{2. Proof that $\widehat{\delta}$ satisfies $\forall\left(u,v\right)\in\mathbb{R}^{2},\,\,\,\,\widehat{\delta}\left(u+v\right)\widehat{\delta}\left(u-v\right)=\left[\widehat{\delta}\left(u\right)\right]^{2}\left[\widehat{\delta}\left(v\right)\right]^{2}.$
$G=\left\{ t\in\mathbb{R}:\,\widehat{\delta}\left(t\right)\neq0\right\} $
is a group. $\widehat{\delta}$ never vanishes. Determination of $\widehat{\delta}$
and $\left|\widehat{\gamma}\right|$}

The transfer theorem allows to state the equality of the Fourier transform
of $\gamma$ and $\overline{\gamma},$
\begin{equation}
\forall t\in\mathbb{R},\,\,\,\,\widehat{\overline{\gamma}}\left(t\right)=\overline{\widehat{\gamma}}\left(t\right).\label{eq:four_tr_gamma_equality_with_conj}
\end{equation}
Therefore, for every real number $t,$
\begin{equation}
\widehat{\delta}\left(t\right)=\widehat{\gamma}\left(t\right)\widehat{\overline{\gamma}}\left(t\right)=\widehat{\gamma}\left(t\right)\overline{\widehat{\gamma}\left(t\right)}=\left|\widehat{\gamma}\left(t\right)\right|^{2}\geqslant0.\label{eq:f_transfor_delta}
\end{equation}
By taking the modulus in the relation $\refpar{eq:tr_gamma_u+v_tr_gamma_u-v},$
we obtain for every real numbers $u$ and $v,$
\begin{equation}
\widehat{\delta}\left(u+v\right)\widehat{\delta}\left(u-v\right)=\left[\widehat{\delta}\left(u\right)\right]^{2}\left[\widehat{\delta}\left(v\right)\right]^{2}.\label{eq:f_tr_delta_u+v_f_tr_delta_u-v}
\end{equation}

It follows that if $u$ and $v$ are such that $\widehat{\delta}\left(u\right)$
and $\widehat{\delta}\left(v\right)$ are different from 0, then the
above relation implies $\widehat{\delta}\left(u+v\right)\neq0$ and
$\widehat{\delta}\left(-u\right)\neq0.$ Moreover as $\widehat{\delta}\left(0\right)=1,$
\textbf{$G$ is a group}.

Additionally, since $\widehat{\delta}$ is continuous, $G$ is open
and thus identical to $\mathbb{R}.$ Thus, $\widehat{\delta}$ never
vanishes and is positive. 

Define, for every real number $t,$
\[
f\left(t\right)=-\ln\widehat{\delta}\left(t\right).
\]
The relation $\refpar{eq:f_tr_delta_u+v_f_tr_delta_u-v}$ yields for
every real numbers $u$ and $v,$
\begin{equation}
f\left(u+v\right)+f\left(u-v\right)=2\left[f\left(u\right)+f\left(v\right)\right].\label{eq:f_u+v_plus_f_u-v}
\end{equation}
Since $f$ is continuous, nonnegative, and satisfies $f\left(0\right)=0,$
it follows that there exists\footnotemark $a>0$ such that, for every
real number $u,$ $f\left(u\right)=au^{2}.$ Therefore $\widehat{\delta}\left(u\right)=\text{e}^{-au^{2}}.$ 

Using the relation $\refpar{eq:f_transfor_delta},$ we conclude that
there exists $a>0$ such that\boxeq{
\begin{equation}
\forall t\in\mathbb{R},\,\,\,\,\left|\widehat{\gamma}\left(t\right)\right|=\text{e}^{-\frac{at^{2}}{2}}.\label{eq:module_tr_f_gamma}
\end{equation}
}

\textbf{3. Proof that $g$ verifies $\forall\left(u,v\right)\in\mathbb{R}^{2},\,\,\,\,g^{2}\left(u+v\right)=g^{2}\left(u\right)g^{2}\left(v\right).$}

Taking into account the modules in the relation $\refpar{eq:tr_gamma_u+v_tr_gamma_u-v},$
it follows, for every real numbers $u$ and $v,$
\begin{equation}
\left|\widehat{\gamma}\left(u+v\right)\right|\left|\widehat{\gamma}\left(u-v\right)\right|=\left|\widehat{\gamma}\left(u\right)\right|^{2}\left|\widehat{\gamma}\left(v\right)\right|^{2}.\label{eq:modules_f_tr_gamma_u+v_u-v}
\end{equation}
Dividing term by term the relations $\refpar{eq:tr_gamma_u+v_tr_gamma_u-v}$
and $\refpar{eq:modules_f_tr_gamma_u+v_u-v},$ we get
\[
\forall\left(u,v\right)\in\mathbb{R}^{2},\,\,\,\,g\left(u+v\right)g\left(u-v\right)=g^{2}\left(u\right).
\]
Exchanging the roles of $u$ and $v,$ we also have
\[
\forall\left(u,v\right)\in\mathbb{R}^{2},\,\,\,\,g\left(u+v\right)g\left(v-u\right)=g^{2}\left(v\right).
\]
Multiplying these two equalities term by term, and using the facts
that for every $t,$ $g\left(-t\right)=\overline{g\left(t\right)},$
and $\left|g\left(t\right)\right|=1,$ we obtain\boxeq{
\begin{equation}
\forall\left(u,v\right)\in\mathbb{R}^{2},\,\,\,\,g^{2}\left(u+v\right)=g^{2}\left(u\right)g^{2}\left(v\right).\label{eq:g_square_u_plus_v}
\end{equation}
}

\textbf{4. Proof that there exists $c$ such that $\forall t,$ $\Phi\left(t\right)=\text{e}^{\text{i}ct}.$ }

If we had $f\left(a\right)=0$ for every $a,$ then both the functions
$\text{Re}\left(\Phi\right)$ and $\text{Im}\left(\Phi\right)$ would
vanish almost everywhere, which is impossible since $\left|\Phi\right|=1.$ 

Choose $a$ such that $f\left(a\right)\neq0.$ We have
\[
f\left(x+a\right)-f\left(x\right)=\intop^{x+a}_{x}\Phi\left(t\right)\text{d}t.
\]
With the change of variables $s=t-x,$ and using the semigroup property
$\refpar{eq:phi_s_plus_t}$ of $\Phi,$ we get
\[
f\left(x+a\right)-f\left(x\right)=\intop^{a}_{0}\Phi\left(x+s\right)\text{d}s=\Phi\left(x\right)\intop^{a}_{0}\Phi\left(s\right)\text{d}s=\Phi\left(x\right)f\left(a\right).
\]
Hence,
\begin{equation}
\Phi\left(x\right)=\dfrac{f\left(x+a\right)-f\left(x\right)}{f\left(a\right)}.\label{eq:Phi_as_growth_toll}
\end{equation}
The continuity of $f$ implies the one of $\Phi.$ Then the function
$f,$ defined as function of the upper-bound of the integral of $\Phi$
is differentiable, and, by $\refpar{eq:Phi_as_growth_toll},$ the
function $\Phi$ is differentiable as well. 

Differentiating with respect to $s$ the identity $\Phi\left(s+t\right)=\Phi\left(s\right)\Phi\left(t\right),$
we obtain, for every $s$ and $t,$
\[
\Phi^{\prime}\left(s+t\right)=\Phi^{\prime}\left(s\right)\Phi\left(t\right),
\]
and thus, for every $t\in\mathbb{R},$
\begin{equation}
\Phi^{\prime}\left(t\right)=\Phi^{\prime}\left(0\right)\Phi\left(t\right).\label{eq:deriv_Phi}
\end{equation}

\begin{itemize}
\item If $\Phi^{\prime}\left(0\right)=0,$ then, for every $t\in\mathbb{R},$
$\Phi^{\prime}\left(t\right)=0,$ hence $\Phi$ is constant. Since
$\left|\Phi\right|=1,$ $\Phi$ is nonzero, and by $\refpar{eq:phi_s_plus_t},$
we must have $\Phi\left(t\right)=1$ for every $t\in\mathbb{R}$ so
the choice $c=0$ works.
\item If $\Phi^{\prime}\left(0\right)\neq0,$ it follows by $\refpar{eq:deriv_Phi}$
that $\Phi\left(0\right)=1.$ \\
Set $c=-\text{i}\Phi^{\prime}\left(0\right).$ Then $\Phi^{\prime}\left(t\right)=\text{i}c\Phi\left(t\right),$
and therefore
\[
\dfrac{\text{d}}{\text{d}t}\left(\Phi\left(t\right)\text{e}^{-\text{i}ct}\right)=\Phi^{\prime}\left(t\right)\text{e}^{-\text{i}ct}-\text{i}c\Phi\left(t\right)\text{e}^{-\text{i}ct}=0.
\]
It follows that
\[
\Phi\left(t\right)\text{e}^{-\text{i}ct}=\Phi\left(0\right)=1,
\]
so, for every $t\in\mathbb{R},$ $\Phi\left(t\right)=\text{e}^{\text{i}ct}.$
\\
Finally, since $\left|\Phi\left(1\right)\right|=\left|\text{e}^{\text{i}c}\right|=1,$
$c$ must be a \textbf{real number}.
\end{itemize}
\textbf{5. Existence of $m$ and $a>0$ such that $\forall t\in\mathbb{R},$
$\widehat{\gamma}\left(t\right)=\text{e}^{\text{i}mt-a\frac{t^{2}}{2}}.$
Gaussianity of $X$ and $Y$}

By the relation $\refpar{eq:g_square_u_plus_v}$ and the definition
of $g,$ we may apply the result of the previous question to the function
$g^{2}.$ Hence, there exists a real number $m$ such that, for every
real number $t,$ 
\[
g^{2}\left(t\right)=\text{e}^{\text{i}2mt}.
\]
Since $g\left(0\right)=1,$ the continuity of $g$ implies that $g\left(t\right)=\text{e}^{\text{i}mt}.$
It then follows from $\refpar{eq:module_tr_f_gamma}$ that, for every
$t\in\mathbb{R},$\boxeq{
\[
\widehat{\gamma}\left(t\right)=\text{e}^{\text{i}mt-a\frac{t^{2}}{2}},
\]
}where $a>0.$ In other words, $\gamma$ is the Gaussian probability
$\mathscr{N}_{\mathbb{R}}\left(m,a\right).$ Therefore, the random
variable $X+Y$ is Gaussian. One shows similarly that $X-Y$ is Gaussian.
Since these two random variables are independent, the random variable
$\left(X+Y,X-Y\right)$ is Gaussian. Finally, the random variables
$X$ and $Y$ are Gaussian, as linear transformations of the Gaussian
random variable $\left(X+Y,X-Y\right).$

\textbf{6. Generalization to random variables with values in $\mathbb{R}^{d}$}

To extend the result to $\mathbb{R}^{d},$ it suffices to apply, for
every $x$ and $y$ of $\mathbb{R}^{d},$ the previous result to the
real-valued random variables $\left\langle X,x\right\rangle $ and
$\left\langle Y,y\right\rangle .$

\end{solution}

\footnotetext{We use the classical argument on the integers, then
on the rationnals, and then by continuity, we extend to the real numbers.}

\begin{solution}{}{solexercise14.5}

\textbf{1. Law of $X.$ Expression of $M_{n}$ and $\Sigma_{n}.$
Independence of $M_{n}$ and $\Sigma_{n}.$ Laws of $M_{n}$ and $n\Sigma_{n}$
when $m=0$ and $\sigma=1$}

The random variables $X_{i}$ are \textbf{independent} and have Gaussian
law $\mathscr{N}_{\mathbb{R}}\left(m,\sigma^{2}\right).$ The random
variable $X,$ taking values in $\mathbb{R}^{n},$ is therefore also
Gaussian, with expectation $\left(m,m,\cdots,m\right)$ and covariance
matrix $\sigma^{2}\text{I}$ where $\text{I}$ is the identity matrix
of $\mathbb{R}^{n}.$ It follows from the definition of $C$ that\boxeq{
\[
M_{n}=\dfrac{1}{\sqrt{n}}\left(CX\right)_{1},
\]
}and thus that
\[
n\Sigma_{n}=\left\Vert X\right\Vert ^{2}-\left|\left(CX\right)_{1}\right|^{2}.
\]
Since the matrix $C$ is orthogonal, it preserves the norm, which
yields\boxeq{
\[
n\Sigma_{n}=\sum^{n}_{i=2}\left[\left(CX\right)_{i}\right]^{2}.
\]
}

The random variable $CX,$ being a linear transformation of the Gaussian
random variable $X,$ is itself Gaussian, with expectation $\left(Cm,Cm,\cdots,Cm\right)$
and the covariance matrix $\sigma^{2}C\text{I}C^{\ast}=\sigma^{2}\text{I}$---since
the matrix $C$ is orthogonal. It follows that the components of $CX$
are independent, and therefore the random variables $M_{n}$ and $\Sigma_{n}$
are also independent.

In the case where $m=0$ and $\sigma=1,$ the random variable $M_{n}$
follows the law $\mathscr{N}_{\mathbb{R}}\left(0,\dfrac{1}{n}\right).$
Moreover, since $n\Sigma_{n}$ is a sum of $n-1$ squares of independent
reduced centered Gaussian random variables, the random variable $n\Sigma_{n}$
follows the chi-squared law with $n-1$ degrees of freedom.

\textbf{2. }The random variables $M_{n}$ and $\Sigma_{n}$ being
independent, the variables $S_{n}$ and $V_{n}$ are also independent.
We first assume that the random variables $X_{i}$ are centered.

\textbf{(a) Expectation $\mathbb{E}\left(V_{n}\right)$ as a function
of $\sigma.$}

We have
\[
V_{n}=\sum^{n}_{i=1}X^{2}_{i}-\dfrac{S^{2}_{n}}{n},\,\,\,\,\text{and thus,}\,\,\,\,\mathbb{E}\left(V_{n}\right)=\sum^{n}_{i=1}\mathbb{E}\left(X^{2}_{i}\right)-\dfrac{\mathbb{E}\left(S^{2}_{n}\right)}{n}.
\]
But since 
\[
\mathbb{E}\left(S_{n}\right)=\sum^{n}_{i=1}\mathbb{E}\left(X_{i}\right)=0,
\]
 we have, by independence of the $X_{i},$
\[
\mathbb{E}\left(S^{2}_{n}\right)=\sigma^{2}_{S_{n}}=\sum^{n}_{i=1}\sigma^{2}_{X_{i}}=n\sigma^{2}.
\]

Finally, since the $X_{i}$ are centered, $\sigma^{2}_{X_{i}}=\mathbb{E}\left(X^{2}_{i}\right),$
which yields\boxeq{
\[
\mathbb{E}\left(V_{n}\right)=\left(n-1\right)\sigma^{2}.
\]
}

\textbf{(b) $\left(u,v\right)\mapsto\mathbb{E}\left(\text{e}^{\text{i}\left(uS_{n}+vV_{n}\right)}\right)$
is differentiable. Justification of $\forall u\in\mathbb{R},\,\,\,\,\mathbb{E}\left(V_{n}\text{e}^{\text{i}\left(uS_{n}\right)}\right)=\left[\varphi\left(u\right)\right]^{n}\mathbb{E}\left(V_{n}\right)$}

The random variables $S_{n}$ and $V_{n}$ are independent. Thus,
for every $\left(u,v\right)\in\mathbb{R}^{2},$
\[
\mathbb{E}\left(\text{e}^{\text{i}\left(uS_{n}+vV_{n}\right)}\right)=\varphi_{S_{n}}\left(u\right)\varphi_{V_{n}}\left(v\right),
\]
where $\varphi_{S_{n}}$ and $\varphi_{V_{n}}$ are the characteristic
functions of $S_{n}$ and $V_{n}.$ Since these random variables have
a moment of order one, their characteristic functions are differentiable.
Consequentely, the function $\left(u,v\right)\mapsto\mathbb{E}\left(\text{e}^{\text{i}\left(uM_{n}+v\Sigma_{n}\right)}\right)$
is differentiable.

Moreover, the random variables $S_{n}$ and $V_{n}$ being independent,
the random variables $V_{n}$ and $\text{e}^{\text{i}\left(uS_{n}\right)}$
are also independent. Since they are integrable, for every real number
$u,$
\[
\mathbb{E}\left(V_{n}\text{e}^{\text{i}\left(uS_{n}\right)}\right)=\mathbb{E}\left(V_{n}\right)\mathbb{E}\left(\text{e}^{\text{i}\left(uS_{n}\right)}\right)=\varphi_{S_{n}}\left(u\right)\mathbb{E}\left(V_{n}\right).
\]
Finally, since the $X_{i}$ are independent, it then follows that\boxeq{
\begin{equation}
\forall u\in\mathbb{R},\,\,\,\,\mathbb{E}\left(V_{n}\text{e}^{i\left(uS_{n}\right)}\right)=\left[\varphi\left(u\right)\right]^{n}\mathbb{E}\left(V_{n}\right).\label{eq:expectation_V_neiuS_n}
\end{equation}
}

\begin{remark}{}{}We could also have already exploited the differentiability
established above and observed that
\[
\left[\dfrac{\partial}{\partial v}\mathbb{E}\left(\text{e}^{i\left(uS_{n}+vV_{n}\right)}\right)\right]_{v=0}=\varphi_{S_{n}}\left(u\right)\varphi^{\prime}_{V_{n}}\left(0\right).
\]
Moreover, since $\left|\text{e}^{\text{i}\left(uS_{n}+vV_{n}\right)}\right|\leqslant1,$
\[
\left|\dfrac{\partial}{\partial v}\text{e}^{\text{i}\left(uS_{n}+vV_{n}\right)}\right|=\left|V_{n}\text{e}^{\text{i}\left(uS_{n}+vV_{n}\right)}\right|\leqslant\left|V_{n}\right|.
\]
Because $V_{n}$ is $P-$integrable, the theorem on differentiation
under the integral sign for a parameter-dependent integral, gives
\[
\dfrac{\partial}{\partial v}\mathbb{E}\left(\text{e}^{i\left(uS_{n}+vV_{n}\right)}\right)=\text{i}\mathbb{E}\left(V_{n}\text{e}^{\text{i}\left(uS_{n}+vV_{n}\right)}\right).
\]

Since
\[
\varphi^{\prime}_{V_{n}}\left(0\right)=\text{i}\mathbb{E}\left(V_{n}\right),
\]
and using the independence of the $X_{i},$ we recover the relation
$\refpar{eq:expectation_V_neiuS_n}$.

\end{remark}

\textbf{(c) Computation of $\mathbb{E}\left(V_{n}\text{e}^{\text{i}\left(uS_{n}\right)}\right)$}

Since $V_{n}=\sum^{n}_{k=1}X^{2}_{k}-\dfrac{S^{2}_{n}}{n},$ we have
\begin{equation}
\mathbb{E}\left(V_{n}\text{e}^{\text{i}\left(uS_{n}\right)}\right)=\sum^{n}_{k=1}\mathbb{E}\left(X^{2}_{k}\text{e}^{\text{i}\left(uS_{n}\right)}\right)-\dfrac{1}{n}\mathbb{E}\left(S^{2}_{n}\text{e}^{\text{i}\left(uS_{n}\right)}\right).\label{eq:expect_v_n_eiusn}
\end{equation}
By independence of the $X_{k},$ for every fixed $k\leqslant n,$
\[
\mathbb{E}\left(X^{2}_{k}\text{e}^{\text{i}uS_{n}}\right)=\left[\prod_{\substack{l\neq k\\
1\leqslant l\leqslant n
}
}\mathbb{E}\left(\text{e}^{\text{i}uX_{l}}\right)\right]\mathbb{E}\left(X^{2}_{k}\text{e}^{\text{i}uX_{k}}\right).
\]
Hence, since the random variables $X_{l}$ are independent,
\begin{equation}
\mathbb{E}\left(X^{2}_{k}\text{e}^{\text{i}uS_{n}}\right)=\left[\varphi\left(u\right)\right]^{n-1}\mathbb{E}\left(X^{2}_{k}\text{e}^{\text{i}uX_{k}}\right).\label{eq:exp_x_ksqu_e_iuS_n}
\end{equation}

Since the random variables $X_{k}$ and $S_{n}$ admit a second-order
moment, their characteristic function are twice differentiable, and
\begin{equation}
\mathbb{E}\left(X^{2}_{k}\text{e}^{\text{i}uX_{k}}\right)=-\varphi^{\prime\prime}_{X_{k}}\left(u\right)\,\,\,\,\text{and}\,\,\,\,\mathbb{E}\left(S^{2}_{n}\text{e}^{\text{i}uS_{n}}\right)=-\varphi^{\prime\prime}_{S_{n}}\left(u\right).\label{eq:moment_order_two_x_k_and_s_n}
\end{equation}

Because the random variables $X_{k}$ are independent and follow the
same law, 
\[
\varphi_{S_{n}}\left(u\right)=\left[\varphi\left(u\right)\right]^{n}.
\]
It then follows by $\refpar{eq:expect_v_n_eiusn},$ $\refpar{eq:exp_x_ksqu_e_iuS_n}$
and $\refpar{eq:moment_order_two_x_k_and_s_n}$ that 
\begin{align*}
\mathbb{E}\left(V_{n}\text{e}^{\text{i}uS_{n}}\right) & =-n\left(\varphi\left(u\right)\right)^{n-1}\varphi^{\prime\prime}\left(u\right)\\
 & +\dfrac{1}{n}\left[n\left(n-1\right)\left[\varphi\left(u\right)\right]^{n-2}\left[\varphi^{\prime}\left(u\right)\right]^{2}+n\text{\ensuremath{\left[\varphi\left(u\right)\right]}}^{n-1}\varphi^{\prime\prime}\left(u\right)\right].
\end{align*}
Hence,
\begin{equation}
\mathbb{E}\left(V_{n}\text{e}^{\text{i}uS_{n}}\right)=-\left(n-1\right)\left(\varphi\left(u\right)\right)^{n-1}\varphi^{\prime\prime}\left(u\right)+\left(n-1\right)\left[\varphi\left(u\right)\right]^{n-2}\left[\varphi^{\prime}\left(u\right)\right]^{2}.\label{eq:exp_v_n_exp_iuS_n_2}
\end{equation}

\textbf{(d) $\varphi$ is solution of $\dfrac{\varphi^{\prime\prime}}{\varphi}-\left(\dfrac{\varphi^{\prime}}{\varphi}\right)^{2}=-\sigma^{2}.$
Proof that $\mu$ is the Gaussian law $\mathscr{N}_{\mathbb{R}}\left(0,\sigma^{2}\right)$}

First note that, since $\varphi$ is continuous and $\varphi\left(0\right)=1,$
the set $\varphi^{-1}\left(\left\{ 0\right\} ^{c}\right)$ is an open
neighborhood of $0.$ 

Using the equalities $\refpar{eq:expectation_V_neiuS_n}$ and $\refpar{eq:exp_v_n_exp_iuS_n_2},$
together with the value of $\mathbb{E}\left(V_{n}\right),$ we see
that on the open set $\varphi^{-1}\left(\left\{ 0\right\} ^{c}\right),$
the function $\varphi$ satisfies the differential equation
\[
\dfrac{\psi^{\prime\prime}\left(u\right)}{\psi\left(u\right)}-\left(\dfrac{\psi^{\prime}\left(u\right)}{\psi\left(u\right)}\right)^{2}=-\sigma^{2}.
\]
Moreover, since $\varphi\left(0\right)=1$ and $\varphi^{\prime}\left(0\right)=0$---since
the $X_{i}$ are centered---there exists $a_{1}>0$ such that, for
every $u\in\left[-a_{1},a_{1}\right],$
\begin{equation}
\varphi\left(u\right)=\text{e}^{-\frac{\sigma^{2}u^{2}}{2}}.\label{eq:varphi_exp}
\end{equation}
In particular, $\varphi\left(a_{1}\right)\neq0.$ By continuity, $\varphi$
does not vanish on some interval $\left[-a_{1},a_{2}\right],$ with
$a_{2}>a_{1},$ and the equality $\refpar{eq:varphi_exp}$ remains
valid there. 

By induction, we prove that similarly there exists a strictly increasing
sequence of real numbers $a_{n}>0$ such that the equality $\refpar{eq:varphi_exp}$
remains valid on the interval $\left[-a_{1},a_{n}\right].$ If this
sequence were bounded, it would converge to some real number $a>0.$
Then, for every $n\in\mathbb{N},$ 
\[
\varphi\left(a_{n}\right)=\text{e}^{-\frac{\sigma^{2}a^{2}_{n}}{2}}
\]
and, by continuity of $\varphi$ and of the exponential,
\[
\varphi\left(a\right)=\text{e}^{-\frac{\sigma^{2}a^{2}}{2}}>0
\]
which yields a contradiction. 

Therefore, the equality $\refpar{eq:varphi_exp}$ holds on $\left[-a_{1},+\infty\right[,$
and hence, on all $\mathbb{R},$ since, for every real number $u,$
\[
\varphi\left(-u\right)=\overline{\varphi\left(u\right)}.
\]
Thus, \textbf{$\mu$ follows the Gaussian law $\mathscr{N}_{\mathbb{R}}\left(0,\sigma^{2}\right).$}

\textbf{(e) Case where the random variables $X_{i}$ are no longer
centered}

If the random variables $X_{i}$ have expectation $m,$ consider the
centered random variables $\mathring{X_{i}}=X_{i}-m,$ which are still
independent. A straightforward computation gives
\[
\mathring{S_{n}}=\sum^{n}_{i=1}\mathring{X}_{i}=S_{n}-nm\,\,\,\,\text{and}\,\,\,\,\mathring{V}_{n}=\sum^{n}_{i=1}\mathring{X}^{2}_{i}-\dfrac{\mathring{S}^{2}_{n}}{n}=V_{n}.
\]
If the random variables $S_{n}$ and $V_{n}$ are independent, then
the random variables $\mathring{S}_{n}$ and $\mathring{V}_{n}$ are
also independent. By the previous question, the random variables $\mathring{X}_{i}$
follow the law $\mathscr{N}_{\mathbb{R}}\left(0,\sigma^{2}\right).$
It follows that the measure \textbf{$\mu$ }follows the\textbf{ Gaussian
law $\mathscr{N}_{\mathbb{R}}\left(m,\sigma^{2}\right).$}

\end{solution}

\begin{solution}{}{solexercise14.6}

Since the random variables $S$ and $V$ are independent and Gaussian,
the random variable $\left(S,V\right)$ is Gaussian, and thus any
linear transformation, in particular $\left(aR+bS,R\right)$ is also
Gaussian. Hence, for $aR+bS$ and $R$ to be independent, it is necessary
and sufficient that 
\[
\text{cov}\left(aR+bS,R\right)=0.
\]
By bilinearity of the covariance and using the independence of the
random variables $S$ and $V,$ we obtain
\[
\text{cov}\left(aR+bS,R\right)=a\sigma^{2}_{R}+b\text{cov}\left(S,R\right)=a\left(t^{2}\sigma^{2}+t\right)+bt\sigma^{2}.
\]

We choose $a$ and $b$ nonzero such that
\[
a=-\dfrac{b\sigma^{2}}{t\sigma^{2}+1},
\]
so that the random variables $aR+bS$ and $R$ are independent. 

We then have
\[
\mathbb{E}^{\sigma\left(R\right)}\left(aR+bS\right)=\mathbb{E}\left(aR+bS\right)=a\mathbb{E}\left(R\right)+b\mathbb{E}\left(S\right),
\]
and also
\[
\mathbb{E}^{\sigma\left(R\right)}\left(aR+bS\right)=aR+b\mathbb{E}^{\sigma\left(R\right)}\left(S\right).
\]
Therefore,
\[
b\mathbb{E}^{\sigma\left(R\right)}\left(S\right)=-a\left(R-\mathbb{E}\left(R\right)\right)+b\mathbb{E}\left(S\right).
\]

Taking into account the choice of $a$ and $b,$ we obtain
\[
\mathbb{E}^{\sigma\left(R\right)}\left(S\right)=\dfrac{\sigma^{2}}{t\sigma^{2}+1}\left(R-tm\right)+m,
\]
and hence\boxeq{
\[
\mathbb{E}^{\sigma\left(R\right)}\left(S\right)=\dfrac{m+\sigma^{2}R}{t\sigma^{2}+1}.
\]
}

\end{solution}

\begin{solution}{}{solexercise14.7}

Since the random variables $S$ and $W_{t_{1}},W_{t_{2}},\cdots,W_{t_{n}}$
are independent and Gaussian, the random variable $\left(S,W\right)$
is Gaussian and thus so is its linear transform $\left(\left\langle u,R\right\rangle +bS,R\right).$
Hence, for $\left\langle u,R\right\rangle +bS$ and $R$ to be independent,
it is necessary and sufficient that the covariance matrix of $\left(\left\langle u,R\right\rangle +bS,R\right)$
is zero, which can also be written
\[
\forall j\in\left\llbracket 1,n\right\rrbracket ,\,\,\,\,\text{cov}\left(\left\langle u,R\right\rangle +bS,R_{j}\right)=0,
\]
that is,
\begin{equation}
\forall j\in\left\llbracket 1,n\right\rrbracket ,\,\,\,\,\text{cov}\left(\left\langle u,R\right\rangle ,R_{j}\right)+b\text{cov}\left(S,R_{j}\right)=0.\label{eq:cond_covur+bsandrzero}
\end{equation}
By independence of $S$ and $W_{t_{j}},$
\[
\text{cov}\left(S,R_{j}\right)=\text{cov}\left(S,t_{j}S+W_{t_{j}}\right)=t_{j}\sigma^{2}_{S}.
\]
Moreover,
\[
\text{cov}\left(\left\langle u,R\right\rangle ,R_{j}\right)=\left\langle u,t\right\rangle \text{cov}\left(S,R_{j}\right)+\text{cov}\left(\left\langle u,W\right\rangle ,R_{j}\right).
\]
Hence,
\[
\text{cov}\left(\left\langle u,R\right\rangle ,R_{j}\right)=\left\langle u,t\right\rangle \left[t_{j}\sigma^{2}_{S}\right]+\left[t_{j}\text{cov}\left(\left\langle u,W\right\rangle ,S\right)+\text{cov}\left(\left\langle u,W\right\rangle ,W_{t_{j}}\right)\right].
\]
Since the random variables $\left\langle u,W\right\rangle $ and $S$
are independent,
\[
\text{cov}\left(\left\langle u,R\right\rangle ,R_{j}\right)=\left\langle u,t\right\rangle \left[t_{j}\sigma^{2}_{S}\right]+\text{cov}\left(\left\langle u,W\right\rangle ,W_{t_{j}}\right).
\]
Since the components of the random variable $W$ are independent,
\[
\text{cov}\left(\left\langle u,W\right\rangle ,W_{t_{j}}\right)=u_{j}\sigma^{2}_{W_{t_{j}}}=u_{j}t_{j}.
\]
Thus, the independence condition $\refpar{eq:cond_covur+bsandrzero}$
can be written, after dividing by $t_{j},$ as
\begin{equation}
\forall j\in\left\llbracket 1,n\right\rrbracket ,\,\,\,\,\left\langle u,t\right\rangle \sigma^{2}+u_{j}+b\sigma^{2}=0.\label{eq:cond_covur+bsandrzero-1}
\end{equation}
We take for $u$ the vector $\boldsymbol{1}$ whose components are
all equal to 1 and choose $b$ such that
\begin{equation}
\left\langle \boldsymbol{1},t\right\rangle \sigma^{2}+1+b\sigma^{2}=0.\label{eq:cond_cov_u_zero_when_u_1}
\end{equation}
With this choice, the random variables $\left\langle \boldsymbol{1},R\right\rangle +bS$
and $R$ are independent. 

We then have
\[
\mathbb{E}^{\sigma\left(R\right)}\left(\left\langle \boldsymbol{1},R\right\rangle +bS\right)=\mathbb{E}\left(\left\langle \boldsymbol{1},R\right\rangle +bS\right)=\mathbb{E}\left(\left\langle \boldsymbol{1},R\right\rangle \right)+b\mathbb{E}\left(S\right),
\]
and also
\[
\mathbb{E}^{\sigma\left(R\right)}\left(\left\langle \boldsymbol{1},R\right\rangle +bS\right)=\left\langle \boldsymbol{1},R\right\rangle +b\mathbb{E}^{\sigma\left(R\right)}\left(S\right),
\]
Comparing the right-hand sides of these equalities yields
\[
b\mathbb{E}^{\sigma\left(R\right)}\left(S\right)=-\left\langle \boldsymbol{1},R\right\rangle +\mathbb{E}\left(\left\langle \boldsymbol{1},R\right\rangle \right)+bm.
\]
Using the value of $b$ given by the equality $\refpar{eq:cond_cov_u_zero_when_u_1}$
and the fact that
\[
\mathbb{E}\left(\left\langle \boldsymbol{1},R\right\rangle \right)=\left\langle \boldsymbol{1},\mathbb{E}\left(R\right)\right\rangle =m\left\langle \boldsymbol{1},t\right\rangle .
\]

It follows that 
\[
\mathbb{E}^{\sigma\left(R\right)}\left(S\right)=\dfrac{\sigma^{2}}{\left\langle \boldsymbol{1},t\right\rangle \sigma^{2}+1}\left\langle \boldsymbol{1},R\right\rangle +m\cdot\dfrac{1}{\left\langle \boldsymbol{1},t\right\rangle \sigma^{2}+1},
\]
which can also be written as\boxeq{
\[
\mathbb{E}^{\sigma\left(R\right)}\left(S\right)=\dfrac{m+\sigma^{2}\sum^{n}_{j=1}R_{j}}{1+\sigma^{2}\sum^{n}_{j=1}t_{j}}.
\]
}

\end{solution}

\begin{solution}{}{solexercise14.8}

\textbf{1. $U$ and $V$ are independent. Identification of their
laws.}

Let $\left(e_{i}\right)_{i=1,\cdots,d}$ be an orthonormal basis of
$E$ whose first vector is $a.$ We then have
\[
U=\left\langle X,e_{1}\right\rangle ,\,\,\,\,\text{and}\,\,\,\,V=\sum^{d}_{i=2}\left\langle X,e_{i}\right\rangle ^{2}
\]
and the random variables $\left\langle X,e_{i}\right\rangle ,\,i\in\left\llbracket 1,d\right\rrbracket ,$
are independent. The independence of $U$ and $V$ follows. 

Since the law of $U$ is Gaussian,
\[
\mathbb{E}\left(U\right)=\left\langle E\left(X\right),e_{1}\right\rangle =0\,\,\,\,\text{and\,\,\,\,\ensuremath{\sigma}}^{2}_{\left\langle X,e_{1}\right\rangle }=\left\langle \Lambda_{X}a,a\right\rangle =\left\Vert a\right\Vert ^{2}=1.
\]
\textbf{Thus $U$ has law $\mathscr{N}_{\mathbb{R}}\left(0,1\right).$ }

Similarly, the random variables $\left\langle X,e_{i}\right\rangle ,\,i\in\left\llbracket 2,d\right\rrbracket ,$
are independent and have law $\mathscr{N}_{\mathbb{R}}\left(0,1\right).$
Therefore \textbf{$V$ follows the chi-squared law with $d-1$ degrees
of freedom.}

\textbf{2. Law of $\left\Vert Y\right\Vert ^{2}$ as convolution of
the chi-squared law at $d-1$ degrees of freedom and the law of the
square of a Gaussian real-valued random variable of law $\mathscr{N}_{\mathbb{R}}\left(\left\Vert m\right\Vert ,\text{I}\right).$}

We have
\[
\left\Vert Y\right\Vert ^{2}=\left\Vert Y-m\right\Vert ^{2}+2\left\langle Y-m,m\right\rangle +\left\Vert m\right\Vert ^{2},
\]
which can also be written as
\[
\left\Vert Y\right\Vert ^{2}=\left[\left\Vert Y-m\right\Vert ^{2}-\left\langle Y-m,\dfrac{m}{\left\Vert m\right\Vert }\right\rangle ^{2}\right]+\left[\left\langle Y-m,\dfrac{m}{\left\Vert m\right\Vert }\right\rangle +\left\Vert m\right\Vert \right]^{2}.
\]
The random variable $Y-m$ follows the law $\mathscr{N}_{E}\left(0,\text{I}\right).$
It follows from the previous question that the random variables $\left\Vert Y-m\right\Vert ^{2}-\left\langle Y-m,\dfrac{m}{\left\Vert m\right\Vert }\right\rangle ^{2}$
and $\left\langle Y-m,\dfrac{m}{\left\Vert m\right\Vert }\right\rangle +\left\Vert m\right\Vert $
are independent, with respective laws: the chi-squared law with $d-1$
degrees of freedom, and the Gaussian law $\mathscr{N}_{\mathbb{R}}\left(\left\Vert m\right\Vert ,1\right).$
This shows that\textbf{ the law of $\left\Vert Y\right\Vert ^{2}$
is the convolution of a chi-squared law with $d-1$ degrees of freedom
and the law of the square of a Gaussian real-valued random variable
with law $\mathscr{N}_{\mathbb{R}}\left(\left\Vert m\right\Vert ,1\right).$}

\end{solution}

\begin{solution}{}{solexercise14.9}

\textbf{1. $\left\Vert Y\right\Vert ^{2}$ and $\dfrac{Y}{\left\Vert Y\right\Vert }$
are independent. Specification of their laws in the case $E=\mathbb{R}$}

We first consider the case $E=\mathbb{R}$ and let $f$ and $g$ be
nonnegative measurable functions from $\mathbb{R}$ into itself. By
the transfer theorem,
\[
\mathbb{E}\left(f\left(\left\Vert Y\right\Vert ^{2}\right)g\left(\dfrac{Y}{\left\Vert Y\right\Vert }\right)\right)=\intop_{\mathbb{R}^{n}}f\left(\left\Vert y\right\Vert ^{2}\right)g\left(\dfrac{y}{\left\Vert y\right\Vert }\right)\dfrac{1}{\left(2\pi\right)^{\frac{n}{2}}}\text{e}^{-\frac{\left\Vert y\right\Vert ^{2}}{2}}\text{d}y.
\]
We make the change of variables to spherical coordinates defined by
\[
\left\{ \begin{array}{l}
y_{1}=\rho\cos\varphi_{1}\\
y_{2}=\rho\sin\varphi_{1}\cos\varphi_{2}\\
\vdots\\
y_{n-1}=\rho\sin\varphi_{1}\cdots\sin\varphi_{n-2}\cos\varphi_{n-1}\\
y_{n}=\rho\sin\varphi_{1}\cdots\sin\varphi_{n-2}\sin\varphi_{n-1},
\end{array}\right.
\]
which defines a diffeomorphism from $\mathbb{R}^{n}\backslash\left(\bigcup^{n}_{i=1}D_{i}\right)$
to $\left]0,+\infty\right[\times\left]0,\pi\right[^{n-2}\times\left]0,2\pi\right[,$
where $D_{i}$ is the line generated by the $i-$th vector of the
canonical basis of $\mathbb{R}^{n}.$ 

The Jacobian of the transformation is
\[
J_{\left(\rho,\varphi_{1},\varphi_{2},\cdots,\varphi_{n-1}\right)}=\rho^{n-1}\prod^{n-2}_{i=1}\left(\sin\varphi_{i}\right)^{n-i-1}.
\]
The change of variables and an application of the Fubini theorem lead
to
\begin{equation}
\mathbb{E}\left(f\left(\left\Vert Y\right\Vert ^{2}\right)g\left(\dfrac{Y}{\left\Vert Y\right\Vert }\right)\right)=I_{1}\left(f\right)I_{2}\left(g\right)\label{eq:exp_f_norm_Y_sqg_Yover_normY}
\end{equation}
 where we set
\[
I_{1}\left(f\right)=\intop_{\left]0,+\infty\right[}\dfrac{1}{\left(2\pi\right)^{\frac{n}{2}}}\rho^{n-1}f\left(\rho^{2}\right)\text{e}^{-\frac{\rho^{2}}{2}}\text{d}\rho,
\]
\begin{equation}
I_{2}\left(g\right)=\intop_{\left]0,\pi\right[^{n-2}\times\left]0,2\pi\right[}\prod^{n-2}_{i=1}\left(\sin\varphi_{i}\right)^{n-i-1}g\left(\Phi\left(\varphi_{1},\varphi_{2},\cdots,\varphi_{n-1}\right)\right)\text{d}\left(\varphi_{1},\varphi_{2},\cdots,\varphi_{n-1}\right),\label{eq:I_2_exp}
\end{equation}
and
\[
\Phi\left(\varphi_{1},\varphi_{2},\cdots,\varphi_{n-1}\right)=\left(\cos\varphi_{1},\sin\varphi_{1}\cos\varphi_{2},\cdots,\sin\varphi_{1}\cdots\sin\varphi_{n-2}\sin\varphi_{n-1}\right).
\]
In particular, we obtain
\[
\mathbb{E}\left(f\left(\left\Vert Y\right\Vert ^{2}\right)\right)=I_{1}\left(f\right)I_{2}\left(\boldsymbol{1}\right)\,\,\,\,\text{and}\,\,\,\,\mathbb{E}\left(g\left(\dfrac{Y}{\left\Vert Y\right\Vert }\right)\right)=I_{1}\left(\boldsymbol{1}\right)I_{2}\left(g\right),
\]
and
\[
\mathbb{E}\left(\boldsymbol{1}\right)=I_{1}\left(\boldsymbol{1}\right)I_{2}\left(\boldsymbol{1}\right)=\boldsymbol{1}.
\]

It then follows from $\refpar{eq:exp_f_norm_Y_sqg_Yover_normY}$ that,
for every nonnegative measurable functions $f$ and $g,$
\begin{equation}
\mathbb{E}\left(f\left(\left\Vert Y\right\Vert ^{2}\right)g\left(\dfrac{Y}{\left\Vert Y\right\Vert }\right)\right)=\mathbb{E}\left(f\left(\left\Vert Y\right\Vert ^{2}\right)\right)\mathbb{E}\left(g\left(\dfrac{Y}{\left\Vert Y\right\Vert }\right)\right)\label{eq:exp_f_squ_norm_Y_g_Y_over_normY}
\end{equation}
which is necessary and sufficient to ensure\textbf{ the independence
of the random variables $\left\Vert Y\right\Vert ^{2}$ and $\dfrac{Y}{\left\Vert Y\right\Vert }.$} 

\textbf{If now $E$ is an arbitrary Euclidean space, }let $\left(e_{i}\right)_{i\in\left\llbracket 1,d\right\rrbracket }$
be an orthonormal basis of $E.$ The real-valued random variables
$Z_{i,j}=\left\langle Y_{i},e_{j}\right\rangle ,$ where $i\in\left\llbracket 1,n\right\rrbracket $
and $j\in\left\llbracket 1,d\right\rrbracket $ are independent with
Gaussian law $\mathscr{N}_{\mathbb{R}}\left(0,1\right),$ since they
are linear transforms of the Gaussian random variable $Y$ taking
values in $E^{n},$ and since, for two different couples $\left(i,j\right)$
and $\left(k,l\right),$ 
\[
\text{cov}\left(\left\langle Y_{i},e_{j}\right\rangle ,\left\langle Y_{k},e_{l}\right\rangle \right)=0.
\]
 Moreover,
\[
\left\Vert Y\right\Vert ^{2}=\sum^{n}_{i=1}\sum^{d}_{j=1}Z^{2}_{i,j}.
\]

Let $Z$ be the random variable taking values in $E^{nd}$ defined
by
\[
Z=\left(Z_{1,1},\cdots,Z_{1,d},Z_{2,1},\cdots,Z_{2,d},\cdots,Z_{n,1},\cdots,Z_{n,d}\right).
\]
Then $\left\Vert Z\right\Vert ^{2}=\left\Vert Y\right\Vert ^{2},$
and by the previously established property, the random variables $\left\Vert Z\right\Vert ^{2}$
and $\dfrac{Z}{\left\Vert Z\right\Vert }$ are independent. It follows
that the random variables $\left\Vert Y\right\Vert ^{2}$ and $\dfrac{1}{\left\Vert Z\right\Vert }\left(\sum^{d}_{j=1}Z_{1,j}e_{j},\cdots,\sum^{d}_{j=1}Z_{n,j}e_{j}\right)$
are also independent. Since
\[
\dfrac{1}{\left\Vert Z\right\Vert }\left(\sum^{d}_{j=1}Z_{1,j}e_{j},\cdots,\sum^{d}_{j=1}Z_{n,j}e_{j}\right)=\dfrac{1}{\left\Vert Y\right\Vert }\left(Y_{1},\cdots,Y_{n}\right)=\dfrac{Y}{\left\Vert Y\right\Vert },
\]
\textbf{the random variables $\left\Vert Y\right\Vert ^{2}$ and $\dfrac{Y}{\left\Vert Y\right\Vert }$
are independent}.

If $E=\mathbb{R},$ \textbf{the law of $\left\Vert Y\right\Vert ^{2}$
is the chi-squared law with $n$ degrees of freedom}, since it is
the sum of $n$ independent squares of random variables with Gaussian
laws $\mathscr{N}_{\mathbb{R}}\left(0,1\right).$

We have shown that for every nonnegative measurable function $g,$
\[
\mathbb{E}\left(g\left(\dfrac{Y}{\left\Vert Y\right\Vert }\right)\right)=I_{1}\left(\boldsymbol{1}\right)I_{2}\left(g\right).
\]
Let $S_{n}$ be the sphere of $\mathbb{R}^{n}$ of center 0 and of
radius 1. Let $\mu$ be the image measure of the measure 
\[
I_{1}\left(\boldsymbol{1}\right)\left[\prod^{n-2}_{i=1}\left(\sin\varphi_{i}\right)^{n-i-1}\right]\text{d}\left(\varphi_{1},\varphi_{2},\cdots,\varphi_{n-1}\right)
\]
on $\left]0,\pi\right[^{n-2}\times\left]0,2\pi\right[$ by the function
$\Phi.$ The measure $\mu$ can be called the \textbf{uniform probability
on $S_{n}$\index{uniform probability} }and
\[
\mathbb{E}\left(g\left(\dfrac{Y}{\left\Vert Y\right\Vert }\right)\right)=\intop_{S_{n}}g\left(x\right)\text{d}\mu\left(x\right).
\]
Hence, \textbf{the law of $\dfrac{Y}{\left\Vert Y\right\Vert }$ is
the uniform law}---in the previous sense---\textbf{on $S_{n}.$} 

\textbf{2. }Since the random variables $X_{1},X_{2},\cdots,X_{n}$
are independent with Gaussian law $\mathscr{N}_{\mathbb{R}}\left(m,\sigma^{2}\right),$
the random variable $X$ is Gaussian with law $\mathscr{N}_{\mathbb{R}^{n}}\left(me,\sigma^{2}\boldsymbol{1}_{\mathbb{R}^{n}}\right).$
The random variable $\left(M,X^{\prime}\right),$ being a linear transform
of $X,$ is therefore Gaussian on $\mathbb{R}^{n+1}.$

\textbf{(a) $M$ and $X^{\prime}$ are independent}

Hence, for $M$ and $X^{\prime}$ to be independent, it is necessary
and sufficient that their cross-covariance operator $\Lambda_{M,X^{\prime}}$
is zero. Note that $X^{\prime}$ is centered since
\[
X^{\prime}=X-Me\,\,\,\,\text{and}\,\,\,\,M=\dfrac{1}{n}\left\langle X,e\right\rangle .
\]
Indeed,
\[
\mathbb{E}\left(X^{\prime}\right)=\mathbb{E}\left(X\right)-\mathbb{E}\left(\dfrac{1}{n}\left\langle X,e\right\rangle e\right)=me-\dfrac{1}{n}\left\langle me,e\right\rangle e=0.
\]

For every $u\in\mathbb{R}^{n},$
\[
\Lambda_{M,X^{\prime}}u=\text{cov}\left(M,\left\langle X^{\prime},u\right\rangle \right)=\text{cov}\left(M,\left\langle X,u\right\rangle -M\left\langle e,u\right\rangle \right).
\]
Hence
\[
\Lambda_{M,X^{\prime}}u=\text{cov}\left(\dfrac{1}{n}\left\langle X,e\right\rangle ,\left\langle X,u\right\rangle \right)-\left\langle e,u\right\rangle \sigma^{2}_{M}=\dfrac{1}{n}\left\langle \Lambda_{X}e,u\right\rangle -\left\langle e,u\right\rangle \sigma^{2}_{M}.
\]
Since $\Lambda_{X}=\sigma^{2}\boldsymbol{1}_{\mathbb{R}^{n}},$ 
\begin{equation}
\sigma^{2}_{M}=\dfrac{1}{n^{2}}\sigma^{2}_{\left\langle X,e\right\rangle }=\dfrac{1}{n^{2}}\left\langle \Lambda_{X}e,e\right\rangle =\dfrac{\sigma^{2}\left\Vert e\right\Vert ^{2}}{n^{2}}=\dfrac{\sigma^{2}}{n},\label{eq:square_sigma_M}
\end{equation}
and therefore
\[
\Lambda_{M,X^{\prime}}u=\dfrac{\sigma^{2}}{n}\left\langle e,u\right\rangle -\dfrac{\sigma^{2}}{n}\left\langle e,u\right\rangle =0.
\]
Hence, $\Lambda_{M,X^{\prime}}=0$ and \textbf{the random variables
$M$ and $X^{\prime}$ are independent.}

\textbf{(b) Computation of the auto-covariance operator of $X^{\prime}$}

Since $X=X^{\prime}+Me$ and $M$ and $X^{\prime}$ are independent,
we have
\[
\Lambda_{X}=\Lambda_{X^{\prime}}+\Lambda_{Me}.
\]
For every $u$ and $v$ of $\mathbb{R}^{n},$
\[
\left\langle \Lambda_{Me}u,v\right\rangle =\text{cov}\left(\left\langle Me,u\right\rangle ,\left\langle Me,v\right\rangle \right)=\left\langle e,u\right\rangle \left\langle e,v\right\rangle \sigma^{2}_{M}.
\]
Using $\refpar{eq:square_sigma_M},$ we obtain\boxeq{
\[
\Lambda_{X^{\prime}}=\sigma^{2}\left(\boldsymbol{1}_{\mathbb{R}^{n}}-\dfrac{1}{n}ee^{\ast}\right),
\]
 }where $ee^{\ast}$ is the endomorphism defined, for every $u$
and $v$ of $\mathbb{R}^{n},$ by
\[
\left\langle ee^{\ast}u,v\right\rangle =\left\langle e,u\right\rangle \left\langle e,v\right\rangle ,
\]
and its matrix representation in the canonical basis is the Kronecker
production of $e$ with itself.

\textbf{(c) Existence of an isometry $B$ from $\mathbb{R}^{n-1}$
to $H$ and a random variable $U$ on $\mathbb{R}^{n-1}$ with law
$\mathscr{N}_{\mathbb{R}^{n-1}}\left(0,\sigma^{2}\boldsymbol{1}_{\mathbb{R}^{n-1}}\right)$
such that $X^{\prime}=BU$ $P-$almost surely}

We have $\text{Ker}\Lambda_{X^{\prime}}=\mathbb{R}e$ and since $\Lambda_{X^{\prime}}$
is self-adjoint, 
\[
\text{Im}\left(\Lambda_{X^{\prime}}\right)=\left(\text{Ker}\Lambda_{X^{\prime}}\right)^{\perp}=H.
\]
Hence, $X^{\prime}$ takes its values $P-$almost surely in $H$---which
has dimension $n-1.$ Let $B$ be an isometry from $\mathbb{R}^{n-1}$
to $H.$ Denote by $i_{H}$ the canonical injection from $H$ into
$\mathbb{R}^{n}.$ Let $U$ be the random variable taking values in
$\mathbb{R}^{n-1},$ defined by
\[
U=B^{\ast}i^{\ast}_{H}X^{\prime}=\left(i_{H}B\right)^{\ast}X^{\prime}.
\]
Then $U$ follows a Gaussian law 
\[
\mathscr{N}_{\mathbb{R}^{n-1}}\left(\mathbb{E}\left(\left(i_{H}B\right)^{\ast}X^{\prime}\right),\left(i_{H}B\right)^{\ast}\Lambda_{X^{\prime}}\left(i_{H}B\right)\right)=\mathscr{N}_{\mathbb{R}^{n-1}}\left(0,\sigma^{2}\boldsymbol{1}_{\mathbb{R}^{n-1}}\right),
\]
since $X^{\prime}$ is centered and
\[
B^{\ast}i^{\ast}_{H}\Lambda_{X^{\prime}}i_{H}B=\sigma^{2}\boldsymbol{1}_{\mathbb{R}^{n-1}}.
\]
Finally, since $BB^{\ast}=\boldsymbol{1}_{H}$ and $X^{\prime}$ takes
values $P-$almost surely in $H,$ we have $P-$almost surely 
\[
X^{\prime}=BU.
\]

\textbf{(d) $M,V$ and $Z$ are independent}

It follows from the first question that the random variables $\left\Vert U\right\Vert ^{2}$
and $\dfrac{U}{\left\Vert U\right\Vert }$ are independent. Therefore,
the \textbf{random variables $V=\left\Vert X^{\prime}\right\Vert ^{2}$
and $Z=\dfrac{X^{\prime}}{\left\Vert X^{\prime}\right\Vert }$ are
also independent.} Since the random variables $V$ and $Z$ are $\sigma\left(X^{\prime}\right)-$measurable
and that $M$ and $X^{\prime}$ are independent, it follows that the
\textbf{random variables $M,\,V$ and $Z$ are independent}.

\textbf{(e) Laws of $M$ and $\dfrac{1}{\sigma^{2}}V.$}

Since
\[
\mathbb{E}\left(M\right)=\dfrac{1}{n}\left\langle \mathbb{E}\left(X\right),e\right\rangle =\dfrac{m}{n}\left\Vert e\right\Vert ^{2}=m,
\]
it follows from $\refpar{eq:square_sigma_M}$ that $M$ follows the
law $\mathscr{N}_{\mathbb{R}}\left(m,\dfrac{\sigma^{2}}{n}\right).$ 

Moreover,
\[
\dfrac{V}{\sigma^{2}}=\left\Vert \dfrac{X^{\prime}}{\sigma}\right\Vert ^{2}=\left\Vert \dfrac{U}{\sigma}\right\Vert ^{2}\,\,\,\,P-\text{almost surely.}
\]
Hence, $\dfrac{U}{\sigma}$ follows the law $\mathscr{N}_{\mathbb{R}^{n-1}}\left(0,\boldsymbol{1}_{\mathbb{R}^{n-1}}\right),$
and $\dfrac{V}{\sigma^{2}}$ follows the chi-squared law with $n-1$
degrees of freedom.

\end{solution}

\begin{remark}{}{}

The results of this last exercise lead to the \textbf{Student test}\index{Student test}.
This is a parametric test.

Consider a real-valued random variable $X$ with law $\mathscr{N}_{\mathbb{R}}\left(m,\sigma^{2}\right)$
where the parameters are unknown. We wish to test the hypothesis that
$m$ is lower or equal to a given value $m_{0}$ based on a sample
$\left(x_{1},x_{2},\cdots,x_{n}\right).$ 

Let $\left(X_{1},X_{2},\cdots,X_{n}\right)$ be an empirical sample
of $X,$ that is, $n$ independent random variables with the same
law as $X.$ Introduce the normalized centered random variables 
\[
\mathring{X_{i}}=\dfrac{X_{i}-m}{\sigma}
\]
which have law $\mathscr{N}_{\mathbb{R}}\left(0,1\right),$ and the
empirical moments associated to this sample, that is
\[
\begin{array}{ccc}
M_{n}=\dfrac{1}{n}\sum^{n}_{j=1}X_{j} & \,\,\,\,\,\,\,\, & \Sigma^{2}_{n}=\dfrac{1}{n-1}\sum^{n}_{j=1}\left(X_{j}-M_{n}\right)^{2}\\
\mathring{M}_{n}=\dfrac{1}{n}\sum^{n}_{j=1}\mathring{X}_{j} & \,\,\,\,\,\,\,\, & \left(\mathring{\Sigma}_{n}\right)^{2}=\dfrac{1}{n-1}\sum^{n}_{j=1}\left(\mathring{X}_{j}-\mathring{M}_{n}\right)^{2}.
\end{array}
\]

We have
\[
M_{n}=\sigma\mathring{M}_{n}+m\,\,\,\,\text{and}\,\,\,\,\Sigma^{2}_{n}=\sigma^{2}\left(\left(\mathring{\Sigma}_{n}\right)^{2}\right).
\]

Thus the random variables $\mathring{M}_{n}$ and $\left(\mathring{\Sigma}_{n}\right)^{2}$
are independent. As shown above, the same holds for the random variables
$M_{n}$ and $\Sigma^{2}_{n}.$ Moreover $\sqrt{n}\mathring{M_{n}}$
follows the law $\mathscr{N}_{\mathbb{R}}\left(0,1\right)$ and $\left(n-1\right)\left(\mathring{\Sigma}_{n}\right)^{2}$
follows the chi-squared law $\chi^{2}_{n-1}.$ 

Then the random variable 
\[
T_{n}=\dfrac{\sqrt{n}\mathring{M}_{n}}{\mathring{\Sigma}_{n}}=\sqrt{n-1}\dfrac{\sqrt{n}\mathring{M_{n}}}{\sqrt{n-1}\mathring{\Sigma}_{n}},
\]
follows the Student law with parameter $n-1$---see the exercise
of Chapter \ref{chap:PartIIChap10} on the Student laws. This law
is tabulated. Note that
\[
T_{n}=\sqrt{n}\dfrac{M_{n}-m}{\Sigma_{n}}.
\]
For a given threshold $\alpha,$ one determines from the table the
value $t_{n-1,\alpha}$ such that 
\[
P\left(T_{n}\leqslant t_{n-1,\alpha}\right)=1-\alpha.
\]
Now 
\[
T_{n}\leqslant t_{n-1,\alpha}\,\,\,\,\Leftrightarrow\,\,\,\,M_{n}\in\left]-\infty,t_{n-1,\alpha}\dfrac{\Sigma_{n}}{\sqrt{n}}+m\right].
\]
Under the hypothesis that the true---unknown---value of $m$ satisfies
$m\leqslant m_{0},$ we then have 
\[
M_{n}\in\left]-\infty,t_{n-1,\alpha}\dfrac{\Sigma_{n}}{\sqrt{n}}+m_{0}\right]
\]
with probability at least $1-\alpha.$ 

The Student test consists in accepting the hypothesis, with probability
of error at most $\alpha$, if the observed sample satisfies 
\[
\overline{x}\in\left]-\infty,t_{n-1,\alpha}\dfrac{s}{\sqrt{n}}+m_{0}\right],
\]
where
\[
\overline{x}=\dfrac{1}{n}\sum^{n}_{j=1}x_{j},\,\,\,\,s^{2}=\dfrac{1}{n-1}\sum^{n}_{j=1}\left(x_{j}-\overline{x}\right)^{2}.
\]
In this case, we say that a \textbf{confidence interval\index{confidence interval}}
at the threshold or level $\alpha$ has been determined.

\end{remark}

\chapter{Convergences of Measures and Convergence in Law}\label{chap:PartIIChapConvMeasuresandInLaw15}

\begin{objective}{}{}

Chapter \ref{chap:PartIIChapConvMeasuresandInLaw15} is devoted to
the study of convergence of measures and convergence in law.
\begin{itemize}
\item Section \ref{sec:Convergence-of-Bounded} studies the convergence
of bounded measures on $\mathbb{R}^{d}.$ It begins by introducing
the vague, weak and narrow topologies, as well as the notions of vague,
weak and narrow convergence for a sequence of bounded measures. These
topologies are then compared. A necessary and sufficient condition
for narrow convergence is established. The concept of $\mu-$continuity
for Borel sets is introduced before presenting practical criteria
for narrow convergence. The notion of a tight sequence is then defined,
followed by a sufficient condition ensuring tightness. The section
concludes with a characterization of narrow convergence of a sequence
of measures in terms of Fourier transforms, via the fundamental Lévy
theorem.
\item Section \ref{sec:Convergence-in-Law} begins by defining the convergence
in law toward a random variable and toward a probability measure.
The Lévy theorem is then formulated in terms of convergence in law.
A comparison between convergence in law and convergence in probability
is provided through a proposition, together with a partial converse
in the case of $P-$almost surely constant random variables. The Scheffé
lemma is then stated, giving a sufficient condition for convergence
in law when the random variables admit a density. A criterion for
convergence in law of a sequence of random variables is presented,
and the effect of the convergence in law on the cumulative distribution
function is analyzed. The section concludes with the study of rare
events, also known as the Poisson theorem.
\item Section \ref{sec:Central-Limit-Theorem} is devoted to the central
limit theorem, as well as to an important application: the Pearson
theorem, which forms the theoretical foundation of the chi-squared
test, is subsequently developed in detail.
\item Section \ref{sec:Estimation} briefly addresses the problem of estimating
the law of a random variable. It presents a method for constructing
an estimator via the method of maximum likelihood.
\end{itemize}
\end{objective}

\textbf{To simplify the exposition, we restrict ourselves to the study
of measures on $\mathbb{R}^{d}.$ All the results that follow extend
to the case where $E$ is a metric space that is locally compact and
countable at infinity}{\bfseries\footnote{We say that a locally compact space $E$ is \textbf{\index{countable at infinity}countable
at infinity} if there exists a sequence of compact sets $\left(K_{n}\right)_{n\in\mathbb{N}}$
such that $K_{n}\subset\mathring{K}_{n+1}$ for every $n\in\mathbb{N},$
and such that $\bigcup_{n\in\mathbb{N}}K_{n}=E.$ The open and closed
subsets of $\mathbb{R}^{d}$ are locally compact and countable at
infinity---for a bounded open subset of $\mathbb{R}^{d},$ we can
take $K_{n}$ to be the set of points at distance less than or equal
to $\dfrac{1}{n}$ from the boundary of the open set. If $E$ is compactified
by adjoining a point at infinity---Alexandrov compactification---this
is equivalent to saying that the point at infinity has a countable
basis of neighborhoods.}}\textbf{. In particular this applies when $E$ is a compact subset,
an open subset or a closed subset of $\mathbb{R}^{d},$ or when the
vector space is finite-dimensional.}

\textbf{The fundamental reference for questions concerning convergence
of measures is the book by Billingsley \cite{billingsley2013convergence}.
The necessary topological background can be found, for example in
the book by J. Dieudonné, ``Foundation of the modern analysis''
\cite{dieudonne2011foundations}}{\bfseries\footnote{Dieudonné J. (1965) ``Fondements de l'analyse moderne'', Cahiers
Scientifiques, Book XXVIII, Gauthier-Villars Editor. Tr.N. We give
here the reference to the English translation from 2011, the author
was refering to the French original version.}}\textbf{.}

\section{Convergence of Bounded Measures on $\mathbb{R}^{d}$}\label{sec:Convergence-of-Bounded}

Let $\mathscr{M}$ denote the set of bounded nonnegative measures
on $\mathbb{R}^{d},$ equipped with its Borel $\sigma-$algebra. For
$b>0,$ define $\mathscr{M}\left(b\right)$ the subset of measures
$\mu$ whose mass is less than or equal to $b;$ that is, such that
\[
\mu\left(\mathbb{R}^{d}\right)\leqslant b.
\]
 We denote by $\mathscr{M}^{1}$ the set of probability measures on
$\mathbb{R}^{d}.$

We now introduce three vector spaces of continuous real-valued functions
on $\mathbb{R}^{d}:$
\begin{itemize}
\item $\mathscr{C}_{\mathscr{K}}\left(\mathbb{R}^{d}\right),$ the space
of continuous functions with compact support.
\item $\mathscr{C}_{0}\left(\mathbb{R}^{d}\right),$ the space of continuous
functions that tend to 0 at infinity\footnote{In a locally compact space $E,$ we say that a real-valued function
$f$ is said to tend to 0 at infinity if, for every $\epsilon>0,$
there exists a compact set $K$ such that 
\[
\sup_{x\in K^{c}}\left|f\left(x\right)\right|\leqslant\epsilon.
\]
If, moreover, $E$ is countable at infinity, it is enough that for
every sequence $\left(x_{n}\right)_{n\in\mathbb{N}}$ tending to infinity,
the sequence $\left(f\left(x_{n}\right)\right)_{n\in\mathbb{N}}$
converges to 0.

By definition, a sequence $\left(x_{n}\right)_{n\in\mathbb{N}}$ tends
to infinity if, for every compact set $K,$ the sequence eventually
lies in $K^{c}$ from some index onward.}.
\item $\mathscr{C}_{b}\left(\mathbb{R}^{d}\right),$ the space of bounded
continuous functions.
\end{itemize}
These spaces satisfy the inclusions
\[
\mathscr{C}_{\mathscr{K}}\left(\mathbb{R}^{d}\right)\subset\mathscr{C}_{0}\left(\mathbb{R}^{d}\right)\subset\mathscr{C}_{b}\left(\mathbb{R}^{d}\right).
\]
Equipped with the norm, 
\[
\left\Vert f\right\Vert =\sup_{x\in\mathbb{R}^{d}}\left|f\left(x\right)\right|,
\]
 the space $\mathscr{C}_{b}\left(\mathbb{R}^{d}\right)$ is a Banach
space. 

Moreover:
\begin{itemize}
\item $\mathscr{C}_{0}\left(\mathbb{R}^{d}\right)$ is a closed subspace
of $\mathscr{C}_{b}\left(\mathbb{R}^{d}\right),$ 
\item $\mathscr{C}_{\mathscr{K}}\left(\mathbb{R}^{d}\right)$ is dense in
$\mathscr{C}_{0}\left(\mathbb{R}^{d}\right).$ 
\item $\mathscr{C}_{0}\left(\mathbb{R}^{d}\right)$ is \textbf{separable}\index{separable}\footnote{(a) A metric space is \textbf{separable\mindex{metric space!separable}}
if it has a dense countable subset.

(b) A subset $H$ of a normed vector space is said \textbf{total\mindex{normed vector space! total}}
if the vector subspace generated by $H$---that is the set of finite
linear combination of elements of $H$---is dense in $E.$

(c) A normed vector space that possesses a total countable subset
$H$ is separable---consider the finite linear combinations of elements
of $H$ with real-valued coefficients.

(d) There exists in $\mathscr{C}_{\mathscr{K}}\left(\mathbb{R}^{d}\right)$
and in $\mathscr{C}_{0}\left(\mathbb{R}^{d}\right)$ countable total
sets, which is not the case in $\mathscr{C}_{b}\left(\mathbb{R}^{d}\right).$}, 
\item Whereas $\mathscr{C}_{b}\left(\mathbb{R}^{d}\right)$ is not separable.
\end{itemize}
\begin{definition}{Vague, Weak and Narrow Topologies. Vaguely, Weakly and Narrowly Convergence}{}

On $\mathscr{M},$ we define the \textbf{vague topology}\index{vague topology}\mindex{topology!vague},
the \textbf{weak topology\index{weak topology}}\mindex{topology!weak}
and the \textbf{narrow topology\index{narrow topology}}\mindex{topology!narrow}
as the least fine\footnotemark topologies that make the mappings
\[
\mu\mapsto\intop f\text{d}\mu
\]
from $\mathscr{M}$ to $\mathbb{R}$ continuous:
\begin{itemize}
\item For every $f\in\mathscr{C}_{\mathscr{K}}\left(\mathbb{R}^{d}\right)$
in the case of the vague topology,
\item For every $f\in\mathscr{C}_{0}\left(\mathbb{R}^{d}\right)$ in the
case of the weak topology,
\item For every $f\in\mathscr{C}_{b}\left(\mathbb{R}^{d}\right)$ in the
case of the narrow topology.
\end{itemize}
In particular, a sequence $\left(\mu_{n}\right)_{n\in\mathbb{N}}$
of bounded measures converges to a measure $\mu,$
\begin{itemize}
\item \textbf{Vaguely\index{vague convergence}}\mindex{convergence!vague}
if 
\[
\lim_{n\to+\infty}\intop f\text{d}\mu_{n}=\intop f\text{d}\mu,\,\,\,\,\text{for every }f\in\mathscr{C}_{\mathscr{K}}\left(\mathbb{R}^{d}\right).
\]
\item \textbf{\index{weak convergence}}\mindex{convergence!weak}\textbf{Weakly}
if
\[
\lim_{n\to+\infty}\intop f\text{d}\mu_{n}=\intop f\text{d}\mu,\,\,\,\,\text{for every }f\in\mathscr{C}_{0}\left(\mathbb{R}^{d}\right).
\]
\item \textbf{\index{narrow convergence}}\mindex{convergence!narrow}\textbf{Narrowly}
if 
\[
\lim_{n\to+\infty}\intop f\text{d}\mu_{n}=\intop f\text{d}\mu,\,\,\,\,\text{for every }f\in\mathscr{C}_{b}\left(\mathbb{R}^{d}\right).
\]
\end{itemize}
\end{definition}

\footnotetext{(a) These are ``initial'' topologies.

(b) A topology $\mathscr{T}_{1}$ on a set $X$ is said to be \textbf{\index{less fine topology}less
fine} than a topology $\mathscr{T}_{2}$ on $X$ if every $\mathscr{T}_{1}-$open
set is also $\mathscr{T}_{2}-$open; in other words, if $\mathscr{T}_{1}$
contains fewer open sets than $\mathscr{T}_{2}.$

Equivalently, this means that the identity map 
\[
\left(X,\mathscr{T}_{2}\right)\to\left(X,\mathscr{T}_{1}\right)
\]
 is continuous.

}

\begin{remark}{}{}

A neighborhood basis of a measure $\mu,$ for any of these topologies
is given by the sets
\[
V_{\epsilon,f_{1},\cdots,f_{n}}\left(\mu\right)=\left\{ \nu\in\mathscr{M}:\,\sup_{1\leqslant i\leqslant n}\left|\intop f_{i}\text{d}\mu-\intop f_{i}\text{d}\nu\right|\leqslant\epsilon\right\} ,
\]
where $\epsilon>0$ and the functions $f_{i}$ belong respectively
to $\mathscr{C}_{\mathscr{K}}\left(\mathbb{R}^{d}\right),\,\mathscr{C}_{0}\left(\mathbb{R}^{d}\right)$
or $\mathscr{C}_{b}\left(\mathbb{R}^{d}\right)$ depending on whether
we consider the vague, weak or narrow topology.

It is immediate that the vague topology is less fine than the weak
topology, which itself is less fine than the narrow topology. Consequently,
for a sequence
\begin{itemize}
\item \textbf{Narrow convergence implies weak convergence, }
\item \textbf{Weak convergence implies vague convergence}.
\end{itemize}
Moreover, the weak topology on $\mathscr{M}\left(1\right)$ is strictly
less fine than the narrow topology, as shown by the following example. 

Let $x$ be a nonzero vector of $\mathbb{R}^{d},$ and define $\mu_{n}=\delta_{nx}.$
The sequence $\left(\mu_{n}\right)_{n\in\mathbb{N}}$ converges weakly
to the zero measure $\mu.$ Indeed, for every $f\in\mathscr{C}_{0}\left(\mathbb{R}^{d}\right),$
\[
\lim_{n\to+\infty}\intop f\text{d}\mu_{n}=\lim_{n\to+\infty}f\left(nx\right)=0.
\]
However, the sequence does not converge narrowly to $\mu.$ Indeed,
\[
\lim_{n\to+\infty}\intop1\text{d}\mu_{n}=1\,\,\,\,\text{while}\,\,\,\,\intop1\text{d}\mu=0.
\]
This example also shows that \textbf{$\mathscr{M}^{1}$ is not weakly
closed: a weak limit of probability measures need not to be a probability
measure.}

\end{remark}

We now compare the topologies on $\mathscr{M}\left(b\right)$---and
consequently, the associated notions of convergence for sequences
of measures with mass at most $b$---and examing some of their fundamental
properties.

\begin{proposition}{Comparison of Topologies}{}

(a) On $\mathscr{M}\left(b\right)$ the \textbf{vague} and \textbf{weak}
topologies coincide with the least fine topology that makes the functions
\[
\mu\mapsto\intop f\text{d}\mu
\]
continuous, where $f$ ranges over a \textbf{total} subset $\mathscr{H}$
in $\mathscr{C}_{\mathscr{K}}\left(\mathbb{R}^{d}\right)$ in the
vague case, or in $\mathscr{C}_{0}\left(\mathbb{R}^{d}\right)$ in
the weak case.

(b) On $\mathscr{M}^{1},$ the three topologies coincide.

(c) The space $\mathscr{M}\left(b\right)$ is metrizable and \textbf{compact}
for the \textbf{weak topology}.

\end{proposition}

\begin{figure}[t]
\begin{center}\includegraphics[width=0.4\textwidth]{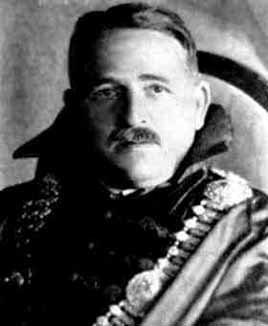}

{\tiny Public Domain}\end{center}

\caption{\textbf{\protect\href{https://en.wikipedia.org/wiki/Frigyes_Riesz}{Frigyes Riesz}}
(1880-1956)}\sindex[fam]{Riesz, Frigyes}
\end{figure}

\begin{proof}{}{}

First observe that if the functions 
\[
\mu\mapsto\intop f\text{d}\mu
\]
are continuous for every $f$ ranges over the set $\mathscr{H},$
then they are also continuous when $f$ ranges over the vector space
$\widetilde{\mathscr{H}}$ generated by $\mathscr{H}.$

(a) Let $f$ be an arbitrary function in $\mathscr{C}_{\mathscr{K}}\left(\mathbb{R}^{d}\right)$
or $\mathscr{C}_{0}\left(\mathbb{R}^{d}\right).$ Let $\epsilon>0.$
There exists a function $g$ in $\widetilde{\mathscr{H}}$ such that
\[
\left\Vert f-g\right\Vert \leqslant\dfrac{\epsilon}{4b}.
\]
By the triangle inequality,
\[
\left|\intop f\text{d}\mu-\intop f\text{d}\nu\right|\leqslant\left|\intop f\text{d}\mu-\intop g\text{d}\mu\right|+\left|\intop g\text{d}\mu-\intop g\text{d}\nu\right|+\left|\intop g\text{d}\nu-\intop f\text{d}\nu\right|.
\]
Hence,
\[
\left|\intop f\text{d}\mu-\intop f\text{d}\nu\right|\leqslant2b\left\Vert f-g\right\Vert +\left|\intop g\text{d}\mu-\intop g\text{d}\nu\right|.
\]
It follows that whenever $\nu\in V_{\frac{\epsilon}{2},g}\left(\mu\right),$
we have 
\[
\left|\intop f\text{d}\mu-\intop f\text{d}\nu\right|\leqslant\epsilon.
\]
In other words, 
\[
V_{\frac{\epsilon}{2},g}\left(\mu\right)\subset V_{\epsilon,f}\left(\mu\right),
\]
which shows that $V_{\epsilon,f}\left(\mu\right)$ is a neighborhood
of $\mu$ for the initial topology associated with $\widetilde{\mathscr{H}}.$ 

Since this topology is clearly no finer than the vague and weak topologies,
the result follows.

\begin{remark}{}{}

For example, the countable set $\mathscr{H}$ consisting of functions
of the form $x^{n}\text{e}^{-x^{2}},\,n\in\mathbb{N},$ is dense in
$\mathscr{C}_{0}\left(\mathbb{R}\right).$ 

By contrast, the analogous statement for the narrow topology is not
meaningful: a total set $\mathscr{H}$ in $\mathscr{C}_{b}\left(\mathbb{R}^{d}\right)$
cannot be countable and must necessarily be very large. 

\end{remark}

(b) We now show that the weak and narrow topologies coincide on $\mathscr{M}^{1}.$ 

To this end, it is enough to prove that for any $P\in\mathscr{M}^{1},$
every neighborhood of $P$ of the type $V_{\epsilon,f}\left(P\right)$
for the narrow topology---where $f\in\mathscr{C}_{b}\left(\mathbb{R}^{d}\right)$
and $\epsilon>0$ are arbitrary---is in fact a weak neighborhood.

Let $f\in\mathscr{C}_{b}\left(\mathbb{R}^{d}\right)$ and $\epsilon>0.$
Let $\left(h_{p}\right)_{p\in\mathbb{N}}$ be a sequence of nonnegative
functions in $\mathscr{C}_{0}\left(\mathbb{R}^{d}\right),$ such that
$h_{p}$ converges pointwise to 1 by nondecreasing. For each integer
$p,$ $fh_{p}\in\mathscr{C}_{0}\left(\mathbb{R}^{d}\right),$ and
for every $Q\in\mathscr{M}^{1},$
\[
\left|\intop f\text{d}P-\intop f\text{d}Q\right|\leqslant\left|\intop\left(f-fh_{p}\right)\text{d}P\right|+\left|\intop fh_{p}\text{d}P-\intop fh_{p}\text{d}Q\right|+\left|\intop\left(fh_{p}-f\right)\text{d}Q\right|.
\]
 Hence,
\[
\left|\intop f\text{d}P-\intop f\text{d}Q\right|\leqslant\left\Vert f\right\Vert \intop\left(1-h_{p}\right)\text{d}P+\left|\intop fh_{p}\text{d}P-\intop fh_{p}\text{d}Q\right|+\left\Vert f\right\Vert \intop\left(1-h_{p}\right)\text{d}Q
\]
Since $P$ and $Q$ are probabilities,
\[
\intop\left(1-h_{p}\right)\text{d}P+\intop\left(1-h_{p}\right)\text{d}Q=2\left(1-\intop h_{p}\text{d}P\right)+\intop h_{p}\text{d}P-\intop h_{p}\text{d}Q,
\]
and therefore, a fortiori,
\[
0\leqslant\intop\left(1-h_{p}\right)\text{d}P+\intop\left(1-h_{p}\right)\text{d}Q=2\left(1-\intop h_{p}\text{d}P\right)+\left|\intop h_{p}\text{d}P-\intop h_{p}\text{d}Q\right|.
\]
It follows that
\[
\left|\intop f\text{d}P-\intop f\text{d}Q\right|\leqslant\left\Vert f\right\Vert \left[2\left(1-\intop h_{p}\text{d}P\right)+\left|\intop h_{p}\text{d}P-\intop h_{p}\text{d}Q\right|\right]+\left|\intop fh_{p}\text{d}P-\intop fh_{p}\text{d}Q\right|.
\]
By decreasing monotone convergence---$P$ is a probability---,
\[
\lim_{p\to+\infty}\intop\left(1-h_{p}\right)\text{d}P=0.
\]
We therefore choose $p$ such that 
\[
0\leqslant\intop\left(1-h_{p}\right)\text{d}P\leqslant\dfrac{\epsilon}{4\left\Vert f\right\Vert }.
\]
Then, for every 
\[
Q\in\mathscr{M}^{1}\cap V_{\frac{\epsilon}{4\left\Vert f\right\Vert },h_{p}}\left(P\right)\cap V_{\frac{\epsilon}{2},fh_{p}}\left(P\right),
\]
we obtain
\[
\left|\intop f\text{d}P-\intop f\text{d}Q\right|\leqslant\epsilon,
\]
and therefore 
\[
Q\in\mathscr{M}^{1}\cap V_{\epsilon,f}\left(P\right).
\]
Thus, we have shown that
\[
\mathscr{M}^{1}\cap V_{\frac{\epsilon}{4\left\Vert f\right\Vert },h_{p}}\left(P\right)\cap V_{\frac{\epsilon}{2},fh_{p}}\left(P\right)\subset\mathscr{M}^{1}\cap V_{\epsilon,f}\left(P\right).
\]
Since $h_{p}$ and $fh_{p}$ belong to $\mathscr{C}_{0}\left(\mathbb{R}^{d}\right),$
this proves that 
\[
\mathscr{M}^{1}\cap V_{\epsilon,f}\left(P\right)
\]
is a neighborhood of $P$ for the weak topology on $\mathscr{M}^{1}.$ 

Because the weak topology is coarser than the narrow topology, we
conclude that the two topologies coincide on $\mathscr{M}^{1}.$

(c) Choose a sequence $\left(f_{n}\right)_{n\in\mathbb{N}}$ of elements
of $\mathscr{C}_{\mathscr{K}}\left(\mathbb{R}^{d}\right)$ that is
less dense in $\mathscr{C}_{0}\left(\mathbb{R}^{d}\right).$ Define
the distance $d$ on $\mathscr{M}\left(b\right)$ by
\[
d\left(\mu,\nu\right)=\sum^{+\infty}_{n=1}\dfrac{1}{2^{n}\left\Vert f_{n}\right\Vert }\left|\intop f_{n}\text{d}\mu-\intop f_{n}\text{d}\nu\right|.
\]
This is indeed a distance: $d\left(\mu,\nu\right)$ is always finite,
and symmetry and the triangle inequality are immediate. Moreover,
if $d\left(\mu,\nu\right)=0,$ then for every $n\in\mathbb{N},$ 
\[
\intop f_{n}\text{d}\mu=\intop f_{n}\text{d}\nu.
\]
By density, it follows that for every $f\in\mathscr{C}_{\mathscr{K}}\left(\mathbb{R}^{d}\right),$
\[
\intop f\text{d}\mu=\intop f\text{d}\nu,
\]
and therefore $\mu=\nu$---see Chapter \ref{chap:PartIIChap9}, Corollary
$\ref{co:equality_radon_measures}.$ By $\left(a\right),$ taking
for $\mathscr{H}$ the total set constituted by the functions $f_{n},\,n\in\mathbb{N},$
the topology induced by $d$ coïncides with the weak topology.

\textbf{To show that $\mathscr{M}\left(b\right)$ is weakly compact},
it is therefore enough to prove that from every sequence, one can
extract a convergent subsequence. Let then, for every $p,$ $\mu_{p}\in\mathscr{M}\left(b\right).$
We use a \textbf{diagonal extraction} process to obtain from the sequence
$\left(\mu_{p}\right)_{p\in\mathbb{N}}$ a convergent subsequence. 

Fix $n\in\mathbb{N},$ the real number sequence 
\[
\left(\intop f_{n}\text{d}\mu_{p}\right)_{p\in\mathbb{N}}
\]
is bounded by $b\left\Vert f_{n}\right\Vert .$ Hence, we can then
extract a convergent subsequence from the sequence $\left(\intop f_{1}\text{d}\mu_{p}\right)_{p\in\mathbb{N}}.$
Let $\varphi_{1}$ denote the nondecreasing injection from $\mathbb{N}$
to itself that defines the extracted sequence.

Applying the same process, we can extract a convergent subsequence
of the sequence 
\[
\left(\intop f_{2}\text{d}\mu_{\varphi_{1}\left(p\right)}\right)_{p\in\mathbb{N}}.
\]
and we denote by $\varphi_{2}$ the nondecreasing injection from $\mathbb{N}$
to itself that defines the extracted sequence. Then both sequences
$\left(\intop f_{1}\text{d}\mu_{\varphi_{2}\left(p\right)}\right)_{p\in\mathbb{N}}$
and $\left(\intop f_{2}\text{d}\mu_{\varphi_{2}\left(p\right)}\right)_{p\in\mathbb{N}}$
are convergent. 

Proceeding by induction, we similarly construct for each integer $k\geqslant1,$
the sequence $\left(\mu_{\varphi_{k}\left(p\right)}\right)_{p\in\mathbb{N}},$
subsequence of $\left(\mu_{\varphi_{k-1}\left(p\right)}\right)_{p\in\mathbb{N}}$
such that for every $i\leqslant k,$ the sequence 
\[
\left(\intop f_{i}\text{d}\mu_{\varphi_{k}\left(p\right)}\right)_{p\in\mathbb{N}}
\]
converges. 

Then, for every integer $k\in\mathbb{N},$ the sequence $\left(\intop f_{k}\text{d}\mu_{\varphi_{p}\left(p\right)}\right)_{p\in\mathbb{N}}$
converges: from a certain rank $k,$ it is a subsequence of the convergent
sequence $\left(\intop f_{k}\text{d}\mu_{\varphi_{k}\left(p\right)}\right)_{p\in\mathbb{N}}.$ 

By density, it follows that, for every $f\in\mathscr{C}_{0}\left(\mathbb{R}^{d}\right),$
the sequence 
\[
\left(\intop f\text{d}\mu_{\varphi_{k}\left(k\right)}\right)_{k\in\mathbb{N}}
\]
is convergent of limit, denoted by $\Psi\left(f\right).$ The function
$\Psi$ is a nonnegative linear form on $\mathscr{C}_{\mathscr{K}}\left(\mathbb{R}^{d}\right).$
By the \textbf{\index{Riesz theorem}Riesz\sindex[fam]{Riesz, Frigyes}
representation }theorem\footnotemark\footnotemark---see for instance,
\cite{metivier2014stochastic} (p.156 in the 1980 edition)---there
exists a unique measure $\mu$ such that $\Psi\left(f\right)=\intop f\text{d}\mu$
for every $f\in\mathscr{C}_{\mathscr{K}}\left(\mathbb{R}^{d}\right).$ 

Using again the density of $\mathscr{C}_{\mathscr{K}}\left(\mathbb{R}^{d}\right)$
in $\mathscr{C}_{0}\left(\mathbb{R}^{d}\right),$ we conclude that
the subsequence $\left(\mu_{\varphi\left(k\right)}\right)_{k\in\mathbb{N}}$
converges weakly to $\mu.$ 

It remains to check that $\mu\in\mathscr{M}\left(b\right).$ Let $\left(h_{p}\right)_{p\in\mathbb{N}}$
be a sequence of nonnegative functions in $\mathscr{C}_{\mathscr{K}}\left(\mathbb{R}^{d}\right)$
that converges pointwise to $1$ by nondecreasing. For every $p,$
\[
0\leqslant\lim_{k\to+\infty}\intop h_{p}\text{d}\mu_{\varphi_{k}\left(k\right)}=\intop h_{p}\text{d}\mu\leqslant b.
\]
By monotone convergence,
\[
\lim_{p\to+\infty}\intop h_{p}\text{d}\mu=\mu\left(\mathbb{R}^{d}\right)\leqslant b.
\]

\end{proof}

\addtocounter{footnote}{-1}

\footnotetext{\textbf{\href{https://en.wikipedia.org/wiki/Frigyes_Riesz}{Frigyes Riesz}\sindex[fam]{Riesz, Frigyes}}
(1880-1956) was a Hungarian mathematician who made foundational contributions
to functional analysis. His work has had significant impact on the
development of modern analysis and theoretical physics.}

\stepcounter{footnote}

\footnotetext{\textbf{Riesz representation theorem\index{Riesz theorem}.}
Let $\Phi$ be a nonnegative linear form on $\mathscr{C}_{\mathscr{K}}\left(\mathbb{R}^{d}\right).$
There exists a unique measure $\mu$ on $\mathbb{R}^{d},$ equipped
with its Borel $\sigma-$algebra, representing $\Phi,$ in the sense
that for every $f\in\mathscr{C}_{\mathscr{K}}\left(\mathbb{R}^{d}\right),$
\[
\Phi\left(f\right)=\intop_{\mathbb{R}^{d}}f\text{d}\mu.
\]
This measure has the following properties:

(i) $\mu$ is finite on every compact---that is, $\mu$ is a Radon
measure.

(ii) For every $B\in\mathscr{B}_{\mathbb{R}^{d}},$
\begin{gather*}
\mu\left(B\right)=\inf\left\{ \mu\left(O\right):\,O\,\text{open,}\,\,O\supset B\right\} ,
\end{gather*}
and
\[
\mu\left(B\right)=\sup\left\{ \mu\left(K\right):\,K\,\text{compact,\,}\,K\subset B\right\} .
\]
}

\begin{denotation}{}{}We traditionally write $\mu_{n}\Rightarrow\mu$
to denote that the sequence $\left(\mu_{n}\right)_{n\in\mathbb{N}}$
converges narrowly to $\mu.$

\end{denotation}

\begin{remark}{}{}

As a consequence of the previous proposition, a sequence $\left(\mu_{n}\right)_{n\in\mathbb{N}}$
in $\mathscr{M}\left(b\right)$ converges weakly to $\mu$ if and
only if the sequence $\left(\intop f\text{d}\mu_{n}\right)_{n\in\mathbb{N}}$
converges to $\intop f\text{d}\mu$ for every $f$ ranging over a
total subset $\mathscr{H}$ of $\mathscr{C}_{0}\left(\mathbb{R}^{d}\right),$
for instance $\mathscr{C}_{\mathscr{K}}\left(\mathbb{R}^{d}\right)$
itself. 

Moreover, if the measures $\mu_{n}$ and $\mu$ are probability measures,
then the sequence $\left(\mu_{n}\right)_{n\in\mathbb{N}}$ converges
narrowly to $\mu$ if and only if the sequence $\left(\intop f\text{d}\mu_{n}\right)_{n\in\mathbb{N}}$
converges to $\intop f\text{d}\mu$ for every $f$ ranging over the
total subset $\mathscr{H}$ of $\mathscr{C}_{\mathscr{K}}\left(\mathbb{R}^{d}\right).$

\end{remark}

\begin{proposition}{Necessary and Sufficient Condition for Narrow Convergence}{narrow_cv_nsc}

Let $\mu_{n},\,n\in\mathbb{N},$ and $\mu$ be measures of $\mathscr{M}\left(b\right)$
such that the sequence $\left(\mu_{n}\right)_{n\in\mathbb{N}}$ converges
vaguely---or weakely---to $\mu.$ 

Then the sequence $\left(\mu_{n}\right)_{n\in\mathbb{N}}$ converges
narrowly to $\mu$ if and only if 
\[
\lim_{n\to+\infty}\mu_{n}\left(\mathbb{R}^{d}\right)=\mu\left(\mathbb{R}^{d}\right).
\]

\end{proposition}

\begin{proof}{}{}

The condition is clearly necessary. 

For sufficiency, let $f\in\mathscr{C}_{b}\left(\mathbb{R}^{d}\right)$
and $\epsilon>0$ be arbitrary. Since the measure $\mu$ has finite
mass, we can choose $\varphi\in\mathscr{C}^{+}_{\mathscr{K}}\left(\mathbb{R}^{d}\right)$
such that
\[
\left\Vert \varphi\right\Vert \leqslant1\,\,\,\,\text{and}\,\,\,\,0\leqslant\intop\left(1-\varphi\right)\text{d}\mu\leqslant\dfrac{\epsilon}{8\left\Vert f\right\Vert }.
\]
Because $f\varphi\in\mathscr{C}_{\mathscr{K}}\left(\mathbb{R}^{d}\right)$
and $\mu_{n}$ converges vaguely---or weakly---to $\mu,$ there
exists an integer $N_{1}$ such that, for every $n\geqslant N_{1},$
\[
\left|\intop f\varphi\text{d}\mu_{n}-\intop f\varphi\text{d}\mu\right|\leqslant\dfrac{\epsilon}{2}.
\]
Moreover,
\[
\intop\left(1-\varphi\right)\text{d}\mu_{n}=\mu_{n}\left(\mathbb{R}^{d}\right)-\intop\varphi\text{d}\mu_{n}.
\]
It follows from the hypothesis that the sequence with general term
$\intop\left(1-\varphi\right)\text{d}\mu_{n}$ converges to $\mu\left(\mathbb{R}^{d}\right)-\intop\varphi\text{d}\mu=\intop\left(1-\varphi\right)\text{d}\mu.$ 

Thus, there exists an integer $N_{2}$ such that, for every $n\geqslant N_{2},$
\[
0\leqslant\intop\left(1-\varphi\right)\text{d}\mu_{n}\leqslant\dfrac{\epsilon}{4\left\Vert f\right\Vert }.
\]
Let $N=\max\left(N_{1},N_{2}\right).$ For every $n\geqslant N,$
the triangle inequality gives
\[
\left|\intop f\text{d}\mu_{n}-\intop f\text{d}\mu\right|\leqslant\left|\intop\left(f-f\varphi\right)\text{d}\mu_{n}\right|+\left|\intop f\varphi\text{d}\mu_{n}-\intop f\varphi\text{d}\mu\right|+\left|\intop\left(f\varphi-f\right)\text{d}\mu\right|.
\]
Hence, a fortiori, 
\[
\left|\intop f\text{d}\mu_{n}-\intop f\text{d}\mu\right|\leqslant\left\Vert f\right\Vert \intop\left(1-\varphi\right)\text{d}\mu_{n}+\left|\intop f\varphi\text{d}\mu_{n}-\intop f\varphi\text{d}\mu\right|+\left\Vert f\right\Vert \intop\left(1-\varphi\right)\text{d}\mu,
\]
and, consequentely
\[
\left|\intop f\text{d}\mu_{n}-\intop f\text{d}\mu\right|\leqslant\epsilon,
\]
which shows that the sequence $\left(\mu_{n}\right)_{n\in\mathbb{N}}$
converges narrowly to $\mu.$ 

\end{proof}

\begin{remark}{}{}

On the measurable space $\left(\mathbb{R},\mathscr{B}_{\mathbb{R}}\right),$
the sequence of Dirac measures $\delta_{\frac{1}{n}}$ converges narrowly
to the Dirac measure $\delta_{0}.$ Nevertheless, for every $n\in\mathbb{N}^{\ast},$
$\delta_{\frac{1}{n}}\left(\left\{ 0\right\} \right)=0,$and therefore
\[
\lim_{n\to+\infty}\delta_{\frac{1}{n}}\left(\left\{ 0\right\} \right)\neq\delta_{0}\left(\left\{ 0\right\} \right).
\]
Thus, narrow convergence of a sequence of measures $\left(\mu_{n}\right)_{n\in\mathbb{N}}$
to a measure $\mu$ does not imply that, for every Borel $A,$ the
sequence $\left(\mu_{n}\left(A\right)\right)_{n\in\mathbb{N}}$ converges.
The following Proposition $\ref{pr:narrow_cv_criterion}$ gives an
answer to this convergence problem.

\end{remark}

\begin{definition}{$\mu-$continuity of a Borelian}{}

Let $\mu$ be a measure on $\mathbb{R}^{d}.$ A Borel set $A$ is
said to be $\mu-$continuous if 
\[
\mu\left(\partial A\right)=0,
\]
 where $\partial A$ denotes the boundary of $A.$

\end{definition}

\begin{proposition}{Narrow Convergence Criteria}{narrow_cv_criterion}

Let $\mu_{n},\,n\in\mathbb{N},$ and $\mu$ be measures of $\mathscr{M}\left(b\right).$
The following assertions are equivalent:

(i) The sequence $\left(\mu_{n}\right)_{n\in\mathbb{N}}$ converges
narrowly to $\mu.$

(ii) For every \textbf{closed set} $F\subset\mathbb{R}^{d},$ 
\[
\limsup_{n\to+\infty}\mu_{n}\left(F\right)\leqslant\mu\left(F\right),
\]
and moreover,
\[
\lim_{n\to+\infty}\mu_{n}\left(\mathbb{R}^{d}\right)=\mu\left(\mathbb{R}^{d}\right).
\]

(iii) For every \textbf{open set} $O\subset\mathbb{R}^{d},$ 
\[
\liminf_{n\to+\infty}\mu_{n}\left(O\right)\geqslant\mu\left(O\right),
\]
and moreover,
\[
\lim_{n\to+\infty}\mu_{n}\left(\mathbb{R}^{d}\right)=\mu\left(\mathbb{R}^{d}\right).
\]

(iv) For every \textbf{Borel} set $A$ that is\textbf{ $\mu-$continuous},
\[
\lim_{n\to+\infty}\mu_{n}\left(A\right)=\mu\left(A\right).
\]

\end{proposition}

\begin{proof}{}{}We prove it using the following schema:
\[
\left(\text{i}\right)\Leftrightarrow\left(\text{ii}\right)\Leftrightarrow\left(\text{iii}\right)\Rightarrow\left(\text{iv}\right)\Rightarrow\left(\text{ii}\right).
\]

\begin{itemize}
\item \textbf{$\left(\text{i}\right)\Rightarrow\left(\text{ii}\right)$}\\
For every $j\in\mathbb{N}^{\ast},$ define $\varphi_{j}\in\mathscr{C}_{b}\left(\mathbb{R}\right)$
by
\[
\varphi_{j}\left(u\right)=\begin{cases}
1, & \text{if }u\leqslant0,\\
1-ju & \text{if }0<u<\dfrac{1}{j},\\
0, & \text{if }u\geqslant\dfrac{1}{j}.
\end{cases}
\]
The sequence $\left(\varphi_{j}\right)_{j\in\mathbb{N}^{\ast}}$ is
nonincreasing and converges pointwise to $\boldsymbol{1}_{\left]-\infty,0\right]}.$
Let $F$ be a closed subset of $\mathbb{R}^{d}.$ Define the function
$f_{j}\in\mathscr{C}_{b}\left(\mathbb{R}^{d}\right)$ by setting,
for $x\in\mathbb{R}^{d},$ $f_{j}\left(x\right)=\varphi_{j}\left(d\left(x,F\right)\right).$
Then, for every $x\in\mathbb{R}^{d},$ 
\[
\lim_{j\to+\infty}\searrow f_{j}\left(x\right)=\boldsymbol{1}_{\left]-\infty,0\right]}\left(d\left(x,F\right)\right)=\boldsymbol{1}_{F}\left(x\right).
\]
The functions $f_{j}$ being bounded by 1 and the measure $\mu$ being
bounded, we may apply the decreasing monotone convergence theorem,
and obtain
\[
\lim_{j\to+\infty}\intop f_{j}\text{d\ensuremath{\mu}}=\mu\left(F\right).
\]
Let $\epsilon>0.$ Then there exists $j_{0}$ such that
\[
\mu\left(F\right)\leqslant\intop f_{j_{0}}\text{d\ensuremath{\mu}}\leqslant\mu\left(F\right)+\epsilon.
\]
Since, by hypothesis, the sequence $\left(\mu_{n}\right)_{n\in\mathbb{N}}$
narrowly converges to $\mu,$ we then have
\[
\lim_{n\to+\infty}\intop f_{j_{0}}\text{d}\mu_{n}=\intop f_{j_{0}}\text{d}\mu\leqslant\mu\left(F\right)+\epsilon.
\]
Moreover, for every $n\in\mathbb{N},$ 
\[
\mu_{n}\left(F\right)\leqslant\intop f_{j_{0}}\text{d}\mu_{n}.
\]
Hence,
\[
\limsup_{n\to+\infty}\mu_{n}\left(F\right)\leqslant\limsup_{n\to+\infty}\intop f_{j_{0}}\text{d}\mu_{n}\leqslant\mu\left(F\right)+\epsilon.
\]
Since $\epsilon>0$ is arbitrary, this proves that 
\[
\limsup_{n\to+\infty}\mu_{n}\left(F\right)\leqslant\mu\left(F\right).
\]
\item \textbf{$\left(\text{ii}\right)\Rightarrow\left(\text{iii}\right)$}\\
Let $O$ be an open set. Since $\mu$ is bounded, we can pass to complements.
For every $n\in\mathbb{N},$
\[
\mu_{n}\left(O\right)=\mu_{n}\left(\mathbb{R}^{d}\right)-\mu_{n}\left(O^{c}\right).
\]
Taking lower limits and using the assumptions in (ii), we get
\[
\liminf_{n\to+\infty}\mu_{n}\left(O\right)=\lim_{n\to+\infty}\mu_{n}\left(\mathbb{R}^{d}\right)-\limsup_{n\to+\infty}\mu_{n}\left(O^{c}\right)\geqslant\mu\left(\mathbb{R}^{d}\right)-\mu\left(O^{c}\right)=\mu\left(O\right).
\]
\item \textbf{$\left(\text{iii}\right)\Rightarrow\left(\text{ii}\right)$}\\
The argument is the same, again by taking complements.
\item \textbf{$\left(\text{iii}\right)\Rightarrow\left(\text{iv}\right)$}\\
Let $A$ be a Borel set that is $\mu-$continuous. Since $\left(\text{ii}\right)$
and $\left(\text{iii}\right)$ are equivalent, and using the nondecreasing
of measures together with the inclusions 
\[
\mathring{A}\subset A\subset\overline{A},
\]
we obtain
\[
\mu\left(\mathring{A}\right)\leqslant\liminf_{n\to+\infty}\mu_{n}\left(\mathring{A}\right)\leqslant\liminf_{n\to+\infty}\mu_{n}\left(A\right)\leqslant\limsup_{n\to+\infty}\mu_{n}\left(A\right)\leqslant\limsup_{n\to+\infty}\mu_{n}\left(\overline{A}\right)\leqslant\mu\left(\overline{A}\right).
\]
Now, since $\partial A=\overline{A}\backslash\mathring{A}$ and $\mu\left(\partial A\right)=0,$
we have $\mu\left(\overline{A}\right)=\mu\left(\mathring{A}\right)=\mu\left(A\right),$
which proves that
\[
\liminf_{n\to+\infty}\mu_{n}\left(A\right)=\limsup_{n\to+\infty}\mu_{n}\left(A\right)=\mu\left(A\right),
\]
that is the sequence with general term $\mu_{n}\left(A\right)$ converges
to $\mu\left(A\right).$
\item \textbf{$\left(\text{ii}\right)\Rightarrow\left(\text{i}\right)$
and thus $\left(\text{iii}\right)\Rightarrow\left(\text{i}\right)$}\\
We first establish an identity analogous to the usual formula expressing
the expectation of a nonnegative random variable in terms of the tail
of its cumulative distribution function: if $X$ is a nonnegative
random variable, then
\[
\intop_{\Omega}X\text{d}P=\intop_{\mathbb{R}}P\left(X>x\right)\text{d}\lambda\left(x\right)=\intop_{\mathbb{R}}P\left(X\geqslant x\right)\text{d}\lambda\left(x\right).
\]
By the Fubini theorem, for every $f\in\mathscr{C}^{+}_{b}\left(\mathbb{R}\right)$
and every $\mu\in\mathscr{M}\left(b\right),$ if $\lambda$ denotes
the Lebesgue measure on $\mathbb{R},$ then
\begin{align*}
\intop_{\left[0,\left\Vert f\right\Vert \right]}\mu\left(f\geqslant u\right)\text{d}\lambda\left(u\right) & =\intop_{\left[0,\left\Vert f\right\Vert \right]}\left[\intop_{\mathbb{R}^{d}}\boldsymbol{1}_{\left(f\left(x\right)\geqslant u\right)}\text{d}\mu\left(x\right)\right]\text{d}\lambda\left(u\right)\\
 & =\intop_{\mathbb{R}^{d}}\left[\intop_{\left[0,\left\Vert f\right\Vert \right]}\boldsymbol{1}_{\left(\left[0,\left\Vert f\right\Vert \right]\right)}\boldsymbol{1}_{\left(f\left(x\right)\geqslant u\right)}\text{d}\lambda\left(u\right)\right]\text{d}\mu\left(x\right).
\end{align*}
Hence,
\begin{equation}
\intop_{\left[0,\left\Vert f\right\Vert \right]}\mu\left(f\geqslant u\right)\text{d}\lambda\left(u\right)=\intop_{\mathbb{R}^{d}}f\text{d}\mu.\label{eq:int_0_norm_f_mufoverudlambau}
\end{equation}
Similarly,
\begin{equation}
\intop_{\left[0,\left\Vert f\right\Vert \right]}\mu\left(f>u\right)\text{d}\lambda\left(u\right)=\intop_{\mathbb{R}^{d}}f\text{d}\mu.\label{eq:int_0_norm_f_mufoverudlambau-1}
\end{equation}
Now let $f\in\mathscr{C}^{+}_{b}\left(\mathbb{R}^{d}\right).$ Since
the measure $\lambda\mid_{\left[0,\left\Vert f\right\Vert \right]}$
is bounded, the Fatou-Lebesgue lemma yields
\begin{align*}
\limsup_{n\to+\infty}\intop_{\mathbb{R}^{d}}f\text{d}\mu_{n} & =\limsup_{n\to+\infty}\intop_{\left[0,\left\Vert f\right\Vert \right]}\mu_{n}\left(f\geqslant u\right)\text{d}\lambda\left(u\right)\\
 & \leqslant\intop_{\left[0,\left\Vert f\right\Vert \right]}\limsup_{n\to+\infty}\mu_{n}\left(f\geqslant u\right)\text{d}\lambda\left(u\right).
\end{align*}
Since $\left(f\geqslant u\right)$ is \textbf{closed} and since we
suppose the assertion $\left(ii\right)$ to be true, we obtain
\[
\limsup_{n\to+\infty}\intop_{\mathbb{R}^{d}}f\text{d}\mu_{n}\leqslant\intop_{\left[0,\left\Vert f\right\Vert \right]}\mu\left(f\geqslant u\right)\text{d}\lambda\left(u\right)=\intop_{\mathbb{R}^{d}}f\text{d}\mu.
\]
Similarly,
\begin{align*}
\liminf_{n\to+\infty}\intop_{\mathbb{R}^{d}}f\text{d}\mu_{n} & =\liminf_{n\to+\infty}\intop_{\left[0,\left\Vert f\right\Vert \right]}\mu_{n}\left(f>u\right)\text{d}\lambda\left(u\right)\\
 & \geqslant\intop_{\left[0,\left\Vert f\right\Vert \right]}\liminf_{n\to+\infty}\mu_{n}\left(f>u\right)\text{d}\lambda\left(u\right).
\end{align*}
Since $\left(f>u\right)$ is an open set, and assuming the assertion
$\left(\text{iii}\right),$ we obtain
\[
\liminf_{n\to+\infty}\intop_{\mathbb{R}^{d}}f\text{d}\mu_{n}\geqslant\intop_{\left[0,\left\Vert f\right\Vert \right]}\mu\left(f>u\right)\text{d}\lambda\left(u\right)=\intop_{\mathbb{R}^{d}}f\text{d}\mu.
\]
Hence, we proved that for every $f\in\mathscr{C}^{+}_{b}\left(\mathbb{R}\right),$
\[
\lim_{n\to+\infty}\intop_{\mathbb{R}^{d}}f\text{d}\mu_{n}=\intop_{\mathbb{R}^{d}}f\text{d}\mu.
\]
By linearity, the same holds for every $f\in\mathscr{C}_{b}\left(\mathbb{R}\right)$
of arbitrary sign---it suffices to apply the previous result to the
nonnegative function $\left\Vert f\right\Vert -f.$
\item \textbf{$\left(\text{iv}\right)\Rightarrow\left(\text{ii}\right)$}\\
Let $F$ be an arbitrary closed set. For each $\epsilon>0,$ define
\[
F_{\epsilon}=\left\{ x\in\mathbb{R}^{d}:\,d\left(x,F\right)\leqslant\epsilon\right\} .
\]
Consider the function $\Phi$ from $\left[0,1\right]$ to $\left[0,b\right]$
defined by $\Phi\left(\epsilon\right)=\mu\left(F_{\epsilon}\right).$\\
It is nondecreasing and bounded, hence, it admits a set $I$ of discontinuity
points at most countable. \\
For every $\epsilon\in\left[0,1\right]\backslash I,$ $F_{\epsilon}$
being closed,
\[
\partial F_{\epsilon}=F_{\epsilon}\backslash\left[\bigcup_{n\in\mathbb{N}^{\ast}}F_{\epsilon-\frac{1}{n}}\right].
\]
Therefore,
\[
\mu\left(\partial F_{\epsilon}\right)=\mu\left(F_{\epsilon}\right)-\lim_{n\to+\infty}\mu\left(F_{\epsilon-\frac{1}{n}}\right).
\]
Because $\epsilon$ is a continuity point of $\Phi,$ it follows that
\[
\mu\left(\partial F_{\epsilon}\right)=0.
\]
Hence, we can choose a sequence $\left(\epsilon_{k}\right)_{k\in\mathbb{N}^{\ast}}$
nonincreasing to 0 such that, for every $k\in\mathbb{N}^{\ast},$
\[
\mu\left(\partial F_{\epsilon_{k}}\right)=0.
\]
By hypothesis, for every $k\in\mathbb{N}^{\ast},$
\[
\limsup_{n\to+\infty}\mu_{n}\left(F\right)\leqslant\limsup_{n\to+\infty}\mu_{n}\left(F_{\epsilon_{k}}\right)=\lim_{n\to+\infty}\mu_{n}\left(F_{\epsilon_{k}}\right)=\mu\left(F_{\epsilon_{k}}\right).
\]
As $F=\bigcap_{k\in\mathbb{N}^{\ast}}F_{\epsilon_{k}},$ and the sequence
of sets $F_{\epsilon_{k}}$ is nonincreasing, we have 
\[
\lim_{k\to+\infty}\mu\left(F_{\epsilon_{k}}\right)=\mu\left(F\right),
\]
which yields 
\[
\limsup_{n\to+\infty}\mu_{n}\left(F\right)\leqslant\mu\left(F\right).
\]
Last, since $\mathbb{R}^{d}$ is a set of $\mu-$continuity, then
\[
\limsup_{n\to+\infty}\mu_{n}\left(\mathbb{R}^{d}\right)=\mu\left(\mathbb{R}^{d}\right).
\]
\end{itemize}
\end{proof}

\begin{remark}{}{}

On the measurable space $\left(\mathbb{R},\mathscr{B}_{\mathbb{R}}\right),$
consider, for every $n\in\mathbb{N}^{\ast},$ the measure 
\[
\mu_{n}=\dfrac{1}{n}\sum^{n}_{j=1}\delta_{\frac{j}{n}}.
\]
This sequence of measure converges narrowly to $\boldsymbol{1}_{\left[0,1\right]}\cdot\lambda,$
where $\lambda$ denotes the Lebesgue measure on $\mathbb{R}.$ Indeed,
for every $f\in\mathscr{C}_{b}\left(\mathbb{R}\right),$
\[
\lim_{n\to+\infty}\intop_{\mathbb{R}}f\text{d}\mu_{n}=\lim_{n\to+\infty}\dfrac{1}{n}\sum^{n}_{j=1}f\left(\dfrac{j}{n}\right)=\intop^{1}_{0}f\left(x\right)\text{d}x=\intop_{\mathbb{R}}f\boldsymbol{1}_{\left[0,1\right]}\text{d}\lambda.
\]
The sum $\dfrac{1}{n}\sum^{n}_{j=1}f\left(\dfrac{j}{n}\right)$ is
precisely the Riemann sum of $f$ associated with the subdivision
of $\left[0,1\right]$ given by the points $\dfrac{j}{n}.$ 

Nevertheless, for every $n\in\mathbb{N},$
\[
\mu_{n}\left(\mathbb{Q}\cap\left[0,1\right]\right)=1,
\]
so 
\[
\lim_{n\to+\infty}\mu_{n}\left(\mathbb{Q}\cap\left[0,1\right]\right)=1,
\]
whereas 
\[
\boldsymbol{1}_{\left[0,1\right]}\cdot\lambda\left(\mathbb{Q}\cap\left[0,1\right]\right)=0.
\]
 This does not contradict the previous proposition, because $\mathbb{Q}\cap\left[0,1\right]$
is not a set of $\boldsymbol{1}_{\left[0,1\right]}\cdot\lambda-$continuity.
Indeed 
\[
\partial\left(\mathbb{Q}\cap\left[0,1\right]\right)=\left[0,1\right]
\]
 and therefore
\[
\boldsymbol{1}_{\left[0,1\right]}\cdot\lambda\left(\partial\left(\mathbb{Q}\cap\left[0,1\right]\right)\right)=1.
\]

\end{remark}

\begin{definition}{Tight Sequence}{}

A \textbf{sequence} $\left(\mu_{n}\right)_{n\in\mathbb{N}}$ of \textbf{measures}
in $\mathscr{M}\left(b\right)$ is said to be \textbf{\index{tight measures sequence}tight}
if, for every $\epsilon>0,$ there exists a compact $K$ on $\mathbb{R}^{d}$
such that 
\[
\sup_{n\in\mathbb{N}}\mu_{n}\left(K^{c}\right)\leqslant\epsilon.
\]

\end{definition}

\begin{corollary}{Sufficient Condition To Obtain a Sequence of Tight Sequence}{sc_tight_seq}

If the sequence $\left(\mu_{n}\right)_{n\in\mathbb{N}}$ of measures
in $\mathscr{M}\left(b\right)$ converges narrowly to some $\mu\in\mathscr{M}\left(b\right),$
then the sequence is \textbf{tight}.

\end{corollary}

\begin{proof}{}{}

Let $\epsilon>0.$ Choose an open ball $O$ such that 
\[
\mu\left(O\right)\geqslant\mu\left(\mathbb{R}^{d}\right)-\dfrac{\epsilon}{2}.
\]
Since 
\[
\lim_{n\to+\infty}\mu_{n}\left(\mathbb{R}^{d}\right)=\mu\left(\mathbb{R}^{d}\right),
\]
there exists an integer $N_{1}$ such that for every $n\geqslant N_{1},$
\[
\mu_{n}\left(\mathbb{R}^{d}\right)\leqslant\mu\left(\mathbb{R}^{d}\right)+\dfrac{\epsilon}{2}.
\]
Moreover, by Proposition $\ref{pr:narrow_cv_criterion},$ 
\[
\liminf_{n\to+\infty}\mu_{n}\left(O\right)\geqslant\mu\left(O\right)\geqslant\mu\left(\mathbb{R}^{d}\right)-\dfrac{\epsilon}{2}.
\]

Hence, there exists an integer $N_{2}$ such that, for every $n\geqslant N_{2},$
\[
\mu_{n}\left(O\right)\geqslant\mu\left(\mathbb{R}^{d}\right)-\dfrac{\epsilon}{2}.
\]
Let $N=\max\left(N_{1},N_{2}\right).$ Then, for every $n\geqslant N,$
\[
\mu_{n}\left(O^{c}\right)=\mu_{n}\left(\mathbb{R}^{d}\right)-\mu_{n}\left(O\right)\leqslant\left[\mu\left(\mathbb{R}^{d}\right)+\dfrac{\epsilon}{2}\right]-\mu\left(\mathbb{R}^{d}\right)+\dfrac{\epsilon}{2}=\epsilon.
\]
Therefore,
\begin{equation}
\sup_{n\geqslant N}\mu_{n}\left(O^{c}\right)\leqslant\epsilon.\label{eq:sup_mu_n_Oc}
\end{equation}
It remains to choose a compact $K$ containing $O$ such that $\mu_{n}\left(K^{c}\right)\leqslant\epsilon$
as soon as $0\leqslant n\leqslant N.$ This is possible because, for
each fixed $n,$ in a finite number, we have
\[
\lim_{p\to+\infty}\mu_{n}\left(B_{f}\left(0,p\right)^{c}\right)=0,
\]
where $B_{f}\left(0,p\right)$ denotes the closed ball centered at
0 with radius $p.$ 

Finally, since $K\supset O,$ we have $K^{c}\subset O^{c},$ so combining
with $\refpar{eq:sup_mu_n_Oc}$ gives
\[
\sup_{n\in\mathbb{N}}\mu_{n}\left(K^{c}\right)\leqslant\epsilon.
\]

\end{proof}

The famous and fundamental Lévy theorem characterizes narrow convergence
of measure sequences in terms of their Fourier transforms.

\begin{theorem}{Lévy Theorem}{levy_th_measures}

Let, for each $n\in\mathbb{N},$ consider a measure $\mu_{n}\in\mathscr{M}\left(b\right).$

(a) If the sequence $\left(\mu_{n}\right)_{n\in\mathbb{N}}$ converges
narrowly to $\mu,$ then the sequence $\left(\widehat{\mu}_{n}\right)_{n\in\mathbb{N}}$
of the Fourier transforms of $\mu_{n}$ converges pointwise to $\widehat{\mu},$
the Fourier transform of $\mu.$

(b) Conversely, if the sequence $\left(\widehat{\mu}_{n}\right)_{n\in\mathbb{N}}$
of the Fourier transforms of $\mu_{n}$ converges pointwise to a function
$\varphi$ that is continuous at 0, then there exists a unique measure
$\mu\in\mathscr{M}\left(b\right)$ such that $\varphi=\widehat{\mu}.$
Moreover, the sequence $\left(\mu_{n}\right)_{n\in\mathbb{N}}$ converges
narrowly to $\mu.$

(c) In each of the above situations, the convergence of the sequence
$\left(\widehat{\mu}_{n}\right)_{n\in\mathbb{N}}$ is uniform on every
compact of $\mathbb{R}^{d}.$

\end{theorem}

\begin{proof}{Beginning of the Proof}{}

First, observe that if a sequence $\left(\mu_{n}\right)_{n\in\mathbb{N}}$
of bounded measures on $\mathbb{R}^{d}$ converges weakly---respectively
narrowly---to a measure $\mu,$ then 
\[
\lim_{n\to+\infty}\intop f\text{d}\mu_{n}=\intop f\text{d}\mu
\]
for every continuous function $f$ on $\mathbb{R}^{d}$ that vanishes
at infinity---respectively for every bounded continuous function
on $\mathbb{R}^{d}$---and taking complex values. Indeed, it suffices
to observe that the convergence occurs for $\text{Re}\left(f\right)$
and $\text{Im\ensuremath{\left(f\right)}}.$ In this proof, the notation
$\mathscr{C}_{0}\left(\mathbb{R}^{d}\right)$ denotes the space of
complex-valued continuous function on $\mathbb{R}^{d}$ and taking
complex values.

(a) For every $t\in\mathbb{R}^{d},$ the function $x\mapsto\text{e}^{\text{i}\left\langle x,t\right\rangle }$
is continuous and bounded, and thus the sequence with general term
$\widehat{\mu}_{n}\left(t\right)$ converges to $\widehat{\mu}\left(t\right).$

(b) We start by proving that the sequence $\left(\mu_{n}\right)_{n\in\mathbb{N}}$
is \textbf{weakly convergent}. Since $\mathscr{M}\left(b\right)$
is metrizable and compact for the weak topology, for $\left(\mu_{n}\right)_{n\in\mathbb{N}}$
to be weakly convergent, it must and it is enough that this sequence
admits at most one weak adherence value.

Indeed, in a compact metric space, every sequence has at least one
adherent point, and a sequence that admits only one adherent point
converges to this point.

Let $\mu$ be a weak adherence value of the sequence $\left(\mu_{n}\right)_{n\in\mathbb{N}}$
and let $\left(\mu_{\psi\left(n\right)}\right)_{n\in\mathbb{N}}$
be a subsequence that converges weakly to $\mu$---where $\psi$
denotes the injection from $\mathbb{N}$ to $\mathbb{N}$ defining
the subsequence. 

We are going to show that $\mu_{\psi\left(n\right)}$ tends narrowly
to $\mu$ as $n$ tends to infinity. By (a), this will ensure the
sequence with general term $\widehat{\mu_{\psi\left(n\right)}}$ to
converge simply to $\widehat{\mu}.$ Since, by assumption, $\widehat{\mu_{n}}$
pointwise tends to $\varphi$ when $n$ tends to the infinity, it
is the same for every subsequence, and we will have $\widehat{\mu}=\varphi.$

The uniqueness of the weak adherence value $\mu$ then follows from
the injectivity of the Fourier transform, and we will have proved
the weak convergence of the sequence $\left(\mu_{n}\right)_{n\in\mathbb{N}}$
to $\mu.$ 

Thus, it remains to show that the sequence $\left(\mu_{\psi\left(n\right)}\right)_{n\in\mathbb{N}}$
converges narrowly to $\mu.$ Since we already have a weak convergence,
by Proposition $\refpar{pr:narrow_cv_nsc},$ it is enough to prove
that
\[
\lim_{n\to+\infty}\mu_{\psi\left(n\right)}\left(\mathbb{R}^{d}\right)=\mu\left(\mathbb{R}^{d}\right).
\]
By assumption,
\[
\lim_{n\to+\infty}\mu_{\psi\left(n\right)}\left(\mathbb{R}^{d}\right)=\lim_{n\to+\infty}\mu_{n}\left(\mathbb{R}^{d}\right)=\lim_{n\to+\infty}\widehat{\mu_{n}}\left(0\right)=\varphi\left(0\right).
\]

Since 
\[
\mu\left(\mathbb{R}^{d}\right)=\widehat{\mu}\left(0\right),
\]
it remains to show that 
\[
\widehat{\mu}\left(0\right)=\varphi\left(0\right).
\]
 To this end, observe first that for $\epsilon>0,$
\begin{equation}
\lim_{n\to+\infty}\intop_{\left[0,\epsilon\right]^{d}}\widehat{\mu_{n}}\left(t\right)\text{d}t=\intop_{\left[0,\epsilon\right]^{d}}\varphi\left(t\right)\text{d}t.\label{eq:lim_int_transform_mu_n}
\end{equation}

Indeed, since $\widehat{\mu_{n}}$ converges pointwise to $\varphi$
when $n$ tends to infinity, and the functions $\widehat{\mu_{n}}$
are bounded in modulus by $b,$ this follows from the dominated convergence
theorem.

We then use the following lemma.

\end{proof}

\begin{lemma}{}{}

Let $\epsilon>0.$ There exists a function $f_{\epsilon}\in\mathscr{C}_{0}\left(\mathbb{R}^{d}\right)$
such that for every bounded measure $\nu$ on $\mathbb{R}^{d},$
\begin{equation}
\intop_{\left[0,\epsilon\right]^{d}}\widehat{\nu}\left(t\right)\text{d}t=\intop_{\mathbb{R}^{d}}f_{\epsilon}\text{d}\nu.\label{eq:int_transform_bounded_meas_ex_f}
\end{equation}

\end{lemma}

\begin{proof}{Lemma Proof}{}

By the Fubini theorem,
\begin{align*}
\intop_{\left[0,\epsilon\right]^{d}}\widehat{\nu}\left(t\right)\text{d}t & =\intop_{\left[0,\epsilon\right]^{d}}\left[\intop_{\mathbb{R}^{d}}\text{e}^{\text{i}\left\langle x,t\right\rangle }\text{d}\nu\left(x\right)\right]\text{d}t\\
 & =\intop_{\mathbb{R}^{d}}\left[\intop_{\left[0,\epsilon\right]^{d}}\text{e}^{\text{i}\left\langle x,t\right\rangle }\text{d}t\right]\text{d}\nu\left(x\right).
\end{align*}
Again, by the Fubini theorem,
\[
\intop_{\left[0,\epsilon\right]^{d}}\text{e}^{\text{i}\left\langle x,t\right\rangle }\text{d}t=\prod^{d}_{j=1}\left[\intop^{\epsilon}_{0}\text{e}^{\text{i}x_{j}t_{j}}\text{d}t_{j}\right].
\]
We thus obtain $\refpar{eq:int_transform_bounded_meas_ex_f}$ by defining,
for $u\in\mathbb{R},$
\[
g_{\epsilon}\left(u\right)=\intop_{\left[0,\epsilon\right]}\text{e}^{\text{i}ut}\text{d}t=\begin{cases}
\dfrac{\text{e}^{\text{i}\epsilon u}-1}{\text{i}u}, & \text{if }u\neq0,\\
\epsilon, & \text{if }u=0,
\end{cases}
\]
and, for $x\in\mathbb{R}^{d},$ 
\[
f_{\epsilon}\left(x\right)=\prod^{d}_{j=1}g_{\epsilon}\left(x_{j}\right).
\]
It is then clear\footnotemark that $f_{\epsilon}\in\mathscr{C}_{0}\left(\mathbb{R}^{d}\right).$

\end{proof}

\footnotetext{The relation $\refpar{eq:int_transform_bounded_meas_ex_f}$
is a particular case of the relation
\[
\intop\widehat{\mu}\text{d}\nu=\intop\widehat{\nu}\text{d}\mu
\]
that holds for arbitrary bounded measures $\mu,\nu$ on $\mathbb{R}^{d}.$ 

When $\mu$ is a measure of density $h$ with respect to the Lebesgue
measure---$h\in\text{L}^{1}\left(\mathbb{R}^{d}\right)$---we write
$\widehat{h}=\widehat{\mu}$ and say that $\widehat{h}$ is the Fourier
transform of the function $h.$ In this case, 
\[
\intop\widehat{h}\text{d}\nu=\intop\widehat{\nu}\left(t\right)h\left(t\right)\text{d}t.
\]
We obtain $\refpar{pr:narrow_cv_criterion}$ by taking $h=\boldsymbol{1}_{\left[0,\epsilon\right]^{d}}\cdot\widehat{h}=f_{\epsilon}.$
The fact that $\widehat{h}\in\mathscr{C}_{0}\left(\mathbb{R}^{d}\right)$
is a general result: the Riemann-Lebesgue lemma.}

\begin{proof}{End of the Proof of (b)}{}

Since the sequence with general term $\mu_{\psi\left(n\right)}$ converges
weakly to $\mu$ and since $f_{\epsilon}\in\mathscr{C}_{0}\left(\mathbb{R}^{d}\right),$
\[
\lim_{n\to+\infty}\intop f_{\epsilon}\text{d}\mu_{\psi\left(n\right)}=\intop f_{\epsilon}\text{d}\mu.
\]
Hence, by the previous lemma,
\[
\lim_{n\to+\infty}\intop_{\left[0,\epsilon\right]^{d}}\widehat{\mu_{\psi\left(n\right)}}\left(t\right)\text{d}t=\intop_{\left[0,\epsilon\right]^{d}}\widehat{\mu}\left(t\right)\text{d}t.
\]

By $\refpar{eq:lim_int_transform_mu_n}$ applied to the subsequence
$\left(\mu_{\psi\left(n\right)}\right)_{n\in\mathbb{N}},$ 
\[
\dfrac{1}{\epsilon^{d}}\intop_{\left[0,\epsilon\right]^{d}}\widehat{\mu}\left(t\right)\text{d}t=\dfrac{1}{\epsilon^{d}}\intop_{\left[0,\epsilon\right]^{d}}\varphi\left(t\right)\text{d}t.
\]
By continuity of $\widehat{\mu}$ and of $\varphi,$ letting $\epsilon$
tends to 0 in the two terms of the above equality yields
\begin{equation}
\widehat{\mu}\left(0\right)=\varphi\left(0\right).\label{eq:four_transf_in_zero}
\end{equation}

We have thus proved the weak convergence of the sequence $\left(\mu_{n}\right)_{n\in\mathbb{N}}.$
Finally, by $\refpar{eq:four_transf_in_zero}$ and the pointwise convergence
of the sequence $\left(\widehat{\mu}_{n}\right)_{n\in\mathbb{N}}$
to $\varphi,$ we obtain
\[
\lim_{n\to+\infty}\widehat{\mu_{n}}\left(0\right)=\widehat{\mu}\left(0\right).
\]
That is,
\[
\lim_{n\to+\infty}\mu_{n}\left(\mathbb{R}^{d}\right)=\mu\left(\mathbb{R}^{d}\right).
\]
This completes the proof of the \textbf{narrow} convergence of the
sequence $\left(\mu_{n}\right)_{n\in\mathbb{N}}$ to its weak limit
$\mu.$

(c) By Corollary $\ref{co:sc_tight_seq},$ the sequence $\left(\mu_{n}\right)_{n\in\mathbb{N}},$
which converges narrowly, is tight. Let $\epsilon>0.$ Choose a compact
set $K_{\epsilon}$ such that
\[
\sup_{n\in\mathbb{N}}\mu_{n}\left(K^{c}_{\epsilon}\right)\leqslant\dfrac{\epsilon}{3}.
\]
For every $n\in\mathbb{N},$ and every $t,t^{\prime}$ of $\mathbb{R}^{d},$
\begin{align*}
\left|\widehat{\mu_{n}}\left(t\right)-\widehat{\mu_{n}}\left(t^{\prime}\right)\right| & =\left|\intop_{\mathbb{R}^{d}}\left[\text{e}^{\text{i}\left\langle x,t\right\rangle }-\text{e}^{\text{i}\left\langle x,t^{\prime}\right\rangle }\right]\text{d}\mu_{n}\left(x\right)\right|\\
 & \leqslant\intop_{K_{\epsilon}}\left|\text{e}^{\text{i}\left\langle x,t\right\rangle }-\text{e}^{\text{i}\left\langle x,t^{\prime}\right\rangle }\right|\text{d}\mu_{n}\left(x\right)+2\mu_{n}\left(K^{c}_{\epsilon}\right),
\end{align*}
Hence, by the \textbf{finite increment inequality},
\[
\left|\widehat{\mu_{n}}\left(t\right)-\widehat{\mu_{n}}\left(t^{\prime}\right)\right|\leqslant\intop_{K_{\epsilon}}\left|\left\langle x,t-t^{\prime}\right\rangle \right|\text{d}\mu_{n}\left(x\right)+2\mu_{n}\left(K^{c}_{\epsilon}\right).
\]
Therefore, for every $t,t^{\prime}$ of $\mathbb{R}^{d},$
\[
\sup_{n\in\mathbb{N}}\left|\widehat{\mu_{n}}\left(t\right)-\widehat{\mu_{n}}\left(t^{\prime}\right)\right|\leqslant\left\Vert t-t^{\prime}\right\Vert \left[b\sup_{x\in K_{\epsilon}}\left\Vert x\right\Vert \right]+\dfrac{2\epsilon}{3},
\]
and, for every $t,t^{\prime}$ of $\mathbb{R}^{d},$ such that $\left\Vert t-t^{\prime}\right\Vert \leqslant\dfrac{\epsilon}{3b\sup_{x\in K_{\epsilon}}\left\Vert x\right\Vert },$
\[
\sup_{n\in\mathbb{N}}\left|\widehat{\mu_{n}}\left(t\right)-\widehat{\mu_{n}}\left(t^{\prime}\right)\right|\leqslant\epsilon.
\]

Thus, the sequence of functions $\widehat{\mu_{n}}$ is equicontinuous---uniformly
in $t.$ Since it converges pointwise, it converges uniformly on every
compact subset of $\mathbb{R}^{d}.$

\end{proof}

\section{Convergence in Law}\label{sec:Convergence-in-Law}

\textbf{For the sake of simplicity, we restrict our attention to random
variables taking values in $\mathbb{R}^{d}.$ All the statements below
remain valid when the random variables take values in a metric space
$E$ that is locally compact and countable at infinity.}

\begin{definition}{Convergence In Law To A Random Variable}{}

For every $n\in\mathbb{N},$ let $X_{n}$ be a random variable defined
on a probabilized space $\left(\Omega_{n},\mathscr{A}_{n},P^{n}\right),$
taking values in $\mathbb{R}^{d},$ and let $X$ be a \textbf{random
variable} defined on a probabilized space $\left(\Omega,\mathscr{A},P\right),$
also taking values in $\mathbb{R}^{d}.$ 

The sequence of random variables $\left(X_{n}\right)_{n\in\mathbb{N}}$
is said to \textbf{converge in law\index{convergence in law}} to
$X$ if the sequence $\left(P^{n}_{X_{n}}\right)_{n\in\mathbb{N}}$
of the laws of $X_{n}$ converges narrowly to the law $P_{X}$ of
$X.$ This is denoted by
\[
X_{n}\stackrel[n\to+\infty]{\mathscr{L}}{\longrightarrow}X.
\]

\end{definition}

\begin{remark}{}{}

This notion of convergence does not concern the random variables as
functions, but rather their laws. In particular, it provides a framework
for \textbf{approximating the laws\index{approximation of laws}}
of random variables. It is important to note that the random variables
need not necessarily to be defined on the same probabilized space.
In practice, there is not always a natural limiting random variable
$X.$ This leads to a second definition of convergence in law for
a sequence of random variables, which is often used alongside the
previous definition. 

\end{remark}

\begin{definition}{Convergence in Law to a Probability}{}

For each $n\in\mathbb{N},$ let $X_{n}$ be a random variable defined
on a probabilized space $\left(\Omega_{n},\mathscr{A}_{n},P^{n}\right),$
taking values in $\mathbb{R}^{d},$ and let $\mu$ be a \textbf{probability}
on $\mathbb{R}^{d}.$ The sequence of random variables $\left(X_{n}\right)_{n\in\mathbb{N}}$
is said to \textbf{converge in law\index{convergence in law} to $\mu$}
if the sequence $\left(P^{n}_{X_{n}}\right)_{n\in\mathbb{N}}$ of
the laws of $X_{n}$ converges narrowly to the law $\mu.$ This is
denoted
\[
X_{n}\stackrel[n\to+\infty]{\mathscr{L}}{\longrightarrow}\mu.
\]

\end{definition}

It is worth to note that in this case, the mathematical object on
the two sides of the arrow are of \textbf{a different nature}.

\begin{example}{}{}

We will prove later that if, for every $n\in\mathbb{N},$ $X_{n}$
is a real-valued random variable following the binomial law $\mathscr{B}\left(n,\dfrac{\lambda}{n}\right),$
where $\lambda>0,$ then the sequence $\left(X_{n}\right)_{n\in\mathbb{N}}$
converges in law to the Poisson law $\mathscr{P}\left(\lambda\right).$

\end{example}

Conceptually, there is therefore nothing new compared to the narrow
convergence of a sequence of probability measures, and the criteria
for the convergence in law are precisely those for the narrow convergence
of a sequence of \textbf{probability}. Nevertheless, we now state
the Levy theorem in the framework of convergence in law.

\begin{theorem}{Lévy Theorem. Convergence in Law}{levy_th_convergence_law}

Let, for every $n\in\mathbb{N},$ $X_{n}$ be a random variable defined
on the probabilized space $\left(\Omega_{n},\mathscr{A}_{n},P^{n}\right),$
taking values in $\mathbb{R}^{d},$ with characteristic function $\varphi_{X_{n}}.$

(a) If the sequence of random variables $\left(X_{n}\right)_{n\in\mathbb{N}}$
converges in law to $X,$ where $X$ is the random variable defined
on the probabilized space $\left(\Omega,\mathscr{A},P\right),$ taking
values in $\mathbb{R}^{d},$ then the sequence $\left(\varphi_{X_{n}}\right)_{n\in\mathbb{N}}$
of characteristic functions converges pointwise---and even uniformly
on every compact of $\mathbb{R}^{d}$---to the characteristic function
$\varphi_{X}$ of $X.$

(b) Conversely, if the sequence $\left(\varphi_{X_{n}}\right)_{n\in\mathbb{N}}$
of characteristic functions converges pointwise to a function $\varphi$
that is \textbf{continuous at 0}, then $\varphi$ is the Fourier transform
of a \textbf{probability measure} $\mu$ on $\mathbb{R}^{d},$ and
the sequence of random variables $\left(X_{n}\right)_{n\in\mathbb{N}}$
\textbf{converges in law to $\mu.$}

Moreover, there exists a non-unique random variable $X,$ defined
on a probabilized space $\left(\Omega,\mathscr{A},P\right),$ and
taking values in $\mathbb{R}^{d},$ such that the sequence of random
variables $\left(X_{n}\right)_{n\in\mathbb{N}}$ \textbf{converges
in law to $X.$} 

\end{theorem}

\begin{proof}{}{}This is simply a reformulation of the Lévy theorem
for the narrow convergence of bounded measures, recalling that $\varphi_{X_{n}}$
is, by definition, the Fourier transform of the law of $X_{n}.$ 

Only the last point in the converse requires a clarification. By the
Lévy theorem---Theorem $\ref{th:levy_th_measures}$---the sequence
$\left(X_{n}\right)_{n\in\mathbb{N}}$ converges in law to the probability
$\mu$ such that $\widehat{\mu}=\varphi$---$\mu$ is then a probability,
since 
\[
\lim_{n\to+\infty}\varphi_{X_{n}}\left(0\right)\equiv1=\varphi\left(0\right)=\widehat{\mu}\left(0\right).
\]
We then consider the identity function $X$ from $\mathbb{R}^{d}$
onto itself. This is a random variable defined on the probabilized
space $\left(\mathbb{R}^{d},\mathscr{B}_{\mathbb{R}^{d}},\mu\right)$
taking values in $\mathbb{R}^{d},$ with law $\mu,$ and such that
the sequence of random variables $\left(X_{n}\right)_{n\in\mathbb{N}}$
converges in law to $X.$

\end{proof}

\begin{example}{}{}

Let, for each $n\in\mathbb{N},$ $X_{n}$ be a random variable defined
on the probabilized space $\left(\Omega_{n},\mathscr{A}_{n},P^{n}\right),$
taking values on $\mathbb{R}^{d},$ and let $X$ be a random variable
defined on the probabilized space $\left(\Omega,\mathscr{A},P\right),$
taking values in $\mathbb{R}^{d}.$ Show the equivalence
\[
X_{n}\stackrel[n\to+\infty]{\mathscr{L}}{\longrightarrow}X\,\,\,\,\Longleftrightarrow\,\,\,\,\forall t\in\mathbb{R}^{d},\,\,\,\,\left\langle X_{n},t\right\rangle \stackrel[n\to+\infty]{\mathscr{L}}{\longrightarrow}\left\langle X,t\right\rangle .
\]

\end{example}

\begin{solutionexample}{}{}

Indeed, for every $t\in\mathbb{R}^{d}$ and for every real number
$\alpha,$
\[
\varphi_{X_{n}}\left(\alpha t\right)=\varphi_{\left\langle X_{n},t\right\rangle }\left(\alpha\right)\,\,\,\,\text{and}\,\,\,\,\varphi_{X}\left(\alpha t\right)=\varphi_{\left\langle X,t\right\rangle }\left(\alpha\right).
\]
By Theorem $\ref{th:levy_th_convergence_law},$ we conclude.

\end{solutionexample}

\begin{remark}{}{}

With the same notation as in Theorem $\ref{th:levy_th_convergence_law},$
it is clear that if $f$ is a continuous function from $\mathbb{R}^{d}$
to $\mathbb{R}^{k}$ and if the sequence $\left(X_{n}\right)_{n\in\mathbb{N}}$
converges in law to $X,$ then the sequence $\left(f\left(X_{n}\right)\right)_{n\in\mathbb{N}}$
converges in law to $f\left(X\right).$ Indeed, for every function
$g\in\mathscr{C}_{b}\left(\mathbb{R}^{k}\right),$ by the transfer
theorem
\[
\intop_{\mathbb{R}^{k}}g\text{d}P^{n}_{f\left(X_{n}\right)}=\intop_{\mathbb{R}^{d}}f\circ g\text{d}P^{n}_{X_{n}}\,\,\,\,\text{and}\,\,\,\,\intop_{\mathbb{R}^{k}}g\text{d}P_{X}=\intop_{\mathbb{R}^{d}}f\circ g\text{d}P_{X}.
\]

\end{remark}

The following proposition extends the class of functions for which
the property stated in the previous remark holds.

\begin{proposition}{}{}

Let, for each $n\in\mathbb{N},$ $X_{n}$ be a random variable defined
on the probabilized space $\left(\Omega_{n},\mathscr{A}_{n},P^{n}\right),$
taking values in $\mathbb{R}^{d},$ and let $X$ be a random variable
defined on a probabilized space $\left(\Omega,\mathscr{A},P\right),$
taking values also in $\mathbb{R}^{d}.$

Let $f$ be a Borel function from $\mathbb{R}^{d}$ to $\mathbb{R}^{k}$
such that $f$ is $P_{X}-$almost surely continuous.

If the sequence $\left(X_{n}\right)_{n\in\mathbb{N}}$ converges in
law to $X,$ then the sequence $\left(f\left(X_{n}\right)\right)_{n\in\mathbb{N}}$
converges in law to $f\left(X\right).$

\end{proposition}

\begin{proof}{}{}

Let $C\in\mathscr{B}_{\mathbb{R}^{d}}$ be such that $P_{X}\left(C\right)=1,$
and such that $f$ is continuous on $C.$ Let $F$ be an arbitrary
closed subset on $\mathbb{R}^{k}.$ Then, for every $n\in\mathbb{N},$
\[
P^{n}_{f\left(X_{n}\right)}\left(F\right)=P^{n}_{X_{n}}\left(f^{-1}\left(F\right)\right)\leqslant P^{n}_{X_{n}}\left(\overline{f^{-1}\left(F\right)}\right),
\]
and thus, by Proposition $\ref{pr:narrow_cv_criterion},$
\begin{equation}
\limsup_{n\to+\infty}P^{n}_{f\left(X_{n}\right)}\left(F\right)\leqslant\limsup_{n\to+\infty}P^{n}_{X_{n}}\left(\overline{f^{-1}\left(F\right)}\right)\leqslant P_{X}\left(\overline{f^{-1}\left(F\right)}\right).\label{eq:upperbound_limsup_n_PfX_nF}
\end{equation}
Moreover, we have the inclusions
\[
f^{-1}\left(F\right)\subset\overline{f^{-1}\left(F\right)}\subset C^{c}\cup f^{-1}\left(F\right).
\]
Indeed, let $x\in\overline{f^{-1}\left(F\right)}.$ 
\begin{itemize}
\item If $x\in C^{c},$ then $x\in C^{c}\cup f^{-1}\left(F\right).$
\item If $x\in C,$ since $x\in\overline{f^{-1}\left(F\right)},$ there
exists a sequence $\left(x_{n}\right)_{n\in\mathbb{N}}$ of points
in $f^{-1}\left(F\right)$---i.e. such that $f\left(x_{n}\right)\in F$---which
converges to $x,$ a point of continuity of $f.$ Then 
\[
\lim_{n\to+\infty}f\left(x_{n}\right)=f\left(x\right),
\]
and since $F$ is closed, $f\left(x\right)\in F.$ Hence, $x\in f^{-1}\left(F\right),$
and so $x\in C^{c}\cup f^{-1}\left(F\right).$ 
\end{itemize}
Since $P_{X}\left(C^{c}\right)=0,$ it follows that
\[
P_{X}\left(\overline{f^{-1}\left(F\right)}\right)=P_{X}\left(f^{-1}\left(F\right)\right)=P_{f\left(X\right)}\left(F\right),
\]
and, by substituting into the inequality $\refpar{eq:upperbound_limsup_n_PfX_nF},$
\[
\limsup_{n\to+\infty}P^{n}_{f\left(X_{n}\right)}\left(F\right)\leqslant P_{f\left(X\right)}\left(F\right),
\]
which shows that the sequence $\left(f\left(X_{n}\right)\right)_{n\in\mathbb{N}}$
converges in law to $f\left(X\right).$

\end{proof}

\begin{example}{}{}

Let, for each $n\in\mathbb{N},$ $\left(X_{n},Y_{n}\right)$ be a
couple of random variables on the probabilized space $\left(\Omega,\mathscr{A},P\right),$
and $\left(X,Y\right)$ taking values in $\mathbb{R}^{2},$ such that
$\left(X_{n},Y_{n}\right)\stackrel[n\to+\infty]{\mathscr{L}}{\longrightarrow}\left(X,Y\right).$

Then, for instance,
\[
X_{n}+Y_{n}\stackrel[n\to+\infty]{\mathscr{L}}{\longrightarrow}X+Y\,\,\,\,\text{and}\,\,\,\,X_{n}Y_{n}\stackrel[n\to+\infty]{\mathscr{L}}{\longrightarrow}XY.
\]

Assume, for simplicity, that $Y_{n}\neq0$ everywhere, for every $n\in\mathbb{N}.$
If $P_{Y}\left(\left\{ 0\right\} \right)=0$---that is $P_{\left(X,Y\right)}\left(\mathbb{R}\times\left\{ 0\right\} =0\right)$---then
\[
\dfrac{X_{n}}{Y_{n}}\stackrel[n\to+\infty]{\mathscr{L}}{\longrightarrow}\dfrac{X}{Y}.
\]

\end{example}

\textbf{We now compare, when meaningful, the convergence in law and
the convergence in probability for a sequence of random variables.}

\begin{proposition}{Convergence in Probability Implies Convergence in Law}{}

If a sequence $\left(X_{n}\right)_{n\in\mathbb{N}}$ of random variables,
defined on the \textbf{same} probabilized space $\left(\Omega,\mathscr{A},P\right)$
and taking values in $\mathbb{R}^{d},$ converges in probability to
a random variable $X$---defined on $\left(\Omega,\mathscr{A},P\right)$
and taking values in $\mathbb{R}^{d}$---then it also converges in
law to $X.$

\end{proposition}

\begin{proof}{}{}

Let $f\in\mathscr{C}_{b}\left(\mathbb{R}^{d}\right).$ For every $\epsilon>0,$
by the transfer theorem,
\begin{align*}
\left|\intop_{\mathbb{R}^{d}}f\text{d}P_{X_{n}}-\intop_{\mathbb{R}^{d}}f\text{d}P_{X}\right| & =\left|\mathbb{E}\left(f\left(X_{n}\right)\right)-\mathbb{E}\left(f\left(X\right)\right)\right|\\
 & \leqslant\epsilon+2\left\Vert f\right\Vert P\left(\left|f\left(X_{n}\right)-f\left(X\right)\right|>\epsilon\right).
\end{align*}
Since $f$ is continous, the sequence $\left(f\left(X_{n}\right)\right)_{n\in\mathbb{N}}$
converges in probability to $f\left(X\right),$ and it follows that
\[
0\leqslant\limsup_{n\to+\infty}\left|\intop_{\mathbb{R}^{d}}f\text{d}P_{X_{n}}-\intop_{\mathbb{R}^{d}}f\text{d}P_{X}\right|\leqslant\epsilon,
\]
which, by the arbitrariness of $\epsilon,$ proves that
\[
\lim_{n\to+\infty}\int_{\mathbb{R}^{d}}f\text{d}P_{X_{n}}=\int_{\mathbb{R}^{d}}f\text{d}P_{X}.
\]

\end{proof}

\begin{remark}{The Converse is False!}{}

The converse is false and, as shown by the following counter-example,
it even does not hold for a stationary sequence.

\end{remark}

\begin{counterexample}{}{}

Consider on the probabilized space \preds a Bernoulli random variable
$X$ with parameter $\dfrac{1}{2}$ and set, for every $n\in\mathbb{N},$
$X_{n}=X.$

Prove that the sequence converges in law without converging in probability
to $X.$

\end{counterexample}

\begin{solutioncounterexample}{}{}

It is clear that $X_{n}\stackrel[n\to+\infty]{\mathscr{L}}{\longrightarrow}X.$

The random variable $Y=1-X$ also follows a Bernoulli law with parameter
$\dfrac{1}{2}.$ Hence $X_{n}\stackrel[n\to+\infty]{\mathscr{L}}{\longrightarrow}Y.$
Nonetheless, since $P-$almost surely
\[
\left|X_{n}-Y\right|=\left|2X-1\right|=1,
\]
for every $\epsilon\in\left]0,1\right[,$ 
\[
P\left(\left|X_{n}-Y\right|>\epsilon\right)=1
\]
and therefore the sequence $\left(X_{n}\right)_{n\in\mathbb{N}}$
does not converge in probability to $Y.$

\end{solutioncounterexample}

Nonetheless, we have a partial converse.

\begin{proposition}{Partial Converse For a $P-$almost surely Constant Random Variable}{}

If a sequence $\left(X_{n}\right)_{n\in\mathbb{N}}$ of random variables,
defined on the same probabilized space $\left(\Omega,\mathscr{A},P\right)$
and taking values in $\mathbb{R}^{d},$ converges $P-$almost surely
in law to a constant random variable $a,$ then it also converges
in probability to $a.$

\end{proposition}

\begin{proof}{}{}

For every $\epsilon>0,$
\[
\delta_{a}\left(\partial B_{f}\left(a,\epsilon\right)\right)=0,
\]
where $B_{f}\left(a,\epsilon\right)$ is the closed ball with center
$a$ and radius $\epsilon.$

Hence, by Proposition $\ref{pr:narrow_cv_criterion},$
\[
\lim_{n\to+\infty}P_{X_{n}}\left(B_{f}\left(a,\epsilon\right)\right)=\delta_{a}\left(B_{f}\left(a,\epsilon\right)\right)=1.
\]
It follows that
\[
\lim_{n\to+\infty}P\left(\left\Vert X_{n}-a\right\Vert >\epsilon\right)=0.
\]

\end{proof}

The following Scheffé lemma states a sufficient condition of convergence
in law in the case where the random variables admit a density.

\begin{lemma}{Scheffé Lemma}{scheffe_lem}

Let, for each $n\in\mathbb{N},$ $X_{n}$ be a random variable defined
on the probabilized space $\left(\Omega_{n},\mathscr{A}_{n},P^{n}\right),$
taking values in $\mathbb{R}^{d}$ and admitting a density $f_{X_{n}}.$
If the sequence $\left(f_{X_{n}}\right)_{n\in\mathbb{N}}$ converges
$\lambda_{d}-$almost everywhere to a function $f$ such that 
\[
\intop_{\mathbb{R}^{d}}f\text{d}\lambda_{d}=1,
\]
then the sequence $\left(X_{n}\right)_{n\in\mathbb{N}}$ converges
in law to the law $f\cdot\lambda_{d}.$

Moreover,
\[
\lim_{n\to+\infty}\sup_{A\in\mathscr{B}_{\mathbb{R}^{d}}}\left|P_{X_{n}}\left(A\right)-\intop_{A}f\text{d}\lambda_{d}\right|=0.
\]

\end{lemma}

\begin{proof}{}{}

For every $A\in\mathscr{B}_{\mathbb{R}^{d}},$
\[
\left|P_{X_{n}}\left(A\right)-\intop_{A}f\text{d}\lambda_{d}\right|=\left|\intop_{A}\left(f_{X_{n}}-f\right)\text{d}\lambda_{d}\right|\leqslant\left|\intop_{\mathbb{R}^{d}}\left(f_{X_{n}}-f\right)\text{d}\lambda_{d}\right|.
\]
Thus,
\begin{equation}
\sup_{A\in\mathscr{B}_{\mathbb{R}^{d}}}\left|P_{X_{n}}\left(A\right)-\intop_{A}f\text{d}\lambda_{d}\right|\leqslant\intop_{\mathbb{R}^{d}}\left|f_{X_{n}}-f\right|\text{d}\lambda_{d}.\label{eq:upperbound_sup_diff_P_X_n_andf}
\end{equation}

Recall the useful identity, for every $a,\,b\in\mathbb{R},$
\[
\left|a-b\right|=a+b-2\min\left(a,b\right).
\]
Hence,
\[
\intop_{\mathbb{R}^{d}}\left|f_{X_{n}}-f\right|\text{d}\lambda_{d}=\intop_{\mathbb{R}^{d}}f_{X_{n}}\text{d}\lambda_{d}+\intop_{\mathbb{R}^{d}}f\text{d}\lambda_{d}-2\intop_{\mathbb{R}^{d}}\min\left(f_{X_{n}},f\right)\text{d}\lambda_{d}.
\]
Taking into account that $f_{X_{n}}$ and $f$ are probability densities,
we obtain
\[
\intop_{\mathbb{R}^{d}}\left|f_{X_{n}}-f\right|\text{d}\lambda_{d}=2-2\intop_{\mathbb{R}^{d}}\min\left(f_{X_{n}},f\right)\text{d}\lambda_{d}.
\]
Since, for every $n\in\mathbb{N},$
\[
0\leqslant\min\left(f_{X_{n}},f\right)\leqslant f,
\]
and since $\lambda_{d}-$almost everywhere
\[
\lim_{n\to+\infty}\min\left(f_{X_{n}},f\right)=f,
\]
it follows from the dominated convergence theorem that
\[
\lim_{n\to+\infty}\intop_{\mathbb{R}^{d}}\left|f_{X_{n}}-f\right|\text{d}\lambda_{d}=2-2\lim_{n\to+\infty}\intop_{\mathbb{R}^{d}}\min\left(f_{X_{n}},f\right)\text{d}\lambda_{d}=2-2\intop_{\mathbb{R}^{d}}f\text{d}\lambda_{d}=0.
\]
Taking into account $\refpar{eq:upperbound_sup_diff_P_X_n_andf},$
this yields the result.

\end{proof}

The following proposition gives a \textbf{criterion for convergence
in law for a sequence of discrete random variables} taking values
in $\mathbb{Z}.$

\begin{proposition}{Convergence in Law Criterion For Discrete Random Variables Sequence}{}

Let $X_{n},\,n\in\mathbb{N}$ and $X$ be random variables defined
on a probabilized space \preds, taking values in $\mathbb{Z}.$ The
following equivalence holds
\[
X_{n}\stackrel[n\to+\infty]{\mathscr{L}}{\longrightarrow}X\,\,\,\,\Longleftrightarrow\,\,\,\,\forall r\in\mathbb{Z},\,\lim_{n\to+\infty}P\left(X_{n}=r\right)=P\left(X=r\right).
\]

\end{proposition}

\begin{proof}{}{}

Assume that the sequence $\left(X_{n}\right)_{n\in\mathbb{N}}$ converges
in law to $X.$ Fix $r.$ Choose $f\in\mathscr{C}_{\mathscr{K}}\left(\mathbb{R}\right)$
with support contained in the interval $\left]r-\dfrac{1}{2},r+\dfrac{1}{2}\right[$
and such that $f\left(r\right)\neq0.$

Since
\[
\intop_{\mathbb{R}}f\text{d}P_{X_{n}}=f\left(r\right)P\left(X_{n}=r\right)\,\,\,\,\text{and}\,\,\,\,\intop_{\mathbb{R}}f\text{d}P_{X}=f\left(r\right)P\left(X=r\right),
\]
and since
\[
\lim_{n\to+\infty}\intop_{\mathbb{R}}f\text{d}P_{X_{n}}=\intop_{\mathbb{R}}f\text{d}P_{X},
\]
it follows that
\[
\lim_{n\to+\infty}P\left(X_{n}=r\right)=P\left(X=r\right).
\]
Conversely, let $f\in\mathscr{C}_{\mathscr{K}}\left(\mathbb{R}\right)$
have compact support $K.$ Then
\[
\intop_{\mathbb{R}}f\text{d}P_{X_{n}}=\sum_{r\in K}f\left(r\right)P\left(X_{n}=r\right),
\]
where the sum has only finitely many terms. 

If, for every $r\in\mathbb{Z},$
\[
\lim_{n\to+\infty}P\left(X_{n}=r\right)=P\left(X=r\right),
\]
it follows that
\[
\lim_{n\to+\infty}\intop_{\mathbb{R}}f\text{d}P_{X_{n}}=\intop_{\mathbb{R}}f\text{d}P_{X},
\]
which proves the narrow convergence of the sequence of probabilities
$P_{X_{n}}$ to the \textbf{probability} $P_{X}.$

\end{proof}

Historically, the convergence in law was defined in terms of the convergence
of cumulative distribution function sequences. However, as the following
proposition shows, this definition is not entirely straightforward.

\begin{proposition}{Convergence in Law and Cumulative Distribution Function}{cv_law_dist_function}

Let, for each $n\in\mathbb{N},$ $X_{n}$ be a random variable defined
on the probabilized space $\left(\Omega_{n},\mathscr{A}_{n},P^{n}\right),$
with a cumulative distribution function $F_{X}$ and let $X$ be a
real-valued random variable defined on a probabilized space $\left(\Omega,\mathscr{A},P\right),$
with cumulative distribution function $F_{X}.$ 

The sequence $\left(X_{n}\right)_{n\in\mathbb{N}}$ converges in law
to $X$ if and only if the sequence $\left(F_{X_{n}}\left(x\right)\right)_{n\in\mathbb{N}}$
converges to $F_{X}\left(x\right)$ \textbf{in every point $x$ where}
$F_{X}$ \textbf{is continuous}.

\end{proposition}

\begin{proof}{}{}

Suppose that the sequence $\left(X_{n}\right)_{n\in\mathbb{N}}$ converges
in law to $X.$ Let $x$ be a point of continuity of $F_{X}.$ Since
$\partial\left(\left]-\infty,x\right]\right)=\left\{ x\right\} $
and that 
\[
P_{X}\left(\left\{ x\right\} \right)=F_{X}\left(x\right)-F_{X}\left(x-0\right)=0,
\]
the half-line $\left]-\infty,x\right]$ is a $P_{X}-$continuity set.
Hence, by Proposition $\ref{pr:narrow_cv_criterion},$
\[
\lim_{n\to+\infty}P_{X_{n}}\left(\left]-\infty,x\right]\right)=P_{X}\left(\left]-\infty,x\right]\right)=F_{X}\left(x\right),
\]
which proves the necessary condition.

Conversely, suppose that the sequence $\left(F_{X_{n}}\left(x\right)\right)_{n\in\mathbb{N}}$
converges to $F_{X}\left(x\right)$ at every point $x$ where $F_{X}$
is continuous. Let $f\in\mathscr{C}_{0}\left(\mathbb{R}\right)$ and
$\epsilon>0.$ Since the set of discontinuity points of $F_{X}$ is
countable---possibly empty---there exists a step function of the
form
\[
g=\sum^{k}_{j=1}\alpha_{j}\boldsymbol{1}_{\left]a_{j},b_{j}\right]},
\]
with $a_{j}<b_{j}\leqslant a_{j+1}<b_{j+1},$ and where the $a_{j}$
and $b_{j}$ are continuity points of $F_{X},$ such that $\left\Vert f-g\right\Vert \leqslant\epsilon.$ 

Then, by hypothesis,
\[
\intop_{\mathbb{R}}g\text{d}P_{X_{n}}=\sum^{k}_{j=1}\alpha_{j}\left(F_{X_{n}}\left(b_{j}\right)-F_{X_{n}}\left(a_{j}\right)\right)\underset{n\to+\infty}{\longrightarrow}\sum^{k}_{j=1}\alpha_{j}\left(F_{X}\left(b_{j}\right)-F_{X}\left(a_{j}\right)\right).
\]
That is
\[
\lim_{n\to+\infty}\intop_{\mathbb{R}}g\text{d}P_{X_{n}}=\intop_{\mathbb{R}}g\text{d}P_{X}.
\]
By the triangle inequality,
\[
\left|\intop_{\mathbb{R}}f\text{d}P_{X_{n}}-\intop_{\mathbb{R}}f\text{d}P_{X}\right|\leqslant\left|\intop_{\mathbb{R}}\left(f-g\right)\text{d}P_{X_{n}}\right|+\left|\intop_{\mathbb{R}}g\text{d}P_{X_{n}}-\intop_{\mathbb{R}}g\text{d}P_{X}\right|+\left|\intop_{\mathbb{R}}\left(g-f\right)\text{d}P_{X}\right|.
\]
Hence, 
\[
\left|\intop_{\mathbb{R}}f\text{d}P_{X_{n}}-\intop_{\mathbb{R}}f\text{d}P_{X}\right|\leqslant2\left\Vert f-g\right\Vert +\left|\intop_{\mathbb{R}}g\text{d}P_{X_{n}}-\intop_{\mathbb{R}}g\text{d}P_{X}\right|.
\]
Thus,
\[
0\leqslant\limsup_{n\to+\infty}\left|\intop_{\mathbb{R}}f\text{d}P_{X_{n}}-\intop_{\mathbb{R}}f\text{d}P_{X}\right|\leqslant2\epsilon,
\]
and, since $\epsilon$ is arbitrary, we conclude that 
\[
\lim_{n\to+\infty}\intop_{\mathbb{R}}f\text{d}P_{X_{n}}=\intop_{\mathbb{R}}f\text{d}P_{X}.
\]
Hence, the sequence of probabilities $P_{X_{n}}$ converges narrowly
to $P_{X}.$

\end{proof}

\begin{remark}{}{}

As the following example illustrates, one cannot expect the simple
convergence---pointwise everywhere---of the sequence of cumulative
distribution functions $F_{X_{n}}.$ If $X_{n}=\dfrac{1}{n}$ and
if $X=0,$ then $X_{n}\stackrel[n\to+\infty]{\mathscr{L}}{\longrightarrow}X.$
However, for every $n\in\mathbb{N},$ we have $F_{X_{n}}\left(0\right)=0$
whereas $F_{X}\left(0\right)=1.$ A similar statements holds for random
variables taking values in $\mathbb{R}^{d},$ although it is of limited
practical use.

\end{remark}

\begin{example}{The Converse of the Scheffé Lemma is False}{}

Let, for each $n\in\mathbb{N}^{\ast},$ $X_{n}$ be a real-valued
random variable defined on a probabilized space \preds, admitting
a density $f_{X_{n}}$ defined, for every real number $x,$ by
\[
f_{X_{n}}\left(x\right)=\boldsymbol{1}_{\left]0,1\right]}\left(x\right)\left(1-\cos\left(2\pi nx\right)\right).
\]

The sequence $\left(f_{X_{n}}\right)_{n\in\mathbb{N}^{\ast}}$ does
not converge $\lambda-$almost everywhere---it diverges at every
point of $\left]0,1\right[,$ and converges elsewhere. Nevertheless,
the sequence $\left(X_{n}\right)_{n\in\mathbb{N}^{\ast}}$ converges
in law to the law $\boldsymbol{1}_{\left]0,1\right]}\cdot\lambda.$ 

Indeed, the cumulative distribution function of $X_{n}$ is given
by
\[
F_{X_{n}}\left(x\right)=\begin{cases}
0, & \text{if }x<0,\\
x-\dfrac{\sin\left(2\pi nx\right)}{2\pi n}, & \text{if }0<x\leqslant1,\\
1, & \text{if\,}x>1,
\end{cases}
\]
so that
\[
\lim_{n\to+\infty}F_{X_{n}}\left(x\right)=\begin{cases}
0, & \text{if }x<0,\\
x, & \text{if }0<x\leqslant1,\\
1, & \text{if\,}x>1.
\end{cases}
\]
Hence, by Proposition $\ref{pr:cv_law_dist_function},$ this proves
that the sequence $\left(X_{n}\right)_{n\in\mathbb{N^{\ast}}}$ converges
in law to the probability $\boldsymbol{1}_{\left]0,1\right]}\cdot\lambda_{d},$
that is, the uniform law on $\left]0,1\right].$

\end{example}

We now state two convergence-in-law theorems related to the Poisson
law. For their proofs, we use the classical lemma below.

\begin{lemma}{Complex Exponential as Limit of $\left(1+\dfrac{z}{n}\right)^{n}$}{cpl_exp_lim_power}

For every complex number $z$ and for every $n\in\mathbb{N}^{\ast},$
\begin{equation}
\left|\text{e}^{z}-\left(1+\dfrac{z}{n}\right)^{n}\right|\leqslant\text{e}^{\left|z\right|}-\left(1+\dfrac{\left|z\right|}{n}\right)^{n}.\label{eq:ineq_exp_z_minus_power_n_1pluszovern}
\end{equation}
It follows that, for every $z\in\mathbb{C},$ the sequence with general
term $\left(1+\dfrac{z}{n}\right)^{n}$ converges and that
\[
\lim_{n\to+\infty}\left(1+\dfrac{z}{n}\right)^{n}=\text{e}^{z}.
\]

\end{lemma}

\begin{proof}{}{}

The binom formula yields, for every $z\in\mathbb{C},$
\[
\text{e}^{z}-\left(1+\dfrac{z}{n}\right)^{n}=\sum^{+\infty}_{j=0}\dfrac{z^{j}}{j!}-\sum^{n}_{j=0}\binom{n}{j}\dfrac{z^{j}}{n^{j}}.
\]
Since 
\[
\dfrac{\binom{n}{j}}{n^{j}}=\dfrac{1}{j!}\prod^{j-1}_{k=0}\left(1-\dfrac{k}{n}\right),
\]
it follows that
\begin{equation}
\text{e}^{z}-\left(1+\dfrac{z}{n}\right)^{n}=\sum^{+\infty}_{j=n+1}\dfrac{z^{j}}{j!}+\sum^{n}_{j=0}\dfrac{z^{j}}{j!}\left[1-\prod^{j-1}_{k=0}\left(1-\dfrac{k}{n}\right)\right].\label{eq:diff_e_z_with_sum}
\end{equation}

Since $1-\prod^{j-1}_{k=0}\left(1-\dfrac{k}{n}\right)\geqslant0,$
we obtain
\[
\left|\text{e}^{z}-\left(1+\dfrac{z}{n}\right)^{n}\right|\leqslant\sum^{+\infty}_{j=n+1}\dfrac{\left|z\right|^{j}}{j!}+\sum^{n}_{j=0}\dfrac{\left|z\right|^{j}}{j!}\left[1-\prod^{j-1}_{k=0}\left(1-\dfrac{k}{n}\right)\right],
\]
which yields $\refpar{eq:ineq_exp_z_minus_power_n_1pluszovern},$
by using $\refpar{eq:diff_e_z_with_sum}$ with $\left|z\right|.$

Finally, since $\ln\left(1+\dfrac{\left|z\right|}{n}\right)^{n}=\left|z\right|+o\left(1\right),$
\[
\lim_{n\to+\infty}\left(1+\dfrac{\left|z\right|}{n}\right)^{n}=\text{e}^{\left|z\right|},
\]
which implies, by $\refpar{eq:ineq_exp_z_minus_power_n_1pluszovern},$
\[
\lim_{n\to+\infty}\left(1+\dfrac{z}{n}\right)^{n}=\text{e}^{z}.
\]

\end{proof}

We now prove the \textbf{Poisson theorem}---an elementary version
was already stated in Part \ref{part:Introduction-to-Probability}---using
the Lévy theorem.

\begin{theorem}{Poisson Theorem}{Poisson_theorem_first}

Let, for each $n\in\mathbb{N}^{\ast},$ $X_{n}$ be a random variable
with binomial law $\mathscr{B}\left(n,p_{n}\right).$ Suppose that
$\lim_{n\to+\infty}p_{n}=\lambda,$ where $\lambda>0.$

Then, the sequence $\left(X_{n}\right)_{n\in\mathbb{N}^{\ast}}$ converges
in law to the Poisson law $\mathscr{P}\left(\lambda\right).$

\end{theorem}

\begin{proof}{}{}

For every $n\in\mathbb{N}^{\ast},$ the characteristic function of
$X_{n}$ is given by
\[
\forall t\in\mathbb{R},\,\,\,\,\varphi_{X_{n}}\left(t\right)=\left[p_{n}\text{e}^{\text{i}t}+\left(1-p_{n}\right)\right]^{n}=\left[1+p_{n}\left(\text{e}^{\text{i}t}-1\right)\right]^{n}.
\]
By $\refpar{eq:ineq_exp_z_minus_power_n_1pluszovern},$ it follows
that, for every $z\in\mathbb{C},$
\[
\left|\text{e}^{np_{n}z}-\left(1+p_{n}z\right)^{n}\right|\leqslant\text{e}^{np_{n}\left|z\right|}-\left(1+p_{n}\left|z\right|\right)^{n}.
\]
Since, by hypothesis,
\[
\ln\left(1+p_{n}\left|z\right|\right)^{n}=n\left[\dfrac{\lambda\left|z\right|}{n}+o\left(\dfrac{1}{n}\right)\right]\underset{n\to+\infty}{\longrightarrow}\lambda\left|z\right|,
\]
we obtain
\[
\lim_{n\to+\infty}\left[\text{e}^{np_{n}\left|z\right|}-\left(1+p_{n}\left|z\right|\right)^{n}\right]=0,
\]
and thus
\[
\lim_{n\to+\infty}\left(1+p_{n}z\right)^{n}=\text{e}^{\lambda z}.
\]

Taking $z=\text{e}^{\text{i}t}-1,$ it follows that
\[
\forall t\in\mathbb{R},\,\,\,\,\lim_{n\to+\infty}\varphi_{X_{n}}\left(t\right)=\text{e}^{\lambda\left(\text{e}^{\text{i}t}-1\right)},
\]
 which, by the Levy theorem, proves the result.

\end{proof}

\textbf{This first Poisson theorem is extended as follows.}

\begin{theorem}{Rare Events Theorem. Poisson Theorem}{Poisson_law_second}

Let, for each $n\in\mathbb{N}^{\ast},$ $\left\{ A_{n,j}:\,1\leqslant j\leqslant M_{n}\right\} $
be a finite family of independent events defined on the probabilized
space $\left(\Omega,\mathscr{A},P\right).$ Denote $P\left(A_{n,j}\right)=p_{n,j}$
and
\[
S_{n}=\sum^{M_{n}}_{j=1}\boldsymbol{1}_{A_{n,j}}.
\]
Assume that the sequence with general term $M_{n}$ tends by increasing
to $+\infty,$ that
\begin{equation}
\max_{1\leqslant j\leqslant M_{n}}p_{n,j}\underset{n\to+\infty}{\longrightarrow}0\,\,\,\,\text{and}\,\,\,\,\sum^{M_{n}}_{j=1}p_{n,j}\underset{n\to+\infty}{\longrightarrow}\lambda,\label{eq:M_n_assumptions}
\end{equation}
where $\lambda>0.$ 

Then the sequence $\left(S_{n}\right)_{n\in\mathbb{N}^{\ast}}$ converges
in law to the Poisson law $\mathscr{P}\left(\lambda\right).$

\end{theorem}

\begin{proof}{}{}

We again use the Lévy theorem. By independence of the events $A_{n,j},\,1\leqslant j\leqslant M_{n},$
for every $t\in\mathbb{R},$
\begin{align*}
\varphi_{S_{n}}\left(t\right) & =\prod^{M_{n}}_{j=1}\varphi_{\boldsymbol{1}_{A_{n,j}}}\left(t\right)=\prod^{M_{n}}_{j=1}\left[p_{n,j}\text{e}^{\text{i}t}+\left(1-p_{n,j}\right)\right]=\prod^{M_{n}}_{j=1}\left[1+p_{n,j}\left(\text{e}^{\text{i}t}-1\right)\right].
\end{align*}
Let $\text{Log}$ denote the principal branch of the complex logarithm.
By the Taylor formula with integral remainder, for every $z$ such
that $\left|z\right|<1,$ 
\[
\text{Log}\left(1+z\right)=z-z^{2}\intop^{1}_{0}\left(1-u\right)\dfrac{1}{\left(1+uz\right)^{2}}\text{d}u.
\]
Set $z=\text{e}^{\text{i}t}-1.$ 

Since $\max_{1\leqslant j\leqslant M_{n}}p_{n,j}\underset{n\to+\infty}{\longrightarrow}0,$
there exists $N\in\mathbb{N},$ such that for every $n\geqslant N,$
\[
\max_{1\leqslant j\leqslant M_{n}}\left|p_{n,j}z\right|<\dfrac{1}{2}.
\]
For every $n\geqslant N,$
\[
\text{Log}\varphi_{S_{n}}\left(t\right)=z\sum^{M_{n}}_{j=1}p_{n,j}-z^{2}\sum^{M_{n}}_{j=1}p^{2}_{n,j}\intop^{1}_{0}\dfrac{1-u}{\left(1+up_{n,j}z\right)^{2}}\text{d}u.
\]

By the triangle inequality, for every $n\geqslant N$ and every $u\in\left[0,1\right],$
\[
\left|1+up_{n,j}z\right|\geqslant1-p_{n,j}\left|z\right|\geqslant\dfrac{1}{2}.
\]
Hence, for every $n\geqslant N,$
\[
\left|\sum^{M_{n}}_{j=1}p^{2}_{n,j}\intop^{1}_{0}\dfrac{1-u}{\left(1+up_{n,j}z\right)^{2}}\text{d}u\right|\leqslant2\left[\max_{1\leqslant j\leqslant M_{n}}p_{n,j}\right]\left[\sum^{M_{n}}_{j=1}p_{n,j}\right].
\]
By the hypothesis, it follows that $\lim_{n\to+\infty}\text{Log}\varphi_{S_{n}}\left(t\right)=\lambda z.$
That is,
\[
\forall t\in\mathbb{R},\,\,\,\,\lim_{n\to+\infty}\varphi_{S_{n}}\left(t\right)=\text{e}^{\lambda\left(\text{e}^{\text{i}t}-1\right)},
\]
which, by the Lévy theorem, proves the result.

\end{proof}

\begin{remark}{}{}

The name of Theorem $\ref{th:Poisson_law_second}$ comes from the
fact that it shows that a random phenomenon can be represented as
a superposition of rare---that is, events with ``small'' probability,
in the sense of conditions $\refpar{eq:M_n_assumptions}$---and independent
events, follows a Poisson law. 

Moreover, this theorem extends the Poisson theorem $\ref{th:Poisson_theorem_first}.$
Indeed, suppose---with the notations of Theorem $\ref{th:Poisson_law_second}$---that
$M_{n}=n,$ and that, for each $n\in\mathbb{N}^{\ast},$ the family
of independent events $\left\{ A_{n,j}:\,1\leqslant j\leqslant n\right\} $
is such that $P\left(A_{n,j}\right)=p_{n},$ independently of $j$
with $1\leqslant j\leqslant n,$ where the sequence $\left(p_{n}\right)_{n\in\mathbb{N}}$
additionally satisfies the condition $\lim_{n\to+\infty}np_{n}=\lambda,$
with $\lambda>0.$

Then, the random variable $S_{n}$ follows the binomial law $\mathscr{B}\left(n,p_{n}\right)$
and the conditions $\refpar{eq:M_n_assumptions}$ are clearly satisfied,
since
\[
p_{n}\equiv\max_{1\leqslant j\leqslant n}p_{n,j}\underset{n\to+\infty}{\longrightarrow}0\,\,\,\,\text{and that}\,\,\,\,\sum^{n}_{j=1}p_{n,j}=np_{n}\underset{n\to+\infty}{\longrightarrow}\lambda.
\]
Both theorems therefore assert that the sequence $\left(S_{n}\right)_{n\in\mathbb{N}^{\ast}}$
converges in law to the Poisson law $\mathscr{P}\left(\lambda\right).$

\end{remark}

\section{Central Limit Theorem}\label{sec:Central-Limit-Theorem}

The \textbf{\index{central limit theorem}central limit theorem} of
probability theory states that, under fairly general conditions, the
law of the sum of a large number of independent random variable is
close to a normal law. There exists numerous versions of this theorem---particularly
those involving conditions of \textbf{\sindex[fam]{Lindeberg, Jarl}Lindeberg}{\bfseries\footnote{\textbf{\href{https://en.wikipedia.org/wiki/Jarl_Waldemar_Lindeberg}{Jarl Waldemar Lindeberg}}\sindex[fam]{Lindeberg, Jarl}
(1876-1932) was a Finnish mathematician. His main known work is on
the central limit theorem.}}\footnote{For a central limit theory with a Lindeberg-type condition, see for
instance \cite{renyiprobatheory,renyi2007foundations} respectively
p.443 and p.223.} type. 

\begin{figure}[t]
\begin{center}\includegraphics[width=0.4\textwidth]{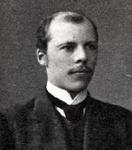}

{\tiny Source: \href{https://expo.oscapps.jyu.fi/s/minerva/item/57144}{Minerva}
Public Domain}\end{center}

\caption{\textbf{\protect\href{https://en.wikipedia.org/wiki/Jarl_Waldemar_Lindeberg}{Jarl Waldemar Lindeberg}}
(1876-1932)}\sindex[fam]{Lindeberg, Jarl}
\end{figure}

We will present only an elementary version. Modern proofs of the various
forms of the theorem rely on Lévy theorem and consist of performing
an asymptotic expansion of the characteristic function of the sum---centered
and reduced---of $n$ independent random variables.

By Theorem $\ref{th:der_int_param}$ in Chapter \ref{chap:PartIIChap13},
if a random variable admits a moment of order $2k,$ then its characteristic
function admits a Taylor expansion of order $2k.$ The following lemma
gives---in the case of a Taylor expansion of order 2---an upper-bound
of the remainder that might sometimes be useful.

\begin{lemma}{Taylor Expansion Remainder of a Characteristic Function Upper-Bound}{u_bound_rem_tayl_exp_ch_f}

If the real-valued random variable $X$ admits a second-order moment,
then its characteristic function $\varphi_{X}$ admits a second-order
Taylor expansion at 0 given, for every real number $t,$ by
\[
\varphi_{X}\left(t\right)=1+\text{i}t\mathbb{E}\left(X\right)-\dfrac{t^{2}}{2}\mathbb{E}\left(X^{2}\right)+o\left(t^{2}\right).
\]
More precisely, for every real number $t,$ we have the inequality,
\begin{equation}
\left|\varphi_{X}\left(t\right)-\left(1+\text{i}t\mathbb{E}\left(X\right)-\dfrac{t^{2}}{2}\mathbb{E}\left(X^{2}\right)\right)\right|\leqslant t^{2}\mathbb{E}\left(\min\left(X^{2},\left|t\right|\dfrac{\left|X\right|^{3}}{6}\right)\right).\label{eq:upper_boud_remainder_char_fct_with_order_2_taylor_exp}
\end{equation}

\end{lemma}

\begin{proof}{}{}

The second-order Taylor expansion with integral remainder yields,
for every real number $x,$
\[
\text{e}^{\text{i}x}=1+\text{i}x-x^{2}\intop^{1}_{0}\left(1-u\right)\text{e}^{\text{i}ux}\text{d}u.
\]
Since $\intop^{1}_{0}\left(1-u\right)\text{d}u=\dfrac{1}{2},$ we
obtain
\[
\text{e}^{\text{i}x}-\left(1+\text{i}x-\dfrac{x^{2}}{2}\right)=-x^{2}\intop^{1}_{0}\left(1-u\right)\left[\text{e}^{\text{i}ux}-1\right]\text{d}u.
\]
It follows that
\[
\left|\text{e}^{\text{i}x}-\left(1+\text{i}x-\dfrac{x^{2}}{2}\right)\right|\leqslant x^{2}.
\]
The third-order Taylor expansion gives, for every real number $x,$
\[
\text{e}^{\text{i}x}=1+\text{i}x-\dfrac{x^{2}}{2}-\text{i}\dfrac{x^{3}}{3}\intop^{1}_{0}\left(1-u\right)^{2}\text{e}^{\text{i}ux}\text{d}u.
\]
This implies that
\[
\left|\text{e}^{\text{i}x}-\left(1+\text{i}x-\dfrac{x^{2}}{2}\right)\right|\leqslant\dfrac{\left|x\right|^{3}}{6}.
\]
Hence, for every real number $x,$
\[
\left|\text{e}^{\text{i}x}-\left(1+\text{i}x-\dfrac{x^{2}}{2}\right)\right|\leqslant\min\left(x^{2},\dfrac{\left|x\right|^{3}}{6}\right).
\]

The upper-bound $\refpar{eq:upper_boud_remainder_char_fct_with_order_2_taylor_exp}$
follows immediately. Moreover, by dominated convergence---consider
an arbitrary sequence tending to 0---we have
\[
\lim_{t\to0}\mathbb{E}\left(\min\left(X^{2},\left|t\right|\dfrac{\left|X\right|^{3}}{6}\right)\right)=0,
\]
which completes the proof.

\end{proof}

\begin{remark}{}{}

If this sharp upper-bound for the remainder is not required, it suffices,
in order to establish the Taylor expansion, to apply the Taylor-Young
formula to $\varphi_{X},$ which is twice differentiable in this setting---see
Theorem $\ref{th:der_int_param}$, Chapter \ref{chap:PartIIChap13}.

\end{remark}

\begin{theorem}{Central Limit Theorem}{}

Let $\left(X_{n}\right)_{n\in\mathbb{N}^{\ast}}$ be a sequence of
random variables defined on the same probabilized space $\left(\Omega,\mathscr{A},P\right),$
taking values in $\mathbb{R}^{d},$ independent, with same law, and
admitting a second-order\footnotemark moment. The sequence with general
term $Y_{n},$ defined for every $n\in\mathbb{N}^{\ast}$ by
\[
Y_{n}=\dfrac{1}{\sqrt{n}}\sum^{n}_{j=1}\left(X_{j}-\mathbb{E}\left(X_{j}\right)\right),
\]
converges in law to the Gaussian law $\mathscr{N}_{\mathbb{R}^{d}}\left(0,C_{X_{1}}\right),$
where $C_{X_{1}}$ is the covariance matrix of the $X_{j}.$

In particular, if $d=1$ and if $Z_{n}=\dfrac{Y_{n}}{\sigma_{n}},$
where $\sigma_{X_{1}}$ is the standard-deviation of $X_{j},$ then
the sequence of cumulative distribution functions $F_{Z_{n}}$ of
the $Z_{n}$ converge pointwise to $\Phi,$ the cumulative distribution
function of the law $\mathscr{N}_{\mathbb{R}}\left(0,1\right),$ given,
for every real number $z,$ by
\[
\Phi\left(z\right)=\intop^{z}_{-\infty}\dfrac{1}{\sqrt{2\pi}}\text{e}^{-\frac{x^{2}}{2}}\text{d}x.
\]

\end{theorem}

\footnotetext{That is square-integrable---integrable in norm.

}

\begin{proof}{}{}

Since the random variables $X_{j}$ are independent and with same
law, the characteristic function of $Y_{n}$ is, for every $t\in\mathbb{R}^{d},$
given by
\[
\varphi_{Y_{n}}\left(t\right)=\prod^{n}_{j=1}\varphi_{\left(X_{j}-\mathbb{E}\left(X_{j}\right)\right)}\left(\dfrac{t}{\sqrt{n}}\right)=\left[\varphi_{\left\langle X_{1}-\mathbb{E}\left(X_{1}\right),t\right\rangle }\left(\dfrac{t}{\sqrt{n}}\right)\right]^{n}.
\]

Lemma $\ref{lm:u_bound_rem_tayl_exp_ch_f},$ applied to the centered
real-valued random variable $\left\langle X_{1}-\mathbb{E}\left(X_{1}\right),t\right\rangle ,$
yields the asymptotic expansion
\[
\varphi_{Y_{n}}\left(t\right)=\left[1-\dfrac{1}{2n}\mathbb{E}\left(\left\langle X_{1}-\mathbb{E}\left(X_{1}\right),t\right\rangle ^{2}\right)+o\left(\dfrac{1}{n}\right)\right]^{n}.
\]
Lemma $\ref{lm:cpl_exp_lim_power}$ ensures that the sequence with
general term $\varphi_{Y_{n}}\left(t\right)$ converges and that
\[
\lim_{n\to+\infty}\varphi_{Y_{n}}\left(t\right)=\text{e}^{-\frac{1}{2}\mathbb{E}\left(\left\langle X_{1}-\mathbb{E}\left(X_{1}\right),t\right\rangle ^{2}\right)}.
\]
Since $\mathbb{E}\left(\left\langle X_{1}-\mathbb{E}\left(X_{1}\right),t\right\rangle ^{2}\right)=\left\langle C_{X_{1}}t,t\right\rangle ,$
\[
\lim_{n\to+\infty}\varphi_{Y_{n}}\left(t\right)=\text{e}^{-\frac{1}{2}\left\langle C_{X_{1}}t,t\right\rangle }.
\]
The Lévy theorem implies that
\[
Y_{n}\stackrel[n\to+\infty]{\mathscr{L}}{\longrightarrow}\mathscr{N}_{\mathbb{R}^{d}}\left(0,C_{X_{1}}\right).
\]
If $d=1,$ then for every real number $t,$
\[
\varphi_{Z_{n}}\left(t\right)=\varphi_{Y_{n}}\left(\dfrac{t}{\sigma_{X_{1}}}\right)\,\,\,\,\text{and\,thus}\,\,\,\,\lim_{n\to+\infty}\varphi_{Z_{n}}\left(t\right)=\text{e}^{-\frac{t^{2}}{2}}.
\]
Hence, the sequence of probabilities $P_{Z_{n}}$ converges narrowly
to the law $\mathscr{N}_{\mathbb{R}}\left(0,1\right).$ Since $\Phi$
is continuous, this is equivalent to the pointwise convergence of
the sequence of cumulative distribution functions of the $Z_{n}$
to the cumulative distribution function $\Phi$ of the limiting law.

\end{proof}

\begin{remark}{}{}

The cumulative distribution function $\Phi$ of the law $\mathscr{N}_{\mathbb{R}}\left(0,1\right)$
is tabulated\footnotemark. We recall three commonly nondecreasing
values of the function $\Phi:$ they were already given in Part \ref{part:Introduction-to-Probability}.

\boxeq{
\begin{gather*}
\Phi\left(1.64\right)-\Phi\left(-1.64\right)\approx0.9\\
\Phi\left(1.96\right)-\Phi\left(-1.96\right)\approx0.95\\
\Phi\left(3.09\right)-\Phi\left(-3.09\right)\approx0.99
\end{gather*}
}

\end{remark}

\footnotetext{A table of the cumulative distribution function $\Phi$
of the centered reduced Gauss law is given for instance in \cite{renyiprobatheory,renyi2007foundations}.

}

\begin{figure}[t]
\begin{center}\includegraphics[width=0.4\textwidth]{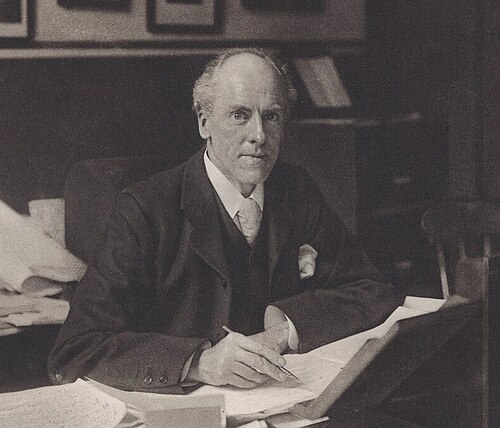}

{\tiny Public Domain}\end{center}

\caption{\textbf{\protect\href{https://en.wikipedia.org/wiki/Karl_Pearson}{Karl Pearson}}
(1857-1936)}\sindex[fam]{Pearson, Karl}
\end{figure}

An application of the central limit theorem is the proof of the \textbf{Karl
Pearson}{\bfseries\footnote{\textbf{\href{https://en.wikipedia.org/wiki/Karl_Pearson}{Karl Pearson}}
\textbf{\sindex[fam]{Pearson, Karl}} (1857-1936) was an English mathematician
and biostatistician. He write foundational works of the discipline
of mathematical statistics, creating the first deparment of mathematical
statistics at University College London in 1911. }}\textbf{\sindex[fam]{Pearson, Karl} theorem}\index{Karl Pearson theorem},
which forms the foundation of the \textbf{chi-squared test}.

\begin{theorem}{Karl Pearson Theorem}{}

Fix $k\in\mathbb{N}^{\ast}.$ For each $n\in\mathbb{N}^{\ast},$ let
$\left(A^{n}_{j}\right)_{1\leqslant j\leqslant k}$ be a partition
of $\Omega$ into $\mathscr{A}-$measurable sets. Assume that these
partitions are independent, in the sense that the families, indexed
by $n,$ consisting of the elements of the partitions, are independent.
Moreover, suppose that
\[
\forall n\in\mathbb{N}^{\ast},\,\,\,\,\forall j\in\left\llbracket 1,k\right\rrbracket ,\,\,\,\,P\left(A^{n}_{j}\right)=p_{j},
\]
where $p_{j}>0$ and $\sum^{k}_{j=1}p_{j}=1.$ 

For each $j\in\left\llbracket 1,k\right\rrbracket ,$ define the real-valued
random variables
\[
N^{n}_{j}=\sum^{n}_{l=1}\boldsymbol{1}_{A^{l}_{j}}.
\]
as well as the random variable
\[
\chi^{2}_{k,n}=\sum^{k}_{j=1}\dfrac{\left(N^{n}_{j}-np_{j}\right)^{2}}{np_{j}}\equiv n\sum^{k}_{j=1}\dfrac{\left(\frac{1}{n}N^{n}_{j}-p_{j}\right)^{2}}{p_{j}}.
\]
Then the sequence of laws $P_{\chi^{2}_{k,n}}$ converges narrowly
to the \textbf{chi-squared law }$\chi^{2}_{k-1}$ \textbf{with $k-1$
degrees of freedom}. That is, the sequence of random variables $\chi^{2}_{k,n}$
converges in law to the \textbf{chi-squared law with $k-1$ degrees
of freedom}.

\end{theorem}

\begin{proof}{}{}

For every $n\in\mathbb{N}^{\ast},$ define the random variables $X^{n}$
and $N^{n},$ taking values in $\mathbb{R}^{k},$ by
\[
X^{n}=\left(\begin{array}{c}
\boldsymbol{1}_{A^{n}_{1}}\\
\vdots\\
\boldsymbol{1}_{A^{n}_{k}}
\end{array}\right)\,\,\,\,\text{and}\,\,\,\,N^{n}=\sum^{n}_{j=1}X^{j}\equiv\left(\begin{array}{c}
N^{n}_{1}\\
\vdots\\
N^{n}_{k}
\end{array}\right).
\]
Recall that the law of $N^{n}$ is the multinomial law $\mathscr{M}\left(n,p_{1},p_{2},\cdots,p_{k-1}\right),$
and that the random variables $X^{n}$ have the same law, of expectation
and covariance matrix $C_{X^{n}}$ given by
\[
\mathbb{E}\left(X^{n}\right)=\left(\begin{array}{c}
p_{1}\\
\vdots\\
p_{k}
\end{array}\right)\,\,\,\,\text{and}\,\,\,\,\left(C_{X^{n}}\right)_{i,j}=\begin{cases}
p_{i}\left(1-p_{i}\right), & \text{if\,}i=j,\\
-p_{i}p_{j}, & \text{if}\,i\neq j,
\end{cases}
\]
which can be written, denoting $p$ the vector of components $p_{j},\,j\in\left\llbracket 1,k\right\rrbracket ,$
\[
\mathbb{E}\left(X^{n}\right)=p\,\,\,\,\text{and}\,\,\,\,C_{X^{n}}=\left(\begin{array}{cccc}
p_{1} & 0 & \cdots & 0\\
0 & \ddots &  & \vdots\\
\vdots &  & \ddots & 0\\
0 & \cdots & 0 & p_{k}
\end{array}\right)-pp^{\ast}.
\]
Since the random variables $X^{j}$ are independent, the central limit
theorem implies that the sequence with general term $Y_{n},$ defined
for every $n\in\mathbb{N}^{\ast}$ by
\[
Y_{n}=\dfrac{1}{\sqrt{n}}\left[\sum^{n}_{j=1}X_{j}-np\right],
\]
converges in law to the Gaussian law $\mathscr{N}_{\mathbb{R}^{d}}\left(0,C_{X_{1}}\right).$ 

By the Lévy theorem, this is equivalent to 
\begin{equation}
\forall t\in\mathbb{R}^{k},\,\,\,\,\lim_{n\to+\infty}\varphi_{Y_{n}}\left(t\right)=\text{e}^{-\frac{\left\langle C_{X_{1}}t,t\right\rangle }{2}}.\label{eq:lim_charact_f_Y_n}
\end{equation}
Let $M$ be the diagonal matrix defined by
\[
\forall j\in\left\llbracket 1,k\right\rrbracket ,\,\,\,\,M_{j,j}=\dfrac{1}{p_{j}},\,\,\,\,\text{and}\,\,\,\,M_{i,j}=0,\,\,\,\,\text{if}\,i\neq j.
\]
Then
\[
\chi^{2}_{k,n}=\left\langle MY_{n},Y_{n}\right\rangle =\left\Vert M^{\frac{1}{2}}Y_{n}\right\Vert ^{2}.
\]
Since, for every $t\in\mathbb{R}^{k},$
\[
\varphi_{M^{\frac{1}{2}}Y_{n}}\left(t\right)=\varphi_{Y_{n}}\left(M^{\frac{1}{2}}t\right),
\]
it follows, by $\refpar{eq:lim_charact_f_Y_n},$ that
\[
\forall t\in\mathbb{R}^{k},\,\,\,\,\lim_{n\to+\infty}\varphi_{M^{\frac{1}{2}}Y_{n}}\left(t\right)=\text{e}^{-\frac{\left\langle M^{\frac{1}{2}}C_{X_{1}}M^{\frac{1}{2}}t,t\right\rangle }{2}}.
\]
Using again the Lévy theorem, the sequence with general term $M^{\frac{1}{2}}Y_{n}$
converges in law to the Gaussian law $\mathscr{N}_{\mathbb{R}^{d}}\left(0,M^{\frac{1}{2}}C_{X_{1}}M^{\frac{1}{2}}\right).$
But,
\[
M^{\frac{1}{2}}C_{X_{1}}M^{\frac{1}{2}}=M^{\frac{1}{2}}\left[M^{-1}-pp^{\ast}\right]M^{\frac{1}{2}}=I-\left(M^{\frac{1}{2}}p\right)\left(M^{\frac{1}{2}}p\right)^{\ast}
\]
and 
\[
\left\Vert M^{\frac{1}{2}}p\right\Vert ^{2}=\sum^{k}_{j=1}p_{j}=1.
\]
That is, $M^{\frac{1}{2}}p$ is a unit vector. Choosing the orthogonal
transformation $O,$ such that $O\left(M^{\frac{1}{2}}p\right)=e_{1},$
we then have
\[
O\left[M^{\frac{1}{2}}C_{X_{1}}M^{\frac{1}{2}}\right]O^{\ast}=I-\left(e_{1}\right)\left(e_{1}\right)^{\ast}=\left(\begin{array}{cc}
0 & 0\\
0 & I_{\mathbb{R}^{k-1}}
\end{array}\right).
\]
Hence, still by the Lévy theorem, the sequence with general term $OM^{\frac{1}{2}}Y_{n}$
converges in law to the Gaussian law $\mathscr{N}_{\mathbb{R}^{d}}\left(0,O\left[M^{\frac{1}{2}}C_{X_{1}}M^{\frac{1}{2}}\right]O^{\ast}\right)=\mathscr{N}_{\mathbb{R}^{d}}\left(0,I-\left(e_{1}\right)\left(e_{1}\right)^{\ast}\right).$

Since $O$ is orthogonal, 
\[
\chi^{2}_{k,n}=\left\Vert OM^{\frac{1}{2}}Y_{n}\right\Vert ^{2}.
\]
It follows that, if $U$ is a random variable of law $\mathscr{N}_{\mathbb{R}^{d}}\left(0,I-\left(e_{1}\right)\left(e_{1}\right)^{\ast}\right),$
then the sequence of law of $\chi^{2}_{k,n}$ converges narrowly to
the law of $\left\Vert U\right\Vert ^{2},$ that is the chi-squared
law with $k-1$ degrees of freedom. 

\end{proof}

\begin{remark}{}{}

Let $q=\left(q_{1},\cdots,q_{k}\right)^{\ast}$ be a vector of $\mathbb{R}^{k}$
distinct from $p.$ Define
\[
\kappa^{2}_{k,n}\equiv n\sum^{k}_{j=1}\dfrac{\left(\frac{1}{n}N^{n}_{j}-q_{j}\right)^{2}}{q_{j}}.
\]
By the strong law of large numbers, for every $j\in\left\llbracket 1,k\right\rrbracket ,$
\[
\lim_{n\to+\infty}\dfrac{1}{n}N^{n}_{j}=p_{j}\,\,\,\,\,P-\text{almost surely.}
\]

Since $q\neq p,$ the sequence with general term
\[
\sum^{k}_{j=1}\dfrac{\left(\frac{1}{n}N^{n}_{j}-q_{j}\right)^{2}}{q_{j}}
\]
converges $P-$almost surely to a constant $\alpha>0.$ Consequently,
the sequence with general term $\kappa^{2}_{k,n}$ converges $P-$almost
surely to $+\infty$ as $n$ tends to infinity.

\end{remark}

\begin{example}{At the Basis of the Chi-Squared Test}{ex_at_basis_chi_two}

Let $\left(X_{n}\right)_{n\in\mathbb{N}^{\ast}}$ be a sequence of
independent random variables taking values in $\mathbb{R}^{d},$ all
with same law $\mu.$ Let $\left(D_{j}\right)_{1\leqslant j\leqslant k}$
be a Borel partition of $\mathbb{R}^{d}$ such that $\mu\left(D_{j}\right)=p_{j}>0,$
for every $j\in\left\llbracket 1,k\right\rrbracket .$ 

For each $n\in\mathbb{N},$ the $A^{n}_{j}=X^{-1}_{n}\left(D_{j}\right),\,j\in\left\llbracket 1,k\right\rrbracket ,$
form a partition of $\Omega;$ additionally these partitions are independent
as $n$ varies. For every $n\in\mathbb{N}^{\ast},$
\[
P\left(A^{n}_{j}\right)=\mu\left(D_{j}\right)=p_{j}.
\]
For each $j\in\left\llbracket 1,k\right\rrbracket ,$ define the real-valued
random variables
\[
N^{n}_{j}=\sum^{n}_{l=1}\boldsymbol{1}_{D_{j}}\left(X_{l}\right),
\]
Define the random variables
\[
\chi^{2}_{k,n}=n\sum^{k}_{j=1}\dfrac{\left(\frac{1}{n}N^{n}_{j}-p_{j}\right)^{2}}{p_{j}}\,\,\,\,\text{and}\,\,\,\,\kappa^{2}_{k,n}\equiv n\sum^{k}_{j=1}\dfrac{\left(\frac{1}{n}N^{n}_{j}-q_{j}\right)^{2}}{q_{j}}.
\]
Then the sequence of laws $P_{\chi^{2}_{k,n}}$ converges narrowly
to the chi-squared law $\chi^{2}_{k-1}$ with $k-1$ degrees of freedom.
Moreover, for every vector $q\neq p,$ the sequence with general term
$\kappa^{2}_{k,n}$ converges to $+\infty$ $P-$almost surely.

\end{example}

This example is at the foundation of the chi-squared test. We now
present the underlying problem, followed by an example of its practical
implementation.

\paragraph*{The problem.}

After modelling a random phenomenon probabilistically, we consider
a random variable $X$ taking values in $\mathbb{R}^{d},$ representing
a vector quantity associated with this phenomenon. The law $\mu$
of $X$ is \textbf{unknown} to the experimenter. However, after computations
and reasonning, the experimenter is led to formulate a \textbf{hypothesis}
about this law. 

The goal is to \textbf{test} the hypothesis $H$ that $X$ has law
$\mu,$ given a \textbf{sample} of size $n,$ 
\[
\uuline{x_{n}}=\left(x_{1},x_{2},\cdots,x_{n}\right),
\]
obtained by observing $n$ \textbf{independent} realizations of the
phenomenon. This sample is assumed to be the outcome of $n$ random
variables $X_{1},X_{2},\cdots,X_{n},$ independent and with the same
law\footnote{Recall that $\uuline{X_{n}}=\left(X_{1},X_{2},\cdots,X_{n}\right)$
is called the \textbf{empirical sample} of size $n$ associated with
the random variable $X.$} as $X.$

We keep the notations of Example $\ref{ex:ex_at_basis_chi_two},$
where the sets $D_{j}$ are called \textbf{classes}. 

\subsubsection*{Empirical and Theoretical Frequencies}

For each vector $\uuline{x_{n}}$ of $\mathbb{R}^{d_{n}},$ define
\boxeq{
\[
f_{j}\left(\uuline{x_{n}}\right)=\dfrac{1}{n}\sum^{n}_{l=1}\boldsymbol{1}_{D_{j}}\left(x_{l}\right),
\]
}the frequency of the number of points $x_{l}$ that fall into the
class $D_{j}.$ The \textbf{observed number} of points $x_{l}$ in
$D_{j}$ for the sample $\uuline{x_{n}}$ is then $n_{j}\left(\uuline{x_{n}}\right)=nf_{j}\left(\uuline{x_{n}}\right),$
while the \textbf{theoretical number} in $D_{j}$ under hypothesis
$H$ is $np_{j}.$ The \textbf{chi-squared distance} between observed
and theoretical counts is\boxeq{
\[
\Delta\left(\uuline{x_{n}}\right)=\sum^{k}_{j=1}\dfrac{\left(n_{j}\left(\uuline{x_{n}}\right)-np_{j}\right)^{2}}{np_{j}}\equiv n\sum^{k}_{j=1}\dfrac{\left(f_{j}\left(\uuline{x_{n}}\right)-p_{j}\right)^{2}}{p_{j}}.
\]
}Thus, $\Delta\left(\uuline{x_{n}}\right)$ is a realization of the
random variable $\chi^{2}_{k,n}.$

\subsubsection*{Rejection Region}

For every real number $c>0,$ define the set \boxeq{
\[
R^{n}_{c}=\left\{ \uuline{x_{n}}=\left(x_{1},x_{2},\cdots,x_{n}\right):\,\Delta\left(\uuline{x_{n}}\right)>c\right\} ,
\]
}called the \textbf{rejection region} of the hypothesis $H.$ 

From the previous example,
\[
P\left(\uuline{X_{n}}\in R^{n}_{c}\right)=P\left(\Delta\left(\uuline{X_{n}}\right)>c\right)\underset{n\to+\infty}{\longrightarrow}\chi^{2}_{k-1}\left(\left]c,+\infty\right[\right),
\]
independently of the law $\mu$ of $X.$ 

\subsubsection*{Decision Rule and Error Risk}

This leads to adopt the following \textbf{chi-squared test} \textbf{rule}: 
\begin{itemize}
\item Accept the hypothesis $H$ if
\[
\Delta\left(\uuline{x_{n}}\right)\leqslant c,
\]
\item Reject $H$ otherwise. 
\end{itemize}
The \textbf{error risk} corresponds to the probability of rejecting
$H$ when $H$ is true. For sufficiently large $n,$ it is approximately
equal to 
\[
\chi^{2}_{k-1}\left(\left]c,+\infty\right[\right).
\]
In practice, this approximation is considered reliable as soon as
$np_{j}\geqslant5,$ for every $j\in\left\llbracket 1,k\right\rrbracket .$

\subsubsection*{Choice of the Critical Value}

For a prescribed error level $\alpha$---most often 0.1 or 0.05---one
determines the real number $c_{\alpha}$ such that 
\[
\chi^{2}_{k-1}\left(\left]c,+\infty\right[\right)=\alpha.
\]
The hypothesis $H$ is then accepted or rejected, given the observed
sample $\uuline{x_{n}}$ depending on whether 
\[
\Delta\left(\uuline{x_{n}}\right)\leqslant c_{\alpha}
\]
 or not.

\begin{example}{}{}

We throw a dice $n$ times. We obtain $n_{j}$ occurrences of the
number $j,$ for $j\in\left\llbracket 1,6\right\rrbracket .$ We wonder
whether the dice is fair, given the following samples:
\[
\begin{array}{lllllll}
n=60 & n_{1}=11 & n_{2}=8 & n_{3}=12 & n_{4}=9 & n_{5}=8 & n_{6}=12\\
n=600 & n_{1}=110 & n_{2}=80 & n_{3}=120 & n_{4}=90 & n_{5}=80 & n_{6}=120.
\end{array}
\]

\end{example}

\begin{solutionexample}{}{}

With the dice throws, we associate a random variable $X$ with law
$\mu$ supported on the set $\left\{ 1,2,\cdots,6\right\} .$ The
classes are the singletons $\left\{ j\right\} ,\,j\in\left\llbracket 1,6\right\rrbracket .$
The hypothesis $H$ is that the law $\mu$ is uniform. We choose $\alpha=0.05,$
so that the table gives $c_{\alpha}=11.1,$ that is,
\[
\chi^{2}_{5}\left(\left]11.1,+\infty\right[\right)=0.05.
\]

We then compute $\Delta\left(\uuline{x_{n}}\right)$ for these two
samples.
\begin{itemize}
\item For $n=60,$
\[
\Delta\left(\uuline{x_{n}}\right)=\dfrac{1}{\frac{60}{6}}\left[\substack{\left(11-10\right)^{2}+\left(8-10\right)^{2}+\left(12-10\right)^{2}+\left(9-10\right)^{2}+\left(8-10\right)^{2}+\left(11-10\right)^{2}}
\right],
\]
hence
\[
\Delta\left(\uuline{x_{n}}\right)=1.8\leqslant11.1,
\]
 which leads us to accept $H,$ given $\uuline{x_{60}}.$
\item For $n=600,$
\[
\Delta\left(\uuline{x_{n}}\right)=\dfrac{1}{\frac{600}{6}}\left[\substack{\left(110-100\right)^{2}+\left(80-100\right)^{2}+\left(120-100\right)^{2}+\left(90-100\right)^{2}+\left(80-100\right)^{2}+\left(110-100\right)^{2}}
\right],
\]
hence
\[
\Delta\left(\uuline{x_{n}}\right)=18\geqslant11.1,
\]
 which leads to reject $H,$ given $\uuline{x_{600}}.$
\end{itemize}
\end{solutionexample}

\section{Estimation}\label{sec:Estimation}

We now briefly present the problem of estimating the law of a random
variable and describe a construction method of the estimator, known
as the \textbf{maximum likelihood estimator}. Under suitable regularity
hypothesis on the density, this estimator possesses important \textbf{asymptotic
properties}. For a detailed treatment of estimation theory, the reader
may consult the books \cite{fourgeaud1967statistique} and \cite{dacuna1993probabilites}. 

\paragraph*{The Problem of Parametric Estimation}

After modelling a random phenomenon probabilistically, we consider
a random variable $X$ taking values in $\mathbb{R},$ intended to
represent a real-valued quantity associated with the phenomenon. The
law $\mu$ of $X$ is \textbf{unknown} to the experimenter. However,
after computations and reasonning, the experimenter is led to suppose
that this law belongs to a family of laws depending on a parameter
$\theta\in\Theta,$ where $\Theta$ is an open subset of $\mathbb{R}^{p}.$ 

The objective is to \textbf{estimate} the true value $\theta_{0}$
of the parameter based on a \textbf{sample} of size $n,$ 
\[
\uuline{x_{n}}=\left(x_{1},x_{2},\cdots,x_{n}\right),
\]
obtained by observing $n$ independent realization of the phenomenon.
This sample is supposed to be the realization of $n$ random variables
$X_{1},X_{2},\cdots,X_{n},$ independent and with the same law as
$X.$

\paragraph*{Statistic Modeling of the Problem}

We consider a \textbf{statistic structure}, that is, a family of probabilized
spaces $\left[\left(\Omega,\mathscr{A},P_{\theta}\right)\right]_{\theta\in\Theta},$
where $\Theta$ is an open subset of $\mathbb{R}^{p}.$ On each we
define a real-valued random variable $X$ and an infinite \textbf{sample},
that is, a sequence $\left(X_{n}\right)_{n\in\mathbb{N}^{\ast}}$
of $P_{\theta}-$independent random variables, all having the same
law $\mu_{\theta}$ as $X$---the image measure of $P_{\theta}$
by $X$---for every $\theta\in\Theta.$ 

Assume that the function $\theta\mapsto P_{\theta}$ is injective.
Let $g$ be a function from $\Theta$ into an open subset $\Theta^{\prime}$
of $\mathbb{R}^{k}$ with $k\leqslant p.$ An estimator of $g\left(\theta\right),$
based on a sample of size $n,$ is a random variable---called a \textbf{statistic\index{statistic}}
by statisticians---of the form
\[
T_{n}=\varphi_{n}\left(\uuline{X_{n}}\right),
\]
where $\varphi_{n}$ is a measurable function from $\mathbb{R}^{n}$
to $\Theta^{\prime}.$ 

This estimator is said to be \textbf{\index{unbiased}unbiased} if
\[
E_{\theta}\left(T_{n}\right)=g\left(\theta\right),
\]
where $E_{\theta}$ denotes\textbf{ the integration with respect to
the probability $P_{\theta}.$} 

When the sample is infinite, the sequence $T=\left(T_{n}\right)_{n\in\mathbb{N}^{\ast}}$
is called an \textbf{estimator\index{estimator}} of $g\left(\theta\right).$
It is of interest if it is consistent---in probability or almost
surely---that is, if the sequence $\left(T_{n}\right)_{n\in\mathbb{N}^{\ast}}$
converges in $P_{\theta}-$probability, respectively $P_{\theta}-$almost
surely.

The \textbf{\mindex{maximum likelihood!method}method of maximum likelihood}
often provides such estimators. Its foundation is only empirical,
based on the following experiment. Suppose we draw a ball at random
from one of two urns $U_{1}$ and $U_{2},$ without knowing which
urn was chosen but knowing their composition. For instance, $U_{1}$
contains one red ball and nine black balls, while $U_{2}$ contains
nine red balls and one black ball. If a red ball is drawn, we are
naturally let to think that the draw was made from $U_{2},$ since
this urn assigns a higher probability to drawing a red ball.

In the following, we make the hypothesis that there exists a $\sigma-$finite
measure $\mu$ on $\mathbb{R}$---often the Lebesgue measure---such
that, for every $\theta\in\Theta,$ the law $\mu_{\theta}$ of $X$
has density the function $f\left(\cdot,\theta\right)$ with respect
to $\mu.$ Then the random variable $\uuline{X_{n}}$ has density,
with respect to the product measure $\mu^{\otimes n},$ the function
$L_{n},$ called the \textbf{likelihood function}\mindex{likelihood!function},
relatively to the sample of size $n,$ defined for every $\uuline{x_{n}}\in\mathbb{R}^{n}$
by \boxeq{
\[
L_{n}\left(\uuline{x_{n}},\theta\right)=\prod^{n}_{j=1}f\left(x_{j},\theta\right).
\]
}

An estimator $T_{n}$ of $\theta$ is called a \textbf{\mindex{maximum likelihood!estimator}maximum
likelihood estimator} of $\theta$ if it can be written in the form
\boxeq{
\[
T_{n}=\widehat{\varphi}_{n}\left(X_{n}\right),
\]
}where $\widehat{\varphi}_{n}$ is a measurable function satisfying
the condition\boxeq{
\begin{equation}
\forall\uuline{x_{n}}\in\mathbb{R}^{n},\,\,\,\,L_{n}\left(\uuline{x_{n}},\widehat{\varphi}_{n}\left(\uuline{x_{n}}\right)\right)=\sup_{\theta\in\Theta}L_{n}\left(\uuline{x_{n}},\theta\right).\label{eq:likelihood_f_x_n_varphi}
\end{equation}
}If, moreover $f\left(x,\cdot\right)$ is differentiable, then $\widehat{\varphi}_{n}\left(\uuline{x_{n}}\right)$
is solution of the \textbf{likelihood equation}\mindex{likelihood!equation}\boxeq{
\begin{equation}
\forall\uuline{x_{n}}\in\mathbb{R}^{n},\,\,\,\,\dfrac{\partial}{\partial\theta}L_{n}\left(\uuline{x_{n}},\widehat{\varphi}_{n}\left(\uuline{x_{n}}\right)\right)=0.\label{eq:likelihood_eq}
\end{equation}
}

It should be noted that $\widehat{\varphi}_{n}\left(\uuline{x_{n}}\right)$
is, a priori, just a stationary point; further analysis is required
to ensure the existence of the maximum.

In the case where $f\left(x,\theta\right)>0$ for every $\left(x,\theta\right),$
the equation $\refpar{eq:likelihood_f_x_n_varphi}$ is equivalent
to the more convenient equation obtained by taking logarithms. Then
$\widehat{\varphi}_{n}$ is solution of the equation\boxeq{
\[
\forall\uuline{x_{n}}\in\mathbb{R}^{n},\,\,\,\,\ln L_{n}\left(\uuline{x_{n}},\widehat{\varphi}_{n}\left(\uuline{x_{n}}\right)\right)=\sup_{\theta\in\Theta}\ln L_{n}\left(\uuline{x_{n}},\theta\right).
\]
}

If moreover $f\left(x,\cdot\right)$ is differentiable, $\widehat{\varphi}_{n}\left(\uuline{x_{n}}\right)$
is solution of the \textbf{log-likelihood equation}\index{log-likelihood equation}\boxeq{
\[
\forall\uuline{x_{n}}\in\mathbb{R}^{n},\,\,\,\,\dfrac{\partial}{\partial\theta}\ln L_{n}\left(\uuline{x_{n}},\widehat{\varphi}_{n}\left(\uuline{x_{n}}\right)\right)=0.
\]
}Such estimators, when they exist, are generally not unique. They
are often consistant and possess \textbf{asymptotic normality} properties.
The problem of existence is a maximization problem. The conditions
ensuring existence are frequently of a differentiability type, although
this is not always the case. We now present an example of each situation,
using the previous notation.

\paragraph*{Gaussian Case}

$\theta=\mathbb{R}\times\mathbb{R}^{+\ast},$ and for every $\theta=\left(m,\sigma^{2}\right),$
$\mu_{\theta}=\mathscr{N}_{\mathbb{R}}\left(m,\sigma^{2}\right).$

For every $\uuline{x_{n}}\in\mathbb{R}^{n},$
\[
\ln L_{n}\left(\uuline{x_{n}},\theta\right)=\sum^{n}_{j=1}\ln f\left(x_{j},\theta\right)=-\dfrac{n}{2}\left(\text{\ensuremath{\ln}2\ensuremath{\pi}-\ensuremath{\ln\sigma^{2}}}\right)-\dfrac{1}{2\sigma^{2}}\sum^{n}_{j=1}\left(x_{j}-m\right)^{2}.
\]
Hence,
\[
\dfrac{\partial}{\partial m}\ln L_{n}\left(\uuline{x_{n}},\left(m,\sigma^{2}\right)\right)=\dfrac{1}{\sigma^{2}}\sum^{n}_{j=1}\left(x_{j}-m\right),
\]
\[
\dfrac{\partial}{\partial\sigma^{2}}\ln L_{n}\left(\uuline{x_{n}},\left(m,\sigma^{2}\right)\right)=-\dfrac{n}{2}\dfrac{1}{\sigma^{2}}+\dfrac{1}{2\sigma^{4}}\sum^{n}_{j=1}\left(x_{j}-m\right)^{2}.
\]

The solutions of the \textbf{maximum likelihood equation} have to
be sought among the stationary points, that is, the solutions of the
log-likelihood equations
\[
\dfrac{\partial}{\partial m}\ln L_{n}\left(\uuline{x_{n}},\left(m,\sigma^{2}\right)\right)=0\,\,\,\,\text{and}\,\,\,\,\dfrac{\partial}{\partial\sigma^{2}}\ln L_{n}\left(\uuline{x_{n}},\left(m,\sigma^{2}\right)\right)=0.
\]
In the present case,
\[
\widehat{m_{n}}=\dfrac{1}{n}\sum^{n}_{j=1}x_{j}\,\,\,\,\text{and}\,\,\,\,\widehat{\sigma^{2}_{n}}=\dfrac{1}{n}\sum^{n}_{j=1}\left(x_{j}-\widehat{m_{n}}\right)^{2}.
\]
It remains to verify that this stationary point corresponds to a \textbf{maximum}.
One could examine the second derivatives, but here this can be seen
directly from the following computation:
\begin{align*}
 & \ln L_{n}\left(\uuline{x_{n}},\left(\widehat{m_{n}},\widehat{\sigma^{2}_{n}}\right)\right)-\ln L_{n}\left(\uuline{x_{n}},\left(m,\sigma^{2}\right)\right)\\
 & \qquad=-\dfrac{n}{2}\ln\left(\dfrac{\widehat{\sigma^{2}_{n}}}{\sigma^{2}}\right)-\dfrac{1}{2\widehat{\sigma^{2}_{n}}}\sum^{n}_{j=1}\left(x_{j}-\widehat{m_{n}}\right)^{2}+\dfrac{1}{2\sigma^{2}}\sum^{n}_{j=1}\left(x_{j}-m\right)^{2}\\
 & \qquad=\dfrac{n}{2}\left[\left(\dfrac{\widehat{\sigma^{2}_{n}}}{\sigma^{2}}-1\right)-\ln\left(\dfrac{\widehat{\sigma^{2}_{n}}}{\sigma^{2}}\right)\right]+n\dfrac{\left(\widehat{m_{n}}-m\right)^{2}}{2\sigma^{2}}\geqslant0.
\end{align*}
Since we have, for every $x>0,$ $x-1>\ln x.$

Hence, there exists a unique maximum likelihood estimator $\left(M_{n},\Sigma_{n}\right),$
where 
\[
M_{n}=\dfrac{1}{n}\sum^{n}_{j=1}x_{j}\,\,\,\,\text{and}\,\,\,\,\Sigma_{n}=\dfrac{1}{n}\sum^{n}_{j=1}\left(X_{j}-M_{n}\right)^{2}
\]
 are the empirical mean and variance of the sample.

\begin{remark}{}{}

We have 
\[
E_{\theta}\left(\Sigma_{n}\right)=\dfrac{n-1}{n}\sigma^{2}.
\]
 We say that $\Sigma_{n}$ is a \textbf{biased} \textbf{estimator\mindex{variance!biased estimator}}
of $\sigma^{2}.$

\end{remark}

\paragraph*{Uniform Case}

$\Theta=\mathbb{R}^{+\ast},$ and for every $\theta>0,$ $\mu_{\theta}$
is the uniform law $\mathscr{U}\left(\left[0,\theta\right]\right)$
on the interval $\left[0,\theta\right].$ For every $\uuline{x_{n}}\in\mathbb{R}^{n},$
\[
L_{n}\left(\uuline{x_{n}},\theta\right)=\dfrac{1}{\theta^{n}}\prod^{n}_{j=1}\boldsymbol{1}_{\left[0,\theta\right]}\left(x_{j}\right).
\]
The function $L_{n}\left(\uuline{x_{n}},\cdot\right)$ is not differentiable.
However, $L_{n}\left(\uuline{x_{n}},\cdot\right)$ is maximal when
$\theta$ is as small as possible, that is, at 
\[
\widehat{\varphi}_{n}\left(\uuline{x_{n}}\right)=\max_{1\leqslant j\leqslant n}\left(\uuline{x_{j}}\right).
\]
In this case, there exists a unique maximum likelihood estimator of
$\theta,$ 
\[
T_{n}=\max_{1\leqslant j\leqslant n}\left(\uuline{X_{j}}\right).
\]

In fact, the asymptotic properties are properties of solutions of
the likelihood equation---which provide stationary points---and
not necessarily properties of the solutions that achieve the maximum
of $L_{n}\left(\uuline{x_{n}},\cdot\right).$ For instance, the following
theorem states, under very strong hypotheses, a result of asymptotic
normality in the case where the parameter is a real number. This theorem
can be generalized to the case of a multidimensional parameter.

\begin{theorem}{Asympotic Normality For Real Parameter}{}

Let $\theta_{0}\in\Theta$ be the true value of the parameter. Assume
that $f\left(x,\theta\right)>0$ for every $\left(x,\theta\right)\in\mathbb{R}\times\Theta,$
where $\Theta$ is an open subset of $\mathbb{R}.$ Assume that, for
every real number $x,$ $f\left(x,\cdot\right)$ is twice continuously
differentiable, and that the function $\dfrac{\partial^{2}}{\partial\theta^{2}}\ln f\left(x,\cdot\right)$
is continuous in $\theta,$ uniformly in $x.$ 

Moreover, assume that there exists a $\mu-$integrable function $g$
and such that
\begin{equation}
\forall\left(x,\theta\right)\in\mathbb{R}\times\Theta,\,\,\,\,\left|\dfrac{\partial^{2}}{\partial\theta^{2}}f\left(x,\theta\right)\right|\leqslant g\left(x\right),\label{eq:upper_bound_second_part_deriv_f}
\end{equation}
and a function $h$ such that the function $x\mapsto h\left(x\right)f\left(x,\theta_{0}\right)$
is $\mu-$integrable and
\begin{equation}
\forall\left(x,\theta\right)\in\mathbb{R}\times\Theta,\,\,\,\,\left|\dfrac{\partial^{2}}{\partial\theta^{2}}\ln f\left(x,\theta\right)\right|\leqslant h\left(x\right).\label{eq:upper_bound_second_part_deriv_lnf}
\end{equation}

Define
\[
I\left(\theta_{0}\right)=-\intop_{\mathbb{R}}\left(\dfrac{\partial^{2}}{\partial\theta^{2}}\ln f\left(x,\theta_{0}\right)\right)f\left(x,\theta_{0}\right)\text{d}\mu\left(x\right).
\]
The quantity $I\left(\theta_{0}\right)$ is called the \textbf{\index{Fisher information}Fisher
information}. We have $0\leqslant I\left(\theta_{0}\right)<+\infty$
and we suppose that $I\left(\theta_{0}\right)>0.$

Let $\left(\widehat{\varphi}_{n}\right)_{n\in\mathbb{N}^{\ast}}$
be a sequence of solutions of the likelihood equation $\refpar{eq:likelihood_eq},$
and define $T_{n}=\widehat{\varphi}_{n}\left(\uuline{X_{n}}\right).$
If this sequence of estimators $\left(T_{n}\right)_{n\in\mathbb{N}^{\ast}}$
converges $P_{\theta_{0}}-$almost surely to $\theta_{0},$ then the
sequence with general term $Y_{n}=\sqrt{nI\left(\theta_{0}\right)}\left(T_{n}-\theta_{0}\right)$
converges in law to the Gaussian law $\mathscr{N}_{\mathbb{R}}\left(0,1\right).$

\end{theorem}

\begin{proof}{}{}

Let $\varphi$ be the function defined, for every $\left(x,\theta\right),$
by 
\[
\varphi\left(x,\theta\right)=\dfrac{\partial}{\partial\theta}\ln f\left(x,\theta\right)
\]
and let $U_{\theta}$ be the random variable $U_{\theta}=\varphi\left(X,\theta\right).$ 

The condition $\refpar{eq:upper_bound_second_part_deriv_lnf}$ implies
that
\begin{equation}
\intop_{\mathbb{R}}\left|\left[\dfrac{\partial^{2}}{\partial\theta^{2}}\ln f\left(x,\theta\right)\right]_{\theta=\theta_{0}}\right|f\left(x,\theta_{0}\right)\text{d}\mu\left(x\right)<+\infty.\label{eq:reg_model}
\end{equation}
In this case, the model is said to be \textbf{regular\mindex{regular!model}\mindex{model!regular}}
in $\theta_{0}.$ Then the transfer theorem and $\refpar{eq:reg_model}$
ensure that $U_{\theta_{0}}$ is square $P_{\theta_{0}}-$integrable
and that $E_{\theta_{0}}\left(U_{\theta_{0}}\right)=0.$

Indeed,
\begin{align*}
\intop_{\Omega}\left(U_{\theta_{0}}\right)^{2}\text{d}P_{\theta_{0}} & =\intop_{\mathbb{R}}\left(\left[\dfrac{\partial}{\partial\theta}\ln f\left(x,\theta\right)\right]_{\theta=\theta_{0}}\right)^{2}f\left(x,\theta_{0}\right)\text{d}\mu\left(x\right)\\
 & =\intop_{\mathbb{R}}\left[\dfrac{\left[\dfrac{\partial}{\partial\theta}f\left(x,\theta\right)\right]_{\theta=\theta_{0}}}{f\left(x,\theta_{0}\right)}\right]^{2}f\left(x,\theta_{0}\right)\text{d}\mu\left(x\right).
\end{align*}
Since
\[
\dfrac{\partial^{2}}{\partial\theta^{2}}\ln f\left(x,\theta\right)=\dfrac{\dfrac{\partial^{2}}{\partial\theta^{2}}f\left(x,\theta\right)f\left(x,\theta\right)-\left[\dfrac{\partial}{\partial\theta}f\left(x,\theta\right)\right]^{2}}{\left[f\left(x,\theta\right)\right]^{2}},
\]
it follows, taking into account the conditions $\refpar{eq:upper_bound_second_part_deriv_f}$
and $\refpar{eq:upper_bound_second_part_deriv_lnf},$
\begin{multline}
\intop_{\Omega}\left(U_{\theta_{0}}\right)^{2}\text{d}P_{\theta_{0}}\\
=-\intop_{\mathbb{R}}\left[\dfrac{\partial^{2}}{\partial\theta^{2}}\ln f\left(x,\theta\right)\right]_{\theta=\theta_{0}}f\left(x,\theta_{0}\right)\text{d}\mu\left(x\right)+\intop_{\mathbb{R}}\left[\dfrac{\partial^{2}}{\partial\theta^{2}}f\left(x,\theta\right)\right]_{\theta=\theta_{0}}\text{d}\mu\left(x\right)<+\infty.\label{eq:int_u_thetao_sq}
\end{multline}
Hence
\[
\intop_{\Omega}\left|U_{\theta_{0}}\right|\text{d}P_{\theta_{0}}<+\infty,
\]
 and, by the transfer theorem,
\[
\intop_{\Omega}\left|U_{\theta_{0}}\right|\text{d}P_{\theta_{0}}=\intop_{\mathbb{R}}\left|\left[\dfrac{\partial}{\partial\theta}\ln f\left(x,\theta\right)\right]_{\theta=\theta_{0}}\right|f\left(x,\theta_{0}\right)\text{d}\mu\left(x\right)<+\infty.
\]
Moreover
\[
E_{\theta_{0}}\left(U_{\theta_{0}}\right)=\intop_{\mathbb{R}}\dfrac{\left[\dfrac{\partial}{\partial\theta}f\left(x,\theta\right)\right]_{\theta=\theta_{0}}}{f\left(x,\theta_{0}\right)}f\left(x,\theta_{0}\right)\text{d}\mu\left(x\right)=\intop_{\mathbb{R}}\dfrac{\partial}{\partial\theta}f\left(x,\theta_{0}\right)\text{d}\mu\left(x\right).
\]
Let $K$ be a compact neighborhood of $\theta_{0}$ contained in $\theta.$
The finite increments theorem and the condition $\refpar{eq:upper_bound_second_part_deriv_f}$
ensure that, for every $\theta\in K,$
\[
\left|\dfrac{\partial}{\partial\theta}f\left(x,\theta\right)\right|\leqslant\left|\dfrac{\partial}{\partial\theta}f\left(x,\theta_{0}\right)\right|+cg\left(x\right),
\]
where $c>0$ is a constant depending on $K.$ 

We may therefore apply the differentiation theorem for an integral
with parameter, which yields
\[
E_{\theta_{0}}\left(U_{\theta_{0}}\right)=\left[\dfrac{\partial}{\partial\theta}\intop_{\mathbb{R}}f\left(x,\theta\right)\text{d}\mu\left(x\right)\right]_{\theta=\theta_{0}}.
\]
Since, for every $\theta,$ 
\[
\intop_{\mathbb{R}}f\left(x,\theta\right)\text{d}\mu\left(x\right)=1,
\]
it follows that $E_{\theta_{0}}\left(U_{\theta_{0}}\right)=0.$

Similarly, by the condition $\refpar{eq:upper_bound_second_part_deriv_f},$
\[
\intop_{\mathbb{R}}\dfrac{\partial^{2}}{\partial\theta^{2}}f\left(x,\theta\right)\text{d}\mu\left(x\right)=\dfrac{\partial^{2}}{\partial\theta^{2}}\intop_{\mathbb{R}}f\left(x,\theta\right)\text{d}\mu\left(x\right)=0.
\]

It follows that, substituting into $\refpar{eq:int_u_thetao_sq},$
\[
E_{\theta_{0}}\left(U^{2}_{\theta_{0}}\right)=-\intop_{\mathbb{R}}\left[\dfrac{\partial^{2}}{\partial\theta^{2}}\ln f\left(x,\theta\right)\right]_{\theta=\theta_{0}}f\left(x,\theta_{0}\right)\text{d}\mu\left(x\right)=I\left(\theta_{0}\right),
\]
which shows, in particular, that $0\leqslant I\left(\theta_{0}\right)<+\infty.$

For every $n\in\mathbb{N}^{\ast},$ let $\widehat{\varphi}_{n}\left(\uuline{x_{n}}\right)$
be a solution of the log-likelihood equation
\[
\dfrac{\partial}{\partial\theta}\ln L_{n}\left(\uuline{x_{n}},\theta\right)=\sum^{n}_{j=1}\varphi\left(x_{j},\theta\right)=0.
\]
The first-order Taylor formula with integral remainder, applied at
$\widehat{\varphi}_{n}\left(\uuline{x_{n}}\right),$ gives
\[
\sum^{n}_{j=1}\varphi\left(x_{j},\theta_{0}\right)=\left(\theta_{0}-\widehat{\varphi}_{n}\left(\uuline{x_{n}}\right)\right)\intop^{1}_{0}\sum^{n}_{j=1}\dfrac{\partial}{\partial\theta}\varphi\left(x_{j},\widehat{\varphi}_{n}\left(\uuline{x_{n}}\right)+t\left(\theta_{0}-\widehat{\varphi}_{n}\left(\uuline{x_{n}}\right)\right)\right)\text{d}t,
\]
This yields
\[
\sum^{n}_{j=1}\varphi\left(X_{j},\theta_{0}\right)=\left(\theta_{0}-T_{n}\right)\intop^{1}_{0}\sum^{n}_{j=1}\dfrac{\partial}{\partial\theta}\varphi\left(X_{j},T_{n}+t\left(\theta_{0}-T_{n}\right)\right)\text{d}t,
\]
and therefore
\begin{equation}
\dfrac{1}{\sqrt{n}}\sum^{n}_{j=1}\varphi\left(X_{j},\theta_{0}\right)=\sqrt{n}\left(\theta_{0}-T_{n}\right)\intop^{1}_{0}\dfrac{1}{n}\sum^{n}_{j=1}\dfrac{\partial}{\partial\theta}\varphi\left(X_{j},T_{n}+t\left(\theta_{0}-T_{n}\right)\right)\text{d}t.\label{eq:eq15_20}
\end{equation}
The random variables $\varphi\left(X_{j},\theta_{0}\right)$ are independent
with same law under $P_{\theta_{0}},$ as the law of $U_{\theta_{0}}.$
In particular, they admit a second-order moment. The central limit
theorem then gives
\begin{equation}
\dfrac{1}{\sqrt{n}}\sum^{n}_{j=1}\varphi\left(X_{j},\theta_{0}\right)\stackrel[n\to+\infty]{\mathscr{L}}{\longrightarrow}\mathscr{N}_{\mathbb{R}}\left(0,I\left(\theta_{0}\right)\right).\label{eq:15.21}
\end{equation}

It remains to study the sequence with general term
\[
\intop^{1}_{0}\dfrac{1}{n}\sum^{n}_{j=1}\dfrac{\partial}{\partial\theta}\varphi\left(X_{j},T_{n}+t\left(\theta_{0}-T_{n}\right)\right)\text{d}t.
\]
By the condition $\refpar{eq:upper_bound_second_part_deriv_lnf},$
the random variables $\dfrac{\partial}{\partial\theta}\varphi\left(X_{j},\theta_{0}\right)$
admit an expectation under $P_{\theta_{0}},$ equal to $-I\left(\theta_{0}\right)$.
Moreover, they are independent. Hence, by the strong law of large
numbers,
\begin{equation}
\lim_{n\to+\infty}\dfrac{1}{n}\sum^{n}_{j=1}\dfrac{\partial}{\partial\theta}\varphi\left(X_{j},\theta_{0}\right)=-I\left(\theta_{0}\right)\,\,\,\,P_{\theta_{0}}-\text{almost surely.}\label{eq:lim15_22}
\end{equation}
Now, define
\[
A_{n}\left(t\right)=\dfrac{1}{n}\sum^{n}_{j=1}\left[\dfrac{\partial}{\partial\theta}\varphi\left(X_{j},T_{n}+t\left(\theta_{0}-T_{n}\right)\right)-\dfrac{\partial}{\partial\theta}\varphi\left(X_{j},\theta_{0}\right)\right].
\]
We will prove that $P_{\theta_{0}}-$almost surely, for every $t\in\left[0,1\right],$
$\lim_{n\to+\infty}A_{n}\left(t\right)=0$---with attention to the
order in which the assertions are stated.

Let $\epsilon>0.$ Since the function $\theta\mapsto\dfrac{\partial}{\partial\theta}\varphi\left(x,\theta\right)$
is continuous, uniformly in $x,$ there exists an interval $V,$ centered
at $\theta_{0},$ and contained in $\Theta,$ such that, for every
$\theta\in V,$
\begin{equation}
\sup_{x\in\mathbb{R}}\left|\dfrac{\partial}{\partial\theta}\varphi\left(x,\theta\right)-\dfrac{\partial}{\partial\theta}\varphi\left(x,\theta_{0}\right)\right|\leqslant\epsilon.\label{eq:15_23}
\end{equation}

By hypothesis, there exists $N\in\mathscr{A}$ such that $P_{\theta_{0}}\left(N\right)=0$
and such that, for every $\omega\notin N,$
\[
\lim_{n\to+\infty}T_{n}\left(\omega\right)=\theta_{0}.
\]
Fix such an $\omega,$ and let $K\left(\omega\right)$ be such that,
for every $n\geqslant K\left(\omega\right),$ $T_{n}\left(\omega\right)\in V.$
Then, for $0<t<1,$
\begin{align*}
\left|A_{n}\left(t\right)\left(\omega\right)\right| & \leqslant\dfrac{1}{n}\sum^{K\left(\omega\right)}_{j=1}\left|\dfrac{\partial}{\partial\theta}\varphi\left(X_{j}\left(\omega\right),T_{n}\left(\omega\right)+t\left(\theta_{0}-T_{n}\left(\omega\right)\right)\right)-\dfrac{\partial}{\partial\theta}\varphi\left(X_{j}\left(\omega\right),\theta_{0}\right)\right|\\
 & \qquad+\dfrac{1}{n}\sum^{n}_{j=K\left(\omega\right)+1}\left|\dfrac{\partial}{\partial\theta}\varphi\left(X_{j}\left(\omega\right),T_{n}\left(\omega\right)+t\left(\theta_{0}-T_{n}\left(\omega\right)\right)\right)-\dfrac{\partial}{\partial\theta}\varphi\left(X_{j}\left(\omega\right),\theta_{0}\right)\right|.
\end{align*}
Thus, taking $\refpar{eq:15_23}$ into account,
\[
\left|A_{n}\left(t\right)\left(\omega\right)\right|\leqslant\dfrac{2}{n}\sum^{K\left(\omega\right)}_{j=1}g\left(X_{j}\left(\omega\right)\right)+\epsilon.
\]
It follows that 
\[
\limsup_{n\to+\infty}\left|A_{n}\left(t\right)\left(\omega\right)\right|\leqslant\epsilon,
\]
and, since $\epsilon$ is arbitrary, this proves that, for every $t\in\left[0,1\right],$
\[
\lim_{n\to+\infty}A_{n}\left(t\right)\left(\omega\right)=0.
\]
 Since we have
\[
\left|A_{n}\left(t\right)\left(\omega\right)\right|\leqslant2g\left(x\right),
\]
it follows by the dominated convergence theorem that the sequence
with general term
\[
\intop^{1}_{0}\dfrac{1}{n}\sum^{n}_{j=1}\left[\dfrac{\partial}{\partial\theta}\varphi\left(X_{j}\left(\omega\right),T_{n}\left(\omega\right)+t\left(\theta_{0}-T_{n}\left(\omega\right)\right)\right)-\dfrac{\partial}{\partial\theta}\varphi\left(X_{j}\left(\omega\right),\theta_{0}\right)\right]\text{d}t
\]
converges to 0. Since this holds for every $\omega\notin\mathbb{N},$
it follows from $\refpar{eq:lim15_22}$ that
\begin{align*}
\lim_{n\to+\infty}\intop^{1}_{0}\dfrac{1}{n}\sum^{n}_{j=1}\dfrac{\partial}{\partial\theta}\varphi\left(X_{j},T_{n}+t\left(\theta_{0}-T_{n}\right)\right)\text{d}t & =\lim_{n\to+\infty}\dfrac{1}{n}\sum^{n}_{j=1}\dfrac{\partial}{\partial\theta}\varphi\left(X_{j},\theta_{0}\right)\\
 & =-I\left(\theta_{0}\right)\,\,\,\,P_{\theta_{0}}-\text{almost surely.}
\end{align*}
Since $I\left(\theta_{0}\right)\neq0,$
\[
\lim_{n\to+\infty}\left[\dfrac{1}{n}\sum^{n}_{j=1}\dfrac{\partial}{\partial\theta}\varphi\left(X_{j},\theta_{0}\right)\right]^{-1}=-\left[I\left(\theta_{0}\right)\right]^{-1}\,\,\,\,P_{\theta_{0}}-\text{almost surely,}
\]
and this convergence also holds in law. 

Finally, the equality $\refpar{eq:eq15_20}$ and the convergence in
law established in $\refpar{eq:lim15_22}$ implies that, by the Slutsky
lemma---see Exercise $\ref{exo:exercise15.8}$---that
\[
\sqrt{nI\left(\theta_{0}\right)}\left(T_{n}-\theta_{0}\right)\stackrel[n\to+\infty]{\mathscr{L}}{\longrightarrow}\mathscr{N}_{\mathbb{R}}\left(0,1\right).
\]

\end{proof}

\section*{Exercises}

\addcontentsline{toc}{section}{Exercises}

\begin{exercise}{Narrow Convergence of a Sequence of Probabilities Carried by $\mathbb{Z}$}{exercise15.1}

For every $n\in\mathbb{N},$ let $\mu_{n}$ be a probability measure
on the measurable space $\left(\mathbb{R},\mathscr{B}_{\mathbb{R}}\right),$
carried by $\mathbb{Z}.$ Thus $\mu_{n}$ is of the form 
\[
\mu_{n}=\sum_{r\in\mathbb{Z}}a^{n}_{r}\delta_{r},
\]
where for every $r\in\mathbb{Z},$ $a^{n}_{r}\geqslant0.$ 

Prove that the sequence $\left(\mu_{n}\right)_{n\in\mathbb{N}}$ converges
to a measure $\mu$ if and only, for every $r\in\mathbb{Z},$ the
sequence $\left(a^{n}_{r}\right)_{n\in\mathbb{N}}$ converges to a
real number $a_{r}\geqslant0,$ and we have $\sum_{r\in\mathbb{Z}}a_{r}=1$
and $\mu=\sum_{r\in\mathbb{Z}}a_{r}\delta_{r}$---$\mu$ is carried
by $\mathbb{Z}.$

\end{exercise}

\begin{exercise}{Binomial Approximation of the Hypergeometric Law}{exercise15.2}

Let $j\in\mathbb{N}^{\ast}.$ Consider a finite set $U^{j}=U^{j}_{1}\uplus U^{j}_{2}$
partitioned into two nonempty subsets $U^{j}_{1}$ and $U^{j}_{2}.$ 

Denote $\left|U^{j}\right|=r^{j}$ and $\left|U^{j}_{1}\right|=r^{1}_{j}$---and
thus $\left|U_{2}\right|=r^{j}-r^{1}_{j}\geqslant1.$ 

Let $n$ be an integer such that $1\leqslant n<r^{j}.$ We extract
``at random'', that is, in an uniform manner, $n$ elements of $U^{j}.$
Determine, for every $k$ such that $0\leqslant k\leqslant n,$ the
probability to obtain exactly $k$ elements of $U^{j}_{1}$---and
thus $n-k$ elements of $U^{j}_{2}.$

Moreover, suppose that the two sequences of integers $\left(r^{1}_{j}\right)_{j\in\mathbb{N}^{\ast}}$
and $\left(r_{j}\right)_{j\in\mathbb{N}^{\ast}}$ are nondecreasing
and tend to infinity with $j$ such that $\dfrac{r^{1}_{j}}{r_{j}}\underset{j\to+\infty}{\longrightarrow}p,$
where $p\in\left]0,1\right[.$ 

Let us consider a couple of integers such that $0\leqslant k\leqslant n.$
Prove that there exists an integer $j_{0}$ such that, for every $j\geqslant j_{0},$
$n\leqslant r^{j}-r^{1}_{j}$ and $n\leqslant r^{1}_{j}.$

For $j\geqslant j_{0},$ define
\[
P_{k,n}\left(r^{1}_{j},r_{j}\right)=\dfrac{\binom{r^{i}_{j}}{k}\binom{r^{j}-r^{1}_{j}}{n-k}}{\binom{r^{j}}{n}}.
\]
Prove that
\[
P_{k,n}\left(r^{1}_{j},r_{j}\right)\underset{j\to+\infty}{\longrightarrow}\binom{n}{k}p^{k}\left(1-p\right)^{n-k}.
\]

Interpret this result in terms of narrow convergence of a sequence
of probabilities.

\textit{Hint: use the previous exercise.}

\end{exercise}

\begin{exercise}{Geometric and Exponential Laws}{exercise15.3}

Let $X$ be a nonnegative real-valued random variable. For every $a>0,$
define the random variables $V_{a}=\left\lfloor \dfrac{X}{a}\right\rfloor $
and $X_{a}=a\left\lfloor \dfrac{X}{a}\right\rfloor .$

1. If $X$ follows the exponential law $\text{e}^{\lambda},$ where
$\lambda>0,$ determine the law of $V_{a}.$

2. Suppose that for every $a>0,$ the law of $V_{a}$ is the geometric
law on $\mathbb{N}$ with parameter $1-\text{e}^{-\lambda a}.$ Compute
the cumulative distribution function of $X_{a}.$ Study the narrow
convergence of the family of laws of $X_{a}$ when $a$ tends to $0.$ 

\end{exercise}

\begin{exercise}{Narrow Convergence of a Sequence of Gaussian Probabilities on $\mathbb{R}$}{exercise15.4}

On the measurable space $\left(\mathbb{R},\mathscr{B}_{\mathbb{R}}\right),$
consider the sequence of Gaussian measures $\mu_{n}=f_{n}\cdot\lambda,\,n\in\mathbb{N},$
where $\lambda$ is the Lebesgue measure and $f_{n}$ is the density
defined, for every real number $x,$ by
\[
f_{n}\left(x\right)=\dfrac{1}{\sigma_{n}\sqrt{2\pi}}\text{e}^{-\frac{\left(x-m_{n}\right)^{2}}{2\sigma^{2}_{n}}},
\]
where $m_{n}$ is an arbitrary real number and $\sigma_{n}$ is a
positive real number.

1. If the sequences $\left(m_{n}\right)_{n\in\mathbb{N}}$ and $\left(\sigma_{n}\right)_{n\in\mathbb{N}}$
converge respectively to $m$ and $\sigma,$ study the narrow convergence
of the sequence $\left(\mu_{n}\right)_{n\in\mathbb{N}}$ directly
from the definition of this notion of convergence. What does the Scheffé
lemma add in the case where $\sigma>0?$

2. If the sequence $\left(m_{n}\right)_{n\in\mathbb{N}}$ is bounded
and if the sequence $\left(\sigma_{n}\right)_{n\in\mathbb{N}}$ converges
to $+\infty$ as $n$ tends to infinity, study the weak and narrow
convergences of the sequence $\left(\mu_{n}\right)_{n\in\mathbb{N}}.$

\end{exercise}

\begin{exercise}{Gaussian Random Variables and Convergence in Law}{exercise15.5}

\textit{Hint: You may use the results from the previous exercise.}

Let $\left\{ X_{0};Z_{n},\,n\in\mathbb{N}\right\} $ be a family of
independent Gaussian real-valued random variables defined on a probabilized
space $\left(\Omega,\mathscr{A},P\right).$ Assume that the $Z_{n}$
follow all the same law $\mathscr{N}_{\mathbb{R}}\left(0,\sigma^{2}\right),$
where $\sigma>0.$ Let $\rho$ be a nonzero real number. For every
$n\in\mathbb{N}^{\ast},$ define the random variable
\[
X_{n}=\rho X_{n-1}+Z_{n}.
\]

1. Prove that $X_{n}$ admits a second-order moment, and compute its
expectation and variance.

2. Study the convergence in law of the sequence $\left(X_{n}\right)_{n\in\mathbb{N}}.$

\end{exercise}

\begin{exercise}{Convergence in Law}{exercise15.6}

On the probabilized space $\left(\Omega,\mathscr{A},P\right),$ let,
for each $n\in\mathbb{N}^{\ast},$ $X_{n}$ and $Y_{n}$ be random
variables. Assume that the $X_{n},\,n\in\mathbb{N}^{\ast},$ all follow
the Gaussian law $\mathscr{N}_{\mathbb{R}}\left(0,1\right)$ and that
the law of $Y_{n}$ is 
\[
P_{Y_{n}}=\left(1-\dfrac{1}{n}\right)\delta_{1}+\dfrac{1}{n}\delta_{0}.
\]
Study the convergence in law of the sequence $\left(X_{n}Y_{n}\right)_{n\in\mathbb{N}^{\ast}}.$

\end{exercise}

\begin{exercise}{Convergence in Law of a Sequence of Random Variables Taking Values in $\mathbb{R}^2$ and of the Sequence of Its Marginals}{exercise15.7}

On the same probabilized space \preds, let $\left(X_{n}\right)_{n\in\mathbb{N}}$
and $\left(Y_{n}\right)_{n\in\mathbb{N}}$ be two sequences of real-valued
random variables that converge in law respectively to the random variables
$X$ and $Y,$ assumed to be independent.

1. Assume that, for every $n\in\mathbb{N},$ the random variables
$X_{n}$ and $Y_{n}$ are independent. Prove that the sequence of
random variables $\left(X_{n},Y_{n}\right),\,n\in\mathbb{N},$ converges
in law to $\left(X,Y\right).$ In particular, deduce that the sequence
of random variables $X_{n}+Y_{n},\,n\in\mathbb{N},$ converges in
law to $X+Y.$

2. \textbf{We now study a counter-example, showing what may fail if
we drop the assumption ``for every $n\in\mathbb{N},$ $X_{n}$ and
$Y_{n}$ are independent''.} Let $X$ and $Y$ be two independent
real-valued random variables following the same Bernoulli law $\dfrac{\delta_{0}+\delta_{1}}{2}.$ 

For every $n\in\mathbb{N}^{\ast},$ define
\[
X_{n}=X+\dfrac{1}{n}\,\,\,\,\text{and}\,\,\,\,Y_{n}=\left(1-X\right)-\dfrac{1}{n}.
\]

Study the convergence in law of the three sequences $\left(X_{n}\right)_{n\in\mathbb{N}^{\ast}},$
$\left(Y_{n}\right)_{n\in\mathbb{N}^{\ast}}$ and $\left(X_{n}+Y_{n}\right)_{n\in\mathbb{N}^{\ast}}.$
Conclude that the sequence of random variables $\left(X_{n},Y_{n}\right)_{n\in\mathbb{N}^{\ast}}$
does not converge in law to $\left(X,Y\right).$

\end{exercise}

The Slutsky lemma gives an alternative assumption to independence
that guarantees the convergence in law property studied in the first
question of the previous exercise.

\begin{exercise}{Slutsky Lemma}{exercise15.8}

On the same probabilized space \preds, let $\left(X_{n}\right)_{n\in\mathbb{N}}$
and $\left(Y_{n}\right)_{n\in\mathbb{N}}$ be two sequences of real-valued
random variables that converge in law respectively to the independent
random variables $X$ and to a constant $y_{0}.$

1. (a) Prove that the sequence of random variables $\left(X_{n},Y_{n}\right),\,n\in\mathbb{N}$
converges in law to $\left(X,y_{0}\right)$

\textit{Hint: You may admit that the set
\[
\mathscr{H}=\left\{ \left(x,y\right)\mapsto f\left(x\right)g\left(y\right):\,f,g\in\mathscr{C}_{0}\left(\mathbb{R}\right)\right\} 
\]
is total in $\mathscr{C}_{0}\left(\mathbb{R}^{2}\right)$ or, alternatively
use the Lévy theorem. }

(b) Deduce in particular that the sequence of random variables $X_{n}+Y_{n},\,n\in\mathbb{N},$
converges in law to $X+y_{0}.$

2. Prove that if the sequence $\left(X_{n}\right)_{n\in\mathbb{N}}$
converges in law to a random variable $X,$ and if the sequence $\left(X_{n}-Y_{n}\right)_{n\in\mathbb{N}}$
converges in probability to 0, then the sequence $\left(Y_{n}\right)_{n\in\mathbb{N}}$
converges in law to $X.$

\end{exercise}

\begin{exercise}{Decimal Expansion. Convergence in Law and Lévy Theorem}{exercise15.9}

Let $\left(X_{n}\right)_{n\in\mathbb{N}}$ be a sequence of independent
real-valued random variables defined on the same probabilized space
$\left(\Omega,\mathscr{A},P\right).$ Each $X_{n}$ follows the uniform
law on the set of integers $\left\llbracket 0,9\right\rrbracket .$
For every $n\in\mathbb{N},$ define the random variable 
\[
Y_{n}=\sum^{n}_{j=0}\dfrac{X_{j}}{10^{j}}.
\]

Prove that the sequence $\left(Y_{n}\right)_{n\in\mathbb{N}}$ converges
$P-$almost surely to a random variable $Y,$ and determine its law.

\end{exercise}

\begin{exercise}{Convergence in Law and Cumulative Distribution Functions}{exercise15.10}

Let $\left(X_{n}\right)_{n\in\mathbb{N}^{\ast}}$ be a sequence of
independent real-valued random variables defined on the same probabilized
space $\left(\Omega,\mathscr{A},P\right),$ following the same law,
with cumulative distribution function $F.$ Define for $n\in\mathbb{N}^{\ast}$
the random variables $I_{n}$ and $M_{n}$ by
\[
I_{n}=\min_{1\leqslant i\leqslant n}X_{i}\,\,\,\,\text{and}\,\,\,\,M_{n}=\max_{1\leqslant i\leqslant n}X_{i}.
\]

1. Study the convergence in law of the sequences $\left(I_{n}\right)_{n\in\mathbb{N}^{\ast}}$
and $\left(M_{n}\right)_{n\in\mathbb{N}^{\ast}}.$

2. Suppose that the $X_{n},\,n\in\mathbb{N}^{\ast}$ follow the same
exponential law $\text{e}^{\lambda}$ where $\lambda>0.$ For every
$n\in\mathbb{N}^{\ast},$ define $Z_{n}=\dfrac{M_{n}}{\ln n}.$ Study
the convergence in law of the sequence $\left(Z_{n}\right)_{n\in\mathbb{N}^{\ast}}.$

\end{exercise}

\begin{exercise}{Integrale Inequality for the Real Part of a Characteristic Function. Convergence in Law of a Series of Independent Random Variables---Lévy Theorem}{exercise15.11}

All random variables are defined on the same probabilized space $\left(\Omega,\mathscr{A},P\right).$

Let $X$ be a real-valued random variable with characteristic function
$\varphi_{X}.$ Let $g$ be the real-valued function defined on $\mathbb{R}$
by
\[
g\left(x\right)=\begin{cases}
1-\frac{\sin x}{x}, & \text{if }x\neq0,\\
0, & \text{if }x=0.
\end{cases}
\]

1. Verify that $g\in\mathscr{C}_{b}\left(\mathbb{R}\right)$ is nonnegative,
and that $g\left(x\right)=0$ if and only if $x=0.$

Let $\delta>0.$ Prove the equality
\begin{equation}
\dfrac{1}{\delta}\intop^{\delta}_{0}\left(1-\text{Re}\left(\varphi_{X}\left(t\right)\right)\right)\text{d}t=\intop_{\Omega}g\left(\delta X\right)\text{d}P.\label{eq:int_equality_re_varphi_X_g_deltaX}
\end{equation}

2. For every $\epsilon>0,$ define $I_{\epsilon}=\inf_{\left|x\right|>\epsilon}g\left(x\right)>0.$ 

Prove that
\begin{equation}
P\left(\left|X\right|>\epsilon\right)\leqslant\dfrac{1}{I_{\epsilon}\delta}\intop^{\delta}_{0}\left(1-\text{Re}\left(\varphi_{X}\left(t\right)\right)\right)\text{d}t=\dfrac{1}{2I_{\epsilon}\delta}\intop^{\delta}_{-\delta}\left(1-\varphi_{X}\left(t\right)\right)\text{d}t.\label{eq:P_abs_X_over_eps_upper_bound}
\end{equation}

Let $\left(X_{n}\right)_{n\in\mathbb{N}^{\ast}}$ be a sequence of
real-valued random variables. Set $S_{n}=\sum^{n}_{j=1}X_{j}.$

3. Prove that the sequence $\left(X_{n}\right)_{n\in\mathbb{N}^{\ast}}$
converges in law to 0---and hence in probability to 0---if and only
if there exists $\delta>0$ such that the sequence $\left(\varphi_{X_{n}}\left(t\right)\right)_{n\in\mathbb{N}^{\ast}}$
converges to $1$ for every $t\in\left[-\delta,\delta\right].$

4. \textbf{Suppose that the random variables $X_{n},\,n\in\mathbb{N}^{\ast},$
are independent}. Prove that the sequence $\left(S_{n}\right)_{n\in\mathbb{N}^{\ast}}$
converges in law if and only if it converges in probability---\textbf{Lévy}
theorem.

\end{exercise}

\begin{exercise}{Gaussian Random Variables. Conditional Laws. Characteristic Functions and Convergence in Law}{exercise15.12}

We denote by $\uuline{x_{n}}$ a vector $\left(x_{1},x_{2},\cdots,x_{n}\right)$
in $\mathbb{R}^{n}.$ Let $\left(X_{n}\right)_{n\in\mathbb{N}^{\ast}}$
be a sequence of real-valued random variables defined on the same
probabilized space $\left(\Omega,\mathscr{A},P\right).$ Assume that
$X_{1}$ follows a Gaussian law $\mathscr{N}_{\mathbb{R}}\left(0,1\right),$
and that, for every $n>1,$ a conditional law $P^{\uuline{X_{n}}=\uuline{x_{n}}}_{X_{n+1}}$
of $X_{n+1}$ given $\uuline{X_{n}}$ is, for every $\uuline{x_{n}}\in\mathbb{R}^{n},$
the Gaussian law $\mathscr{N}_{\mathbb{R}}\left(x_{n},1\right).$

1. What is the law of $\left(X_{1},X_{2}\right)?$ Find, up to a multiplicative
factor, a linear combination of $X_{1}$ and $X_{2}$ that is independent
of $X_{1}.$

2. Let $\mathscr{B}_{n}$ be the $\sigma-$algebra generated by $\uuline{X_{n}}.$
Compute the conditional expectations $\mathbb{E}^{\mathscr{B}_{n}}\left(X_{n+1}\right)$
and $\mathbb{E}^{\mathscr{B}_{n}}\left(X^{2}_{n+1}\right).$ Deduce
from this the expectation and variance of $X_{n}.$ Prove that the
sequence $\left(X_{n}\right)_{n\in\mathbb{N}^{\ast}}$ does not converge
in $\text{L}^{2}.$

3. Justify the existence of a density $f_{\uuline{X_{n}}}$ for the
random variable $\uuline{X_{n}}$ and compute it---start with the
case $n=3.$ What is the characteristic function of $\uuline{X_{n}}?$

4. Let $j<k.$ What is the law of the random variable $\left(X_{j},X_{k}\right)?$
What is the correlation coefficient of $X_{j}$ and $X_{k}?$ Study
the convergence in law of the sequence of random variables $\left(X_{j},\dfrac{X_{k}}{\sqrt{k}}\right)_{k\in\mathbb{N}^{\ast}}.$
What can be said about the limiting law?

5. For every $n\in\mathbb{N}^{\ast},$ consider the random variable
\[
Z_{n}=\dfrac{1}{n\sqrt{n}}\sum^{n}_{j=1}X_{j}.
\]
Study the convergence in law of the sequence of random variables $\left(Z_{n}\right)_{n\in\mathbb{N}^{\ast}}.$

\end{exercise}

\section*{Solutions of Exercises}

\addcontentsline{toc}{section}{Solutions of Exercises}

\begin{solution}{}{solexercise15.1}

Assume that the sequence $\left(\mu_{n}\right)_{n\in\mathbb{N}}$
converges narrowly to a probability $\mu.$ Then, for every open interval
$\left]r-1,r\right[,$ where $r\in\mathbb{Z},$
\[
0\leqslant\mu\left(\left]r-1,r\right[\right)\leqslant\liminf_{n\to+\infty}\mu_{n}\left(\left]r-1,r\right[\right).
\]

Since, for every $n\in\mathbb{N},$ $\mu_{n}\left(\left]r-1,r\right[\right)=0,$
it follows that 
\[
\mu\left(\left]r-1,r\right[\right)=0.
\]
The probability $\mu$ is thus supported by a subset of $\mathbb{Z},$
and is of the form $\mu=\sum_{r\in\mathbb{Z}}a_{r}\delta_{r},$ where
$a_{r}\geqslant0.$ 

For every function $f\in\mathscr{C}_{\mathscr{K}}\left(\mathbb{R}\right)$
with support in the interval $\left]r-\frac{1}{2},r+\frac{1}{2}\right[$
such that $f\left(r\right)\neq0,$ then
\[
\intop_{\mathbb{R}}f\text{d}\mu_{n}=f\left(r\right)\mu_{n}\left(\left\{ r\right\} \right)=f\left(r\right)a^{n}_{r}\,\,\,\,\text{and}\,\,\,\,\intop_{\mathbb{R}}f\text{d}\mu=f\left(r\right)\mu\left(\left\{ r\right\} \right).
\]
Since 
\[
\lim_{n\to+\infty}\intop_{\mathbb{R}}f\text{d}\mu_{n}=\intop_{\mathbb{R}}f\text{d}\mu,
\]
 it follows that
\[
\lim_{n\to+\infty}a^{n}_{r}=\mu\left(\left\{ r\right\} \right)\geqslant0.
\]
Finally, $\mu$ being by hypothesis a probability, we have
\[
\sum_{r\in\mathbb{Z}}\mu\left(\left\{ r\right\} \right)=\sum_{r\in\mathbb{Z}}a_{r}=1.
\]

Conversely, suppose that, for every $r\in\mathbb{Z},$ the sequence
$\left(a^{n}_{r}\right)_{n\in\mathbb{N}}$ converges to a real number
$a_{r}\geqslant0,$ and that 
\[
\sum_{r\in\mathbb{Z}}a_{r}=1\,\,\,\,\text{and}\,\,\,\,\mu=\sum_{r\in\mathbb{Z}}a_{r}\delta_{r}.
\]
For every $f\in\mathscr{C}_{\mathscr{K}}\left(\mathbb{R}\right)$
with compact support $K,$
\[
\intop_{\mathbb{R}}f\text{d}\mu_{n}=\sum_{r\in K}f\left(r\right)a^{n}_{r}\,\,\,\,\text{and}\,\,\,\,\intop_{\mathbb{R}}f\text{d}\mu_{n}=\sum_{r\in K}f\left(r\right)a_{r},
\]
the sums having only a finite number of terms. 

It follows that
\[
\lim_{n\to+\infty}\intop_{\mathbb{R}}f\text{d}\mu_{n}=\intop_{\mathbb{R}}f\text{d}\mu,
\]
which shows that the vague convergence, and thus the narrow convergence
of the sequence of probabilities $\mu_{n}$ to the probability $\mu.$ 

\end{solution}

\begin{solution}{}{solexercise15.2}

An outcome is a subset of $U^{j}$ with $n$ elements. We choose for
the set of outcomes
\[
\Omega^{j}=\left\{ A\in\mathscr{P}\left(U^{j}\right):\,\left|A\right|=n\right\} .
\]
The event under study is the subset $A_{k}$ of $\Omega^{j},$
\[
A_{k}=\left\{ A\in\Omega^{j}:\,\left|A\cap U^{j}_{1}\right|=k\right\} .
\]
On the probabilizable space $\left(\Omega^{j},\mathscr{A}^{j}\right),$
where $\mathscr{A}^{j}=P\left(\Omega^{j}\right),$ we consider the
uniform probability $P^{j}$---this translates the usual meaning
of the expression ``at random''. We are therefore looking for the
probability $P^{j}\left(A_{k}\right).$ 

The set $A_{k}$ is empty if and only if $r^{1}_{j}<k\leqslant n$
or if $0\leqslant k<n-\left(r^{j}-r^{1}_{j}\right).$

Otherwise, that is, if 
\[
\max\left(0,n-\left(r^{j}-r^{1}_{j}\right)\right)\leqslant k\leqslant\min\left(n,r^{1}_{j}\right),
\]
then
\[
\left|A_{k}\right|=\binom{r^{1}_{j}}{k}\binom{r^{j}-r^{1}_{j}}{n-k}.
\]
Moreover,
\[
\left|\Omega^{j}\right|=\binom{r^{j}}{n}.
\]
It follows that\boxeq{
\[
P^{j}\left(A_{k}\right)=\dfrac{\binom{r^{1}_{j}}{k}\binom{r^{j}-r^{1}_{j}}{n-k}}{\binom{r^{j}}{n}}.
\]
}

The events $A_{k},\,0\leqslant k\leqslant n$ form a partition of
$\Omega^{j}.$ The measure
\[
\mu^{j}=\sum^{n}_{k=0}P^{J}\left(A_{k}\right)\delta_{k}=\sum^{\min\left(n,r^{1}_{j}\right)}_{k=\max\left(0,n-\left(r^{j}-r^{1}_{j}\right)\right)}\dfrac{\binom{r^{1}_{j}}{k}\binom{r^{j}-r^{1}_{j}}{n-k}}{\binom{r^{j}}{n}}\delta_{k}
\]
is therefore a probability measure, called \textbf{hypergeometric
law}\index{hypergeometric law}\mindex{law!hypergeometric}.

Let $\epsilon>0$ such that $\epsilon<\min\left(p,1-p\right).$ There
exists $j_{1}$ such that, for every $j\geqslant j_{1},$ 
\[
p-\epsilon\leqslant\dfrac{r^{1}_{j}}{r^{j}}\leqslant p+\epsilon,
\]
 which implies 
\[
r^{j}\left(1-p-\epsilon\right)\leqslant r_{j}-r^{1}_{j}\leqslant r^{j}\left(1-p+\epsilon\right).
\]
Since the sequences $\left(r^{1}_{j}\right)_{j\in\mathbb{N}^{\ast}}$
and $\left(r_{j}\right)_{j\in\mathbb{N}^{\ast}}$ tend to infinity
with $j,$ then there exists $j_{0}$ such that, for every $j\geqslant j_{0},$
we have $n\leqslant r^{j}-r^{1}_{j}$ and $n\leqslant r^{1}_{j}.$
For such $j,$
\begin{equation}
\mu^{j}=\sum^{n}_{k=0}P^{j}\left(A_{k}\right)\delta_{k}=\sum^{n}_{k=0}P_{k,n}\left(r^{1}_{j},r_{j}\right)\delta_{k}.\label{eq:15_24}
\end{equation}
After simplification of the binomial coefficients,
\[
P_{k,n}\left(r^{1}_{j},r_{j}\right)=\binom{n}{k}\prod^{k-1}_{l=0}\left(\dfrac{r^{1}_{j}-l}{r_{j}-l}\right)\prod^{\left(n-k\right)-1}_{l=0}\left(\dfrac{r_{j}-r^{1}_{j}-l}{r_{j}-l}\right).
\]
Hence,
\[
P_{k,n}\left(r^{1}_{j},r_{j}\right)=\binom{n}{k}\prod^{k-1}_{l=0}\left(\dfrac{\dfrac{r^{1}_{j}}{r^{j}}-\dfrac{l}{r^{j}}}{1-\dfrac{l}{r^{j}}}\right)\prod^{\left(n-k\right)-1}_{l=0}\left(\dfrac{1-\dfrac{r^{1}_{j}}{r^{j}}-\dfrac{l}{r^{j}}}{1-\dfrac{l}{r^{j}}}\right).
\]
By the hypothesis that have been made, this shows the convergence\boxeq{
\[
P_{k,n}\left(r^{1}_{j},r_{j}\right)\underset{j\to+\infty}{\longrightarrow}\binom{n}{k}p^{k}\left(1-p\right)^{n-k}.
\]

}

The equality $\refpar{eq:15_24}$ and the previous exercise show the
narrow convergenece and the sequence of probabilities $\mu^{j}$ to
the binomial law 
\[
\mathscr{B}\left(n,p\right)=\sum^{n}_{k=0}\binom{n}{k}p^{k}\left(1-p\right)^{n-k}\delta_{k}.
\]

\end{solution}

\begin{solution}{}{solexercise15.3}

\textbf{1. Determination of the law of $V_{a}.$}

The random variable $V_{a}$ takes values in $\mathbb{N}.$ Assume
that $X$ follows an exponential law $\text{e}^{\lambda}.$ Then,
for every $n\in\mathbb{N},$ 
\[
P\left(V_{a}=n\right)=P\left(na\leqslant X<\left(n+1\right)a\right)=\intop^{\left(n+1\right)a}_{na}\lambda\text{e}^{-\lambda x}\text{d}x.
\]
Computing the integral gives\boxeq{
\[
P\left(V_{a}=n\right)=\text{e}^{-\lambda na}\left(1-\text{e}^{-\lambda a}\right).
\]
}Hence, \textbf{the law of $V_{a}$ is the geometric law on $\mathbb{N}$
with parameter $1-\text{e}^{-\lambda a}.$} 

\textbf{2. Cumulative distribution function of $X_{a}.$ Narrow convergence
when $a$ tends to $0$}

Conversely, assume that for every $a>0,$ the law of $V_{a}$ is the
geometric law on $\mathbb{N}$ with parameter $1-\text{e}^{-\lambda a}.$ 

Since $X_{a}$ takes values in $a\mathbb{N},$ for every $n\in\mathbb{N},$
\[
P\left(X_{a}=na\right)=P\left(V_{a}=n\right)=\text{e}^{-\lambda na}\left[1-\text{e}^{-\lambda a}\right].
\]
Then, for every real number $x\geqslant0,$
\[
P\left(X_{a}>x\right)=\sum_{n:na>x}\text{e}^{-\lambda na}\left[1-\text{e}^{-\lambda a}\right].
\]
Define
\[
n_{0}\left(x\right)=\inf\left\{ n\in\mathbb{N}:\,na>x\right\} =\left\lfloor \dfrac{x}{a}\right\rfloor +1.
\]
Thus,
\[
P\left(X_{a}>x\right)=\left[1-\text{e}^{-\lambda a}\right]\sum^{+\infty}_{n=n_{0}\left(x\right)}\text{e}^{-\lambda na}.
\]
Hence, 
\[
P\left(X_{a}>x\right)=\text{e}^{-\lambda n_{0}\left(x\right)a}=\text{e}^{-\lambda a}\text{e}^{-\lambda a\left\lfloor \frac{x}{a}\right\rfloor }.
\]
Moreover, for every $x<0,$
\[
P\left(X_{a}>x\right)=1.
\]
 Therefore, the cumulative distribution function $F_{X_{a}}$ of $X_{a}$
is given by
\[
F_{X_{a}}\left(x\right)=\begin{cases}
0, & \text{if }x<0,\\
1-\text{e}^{-\lambda a}\text{e}^{-\lambda a\left\lfloor \frac{x}{a}\right\rfloor }, & \text{if }x\geqslant0.
\end{cases}
\]
Now, observe that for every $x\geqslant0,$ 
\[
x-a\leqslant a\left\lfloor \dfrac{x}{a}\right\rfloor \leqslant x.
\]
It follows that 
\[
\lim_{a\to0}a\left\lfloor \dfrac{x}{a}\right\rfloor =x.
\]
Consequently,
\[
\lim_{a\to0}F_{X_{a}}\left(x\right)=\begin{cases}
0, & \text{if }x<0,\\
1-\text{e}^{-\lambda x}, & \text{if }x\geqslant0.
\end{cases}
\]
This proves that \textbf{the family of laws of the random variables
$X_{a}$ converges narrowly to the exponential law $\exp\left(\lambda\right)$
as $a$ tends to 0.} We also say that the family of random variables
$X_{a}$ converges in law to $\exp\left(\lambda\right)$ when $a$
tends to $0.$

\end{solution}

\begin{solution}{}{solexercise15.4}

For every $f\in\mathscr{C}_{b}\left(\mathbb{R}\right),$
\[
\intop_{\mathbb{R}}f\text{d}\mu_{n}=\intop_{\mathbb{R}}f\left(x\right)\dfrac{1}{\sigma_{n}\sqrt{2\pi}}\text{e}^{-\frac{\left(x-m_{n}\right)^{2}}{2\sigma^{2}_{n}}}\text{d}x.
\]
With the change of variables $y=\dfrac{x-m_{n}}{\sigma_{n}},$ we
obtain\boxeq{
\begin{equation}
\intop_{\mathbb{R}}f\text{d}\mu_{n}=\intop_{\mathbb{R}}f\left(y\sigma_{n}+m_{n}\right)\dfrac{1}{\sqrt{2\pi}}\text{e}^{-\frac{y^{2}}{2}}\text{d}y.\label{eq:15_25}
\end{equation}
}

\textbf{1. Narrow convergence of $\left(\mu_{n}\right)_{n\in\mathbb{N}}$
when $\left(m_{n}\right)_{n\in\mathbb{N}}$ and $\left(\sigma_{n}\right)_{n\in\mathbb{N}}$
converge}

Assume that the sequence $\left(m_{n}\right)_{n\in\mathbb{N}}$ and
$\left(\sigma_{n}\right)_{n\in\mathbb{N}}$ converge. 

Since $f$ is countinous,
\[
\lim_{n\to+\infty}f\left(y\sigma_{n}+m_{n}\right)=f\left(y\sigma+m\right).
\]
Moreover, for every $n\in\mathbb{N},$
\[
\left|f\left(y\sigma_{n}+m_{n}\right)\right|\leqslant\left\Vert f\right\Vert _{\infty},
\]
and the function $y\mapsto f\left(y\sigma_{n}+m_{n}\right)$ is integrable
with respect to the Gaussian probability $\mathscr{N}_{\mathbb{R}}\left(0,1\right).$
Hence, by the dominated convergence theorem,
\begin{equation}
\lim_{n\to+\infty}\intop_{\mathbb{R}}f\text{d}\mu_{n}=\intop_{\mathbb{R}}f\left(y\sigma+m\right)\dfrac{1}{\sqrt{2\pi}}\text{e}^{-\frac{y^{2}}{2}}\text{d}y.\label{eq:15_26}
\end{equation}

\textbf{If $\sigma>0,$} we perform the change of variables defined
by $x=y\sigma+m,$ which yields
\begin{equation}
\lim_{n\to+\infty}\intop_{\mathbb{R}}f\text{d}\mu_{n}=\intop_{\mathbb{R}}f\left(x\right)\dfrac{1}{\sqrt{2\pi}}\text{e}^{-\frac{\left(x-m\right)^{2}}{2\sigma^{2}}}\text{d}x.\label{eq:15_26-1}
\end{equation}
Thus, we have proved \textbf{the narrow convergence of the sequence
$\left(\mu_{n}\right)_{n\in\mathbb{N}}$ to the Gaussian probability
$\mathscr{N}_{\mathbb{R}}\left(m,\sigma^{2}\right).$}

Moreover, since for every real number $x,$
\[
\lim_{n\to+\infty}f_{n}\left(x\right)=\dfrac{1}{\sqrt{2\pi}}\text{e}^{-\frac{\left(x-m\right)^{2}}{2\sigma^{2}}},
\]
the Scheffé lemma applies and yields uniform convergence on the Borel
sets, that is, the sequence with general term
\[
\sup_{A\in\mathscr{B}_{\mathbb{R}}}\left|\intop_{A}\dfrac{1}{\sqrt{2\pi}}\text{e}^{-\frac{\left(x-m_{n}\right)^{2}}{2\sigma^{2}_{n}}}\text{d}x-\intop_{A}\dfrac{1}{\sqrt{2\pi}}\text{e}^{-\frac{\left(x-m\right)^{2}}{2\sigma^{2}}}\text{d}x\right|
\]
converges to 0.

\textbf{If $\sigma=0,$} since 
\[
\intop_{\mathbb{R}}\dfrac{1}{\sqrt{2\pi}}\text{e}^{-\frac{y^{2}}{2}}\text{d}y=1,
\]
the relation $\refpar{eq:15_26}$ gives
\[
\lim_{n\to+\infty}\intop_{\mathbb{R}}f\text{d}\mu_{n}=f\left(m\right)=\intop_{\mathbb{R}}f\text{d}\delta_{m}.
\]
Hence, \textbf{the sequence $\left(\mu_{n}\right)_{n\in\mathbb{N}}$
converges narrowly to the Dirac measure at $m.$} 

\textbf{2. Weak and narrow convergence of $\left(\mu_{n}\right)_{n\in\mathbb{N}}$
when $\left(m_{n}\right)_{n\in\mathbb{N}}$ is bounded and $\left(\sigma_{n}\right)_{n\in\mathbb{N}}$
tends to $+\infty$ with $n$}

For every $f\in\mathscr{C}_{0}\left(\mathbb{R}\right),$ the relation
$\refpar{eq:15_25}$ remains valid. 

Assume that the sequence \textbf{$\left(m_{n}\right)_{n\in\mathbb{N}}$}
is bounded and that the sequence $\left(\sigma_{n}\right)_{n\in\mathbb{N}}$
tends to $+\infty$ with $n.$ 

Then, for every $y\neq0,$ 
\[
\lim_{n\to+\infty}\left|y\sigma_{n}+m_{n}\right|=0,
\]
and therefore 
\[
\lim_{n\to+\infty}f\left(y\sigma_{n}+m_{n}\right)=0.
\]
The dominated convergence theorem implies that 
\[
\lim_{n\to+\infty}\intop_{\mathbb{R}}f\text{d}\mu_{n}=0.
\]
Thus, the sequence $\left(\mu_{n}\right)_{n\in\mathbb{N}}$ converges
weakly to the zero measure. However, there is no narrow convergence,
since 
\[
\lim_{n\to+\infty}\mu_{n}\left(\mathbb{R}\right)=1\,\,\,\,\text{and}\,\,\,\,0\left(\mathbb{R}\right)=0.
\]

\begin{remark}{}{}

Under the same assumptions, for every $f\in\mathscr{C}_{b}\left(\mathbb{R}\right),$
using the change of variables defined by $y=x-m_{n},$
\[
\intop_{\mathbb{R}}f\text{d}\mu_{n}=\intop_{\mathbb{R}}f\left(y+m_{n}\right)\dfrac{1}{\sigma_{n}\sqrt{2\pi}}\text{e}^{-\frac{y^{2}}{2\sigma^{2}_{n}}}\text{d}y.
\]
Moreover,
\[
\lim_{n\to+\infty}f\left(y+m_{n}\right)\dfrac{1}{\sigma_{n}\sqrt{2\pi}}\text{e}^{-\frac{y^{2}}{2\sigma^{2}_{n}}}=0.
\]
Hence, this provides \textbf{an example where the dominated convergence
theorem cannot be applied}.

\end{remark}

\end{solution}

\begin{solution}{}{solexercise15.5}

\textbf{1. Existence of a second-order moment for $X_{n}.$ Expectation
and variance }

The random variable $X_{0}$ admits a second-order moment. Suppose
that this holds for $X_{n}.$ Since $Z_{n+1}$ is Gaussian, it admits
a second-order moment, and therefore so does $X_{n+1}.$ 

By linearity of expectation, for every $n\in\mathbb{N}^{\ast},$
\[
\mathbb{E}\left(X_{n}\right)=\rho\mathbb{E}\left(X_{n-1}\right),
\]
and thus\boxeq{
\[
\mathbb{E}\left(X_{n}\right)=\rho^{n}\mathbb{E}\left(X_{0}\right).
\]
}

The random variable $X_{n-1}$ is a linear function of $\left(X_{0},Z_{1},\cdots,Z_{n-1}\right).$
Since the random variables $X_{0},Z_{1},\cdots,Z_{n-1}$ are independent,
so the random variables $X_{n-1}$ and $Z_{n}$ are. 

Hence, for every $n\in\mathbb{N}^{\ast},$
\[
\sigma^{2}_{X_{n}}=\rho^{2}\sigma^{2}_{X_{n-1}}+\sigma^{2}_{Z_{n}}=\rho^{2}\sigma^{2}_{X_{n-1}}+\sigma^{2}.
\]
Solving this linear recurrence yields\boxeq{
\[
\sigma^{2}_{X_{n}}=\begin{cases}
\rho^{2n}\left(\sigma^{2}_{X_{0}}-\dfrac{\sigma^{2}}{1-\rho^{2}}\right)+\dfrac{\sigma^{2}}{1-\rho^{2}}, & \text{if }\left|\rho\right|\neq1,\\
\sigma^{2}_{X_{0}}+n\sigma^{2}, & \text{if }\left|\rho\right|=1.
\end{cases}
\]
}

\textbf{2. Convergence in law of the sequence $\left(X_{n}\right)_{n\in\mathbb{N}}$}

The random variables $X_{0},Z_{1},\cdots,Z_{n}$ are independent and
Gaussian. Therefore, the vector random variable $\left(X_{0},Z_{1},\cdots,Z_{n}\right)$
is Gaussian. Since the real-valued random variable $X_{n}$ is a linear
function of $\left(X_{0},Z_{1},\cdots,Z_{n}\right),$ it is also Gaussian. 

To study the convergence in law of the sequence $\left(X_{n}\right)_{n\in\mathbb{N}},$
we use the results of the previous exercise.
\begin{itemize}
\item \textbf{Case $\left|\rho\right|<1$} 
\[
\lim_{n\to+\infty}\mathbb{E}\left(X_{n}\right)=0\,\,\,\,\text{and}\,\,\,\,\lim_{n\to+\infty}\sigma^{2}_{X_{n}}=\dfrac{\sigma^{2}}{1-\rho^{2}}.
\]
The sequence $\left(X_{n}\right)_{n\in\mathbb{N}}$ converges in law
to the law $\mathscr{N}_{\mathbb{R}}\left(0,\dfrac{\sigma^{2}}{1-\rho^{2}}\right).$
\item \textbf{Case $\left|\rho\right|\geqslant1$}\\
In this case, 
\[
\lim_{n\to+\infty}\sigma^{2}_{X_{n}}=+\infty.
\]

\begin{itemize}
\item If $\mathbb{E}\left(X_{0}\right)=0$ or if $\left|\rho\right|=1,$
then the sequence $\left(\mathbb{E}\left(X_{n}\right)\right)_{n\in\mathbb{N}}$
is bounded. By the previous exercise, \textbf{the sequence $\left(X_{n}\right)_{n\in\mathbb{N}}$
does not converge in law.}
\item If $\mathbb{E}\left(X_{0}\right)\neq0,$ and $\left|\rho\right|>1,$
then 
\[
\lim_{n\to+\infty}\left|\mathbb{E}\left(X_{n}\right)\right|=+\infty.
\]
This case was not studied in the previous exercise. Set $m_{n}=\mathbb{E}\left(X_{n}\right)$
and $\sigma_{n}=\sigma_{X_{n}}.$ \\
For every $f\in\mathscr{C}_{0}\left(\mathbb{R}\right),$
\[
\intop_{\mathbb{R}}f\text{d}P_{X_{n}}=\intop_{\mathbb{R}}f\left(x\right)\dfrac{1}{\sigma_{n}\sqrt{2\pi}}\text{e}^{-\frac{\left(x-m_{n}\right)^{2}}{2\sigma^{2}_{n}}}\text{d}x.
\]
With the change of variables defined by $y=\dfrac{x-m_{n}}{\sigma_{n}},$
we obtain
\begin{equation}
\intop_{\mathbb{R}}f\text{d}P_{X_{n}}=\intop_{\mathbb{R}}f\left(y\sigma_{n}+m_{n}\right)\dfrac{1}{\sqrt{2\pi}}\text{e}^{-\frac{y^{2}}{2}}\text{d}y.\label{eq:int_f_dP_X_n}
\end{equation}
In this case,
\[
y\sigma_{n}+m_{n}=\rho^{n}\left[y\left(\sigma^{2}_{X_{\Omega}}-\dfrac{\sigma^{2}}{1-\rho^{2}}+\dfrac{\sigma^{2}}{\rho^{2n}\left(1-\rho^{2}\right)}\right)^{\frac{1}{2}}+\mathbb{E}\left(X_{0}\right)\right].
\]
It follows that $\lim_{n\to+\infty}f\left(y\sigma_{n}+m_{n}\right)=0,$
for $\lambda-$almost every $y$ and by the dominated convergence
theorem 
\[
\lim_{n\to+\infty}\intop_{\mathbb{R}}f\text{d}P_{X_{n}}=0.
\]
 That is, \textbf{the sequence $\left(P_{X_{n}}\right)_{n\in\mathbb{N}}$
narrowly converges to the zero measure 0.} Of course, there is no
narrow convergence.
\end{itemize}
\end{itemize}
\begin{remark}{}{}

In summary, \textbf{the sequence $\left(X_{n}\right)_{n\in\mathbb{N}}$
converges in law, if and only if $\left|\rho\right|<1.$} This exercise
could also be treated using the Lévy theorem.

\end{remark}

\end{solution}

\begin{solution}{}{solexercise15.6}

The sequence $\left(Y_{n}\right)_{n\in\mathbb{N}^{\ast}}$ converges
in law to 1, and 
\[
\lim_{n\to+\infty}P\left(Y_{n}=0\right)=0.
\]
 For every $f\in\mathscr{C}_{b}\left(\mathbb{R}\right)$ and every
$n\in\mathbb{N}^{\ast},$
\[
\intop_{\Omega}f\left(X_{n}\right)\text{d}P=\intop_{\mathbb{R}}f\text{d}\mathscr{N}_{\mathbb{R}}\left(0,1\right).
\]
Hence,
\[
\left|\intop_{\Omega}f\left(X_{n}Y_{n}\right)\text{d}P-\intop_{\Omega}f\text{d}\mathscr{N}_{\mathbb{R}}\left(0,1\right)\right|=\left|\intop_{\Omega}f\left(X_{n}Y_{n}\right)\text{d}P-\intop_{\Omega}f\left(X_{n}\right)\text{d}P\right|.
\]
Observe that the sets $\left(Y_{n}=1\right)$ and $\left(Y_{n}=0\right)$
form a partition of $\Omega,$ up to a set of probability zero. Moreover,
on the set $\left(Y_{n}=1\right),$ we have $X_{n}=X_{n}Y_{n}.$ Therefore,
\begin{multline*}
\left|\intop_{\Omega}f\left(X_{n}Y_{n}\right)\text{d}P-\intop_{\Omega}f\text{d}\mathscr{N}_{\mathbb{R}}\left(0,1\right)\right|\\
=\left|\intop_{\Omega}f\left(X_{n}Y_{n}\right)\text{d}P-\intop_{\left(Y_{n}=1\right)}f\left(X_{n}Y_{n}\right)\text{d}P-\intop_{\left(Y_{n}=0\right)}f\left(X_{n}\right)\text{d}P\right|,
\end{multline*}
which also gives
\[
\left|\intop_{\Omega}f\left(X_{n}Y_{n}\right)\text{d}P-\intop_{\Omega}f\text{d}\mathscr{N}_{\mathbb{R}}\left(0,1\right)\right|=\left|\intop_{\left(Y_{n}=1\right)}f\left(X_{n}Y_{n}\right)\text{d}P-\intop_{\left(Y_{n}=0\right)}f\left(X_{n}\right)\text{d}P\right|.
\]
It follows that
\[
\left|\intop_{\Omega}f\left(X_{n}Y_{n}\right)\text{d}P-\intop_{\Omega}f\text{d}\mathscr{N}_{\mathbb{R}}\left(0,1\right)\right|\leqslant2\left\Vert f\right\Vert _{\infty}P\left(Y_{n}=0\right),
\]
which proves that
\[
\lim_{n\to+\infty}\intop_{\Omega}f\left(X_{n}Y_{n}\right)\text{d}P=\intop_{\mathbb{R}}f\text{d}\mathscr{N}_{\mathbb{R}}\left(0,1\right).
\]
Thus, the sequence $\left(X_{n}Y_{n}\right)_{n\in\mathbb{N}}$ converges
in law to the law $\mathscr{N}_{\mathbb{R}}\left(0,1\right).$ 

\end{solution}

\begin{solution}{}{solexercise15.7}

\textbf{1. Convergence in law to $\left(X,Y\right)$ of $\left(X_{n},Y_{n}\right)_{n\in\mathbb{N}}.$
Convergence in law to $X+Y$ of $\left(X_{n}+Y_{n}\right)_{n\in\mathbb{N}}$}

For every $n\in\mathbb{N},$ the random variables $X_{n}$ and $Y_{n}$
are independent. Thus, the characteristic function $\varphi_{\left(X_{n},Y_{n}\right)}$
of $\left(X_{n},Y_{n}\right)$ is given, for every $\left(u,v\right)\in\mathbb{R}^{2},$
by
\[
\varphi_{\left(X_{n},Y_{n}\right)}\left(u,v\right)=\varphi_{X_{n}}\left(u\right)\varphi_{Y_{n}}\left(v\right).
\]
Since the two sequences $\left(X_{n}\right)_{n\in\mathbb{N}}$ and
$\left(Y_{n}\right)_{n\in\mathbb{N}}$ converge in law respectively
to the random variables $X$ and $Y,$ the Lévy theorem ensures that
\[
\lim_{n\to+\infty}\varphi_{X_{n}}\left(u\right)=\varphi_{X}\left(u\right)\,\,\,\,\text{and}\,\,\,\,\lim_{n\to+\infty}\varphi_{Y_{n}}\left(v\right)=\varphi_{Y}\left(v\right),
\]
which implies
\[
\lim_{n\to+\infty}\varphi_{\left(X_{n},Y_{n}\right)}\left(u,b\right)=\varphi_{X}\left(u\right)\varphi_{Y}\left(v\right).
\]

Since the random variables $X$ and $Y$ are independent, we also
have
\[
\lim_{n\to+\infty}\varphi_{\left(X_{n},Y_{n}\right)}\left(u,v\right)=\varphi_{\left(X,Y\right)}\left(u,v\right).
\]

The converse part (b) of the Lévy theorem---Theorem $\ref{th:levy_th_measures}$---then
proves that \textbf{the sequence of random variables $\left(X_{n},Y_{n}\right),\,n\in\mathbb{\mathbb{N}},$
converges in law to $\left(X,Y\right).$} Since the random variable
$X_{n}+Y_{n}$ is a \textbf{continuous} function of $\left(X_{n},Y_{n}\right),$
it follows that \textbf{the sequence of random variables $X_{n}+Y_{n},\,n\in\mathbb{N},$
then converges in law to $X+Y.$}

\textbf{2. Counter example}

The sequences $\left(X_{n}\right)_{n\in\mathbb{N}}$ and $\left(Y_{n}\right)_{n\in\mathbb{N}}$
converge $P-$almost surely, and thus in law, respectively to $X$
and $1-X.$ Since the random variables $X,\,1-X$ and $Y$ have the
same law,
\[
X_{n}\stackrel[n\to+\infty]{\mathscr{L}}{\longrightarrow}X\,\,\,\,\text{and}\,\,\,\,Y_{n}\stackrel[n\to+\infty]{\mathscr{L}}{\longrightarrow}Y.
\]

On the other hand, for every $n\in\mathbb{N}^{\ast},$ $X_{n}+Y_{n}=1.$
Thus, the \textbf{sequence of random variables $X_{n}+Y_{n},n\in\mathbb{N},$
converges in law to $\delta_{1}.$} 

However, since the random variables $X$ and $Y$ are independent,
\[
P_{X+Y}=\dfrac{1}{4}\left(\delta_{0}+\delta_{2}\right)+\dfrac{1}{2}\delta_{1}.
\]
\textbf{the sequence of random variables $X_{n}+Y_{n},n\in\mathbb{N}$
does not converge in law to $X+Y.$} A fortiori, \textbf{the sequence
of random variables $\left(X_{n},Y_{n}\right),\,n\in\mathbb{N},$
does not converge in law to $\left(X,Y\right).$}

\end{solution}

\begin{solution}{}{solexercise15.8}

\textbf{1. (a) Proof that $\left(X_{n},Y_{n}\right),\,n\in\mathbb{N}$
converges in law to $\left(X,y_{0}\right)$ }

Let $f$ and $g$ belong to $\mathscr{C}_{0}\left(\mathbb{R}\right).$ 

Since the sequence $\left(Y_{n}\right)_{n\in\mathbb{N}}$ converges
in law to a constant $y_{0},$ then it also converges in probability
to $y_{0}.$ Because the function $g$ is continuous, then the sequence
$\left(g\left(Y_{n}\right)\right)_{n\in\mathbb{N}}$ converges in
probability to $g\left(y_{0}\right).$ Let $\epsilon>0$ be arbitratry.
Then
\begin{equation}
\lim_{n\to+\infty}P\left(\left|g\left(Y_{n}\right)-g\left(y_{0}\right)\right|>\epsilon\right)=0.\label{eq:lim_P_abs_diff_g_Y_n_and_g_y_0}
\end{equation}
By the transfer theorem, 
\begin{multline*}
\left|\intop_{\mathbb{R}^{2}}f\left(x\right)g\left(y\right)\text{d}P_{\left(X_{n}Y_{n}\right)}\left(x,y\right)-\intop_{\mathbb{R}^{2}}f\left(x\right)g\left(y\right)\text{d}P_{X}\otimes\delta_{y_{0}}\left(x,y\right)\right|\\
=\left|\intop_{\Omega}f\left(X_{n}\right)g\left(Y_{n}\right)\text{d}P-\intop_{\Omega}f\left(X\right)g\left(y_{0}\right)\text{d}P\right|.
\end{multline*}
By the triangle inequality,
\begin{multline*}
\left|\intop_{\mathbb{R}^{2}}f\left(x\right)g\left(y\right)\text{d}P_{\left(X_{n}Y_{n}\right)}\left(x,y\right)-\intop_{\mathbb{R}^{2}}f\left(x\right)g\left(y\right)\text{d}P_{X}\otimes\delta_{y_{0}}\left(x,y\right)\right|\\
\leqslant\left|\intop_{\Omega}f\left(X_{n}\right)g\left(Y_{n}\right)\text{d}P-\intop_{\Omega}f\left(X_{n}\right)g\left(y_{0}\right)\text{d}P\right|\\
+\left|\intop_{\Omega}f\left(X_{n}\right)g\left(y_{0}\right)\text{d}P-\intop_{\Omega}f\left(X\right)g\left(y_{0}\right)\text{d}P\right|,
\end{multline*}
Hence,
\begin{multline*}
\left|\intop_{\mathbb{R}^{2}}f\left(x\right)g\left(y\right)\text{d}P_{\left(X_{n}Y_{n}\right)}\left(x,y\right)-\intop_{\mathbb{R}^{2}}f\left(x\right)g\left(y\right)\text{d}P_{X}\otimes\delta_{y_{0}}\left(x,y\right)\right|\\
\leqslant\left\Vert f\right\Vert _{\infty}\left[\intop_{\Omega}\left|g\left(Y_{n}\right)-g\left(y_{0}\right)\right|\text{d}P\right]+\left|g\left(y_{0}\right)\right|\left|\intop_{\Omega}\left(f\left(X_{n}\right)-f\left(X\right)\right)\text{d}P\right|.
\end{multline*}
Partitioning the first integral over the set $\left(\left|g\left(Y_{n}\right)-g\left(y_{0}\right)\right|>\epsilon\right)$
and its complement yields
\begin{multline*}
\left|\intop_{\mathbb{R}^{2}}f\left(x\right)g\left(y\right)\text{d}P_{\left(X_{n}Y_{n}\right)}\left(x,y\right)-\intop_{\mathbb{R}^{2}}f\left(x\right)g\left(y\right)\text{d}P_{X}\otimes\delta_{y_{0}}\left(x,y\right)\right|\\
\leqslant\left\Vert f\right\Vert _{\infty}\left[\epsilon+\intop_{\left(\left|g\left(Y_{n}\right)-g\left(y_{0}\right)\right|>\epsilon\right)}\left|g\left(Y_{n}\right)-g\left(y_{0}\right)\right|\text{d}P\right]+\left\Vert g\right\Vert _{\infty}\left|\intop_{\Omega}\left(f\left(X_{n}\right)-f\left(X\right)\right)\text{d}P\right|\\
\leqslant\left\Vert f\right\Vert _{\infty}\left[\epsilon+2\left\Vert g\right\Vert _{\infty}P\left(\left|g\left(Y_{n}\right)-g\left(y_{0}\right)\right|>\epsilon\right)\right]+\left\Vert g\right\Vert _{\infty}\left|\intop_{\mathbb{R}}f\text{d}P_{X_{n}}-\intop_{\mathbb{R}}f\text{d}P_{X}\right|.
\end{multline*}
By the convergence in law to $X$ of the sequence $\left(X_{n}\right)_{n\in\mathbb{N}},$
and by the relation $\refpar{eq:lim_P_abs_diff_g_Y_n_and_g_y_0},$
it yields 
\[
\limsup_{n\to+\infty}\left|\intop_{\mathbb{R}^{2}}f\left(x\right)g\left(y\right)\text{d}P_{\left(X_{n}Y_{n}\right)}\left(x,y\right)-\intop_{\mathbb{R}^{2}}f\left(x\right)g\left(y\right)\text{d}P_{X}\otimes\delta_{y_{0}}\left(x,y\right)\right|\leqslant\left\Vert f\right\Vert _{\infty}\epsilon.
\]
Since the right hand side term of this inequality is nonnegative,
the arbitrary of $\epsilon$ then ensures that
\[
\lim_{n\to+\infty}\intop_{\mathbb{R}^{2}}f\left(x\right)g\left(y\right)\text{d}P_{\left(X_{n}Y_{n}\right)}\left(x,y\right)=\intop_{\mathbb{R}^{2}}f\left(x\right)g\left(y\right)\text{d}P_{X}\otimes\delta_{y_{0}}\left(x,y\right).
\]

\textbf{(b) $X_{n}+Y_{n},\,n\in\mathbb{N},$ converges in law to $X+y_{0}.$}

Because the set $\mathscr{H}$ is total in $\mathscr{C}_{0}\left(\mathbb{R}^{2}\right),$
this proves that \textbf{the sequence of random variables $\left(\left(X_{n},Y_{n}\right)\right)_{n\in\mathbb{N}}$
converges in law to $\left(X,y_{0}\right).$} Since addition is a
continuous function, it follows that \textbf{the sequence of random
variables $\left(X_{n}+Y_{n}\right)_{n\in\mathbb{N}}$ converges in
law to $X+y_{0}.$}

An analogous proof can be given using the Lévy theorem. This is not
surprising, since both approaches rely on the same density argument,
which is itself a component of the proof of the Lévy theorem. 

We now present this alternative proof.

Since the sequence $\left(Y_{n}\right)$ converges in law to a constant
$y_{0},$ it converges in probability to $y_{0}.$ Since the function
$y\mapsto\text{e}^{\text{i}vy}$ is continuous, the sequence $\left(\text{e}^{\text{i}vY_{n}}\right)_{n\in\mathbb{N}}$
converges in probability to $\text{e}^{\text{i}vy_{0}}.$ Let $\epsilon>0$
be arbitrary. Then
\begin{equation}
\lim_{n\to+\infty}P\left(\left|\text{e}^{\text{i}vY_{n}}-\text{e}^{\text{i}vy_{0}}\right|>\epsilon\right)=0.\label{eq:modP_e_ivY_n_minus_eiivy_0}
\end{equation}

Let $\varphi_{\left(X_{n},Y_{n}\right)}$ denote the characteristic
function of $\left(X_{n},Y_{n}\right),$ and let $\widehat{P_{X}\otimes\delta_{y_{0}}}$
be the Fourier transform of the probability $P_{X}\otimes\delta_{y_{0}}.$
Then, by the triangle inequality,
\begin{multline*}
\left|\varphi_{\left(X_{n},Y_{n}\right)}\left(u,v\right)-\widehat{P_{X}\otimes\delta_{y_{0}}}\left(u,v\right)\right|\\
\leqslant\left|\intop_{\Omega}\text{e}^{\text{i}uX_{n}}\text{e}^{\text{i}vY_{n}}\text{d}P-\intop_{\Omega}\text{e}^{\text{i}uX_{n}}\text{e}^{\text{i}vy_{0}}\text{d}P\right|+\left|\intop_{\Omega}\text{e}^{\text{i}uX_{n}}\text{e}^{\text{i}vy_{0}}\text{d}P-\intop_{\Omega}\text{e}^{\text{i}uX}\text{e}^{\text{i}vy_{0}}\text{d}P\right|.
\end{multline*}
Thus,
\begin{multline*}
\left|\varphi_{\left(X_{n},Y_{n}\right)}\left(u,v\right)-\widehat{P_{X}\otimes\delta_{y_{0}}}\left(u,v\right)\right|\leqslant\intop_{\Omega}\left|\text{e}^{\text{i}vY_{n}}-\text{e}^{\text{i}vy_{0}}\right|\text{d}P+\left|\intop_{\Omega}\left(\text{e}^{\text{i}uX_{n}}-\text{e}^{\text{i}uX}\right)\text{d}P\right|.
\end{multline*}
Partitioning $\Omega$ into the set $\left(\left|\text{e}^{\text{i}vY_{n}}-\text{e}^{\text{i}vy_{0}}\right|>\epsilon\right)$
and its complement,
\begin{align*}
 & \left|\varphi_{\left(X_{n},Y_{n}\right)}\left(u,v\right)-\widehat{P_{X}\otimes\delta_{y_{0}}}\left(u,v\right)\right|\\
 & \qquad\leqslant\epsilon+\intop_{\left|\text{e}^{\text{i}vY_{n}}-\text{e}^{\text{i}vy_{0}}\right|>\epsilon}\left|\text{e}^{\text{i}vY_{n}}-\text{e}^{\text{i}vy_{0}}\right|\text{d}P+\left|\varphi_{X_{n}}\left(u\right)-\varphi_{X}\left(u\right)\right|\\
 & \qquad\leqslant\epsilon+2P\left(\left|\text{e}^{\text{i}vY_{n}}-\text{e}^{\text{i}vy_{0}}\right|>\epsilon\right)+\left|\varphi_{X_{n}}\left(u\right)-\varphi_{X}\left(u\right)\right|.
\end{align*}
By the convergence in law of $X_{n}$ to $X,$ via the Lévy theorem,
and by the relation $\refpar{eq:modP_e_ivY_n_minus_eiivy_0},$
\[
0\leqslant\limsup_{n\to+\infty}\left|\varphi_{\left(X_{n},Y_{n}\right)}\left(u,v\right)-\widehat{P_{X}\otimes\delta_{y_{0}}}\left(u,v\right)\right|\leqslant\epsilon.
\]
Given the arbitrary of $\epsilon,$ it proves that
\[
\lim_{n\to+\infty}\varphi_{\left(X_{n},Y_{n}\right)}\left(u,v\right)=\widehat{P_{X}\otimes\delta_{y_{0}}}\left(u,v\right).
\]
The converse part of the Lévy theorem then implies that \textbf{the
sequence of random variables $\left[\left(X_{n},Y_{n}\right)\right]_{n\in\mathbb{N}}$
converges in law to $\left(X,y_{0}\right).$}

\textbf{2. $\left(Y_{n}\right)_{n\in\mathbb{N}}$ convergence in law
to $X.$}

Assume that the sequence $\left(X_{n}\right)_{n\in\mathbb{N}}$ converges
\textbf{in law} to a random variable $X$ and that the sequence $\left(X_{n}-Y_{n}\right)_{n\in\mathbb{N}}$
converges \textbf{in probability} to 0. Then the sequence $\left(X_{n}-Y_{n}\right)_{n\in\mathbb{N}}$
converges \textbf{in law} to 0. By the Slutsky lemma previously proved,
\[
\left(X_{n},Y_{n}-X_{n}\right)\stackrel[n\to+\infty]{\mathscr{L}}{\longrightarrow}\left(X,0\right).
\]
Since $Y_{n}=\left(Y_{n}-X_{n}\right)+X_{n},$ it follows that the
sequence $\left(Y_{n}\right)_{n\in\mathbb{N}}$ converges in law to
$X.$

\end{solution}

\begin{solution}{}{solexercise15.9}

We have $P-$almost surely, for every $n\in\mathbb{N},$
\[
0\leqslant\dfrac{X_{n}}{10^{n}}\leqslant\dfrac{1}{10^{n-1}},
\]
which shows that the sequence with general term $\dfrac{X_{n}}{10^{n}}$
is $P-$almost surely convergent. Hence \textbf{the sequence $\left(Y_{n}\right)_{n\in\mathbb{N}}$
converges $P-$almost surely to a random variable $Y.$} In particular,
it also converges in probability and in law. 

We now use the Lévy theorem to identify the law of $Y.$ Since the
random variables $X_{n}$ are independent and with same law---hence,
have the same characteristic function---the characteristic function
of $Y_{n}$ is, for every $t\in\mathbb{R},$
\[
\varphi_{Y_{n}}\left(t\right)=\prod^{n}_{j=0}\varphi_{X_{j}}\left(\dfrac{t}{10^{j}}\right)=\prod^{n}_{j=0}\varphi_{X_{0}}\left(\dfrac{t}{10^{j}}\right).
\]
The characteristic function of $X_{0}$ is, for every $t\in\mathbb{R},$
\[
\varphi_{X_{0}}\left(t\right)=\dfrac{1}{10}\sum^{9}_{j=0}\text{e}^{\text{i}jt}=\begin{cases}
\frac{1}{10}\cdot\frac{1-\text{e}^{\text{i}10t}}{1-\text{e}^{\text{i}t}}, & \text{if }\text{e}^{\text{i}t}\neq1,\\
1, & \text{otherwise.}
\end{cases}
\]

Assume $\text{e}^{\text{i}t}\neq1,$ i.e. $t\notin2\pi\mathbb{Z}.$
Then for every $j\in\mathbb{N},$ $\text{e}^{\text{i}\frac{t}{10^{j}}}\neq1$---indeed,
if for some $j\in\mathbb{N}$ one had $\text{e}^{\text{i}\frac{t}{10^{j}}}=1,$
then necessarily $\text{e}^{\text{i}t}=1.$ Therefore, by simplification,
\[
\varphi_{Y_{n}}\left(t\right)=\dfrac{1}{10^{n+1}}\prod^{n}_{j=0}\dfrac{1-\text{e}^{\text{i}\frac{t}{10^{j-1}}}}{1-\text{e}^{\text{i}\frac{t}{10^{j}}}}=\dfrac{1}{10^{n+1}}\cdot\dfrac{1-\text{e}^{\text{i}10t}}{1-\text{e}^{\text{i}\frac{t}{10^{n}}}}.
\]
In this case,
\[
\varphi_{Y_{n}}\left(t\right)=\dfrac{1-\text{e}^{\text{i}10t}}{10^{n+1}}\left[\dfrac{1}{1-\left(1+\text{i}\frac{t}{10^{n}}+o\left(\frac{t}{10^{n}}\right)\right)}\right].
\]
Thus, for every $t\notin2\pi\mathbb{Z},$
\begin{equation}
\lim_{n\to+\infty}\varphi_{Y_{n}}\left(t\right)=\dfrac{\text{e}^{\text{i}10t}-1}{10\text{i}t}=\widehat{\mathscr{U}\left(\left[0,10\right]\right)}\left(t\right)\label{eq:lim_varphi_Y_n_t}
\end{equation}
where $\widehat{\mathscr{U}\left(\left[0,10\right]\right)}$ denotes
the Fourier transform of the uniform law on the interval $\left[0,10\right].$ 

Since the sequence $\left(Y_{n}\right)_{n\in\mathbb{N}}$ converges
in law to $Y,$ it follows by the Lévy theorem that the sequence $\left(\varphi_{Y_{n}}\right)_{n\in\mathbb{N}}$
converges pointwise to the characteristic function $\varphi_{Y}$
of $Y.$ Hence, the relation $\refpar{eq:lim_varphi_Y_n_t}$ implies
that, for every $t\neq2\pi\mathbb{Z},$ 
\[
\varphi_{Y}\left(t\right)=\widehat{\mathscr{U}\left(\left[0,10\right]\right)}\left(t\right).
\]
Since the functions $\varphi_{Y}$ and $\widehat{\mathscr{U}\left(\left[0,10\right]\right)}$
are continuous, they coincide everywhere. By the Fourier transform
injectivity, this proves that the law of $Y$ is the uniform law on
the interval $\left[0,10\right].$

\end{solution}

\begin{solution}{}{solexercise15.10}

Since the random variables $I_{n}$ and $M_{n}$ are defined through
operations related to the order structure, it is natural to use here
the convergence-in-law criterion in terms of cumulative distribution
functions.

\textbf{1. Convergence in law of the sequences $\left(I_{n}\right)_{n\in\mathbb{N}^{\ast}}$
and $\left(M_{n}\right)_{n\in\mathbb{N}^{\ast}}$}

For every $x\in\mathbb{R},$ using the independence of the random
variables $X_{n},$
\[
P\left(I_{n}>x\right)=P\left(\bigcap_{1\leqslant i\leqslant n}\left(X_{i}>x\right)\right)=\prod_{1\leqslant i\leqslant n}P\left(X_{i}>x\right).
\]

Since the $X_{n}$ have same cumulative distribution function $F,$
the cumulative distribution function of $I_{n}$ is given by, for
every real number $x,$
\[
F_{I_{n}}\left(x\right)=1-\left(1-F\left(x\right)\right)^{n}.
\]
It follows that
\[
\lim_{n\to+\infty}F_{I_{n}}\left(x\right)=\begin{cases}
0, & \text{if }F\left(x\right)=0,\\
1, & \text{if }0<F\left(x\right)\leqslant1.
\end{cases}
\]

\begin{itemize}
\item If $x_{i}=\inf\left\{ x:\,F\left(x\right)>0\right\} >-\infty,$ then
\[
\lim_{n\to+\infty}F_{I_{n}}\left(x\right)=\begin{cases}
0, & \text{if }x<x_{i},\\
1, & \text{if }x>x_{i}.
\end{cases}
\]
This proves that\boxeq{
\[
I_{n}\stackrel[n\to+\infty]{\mathscr{L}}{\longrightarrow}x_{i}.
\]
}
\item If $x_{i}=\inf\left\{ x\in\mathbb{R}:\,F\left(x\right)>0\right\} =-\infty,$
then for every $x\in\mathbb{R},$
\[
\lim_{n\to+\infty}F_{I_{n}}\left(x\right)=1.
\]
 The limiting function is not a cumulative distribution function:
hence, the sequence $\left(I_{n}\right)_{n\in\mathbb{N}^{\ast}}$
does not converge in law. 
\end{itemize}
Similarly, for every $x\in\mathbb{R},$ by taking into account the
independence of random variables $X_{n},$
\[
P\left(M_{n}\leqslant x\right)=P\left(\bigcap_{1\leqslant j\leqslant n}\left(X_{i}\leqslant x\right)\right)=\prod_{1\leqslant j\leqslant n}P\left(X_{i}\leqslant x\right).
\]
Since the $X_{n}$ have the same cumulative distribution function
$F,$ the cumulative distribution function of $M_{n}$ is
\[
F_{M_{n}}\left(x\right)=\left[F\left(x\right)\right]^{n}.
\]
Hence,
\[
\lim_{n\to+\infty}F_{M_{n}}\left(x\right)=\begin{cases}
0, & \text{if }F\left(x\right)<1,\\
1, & \text{if }F\left(x\right)=1.
\end{cases}
\]

\begin{itemize}
\item If $x_{s}=\inf\left\{ x\in\mathbb{R}:\,F\left(x\right)=1\right\} <+\infty,$
then
\[
\lim_{n\to+\infty}F_{M_{n}}\left(x\right)=\begin{cases}
0, & \text{if }x<x_{s},\\
1, & \text{if }x>x_{s}.
\end{cases}
\]
This proves that\boxeq{
\[
M_{n}\stackrel[n\to+\infty]{\mathscr{L}}{\longrightarrow}x_{s}.
\]
}
\item If $x_{s}=\inf\left\{ x\in\mathbb{R}:\,F\left(x\right)=1\right\} =+\infty,$
then for every $x\in\mathbb{R},$
\[
\lim_{n\to+\infty}F_{M_{n}}\left(x\right)=0.
\]
The limiting function is not a cumulative distribution function: there
is no convergence in law of the sequence $\left(M_{n}\right)_{n\in\mathbb{N}^{\ast}}.$
\end{itemize}
\textbf{2. Convergence in law of the sequence $\left(Z_{n}\right)_{n\in\mathbb{N}^{\ast}}$}

In this case, for every $x\in\mathbb{R},$
\[
F\left(x\right)=\begin{cases}
0, & \text{if }x\leqslant0,\\
1-\text{e}^{-\lambda x}, & \text{if }x>0.
\end{cases}
\]
Hence, the cumulative distribution function of $Z_{n}$ is given for
every $x>0$ by
\[
F_{Z_{n}}\left(x\right)=F_{M_{n}}\left(x\ln n\right)=\begin{cases}
0, & \text{if }x\leqslant0,\\
\left[1-\text{e}^{-\lambda x\ln n}\right]^{n}, & \text{if }x>0.
\end{cases}
\]
If $x>0,$ then
\[
\ln F_{Z_{n}}\left(x\right)=n\ln\left[1-\dfrac{1}{n^{\lambda x}}\right]=n\left[-\dfrac{1}{n^{\lambda x}}+o\left(\dfrac{1}{n^{\lambda x}}\right)\right].
\]
Hence,
\[
\lim_{n\to+\infty}\ln F_{Z_{n}}\left(x\right)=\begin{cases}
0, & \text{if }\lambda x>1,\\
-\infty, & \text{if }0<\lambda x<1.
\end{cases}
\]
Moreover, for every $x\leqslant0,$ we have 
\[
\lim_{n\to+\infty}F_{Z_{n}}\left(x\right)=0.
\]
Therefore,
\[
\lim_{n\to+\infty}F_{Z_{n}}\left(x\right)=\begin{cases}
1, & \text{if }x>\frac{1}{\lambda},\\
0, & \text{if }x<\frac{1}{\lambda}.
\end{cases}
\]
Noting that $\frac{1}{\lambda}=\mathbb{E}\left(X_{1}\right),$ this
proves that\boxeq{
\[
Z_{n}\stackrel[n\to+\infty]{\mathscr{L}}{\longrightarrow}\mathbb{E}\left(X_{1}\right).
\]
}

\end{solution}

\begin{solution}{}{solexercise15.11}

\textbf{1. $g\in\mathscr{C}_{b}\left(\mathbb{R}\right)$ is nonnegative.
$g\left(x\right)=0$ if and only if $x=0.$ Proof of $\dfrac{1}{\delta}\intop^{\delta}_{0}\left(1-\text{Re}\left(\varphi_{X}\left(t\right)\right)\right)\text{d}t=\intop_{\Omega}g\left(\delta X\right)\text{d}P$}

Continuity at 0 follows from the fact that 
\[
\lim_{x\to0}\dfrac{\sin x}{x}=1.
\]
Moreover, $g$ is continuous at every other point, is even and nonnegative---since
$\left|\sin x\right|\leqslant\left|x\right|.$ Furthermore, 
\[
\lim_{x\to+\infty}g\left(x\right)=1,
\]
which proves that $g\in\mathscr{C}_{b}\left(\mathbb{R}\right).$ Finally,
for every $x>0,$
\[
xg\left(x\right)=\intop^{x}_{0}\left(1-\sin u\right)\text{d}u>0,
\]
which implies that $g\left(x\right)=0$ if and only if $x=0.$

We have
\[
\dfrac{1}{\delta}\intop^{\delta}_{0}\left(1-\text{Re}\left(\varphi_{X}\left(t\right)\right)\right)\text{d}t=\dfrac{1}{\delta}\intop^{\delta}_{0}\left(1-\mathbb{E}\left(\cos\left(tX\right)\right)\right)\text{d}t.
\]
Since $0\leqslant1-\cos\left(tX\right),$ the Fubini theorem gives
\[
\dfrac{1}{\delta}\intop^{\delta}_{0}\left(1-\text{Re}\left(\varphi_{X}\left(t\right)\right)\right)\text{d}t=\dfrac{1}{\delta}\intop_{\Omega}\left[\intop^{\delta}_{0}\left(1-\cos\left(tX\right)\right)\text{d}t\right]\text{d}P.
\]
It follows that
\begin{align*}
\dfrac{1}{\delta}\intop^{\delta}_{0}\left(1-\text{Re}\left(\varphi_{X}\left(t\right)\right)\right)\text{d}t & =\intop_{\left(X\neq0\right)}\left[\dfrac{1}{\delta}\intop^{\delta}_{0}\left(1-\cos\left(tX\right)\right)\text{d}t\right]\text{d}P\\
 & =\intop_{\left(X\neq0\right)}\left[1-\dfrac{\sin\delta X}{\delta X}\right]\text{d}P,
\end{align*}
which proves $\refpar{eq:int_equality_re_varphi_X_g_deltaX},$ since
$g\left(0\right)=0.$

\textbf{2. Proof of $P\left(\left|X\right|>\epsilon\right)\leqslant\dfrac{1}{I_{\epsilon}\delta}\intop^{\delta}_{0}\left(1-\text{Re}\left(\varphi_{X}\left(t\right)\right)\right)\text{d}t=\dfrac{1}{2I_{\epsilon}\delta}\intop^{\delta}_{-\delta}\left(1-\varphi_{X}\left(t\right)\right)\text{d}t.$}

Since $\varphi_{X}\left(-t\right)=\overline{\varphi_{X}\left(t\right)},$
and after making the change of variable defined by $-t=u,$ we obtain
\[
\intop^{0}_{-\delta}\left(1-\varphi_{X}\left(t\right)\right)\text{d}t=\intop^{\delta}_{0}\left(1-\varphi_{X}\left(-u\right)\right)\text{d}u=\intop^{\delta}_{0}\left(1-\overline{\varphi_{X}\left(u\right)}\right)\text{d}u,
\]
which implies
\[
\intop^{\delta}_{-\delta}\left(1-\varphi_{X}\left(t\right)\right)\text{d}t=2\intop^{\delta}_{0}\left(1-\text{Re}\left(\varphi_{X}\left(t\right)\right)\right)\text{d}t.
\]
This gives the equality in the relation $\refpar{eq:P_abs_X_over_eps_upper_bound}.$
Using $\refpar{eq:int_equality_re_varphi_X_g_deltaX},$ the nonnegativity
of $g$ and the definition of $I_{\epsilon},$ we have
\begin{align*}
\dfrac{1}{\delta}\intop^{\delta}_{0}\left(1-\text{Re}\left(\varphi_{X}\left(t\right)\right)\right)\text{d}t & =\intop_{\Omega}g\left(\delta X\right)\text{d}P\geqslant\intop_{\left(\left|X\right|>\epsilon\right)}g\left(\delta X\right)\text{d}P\geqslant I_{\epsilon}P\left(\left|X\right|>\epsilon\right),
\end{align*}
which completes the proof of $\refpar{eq:P_abs_X_over_eps_upper_bound}.$

\textbf{3. Proof that $\left(X_{n}\right)_{n\in\mathbb{N}^{\ast}}$
converges in law to 0 if and only if there exists $\delta>0$ such
that the sequence $\left(\varphi_{X_{n}}\left(t\right)\right)_{n\in\mathbb{N}^{\ast}}$
converges to $1$ for every $t\in\left[-\delta,\delta\right]$}

If the sequence $\left(X_{n}\right)_{n\in\mathbb{N}^{\ast}}$ converges
in law to $0,$ the Lévy theorem yields pointwise convergence of the
sequence $\left(\varphi_{X_{n}}\left(t\right)\right)_{n\in\mathbb{N}^{\ast}}$
to 1, hence, in particular on any interval $\left[-\delta,\delta\right].$ 

Conversely, suppose there exists $\delta>0$ such that the sequence
$\left(\varphi_{X_{n}}\left(t\right)\right)_{n\in\mathbb{N}^{\ast}}$
converges to $1$ for every $t\in\left[-\delta,\delta\right].$ By
$\refpar{eq:P_abs_X_over_eps_upper_bound},$ for every $\epsilon>0$
and every $n\in\mathbb{N}^{\ast},$
\[
P\left(\left|X_{n}\right|>\epsilon\right)\leqslant\dfrac{1}{2I_{\epsilon}\delta}\intop^{\delta}_{-\delta}\left|1-\varphi_{X_{n}}\left(t\right)\right|\text{d}t.
\]

Since $\lim_{n\to+\infty}\left|1-\varphi_{X_{n}}\left(t\right)\right|=0$
on $\left[-\delta,\delta\right]$ and $\left|1-\varphi_{X_{n}}\left(t\right)\right|\leqslant2,$
the dominated convergence theorem gives 
\[
\lim_{n\to+\infty}P\left(\left|X_{n}\right|>\epsilon\right)=0.
\]
Thus\textbf{ the sequence $\left(X_{n}\right)_{n\in\mathbb{N}^{\ast}}$
converges in probability, and therefore also in law, to 0}.

\textbf{4. Proof that the sequence $\left(S_{n}\right)_{n\in\mathbb{N}^{\ast}}$
converges in law if and only if it converges in probability}

Fix $\delta>0,$ and let $m$ and $n$ be any integers such that $m<n.$
By $\refpar{eq:P_abs_X_over_eps_upper_bound},$ for every $\epsilon>0,$
\begin{equation}
P\left(\left|S_{n}-S_{m}\right|>\epsilon\right)\leqslant\dfrac{1}{2I_{\epsilon}\delta}\intop^{\delta}_{-\delta}\left|1-\varphi_{S_{n}-S_{m}}\left(t\right)\right|\text{d}t.\label{eq:Pabs_S_n-S_m}
\end{equation}
The random variables $S_{m}$ and $S_{n}-S_{m}$ are independent,
hence, for every $t\in\mathbb{R},$
\[
\varphi_{S_{n}}\left(t\right)=\varphi_{S_{m}}\left(t\right)\varphi_{S_{n}-S_{m}}\left(t\right),
\]
and therefore
\[
\varphi_{S_{n}}\left(t\right)-\varphi_{S_{m}}\left(t\right)=\varphi_{S_{m}}\left(t\right)\left[1-\varphi_{S_{n}-S_{m}}\left(t\right)\right].
\]
Since the sequence $\left(S_{n}\right)_{n\in\mathbb{N}^{\ast}}$ converges
in law, it follows by the Lévy theorem that the sequence $\left(\varphi_{S_{n}}\right)_{n\in\mathbb{N}^{\ast}}$
converges pointwise to the function $\varphi$ which is equal to 1
at 0, and this convergence is \textbf{uniform} on $\left[-\delta,\delta\right].$
Hence, there exists $N\in\mathbb{N}$ such that $\left|\varphi_{S_{m}}\left(t\right)\right|\geqslant\dfrac{1}{2}$
as soon as $m\geqslant N.$ 

For $n>m\geqslant N,$
\[
\left|\varphi_{S_{n}}\left(t\right)-\varphi_{S_{m}}\left(t\right)\right|=\left|\varphi_{S_{m}}\left(t\right)\right|\left|1-\varphi_{S_{n}-S_{m}}\left(t\right)\right|\geqslant\dfrac{1}{2}\left|1-\varphi_{S_{n}-S_{m}}\left(t\right)\right|.
\]

Plugging this into $\refpar{eq:Pabs_S_n-S_m}$ yields
\[
P\left(\left|S_{n}-S_{m}\right|>\epsilon\right)\leqslant\dfrac{1}{I_{\epsilon}\delta}\intop^{\delta}_{-\delta}\left|\varphi_{S_{n}}\left(t\right)-\varphi_{S_{m}}\left(t\right)\right|\text{d}t.
\]
 Since 
\[
\lim_{n,m\to+\infty}\left|\varphi_{S_{n}}\left(t\right)-\varphi_{S_{m}}\left(t\right)\right|=0,
\]
and since $\left|\varphi_{S_{n}}\left(t\right)-\varphi_{S_{m}}\left(t\right)\right|\leqslant2,$
a double application of the dominated convergence theorem gives
\[
\lim_{n,m\to+\infty}P\left(\left|S_{n}-S_{m}\right|>\epsilon\right)=0.
\]
Thus the sequence $\left(S_{n}\right)_{n\in\mathbb{N}^{\ast}}$ is
Cauchy for the convergence in probability. Hence, the sequence $\left(S_{n}\right)_{n\in\mathbb{N}^{\ast}}$
converges in probability. Since the convergence in probability implies
the convergence in law, the converse implication follows at well.

\begin{remark}{}{}

As a consequence of the Ottoviani inequality, we proved in Exercise
$\ref{exo:exercise11.10}$ the other part of the Lévy theorem. That
is, for a series of independent random variables, the convergence
in probability and $P-$almost sure are equivalent.

\end{remark}

\end{solution}

\begin{solution}{}{solexercise15.12}

\textbf{1. Law of $\left(X_{1},X_{2}\right).$ Linear combination
of $X_{1}$ and $X_{2}$ that is independent of $X_{1}$}

Since $X_{1}$ has a density, and since a conditional density of $X_{2}$
given $X_{1}$ exists, the random variable $\left(X_{1},X_{2}\right)$
admits a density $f_{\left(X_{1},X_{2}\right)}$ given, for every
$\left(x_{1},x_{2}\right)\in\mathbb{R}^{2},$ by
\begin{align*}
f_{\left(X_{1},X_{2}\right)}\left(x_{1},x_{2}\right) & =f_{X_{1}}\left(x_{1}\right)f^{X_{1}=x_{1}}_{X_{2}}\left(x_{2}\right)\\
 & =\dfrac{1}{\sqrt{2\pi}}\text{e}^{-\frac{\left(x_{2}-x_{1}\right)^{2}}{2}}\dfrac{1}{\sqrt{2\pi}}\text{e}^{-\frac{x^{2}_{1}}{2}},
\end{align*}
which proves that $\left(X_{1},X_{2}\right)$ is Gaussian, with density
given for every $\left(x_{1},x_{2}\right)\in\mathbb{R}^{2}$ by\boxeq{
\[
f_{\left(X_{1},X_{2}\right)}\left(x_{1},x_{2}\right)=\dfrac{1}{2\pi}\text{e}^{-\frac{1}{2}\left(x^{2}_{1}+\left(x_{2}-x_{1}\right)^{2}\right)}.
\]
}

The random variable $\left(Y,Z\right)=\left(aX_{1}+bX_{2},X_{1}\right)$
is Gaussian, as linear transformation of the Gaussian random variable
$\left(X_{1},X_{2}\right).$ Therefore, for $Y$ and $Z$ to be independent,
it must and it is enough that $\text{cov}\left(Y,Z\right)=0,$ or
also, since $Z$ is centered, that $\mathbb{E}\left(YZ\right)=0.$
We have
\[
\mathbb{E}\left(YZ\right)=a\mathbb{E}\left(X^{2}_{1}\right)+b\mathbb{E}\left(X_{1}X_{2}\right)=a+b\intop_{\mathbb{R}}x_{1}m^{X_{1}=x_{1}}_{X_{2}}f_{X_{1}}\left(x_{1}\right)\text{d}x_{1},
\]
where $m^{X_{1}=\cdot}_{X_{2}}$ denotes the conditional mean of $X_{2}$
given $X_{1}.$ 

Hence
\[
\mathbb{E}\left(YZ\right)=a+b\mathbb{E}\left(X^{2}_{1}\right)=a+b.
\]
Thus \textbf{$Y$ and $Z$ are independent if and only if $a+b=0.$}

\textbf{2. Computation of $\mathbb{E}^{\mathscr{B}_{n}}\left(X_{n+1}\right)$
and $\mathbb{E}^{\mathscr{B}_{n}}\left(X^{2}_{n+1}\right).$ Expectation
and variance of $X_{n}.$ Proof of $\left(X_{n}\right)_{n\in\mathbb{N}^{\ast}}$
does not converge in $\text{L}^{2}$}

A representative---or version---of the condition expectation $\mathbb{E}^{\mathscr{B}_{n}}\left(X_{n+1}\right)$
is obtained by composing the conditional mean $m^{\uuline{X_{n}}=\uuline{x_{n}}}_{X_{n+1}}$
with $\uuline{X_{n}},$ which yields
\[
\mathbb{E}^{\mathscr{B}_{n}}\left(X_{n+1}\right)=X_{n}.
\]
We then say that the sequence $\left(X_{n}\right)_{n\in\mathbb{N}}$
is a \textbf{martingale}\footnotemark\textbf{\index{martingale}}
with respect to the nondecreasing sequence---for inclusion---of
sub-\salg $\left(\mathscr{B}_{n}\right)_{n\in\mathbb{N}},$ called
a \textbf{\index{filtration}filtration}.

Similarly, by the conditional transfer theorem,
\[
\mathbb{E}^{\mathscr{B}_{n}}\left(X^{2}_{n+1}\right)=\left[\intop_{\mathbb{R}}u^{2}\text{d}P^{\uuline{X_{n}}=\cdot}_{X_{n+1}}\left(u\right)\right]\circ\uuline{X_{n}}=1+X^{2}_{n}.
\]

It follows that
\[
\mathbb{E}\left(X_{n+1}\right)=\mathbb{E}\left(\mathbb{E}^{\mathscr{B}_{n}}\left(X_{n+1}\right)\right)=\mathbb{E}\left(X_{n}\right),
\]
and therefore\boxeq{
\[
\mathbb{E}\left(X_{n}\right)=\mathbb{E}\left(X_{1}\right)=0.
\]
}

Likewise,
\[
\mathbb{E}\left(X^{2}_{n}\right)=\mathbb{E}\left(\mathbb{E}^{\mathscr{B}_{n}}\left(X^{2}_{n+1}\right)\right)=1+\mathbb{E}\left(X^{2}_{n}\right),
\]
and since $\mathbb{E}\left(X^{2}_{1}\right)=1,$\boxeq{
\[
\mathbb{E}\left(X^{2}_{n}\right)=n.
\]
}

\textbf{The sequence $\left(X_{n}\right)_{n\in\mathbb{N}^{\ast}}$
is not bounded in $\text{L}^{2},$ and therefore does not converge
in $\text{L}^{2}.$}

\textbf{3. Existence of a density $f_{\uuline{X_{n}}}$ for $\uuline{X_{n}.}$
Computation of $f_{\uuline{X_{n}}}.$ Characteristic function of $\uuline{X_{n}}$}

The same reasonning as in the first question shows that the random
variable $\left(X_{1},X_{2},X_{3}\right)$ admits a density $f_{\left(X_{1},X_{2},X_{3}\right)}$
given, for every $\left(x_{1},x_{2},x_{3}\right)\in\mathbb{R}^{3},$
by
\[
f_{\left(X_{1},X_{2},X_{3}\right)}\left(x_{1},x_{2},x_{3}\right)=f_{X_{1}}\left(x_{1}\right)f^{X_{1}=x_{1}}_{X_{2}}\left(x_{2}\right)f^{\uuline{X_{2}}=\uuline{x_{2}}}_{X_{3}}\left(x_{3}\right),
\]
which proves that $\left(X_{1},X_{2},X_{3}\right)$ is Gaussian with
density given, for every $\left(x_{1},x_{2},x_{3}\right)\in\mathbb{R}^{3},$
by\boxeq{
\[
f_{\left(X_{1},X_{2},X_{3}\right)}\left(x_{1},x_{2},x_{3}\right)=\dfrac{1}{\left(2\pi\right)^{\frac{3}{2}}}\text{e}^{-\frac{1}{2}\left(x^{2}_{1}+\left(x_{2}-x_{1}\right)^{2}+\left(x_{3}-x_{2}\right)^{2}\right)}.
\]
}

Note that
\[
x^{2}_{1}+\left(x_{2}-x_{1}\right)^{2}+\left(x_{3}-x_{2}\right)^{2}=\left\langle A\uuline{x_{3}},\uuline{x_{3}}\right\rangle ,
\]
where
\[
A=\left(\begin{array}{ccc}
2 & -1 & 0\\
-1 & 2 & -1\\
0 & -1 & 1
\end{array}\right).
\]
A straightforward computation gives
\[
A^{-1}=\left(\begin{array}{ccc}
1 & 1 & 1\\
1 & 2 & 2\\
1 & 2 & 3
\end{array}\right).
\]
The same reasonning then shows that the random variable $\uuline{X_{n}}$
admits a density $f_{\uuline{X_{n}}}$ given, for every $\uuline{x_{n}}\in\mathbb{R}^{n},$
by
\[
f_{\uuline{X_{n}}}\left(\uuline{x_{n}}\right)=f_{X_{1}}\left(x_{1}\right)\prod^{n-1}_{j=1}f^{\uuline{X_{j}}=\uuline{x_{j}}}_{X_{j+1}}\left(x_{j+1}\right),
\]
which proves that $\uuline{X_{n}}$ is \textbf{Gaussian} with density
given, for every $\uuline{x_{n}}\in\mathbb{R}^{n},$ by\boxeq{
\[
f_{\uuline{X_{n}}}\left(\uuline{x_{n}}\right)=\dfrac{1}{\left(2\pi\right)^{\frac{n}{2}}}\text{e}^{-\frac{1}{2}\left(x^{2}_{1}+\left(x_{2}-x_{1}\right)^{2}+\cdots+\left(x_{n}-x_{n-1}\right)^{2}\right)}.
\]
}

Moreover,
\[
x^{2}_{1}+\left(x_{2}-x_{1}\right)^{2}+\cdots+\left(x_{n}-x_{n-1}\right)^{2}=\left\langle A_{n}\uuline{x_{n}},\uuline{x_{n}}\right\rangle ,
\]
where
\[
A_{n}=\left(\begin{array}{ccccccc}
2 & -1 & 0 & \cdots & \cdots & \cdots & 0\\
-1 & 2 & -1 & 0 & \cdots & \cdots & 0\\
0 & -1 & 2 & -1 & 0 & \cdots & \vdots\\
\vdots & \ddots & \ddots & \ddots & \ddots & \ddots & \vdots\\
\vdots & \ddots & \ddots & \ddots & \ddots & \ddots & 0\\
0 & \cdots & \cdots & 0 & -1 & 2 & -1\\
0 & \cdots & \cdots & \cdots & 0 & -1 & 1
\end{array}\right).
\]
Inverting $A_{n}$---obtained, for instance, by solving the associated
linear system---yields
\[
A_{n}=\left(\begin{array}{cccccc}
1 & 1 & 1 & \cdots & \cdots & 1\\
1 & 2 & 2 & \cdots & \cdots & 2\\
1 & 2 & 3 & 3 & \cdots & 3\\
\vdots & \vdots & \vdots & \ddots & \ddots & \vdots\\
1 & 2 & 3 & \cdots & n-1 & n-1\\
1 & 2 & 3 & \cdots & n-1 & n
\end{array}\right).
\]
We saw in the previous question that the Gaussian random variable
$\uuline{X_{n}}$ is centered. Its characteristic function is therefore,
for every $\uuline{t_{n}}\in\mathbb{R}^{n},$\boxeq{
\begin{equation}
\varphi_{\uuline{X_{n}}}\left(\uuline{t_{n}}\right)=\text{e}^{-\frac{1}{2}\left\langle A^{-1}_{n}\uuline{t_{n}},\uuline{t_{n}}\right\rangle }.\label{eq:varphi_uu_X_n__uu_t_n}
\end{equation}
}

\textbf{4. Law of $\left(X_{j},X_{k}\right)$ with $j<k.$ Correlation
coefficient of $X_{j}$ and $X_{k}.$ Convergence in law of $\left(X_{j},\dfrac{X_{k}}{\sqrt{k}}\right)_{k\in\mathbb{N}^{\ast}}.$
Limiting law}

Let $j<k.$ The random variable $\left(X_{j},X_{k}\right)$ being
a marginal of $\uuline{X_{n}},$ is still centered Gaussian and its
covariance matrix is
\[
C_{\left(X_{j},X_{k}\right)}=\left(\begin{array}{cc}
j & j\\
j & k
\end{array}\right).
\]

The correlation coefficient is
\[
\rho_{X_{j},X_{k}}=\dfrac{\text{cov}\left(X_{j},X_{k}\right)}{\sigma_{X_{j}}\sigma_{X_{k}}}=\dfrac{j}{\sqrt{jk}}=\sqrt{\dfrac{j}{k}}.
\]

The characteristic function of $\left(X_{j},X_{k}\right)$ is, for
every $\left(u,v\right)\in\mathbb{R}^{2},$
\[
\varphi_{\left(X_{j},X_{k}\right)}\left(u,v\right)=\text{e}^{-\frac{1}{2}\left\langle C_{\left(X_{j},X_{k}\right)}\binom{u}{v},\binom{u}{v}\right\rangle },
\]
hence,
\[
\varphi_{\left(X_{j},X_{k}\right)}\left(u,v\right)=\text{e}^{-\frac{1}{2}\left(ju^{2}+2juv+kv^{2}\right)}.
\]

It follows that
\[
\varphi_{\left(X_{j},\frac{X_{k}}{\sqrt{k}}\right)}\left(u,v\right)=\varphi_{\left(X_{j},X_{k}\right)}\left(u,\frac{v}{\sqrt{k}}\right)=\text{e}^{-\frac{1}{2}\left(ju^{2}+2\frac{j}{\sqrt{k}}uv+v^{2}\right)},
\]
which implies
\[
\lim_{k\to+\infty}\varphi_{\left(X_{j},\frac{X_{k}}{\sqrt{k}}\right)}\left(u,v\right)=\text{e}^{-\frac{1}{2}\left(ju^{2}+v^{2}\right)}.
\]
By the Lévy theorem,
\[
\left(X_{j},\dfrac{X_{k}}{\sqrt{k}}\right)\stackrel[k\to+\infty]{\mathscr{L}}{\longrightarrow}\mathscr{N}_{\mathbb{R}^{2}}\left(0,\left(\begin{array}{cc}
j & 0\\
0 & 1
\end{array}\right)\right).
\]
\textbf{The limiting law is the product law $\mathscr{N}_{\mathbb{R}}\left(0,j\right)\otimes\mathscr{N}_{\mathbb{R}}\left(0,1\right).$}
We say, that for fixed $j,$ \textbf{the random variables $X_{j}$
and $\dfrac{X_{k}}{\sqrt{k}}$ are \mindex{random variables!asymptotically independent}asymptotically
independent}.

\textbf{5. Convergence in law of $\left(Z_{n}\right)_{n\in\mathbb{N}^{\ast}}$}

Let $\uuline{1_{n}}$ be the vector of $\mathbb{R}^{n}$ whose components
are all equal to 1. Then 
\[
Z_{n}=n^{-\frac{3}{2}}\left\langle \uuline{X_{n}},\uuline{1_{n}}\right\rangle .
\]
The characteristic function of $Z_{n}$ is, for every real number
$t,$
\[
\varphi_{Z_{n}}\left(t\right)=\varphi_{\left\langle \uuline{X_{n}},\uuline{1_{n}}\right\rangle }\left(n^{-\frac{3}{2}}t\right)=\varphi_{\uuline{X_{n}}}\left(n^{-\frac{3}{2}}\uuline{1_{n}}\right).
\]
It follows from $\refpar{eq:varphi_uu_X_n__uu_t_n}$ that
\[
\varphi_{Z_{n}}\left(t\right)=\text{e}^{-\frac{t^{2}}{2n^{3}}\left\langle A^{-1}_{n}\uuline{1_{n}},\uuline{1_{n}}\right\rangle }=\text{e}^{-\frac{t^{2}S_{n}}{2n^{3}}},
\]
where $S_{n}$ denotes the sum of the entries of $A^{-1}_{n}.$ To
compute this sum, we sum along the first diagonal, which gives
\begin{align*}
S_{n} & =\left(1+2+\cdots+n\right)\\
 & +2\left[\left(1+2+\cdots+\left(n-1\right)\right)+\left(1+2+\cdots+\left(n-2\right)\right)+\cdots+1\right].
\end{align*}
Hence,
\begin{align*}
S_{n} & =\dfrac{n\left(n+1\right)}{2}+2\left[\dfrac{\left(n-1\right)n}{2}+\dfrac{\left(n-2\right)\left(n-1\right)}{2}+\cdots+1\right]\\
 & =\dfrac{n\left(n+1\right)}{2}+2\left[\binom{n}{2}+\binom{n-1}{2}+\cdots+\binom{3}{2}+\binom{3}{3}\right].
\end{align*}
Now, by the Pascal identity,
\[
\binom{n+1}{3}=\binom{n}{2}+\binom{n-1}{2}+\cdots+\binom{3}{2}+\binom{3}{3},
\]
so that
\[
S_{n}=\dfrac{n\left(n+1\right)\left(2n+1\right)}{6}.
\]
Therefore,
\[
\varphi_{Z_{n}}\left(t\right)=\text{e}^{-\frac{t^{2}}{12}\cdot\frac{\left(n+1\right)\left(2n+1\right)}{n^{2}}}.
\]
It follows that
\[
\lim_{n\to+\infty}\varphi_{Z_{n}}\left(t\right)=\text{e}^{-\frac{t^{2}}{6}},
\]
which, by the Lévy theorem, yields\boxeq{
\[
Z_{n}\stackrel[n\to+\infty]{\mathscr{L}}{\longrightarrow}\mathscr{N}_{\mathbb{R}}\left(0,\dfrac{1}{3}\right).
\]
}

\end{solution}

\footnotetext{The martingale theory is developped in Chapter \ref{chap:Processus-and-Discrete}.}

\chapter{Processes and Discrete Martingales}\label{chap:Processus-and-Discrete}

\begin{objective}{}{}

Chapter \ref{chap:Processus-and-Discrete} is devoted to processes
and discrete martingales.
\begin{itemize}
\item Section \ref{sec:Some-Example-of} introduces several examples of
processes: motion of a particle in a fluid, random walks in $\mathbb{R}^{n},$
renewal processes, and jump processes.
\item Section \ref{sec:Processes-and-Martingales:} provides the definitions
of a stochastic process and a martingale, followed by several illustrative
examples.
\item Section \ref{sec:Stopping-Time} addresses the notion of stopping
time and defines the \salg of events prior to a stopping time. The
properties of stopping time are then studied, including the characterization
of functions measurable with respect to this $\sigma-$algebra, and
the measurability of the stopped process. A comparison between stopped
\salg is also carried out.
\item Section \ref{sec:First-Stopping-Theorem} begins with a characterization
of martingales via bounded stopping times. It then introduces the
martingale stopped at a stoping time, and states the first Doob stopping
theorem.
\item Section \ref{sec:Maximal-Lemma-and} first presents the maximal lemma,
also known as Doob's maximal inequality. It then states the Doob inequality
and establishes a convergence theorem for martingales bounded in $\text{L}^{2}.$
\item Section \ref{sec:Doob-Decomposition} starts by defining predictable
processes and predictable nondecreasing processes, and then states
the Doob decomposition theorem. The predictable nondecreasing process
associated with a martingale is defined, followed by the strong law
of large numbers for martingales in $\text{L}^{2}.$
\item Section \ref{sec:Convergence-of-Integrable} begins with the definition
of the quadratic variation of a process. A lemma is given describing
the limit of the quadratic variation of a bounded martingale in $\text{L}^{1}.$
A martingale convergence theorem follows, together with a corollary
on the $P-$almost sure convergence of a stopped martingale. The $\text{L}^{1}$-convergence
of an integrable is then characterized, and the equi-integrability
of the conditional expectations of a random variable in $\text{L}^{1}$
is studied.
\item Section \ref{sec:Second-Stopping-Theorem} first establishes a theorem
characterizing closed martingales in terms of arbitrary stopping times,
from which the second Doob stopping theorem is deduced.
\item Section \ref{sec:Convergence-of-Sub-} concludes the chapter with
a submartingale convergence theorem and a corollary on the $P-$almost
convergence of a nonnegative submartingale.
\end{itemize}
\end{objective}

\textbf{We begin by introducing, through examples, several notions
related to stochastic processes. We then turn to the study of bounded
martingales in $\text{L}^{2},$ with particular emphasis on almost
sure convergence results.}

\section{Some Examples of Processes}\label{sec:Some-Example-of}

\paragraph*{Movement of a Particle in a Fluid.}

Let $\left(X_{t},V_{t}\right)$ denote the position-velocity pair
of a particle in a fluid at time $t.$ Because the particle undergoes
numerous collisions with other particles, it is natural to model the
phenomenon via the couple $\left(X_{t},V_{t}\right)$ as a random
variable. 

The family $\left\{ \left(X_{t},V_{t}\right)\right\} _{t\in\mathbb{R}^{+}}$
is therefore a \textbf{continuous-time} \textbf{stochastic process}.
We assume that all random variables are defined on the same probabilized
space $\left(\Omega,\mathscr{A},P\right).$ 

From both a probabilistic and physical perspective of the phenomenon,
we are interested in certain quantities. The observable---or measurable---quantities,
at time $t$ are those that depend only on the past \textbf{history}
of the process up to time $t$---that is, they are ``functions''
of the values of $X_{s}$ and $V_{s}$ for $s\leqslant t.$ 

A classical theorem in measure theory states that a random variable
$Y$ is a ``function'' of a random variable $X,$ meaning that 
\[
Y=f\left(X\right),
\]
for some measurable function $f,$ if and only if $Y$ is measurable
with respect to the $\sigma-$algebra $\sigma\left(X\right)$ generated
by $X.$ 

Extending this result to an uncountable family of random variables---here
$X_{s},\,V_{s},\,s\leqslant t$---is not straighforward. In particular,
one may ask what means a measurable function of all the $X_{s},\,V_{s}$
for $s\leqslant t.$ Nevertheless, it is natural to define the observable
quantities at time $t$ as those measurable with respect to the $\sigma-$algebra
\[
\mathscr{A}_{t}=\sigma\left(\left(X_{s},V_{s}\right):\,s\leqslant t\right).
\]
In this sense, the \textbf{history} of the process up to time $t$
is summarized by the $\sigma-$algebra $\mathscr{A}_{t}.$

Let $f\left(x,v\right)$ denote the value of a physical quantity measure
when the particle has position-velocity pair $\left(x,v\right).$ 

If an observer performs measurements at an increasing sequence of
instants $t_{1},t_{2},\dots,t_{n},\dots$ the available information
is the discrete observation process $\left\{ f\left(X_{t_{n}},V_{t_{n}}\right)\right\} _{n\in\mathbb{N}^{\ast}}$
obtained by the observer. The history of this discrete process at
time $t_{n}$ is summarized by the \salg 
\[
\mathscr{B}_{n}=\sigma\left(f\left(X_{t_{i}},V_{t_{i}}\right):\,i\leqslant n\right).
\]
One may study this observation process together with its own history,
that is, the \textbf{filtration\index{filtration}} $\left(\mathscr{B}_{n}\right)_{n\in\mathbb{N}^{\ast}},$
or consider a richer history---for instance, the one of the full
position-velocity process---namely the \textbf{filtration} $\left(\mathscr{A}_{t_{n}}\right)_{n\in\mathbb{N}^{\ast}}.$
Intermediate situations are also possible.

\paragraph*{Random Walks in $\mathbb{R}^{n}$}

Let $X=\left(X_{n}\right)_{n\in\mathbb{N}}$ be a sequence of independent
random variables taking values in $\mathbb{R}^{n},$ such that the
random variables $X_{n},$ for $n\in\mathbb{N}^{\ast},$ follow the
same law. For every $n\in\mathbb{N},$ define $S_{n}=\sum^{n}_{j=0}X_{j}.$
The family of random variables $S=\left(S_{n}\right)_{n\in\mathbb{N}}$
is a discrete-time process called \textbf{random walk on $\mathbb{R}^{n},$}
starting from the point---possibly random---$X_{0}.$

\paragraph*{Renewal Process}

A sequence of random variables $\left(X_{n}\right)_{n\in\mathbb{N}}$
taking values in $\mathbb{R},$ where each $X_{n}$ is nonnegative
defines a \textbf{\index{renewal process}renewal process}. We now
describe the example that motivates this terminology.

Consider a machine operating continuously, with a component that may
fail. Whenever the component fails, it is immediately replaced by
an identical one. The random variable $X_{n}$ represents the \textbf{lifetime}
of the $n-$th component. Setting 
\[
S_{0}=0,\,\,\,\,S_{n}=\sum^{n}_{j=0}X_{j},
\]
the quantity $S_{n}$ represents the time at which the $n-$th renewal
occurs, it is the \textbf{renewal} \textbf{time}\index{renewal time}.

Another classical example of a renewal process arises in \textbf{queuing
theory}\index{queuing theory}. In this context, $S_{n}$ denotes
the \textbf{arrival time} at the service desk of the $n-$th customer.

\paragraph*{Jump Process}

In the preceding setting of a renewal process $S=\left(S_{n}\right)_{n\in\mathbb{N}},$
for any fixed $t\geqslant0,$ we consider the number $N_{t}$ of indices
$n$ such that $S_{n}\leqslant t.$ This number is random; hence,
it is a random variable. In the examples, it represents either the
number of component replacements that have occured before time $t,$
or the number of customers who have arrived at a service desk during
the interval $\left[0,t\right].$ 

The family $\left(N_{t}\right)_{t\in\mathbb{R}^{+}}$ define a continuous-time
process whose outcomes $\omega\mapsto N_{t}\left(\omega\right)$ are
nondecreasing function taking values in $\mathbb{N}.$ Such process
is called a \textbf{\index{jump process}jump process}. In the particular
case where the random variables $X_{n},\,n\in\mathbb{N}^{\ast}$ are
independent and follow the same exponential law, the process $\left(N_{t}\right)_{t\in\mathbb{R}^{+}}$
is a Poisson process---see Exercise $\ref{exo:exercise12.3}.$ Another
example is the number of pulses $N_{t}$ recorded by a Geiger counter
during the time interval $\left[0,t\right].$

\paragraph*{The Index Does Not Always Have a Temporal Interpretation.}

To study the repartition of gaze molecules in a given volume at a
fixed time, one may partition the space into numbered cubes. We then
consider the discrete process $\left(X_{n}\right)_{n\in\mathbb{N}^{\ast}},$
where $X_{n}$ is the random variable representing the number of molecules
contained in the $n-$th cube.

However, it is not necessary to discretize the space. One may instead
define the notion of process indexed by $\mathbb{R}^{3}.$ In the
case of gaze molecules, this amounts to considering, for each Borel
subset $A$ of $\mathbb{R}^{3},$ the random variable $X_{A}$ which
gives the number of gaze molecules located in the region $A.$

\section{Processes and Martingales: Definitions}\label{sec:Processes-and-Martingales:}

For simplicity, we restrict attention to processes taking values in
$\mathbb{R}$ or $\overline{\mathbb{R}},$ the following definitions
extend in a straightforward manner to processes with values in $\mathbb{R}^{n}.$

\begin{definition}{Process. Discrete Process. Filtration. Process Basis. Natural Filtration. Adapted Process.}{}

An indexed \textbf{\index{process}process}---indexed by a partially
order set of indices $I$---is the family $\left(X_{i}\right)_{i\in I}$
of random variables defined on the same probabilized space $\left(\Omega,\mathscr{A},P\right).$ 

A \textbf{\mindex{discrete!process}\mindex{process!discrete}discrete}
\textbf{process} is a process for which the index set $I$ is at most
countable infinite. In general $I$ will be equal to $\mathbb{N},\,\mathbb{N}^{\ast}$
or $\overline{\mathbb{N}}.$

A nondecreasing---with respect to the inclusion---family $\left(\mathscr{A}_{i}\right)_{i\in I}$
of sub-$\sigma$-algebra of $\mathscr{A}$ is called a \textbf{\index{filtration}filtration}. 

The object $\left(\Omega,\mathscr{A},P,\left(\mathscr{A}_{i}\right)_{i\in I}\right)$
is called\textbf{\mindex{basis!process}} \textbf{\mindex{process!basis}process
basis}. 

If $X=\left(X_{i}\right)_{i\in I}$ is a discrete process, its \textbf{natural
filtration\index{natural filtration}\mindex{filtration!natural}}
is the family of sub-$\sigma$-algebra $\left(\mathscr{A}_{i}\right)_{i\in I},$
where
\[
\mathscr{A}_{i}=\sigma\left(X_{j}:\,j\leqslant i\right).
\]
 The discrete process $X=\left(X_{i}\right)_{i\in I}$ is said to
be \textbf{adapted}\mindex{process!discrete!adapted}---implicitly
with respect to the filtration $\left(\mathscr{A}_{i}\right)_{i\in I}$---if,
for every $i\in I,$ $X_{i}$ is $\mathscr{A}_{i}-$measurable.

\end{definition}

\begin{example}{}{}

For a random walk $S=\left(S_{n}\right)_{n\in\mathbb{N}^{\ast}},$
it is straightforward to verify that its natural filtration coincides
with the natural filtration of the increment process $X=\left(X_{n}\right)_{n\in\mathbb{N}^{\ast}}$
from which it is defined.

\end{example}

An important class of processes is that of \textbf{real-valued} \textbf{discrete
martingales}. The study of almost sure convergence of sequences of
random variables is often facilitated by introducing associated martingales,
for which powerful convergence theorems are available.

\begin{definition}{Adapted Process. Integrable Submartingale. Integrable Supermartingale. Integrable Martingale.}{}

Let $\left(\Omega,\mathscr{A},P,\left(\mathscr{A}_{i}\right)_{i\in I}\right)$
be a process basis, where $I$ is countable and partially ordered.
Let $X=\left(X_{i}\right)_{i\in I}$ be an adapted process such that,
for every $i\in I,$ the random variable $X_{i}$ is $P-$integrable---respectively
non-negative. The process $X$ is:
\begin{itemize}
\item An \textbf{integrable submartingale}\index{integrable submartingale}\mindex{submartingale!integrable}---respectively
a \textbf{nonnegative submartingale}\index{nonnegative submartingale}\mindex{submartingale!nonnegative}---if,
for every $i$ and $j$ such that $i\leqslant j,$ \boxeq{
\[
\mathbb{E}^{\mathscr{A}_{i}}\left(X_{j}\right)\geqslant X_{i}.
\]
}
\item An \textbf{integrable supermartingale}\index{integrable supermartingale}\mindex{supermartingale!integrable}---respectively
\textbf{nonnegative supermartingale}\index{nonnegative supermartingale}\mindex{supermartingale!nonnegative}---if,
for every $i$ and $j$ such that $i\leqslant j,$ \boxeq{
\[
\mathbb{E}^{\mathscr{A}_{i}}\left(X_{j}\right)\leqslant X_{i}.
\]
}
\item An \textbf{integrable martingale\index{integrable martingale}}\mindex{martingale!integrable}
if it is both an integrable submartingale and an integrable supermartingale,
which is equivalent, for every $i$ and $j$ such that $i\leqslant j,$
to \boxeq{
\[
\mathbb{E}^{\mathscr{A}_{i}}\left(X_{j}\right)=X_{i}.
\]
}
\end{itemize}
If $I=\mathbb{N}$ or $\mathbb{N}^{\ast},$ we talk of \textbf{discrete
submartingale}\mindex{discrete!submartingale}\mindex{submartingale!discrete},
\textbf{discrete supermartingale}\mindex{discrete!supermartingale}\mindex{supermartingale!discrete}
or \textbf{discrete martingale}\mindex{discrete!martingale}\mindex{martingale!discrete}. 

If $I=\overline{\mathbb{N}}$ and $\mathscr{A}_{\infty}=\bigvee_{n\in\mathbb{N}}\mathscr{A}_{n},$
where $\mathscr{A}_{\infty}$ is the $\sigma-$algebra generated by
the union of the $\sigma-$algebras $\mathscr{A}_{n},n\in\mathbb{N},$
one speaks of \textbf{closed discrete martingale}\index{closed discrete martingale}\mindex{martingale!discrete!closed}. 

With the same definition of $\mathscr{A}_{\infty},$ a discrete integrable
martingale $X=\left(X_{n}\right)_{n\in\mathbb{N}}$ is said to be
\textbf{closable}\mindex{martingale!discrete integrable!closable}
if there exists an $\mathscr{A}_{\infty}-$measurable random variable
$X_{\infty}$ such that for every $n\in\mathbb{N},$ 
\[
X_{n}=\mathbb{E}^{\mathscr{A}_{n}}\left(X_{\infty}\right).
\]
 The extended process $X=\left(X_{n}\right)_{n\in\overline{\mathbb{N}}}$
is then a closed martingale.

An integrable discrete martingale $X=\left(X_{n}\right)_{n\in\mathbb{N}}$
is said \textbf{bounded}\mindex{martingale!discrete integrable!bounded}
in $\text{L}^{1}$ if 
\[
\sup_{n\in\mathbb{N}}\mathbb{E}\left(\left|X_{n}\right|\right)<+\infty.
\]

A discrete martingale $X=\left(X_{n}\right)_{n\in\mathbb{N}}$ is
said to be \textbf{square-integrable}\mindex{martingale!discrete!square-integrable}\textbf{---}or
in $\text{L}^{2}$---if for every $n\in\mathbb{N},$ $X_{n}$ is
square-integrable, that is
\[
\mathbb{E}\left(X^{2}_{n}\right)<+\infty,
\]
and it is \textbf{bounded in $\text{L}^{2}$} if 
\[
\sup_{n\in\mathbb{N}}\mathbb{E}\left(X^{2}_{n}\right)<+\infty.
\]

\end{definition}

\begin{remark}{}{}

A submartingale---respectively supermartingale---is nondecreasing---respectively
nonincreasing---in conditional expectation. A martingale is constant
in conditional expectation.

\end{remark}

Here are some examples of such process.

\begin{example}{}{}

\textbf{1. Renewal Process}

A renewal process is a submartingale with respect to its natural filtration.

\textbf{2. Simple Random Walk on $\mathbb{Z}$}

Consider the random walk in $\mathbb{Z}$ defined by $S_{0}=a$ and
$S_{n}=\sum^{n}_{j=1}X_{j},n\in\mathbb{N}^{\ast},$ where $\left(X_{n}\right)_{n\in\mathbb{N}^{\ast}}$
is a sequence of independent random variables with same law 
\[
p\delta_{1}+\left(1-p\right)\delta_{-1},\,\,\,\,0<p<1.
\]
This walk may represent, for instance, the fortune of a player tossing
a coin: at each toss, the player wins one unit with probability $p$
and loses one unit with probability $1-p.$ 

Let $\left(\mathscr{A}_{n}\right)_{n\in\mathbb{N}}$ be the natural
filtration of the process $S.$ For $n\in\mathbb{N}^{\ast},$ 
\[
\mathscr{A}_{n}=\sigma\left(X_{j}:\,1\leqslant j\leqslant n\right).
\]
 We compute
\[
\mathbb{E}^{\mathscr{A}_{n}}\left(S_{n+1}\right)=S_{n}+\mathbb{\mathbb{E}}^{\mathscr{A}_{n}}\left(X_{n+1}\right).
\]
Since the variables $X_{n}$ are independent,
\[
\mathbb{E}^{\mathscr{A}_{n}}\left(X_{n+1}\right)=\mathbb{E}\left(X_{n+1}\right)=p-\left(1-p\right)=2p-1.
\]
If follows that
\[
\begin{cases}
S\text{ is a submartingale,} & \text{if }p>\frac{1}{2},\\
S\text{ is a supermartingale,} & \text{if }p<\frac{1}{2},\\
S\text{ is a martingale,} & \text{if }p=\frac{1}{2}.
\end{cases}
\]

\end{example}

\begin{remarks}{}{}

1. The process $X$ is a submartingale if and only if the process
$-X$ is a supermartingale.

2. If $X$ and $Y$ are submartingales, then for every nonnegative
real numbers $a$ and $b,$ the process 
\[
aX+bY=\left(aX_{i}+bY_{i}\right)_{i\in I}
\]
is also a submartingale.

3. If $X$ and $Y$ are submartingales, then the process 
\[
X\vee Y=\left(X_{i}\vee Y_{i}\right)_{i\in I}
\]
is a submartingale. 

Similarly, if $X$ and $Y$ are supermartingales, then the process
\[
X\land Y=\left(X_{i}\land Y_{i}\right)_{i\in I}
\]
 is a supermartingale\footnotemark. 

Indeed, the first assertion results from that, if $i\leqslant j,$
\[
\mathbb{E}^{\mathscr{A}_{i}}\left(X_{j}\vee Y_{j}\right)\geqslant\mathbb{E}^{\mathscr{A}_{i}}\left(X_{j}\right)\geqslant X_{i}
\]
and
\[
\mathbb{E}^{\mathscr{A}_{i}}\left(X_{j}\vee Y_{j}\right)\geqslant\mathbb{E}^{\mathscr{A}_{i}}\left(Y_{j}\right)\geqslant Y_{i},
\]
and the second results from that
\[
\mathbb{E}^{\mathscr{A}_{i}}\left(X_{j}\land Y_{j}\right)\leqslant\mathbb{E}^{\mathscr{A}_{i}}\left(X_{j}\right)\leqslant X_{i}
\]
and
\[
\mathbb{E}^{\mathscr{A}_{i}}\left(X_{j}\land Y_{j}\right)\leqslant\mathbb{E}^{\mathscr{A}_{i}}\left(Y_{j}\right)\leqslant Y_{i}.
\]

4. If $X$ is a submartingale---respectively a supermartingale---integrable
and if $i\leqslant j,$ then $\mathbb{E}\left(X_{j}\right)\geqslant\mathbb{E}\left(X_{i}\right)$---respectively
$\mathbb{E}\left(X_{j}\right)\leqslant\mathbb{E}\left(X_{i}\right).$
In particular, if $X$ is an integrable martingale, then we have $\mathbb{E}\left(X_{j}\right)=\mathbb{E}\left(X_{i}\right).$

5. A process $X$ is a \textbf{discrete} submartingale---respectively
supermartingale, or martingale---if and only if, for every $n\in\mathbb{N},$
\[
\mathbb{E}^{\mathscr{A}_{n}}\left(X_{n+1}\right)\geqslant X_{n}
\]
---respectively $\mathbb{E}^{\mathscr{A}_{n}}\left(X_{n+1}\right)\leqslant X_{n}$
or $\mathbb{E}^{\mathscr{A}_{n}}\left(X_{n+1}\right)=X_{n}.$ 

6. Let $X=\left(X_{n}\right)_{n\in\mathbb{N}}$ be a martingale in
$\text{L}^{2}.$ Then the process $X^{2}=\left(X^{2}_{n}\right)_{n\in\mathbb{N}}$
is a submartingale. Consequently, the sequence $\left(\mathbb{E}\left(X^{2}_{n}\right)\right)_{n\in\mathbb{N}}$
is nondecreasing. 

Indeed, since $X_{n}$ is $\mathscr{A}_{n}-$measurable, then, for
every $n\in\mathbb{N},$
\[
\mathbb{E}^{\mathscr{A}_{n}}\left[\left(X_{n+1}-X_{n}\right)^{2}\right]=\mathbb{E}^{\mathscr{A}_{n}}\left(X^{2}_{n+1}\right)+\mathbb{E}^{\mathscr{A}_{n}}\left(X^{2}_{n}\right)-2X_{n}\mathbb{E}^{\mathscr{A}_{n}}\left(X_{n+1}\right)=\mathbb{E}^{\mathscr{A}_{n}}\left(X^{2}_{n+1}\right)-\mathbb{E}^{\mathscr{A}_{n}}\left(X^{2}_{n}\right),
\]
which proves that 
\[
\mathbb{E}^{\mathscr{A}_{n}}\left(X^{2}_{n+1}\right)\geqslant\mathbb{E}^{\mathscr{A}_{n}}\left(X^{2}_{n}\right).
\]
 By integrating, it follows that $\mathbb{E}\left(X^{2}_{n+1}\right)\geqslant\mathbb{E}\left(X^{2}_{n}\right).$

\end{remarks}

\footnotetext{Notation recall: for every real numbers $a$ and $b,$
$a\vee b=\max\left(a,b\right)$ and $a\land b=\min\left(a,b\right),$
which may be read respectively ``$a\,\text{sup\,}b$'' and ``$a\,\text{inf\,}b$''.}

\textbf{In what follows, unless explicit stated otherwise, the processes
introduced are defined on the same process basis $\left(\Omega,\mathscr{A},P,\left(\mathscr{A}_{n}\right)_{n\in\mathbb{N}}\right).$} 

When there is no ambiguity, the adjective \textit{discrete} will be
suppressed.

\begin{example}{}{}

Let $\left(X_{n}\right)_{n\in\mathbb{N}}$ be a sequence of integrable
random variables defined on the probabilized space \preds. For each
$n\in\mathbb{N},$ denote by $\mathscr{B}_{n}$ the \salg $\sigma\left(X_{i}:\,i\leqslant n\right).$
Define the random variables $\left(Y_{n}\right)_{n\in\mathbb{N}}$
by $Y_{0}=0,$ and, for $n\geqslant1,$
\[
Y_{n}=\sum^{n}_{i=1}\left(X_{i}-\mathbb{E}^{\mathscr{B}_{i-1}}\left(X_{i}\right)\right).
\]
Then the process $\left(Y_{n}\right)_{n\in\mathbb{N}}$ is a martingale
on the process basis $\left(\Omega,\mathscr{A},P,\left(\mathscr{B}_{n}\right)_{n\in\mathbb{N}}\right).$ 

One can also says that $\left(Y_{n}\right)_{n\in\mathbb{N}}$ is a
martingale with respect to the filtration $\left(\mathscr{B}_{n}\right)_{n\in\mathbb{N}},$
or more briefly, when there is no ambiguity, a martingale.

\end{example}

An important special case is when the random variables $X_{n}$ are
independent. In that case, for $n\geqslant1,$
\[
Y_{n}=\sum^{n}_{j=1}\left(X_{i}-\mathbb{E}\left(X_{i}\right)\right).
\]

\begin{example}{}{}

Let $U$ be an integrable random variable and let $\left(\mathscr{A}_{n}\right)_{n\in\mathbb{N}}$
be a filtration on the probabilized space \preds. Define $X_{n}=\mathbb{E}^{\mathscr{A}_{n}}\left(U\right).$
Then the process $\left(X_{n}\right)_{n\in\mathbb{N}}$ is a bounded
martingale in $\text{L}^{1}.$

\end{example}

\begin{example}{}{}

Let $X=\left(X_{i}\right)_{i\in I}$ be a nonnegative submartingale,
and let $f$ be a convex, nondecreasing function from $\mathbb{R}^{+}$
onto itself such that $f\left(X_{i}\right)$ is integrable for every
$i\in I.$ By the Jensen inequality, the process $f\left(X\right)=\left(f\left(X_{i}\right)\right)_{i\in I}$
is a nonnegative submartingale. In particular, this holds for the
processes $X^{p}=\left(X^{p}_{i}\right)_{i\in I}$ with $p\geqslant1$
and $X^{+}=\left(X^{+}_{i}\right)_{i\in I}.$

\end{example}

\begin{example}{}{}

Let $X=\left(X_{i}\right)_{i\in I}$ be a martingale. The process
$\left(\left|X_{i}\right|\right)_{i\in I}$ is a submartingale since,
if $i\leqslant j,$ then 
\[
\left|X_{i}\right|=\left|\mathbb{E}^{\mathscr{A}_{i}}\left(X_{j}\right)\right|\leqslant\mathbb{E}^{\mathscr{A}_{i}}\left(\left|X_{j}\right|\right).
\]
More generally, if $f$ is a continuous convex fuction, such that
$f\left(X_{i}\right)$ is integrable for every $i\in I,$ by the Jensen
inequality, the process $f\left(X\right)=\left(f\left(X_{i}\right)\right)_{i\in I}$
is a submartingale. In particular, if $X$ is a martingale $\text{L}^{2},$
then the process $X^{2}=\left(X^{2}_{i}\right)_{i\in I}$ is a submartingale.

\end{example}

\section{Stopping Time}\label{sec:Stopping-Time}

The notion of time used in modeling a stochastic process refers to
the observer's clock. However, the random phenomenon under study does
not necessarily evolve in a simple manner according to this clock.
This leads to the introduction of random times, called \textbf{stopping
time}\index{stopping time}, which reflect the internal clock of the
process.

Let $\left(\Omega,\mathscr{A},P,\left(\mathscr{A}_{n}\right)_{n\in\mathbb{N}}\right)$
be a process basis. Denote $\mathscr{A}_{\infty}=\bigvee_{n\in\mathbb{N}}\mathscr{A}_{n}.$

\begin{definition}{Stopping Time}{}

A function $T$ from $\Omega$ to $\overline{\mathbb{N}}$ is called
a \textbf{stopping time}\index{stopping time}, if, for every $n\in\mathbb{N},$
$\left(T=n\right)\in\mathscr{A}_{n}.$

\end{definition}

\begin{remark}{}{}

If $\left(\mathscr{A}_{n}\right)_{n\in\mathbb{N}}$ is the natural
filtration of a process $X=\left(X_{n}\right)_{n\in\mathbb{N}}$ taking
values in a measurable space $\left(E,\mathscr{E}\right),$ then a
function $T$ from $\Omega$ to $\overline{\mathbb{N}}$ is a stopping
time if, for every $n\in\mathbb{N},$ there exists a measurable function
$f_{n}$ of $\left(E^{n+1},\mathscr{E}^{\otimes\left(n+1\right)}\right)$
taking values 0 or 1 such that
\[
\boldsymbol{1}_{\left(T=n\right)}=f_{n}\left(X_{0},X_{1},\cdots,X_{n}\right).
\]

If $T$ is a stopping time, then we have $\left(T=+\infty\right)\in\mathscr{A}_{\infty}.$
Indeed, $\left(T=+\infty\right)=\bigcup_{n\in\mathbb{N}}\left(T=n\right)$
and, for every $n\in\mathbb{N},$ we have $\left(T=n\right)\in\mathscr{A}_{\infty},$
since $\mathscr{A}_{n}\subset\mathscr{A}_{\infty}.$

\end{remark}

\begin{example}{}{}

Every constant function $T$ from $\Omega$ to $\overline{\mathbb{N}}$
is a stopping time.

\end{example}

\begin{example}{}{}

Let $X=\left(X_{n}\right)_{n\in\mathbb{N}}$ be an adapted process
and let $A$ be a Borel set.

The \textbf{hitting time} $T_{A}$ of $A$ is defined by
\[
T_{A}=\inf\left\{ n\in\mathbb{N}:\,X_{n}\in A\right\} ,
\]
with the convention 
\[
\inf\emptyset=+\infty.
\]

Then $T_{A}$ is a stopping time.

Indeed, by the adaptness of $X$ and by the fact that the sequence
of sub-$\sigma$-algebra of the filtration is nondecreasing,
\[
\left(T_{A}=0\right)=\left(X_{0}\in A\right)\in\mathscr{A}_{0},
\]
and, for every $n\in\mathbb{N}^{\ast},$
\[
\left(T_{A}=n\right)=\left[\bigcap^{n-1}_{k=0}\left(X_{k}\notin A\right)\right]\cap\left(X_{n}\in A\right)\in\mathscr{A}_{n}.
\]

\end{example}

\begin{example}{}{}

The \textbf{\index{last visit time}last visit time} $\tau_{A}$ of
an adapted process $X$ to a Borel set $A$ is not, in general, a
stopping time. The last visit time is defined by
\[
\tau_{A}=\sup\left\{ n\in\mathbb{N}^{\ast}:\,X_{n}\in A\right\} 
\]
with the convention 
\[
\sup\emptyset=0.
\]

Indeed, 
\[
\left(\tau_{A}=0\right)=\bigcap_{n\in\mathbb{N}^{\ast}}\left(X_{n}\notin A\right)\in\mathscr{A}_{\infty},
\]
but, in general $\left(\tau_{A}=0\right)\notin\mathscr{A}_{0}.$ 

Moreover, for every $n\in\mathbb{N}^{\ast},$
\[
\left(\tau_{A}=n\right)=\left(X_{n}\in A\right)\bigcap\left[\bigcap_{k\geqslant n+1}\left(X_{n}\notin A\right)\right]\notin\mathscr{A}_{n}.
\]

\end{example}

\begin{denotation}{}{}

We denote by $\mathscr{T}$---respectively $\mathscr{T}_{b}$---the
set of stopping time---respectively bounded stopping times---relative
to the filtration $\left(\mathscr{A}_{n}\right)_{n\in\mathbb{N}}.$

\end{denotation}

\begin{definition}{$\sigma-$algebra of the Events Preceding a Stopping Time}{}

Let $T$ be a stopping time. The family of events $\mathscr{A}_{T}$
defined by
\[
\mathscr{A}_{T}=\left\{ A\in\mathscr{A}_{\infty}:\,\forall n\in\mathbb{N},\,\,\,\,A\cap\left(T=n\right)\in\mathscr{A}_{n}\right\} 
\]
is a $\sigma-$algebra. It is called the \textbf{\salg of events
preceding $T.$\mindex{sigma-algebra@\salg! of events preceding a stopping time}} 

\end{definition}

\begin{remark}{}{}

The stopping time $T$ is $\mathscr{A}_{T}-$measurable.

\end{remark}

\begin{lemma}{}{}

(a) Let $T$ be a stopping time and let $A\in\mathscr{A}_{\infty}.$ 

Then $A\in\mathscr{A}_{T}$ if and only if\footnotemark, for every
$n\in\mathbb{N},$ $\left(T\leqslant n\right)\in\mathscr{A}_{n}.$

(b) If $T_{1},T_{2},\cdots,T_{k}$ are stopping times, then it is
the same for the functions $\inf_{1\leqslant i\leqslant k}T_{i}$
and $\sup_{1\leqslant i\leqslant k}T_{i}.$ In particular, if $T$
is a stopping time, for every integer $k,$ $T\land k$ is a bounded
stopping time.

\end{lemma}

\footnotetext{This characterization justifies the name \textbf{\salg
of events preceding $T.$}}

\begin{proof}{}{}

(a) Suppose that $T$ is a stopping time and let $A\in\mathscr{A}_{T}.$
If $k\leqslant n,$ then $A\cap\left(T=k\right)\in\mathscr{A}_{k}\subset\mathscr{A}_{n}$
and therefore
\[
A\cap\left(T\leqslant n\right)=A\cap\left[\bigcup_{k\leqslant n}\left(T=k\right)\right]\in\mathscr{A}_{n}.
\]

Conversely, suppose that for every $n\in\mathbb{N},$ $A\cap\left(T\leqslant n\right)\in\mathscr{A}_{n}.$
If $n\geqslant1,$ then 
\[
A\cap\left(T\leqslant n-1\right)\in\mathscr{A}_{n-1}\subset\mathscr{A}_{n},
\]
and consequently
\[
A\cap\left(T=n\right)=\left[A\cap\left(T\leqslant n\right)\right]\backslash\left[A\cap\left(T\leqslant n-1\right)\right]\in\mathscr{A}_{n}.
\]
Moreover,
\[
A\cap\left(T=0\right)=A\cap\left(T\leqslant0\right)\in\mathscr{A}_{0}.
\]

(b) From the previous characterization, by noting that $T$ is a stopping
time if and only if $\Omega\in\mathscr{A}_{T}.$

(c) The $T_{i},\,1\leqslant i\leqslant k$ are stopping times. By
the previous stopping time characterization, for every $n\in\mathbb{N},$
\[
\left(\inf_{1\leqslant i\leqslant k}T_{i}\leqslant n\right)=\bigcup_{1\leqslant i\leqslant k}\left(T_{i}\leqslant n\right)\in\mathscr{A}_{n}
\]
and
\[
\left(\sup_{1\leqslant i\leqslant k}T_{i}\leqslant n\right)=\bigcap_{1\leqslant i\leqslant k}\left(T_{i}\leqslant n\right)\in\mathscr{A}_{n}.
\]
This proves the result still by the same characterization.

\end{proof}

The following lemma characterizes $\mathscr{A}_{T}-$measurable functions
and gives the expression of the conditional expectation of a random
variable with respect to the \salg $\mathscr{A}_{T}.$

\paragraph*{Convention}

A function $X$ defined on a subset $\Omega^{\prime}$ of $\Omega$
is said to be $\mathscr{A}-$measurable if it is measurable with respect
to the trace space $\left(\Omega^{\prime},\Omega^{\prime}\cap\mathscr{A}\right).$

\begin{lemma}{Characterization of $\mathscr{A}_T-$measurable Functions. Conditional Expectation}{}

(a) Let $T$ be a stopping time and let $X:\Omega\to\overline{\mathbb{R}}$
be an $\mathscr{A}_{\infty}-$measurable function.

Then the function $X$ is $\mathscr{A}_{T}-$measurable if and only
if, for every $n\in\mathbb{N},$ its restriction $X_{\left(T=n\right)}$
to the set $\left(T=n\right)$ is $\mathscr{A}_{n}-$measurable.

(b) Let $X$ be a numerical random variable defined on \preds, nonnegative
or integrable. For every $n\in\overline{\mathbb{N}},$
\[
\mathbb{E}^{\mathscr{A}_{T}}\left(X\right)=\mathbb{E}^{\mathscr{A}_{n}}\left(X\right)\,\,\,\,\text{on}\,\left(T=n\right),
\]
that is
\[
\mathbb{E}^{\mathscr{A}_{T}}\left(X\right)=\sum_{n\in\overline{\mathbb{N}}}\boldsymbol{1}_{\left(T=n\right)}\mathbb{E}^{\mathscr{A}_{n}}\left(X\right).
\]

\end{lemma}

\begin{proof}{}{}

(a) First suppose that $X=\boldsymbol{1}_{B},$ where $B\in\mathscr{A}_{\infty}.$
For every $n\in\mathbb{N},$ its restriction $X_{\left(T=n\right)}$
to the set $\left(T=n\right)$ is $\boldsymbol{1}_{B\cap\left(T=n\right)}.$
The announced equivalence then follows directly from the definition
of the \salg $\mathscr{A}_{T}.$

By linearity, the result holds when $X$ is a step function and also
when $X$ is nonnegative---$X$ is then a pointwise limit of a nondecreasing
sequence of step functions. The general case follows by decomposing
$X$ into its nonnegative and nonpositive parts.

(b) Let $X$ be a nonnegative numerical random variable and, for every
$n\in\overline{\mathbb{N}},$ let $Y_{n}$ be a representative of
$\mathbb{E}^{\mathscr{A}_{n}}\left(X\right).$ Each $Y_{n}$ is $\mathscr{A}-$measurable
and nonnegative. Then, the nonnegative random variable $Y=\sum_{n\in\mathbb{N}}\boldsymbol{1}_{\left(T=n\right)}Y_{n}$
is, by the previous property, $\mathscr{A}_{T}-$measurable. Moreover,
since $X$ is nonnegative, for every $A\in\mathscr{A}_{T},$
\[
\intop_{A}X\text{d}P=\sum_{n\in\overline{\mathbb{N}}}\intop_{A\cap\left(T=n\right)}X\text{d}P.
\]
Because $A\cap\left(T=n\right)\in\mathscr{A}_{n},$ and the integrands
are nonnegative,
\[
\intop_{A}X\text{d}P=\sum_{n\in\overline{\mathbb{N}}}\intop_{A\cap\left(T=n\right)}\mathbb{E}^{\mathscr{A}_{n}}\left(X\right)\text{d}P=\intop_{A}\left[\sum_{n\in\overline{\mathbb{N}}}\boldsymbol{1}_{\left(T=n\right)}Y_{n}\right]\text{d}P=\intop_{A}Y\text{d}P.
\]
Hence, the result when $X$ is nonnegative. The case where $X$ is
integrable with arbitrary sign follows by decomposing $X$ into its
nonnegative and nonpositive parts.

\end{proof}

\begin{corollary}{Stopped Process is $\mathscr{A}_T-$measurable}{}

Let $X=\left(X_{n}\right)_{n\in\mathbb{N}}$ be an adapted process
and let $T$ be a finite stopping time. Define the function $X_{T}$
on the set $\left(T<+\infty\right),$ for every $n\in\mathbb{N}$
by 
\[
X_{T}=X_{n}\,\,\,\,\text{on}\,\left(T=n\right).
\]
Then $X_{T}$ is $\mathscr{A}_{T}-$measurable.

Moreover, let $\left(X_{n}\right)_{n\in\overline{\mathbb{N}}}$ be
an adapted process, and let $T$ be a stopping time. Define the function
$X_{T},$ for every $n\in\overline{\mathbb{N}}$, by 
\[
X_{T}=X_{n}\,\,\,\,\text{on}\,\left(T=n\right).
\]
Then $X_{T}$ is $\mathscr{A}_{T}-$measurable.

\end{corollary}

\begin{proposition}{Comparison of Stopped $\sigma-$algebras}{}

Let $S$ and $T$ be stopping times.

(a) The events $\left(S<T\right),\,\left(S=T\right)$ and $\left(S\leqslant T\right)$
belong to $\mathscr{A}_{S}$ and to $\mathscr{A}_{T}.$

(b) If $B\in\mathscr{A}_{S},$ then $B\cap\left(S\leqslant T\right)\in\mathscr{A}_{T}.$

(c) Consequently, if the stopping time $S$ and $T$ are such that
$S\leqslant T,$ then $\mathscr{A}_{S}\subset\mathscr{A}_{T}.$ It
follows that the families of sub-\salg $\left(\mathscr{A}_{T}\right)_{T\in\mathscr{T}_{b}}$
and $\left(\mathscr{A}_{T}\right)_{T\in\mathscr{T}}$ are filtrations.

\end{proposition}

\begin{proof}{}{}

(a) For every $n\in\mathbb{N},$
\[
\left(S<T\right)\cap\left(S=n\right)=\left(n<T\right)\cap\left(S=n\right)\in\mathscr{A}_{n},
\]
since $\left(n<T\right)=\left(T\leqslant n\right)^{c}$ and both $S$
and $T$ are stopping times. Hence, $\left(S<T\right)$ belongs to
$\mathscr{A}_{S}.$ 

Moreover, for every $n\in\mathbb{N}^{\ast},$
\[
\left(S<T\right)\cap\left(S=n\right)=\left(S<n-1\right)\cap\left(T=n\right)\in\mathscr{A}_{n},
\]
since $\left(S<n-1\right)\in\mathscr{A}_{n-1}\subset\mathscr{A}_{n}$
and both $S$ and $T$ are stopping times.

Also,
\[
\left(S<T\right)\cap\left(T=0\right)=\emptyset\in\mathscr{A}_{0}.
\]
It follows that $\left(S<T\right)$ belongs to $\mathscr{A}_{T}.$
Thus, $\left(S<T\right)\in\mathscr{A}_{S}\cap\mathscr{A}_{T}.$

For every $n\in\mathbb{N},$
\[
\left(S=T\right)\cap\left(S=n\right)=\left(S=n\right)\cap\left(T=n\right)\in\mathscr{A}_{n},
\]
which proves that $\left(S=T\right)\in\mathscr{A}_{S}.$ 

Similarly, 
\[
\left(S=T\right)\in\mathscr{A}_{T}.
\]

Taking into account these results, it yields
\[
\left(S\leqslant T\right)=\left(S<T\right)\cap\left(S=T\right)\in\mathscr{A}_{S}\cap\mathscr{A}_{T}.
\]

(b) If $B\in\mathscr{A}_{S},$ then for every $n\in\mathbb{N},$ 
\[
B\cap\left(S\leqslant n\right)\in\mathscr{A}_{n}.
\]
Hence,
\[
\left[B\cap\left(S\leqslant T\right)\right]\cap\left(T=n\right)=\left[B\cap\left(S\leqslant n\right)\right]\cap\left(T=n\right)\in\mathscr{A}_{n},
\]
which proves that 
\[
B\cap\left(S\leqslant T\right)\in\mathscr{A}_{T}.
\]

(c) If the stopping times $S$ and $T$ are such that $S\leqslant T,$
then $\left(S\leqslant T\right)=\Omega,$ and by the previous statement
ensures that $\mathscr{A}_{S}\subset\mathscr{A}_{T}.$

\end{proof}

\begin{remark}{}{}

If $B\in\mathscr{A}_{S},$ then
\[
B\cap\left(S=T\right)\in\mathscr{A}_{T}\,\,\,\,\text{and}\,\,\,\,B\cap\left(S<T\right)\in\mathscr{A}_{T}.
\]
Hence, by the previous proposition,
\[
B\cap\left(S=T\right)=\left[B\cap\left(S\leqslant T\right)\right]\cap\left(S=T\right)\in\mathscr{A}_{T}
\]
and
\[
B\cap\left(S<T\right)=\left[B\cap\left(S\leqslant T\right)\right]\cap\left(S<T\right)\in\mathscr{A}_{T}.
\]

\end{remark}

\section{First Stopping Theorem}\label{sec:First-Stopping-Theorem}

We prove a \textbf{characterization of martingales} in terms of \textbf{bounded}
stopping times. From this, we deduce the first \textbf{Doob} stopping
theorem, which is frequently used.

\begin{theorem}{Characterization of Martingales By Bounded Stopping Times}{ch_mart_bounded_stopping_time}

Let $X=\left(X_{n}\right)_{n\in\mathbb{N}}$ be an adapted process.
The following properties are equivalent:

(i) $X$ is an integrable martingale.

(ii) For every $T\in\mathscr{T}_{b},$ $X_{T}\in\text{L}^{1}\left(\Omega,\mathscr{A}_{T},P\right)$
and $\mathbb{E}\left(X_{T}\right)=\mathbb{E}\left(X_{0}\right).$

(iii) The process $\left(X_{T}\right)_{T\in\mathscr{T}_{b}}$ is a
martingale with respect to the filtration $\left(\mathscr{A}_{T}\right)_{T\in\mathscr{T}_{b}}.$

\end{theorem}

\begin{proof}{}{}
\begin{itemize}
\item \textbf{$\text{(i)}\Rightarrow\text{(ii)}$}\\
Suppose that $X$ is an integrable martingale. If $T\in\mathscr{T}_{b}$
is bounded by $k,$ then
\[
X_{T}=\sum^{k}_{j=0}\boldsymbol{1}_{\left(T=j\right)}X_{j},
\]
and thus $X_{T}\in\text{L}^{1}\left(\Omega,\mathscr{A}_{T},P\right).$\\
Let $A\in\mathscr{A}_{T}.$ We have
\[
A=\bigcup^{k}_{j=0}\left[A\cap\left(T=j\right)\right]
\]
and thus
\[
\intop_{A}X_{T}\text{d}P=\sum^{k}_{j=0}\intop_{A\cap\left(T=j\right)}X_{j}\text{d}P.
\]
Since $X$ is a martingale and, for every $j\in\mathbb{N},$ $A\cap\left(T=j\right)\in\mathscr{A}_{j},$
it follows that
\[
\intop_{A}X_{T}\text{d}P=\sum^{k}_{j=0}\intop_{A\cap\left(T=j\right)}X_{k}\text{d}P=\intop_{A}X_{k}\text{d}P.
\]
Since $X_{T}$ is $\mathscr{A}_{T}-$measurable, we have proved that
$X_{T}=\mathbb{E}^{\mathscr{A}_{T}}\left(X_{k}\right).$ \\
Taking the expectations, and using that $X$ is a martingale, we obtain
\[
\mathbb{E}\left(X_{T}\right)=\mathbb{E}\left(X_{k}\right)=\mathbb{E}\left(X_{0}\right).
\]
\item \textbf{$\text{(ii)}\Rightarrow\text{(iii)}$}\\
Assume that the property (ii) holds. Let $S$ and $T$ be two stopping
times bounded by $k$ such that $S\leqslant T\leqslant k.$ Then,
\[
\mathscr{A}_{S}\subset\mathscr{A}_{T}\subset\mathscr{A}_{k}.
\]
Let $A\in\mathscr{A}_{S}.$ Define the function $R$ by 
\[
R=S\boldsymbol{1}_{A}+k\boldsymbol{1}_{A^{c}}
\]
It is bounded by $k$ and is a stopping time. Indeed,
\[
\left(R=n\right)=\begin{cases}
\left(S=n\right)\cap A\in\mathscr{A}_{n}, & \text{if }n<k,\\
\left[\left(S=k\right)\cap A\right]\cup A^{c}, & \text{if }n=k,\\
\emptyset, & \text{if }n>k.
\end{cases}
\]
When $n=k,$ since $\left(S=k\right)\cap A\in\mathscr{A}_{k}$ and
that $\mathscr{A}_{S}\subset\mathscr{A}_{k},$ we have $\left(R=n\right)\in\mathscr{A}_{k}=\mathscr{A}_{n}.$
Hence, for every $n\in\mathbb{N},$ $\left(R=n\right)\in\mathscr{A}_{n}.$\\
By applying the hypothesis to the bounded stopping times $R$ and
$k,$ it follows that
\[
\mathbb{E}\left(X_{R}\right)=\mathbb{E}\left(X_{k}\right)=\mathbb{E}\left(X_{0}\right)
\]
which yields
\[
\mathbb{E}\left(\boldsymbol{1}_{A}X_{S}+\boldsymbol{1}_{A^{c}}X_{k}\right)=\mathbb{E}\left(X_{k}\right),
\]
or equivalently,
\[
\mathbb{E}\left(\boldsymbol{1}_{A}X_{S}\right)=\mathbb{E}\left(\boldsymbol{1}_{A}X_{k}\right).
\]
Since $X_{S}$ is $\mathscr{A}_{S}-$measurable, it follows that 
\[
X_{S}=\mathbb{E}^{\mathscr{A}_{S}}\left(X_{k}\right).
\]
 Similarly 
\[
X_{T}=\mathbb{E}^{\mathscr{A}_{T}}\left(X_{k}\right).
\]
 Since, $\mathscr{A}_{S}\subset\mathscr{A}_{T},$
\[
\mathbb{E}^{\mathscr{A}_{S}}\left(X_{T}\right)=\mathbb{E}^{\mathscr{A}_{S}}\left(\mathbb{E}^{\mathscr{A}_{T}}\left(X_{k}\right)\right)=\mathbb{E}^{\mathscr{A}_{S}}\left(X_{k}\right)=X_{S},
\]
which shows that the property (iii) holds.
\item \textbf{$\text{(iii)}\Rightarrow\text{(i)}$}\\
It is enough to consider constant stopping times.
\end{itemize}
\end{proof}

\begin{corollary}{Martingale Stopped at a Stopping Time}{}

Let $X=\left(X_{n}\right)_{n\in\mathbb{N}}$ be a martingale. For
every stopping time $T,$ the process 
\[
X^{T}=\left(X_{T\land n}\right)_{n\in\mathbb{N}}
\]
is a martingale. It is called the \textbf{\mindex{martingale!stopped at given stopping time}martingale
stopped} \textbf{at time $T.$}

\end{corollary}

\begin{proof}{}{}

For every bounded stopping time $S,$ we have $X^{T}_{S}=X_{T\land S},$
and since $T\land S$ is a bounded stopping time, Theorem $\ref{th:ch_mart_bounded_stopping_time}$
applied to the martingale $X$ gives
\[
\mathbb{E}\left(X^{T}_{S}\right)=\mathbb{E}\left(X_{T\land S}\right)=\mathbb{E}\left(X_{0}\right)=\mathbb{E}\left(X^{T}_{0}\right).
\]
which implies, again from this theorem and the arbitrary of $S,$
that $X^{T}$ is a martingale.

\end{proof}

\begin{figure}[t]
\begin{center}\includegraphics[width=0.4\textwidth]{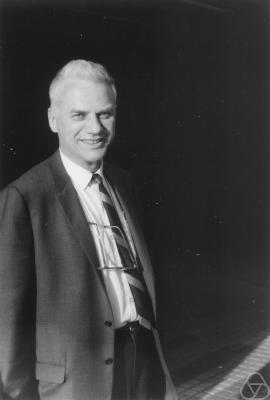}

{\tiny Konrad Jacobs---\href{https://opc.mfo.de/detail?photo_id=919}{https://opc.mfo.de/detail?photo\_id=919}---CC
BY-SA 2.0 de}\end{center}

\caption{\textbf{\protect\href{https://en.wikipedia.org/wiki/Joseph_L._Doob}{Joseph Leo Doob}}
(1910 - 2004)}\sindex[fam]{Doob, Joseph}
\end{figure}

An explicit formulation of the implication $\text{(i)}\Rightarrow\text{(iii)}$
in Theorem $\ref{th:ch_mart_bounded_stopping_time}$ yields the \textbf{\sindex[fam]{Doob, Joseph}Doob}\footnote{\textbf{\href{https://en.wikipedia.org/wiki/Joseph_L._Doob}{Joseph Leo Doob}\sindex[fam]{Doob, Joseph}}
(1910 - 2004) was an American mathematician. He has developed the
theory of maringales. He worked at the department of mathematics at
the University of Illinois.} first stopping theorem.

\begin{theorem}{Doob First Stopping Theorem}{}

Let $X=\left(X_{n}\right)_{n\in\mathbb{N}}$ be a martingale. For
every bounded stopping times $S$ and $T$ such that $S\leqslant T,$
\[
\mathbb{E}^{\mathscr{A}_{S}}\left(X_{T}\right)=X_{S}.
\]

\end{theorem}

\section{Maximal Lemma and Martingales in $\text{L}^{2}$}\label{sec:Maximal-Lemma-and}

We now aim to prove an almost sure convergence theorem for martingales
in $\text{L}^{2}$ that are bounded in $\text{L}^{2}.$ 

Before doing so, we state the maximal lemma, also known as Doob maximal
inequality, for submartingales. This result extends the Kolmogorov
inequality for sums of independent random variables.

\begin{lemma}{Maximal Lemma or Maximal Inequality of Doob}{}

(a) Let $X=\left(X_{n}\right)_{n\in\mathbb{N}}$ be a nonnegative
or integrable submartingale. Then, for every integer $N$ and every
$\epsilon>0,$
\begin{equation}
P\left(\sup_{0\leqslant n\leqslant N}X_{n}>\epsilon\right)\leqslant\dfrac{1}{\epsilon}\left(\intop_{\left(\sup_{0\leqslant n\leqslant N}X_{n}>\epsilon\right)}X_{N}\text{d}P\right).\label{eq:16_1}
\end{equation}
In particular,
\begin{equation}
P\left(\sup_{0\leqslant n\leqslant N}X_{n}>\epsilon\right)\leqslant\dfrac{1}{\epsilon}\mathbb{E}\left(\left|X_{N}\right|\right).\label{eq:16_2}
\end{equation}
Moreover, for every $\epsilon>0,$
\begin{equation}
P\left(\sup_{n\in\mathbb{N}}X_{n}>\epsilon\right)\leqslant\dfrac{1}{\epsilon}\sup_{n\in\mathbb{N}}\mathbb{E}\left(\left|X_{N}\right|\right).\label{eq:16_3}
\end{equation}

(b) In particular, if $X$ is an integrable martingale bounded in
$\text{L}^{1}$---that is such that $\sup_{n\in\mathbb{N}}\mathbb{E}\left(\left|X_{n}\right|\right)<+\infty$---then
the random variable $X^{\ast}=\sup_{n\in\mathbb{N}}\left|X_{n}\right|$
is finite $P-$almost surely.

\end{lemma}

\begin{proof}{}{}

(a) Our goal is to obtain an upper bound for the probability of the
set 
\[
E=\sup_{0\leqslant n\leqslant N}X_{n}>\epsilon.
\]
If $E$ is empty, the inequality is trivial, so we assume that $E\neq\emptyset.$
Let $k$ be the first index at which $X_{k}$ is over $\epsilon.$
To formulate this, define the sets
\[
E_{0}=\left(X_{0}>\epsilon\right)
\]
and, for $1\leqslant k\leqslant N,$ 
\[
E_{k}=\left(X_{k}>\epsilon\right)\cap\left[\bigcap^{k-1}_{i=0}\left(X_{i}\leqslant\epsilon\right)\right].
\]
These sets form a partition of $E.$ Consequently,
\[
\intop_{E}X_{N}\text{d}P=\sum^{N}_{k=0}\intop_{E_{k}}X_{N}\text{d}P.
\]
By definition of $E_{k},$
\[
\intop_{E}X_{N}\text{d}P\geqslant\epsilon\sum^{N}_{k=0}P\left(E_{k}\right)=\epsilon P\left(E\right),
\]
which proves the inequality $\refpar{eq:16_1}.$ 

The inequality $\refpar{eq:16_2}$ follows from it immediately.

Finally, for every $\epsilon>0,$ the sequence of sets $\left(\sup_{0\leqslant n\leqslant N}X_{n}>\epsilon\right)$
is nondecreasing in $N,$ and their union is $\left(\sup_{n\in\mathbb{N}}X_{n}>\epsilon\right).$
Hence,
\[
P\left(\sup_{n\in\mathbb{N}}X_{n}>\epsilon\right)=\lim_{N\to+\infty}P\left(\sup_{0\leqslant n\leqslant N}X_{n}>\epsilon\right)\leqslant\dfrac{1}{\epsilon}\sup_{n\in\mathbb{N}}\mathbb{E}\left(\left|X_{n}\right|\right).
\]

(b) Since the process $\left|X\right|$ is a submartingale, it follows
from the equality $\refpar{eq:16_3}$ that, for every $k\in\mathbb{N}^{\ast},$
\[
P\left(\sup_{n\in\mathbb{N}}\left|X_{n}\right|>k\right)\leqslant\dfrac{1}{k}\sup_{n\in\mathbb{N}}\mathbb{E}\left(\left|X_{n}\right|\right),
\]
Letting $k$ tends to the infinity, yields
\[
P\left(\sup_{n\in\mathbb{N}}\left|X_{n}\right|=+\infty\right)=\lim_{k\to+\infty}P\left(\sup_{n\in\mathbb{N}}\left|X_{n}\right|>k\right)=0.
\]

\end{proof}

As a corollary, we obtain the Doob inequality for martingales bounded
in $\text{L}^{2}.$

\begin{theorem}{Inequality of Doob}{doob_ineq}

Let $X=\left(X_{n}\right)_{n\in\mathbb{N}}$ be a martingale bounded
in $\text{L}^{2}.$ Then the random variable 
\[
X^{\ast}=\sup_{n\in\mathbb{N}}\left|X_{n}\right|
\]
belongs to $\text{L}^{2},$ and we have the Doob inequality
\begin{equation}
\left\Vert X^{\ast}\right\Vert _{\text{L}^{2}}\leqslant2\sup_{n\in\mathbb{N}}\left\Vert X_{n}\right\Vert _{\text{L}^{2}}.\label{eq:Doob_ineq}
\end{equation}

\end{theorem}

\begin{proof}{}{}

Let $M_{n}=\sup_{0\leqslant k\leqslant n}\left|X_{k}\right|.$ Since
$M_{n}\leqslant\sum^{n}_{k=0}\left|X_{k}\right|,$ it follows that
$M_{n}\in\text{L}^{2}.$ Because the process $\left|X\right|$ is
a submartingale nonnegative and integrable, the maximal lemma implies
that, for every $a>0,$
\[
a\mathbb{E}\left(\boldsymbol{1}_{\left(M_{n}>a\right)}\right)\leqslant\mathbb{E}\left(\left|X_{n}\right|\boldsymbol{1}_{\left(M_{n}>a\right)}\right).
\]
Integrating with respect to the Lebesgue measure $\lambda$ on $\mathbb{R}^{+},$
we obtain the inequality
\[
\intop_{\mathbb{R}^{+}}a\mathbb{E}\left(\boldsymbol{1}_{\left(M_{n}>a\right)}\right)\text{d}\lambda\left(a\right)\leqslant\intop_{\mathbb{R}^{+}}\mathbb{E}\left(\left|X_{n}\right|\boldsymbol{1}_{\left(M_{n}>a\right)}\right)\text{d}\lambda\left(a\right).
\]
Hence, by the Fubini theorem and by integration,
\[
\mathbb{E}\left(\intop_{\left[0,M_{n}\right[}a\text{d}\lambda\left(a\right)\right)=\dfrac{1}{2}\mathbb{E}\left(M^{2}_{n}\right)\leqslant\mathbb{E}\left(\left|X_{n}\right|M_{n}\right).
\]
Applying the Schwarz inequality to the right-hand side yields
\[
\dfrac{1}{2}\mathbb{E}\left(M^{2}_{n}\right)\leqslant\left[\mathbb{E}\left(X^{2}_{n}\right)\right]^{\frac{1}{2}}\left[\mathbb{E}\left(M^{2}_{n}\right)\right]^{\frac{1}{2}},
\]
which implies
\[
\left[\mathbb{E}\left(M^{2}_{n}\right)\right]^{\frac{1}{2}}\leqslant2\left[\mathbb{E}\left(X^{2}_{n}\right)\right]^{\frac{1}{2}}.
\]
A fortiori,
\[
\left[\mathbb{E}\left(M^{2}_{n}\right)\right]^{\frac{1}{2}}\leqslant2\sup_{n\in\mathbb{N}}\left\Vert X_{n}\right\Vert _{\text{L}^{2}}.
\]
Since the sequence $\left(M_{n}\right)_{n\in\mathbb{N}}$ is nondecreasing
and converges to $X^{\ast},$ letting $n$ tends to infinity and applying
the Beppo Levi property yields the inequality $\refpar{eq:Doob_ineq}.$ 

\end{proof}

\begin{remark}{}{}

Under the asumptions of Theorem $\ref{th:doob_ineq},$ we obtain the
double inequality
\[
\sup_{n\in\mathbb{N}}\left\Vert X_{n}\right\Vert _{\text{L}^{2}}\leqslant\left\Vert X^{\ast}\right\Vert _{\text{L}^{2}}\leqslant2\sup_{n\in\mathbb{N}}\left\Vert X_{n}\right\Vert _{\text{L}^{2}}.
\]

\end{remark}

As a corollary, we obtain the following convergence theorem for martingales
bounded in $\text{L}^{2}.$

\begin{theorem}{Convergence Theorem in $\textrm{L}^2$}{cv_th_L2}

Let $X=\left(X_{n}\right)_{n\in\mathbb{N}}$ be a martingale bounded
in $\text{L}^{2}.$ 

Then the sequence $\left(X_{n}\right)_{n\in\mathbb{N}}$ converges
$P-$almost surely and in $\text{L}^{2}$ to a random variable $X_{\infty}.$
For every $n\in\mathbb{N},$ 
\[
X_{n}=\mathbb{E}^{\mathscr{A}_{n}}\left(X_{\infty}\right).
\]

If, in addition, the filtration is complete, in the sense where the
\salg $\mathscr{A}_{0}$ contains every $\mathscr{A}-$negligible
sets, then $X_{\infty}$ is $\mathscr{A}_{\infty}-$measurable and
the martingale $X$ is closable.

\end{theorem}

\begin{proof}{}{}

We first prove the $P-$almost sure convergence. Denote by $\left\{ X\longrightarrow\right\} $
the set of $\omega$ for which the sequence $\left(X_{n}\left(\omega\right)\right)_{n\in\mathbb{N}}$
converges in $\mathbb{R}.$ By the Cauchy criterion,
\[
\left\{ X\longrightarrow\right\} =\bigcap_{\epsilon\in\mathbb{Q}^{+}}\bigcup_{N\in\mathbb{N}^{\ast}}\bigcap_{m,n\geqslant N}\left\{ \left|X_{n}-X_{m}\right|\leqslant\epsilon\right\} ,
\]
and therefore
\[
\left\{ X\longrightarrow\right\} ^{c}=\bigcup_{\epsilon\in\mathbb{Q}^{+}}\bigcap_{N\in\mathbb{N}^{\ast}}\bigcup_{m,n\geqslant N}\left\{ \left|X_{n}-X_{m}\right|>\epsilon\right\} .
\]

For $N\in\mathbb{N}^{\ast},$
\[
\bigcup_{m,n\geqslant N}\left\{ \left|X_{n}-X_{m}\right|>\epsilon\right\} \subset\left\{ \sup_{m,n\geqslant N}\left|X_{n}-X_{m}\right|>\epsilon\right\} \subset\left\{ \sup_{n\geqslant N}\left|X_{n}-X_{N}\right|>\dfrac{\epsilon}{3}\right\} .
\]
The last inclusion follows from the contrapositive of
\begin{align*}
\sup_{n\geqslant N}\left|X_{n}-X_{N}\right|\leqslant\dfrac{\epsilon}{3} & \,\,\,\,\Longrightarrow\,\,\,\,\forall m,n\geqslant N,\,\left|X_{m}-X_{N}\right|\leqslant\dfrac{\epsilon}{3}\,\,\,\,\text{and}\,\,\,\,\left|X_{n}-X_{N}\right|\leqslant\dfrac{\epsilon}{3}\\
 & \,\,\,\,\Longrightarrow\,\,\,\,\forall m,n\geqslant N,\,\left|X_{m}-X_{n}\right|\leqslant\dfrac{2\epsilon}{3}<\epsilon\\
 & \,\,\,\,\Longrightarrow\,\,\,\,\sup_{m,n\geqslant N}\left|X_{m}-X_{n}\right|\leqslant\epsilon.
\end{align*}
By the Markov inequality,
\[
P\left(\sup_{n\geqslant N}\left|X_{n}-X_{N}\right|>\dfrac{\epsilon}{3}\right)\leqslant\dfrac{9}{\epsilon^{2}}\left\Vert \sup_{n\geqslant N}\left|X_{n}-X_{N}\right|\right\Vert ^{2}_{\text{L}^{2}}.
\]

Applying the Doob inequality to the martingale $\left(X_{n}-X_{N}\right)_{n\geqslant N}$---or
equivalently to the martingale $\left(Y_{n}\right)_{n\in\mathbb{N}}$
defined by $Y_{n}=0$ if $0<n\leqslant N-1,$ and $Y_{n}=X_{n}-X_{N},$
if $n\geqslant N$---we obtain
\begin{equation}
P\left(\sup_{n\geqslant N}\left|X_{n}-X_{N}\right|>\dfrac{\epsilon}{3}\right)\leqslant\dfrac{36}{\epsilon^{2}}\sup_{n\geqslant N}\left\Vert X_{n}-X_{N}\right\Vert ^{2}_{\text{L}^{2}}.\label{eq:16_5}
\end{equation}

Since $X$ is a martingale in $\text{L}^{2},$
\begin{align*}
\left\Vert X_{n}-X_{N}\right\Vert ^{2}_{\text{L}^{2}} & =\mathbb{E}\left(\left[X_{n}-X_{N}\right]^{2}\right)=\mathbb{E}\left(X^{2}_{n}\right)+\mathbb{E}\left(X^{2}_{N}\right)-2\mathbb{E}\left(X_{n}X_{N}\right)\\
 & =\mathbb{E}\left(X^{2}_{n}\right)+\mathbb{E}\left(X^{2}_{N}\right)-2\mathbb{E}\left(X_{N}\mathbb{E}^{\mathscr{A}_{N}}\left(X_{N}\right)\right)\\
 & =\mathbb{E}\left(X^{2}_{n}\right)-\mathbb{E}\left(X^{2}_{N}\right).
\end{align*}
Hence, the sequence $\left(\mathbb{E}\left(X^{2}_{n}\right)\right)_{n\in\mathbb{N}}$
is nondecreasing and thus converges, since $X$ is a martingale bounded
in $\text{L}^{2}.$ Consequently,
\[
\sup_{n\geqslant N}\left\Vert X_{n}-X_{N}\right\Vert ^{2}_{\text{L}^{2}}=\sup_{n\geqslant N}\mathbb{E}\left(X^{2}_{n}\right)-\mathbb{E}\left(X^{2}_{N}\right)=\lim_{n\to+\infty}\mathbb{E}\left(X^{2}_{n}\right)-\mathbb{E}\left(X^{2}_{N}\right).
\]
Substituting into the inequality $\refpar{eq:16_5},$
\[
P\left(\sup_{n\geqslant N}\left|X_{n}-X_{N}\right|>\dfrac{\epsilon}{3}\right)\leqslant\dfrac{36}{\epsilon^{2}}\left[\lim_{n\to+\infty}\mathbb{E}\left(X^{2}_{n}\right)-\mathbb{E}\left(X^{2}_{N}\right)\right].
\]

Since the sets $\left(\sup_{n\geqslant N}\left|X_{n}-X_{N}\right|>\dfrac{\epsilon}{3}\right)$
form a nonincreasing sequence in $N,$
\[
P\left(\bigcap_{N\in\mathbb{N}^{\ast}}\left[\sup_{n\geqslant N}\left|X_{n}-X_{N}\right|>\dfrac{\epsilon}{3}\right]\right)\leqslant\dfrac{36}{\epsilon^{2}}\lim_{N\to+\infty}\left[\lim_{n\to+\infty}\mathbb{E}\left(X^{2}_{n}\right)-\mathbb{E}\left(X^{2}_{N}\right)\right]=0.
\]
Therefore,
\[
P\left(\left\{ X\longrightarrow\right\} ^{c}\right)\leqslant\sum_{\epsilon\in\mathbb{Q}^{+}}P\left(\bigcap_{N\in\mathbb{N}^{\ast}}\bigcup_{m,n\geqslant N}\left\{ \left|X_{n}-X_{m}\right|>\epsilon\right\} \right)=0.
\]
Hence the sequence $\left(X_{n}\right)_{n\in\mathbb{N}}$ converges
$P-$almost surely to a random variable $X_{\infty}.$

By the Fatou lemma,
\[
\intop_{\Omega}X^{2}_{\infty}\text{d}P\leqslant\liminf_{n\to+\infty}\mathbb{E}\left(X^{2}_{n}\right)\leqslant\sup_{n\in\mathbb{N}}\mathbb{E}\left(X^{2}_{n}\right)<+\infty,
\]
so $X_{\infty}\in\text{L}^{2}.$

Since $\mathbb{E}\left(\left[X_{n}-X_{m}\right]^{2}\right)=\mathbb{E}\left(X^{2}_{n}\right)-\mathbb{E}\left(X^{2}_{m}\right)$
and since the sequence $\left(\mathbb{E}\left(X^{2}_{n}\right)\right)_{n\in\mathbb{N}}$
converges, the sequence $\left(X_{n}\right)_{n\in\mathbb{N}}$ is
Cauchy in $\text{L}^{2}.$ Hence it converges in $\text{L}^{2}$ to
$X_{\infty}.$ Since for every $m$ and $n$ such that $m\geqslant n,$
$X_{n}=\mathbb{E}^{\mathscr{A}_{n}}\left(X_{m}\right),$ by continuity
of the conditional expectation for the norm $\text{L}^{2},$
\[
X_{n}=\lim_{m\to+\infty}\mathbb{E}^{\mathscr{A}_{n}}\left(X_{m}\right)=\mathbb{E}^{\mathscr{A}_{n}}\left(X_{\infty}\right).
\]
In particular, if the filtration is complete, $X_{\infty}$ is $\mathscr{A}_{\infty}-$measurable
and the martingale $X$ is closable.

\end{proof}

Now, we give two examples of $\text{L}^{2}$ martingales: one bounded
in $\text{L}^{2},$ and the other unbounded in $\text{L}^{2}.$

\begin{example}{}{}

Let $\left(a_{n}\right)_{n\in\mathbb{N}}$ be a sequence of real numbers.
Consider the process $S=\left(S_{n}\right)_{n\in\mathbb{N}}$ defined
for every $n\in\mathbb{N},$ by
\[
S_{n}=\sum^{n}_{i=0}a_{j}X_{j},
\]
where $\left(\text{\ensuremath{X_{n}}}\right)_{n\in\mathbb{N}}$ is
a sequence of independent random variables with common law 
\[
\dfrac{1}{2}\delta_{1}+\dfrac{1}{2}\delta_{-1}.
\]

For $n\in\mathbb{N},$ set
\[
\mathscr{A}_{n}=\sigma\left(X_{j}:\,0\leqslant j\leqslant n\right).
\]
Then, the process $S=\left(S_{n}\right)_{n\in\mathbb{N}}$ is a martingale
with respect to the filtration $\left(\mathscr{A}_{n}\right)_{n\in\mathbb{N}}.$
Indeed, 
\[
\mathbb{E}^{\mathscr{A}_{n}}\left(S_{n+1}\right)=S_{n}+a_{n+1}\mathbb{E}^{\mathscr{A}_{n}}\left(X_{n+1}\right),
\]
and, since the $X_{n}$ are independent and centered,
\[
\mathbb{E}^{\mathscr{A}_{n}}\left(X_{n+1}\right)=\mathbb{E}\left(X_{n+1}\right)=0.
\]
It follows that
\[
\mathbb{E}^{\mathscr{A}_{n}}\left(S_{n+1}\right)=S_{n}.
\]
Of course, $S_{n}\in\text{L}^{2},$ and since the $X_{n}$ are independent,
centered, and of variance 1,
\[
\mathbb{E}\left(S^{2}_{n}\right)=\sigma^{2}_{S_{n}}=\sum^{n}_{j=0}a^{2}_{j}.
\]
If we assume that $\sum^{+\infty}_{j=0}a^{2}_{j}<+\infty,$ then the
martingale $S$ is bounded in $\text{L}^{2}$ and therefore converges
$P-$almost surely and in $\text{L}^{2}.$

If we choose all the $a_{n}$ equal to 1, the martingale $S$ is in
$\text{L}^{2},$ but it is not bounded in $\text{L}^{2}.$ We now
prove that, $P-$almost surely, the sequence $\left(S_{n}\right)_{n\in\mathbb{N}}$
does not converge. It is enough to show that
\begin{equation}
P\left(\limsup_{n\to+\infty}\dfrac{S_{n}}{\sqrt{n}}=+\infty\right)=1.\label{eq:16_6}
\end{equation}

For every $c>0,$ by the Fatou lemma for sets,
\[
\limsup_{n\to+\infty}P\left(\dfrac{S_{n}}{\sqrt{n}}>c\right)\leqslant P\left(\limsup_{n\to+\infty}\left(\dfrac{S_{n}}{\sqrt{n}}>c\right)\right).
\]
Therefore, a fortiori
\[
\limsup_{n\to+\infty}P\left(\dfrac{S_{n}}{\sqrt{n}}>c\right)\leqslant P\left(\limsup_{n\to+\infty}\dfrac{S_{n}}{\sqrt{n}}\geqslant c\right).
\]
But, by the central limit theorem,
\[
\limsup_{n\to+\infty}P\left(\dfrac{S_{n}}{\sqrt{n}}>c\right)=\intop^{+\infty}_{c}\dfrac{1}{\sqrt{2\pi}}\text{e}^{-\frac{x^{2}}{2}}\text{d}x>0.
\]
Hence, for every $c>0,$
\[
P\left(\limsup_{n\to+\infty}\dfrac{S_{n}}{\sqrt{n}}\geqslant c\right)>0,
\]
and thus, by the Kolmogorov zero-one law,
\[
P\left(\limsup_{n\to+\infty}\dfrac{S_{n}}{\sqrt{n}}\geqslant c\right)=1,
\]
since the event $\left(\limsup_{n\to+\infty}\dfrac{S_{n}}{\sqrt{n}}\geqslant c\right)$
is asymptotic. The equality $\refpar{eq:16_6}$ follows immediately
by writing, for instance, that
\[
\left(\limsup_{n\to+\infty}\dfrac{S_{n}}{\sqrt{n}}=+\infty\right)=\bigcup_{p\in\mathbb{N}^{\ast}}\left(\limsup_{n\to+\infty}\dfrac{S_{n}}{\sqrt{n}}\geqslant p\right),
\]
which completes the proof.

\end{example}

\section{Doob Decomposition}\label{sec:Doob-Decomposition}

Given a process $X=\left(X_{n}\right)_{n\in\mathbb{N}},$ we associate
its \textbf{increment process\index{increment process}} $\Delta X=\left(\Delta X_{n}\right)_{n\in\mathbb{N}}$
defined by
\[
\Delta X_{0}=X_{0}
\]
and, for every $n\in\mathbb{N}^{\ast},$
\[
\Delta X_{n}=X_{n}-X_{n-1}.
\]
Then, for every $n\in\mathbb{N},$
\[
X_{n}=\sum^{n}_{j=0}\Delta X_{j}
\]

\begin{definition}{Predictable Process. Predictable nondecreasing Process}{}

(a) A process $X=\left(X_{n}\right)_{n\in\mathbb{N}}$ is \textbf{predictable\mindex{process!predictable}}
if $X_{0}$ is $\mathscr{A}_{0}-$measurable, and, for every $n\in\mathbb{N}^{\ast},$
$X_{n}$ is $\mathscr{A}_{n-1}-$measurable.

(b) A process $A=\left(A_{n}\right)_{n\in\mathbb{N}}$ is \textbf{predictable
nondecreasing\mindex{process!predictable!nondecreasing}} if it is
predicatable, if $A_{0}=0$ and if, for every $n\in\mathbb{N}^{\ast},$
\[
0\leqslant A_{n}\leqslant A_{n+1}<+\infty\,\,\,\,P-\text{almost surely}.
\]
We then denote by $A_{\infty}$ the limit in $\overline{\mathbb{R}}^{+}$
of the sequence $\left(A_{n}\right)_{n\in\mathbb{N}}.$

\end{definition}

\begin{theorem}{Doob Decomposition}{doob_decomp}

Let $X=\left(X_{n}\right)_{n\in\mathbb{N}}$ be an integrable submartingale.

(a) There exists a unique integrable martingale $M=\left(M_{n}\right)_{n\in\mathbb{N}}$
and a unique predictable nondecreasing process $A=\left(A_{n}\right)_{n\in\mathbb{N}}$
such that $X=M+A.$

(b) The following equivalence holds\footnotemark
\[
\sup_{n\in\mathbb{N}}\mathbb{E}\left(X^{+}_{n}\right)<+\infty\Longleftrightarrow\sup_{n\in\mathbb{N}}\mathbb{E}\left(\left|M_{n}\right|\right)<+\infty\,\,\,\,\text{and}\,\,\,\,A_{\infty}\in\mathscr{L}^{1}.
\]

\end{theorem}

\footnotetext{Tr.N. Recall that for a random variable $X,$ $X^{+}=\text{max}\left(0,X\right)$
and $X^{-}=\text{min}\left(0,-X\right).$ Here, only the positive
part needs to be controlled, since the negative part is already under
controlled because $X$ is a submartingale.}

\begin{proof}{}{}

\textbf{(a) Existence}

Define $M$ and $A$ through their increments as follows
\[
M_{0}=X_{0}\,\,\,\,\text{and}\,\,\,\,\forall n\in\mathbb{N}^{\ast},\,\,\,\,\Delta M_{n}=X_{n}-\mathbb{E}^{\mathscr{A}_{n-1}}\left(X_{n}\right).
\]
\[
A_{0}=0\,\,\,\,\text{and}\,\,\,\,\forall n\in\mathbb{N}^{\ast},\,\,\,\,\Delta A_{n}=\mathbb{E}^{\mathscr{A}_{n-1}}\left(X_{n}\right)-X_{n-1}.
\]
Then, 
\[
\mathbb{E}^{\mathscr{A}_{n-1}}\left(\Delta M_{n}\right)=0
\]
so $M$ is indeed an integrable martingale. Moreover, since $X$ is
a submartingale, $\Delta A_{n}\geqslant0.$ Finally, by construction,
$X=M+A.$

\textbf{Unicity}

Let $X=M^{\prime}+A^{\prime}$ be another decomposition, where $M^{\prime}$
is an integrable martingale and $A^{\prime}$ is a predictible increasing
process.

Then
\[
\Delta A^{\prime}_{n}=\Delta X_{n}-\Delta M^{\prime}_{n}.
\]
Since $M^{\prime}$ is a martingale and $A^{\prime}$ is a predictible
and increasing process, we obtain
\[
\Delta A^{\prime}_{n}=\mathbb{E}^{\mathscr{A}_{n-1}}\left(\Delta X_{n}\right)=\Delta A_{n}.
\]
Hence $A=A^{\prime},$ and consequently $M=M^{\prime}.$

\textbf{(b) }
\begin{itemize}
\item Suppose that $\sup_{n\in\mathbb{N}}\mathbb{E}\left(\left|M_{n}\right|\right)<+\infty\,\text{and}\,A_{\infty}\in\mathscr{L}^{1}.$\\
We have
\[
X^{+}_{n}=\left(M_{n}+A_{n}\right)^{+}\leqslant M^{+}_{n}+A_{n},
\]
and therefore
\[
\sup_{n\in\mathbb{N}}\mathbb{E}\left(X^{+}_{n}\right)\leqslant\sup_{n\in\mathbb{N}}\mathbb{E}\left(M^{+}_{n}\right)+\mathbb{E}\left(A_{\infty}\right)<+\infty.
\]
\item Conversely, suppose that 
\[
\sup_{n\in\mathbb{N}}\mathbb{E}\left(X^{+}_{n}\right)<+\infty.
\]
Since $M_{n}\leqslant X_{n},$ we have $M^{+}_{n}\leqslant X^{+}_{n}.$
\\
And thus
\[
\sup_{n\in\mathbb{N}}\mathbb{E}\left(M^{+}_{n}\right)<+\infty.
\]
Moreover 
\[
A_{n}=X_{n}-M_{n}\leqslant X^{+}_{n}-M_{n}
\]
 and since $\mathbb{E}\left(M_{n}\right)=\mathbb{E}\left(M_{0}\right),$
it follows that
\[
\mathbb{E}\left(A_{n}\right)\leqslant\sup_{n\in\mathbb{N}}\mathbb{E}\left(X^{+}_{n}\right)-\mathbb{E}\left(M_{0}\right).
\]
As $A$ is nonnegative and nondecreasing, the Beppo Levi lemma yields
\[
\mathbb{E}\left(A_{\infty}\right)\leqslant\sup_{n\in\mathbb{N}}\mathbb{E}\left(X^{+}_{n}\right)-\mathbb{E}\left(M_{0}\right)<+\infty.
\]
It remains to observe that if $M$ is an integrable martingale, then
\[
\sup_{n\in\mathbb{N}}\mathbb{E}\left(M^{+}_{n}\right)<+\infty\Longleftrightarrow\sup_{n\in\mathbb{N}}\mathbb{E}\left(\left|M_{n}\right|\right)<+\infty.
\]
The implication from right to left follows from
\[
M^{+}_{n}\leqslant\left|M_{n}\right|.
\]
Conversely, implication follows from the fact that
\[
\left|M_{n}\right|=2M^{+}_{n}-M_{n}.
\]
Since $M$ is a martingale,
\[
\mathbb{E}\left(\left|M_{n}\right|\right)=2\mathbb{E}\left(M^{+}_{n}\right)-\mathbb{E}\left(M_{0}\right).
\]
\end{itemize}
\end{proof}

\begin{definition}{Predictable Increasing Process of a Martingale}{}

Let $X$ be an $\text{L}^{2}$-martingale. The predictable nondecreasing
process in the Doob decomposition of the integrable submartingale
$X^{2}$ is called the \textbf{predictable nondecreasing process of
the martingale\mindex{martingale!predictable nondecreasing process}}
$X$ and is denoted by $\left\langle X\right\rangle .$ 

It is the unique predictable nondecreasing process such that $X^{2}-\left\langle X\right\rangle $
is a martingale.

\end{definition}

\begin{remark}{}{}

An $\text{L}^{2}$-martingale $X$ is bounded in $\text{L}^{2}$ if
and only if $\left\langle X\right\rangle _{\infty}$ is integrable.
Moreover
\[
\sup_{n\in\mathbb{N}}\left(\mathbb{E}\left(X^{2}_{n}\right)\right)=\mathbb{E}\left(X^{2}_{0}\right)+\mathbb{E}\left(\left\langle X\right\rangle _{\infty}\right).
\]

\end{remark}

We now state a strong law of large numbers for a martingale of $\text{L}^{2}.$

\begin{theorem}{Strong Law of Large Numbers for a $\textrm{L}^2-$Martingale}{st_law_large_nb_L2-marting}

Let $X$ be an $\text{L}^{2}$-martingale. On the set $\left\{ \left\langle X\right\rangle _{\infty}=+\infty\right\} ,$
the sequence with general term $\left\langle X\right\rangle _{n}$
is eventually nonzero from some rank onward, and the sequence with
general term $\dfrac{X_{n}}{\left\langle X\right\rangle _{n}}$ converges
$P-$almost surely to 0.

\end{theorem}

\begin{proof}{}{}

On the set $\left\{ \left\langle X\right\rangle _{\infty}=+\infty\right\} ,$
the sequence with general term $\left\langle X\right\rangle _{\infty}$
tends to $+\infty$ by nondecreasing. Hence it is nonzero from some,
random, rank onward.

Define a process $Y$ through its increment by
\[
Y_{0}=X_{0}\,\,\,\,\text{and}\,\,\,\,\forall n\in\mathbb{N}^{\ast},\,\,\,\,\Delta Y_{n}=\dfrac{\Delta X_{n}}{1+\left\langle X\right\rangle _{n}}.
\]
 Then $Y$ is a martingale bounded in $\text{L}^{2}.$ 

Indeed: 
\begin{itemize}
\item $Y$ is a martingale, because $\left\langle X\right\rangle _{n}$
is $\mathscr{A}_{n-1}-$measurable. For every $n\in\mathbb{N}^{\ast},$
\[
\mathbb{E}^{\mathscr{A}_{n-1}}\left(\Delta Y_{n}\right)=\dfrac{1}{1+\left\langle X\right\rangle _{n}}\mathbb{E}^{\mathscr{A}_{n-1}}\left(\Delta X_{n}\right)=0.
\]
\item It belongs to $\text{L}^{2}$ since 
\[
\left(\Delta Y_{n}\right)^{2}\leqslant\left(\Delta X_{n}\right)^{2}.
\]
\end{itemize}
Moreover, for every $k\in\mathbb{N}^{\ast},$
\[
\mathbb{E}\left(\left[Y_{k}-Y_{k-1}\right]^{2}\right)=\mathbb{E}\left(Y^{2}_{k}\right)-\mathbb{E}\left(Y^{2}_{k-1}\right),
\]
so that
\[
\mathbb{E}\left(Y^{2}_{n}\right)=\mathbb{E}\left(Y^{2}_{0}\right)+\mathbb{E}\left(\sum^{n}_{k=1}\left(\Delta Y_{k}\right)^{2}\right).
\]
Since $\dfrac{1}{\left(1+\left\langle X\right\rangle _{k}\right)^{2}}$
is $\mathscr{A}_{k-1}-$measurable, we obtain
\[
\mathbb{E}\left(\sum^{n}_{k=1}\left(\Delta Y_{k}\right)^{2}\right)=\sum^{n}_{k=1}\mathbb{E}\left(\dfrac{1}{\left(1+\left\langle X\right\rangle _{k}\right)^{2}}\mathbb{E}^{\mathscr{A}_{k-1}}\left(\left(\Delta X_{k}\right)^{2}\right)\right).
\]
By definition of the process $\left\langle X\right\rangle ,$ and
using its nondecreasing monotonicity,
\[
\mathbb{E}\left(\sum^{n}_{k=1}\left(\Delta Y_{k}\right)^{2}\right)=\sum^{n}_{k=1}\mathbb{E}\left(\dfrac{\Delta\left\langle X\right\rangle _{k}}{\left(1+\left\langle X\right\rangle _{k}\right)^{2}}\right)\leqslant\mathbb{E}\left(\sum^{n}_{k=1}\intop^{\left\langle X\right\rangle _{k}}_{\left\langle X\right\rangle _{k-1}}\dfrac{1}{\left(1+x\right)^{2}}\text{d}x\right).
\]
It follows that
\[
\mathbb{E}\left(\sum^{n}_{k=1}\left(\Delta Y_{k}\right)^{2}\right)\leqslant\intop^{+\infty}_{0}\dfrac{1}{\left(1+x\right)^{2}}\text{d}x<+\infty,
\]
and therefore
\[
\sup_{n\in\mathbb{N}}\mathbb{E}\left(Y^{2}_{n}\right)<+\infty.
\]

By Theorem $\ref{th:cv_th_L2},$ the sequence with general term $Y_{n}$
converges $P-$almost surely and in $\text{L}^{2}.$ By the Kronecker
lemma that, on the set $\left\{ \left\langle X\right\rangle _{\infty}=+\infty\right\} ,$
the sequence with general term $\dfrac{1}{1+\left\langle X\right\rangle _{n}}\left(\sum^{n}_{k=1}\Delta X_{k}\right)$
converges to 0, which yields the result.

\end{proof}

\begin{remark}{}{}This theorem is a generalization, in the $\text{L}^{2}$
setting, of strong laws of large numbers for independent random variables.
These classical results can also be recovered from Theorem $\ref{th:st_law_large_nb_L2-marting}.$
Let us illustrate this in a simple case.

Suppose that $\left(X_{n}\right)_{n\in\mathbb{N}}$ is a sequence
of independent and centered random variables, following the same law,
with finite second-order moment. Define the process $S=\left(S_{n}\right)_{n\in\mathbb{N}}$
by
\[
S_{n}=\sum^{n}_{j=0}X_{j}.
\]
For $n\in\mathbb{N},$ set 
\[
\mathscr{A}_{n}=\sigma\left(X_{j}:\,0\leqslant j\leqslant n\right).
\]
Then the process $S=\left(S_{n}\right)_{n\in\mathbb{N}}$ with respect
to the filtration $\left(\mathscr{A}_{n}\right)_{n\in\mathbb{N}}$
is a martingale, and it belongs to $\text{L}^{2}.$ Let us compute
its predictable nondecreasing process $\left\langle S\right\rangle .$
Recall that
\[
\mathbb{E}^{\mathscr{A}_{n}}\left(S^{2}_{n+1}-S^{2}_{n}\right)=\mathbb{E}^{\mathscr{A}_{n}}\left(\left(\Delta S_{n+1}\right)^{2}\right).
\]
Since the $X_{n}$ are independent and follow the same law,
\[
\mathbb{E}^{\mathscr{A}_{n}}\left(\left(\Delta S_{n+1}\right)^{2}\right)=\mathbb{E}^{\mathscr{A}_{n}}\left(X^{2}_{n+1}\right)=\mathbb{E}\left(X^{2}_{n+1}\right)=\sigma^{2},
\]
where $\sigma^{2}$ denotes the common variance of the $X_{n}.$ Hence
\[
\left\langle S\right\rangle _{n}=n\sigma^{2}.
\]
By the martingale strong law,
\[
\dfrac{1}{n}\sum^{n}_{j=1}X_{j}=\sigma^{2}\dfrac{S_{n}-X_{0}}{\left\langle S_{n}\right\rangle }\stackrel[n\to+\infty]{P-\text{a.s.}}{\longrightarrow}0.
\]

\end{remark}

\section{Convergence of Integrable Martingales}\label{sec:Convergence-of-Integrable}

\begin{definition}{Quadratic Variation Process}{}

Given a process $X=\left(X_{n}\right)_{n\in\mathbb{N}},$ we associate
to it its \textbf{quadratic variation process\index{quadratic variation process}}
$\left[X\right]=\left(\left[X\right]_{n}\right)_{n\in\mathbb{N}}$
defined, for every $n\in\mathbb{N},$ by
\[
\left[X\right]_{n}=\sum^{n}_{j=0}\left(\Delta X_{j}\right)^{2}.
\]

We denote by $\left[X\right]_{\infty}$ the limit in $\overline{\mathbb{R}}^{+}$
of the nondecreasing sequence $\left(X_{n}\right)_{n\in\mathbb{N}}.$

\end{definition}

\begin{remark}{}{}

Let $x=\left(x_{n}\right)_{n\in\mathbb{N}}$ be a sequence of real
numbers. With the previous notations, the condition 
\[
\sum^{+\infty}_{j=0}\left|\Delta x_{j}\right|<+\infty
\]
ensures the convergence of the sequence $x.$ By contrast, as the
following example shows, the sequence $x$ may diverge while its quatratic
variation is finite. This happens for the sequence defined by its
increments: $\Delta x_{0}=0$ and, for every $n\in\mathbb{N}^{\ast},$
$\Delta x_{n}=\dfrac{1}{n}.$

Moreover, the sequence $x$ may converge while its quadratic variation
is infinite. This occurs for the sequence defined by its increments:
$\Delta x_{0}=0$ and, for every $n\in\mathbb{N}^{\ast},$ $\Delta x_{n}=\dfrac{\left(-1\right)^{n}}{\sqrt{n}},$
the convergence of the sequence $x$ coming from the alternating series
Leibniz criterion. The following lemma shows that this situation cannot
happen for martingales bounded in $\text{L}^{1}.$

\end{remark}

\begin{lemma}{Quadratic Variation Limit of a Bounded Martingale in $\text{L}^{1}$}{}

If $X$ is a bounded martingale in $\text{L}^{1},$ then
\[
\left[X\right]_{\infty}<+\infty\,\,\,\,P-\text{almost surely.}
\]

\end{lemma}

\begin{proof}{}{}

(a) For every $n\in\mathbb{N}^{\ast},$
\begin{align*}
\left[X\right]_{n} & =X^{2}_{0}+\sum^{n}_{j=1}\left(X^{2}_{j}+X^{2}_{j-1}-2X_{j}X_{j-1}\right)\\
 & =\sum^{n}_{j=0}X^{2}_{j}+\sum^{n-1}_{j=0}X^{2}_{j}-2\sum^{n}_{j=1}X_{j-1}\left(X_{j}-X_{j-1}\right)-2\sum^{n}_{j=1}X^{2}_{j-1},
\end{align*}
which can be written as
\[
\left[X\right]_{n}=X^{2}_{n}-2\sum^{n}_{j=1}X_{j-1}\Delta X_{j}.
\]
It follows that, for every $n\geqslant2,$
\[
\left[X\right]_{n-1}+X^{2}_{n-1}=2X_{n}X_{n-1}-2\sum^{n}_{j=1}X_{j-1}\Delta X_{j}.
\]
One checks that this equality also holds for $n=1.$ Hence, for every
$n\in\mathbb{N}^{\ast},$
\begin{equation}
\left[X\right]_{n-1}\leqslant2X_{n}X_{n-1}-2\sum^{n}_{j=1}X_{j-1}\Delta X_{j}\label{eq:upper_bound_quad_var_X_n-1}
\end{equation}

(b) Let $\lambda>0,$ and let $T_{\lambda}$ be the stopping time
defined by
\[
T_{\lambda}=\inf\left\{ n\in\mathbb{N}:\,\left|X_{n}\right|>\lambda\right\} ,
\]
with the convention $\inf\emptyset=+\infty.$ 

For every integer $k\geqslant2,$ define the bounded stopping time
$S_{k}=T_{k}\land k.$

Set $\left\Vert X\right\Vert _{1}=\sup_{n\in\mathbb{N}}\mathbb{E}\left(\left|X_{n}\right|\right).$ 

We claim that\boxeq{
\begin{equation}
\mathbb{E}\left(\boldsymbol{1}_{\left(S_{k}\geqslant1\right)}\left[X\right]_{S_{k}-1}\right)\leqslant2\lambda\left\Vert X\right\Vert _{1}.\label{eq:bound_exp_quad_var_X_S_k-1}
\end{equation}
}Indeed, on the set $\left(S_{k}\geqslant1\right),$
\begin{equation}
\sum^{S_{k}}_{j=1}X_{j-1}\Delta X_{j}=\sum^{k}_{j=1}X_{j-1}\boldsymbol{1}_{\left(j\leqslant T_{\lambda}\right)}\Delta X_{j}.\label{eq:sumX_j-1incrementX_j}
\end{equation}
Note that if $j\in\mathbb{N}^{\ast},$ then
\[
\left(j\leqslant T_{\lambda}\right)=\left(T_{\lambda}\leqslant j-1\right)^{c}\in\mathscr{A}_{j-1}
\]
and
\[
\left(S_{k}\geqslant1\right)^{c}=\left(T_{\lambda}=0\right)\in\mathscr{A}_{0}.
\]
Hence the random variable $\boldsymbol{1}_{\left(S_{k}\geqslant1\right)}X_{j-1}\boldsymbol{1}_{\left(j\leqslant T_{\lambda}\right)}$
is $\mathscr{A}_{j-1}-$measurable. Moreover, by definition of $T_{\lambda},$
\[
\left|X_{j-1}\right|\boldsymbol{1}_{\left(j\leqslant T_{\lambda}\right)}\leqslant\lambda.
\]
Since $\Delta X_{j}$ is integrable, it follows that 
\[
\boldsymbol{1}_{\left(S_{k}\geqslant1\right)}\sum^{k}_{j=1}X_{j-1}\boldsymbol{1}_{\left(j\leqslant T_{\lambda}\right)}\Delta X_{j}
\]
is integrable as well. Integrating both sides of $\refpar{eq:sumX_j-1incrementX_j}$
over $\left(S_{k}\geqslant1\right),$ we obtain, since $X$ is a martingale,
\begin{equation}
\mathbb{E}\left(\boldsymbol{1}_{\left(S_{k}\geqslant1\right)}\sum^{S_{k}}_{j=1}X_{j-1}\Delta X_{j}\right)=\sum^{k}_{j=1}\mathbb{E}\left(\boldsymbol{1}_{\left(S_{k}\geqslant1\right)}X_{j-1}\boldsymbol{1}_{\left(j\leqslant T_{\lambda}\right)}\mathbb{E}^{\mathscr{A}_{j-1}}\left(\Delta X_{j}\right)\right)=0.\label{eq:eq16_10}
\end{equation}
Still by definition of $T_{\lambda},$ and thus of $S_{k},$
\begin{equation}
\boldsymbol{1}_{\left(S_{k}\geqslant1\right)}\left|X_{S_{k}}X_{S_{k-1}}\right|\leqslant\boldsymbol{1}_{\left(S_{k}\geqslant1\right)}\lambda\left|X_{S_{k}}\right|.\label{eq:eq16_11}
\end{equation}
Since $S_{k}$ is bounded, the Doob first stopping theorem shows that
$X_{S_{k}}$ is integrable and that 
\[
X_{S_{k}}=\mathbb{E}^{\mathscr{A}_{S_{k}}}\left(X_{k}\right).
\]
Hence,
\begin{equation}
\mathbb{E}\left(\left|X_{S_{k}}\right|\right)\leqslant\mathbb{E}\left(\left|X_{k}\right|\right)\leqslant\left\Vert X\right\Vert _{1}.\label{eq:eq16_12}
\end{equation}
Therefore, from $\refpar{eq:eq16_11},$
\begin{equation}
\mathbb{E}\left(\boldsymbol{1}_{\left(S_{k}\geqslant1\right)}\left|X_{S_{k}}X_{S_{k-1}}\right|\right)\leqslant\lambda\left\Vert X\right\Vert _{1}.\label{eq:eq16_13}
\end{equation}
Finally, using $\refpar{eq:upper_bound_quad_var_X_n-1},$
\begin{equation}
\boldsymbol{1}_{\left(S_{k}\geqslant1\right)}\left[X\right]_{S_{k}-1}\leqslant\boldsymbol{1}_{\left(S_{k}\geqslant1\right)}\left[2\left|X_{S_{k}}X_{S_{k}-1}\right|-2\sum^{S_{k}}_{j=1}X_{j-1}\Delta X_{j}\right].\label{eq:eq16_14}
\end{equation}
Integrating both sides of this inequality yields immediately the claimed
inequality $\refpar{eq:bound_exp_quad_var_X_S_k-1},$ using the relations
$\refpar{eq:eq16_10}$ and $\refpar{eq:eq16_13}.$

(c) Since the sequence with general term $S_{k}$ converges by nondecreasing
to $T_{\lambda},$ the sequence with general term $\boldsymbol{1}_{\left(S_{k}\geqslant1\right)}\left[X\right]_{S_{k}-1}$
converges by nondecreasing to $\boldsymbol{1}_{\left(T_{\lambda}\geqslant1\right)}\left[X\right]_{T_{\lambda}-1}$
because $\left[X\right]_{n}$ is a square sum. Taking limits in $\refpar{eq:bound_exp_quad_var_X_S_k-1}$
and applying the Beppo Levi theorem gives
\begin{equation}
\mathbb{E}\left(\boldsymbol{1}_{\left(T_{\lambda}\geqslant1\right)}\left[X\right]_{T_{\lambda}-1}\right)\leqslant2\lambda\left\Vert X\right\Vert _{1}.\label{eq:eq16_15}
\end{equation}

(d) Since $\left(T_{\lambda}<+\infty\right)\subset\left(X^{\ast}>\lambda\right),$
the maximal lemma yields, for every $\alpha>0$ and $\lambda>0,$
\begin{equation}
P\left(\left(\left[X\right]_{\infty}\geqslant\alpha^{2}\right)\cap\left(T_{\lambda}=+\infty\right)\right)\leqslant P\left(T_{\lambda}<+\infty\right)\leqslant\dfrac{\left\Vert X\right\Vert _{1}}{\lambda}.\label{eq:eq16_16}
\end{equation}

Moreover,
\[
P\left(\left(\left[X\right]_{\infty}\geqslant\alpha^{2}\right)\cap\left(T_{\lambda}=+\infty\right)\right)\leqslant P\left(\left(\boldsymbol{1}_{\left(T_{\lambda}\geqslant1\right)}\left[X\right]_{T_{\lambda}-1}\geqslant\alpha^{2}\right)\cap\left(T_{\lambda}=+\infty\right)\right).
\]
Hence,
\begin{equation}
P\left(\left(\left[X\right]_{\infty}\geqslant\alpha^{2}\right)\cap\left(T_{\lambda}=+\infty\right)\right)\leqslant P\left(\boldsymbol{1}_{\left(T_{\lambda}\geqslant1\right)}\left[X\right]_{T_{\lambda}-1}\geqslant\alpha^{2}\right).\label{eq:eq16_17}
\end{equation}
By the Markov inequality and by $\refpar{eq:eq16_15},$ we obtain
\begin{equation}
P\left(\left(\left[X\right]_{\infty}\geqslant\alpha^{2}\right)\cap\left(T_{\lambda}=+\infty\right)\right)\leqslant\dfrac{2\lambda}{\alpha^{2}}\left\Vert X\right\Vert _{1}.\label{eq:eq16_18}
\end{equation}

Taking $\lambda=\alpha$ in the last inequality, and adding term by
term the inequalities $\refpar{eq:eq16_16}$ and $\refpar{eq:eq16_18},$
we get
\[
P\left(\left[X\right]_{\infty}\geqslant\alpha^{2}\right)\leqslant\dfrac{3}{\alpha}\left\Vert X\right\Vert _{1}.
\]
Since $\alpha$ is arbitrary, this implies $P\left(\left[X\right]_{\infty}=+\infty\right)=0,$
which is the desired result.

\end{proof}

The next theorem concerns martingale convergence and was first proved
by D.L. Burkholder. It is a consequence of the convergence theorem
for the bounded martingales in $\text{L}^{2},$ the previous lemma
on quadratic variation, and the maximal lemma. Its proof follows an
article of Louis H.Y. Chen published in the Proceeding of the AMS
in 1981 \cite{chen1981martingale}.

\begin{denotation}{}{}

If $X=\left(X_{n}\right)_{n\in\mathbb{N}}$ is a process, we denote
\[
X^{\ast}=\sup_{n\in\mathbb{N}}\left|X_{n}\right|.
\]

\end{denotation}

\begin{theorem}{Martingale Convergence}{}

Let $M$ and $N$ be two martingales on the same process basis. Assume
that $M$ is bounded in $\text{L}^{1}.$ If the quadratic variation
process of $M$ and $N$ are such that $\left[N\right]\leqslant\left[M\right],$
then the sequence $\left(N_{n}\right)_{n\in\mathbb{N}}$ converges
$P-$almost surely.

In particular, every martingale bounded in $\text{L}^{1}$ converges
$P-$almost surely.

\end{theorem}

\begin{proof}{}{}

Set $\left\Vert M\right\Vert _{1}=\sup_{n\in\mathbb{N}}\mathbb{E}\left(\left|M_{n}\right|\right).$
\begin{itemize}
\item Let $\lambda>0$ et define the stopping time $T_{\lambda}$ by
\[
T_{\lambda}=\inf\left\{ n\in\mathbb{N}:\,\left|M_{n}\right|>\lambda\,\,\,\,\text{or}\,\,\,\,\left[M\right]_{n}>\lambda^{2}\right\} ,
\]
with the convention $\inf\emptyset=+\infty.$ We prove the following
inequality for the increment process of the stopped martingale $N^{T_{\lambda}}.$
\begin{equation}
\mathbb{E}\left(\left(\Delta N^{T_{\lambda}}\right)^{\ast}\right)\leqslant2\lambda+\left\Vert M\right\Vert _{1}<+\infty.\label{eq:eq16_19}
\end{equation}
For every $n\in\mathbb{N},$
\[
\left(\Delta N^{T_{\lambda}}_{n}\right)^{2}\leqslant\left[N^{T_{\lambda}}\right]_{n}=\left[N\right]_{n\land T_{\lambda}}\leqslant\left[M\right]_{n\land T_{\lambda}}=\boldsymbol{1}_{\left(T_{\lambda}>n\right)}\left[M\right]_{n}+\boldsymbol{1}_{\left(T_{\lambda}\leqslant n\right)}\left[M\right]_{T_{\lambda}}.
\]
Hence,
\begin{equation}
\left|\Delta N^{T_{\lambda}}_{n}\right|\leqslant\boldsymbol{1}_{\left(T_{\lambda}>n\right)}\left[M\right]^{\frac{1}{2}}_{n}+\boldsymbol{1}_{\left(T_{\lambda}\leqslant n\right)}\left[M\right]^{\frac{1}{2}}_{T_{\lambda}}\leqslant\lambda\boldsymbol{1}_{\left(T_{\lambda}>n\right)}+\boldsymbol{1}_{\left(T_{\lambda}\leqslant n\right)}\left[M\right]^{\frac{1}{2}}_{T_{\lambda}}.\label{eq:eq16_20}
\end{equation}
By definition of $T_{\lambda},$ on $\left(T_{\lambda}<+\infty\right),$
\[
\left[M\right]^{\frac{1}{2}}_{T_{\lambda}}=\left[\left[M\right]_{T_{\lambda}-1}+\left(\Delta M_{T_{\lambda}}\right)^{2}\right]^{\frac{1}{2}}\leqslant\lambda+\left|\Delta M_{T_{\lambda}}\right|.
\]
Substituting into $\refpar{eq:eq16_20}$ gives
\begin{equation}
\left|\Delta N^{T_{\lambda}}_{n}\right|\leqslant\lambda+\left|\Delta M_{T_{\lambda}}\right|\boldsymbol{1}_{\left(T_{\lambda}\leqslant n\right)}.\label{eq:eq16_21}
\end{equation}
Therefore
\begin{equation}
\left(\Delta N^{T_{\lambda}}_{n}\right)^{\ast}\leqslant\lambda+\left|\Delta M_{T_{\lambda}}\right|\boldsymbol{1}_{\left(T_{\lambda}<+\infty\right)}.\label{eq:eq16_22}
\end{equation}
Moreover, by the triangle inequality, on $\left(T_{\lambda}<+\infty\right),$
\[
\left|\Delta M_{T_{\lambda}}\right|\leqslant\left|M_{T_{\lambda}-1}\right|+\left|M_{T_{\lambda}}\right|\leqslant\lambda+\left|M_{T_{\lambda}}\right|,
\]
so that
\begin{equation}
\left(\Delta N^{T_{\lambda}}_{n}\right)^{\ast}\leqslant2\lambda+\left|M_{T_{\lambda}}\right|\boldsymbol{1}_{\left(T_{\lambda}<+\infty\right)}.\label{eq:eq16_23}
\end{equation}
Hence, by integrating,
\begin{equation}
\mathbb{E}\left(\left(\Delta N^{T_{\lambda}}_{n}\right)^{\ast}\right)\leqslant2\lambda+\mathbb{E}\left(\left|M_{T_{\lambda}}\right|\boldsymbol{1}_{\left(T_{\lambda}<+\infty\right)}\right).\label{eq:eq16_24}
\end{equation}
It remains to bound the right-hand side. Since
\[
\lim_{n\to+\infty}\boldsymbol{1}_{\left(T_{\lambda}<+\infty\right)}\left|M_{T_{\lambda}\land n}\right|=\boldsymbol{1}_{\left(T_{\lambda}<+\infty\right)}\left|M_{T_{\lambda}}\right|.
\]
The Fatou lemma and the inequality $\refpar{eq:eq16_24}$ yield
\[
\mathbb{E}\left(\left(\Delta N^{T_{\lambda}}_{n}\right)^{\ast}\right)\leqslant2\lambda+\liminf_{n\to+\infty}\mathbb{E}\left(\boldsymbol{1}_{\left(T_{\lambda}<+\infty\right)}\left|M_{T_{\lambda}\land n}\right|\right)\leqslant2\lambda+\sup_{n\in\mathbb{N}}\mathbb{E}\left(\left|M_{T_{\lambda}\land n}\right|\right).
\]
But $T_{\lambda}\land n$ is a bounded stopping time, so by Doob first
stopping theorem, 
\[
M_{T_{\lambda}\land n}=\mathbb{E}^{\mathscr{A}_{T_{\lambda}\land n}}\left(M_{n}\right),
\]
and thus
\[
\mathbb{E}\left(\left|M_{T_{\lambda}\land n}\right|\right)\leqslant\mathbb{E}\left(\left|M_{n}\right|\right)\leqslant\left\Vert M\right\Vert _{1},
\]
This proves the inequality $\refpar{eq:eq16_19}.$
\item Let $U=1+\left[N^{T_{\lambda}}\right].$ Define the process $Y^{\lambda},$
as sum normalized by $U$ of the process of $N^{T_{\lambda}}$ increments.
For every $n\in\mathbb{N},$
\[
\Delta Y^{\lambda}_{n}=\dfrac{\Delta N^{T_{\lambda}}_{n}}{U_{n}}.
\]
For every $n\in\mathbb{N}^{\ast},$
\[
\left(\dfrac{\Delta N^{T_{\lambda}}_{n}}{U_{n}}\right)^{2}=\dfrac{\Delta U_{n}}{U^{2}_{n}}\leqslant\intop^{U_{n}}_{U_{n-1}}\dfrac{\text{d}x}{x^{2}}.
\]
Hence,
\[
\sum^{+\infty}_{n=1}\left(\dfrac{\Delta N^{T_{\lambda}}_{n}}{U_{n}}\right)^{2}\leqslant\sum^{+\infty}_{n=1}\intop^{U_{n}}_{U_{n-1}}\dfrac{\text{d}x}{x^{2}}\leqslant\intop^{+\infty}_{1}\dfrac{\text{d}x}{x^{2}}=1,
\]
which implies
\begin{equation}
\sum^{+\infty}_{n=1}\mathbb{E}\left(\left(\Delta Y^{\lambda}_{n}\right)^{2}\right)=\mathbb{E}\left(\sum^{+\infty}_{n=1}\left(\dfrac{\Delta N^{T_{\lambda}}_{n}}{U_{n}}\right)^{2}\right)\leqslant1.\label{eq:eq16_25}
\end{equation}
\item \textbf{We now prove that the sequence $\left(Y^{\lambda}_{n}\right)_{n\in\mathbb{N}}$
converges $P-$almost surely.} Introduce the $\text{L}^{2}$-martingale
$Z$ defined by its increments
\[
\Delta Z_{0}=0
\]
and for $n\in\mathbb{N}^{\ast},$
\[
\Delta Z_{n}=\Delta Y^{\lambda}_{n}-\mathbb{E}^{\mathscr{A}_{n-1}}\left(\Delta Y^{\lambda}_{n}\right).
\]
Then $Z$ is a \textbf{bounded martingale in $\text{L}^{2}.$} Indeed,
by a classical computation of the conditional covariance, for every
$n\geqslant1,$
\begin{align*}
\mathbb{E}\left(\left(\Delta Z_{n}\right)^{2}\right) & =\mathbb{E}\left(\left(\Delta Y^{\lambda}_{n}\right)^{2}\right)-2\mathbb{E}\left(\Delta Y^{\lambda}_{n}\left(\mathbb{E}^{\mathscr{A}_{n-1}}\left(\Delta Y^{\lambda}_{n}\right)\right)\right)+\mathbb{E}\left(\left(\mathbb{E}^{\mathscr{A}_{n-1}}\left(\Delta Y^{\lambda}_{n}\right)\right)^{2}\right)\\
 & =\mathbb{E}\left(\left(\Delta Y^{\lambda}_{n}\right)^{2}\right)-\mathbb{E}\left(\left(\mathbb{E}^{\mathscr{A}_{n-1}}\left(\Delta Y^{\lambda}_{n}\right)\right)^{2}\right)\leqslant\mathbb{E}\left(\left(\Delta Y^{\lambda}_{n}\right)^{2}\right).
\end{align*}
Since $Z$ is an $\text{L}^{2}$-martingale,
\[
\mathbb{E}\left(Z^{2}_{n}\right)-\mathbb{E}\left(Z^{2}_{n-1}\right)=\mathbb{E}\left(\left(\Delta Z_{n}\right)^{2}\right),
\]
so
\[
\mathbb{E}\left(Z^{2}_{n}\right)=\mathbb{E}\left(Z^{2}_{0}\right)+\sum^{n}_{j=1}\mathbb{E}\left(\left(\Delta Z_{j}\right)^{2}\right)\leqslant\mathbb{E}\left(Z^{2}_{0}\right)+\sum^{n}_{j=1}\mathbb{E}\left(\left(\Delta Y^{\lambda}_{j}\right)^{2}\right).
\]
By inequality $\refpar{eq:eq16_25},$
\[
\sup_{n\in\mathbb{N}^{\ast}}\mathbb{E}\left(Z^{2}_{n}\right)\leqslant\mathbb{E}\left(Z^{2}_{0}\right)+1.
\]
Hence, \textbf{the sequence $\left(Z_{n}\right)_{n\in\mathbb{N}}$
converges $P-$almost surely}.\\
Now, we show that the series with general term $\mathbb{E}^{\mathscr{A}_{n-1}}\left(\Delta Y^{\lambda}_{n}\right)$
converges absolutely $P-$almost surely. Since $U_{n-1}$ is $\mathscr{A}_{n-1}-$measurable
and $N^{T_{\lambda}}$ is a martingale,
\[
\mathbb{E}^{\mathscr{A}_{n-1}}\left(\dfrac{\Delta N^{T_{\lambda}}_{n}}{U_{n-1}}\right)=\dfrac{1}{U_{n-1}}\mathbb{E}^{\mathscr{A}_{n-1}}\left(\Delta N^{T_{\lambda}}_{n}\right)=0.
\]
Therefore,
\begin{align*}
\mathbb{E}\left(\sum^{+\infty}_{n=1}\left|\mathbb{E}^{\mathscr{A}_{n-1}}\left(\Delta Y^{\lambda}_{n}\right)\right|\right) & =\mathbb{E}\left(\sum^{+\infty}_{n=1}\left|\mathbb{E}^{\mathscr{A}_{n-1}}\left(\Delta Y^{\lambda}_{n}\right)-\mathbb{E}^{\mathscr{A}_{n-1}}\left(\dfrac{\Delta N^{T_{\lambda}}_{n}}{U_{n-1}}\right)\right|\right)\\
 & =\sum^{+\infty}_{n=1}\mathbb{E}\left(\left|\mathbb{E}^{\mathscr{A}_{n-1}}\left(\Delta N^{T_{\lambda}}_{n}\left[\dfrac{1}{U_{n-1}}-\dfrac{1}{U_{n}}\right]\right)\right|\right),
\end{align*}
Since $\dfrac{1}{U_{n-1}}-\dfrac{1}{U_{n}}$ is nonnegative, we get
the upper-bound
\[
\mathbb{E}\left(\sum^{+\infty}_{n=1}\left|\mathbb{E}^{\mathscr{A}_{n-1}}\left(\Delta Y^{\lambda}_{n}\right)\right|\right)\leqslant\sum^{+\infty}_{n=1}\mathbb{E}\left(\left|\Delta N^{T_{\lambda}}_{n}\right|\left[\dfrac{1}{U_{n-1}}-\dfrac{1}{U_{n}}\right]\right),
\]
and hence
\begin{equation}
\mathbb{E}\left(\sum^{+\infty}_{n=1}\left|\mathbb{E}^{\mathscr{A}_{n-1}}\left(\Delta Y^{\lambda}_{n}\right)\right|\right)\leqslant\mathbb{E}\left(\left(\Delta N^{T_{\lambda}}\right)^{\ast}\sum^{+\infty}_{n=1}\left[\dfrac{1}{U_{n-1}}-\dfrac{1}{U_{n}}\right]\right).\label{eq:eq16_26}
\end{equation}
However, for $n\geqslant1,$
\[
\sum^{n}_{j=1}\left(\dfrac{1}{U_{j-1}}-\dfrac{1}{U_{j}}\right)=\dfrac{1}{U_{0}}-\dfrac{1}{U_{n}}\leqslant1,
\]
which implies, since $U_{n}$ is nonnegative, that $\sum^{+\infty}_{n=1}\left(\dfrac{1}{U_{n-1}}-\dfrac{1}{U_{n}}\right)\leqslant1.$
The inequalities $\refpar{eq:eq16_26}$ and $\refpar{eq:eq16_19}$
then yield
\[
\mathbb{E}\left(\sum^{+\infty}_{n=1}\left|\mathbb{E}^{\mathscr{A}_{n-1}}\left(\Delta Y^{\lambda}_{n}\right)\right|\right)\leqslant\mathbb{E}\left(\left(\Delta N^{T_{\lambda}}\right)^{\ast}\right)<+\infty,
\]
and it follows that \textbf{the series with general term $\mathbb{E}^{\mathscr{A}_{n-1}}\left(\Delta Y^{\lambda}_{n}\right)$
$P-$almost surely absolutely converges}. Finally, for every $n\in\mathbb{N}^{\ast},$
\[
Y^{\lambda}_{n}=Y^{\lambda}_{0}+Z_{n}+\sum^{n}_{j=1}\mathbb{E}^{\mathscr{A}_{j-1}}\left(\Delta Y^{\lambda}_{j}\right).
\]
Since the sequence $\left(Z_{n}\right)_{n\in\mathbb{N}}$ converges,
then $P-$almost surely, \textbf{the sequence $\left(Y^{\lambda}_{n}\right)_{n\in\mathbb{N}}$
converges $P-$almost surely}.
\item Define now the process $Y,$ for every $n\in\mathbb{N},$ by
\[
Y_{n}=\sum^{n}_{j=0}\dfrac{\Delta N_{j}}{U_{j}}.
\]
For every $\lambda>0,$ the process $Y$ coincides on $\left(T_{\lambda}=+\infty\right)$
with the process $Y^{\lambda}.$ Hence the sequence $\left(Y_{n}\right)_{n\in\mathbb{N}}$
converges $P-$almost surely on $\bigcup_{\lambda\in\mathbb{Q}^{+}}\left(T_{\lambda}=+\infty\right).$
But,
\begin{align*}
\bigcup_{\lambda\in\mathbb{Q}^{+}}\left(T_{\lambda}=+\infty\right) & =\bigcup_{\lambda\in\mathbb{Q}^{+}}\left[\left(\left[M\right]_{\infty}\leqslant\lambda^{2}\right)\cap\left(M^{\ast}\leqslant\lambda\right)\right]\\
 & =\left(\left[M\right]_{\infty}\leqslant+\infty\right)\cap\left(M^{\ast}\leqslant+\infty\right).
\end{align*}
Since the martingale $M$ is bounded in $\text{L}^{1},$ the lemma
on quadratic variation and the maximal lemma imply
\[
P\left(\bigcup_{\lambda\in\mathbb{Q}^{+}}\left(T_{\lambda}=+\infty\right)\right)=1.
\]
Therefore, \textbf{the sequence $\left(Y_{n}\right)_{n\in\mathbb{N}}$
converges $P-$almost surely}. Finally, since we have the inclusion
of sets 
\[
\left(\left[M\right]_{\infty}<+\infty\right)\subset\left(\left[N\right]_{\infty}<+\infty\right),
\]
the sequence $\left(U_{n}\right)_{n\in\mathbb{N}}$ converges $P-$almost
surely to a finite limit. A standard result from real analysis then
implies that the sequence $\left(\sum^{n}_{j=0}\Delta N_{j}\right)_{n\in\mathbb{N}}$
converges $P-$almost surely, i.e. the sequence $\left(N_{n}\right)_{n\in\mathbb{N}}$
converges $P-$almost surely. This completes the proof.
\end{itemize}
\end{proof}

\begin{corollary}{}{}

Let $X$ be a bounded martingale in $\text{L}^{1},$ and let $T$
be a stopping time. Then the stopped martingale $X^{T}$ converges
$P-$almost surely.

\end{corollary}

\begin{proof}{}{}

It suffices to prove that the martingale $X^{T}$ is bounded in $\text{L}^{1}.$
However, by the Doob first stopping theorem, for every $n\in\mathbb{N},$
$X_{T\land n}=\mathbb{E}^{\mathscr{A}_{T\land n}}\left(X_{n}\right),$
so
\begin{align*}
\mathbb{E}\left(\left|X^{T}_{n}\right|\right) & =\mathbb{E}\left(\left|X_{T\land n}\right|\right)=\mathbb{E}\left(\left|\mathbb{E}^{\mathscr{A}_{T\land n}}\left(X_{n}\right)\right|\right)\leqslant\mathbb{E}\left(\mathbb{E}^{\mathscr{A}_{T\land n}}\left(\left|X_{n}\right|\right)\right)\\
 & =\mathbb{E}\left(\left|X_{n}\right|\right)\leqslant\sup_{n\in\mathbb{N}}\mathbb{E}\left(\left|X_{n}\right|\right)<+\infty,
\end{align*}
Hence,
\[
\sup_{n\in\mathbb{N}}\mathbb{E}\left(\left|X^{T}_{n}\right|\right)<+\infty.
\]

\end{proof}

\begin{remark}{}{}

The following counter-example shows that a bounded martingale in $\text{L}^{1}$
may possibly not converge in $\text{L}^{1}.$

Let $\left(X_{n}\right)_{n\in\mathbb{N}}$ be a sequence of independent
random variables, all following the same law $\dfrac{1}{2}\left(\delta_{0}+\delta_{2}\right).$
For every $n\in\mathbb{N},$ define $Y_{n}=\prod^{n}_{j=0}X_{j}$
and $\mathscr{A}_{n}=\sigma\left(X_{j}:\,0\leqslant j\leqslant n\right).$ 

We have 
\[
\mathbb{E}^{\mathscr{A}_{n}}\left(Y_{n+1}\right)=Y_{n}\mathbb{E}^{\mathscr{A}_{n}}\left(X_{n+1}\right).
\]

Since $X_{n+1}$ is independent of $\mathscr{A}_{n},$ it follows
that 
\[
\mathbb{E}^{\mathscr{A}_{n}}\left(Y_{n+1}\right)=Y_{n}\mathbb{E}\left(X_{n+1}\right)=Y_{n}.
\]
Hence the process $Y$ is a martingale. Moreover, it is bounded in
$\text{L}^{1}.$ Indeed, by independence, for every $n\in\mathbb{N},$
\[
\mathbb{E}\left(\left|Y_{n}\right|\right)=\prod^{n}_{j=0}\mathbb{E}\left(\left|X_{j}\right|\right)=1.
\]
Thus, the sequence $\left(Y_{n}\right)_{n\in\mathbb{N}}$ converges
$P-$almost surely to a random variable $Y_{\infty}.$ Observe that
$Y_{n}$ takes only the values 0 and $2^{n+1}$ $P-$almost surely,
and that 
\[
P\left(Y_{n}=2^{n+1}\right)=2^{-\left(n+1\right)}.
\]
It follows that the sequence $\left(Y_{n}\right)_{n\in\mathbb{N}}$
converges in probability to 0, and therefore $Y_{\infty}=0.$ However,
the convergence cannot take place in $\text{L}^{1},$ since $\mathbb{E}\left(Y_{n}\right)=1.$
It is however easy to verify directly that this sequence is not equi-integrable. 

\end{remark}

The $\text{L}^{1}$-convergence of an integrable martingale is characterized
by the following proposition.

\begin{proposition}{$\textrm{L}^1$ Convergence of an Integrable Martingale}{L1conv_int_mart}

Let $X=\left(X_{n}\right)_{n\in\mathbb{N}}$ be an integrable martingale.
The sequence $\left(X_{n}\right)_{n\in\mathbb{N}}$ is equi-integrable
if and only if the martingale is closable. In this case, the sequence
$\left(X_{n}\right)_{n\in\mathbb{N}}$ converges $P-$almost surely
and in $\text{L}^{1}$ to an $\mathscr{A}_{\infty}-$measurable random
variable $X_{\infty}.$

\end{proposition}

\begin{proof}{}{}

The sequence is equi-integrable. In particular, it is bounded in $\text{L}^{1},$
hence the martingale $X$ converges $P-$almost surely, and by equi-integrability,
the convergence also holds in $\text{L}^{1}.$ 

Set 
\[
X_{\infty}=\limsup_{n\to+\infty}X_{n}.
\]
Since $\left(X_{n}\right)_{n\in\mathbb{N}}$ is adapted, the limit
$X_{\infty}$ is $\mathscr{A}_{\infty}-$measurable.

Moreover, for $n\leqslant p,$ we have 
\[
X_{n}=\mathbb{E}^{\mathscr{A}_{n}}\left(X_{p}\right).
\]
Passing to the limit in $\text{L}^{1},$ we obtain 
\[
X_{n}=\mathbb{E}^{\mathscr{A}_{n}}\left(X_{\infty}\right).
\]

Conversely, let $X_{\infty}$ be an $\mathscr{A}_{\infty}-$measurable
random variable such that, for every $n\in\mathbb{N},$ 
\[
X_{n}=\mathbb{E}^{\mathscr{A}_{n}}\left(X_{\infty}\right).
\]
 The equi-integrability of the sequence $\left(X_{n}\right)_{n\in\mathbb{N}}$
then follows from the following general lemma.

\end{proof}

\begin{lemma}{Equi-Integrability of the Conditional Expectation of a Random Variable in $\textrm{L}^1$}{eq_int_cond_exp}

Let $X\in\text{L}^{1}\left(\Omega,\mathscr{A},P\right),$ and let
$\left(\mathscr{A}_{i}\right)_{i\in I}$ be a family of sub-$\sigma$-algebra
of $\mathscr{A}.$ Note $X_{i}=\mathbb{E}^{\mathscr{A}_{i}}\left(X\right).$ 

Then the family $\left(X_{i}\right)_{i\in I}$ is equi-integrable.

\end{lemma}

\begin{proof}{}{}

Let $a>0$ be arbitrary. Since 
\[
\left|X_{i}\right|\leqslant\mathbb{E}^{\mathscr{A}_{i}}\left(\left|X\right|\right),
\]
 and since $\left(\left|X_{i}\right|>a\right)\in\mathscr{A}_{i},$
for every $i\in I,$
\[
P\left(\left|X_{i}\right|>a\right)\leqslant\dfrac{1}{a}\intop\left|X_{i}\right|\text{d}P\leqslant\dfrac{1}{a}\intop\left|X\right|\text{d}P.
\]
Moreover,
\begin{equation}
\sup_{i\in I}\intop_{\left(\left|X_{i}\right|>a\right)}\left|X_{i}\right|\text{d}P\leqslant\sup_{i\in I}\intop_{\left(\left|X_{i}\right|>a\right)}\left|X\right|\text{d}P.\label{eq:eq16_27}
\end{equation}
Thus, for every $\eta>0,$ there exists $A\left(\eta\right)>0$ such
that 
\[
\sup_{i\in I}P\left(\left|X_{i}\right|>a\right)\leqslant\eta
\]
whenever $a\geqslant A\left(\eta\right).$ 

On the other hand, for every $\epsilon>0,$ there exists $\eta\left(\epsilon\right)>0$
such that $\intop_{A}\left|X\right|\text{d}P\leqslant\epsilon$ as
soon as $P\left(A\right)\leqslant\eta\left(\epsilon\right).$ 

It then follows from $\refpar{eq:eq16_27}$ that, if $a\geqslant A\left(\eta\left(\epsilon\right)\right),$
then 
\[
\sup_{i\in I}\intop_{\left(\left|X_{i}\right|>a\right)}\left|X_{i}\right|\text{d}P\leqslant\epsilon,
\]
which proves the announced result.

\end{proof}

\section{Second Stopping Theorem}\label{sec:Second-Stopping-Theorem}

We now prove a characterization theorem for closed martingale in terms
of arbitrary stopping times, from which we deduce Doob second stopping
theorem.

\begin{theorem}{Closed Martingale Characterization with Arbitrary Stopping Times}{cl_mart_ch_arb_stopping_time}

Suppose that $\mathscr{A}_{\infty}=\bigvee_{n\in\mathbb{N}}\mathscr{A}_{n}.$
Let $X=\left(X_{n}\right)_{n\in\mathbb{N}}$ be an adapted process.
The following properties are equivalent:

(i) $X$ is a closed martingale.

(ii) For every $T\in\mathscr{T},$ $X\in\text{L}^{1}\left(\Omega,\mathscr{A},P\right)$
and $\mathbb{E}\left(X_{T}\right)=\mathbb{E}\left(X_{0}\right).$

(iii) The process $\left(X_{T}\right)_{T\in\mathscr{T}}$ is a martingale
with respect to the filtration $\left(\mathscr{A}_{T}\right)_{T\in\mathscr{T}}.$

\end{theorem}

\begin{proof}{}{}
\begin{itemize}
\item \textbf{$\text{(i)}\Rightarrow\text{(ii)}$}\\
Assume that $X$ is a closed martingale. For every stopping time $T$
bounded by an integer $k,$ the first Doob stopping theorem states
that $X_{T}=\mathbb{E}^{\mathscr{A}_{T}}\left(X_{k}\right).$ Since
$X$ is closed, then $X_{k}=\mathbb{E}^{\mathscr{A}_{k}}\left(X_{\infty}\right).$
Because $\mathscr{A}_{T}\subset\mathscr{A}_{k},$ we obtain
\[
X_{T}=\mathbb{E}^{\mathscr{A}_{T}}\left(\mathbb{E}^{\mathscr{A}_{k}}\left(X_{\infty}\right)\right)=\mathbb{E}^{\mathscr{A}_{T}}\left(X_{\infty}\right).
\]
By Lemma $\ref{lm:eq_int_cond_exp},$ the family of random variables
\[
\left\{ X_{T}:\,T\in\mathscr{T}_{b}\right\} 
\]
is \textbf{equi-integrable}. \\
Moreover, $\mathbb{E}\left(X_{T}\right)=\mathbb{E}\left(X_{0}\right)$
for every $T\in\mathscr{T}_{b}.$ \\
Now let $T$ be an arbitrary stopping time. The family of random variables
\[
\left\{ X_{T\land n}:\,n\in\mathbb{N}\right\} 
\]
is equi-integrable. The stopped martingale $X^{T}$ converges $P-$almost
surely, and for every $n\in\mathbb{N},$
\[
\mathbb{E}\left(X_{T\land n}\right)=\mathbb{E}\left(X_{0}\right).
\]
Since $T=\lim_{n\to+\infty}T\land n,$ we have on $\left(T<+\infty\right),$
\[
\lim_{n\to+\infty}X_{T\land n}=X_{T}.
\]
Moreover, on $\left(T=+\infty\right),$
\[
X_{T\land n}=X_{n}.
\]
Since, by Proposition $\ref{pr:L1conv_int_mart},$
\[
\lim_{n\to+\infty}X_{n}=X_{\infty}\,\,\,\,P-\text{almost surely}
\]
Hence, $P-$almost surely on $\left(T=+\infty\right),$
\[
\lim_{n\to+\infty}X_{T\land n}=X_{\infty}.
\]
Consequently,
\[
\lim_{n\to+\infty}X_{T\land n}=X_{\infty}\,\,\,\,P-\text{almost surely,}
\]
and by equi-integrability,
\[
\lim_{n\to+\infty}\mathbb{E}\left(X_{T\land n}\right)=\mathbb{E}\left(X_{\infty}\right).
\]
It follows that for every stopping time $T,$
\[
\mathbb{E}\left(X_{T}\right)=\mathbb{E}\left(X_{0}\right).
\]
\item \textbf{$\text{(ii)}\Rightarrow\text{(iii)}$}\\
Let $S\in\mathscr{T}.$ We first prove that $X_{s}=\mathbb{E}^{\mathscr{A}_{s}}\left(X_{\infty}\right).$
\\
Let $A\in\mathscr{A}_{S}.$ The function $R$ defined by
\[
R=S\boldsymbol{1}_{A}+\left(+\infty\right)\boldsymbol{1}_{A^{c}}
\]
is a stopping time. Indeed, for every $n\in\mathbb{N},$
\[
\left(R=n\right)=\left(S=n\right)\cap A\in\mathscr{A}_{n}.
\]
Applying the hypothesis to the stopping time $R$ and $+\infty,$
we obtain 
\[
\mathbb{E}\left(X_{R}\right)=\mathbb{E}\left(X_{\infty}\right).
\]
This yields
\[
\mathbb{E}\left(S\boldsymbol{1}_{A}+\boldsymbol{1}_{A^{c}}X_{\infty}\right)=\mathbb{E}\left(X_{\infty}\right).
\]
Hence,
\[
\mathbb{E}\left(\boldsymbol{1}_{A}X_{S}\right)=\mathbb{E}\left(\boldsymbol{1}_{A}X_{\infty}\right).
\]
Since $X_{S}$ is $\mathscr{A}_{S}-$measurable, we conclude that
$X_{S}=\mathbb{E}^{\mathscr{A}_{S}}\left(X_{\infty}\right).$\\
Now, if $S$ and $T$ are two stopping times such that $S\leqslant T,$
then similarly 
\[
X_{T}=\mathbb{E}^{\mathscr{A}_{T}}\left(X_{\infty}\right).
\]
Because $\mathscr{A}_{S}\subset\mathscr{A}_{T},$ we obtain
\[
\mathbb{E}^{\mathscr{A}_{T}}\left(X_{T}\right)=\mathbb{E}^{\mathscr{A}_{S}}\left(\mathbb{E}^{\mathscr{A}_{T}}\left(X_{\infty}\right)\right)=\mathbb{E}^{\mathscr{A}_{S}}\left(X_{\infty}\right)=X_{S}.
\]
Thus the property (iii) holds.
\item \textbf{$\text{(iii)}\Rightarrow\text{(i)}$}\\
It is enough to consider constant stopping times, possibly equal to
$+\infty.$
\end{itemize}
\end{proof}

An explicit formulation of the implication $\text{(i)}\Rightarrow\text{(iii)}$
in Theorem $\ref{th:cl_mart_ch_arb_stopping_time}$ yields the Doob
second stopping theorem.

\begin{theorem}{Second Stopping Theorem of Doob}{}

Let $X=\left(X_{n}\right)_{n\in\overline{\mathbb{N}}}$ be a closed
martingale. For every stopping time $S$ and $T$ such that $S\leqslant T,$
\[
\mathbb{E}^{\mathscr{A}_{S}}\left(X_{T}\right)=X_{S}.
\]

\end{theorem}

\section{Convergence of Submartingales and Supermartingales }\label{sec:Convergence-of-Sub-}

We deduce the convergence theorems for submartingales and supermartingales
from the convergence theorem of $\text{L}^{1}$ bounded martingales
and from the Doob decomposition in submartingales.

\begin{theorem}{Submartingale Convergence}{}

Let $X=\left(X_{n}\right)_{n\in\mathbb{N}}$ be a submartingale such
that
\[
\sup_{n\in\mathbb{N}}\mathbb{E}\left(X^{+}_{n}\right)<+\infty.
\]
Then the sequence $\left(X_{n}\right)_{n\in\mathbb{N}}$ converges
$P-$almost surely.

\end{theorem}

\begin{proof}{}{}

By Theorem $\ref{th:doob_decomp}$ and the hypothesis, $X$ admits
the Doob decomposition $X=M+A,$ where $M$ is a martingale bounded
in $\text{L}^{1},$ hence converging $P-$almost surely, and $A$
is a nondecreasing predictable process such that $A_{\infty}\in\mathscr{L}^{1},$
and therefore $P-$almost surely. The sequence $\left(A_{n}\right)_{n\in\mathbb{N}}$
and thus also the sequence $\left(X_{n}\right)_{n\in\mathbb{N}}$
converge $P-$almost surely.

\end{proof}

\begin{remark}{}{}

Let $X=\left(X_{n}\right)_{n\in\mathbb{N}}$ be a submartingale with
Doob decomposition $X=M+A.$ The sequence $\left(X_{n}\right)_{n\in\mathbb{N}}$
converges in $\text{L}^{1}$ if and only if the sequence $\left(M_{n}\right)_{n\in\mathbb{N}}$
is equi-integrable and $A_{\infty}\in\mathscr{L}^{1}.$ Indeed, suppose
that the sequence $\left(X_{n}\right)_{n\in\mathbb{N}}$ converges
in $\text{L}^{1}$ to $X_{\infty}.$ Then
\begin{align*}
\mathbb{E}\left(A_{\infty}\right) & =\lim_{n\to+\infty}\nearrow\mathbb{E}\left(A_{n}\right)=\lim_{n\to+\infty}\left[\mathbb{E}\left(X_{n}\right)-\mathbb{E}\left(M_{n}\right)\right]\\
 & =\lim_{n\to+\infty}\left[\mathbb{E}\left(X_{n}\right)-\mathbb{E}\left(M_{0}\right)\right]=\mathbb{E}\left(X_{\infty}\right)-\mathbb{E}\left(X_{0}\right)<+\infty,
\end{align*}
which proves that $A_{\infty}\in\mathscr{L}^{1}.$ 

Since the sequence $\left(A_{n}\right)_{n\in\mathbb{N}}$ converges
by nondecreasing to $A_{\infty},$ it also converges in $\text{L}^{1}.$
It follows that the sequence $\left(M_{n}\right)_{n\in\mathbb{N}}$
converges in $\text{L}^{1},$ and therefore this sequence is equi-integrable.

Conversely, if the martingale $M$ is equi-integrable and $A_{\infty}\in\mathscr{L}^{1},$
then the sequences $\left(M_{n}\right)_{n\in\mathbb{N}}$ and $\left(A_{n}\right)_{n\in\mathbb{N}}$
converges in $\text{L}^{1}$ respectively to $M_{\infty}$ and $A_{\infty}.$
Consequently, the sequence $\left(X_{n}\right)_{n\in\mathbb{N}}$
converges in $\text{L}^{1}$ to $M_{\infty}+A_{\infty}.$

\end{remark}

\begin{corollary}{}{}

Let $X=\left(X_{n}\right)_{n\in\mathbb{N}}$ be a nonnegative submartingale.
Then the sequence $\left(X_{n}\right)_{n\in\mathbb{N}}$ converges
$P-$almost surely to a random variable $X_{\infty}$ taking values
in $\overline{\mathbb{R}}^{+},$ and for every $n\in\mathbb{N},$
\[
X_{n}\geqslant\mathbb{E}^{\mathscr{A}_{n}}\left(X_{\infty}\right).
\]

\end{corollary}

\begin{proof}{}{}

There exists $c>0,$ such that for every $n\in\mathbb{N},$ 
\[
0\leqslant X_{n}\leqslant c.
\]
Then $-X$ is an integrable submartingale satisfying 
\[
\sup_{n\in\mathbb{N}}\mathbb{E}\left(-X_{n}\right)^{+}<+\infty.
\]
By the previous theorem, the sequence $\left(X_{n}\right)_{n\in\mathbb{N}}$
converges $P-$almost surely to a random variable $X_{\infty}.$ Under
this additional hypothesis, the sequence is equi-integrable. Moreover,
for $p\geqslant n,$
\[
X_{n}\geqslant\mathbb{E}^{\mathscr{A}_{n}}\left(X_{p}\right).
\]
Passing to the limit in $\text{L}^{1}$ as $p$ tends to infinity
yields 
\[
X_{n}\geqslant\mathbb{E}^{\mathscr{A}_{n}}\left(X_{\infty}\right).
\]

\textbf{General case} 

For every $q\in\mathbb{Q}^{+},$ the process 
\[
X\land q=\left(X_{n}\land q\right)_{n\in\mathbb{N}}
\]
is a supermartingale being the minimum of two supermartingales. Moreover,
it is nonnegative and bounded by $q.$ From the previous step, the
sequence $\left(X_{n}\land q\right)_{n\in\mathbb{N}}$ converges $P-$almost
surely. By a classical argument, it follows that, $P-$almost surely,
for every $q\in\mathbb{Q}^{+},$ the sequence $\left(X_{n}\land q\right)_{n\in\mathbb{N}}$
converges. By deterministic Lemma $\ref{lm:determ_lemma},$ we conclude
that the sequence $\left(X_{n}\right)_{n\in\mathbb{N}}$ converges
in $\overline{\mathbb{R}}^{+}$ $P-$almost surely. 

Finally, since for $p\geqslant n,$ 
\[
X_{n}\geqslant\mathbb{E}^{\mathscr{A}_{n}}\left(X_{p}\right),
\]
 the conditional Fatou lemma implies 
\[
X_{n}\geqslant\mathbb{E}^{\mathscr{A}_{n}}\left(X_{\infty}\right).
\]

\end{proof}

\begin{lemma}{}{determ_lemma}

Let $\left(x_{n}\right)_{n\in\mathbb{N}}$ be a sequence of nonnegative
real numbers such that, for every $q\in\mathbb{Q}^{+},$ the sequence
$\left(x_{n}\land q\right)_{n\in\mathbb{N}}$ converges. Then the
sequence $\left(x_{n}\right)_{n\in\mathbb{N}}$ converges in $\overline{\mathbb{R}}^{+}.$

\end{lemma}

\begin{proof}{}{}

If 
\[
\sup_{n\in\mathbb{N}}x_{n}<+\infty,
\]
choose $q>\sup_{n\in\mathbb{N}}x_{n}.$ Then for every $n\in\mathbb{N},$
\[
x_{n}\land q=x_{n},
\]
and by hypothesis the sequence $\left(x_{n}\right)_{n\in\mathbb{N}}$
converges in $\mathbb{R}^{+}.$

If instead 
\[
\sup_{n\in\mathbb{N}}x_{n}=+\infty,
\]
then for every $q\in\mathbb{Q}^{+},$ denote by $l_{q}$ the limit
of the sequence $\left(x_{n}\land q\right)_{n\in\mathbb{N}}.$ We
have 
\[
0\leqslant l_{q}\leqslant q.
\]
For every $\epsilon>0,$ there exists $N\left(\epsilon,q\right)$
such that 
\[
l_{q}-\epsilon\leqslant x_{n}\land q\leqslant l_{q}+\epsilon
\]
for $n\geqslant N\left(\epsilon,q\right).$ 

Let $B>0$ be arbitrary, and set $\epsilon=\dfrac{B}{4}.$ Choose
$q>B.$

If we had 
\[
l_{q}\leqslant B-\epsilon,
\]
then for every $n\geqslant N\left(\dfrac{B}{4},q\right),$ 
\[
x_{n}\land q\leqslant B,
\]
and hence $x_{n}\leqslant B.$ This would imply that 
\[
\sup_{n\in\mathbb{N}}x_{n}<+\infty,
\]
and there would be a contradiction. Therefore, 
\[
l_{q}>B-\epsilon.
\]
For $n\geqslant N\left(\dfrac{B}{4},q\right),$ 
\[
x_{n}\geqslant x_{n}\land q\geqslant l_{q}-\epsilon>B-2\epsilon=\dfrac{B}{2}.
\]
Since $B$ is arbitrary, we conclude that
\[
\lim_{n\to+\infty}x_{n}=+\infty.
\]
In all cases, the sequence $\left(x_{n}\right)_{n\in\mathbb{N}}$
converges in $\overline{\mathbb{R}}^{+}.$

\end{proof}

\section*{Exercises}

\addcontentsline{toc}{section}{Exercises}

\begin{exercise}{$\sigma-$algebra of Previous Events to a Stopping Time}{exercise16.1}

Let $\left(\mathscr{A}_{n}\right)_{n\in\mathbb{N}}$ be a filtration
on the probabilized space \preds and define 
\[
\mathscr{A}_{\infty}=\bigvee_{n\in\mathbb{N}}\mathscr{A}_{n},
\]
the \salg generated by the union of the \salgs $\mathscr{A}_{n},\,n\in\mathbb{N}.$ 

Let $S$ and $T$ be two stopping times. Prove that, for every $Y\in\text{L}^{1}\left(\Omega,\mathscr{A},P\right),$
\[
\mathbb{E}^{\mathscr{A}_{S}}\left(\mathbb{E}^{\mathscr{A}_{T}}\left(Y\right)\right)=\mathbb{E}^{\mathscr{A}_{T}}\left(\mathbb{E}^{\mathscr{A}_{S}}\left(Y\right)\right)=\mathbb{E}^{\mathscr{A}_{S\land T}}\left(Y\right).
\]

\end{exercise}

\begin{exercise}{The Gambler's Ruin Problem}{exercise16.2}

A gambler repeatedly tosses---a not necessarily fair---coin. Denote
by $p$ the probability of obtaining heads at a toss. The gambler
receives one euro from the bank if a head appears and gives one euro
to the bank if a tail appears. The gambler's initial fortune is $a\in\mathbb{N}^{\ast}$
euros, and the bank's initial fortune is $b\in\mathbb{N}^{\ast}$
euros. The gambler plays until either they are ruined or the bank
is ruined. 

We model this game as follows. Let $\left(Y_{n}\right)_{n\in\mathbb{N}^{\ast}}$
be a sequence of independent random variables defined on a probabilized
space $\left(\Omega,\mathscr{A},P\right),$ all following the same
law $p\delta_{1}+q\delta_{-1},$ where $q=1-p.$ Let $S_{n}$ denote
the gambler's fortune after $n$ rounds, for a game that never ends.
Define
\[
S_{0}=a\,\,\,\,\text{and}\,\,\,\,S_{n}=a+\sum^{n}_{j=1}Y_{j}.
\]
By setting $Y_{0}=a,$ the natural filtrations $\left(\mathscr{A}_{n}\right)_{n\in\mathbb{N}}$
of the processes $Y$ and $S$ coincide. Let $T$ be the stopping
time of the game, that is
\[
T=\inf\left\{ n\in\mathbb{N}^{\ast}:\,S_{n}=0\,\,\,\,\text{or}\,\,\,\,a+b\right\} .
\]

We address the following questions: 
\begin{itemize}
\item What is the probability $P\left(T<+\infty\right)$ that the game eventually
stops?
\item What is the probability $\rho=P\left(S_{T}=a+b\right)$ that the gambler
wins?
\item What is the expected duration $\mathbb{E}\left(T\right)$ of the game
to stop?
\end{itemize}
1. Determine the nature of the process $S=\left(S_{n}\right)_{n\in\mathbb{N}}$
depending on the values of $p.$

\textbf{2. Study of the case $p\neq q.$ }

Assume that $p>q.$ 

(a) Write the Doob decomposition of the submartingale $S$ and specify
its predictable nondecreasing process $A.$ 

(b) Deduce that $\mathbb{E}\left(T\right)<+\infty.$ 

(c) Determine the value of $P\left(T<+\infty\right)$ and express
$\mathbb{E}\left(T\right)$ as a function of $\rho.$ 

For $s>0,$ define the process $U,$ for every $n\in\mathbb{N},$
by $U_{n}=s^{S_{n}}.$ 

(d) Determine $s$ so that $U$ is a nonconstant martingale. 

(e) Verify that the stopped martingale $U^{T}$ converges $P-$almost
surely and in $\text{L}^{1}$ to $U_{T}.$ 

(f) Deduce from this the value of $\rho,$ and then compute $\mathbb{E}\left(T\right).$

\textbf{3. Study of the case $p=\dfrac{1}{2}.$}

(a) Show that $S$ is a square-integrable martingale and determine
its predictable nondecreasing process $B.$

(b) Deduce from (a) that $\mathbb{E}\left(T\right)<+\infty.$ Determine
the value of $P\left(T<+\infty\right).$

(c) Verify that the stopped martingale $S^{T}$ converges $P-$almost
surely in $\text{L}^{1}$ and $\text{L}^{2}$ to $S_{T}.$ 

(d) Deduce the values of $\mathbb{E}\left(S_{T}\right),$ $\rho$
and $\mathbb{E}\left(T\right).$

\end{exercise}

\begin{exercise}{Game of Heads and Tails With Two Unfair Coins and a Learning Strategy}{exercise16.3}

A gambler has two coins $A$ and $B.$ The probabilities, unknown
of the gambler, to obtain head when tossing the coin $A$---respectively
$B$---is $p^{A}$---respectively $p^{B}.$ The gambler wins one
unit each time a head appears. After each toss, the gambler chooses
which coin to use for the next toss based on the outcomes of the previous
tosses. The gambler strategy is to identify which coin has the largest
probability $p$ of heads in order to maximize their gain. We model
this game as follows.

Let $\left\{ \left(X^{A}_{n}\right)_{n\in\mathbb{N}},\left(X^{B}_{n}\right)_{n\in\mathbb{N}}\right\} $
be a family of independent random variables defined on a probabilized
space \preds. For every $n\in\mathbb{N},$ $X^{A}_{n}$---respectively
$X^{B}_{n}$---follow the Bernoulli laws with parameter $p^{A}$
and $p^{B},$ respectively. Let $\mathscr{A}_{n}$ denote the \salg
$\sigma\left(X^{A}_{j},X^{B}_{j}:\,0\leqslant j\leqslant n\right)$
and let $U=\left(U_{n}\right)_{n\in\mathbb{N}}$ be a process adapted
to the filtration $\left(\mathscr{A}_{n}\right)_{n\in\mathbb{N}}$
and taking values in the set $\left\{ A,B\right\} .$ Define the process
of toss outcomes $\left(X_{n}\right)_{n\in\mathbb{N}}$ by, for every
$n\in\mathbb{N},$
\[
X_{n+1}=\boldsymbol{1}_{\left(U_{n}=A\right)}X^{A}_{n+1}+\boldsymbol{1}_{\left(U_{n}=B\right)}X^{B}_{n+1}=X^{U_{n}}_{n+1}.
\]
Set
\[
G_{0}=0\,\,\,\,\text{and, if }n\geqslant1,\,\,\,\,G_{n}=\sum^{n}_{j=1}X_{j}.
\]

1. Compute the conditional expectation $\mathbb{E}^{\mathscr{A}_{n}}\left(X_{n+1}\right).$

2. Define the process $M$ by $M_{0}=0$ and, for $n\in\mathbb{N}^{\ast},$
by
\[
M_{n}=\sum^{n}_{j=1}\left(X_{j}-p^{U_{j-1}}\right).
\]

(a) Verify that $M$ is a square-integrable martingale and compute
its predictable nondecreasing process. 

(b) Deduce that the sequence with general term 
\[
\dfrac{G_{n}}{n}-\dfrac{1}{n}\sum^{n}_{j=1}p^{U_{j-1}}
\]
converges $P-$almost surely to 0. 

3. For $J\in\left\{ A,B\right\} ,$ define the processes $N^{J},$
$M^{J},$ and $\widetilde{p}^{J}$ by
\[
N^{J}_{n}=\sum^{n}_{j=0}\boldsymbol{1}_{\left(U_{j}=J\right)},\,\,\,\,M^{J}_{0}=0,\,\,\,\,\text{and}\,\,\,\,M^{J}_{n}=\sum^{n}_{j=0}\left[\boldsymbol{1}_{\left(U_{j-1}=J,X_{j}=1\right)}-\boldsymbol{1}_{\left(U_{j-1}=J\right)}p^{J}\right],
\]
and for $n\in\mathbb{N}^{\ast},$
\[
\widetilde{p}^{J}_{n}=\dfrac{1}{N^{J}_{n-1}}\sum^{n-1}_{j=0}\boldsymbol{1}_{\left(U_{j}=J,X_{j+1}=1\right)}
\]
with the convention $\dfrac{0}{0}=0.$

Prove that $M^{J}$ is a square-integrable martingale and compute
its predictable nondecreasing process $\left\langle M^{J}\right\rangle .$
Deduce that, on the set $\left\{ N^{J}_{n}\longrightarrow+\infty\right\} ,$
the sequence with general term $\widetilde{p}^{J}_{n}$ converges
$P-$almost surely to $p^{J}.$

4. Let $\left(v_{n}\right)_{n\in\mathbb{N}}$ be an increasing sequence
of nonnegative integers such that 
\[
\lim_{n\to+\infty}\dfrac{v_{n}}{n}=+\infty.
\]
 The gambler adopts the following process of choice $U=\left(U_{n}\right)_{n\in\mathbb{N}}:$
\begin{gather*}
\text{If }n\notin\left\{ v_{j}:j\in\mathbb{N}\right\} ,\,\,\,\,U_{n}=A\boldsymbol{1}_{\left(\widetilde{p}^{A}_{n}\geqslant\widetilde{p}^{B}_{n}\right)}+B\boldsymbol{1}_{\left(\widetilde{p}^{A}_{n}<\widetilde{p}^{B}_{n}\right)},\\
\text{and}\,\,\,\,U_{v_{2n}}=A,\,\,\,\,U_{V_{2n+1}}=B.
\end{gather*}
That is, they choose the coin that has appeared the most frequently
when $n$ is not in the support of the sequence $v,$ and choose alternatively
the coins $A$ and $B$ along the sequence $v.$ For this strategy,
for each $J\in\left\{ A,B\right\} ,\,N^{J}_{n}\to+\infty.$ 

Assume, for instance, $p^{A}>p^{B}$---which is unknown to the gambler.

(a) Study, for $J\in\left\{ A,B\right\} ,$ the $P-$almost sure convergence
of the sequences of general term $\dfrac{N^{J}_{n}}{n}.$

(b) Then study, for this choice, the $P-$almost sure convergence
of the sequence with general term $\dfrac{G_{n}}{n}.$

\end{exercise}

\begin{exercise}{Game of Heads and Tails with a Gain Depending on Two Consecutive Tosses}{exercise16.4}

Let $\left(Y_{n}\right)_{n\in\mathbb{N}^{\ast}}$ be a sequence of
independent random variables defined on a probabilized space \preds,
all following the same law $p\delta_{1}+q\delta_{-1},$ where $q=1-p.$
Let $S_{n}$ denote the fortune of a gambler after $n$ tosses. Suppose
the gain rule is given by
\[
S_{0}=a\,\,\,\,\text{and}\,\,\,\,S_{n}=a+\sum^{n}_{j=1}Y_{j-1}Y_{j}.
\]
All processes are considered with respect to the natural filtration
$\left(\mathscr{A}_{n}\right)_{n\in\mathbb{N}}$ of the process $Y.$

1. Compute the probability $P\left(S_{n}>S_{n-1}\right)$ and verify
that it is strictly greater than $\dfrac{1}{2}$ whenever $p\neq q.$

2. For $n\in\mathbb{N}^{\ast},$ compute the conditional expectation
$\mathbb{E}^{\mathscr{A}_{n-1}}\left(S_{n}\right).$ What is the nature
of the process $S$ when $p=\dfrac{1}{2}?$\\
Study the convergence of the sequence with general term $\mathbb{E}\left(S_{n}\right).$

3. Let $s>0$ be arbitrary. 

(a) For $n\in\mathbb{N}^{\ast},$ compute the conditional expectation
$\mathbb{E}^{\mathscr{A}_{n-1}}\left(s^{S_{n}}\right).$ 

Set $u=s+\dfrac{1}{s}.$

(b) Prove that the process $\left(\dfrac{s^{S_{n}}}{u^{n}}\right)_{n\in\mathbb{N}}$
is a nonnegative supermartingale.

(c) Study the $P-$almost sure and $\text{L}^{1}$ convergences of
the sequence $\left(\dfrac{s^{S_{n}}}{u^{n}}\right)_{n\in\mathbb{N}}.$

4. (a) Prove that $S$ can be written in a unique manner as the sum
of a square-integrable martingale $W$ and a predictable integrable
process $T$ satisfying $T_{0}=0.$ 

(b) Compute the predictable nondecreasing process $\left\langle W\right\rangle $
associated with the martingale $W.$ 

(c) Study the $P-$almost sure convergence of the sequence $\left(\dfrac{S_{n}}{n}\right)_{n\in\mathbb{N}^{\ast}}$
and deduce from it, in the case where $p\neq q,$ the $P-$almost
sure convergence of the sequence $\left(S_{n}\right)_{n\in\mathbb{N}}.$

\end{exercise}

\begin{exercise}{A Model of Stock Portfolio}{exercise16.5}

The process $S=\left(S_{n}\right)_{n\in\mathbb{N}}$ models the temporal
evolution of a stock price with the help of a probabilized space \preds
on which the sequence of random variables $\left(S_{n}\right)_{n\in\mathbb{N}}$
is defined by
\[
S_{0}=s_{0}>0\,\,\,\,\text{and,}\,\text{if}\,n\in\mathbb{N}^{\ast},\,\,\,\,S_{n}=\left(1+\mu\right)S_{n-1}+\sigma S_{n-1}\epsilon_{n},
\]
where $\left(\epsilon_{n}\right)_{n\in\mathbb{N}^{\ast}}$ is a \textbf{noise}
process consisting in a sequence of independent random variables following
the same law $\dfrac{1}{2}\left(\delta_{0}+\delta_{-1}\right),$ and
where the real parameters $\mu$---called discount rate---and $\sigma$---called
volatility coefficient---satisfy the inequality
\[
\left|\sigma\right|<1+\mu.
\]
Denote $\lambda$ the real number $\lambda=\left[\left(1+\mu\right)^{2}-\sigma^{2}\right]^{\frac{1}{2}}.$
Let $\left(\mathscr{A}_{n}\right)_{n\in\mathbb{N}}$ be the natural
filtration of the process $S.$

1. Let $f$ be the real-valued function defined on $\mathbb{R}^{2}$
by $f\left(x,y\right)=\left(1+\mu\right)x+\sigma xy.$ Prove that,
for every $x\in\mathbb{R}^{+},$
\[
f\left(x,1\right)\geqslant0\,\,\,\,\text{and}\,\,\,\,f\left(x,-1\right)\geqslant0.
\]

2. For every $n\in\mathbb{N}^{\ast},$ compute $\mathbb{E}^{\mathscr{A}_{n-1}}\left(\Delta S_{n}\right).$
Deduce the nature of the process $S,$ depending on the values of
the parameters $\mu$ and $\sigma.$ In the case $\mu<0,$ prove that
the sequence with general term $S_{n}$ converges $P-$almost surely
to a limit that must be determined.

3. Verify that $S_{n}$ is square-integrable and compute $\mathbb{E}\left(S^{2}_{n}\right).$

4. Define the process $Z$ for every $n\in\mathbb{N},$ by $Z_{n}=\ln S_{n}.$
Prove that, depending on the values of $\lambda,$ the process $Z$
is a martingale, a submartingale, or a supermartingale.

Write $Z_{n}$ as a sum of independent random variables and deduce,
according to the values of $\lambda,$ the $P-$almost sure convergence---in
$\overline{\mathbb{R}}$---of the sequence $\left(Z_{n}\right)_{n\in\mathbb{N}}$
to a limit to be specified. How do these results translate for the
sequence $\left(S_{n}\right)_{n\in\mathbb{N}}?$

5. In the special case where $\left|\sigma\right|<\dfrac{1}{\sqrt{2}}$
and $\left(1+\mu\right)^{2}+\sigma^{2}<1,$ show that the previous
results imply that $-S$ is a submartingale. 

Write its Doob decomposition $-S=M+A,$ where $M$ is an integrable
martingale and $A$ a predictable nondecreasing process, nonzero at
zero.

Verify that $M$ is a martingale in $\text{L}^{2}$ and compute its
predictable nondecreasing process $\left\langle M\right\rangle .$ 

Deduce from this the $P-$almost sure convergence of the series with
general term $S^{2}_{n}.$

6. Define the process $W,$ for every $n\in\mathbb{N},$ by $W_{n}=\ln\left(\dfrac{S_{n}}{\lambda^{n}}\right).$
Prove that $W$ is a martingale in $\text{L}^{2}$ and compute its
predictable nondecreasing process $\left\langle W\right\rangle $
as a function of $\delta=-\left[\ln\left(\dfrac{1+\mu+\sigma}{\lambda}\right)\right]\left[\ln\left(\dfrac{1+\mu-\sigma}{\lambda}\right)\right].$
Verify that this number is positive.

Deduce from this the $P-$almost sure convergence of the sequence
$\left(S^{\frac{1}{n}}_{n}\right)_{n\in\mathbb{N}^{\ast}}$ to a limit
to be specified.

7. Define the process $R,$ for every $n\in\mathbb{N},$ by $R_{n}=\lambda^{-\sqrt{n}}S^{\frac{1}{\sqrt{n}}}_{n}.$
Prove that the sequence of laws $P_{R_{n}}$ converges narrowly to
a law admitting a density with respect to Lebesgue measure. Determine
this law.

\end{exercise}

\section*{Solutions of Exercises}

\addcontentsline{toc}{section}{Solutions of Exercises}

\begin{solution}{}{solexercise16.1}

For every $n\in\mathbb{N},$ set $X_{n}=\mathbb{E}^{\mathscr{A}_{n}}\left(Y\right).$
Then the process $\left(X_{n}\right)_{n\in\mathbb{N}}$ is an equi-integrable
martingale, hence it converges $P-$almost surely and in $\text{L}^{1}.$
Setting 
\[
X_{\infty}=\limsup_{n\to+\infty}X_{n},
\]
we have $P-$almost surely $Y=X_{\infty}$ and the process $\left(X_{n}\right)_{n\in\overline{\mathbb{N}}}$
is a closed martingale. By the Doob second stopping theorem, we obtain
\[
X_{T}=\mathbb{E}^{\mathscr{A}_{T}}\left(Y\right),
\]
and for every $n\in\mathbb{N},$
\[
X^{T}_{n}=X_{T\land n}=\mathbb{E}^{\mathscr{A}_{T\land n}}\left(X_{T}\right),
\]
where $X^{T}$ denotes the martingale $X$ stopped at $T.$ 

Fix $n\in\mathbb{N}$ and evaluate $\mathbb{E}^{\mathscr{A}_{T\land n}}\left(X_{T}\right).$
For every $j\in\overline{\mathbb{N}},$ on $\left(T=j\right),$ we
have
\[
\mathbb{E}^{\mathscr{A}_{T\land n}}\left(X_{T}\right)=\mathbb{E}^{\mathscr{A}_{j\land n}}\left(X_{T}\right)=\begin{cases}
\mathbb{E}^{\mathscr{A}_{j}}\left(X_{T}\right), & \text{if }j<n,\\
\mathbb{E}^{\mathscr{A}_{n}}\left(X_{T}\right), & \text{if }j\geqslant n.
\end{cases}
\]
Now, for $j<n,$ since $\left(T=j\right)\in\mathscr{A}_{j},$ since
$\mathscr{A}_{j}\subset\mathscr{A}_{n},$ and since $X_{j}$ is $\mathscr{A}_{j}-$measurable,
and so $\mathscr{A}_{n}-$measurable, then
\[
\boldsymbol{1}_{\left(T=j\right)}\mathbb{E}^{\mathscr{A}_{j}}\left(X_{T}\right)=\mathbb{E}^{\mathscr{A}_{j}}\left(\boldsymbol{1}_{\left(T=j\right)}X_{T}\right)=\boldsymbol{1}_{\left(T=j\right)}X_{j}=\mathbb{E}^{\mathscr{A}_{n}}\left(\boldsymbol{1}_{\left(T=j\right)}X_{T}\right)=\boldsymbol{1}_{\left(T=j\right)}\mathbb{E}^{\mathscr{A}_{n}}\left(X_{T}\right).
\]
It follows that, for every $j\in\overline{\mathbb{N}},$ in all cases
on $\left(T=j\right),$
\[
\mathbb{E}^{\mathscr{A}_{T\land n}}\left(X_{T}\right)=\mathbb{E}^{\mathscr{A}_{n}}\left(X_{T}\right).
\]
Hence, for every $n\in\mathbb{N},$
\begin{equation}
X^{T}_{n}=X_{T\land n}=\mathbb{E}^{\mathscr{A}_{n}}\left(X_{T}\right).\label{eq:X_nT_as_cond_exp_XT_A_n}
\end{equation}
Since $X_{T}=\mathbb{E}^{\mathscr{A}_{T}}\left(Y\right),$ the random
variable $X_{T}$ is integrable and the stopped martingale $X^{T}$
is closable. Apply the Doob second stopping theorem to $X^{T}$ with
the stopping time $S.$ Using $\refpar{eq:X_nT_as_cond_exp_XT_A_n},$
we obtain
\[
X^{T}_{S}=X_{T\land S}=\mathbb{E}^{\mathscr{A}_{S}}\left(X_{T}\right)=\mathbb{E}^{\mathscr{A}_{S}}\left(\mathbb{E}^{\mathscr{A}_{T}}\left(Y\right)\right),
\]
and in particular, 
\[
X_{T\land S}=\mathbb{E}^{\mathscr{A}_{S}}\left(\mathbb{E}^{\mathscr{A}_{T}}\left(Y\right)\right).
\]
 By interchanging the roles of $S$ and $T,$ we also have 
\[
X_{T\land S}=\mathbb{E}^{\mathscr{A}_{T}}\left(\mathbb{E}^{\mathscr{A}_{S}}\left(Y\right)\right),
\]
which implies the equality 
\[
\mathbb{E}^{\mathscr{A}_{S}}\left(\mathbb{E}^{\mathscr{A}_{T}}\left(Y\right)\right)=\mathbb{E}^{\mathscr{A}_{T}}\left(\mathbb{E}^{\mathscr{A}_{S}}\left(Y\right)\right).
\]
Finally, applying the Doob second stopping theorem to the closed martingale
$\left(X_{n}\right)_{n\in\overline{\mathbb{N}}},$ with the stopping
theorem $S\land T,$ we obtain $X_{T\land S}=\mathbb{E}^{\mathscr{A}_{S\land T}}\left(Y\right),$
which proves the last equality.

\end{solution}

\begin{solution}{}{solexercise16.2}

\textbf{1. Nature of $S=\left(S_{n}\right)_{n\in\mathbb{N}}$}

For $n\in\mathbb{N}^{\ast},$ 
\[
\mathbb{E}^{\mathscr{A}_{n-1}}\left(\Delta S_{n}\right)=\mathbb{E}^{\mathscr{A}_{n-1}}\left(Y_{n}\right).
\]
Since the random variables $Y_{n}$ are independent, it follows
\[
\mathbb{E}^{\mathscr{A}_{n-1}}\left(\Delta S_{n}\right)=\mathbb{E}\left(Y_{n}\right)=p-q,
\]
which gives the following classification for the process $S:$
\begin{itemize}
\item If $p>q,$ then  $S$ is a submartingale. 
\item If $p=q=\dfrac{1}{2},$ then $S$ is a martingale.
\item If $p<q,$ then $S$ is a supermartingale.
\end{itemize}
\textbf{2. Study of the case $p>q.$ }

\textbf{(a) Doob decomposition of $S.$ Predictable nondecreasing
process $A$ }

The submartingale $S$ admits the Doob decomposition $S=M+A,$ where
the predictable nondecreasing process $A$ is defined by $A_{0}=0$
and $\Delta A_{n}=\mathbb{E}^{\mathscr{A}_{n-1}}\left(\Delta S_{n}\right),$
which yields\boxeq{
\[
A_{0}=0,\,\,\,\,\text{and, for}\,n\geqslant1,\,A_{n}=n\left(p-q\right).
\]
}

\textbf{(b) Proof that $\mathbb{E}\left(T\right)<+\infty.$ }

Applying the Doob first stopping theorem to the martingale 
\[
M=\left(S_{n}-n\left(p-q\right)\right)_{n\in\mathbb{N}}
\]
and to the bounded stopping time $T\land n$ then gives
\[
a=\mathbb{E}\left(S_{0}\right)=\mathbb{E}\left(S_{T\land n}-T\land n\left(p-q\right)\right).
\]
Hence
\begin{equation}
\left(p-q\right)\mathbb{E}\left(T\land n\right)=\mathbb{E}\left(S_{T\land n}\right)-a.\label{eq:expTlandn}
\end{equation}
By definition of $T,$ for every $n\in\mathbb{N},$ 
\[
0\leqslant S_{T\land n}\leqslant a+b,
\]
so
\[
0\leqslant\left(p-q\right)\mathbb{E}\left(T\land n\right)\leqslant b.
\]
By the Beppo Levi property,\boxeq{
\[
\mathbb{E}\left(T\right)=\lim_{n\to+\infty}\nearrow\mathbb{E}\left(T\land n\right).
\]
}

\textbf{(c) Value of $P\left(T<+\infty\right).$ Expression of $\mathbb{E}\left(T\right)$
as a function of $\rho$ }

It follows that $T$ is integrable, hence in particular 
\[
P\left(T<+\infty\right)=1.
\]

The sequence $\left(S_{T\land n}\right)_{n\in\mathbb{N}}$ converges
$P-$almost surely to $S_{T}.$ Passing to the limit in $\refpar{eq:expTlandn},$
and using dominated convergence gives
\begin{equation}
\left(p-q\right)\mathbb{E}\left(T\right)=\mathbb{E}\left(S_{T}\right)-a.\label{eq:p-qET}
\end{equation}
Since, by definition of $T,$
\[
\mathbb{E}\left(S_{T}\right)=\left(a+b\right)P\left(S_{T}=a+b\right),
\]
we obtain\boxeq{
\begin{equation}
\mathbb{E}\left(T\right)=\dfrac{\left(a+b\right)\rho-a}{p-q}.\label{eq:ET_exp_rho}
\end{equation}
}

\textbf{(d) Value of $s$ such that $U$ is a nonconstant martingale. }

Let $s>0.$ Since $s^{S_{n-1}}$ is $\mathscr{A}_{n-1}-$measurable
and independent of $s^{Y_{n}},$ for every $n\in\mathbb{N}^{\ast},$
\[
\mathbb{E}^{\mathscr{A}_{n-1}}\left(U_{n}\right)=\mathbb{E}^{\mathscr{A}_{n-1}}\left(s^{S_{n-1}}s^{Y_{n}}\right)=s^{S_{n-1}}\mathbb{E}^{\mathscr{A}_{n-1}}\left(s^{Y_{n}}\right)=s^{S_{n-1}}\mathbb{E}\left(s^{Y_{n}}\right).
\]
Thus
\[
\mathbb{E}^{\mathscr{A}_{n-1}}\left(U_{n}\right)=U_{n-1}\left(sp+\dfrac{1}{s}q\right).
\]
We have $sp+\dfrac{1}{s}q=1$ if and only if $s^{2}p-s+q=0,$ equation
which has the trivial root $s=1$---since $p+q=1$---and the other
root is $\dfrac{q}{p}.$ Hence, \textbf{for $s=\dfrac{q}{p},$ $U$
is a nonconstant martingale}.

\textbf{(e) $U^{T}$ converges $P-$almost surely and in $\text{L}^{1}$
to $U_{T}.$ }

By definition of $T$ and since $\dfrac{q}{p}<1,$ for every $n\in\mathbb{N},$
\[
0\leqslant U_{T\land n}\leqslant1.
\]
Thus the stopped martingale $U^{T}$ is equi-integrable and converges
$P-$almost surely and in $\text{L}^{1}$ to $U_{T}.$

\textbf{(f) Values of $\rho$ and $\mathbb{E}\left(T\right).$}

By definition of $T,$ $U_{T}$ takes $P-$almost surely the values
$1$ or $\left(\dfrac{q}{p}\right)^{a+b}.$ Its expectation is thus
equal to
\[
\mathbb{E}\left(U_{T}\right)=P\left(S_{T}=0\right)+\left(\dfrac{q}{p}\right)^{a+b}P\left(S_{T}=a+b\right).
\]
Hence,
\begin{equation}
\mathbb{E}\left(U_{T}\right)=1-\rho+\left(\dfrac{q}{p}\right)^{a+b}\rho.\label{eq:EU_T_f_rho}
\end{equation}
Moreover, by the first stopping theorem, for every $n\in\mathbb{N},$
\[
\mathbb{E}\left(U_{T\land n}\right)=\mathbb{E}\left(U_{0}\right)=\left(\dfrac{q}{p}\right)^{a},
\]
hence by the dominated convergence theorem
\[
\mathbb{E}\left(U_{T}\right)=\lim_{n\to+\infty}\mathbb{E}\left(U_{T\land n}\right)=\left(\dfrac{q}{p}\right)^{a}.
\]
Combining with $\refpar{eq:EU_T_f_rho}$ and $\refpar{eq:ET_exp_rho},$
we obtain\boxeq{
\[
\rho=\dfrac{1-\left(\dfrac{q}{p}\right)^{a}}{1-\left(\dfrac{q}{p}\right)^{a+b}}\,\,\,\,\text{and}\,\,\,\,\mathbb{E}\left(T\right)=\dfrac{1}{p-q}\left[\dfrac{1-\left(\dfrac{q}{p}\right)^{a}}{1-\left(\dfrac{q}{p}\right)^{a+b}}-a\right].
\]
}

\textbf{3. Study of the case $p=q=\dfrac{1}{2}$}

\textbf{(a) $S$ is a square-integrable martingale. Predictable nondecreasing
process $B$}

In this case, the process $S$ is a martingale in $\text{L}^{2},$
since the $Y_{n}$ are bounded. Its predictable nondecreasing process
$B$ is defined by $B_{0}=0$ and its increments given, for $n\geqslant1,$
by $\Delta B_{n}=\mathbb{E}^{\mathscr{A}_{n-1}}\left(\left(\Delta S_{n}\right)^{2}\right).$

Hence, by independence of $\mathscr{A}_{n-1}$ and $Y^{2}_{n},$
\[
\Delta B_{n}=\mathbb{E}^{\mathscr{A}_{n-1}}\left(Y^{2}_{n}\right)=\mathbb{E}\left(Y^{2}_{n}\right)=1.
\]
Thus
\[
B_{0}=0\,\,\,\,\text{and, for }n\geqslant1,\,B_{n}=n.
\]

\textbf{(b) Proof that $\mathbb{E}\left(T\right)<+\infty.$ Value
of $P\left(T<+\infty\right)$}

The Doob first stopping theorem applied to the martingale $\left(S^{2}_{n}-n\right)_{n\in\mathbb{N}}$
and the bounded stopping time $T\land n$ then yields, since $S_{0}=a,$
\begin{equation}
\mathbb{E}\left(S^{2}_{0}\right)=\mathbb{E}\left(S^{2}_{T\land n}-T\land n\right)=a^{2}.\label{eq:EsqS_0}
\end{equation}
Since $S^{2}_{T\land n}\leqslant\left(a+b\right)^{2},$
\[
\mathbb{E}\left(S^{2}_{T\land n}\right)=\mathbb{E}\left(T\land n\right)\leqslant\left(a+b\right)^{2}-a^{2},
\]
By the Beppo Levi property,
\[
\mathbb{E}\left(T\right)=\lim_{n\to+\infty}\mathbb{E}\left(T\land n\right)\leqslant\left(a+b\right)^{2}-a^{2}.
\]

Hence $T$ is integrable, and in particular $P\left(T<+\infty\right)=1.$

\textbf{(c) $S^{T}$ converges $P-$almost surely in $\text{L}^{1}$
and $\text{L}^{2}$ to $S_{T}$ }

The sequence $\left(S_{T\land n}\right)_{n\in\mathbb{N}}$ converges
$P-$almost surely to $S_{T}.$ Since for every $n\in\mathbb{N}^{\ast},$
\[
0\leqslant S_{T\land n}\leqslant a+b,
\]
the dominated convergence theorem implies that there is also $\text{L}^{1}$
and $\text{L}^{2}$ convergence to $S_{T}.$ 

\textbf{(d) Values of $\mathbb{E}\left(S_{T}\right),$ $\rho$ and
$\mathbb{E}\left(T\right)$}

The first stopping theorem applied to the martingale $S$ and the
bounded stopping time $T\land n$ yields
\[
\mathbb{E}\left(S_{T\land n}\right)=\mathbb{E}\left(S_{0}\right)=a,
\]
By $\text{L}^{1}$ convergence\boxeq{
\[
\mathbb{E}\left(S_{T}\right)=a.
\]
}

Since, by definition of $T,$
\[
\mathbb{E}\left(S_{T}\right)=\left(a+b\right)P\left(S_{T}=a+b\right)=\rho\left(a+b\right),
\]
we obtain\boxeq{
\[
\rho=\dfrac{a}{a+b}.
\]
}

By the relation $\refpar{eq:EsqS_0},$ for every $n\in\mathbb{N}^{\ast},$
\begin{equation}
\mathbb{E}\left(T\land n\right)=\mathbb{E}\left(S^{2}_{T\land n}\right)-a^{2}.\label{eq:expTlandn-1}
\end{equation}
Since the sequence $\left(S_{T\land n}\right)_{n\in\mathbb{N}}$ converges
to $S_{T}$ in $\text{L}^{2},$ passing to the limit in $\refpar{eq:expTlandn-1}$
yields
\[
\mathbb{E}\left(T\right)=\mathbb{E}\left(S^{2}_{T}\right)-a^{2}.
\]
By definition of $T,$
\[
\mathbb{E}\left(S^{2}_{T}\right)=\left(a+b\right)^{2}P\left(S_{T}=a+b\right)=\rho\left(a+b\right)^{2}=a\left(a+b\right).
\]
Therefore,\boxeq{
\[
\mathbb{E}\left(T\right)=ab.
\]
}

\end{solution}

\begin{solution}{}{solexercise16.3}

\textbf{1. Computation of $\mathbb{E}^{\mathscr{A}_{n}}\left(X_{n+1}\right)$}

Since the process $U$ is adapted,
\[
\mathbb{E}^{\mathscr{A}_{n}}\left(X_{n+1}\right)=\boldsymbol{1}_{\left(U_{n}=A\right)}\mathbb{E}^{\mathscr{A}_{n}}\left(X^{A}_{n+1}\right)+\boldsymbol{1}_{\left(U_{n}=B\right)}\mathbb{E}^{\mathscr{A}_{n}}\left(X^{B}_{n+1}\right).
\]
By independence of $X^{J}_{n+1}$---for $J\in\left\{ A,B\right\} $---and
$\mathscr{A}_{n},$
\[
\mathbb{E}^{\mathscr{A}_{n}}\left(X_{n+1}\right)=\boldsymbol{1}_{\left(U_{n}=A\right)}\mathbb{E}\left(X^{A}_{n+1}\right)+\boldsymbol{1}_{\left(U_{n}=B\right)}\mathbb{E}\left(X^{B}_{n+1}\right)=\boldsymbol{1}_{\left(U_{n}=A\right)}p^{A}+\boldsymbol{1}_{\left(U_{n}=B\right)}p^{B},
\]
which yields\boxeq{
\[
\mathbb{E}^{\mathscr{A}_{n}}\left(X_{n+1}\right)=p^{U_{n}}.
\]
}

\textbf{2. (a) M is a martingale of integrable square. Computation
of its predictable nondecreasing process}

We have
\[
\mathbb{E}^{\mathscr{A}_{n-1}}\left(\Delta M_{n}\right)=\mathbb{E}^{\mathscr{A}_{n-1}}\left(X_{n}-p^{U_{n-1}}\right)=0,
\]
so $M$ is a martingale. It is square-integrable since $0\leqslant X_{n}\leqslant1$
and $0\leqslant p^{U_{n-1}}\leqslant1.$ 

Its predictable nondecreasing process $\left\langle M\right\rangle $
is defined by $\left\langle M\right\rangle _{0}$ and its increments
are given, for $n\geqslant1,$ by
\begin{align*}
\Delta\left\langle M\right\rangle _{n} & =\mathbb{E}^{\mathscr{A}_{n-1}}\left(\left[X_{n}-p^{U_{n-1}}\right]^{2}\right)\\
 & =\boldsymbol{1}_{\left(U_{n-1}=A\right)}\mathbb{E}^{\mathscr{A}_{n-1}}\left(\left[X^{A}_{n}-p^{A}\right]^{2}\right)+\boldsymbol{1}_{\left(U_{n-1}=B\right)}\mathbb{E}^{\mathscr{A}_{n-1}}\left(\left[X^{B}_{n}-p^{B}\right]^{2}\right).
\end{align*}
Hence, since $U_{n-1}$ is $\mathscr{A}_{n-1}-$measurable and that
$\mathscr{A}_{n-1}$ and $X^{J}_{n}$ are independent,
\begin{align*}
\Delta\left\langle M\right\rangle _{n} & =\boldsymbol{1}_{\left(U_{n-1}=A\right)}\mathbb{E}\left(\left[X^{A}_{n}-p^{A}\right]^{2}\right)+\boldsymbol{1}_{\left(U_{n-1}=B\right)}\mathbb{E}\left(\left[X^{B}_{n}-p^{B}\right]^{2}\right)\\
 & =\boldsymbol{1}_{\left(U_{n-1}=A\right)}p^{A}\left(1-p^{A}\right)+\boldsymbol{1}_{\left(U_{n-1}=B\right)}p^{B}\left(1-p^{B}\right).
\end{align*}
Thus\boxeq{
\[
\left\langle M\right\rangle _{n}=\sum^{n}_{j=1}p^{U_{j-1}}\left(1-p^{U_{j-1}}\right).
\]
}

\textbf{(b) Convergence $P-$almost sure to 0 of the sequence with
general term $\dfrac{G_{n}}{n}-\dfrac{1}{n}\sum^{n}_{j=1}p^{U_{j-1}}.$}

Let $m=\min\left(p^{A}\left(1-p^{A}\right),p^{B}\left(1-p^{B}\right)\right)$
and $s=\max\left(p^{A}\left(1-p^{A}\right),p^{B}\left(1-p^{B}\right)\right).$

Then for every $n\in\mathbb{N}^{\ast},$
\[
0\leqslant nm\leqslant\left\langle M\right\rangle _{n}\leqslant ns.
\]

Hence, the sequence with general term $\left\langle M\right\rangle _{n}$
tends to $+\infty$ with $n.$ By the strong law of large numbers
for $\text{L}^{2}-$martingales,
\[
\dfrac{M_{n}}{\left\langle M\right\rangle _{n}}\stackrel[n\to+\infty]{P-\text{a.s.}}{\longrightarrow}0,
\]
and therefore, using the bounds above,
\[
\dfrac{M_{n}}{n}\stackrel[n\to+\infty]{P-\text{a.s.}}{\longrightarrow}0.
\]
This proves that\boxeq{
\begin{equation}
\dfrac{G_{n}}{n}-\dfrac{1}{n}\sum^{n}_{j=1}p^{U_{j-1}}\stackrel[n\to+\infty]{P-\text{a.s.}}{\longrightarrow}0.\label{eq:limG_novern-avp_U_j-1}
\end{equation}
}

\textbf{3. Martingale $M^{J}.$ Computation of $\left\langle M^{J}\right\rangle .$
$\widetilde{p}^{J}_{n}$ converges $P-$almost surely to $p^{J}$}

Fix $J\in\left\{ A,B\right\} .$ Since the process $U$ is adapted,
\[
\mathbb{E}^{\mathscr{A}_{n-1}}\left(\Delta M^{J}_{n}\right)=\boldsymbol{1}_{\left(U_{n-1}=J\right)}\mathbb{E}^{\mathscr{A}_{n-1}}\left(\boldsymbol{1}_{\left(X^{J}_{n}=1\right)}-p^{J}\right).
\]
By independence,
\[
\mathbb{E}^{\mathscr{A}_{n-1}}\left(\Delta M^{J}_{n}\right)=\boldsymbol{1}_{\left(U_{n-1}=J\right)}\mathbb{E}\left(\boldsymbol{1}_{\left(X^{J}_{n}=1\right)}-p^{J}\right)=0.
\]
The process $M^{J}$ is a martingale, and it is straightforwardly
square-integrable.

Its predictable nondecreasing process $\left\langle M^{J}\right\rangle $
is defined by $\left\langle M^{J}\right\rangle _{0}=0$ and its increments
are given, for $n\geqslant1,$ by
\[
\Delta\left\langle M^{J}\right\rangle _{n}=\mathbb{E}^{\mathscr{A}_{n-1}}\left(\Delta M^{J}_{n}\right)^{2}=\mathbb{E}^{\mathscr{A}_{n-1}}\left(\boldsymbol{1}_{\left(U_{n-1}=J\right)}\left(\boldsymbol{1}_{\left(X^{J}_{n}=1\right)}-p^{J}\right)^{2}\right).
\]
Since $U_{n-1}$ is $\mathscr{A}_{n-1}-$measurable and since $\mathscr{A}_{n-1}$
and $X^{J}_{n}$ are independent
\[
\Delta\left\langle M^{J}\right\rangle _{n}=\boldsymbol{1}_{\left(U_{n-1}=J\right)}\mathbb{E}\left(\left(\boldsymbol{1}_{\left(X^{J}_{n}=1\right)}-p^{J}\right)^{2}\right)=\boldsymbol{1}_{\left(U_{n-1}=J\right)}p^{J}\left(1-p^{J}\right).
\]
Hence, for every $n\in\mathbb{N}^{\ast},$\boxeq{
\[
\left\langle M\right\rangle _{n}=N^{J}_{n-1}p^{J}\left(1-p^{J}\right).
\]
}By the strong law of large numbers for $\text{L}^{2}$-martingales,
on the set $\left\{ N^{J}_{n}\longrightarrow+\infty\right\} ,$
\[
\dfrac{M^{J}_{n}}{N^{J}_{n}}\stackrel[n\to+\infty]{P-\text{a.s.}}{\longrightarrow}0.
\]
Moreover, 
\[
\widetilde{p^{J}_{n}}=\dfrac{M^{J}_{n}}{N^{J}_{n-1}}+p^{J},
\]
so we obtain\boxeq{
\[
\text{on the set }\left\{ N^{J}_{n}\longrightarrow+\infty\right\} ,\,\,\,\,\widetilde{p^{J}_{n}}\stackrel[n\to+\infty]{P-\text{a.s.}}{\longrightarrow}p^{J}.
\]
}

\textbf{4. }

\textbf{(a) $P-$almost sure convergence of $\dfrac{N^{J}_{n}}{n}$}

For this strategy, for each $J\in\left\{ A,B\right\} ,$ $N^{J}_{n}\underset{n\to+\infty}{\longrightarrow}+\infty.$ 

Assume for instance that $p^{A}>p^{B}.$ 

Let $\omega$ be such that the sequence with general term $\widetilde{p^{J}_{n}}\left(\omega\right)$
converges to $p^{J}$ for $J=A$ and $B.$ 

Then there exists an integer $N\left(\omega\right)$ such that for
every $n\geqslant N\left(\omega\right),$ 
\[
\widetilde{p^{A}_{n}}\left(\omega\right)\geqslant\widetilde{p^{B}_{n}}\left(\omega\right)
\]
Hence, 
\[
\lim_{n\to+\infty}\dfrac{1}{n}\sum^{n}_{j=0}\boldsymbol{1}_{\left(\widetilde{p^{A}_{n}}\left(\omega\right)\geqslant\widetilde{p^{B}_{n}}\left(\omega\right)\right)}\left(\omega\right)=1.
\]
By definition of the process of choice $U,$ for every $n>N\left(\omega\right),$
\[
N^{A}_{n}\left(\omega\right)-N^{A}_{N\left(\omega\right)}\left(\omega\right)=\sum^{n}_{j=N\left(\omega\right)+1}\boldsymbol{1}_{\left(\widetilde{p^{A}_{n}}\left(\omega\right)\geqslant\widetilde{p^{B}_{n}}\left(\omega\right)\right)}\left(\omega\right)=\sum^{+\infty}_{j=1}\boldsymbol{1}_{\left(N\left(\omega\right)<v_{2j}\leqslant n\right)}.
\]
Since $\dfrac{v_{n}}{n}\underset{n\to+\infty}{\longrightarrow}+\infty,$
there exists an integer $N^{\prime}$ such that $v_{n}>n$ for $n\geqslant N^{\prime}.$
Therefore,
\[
L_{n}\equiv\sum^{n}_{j=1}\boldsymbol{1}_{\left(v_{2j}\leqslant n\right)}\leqslant\#\left\{ j:\,v_{2j}\leqslant N^{\prime}\right\} ,
\]
which implies
\[
\dfrac{L_{n}}{n}\underset{n\to+\infty}{\longrightarrow}0.
\]

Hence
\[
\dfrac{N^{A}_{n}\left(\omega\right)}{n}\underset{n\to+\infty}{\longrightarrow}1.
\]

By a similar argument, 
\[
\dfrac{N^{B}_{n}\left(\omega\right)}{n}\underset{n\to+\infty}{\longrightarrow}0.
\]
Thus,\boxeq{
\begin{equation}
\dfrac{N^{A}_{n}\left(\omega\right)}{n}\overset{P-\text{a.s.}}{\underset{n\to+\infty}{\longrightarrow}}1\,\,\,\,\text{and}\,\,\,\,\dfrac{N^{B}_{n}\left(\omega\right)}{n}\overset{P-\text{a.s.}}{\underset{n\to+\infty}{\longrightarrow}}0.\label{eq:P_as_cvNAnomandNB}
\end{equation}
}

\textbf{(b) $P-$almost sure convergence of $\dfrac{G_{n}}{n}$}

We have
\[
\sum^{n}_{j=1}p^{U_{j-1}}=p^{A}\sum^{n-1}_{j=0}\boldsymbol{1}_{\left(U_{j}=A\right)}+p^{B}\sum^{n-1}_{j=0}\boldsymbol{1}_{\left(U_{j}=B\right)}=p^{A}N^{A}_{n-1}+p^{B}N^{B}_{n-1}.
\]
Using $\refpar{eq:limG_novern-avp_U_j-1}$ and $\refpar{eq:P_as_cvNAnomandNB},$
we obtain\boxeq{
\[
\dfrac{G_{n}}{n}\stackrel[n\to+\infty]{P-\text{a.s.}}{\longrightarrow}p^{A}.
\]
}

\end{solution}

\begin{solution}{}{solexercise16.4}

\textbf{1. Computation of $P\left(S_{n}>S_{n-1}\right).$ Proof of
$P\left(S_{n}>S_{n-1}\right)>\dfrac{1}{2}$ when $p\neq q$}

Since 
\[
\left(S_{n}>S_{n-1}\right)=\left(Y_{n-1}Y_{n}=1\right)
\]
and since the $Y_{n}$ take $P-$almost surely the values $\pm1,$
we have
\[
P\left(S_{n}>S_{n-1}\right)=P\left(\left(Y_{n-1}=1\right)\cap\left(Y_{n}=1\right)\right)+P\left(\left(Y_{n-1}=-1\right)\cap\left(Y_{n}=-1\right)\right).
\]
By independence of $Y_{n-1}$ and $Y_{n},$\boxeq{
\[
P\left(S_{n}>S_{n-1}\right)=p^{2}+q^{2}.
\]
}

Using $p+q=1,$ we have
\[
p^{2}+q^{2}=2p^{2}-2p+1\equiv h\left(p\right).
\]
Then, $h^{\prime}\left(p\right)=2\left(2p-1\right)$ and $h^{\prime\prime}\left(p\right)=4,$
so $h$ admits a minimum at $p=\dfrac{1}{2},$ and $h\left(\dfrac{1}{2}\right)=\dfrac{1}{2}.$
Therefore, if $p\neq q,$\boxeq{
\[
P\left(S_{n}>S_{n-1}\right)>\dfrac{1}{2}.
\]
}

\textbf{2. Computation of $\mathbb{E}^{\mathscr{A}_{n-1}}\left(S_{n}\right).$
Nature of $S$ when $p=\dfrac{1}{2}.$ Convergence of $\left(\mathbb{E}\left(S_{n}\right)\right)_{n\in\mathbb{N}}$}

For $n\in\mathbb{N}^{\ast},$
\[
S_{n}=S_{n-1}+Y_{n-1}Y_{n}.
\]
Since the random variables $S_{n-1}$ and $Y_{n-1}$ are $\mathscr{A}_{n-1}-$measurable,
\[
\mathbb{E}^{\mathscr{A}_{n-1}}\left(S_{n}\right)=S_{n-1}+Y_{n-1}\mathbb{E}^{\mathscr{A}_{n-1}}\left(Y_{n}\right).
\]
Since the \salg $\mathscr{A}_{n-1}$ and $\sigma\left(Y_{n}\right)$
are independent, 
\[
\mathbb{E}^{\mathscr{A}_{n-1}}\left(Y_{n}\right)=\mathbb{E}\left(Y_{n}\right),
\]
hence
\begin{equation}
\mathbb{E}^{\mathscr{A}_{n-1}}\left(S_{n}\right)=S_{n-1}+\left(p-q\right)Y_{n-1}.\label{eq:cond_exp_A_n-1S_n}
\end{equation}
In particular, \textbf{if $p=q,$ }then\textbf{ $S$ is a martingale}.

It follows by $\refpar{eq:cond_exp_A_n-1S_n}$ that
\[
\mathbb{E}\left(\Delta S_{n}\right)=\mathbb{E}\left(\mathbb{E}^{\mathscr{A}_{n-1}}\left(\Delta S_{n}\right)\right)=\left(p-q\right)\mathbb{E}\left(Y_{n-1}\right)=\left(p-q\right)^{2}.
\]
Therefore,
\[
\mathbb{E}\left(S_{n}\right)=a+n\left(p-q\right)^{2}.
\]
Hence, if $p\neq q,$ 
\[
\lim_{n\to+\infty}\mathbb{E}\left(S_{n}\right)=+\infty,
\]
whereas if $p=q,$ the sequence is constant.

\textbf{3. (a) Computation of $\mathbb{E}^{\mathscr{A}_{n-1}}\left(s^{S_{n}}\right)$ }

For $n\in\mathbb{N}^{\ast},$ $s^{S_{n}}=s^{S_{n-1}}s^{Y_{n-1}Y_{n}}.$
Since the random variable $s^{S_{n-1}}$ is $\mathscr{A}_{n-1}-$measurable,
\[
\mathbb{E}^{\mathscr{A}_{n-1}}\left(s^{S_{n}}\right)=s^{S_{n-1}}\mathbb{E}^{\mathscr{A}_{n-1}}\left(s^{Y_{n-1}Y_{n}}\right).
\]
Now, for every $\left(y_{0},y_{1},\cdots,y_{n-1}\right)\in\left\{ -1,1\right\} ^{n},$
we have the relations on the conditional means
\begin{multline*}
m^{\left(Y_{0},Y_{1},\cdots,Y_{n-1}\right)=\left(y_{0},y_{1},\cdots,y_{n-1}\right)}\left(s^{Y_{n-1}Y_{n}}\right)=m^{\left(Y_{0},Y_{1},\cdots,Y_{n-1}\right)=\left(y_{0},y_{1},\cdots,y_{n-1}\right)}\left(s^{y_{n-1}Y_{n}}\right),
\end{multline*}
which yields, by independence of $\left(Y_{0},Y_{1},\cdots,Y_{n-2}\right)$
and $Y_{n-1}$ 
\[
m^{\left(Y_{0},Y_{1},\cdots,Y_{n-1}\right)=\left(y_{0},y_{1},\cdots,y_{n-1}\right)}\left(s^{Y_{n-1}Y_{n}}\right)=\mathbb{E}\left(s^{y_{n-1}Y_{n}}\right)=ps^{y_{n-1}}+qs^{-y_{n-1}}.
\]
Since the \salg $\mathscr{A}_{n-1}$ is generated by $\left(Y_{0},Y_{1},\cdots,Y_{n-1}\right),$
we obtain
\[
\mathbb{E}^{\mathscr{A}_{n-1}}\left(s^{Y_{n-1}Y_{n}}\right)=ps^{Y_{n-1}}+qs^{-Y_{n-1}}.
\]
Hence,\boxeq{
\[
\mathbb{E}^{\mathscr{A}_{n-1}}\left(s^{S_{n}}\right)=s^{S_{n-1}}\left[ps^{Y_{n-1}}+qs^{-Y_{n-1}}\right].
\]
}

\textbf{(b) $\left(\dfrac{s^{S_{n}}}{u^{n}}\right)_{n\in\mathbb{N}}$
is a nonnegative supermartingale}

We can write
\[
\mathbb{E}^{\mathscr{A}_{n-1}}\left(s^{S_{n}}\right)=s^{S_{n-1}}\left[\boldsymbol{1}_{\left(Y_{n-1}=1\right)}\left(ps+\dfrac{q}{s}\right)+\boldsymbol{1}_{\left(Y_{n-1}=-1\right)}\left(\dfrac{p}{s}+qs\right)\right],
\]
and therefore
\[
\mathbb{E}^{\mathscr{A}_{n-1}}s^{S_{n}}\leqslant s^{S_{n-1}}\left(s+\dfrac{1}{s}\right).
\]
Set $u=s+\dfrac{1}{s}.$ Dividing by $u^{n},$
\[
\mathbb{E}^{\mathscr{A}_{n-1}}\left(\dfrac{s^{S_{n}}}{u^{n}}\right)\leqslant\dfrac{s^{S_{n-1}}}{u^{n-1}}.
\]
Thus the process $\left(\dfrac{s^{S_{n}}}{u^{n}}\right)_{n\in\mathbb{N}}$
is a nonnegative supermartingale.

\textbf{(c) $P-$almost sure and $\text{L}^{1}$ convergence of $\left(\dfrac{s^{S_{n}}}{u^{n}}\right)_{n\in\mathbb{N}}$}

Hence, we already know that it converges $P-$almost surely. In fact,
we are going to see that this sequence is bounded. Indeed, for every
$n\in\mathbb{N}^{\ast},$ $P-$almost surely
\[
a-n\leqslant S_{n}\leqslant a+n.
\]

\begin{itemize}
\item If $0<s<1,$ then $0\leqslant s^{S_{n}}\leqslant s^{a-n},$ so
\[
0\leqslant\dfrac{s^{S_{n}}}{u^{n}}\leqslant\dfrac{s^{a}}{\left(us\right)^{n}}\leqslant s^{a},\,\,\,\,\text{since}\,\,\,\,us=1+s^{2}>1.
\]
\item If $s\geqslant1,$ similarly $0\leqslant s^{S_{n}}\leqslant s^{a+n},$
so
\[
0\leqslant\dfrac{s^{S_{n}}}{u^{n}}\leqslant s^{a}\left(\dfrac{s}{u}\right)^{n}\leqslant s^{a},\,\,\,\,\text{since}\,\,\,\,\dfrac{u}{s}=1+\dfrac{1}{s^{2}}>1.
\]
\end{itemize}
Hence, for $s>0,$
\[
\lim_{n\to+\infty}\dfrac{s^{S_{n}}}{u^{n}}=0.
\]
Moreover, for every $s>0,$ and every $n\in\mathbb{N}^{\ast},$
\[
0\leqslant\dfrac{s^{S_{n}}}{u^{n}}\leqslant s^{a},
\]
which implies that the sequence is equi-integrable and thus also converges
to 0 in $\text{L}^{1}.$ 

\textbf{4. (a) $S=W+T$}

Let $W$ be the process defined by $W_{0}=S_{0}=a$ and the increments
given, for $n\geqslant1,$ by
\[
\Delta W_{n}=\Delta S_{n}-\mathbb{E}^{\mathscr{A}_{n-1}}\left(\Delta S_{n}\right)=Y_{n-1}Y_{n}-\left(p-q\right)Y_{n-1}.
\]
Equivalently,
\[
W_{0}=a,\,\,\,\,\Delta W_{n}=Y_{n-1}\left[Y_{n}-\left(p-q\right)\right].
\]
By construction, $W$ is a martingale. It is in $\text{L}^{2},$ since
the increments $\Delta W_{n}$ are $P-$almost surely bounded. Let
$T$ be the process defined by $T_{0}=0$ and its increments given,
for $n\geqslant1,$ by
\[
\Delta T_{n}=\mathbb{E}^{\mathscr{A}_{n-1}}\left(\Delta S_{n}\right)=\left(p-q\right)Y_{n-1}.
\]
Then
\[
S=W+T
\]
and $T$ is a predictable and integrable process.

Let us show its uniqueness. Suppose that $S=W^{\prime}+T^{\prime}$
is another decomposition of the same type. Then, for every $n\in\mathbb{N},$
\[
\Delta\left(W-W^{\prime}\right)_{n}=\Delta\left(T-T^{\prime}\right)_{n}.
\]
Taking conditional expectations and taking into account the processes
properties, we obtain
\[
0=\mathbb{E}^{\mathscr{A}_{n-1}}\left(\Delta\left(W-W^{\prime}\right)_{n}\right)=\Delta\left(T-T^{\prime}\right)_{n},
\]
which proves that $T=T^{\prime},$ and consequently $W=W^{\prime}.$
Therefore the decomposition is unique.

\textbf{(b) Computation of $\left\langle W\right\rangle $ }

The predictable nondecreasing process of $W,$ $\left\langle W\right\rangle ,$
is determined by its increments. Since $Y^{2}_{n-1}=1$ $P-$almost
surely, they are given, for every $n\in\mathbb{N},$ by
\[
\Delta\left\langle W\right\rangle _{n}=\mathbb{E}^{\mathscr{A}_{n-1}}\left(\left[\Delta W_{n}\right]^{2}\right)=\mathbb{E}^{\mathscr{A}_{n-1}}\left(\left[Y_{n}-\left(p-q\right)\right]^{2}\right).
\]
By independence of $\mathscr{A}_{n-1}$ and $\sigma\left(Y_{n}\right),$
\[
\Delta\left\langle W\right\rangle _{n}=\mathbb{E}\left(\left[Y_{n}-\left(p-q\right)\right]^{2}\right)=\sigma^{2}_{Y_{n}}=\mathbb{E}\left(Y^{2}_{n}\right)-\left(\mathbb{E}\left(Y_{n}\right)\right)^{2}=1-\left(p-q\right)^{2}=4pq.
\]
Thus,\boxeq{
\[
\left\langle W\right\rangle _{n}=4pqn.
\]
}

\textbf{(c) $P-$almost sure sonvergence of $\left(\dfrac{S_{n}}{n}\right)_{n\in\mathbb{N}^{\ast}},$
and $\left(S_{n}\right)_{n\in\mathbb{N}}$ when $p\neq q$}

The strong law of large numbers for $\text{L}^{2}$-martingales yields
\[
\dfrac{W_{n}}{\left\langle W_{n}\right\rangle }\stackrel[n\to+\infty]{P-\text{a.s.}}{\longrightarrow}0.
\]
Moreover, 
\[
\dfrac{S_{n}}{n}=\dfrac{W_{n}}{n}+\dfrac{T_{n}}{n}=\dfrac{W_{n}}{n}+\left(p-q\right)\dfrac{1}{n}\sum^{n}_{j=1}Y_{j-1}.
\]
By the strong law of large numbers for the sequence of independent
random variables, with same law, yields
\[
\dfrac{1}{n}\sum^{n}_{j=1}Y_{j-1}\stackrel[n\to+\infty]{P-\text{a.s.}}{\longrightarrow}\mathbb{E}\left(Y_{0}\right).
\]
 Since $\mathbb{E}\left(Y_{0}\right)=p-q,$ it follows\boxeq{
\[
\dfrac{S_{n}}{n}\stackrel[n\to+\infty]{P-\text{a.s.}}{\longrightarrow}\left(p-q\right)^{2}.
\]
}If $p\neq q,$ this limit is positive, and thus it follows that
the sequence with general term $S_{n}$ tends to $+\infty$ $P-$almost
surely as $n$ tends to infinity.

\end{solution}

\begin{solution}{}{solexercise16.5}

\textbf{1. Proof that $f\left(x,1\right)\geqslant0$ and $f\left(x,-1\right)\geqslant0$
for $x\in\mathbb{R}^{+}$}

We have
\[
f\left(x,1\right)=\left(1+\mu+\sigma\right)x
\]
and
\[
f\left(x,-1\right)=\left(1+\mu-\sigma\right)x.
\]

The inequalities $-\sigma\leqslant\left|\sigma\right|<1+\mu$ and
$\sigma\leqslant\left|\sigma\right|<1+\mu$ imply that, for every
$x\in\mathbb{R}^{+},$ $f\left(x,1\right)\geqslant0$ and $f\left(x,-1\right)\geqslant0.$ 

It follows that $S$ is a process taking values in $\mathbb{R}^{+}$
$P-$almost surely.

\textbf{2. Computation of $\mathbb{E}^{\mathscr{A}_{n-1}}\left(\Delta S_{n}\right).$
Nature of $S,$ depending on $\mu$ and $\sigma.$ When $\mu<0,$
convergence $P-$almost sure of $S_{n}$ and limit }

For every $n\in\mathbb{N}^{\ast},$ 
\[
\Delta S_{n}=S_{n-1}\left(\mu+\sigma\epsilon_{n}\right)
\]
Since $S_{n-1}$ is $\mathscr{A}_{n-1}-$measurable and since the
\salgs $\mathscr{A}_{n-1}$ and $\sigma\left(\epsilon_{n}\right)$
are independent,
\[
\mathbb{E}^{\mathscr{A}_{n-1}}\left(\Delta S_{n}\right)=S_{n-1}\mathbb{E}^{\mathscr{A}_{n-1}}\left(\mu+\sigma\epsilon_{n}\right)=S_{n-1}\mathbb{E}\left(\mu+\sigma\epsilon_{n}\right).
\]

Hence,\boxeq{
\begin{equation}
\mathbb{E}^{\mathscr{A}_{n-1}}\left(\Delta S_{n}\right)=\mu S_{n-1}.\label{eq:E_A_n-1_Delta_S_nismuS_n-1}
\end{equation}
}

Since $S_{n-1}\geqslant0$ $P-$almost surely, it follows that\boxeq{
\[
\begin{cases}
S\text{ is a submartingale,} & \text{if }\mu>0.\\
S\text{ is a martingale,} & \text{if }\mu=0.\\
S\text{ is a supermartingale,} & \text{if }\mu<0.
\end{cases}
\]
}

In the case where $\mu<0,$ $S$ is a nonnegative supermartingale.
Therefore the sequence with general term $S_{n}$ converges $P-$almost
surely in $\overline{\mathbb{R}}^{+}.$ Moreover, by $\refpar{eq:E_A_n-1_Delta_S_nismuS_n-1},$
\[
\mathbb{E}\left(S_{n}\right)=\left(1+\mu\right)\mathbb{E}\left(S_{n-1}\right)=\left(1+\mu\right)^{n}s_{0}.
\]
Since in this case, $0<1+\mu<1,$
\[
\lim_{n\to+\infty}\mathbb{E}\left(S_{n}\right)=0.
\]
Since the $S_{n}$ are nonnegative, the sequence with general term
$S_{n}$ converges in $\text{L}^{1}.$ Thus, we also have convergence
$P-$almost sure to 0 when $n$ tends to infinity.

\textbf{3. Verify that $S_{n}$ is square-integrable and computation
of $\mathbb{E}\left(S^{2}_{n}\right)$}

Since the random variables $S_{n-1}$ and $\epsilon_{n}$ are independent,
\[
\mathbb{E}\left(S^{2}_{n}\right)=\mathbb{E}\left(S^{2}_{n-1}\right)\mathbb{E}\left(\left[1+\mu+\sigma\epsilon_{n}\right]^{2}\right).
\]
Now
\[
\mathbb{E}\left(\left[1+\mu+\sigma\epsilon_{n}\right]^{2}\right)=\dfrac{1}{2}\left[\left(1+\mu+\sigma\right)^{2}+\left(1+\mu-\sigma\right)^{2}\right]=\left(1+\mu\right)^{2}+\sigma^{2}.
\]
Hence, for every $n\in\mathbb{N}^{\ast},$\boxeq{
\begin{equation}
\mathbb{E}\left(S^{2}_{n}\right)=s^{2}_{0}\left[\left(1+\mu\right)^{2}+\sigma^{2}\right]^{n}.\label{eq:EsqS_nfs_0muandsig}
\end{equation}
}

\textbf{4. Nature of $Z$ depending on $\lambda.$ Representation
of $Z_{n}.$ $P-$almost sure convergence of $\left(Z_{n}\right)_{n\in\mathbb{N}}$
and consequences for $\left(S_{n}\right)_{n\in\mathbb{N}}$}

We have $Z_{0}=\ln S_{0}=\ln s_{0}$ and, for every $n\in\mathbb{N}^{\ast},$
\[
Z_{n}=Z_{n-1}+\ln\left(1+\mu+\sigma\epsilon_{n}\right).
\]
By independence of the \salgs $\mathscr{A}_{n-1}$ and $\sigma\left(\epsilon_{n}\right),$
\[
\mathbb{E}^{\mathscr{A}_{n-1}}\left(\Delta S_{n}\right)=\mathbb{E}\left(\ln\left(1+\mu+\sigma\epsilon_{n}\right)\right)=\dfrac{1}{2}\left[\ln\left(1+\mu+\sigma\right)+\ln\left(1+\mu-\sigma\right)\right].
\]
Thus,
\[
\mathbb{E}^{\mathscr{A}_{n-1}}\left(\Delta S_{n}\right)=\dfrac{1}{2}\ln\left(\left[1+\mu\right]^{2}-\sigma^{2}\right),
\]
that is,\boxeq{
\begin{equation}
\mathbb{E}^{\mathscr{A}_{n-1}}\left(\Delta S_{n}\right)=\ln\lambda.\label{eq:E_A_n-1_DeltaS_n_is_lnlamda}
\end{equation}
}

Hence, depending on whether $\lambda$ is greater than, equal to,
or less than 1, $Z$ is respectively a submartingale, a martingale
or a supermartingale.

Moreover, for every $n\in\mathbb{N}^{\ast},$
\begin{equation}
Z_{n}=Z_{0}+\sum^{n}_{j=1}\ln\left(1+\mu+\sigma\epsilon_{j}\right).\label{eq:Z_n_expr}
\end{equation}
Since the random variables $\ln\left(1+\mu+\sigma\epsilon_{n}\right)$
are independent, following the same law, and integrable, the strong
law of large numbers yields
\[
\dfrac{1}{n}\sum^{n}_{j=1}\ln\left(1+\mu+\sigma\epsilon_{j}\right)\stackrel[n\to+\infty]{P-\text{a.s.}}{\longrightarrow}\mathbb{E}\left(\ln\left(1+\mu+\sigma\epsilon_{1}\right)\right),
\]
so 
\[
\dfrac{Z_{n}}{n}\stackrel[n\to+\infty]{P-\text{a.s.}}{\longrightarrow}\ln\lambda.
\]
Thus\boxeq{
\[
\left\{ \begin{array}{cccc}
\text{If }\lambda>1, & Z_{n}\stackrel[n\to+\infty]{P-\text{a.s.}}{\longrightarrow}+\infty, & \text{and} & S_{n}\stackrel[n\to+\infty]{P-\text{a.s.}}{\longrightarrow}+\infty.\\
\text{If }\lambda=1, & Z_{n}\stackrel[n\to+\infty]{P-\text{a.s.}}{\longrightarrow}0, & \text{and} & S^{\frac{1}{n}}_{n}\stackrel[n\to+\infty]{P-\text{a.s.}}{\longrightarrow}1.\\
\text{If }\lambda<1, & Z_{n}\stackrel[n\to+\infty]{P-\text{a.s.}}{\longrightarrow}-\infty, & \text{and} & S_{n}\stackrel[n\to+\infty]{P-\text{a.s.}}{\longrightarrow}0.
\end{array}\right.
\]
}

\textbf{5. $M$ is a $\text{L}^{2}$-martingale. Computation of $\left\langle M\right\rangle .$
$P-$almost sure convergence of $S^{2}_{n}$}

Since $\left|\sigma\right|<1+\mu,$ the condition $\left(1+\mu\right)^{2}+\sigma^{2}<1$
is realized as soon as $\left|\sigma\right|<\dfrac{1}{\sqrt{2}}.$
In this case, necessarily $\mu<0,$ and thus $-S$ is an integrable
submartingale.

Let $-S=M+A$ be its Doob decomposition, where $M$ is an integrable
martingale and $A$ is a predictable, nondecreasing process, null
at zero. By $\refpar{eq:E_A_n-1_Delta_S_nismuS_n-1},$
\[
\Delta A_{n}=\mathbb{E}^{\mathscr{A}_{n-1}}\left(\Delta\left(-S\right)_{n}\right)=-\mu S_{n-1}.
\]
Hence\boxeq{
\[
\begin{cases}
A_{0}=0,\,\text{and\,\ensuremath{\text{if}}\,}n\ensuremath{\in}\mathbb{N}^{\ast}, & A_{n}=-\mu\sum^{n-1}_{j=0}S_{j}.\\
M_{0}=s_{0},\,\text{and\,\ensuremath{\text{if}}\,}n\ensuremath{\in}\mathbb{N}^{\ast}, & M_{n}=-S_{n}+\mu\sum^{n-1}_{j=0}S_{j}.
\end{cases}
\]
}Therefore, the process $M$ is a martingale in $\text{L}^{2},$
and its predictable nondecreasing process $\left\langle M\right\rangle $
is given by its increments
\begin{align*}
\Delta\left\langle M\right\rangle _{n} & =\mathbb{E}^{\mathscr{A}_{n-1}}\left(\left(\Delta M_{n}\right)^{2}\right)=\mathbb{E}^{\mathscr{A}_{n-1}}\left(\left[\Delta S_{n}+\Delta A_{n}\right]^{2}\right)\\
 & =\mathbb{E}^{\mathscr{A}_{n-1}}\left(\left(S_{n-1}\left(\mu+\sigma\epsilon_{n}\right)-\mu S_{n-1}\right)^{2}\right)=\sigma^{2}S^{2}_{n-1}.
\end{align*}
Therefore,\boxeq{
\[
\left\langle M\right\rangle _{0}=0\,\,\,\,\text{and,\,if\,}n\in\mathbb{N}^{\ast},\,\,\,\,\left\langle M\right\rangle _{n}=\sigma^{2}\sum^{n-1}_{j=0}S^{2}_{j}.
\]
}

Then
\[
\mathbb{E}\left(\left\langle M\right\rangle _{n}\right)=\sigma^{2}\sum^{n-1}_{j=0}\mathbb{E}\left(S^{2}_{j}\right),
\]
and by $\refpar{eq:E_A_n-1_DeltaS_n_is_lnlamda},$
\[
\mathbb{E}\left(\left\langle M\right\rangle _{n}\right)=\sigma^{2}s^{2}_{0}\dfrac{1-\left[\left(1+\mu\right)^{2}+\sigma^{2}\right]^{n}}{1-s^{2}_{0}\left[\left(1+\mu\right)^{2}+\sigma^{2}\right]}.
\]
Since, by hypothesis $\left(1+\mu\right)^{2}+\sigma^{2}<1,$ the sequence
with general term $\mathbb{E}\left(\left\langle M\right\rangle _{n}\right)$
converges, hence
\[
\sum^{+\infty}_{j=0}\mathbb{E}\left(S^{2}_{j}\right)=\mathbb{E}\left(\sum^{+\infty}_{j=0}S^{2}_{j}\right)<+\infty.
\]
It follows that the series with general term $S^{2}_{n}$ converges
$P-$almost surely.

\textbf{6. $W$ is an $\text{L}^{2}$-martingale. Computation of $\left\langle W\right\rangle ,$
and $P-$almost sure convergence of $\left(S^{\frac{1}{n}}_{n}\right)_{n\in\mathbb{N}^{\ast}}$}

We have
\[
\Delta W_{n}=\ln\left(\dfrac{S_{n}}{\lambda S_{n-1}}\right)=\Delta Z_{n}-\ln\lambda.
\]
Using $\refpar{eq:E_A_n-1_DeltaS_n_is_lnlamda},$ we get 
\[
\mathbb{E}^{\mathscr{A}_{n-1}}\left(\Delta W_{n}\right)=0,
\]
so $W$ is a martingale. It is in $\text{L}^{2}$ since the random
variables $Z_{n}$ are in $\text{L}^{2}.$ Its predictable nondecreasing
process $\left\langle W\right\rangle $ satisfies, using the independence
of \salgs $\mathscr{A}_{n-1}$ and $\sigma\left(\epsilon_{n}\right),$
\[
\Delta\left\langle W\right\rangle _{n}=\mathbb{E}^{\mathscr{A}_{n-1}}\left[\left(\Delta W_{n}\right)^{2}\right]=\mathbb{E}^{\mathscr{A}_{n-1}}\left[\left(\ln\dfrac{1+\mu+\sigma\epsilon_{n}}{\lambda}\right)^{2}\right]=\mathbb{E}\left[\left(\ln\dfrac{1+\mu+\sigma\epsilon_{n}}{\lambda}\right)^{2}\right].
\]
Hence,
\begin{align*}
\Delta\left\langle W\right\rangle _{n} & =\dfrac{1}{2}\left[\left(\ln\dfrac{1+\mu+\sigma}{\lambda}\right)^{2}+\left(\ln\dfrac{1+\mu-\sigma}{\lambda}\right)^{2}\right]\\
 & =\dfrac{1}{2}\left[\left(\ln\dfrac{1+\mu+\sigma}{\lambda}+\ln\dfrac{1+\mu-\sigma}{\lambda}\right)^{2}-2\ln\dfrac{1+\mu+\sigma}{\lambda}\ln\dfrac{1+\mu-\sigma}{\lambda}\right].
\end{align*}
Since
\[
\ln\dfrac{1+\mu+\sigma}{\lambda}+\ln\dfrac{1+\mu-\sigma}{\lambda}=\ln\dfrac{\left(1+\mu\right)^{2}-\sigma^{2}}{\lambda^{2}}=0,
\]
it follows that $\Delta\left\langle W\right\rangle _{n}=\delta$---which
proves that $\delta>0$---and thus that, for every $n\in\mathbb{N},$
$\left\langle W\right\rangle _{n}=n\delta.$

The strong law of large numbers for the martingales in $\text{L}^{2}$
then ensures that $\dfrac{W_{n}}{\left\langle W\right\rangle _{n}}\stackrel[n\to+\infty]{P-\text{a.s.}}{\longrightarrow}0,$
hence also
\[
\dfrac{\ln S_{n}-n\ln\lambda}{n}\stackrel[n\to+\infty]{P-\text{a.s.}}{\longrightarrow}0,
\]
which implies that\boxeq{
\[
S^{\frac{1}{n}}_{n}\stackrel[n\to+\infty]{P-\text{a.s.}}{\longrightarrow}\lambda.
\]
}

\textbf{7. Narrow convergence of $P_{R_{n}}.$ Determination of the
limit}

Since 
\[
\ln R_{n}=\dfrac{1}{\sqrt{n}}Z_{n}-\sqrt{n}\ln\lambda,
\]
and using $\refpar{eq:Z_n_expr},$
\[
\ln R_{n}=\dfrac{1}{\sqrt{n}}\left[Z_{0}+\sum^{n}_{j=1}\left(\ln\left[1+\mu+\sigma\epsilon_{j}\right]-\ln\lambda\right)\right].
\]
Moreover, the random variables $\ln\left[1+\mu+\sigma\epsilon_{j}\right]$
are independent, following the same law, and have a finite second
order moment. Their expectation is
\[
\mathbb{E}\left(\ln\left[1+\mu+\sigma\epsilon_{n}\right]\right)=\dfrac{1}{2}\left[\ln\left(1+\mu+\sigma\right)+\ln\left(1+\mu-\sigma\right)\right]=\ln\lambda.
\]
and their second order moment is
\[
\mathbb{E}\left(\left(\ln\left[1+\mu+\sigma\epsilon_{n}\right]\right)^{2}\right)=\dfrac{1}{2}\left[\left(\ln\left(1+\mu+\sigma\right)\right)^{2}+\left(\ln\left(1+\mu-\sigma\right)\right)^{2}\right].
\]
Hence, their variance is
\begin{multline*}
\sigma^{2}_{\ln\left[1+\mu+\sigma\epsilon_{n}\right]}=\dfrac{1}{2}\left[\left(\ln\left(1+\mu+\sigma\right)\right)^{2}+\left(\ln\left(1+\mu-\sigma\right)\right)^{2}\right]\\
-\dfrac{1}{4}\left[\ln\left(1+\mu+\sigma\right)+\ln\left(1+\mu-\sigma\right)\right]^{2}.
\end{multline*}
Hence,
\begin{multline*}
\sigma^{2}_{\ln\left[1+\mu+\sigma\epsilon_{n}\right]}=\dfrac{1}{2}\left[\left(\ln\left(1+\mu+\sigma\right)\right)^{2}+\left(\ln\left(1+\mu-\sigma\right)\right)^{2}\right.\\
\left.-2\ln\left(1+\mu+\sigma\right)\ln\left(1+\mu-\sigma\right)\right].
\end{multline*}

Equivalently,
\[
\sigma^{2}_{\ln\left[1+\mu+\sigma\epsilon_{n}\right]}=\dfrac{1}{4}\left[\ln\dfrac{1+\mu+\sigma}{1+\mu-\sigma}\right]^{2}=\rho^{2}.
\]
By the central limit theorem,
\[
\dfrac{1}{\sqrt{n}\left|\rho\right|}\left[\sum^{n}_{j=1}\left(\ln\left(1+\mu+\sigma\epsilon_{j}\right)-\ln\lambda\right)\right]\stackrel[n\to+\infty]{\mathscr{L}}{\longrightarrow}\mathscr{N}_{\mathbb{R}}\left(0,1\right).
\]
Thus, the sequence of laws $P_{\ln R_{n}}$ converges narrowly to
the law $\mathscr{N}_{\mathbb{R}}\left(0,\rho^{2}\right).$ Therefore,
for every $f\in\mathscr{C}_{b}\left(\mathbb{R}\right),$ 
\[
\lim_{n\to+\infty}\intop_{\mathbb{R}}f\left(R_{n}\right)\text{d}P=\intop_{\mathbb{R}}f\left(\exp\left(\ln R_{n}\right)\right)\text{d}P=\intop_{\mathbb{R}}f\left(\text{e}^{x}\right)\dfrac{1}{\left|\rho\right|\sqrt{2\pi}}\text{e}^{-\frac{x^{2}}{2\rho^{2}}}\text{d}x,
\]
since $f\circ\exp\in\mathscr{C}_{b}\left(\mathbb{R}\right).$ By making
the change of variables from $\mathbb{R}$ onto $\mathbb{R}^{+\ast}$
defined by $y=\text{e}^{x},$ we obtain
\[
\lim_{n\to+\infty}\intop_{\mathbb{R}}f\left(R_{n}\right)\text{d}P=\intop_{\mathbb{R}^{+\ast}}f\left(y\right)\dfrac{1}{\left|\rho\right|\sqrt{2\pi}}\dfrac{1}{y}\text{e}^{-\frac{\left(\ln y\right)^{2}}{2\rho^{2}}}\text{d}y.
\]
Hence the sequence of laws $P_{R_{n}}$ converges narrowly to the
probability with density with respect to the Lebesgue measure, the
function
\[
y\mapsto\boldsymbol{1}_{\mathbb{R}^{+\ast}}\left(y\right)\dfrac{1}{\left|\rho\right|\sqrt{2\pi}}\dfrac{1}{y}\text{e}^{-\frac{\left(\ln y\right)^{2}}{2\rho^{2}}}.
\]
This is the \textbf{Log-normal law\index{Log-normal law}\mindex{law!Log-normal}}
with parameters $0$ and $\rho^{2}=\dfrac{1}{4}\left[\ln\dfrac{1+\mu+\sigma}{1+\mu-\sigma}\right]^{2}.$

\end{solution}

\chapter{Markov Chain}\label{chap:PartIIChap17}

\begin{objective}{}{}

Chapter \ref{chap:PartIIChap17} is devoted to the study of Markov
chains, an important class of discrete-time stochastic processes.
\begin{itemize}
\item Section \ref{sec:Introduction} presents guiding examples used throughout
the chapter to illustrate Markov chains. The first example concerns
the Bernoulli-Laplace gaze diffusion model, the second is the Ehrenfest
heat diffusion model between two insulated bodies, the third is the
P\'olya contagion model, and the last example is a particular case
of random walks.
\item Section \ref{sec:Conditional-Independence} focuses on conditional
independence of $\sigma-$algebra. It starts with the definition and
then gives a characterization of the conditional independence of two
\salgs with respect to a third one. A theorem is then stated concerning
the \salg generated by two \salgs that are conditionally independent
with respect to a third one.
\item Section \ref{sec:Markov-Chains:-General} contains three subsections.
\begin{itemize}
\item The first subsection concerns the Markov property and transition matrices.
It starts by defining Markov chains, as well as random variables with
the Markov property. A characterization of Markov chains is stated
as a theorem. An example then focuses on an autoregressive process.
To facilitate the study of Markov chains, the definition of the transition---or
stochastic---matrix is given, before stating the Chapman-Kolmogorov
equation. Homogeneous Markov chains are then introduced, together
with a proposition providing a criterion---via the conditional expectation---to
prove that a Markov chain is homogeneous, illustrated by examples.
A version of the Chapman-Kolmogorov equation is then given for homogeneous
Markov chains, as well as a necessary and sufficient condition to
be such a chain.
\item The second subsection first focuses on the simple Markov property
and finite-dimensional laws, before characterizing homogeneous Markov
chains with respect to the natural filtration.
\item The third subsection focuses only on homogeneous Markov chains. It
starts by defining chains with a given initial law. A first result
explains how we can start a chain with any initial law, when we know
how to start a chain from any point. This is followed by a result
on changing the starting point of a homogeneous Markov chain. The
single Markov property is then restated, before stating the strong
Markov property.
\end{itemize}
\item Section \ref{sec:Visits-with-one} deals with visits to a fixed state.
It starts with the study of the sequence of hitting times of a given
point, and continues with the law of the number of visits to that
point, and of the first hitting time. The potential matrix is then
defined, before introducing the potential of a path; in particular,
the potential matrix is expressed in terms of the transition matrix.
\item Section \ref{sec:State-Classification} concerns state classification.
First, it defines reachability---how a state leads to another---and
establishes the transitivity of this relation. Communication states
are then introduced, and communication is shown to be an equivalence
relation; its equivalence classes are called the communication classes.
The concept of an irreducible chain derives. Periodic and aperiodic
states are then defined, together with results on the periodicity
of communication classes. \\
Second, recurrence is studied: recurrent, null recurrent, positive
recurrent, and transient states are defined, followed by results on
the probability of multiple visits, the probability of infinitely
many visits, and the probability of returning to the same state. A
state-classification theorem is then stated. \\
Third, asymptotic behaviour is examined, depending on the type of
state targeted by the trajectory: from transient states to recurrent
states, via aperiodic states. The nature of states within a communication
class and the nature of the class are then discussed. Closed communication
classes are defined, their possible types are listed, their possible
types are listed, and inessential states are introduced, and characterized
as transient states. \\
Fourth, an analytic recurrence criterion is examined, giving necessary
and sufficent conditions for a chain to be recurrent.
\item Section \ref{sec:Computation-of-the} begins with the computation
of the potential matrix, depending on the nature of the targeted state.
The computation of the probability that the first hitting time is
finite is then considered.
\item Section \ref{sec:Invariant-Measures} introduces invariant measures,
linked to the asymptotic behaviour of a homogeneous Markov chain.
The existence and uniqueness of an invariant probability are then
discussed. A positive recurrence criterion is given that characterizes
when a homogeneous Markov chain admits a unique invariant probability,
illustrated by the Ehrenfest diffusion model. A necessary and sufficient
condition for the existence of a limiting probability is then examined,
with an illustration on a genetic model.
\item Section \ref{sec:Strong-Law-of} begins with the Chacon-Orstein theorem,
before stating a strong law of large numbers for homogeneous Markov
chains. It is then applied to the estimation of the transition matrix.
\end{itemize}
\end{objective}

\textbf{Markov chains form an important class of discrete-time stochastic
processes. They model temporal random phenomena for which the probabilistic
evolution at a given time depends only on the current state of the
system, and not on the entire past history. In other words, Markov
chains describe memoryless phenomena. In this chapter, we restrict
attention to Markov chains with a countable state space.}

\section{Introduction}\label{sec:Introduction}

In Part \ref{part:Introduction-to-Probability}, Section \ref{sec:Evolutive-phenomena-modeling},
we proved that an evolutive phenomenon, finite both in time and in
space may be modeled as a Markovian process, that is, as a process
without memory beyond the present state. 

The purpose of this introduction is to illustrate this elementary
formalization by the study of a classical historical model, and to
emphasize the mathematical difficulties that arise in its axiomatization.
This discussion naturally leads to the definition of homogeneous Markov
chains given in $\ref{pr:cns_hom_Markov_chain}$.

We begin with the gaze diffusion model known as the \textbf{Bernoulli-Laplace
model\index{Bernoulli-Laplace model}}.

\begin{example}{Bernoulli-Laplace Gaze Diffusion Model}{}

Two urns, labelled 1 and 2, each contains $m$ balls. Among the total
of $2m$ balls, $r,$ with $1\leqslant r\leqslant m,$ are red and
$2m-r$ are blue. Time scale is chosen discrete and, after reindexing,
is identified with $\mathbb{N}.$ At each time step, one ball is drawn
independently at random from each urn, and the two selected balls
are exchanged, i.e. each ball is placed into the other urn.

We choose to represent the state of the system at each instant $n,$
$n\geqslant1,$ by the number $X_{n}$ of red balls in Urn 1 after
the $n-$th exchange of the balls in the urns. The initial state is
denoted $X_{0}.$ The set of possible states is the finite integer
interval $E=\left\llbracket 0,r\right\rrbracket .$

Thus, the state $X_{n}$ may be viewed as a random point moving on
$E.$ From one time step to the next, this point can either remain
at its current position or move to one of its nearest neighbours.
It is customary to represent this evolution by a graph whose vertices
are the elements of $E,$ with edges connecting each state to itself
and to its adjacent states, as represented in the figure just below.

\begin{center}
\begin{tikzpicture}[
  scale=0.75,
  transform shape,
  >=Stealth,
  every node/.style={
    circle,
    draw=black,
    minimum size=10mm,
    inner sep=1pt
  },
  edge/.style={
    ->,
    draw=black,
    line width=0.75pt
  }
]

% Nodes (rectangular layout)
\node (a0) at (0,0) {0};
\node (a1) at (2,0) {1};
\node[draw=none] (a2) at (4,0) {};
\node[draw=none] at (4.25,0)  {\(\cdot\)};
\node[draw=none] at (4.5,0)  {\(\cdot\)};
\node[draw=none] at (4.75,0)  {\(\cdot\)};
\node[draw=none] (a3) at (5.0,0) {};
\node (a4) at (7,0) {$i-1$};
\node (a5) at (9,0) {$i$};
\node (a6) at (11,0) {$i+1$};
\node[draw=none] (a7) at (13.0,0) {};
\node[draw=none] at (13.25,0)  {\(\cdot\)};
\node[draw=none] at (13.5,0)  {\(\cdot\)};
\node[draw=none] at (13.75,0)  {\(\cdot\)};
\node[draw=none] (a8) at (14.0,0) {};
\node (a9) at (16.0,0) {$r$};

% Transitions
\draw[edge] (a0) to[bend right=20] (a1);
\draw[edge] (a1) to[bend right=20] (a0);
\draw[edge] (a1) to[bend right=20] (a2);
\draw[edge] (a2) to[bend right=20] (a1);
\draw[edge] (a3) to[bend right=20] (a4);
\draw[edge] (a4) to[bend right=20] (a3);
\draw[edge] (a4) to[bend right=20] (a5);
\draw[edge] (a5) to[bend right=20] (a4);
\draw[edge] (a5) to[bend right=20] (a6);
\draw[edge] (a6) to[bend right=20] (a5);
\draw[edge] (a6) to[bend right=20] (a7);
\draw[edge] (a7) to[bend right=20] (a6);
\draw[edge] (a8) to[bend right=20] (a9);
\draw[edge] (a9) to[bend right=20] (a8);

%Self-loops
\draw[edge] (a0) to[out=120,in=60,looseness=6] (a0);
\draw[edge] (a1) to[out=120,in=60,looseness=6] (a1);
\draw[edge] (a4) to[out=120,in=60,looseness=6] (a4);
\draw[edge] (a5) to[out=120,in=60,looseness=6] (a5);
\draw[edge] (a6) to[out=120,in=60,looseness=6] (a6);
\draw[edge] (a9) to[out=120,in=60,looseness=6] (a9);
\end{tikzpicture}
\end{center}

It is intuitively clear that the \textbf{process} $\left(X_{n}\right)_{n\in\mathbb{N}^{\ast}}$
is \textbf{Markovian}\index{Markovian process}; that is, for every
$n\geqslant0,$ and every $\left(n+2\right)-$uple states $x_{0},x_{1},\cdots,x_{n+1},$
it verifies 
\begin{equation}
P\left(X_{n+1}=x_{n+1}\mid X_{0}=x_{0},\cdots,X_{n}=x_{n}\right)=P\left(X_{n+1}=x_{n+1}\mid X_{n}=x_{n}\right),\label{eq:cond_P_markovian}
\end{equation}
these conditional probabilities being understood in the elementary
sense, that is when the conditioned event has intuitively nonzero
probability. The evolution of the process $\left(X_{n}\right)_{n\in\mathbb{N}^{\ast}}$
is therefore refined by determining these conditional probabilities,
called \textbf{\index{transition probability}transition probabilities}.

Assume that $X_{n}=i.$ Before the $n+1$-th drawing, urn 1 contains
$i$ red balls and $m-i$ blue balls, while urn 2 contains $r-i$
red balls and $m-\left(r-i\right)$ blue balls. 

For $j=1,2,$ denote by $R^{n+1}_{j}$ and $B^{n+1}_{j},$ the events
that the drawn ball at the $n+1$-th drawing from urn $j$ is respectively
red or blue. The independence and the uniformity of the drawings lead
to the tables below.
\begin{itemize}
\item If $1\leqslant i\leqslant r-1,$\\
\begin{center}%
\begin{tabular}{|c|c|c|}
\hline 
\begin{cellvarwidth}[t]
\centering
Configurations\\
of the $n+1$-th drawing
\end{cellvarwidth} &
\begin{cellvarwidth}[t]
\centering
Transition from one\\
state to another
\end{cellvarwidth} &
\begin{cellvarwidth}[t]
\centering
Probability of\\
transition
\end{cellvarwidth}\tabularnewline
\hline 
$R^{n+1}_{1}R^{n+1}_{2}$ &
$i\rightarrow i$ &
$\dfrac{i}{m}\cdot\dfrac{r-i}{m}$\tabularnewline
\hline 
$R^{n+1}_{1}B^{n+1}_{2}$ &
$i\rightarrow i-1$ &
$\dfrac{i}{m}\cdot\dfrac{m-\left(r-i\right)}{m}$\tabularnewline
\hline 
$B^{n+1}_{1}R^{n+1}_{2}$ &
$i\rightarrow i+1$ &
$\dfrac{m-i}{m}\cdot\dfrac{r-i}{m}$\tabularnewline
\hline 
$B^{n+1}_{1}B^{n+1}_{2}$ &
$i\rightarrow i$ &
$\dfrac{m-i}{m}\cdot\dfrac{m-\left(r-i\right)}{m}$\tabularnewline
\hline 
\end{tabular}\end{center}
\item If $i=0,$\\
\begin{center}%
\begin{tabular}{|c|c|c|}
\hline 
\begin{cellvarwidth}[t]
\centering
Configurations\\
of the $n+1$-th drawing
\end{cellvarwidth} &
\begin{cellvarwidth}[t]
\centering
Transition from one\\
state to another
\end{cellvarwidth} &
\begin{cellvarwidth}[t]
\centering
Probability of\\
transition
\end{cellvarwidth}\tabularnewline
\hline 
$B^{n+1}_{1}R^{n+1}_{2}$ &
$0\rightarrow1$ &
$\dfrac{r}{m}$\tabularnewline
\hline 
$B^{n+1}_{1}B^{n+1}_{2}$ &
$0\rightarrow0$ &
$\dfrac{m-r}{m}$\tabularnewline
\hline 
\end{tabular}\end{center}
\item If $i=r,$\\
\begin{center}%
\begin{tabular}{|c|c|c|}
\hline 
\begin{cellvarwidth}[t]
\centering
Configurations\\
of the $n+1$-th drawing
\end{cellvarwidth} &
\begin{cellvarwidth}[t]
\centering
Transition from one\\
state to another
\end{cellvarwidth} &
\begin{cellvarwidth}[t]
\centering
Probability of\\
transition
\end{cellvarwidth}\tabularnewline
\hline 
$R^{n+1}_{1}B^{n+1}_{2}$ &
$r\rightarrow r-1$ &
$\dfrac{r}{m}$\tabularnewline
\hline 
$B^{n+1}_{1}B^{n+1}_{2}$ &
$r\rightarrow r$ &
$\dfrac{m-r}{m}$\tabularnewline
\hline 
\end{tabular}\end{center}
\end{itemize}
In all three cases, we observe that the probabilities of passage from
one state to the other can be written in a unique form. The conditional
probabilities sought, called \textbf{transition probabilities}\index{transition probability},
are therefore given, for every $i$ such that $0\leqslant i\leqslant r,$
by
\[
\left\{ \begin{array}{lll}
P\left(X_{n+1}=i\mid X_{n}=i\right) & = & \dfrac{i}{m}\cdot\dfrac{r-i}{m}+\dfrac{m-i}{m}\cdot\dfrac{m-\left(r-i\right)}{m},\\
P\left(X_{n+1}=i-1\mid X_{n}=i\right) & = & \dfrac{i}{m}\cdot\dfrac{m-\left(r-i\right)}{m},\\
P\left(X_{n+1}=i+1\mid X_{n}=i\right) & = & \dfrac{m-i}{m}\cdot\dfrac{r-i}{m}.
\end{array}\right.
\]
It is customary to gather these transition probabilities into a matrix
$M=\left(M_{ij}\right),$ called the \textbf{transition matrix}\index{transition matrix},
whose general term is defined by\boxeq{
\[
M_{ij}=P\left(X_{n+1}=j\mid X_{n}=i\right).
\]
}

We shall see in Proposition $\ref{pr:P_Markov_transition_mat}$ once
the transition matrix and the law of the initial random variable $X_{0}$
are given, they completely determine the law of the random variable
$\left(X_{0},X_{1},\cdots,X_{n}\right),$ and, hence fully determine
the probabilistic behaviour of the process $X.$ 

In particular, after developing the theory of Markov chains, it will
be possible to describe rapidly the qualitative \textbf{asymptotic
behaviour\mindex{asymptotic!behaviour}} of this process, and to determine
the \textbf{limits of the probabilities\index{limits of the probabilities}}
of being in any given state as time tends to infinity. This problem
was historically solved by \textbf{Bernoulli\sindex[fam]{Bernoulli, Daniel}}
and \textbf{\sindex[fam]{Laplace, Pierre-Simon}Laplace}.

Let us also briefly mention two other classical models. The first
concerns \textbf{heat exchange between two isolated bodies}, and is
known as the \textbf{\sindex[fam]{Ehrenfest, Paul}Ehrenfest\footnotemark
model}\index{Ehrenfest model}. The second, due to \textbf{P\'olya}\sindex[fam]{Pólya, George@P\'olya, George}\mindex{Polya model@P\'olya model},
models the \textbf{propagation of contagious diseases}. These physical
phenomena are likewise represented by urn models involving ball drawings.
They are described in detail in Feller's book \cite{feller1958introduction},
and in many more recent references. These models will serve throughout
the chapter as illustrative examples of the concepts and results that
are introduced.

\end{example}

\textbf{\footnotetext{\href{https://fr.wikipedia.org/wiki/Paul_Ehrenfest}{Paul Ehrenfest}\sindex[fam]{Ehrenfest, Paul}}
(1880 - 1933) was an Austrian theoretical physicist. He did major
contributions within statisical mechanic and quantum mechanic. He
taught at Leyde University.}

\begin{figure}[t]
\begin{center}\includegraphics[width=0.4\textwidth]{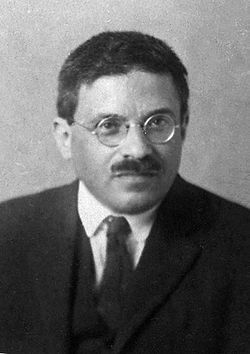}

{\tiny Unknown author---\href{http://www.sil.si.edu/digitalcollections/hst/scientific-identity/cf/by_name_display_results.cfm?scientist=Ehrenfest,\%2520Paul}{Dibner Library}---Public
domain}\end{center}

\caption{\textbf{\protect\href{https://fr.wikipedia.org/wiki/Paul_Ehrenfest}{Paul Ehrenfest}}
(1880 - 1933)}\sindex[fam]{Ehrenfest, Paul}
\end{figure}

\begin{example}{Heat Diffusion Model Between Two Insulated Bodies of Ehrenfest\footnotemark}{ehrenfest_model}

Two receptacles, labelled 1 and 2, contain a total of $m$ particles
which can diffuse from one receptacle to the other. The particle diffusion
phenomenon can be envisioned as, at each unit of time, the choice
at random of a particle in one receptacle and its transfer in the
other one. Repeating in a similar way these choices and transfers,
we are interested in the \textbf{distribution} of the particles in
each receptacle, considered as urn, after $n$ steps.

The modelling under the form of draws of balls from an urn is then
the following. We identify the first receptacle particles to red balls
in quantity $r$, and the ones of the second receptacle to blue balls,
in quantity $m-r.$ To each draw, the drawn ball is replaced by a
ball of the opposite color. The chosen scale of time is discrete,
and after reindexing, it corresponds to $\mathbb{N}.$ The \textbf{state
at the instant $n$} is the \textbf{number of red balls} contained
in the urn.

\end{example}

\footnotetext{\cite{ehrenfest1907zwei}

}

\begin{figure}[t]
\begin{center}\includegraphics[width=0.4\textwidth]{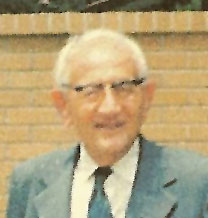}

{\tiny Thane Plambeck---\href{https://www.flickr.com/photos/thane}{https://www.flickr.com/photos/thane}---CC
BY 2.0}\end{center}

\caption{\textbf{\protect\href{https://en.wikipedia.org/wiki/George_P\%25C3\%25B3lya}{George Polya}}
(1887 - 1985)}\sindex[fam]{Pólya, George@P\'olya, George}
\end{figure}

\begin{example}{Diffusion Model of Contagious Diseases of P\'olya}{polya_diff_mod}

This model, introduced by P\'olya\sindex[fam]{Pólya, George@P\'olya, George}\footnotemark
describes the spread of contagious diseases. It reflects how the probability
of individual contagion increases or decreases with each appearance
or disappearance of a new case within a population.

The modelling, in terms of ball draws from an urn is as follows. Healthy
individuals are represented by red balls, in number $r,$ and sick
individuals by blue balls, in number $b.$ After each random draw,
the selected ball is returned to the urn with $c$ additional balls
of the same color. 

The \textbf{state at time $n$} is defined as the \textbf{proportion
$Y_{n}$ }of blue balls contained in the urn \textbf{after} the $n-$th
draw and after the $c$ balls have been added.

If this process evolves indefinitely, the set $E$ of possible states
is no longer finite, but countably infinite, and a priori \textbf{contains}
the set of rational numbers in the interval $\left[0,1\right].$ It
therefore becomes more difficult to determine, in an elementary manner,
whether the process $\left(Y_{n}\right)_{n\in\mathbb{N}^{\ast}}$
is Markovian, via a relation of the form $\refpar{eq:cond_P_markovian}.$

Indeed, constructing a model first requires specifying the set $E$
of possible states. Once this is done, and before any analysis of
the model, one cannot determine whether, for every $\left(n+1\right)-$uple
of states $y_{0},y_{1},\cdots,y_{n},$ the probability of the conditioning
event
\[
P\left(Y_{0}=y_{0},\cdots,Y_{n}=y_{n}\right)
\]
vanishes or not. 

This creates a difficulty in defining conditional probability in the
elementary sense and consequently calls into question the relevance
of defining the Markov property through relations of the type $\refpar{eq:cond_P_markovian}$---namely
the idea that the process is memoryless except for its present state.
We are therefore led to adopt a definition that is more mathematically
appropriate. Under this refined definition, the P\'olya process $\left(Y_{n}\right)_{n\in\mathbb{N}}$
is indeed \textbf{Markovian}, although this is not entirely intuitive.

\end{example}

\footnotetext{\textbf{\href{https://en.wikipedia.org/wiki/George_P\%25C3\%25B3lya}{George Polya}}\sindex[fam]{Pólya, George@P\'olya, George}
(1887 - 1985) was an Hungarian-American mathematician. He taught at
ETH Zürich from 1914 to 1940 and at Stanford University from 1940
to 1953. He made several key contributions to combinatorics, number
theory, numerical analysis and probability theory.}

Another example, which is a particular case of \textbf{random walks}---these
will be presented later---once again highlights the same issue in
defining the \textbf{Markov property}. In this case, however the Markov
property is intuitively clear, provided one understands it as meaning
that \textbf{the probabilistic evolution of the process after each
instant $n$ does not depend on the past beyond the present state
at the time $n.$} This example also emphasizes the importance of
the initial law.

\begin{example}{}{}

Let $\left(X_{n}\right)_{n\in\mathbb{N}}$ be a sequence of independent
random variables defined on a probabilized space \preds, taking values
in $\mathbb{Z}.$ Suppose that, for every $n\in\mathbb{N}^{\ast},$
the random variables $X_{n}$ follow the same law $p\delta_{-2}+q\delta_{2}$
with $p+q=1$ and $p,q>0.$ 

We define
\[
S_{n}=\sum^{n}_{j=0}X_{j}\,\,\,\,\text{and}\,\,\,\,Y_{n}=\sum^{n}_{j=0}X^{2}_{j}.
\]
Hence, $S_{n+1}=S_{n}+X_{n+1}$ and $Y_{n+1}=Y_{n}+X^{2}_{n+1}.$
The quantity $S_{n}$ can be interpreted as the \textbf{position of
a particle} that, at each unit of time, jumps from one integer to
another. The process $\left(S_{n}\right)_{n\in\mathbb{N}}$ is therefore
a \textbf{random walk}.

If $X_{0}=0,$ then $S_{n}$ takes values in $2\mathbb{Z}.$ A straightforward
convolution argument shows that every even integers between $-2n$
and $2n,$ and only those integers, is visited by $S_{n}$ with nonzero
probability. 

However, if $X_{0}$ follows the law $\dfrac{1}{2}\left(\delta_{0}+\delta_{1}\right),$
then $S_{n}$ takes values in $\mathbb{Z}.$ Some trajectories of
process $\left(S_{n}\right)_{n\in\mathbb{N}}$ are contained in $2\mathbb{Z},$
while others are contained in $2\mathbb{Z}+1.$ In any case, for every
$\left(n+1\right)-$uple $x_{0},x_{1},\cdots,x_{n}$ of consecutive
elements of the state space $\mathbb{Z},$ the probability of the
conditioning event $P\left(X_{0}=x_{0},\cdots,X_{n}=x_{n}\right)$
appearing in a relation of the form $\refpar{eq:cond_P_markovian}$
is equal to zero. Thus, once again, defining the Markov property through
such relations fails in this setting.

Finally, the probabilistic evolution of the process $\left(Y_{n}\right)_{n\in\mathbb{N}}$
at any time $n$ does not depend on the state of the process at that
instant. In this sense, the process $\left(Y_{n}\right)_{n\in\mathbb{N}}$
has the Markov property. Moreover, the probabilistic \textbf{``history}''
\textbf{at time} $n$ can be observed in more or less exhaustive way,
depending on whether one has access to the values of the $X_{j},$
or only to those of $X^{2}_{j},\,0\leqslant j\leqslant n.$ Nevertheless,
this history influences the probability that $Y_{n+1}$ takes a given
value at time $n,$ only through the knowledge of $Y_{n}$ itself.
To account for such situations, \textbf{Markov chains will be defined
with respect to filtrations}.

\end{example}

We now proceed to define the Markov property using the notion of \textbf{conditional
independence of \salg,} a concept particular useful in several contexts,
in the following section.

\section{Conditional Independence}\label{sec:Conditional-Independence}

\begin{denotations}{}{}

In this chapter, we denote $\mathbb{E}^{X=\cdot}\left(Y\right)$ the
\textbf{conditional mean\index{conditional mean}} function of the
random variable $Y$ conditionally to the random variable $X.$ This
should not be confused\footnotemark with the \textbf{conditional
expectation} of $Y$ with respect to the \salg $\sigma\left(X\right)$
generated by $X,$ denoted $\mathbb{E}^{\sigma\left(X\right)}\left(Y\right)$
or $\mathbb{E}\left(Y\mid X\right),$ depending on the context. This
latter is a \textbf{class of random variables}; as usual, we denote
both a representative and its equivalence class by the same symbol.

\end{denotations}

We begin with an example that clarifies the notion of conditional
dependency, which we will define and study in full generality immediately
afterwards.

Let $X,Y,Z$ be three independent real-valued random variables defined
on a probabilized space \preds. Define
\[
U=X+Y\,\,\,\,\text{and}\,\,\,\,V=XZ.
\]

Let $f$ and $g$ be arbitrary bounded Borel functions. Since the
random variables $X$ and $\left(Y,Z\right)$ are independent, the
conditional transfer theorem---Theorem $\ref{th:cond_transfer}$---implies
that, the \textbf{conditional mean} of $f\left(U\right)g\left(V\right)$
conditionally to the random variable $X,$ verifies, for $P_{X}-$almost
every $x,$
\[
\mathbb{E}^{X=x}\left[f\left(U\right)g\left(V\right)\right]=\mathbb{E}^{X=x}\left[f\left(x+Y\right)g\left(xZ\right)\right]=\mathbb{E}\left[f\left(x+Y\right)g\left(xZ\right)\right].
\]
Taking into account the independence of the random variables $Y$
and $Z,$ we obtain
\[
\mathbb{E}^{X=x}\left[f\left(U\right)g\left(V\right)\right]=\mathbb{E}\left[f\left(x+Y\right)\right]\mathbb{E}\left[g\left(xZ\right)\right].
\]
In particular, by succesively taking for $f$ and $g$ the constant
function equal to 1, we deduce
\[
\mathbb{E}^{X=x}\left[f\left(U\right)\right]=\mathbb{E}\left[f\left(x+Y\right)\right]\,\,\,\,\text{and}\,\,\,\,\mathbb{E}^{X=x}\left[g\left(V\right)\right]=\mathbb{E}\left[g\left(xZ\right)\right].
\]
Hence, for $P_{X}-$almost every $x,$\boxeq{
\[
\mathbb{E}^{X=x}\left[f\left(U\right)g\left(V\right)\right]=\mathbb{E}^{X=x}\left[f\left(U\right)\right]\mathbb{E}^{X=x}\left[g\left(V\right)\right].
\]
}

We therefore obtain the equality of conditional expectations\boxeq{
\[
\mathbb{E}^{\sigma\left(X\right)}\left[f\left(U\right)g\left(V\right)\right]=\mathbb{E}^{\sigma\left(X\right)}\left[f\left(U\right)\right]\mathbb{E}^{\sigma\left(X\right)}\left[g\left(V\right)\right].
\]
}

These two relations express, in equivalent forms, the conditional
independence of the random variables $U$ and $V$ given the \textbf{random
variable $X;$} equivalent, the \textbf{independence} of the generated
\salg by the random variables $U$ and $V,$ \textbf{conditionally}
to the \salg $\sigma\left(X\right)$ generated by the random variable
$X.$

We now give the general definition of \textbf{conditional independence
of \salgs,} a fundamental notion used later to define the \textbf{Markov
property}. We then study some of its properties.

All \salgs considered below are \textbf{sub-$\sigma$-algebras}.
defined on the same probabilized space \preds.

\begin{definition}{}{}

Let $\mathscr{A}_{i},\,i=1,2,3$ be three $\sigma-$algebras. The
\salgs $\mathscr{A}_{1}$ and $\mathscr{A}_{3}$ are said to be conditionally
independent with respect to $\mathscr{A}_{2}$ if, for $i=1,3$ and
for every bounded $\mathscr{A}_{i}-$measurable real-valued random
variables $Y_{i},$ which we denote $Y_{i}\in b\mathscr{A}_{i},$
the following holds:
\begin{equation}
\mathbb{E}^{\mathscr{A}_{2}}\left(Y_{1}Y_{3}\right)=\mathbb{E}^{\mathscr{A}_{2}}\left(Y_{1}\right)\mathbb{E}^{\mathscr{A}_{2}}\left(Y_{3}\right).\label{eq:17_2}
\end{equation}

In particular, if $\mathscr{A}_{2}$ is the $\sigma-$algebra $\sigma\left(X\right)$
generated by a random variable $X,$ we simply say that the \salgs
$\mathscr{A}_{1}$ and $\mathscr{A}_{3}$ are conditionally independent
with respect to $X.$ The relation $\refpar{eq:17_2}$ may then be
written as
\begin{equation}
\mathbb{E}\left(Y_{1}Y_{3}\mid X\right)=\mathbb{E}^{\mathscr{A}_{2}}\left(Y_{1}\mid X\right)\mathbb{E}^{\mathscr{A}_{2}}\left(Y_{3}\mid X\right).\label{eq:17_3}
\end{equation}

\end{definition}

\begin{remark}{}{}

Conditional independence of the \salgs $\mathscr{A}_{1}$ and $\mathscr{A}_{3}$
with respect to $\mathscr{A}_{2}$ does not imply their independence.
However, if $\mathscr{A}_{2}$ is the trivial \salg $\left\{ \Omega,\emptyset\right\} ,$
then conditional independence of the \salgs $\mathscr{A}_{1}$ and
$\mathscr{A}_{3}$ with respect to $\mathscr{A}_{2}$ is equivalent
to ordinary independence.

\end{remark}

We denote by $P^{\mathscr{B}}\left(A\right)$ the conditional probability
of an event $A$ given the \salg $\mathscr{B}$ defined as $\mathbb{E}^{\mathscr{B}}\left(\boldsymbol{1}_{A}\right).$
If $\mathscr{A}_{2}$ is the \salg $\sigma\left(X\right)$ generated
by a random variable $X,$ we simply write $P\left(A\mid X\right).$

\begin{lemma}{}{}

The \salg $\mathscr{A}_{1}$ and $\mathscr{A}_{3}$ are conditionally
independent with respect to $\mathscr{A}_{2}$ if and only if, for
every events $A_{1}\in\mathscr{A}_{1}$ and $A_{3}\in\mathscr{A}_{3},$
\begin{equation}
P^{\mathscr{A}_{2}}\left(A_{1}\cap A_{3}\right)=P^{\mathscr{A}_{2}}\left(A_{1}\right)P^{\mathscr{A}_{2}}\left(A_{3}\right).\label{eq:cond_prob_and_indep}
\end{equation}

In particular, if $\mathscr{A}_{2}$ is the \salg $\sigma\left(X\right)$
generated by a random variable $X,$ the \salgs $\mathscr{A}_{1}$
and $\mathscr{A}_{3}$ are conditionally independent with respect
to $X$ if and only if for every events $A_{1}\in\mathscr{A}_{1}$
and $A_{3}\in\mathscr{A}_{3},$
\begin{equation}
P\left(A_{1}\cap A_{3}\mid X\right)=P\left(A_{1}\mid X\right)P\left(A_{3}\mid X\right).\label{eq:cond_prod_and_indep_with_sigmaX}
\end{equation}

\end{lemma}

\begin{proof}{}{}

The necessary condition is straightforward. Conversely, from relation
$\refpar{eq:cond_prob_and_indep},$ we deduce that $\refpar{eq:17_2}$
holds for every $\mathscr{A}_{i}-$measurable step random variables
$Y_{i},$ with $i=1,3.$ The general case follows via the usual process
of integration. 

\end{proof}

The previous lemma is mainly technical. The following theorem, however,
is fundamental for the study of Markov chains. 

\begin{theorem}{Generated \salg of Two \salgs and Independence}{gen_salg_two_salgs_ind}

Let $\mathscr{A}_{1,2}$ ---or also $\mathscr{A}_{1}\vee\mathscr{A}_{2}$---be
the \salg generated by $\mathscr{A}_{1}$ and $\mathscr{A}_{2}.$
The \salgs $\mathscr{A}_{1}$ and $\mathscr{A}_{3}$ are conditionally
independent with respect to $\mathscr{A}_{2}$ if and only if, for
every $Y_{3}\in\mathscr{L}^{1}\left(\Omega,\mathscr{A}_{3},P\right),$
\begin{equation}
\mathbb{E}^{\mathscr{A}_{1,2}}\left(Y_{3}\right)=\mathbb{E}^{\mathscr{A}_{2}}\left(Y_{3}\right)\label{eq:cond_ind_relatively_to_salgA_2}
\end{equation}

\end{theorem}

\begin{remark}{}{}

In fact, by standard arguments of integration, conditional independence
holds if and only if relation $\refpar{eq:cond_ind_relatively_to_salgA_2}$
is satisfied for every $Y_{3}\in b\mathscr{A}_{3}.$

\end{remark}

\begin{proof}{}{}

\textbf{Necessity}

We use the principle of extension by measurability. We easily verify
that the family of events
\[
\mathscr{S}=\left\{ A\in\mathscr{A}_{1,2}:\,\intop_{A}Y_{3}\text{d}P=\intop_{A}\mathbb{E}^{\mathscr{A}_{2}}\left(Y_{3}\right)\text{d}P\,\,\,\,\forall Y_{3}\in\mathscr{L}^{1}\left(\Omega,\mathscr{A}_{3},P\right)\right\} 
\]
is a $\lambda-$system. We are going to prove that it contains a $\pi-$system
$\mathscr{C}$ which generates $\mathscr{A}_{1,2}$ defined by
\[
\mathscr{C}=\left\{ A_{1}\cap A_{2}\in\mathscr{A}_{1,2}:\,A_{1}\in\mathscr{A}_{1}\,\,\,\,\text{and}\,\,\,\,A_{2}\in\mathscr{A}_{2}\right\} .
\]
This will show that $\mathscr{S}$ contains $\mathscr{A}_{1,2}$ and
thus that $\refpar{eq:cond_ind_relatively_to_salgA_2}$ holds for
every $Y_{3}\in\mathscr{L}^{1}\left(\Omega,\mathscr{A}_{3},P\right).$

Let $A_{1}\in\mathscr{A}_{1}$ and $A_{2}\in\mathscr{A}_{2}$ be arbitrary.
Using the measurability of concerned random variables with respect
to the \salgs $\mathscr{A}_{1,2}$ and $\mathscr{A}_{2},$ we obtain
\[
\mathbb{E}\left(\boldsymbol{1}_{A_{1}}\boldsymbol{1}_{A_{2}}\mathbb{E}^{\mathscr{A}_{1,2}}\left(Y_{3}\right)\right)=\mathbb{E}\left(\boldsymbol{1}_{A_{1}}\boldsymbol{1}_{A_{2}}Y_{3}\right)=\mathbb{E}\left(\boldsymbol{1}_{A_{2}}\mathbb{E}^{\mathscr{A}_{2}}\left(\boldsymbol{1}_{A_{1}}Y_{3}\right)\right),
\]
Using $\refpar{eq:17_2},$
\[
\mathbb{E}\left(\boldsymbol{1}_{A_{1}}\boldsymbol{1}_{A_{2}}\mathbb{E}^{\mathscr{A}_{1,2}}\left(Y_{3}\right)\right)=\mathbb{E}\left(\boldsymbol{1}_{A_{2}}\mathbb{E}^{\mathscr{A}_{2}}\left(\boldsymbol{1}_{A_{1}}\right)\mathbb{E}^{\mathscr{A}_{2}}\left(Y_{3}\right)\right).
\]
Since $\boldsymbol{1}_{A_{2}}\mathbb{E}^{\mathscr{A}_{2}}\left(Y_{3}\right)$
is $\mathscr{A}_{2}-$measurable,
\[
\mathbb{E}\left(\boldsymbol{1}_{A_{1}}\boldsymbol{1}_{A_{2}}\mathbb{E}^{\mathscr{A}_{1,2}}\left(Y_{3}\right)\right)=\mathbb{E}\left(\boldsymbol{1}_{A_{1}}\boldsymbol{1}_{A_{2}}\mathbb{E}^{\mathscr{A}_{2}}\left(Y_{3}\right)\right).
\]
Thus the necessary condition holds.

\textbf{Sufficiency}

Let, for $i=1,3,$ $Y_{i}\in b\mathscr{A}_{i}$ be arbitrary. Taking
into account the inclusion of \salgs $\mathscr{A}_{2}\subset\mathscr{A}_{1,2},$
and the $\mathscr{A}_{1,2}-$measurability of $Y_{1},$
\[
\mathbb{E}^{\mathscr{A}_{2}}\left(Y_{1}Y_{3}\right)=\mathbb{E}^{\mathscr{A}_{2}}\left(\mathbb{E}^{\mathscr{A}_{1,2}}\left(Y_{1}Y_{3}\right)\right)=\mathbb{E}^{\mathscr{A}_{2}}\left(Y_{1}\mathbb{E}^{\mathscr{A}_{1,2}}\left(Y_{3}\right)\right).
\]
By hypothesis $\refpar{eq:cond_ind_relatively_to_salgA_2},$
\[
\mathbb{E}^{\mathscr{A}_{2}}\left(Y_{1}Y_{3}\right)=\mathbb{E}^{\mathscr{A}_{2}}\left(Y_{1}\mathbb{E}^{\mathscr{A}_{2}}\left(Y_{3}\right)\right),
\]
and therefore
\[
\mathbb{E}^{\mathscr{A}_{2}}\left(Y_{1}Y_{3}\right)=\mathbb{E}^{\mathscr{A}_{2}}\left(Y_{1}\right)\mathbb{E}^{\mathscr{A}_{2}}\left(Y_{3}\right).
\]

\end{proof}

Using again the extension by measurability theorem, one can solve
the following example, let as exercise to the reader.

\begin{example}{}{}

Let $\mathscr{A}_{i},$ $i=1,2,3$ and $\mathscr{B}_{3}$ be four
\salgs. Suppose that $\mathscr{A}_{3}=\mathscr{B}_{3}\vee\mathscr{A}_{2}.$
The \salgs $\mathscr{A}_{1}$ and $\mathscr{A}_{3}$ are conditionally
independent with respect to $\mathscr{A}_{2}$ if and only if the
\salgs $\mathscr{A}_{1}$ and $\mathscr{B}_{3}$ are conditionally
independent with respect to $\mathscr{A}_{2}.$

\end{example}

\section{Markov Chains: General Properties}\label{sec:Markov-Chains:-General}

In this chapter, $E$ is a \textbf{countable} set---finite or infinite---equipped
with the \salg of its subsets, denoted $\mathscr{E}.$ Unless explicitly
stated otherwise, all processes are defined on the probabilized space
\preds.

\subsection{Markov Property. Transition Matrices.}

\begin{definition}{Wide Future. Strict Future. Present $\sigma-$algebra. Markov Chain. Markov Property}{}

Let $X=\left(X_{n}\right)_{n\in\mathbb{N}}$ be a process taking values
in $\left(E,\mathscr{E}\right)$\footnotemark. 

For every $n\in\mathbb{N},$ the $\sigma-$algebra 
\[
\mathscr{F}^{l}_{n}=\sigma\left(X_{j}:j\geqslant n\right)
\]
is called the \salg of \textbf{wide future\mindex{sigma-algebra @ \salg ! wide future}}---or\textbf{
\mindex{sigma-algebra @ \salg !non-strict future}non-strict future}---of
process $X$ after time $n.$

For every $n\in\mathbb{N},$ the $\sigma-$algebra 
\[
\mathscr{F}^{s}_{n}=\sigma\left(X_{j}:j>n\right)
\]
is called the \salg of \textbf{strict future\mindex{sigma-algebra @ \salg ! strict future}}
of process $X$ after time $n.$

For every $n\in\mathbb{N},$ the \salg $\mathscr{P}_{n}=\sigma\left(X_{n}\right)$
is called the \salg of the\textbf{ \mindex{sigma-algebra @ \salg ! present}present}
of the process $X$ at time $n.$

The process $X$ is called a \textbf{\index{Markov chain}Markov chain}
with respect to the filtration $\left(\mathscr{A}_{n}\right)_{n\in\mathbb{N}}$
if it satisfies the following two conditions:

(i) The process $X=\left(X_{n}\right)_{n\in\mathbb{N}}$ is adapted
to the filtration $\left(\mathscr{A}_{n}\right)_{n\in\mathbb{N}}.$

(ii) The process $X$ has the \textbf{Markov property}\index{Markov property},
that is, for every $n\in\mathbb{N},$ the \salg of the \textbf{past\mindex{sigma-algebra @ \salg ! past}}
$\mathscr{A}_{n}$ and the \salg $\mathscr{F}^{l}_{n}$ of the \textbf{wide
future} at time $n$ are conditionally independent with respect to
the \salg $\mathscr{P}_{n}=\sigma\left(X_{n}\right)$ of the \textbf{present}
at time $n.$

We also say, in this case, that $X$ is a Markov chain on the underlying
filtered probabilized space $\left(\Omega,\mathscr{A},\left(\mathscr{A}_{n}\right)_{n\in\mathbb{N}},P\right).$

\end{definition}

\footnotetext{In this context, we simply say ``taking values in
$E$''.}

\begin{remark}{}{}

If $X$ is a Markov chain with respect to the filtration $\left(\mathscr{A}_{n}\right)_{n\in\mathbb{N}},$
then it is also a Markov chain with respect to its natural filtration
$\left(\mathscr{B}_{n}\right)_{n\in\mathbb{N}}.$ When the reference
filtration is the natural filtration, we simply speak of a Markov
chain, without specifying the filtration.

\end{remark}

By the previous example, in the definition of a Markov chain, we may
replace the \salg of the wide future by that of the strict future.

\begin{example}{}{}

Let $X=\left(X_{n}\right)_{n\in\mathbb{N}}$ be a sequence of independent
random variables taking values in $\mathbb{Z}.$ For every $n\in\mathbb{N},$
set
\[
S_{n}=\sum^{n}_{j=0}X_{j}.
\]
The natural filtrations of the processes $X$ and $S=\left(S_{n}\right)_{n\in\mathbb{N}}$
are the same, and $S$ is a Markov chain, called \textbf{random walk
on }$\mathbb{Z}.$

\end{example}

We now give a characterization of Markov chains.

\begin{theorem}{Markov Chains Characterization}{mark_ch_charact}

Let $X=\left(X_{n}\right)_{n\in\mathbb{N}}$ be a process taking values
in $E,$ adapted to the filtration $\left(\mathscr{A}_{n}\right)_{n\in\mathbb{N}}.$
The following properties are equivalent:

(i) The process $X$ is a Markov chain with respect to the filtration
$\left(\mathscr{A}_{n}\right)_{n\in\mathbb{N}}.$

(ii) For every $n\in\mathbb{N},$ and for every $Y\in b\mathscr{F}^{l}_{n},$
\begin{equation}
\mathbb{E}^{\mathscr{A}_{n}}\left(Y\right)=\mathbb{E}^{\mathscr{P}_{n}}\left(Y\right).\label{eq:cond_exp_past_futur_event}
\end{equation}

(iii) For every $n\in\mathbb{N},$ and for every $f\in b\mathscr{E}$---bounded
function on $E$---,
\begin{equation}
\mathbb{E}^{\mathscr{A}_{n}}\left(f\left(X_{n+1}\right)\right)=\mathbb{E}^{\mathscr{P}_{n}}\left(f\left(X_{n+1}\right)\right).\label{eq:cond_exp_past_f_next_event}
\end{equation}

(iv) For every $n\in\mathbb{N}$ and $m\in\mathbb{N},$ such that
$n\leqslant m,$ and for every $f\in b\mathscr{E},$
\begin{equation}
\mathbb{E}^{\mathscr{A}_{n}}\left(f\left(X_{m}\right)\right)=\mathbb{E}^{\mathscr{P}_{n}}\left(f\left(X_{m}\right)\right).\label{eq:cond_exp_past_f_coming_events}
\end{equation}
In particular, the process $X$ is a Markov chain with respect to
its natural filtration $\left(\mathscr{B}_{n}\right)_{n\in\mathbb{N}}$
if and only if, for every $n\in\mathbb{N},$ for every finite nondecreasing
sequence of integers such that $n_{1}\leqslant n_{2}\leqslant\cdots\leqslant n_{k}\leqslant n$
and for every bounded function $f$ on $E,$
\begin{equation}
\mathbb{E}^{\sigma\left(X_{n_{1}},\cdots,X_{n_{k}}\right)}\left(f\left(X_{n}\right)\right)=\mathbb{E}^{\sigma\left(X_{n_{k}}\right)}\left(f\left(X_{n}\right)\right).\label{eq:the_past_does_not_give_info}
\end{equation}

\end{theorem}

\begin{proof}{}{}

The equivalence between (i) and (ii) follows from Theorem $\ref{th:gen_salg_two_salgs_ind}.$
The implications $\text{(iv)}\Rightarrow\text{(iii)},$ $\text{(ii)}\Rightarrow\text{(iii)},$
and $\text{(ii)}\Rightarrow\text{(iv)}$ are immediate.

It remains to prove the implication $\text{(iii)}\Rightarrow\text{(ii)}.$
Suppose that (iii) holds. By the previous remark, it suffices to prove
that the equality $\refpar{eq:cond_exp_past_futur_event}$ holds for
every $Y\in b\mathscr{F}^{s}_{n}.$ Fix $n\in\mathbb{N},$ and consider
\[
\mathscr{H}=\left\{ Y\in b\mathscr{F}^{s}_{n}:\,\mathbb{E}^{\mathscr{A}_{n}}\left(Y\right)=\mathbb{E}^{\mathscr{P}_{n}}\left(Y\right)\right\} .
\]

The set $\mathscr{H}$ is a vector space containing the constants
and stable under bounded monotone limits. By the extension by measurability
theorem, it is therefore enough to prove that $\mathscr{H}$ contains
the indicator functions of the elements of the $\pi-$system $\mathscr{C}$
generating the \salg $\mathscr{F}^{s}_{n},$ namely
\[
\mathscr{C}=\left\{ \bigcap^{k}_{i=1}X^{-1}_{n+i}\left(E_{i}\right):\,k\in\mathbb{N}^{\ast},\,E_{i}\in\mathscr{E}\right\} .
\]

Let $Y=\boldsymbol{1}_{\bigcap^{k}_{i=1}X^{-1}_{n+i}\left(E_{i}\right)}.$ 

If $k=1,$ the result is immediate. 

Assume now that $k\geqslant2.$ We apply the classical backward in
time induction method used in the study of Markov chains. Since $\mathscr{A}_{n}\subset\mathscr{A}_{n+k-1}$
and that $\prod^{k-1}_{i=1}\boldsymbol{1}_{E_{i}}\circ X_{n+i}$ is
$\mathscr{A}_{n+k-1}-$measurable,
\[
\mathbb{E}^{\mathscr{A}_{n}}\left(Y\right)=\mathbb{E}^{\mathscr{A}_{n}}\left(\prod^{k-1}_{i=1}\boldsymbol{1}_{E_{i}}\circ X_{n+i}\right)=\mathbb{E}^{\mathscr{A}_{n}}\left(\left[\prod^{k-1}_{i=1}\boldsymbol{1}_{E_{i}}\circ X_{n+i}\right]\mathbb{E}^{\mathscr{A}_{n+k-1}}\left(\boldsymbol{1}_{E_{k}}\circ X_{n+k}\right)\right),
\]
and thus, by the hypothesis (iii),
\[
\mathbb{E}^{\mathscr{A}_{n}}\left(Y\right)=\mathbb{E}^{\mathscr{A}_{n}}\left(\left[\prod^{k-1}_{i=1}\boldsymbol{1}_{E_{i}}\circ X_{n+i}\right]\mathbb{E}^{\mathscr{P}_{n+k-1}}\left(\boldsymbol{1}_{E_{k}}\circ X_{n+k}\right)\right).
\]
However $\mathbb{E}^{\mathscr{P}_{n+k-1}}\left(\boldsymbol{1}_{E_{k}}\circ X_{n+k}\right)$
is $\mathscr{P}_{n+k-1}-$measurable and bounded; hence there exists
a Borel real-valued function $g_{1}$---trivially measurable---on
$E$ such that
\[
\mathbb{E}^{\mathscr{P}_{n+k-1}}\left(\boldsymbol{1}_{E_{k}}\circ X_{n+k}\right)=g_{1}\left(X_{n+k-1}\right).
\]
Therefore,
\[
\mathbb{E}^{\mathscr{A}_{n}}\left(Y\right)=\mathbb{E}^{\mathscr{A}_{n}}\left(\left[\prod^{k-2}_{i=1}\boldsymbol{1}_{E_{i}}\circ X_{n+i}\right]g_{1}\left(X_{n+k-1}\right)\boldsymbol{1}_{E_{k-1}}\circ X_{n+k-1}\right).
\]
By conditionning successively with respect to the \salgs $\mathscr{A}_{n+k-2},\cdots,\mathscr{A}_{n+1}$
and applying the hypothesis $\refpar{eq:cond_exp_past_f_next_event}$
at each step, and by arguing the same measurability arguments, we
obtain inductively the existence of bounded real-valued functions
$g_{2},\cdots,g_{k-1}$ on $E$ such that
\begin{align*}
\mathbb{E}^{\mathscr{A}_{n}}\left(Y\right) & =\mathbb{E}^{\mathscr{A}_{n}}\left(\left[\prod^{k-2}_{i=1}\boldsymbol{1}_{E_{i}}\circ X_{n+i}\right]\mathbb{E}^{\mathscr{A}_{n+k-2}}\left(g_{1}\left(X_{n+k-1}\right)\boldsymbol{1}_{E_{k-1}}\circ X_{n+k-1}\right)\right)\\
 & =\mathbb{E}^{\mathscr{A}_{n}}\left(\left[\prod^{k-2}_{i=1}\boldsymbol{1}_{E_{i}}\circ X_{n+i}\right]\mathbb{E}^{\mathscr{P}_{n+k-2}}\left(g_{1}\left(X_{n+k-1}\right)\boldsymbol{1}_{E_{k-1}}\circ X_{n+k-1}\right)\right)\\
 & =\mathbb{E}^{\mathscr{A}_{n}}\left(\left[\prod^{k-2}_{i=1}\boldsymbol{1}_{E_{i}}\circ X_{n+i}\right]g_{2}\left(X_{n+k-2}\right)\right)\\
 & =\cdots\\
 & =\mathbb{E}^{\mathscr{A}_{n}}\left(g_{k-1}\left(X_{n+1}\right)\right)\\
 & =\mathbb{E}^{\mathscr{P}_{n}}\left(g_{k-1}\left(X_{n+1}\right)\right).
\end{align*}
Hence $\mathbb{E}^{\mathscr{A}_{n}}\left(Y\right)$ is $\mathscr{P}_{n}-$measurable.
Since $\mathscr{P}_{n}\subset\mathscr{A}_{n},$ it follows that $\mathbb{E}^{\mathscr{A}_{n}}\left(Y\right)=\mathbb{E}^{\mathscr{P}_{n}}\left(Y\right),$
which proves (ii).

In particular, if the process $X$ is a Markov chain with respect
to its natural filtration $\left(\mathscr{B}_{n}\right)_{n\in\mathbb{N}},$
then, for every $n\in\mathbb{N},$ for every finite nondecreasing
sequence of integers such that $n_{1}\leqslant n_{2}\leqslant\cdots\leqslant n_{k}\leqslant n$
and for every bounded function $f$ on $E,$
\[
\sigma\left(X_{n_{k}}\right)\subset\sigma\left(X_{1},\cdots,X_{n_{k}}\right)\subset\mathscr{B}_{n_{k}},
\]
and therefore
\[
\mathbb{E}^{\sigma\left(X_{1},\cdots,X_{n_{k}}\right)}\left(f\left(X_{n}\right)\right)=\mathbb{E}^{\sigma\left(X_{1},\cdots,X_{n_{k}}\right)}\left(\mathbb{E}^{\mathscr{B}_{n_{k}}}\left(f\left(X_{n}\right)\right)\right)=\mathbb{E}^{\sigma\left(X_{n_{k}}\right)}\left(f\left(X_{n}\right)\right).
\]
This proves the equality $\refpar{eq:the_past_does_not_give_info}.$
Conversely, if this last property holds, then by choosing the sequence
of consecutive integers up to $n,$ we obtain that the equality $\refpar{eq:cond_exp_past_f_next_event}$
is satisfied, and thus that $X$ is a Markov chain.

\end{proof}

We now present a common case arising in applications---which can,
of course, be generalized in several directions.

\begin{example}{Auto-Regressive Process}{auto-reg_proc}

Let $E$ be a countable set, equipped with the \salg of its subsets,
and let $g$ be a measurable mapping from $E\times\mathbb{R}$ to
$E.$ Consider a family of independent random variables defined on
the same probabilized space $\left(\Omega,\mathscr{A},P\right).$
One of these random variable, $X_{0},$ takes values in $E,$ the
others form a sequence of real-valued random variables $\left(U_{n}\right)_{n\in\mathbb{N}^{\ast}}$
with common law $\mu.$ In particular, if $\mu$ is the uniform law
on the interval $\left[0,1\right],$ the variables $U_{n}$ can model
random number draws produced by a \textbf{\mindex{generator!random number}\index{random number generator}random
number} \textbf{generator} during a \textbf{simulation}. We construct
the sequence $\left(X_{n}\right)_{n\in\mathbb{N}}$ by setting
\[
X_{n+1}=g\left(X_{n},U_{n+1}\right).
\]
The process $X=\left(X_{n}\right)_{n\in\mathbb{N}}$ is called \textbf{\index{auto-regressive}auto-regressive}:
the randomness injected at each step is explicit.

Let $\mathscr{A}_{0}$ be the \salg $\sigma\left(X_{0}\right)$ and,
for $n\geqslant1,$ let $\mathscr{A}_{n}$ be the \salg $\sigma\left(X_{0},U_{1},\cdots,U_{n}\right),$
which naturally summarizes the information from the past up to time
$n.$ 

Show that the process $X$ is a---homogeneous\footnotemark---Markov
chain with respect to the filtration $\left(\mathscr{A}_{n}\right)_{n\in\mathbb{N}}.$

\end{example}

\footnotetext{The definition of a homogeneous Markov chain is given
later.}

\begin{solutionexample}{}{}

By an easy induction, we see that $X_{n}$ is $\mathscr{A}_{n}-$measurable.
Hence $X$ is adapted to the filtration $\left(\mathscr{A}_{n}\right)_{n\in\mathbb{N}}.$ 

Moreover, for every $f\in bE,$
\[
\mathbb{E}^{\mathscr{A}_{n}}\left(f\left(X_{n+1}\right)\right)=\mathbb{E}^{\mathscr{A}_{n}}\left(f\left(g\left(X_{n},U_{n+1}\right)\right)\right).
\]
Since $U_{n+1}$ and $\mathscr{A}_{n}$ are independent, Proposition
$\ref{pr:cond_exp_exist_indep}$ yields
\begin{equation}
\mathbb{E}^{\mathscr{A}_{n}}\left(f\left(X_{n+1}\right)\right)=\widehat{f}\left(X_{n}\right)\,\,\,\,P-\text{almost surely,}\label{eq:cond_exp_f_X_nplus1_kn_A_n}
\end{equation}
where the function $\widehat{f}$ is defined, for every $x\in E,$
by
\begin{equation}
\widehat{f}\left(x\right)=\mathbb{E}\left(f\left(x,U_{n+1}\right)\right)=\intop_{\mathbb{R}}f\left(x,u\right)\text{d}\mu\left(u\right).\label{eq:def_f_transform_via_exp}
\end{equation}
The equality $\refpar{eq:cond_exp_f_X_nplus1_kn_A_n}$ implies the
equality $\refpar{eq:cond_exp_past_f_next_event},$ namely
\[
\mathbb{E}^{\mathscr{A}_{n}}\left(f\left(X_{n+1}\right)\right)=\mathbb{E}^{\mathscr{P}_{n}}\left(f\left(X_{n+1}\right)\right),
\]
which proves that $X$ is a Markov chain with respect to the filtration
$\left(\mathscr{A}_{n}\right)_{n\in\mathbb{N}},$ and therefore also
a Markov chain with respect to the natural filtration.

\end{solutionexample}

\begin{remark}{}{}

Since $E$ is countable, equality $\refpar{eq:the_past_does_not_give_info}$
is equivalent to the equality of conditional probabilities: for $P_{\left(X_{n_{1}},\cdots,X_{n_{k}}\right)}-$almost
every $\left(x_{n_{1}},\cdots,x_{n_{k}}\right)\in E^{k}$ and every
$x_{n}\in E,$ 
\begin{equation}
P^{\left(X_{n_{1}}=x_{n_{1}},\cdots,X_{n_{k}}=x_{n_{k}}\right)}\left(X_{n}=x_{n}\right)=P^{\left(X_{n_{k}}=x_{n_{k}}\right)}\left(X_{n}=x_{n}\right).\label{eq:cond_prob_eq}
\end{equation}
Hence $X$ is a Markov chain---with respect to the natural filtration---if
and only if, for every $n,k\in\mathbb{N}$ and every finite nondecreasing
sequence of integers such that $n_{1}\leqslant n_{2}\leqslant\cdots\leqslant n_{k}\leqslant n,$
the equality $\refpar{eq:cond_prob_eq}$ holds for $P_{\left(X_{n_{1}},\cdots,X_{n_{k}}\right)}-$almost
every $\left(x_{n_{1}},\cdots,x_{n_{k}}\right)\in E^{k}$ and for
$x_{n}\in E.$

\end{remark}

For integers $n$ and $m$ such that $n\leqslant m,$ and for every
points $x,y$ of the state space $E$ for which it is defined, the
conditional probability
\[
M_{n,m}\left(x,y\right)=P^{\left(X_{n}=x\right)}\left(X_{m}=y\right)
\]
is called the \textbf{\index{transition probability}transition probability}
from state $x$ at time $n$ to state $y$ at time $m.$ We then have
\[
\sum_{y\in E}P^{\left(X_{n}=x\right)}\left(X_{m}=y\right)=1.
\]
In order to handle Markov chains efficiently, we introduce the following
definitions. 

\begin{definition}{Transition Matrix. Stochastic Matrix}{}

A family $A$ of nonnegative real numbers, indexed by $E\times E$
is a \textbf{transition matrix}\index{transition matrix}\footnotemark---or
a \textbf{stochastic matrix}\index{stochastic matrix}---if, for
every $x\in E,$ 
\[
\sum_{y\in E}A\left(x,y\right)=1.
\]

\end{definition}

\footnotetext{If $E$ is finite, this is an ordinary square matrix.
If $E$ is infinite, it is a generalized matrix.

}

\begin{denotation}{}{denot_bounded_f}

We denote by $bE$ the set of bounded functions on $E.$ If $A$ is
a matrix indexed on $E\times E$ with nonnegative entries, the mapping
$A\left(x,\cdot\right)$ defines a measure $\mu_{x}$ on $E.$ If
$f$ is nonnegative or $\mu_{x}-$integrable, we denote by $A\left(x,f\right),$
or equivalently $Af\left(x\right)$ its integral with respect to $\mu_{x},$
namely\boxeq{
\begin{equation}
A\left(x,f\right)=\sum_{y\in E}A\left(x,y\right)f\left(y\right).\label{eq:Ax-f}
\end{equation}
}If moreover $A$ is such that each measure $\mu_{x}$ has mass less
than or equal to 1, and if $f$ is bounded, then the function $A\left(\cdot,f\right)$
is also bounded. This functional point of view extends naturally to
arbitrary state spaces. However, since $E$ is countable here, it
is often convenient for explicit computation---especially when $E$
is finite---to adopt a vectorial viewpoint. We identify a function
$f$ with the ``column vector'' $f\left(y\right)_{y\in E}.$ The
``column vector'' $A\left(\cdot,f\right)$ then has vector components
$A\left(x,f\right),$ as given by the equality $\refpar{eq:Ax-f}.$
Hence, with these identifications, we simply write
\[
A\left(\cdot,f\right)=Af.
\]

\end{denotation}

\begin{definition}{Family of Transition Matrices}{fam_tr_mat}

A Markov chain $X,$ adapted to the filtration $\left(\mathscr{A}_{n}\right)_{n\in\mathbb{N}},$
is said to admit a \textbf{family of transition matrices\mindex{family!transition matrices}}
$\left(M_{n,m}\right)_{n\leqslant m}$ if, for every $n$ and $m$
such that $n<m,$ $M_{n,m}$ is a transition matrix and for every
bounded nonnegative function $f$ on $E,$ we have $P-$almost surely\boxeq{
\begin{equation}
\mathbb{E}^{\mathscr{A}_{n}}\left(f\left(X_{m}\right)\right)=M_{n,m}\left(X_{n},f\right).\label{eq:esp_cond_f_x_m_kn_A_n_as_trans_mat}
\end{equation}
}

\end{definition}

\begin{remark}{}{}

For $P_{X_{n}}-$almost every $x\in E$ and every $y\in E,$ we have
\[
M_{n,m}\left(x,y\right)=P^{\left(X_{n}=x\right)}\left(X_{m}=y\right).
\]
Moreover, if $0\leqslant n_{1}\leqslant n_{2}\leqslant\cdots\leqslant n_{k}=n\leqslant m,$
then for $P_{\left(X_{0},X_{n_{1}},\cdots,X_{n}\right)}-$almost every
$\left(x_{0},x_{n_{1}},\cdots,x_{n}\right)\in E^{k+1}$ and every
$x_{n+1}\in E,$ 
\begin{align*}
P^{\left(X_{0}=x_{0},X_{n_{1}}=x_{n_{1}},\cdots,X_{n}=x_{n}\right)}\left(X_{n+1}=x_{n+1}\right) & =P^{\left(X_{n}=x_{n}\right)}\left(X_{n+1}=x_{n+1}\right)\\
 & =M_{n,n+1}\left(x_{n},x_{n+1}\right).
\end{align*}
Indeed, since $\sigma\left(X_{n}\right)\subset\sigma\left(X_{0},X_{n_{1}},\cdots,X_{n}\right)\subset\mathscr{A}_{n},$
taking for $f$ the function $\boldsymbol{1}_{\left\{ y\right\} },$
we obtain
\begin{align*}
P^{\sigma\left(X_{0},X_{n_{1}},\cdots,X_{n}\right)\left(X_{m}=y\right)} & =\mathbb{E}^{\sigma\left(X_{0},X_{n_{1}},\cdots,X_{n}\right)}\left(\boldsymbol{1}_{\left(X_{m}=y\right)}\right)\\
 & =\mathbb{E}^{\sigma\left(X_{0},X_{n_{1}},\cdots,X_{n}\right)}\left(\mathbb{E}^{\mathscr{A}_{n}}\left(\boldsymbol{1}_{\left(X_{m}=y\right)}\right)\right)\\
 & =M_{n,m}\left(X_{n},y\right).
\end{align*}
It follows that
\[
P^{\sigma\left(X_{n}\right)}\left(X_{m}=y\right)=\mathbb{E}^{\sigma\left(X_{n}\right)}\left(\mathbb{E}^{\sigma\left(X_{0},X_{n_{1}},\cdots,X_{n}\right)}\left(\boldsymbol{1}_{\left(X_{m}=y\right)}\right)\right)=M_{n,m}\left(X_{n},y\right).
\]
\end{remark}

\begin{remark}{}{}

A transition matrix $A$ on $E$ is such that, for every $x\in E,$
the mapping $A\left(x,\cdot\right)$ defines a germ of probability
on the countable set $E.$ Identifying this germ with the generated
probability, one may view $A$ as a \textbf{probability}---or \textbf{kernel}---\textbf{of}
\textbf{transition\index{kernel of transition}\mindex{probability!of transition}}
on $E:$ it is a regular version of the conditional law of $X_{m}$
given $\left(X_{0},X_{n_{1}},\cdots,X_{n}\right).$ This is the viewpoint
that extends naturally to Markov chains with a general state space.

\end{remark}

\begin{figure}[t]
\begin{center}\includegraphics[width=0.4\textwidth]{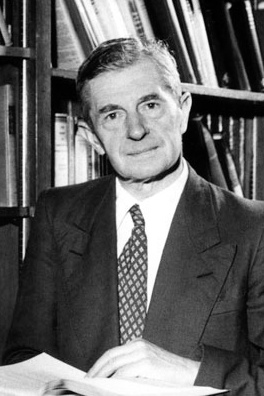}

{\tiny NOAA---\href{http://celebrating200years.noaa.gov/magazine/igy/sydneychapman1_500.html}{NOAA Magazine}---Public
Domain}\end{center}

\caption{\textbf{\protect\href{https://en.wikipedia.org/wiki/Sydney_Chapman_(mathematician)}{Sydney Chapman}}
(1888 - 1970)}\sindex[fam]{Chapman, Sydney}
\end{figure}

\begin{proposition}{Chapman-Kolmogorov Equality}{chapman_kolm_eq}

Let $X$ be a Markov chain adapted to the filtration $\left(\mathscr{A}_{n}\right)_{n\in\mathbb{N}},$
with transition matrix family $\left(M_{n,m}\right)_{n\leqslant m}.$
For every instants $n,r,m$ such that $n<r<m,$ the \textbf{Chapman\footnotemark\sindex[fam]{Chapman, Sydney}-Kolmogorov
relation\index{Chapman-Kolmogorov relation}} holds:\boxeq{
\[
\forall y\in E,\,\,\,\,M_{n,m}\left(X_{n},y\right)=\sum_{z\in E}M_{n,r}\left(X_{n},z\right)M_{r,m}\left(z,y\right).
\]
}Equivalently, in matrix form,\boxeq{
\begin{equation}
M_{n,m}\left(X_{n},\cdot\right)=M_{n,r}M_{r,m}\left(X_{n},\cdot\right).\label{eq:matrix_form_Chapman_Kolmogorov}
\end{equation}
} 

\end{proposition}

\footnotetext{\textbf{\href{https://en.wikipedia.org/wiki/Sydney_Chapman_(mathematician)}{Sydney Chapman}\sindex[fam]{Chapman, Sydney}}
(1888 - 1970) was a British mathematician and geophysicist. He contributed
to kinetic theory of gazes, solar-terrestrial physics and worked on
the Earth's ozone layer.}

\begin{proof}{}{}

The events $\left(X_{r}=z\right),\,z\in E,$ form a complete system
of constituents. Taking for $f$ the function $\boldsymbol{1}_{\left\{ y\right\} }$
in the inequality $\refpar{eq:esp_cond_f_x_m_kn_A_n_as_trans_mat},$
we get
\begin{align*}
M_{n,m}\left(X_{n},y\right) & =\mathbb{E}^{\mathscr{A}_{n}}\left(\boldsymbol{1}_{\left(X_{m}=y\right)}\right)=\mathbb{E}^{\mathscr{A}_{n}}\left(\sum_{z\in E}\boldsymbol{1}_{\left(X_{r}=z\right)}\boldsymbol{1}_{\left(X_{m}=y\right)}\right)\\
 & =\sum_{z\in E}\mathbb{E}^{\mathscr{A}_{n}}\left(\boldsymbol{1}_{\left(X_{r}=z\right)}\boldsymbol{1}_{\left(X_{m}=y\right)}\right).
\end{align*}
Hence, since $\left(X_{r}=z\right)\in\mathscr{A}_{r},$
\begin{align*}
M_{n,m}\left(X_{n},y\right) & =\sum_{z\in E}\mathbb{E}^{\mathscr{A}_{n}}\left(\boldsymbol{1}_{\left(X_{r}=z\right)}\mathbb{E}^{\mathscr{A}_{r}}\left(\boldsymbol{1}_{\left(X_{m}=y\right)}\right)\right)\\
 & =\sum_{z\in E}\mathbb{E}^{\mathscr{A}_{n}}\left(\boldsymbol{1}_{\left(X_{r}=z\right)}M_{r,m}\left(X_{r},y\right)\right)\\
 & =\sum_{z\in E}\mathbb{E}^{\mathscr{A}_{n}}\left(\boldsymbol{1}_{\left(X_{r}=z\right)}M_{r,m}\left(z,y\right)\right)\\
 & =\sum_{z\in E}\mathbb{E}^{\mathscr{A}_{n}}\left(\boldsymbol{1}_{\left(X_{r}=z\right)}\right)M_{r,m}\left(z,y\right)\\
 & =\sum_{z\in E}M_{n,r}\left(X_{n},z\right)M_{r,m}\left(z,y\right).
\end{align*}

\end{proof}

\begin{remark}{}{}

If $X$ is merely a Markov chain with respect to its natural filtration,
one may give the following \textbf{heuristic proof}\footnotemark
of the Chapman-Kolmogorov equality. By the formula of total probability,
and the Markov property,
\[
P\left(X_{m}=y\mid X_{n}=x\right)=\sum_{z\in E}P\left(X_{m}=y\mid X_{r}=z\right)P\left(X_{r}=z\mid X_{n}=x\right).
\]

\end{remark}

\footnotetext{In the sense that we ignore the possible issue of division
by 0.}

A physical system whose state depends on time is said to be \textbf{\index{conservative}conservative}
if, for every $t,$ its law of transition from a state $x$ at time
$s$ to a state $y$ at time $s+t$ does not depend on $s.$ When
such a system is modelled by a Markov chain, this translates into
temporal homogeneity of the conditional laws, and therefore of the
family of transition matrices.

Thus, a Markov chain $X$ with transition matrix family $\left(M_{n,m}\right)_{n\leqslant m}$
is \textbf{temporally homogeneous} if there exists a sequence of transition
matrices $\left(M_{\left(n\right)}\right)_{n\in\mathbb{N}}$ indexed
by $E$ such that, for every integers $n$ and $m,$ 
\[
M_{n,n+m}=M_{\left(m\right)}.
\]
It then follows from $\refpar{eq:matrix_form_Chapman_Kolmogorov}$,
that for every integers $n$ and $m,$
\[
M_{n,n+m}\left(X,\cdot\right)=M^{m}_{\left(1\right)}\left(X_{n},\cdot\right).
\]
In particular, the matrix $M_{\left(1\right)}$ is denoted $M$ and
is called the \textbf{\index{transition matrix}transition matrix}
of the homogeneous Markov chain. Then,
\[
\mathbb{E}^{\mathscr{A}_{n}}\left(f\left(X_{n+m}\right)\right)=M^{m}\left(X_{n},f\right),
\]
where $M^{m}$ denotes the $m-$th power of the matrix $M.$

In particular, for $P_{X_{0}}-$almost every $x\in E$ and every $y\in E,$
\[
M\left(x,y\right)=P^{\left(X_{0}=x\right)}\left(X_{1}=y\right),
\]
and for $P_{X_{n}}-$almost every $x\in E$ and every $y\in E,$\boxeq{
\[
M\left(x,y\right)=P^{\left(X_{n}=x\right)}\left(X_{n+1}=y\right).
\]
}

We are thus led to the following definition.

\begin{definition}{}{}

A \textbf{Markov chain} $X$---with respect to the filtration $\left(\mathscr{A}_{n}\right)_{n\in\mathbb{N}}$---taking
values in $E,$ is \textbf{homogeneous\mindex{Markov chain!homogeneous}}
with \textbf{transition matrix} $M$ if, for every integers $n$ and
$m$ such that $0\leqslant n<m,$ and for every $f\in bE,$ \boxeq{
\begin{equation}
\mathbb{E}^{\mathscr{A}_{n}}\left(f\left(X_{m}\right)\right)=M^{m-n}\left(X_{n},f\right).\label{eq:hom_mark_chain_def}
\end{equation}
}

\end{definition}

The following proposition allows to \textbf{prove that a process is
a homogeneous Markov chain with transition matrix $M.$}

\begin{proposition}{}{}

The process $X,$ adapted to the filtration $\left(\mathscr{A}_{n}\right)_{n\in\mathbb{N}},$
is a homogeneous Markov chain with transition matrix $M$ if and only
if, for every integer $n\in\mathbb{N}$ and every function $f\in bE,$\boxeq{
\begin{equation}
\mathbb{E}^{\mathscr{A}_{n}}\left(f\left(X_{n+1}\right)\right)=M\left(X_{n},f\right).\label{eq:cns_hom_mark_chain_and_cond_exp}
\end{equation}
}

\end{proposition}

\begin{proof}{}{}

The necessity is immediate. Conversely, suppose that the relation
$\refpar{eq:cns_hom_mark_chain_and_cond_exp}$ holds for every integer
$n\in\mathbb{N}.$ Let $n,m\in\mathbb{N}$ such that $0\leqslant n<m.$
It holds
\begin{align*}
\mathbb{E}^{\mathscr{A}_{n}}\left(f\left(X_{m}\right)\right) & =\mathbb{E}^{\mathscr{A}_{n}}\left(\mathbb{E}^{\mathscr{A}_{m-1}}\left[f\left(X_{m}\right)\right]\right)\\
 & =\mathbb{E}^{\mathscr{A}_{n}}\left(M\left(X_{m-1},f\right)\right)\\
 & =\mathbb{E}^{\mathscr{A}_{n}}\left(\mathbb{E}^{\mathscr{A}_{m-2}}\left[M\left(X_{m-1},f\right)\right]\right)\\
 & =\mathbb{E}^{\mathscr{A}_{n}}\left(M\left(X_{m-2},M\left(\cdot,f\right)\right)\right)\\
 & =\mathbb{E}^{\mathscr{A}_{n}}\left(M^{2}\left(X_{m-2},f\right)\right).
\end{align*}
A straighforward induction then yields the relation $\refpar{eq:hom_mark_chain_def}.$

\end{proof}

Hence, the \textbf{auto-regressive\mindex{process!auto-regressive}}
process introduced in Example $\ref{ex:auto-reg_proc}$ is a \textbf{homogeneous
Markov chain} with transition matrix $M,$ determined, for every bounded
function $f,$ by\boxeq{
\[
M\left(x,f\right)=\widetilde{f}\left(x\right)=\intop_{\mathbb{R}}f\left(x,u\right)\text{d}\mu\left(u\right).
\]
}

We now give an example of a homogeneous Markov chain with respect
to a \textbf{filtration that is not its natural filtration}.

\begin{example}{Conditional Random Walk}{}

Consider a family of independent real-valued random variables defined
on the same probabilized space \preds. One of these random variables,
$\Theta,$ takes values in the interval $\left[0,1\right]$ and has
law $\mu.$ The others form a sequence of real-valued random variables
$\left(U_{n}\right)_{n\in\mathbb{N}}$ with uniform law on the interval
$\left[0,1\right],$ modelling random number draws produced by a random
number generator. 

Construct the sequence $\left(X_{n}\right)_{n\in\mathbb{N}}$ by setting
\[
X_{n}=\boldsymbol{1}_{\left(U_{n}\leqslant\Theta\right)}-\boldsymbol{1}_{\left(U_{n}>\Theta\right)},
\]
that is, conditionally on the value $\theta$ of $\Theta$ drawn beforehand,
the random variable $X_{n}$ follows the law $\theta\delta_{1}+\left(1-\theta\right)\delta_{-1}.$ 

Define, for every $n\in\mathbb{N},$ the random variables
\[
S_{n}=\sum^{n}_{j=0}X_{j}\,\,\,\,\text{and}\,\,\,\,Y_{n}=\left(\Theta,S_{n}\right),
\]
and denote $\mathscr{A}_{n}=\sigma\left(\Theta,U_{0},U_{1},\cdots,U_{n}\right)$
the \salg generated by $\Theta,\,U_{0},\,U_{1},\cdots,\,U_{n}.$
The \salg $\mathscr{A}_{n}$ naturally summarizes the information
on the past up to time $n.$ 

To remain in the framework of chains with countable state space, assume
that the law $\mu$ has its support in a countable subset $E$ of
$\left[0,1\right].$ For a given value $\theta,$ the sequence $\left(S_{n}\right)_{n\in\mathbb{N}}$
behaves like a random walk on $\mathbb{Z},$ and the value of $\theta$
must remain constantly in memory in order to follow the walk. Hence,
the process $\left(S_{n}\right)_{n\in\mathbb{N}}$ has a history that
depends permanently on the initial time. 

Nonetheless, we will show that the process $Y=\left(Y_{n}\right)_{n\in\mathbb{N}}$
is a homogeneous Markov chain taking values in $E\times\mathbb{Z},$
adapted to the filtration $\left(\mathscr{A}_{n}\right)_{n\in\mathbb{N}}.$

\end{example}

\begin{solutionexample}{}{}

For this, we compute, for every bounded function $f$ on $E\times\mathbb{Z},$
\[
\mathbb{E}^{\mathscr{A}_{n}}\left(f\left(Y_{n+1}\right)\right)=\mathbb{E}^{\mathscr{A}_{n}}\left(f\left(\Theta,S_{n}+X_{n+1}\right)\right).
\]

Define the function $h$ on $\left[0,1\right]^{2}$ by
\[
h\left(\theta,u\right)=\boldsymbol{1}_{\left(u\leqslant\theta\right)}-\boldsymbol{1}_{\left(u>\theta\right)}.
\]
For $P_{\left(\Theta,U_{0},U_{1},\cdots,U_{n}\right)}-$almost every
$\left(\theta,u_{0},u_{1},\cdots,u_{n}\right),$
\begin{multline*}
\mathbb{E}^{\sigma\left(\Theta,U_{0},U_{1},\cdots,U_{n}\right)=\left(\theta,u_{0},u_{1},\cdots,u_{n}\right)}\left[f\left(Y_{n+1}\right)\right]\\
=\mathbb{E}^{\sigma\left(\Theta,U_{0},U_{1},\cdots,U_{n}\right)=\left(\theta,u_{0},u_{1},\cdots,u_{n}\right)}\left[f\left(\theta,\sum^{n}_{j=0}h\left(\theta,u_{j}\right)+h\left(\theta,U_{n+1}\right)\right)\right].
\end{multline*}
By independence of the random variables $\Theta,U_{0},U_{1},\cdots,U_{n},$
\[
\mathbb{E}^{\sigma\left(\Theta,U_{0},U_{1},\cdots,U_{n}\right)=\left(\theta,u_{0},u_{1},\cdots,u_{n}\right)}\left[f\left(Y_{n+1}\right)\right]=\mathbb{E}\left(f\left(\theta,\sum^{n}_{j=0}h\left(\theta,u_{j}\right)+h\left(\theta,U_{n+1}\right)\right)\right).
\]
Since $U_{n+1}$ has uniform law on $\left[0,1\right],$
\begin{multline*}
\mathbb{E}^{\sigma\left(\Theta,U_{0},U_{1},\cdots,U_{n}\right)=\left(\theta,u_{0},u_{1},\cdots,u_{n}\right)}\left[f\left(Y_{n+1}\right)\right]\\
=\theta f\left(\theta,\sum^{n}_{j=0}h\left(\theta,u_{j}\right)+1\right)+\left(1-\theta\right)f\left(\theta,\sum^{n}_{j=0}h\left(\theta,u_{j}\right)-1\right).
\end{multline*}
Thus, $P-$almost surely,
\[
\mathbb{E}^{\mathscr{A}_{n}}\left[f\left(Y_{n+1}\right)\right]=\Theta f\left(\Theta,\sum^{n}_{j=0}h\left(\Theta,U_{j}\right)+1\right)+\left(1-\Theta\right)f\left(\Theta,\sum^{n}_{j=0}h\left(\Theta,U_{j}\right)-1\right).
\]
Hence, $P-$almost surely,
\[
\mathbb{E}^{\mathscr{A}_{n}}\left[f\left(Y_{n+1}\right)\right]=\Theta f\left(\Theta,S_{n}+1\right)+\left(1-\Theta\right)f\left(\Theta,S_{n}-1\right).
\]

Define the transition matrix $M$ on $E\times\mathbb{Z}$ by
\[
M\left(\left(\theta,s\right),f\right)=\theta f\left(\theta,s+1\right)+\left(1-\theta\right)f\left(\theta,s-1\right).
\]
Then, $P-$almost surely,
\[
\mathbb{E}^{\mathscr{A}_{n}}\left[f\left(Y_{n+1}\right)\right]=M\left(Y_{n},f\right)
\]
which proves that the process $Y$ is a homogeneous Markov chain with
transition matrix $M$ with respect to the filtration $\left(\mathscr{A}_{n}\right)_{n\in\mathbb{N}}.$

\end{solutionexample}

Nonetheless, now we give an example of \textbf{non-homogeneous} Markov
chain.

\begin{example}{}{polya_diff_mod_2}Returning to Example $\ref{ex:polya_diff_mod}$
on \textbf{P\'olya contagious disease diffusion model}---under its
urn-drawing formulation---let us prove that the process of \textbf{proportions}
$Y_{n}$ of blue balls in the urn after the $n$-th draw and after
adding the $c$ balls is a \textbf{non-homogeneous Markov chain},
and also a \textbf{martingale}.

\end{example}

\begin{solutionexample}{}{}

Consider a sequence $\left(X_{n}\right)_{n\in\mathbb{N}^{\ast}}$
of random variables defined on a probabilized space $\left(\Omega,\mathscr{A},P\right),$
taking values in $\left\{ 0,1\right\} $---$X_{n}$ takes the value
0 when the drawn ball is red, and 1 if it is blue. Denote $k_{n}=b+r+nc.$
The number $B_{n}$ and the proportion $Y_{n}$ of blue balls in the
urn after the $n-$th draw are respectively
\[
B_{n}=b+c\sum^{n}_{j=1}X_{j}\,\,\,\,\text{and}\,\,\,\,Y_{n}=\dfrac{B_{n}}{k_{n}}.
\]

We have
\[
P\left(X_{1}=1\right)=\dfrac{b}{b+r}\,\,\,\,\text{and}\,\,\,\,P\left(X_{1}=0\right)=\dfrac{r}{b+r}.
\]
Moreover, since the draws are uniform, for every $n\geqslant2$ and
every $\left(x_{1},\cdots,x_{n}\right)\in\left\{ 0,1\right\} ^{n},$
\begin{equation}
\left\{ \begin{array}{l}
P^{\left(X_{1},\cdots,X_{n}\right)=\left(x_{1},\cdots,x_{n}\right)}\left(X_{n+1}=1\right)=\dfrac{b+c\sum^{n}_{j=1}x_{j}}{k_{n}},\\
P^{\left(X_{1},\cdots,X_{n}\right)=\left(x_{1},\cdots,x_{n}\right)}\left(X_{n+1}=0\right)=\dfrac{r+c\left(n-\sum^{n}_{j=1}x_{j}\right)}{k_{n}}.
\end{array}\right.\label{eq:P_next_knowing_current_red_and_blue}
\end{equation}
The process $Y=\left(Y_{n}\right)_{n\in\mathbb{N}^{\ast}}$ takes
values in
\[
E=\bigcup_{j\in\mathbb{N}^{\ast}}\left\{ \dfrac{j}{k_{n}}:\,b\leqslant j\leqslant k_{n}\right\} ,
\]
a countable infinite subset of the interval $\left[0,1\right].$ 

First, note that the \salgs $\mathscr{A}_{n}=\sigma\left(X_{j}:\,1\leqslant j\leqslant n\right)$
and $\sigma\left(Y_{j}:\,1\leqslant j\leqslant n\right)$ coincide.

Indeed, the mapping $F_{n}$ from $\mathbb{R}^{n}$ onto itself defined
by
\[
F_{n}\left(x_{1},\cdots,x_{n}\right)=\left(y_{1},\cdots,y_{n}\right)\,\,\,\,\text{where for}\,l\in\left\llbracket 1,n\right\rrbracket ,\,\,\,\,y_{l}=\dfrac{b+c\sum^{l}_{j=1}x_{j}}{k_{l}},
\]
is a bijection, and that $\left(Y_{1},\cdots,Y_{n}\right)=F_{n}\left(X_{1},\cdots,X_{n}\right).$

It then follows from $\refpar{eq:P_next_knowing_current_red_and_blue}$
that
\begin{equation}
P^{\left(Y_{1},\cdots,Y_{n}\right)=\left(y_{1},\cdots,y_{n}\right)}\left(X_{n+1}=1\right)=P^{\left(X_{1},\cdots,X_{n}\right)=F^{-1}_{n}\left(y_{1},\cdots,y_{n}\right)}\left(X_{n+1}=1\right)=y_{n},\label{eq:PcondY1Yniny1ynXnp1is1}
\end{equation}
which yields
\begin{equation}
P^{\sigma\left(Y_{1},\cdots,Y_{n}\right)}\left(X_{n+1}=1\right)=Y_{n}\,\,\,\,\text{and}\,\,\,\,P^{\sigma\left(Y_{1},\cdots,Y_{n}\right)}\left(X_{n+1}=0\right)=1-Y_{n}.\label{eq:cond_pr_X_nplus1_kn_gen_sig_alg_Y_1Y_n}
\end{equation}
A direct computation gives
\begin{equation}
k_{n+1}=k_{n}+c\,\,\,\,\text{and}\,\,\,\,Y_{n+1}=\dfrac{k_{n}Y_{n}+cX_{n+1}}{k_{n+1}}.\label{eq:k_nplus1_and_Ynplus1}
\end{equation}
Hence, for every function $f\in bE,$
\[
\mathbb{E}^{\mathscr{A}_{n}}\left(f\left(Y_{n+1}\right)\right)=\mathbb{E}^{\mathscr{A}_{n}}\left(\boldsymbol{1}_{\left(X_{n+1}=1\right)}f\left(\dfrac{k_{n}Y_{n}+c}{k_{n+1}}\right)\right)+\mathbb{E}^{\mathscr{A}_{n}}\left(\boldsymbol{1}_{\left(X_{n+1}=0\right)}f\left(\dfrac{k_{n}Y_{n}}{k_{n+1}}\right)\right),
\]
which yields, by using $\refpar{eq:PcondY1Yniny1ynXnp1is1},$
\[
\mathbb{E}^{\mathscr{A}_{n}}\left(f\left(Y_{n+1}\right)\right)=f\left(\dfrac{k_{n}Y_{n}+c}{k_{n+1}}\right)Y_{n}+f\left(\dfrac{k_{n}Y_{n}}{k_{n+1}}\right)\left(1-Y_{n}\right).
\]
Define, for every $y\in E$ and every $f\in bE,$
\[
M_{n}\left(y,f\right)=f\left(\dfrac{k_{n}y+c}{k_{n+1}}\right)y+f\left(\dfrac{k_{n}y}{k_{n+1}}\right)\left(1-y\right).
\]
Then, for every $n\in\mathbb{N}^{\ast},$
\begin{equation}
\mathbb{E}^{\mathscr{A}_{n}}\left(f\left(Y_{n+1}\right)\right)=M_{n}\left(Y_{n},f\right).\label{eq:cond_exp_as_trans_mat}
\end{equation}
The matrix $M_{n}$ is indeed a transition matrix on $E,$ since
\[
\sum_{z\in E}M_{n}\left(y,z\right)=\sum_{z\in E}\boldsymbol{1}_{\left(\frac{k_{n}y+c}{k_{n+1}}=z\right)}y+\sum_{z\in E}\boldsymbol{1}_{\left(\frac{k_{n}y}{k_{n+1}}=z\right)}\left(1-y\right)=y+\left(1-y\right)=1.
\]
Therefore, the process $Y$ is a \textbf{non-homogeneous} Markov chain
with family of transition matrices $\left(M_{n}\right)_{n\in\mathbb{N}^{\ast}}.$ 

Moreover, it is a \textbf{martingale}. Indeed, by taking for $f,$
in the equality $\refpar{eq:cond_exp_as_trans_mat},$ the identity
map on $E$---which is bounded---we obtain
\[
\mathbb{E}^{\mathscr{A}_{n}}\left(Y_{n+1}\right)=\dfrac{k_{n}Y_{n}+c}{k_{n+1}}Y_{n}+\dfrac{k_{n}Y_{n}}{k_{n+1}}\left(1-Y_{n}\right)=Y_{n}\dfrac{k_{n}+c}{k_{n+1}}=Y_{n}.
\]

It then follows from the martingale property that
\[
\mathbb{E}\left(Y_{n+1}\right)=\mathbb{E}\left(Y_{1}\right)=\mathbb{E}\left(\dfrac{b+cX_{1}}{b+r+c}\right)=\dfrac{1}{b+r+c}\left(b+c\dfrac{b}{b+r}\right),
\]
and therefore\boxeq{
\[
\mathbb{E}\left(Y_{n+1}\right)=\dfrac{b}{b+r}.
\]
}

The random variable $X_{n+1}$ takes only the values 0 and 1. Its
law is therefore determined by its expectation, which we compute using
the relations $\refpar{eq:cond_pr_X_nplus1_kn_gen_sig_alg_Y_1Y_n}.$
We have
\[
P\left(X_{n+1}=1\right)=\mathbb{E}\left(X_{n+1}\right)=\mathbb{E}\left(P^{\sigma\left(Y_{1},\cdots,Y_{n}\right)}\left(X_{n+1}=1\right)\right)=\mathbb{E}\left(Y_{n}\right),
\]
and thus\boxeq{
\[
P\left(X_{n+1}=1\right)=\dfrac{b}{b+r}.
\]
}

We have therefore shown that the law of $X_{n},$ for $n\geqslant1,$
is independent of $n$ and of $c,$ which is, a priori, neither obvious
nor intuitive.

The martingale $Y=\left(Y_{n}\right)_{n\in\mathbb{N}^{\ast}}$ is
bounded; hence it converges $P-$almost surely and in every $\text{L}^{p}$
to a random variable $Y_{\infty}.$ The law of $Y_{\infty}$ is the
beta law $\beta\left(\dfrac{b}{c},\dfrac{r}{c}\right)$ of the first
kind on $\left[0,1\right].$ A proof is proposed in Exercise $\ref{exo:exercise17.12}.$

\end{solutionexample}

\begin{remark}{}{}

Given a homogeneous Markov chain with transition matrix $M,$ for
every $n\in\mathbb{N},$ for $P_{X_{n}}-$almost every $x\in E$ and
for every $y\in E,$ we have\boxeq{
\[
M\left(x,y\right)=P^{\left(X_{n}=x\right)}\left(X_{n+1}=y\right).
\]
}

\end{remark}

We can now reformulate Proposition $\ref{pr:chapman_kolm_eq}$ for
homogeneous Markov chains.

\begin{proposition}{Chapman-Kolmogorov Equality for Homogeneous Markov Chains}{}

Let $X$ be a homogeneous Markov chain with transition matrix $M.$
For every instants $n,\,r,\,n+m$ such that $0\leqslant n<r<n+m,$
the Chapman-Kolmogorov relation holds
\[
\forall y\in E,\,\,\,\,M^{m}\left(X_{n},y\right)=\sum_{z\in E}M^{r-n}\left(X_{n},z\right)M^{n+m-r}\left(z,y\right),
\]
which can be written in matrix form as
\[
M^{m}\left(X_{n},\cdot\right)=M^{r-n}M^{n+m-r}\left(X_{n},\cdot\right).
\]
In particular, for $P_{\left(X_{0},\cdots,X_{n}\right)}-$almost every
$\left(x_{0},\cdots,x_{n}\right)$ and every $y\in E,$
\begin{equation}
P^{\left(X_{0}=x_{0},\cdots,X_{n}=x_{n}\right)}\left(X_{n+m}=y\right)=P^{\left(X_{n}=x\right)}\left(X_{n+m}=y\right)=M^{m}\left(x,y\right).\label{eq:PcondX_nplm_y_kn_prev_steps}
\end{equation}

\end{proposition}

\begin{proposition}{Sufficient and Necessary Condition to be a Homogeneous Markov Chain}{}

A process $X$ taking values in $E$ is a homogeneous Markov chain---with
respect to its natural filtration---with transition matrix $M$ if
and only if, for $P_{\left(X_{0},\cdots,X_{n}\right)}-$almost every
$\left(x_{0},\cdots,x_{n-1},x\right)$ and for every $y\in E,$
\begin{equation}
P^{\left(X_{0}=x_{0},\cdots,X_{n}=x\right)}\left(X_{n+m}=y\right)=P^{\left(X_{n}=x\right)}\left(X_{n+m}=y\right)=M^{m}\left(x,y\right).\label{eq:PcondX_nplm_y_kn_prev_steps-1}
\end{equation}

\end{proposition}

\begin{proof}{}{}

The necessity of the condition is immediate.

Conversely, let $0\leqslant n<m.$ Then, for every function $f$ on
$E,$ 
\begin{multline*}
\mathbb{E}^{\left(X_{0}=x_{0},\cdots,X_{n}=x\right)}\left(f\left(X_{m}\right)\right)=\sum_{y\in E}\left[P^{\left(X_{0}=x_{0},\cdots,X_{n}=x\right)}\left(X_{m}=y\right)f\left(y\right)\right]=M^{m-n}\left(x,f\right),
\end{multline*}
Since
\[
\mathbb{E}^{\mathscr{B}_{n}}\left(f\left(X_{m}\right)\right)=\mathbb{E}^{\left(X_{0}=\cdot,\cdots,X_{n}=\cdot\right)}\left(f\left(X_{m}\right)\right)\circ\left(X_{0},\cdots,X_{n}\right),
\]
we obtain
\[
\mathbb{E}^{\mathscr{B}_{n}}\left(f\left(X_{m}\right)\right)=M^{m-n}\left(X_{n},f\right).
\]
This proves that $X$ is a homogeneous Markov chain with transition
matrix $M.$

\end{proof}

\begin{remark}{}{}

If $X$ is a homogeneous Markov chain with transition matrix $M$
and if $f$ is a nonnegative or bounded function such that $M\left(\cdot,f\right)=f$---such
a function is said \textbf{harmonic}\index{harmonic function}\mindex{function!harmonic}---that
is, in vector form, if $f$ is a right \textbf{eigenvector} of $M$
associated with the eigenvalue 1, then
\[
\mathbb{E}^{\mathscr{A}_{n}}\left(f\left(X_{n+1}\right)\right)=M\left(X_{n},f\right)=f\left(X_{n}\right).
\]
Therefore, the process $\left(f\left(X_{n}\right)\right)_{n\in\mathbb{N}}$
is a \textbf{martingale}.

\end{remark}

\subsection{Simple Markov Property. Finite-Dimensional Laws}

We now generalize the formula $\refpar{eq:hom_mark_chain_def}$ to
the case of an arbitrary functional of the future\footnote{By a \textbf{\mindex{function!future}function of the future}, we
mean a random variable measurable with respect to the \salg of the
\textbf{wide future} after the instant $n.$ Such a functional can
be written in the form $f\left(X_{n},X_{n+1},\cdots\right),$ where
$f$ is a measurable function defined on $E^{\mathbb{N}}.$ An example
is the first hitting time of the process $X$ in a given state after
time $n.$} of a homogeneous Markov chain after time $n.$

The property thus obtained is called the \textbf{simple}\mindex{Markov property!simple}---or
\textbf{weak}\mindex{Markov property!weak}---\textbf{Markov property}.
It expresses the time homogeneity of the chain $X,$ and says that,
at every time $n,$ the conditional expectation of any functional
of the future of the chain, given the present state at time $n,$
is equal to the mean value of the same functional evaluated along
the entire trajectory---that is from the time 0---of a Markov chain
with the same transition matrix as $X,$ but started at time 0 from
the state $X_{n}.$ 

To formulate this property precisely, we introduce the \textbf{shift
operators\index{shift operators}} $\theta_{n},$ $n\in\mathbb{N},$
defined on the space $E^{\mathbb{N}}$ of sequences taking values
in $E.$ For every $y\in E^{\mathbb{N}},$
\[
\theta_{n}\left(y\right)=\left(y_{n+p}\right)_{p\in\mathbb{N}}.
\]
Thus, the sequence with general term $\theta_{n}\left(y\right)$ is
the sequence obtained from $y$ removing the first $n$ terms $y_{0},\cdots,y_{n-1}.$
A \textbf{functional of the future\index{functional of the future}}
of a Markov chain $X$ \textbf{after time $n$} can therefore be written
in the form $f\left(\theta_{n}\left(X\right)\right),$ where $f$
is a function defined on $E^{\mathbb{N}}.$

\begin{proposition}{Simple Markov Property}{simple_markov_prop}

Let $X$ be a homogeneous Markov chain with transition matrix $M.$
For every time $n\in\mathbb{N}$ and every $\mathscr{E}^{\otimes\mathbb{N}}-$measurable
function $f$ on $E^{\mathbb{N}},$ nonnegative or bounded,
\begin{equation}
\mathbb{E}^{\mathscr{A}_{n}}\left(f\left(\theta_{n}\left(X\right)\right)\right)=g\left(X_{n}\right),\label{eq:simple_mark_prop}
\end{equation}
where $g$ is a measurable function on $E$ defined, for every $x$
such that $P\left(X_{0}=x\right)>0,$ by
\[
g\left(x\right)=\mathbb{E}^{X=x}\left(f\left(X\right)\right).
\]
In particular,
\begin{equation}
g\left(X_{0}\right)=\mathbb{E}^{\mathscr{A}_{0}}\left(f\left(X\right)\right).\label{eq:cond_g_simple_mark_prop}
\end{equation}

\end{proposition}

\begin{proof}{}{}

We now present two proofs of this property. The first is heuristic,
and ignores issues related to events of zero probability. Time is
there used in its natural sense\footnotemark. The second proof is
rigorous but more formal. Time is then used in a backward manner\footnotemark.
\begin{itemize}
\item \textbf{Heuristic proof}\\
By an extension by measurability argument---see Theorem $\ref{th:meas_ext_pple},$
it suffices to prove $\refpar{eq:simple_mark_prop}$ for functionals
$f$ of the form $y\mapsto f_{k}\left(y_{0},\cdots,y_{k}\right),$
where $k\in\mathbb{N}$ and $f_{k}$ is a function on $E^{k+1}.$
For such a functional $f,$ the equality $\refpar{eq:cond_exp_past_futur_event}$
of Theorem $\ref{th:mark_ch_charact}$ allows to write that
\begin{equation}
\mathbb{E}^{\mathscr{A}_{n}}\left(f\left(\theta_{n}\left(X\right)\right)\right)=\mathbb{E}^{\mathscr{P}_{n}}\left(f\left(\theta_{n}\left(X\right)\right)\right)=\mathbb{E}^{\mathscr{P}_{n}}\left(f_{k}\left(X_{n},\cdots,X_{n+k}\right)\right).\label{eq:cond_exp_f_theta_n_X_kn_A_n}
\end{equation}
Now, for $P_{X_{n}}-$almost every $x\in E,$
\begin{multline*}
\mathbb{E}^{\left(X_{n}=x\right)}\left(f_{k}\left(X_{n},\cdots,X_{n+k}\right)\right)=\sum_{\left(y_{1},\cdots,y_{k}\right)\in E^{k}}f_{k}\left(x,y_{1},\cdots,y_{k}\right)\\
\times P^{\left(X_{n}=x\right)}\left(X_{n}=x,X_{n+1}=y_{1},\cdots,X_{n+k}=y_{k}\right).
\end{multline*}
By the chain rule for conditional probabilities and the Markov property,
for every $n\in\mathbb{N},$
\begin{align*}
 & P^{\left(X_{n}=x\right)}\left(X_{n}=x,X_{n+1}=y_{1},\cdots,X_{n+k}=y_{k}\right)\\
 & \qquad=P^{\left(X_{n}=x\right)}\left(X_{n+1}=y_{1}\right)P^{\left(X_{n}=x,X_{n+1}=y_{1}\right)}\left(X_{n+2}=y_{2}\right)\times\cdots\\
 & \qquad\qquad\cdots\times P^{\left(X_{n}=x,\cdots,X_{n+k-1}=y_{k-1}\right)}\left(X_{n+k}=y_{k}\right)\\
 & \qquad=P^{\left(X_{n}=x\right)}\left(X_{n+1}=y_{1}\right)P^{\left(X_{n+1}=y_{1}\right)}\left(X_{n+2}=y_{2}\right)\times\cdots\\
 & \qquad\qquad\cdots\times P^{\left(X_{n+k-1}=y_{k-1}\right)}\left(X_{n+k}=y_{k}\right),
\end{align*}
Hence,
\begin{multline}
P^{\left(X_{n}=x\right)}\left(X_{n}=x,X_{n+1}=y_{1},\cdots,X_{n+k}=y_{k}\right)=M\left(x,y_{1}\right)M\left(y_{1},y_{2}\right)\\
\cdots\times M\left(y_{k-1},y_{k}\right).\label{eq:P_X_nisxtoX_npkisyk_kn_X_nisx}
\end{multline}
It then follows that, for every $n\in\mathbb{N},$
\begin{multline*}
\mathbb{E}^{\left(X_{n}=x\right)}\left(f_{k}\left(X_{n},\cdots,X_{n+k}\right)\right)=\sum_{\left(y_{1},\cdots,y_{k}\right)\in E^{k}}f_{k}\left(x,y_{1},\cdots,y_{k}\right)\\
\times M\left(x,y_{1}\right)M\left(y_{1},y_{2}\right)\cdots\times M\left(y_{k-1},y_{k}\right).
\end{multline*}
\\
Defining the function $g$ on $E$ by
\[
g\left(x\right)=\sum_{\left(y_{1},\cdots,y_{k}\right)\in E^{k}}f_{k}\left(x,y_{1},\cdots,y_{k}\right)M\left(x,y_{1}\right)M\left(y_{1},y_{2}\right)\cdots\times M\left(y_{k-1},y_{k}\right).
\]
For $n=0,$ we have
\[
g\left(x\right)=\mathbb{E}^{\left(X_{0}=x\right)}\left(f_{k}\left(X_{0},\cdots,X_{k}\right)\right)=\mathbb{E}^{\left(X_{0}=x\right)}\left(f\left(X\right)\right).
\]
Therefore, by $\refpar{eq:cond_exp_f_theta_n_X_kn_A_n},$ 
\[
\mathbb{E}^{\mathscr{A}_{n}}\left(f\left(\theta_{n}\left(X\right)\right)\right)=g\left(X_{n}\right),
\]
which proves the result heuristically.
\item \textbf{Rigorous proof}\\
Again, by an extension by measurability argument, it suffices to prove
$\refpar{eq:simple_mark_prop}$ for a functional $f$ of the form
$y\mapsto\prod^{k}_{j=0}f_{j}\left(y_{j}\right),$ where $k\in\mathbb{N}$
and where the $f_{j}$ are functions on $E.$ \\
Using conditioning with respect to the \salg $\mathscr{A}_{n+k-1},$
the adaptation of the process $X,$ and equality $\refpar{eq:cond_exp_f_theta_n_X_kn_A_n},$
\begin{align*}
\mathbb{E}^{\mathscr{A}_{n}}\left(f\left(\theta_{n}\left(X\right)\right)\right) & =\mathbb{E}^{\mathscr{A}_{n}}\left(\prod^{k}_{j=0}f_{j}\left(X_{n+j}\right)\right)\\
 & =\mathbb{E}^{\mathscr{A}_{n}}\left(\prod^{k-1}_{j=0}f_{j}\left(X_{n+j}\right)\mathbb{E}^{\mathscr{A}_{n+k-1}}\left(f_{k}\left(X_{n+k}\right)\right)\right).
\end{align*}
Thus,
\[
\mathbb{E}^{\mathscr{A}_{n}}\left(f\left(\theta_{n}\left(X\right)\right)\right)=\mathbb{E}^{\mathscr{A}_{n}}\left(\prod^{k-1}_{j=0}f_{j}\left(X_{n+j}\right)M\left(X_{n+k-1},f_{k}\right)\right).
\]
Conditioning with respect to the \salg $\mathscr{A}_{n+k-2},$ and
applying the same arguments as before, we obtain
\begin{align*}
\mathbb{E}^{\mathscr{A}_{n}}\left(f\left(\theta_{n}\left(X\right)\right)\right) & =\mathbb{E}^{\mathscr{A}_{n}}\left(\prod^{k-2}_{j=0}f_{j}\left(X_{n+j}\right)\mathbb{E}^{\mathscr{A}_{n+k-2}}\left[M\left(X_{n+k-1},f_{k}\right)f_{k-1}\left(X_{n+k-1}\right)\right]\right)\\
 & =\mathbb{E}^{\mathscr{A}_{n}}\left(\prod^{k-2}_{j=0}f_{j}\left(X_{n+j}\right)M\left(X_{n+k-2},f_{k-1}M\left(\cdot,f_{k}\right)\right)\right).
\end{align*}
Now, for every $x\in E,$
\begin{align*}
M\left(x,f_{k-1}M\left(\cdot,f_{k}\right)\right) & =\sum_{y_{k-1}\in E}M\left(x,y_{k-1}\right)f_{k-1}\left(y_{k-1}\right)M\left(y_{k-1},f_{k}\right)\\
 & =\sum_{y_{k-1}\in E}M\left(x,y_{k-1}\right)f_{k-1}\left(y_{k-1}\right)\left[\sum_{y_{k}\in E}M\left(y_{k-1},y_{k}\right)f_{k}\left(y_{k}\right)\right]\\
 & =\sum_{y_{k-1},y_{k}\in E^{2}}M\left(x,y_{k-1}\right)M\left(y_{k-1},y_{k}\right)f_{k-1}\left(y_{k-1}\right)f_{k}\left(y_{k}\right).
\end{align*}
Expanding the last term,
\begin{multline*}
\mathbb{E}^{\mathscr{A}_{n}}\left(f\left(\theta_{n}\left(X\right)\right)\right)=\mathbb{E}^{\mathscr{A}_{n}}\left[\prod^{k-2}_{j=0}f_{j}\left(X_{n+j}\right)\right.\\
\left.\times\sum_{y_{k-1},y_{k}\in E^{2}}M\left(X_{n+k-2},y_{k-1}\right)M\left(y_{k-1},y_{k}\right)f_{k-1}\left(y_{k-1}\right)f_{k}\left(y_{k}\right)\right].
\end{multline*}
Iterating this procedure by induction yields
\[
\mathbb{E}^{\mathscr{A}_{n}}\left(f\left(\theta_{n}\left(X\right)\right)\right)=\sum_{\left(y_{1},\cdots,y_{k}\right)\in E^{k}}M\left(X_{n},y_{1}\right)M\left(y_{1},y_{2}\right)\cdots M\left(y_{k-1},y_{k}\right)f\left(X_{n},y_{1},\cdots,y_{k}\right).
\]
Defining the function $g$ on $E$ by
\[
g\left(x\right)=\sum_{\left(y_{1},\cdots,y_{k}\right)\in E^{k}}M\left(x,y_{1}\right)M\left(y_{1},y_{2}\right)\cdots M\left(y_{k-1},y_{k}\right)f\left(x,y_{1},\cdots,y_{k}\right).
\]
Then
\[
\mathbb{E}^{\mathscr{A}_{n}}\left(f\left(\theta_{n}\left(X\right)\right)\right)=g\left(X_{n}\right).
\]
In particular, this computation shows that
\[
g\left(X_{0}\right)=\mathbb{E}^{\mathscr{A}_{0}}\left(f\left(X\right)\right),
\]
which completes the rigorous proof.
\end{itemize}
\end{proof}

\addtocounter{footnote}{-1}

\footnotetext{Tr.N: This corresponds to the forward direction, reasoning
on the process by conditioning on the past states and analyzing future
states.

}

\stepcounter{footnote}

\footnotetext{Tr.N.: This corresponds to the backward direction,
reasoning on the process by conditioning on future states and analyzing
past states.}

The \textbf{finite-dimensional laws}---that is, the laws of any vector
of components the chain states at finitely many time instants---\textbf{conditional
on the initial state $X_{0}$} of a homogeneous\textbf{ Markov chain
with transition matrix} $M,$ are entirely determined by the matrix
$M$.

\begin{proposition}{Finite-Dimensional Laws Conditional on an Initial State of a Homogeneous Markov Chain}{P_Markov_transition_mat}

Let $X$ be a homogeneous Markov chain with transition matrix $M.$
For every $x\in E$ such that $P\left(X_{0}=x\right)>0,$ denote by
$P_{x}$ the conditional probability $P^{\left(X_{0}=x\right)}.$ 

For every increasing sequence of times $s_{1},s_{2},\cdots,s_{k},$
for every $x$ with $P\left(X_{0}=x\right)>0,$ and for every $\left(x_{1},x_{2},\cdots,x_{k}\right)\in E^{k},$
\begin{multline}
P_{x}\left(X_{s_{1}}=x_{1},X_{s_{2}}=x_{2},\cdots,X_{s_{k}}=x_{k}\right)=M^{s_{1}}\left(x,x_{1}\right)M^{s_{2}-s_{1}}\left(x_{1},x_{2}\right)\\
\cdots M^{s_{k}-s_{k-1}}\left(x_{k-1},x_{k}\right).\label{eq:P_xofgivenstate}
\end{multline}
It follows that
\begin{multline}
P\left(X_{s_{1}}=x_{1},X_{s_{2}}=x_{2},\cdots,X_{s_{k}}=x_{k}\right)=\sum_{x\in E}P\left(X_{0}=x\right)M^{s_{1}}\left(x,x_{1}\right)\\
M^{s_{2}-s_{1}}\left(x_{1},x_{2}\right)\cdots M^{s_{k}-s_{k-1}}\left(x_{k-1},x_{k}\right).\label{eq:P_ofgivenstate}
\end{multline}

\end{proposition}

\begin{proof}{}{}

As in Proposition $\ref{pr:simple_markov_prop},$ we now present two
proofs, one heuristic the other rigorous.
\begin{itemize}
\item \textbf{Heuristic proof}\\
By the chain rule for conditional probabilities and the Markov property,
for every $k\in\mathbb{N},$
\begin{align*}
 & P_{x}\left(X_{s_{1}}=x_{1},X_{s_{2}}=x_{2},\cdots,X_{s_{k}}=x_{k}\right)\\
 & \qquad=P_{x}\left(X_{s_{1}}=x_{1}\right)P^{\left(X_{s_{1}}=x_{1}\right)}\left(X_{s_{2}}=x_{2}\right)\cdots\times P^{\left(X_{s_{1}}=x_{1},\cdots,X_{s_{k-1}}=x_{k-1}\right)}\left(X_{s_{k}}=x_{k}\right)\\
 & \qquad=P_{x}\left(X_{s_{1}}=x_{1}\right)P^{\left(X_{s_{1}}=x_{1}\right)}\left(X_{s_{2}}=x_{2}\right)\cdots\times P^{\left(X_{s_{k-1}}=x_{k-1}\right)}\left(X_{s_{k}}=x_{k}\right),
\end{align*}
which proves the equality $\refpar{eq:P_xofgivenstate},$ by the Chapman-Kolmogorov
equality. It is then enough to apply the total probabilities formula
to obtain the equality $\refpar{eq:P_ofgivenstate}$.
\item \textbf{Rigorous proof}\\
Conditioning with respect to the \salg $\mathscr{A}_{s_{k-1}},$
using the adaptedness of the process $X,$ and applying the equality
$\refpar{eq:hom_mark_chain_def},$ we obtain
\begin{align*}
 & P_{x}\left(X_{s_{1}}=x_{1},X_{s_{2}}=x_{2},\cdots,X_{s_{k}}=x_{k}\right)\\
 & \qquad=\mathbb{E}_{x}\left[\boldsymbol{1}_{\left(X_{s_{1}}=x_{1},X_{s_{2}}=x_{2},\cdots,X_{s_{k-1}}=x_{k-1}\right)}\mathbb{E}^{\mathscr{A}_{s_{k-1}}}_{x}\left(\boldsymbol{1}_{\left(X_{s_{k}}=x_{k}\right)}\right)\right]\\
 & \qquad=\mathbb{E}_{x}\left[\boldsymbol{1}_{\left(X_{s_{1}}=x_{1},X_{s_{2}}=x_{2},\cdots,X_{s_{k-1}}=x_{k-1}\right)}M^{s_{k}-s_{k-1}}\left(X_{s_{k-1}},x_{k}\right)\right]\\
 & \qquad=\mathbb{E}_{x}\left[\boldsymbol{1}_{\left(X_{s_{1}}=x_{1},X_{s_{2}}=x_{2},\cdots,X_{s_{k-1}}=x_{k-1}\right)}\right]M^{s_{k}-s_{k-1}}\left(x_{k-1},x_{k}\right).
\end{align*}
Iterating,
\begin{align*}
 & P_{x}\left(X_{s_{1}}=x_{1},X_{s_{2}}=x_{2},\cdots,X_{s_{k}}=x_{k}\right)\\
 & \,\,\,\,=\mathbb{E}_{x}\left[\boldsymbol{1}_{\left(X_{s_{1}}=x_{1},X_{s_{2}}=x_{2},\cdots,X_{s_{k-1}}=x_{k-2}\right)}\right]M^{s_{k-1}-s_{k-2}}\left(x_{k-2},x_{k-1}\right)M^{s_{k}-s_{k-1}}\left(x_{k-1},x_{k}\right)\\
 & \,\,\,\,=\cdots\\
 & \,\,\,\,=\mathbb{E}_{x}\left[\mathbb{E}^{\mathscr{A}_{0}}\left(\boldsymbol{1}_{\left(X_{s_{1}}=x_{1}\right)}\right)\right]M^{s_{2}-s_{1}}\left(x_{1},x_{2}\right)\cdots M^{s_{k-1}-s_{k-2}}\left(x_{k-2},x_{k-1}\right)M^{s_{k}-s_{k-1}}\left(x_{k-1},x_{k}\right).\\
 & \,\,\,\,=\mathbb{E}_{x}\left[M^{s_{1}}\left(X_{0},x_{1}\right)\right]M^{s_{2}-s_{1}}\left(x_{1},x_{2}\right)\cdots M^{s_{k-1}-s_{k-2}}\left(x_{k-2},x_{k-1}\right)M^{s_{k}-s_{k-1}}\left(x_{k-1},x_{k}\right).
\end{align*}
which proves the equality $\refpar{eq:P_xofgivenstate},$ the equality
$\refpar{eq:P_ofgivenstate}$ follows immediately.
\end{itemize}
\end{proof}

We now obtain a characterization of homogeneous Markov chains with
respect to the natural filtration.

\begin{proposition}{Characterization of Homogeneous Markov Chains with respect to the Natural Filtration}{cns_hom_Markov_chain}

A process $X=\left(X_{n}\right)_{n\in\mathbb{N}},$ defined on the
probabilized space \preds is a homogeneous Markov chain with \textbf{initial
law $\mu$}---that is such that, for every $x\in E,$ $P\left(X_{0}=x\right)=\mu\left(x\right)$---and
transition matrix $M$ if and only if, for every $k\in\mathbb{N}^{\ast}$
and for every $x_{0},x_{1},\cdots,x_{k}\in E,$
\begin{equation}
P\left(X_{0}=x_{0},X_{1}=x_{1},\cdots,X_{k}=x_{k}\right)=\mu\left(x\right)M\left(x_{0},x_{1}\right)M\left(x_{1},x_{2}\right)\cdots M\left(x_{k-1},x_{k}\right).\label{eq:P_path_and_trans_mat}
\end{equation}

\end{proposition}

\begin{proof}{}{}

The necessary condition follows by a straightforward adaptation of
the proof of Proposition $\ref{pr:P_Markov_transition_mat}.$ 

Conversely, suppose that $\refpar{eq:P_path_and_trans_mat}$ holds,
and let $f\in bE$ be arbitrary. For every $x_{0},x_{1},\cdots,x_{k}\in E,$
\begin{align*}
 & \intop_{\left(X_{0}=x_{0},X_{1}=x_{1},\cdots,X_{k}=x_{k}\right)}f\left(X_{k+1}\right)\text{d}P\\
 & \hfill=\sum_{x_{k+1}}f\left(x_{k+1}\right)P\left(X_{0}=x_{0},X_{1}=x_{1},\cdots,X_{k+1}=x_{k+1}\right)\\
 & \hfill=\mu\left(x\right)M\left(x_{0},x_{1}\right)M\left(x_{1},x_{2}\right)\cdots M\left(x_{k-1},x_{k}\right)\sum_{x_{k+1}}f\left(x_{k+1}\right)M\left(x_{k},x_{k+1}\right)\\
 & \hfill=\mu\left(x\right)M\left(x_{0},x_{1}\right)M\left(x_{1},x_{2}\right)\cdots M\left(x_{k-1},x_{k}\right)M\left(x_{k},f\right)\\
 & \hfill=\intop_{\left(X_{0}=x_{0},X_{1}=x_{1},\cdots,X_{k}=x_{k}\right)}M\left(X_{k},f\right)\text{d}P.
\end{align*}
With the previous notations, this proves that
\[
\mathbb{E}^{\mathscr{B}_{k}}\left(f\left(X_{k+1}\right)\right)=M\left(X_{k},f\right),
\]
and thus that $X$ is a homogeneous Markov chain with transition matrix
$M.$ Of course, its initial law is $\mu.$

\end{proof}

\begin{remark}{}{}

By the Carathéodory extension theorem, given a probability $\mu$
on $E$ and a stochastic matrix $M$ indexed on $E\times E,$ there
exists a unique probability $P_{\mu}$ on the probabilizable space
$\left(E^{\mathbb{N}},\mathscr{E}^{\otimes\mathbb{N}}\right)$ such
that the coordinate process $X=\left(X_{n}\right)_{n\in\mathbb{N}}$
is a homogeneous Markov chain with initial law $\mu$ and such that,
for every $k\in\mathbb{N}^{\ast}$ and every $x_{0},x_{1},\cdots,x_{k}\in E,$
equality $\refpar{eq:P_path_and_trans_mat}.$ This chain is called
the \textbf{canonical homogeneous Markov chain\mindex{Markov chain!homogeneous!canonical}}
with initial law $\mu$ and transition matrix $M.$ This result is
a particular case of the Ionescu-Tulcea theorem \cite{tulcea1949mesures}.

\end{remark}

\subsection{Starting Law. Strong Markov Property}

\textbf{Henceforth, we restrict attention to homogeneous Markov chains.}
These processes satisfy the \textbf{strong Markov property}, which
is the analogue of the simple---or weak---Markov property where
deterministic times are replaced by \textbf{stopping times}. This
fundamental property will be essential for establishing the main structural
results of Markov chains. 

Before stating it properly, we first clarify how to construct a homogeneous
Markov chain with a \textbf{prescribed initial law}, which motivates
the following definition.

\begin{definition}{Homogeneous Markov Chain with Given Initial Law}{}

A process $X=\left(X_{n}\right)_{n\in\mathbb{N}}$ taking values in
$E$ is called a homogenous Markov chain on the process basis $\left(\Omega,\mathscr{A},\left(\mathscr{A}_{n}\right)_{n\in\mathbb{N}},P\right)$
with inital law $\nu$ and transition matrix $M$ if

(a) $P_{X_{0}}=\nu$

(b) $X$ is a homogeneous Markov chain with transition matrix $M$
on the underlying filtered probabilized space $\left(\Omega,\mathscr{A},\left(\mathscr{A}_{n}\right)_{n\in\mathbb{N}},P\right).$

\end{definition}

\begin{remark}{}{}

With the preceding notations, if for every $x\in E,$ the process
$X$ is a homogeneous Markov chain on the underlying filtered probabilized
space $\left(\Omega,\mathscr{A},\left(\mathscr{A}_{n}\right)_{n\in\mathbb{N}},P\right)$
with initial law $\delta_{x}$---thus $P\left(X_{0}=x\right)=1$---and
of transition matrix, then $P_{x}=P.$

\end{remark}

The next proposition shows that, if the chain can be started from
any point $x,$ then it can also be started from any initial law $\nu.$

\begin{proposition}{}{}

Suppose that for every $x\in S\subset E,$ the process $X$ is a homogeneous
Markov chain on the underlying filtered probabilized space $\left(\Omega,\mathscr{A},\left(\mathscr{A}_{n}\right)_{n\in\mathbb{N}},P_{x}\right)$
with initial law $\delta_{x}$ and transition matrix $M.$ Let $\nu$
be a probability on $E$ such that $\nu\left(S\right)=1.$ Define
a set function $P_{\nu}$ on $\mathscr{A},$ for every $A\in\mathscr{A},$
by
\[
P_{\nu}\left(A\right)=\sum_{x\in S}\nu\left(x\right)P_{x}\left(A\right)
\]
Then $P_{\nu}$ is a probability on $\left(\Omega,\mathscr{A}\right),$
and $X$ is a homogeneous Markov chain on the underlying filtered
probabilized space $\left(\Omega,\mathscr{A},\left(\mathscr{A}_{n}\right)_{n\in\mathbb{N}},P_{\nu}\right)$
with initial law $\nu$ and transition matrix $M.$

\end{proposition}

\begin{proof}{}{}

Clearly, $P_{\nu}\left(\emptyset\right)=0.$ 

Moreover, $P_{\nu}$ is $\sigma-$additive. Let $\left(A_{n}\right)_{n\in\mathbb{N}}$
be a sequence of disjoint events. Since each $P_{x}$ is a probability
measure,
\[
P_{\nu}\left(\biguplus_{n\in\mathbb{N}}A_{n}\right)=\sum_{x\in S}\nu\left(x\right)\left[\sum_{n\in\mathbb{N}}P_{x}\left(A_{n}\right)\right].
\]
As all terms are nonnegative, we may interchange the sums,
\[
P_{\nu}\left(\biguplus_{n\in\mathbb{N}}A_{n}\right)=\sum_{n\in\mathbb{N}}\left[\sum_{x\in S}\nu\left(x\right)P_{x}\left(A_{n}\right)\right]=\sum_{n\in\mathbb{N}}P_{\nu}\left(A_{n}\right).
\]
Moreover, by standard integration properties, for every nonnegative
or bounded random variable $Y,$
\[
\mathbb{E}_{\nu}\left(Y\right)=\sum_{x\in S}\nu\left(x\right)\mathbb{E}_{x}\left(Y\right),
\]
where $\mathbb{E}_{\nu}$---respectively $\mathbb{E}_{x}$---denotes
expectation with respect to $P_{\nu}$---respectively $P_{x}.$

Let $m,n$ be two nonnegative integers such that $n<m$ and let $f\in bE.$
For every $A\in\mathscr{A}_{n},$
\begin{align*}
\mathbb{E}_{\nu}\left(\boldsymbol{1}_{A}f\left(X_{m}\right)\right) & =\sum_{x\in S}\nu\left(x\right)\mathbb{E}_{x}\left[\boldsymbol{1}_{A}f\left(X_{m}\right)\right]\\
 & =\sum_{x\in S}\nu\left(x\right)\mathbb{E}_{x}\left[\boldsymbol{1}_{A}\mathbb{E}^{\mathscr{A}_{n}}_{x}\left(f\left(X_{m}\right)\right)\right].
\end{align*}
Thus, applying the equality $\refpar{eq:hom_mark_chain_def}$ to the
homogeneous Markov $P_{x}-$chains,
\begin{align*}
\mathbb{E}_{\nu}\left(\boldsymbol{1}_{A}f\left(X_{m}\right)\right) & =\sum_{x\in S}\nu\left(x\right)\mathbb{E}_{x}\left[\boldsymbol{1}_{A}M^{m-n}\left(X_{n},f\right)\right]\\
 & =\mathbb{E}_{\nu}\left[\boldsymbol{1}_{A}M^{m-n}\left(X_{n},f\right)\right].
\end{align*}
Since $M^{m-n}\left(X_{n},f\right)$ is $\mathscr{A}_{n}-$measurable,
\begin{equation}
\mathbb{E}^{\mathscr{A}_{n}}_{\nu}\left[f\left(X_{m}\right)\right]=M^{m-n}\left(X_{n},f\right).\label{eq:cond_exp_nu_fX_m_kn_A_n}
\end{equation}
Thus $X$ is a homogeneous Markov chain on the underlying filtered
probabilized space $\left(\Omega,\mathscr{A},\left(\mathscr{A}_{n}\right)_{n\in\mathbb{N}},P_{\nu}\right).$ 

Finally, for every $B\in\mathscr{E},$
\[
P_{\nu}\left(X_{0}\in B\right)=\sum_{x\in S}\nu\left(x\right)P_{x}\left(X_{0}\in B\right)=\sum_{x\in S}\nu\left(x\right)\delta_{x}\left(B\right)=\nu\left(S\cap B\right)=\nu\left(B\right).
\]
Hence $X$ has initial law $\nu$ for $P_{\nu}.$ 

\end{proof}

To complete this result, we now show that a homogeneous Markov chain
on the underlying filtered probabilized space $\left(\Omega,\mathscr{A},\left(\mathscr{A}_{n}\right)_{n\in\mathbb{N}},P\right)$
can be resumed from $P_{X_{0}}-$almost every point of $E.$ 

\begin{proposition}{Change of Starting Point of a Homogeneous Markov Chain}{}

Let $X$ be a homogeneous Markov chain with transition matrix $M$
on the underlying filtered probabilized space $\left(\Omega,\mathscr{A},\left(\mathscr{A}_{n}\right)_{n\in\mathbb{N}},P\right).$
For every $x\in E$ such that $P\left(X_{0}=x\right)>0,$ define the
probability $P_{x}=P\left(\cdot\mid X_{0}=x\right).$ Then $X$ is
a homogeneous Markov chain on the underlying filtered probabilized
space $\left(\Omega,\mathscr{A},\left(\mathscr{A}_{n}\right)_{n\in\mathbb{N}},P_{x}\right),$
with initial law $\delta_{x}$ and transition matrix $M.$

\end{proposition}

\begin{proof}{}{}

By definition,
\[
P_{x}\left(X_{0}=x\right)=P^{\left(X_{0}=x\right)}\left(X_{0}=x\right)=1
\]
and therefore $X_{0}\left(P_{x}\right)=\delta_{x}.$ 

Moreover, for every nonnegative or bounded random variable $Y,$
\[
\mathbb{E}_{x}\left(Y\right)=\dfrac{1}{P\left(X_{0}=x\right)}\mathbb{E}\left(\boldsymbol{1}_{\left(X_{0}=x\right)}Y\right).
\]

Let $m,n\in\mathbb{N}$ with $n<m,$ and let $f\in bE.$ For every
$A\in\mathscr{A}_{n},$ since $\left(X_{0}=x\right)\cap A\in\mathscr{A}_{n},$
\begin{align*}
\mathbb{E}_{x}\left(\boldsymbol{1}_{A}f\left(X_{m}\right)\right) & =\dfrac{1}{P\left(X_{0}=x\right)}\mathbb{E}\left(\boldsymbol{1}_{\left(X_{0}=x\right)}\boldsymbol{1}_{A}f\left(X_{m}\right)\right)\\
 & =\dfrac{1}{P\left(X_{0}=x\right)}\mathbb{E}\left(\boldsymbol{1}_{\left(X_{0}=x\right)}\boldsymbol{1}_{A}\mathbb{E}^{\mathscr{A}_{n}}\left[f\left(X_{m}\right)\right]\right)\\
 & =\dfrac{1}{P\left(X_{0}=x\right)}\mathbb{E}\left(\boldsymbol{1}_{\left(X_{0}=x\right)}\boldsymbol{1}_{A}M^{m-n}\left(X_{n},f\right)\right)\\
 & =\mathbb{E}_{x}\left[\boldsymbol{1}_{A}M^{m-n}\left(X_{n},f\right)\right].
\end{align*}
Since $M^{m-n}\left(X_{n},f\right)$ is $\mathscr{A}_{n}-$measurable,
it follows that
\[
\mathbb{E}^{\mathscr{A}_{n}}_{x}\left(f\left(X_{m}\right)\right)=M^{m-n}\left(X_{n},f\right),
\]
which proves that $X$ is a homogeneous Markov chain on the underlying
filtered probabilized space $\left(\Omega,\mathscr{A},\left(\mathscr{A}_{n}\right)_{n\in\mathbb{N}},P_{x}\right).$

\end{proof}

\begin{remark}{}{}

Thus, the last two propositions show that, starting by a homogeneous
Markov chain with transition matrix $M$ on the underlying filtered
probabilized space $\left(\Omega,\mathscr{A},\left(\mathscr{A}_{n}\right)_{n\in\mathbb{N}},P\right),$
we can construct, for every probability $\nu$ on $E$ with the same
support as $P_{X_{0}},$ a probability $P_{\nu}$ on $\left(\Omega,\mathscr{A}\right)$
such that $X$ is a homogeneous Markov chain on the underlying filtered
probabilized space $\left(\Omega,\mathscr{A},\left(\mathscr{A}_{n}\right)_{n\in\mathbb{N}},P_{\nu}\right)$
with initial law $\nu$ and transition matrix $M.$

\end{remark}

With the preceding notations, we can now restate the simple Markov
property stated in Proposition $\ref{pr:simple_markov_prop}$

\begin{proposition}{Simple Markov Property}{simple_markov_prop_2}

Let $X$ be a process such that, for every $x\in E,$ it is a homogeneous
Markov chain on the underlying filtered probabilized space $\left(\Omega,\mathscr{A},\left(\mathscr{A}_{n}\right)_{n\in\mathbb{N}},P_{x}\right)$
with initial law $\delta_{x}$ and transition matrix $M.$ For every
$\mathscr{E}^{\otimes\mathbb{N}}-$measurable, nonnegative or bounded
function $f$ on $E^{\mathbb{N}},$ \boxeq{
\begin{equation}
\mathbb{E}^{\mathscr{A}_{n}}_{x}\left[f\left(\theta_{n}\left(X\right)\right)\right]=\mathbb{E}_{X_{n}}\left(f\left(X\right)\right),\label{eq:simple_marko_prop_2}
\end{equation}
}where $\theta_{n}$ is the shift operator at $n\in\mathbb{N}.$

\end{proposition}

\begin{remark}{}{}

It should be clearly understood that $\mathbb{E}_{X_{n}}\left(f\left(X\right)\right)$
denotes the value at $X_{n}$ of the mapping $x\mapsto\mathbb{E}_{x}\left(f\left(X\right)\right).$

\end{remark}

\begin{proposition}{Strong Markov Property}{strong_markov_prop}

Let $X$ be a process such that, for every $x\in E,$ it is a homogeneous
Markov chain on the underlying filtered probabilized space $\left(\Omega,\mathscr{A},\left(\mathscr{A}_{n}\right)_{n\in\mathbb{N}},P_{x}\right)$
with initial law $\delta_{x}$ and transition matrix $M.$ Then $X$
has the strong Markov property: for every $\mathscr{E}^{\otimes\mathbb{N}}-$measurable,
nonnegative or bounded function $f$ on $E^{\mathbb{N}},$ for every
stopping time $T$ and for every $x\in E,$\boxeq{
\begin{equation}
\mathbb{E}^{\mathscr{A}_{T}}_{x}\left(\boldsymbol{1}_{\left(T<+\infty\right)}f\left(\theta_{T}\left(X\right)\right)\right)=\boldsymbol{1}_{\left(T<+\infty\right)}\mathbb{E}_{X_{T}}\left(f\left(X\right)\right).\label{eq:strong_markov_prop}
\end{equation}
}

\end{proposition}

\begin{proof}{}{}

For every $n\in\overline{\mathbb{N}},$
\[
\boldsymbol{1}_{\left(T=n\right)}\mathbb{E}^{\mathscr{A}_{T}}_{x}\left[\boldsymbol{1}_{\left(T<+\infty\right)}f\left(\theta_{T}\left(X\right)\right)\right]=\boldsymbol{1}_{\left(T=n\right)}\mathbb{E}^{\mathscr{A}_{n}}_{x}\left[\boldsymbol{1}_{\left(T<+\infty\right)}f\left(\theta_{T}\left(X\right)\right)\right].
\]
Since $\left(T=n\right)\in\mathscr{A}_{n},$
\[
\boldsymbol{1}_{\left(T=n\right)}\mathbb{E}^{\mathscr{A}_{T}}_{x}\left[\boldsymbol{1}_{\left(T<+\infty\right)}f\left(\theta_{T}\left(X\right)\right)\right]=\mathbb{E}^{\mathscr{A}_{n}}_{x}\left[\boldsymbol{1}_{\left(T=n\right)}\boldsymbol{1}_{\left(T<+\infty\right)}f\left(\theta_{T}\left(X\right)\right)\right].
\]
Hence,
\[
\boldsymbol{1}_{\left(T=n\right)}\mathbb{E}^{\mathscr{A}_{T}}_{x}\left[\boldsymbol{1}_{\left(T<+\infty\right)}f\left(\theta_{T}\left(X\right)\right)\right]=\begin{cases}
\mathbb{E}^{\mathscr{A}_{n}}_{x}\left[\boldsymbol{1}_{\left(T=n\right)}f\left(\theta_{n}\left(X\right)\right)\right], & \text{if\,}n\in\mathbb{N},\\
0, & \text{if\,}n=+\infty.
\end{cases}
\]
Now, let $n\in\mathbb{N}.$ By the simple Markov property,
\begin{align*}
\boldsymbol{1}_{\left(T=n\right)}\mathbb{E}^{\mathscr{A}_{T}}_{x}\left[\boldsymbol{1}_{\left(T<+\infty\right)}f\left(\theta_{T}\left(X\right)\right)\right] & =\boldsymbol{1}_{\left(T=n\right)}\mathbb{E}^{\mathscr{A}_{n}}_{x}\left[f\left(\theta_{n}\left(X\right)\right)\right]\\
 & =\boldsymbol{1}_{\left(T=n\right)}\mathbb{E}_{X_{n}}\left[f\left(X\right)\right]\\
 & =\boldsymbol{1}_{\left(T=n\right)}\mathbb{E}_{X_{T}}\left[f\left(X\right)\right].
\end{align*}
Since the term corresponding to $n=+\infty$ is equal to zero,
\[
\mathbb{E}^{\mathscr{A}_{T}}_{x}\left[\boldsymbol{1}_{\left(T<+\infty\right)}f\left(\theta_{T}\left(X\right)\right)\right]=\sum_{n\in\mathbb{N}}\boldsymbol{1}_{\left(T=n\right)}\mathbb{E}_{X_{T}}\left[f\left(X\right)\right]=\boldsymbol{1}_{\left(T<+\infty\right)}\mathbb{E}_{X_{T}}\left[f\left(X\right)\right].
\]

\end{proof}

\begin{corollary}{}{}

Under the same assumptions as in Proposition $\refpar{pr:strong_markov_prop},$
let $T$ be a finite stopping time. Define the process $Y$ and the
filtration $\left(\mathscr{B}_{n}\right)_{n\in\mathbb{N}}$ by, for
every $n\in\mathbb{N},$ $Y_{n}=X_{T+n}$ and $\mathscr{B}_{n}=\mathscr{A}_{T+n}.$

Then, for every $x\in E,$ the process $Y$ is a homogeneous Markov
chain with transition matrix $M$ on the underlying filtered probabilized
space $\left(\Omega,\mathscr{A},\left(\mathscr{B}_{n}\right)_{n\in\mathbb{N}},P_{x}\right).$

\end{corollary}

\begin{proof}{}{}

For every $f\in bE$ and for every $m,n\in\mathbb{N}$ such that $m<n,$
\[
\mathbb{E}^{\mathscr{B}_{m}}_{x}\left(f\left(Y_{n}\right)\right)=\mathbb{E}^{\mathscr{A}_{T+m}}_{x}\left[f\left(\left[\theta_{T+m}\left(X\right)\right]_{n-m}\right)\right].
\]
By the strong Markov property applied at the stopping time $T+m,$
\[
\mathbb{E}^{\mathscr{B}_{m}}_{x}\left(f\left(Y_{n}\right)\right)=E_{X_{T+m}}\left[f\left(X_{n-m}\right)\right].
\]
Using $\refpar{eq:hom_mark_chain_def},$ we then obtain
\[
\mathbb{E}^{\mathscr{B}_{m}}_{x}\left(f\left(Y_{n}\right)\right)=M^{n-m}\left(X_{T+m},f\right)=M^{n-m}\left(Y_{m},f\right),
\]
which proves the result.

\end{proof}

\section{Visits to a Fixed State}\label{sec:Visits-with-one}

In this section, we consider a process $X=\left(Y_{n}\right)_{n\in\mathbb{N}}$
such that, for every $x\in E,$ it is a homogeneous Markov chain on
the underlying filtered probabilized space $\left(\Omega,\mathscr{A},\left(\mathscr{A}_{n}\right)_{n\in\mathbb{N}},P_{x}\right)$
with initial law $\delta_{x}$ and transition matrix $M.$ 

We study the entrance times of trajectories of $X$ into a subset
$B$ of $E$ and more especially, when $B$ is a singleton, the hitting
times of some points of $E.$ This leads to a classification of the
states of $E,$ depending on how the chain behaves with respect to
them.

\begin{denotations}{}{}

Let $B\subset E.$ With the convention $\inf\emptyset=+\infty,$ define
\[
T_{B}=\inf\left\{ n\in\mathbb{N}^{\ast}:X_{n}\in B\right\} \,\,\,\,\text{and}\,\,\,\,N_{B}=\sum_{j\in\mathbb{N}}\boldsymbol{1}_{\left(X_{j}\in B\right)}.
\]
These are, respectively, the \textbf{\index{first entrance time}first
entrance time} into $B$ after time 1, and the \textbf{total time
spent\index{total time spent}} in B by the chain during the whole
lifetime of the process. 

In particular, if $B=\left\{ y\right\} ,$ with $y\in E,$ we simply
write $T_{y}$ and $N_{y}.$ 

Similarly, define the functionals $\tau_{B},\,n_{B},\,\tau_{y}$ and
$n_{y}$ on $E^{\mathbb{N}},$ for every $u\in E^{\mathbb{N}},$ by
\[
\tau_{B}\left(u\right)=\inf\left\{ n\in\mathbb{N}^{\ast}:u_{n}\in B\right\} \,\,\,\,\text{and}\,\,\,\,n_{B}\left(u\right)=\sum_{j\in\mathbb{N}}\boldsymbol{1}_{\left(u_{j}\in B\right)}.
\]
We denote $\tau_{y}=\tau_{\left\{ y\right\} }$ and $n_{y}=n_{\left\{ y\right\} }.$

\end{denotations}

The following lemma will be used repeatedly. It is particlarly convenient
for applying the strong Markov property.

\begin{lemma}{}{}

With the preceding notation, we have $T_{B}=\tau_{B}\left(X\right),$
and, for every $p\in\mathbb{N}^{\ast},$\boxeq{
\begin{equation}
\text{on }\left(T_{B}>p\right),\,\,\,\,T_{B}=p+\tau_{B}\left[\theta_{p}\left(X\right)\right].\label{eq:T_B_on_T_Boverp}
\end{equation}
}

In particular\boxeq{
\begin{equation}
\text{on }\left(T_{y}>p\right),\,\,\,\,T_{y}=p+\tau_{y}\left[\theta_{p}\left(X\right)\right].\label{eq:T_y_on_T_yoverp}
\end{equation}
}

As a consequence, for every stopping time $T,$\boxeq{
\begin{equation}
\text{on }\left(T_{B}>T\right)\cap\left(T<+\infty\right),\,\,\,\,T_{B}=T+\tau_{B}\left[\theta_{T}\left(X\right)\right].\label{eq:T_B_on_T_BoverTand Tlessinf}
\end{equation}
}

In particular\boxeq{
\begin{equation}
\text{on }\left(T_{y}>T\right)\cap\left(T<+\infty\right),\,\,\,\,T_{y}=T+\tau_{y}\left[\theta_{T}\left(X\right)\right].\label{eq:T_y_on_T_yoverTandTlessinf}
\end{equation}
}

\end{lemma}

\begin{proof}{}{}

It suffices to observe that, on $\left(T_{B}>p\right),$
\begin{align*}
\tau_{B}\left[\theta_{p}\left(X\right)\right] & =\inf\left\{ n\in\mathbb{N}^{\ast}:\,X_{n+p}\in B\right\} \\
 & =\inf\left\{ n\geqslant p+1:\,X_{n}\in B\right\} -p=T_{B}-p.
\end{align*}
To obtain $\refpar{eq:T_B_on_T_BoverTand Tlessinf}$, it is then enough
to apply $\refpar{eq:T_B_on_T_Boverp}$ on the sets $\left(T=p\right)\cap\left(T_{y}>p\right),\,p\in\mathbb{N}^{\ast}.$

\end{proof}

\begin{denotations}{}{}

We also define, by induction, the sequence of successive entrance
times into $B$ by
\[
T^{0}_{B}=0,\,\,\,\,T^{1}_{B}=T_{B},\,\,\,\,T^{p+1}_{B}=\inf\left\{ n>T^{p}_{B}:\,X_{n}\in B\right\} .
\]

In particular, if $B=\left\{ y\right\} ,$ these times are denoted
simply $T^{p}_{y}$---$T^{p}_{y}$ is the \textbf{time of the $p-$th
hitting time} in $y$. We then prove the same relationship\boxeq{
\begin{equation}
\text{on }\left(T_{B}<+\infty\right),\,\,\,\,T^{p+1}_{B}=T^{p}_{B}+\tau_{B}\left[\theta_{T^{p}_{B}}\left(X\right)\right].\label{eq:T_B_pp1_on_T_B_p}
\end{equation}
}

In particular\boxeq{
\begin{equation}
\text{on }\left(T^{p}_{y}<+\infty\right),\,\,\,\,T^{p+1}_{y}=T^{p}_{y}+\tau_{y}\left[\theta_{T^{p}_{y}}\left(X\right)\right].\label{eq:T_pp1_y_on_T_y_p}
\end{equation}
}

\end{denotations}

\subsection{Study of the Sequence of Hitting Times at a Point}

\begin{proposition}{Hitting Times and Successive Visits to a State}{hitt_time_state}

With the previous notations, for every $p\in\mathbb{N}^{\ast},$ $T^{p}_{B}$
is a stopping time. For every $x,y\in E,$ the sequence $\left(T^{p}_{y}\right)_{p\in\mathbb{N}}$
is a homogeneous Markov chain taking values in $\overline{\mathbb{N}^{\ast}}$
on the underlying filtered probabilized space $\left(\Omega,\mathscr{A},\left(\mathscr{A}_{T^{p}_{y}}\right)_{p\in\mathbb{N}^{\ast}},P_{x}\right).$

If $P_{y}\left(T^{1}_{y}<+\infty\right)=1$---that is if the chain
starting from $y$ returns to $y$ in finite time $P_{y}-$almost
surely---then, for every $p\in\mathbb{N}^{\ast},$ $P_{y}\left(T^{p}_{y}<+\infty\right)=1,$
and consequentely $P_{y}\left(N_{y}=+\infty\right)=1.$ Moreover,
the \textbf{sequence $\left(T^{p+1}_{y}-T^{p}_{y}\right)_{p\in\mathbb{N}^{\ast}}$
of inter-visit times at $y$} is a sequence of random variables---defined
and finite $P_{y}-$almost surely--- \textbf{$P_{y}-$independent,
with the same law}---under $P_{y}$---\textbf{than the one of $T^{1}_{y}.$}

\end{proposition}

\begin{proof}{}{}
\begin{itemize}
\item We already know that $T^{1}_{B}$ is a stopping time. Let $p\geqslant2.$
Since $X$ is adapted, for every $n\in\mathbb{N}^{\ast},$
\[
\left(T^{p}_{B}\leqslant n\right)=\left(\sum^{n}_{j=1}\boldsymbol{1}_{B}\left(X_{j}\right)\geqslant p\right)\in\mathscr{A}_{n}.
\]
Hence $T^{p}_{B}$ is a stopping time.
\item For every bounded function $f$ on $\overline{\mathbb{N}^{\ast}}$
and for every $x\in E,$ we have, by the equality $\refpar{eq:T_B_on_T_BoverTand Tlessinf},$
\begin{multline}
\mathbb{E}^{\mathscr{A}_{T^{p}_{y}}}_{x}\left[f\left(T^{p+1}_{y}\right)\right]=\mathbb{E}^{\mathscr{A}_{T^{p}_{y}}}_{x}\left[\boldsymbol{1}_{\left(T^{p}_{y}<+\infty\right)}f\left(T^{p}_{y}+\tau_{y}\left[\theta_{T^{p}_{y}}\left(X\right)\right]\right)\right]\\
+\mathbb{E}^{\mathscr{A}_{T^{p}_{y}}}_{x}\left[\boldsymbol{1}_{\left(T^{p}_{y}=+\infty\right)}f\left(+\infty\right)\right].\label{eq:cond_exp_x_fofTpp1_y_know_A_t_y_p}
\end{multline}
Hence, since $\left(T^{p}_{y}=+\infty\right)\in\mathscr{A}_{T^{p}_{y}}$
and since $\left(T^{p}_{y}=i\right)\in\mathscr{A}_{i},$
\begin{multline*}
\mathbb{E}^{\mathscr{A}_{T^{p}_{y}}}_{x}\left[f\left(T^{p+1}_{y}\right)\right]=\sum_{i\in\mathbb{N}^{\ast}}\boldsymbol{1}_{\left(T^{p}_{y}=i\right)}\mathbb{E}^{\mathscr{A}_{i}}_{x}\left[f\left(i+\tau_{y}\left[\theta_{i}\left(X\right)\right]\right)\right]+\boldsymbol{1}_{\left(T^{p}_{y}=+\infty\right)}f\left(+\infty\right).
\end{multline*}
Applying the simple Markov property at time $i$---since the strong
Markov property cannot be applied here---we have
\begin{multline*}
\mathbb{E}^{\mathscr{A}_{T^{p}_{y}}}_{x}\left[f\left(T^{p+1}_{y}\right)\right]=\sum_{i\in\mathbb{N}^{\ast}}\boldsymbol{1}_{\left(T^{p}_{y}=i\right)}\mathbb{E}_{x}\left[f\left(i+\tau_{y}\left(X\right)\right)\right]+\boldsymbol{1}_{\left(T^{p}_{y}=+\infty\right)}f\left(+\infty\right).
\end{multline*}
Since $\tau_{y}\left(X\right)=T^{1}_{y},$ and since, by definition
of $T^{p}_{y},$ on $\left(T^{p}_{y}=i\right),$ $X_{i}=X_{T^{p}_{y}}=y,$\boxeq{
\[
\mathbb{E}^{\mathscr{A}_{T^{p}_{y}}}_{x}\left[f\left(T^{p+1}_{y}\right)\right]=\sum_{i\in\mathbb{N}^{\ast}}\boldsymbol{1}_{\left(T^{p}_{y}=i\right)}\mathbb{E}_{y}\left[f\left(i+T^{1}_{y}\right)\right]+\boldsymbol{1}_{\left(T^{p}_{y}=+\infty\right)}f\left(+\infty\right).
\]
}Defining the transition probability $N$ on $\overline{\mathbb{N}^{\ast}}$
by
\begin{equation}
N\left(i,f\right)=\begin{cases}
\mathbb{E}_{y}\left[f\left(i+T^{1}_{y}\right)\right], & \text{if }i\in\mathbb{N}^{\ast},\\
f\left(+\infty\right), & \text{if }i=+\infty,
\end{cases}\label{eq:trans_proba_N_i_f}
\end{equation}
Then
\[
\mathbb{E}^{\mathscr{A}_{T^{p}_{y}}}_{x}\left[f\left(T^{p+1}_{y}\right)\right]=\sum_{i\in\mathbb{N}^{\ast}}\boldsymbol{1}_{\left(T^{p}_{y}=i\right)}N\left(i,f\right)+\boldsymbol{1}_{\left(T^{p}_{y}=+\infty\right)}N\left(+\infty,f\right).
\]
Hence, also\boxeq{
\begin{equation}
\mathbb{E}^{\mathscr{A}_{T^{p}_{y}}}_{x}\left[f\left(T^{p+1}_{y}\right)\right]=N\left(T^{p}_{y},f\right).\label{eq:cond_exp_x_f_T_y_pp1_knA_T_y_p_withN}
\end{equation}
}This proves that the process $\left(T^{p+1}_{y}\right)_{p\in\mathbb{N}}$
is a homogeneous Markov chain, taking values in $\overline{\mathbb{N}^{\ast}}$
of transition matrix $N$ given, for $i,j\in\overline{\mathbb{N}^{\ast}},$
by---taking $f=\boldsymbol{1}_{\left\{ j\right\} }$---
\[
N\left(i,j\right)=\begin{cases}
P_{y}\left(T^{1}_{y}=j-i\right), & \text{if }i,j\in\mathbb{N}^{\ast}\text{ and }j-i\geqslant1,\\
0, & \text{if }j\leqslant i,\\
P_{y}\left(T^{1}_{y}=+\infty\right), & \text{if }i\in\mathbb{N}^{\ast}\text{ and }j=+\infty,\\
1 & \text{if }i=j=+\infty.
\end{cases}
\]
\item Taking $f=\boldsymbol{1}_{\mathbb{N}^{\ast}}$ in the equality $\refpar{eq:trans_proba_N_i_f},$
we obtain
\[
N\left(i,\boldsymbol{1}_{\mathbb{N}^{\ast}}\right)=\begin{cases}
P_{y}\left(T^{1}_{y}<+\infty\right), & \text{if }i\in\mathbb{N}^{\ast},\\
0, & \text{if }i=+\infty,
\end{cases}
\]
which leads, by substituting in $\refpar{eq:cond_exp_x_f_T_y_pp1_knA_T_y_p_withN},$
to the equality
\[
\mathbb{E}^{\mathscr{A}_{T^{p}_{y}}}_{x}\left[\boldsymbol{1}_{\left(T^{p+1}_{y}<+\infty\right)}\right]=\boldsymbol{1}_{\left(T^{p}_{y}<+\infty\right)}P_{y}\left(T^{1}_{y}<+\infty\right).
\]
It follows that, by taking the $\mathbb{E}_{x}-$mean in each member
of the previous equality, that\boxeq{
\[
P_{x}\left(T^{p+1}_{y}<+\infty\right)=P_{x}\left(T^{p}_{y}<+\infty\right)P_{y}\left(T^{1}_{y}<+\infty\right).
\]
}In particular, for $x=y,$
\[
P_{y}\left(T^{p+1}_{y}<+\infty\right)=P_{y}\left(T^{p}_{y}<+\infty\right)P_{y}\left(T^{1}_{y}<+\infty\right).
\]
If $P_{y}\left(T^{p+1}_{y}<+\infty\right)=1,$ we then have, for every
$p\in\mathbb{N}^{\ast},$ $P_{y}\left(T^{p}_{y}<+\infty\right)=1.$
Moreover, since the sequence of events $\left(T^{p}_{y}<+\infty\right)$
is nonincreasing and that $\left(N_{y}=+\infty\right)=\bigcap_{p\in\mathbb{N}^{\ast}}\left(T^{p}_{y}<+\infty\right),$\boxeq{
\[
P_{y}\left(N_{y}=+\infty\right)=\lim_{p\to+\infty}P_{y}\left(T^{p}_{y}<+\infty\right)=1.
\]
}Last, for every subset $D$ of $\overline{\mathbb{N}^{\ast}},$
under this hypothesis, we have, by the equality $\refpar{eq:T_pp1_y_on_T_y_p},$
\[
\mathbb{E}^{\mathscr{A}_{T^{p}_{y}}}_{y}\left[\boldsymbol{1}_{D}\left(T^{p+1}_{y}-T^{p}_{y}\right)\right]=\mathbb{E}^{\mathscr{A}_{T^{p}_{y}}}_{y}\left[\boldsymbol{1}_{\left(T^{p}_{y}<+\infty\right)}\boldsymbol{1}_{D}\left[\tau_{y}\left[\theta_{T^{p}_{y}}\left(X\right)\right]\right]\right].
\]
Hence, by the strong Markov property,
\begin{align*}
\mathbb{E}^{\mathscr{A}_{T^{p}_{y}}}_{y}\left[\boldsymbol{1}_{D}\left(T^{p+1}_{y}-T^{p}_{y}\right)\right] & =\boldsymbol{1}_{\left(T^{p}_{y}<+\infty\right)}\mathbb{E}_{X^{p}_{y}}\left[\boldsymbol{1}_{D}\left(\tau_{y}\left(X\right)\right)\right]\\
 & =\boldsymbol{1}_{\left(T^{p}_{y}<+\infty\right)}\mathbb{E}_{y}\left[\boldsymbol{1}_{D}\left(T^{1}_{y}\right)\right].
\end{align*}
Thus,
\[
\mathbb{E}^{\mathscr{A}_{T^{p}_{y}}}_{y}\left[\boldsymbol{1}_{D}\left(T^{p+1}_{y}-T^{p}_{y}\right)\right]=P_{y}\left(T^{1}_{y}\in D\right).
\]
It follows first that the \salgs $\mathscr{A}_{T^{p}_{y}}$ and $\sigma\left(T^{p+1}_{y}-T^{p}_{y}\right)$
are $P_{y}-$independent; since this holds for every $p,$ we can
easily deduce that the random variables $T^{p+1}_{y}-T^{p}_{y}$ are
$P_{y}-$independent. Moreover, taking the $\mathbb{E}_{y}-$mean
of both sides of the previous equality yields
\[
P_{y}\left[\left(T^{p+1}_{y}-T^{p}_{y}\right)\in D\right]=P_{y}\left(T^{1}_{y}\in D\right).
\]
Since $D$ is arbitrary, this proves that $T^{p+1}_{y}-T^{p}_{y}$
and $T^{1}_{y}$ follow the same law under $P_{y}.$
\end{itemize}
\end{proof}

\subsection{Law of the Number of Visits to a Point and First Hitting Time of
This Point}

We now explain how the number of visits to a point $y$ and the first
hitting time in this point relates.The chain starts from $x$ at time
0.

\begin{proposition}{Number of Visits to a Point During the Whole Process Lifetime}{mb_vis_at_point_proc_lif}

The law of the number $N_{y}$ of visits in $y$ during the whole
lifetime of a process is given---with the convention $0^{0}=1$---by
\begin{itemize}
\item If $x\neq y,$
\[
P_{x}\left(N_{y}=m\right)=\begin{cases}
P_{x}\left(T^{1}_{y}<+\infty\right)P_{y}\left(T^{1}_{y}=+\infty\right)\left[P_{y}\left(T^{1}_{y}<+\infty\right)\right]^{m-1}, & \text{if\,}m\in\mathbb{N}^{\ast},\\
P_{x}\left(T^{1}_{y}=+\infty\right), & \text{if\,}m=0.
\end{cases}
\]
\item If $x=y,$
\[
P_{x}\left(N_{y}=m\right)=P_{y}\left(T^{1}_{y}=+\infty\right)\left[P_{y}\left(T^{1}_{y}<+\infty\right)\right]^{m-1}\,\,\,\,\text{if }m\in\mathbb{N}^{\ast}.
\]
\end{itemize}
That is, if $0<P_{y}\left(T^{1}_{y}<+\infty\right)<1,$ then \textbf{the
law of $N_{y}$ under $P_{y}$ is the geometric law on $\mathbb{N}^{\ast}$
with parameter $P_{y}\left(T^{1}_{y}=+\infty\right).$}

\end{proposition}

\begin{proof}{}{}

The event $\left(N_{y}=m\right)=\left(\sum_{j\in\mathbb{N}}\boldsymbol{1}_{\left(X_{j}=y\right)}=m\right)$
is the set of trajectories that visit $y$ exactly $m$ times starting
from time 0. Now, for every $m\in\mathbb{N}^{\ast},$
\begin{multline*}
\left(\sum_{j\in\mathbb{N}^{\ast}}\boldsymbol{1}_{\left(X_{j}=y\right)}=m\right)=\left(T^{1}_{y}<+\infty\right)\cap\left(T^{2}_{y}<+\infty\right)\cap\cdots\\
\cdots\cap\left(T^{m}_{y}<+\infty\right)\cap\left(T^{m+1}_{y}-T^{m}_{y}=+\infty\right),
\end{multline*}
which, by equality $\refpar{eq:T_pp1_y_on_T_y_p},$ becomes
\begin{multline*}
\left(\sum_{j\in\mathbb{N}^{\ast}}\boldsymbol{1}_{\left(X_{j}=y\right)}=m\right)=\left(T^{1}_{y}<+\infty\right)\cap\left(T^{2}_{y}<+\infty\right)\cap\cdots\\
\cdots\cap\left(T^{m}_{y}<+\infty\right)\cap\left(\tau_{y}\left(\theta_{T^{m}_{y}}\left(X\right)\right)=+\infty\right).
\end{multline*}
By integrating with respect to $P_{x}$ and by noting that
\[
\left(T^{1}_{y}<+\infty\right)\cap\left(T^{2}_{y}<+\infty\right)\cap\cdots\cap\left(T^{m}_{y}<+\infty\right)\in\mathscr{A}_{T^{m}_{y}},
\]
we obtain, by conditioning with respect to the \salg $\mathscr{A}_{T^{m}_{y}},$
\begin{multline*}
P_{x}\left(\sum_{j\in\mathbb{N}^{\ast}}\boldsymbol{1}_{\left(X_{j}=y\right)}=m\right)=\mathbb{E}_{x}\left[\boldsymbol{1}_{\left(T^{1}_{y}<+\infty\right)\cap\left(T^{2}_{y}<+\infty\right)\cap\cdots\cap\left(T^{m}_{y}<+\infty\right)}\right.\\
\left.\times\mathbb{E}^{\mathscr{A}_{T^{m}_{y}}}_{x}\left(\boldsymbol{1}_{\left(T^{m}_{y}<+\infty\right)}\boldsymbol{1}_{\left(\tau_{y}\left(\theta_{T^{m}_{y}}\left(X\right)\right)=+\infty\right)}\right)\right].
\end{multline*}
Now, by application of the strong Markov property---computation already
seen---,
\begin{align*}
\mathbb{E}^{\mathscr{A}_{T^{m}_{y}}}_{x}\left(\boldsymbol{1}_{\left(T^{m}_{y}<+\infty\right)}\boldsymbol{1}_{\left(\tau_{y}\left(\theta_{T^{m}_{y}}\left(X\right)\right)=+\infty\right)}\right) & =\boldsymbol{1}_{\left(T^{m}_{y}<+\infty\right)}\mathbb{E}_{X_{T^{m}_{y}}}\left(\boldsymbol{1}_{\left(\tau_{y}\left(X\right)=+\infty\right)}\right)\\
 & =\boldsymbol{1}_{\left(T^{m}_{y}<+\infty\right)}P_{y}\left(T^{1}_{y}=+\infty\right),
\end{align*}
which, by substituting in the previous equality, yields
\begin{multline}
P_{x}\left(\sum_{j\in\mathbb{N}^{\ast}}\boldsymbol{1}_{\left(X_{j}=y\right)}=m\right)=\mathbb{E}_{x}\left[\boldsymbol{1}_{\left(T^{1}_{y}<+\infty\right)\cap\left(T^{2}_{y}<+\infty\right)\cap\cdots\cap\left(T^{m}_{y}<+\infty\right)}\right]\\
\times P_{y}\left(T^{1}_{y}=+\infty\right).\label{eq:P_xnb_visit_m_in_yafter_0}
\end{multline}
By taking into account the equality
\begin{align*}
\left(T^{m-1}_{y}<+\infty\right)\cap\left(T^{m}_{y}<+\infty\right) & =\left(T^{m-1}_{y}<+\infty\right)\cap\left(T^{m}_{y}-T^{m-1}_{y}<+\infty\right)\\
 & =\left(T^{m-1}_{y}<+\infty\right)\cap\left(\tau_{y}\left(\theta_{T^{m-1}_{y}}\left(X\right)\right)<+\infty\right),
\end{align*}
we obtain, by the same method of conditioning with respect to the
\salg $\mathscr{A}_{T^{m-1}_{y}},$ and then by application of the
strong Markov property,
\begin{align*}
 & \mathbb{E}_{x}\left[\boldsymbol{1}_{\left(T^{1}_{y}<+\infty\right)\cap\left(T^{2}_{y}<+\infty\right)\cap\cdots\cap\left(T^{m}_{y}<+\infty\right)}\right]\\
 & \qquad=\mathbb{E}_{x}\left[\boldsymbol{1}_{\left(T^{1}_{y}<+\infty\right)\cap\left(T^{2}_{y}<+\infty\right)\cap\cdots\cap\left(T^{m-1}_{y}<+\infty\right)}\mathbb{E}^{\mathscr{A}_{T^{m-1}_{y}}}_{x}\left(\boldsymbol{1}_{\left(\tau_{y}\left(\theta_{T^{m-1}_{y}}\left(X\right)\right)<+\infty\right)}\right)\right]\\
 & \qquad=\mathbb{E}_{x}\left[\boldsymbol{1}_{\left(T^{1}_{y}<+\infty\right)\cap\left(T^{2}_{y}<+\infty\right)\cap\cdots\cap\left(T^{m-1}_{y}<+\infty\right)}\mathbb{E}_{X_{T^{m-1}_{y}}}\left(\boldsymbol{1}_{\left(\tau_{y}\left(X\right)<+\infty\right)}\right)\right]\\
 & \qquad=\mathbb{E}_{x}\left[\boldsymbol{1}_{\left(T^{1}_{y}<+\infty\right)\cap\left(T^{2}_{y}<+\infty\right)\cap\cdots\cap\left(T^{m-1}_{y}<+\infty\right)}\right]P_{y}\left(T^{1}_{y}<+\infty\right),
\end{align*}
which, by substituting in $\refpar{eq:P_xnb_visit_m_in_yafter_0},$
yields
\begin{multline*}
P_{x}\left(\sum_{j\in\mathbb{N}^{\ast}}\boldsymbol{1}_{\left(X_{j}=y\right)}=m\right)=\mathbb{E}_{x}\left[\boldsymbol{1}_{\left(T^{1}_{y}<+\infty\right)\cap\left(T^{2}_{y}<+\infty\right)\cap\cdots\cap\left(T^{m-1}_{y}<+\infty\right)}\right]\\
\times P_{y}\left(T^{1}_{y}<+\infty\right)P_{y}\left(T^{1}_{y}=+\infty\right).
\end{multline*}
By backward induction, and by the same method, we then obtain the
equality
\begin{multline}
P_{x}\left(\sum_{j\in\mathbb{N}^{\ast}}\boldsymbol{1}_{\left(X_{j}=y\right)}=m\right)=P_{x}\left(T^{1}_{y}<+\infty\right)P_{y}\left(T^{1}_{y}=+\infty\right)\left[P_{y}\left(T^{1}_{y}<+\infty\right)\right]^{m-1}.\label{eq:P_x_to_visite_m_times_y}
\end{multline}

\begin{itemize}
\item If $x\neq y$ and $m\in\mathbb{N^{\ast}},$ then
\[
P_{x}\left(N=m\right)=P_{x}\left(\sum_{j\in\mathbb{N}}\boldsymbol{1}_{\left(X_{j}=y\right)}=m\right)=P_{x}\left(\sum_{j\in\mathbb{N}^{\ast}}\boldsymbol{1}_{\left(X_{j}=y\right)}=m\right),
\]
and the equality $\refpar{eq:P_x_to_visite_m_times_y}$ yields the
announced formula.
\item If $x\neq y$ and $m=0,$ then
\[
P_{x}\left(N_{y}=0\right)=P_{x}\left(T^{1}_{y}=+\infty\right).
\]
\item Finally, if $x=y,$ then $P_{y}\left(N_{y}=0\right)=0$ and, if $m\in\mathbb{N}^{\ast},$
\[
P_{y}\left(N_{y}=m\right)=P_{y}\left(\sum_{j\in\mathbb{N}}\boldsymbol{1}_{\left(X_{j}=y\right)}=m\right)=P_{x}\left(\sum_{j\in\mathbb{N}^{\ast}}\boldsymbol{1}_{\left(X_{j}=y\right)}=m-1\right),
\]
and the equality $\refpar{eq:P_x_to_visite_m_times_y}$ still yields
the announced result.
\end{itemize}
\end{proof}

We now study the time $T^{1}_{y}.$

\begin{denotations}{}{}

For $k\in\mathbb{N}^{\ast}$ and for every $x,y\in E,$ set
\[
F_{k}\left(x,y\right)=P_{x}\left(T^{1}_{y}=k\right)\,\,\,\,\text{and}\,\,\,\,F\left(x,y\right)=P_{x}\left(T^{1}_{y}<+\infty\right).
\]

\end{denotations}

\begin{proposition}{}{first_hitting_point}

The sequence of matrices $F_{k}$ is solution of the following iterative
system, for every $x,y\in E,$
\begin{equation}
\left\{ \begin{array}{lc}
F_{1}\left(x,y\right)=M\left(x,y\right),\\
F_{k}\left(x,y\right)=\sum_{z\in E\backslash\left\{ y\right\} }M\left(x,z\right)F_{k-1}\left(z,y\right), & \text{if\,\ensuremath{k\geqslant2}}
\end{array}\right.\label{eq:F_k(x_y)}
\end{equation}
Therefore, the matrix $F$ is solution of the matrix equation determined,
for every $x,y\in E,$ by
\begin{equation}
F\left(x,y\right)=M\left(x,y\right)+\sum_{z\in E\backslash\left\{ y\right\} }M\left(x,z\right)F\left(z,y\right)\label{eq:F_as_fixed_point}
\end{equation}

\end{proposition}

\begin{proof}{}{}
\begin{itemize}
\item We have $F_{1}\left(x,y\right)=P_{x}\left(T^{1}_{y}=1\right)=P_{x}\left(X_{1}=y\right)=M\left(x,y\right).$
\item Let $k\geqslant2.$ On $\left(T^{1}_{y}>1\right),$ we have 
\[
T^{1}_{y}=1+\tau_{y}\left(\theta_{1}\left(X\right)\right).
\]
Conditioning with respect to the \salg $\mathscr{A}_{1}$ and applying
the simple Markov property---note that $X_{1}$ is $\mathscr{A}_{1}-$measurable---we
obtain
\begin{align*}
F_{k}\left(x,y\right) & =\mathbb{E}_{x}\left[\boldsymbol{1}_{\left(X_{1}\neq y\right)}\mathbb{E}^{\mathscr{A}_{1}}_{x}\left(\boldsymbol{1}_{\left(1+\tau_{y}\left(\theta_{1}\left(X\right)\right)=k\right)}\right)\right]\\
 & =\mathbb{E}_{x}\left[\boldsymbol{1}_{\left(X_{1}\neq y\right)}\mathbb{E}_{X_{1}}\left(\boldsymbol{1}_{\left(1+\tau_{y}\left(X\right)=k\right)}\right)\right].
\end{align*}
Hence,
\[
F_{k}\left(x,y\right)=\mathbb{E}_{x}\left[\boldsymbol{1}_{\left(X_{1}\neq y\right)}P_{X_{1}}\left(T^{1}_{y}=k-1\right)\right]=\sum_{z\in E\backslash\left\{ y\right\} }P_{x}\left(X_{1}=z\right)P_{z}\left(T^{1}_{y}=k-1\right),
\]
which proves $\refpar{eq:F_k(x_y)}.$
\item Finally,
\[
F\left(x,y\right)=P_{x}\left(T^{1}_{y}<+\infty\right)=\sum_{k\in\mathbb{N}^{\ast}}P_{x}\left(T^{1}_{y}=k\right)=\sum_{k\in\mathbb{N}^{\ast}}F_{k}\left(x,y\right).
\]
By $\refpar{eq:F_k(x_y)},$ we get
\[
F\left(x,y\right)=M\left(x,y\right)+\sum_{k\geqslant2}\left[\sum_{z\in E\backslash\left\{ y\right\} }M\left(x,z\right)F_{k-1}\left(z,y\right)\right].
\]
Permuting the sums---all terms are nonnegative---yields $\refpar{eq:F_as_fixed_point}.$ 
\end{itemize}
\end{proof}

We now reformulate the results of Proposition $\ref{pr:mb_vis_at_point_proc_lif}$
using the matrix $F,$ and we state, without proof, some immediate
consequences.

\begin{proposition}{}{reformulation_m_visit_withF}

With the previous notation,
\begin{itemize}
\item If $x\neq y,$
\[
P_{x}\left(N_{y}=m\right)=\begin{cases}
1-F\left(x,y\right), & \text{if\,}m=0,\\
F\left(x,y\right)\left[1-F\left(y,y\right)\right]\left[F\left(y,y\right)\right]^{m-1}, & \text{if\,}m\in\mathbb{N}^{\ast}.
\end{cases}
\]
\item If $x=y$---with the convention $0^{0}=1$---for $m\in\mathbb{N}^{\ast},$
\[
P_{y}\left(N_{y}=m\right)=\left[1-F\left(y,y\right)\right]\left[F\left(y,y\right)\right]^{m-1}.
\]
\item We have the following alternative
\[
P_{y}\left(N_{y}=+\infty\right)=\begin{cases}
1, & \text{if }F\left(y,y\right)<1,\\
0, & \text{if }F\left(y,y\right)=1.
\end{cases}
\]

\begin{itemize}
\item If $F\left(y,y\right)=1,$ then $P_{y}\left(N_{y}=+\infty\right)=1,$
hence $E_{y}\left(N_{y}\right)=+\infty.$
\item If $0<F\left(y,y\right)<1,$ then under $P_{y},$ $N_{y}$ follows
the geometric law on $\mathbb{N}^{\ast}$ with parameter $1-F\left(y,y\right).$
\item If $F\left(y,y\right)=0,$ then $P_{y}-$almost surely $N_{y}=1.$
\end{itemize}
In particular, \textbf{the mean number\index{mean number of visits}
$\mathbb{E}_{y}\left(N_{y}\right)$ of visits to $y$ for the chain
starting from $y$ at time 0} is
\[
\mathbb{E}_{y}\left(N_{y}\right)=\begin{cases}
\dfrac{1}{1-F\left(y,y\right)}, & \text{if }F\left(y,y\right)<1,\\
+\infty, & \text{if }F\left(y,y\right)=1.
\end{cases}
\]

\end{itemize}
\end{proposition}

\begin{definition}{Potential Matrix}{}

The matrix $R$---with values in $\overline{\mathbb{N}}$---defined,
for every $x,y\in E,$ by $R\left(x,y\right)=\mathbb{E}_{x}\left(N_{y}\right),$
that is the mean number of visits to $y$ by the chain starting at
$x$ at time 0, is called the \textbf{\index{potential matrix}potential
matrix} of the chain.

\end{definition}

From Proposition $\ref{pr:reformulation_m_visit_withF},$ we deduce
the following corollary.

\begin{corollary}{Potential of A Path}{potent_path}

With the conventions $\dfrac{1}{0}=+\infty$ and $0\cdot\infty=0,$
the potential matrix satisfies\boxeq{
\begin{equation}
R\left(x,y\right)=\begin{cases}
\dfrac{1}{1-F\left(y,y\right)}, & \text{if }x=y,\\
F\left(x,y\right)R\left(y,y\right), & \text{if }x\neq y.
\end{cases}\label{eq:potential_mat_pr}
\end{equation}
}

\end{corollary}

\begin{remark}{}{}

In practice, it is often easier to compute $R$ first---we shall
give a computation method below---and then deduce $F.$ The following
proposition proves that $R$ satisfies a matrix equation. In particular,
when $E$ is finite, this allows one to compute $R$ once its infinite
entries have been identified.

\end{remark}

\begin{proposition}{}{}

The potential matrix $R$ satisfies\boxeq{
\[
R=\sum^{+\infty}_{n=0}M^{n},
\]
} in the sense that for every $x,y\in E,$ we have the equality,
in $\overline{\mathbb{R}}^{+},$ 
\[
R\left(x,y\right)=\sum^{+\infty}_{n=0}M^{n}\left(x,y\right).
\]
 Moreover, $R$ is solution of the matrix equation\boxeq{
\begin{equation}
R\left(I-M\right)=\left(I-M\right)R=I,\label{eq:R_sol_mat_aq}
\end{equation}
}where $I$ is the identity matrix---that is $I\left(x,y\right)=1$
if $x=y,$ 0 otherwise.

In particular, if $E$ is finite and if $R$ has only finite entries,
then $1-M$ is invertible and $R=\left(1-M\right)^{-1}.$

\end{proposition}

\begin{proof}{}{}

By monotonic convergence,
\[
R\left(x,y\right)=\sum^{+\infty}_{n=0}\mathbb{E}_{x}\left(\boldsymbol{1}_{\left(X_{n}=y\right)}\right)=\sum^{+\infty}_{n=0}M^{n}\left(x=y\right).
\]
Hence,
\[
RM=MR=\sum^{+\infty}_{n=0}M^{n}=R-1,
\]
which yields the equality $\refpar{eq:R_sol_mat_aq}.$

\end{proof}

\section{State Classification}\label{sec:State-Classification}

In this section, we consider a process $X$ such that, for every $x\in E,$
it is a homogeneous Markov chain on the underlying filtered probabilized
space $\left(\Omega,\mathscr{A},\left(\mathscr{A}_{n}\right)_{n\in\mathbb{N}},P_{x}\right),$
with initial law $\delta_{x}$ and transition matrix $M.$ We now
classify the states of $E$ according to how frequently they are visited
by the trajectories of $X.$

\subsection{Communication. Periodicity}

\begin{definition}{}{}

Let $B\subset E.$ We say that a point $x\in E$ \textbf{leads to
$B$\mindex{point!leads to}\mindex{state!leads to}} if $P_{x}\left(T_{B}<+\infty\right)>0.$
We denote this relation by $x\to B.$ In particular, if $B=\left\{ y\right\} ,$
with $y\in E,$ we say that $x$ \textbf{leads to $y$} and we write
$x\to y.$

\end{definition}

\begin{proposition}{}{}The \textbf{reachability relation\index{reachability relation}}
$x\to y$ is transitive. Moreover $x$ leads to $y$ if and only if
there exists $n\in\mathbb{N}^{\ast}$ such that $M^{n}\left(x,y\right)>0.$

\end{proposition}

\begin{proof}{}{}
\begin{itemize}
\item Suppose that $x\to y$ and $y\to z.$ The set of trajectories that
visit $z$ after having visited $y$ is contained in the set of trajectories
that visit $z,$ which gives the inclusion of events
\[
\left(T_{y}<+\infty\right)\cap\left(\tau_{z}\left(\theta_{T_{y}}\left(X\right)<+\infty\right)\right)\subset\left(T_{z}<+\infty\right).
\]
Therefore
\[
P_{x}\left[\left(T_{y}<+\infty\right)\cap\left(\tau_{z}\left(\theta_{T_{y}}\left(X\right)<+\infty\right)\right)\right]\leqslant P_{x}\left(T_{z}<+\infty\right).
\]
Conditioning with respect to $\mathscr{A}_{T_{y}}$ and using the
strong Markov property---computation now classical---lead to the
following equalities
\begin{align*}
 & P_{x}\left[\left(T_{y}<+\infty\right)\cap\left(\tau_{z}\left(\theta_{T_{y}}\left(X\right)<+\infty\right)\right)\right]\\
 & \qquad\qquad=\mathbb{E}_{x}\left[\mathbb{E}^{\mathscr{A}_{T_{y}}}_{x}\left(\boldsymbol{1}_{\left(T_{y}<+\infty\right)}\boldsymbol{1}_{\left(\tau_{z}\left[\theta_{T_{y}}\left(X\right)\right]<+\infty\right)}\right)\right]\\
 & \qquad\qquad=\mathbb{E}_{x}\left[\boldsymbol{1}_{\left(T_{y}<+\infty\right)}\mathbb{E}_{X_{T_{y}}}\left(\boldsymbol{1}_{\left(T_{z}<+\infty\right)}\right)\right]\\
 & \qquad\qquad=P_{x}\left(T_{y}<+\infty\right)P_{y}\left(T_{z}<+\infty\right),
\end{align*}
which proves that
\[
0<P_{x}\left(T_{y}<+\infty\right)P_{y}\left(T_{z}<+\infty\right)\leqslant P_{x}\left(T_{z}<+\infty\right),
\]
and thus that $x\to z.$
\item If $x\to y,$ then, since $\left(T_{y}<+\infty\right)=\bigcup_{n\in\mathbb{N}^{\ast}}\left(X_{n}=y\right),$
\[
0<P_{x}\left(T_{y}<+\infty\right)\leqslant\sum_{n\in\mathbb{N}^{\ast}}P_{x}\left(X_{n}=y\right)=\sum_{n\in\mathbb{N}^{\ast}}M^{n}\left(x,y\right),
\]
which proves there exists $n\in\mathbb{N}^{\ast}$ such that $M^{n}\left(x,y\right)>0.$\\
Conversely, if for some $n.$
\[
0<M^{n}\left(x,y\right)=P_{x}\left(X_{n}=y\right)\leqslant P_{x}\left(T_{y}<+\infty\right),
\]
and thus $x$ leads to $y.$
\end{itemize}
\end{proof}

We deduce from this transitive relation, an equivalence relation by
symmetrizing.

\begin{definition}{Communitation. Communication Classes. Irreducible Chain}{}

We say that $x$ \textbf{communicates with\index{communicates with}\mindex{state!communicates with}\mindex{point!communicates with}}
$y$ if $x$ leads to $y$ and $y$ leads to $x,$ or if $x$ and
$y$ coincide. We denote this relation $x\leftrightarrow y.$

The \textbf{relation of communication\index{relation of communication}}
is an equivalence relation. Its equivalence classes are called \textbf{communication
classes\index{communication classes}} or \textbf{irreducible classes}\index{irreducible classes}.
In particular, if there is only one communication class---that is,
if every pair of points communicates---the \textbf{chain} is called
\textbf{irreducible}\index{irreducible chain}.

\end{definition}

It is customary to associate with a Markov chain of transition matrix
$M$ a directed \textbf{\index{graph}graph} whose vertices are the
points of $E,$ vertices being linked it they effectively communicate,
by arrows that indicate the sense of communication.

For instance, the following graph

\begin{center}
\begin{tikzpicture}[
  scale=0.75,
  transform shape,
  >=Stealth,
  every node/.style={
    circle,
    draw=black,
    fill=white,
    minimum size=10mm,
    inner sep=0pt
  },
  every edge/.style={
    ->,
    draw=black,
    line width=0.75pt
  }
]

  % Nodes
  \node (x) at (0,0) {$x$};
  \node (y) at (4,0) {$y$};
  \node (z) at (8,0) {$z$};
  \node (t) at (12,0) {$t$};

  % Edges (drawn normally, visible)
  \draw (x) edge[loop left] (x);
  \draw (x) edge[bend left=25] (y);
  \draw (y) edge[bend left=25] (z);
  \draw (z) edge[bend left=25] (y);
  \draw (z) edge[bend left=25] (t);
  \draw (t) edge[loop right] (t);

\end{tikzpicture}
\end{center}summarizes the fact that $x$ leads to itself, $x$ leads to $y,$
$y$ communicates with $z,$ $z$ leads to $t$ and $t$ leads to
itself. The classes of communication are $\left\{ x\right\} ,$$\left\{ y,z\right\} $
and $\left\{ t\right\} .$ On this graph, we could also mention the
transition probabilities from a state $x$ to a state $y$ from time
0 to time 1, that is the probabilities $M\left(x,y\right),$ but it
has less interest.

\begin{example}{}{}

1. Let $E=\left\{ 1,2,3,4,5\right\} $ and let $M$ be the matrix
\[
\left(\begin{array}{ccccc}
\dfrac{2}{5} & \dfrac{3}{5} & 0 & 0 & 0\\
\dfrac{1}{4} & \dfrac{3}{4} & 0 & 0 & 0\\
0 & 0 & 0 & 1 & 0\\
0 & 0 & \dfrac{1}{3} & 0 & \dfrac{2}{3}\\
0 & 0 & 0 & 1 & 0
\end{array}\right).
\]

The associated graph is

\begin{center}\begin{tikzpicture}[
  scale=0.75,
  transform shape,
  >=Stealth,
  every node/.style={
    circle,
    draw=black,
    minimum size=10mm,
    inner sep=0pt
  },
  edge/.style={
    ->,
    draw=black,
    line width=0.75pt
  }
]

% Nodes
\node (1) at (0,0) {1};
\node (2) at (3,0) {2};

\node (3) at (7,0) {3};
\node (4) at (10,0) {4};
\node (5) at (13,0) {5};

% Circular self-loops (same width as other edges)
\draw[edge] (1) to[out=210,in=150,looseness=8] (1);
\draw[edge] (2) to[out=30,in=330,looseness=8] (2);

% Between 1 and 2 (bent to separate arrows)
\draw[edge] (1) to[bend left=25] (2);
\draw[edge] (2) to[bend left=25] (1);

% Class {3,4,5} with bending for opposite directions
\draw[edge] (3) to[bend left=25] (4);
\draw[edge] (4) to[bend left=25] (3);

\draw[edge] (4) to[bend left=25] (5);
\draw[edge] (5) to[bend left=25] (4);

\end{tikzpicture}
\end{center}

There are exactly two communication classes $E_{1}=\left\{ 1,2\right\} $
and $E_{2}=\left\{ 3,4,5\right\} .$ Thus, two underlying homogeneous
Markov chains are naturally identified, taking values respectively
in $E_{1}$ and $E_{2}.$ Their transition matrices are the corresponding
sub-matrices of $M,$ namely $M_{1}$ and $M_{2},$ respectively given
by
\[
M_{1}=\left(\begin{array}{cc}
\dfrac{2}{5} & \dfrac{3}{5}\\
\dfrac{1}{4} & \dfrac{3}{4}
\end{array}\right)\,\,\,\,\text{and}\,\,\,\,M_{2}=\left(\begin{array}{ccc}
0 & 1 & 0\\
\dfrac{1}{3} & 0 & \dfrac{2}{3}\\
0 & 1 & 0
\end{array}\right).
\]
2. Let $E=\left\{ 1,2,3,4\right\} $ and let $M$ be the matrix
\[
\left(\begin{array}{cccc}
0 & 0 & \dfrac{1}{3} & \dfrac{2}{3}\\
0 & 0 & \dfrac{1}{2} & \dfrac{1}{2}\\
\dfrac{1}{4} & \dfrac{3}{4} & 0 & 0\\
\dfrac{1}{2} & \dfrac{1}{2} & 0 & 0
\end{array}\right).
\]
The associated graph is

\begin{center}\begin{tikzpicture}[
  scale=0.75,
  transform shape,
  >=Stealth,
  every node/.style={
    circle,
    draw=black,
    minimum size=10mm,
    inner sep=0pt
  },
  edge/.style={
    ->,
    draw=black,
    line width=0.75pt
  }
]

% Nodes (two columns)
\node (1) at (0,1.5) {1};
\node (2) at (0,-1.5) {2};

\node (3) at (4,1.5) {3};
\node (4) at (4,-1.5) {4};

% Straight arrows
\draw[edge] (1) -- (3);
\draw[edge] (1) -- (4);

\draw[edge] (2) -- (3);
\draw[edge] (2) -- (4);

\draw[edge] (3) -- (1);
\draw[edge] (3) -- (2);

\draw[edge] (4) -- (1);
\draw[edge] (4) -- (2);

\end{tikzpicture}
\end{center}

We see that every state communicates with every other: \textbf{the
chain is irreducible}. Nevertheless two sub-classes appear $C_{1}=\left\{ 1,2\right\} $
and $C_{2}=\left\{ 3,4\right\} $ such that if $X_{n}\in C_{1}$ then
$X_{n+1}\in C_{2}$ and if $X_{n}\in C_{2}$ then $X_{n+1}\in C_{1}.$
These are called \textbf{cyclic classes}\index{cyclic class}. This
naturally leads to the definition of the period of a state. 

\end{example}

\begin{definition}{Period of a State. Aperiodic State}{}

Let $x\in E.$ The largest integer $d$ such that we have the inclusion
\[
\left\{ n\in\mathbb{N}^{\ast}:M^{n}\left(x,x\right)>0\right\} \subset d\mathbb{N}^{\ast}
\]
is called the \textbf{period\index{period}} of $x,$ and is denoted
$d\left(x\right).$ This is the greatest common divisor (GCD) of the
set $\left\{ n\in\mathbb{N}^{\ast}:M^{n}\left(x,x\right)>0\right\} .$
If this set is empty, we set $d\left(x\right)=0.$ If $d\left(x\right)=1,$
we say that $x$ is \textbf{\index{aperiodic}aperiodic}.

\end{definition}

\begin{proposition}{Class Period. Aperiodic Class. Aperiodic Markov Chain}{class_per_aper}

Let $C$ be a communication class. Every elements of $C$ have the
same period, denoted $d\left(C\right),$ and called \textbf{period
of the class\mindex{periode!class}\mindex{class!periode} $C.$}
If $d\left(C\right)=1,$ the \textbf{class} is said to be\textbf{\mindex{aperiodic!class}\mindex{class!aperiodic}
aperiodic}. An irreducible homogeneous Markov \textbf{chain} with
at least one aperiodic state is said \textbf{\mindex{Markov chain!irreducible homogeneous!aperiodic}aperiodic}. 

\end{proposition}

\begin{proof}{}{}

Let $x,y\in C.$ Since $x$ and $y$ communicate, there exist $k$
and $l\in\mathbb{N}^{\ast}$ such that $M^{k}\left(x,y\right)>0$
and $M^{l}\left(y,x\right)>0.$ It follows that
\[
M^{k+l}\left(x,x\right)\geqslant M^{k}\left(x,y\right)M^{l}\left(y,x\right)>0.
\]
Hence, from the one hand, $d\left(x\right)\geqslant1,$ and on the
other hand that $k+l\equiv0\left[d\left(x\right)\right].$ 

Note that for every $n$ that is not a multiple of $d\left(x\right),$
this is the same for $n+k+l$ and therefore $M^{n+k+l}\left(x,x\right)=0.$
We deduce that
\[
0=M^{n+k+l}\left(x,x\right)\geqslant M^{k}\left(x,y\right)M^{n}\left(y,y\right)M^{l}\left(y,x\right)\geqslant0
\]
and hence, that $M^{n}\left(y,y\right)=0.$ By contraposition, we
just proved that if $M^{n}\left(y,y\right)>0,$ $d\left(y\right)$
is a multiple of $d\left(x\right),$ and thus that $d\left(y\right)\geqslant d\left(x\right).$
By symmetry, we also obtain $d\left(x\right)\geqslant d\left(y\right),$
which proves the equality $d\left(x\right)=d\left(y\right).$

\end{proof}

Let $C$ be a communication class of period $d>1,$ and let $x_{0}\in C.$
Every point $x\in C$ communicates with $x_{0}.$ Let $k\in\mathbb{N}^{\ast}$
be the smallest integer such that $M^{k}\left(x,x_{0}\right)>0.$
For every $n\in\mathbb{N}^{\ast},$
\[
M^{n+k}\left(x_{0},x_{0}\right)\geqslant M^{n}\left(x_{0},x\right)M^{k}\left(x,x_{0}\right),
\]
which proves that, for every $n\in\mathbb{N}^{\ast}$ such that $M^{n}\left(x_{0},x\right)>0,$
we have $M^{n+k}\left(x_{0},x_{0}\right)>0.$ Since $x_{0}$ has period
$d,$ it follows that $n+k\equiv0\left[d\right].$ Hence, there exists
a unique integer $j\in\left\llbracket 0,d-1\right\rrbracket $---namely,
the remainder of the Euclidean division of $-k$ by $d$---such that
\[
M^{n}\left(x_{0},x\right)>0\Longrightarrow n\equiv j\left[d\right].
\]
We then define the \textbf{cyclic classes} $C_{j},\,j\in\left\llbracket 0,d-1\right\rrbracket ,$
of $C$ by
\[
C_{j}=\left\{ y\in C:\,M^{n}\left(x_{0},y\right)>0\Rightarrow n\equiv j\left[d\right]\right\} .
\]
Equivalently, $y\in C_{j}$ if and only if
\[
\left\{ n\in\mathbb{N}^{\ast}:\,M^{n}\left(x_{0},y\right)>0\right\} \subset j+d\mathbb{N}^{\ast}.
\]
The sets $C_{j},\,j\in\left\llbracket 0,d-1\right\rrbracket ,$ form
a partition of $C.$ Moreover, if $x\in C_{j}$ and $y$ are such
that $M\left(x,y\right)>0,$ then $y\in C_{j+1\left[d\right]}.$ Indeed,
let $n$ be such that $M^{n}\left(x_{0},y\right)>0.$ Then $n$ is
congruent to $j$ modulo $d,$ so $n+1\equiv j+1\left[d\right],$
and by what precedes this implies $y\in C_{j+1}.$

For every $n\in\mathbb{N}^{\ast},$ $X_{j+nd}\in C_{j}$ $P_{x_{0}}-$almost
surely, and the subchain $\left(X_{j+nd}\right)_{n\in\mathbb{N}^{\ast}},$
starting from $x_{0}\in C_{0}$ at time 0, is a homogeneous Markov
chain with state space $C_{j},$ transition matrix $\left(M^{d}\left(x,y\right)\right)_{x,y\in C_{j}},$
and it is irreducible and aperiodic.

To illustrate these notions, we return to the Ehrenfest diffusion
model.

\begin{example}{Ehrenfest Heat Diffusion Model, Continued}{ehrenfest_model_2}

We consider the Ehrenfest model, described in its urn formulation---see
Example $\ref{ex:ehrenfest_model}$ whose notation we keep. 

Let $X_{n}$ denote the number of red balls contained in the urn at
time $n.$ Since each draw is uniform, the process $X=\left(X_{n}\right)_{n\in\mathbb{N}}$
is a homogeneous Markov chain taking values in the integer interval
$\left\llbracket 0,m\right\rrbracket ,$ with transition matrix $M$
defined as follows.
\begin{itemize}
\item For $1\leqslant k\leqslant m-1,$
\begin{itemize}
\item $M\left(k,k+1\right)=1-\dfrac{k}{m},$ 
\item $M\left(k,k-1\right)=\dfrac{k}{m},$ 
\item and $M\left(k,l\right)=0$ if $l\neq k-1$ or $k+1.$
\end{itemize}
\item At the boundaries,
\begin{equation}
M\left(0,1\right)=1,\,\,\,\,M\left(m,m-1\right)=1.\label{eq:ehrenfest_border_cond}
\end{equation}
\end{itemize}
Equivalently, the transition matrix $M$ is given by, for every $k\in\left\llbracket 0,m\right\rrbracket ,$
\[
M\left(k,k+1\right)=p_{k},\,\,\,\,M\left(k,k-1\right)=q_{k},
\]
where
\[
p_{k}=1-\dfrac{k}{m},\,\,\,\,q_{k}=\dfrac{k}{m}.
\]

This is clear that every states of $E$ communicates. The chain is
thus irreducible. Morever, 0 is aperiodic. This is thus the same for
every states of $E.$ Hence, the \textbf{Ehrenfest chain is aperiodic
irreducible}.

\end{example}

\begin{remark}{}{}

This model is a particular case of \textbf{birth-death processes\index{birth-death process}}\mindex{process!birth-death}---see
Exercise $\ref{exo:exercise17.4}.$ Here the states $0$ and $m$
are \textbf{reflecting} \textbf{barriers}\mindex{state!barrier!reflecting}\index{reflecting barriers},
that is, verify the conditions $\refpar{eq:ehrenfest_border_cond}.$

\end{remark}

\subsection{Recurrence}

\begin{definition}{Recurrent, Null Recurrent,  Positive Recurrent, Transitory State}{}

A state $x$ is
\begin{itemize}
\item \textbf{Recurrent\mindex{state!recurrent}} if $P_{x}\left(T^{1}_{x}<+\infty\right)=1.$
\item \textbf{\mindex{state!null recurrent}Null recurrent} if it is recurrent
and if $\mathbb{E}_{x}\left(T^{1}_{x}\right)=+\infty.$
\item \textbf{Positive recurrent\mindex{state!positive recurrent}} if it
is recurrent and if $\mathbb{E}_{x}\left(T^{1}_{x}\right)<+\infty.$
\item \textbf{Transient\mindex{state!transient}} if it is not recurrent,
that is if
\[
P_{x}\left(T^{1}_{x}<+\infty\right)<1.
\]
\end{itemize}
\end{definition}

\begin{remark}{}{}

A state $x$ is null recurrent if the chain, starting from $x,$ returns
to $x$ almost surely in finite time, but ``very slowly''. The reason
for this terminology will become clear when we study the problem of
existence of invariant probabilities---see Theorem $\ref{th:positive_recurrence_crit}$.

\end{remark}

The following lemma is preliminary to the state classification theorem.

\begin{lemma}{Probability of Multiple Visits. Infinite Number of Visits. Potential To Return}{rec_vis_and_proba}

For every $x\in E,$

\textbf{(a) Multiple returns}

For every $p\in\mathbb{N}^{\ast},$ 
\[
P_{x}\left(T^{p}_{x}<+\infty\right)=\left[P_{x}\left(T^{1}_{x}<+\infty\right)\right]^{p}.
\]

\textbf{(b) Infinite number of visits}

The set $R_{x}$ of trajectories that visit $x$ an infinite number
of times, is defined by
\[
R_{x}=\limsup_{n\to+\infty}\left(X_{n}=x\right)=\left(N_{x}=+\infty\right),
\]
and is equal to $\bigcap_{p\in\mathbb{N}^{\ast}}\left(T^{p}_{x}<+\infty\right)$
and verifies
\[
P_{x}\left(R_{x}\right)=\lim_{p\to+\infty}\searrow\left[P_{x}\left(T^{1}_{x}<+\infty\right)\right]^{p}.
\]

\textbf{(c) Potential on the diagonal}

The potential $R\left(x,x\right),$ that is the mean number of visits
in $x$ when the chain starts from $x$ at time 0, is
\[
R\left(x,x\right)=\sum^{+\infty}_{p=0}\left[P_{x}\left(T^{1}_{x}<+\infty\right)\right]^{p}.
\]

\end{lemma}

\begin{proof}{}{}

(a) This is a corollary of Proposition $\ref{pr:hitt_time_state}.$
We nonetheless give a direct proof. Since
\[
\left(T^{p+1}_{x}<+\infty\right)\subset\left(T^{p}_{x}<+\infty\right),
\]
the relation $\refpar{eq:T_pp1_y_on_T_y_p}$ allows to write---by
conditioning with respect to $\mathscr{A}_{T^{p}_{x}},$ and by applying
the strong Markov property---the sequence of equalities
\begin{align*}
P_{x}\left(T^{p+1}_{x}<+\infty\right) & =\mathbb{E}_{x}\left[\mathbb{E}^{\mathscr{A}_{T^{p}_{x}}}_{x}\left(\boldsymbol{1}_{\left(T^{p}_{x}<+\infty\right)}\boldsymbol{1}_{\left(\tau^{1}_{x}\left[\theta_{T^{p}_{x}}\left(X\right)\right]<+\infty\right)}\right)\right]\\
 & =\mathbb{E}_{x}\left[\boldsymbol{1}_{\left(T^{p}_{x}<+\infty\right)}\mathbb{E}_{X_{T^{p}_{x}}}\left(\boldsymbol{1}_{\left(\tau^{1}_{x}\left(X\right)<+\infty\right)}\right)\right].
\end{align*}
Hence, since $\tau^{1}_{x}\left(X\right)=T^{1}_{x}$ and since $X_{T^{p}_{x}}=x,$
\[
P_{x}\left(T^{p+1}_{x}<+\infty\right)=P_{x}\left(T^{p}_{x}<+\infty\right)P_{x}\left(T^{1}_{x}<+\infty\right),
\]
which yields the result by iteration.

(b) By definition of $R_{x}$ and of visiting times in $x,$ we have
the equality 
\[
R_{x}=\bigcap_{p\in\mathbb{N}^{\ast}}\left(T^{p}_{x}<+\infty\right).
\]
It is then enough to note that the sequence of sets $\left(T^{p}_{x}<+\infty\right)$
decreases to obtain the result.

(c) By definition,
\[
R\left(x,x\right)=\sum_{n\in\mathbb{N}}\mathbb{E}_{x}\left(\boldsymbol{1}_{\left(X_{n}=x\right)}\right).
\]
Set $T^{0}_{x}=0.$ The sequence of times $T^{p}_{x}$ is increasing.
Moreover, for every $p,$ we have $p\leqslant T^{p}_{x},$ which implies
that $\lim_{p\to+\infty}T^{p}_{x}=+\infty.$ We partition $\mathbb{N}$
into random intervals $\left[T^{p}_{x},T^{p+1}_{x}\right[,\,p\in\mathbb{N},$
and write
\[
R\left(x,x\right)=\sum_{p\in\mathbb{N}}\mathbb{E}_{x}\left(\sum_{n\in\left[T^{p}_{x},T^{p+1}_{x}\right[}\boldsymbol{1}_{\left(X_{n}=x\right)}\right).
\]
By noting that the interval $\left[T^{p}_{x},T^{p+1}_{x}\right[$
is empty as soon as $T^{p}_{x}=+\infty,$ that by definition of the
times $T^{p}_{x},$ we have $\boldsymbol{1}_{\left(X_{n}=x\right)}=0$
for every $n\in\left]T^{p}_{x},T^{p+1}_{x}\right[,$ and that $\boldsymbol{1}_{\left(X_{T^{p}_{x}}=x\right)}=1$
on $\left(T^{p}_{x}<+\infty\right),$ we obtain
\begin{align*}
R\left(x,x\right) & =\sum_{p\in\mathbb{N}}\mathbb{E}_{x}\left(\boldsymbol{1}_{\left(T^{p}_{x}<+\infty\right)}\right)=\sum_{p\in\mathbb{N}}P_{x}\left(T^{p}_{x}<+\infty\right)\\
 & =\sum_{p\in\mathbb{N}}\left[P_{x}\left(T^{1}_{x}<+\infty\right)\right]^{p}.
\end{align*}

\end{proof}

\begin{theorem}{State Classification}{state_classification}We have
the following alternative:

1. $x$ is recurrent; in this case $P_{x}\left(R_{x}\right)=1$ and
$R\left(x,x\right)=+\infty.$

2. $x$ is transient; in this case $P_{x}\left(R_{x}\right)=0$ and
$R\left(x,x\right)<+\infty.$

Moreover, if $x$ is recurrent and if $x$ leads to $y,$ then $y$
leads to $x,$ $y$ is recurrent, and $P_{y}\left(T^{1}_{x}<+\infty\right)=1.$

\end{theorem}

\begin{proof}{}{}

The alternative follows immediately from Lemma $\ref{lm:rec_vis_and_proba}.$
We now prove the last statement. Suppose that $x$ is recurrent and
$x$ leads to $y.$ If the chain visits $y$ after visiting $x,$
and then never returns to $x,$ this implies that the the total number
of visits to $x$ is finite. Hence, we have the inclusion of sets
\begin{equation}
\left(T^{1}_{x}<+\infty\right)\cap\left(\tau^{1}_{x}\left[\theta_{T^{1}_{x}}\left(X\right)\right]<+\infty\right)\cap\left(\tau^{1}_{x}\left[\theta_{\tau^{1}_{y}\left[\theta_{T^{1}_{x}}\left(X\right)\right]}\left(X\right)\right]=+\infty\right)\subset R^{c}_{x}\label{eq:subset_inclusion_chain_go_finite_time_inx}
\end{equation}
The two last sets of the first term makes appear the functionals of
the future of the process after the time $T^{1}_{x}.$ Hence, conditioning
with respect to $\mathscr{A}_{T^{1}_{x}},$ applying the strong Markov
property and taking into account the fact that $X_{T^{1}_{x}}=x,$
gives
\begin{align*}
 & P_{x}\left[\left(T^{1}_{x}<+\infty\right)\cap\left(\tau^{1}_{x}\left[\theta_{T^{1}_{x}}\left(X\right)\right]<+\infty\right)\cap\left(\tau^{1}_{x}\left[\theta_{\tau^{1}_{y}\left[\theta_{T^{1}_{x}}\left(X\right)\right]}\left(X\right)\right]=+\infty\right)\right]\\
 & \qquad\qquad=\mathbb{E}_{x}\left[\boldsymbol{1}_{\left(T^{1}_{x}<+\infty\right)}\mathbb{E}_{x}\left(\boldsymbol{1}_{\left(T^{1}_{y}<+\infty\right)}\boldsymbol{1}_{\left(\tau^{1}_{x}\left[\theta_{T^{1}_{y}}\left(X\right)\right]=+\infty\right)}\right)\right]\\
 & \qquad\qquad=P_{x}\left(T^{1}_{x}<+\infty\right)\mathbb{E}_{x}\left(\boldsymbol{1}_{\left(T^{1}_{y}<+\infty\right)}\boldsymbol{1}_{\left(\tau^{1}_{x}\left[\theta_{T^{1}_{y}}\left(X\right)\right]=+\infty\right)}\right).
\end{align*}
By now conditioning with respect to $\mathscr{A}_{T^{1}_{y}},$ and
using the strong Markov property together with $X_{T^{1}_{y}}=y,$
yields
\begin{align*}
 & P_{x}\left[\left(T^{1}_{x}<+\infty\right)\cap\left(\tau^{1}_{x}\left[\theta_{T^{1}_{x}}\left(X\right)\right]<+\infty\right)\cap\left(\tau^{1}_{x}\left[\theta_{\tau^{1}_{y}\left[\theta_{T^{1}_{x}}\left(X\right)\right]}\left(X\right)\right]=+\infty\right)\right]\\
 & \qquad\qquad=P_{x}\left(T^{1}_{x}<+\infty\right)\mathbb{E}_{x}\left(\boldsymbol{1}_{\left(T^{1}_{y}<+\infty\right)}\mathbb{E}_{y}\left[\boldsymbol{1}_{\left(\tau^{1}_{x}\left(X\right)=+\infty\right)}\right]\right)\\
 & \qquad\qquad=P_{x}\left(T^{1}_{x}<+\infty\right)P_{x}\left(T^{1}_{y}<+\infty\right)P_{y}\left(T^{1}_{x}=+\infty\right).
\end{align*}
It then follows by the inclusion $\refpar{eq:subset_inclusion_chain_go_finite_time_inx}$
that
\begin{equation}
P_{x}\left(T^{1}_{x}<+\infty\right)P_{x}\left(T^{1}_{y}<+\infty\right)P_{y}\left(T^{1}_{x}=+\infty\right)\leqslant P_{x}\left(R^{c}_{x}\right).\label{eq:upper_bound_prod_proba}
\end{equation}

Since $x$ is recurrent, we have, as shown above, $P_{x}\left(R^{c}_{x}\right)=0.$
Moreover, because $x$ leads to $y,$ we have $P_{x}\left(T^{1}_{y}<+\infty\right)>0.$
It then follows from $\refpar{eq:upper_bound_prod_proba}$ that $P_{y}\left(T^{1}_{x}=+\infty\right)=0,$
or also that $P_{y}\left(T^{1}_{x}<+\infty\right)=1.$ In particular,
$y$ leads to $x.$

Finally, we prove that $y$ is recurrent. Indeed, since $x$ and $y$
communicate, there exist $i,j\in\mathbb{N}^{\ast}$ such that $M^{i}\left(x,y\right)>0$
and $M^{j}\left(y,x\right)>0.$ Moreover, for every $n\in\mathbb{N}^{\ast},$
\[
M^{n+i+j}\left(y,y\right)\geqslant M^{j}\left(y,x\right)M^{n}\left(x,x\right)M^{i}\left(x,y\right),
\]
and, since $x$ is recurrent, we also have $R\left(x,x\right)=\sum_{n\in\mathbb{N}}M^{n}\left(x,x\right)=+\infty.$
Therefore $\sum_{n\in\mathbb{N}}M^{n+i+j}\left(y,y\right)=+\infty,$
and hence $R\left(y,y\right)=+\infty,$ which proves that $y$ is
recurrent. 

\end{proof}

\subsection{Asymptotic Behaviour and Classification}

The law of $X_{n}$ is given by 
\[
P_{x}\left(X_{n}=y\right)=M^{n}\left(x,y\right).
\]
In practice, this quantity is often difficult to compute explicitly
when the matrix $M$ is large or not sufficiently sparse. It is therefore
key to study its asymptotic behaviour.

\begin{proposition}{Paths Targeting a Transient State and Asymptotic Behaviour}{trans_asympt}

If $y$ is transient, then for every $x\in E,$ 
\[
R\left(x,y\right)<+\infty\,\,\,\,\text{and}\,\,\,\,\lim_{n\to+\infty}M^{n}\left(x,y\right)=0.
\]

\end{proposition}

\begin{proof}{}{}

Recall---Corollary $\ref{co:potent_path}$---that
\[
R\left(x,y\right)=\begin{cases}
\dfrac{1}{1-F\left(y,y\right)}, & \text{if }x=y,\\
F\left(x,y\right)R\left(y,y\right), & \text{if }x\neq y.
\end{cases}
\]

Since the state $y$ is transient, $R\left(y,y\right)<+\infty,$ and
thus also $R\left(x,y\right)<+\infty.$ Since $R\left(x,y\right)=\sum^{+\infty}_{n=0}M^{n}\left(x,y\right),$
the general term of this convergent series tends to zero when $n$
tends to infinity. 

\end{proof}

To establish the next result, we rely on the following analytical
lemma, which we state below without proof. The result is nontrivial.
A proof may be found in the book of W. Feller, p.306, \cite{feller1958introduction}.

\begin{lemma}{}{analysis_th_seq_nb_sum_1}

Let $\left(f_{j}\right)_{j\in\mathbb{N}^{\ast}}$ be a sequence of
nonnegative real numbers such that $\sum_{j\in\mathbb{N}^{\ast}}f_{j}=1$
and $\text{GCD}\left\{ j\in\mathbb{N}^{\ast}:f_{j}>0\right\} =1.$
Let $\left(u_{n}\right)_{n\in\mathbb{N}}$ be a sequence of real numbers
satisfying $u_{0}=1$ and 
\[
\forall n\in\mathbb{N}^{\ast},\,\,\,\,u_{n}=\sum^{n}_{j=1}f_{j}u_{n-j}.
\]
The sequence $\left(u_{n}\right)_{n\in\mathbb{N}}$ converges and
\[
\lim_{n\to+\infty}u_{n}=\dfrac{1}{\sum^{+\infty}_{j=1}jf_{j}}.
\]

\end{lemma}

\begin{proposition}{Paths Targetting A Recurrent Aperiodic State and Asymptotic Behaviour}{rec_aper_st}

If $y$ is \textbf{recurrent and aperiodic}, then for every $x\in E,$
the sequence with general term $M^{n}\left(x,y\right)$ converges
and
\begin{equation}
\lim_{n\to+\infty}M^{n}\left(x,y\right)=\dfrac{F\left(x,y\right)}{\mathbb{E}_{y}\left(T^{1}_{y}\right)},\label{eq:lim_of_Mn_x_y_y_rec_aper}
\end{equation}
with the convention $\dfrac{1}{\infty}=0.$

\end{proposition}

\begin{proof}{}{}

Since $\left(X_{n}=y\right)\subset\left(T^{1}_{y}\leqslant n\right),$
\begin{align*}
M^{n}\left(x,y\right) & =\mathbb{E}_{x}\left(\boldsymbol{1}_{\left(X_{n}=y\right)}\boldsymbol{1}_{\left(T^{1}_{y}\leqslant n\right)}\right)=\sum^{n}_{j=1}\mathbb{E}_{x}\left(\boldsymbol{1}_{\left(X_{n}=y\right)}\boldsymbol{1}_{\left(T^{1}_{y}=j\right)}\right)\\
 & =\sum^{n}_{j=1}\mathbb{E}_{x}\left(\boldsymbol{1}_{\left(T^{1}_{y}=j\right)}\boldsymbol{1}_{\left(X_{T^{1}_{y}+n-j}=y\right)}\right).
\end{align*}
Conditioning with respect to $\mathscr{A}^{T^{1}_{y}}$ and applying
the strong Markov property yields
\begin{align*}
M^{n}\left(x,y\right) & =\sum^{n}_{j=1}\mathbb{E}_{x}\left(\boldsymbol{1}_{\left(T^{1}_{y}=j\right)}\mathbb{E}_{X_{T^{1}_{y}}}\left[\boldsymbol{1}_{\left(X_{n-j}=y\right)}\right]\right)\\
 & =\sum^{n}_{j=1}\mathbb{E}_{x}\left[\boldsymbol{1}_{\left(T^{1}_{y}=j\right)}\right]\mathbb{E}_{y}\left[\boldsymbol{1}_{\left(X_{n-j}=y\right)}\right].
\end{align*}
That is
\begin{equation}
M^{n}\left(x,y\right)=\sum^{n}_{j=1}F_{j}\left(x,y\right)M^{n-j}\left(y,y\right).\label{eq:Mnxy_as_sum_F_jxyMn-j}
\end{equation}

\begin{itemize}
\item Suppose $x=y.$ Apply Lemma $\ref{lm:analysis_th_seq_nb_sum_1}$ with
$f_{j}=F_{j}\left(y,y\right)$ and $u_{n}=M^{n}\left(y,y\right).$
Then the relation $\refpar{eq:Mnxy_as_sum_F_jxyMn-j}$ becomes, for
every $n\in\mathbb{N}^{\ast},$ 
\[
u_{n}=\sum^{n}_{j=1}f_{j}u_{n-j}.
\]
Let 
\[
d=\text{GCD\ensuremath{\left\{  j\in\mathbb{N}^{\ast}:\,f_{j}>0\right\} } }.
\]
We prove, by induction with the help of this last relation, that 
\[
\left\{ n\in\mathbb{N}^{\ast}:\,u_{n}>0\right\} \subset d\mathbb{N}^{\ast}.
\]
Since $y$ is aperiodic, we have $d=1.$ Moreover, because $y$ is
recurrent,
\[
\sum_{j\in\mathbb{N}^{\ast}}F_{j}\left(y,y\right)=P_{y}\left(T^{1}_{y}<+\infty\right)=1.
\]
By taking into account the equalities $F\left(y,y\right)=1$ and 
\[
\sum^{+\infty}_{j=1}jf_{j}=\sum^{+\infty}_{j=1}jP_{y}\left(T^{1}_{y}=j\right)=\mathbb{E}_{y}\left(T^{1}_{y}\right),
\]
Lemma $\ref{lm:analysis_th_seq_nb_sum_1}$ establishes the equality
$\refpar{eq:lim_of_Mn_x_y_y_rec_aper}.$
\item Suppose $x\neq y.$ The equality $\refpar{eq:Mnxy_as_sum_F_jxyMn-j}$
can be written
\begin{equation}
M^{n}\left(x,y\right)=\sum^{+\infty}_{j=1}\left[\boldsymbol{1}_{\left(j\leqslant n\right)}M^{n-j}\left(y,y\right)\right]F_{j}\left(x,y\right).\label{eq:expMnxywhen_x_neq_y}
\end{equation}
Interpret this sum as the integral of the function
\[
j\mapsto\boldsymbol{1}_{\left(j\leqslant n\right)}M^{n-j}\left(y,y\right)
\]
with respect to the measure $\sum^{+\infty}_{j=1}F_{j}\left(x,y\right)\delta_{j}$
whose finite total mass is
\[
\sum^{+\infty}_{j=1}F_{j}\left(x,y\right)=P_{x}\left(T^{1}_{y}<+\infty\right)\leqslant1.
\]
For every $n\in\mathbb{N}^{\ast},$ 
\[
0\leqslant\boldsymbol{1}_{\left(j\leqslant n\right)}M^{n-j}\left(y,y\right)\leqslant1,
\]
so the dominated convergence theorem ensures the convergence of the
sequence with general term $M^{n}\left(x,y\right).$ Using the first
part,
\begin{align*}
\lim_{n\to+\infty}M^{n}\left(x,y\right) & =\sum^{+\infty}_{j=1}\left[\lim_{n\to+\infty}\boldsymbol{1}_{\left(j\leqslant n\right)}M^{n-j}\left(y,y\right)\right]F_{j}\left(x,y\right)\\
 & =\sum^{+\infty}_{j=1}\dfrac{F_{j}\left(x,y\right)}{\mathbb{E}_{y}\left(T^{1}_{y}\right)}=\dfrac{F\left(x,y\right)}{\mathbb{E}_{y}\left(T^{1}_{y}\right)}.
\end{align*}
\end{itemize}
\end{proof}

The following proposition describes the \textbf{asymptotic behaviour}
of the chain with respect to a \textbf{periodic recurrent class}.

\begin{proposition}{Path Targetting a Periodic Recurrent State and Asymptotic Behaviour}{per_rec_state}

Let $y$ be a \textbf{periodic recurrent} state with period $d>1.$

(a) If $x$ communicates with $y,$ and if $x\in C_{r}$ and $y\in C_{r+a}$
where $C_{j},\,j\in\left\llbracket 0,d-1\right\rrbracket $ are the
cyclic classes of $C,$ then the sequence with general term $M^{nd+a}\left(x,y\right)$
converges and\boxeq{
\begin{equation}
\lim_{n\to+\infty}M^{nd+a}\left(x,y\right)=\dfrac{d}{\mathbb{E}_{y}\left(T^{1}_{y}\right)}.\label{eq:limMpowndpa_xy}
\end{equation}
}

(b) If $x$ is arbitrary, then for every $a\in\left\llbracket 0,d-1\right\rrbracket ,$\boxeq{
\begin{equation}
\lim_{n\to+\infty}M^{nd+a}\left(x,y\right)=\left[\sum^{+\infty}_{j=1}F_{jd+a}\left(x,y\right)\right]\dfrac{d}{\mathbb{E}_{y}\left(T^{1}_{y}\right)}.\label{eq:limMpowndpa_xy_general}
\end{equation}
}

\end{proposition}

\begin{proof}{}{}

(a) First consider the case $a=0.$ Then $y$ is recurrent and aperiodic
for the homogeneous Markov chain $\left(X_{nd}\right)_{n\in\mathbb{N}}$
of transition matrix $M^{d}.$ By Proposition $\ref{pr:rec_aper_st},$
\[
\lim_{n\to+\infty}M^{nd}\left(x,y\right)=\dfrac{1}{\mathbb{E}_{y}\left(S^{1}_{y}\right)},
\]
where 
\[
S^{1}_{y}=\inf\left\{ n\in\mathbb{N^{\ast}}:\,X_{nd}=y\right\} .
\]
 Since $P_{y}\left(S^{1}_{y}=k\right)=P_{y}\left(T^{1}_{y}=kd\right),$
we obtain 
\[
\mathbb{E}_{y}\left(S^{1}_{y}\right)=\dfrac{1}{d}\mathbb{E}_{y}\left(T^{1}_{y}\right).
\]
This proves $\refpar{eq:limMpowndpa_xy}$ in this case.

Now suppose the result holds up to the order $a<d-1$ and let us prove
it at the order $a+1.$ We have
\[
M^{nd+a+1}\left(x,y\right)=\sum_{z\in E}M\left(x,z\right)M^{nd+a}\left(z,y\right).
\]

By the induction hypothesis and the dominated convergence theorem,
\[
\lim_{n\to+\infty}M^{nd+a+1}\left(x,y\right)=\sum_{z\in E}\left[M\left(x,z\right)\dfrac{d}{\mathbb{E}_{y}\left(T^{1}_{y}\right)}\right]=\dfrac{d}{\mathbb{E}_{y}\left(T^{1}_{y}\right)}.
\]

(b) By the equality $\refpar{eq:Mnxy_as_sum_F_jxyMn-j},$
\[
M^{nd+a}\left(x,y\right)=\sum^{nd+a}_{j=1}F_{j}\left(x,y\right)M^{nd+a-j}\left(y,y\right).
\]
Since $y$ has for period $d,$ $M^{nd+a-j}\left(y,y\right)=0,$ unless
if $a-j\in d\mathbb{Z}.$ Thus,
\begin{align*}
M^{nd+a}\left(x,y\right) & =\sum^{n}_{k=0}F_{kd+a}\left(x,y\right)M^{\left(n-k\right)d}\left(y,y\right)\\
 & =\sum^{+\infty}_{k=0}\left[\boldsymbol{1}_{\left(k\leqslant n\right)}M^{\left(n-k\right)d}\left(y,y\right)\right]F_{kd+a}\left(x,y\right).
\end{align*}
Since
\[
\sum^{+\infty}_{k=0}F_{kd+a}\left(x,y\right)=\sum^{+\infty}_{k=0}P_{x}\left(T^{1}_{y}=kd+a\right)\leqslant P_{x}\left(T^{1}_{y}<+\infty\right)\leqslant1,
\]
we obtain $\refpar{eq:limMpowndpa_xy_general},$ by application of
$\refpar{eq:limMpowndpa_xy}$ together with the dominated convergence
theorem.

\end{proof}

We now deduce two corollaries from this proposition that help to determine
the \textbf{nature of communication classes}.

\begin{corollary}{Nature of States in a Class. Nature of a Class}{}

All \textbf{states} of a \textbf{communication class} $C$ share the
\textbf{same nature}, they are either transient, positive recurrent,
null recurrent, and they are either all periodic with the same period
or all aperiodic. The \textbf{nature of the class\index{nature of the class}\mindex{class!nature}}
$C$ is therefore, by definition, the common nature of any of its
states.

\end{corollary}

\begin{proof}{}{}

Let $x,y\in C.$
\begin{itemize}
\item Since $x$ leads to $y,$ if $y$ is transient, then $x$ must also
be transient---otherwise, by Theorem $\ref{th:state_classification}$
it would be recurrent.
\item If $y$ is recurrent and aperiodic, then the same holds for $x,$
since $x$ and $y$ communicate---see Theorem $\ref{th:state_classification}$
and Proposition $\ref{pr:class_per_aper}.$ Moreover, if $y$ is null
recurrent, Propositions $\ref{pr:rec_aper_st}$ and $\ref{pr:per_rec_state}$
imply that 
\[
\lim_{n\to+\infty}M^{n}\left(y,y\right)=0.
\]
Since $x$ and $y$ communicate, there exist $k,l\in\mathbb{N}^{\ast}$
such that $M^{k}\left(x,y\right)>0$ and $M^{l}\left(y,x\right)>0.$
\\
For every $n\in\mathbb{N}^{\ast},$
\[
M^{n+k+l}\left(y,y\right)\geqslant M^{l}\left(y,x\right)M^{n}\left(x,x\right)M^{k}\left(x,y\right).
\]
It follows that 
\[
\lim_{n\to+\infty}M^{n}\left(x,x\right)=0,
\]
which, again by Proposition $\ref{pr:rec_aper_st}$ proves that $x$
is null recurrent.
\item If $y$ is positive recurrent and aperiodic, then the same holds for
$x;$ otherwise, by the previous assertion, $y$ would be null recurrent
and aperiodic.
\item The remainder of the statement follows directly from Proposition $\ref{pr:class_per_aper}.$
\end{itemize}
\end{proof}

\begin{definition}{Closed Communication Class. Absorbing State}{}

A communication class $C$ is said to be \textbf{\mindex{class!closed}closed}---or
\textbf{\mindex{class!absorbing}absorbing}---if, for every $x\in C,$
\[
P_{x}\left(T_{C^{c}}<+\infty\right)=0
\]
In other words, once the chain enters such a class, it never leaves
it. In particular, if the closed class $C$ consists of a single point,
this point is called an \textbf{\mindex{state!absorbing}absorbing
state}: once the chain reaches it, it remains there forever.

\end{definition}

\begin{corollary}{Closed Communication Class and Possible States}{closed_com_states}

Let $C$ be a \textbf{closed} communication class. If $C$ has \textbf{finite}
cardinality, then it contains neither transient state nor null recurrent
states.

In particular, an \textbf{irreducible} and \textbf{finite} homogeneous
Markov chain contains only \textbf{positive recurrent} points.

\end{corollary}

\begin{proof}{}{}

Assume that all states of $C$ are either transient or null recurrent.
Then, by Propositions $\ref{pr:rec_aper_st}$ and $\ref{pr:per_rec_state},$
possibly after extracting a subsequence, we would have $\lim_{n\to+\infty}M^{n}\left(x,y\right)=0$
for every $x,y\in C.$ Since $C$ is finite, this would imply
\[
\lim_{n\to+\infty}\sum_{y\in C}M^{n}\left(x,y\right)=0.
\]
However, this contradicts the fact that $C$ is closed. Indeed, for
every $n\in\mathbb{N}^{\ast},$
\[
0=P_{x}\left(T_{C^{c}}<+\infty\right)=P_{x}\left(\bigcup_{k\in\mathbb{N}^{\ast}}\left(X_{k}\notin C\right)\right)\geqslant P_{x}\left(X_{n}\notin C\right)=1-\sum_{y\in C}M^{n}\left(x,y\right).
\]
Passing to the limit as $n$ tends to infinity would impose $0\geqslant1,$
which is a contradiction.

\end{proof}

\begin{definition}{Inessential State}{}

A state $x$ is \textbf{\mindex{state!inessential}inessential} if
there exists $m\in\mathbb{N}^{\ast}$ and $y\neq x$ such that $M^{m}\left(x,y\right)>0$
and, for every $n>m$ and every $z\in E,$ $M^{n}\left(z,x\right)=0.$
That is, an inessential state is a state in which with a positive
probability the chain leaves it after a finite---non random---number
of steps without ever coming back.

\end{definition}

\begin{example}{}{}

For instance, for the homogeneous Markov chain, taking values in $E=\left\{ 1,2,3,4\right\} ,$
and with transition matrix $M,$
\[
M=\left(\begin{array}{cccc}
0 & 1 & 0 & 0\\
0 & 0 & \dfrac{1}{2} & \dfrac{1}{2}\\
0 & \dfrac{1}{2} & 0 & \dfrac{1}{2}\\
0 & \dfrac{1}{2} & \dfrac{1}{2} & 0
\end{array}\right).
\]
The associated graph is

\begin{center}\begin{tikzpicture}[
  >=Stealth,
  every node/.style={
    circle,
    draw=black,
    minimum size=10mm,
    inner sep=0pt
  },
  edge/.style={
    ->,
    draw=black,
    line width=0.75pt
  }
]

% Nodes: 1 faces 2; (2,3,4) form a triangle
\node (1) at (-3,0) {1};
\node (2) at (0,0) {2};
\node (3) at (2,1.6) {3};
\node (4) at (2,-1.6) {4};

% Transition 1 -> 2
\draw[edge] (1) -- (2);

% Bidirectional edges among {2,3,4} (use bending to separate directions)
\draw[edge] (2) to[bend left=20] (3);
\draw[edge] (3) to[bend left=20] (2);

\draw[edge] (2) to[bend right=20] (4);
\draw[edge] (4) to[bend right=20] (2);

\draw[edge] (3) to[bend left=20] (4);
\draw[edge] (4) to[bend left=20] (3);

\end{tikzpicture}\end{center}

We see on this graph that $1$ is an inessential state. The other
states all communicate. The chain restricted to the state space $\left\{ 2,3,4\right\} $
is therefore irreducible, positive recurrent, and has period 2.

\end{example}

\begin{proposition}{Inessential States Are Transient}{}

An inessential state is transient

\end{proposition}

\begin{proof}{}{}

Let $x$ be an inessential state. Let $m\in\mathbb{N}^{\ast}$ and
let $y\neq x$ be such that $M^{m}\left(x,y\right)>0,$ and such that,
for every $n>m,$ and every $z\in E,$ $M^{n}\left(z,x\right)=0.$
Then, for every $n>m,$
\[
M^{n+m}\left(x,x\right)=\sum_{z\in E}M^{m}\left(x,z\right)M^{n}\left(z,x\right)=0.
\]
It follows that $R\left(x,x\right)<+\infty,$ and therefore $x$ is
transient. 

\end{proof}

In summary, to study the behaviour of a homogeneous Markov chain,
we first identify the inessential states---which are then transient---and
then the communication classes of the essential states---that is,
those that are not inessential. We then refine the nature of each
class using the characterization established above. 

The following example highlights one of the many links between martingales
and Markov chains.

\begin{example}{Markov Chain and Martingale. An Example of a Genetic Model}{}

Let $X=\left(X_{n}\right)_{n\in\mathbb{N}}$ be a homogeneous Markov
chain taking values in the set of integers $E=\left\llbracket 0,N\right\rrbracket ,$
with transition matrix $M.$ Assume that the process $X$ is also
a martingale with respect to its natural filtration $\left(\mathscr{A}_{n}\right)_{n\in\mathbb{N}},$
for every probability $P_{x}$ under which the chain starts from $x$
at time 0. Assume moreover that \textbf{the boundary points $0$ and
$N$ are absorbing}. 

Define the first hitting times in $x$ by
\[
\tau_{x}=\inf\left\{ n\in\mathbb{N}^{\ast}:X_{n}\right\} =x\,\,\,\,\text{with}\,\,\,\,\inf\emptyset=+\infty.
\]
For every $x\in E,$ prove that\boxeq{
\begin{equation}
P_{x}\left(\tau_{N}<\tau_{0}\right)=\dfrac{x}{N}.\label{eq:P_x_TauNlessTau0}
\end{equation}
}

\end{example}

\begin{solutionexample}{}{}

For every $x\in E$ and for $f\in bE,$
\[
\mathbb{E}^{\mathscr{A}_{n}}_{x}\left(f\left(X_{n+1}\right)\right)=M\left(X_{n},f\right)=\sum_{y\in E}f\left(y\right)M\left(X_{n},y\right).
\]
In particular, taking for $f$ the identity map on $E$---which is
bounded---we obtain
\[
\mathbb{E}^{\mathscr{A}_{0}}_{x}\left(X_{1}\right)=\sum_{y\in E}yM\left(X_{0},y\right).
\]

Since $X$ is a martingale, taking expectations under $P_{x}$ yields
\begin{equation}
\mathbb{E}_{x}\left(X_{0}\right)=\mathbb{E}_{x}\left(X_{1}\right)=\sum_{y\in E}y\mathbb{E}_{x}\left(M\left(X_{0},y\right)\right)=\sum_{y\in E}yM\left(x,y\right).\label{eq:E_x_X0_as_sum_yMxy}
\end{equation}

Taking $x=0$ in $\refpar{eq:E_x_X0_as_sum_yMxy},$ yields
\[
0=\mathbb{E}_{0}\left(X_{0}\right)=\sum_{y\in E}yM\left(0,y\right)=\sum^{N}_{y=1}yM\left(0,y\right).
\]
This implies, since all the terms are nonnegative, that $M\left(0,y\right)=0$
for every $y\in\left\llbracket 1,N\right\rrbracket ,$ hence $M\left(0,0\right)=1.$
That is, the state $0$ is absorbing.

Similarly, taking $x=N$ in $\refpar{eq:E_x_X0_as_sum_yMxy},$ 
\[
N=\mathbb{E}_{N}\left(X_{0}\right)=\sum_{y\in E}yM\left(N,y\right)=\sum^{N-1}_{y=1}yM\left(N,y\right)+NM\left(N,N\right).
\]
Since all terms are nonnegative and 
\[
\sum^{N}_{y=0}M\left(N,y\right)=1,
\]
we must have $M\left(N,y\right)=0$ for every $y\in\left\llbracket 0,N-1\right\rrbracket ,$
and therefore $M\left(N,N\right)=1.$ That is the state $N$ is absorbing.

Since $E$ is bounded, the martingale $X$ is equi-integrable. The
second stopping theorem applied to the stopping time $\tau_{0}\land\tau_{N}$
therefore gives
\[
\mathbb{E}_{x}\left(X_{0}\right)=\mathbb{E}_{x}\left(X_{\tau_{0}\land\tau_{N}}\right).
\]

In particular,
\[
x=\mathbb{E}_{x}\left(X_{0}\right)=\mathbb{E}_{x}\left(\boldsymbol{1}_{\left(\tau_{N}<\tau_{0}\right)}\cdot N+\boldsymbol{1}_{\left(\tau_{N}>\tau_{0}\right)}\cdot0\right)=NP_{x}\left(\tau_{N}<\tau_{0}\right),
\]
which proves $\refpar{eq:P_x_TauNlessTau0}.$

\end{solutionexample}

\subparagraph*{Application to a Genetic Model\protect\footnote{This problem was studied by R.A. Fisher and S. Wright and formulated
in terms of Markov chain by G. Malécot\cite{malecot1944}.}}

We consider a population that reproduces while keeping the same size
$N.$ Each individual in a generation carries two alleles\footnote{For some genetic background, see Example $\ref{ex:genetic_model}$
later in this chapter.} of type $G$ or $g.$ Therefore, at a given generation the total
number of alleles is $2N.$ 

We assume that mating in a given generation is independent and uniform
among individuals of the population---random mating. Let $X_{n}$
be the number of alleles of type $F$ present in the population at
the $n-$th generation. Then the process $X=\left(X_{n}\right)_{n\in\mathbb{N}}$
is a homogeneous Markov chain with transition matrix $M$ generating
the probabilities---we identify the germ and the corresponding probability---defined
by
\[
M\left(i,\cdot\right)=\begin{cases}
\mathscr{B}\left(2N,\dfrac{i}{2N}\right), & \text{if }1\leqslant i\leqslant2N-1,\\
\delta_{0}, & \text{if }i=0,\\
\delta_{2N}, & \text{if }i=2N.
\end{cases}
\]
Hence,
\[
\mathbb{E}^{X_{n}=i}\left(X_{n+1}\right)=\begin{cases}
2N\cdot\dfrac{i}{2N}=i, & \text{if }1\leqslant i\leqslant2N-1,\\
0, & \text{if }i=0,\\
2N, & \text{if }i=2N.
\end{cases}
\]
Letting $E=\left\llbracket 0,2N\right\rrbracket $ be the state space,
this shows that $\mathbb{E}^{\sigma\left(X_{n}\right)}_{x}\left(X_{n+1}\right)=X_{n}$
for every $x\in E.$ With the previous notation, since $X$ is a Markov
chain and $E$ is finite, we then have, for every $x\in E,$
\[
\mathbb{E}^{\mathscr{A}_{n}}_{x}\left(X_{n+1}\right)=X_{n}.
\]
It follows that the points $0$ and $2N$ are absorbing points. \textbf{In
other words, in the long run, only one allele type will remain}.

\subsection{Analytic Criteria of Induction}

We first provide a way to compute the probability that the chain remains
forever in a subset $A\subset E.$ Let $Q$ be the restriction of
$M$ to $A,$ that is, the matrix indexed by $A\times A$ defined,
for every $x,y\in A,$ by 
\[
Q\left(x,y\right)=M\left(x,y\right).
\]
 For every $n\geqslant2,$
\begin{align*}
Q^{n}\left(x,y\right) & =\sum_{x_{1}\in A}\sum_{x_{2}\in A}\cdots\sum_{x_{n-1}\in A}Q\left(x,x_{1}\right)Q\left(x_{1},x_{2}\right)\cdots Q\left(x_{n-1},y\right)\\
 & =P_{x}\left(X_{1}\in A,\cdots,X_{n-1}\in A,X_{n}=y\right),
\end{align*}
so that
\[
P_{x}\left(X_{1}\in A,\cdots,X_{n-1}\in A,X_{n}\in A\right)=\sum_{y\in A}Q^{n}\left(x,y\right)\equiv Q^{n}\left(x,A\right).
\]
For $x\in A,$ define 
\[
f_{n}\left(x\right)=P_{x}\left[\bigcap^{n}_{j=1}\left(X_{j}\in A\right)\right]=\sum_{y\in A}Q^{n}\left(x,y\right).
\]
The sequence with general term $f_{n}\left(x\right)$ is nonincreasing
and converges to 
\[
f\left(x\right)=P_{x}\left[\bigcap_{j\in\mathbb{N}^{\ast}}\left(X_{j}\in A\right)\right].
\]

\begin{proposition}{Maximal Solution To The Fix Point Equation}{max_sol_fix_point_eq}

The function $f,$ defined on $A,$ is the maximal solution of the
system
\begin{equation}
h=Qh\,\,\,\,\,\,\,\,0\leqslant h\leqslant1,\label{eq:h_sol_max_syst_heqQh}
\end{equation}
where, as previously\footnotemark, we denote $h\left(x\right)=\sum_{y\in A}Q\left(x,y\right)h\left(y\right).$
Moreover, either $f=0,$ or $\sup_{x\in A}f\left(x\right)=1.$

\end{proposition}

\footnotetext{For the vector interpretation of functions, see the
notation introduced after Definition $\refpar{df:fam_tr_mat}.$ In
particular, we shall here use the notation $Qf,$ which is better
suited to the vector computations, instead of $Q\left(\cdot,f\right).$}

\begin{proof}{}{}

For every $x\in A,$ by associativity and commutativity of nonnegative
term sums,
\begin{align*}
f_{n+1}\left(x\right) & =\sum_{y\in A}\left[\sum_{z\in A}Q\left(x,z\right)Q^{n}\left(z,y\right)\right]=\sum_{z\in A}Q\left(x,z\right)\left[\sum_{y\in A}Q^{n}\left(z,y\right)\right]\\
 & =Qf_{n}\left(x\right).
\end{align*}
Thus $f_{n+1}=Qf_{n}.$ 

By the dominated convergence theorem, passing to the limit yields
$f=Qf.$ Clearly $0\leqslant f\leqslant1.$ So, $f$ solves the system
$\refpar{eq:h_sol_max_syst_heqQh}.$ 

We now prove its maximality. 

Let $h$ be another solution of the system. Denote by 1 the constant
function equal to 1. For every $n\in\mathbb{N}^{\ast},$
\[
h=Q^{n}h\leqslant Q^{n}1=f_{n}.
\]
Passing to the limit gives, $h\leqslant f.$

Assume now that $f$ is not null, and let $c=\sup_{x\in A}f\left(x\right).$
For every $n\in\mathbb{N}^{\ast},$
\[
f=Q^{n}f\leqslant Q^{n}c=cf_{n}.
\]
Passing to the limit yields, $f\leqslant cf,$ which implies $c\geqslant1.$
Since moreover $c\leqslant1$---as $0\leqslant f\leqslant1$---hence
$c=1.$

\end{proof}

\begin{corollary}{Necessary and Sufficient Condition for a Chain to be Recurrent}{chain_recurrent}

Let $X$ be an irreducible homogeneous Markov chain with transition
matrix $M.$ Fix $x_{0}\in E,$ and let $Q$ be the restriction of
$M$ to $E\backslash\left\{ x_{0}\right\} .$ The chain $X$ is recurrent
if and only if the system
\begin{equation}
h=Qh\,\,\,\,\,\,\,\,0\leqslant h\leqslant1\label{eq:system_hisQh}
\end{equation}
admits the unique solution $h=0.$

\end{corollary}

\begin{proof}{}{}

Assume that $0$ is the unique solution of the system $\refpar{eq:system_hisQh}.$
Since the chain is irreducible, all states have the same nature. Moreover,
by denoting $A=E\backslash\left\{ x_{0}\right\} ,$ there exists $y\in A$
such that $x_{0}$ leads to $y.$ 

By Proposition $\ref{pr:max_sol_fix_point_eq},$ the function $f$
associated with this subset $A$ is the maximal solution of the system
$\refpar{eq:system_hisQh},$ which implies, by hypothesis, that 
\[
P_{y}\left[\bigcap_{j\in\mathbb{N}^{\ast}}\left(X_{j}\in A\right)\right]=0,
\]
hence, also that $P_{y}\left(T_{x_{0}}<+\infty\right)=1.$ We show
that this imply that $x_{0}$ is recurrent.

By the simple Markov property, for every $x\in A$ and for every $n,p\in\mathbb{N}^{\ast},$
\[
P_{x}\left[\bigcap^{n+p}_{i=p}\left(X_{i}\in A\right)\right]=\mathbb{\mathbb{E}}_{x}\left[\boldsymbol{1}_{\left(X_{p}\in A\right)}Q^{n}\left(X_{p},A\right)\right]=\sum_{z\in A}M^{p}\left(x,z\right)Q^{n}\left(z,A\right).
\]
Passing to the limit in $n$---the left-hand side is nonincreasing
in $n$ and the dominated convergence theorem applies on the right---and
since $f=0,$ 
\[
P_{x}\left[\bigcap^{n+p}_{i=p}\left(X_{i}\in A\right)\right]=\sum_{z\in A}M^{p}\left(x,z\right)\left[\lim_{n\to+\infty}Q^{n}\left(z,A\right)\right]=\sum_{z\in A}M^{p}\left(x,z\right)=0.
\]
Thus 
\[
P_{x}\left(\liminf_{n\to+\infty}\left(X_{n}\in A\right)\right)=0,
\]
and therefore, for every $x\in A.$
\begin{equation}
P_{x}\left(N_{x_{0}}=+\infty\right)\geqslant P_{x}\left(\limsup_{n\to+\infty}\left(X_{n}=x_{0}\right)\right)=1.\label{eq:P_xN_x0plusinf}
\end{equation}
Saying that the chain reaches $y$ in finite time and visits $x_{0}$
infinitely many times is equivalent to saying that it reaches $y$
in finite time and, after the visit to $y,$ visits $x_{0}$ infinitely
many times. This can be written as
\[
\left(N_{x_{0}}=+\infty\right)\bigcap\left(T^{1}_{y}<+\infty\right)=\left(T^{1}_{y}<+\infty\right)\cap\left(n_{x_{0}}\left[\theta_{T^{1}_{y}}\left(X\right)\right]=+\infty\right).
\]
Conditioning with respect to $\mathscr{A}_{T^{1}_{y}}$ and applying
the strong Markov property, we obtain
\begin{align*}
 & P_{x_{0}}\left[\left(N_{x_{0}}=+\infty\right)\bigcap\left(T^{1}_{y}<+\infty\right)\right]\\
 & \qquad\qquad=\mathbb{E}_{x_{0}}\left[\boldsymbol{1}_{\left(T^{1}_{y}<+\infty\right)}\mathbb{E}^{\mathscr{A}_{T^{1}_{y}}}_{x_{0}}\left(\boldsymbol{1}_{\left(n_{x_{0}}\left[\theta_{T^{1}_{y}}\left(X\right)=+\infty\right]\right)}\right)\right]\\
 & \qquad\qquad=\mathbb{E}_{x_{0}}\left[\boldsymbol{1}_{\left(T^{1}_{y}<+\infty\right)}\mathbb{E}_{y}\left(\boldsymbol{1}_{\left(n_{x_{0}}\left(X\right)=+\infty\right)}\right)\right].
\end{align*}
Thus
\[
P_{x_{0}}\left[\left(N_{x_{0}}=+\infty\right)\bigcap\left(T^{1}_{y}<+\infty\right)\right]=P_{x_{0}}\left(T^{1}_{y}<+\infty\right)P_{y}\left(N_{x_{0}}=+\infty\right).
\]
Since, moreover, $x_{0}$ leads to $y,$ it follows from this equality
and $\refpar{eq:P_xN_x0plusinf}$ that
\[
P_{x_{0}}\left(N_{x_{0}}=+\infty\right)\geqslant P_{x_{0}}\left(T^{1}_{y}<+\infty\right)P_{y}\left(N_{x_{0}}=+\infty\right)=P_{x_{0}}\left(T^{1}_{y}<+\infty\right)>0.
\]
But $P_{x_{0}}\left(N_{x_{0}}=+\infty\right)$ can take only the values
0 or 1. Hence, 
\[
P_{x_{0}}\left(N_{x_{0}}=+\infty\right)=1,
\]
 which proves that $x_{0}$ is recurrent.

Conversely, if $x_{0}$ is recurrent, then for every $z\in A,$
\[
P_{z}\left(T^{1}_{x_{0}}<+\infty\right)=1
\]
and therefore
\[
f\left(z\right)=P_{z}\left(\bigcap_{j\in\mathbb{N}^{\ast}}\left(X_{j}\in A\right)\right)=0.
\]
Since Proposition $\ref{pr:max_sol_fix_point_eq}$ states that $f$
is the maximal solution of the system $\refpar{eq:system_hisQh},$
it follows that $f$ is the unique solution of this system.

\end{proof}

\begin{remark}{}{}

Although Corollary $\ref{co:chain_recurrent}$ may appear restrictive,
since it is stated for an irreducible chain, it is in fact of general
use for determining whether a communication class $C$ is recurrent.
Indeed, it suffices to apply the corollary to the chain restricted
to the class $C,$ which is itself irreducible. Naturally, this is
of interest only when $E$ is infinite.

\end{remark}

\section{Computation of the Potential Matrix and of $P_{x}\left(T^{1}_{y}<+\infty\right)$}\label{sec:Computation-of-the}

\subsection{Computation of the Potential Matrix}

Let $x,y\in E.$
\begin{itemize}
\item \textbf{If $y$ is recurrent}, it follows by the equality $\refpar{eq:potential_mat_pr}$
from Corollary $\ref{co:potent_path}$
\[
R\left(x,y\right)=\begin{cases}
0, & \text{if }F\left(x,y\right)=0,\\
+\infty, & \text{if }F\left(x,y\right)>0.
\end{cases}
\]
\item \textbf{If $y$ is transient},
\begin{itemize}
\item \textbf{If $x$ is recurrent}, then $x$ does not lead to $y.$ Consequently,
$F\left(x,y\right)=0$ and $R\left(x,y\right)=0.$
\item \textbf{If $x$ is transient}; denote by $D$ the set of transient
states, and by $Q$ and $S$ the restrictions to $D\times D$ of the
matrices $M$ and $R,$ respectively. To treat this case, we focus
on the global computation of $S.$ After possibly reindexing the elements
of $E,$ listing first the recurrent states, the matrix $M$ has the
block structure
\[
M=\left(\begin{array}{cc}
K & 0\\
L & Q
\end{array}\right).
\]
Then,
\[
M^{n}=\left(\begin{array}{cc}
K^{n} & 0\\
L_{n} & Q^{n}
\end{array}\right)
\]
where $L_{n}$ is not a power of $L.$ \\
Consequently,
\[
R=\sum^{+\infty}_{n=0}M^{n}=\left(\begin{array}{cc}
\sum^{+\infty}_{n=0}K^{n} & 0\\
\sum^{+\infty}_{n=0}L_{n} & \sum^{+\infty}_{k=0}Q^{n}
\end{array}\right).
\]
Hence,\boxeq{
\[
S=\sum^{+\infty}_{n=0}Q^{n}.
\]
}Denoting by $I$ the identity matrix on $D,$ we have
\[
SQ=QS=S-I,
\]
and therefore,\boxeq{
\[
\left(I-Q\right)S=S\left(I-Q\right)=I.
\]
}In particular, \textbf{if the set $D$ of transient states is finite},\boxeq{
\[
S=\left(I-Q\right)^{-1}
\]
}
\end{itemize}
\end{itemize}

\subsection{Computation of $F\left(x,y\right)\equiv P_{x}\left(T^{1}_{y}<+\infty\right)$}
\begin{itemize}
\item \textbf{If $x$ and $y$ are recurrent},
\begin{itemize}
\item If they belong to the same communication class, then $F\left(x,y\right)=1,$
\item Otherwise, $F\left(x,y\right)=0.$
\end{itemize}
\item \textbf{If $x$ is recurrent and $y$ is transient}, then $F\left(x,y\right)=0.$
\item \textbf{If $x$ and $y$ are transient}, it follows by the equality
$\refpar{eq:potential_mat_pr}$ from Corollary $\ref{co:potent_path}$
that
\[
F\left(x,y\right)=\begin{cases}
\dfrac{R\left(x,y\right)}{R\left(y,y\right)}, & \text{if }x\neq y,\\
1-\dfrac{1}{R\left(y,y\right)}, & \text{if }x=y.
\end{cases}
\]
\item \textbf{If $x$ is transient and $y$ is recurrent}, the answer is
given by the following propositions.
\end{itemize}
\begin{proposition}{}{}

Let $C$ be a recurrent communication class. For every transient state
$x$ and every $y\in C,$
\[
F\left(x,y\right)=P_{x}\left(T_{C}<+\infty\right).
\]

\end{proposition}

\begin{proof}{}{}

Since $y\in C,$
\[
P_{x}\left(T^{1}_{y}<+\infty\right)\leqslant P_{x}\left(T_{C}<+\infty\right).
\]

Conversely,
\[
P_{x}\left(T_{C}<+\infty\right)=\mathbb{E}_{x}\left[\boldsymbol{1}_{\left(T^{1}_{y}<+\infty\right)}\boldsymbol{1}_{\left(T_{C}<+\infty\right)}\right]+\mathbb{E}_{x}\left[\boldsymbol{1}_{\left(T^{1}_{y}=+\infty\right)}\boldsymbol{1}_{\left(T_{C}<+\infty\right)}\right].
\]
Observe that
\[
\left(T^{1}_{y}=+\infty\right)\cap\left(T_{C}<+\infty\right)\subset\left(T_{C}<+\infty\right)\cap\left(\tau_{y}\left[\theta_{T_{C}}\left(X\right)\right]=+\infty\right).
\]
Conditioning with respect to $\mathscr{A}_{T_{C}},$ and applying
the strong Markov property, we obtain
\[
P_{x}\left(T_{C}<+\infty\right)\leqslant\mathbb{E}_{x}\left[\boldsymbol{1}_{\left(T^{1}_{y}<+\infty\right)}\right]+\mathbb{E}_{x}\left[\boldsymbol{1}_{\left(T_{C}<+\infty\right)}\mathbb{E}_{X_{T_{C}}}\left(\boldsymbol{1}_{\left(T^{1}_{y}=+\infty\right)}\right)\right].
\]
Since $C$ is recurrent and $y\in C,$ 
\[
\mathbb{E}_{X_{T_{C}}}\left(\boldsymbol{1}_{\left(T^{1}_{y}=+\infty\right)}\right)=0,
\]
which yields the other inequality
\[
P_{x}\left(T_{C}<+\infty\right)\leqslant P_{x}\left(T^{1}_{y}<+\infty\right).
\]

This proves the equality.

\end{proof}

Denote by $D$ the set of transient states and by $\left(C_{j}\right)_{j\in J}$
the family of recurrent communication classes. We now describe a method
to compute $P_{x}\left(T_{C}<+\infty\right),$ for every $x\in D$
and for every $y\in J.$ 

This computation is fundamental, since a chain starting from a state
$x\in D$ is eventually absorbed into a unique recurrent class $C_{j}.$

To this end, define a process $Y=\left(Y_{n}\right)_{n\in\mathbb{N}},$
on the same underlying filtration than the process $X,$ taking values
in\footnote{By a slight abuse of notation, we denote $D\cup J$ the set of every
elements of $D$ to which we adjunct the elements of $J.$} $D\cup J$ by
\[
Y_{n}=\begin{cases}
X_{n}, & \text{if }X_{n}\in D,\\
j, & \text{if }X_{n}\in C_{j},\,\,\,\,j\in J.
\end{cases}
\]

\begin{proposition}{}{}

The process $Y$ is a homogeneous Markov chain with transition matrix
$\widehat{M}$ given by
\[
\begin{cases}
\widehat{M}\left(x,y\right)=M\left(x,y\right), & \text{if }x,y\in D,\\
\widehat{M}\left(x,j\right)=\sum_{z\in C_{j}}M\left(x,z\right)\equiv b_{j}\left(x\right), & \text{if }x\in D,\,j\in J,\\
\widehat{M}\left(i,j\right)=\delta_{i,j} & \text{if }i,j\in J,\\
\widehat{M}\left(i,x\right)=0 & \text{if }i\in J,\,x\in D.
\end{cases}
\]
We use the same notations than in the previous section. Let $Q$ be
the restriction to $D\times D$ of the matrix $M$ and $S=\sum^{+\infty}_{n=0}Q^{n}.$
Define the matrix $B,$ indexed on $D\times J$ by
\[
\forall\left(x,j\right)\in D\times J,\,\,\,\,B\left(x,j\right)=b_{j}\left(x\right).
\]
Then, for every $\left(x,j\right)\in D\times J,$\boxeq{
\begin{equation}
P_{x}\left(T_{C_{j}}<+\infty\right)=\left(SB\right)\left(x,j\right).\label{eq:P_x_T_C_jless+infxtransientjrecurrent}
\end{equation}
}

We recall, that in particular, if $D$ is finite, 
\[
S=\left(1-Q\right)^{-1}.
\]

\end{proposition}

\begin{proof}{}{}

Define $g$ the mapping from $E$ to $D\cup J$ by
\[
g\left(x\right)=\begin{cases}
x, & \text{if }x\in D,\\
j, & \text{if }x\in C_{j},\,\,\,\,j\in J.
\end{cases}
\]
Then $Y_{n}=g\left(X_{n}\right).$ 

For every bounded function $f$ defined on $D\cup J,$
\begin{equation}
\mathbb{E}^{\mathscr{A}_{n}}_{x}\left[f\left(Y_{n+1}\right)\right]=\mathbb{E}^{\mathscr{A}_{n}}_{x}\left[f\circ g\left(X_{n+1}\right)\right]=M\left(X_{n},f\circ g\right).\label{eq:exp_x_f_Y_npp1_kn_A_n}
\end{equation}
Since
\[
f\circ g=f\boldsymbol{1}_{D}+\sum_{j\in J}f\left(j\right)\boldsymbol{1}_{C_{j}},
\]
we have, for every $x\in E,$ by linearity of $M\left(x,\cdot\right).$
\begin{align*}
M\left(x,f\circ g\right) & =M\left(x,f\boldsymbol{1}_{D}\right)+\sum_{j\in J}f\left(j\right)M\left(x,\boldsymbol{1}_{C_{j}}\right)\\
 & =\left[\boldsymbol{1}_{D}\left(x\right)M\left(x,f\boldsymbol{1}_{D}\right)+\sum_{j\in J}\boldsymbol{1}_{C_{j}}\left(x\right)M\left(x,f\boldsymbol{1}_{D}\right)\right]\\
 & \qquad\qquad+\left[\sum_{j\in J}f\left(j\right)\boldsymbol{1}_{D}\left(x\right)M\left(x,\boldsymbol{1}_{C_{j}}\right)+\sum_{i\in J}\boldsymbol{1}_{C_{i}}\left(x\right)\left[\sum_{j\in J}f\left(j\right)M\left(x,\boldsymbol{1}_{C_{j}}\right)\right]\right]
\end{align*}
by noting that, if $x\in D,$ then 
\[
M\left(x,\boldsymbol{1}_{C_{j}}\right)=b_{j}\left(x\right)=\widehat{M}\left(x,j\right)
\]
and that, for every $x\in E,$
\[
\boldsymbol{1}_{C_{j}}\left(x\right)M\left(x,f\boldsymbol{1}_{D}\right)=0,
\]
and
\[
\boldsymbol{1}_{C_{i}}\left(x\right)M\left(x,\boldsymbol{1}_{C_{j}}\right)=\delta_{ij}\boldsymbol{1}_{C_{i}}\left(x\right)=\widehat{M}\left(i,j\right)\boldsymbol{1}_{C_{i}}\left(x\right),
\]
it follows
\begin{align*}
M\left(x,f\circ g\right) & =\boldsymbol{1}_{D}\left(x\right)\left[\widehat{M}\left(x,f\boldsymbol{1}_{D}\right)+\sum_{j\in J}f\left(j\right)\widehat{M}\left(x,j\right)\right]+\sum_{i\in J}\boldsymbol{1}_{C_{i}}\left(x\right)\left[\sum_{j\in J}f\left(j\right)\widehat{M}\left(i,j\right)\right]\\
 & =\boldsymbol{1}_{D}\left(x\right)\widehat{M}\left(x,f\right)+\sum_{i\in J}\boldsymbol{1}_{C_{i}}\widehat{M}\left(i,f\right).
\end{align*}

Thus,
\[
M\left(Y_{n},f\circ g\right)=\boldsymbol{1}_{D}\left(Y_{n}\right)\widehat{M}\left(Y_{n},f\right)+\sum_{i\in J}\boldsymbol{1}_{i}\left(Y_{n}\right)\widehat{M}\left(i,f\right)=\widehat{M}\left(Y_{n},f\right).
\]
Substituting into the equality $\refpar{eq:exp_x_f_Y_npp1_kn_A_n}$
gives
\[
\mathbb{E}^{\mathscr{A}_{n}}_{x}\left[f\left(Y_{n+1}\right)\right]=\widehat{M}\left(Y_{n},f\right).
\]
This proves that $Y$ is a homogeneous Markov chain with transition
matrix $\widehat{M}.$

Since
\[
T_{C_{j}}=\inf\left\{ n\in\mathbb{N}^{\ast}:\,X_{n}\in C_{j}\right\} =\inf\left\{ n\in\mathbb{N}^{\ast}:\,Y_{n}=j\right\} \equiv\widehat{T}_{j},
\]
for every $x\in D,$
\[
P_{x}\left(T_{C_{j}}<+\infty\right)=P_{x}\left(\widehat{T}_{j}<+\infty\right)=P_{x}\left[\bigcup_{n\in\mathbb{N}^{\ast}}\left(Y_{n}=j\right)\right].
\]
However, the class $C_{j}$ is recurrent, hence, $P_{x}-$almost surely
\[
\left(Y_{n}=j\right)\subset\left(Y_{n+1}=j\right)
\]
and therefore
\[
P_{x}\left(T_{C_{j}}<+\infty\right)=\lim_{n\to+\infty}\nearrow P_{x}\left(Y_{n}=j\right)=\lim_{n\to+\infty}\widehat{M}^{n}\left(x,j\right).
\]
Now, the matrix $\widehat{M}$ has a block structure
\[
\widehat{M}=\left(\begin{array}{cc}
Q & B\\
0 & I
\end{array}\right).
\]
Thus, for every $n\in\mathbb{N}^{\ast},$
\[
\widehat{M}^{n}=\left(\begin{array}{cc}
Q^{n} & B_{n}\\
0 & I
\end{array}\right),
\]
with
\[
B_{n}=\left(1+Q+Q^{2}+\cdots+Q^{n}\right)B.
\]
Passing to the limit yields
\[
\lim_{n\to+\infty}B_{n}=\left(\sum^{+\infty}_{n=0}Q^{n}\right)B=SB,
\]
which gives the announced result.

\end{proof}

\section{Invariant Measures}\label{sec:Invariant-Measures}

The notion of an invariant measure for a homogeneous Markov chain
with transition matrix $M$ is closely related to its asymptotic behaviour.

Henceforth, since $E$ is countable, we identify a measure $\nu$
on $E$ and its germ. Taking into account the duality between functions
and measures, and maintaining the previously adopted vector viewpoint,
we identify the measure $\nu$ with the row vector $\left(\nu\left(x\right)\right)_{x\in E}.$

\begin{definition}{Invariant Measure. Invariant Probability}{}

Let $M$ be a transition matrix on $E.$ To each measure $\nu$ on
$E,$ we associate the measure $\nu M,$ defined for every $y\in E,$
by
\begin{equation}
\nu M\left(y\right)=\sum_{x\in E}\nu\left(x\right)M\left(x,y\right).\label{eq:def_nuM}
\end{equation}
We say that $\nu$ is an \textbf{invariant measure}\index{invariant measure}---with
respect to the transition matrix $M$---if $\nu M=\nu,$ that is,
under the above identifications, if $\nu$ is a left eigenvector of
$M$ associated with the eigenvalue $1.$ An invariant measure that
is a probability measure is called an \textbf{invariant probability}.

Let $X$ be a homogeneous Markov chain with transition matrix $M$
and let $\nu$ be an invariant measure with respect to $M.$ $\nu$
is called the \textbf{invariant measure of the chain}\index{invariant measure of the chain}.

\end{definition}

\begin{remark}{}{}

If $\nu$ is an invariant measure, then for every $a\geqslant0,$
the measure $a\nu$ is also invariant. Moreover, if $\nu_{1}$ and
$\nu_{2}$ are invariant probabilities, then any convex combination
of $\nu_{1}$ and $\nu_{2}$ is again an invariant probability. 

Consequently, the existence of two distinct invariant probabilities
implies the existence of infinitely many invariant probabilities.

\end{remark}

\begin{proposition}{}{}

Let $X$ be a homogeneous Markov chain with transition matrix $M$
admitting an invariant measure $\nu.$ If $X$ has initial law $\nu,$
then for every $n\in\mathbb{N}^{\ast},$ the random variable $X_{n}$
also follows the law $\nu.$

\end{proposition}

\begin{proof}{}{}

For every $y\in E,$ by denoting $E^{\prime}=\left\{ x\in E:\nu\left(x\right)\neq0\right\} ,$
\[
P\left(X_{n}=y\right)=\sum_{x\in E^{\prime}}P\left(X_{0}=x\right)P^{\left(X_{0}=x\right)}\left(X_{n}=y\right)=\sum_{x\in E}\nu\left(x\right)M^{n}\left(x,y\right).
\]
Then,
\begin{align*}
\sum_{x\in E}\nu\left(x\right)M^{n}\left(x,y\right) & =\sum_{x\in E}\nu\left(x\right)\left[\sum_{z\in E}M\left(x,z\right)M^{n-1}\left(z,y\right)\right]\\
 & =\sum_{z\in E}\left[\sum_{x\in E}\nu\left(x\right)M\left(x,z\right)\right]M^{n-1}\left(z,y\right).
\end{align*}
Now, since $\nu$ is invariant
\[
\sum_{x\in E}\nu\left(x\right)M^{n}\left(x,y\right)=\sum_{z\in E}\nu\left(z\right)M^{n-1}\left(z,y\right).
\]
Thus, for every $n\in\mathbb{N}^{\ast},$ $\nu M^{n}=\nu M^{n-1},$
and by induction, $\nu M^{n}=\nu.$ Hence, for every $y\in E,$ 
\[
P\left(X_{n}=y\right)=\nu\left(y\right).
\]

\end{proof}

We now study the existence and uniqueness of an invariant probability,
as well as its relation with the existence of a limiting probability.
From this, we deduce a criterion for positive recurrence of a homogeneous
Markov chain in terms of invariant probabilities.

\begin{proposition}{Existence And Uniqueness of an Invariant Probability}{ex_unic_inv_prob}

Let $M$ be a transition matrix such that, for every $x,y\in E,$
the sequence with general term $M^{n}\left(x,y\right)$ converges
to a limit $\pi\left(y\right)$ independent of $x.$ Then:

(a) The measure $\pi$ is invariant and mass less than or equal to
1. That is
\[
\pi M=\pi\,\,\,\,\text{and}\,\,\,\,\sum_{y\in E}\pi\left(y\right)\leqslant1.
\]

(b) Either $\pi=0$---that is $\pi\left(y\right)=0$ for every $y\in E$---or
$\pi$ is an invariant probability.

(c) If $\pi=0,$ then there exists no invariant probability for $M.$
If $\pi$ is an invariant probability, then $\pi$ is the unique invariant
probability for $M.$

\end{proposition}

\begin{proof}{}{}

(a) By the Fatou lemma, 
\[
\sum_{y\in E}\pi\left(y\right)=\sum_{y\in E}\lim_{n\to+\infty}M^{n}\left(x,y\right)\leqslant\liminf_{n\to+\infty}\sum_{y\in E}M^{n}\left(x,y\right)=1.
\]

Moreover, for every $y\in E,$ for every $x\in E,$
\[
\pi\left(x\right)=\lim_{n\to+\infty}M^{n}\left(y,x\right).
\]
Hence,
\begin{align*}
\pi M\left(y\right) & =\sum_{x\in E}\left[\lim_{n\to+\infty}M^{n}\left(y,x\right)\right]M\left(x,y\right)\\
 & \leqslant\liminf_{n\to+\infty}\sum_{y\in E}M^{n}\left(y,x\right)M\left(x,y\right)=\liminf_{n\to+\infty}M^{n+1}\left(y,y\right)=\pi\left(y\right).
\end{align*}
Thus $\pi M\leqslant\pi.$ Suppose that the equality does not hold.
Then it would exist $y_{0}$ such that $\pi M\left(y_{0}\right)<\pi\left(y_{0}\right)$
then---since $\pi M\leqslant\pi$---
\[
\sum_{y\in E}\pi\left(y\right)>\sum_{y\in E}\left[\sum_{x\in E}\pi\left(x\right)M\left(x,y\right)\right]=\sum_{x\in E}\pi\left(x\right)\left[\sum_{y\in E}M\left(x,y\right)\right]=\sum_{x\in E}\pi\left(x\right),
\]
which is impossible. Hence, 
\[
\pi M=\pi.
\]
 Therefore, $\pi$ is an invariant measure with mass less than or
equal to 1.

(b) If $\nu$ is an invariant measure of mass less than or equal to
1, then for every $n\in\mathbb{N}^{\ast},$ $\nu M^{n}=\nu.$ Thus,
for every $y\in E,$
\[
\nu\left(y\right)=\lim_{n\to+\infty}\nu M^{n}\left(y\right)=\lim_{n\to+\infty}\sum_{x\in E}\nu\left(x\right)M^{n}\left(x,y\right).
\]
By the dominated convergence theorem and the definition of $\pi,$
\begin{equation}
\nu\left(y\right)=\sum_{x\in E}\nu\left(x\right)\left[\lim_{n\to+\infty}M^{n}\left(x,y\right)\right]=\sum_{x\in E}\nu\left(x\right)\pi\left(y\right)\label{eq:nu_y_dev}
\end{equation}
and thus,
\begin{equation}
\nu\left(y\right)=\pi\left(y\right)\left[\sum_{x\in E}\nu\left(x\right)\right].\label{eq:nu_y_fact}
\end{equation}
In particular, if $\nu=\pi,$ then
\[
\pi\left(y\right)=\pi\left(y\right)\left[\sum_{x\in E}\pi\left(x\right)\right],
\]
which implies the stated result.

(c) If $\pi=0,$ the previous equality $\refpar{eq:nu_y_fact}$ implies
that $\nu=0.$ Hence, no invariant probability exists for $M.$ Still
by the same equality, if $\nu$ is an invariant probability, then
it coincides with $\pi.$

\end{proof}

The following theorem makes it possible, in the case where $E$ is
infinite, to determine whether a communication class is positive recurrent
and to compute the mean return time to a given state.

\begin{theorem}{Positive Recurrence Criterion}{positive_recurrence_crit}

An homogeneous Markov chain with transition matrix $M$ admits a \textbf{unique
invariant probability $\pi$} if and only if it has exactly one positive
recurrent communication class.

In this case,
\begin{equation}
\pi\left(x\right)=\begin{cases}
\dfrac{1}{\mathbb{E}_{x}\left(T^{1}_{x}\right)}, & \text{if }x\in C,\\
0, & \text{otherwise.}
\end{cases}\label{eq:exp_pi_one_positive_rec_com_cl}
\end{equation}

\end{theorem}

\begin{proof}{}{}

\textbf{1. No positive recurrent class}

If $X$ does not have any positive recurrent class, then all states
of $X$ are either transient or null recurrent. In both cases, it
follows from Propositions $\ref{pr:rec_aper_st},$ $\ref{pr:trans_asympt}$
and $\ref{pr:per_rec_state}$ that, for every $x,y\in E,$
\[
\lim_{n\to+\infty}M^{n}\left(x,y\right)=0.
\]

Proposition $\ref{pr:ex_unic_inv_prob}$ then implies that no invariant
probability exists.

\textbf{2. One positive recurrent class}

Assume $X$ admits only one positive recurrent class $C.$ There are
two cases.
\begin{itemize}
\item \textbf{Case 1: $C$ is aperiodic}\\
By Proposition $\ref{pr:rec_aper_st},$ for every $x,y\in C,$
\[
\lim_{n\to+\infty}M^{n}\left(x,y\right)=\dfrac{1}{\mathbb{E}_{y}\left(T^{1}_{y}\right)}=\pi\left(y\right)>0,
\]
and this limit is independent of $x\in C.$ \\
Since the restriction of $M$ to $C\times C$ is in this case a transition
matrix, Proposition $\ref{pr:ex_unic_inv_prob}$ applied to this matrix
shows that $\pi_{\mid_{C}}$ is the unique invariant probability with
respect to the matrix $M_{\mid_{C\times C}}.$ In particular, it follows
that the measure $\pi$ is a probability on $E$ invariant for $M.$
Indeed,
\begin{itemize}
\item If $y\in C,$ since $\pi\left(x\right)=0$ for every $x\notin C,$
then
\begin{align*}
\pi\left(y\right) & =\pi_{\mid_{C}}\left(y\right)=\sum_{x\in C}\pi_{\mid_{C}}\left(x\right)M_{\mid_{C\times C}}\left(x,y\right)\\
 & =\sum_{x\in C}\pi\left(x\right)M\left(x,y\right)=\sum_{x\in E}\pi\left(x\right)M\left(x,y\right)=\pi M\left(y\right).
\end{align*}
\item If $y\notin C,$ since $C$ is a closed class, for every $x\in C,$
$M\left(x,y\right)=0,$ then
\[
\sum_{x\in E}\pi\left(x\right)M\left(x,y\right)=\sum_{x\in C^{c}}\pi\left(x\right)M\left(x,y\right)=0=\pi\left(y\right).
\]
\end{itemize}
Thus the\textbf{ existence} of an invariant probability for $M$ is
proved.

Let us prove the \textbf{uniqueness}. If $\nu$ is an invariant probability
on $E,$ then it follows by $\refpar{eq:nu_y_dev}$ that, for every
$y\in C,$
\[
\nu\left(y\right)=\pi\left(y\right)\sum_{x\in C}\nu\left(x\right)+r\left(y\right),
\]
where
\[
r\left(y\right)=\sum_{x\in C^{c}}\nu\left(x\right)\left[\lim_{n\to+\infty}M^{n}\left(x,y\right)\right].
\]
For every $y\in C,$
\[
\nu\left(y\right)=\pi\left(y\right)\nu\left(C\right)+r\left(y\right).
\]
Summing over $y\in C$ gives
\[
\nu\left(C\right)=\pi\left(C\right)\nu\left(C\right)+r\left(C\right).
\]
Since $\pi$ is a probability on $C,$ then $r\left(C\right)=0.$
Hence, for every $y\in C,$
\begin{equation}
\nu\left(y\right)=\pi\left(y\right)\nu\left(C\right).\label{eq:nu_yispi_ynu_c}
\end{equation}
If $y\notin C,$ then $y$ is transient or null recurrent and, for
every $x\in E,$ 
\[
\lim_{n\to+\infty}M^{n}\left(x,y\right)=0,
\]
so, by $\refpar{eq:nu_y_dev},$ $\nu\left(y\right)=\pi\left(y\right)=0.$ 

Hence, 
\[
\nu\left(C\right)=1,
\]
and it follows from $\refpar{eq:nu_yispi_ynu_c}$ that for every $y\in C,$
\[
\nu\left(y\right)=\pi\left(y\right).
\]
 The probability $\pi$ is thus the unique invariant probability.
\item \textbf{Case 2: $C$ is periodic of period $d$}\\
If $C$ is periodic of period $d,$ let $C_{k},\,k\in\left\llbracket 0,d-1\right\rrbracket $
be the cyclic classes of $C,$ indexed similarly to Proposition $\ref{pr:per_rec_state}.$
Recall that, for every $k\in\left\llbracket 0,d-1\right\rrbracket $
and for every $x,y\in C_{k},$ 
\[
\lim_{n\to+\infty}M^{nd}\left(x,y\right)=\dfrac{d}{m\left(y\right)},
\]
where $m\left(y\right)=\mathbb{E}_{y}\left(T^{1}_{y}\right).$ Let
us prove that the measure $\pi$ defined by $\refpar{eq:exp_pi_one_positive_rec_com_cl},$
which is also
\[
\pi\left(x\right)=\begin{cases}
\dfrac{1}{m\left(x\right)}, & \text{if }x\in C,\\
0, & \text{otherwise,}
\end{cases}
\]
is an invariant probability.
\begin{itemize}
\item If $x\in C,$ then
\[
M^{nd}\left(x,x\right)=\sum_{y\in C}M^{nd-1}\left(x,y\right)M\left(y,x\right),
\]
and, by the Fatou lemma
\begin{align*}
\dfrac{d}{m\left(x\right)} & =\lim_{n\to+\infty}M^{nd}\left(x,x\right)\\
 & \geqslant\sum_{y\in C}\liminf_{n\to+\infty}M^{nd-1}\left(x,y\right)M\left(y,x\right)\\
 & \qquad\qquad=\sum^{d-1}_{k=0}\left[\sum_{y\in C_{k}}\liminf_{n\to+\infty}M^{nd-1}\left(x,y\right)M\left(y,x\right)\right].
\end{align*}
Hence, if $x\in C_{k_{0}}$ and if $k_{1}=k_{0}-1\left[d\right]$
with $k_{1}\in\left\llbracket 0,d-1\right\rrbracket ,$
\[
\dfrac{d}{m\left(x\right)}\geqslant\sum_{y\in C_{k_{1}}}\liminf_{n\to+\infty}M^{nd-1}\left(x,y\right)M\left(y,x\right)=\sum_{y\in C_{k_{1}}}\dfrac{d}{m\left(y\right)}M\left(y,x\right),
\]
which gives, since $M\left(y,x\right)=0$ if $y\notin C_{k_{1}},$
\[
\dfrac{d}{m\left(x\right)}\geqslant\sum^{d-1}_{k=0}\left[\sum_{y\in C_{k_{1}}}\dfrac{d}{m\left(y\right)}M\left(y,x\right)\right]=\sum_{y\in C}\dfrac{d}{m\left(y\right)}M\left(y,x\right).
\]
Hence,
\begin{equation}
\dfrac{1}{m\left(x\right)}\geqslant\sum_{y\in C}\dfrac{1}{m\left(y\right)}M\left(y,x\right).\label{eq:lower_boud_invmx}
\end{equation}
Let us prove that, in fact, for every $x\in C,$ equality holds in
the inequality $\refpar{eq:lower_boud_invmx}.$ Suppose, on the contrary,
that equality fails. Then there would exist $x_{0}$ such that
\[
\dfrac{1}{m\left(x_{0}\right)}>\sum_{y\in C}\dfrac{1}{m\left(y\right)}M\left(y,x_{0}\right).
\]
Summing the inequality over $x\in C,$ we would obtain
\[
\sum_{x\in C}\dfrac{1}{m\left(x\right)}>\sum_{x\in C}\left[\sum_{y\in C}\dfrac{1}{m\left(y\right)}M\left(y,x\right)\right]=\sum_{y\in C}\left[\dfrac{1}{m\left(y\right)}\sum_{x\in C}M\left(y,x\right)\right]=\sum_{y\in C}\dfrac{1}{m\left(y\right)},
\]
which is absurd.\\
Since $\pi$ is supported on $C,$ it follows that, for every $x\in C,$
\[
\dfrac{1}{m\left(x\right)}=\sum_{y\in C}\dfrac{1}{m\left(y\right)}M\left(y,x\right)=\sum_{y\in E}\pi\left(y\right)M\left(y,x\right)=\pi M\left(x\right).
\]
\item If $x\notin C.$ For every $y\notin C,$ we have $\pi\left(y\right)=0.$
Since $C$ is a closed class, we also have $M\left(y,x\right)=0$
for every $y\in C.$ Hence, for every $y\in E,$ 
\[
\pi\left(y\right)M\left(y,x\right)=0,
\]
and therefore
\[
\pi M\left(x\right)=\pi\left(x\right)=0
\]
This proves that $\pi$ is an invariant measure.\\
It remains to prove that $\pi$ is a probability. Yet, since the chain
restricted to $C_{k}$ with transition matrix $M^{d}_{\mid C_{k}\times C_{k}}$
is aperiodic, the previous argument implies that $d\pi_{\mid C_{k}}$
is the unique invariant probability for this restricted chain. It
follows that $\pi\left(C\right)=1.$
\end{itemize}
\item \textbf{Uniqueness}\\
We now prove the \textbf{uniqueness}. Let $\nu$ be an invariant probability
on $E.$ Then, for every $n\in\mathbb{N}^{\ast}$ and every $x\in C,$
\[
\nu\left(x\right)=\sum^{d-1}_{k=0}\sum_{y\in C_{k}}\nu\left(y\right)M^{nd}\left(y,x\right)+r_{n}\left(x\right),
\]
where
\[
r_{n}\left(x\right)=\sum_{y\notin C}\nu\left(y\right)M^{nd}\left(y,x\right).
\]
If $x\in C_{k_{0}},$ then
\[
\lim_{n\to+\infty}M^{nd}\left(y,x\right)=\begin{cases}
\dfrac{d}{m\left(x\right)}, & \text{if }y\in C_{k_{0}},\\
0, & \text{otherwise.}
\end{cases}
\]
Therefore, by applying the dominated convergence theorem to each sum
over $C_{k},$ we obtain the convergence of the sequence with general
term $r_{n}\left(x\right)$ and
\[
\nu\left(x\right)=\sum_{y\in C_{k_{0}}}\left[\nu\left(y\right)\dfrac{d}{m\left(x\right)}\right]+\lim_{n\to+\infty}r_{n}\left(x\right).
\]
Hence,
\begin{equation}
\nu\left(x\right)=d\pi_{\mid_{C_{k_{0}}}}\left(x\right)\nu\left(C_{k_{0}}\right)+\lim_{n\to+\infty}r_{n}\left(x\right).\label{eq:eq17_72}
\end{equation}
Summing over $x\in C_{k_{0}}$ yields
\[
\nu\left(C_{k_{0}}\right)=d\pi_{\mid_{C_{k_{0}}}}\left(C_{k_{0}}\right)\nu\left(C_{k_{0}}\right)+\sum_{x\in C_{k_{0}}}\lim_{n\to+\infty}r_{n}\left(x\right).
\]
Since 
\[
d\pi_{\mid_{C_{k_{0}}}}\left(C_{k_{0}}\right)=1,
\]
it follows that
\[
\sum_{x\in C_{k_{0}}}\lim_{n\to+\infty}r_{n}\left(x\right)=0
\]
and therefore, for every $x\in C_{k_{0}},$ 
\[
\lim_{n\to+\infty}r_{n}\left(x\right)=0.
\]
Using $\refpar{eq:eq17_72},$ we conclude that for every $k_{0}$
such that $0\leqslant k_{0}\leqslant d-1$ and for every $x\in C_{k_{0}},$
\begin{equation}
\nu\left(x\right)=d\pi_{\mid_{C_{k_{0}}}}\left(x\right)\nu\left(C_{k_{0}}\right).\label{eq:eq17_73}
\end{equation}
Moreover, if $y\notin C,$ then $y$ is not positive recurrent and,
by the same reasonning, 
\[
\nu\left(y\right)=\sum_{x\in E}\nu\left(x\right)\left[\lim_{n\to+\infty}M^{n}\left(x,y\right)\right]=0.
\]
This proves that $\nu$ is supported on $C.$\\
Lastly, since $\nu$ is invariant, let $x\in C_{k_{0}}$ and set $k_{1}=k_{0}-1\left[d\right]$
with $0\leqslant k_{1}\leqslant d-1.$ By $\refpar{eq:eq17_73},$
\[
\nu\left(x\right)=d\pi_{\mid_{C_{k_{0}}}}\left(x\right)\nu\left(C_{k_{0}}\right)=\sum^{d-1}_{k=0}\sum_{y\in C_{k}}\nu\left(y\right)M\left(y,x\right)=\sum_{y\in C_{k_{1}}}\nu\left(y\right)M\left(y,x\right).
\]
Using again $\refpar{eq:eq17_73},$ we obtain
\[
\nu\left(x\right)=\sum_{y\in C_{k_{1}}}\left[d\pi_{\mid_{C_{k_{1}}}}\left(y\right)\nu\left(C_{k_{1}}\right)\right]M\left(y,x\right)=d\nu\left(C_{k_{1}}\right)\sum_{y\in C_{k_{1}}}\dfrac{1}{m\left(y\right)}M\left(y,x\right).
\]
Since $x\in C_{k_{0}}$ and $\pi$ is supported on $C,$ this can
also be written as
\[
\nu\left(x\right)=d\nu\left(C_{k_{1}}\right)\sum_{y\in C}\pi\left(y\right)M\left(y,x\right)=d\nu\left(C_{k_{1}}\right)\pi M\left(x\right).
\]
Taking into account the invariance of $\pi,$ we conclude that
\[
\nu\left(x\right)=d\nu\left(C_{k_{1}}\right)\pi\left(x\right).
\]
Summing over $x\in C_{k_{0}}$ yields
\[
\nu\left(C_{k_{0}}\right)=d\nu\left(C_{k_{1}}\right)\pi\left(C_{k_{0}}\right).
\]
Since $d\pi_{\mid_{C_{k_{0}}}}$ is a probability, we then proved
that for every $k_{0}$ such that $0\leqslant k_{0}\leqslant d-1,$
\[
\nu\left(C_{k_{0}}\right)=\nu\left(C_{k_{1}}\right).
\]
It follows that $\nu=\pi.$
\end{itemize}
3. \textbf{Suppose that $X$ has $N,$ positive recurrent classes,
with $N\geqslant2,$ say} $C^{1},C^{2},\cdots,C^{N}.$ 

Define, for every nonnegative real numbers $a_{j},\,j\in\left\llbracket 1,N\right\rrbracket $
such that $\sum^{N}_{j=1}a_{j}=1,$ the measure $\mu$ by
\[
\mu\left(x\right)=\begin{cases}
\dfrac{a_{j}}{m\left(x\right)}, & \text{if }x\in C^{j},\,j\in\left\llbracket 1,N\right\rrbracket ,\\
0, & \text{if }x\notin\bigcup^{N}_{j=1}C^{j}.
\end{cases}
\]
For any $j_{0}\in\left\llbracket 1,N\right\rrbracket $ and for every
$x\in C^{j_{0}},$ since the classes $C^{1},C^{2},\cdots,C^{N}$ are
closed,
\begin{align}
\mu M\left(x\right) & =\sum_{y\in E}\mu\left(y\right)M\left(y,x\right)=\sum^{N}_{j=1}a_{j}\left[\sum_{y\in C^{j}}\dfrac{1}{m\left(y\right)}M\left(y,x\right)\right]\nonumber \\
 & =a_{j_{0}}\sum_{y\in C^{j_{0}}}\dfrac{1}{m\left(y\right)}M\left(y,x\right).\label{eq:eq17_74}
\end{align}
Since the chain restricted to the closed class $C^{j_{0}}$ is, by
nature, irreducible, the result proved in Point 2 ensures that the
measure defined on $C^{j_{0}}$ by $\mu_{j_{0}}\left(y\right)=\dfrac{1}{m\left(y\right)}$
for every $y\in C^{j_{0}}$ is the unique invariant probability of
the restricted chain. 

Substituting in $\refpar{eq:eq17_74},$ we obtain the equality
\[
\mu M\left(x\right)=a_{j_{0}}\mu_{j_{0}}\left(x\right)=a_{j_{0}}\dfrac{1}{m\left(x\right)},
\]
which proves that, for every $x\in C^{j_{0}},$ 
\[
\mu M\left(x\right)=\mu\left(x\right).
\]

Moreover, if $x\notin\bigcup^{N}_{j=1}C^{j},$ then for for every
$j\in\left\llbracket 1,N\right\rrbracket ,$ for every $y\in C^{j},$
$M\left(y,x\right)=0,$ which implies that $\mu M\left(x\right)=0.$
Since also $\mu\left(x\right)=0,$ we still have 
\[
\mu M\left(x\right)=\mu\left(x\right)
\]

Thus $\mu$ is an invariant measure. In fact, it is a probability.
Indeed, since for every $j$ each $\mu_{j}$ is a probability, 
\[
\sum_{x\in C^{j}}\dfrac{1}{m\left(x\right)}=1,
\]
and thus
\[
\mu\left(E\right)=\sum^{N}_{j=1}a_{j}\left[\sum_{x\in C^{j}}\dfrac{1}{m\left(x\right)}\right]=\sum^{N}_{j=1}a_{j}=1.
\]
Hence, we have proved that, in this case, there exists an uncountable
infinity of invariant probabilities.

\end{proof}

We return to the diffusion model of \textbf{Ehrenfest} and determine
the \textbf{invariant probability} of the associated chain.

\begin{example}{Ehrenfest Heat Diffusion Model---Continued}{}

We return to the Ehrenfest model of heat diffusion, described in terms
of drawing balls from an urn---see Examples $\ref{ex:ehrenfest_model}$
and $\ref{ex:ehrenfest_model_2},$ whose notation we keep. Recall
that $X_{n}$ represents the number of red balls in the urn at time
$n,$ and that the process $X=\left(X_{n}\right)_{n\in\mathbb{N}}$
is a homogeneous Markov chain taking values in the integer interval
$E=\left\llbracket 0,m\right\rrbracket ,$ whose transition matrix
$M$ is given, for every $k\in E,$ by
\[
M\left(k,k+1\right)=p_{k},\,\,\,\,M\left(k,k-1\right)=q_{k},
\]
where
\[
p_{k}=1-\dfrac{k}{m},\,\,\,\,q_{k}=\dfrac{k}{m}.
\]
This chain is irreducible, aperiodic, positive recurrent---the latter
follows from Corollary $\ref{co:closed_com_states}.$ Therefore, by
Theorem $\ref{th:positive_recurrence_crit},$ $X$ admits a unique
invariant probability. We now aim to determine it.

\end{example}

\begin{solutionexample}{}{}

We first look for an invariant measure $\mu.$ Such a measure is solution
of the equation system
\[
\mu\left(y\right)=\sum_{x\in E}\mu\left(x\right)M\left(x,y\right).
\]
Writting $\mu_{k}$ for $\mu\left(k\right),$ this becomes
\begin{equation}
\begin{cases}
\mu_{k}=\mu_{k-1}p_{k-1}+\mu_{k+1}q_{k+1}, & \text{if }1\leqslant k\leqslant m-1,\\
\mu_{0}=\mu_{1}q_{1},\\
\mu_{m}=\mu_{m-1}p_{m-1}.
\end{cases}\label{eq:def_mu_k}
\end{equation}
Since $p_{k+1}+q_{k+1}=1,$ for $1\leqslant k\leqslant m-1,$ we obtain
\begin{equation}
\mu_{k+1}-\mu_{k}=\mu_{k+1}p_{k+1}-\mu_{k-1}p_{k-1,}\label{eq:eq17_76}
\end{equation}
and therefore
\[
\sum^{m-1}_{j=k}\left(\mu_{j+1}-\mu_{j}\right)=\sum^{m-1}_{j=k}\left(\mu_{j+1}p_{j+1}-\mu_{j-1}p_{j-1}\right).
\]
Equivalently, after a change of indices,
\[
\mu_{m}-\mu_{k}=\sum^{m}_{j=k+1}\left(\mu_{j}p_{j}\right)-\sum^{m-2}_{j=k-1}\left(\mu_{j}p_{j}\right).
\]
Hence,
\begin{equation}
\mu_{m}-\mu_{k}=\mu_{m}p_{m}+\mu_{m-1}p_{m-1}-\mu_{k-1}p_{k-1}-\mu_{k}p_{k}.\label{eq:eq17_77}
\end{equation}
Now, since
\[
p_{k}+q_{k}=1,\,\,\,\,p_{m}=0,\,\,\,\,\text{and}\,\,\,\,p_{m-1}=\dfrac{1}{m},
\]
we get
\[
\mu_{m}-\dfrac{1}{m}\mu_{m-1}=\mu_{m-1}=\mu_{k}q_{k}-\mu_{k-1}p_{k-1}.
\]
By $\refpar{eq:def_mu_k}$,
\[
\mu_{m}-\dfrac{1}{m}\mu_{m-1}=0,
\]
so for $1\leqslant k\leqslant m-1,$
\begin{equation}
\mu_{k}=\dfrac{p_{k-1}}{q_{k}}\mu_{k-1}.\label{eq:eq16_78}
\end{equation}
Iterating,
\begin{equation}
\mu_{k}=\dfrac{p_{k-1}p_{k-2}\cdots p_{1}p_{0}}{q_{k}q_{k-1}\cdots q_{2}q_{1}}\mu_{0}.\label{eq:eq17_79}
\end{equation}
Moreover,
\begin{align*}
\dfrac{p_{k-1}p_{k-2}\cdots p_{1}p_{0}}{q_{k}q_{k-1}\cdots q_{2}q_{1}} & =\dfrac{\frac{m-k+1}{m}\cdot\frac{m-k+2}{m}\cdots\frac{m-1}{m}\cdot\frac{m}{m}}{\frac{k}{m}\cdot\frac{k-1}{m}\cdots\frac{2}{m}\cdot\frac{1}{m}}\\
 & =\dfrac{m\left(m-1\right)\cdots\left(m-k+1\right)}{k!}=\binom{m}{k}.
\end{align*}
Thus, for every $k\in\left\llbracket 1,m-1\right\rrbracket ,$
\begin{equation}
\mu_{k}=\binom{m}{k}\mu_{0}.\label{eq:eq17_80}
\end{equation}
Moreover, by $\refpar{eq:eq17_80}$ and $\refpar{eq:def_mu_k},$
\[
\mu_{m}=\dfrac{1}{m}\mu_{m-1}=\dfrac{1}{m}m\mu_{0}=\mu_{0},
\]
so also
\[
\mu_{m}=\binom{m}{m}\mu_{0}.
\]
Hence, any invariant \textbf{measure} $\mu$ is determined, for $1\leqslant k\leqslant m,$
by
\begin{equation}
\mu_{k}=\binom{m}{k}\mu_{0}.\label{eq:eq17_81}
\end{equation}
There is therefore a unique invariant \textbf{probability} $\mu,$
determined by the equivalent relations
\[
\mu_{0}+\sum^{m}_{k=1}\binom{m}{k}\mu_{0}=1\Longleftrightarrow\left[\sum^{m}_{k=0}\binom{m}{k}\right]\mu_{0}=1,
\]
which yields
\[
\mu_{0}=\dfrac{1}{2^{m}}.
\]
The \textbf{invariant probability} $\mu$ is thus given, for every
$k\in\left\llbracket 0,m\right\rrbracket $ by\boxeq{
\[
\mu_{k}=\dfrac{\binom{m}{k}}{2^{m}}.
\]
}

In other words, the invariant probability $\mu$ is the binomial law
$\mathscr{B}\left(m,\dfrac{1}{2}\right).$

Hence, the Ehrenfest chain is \textbf{irreducible, aperiodic, positive
recurrent}, and it admits, by Theorem $\ref{th:prob_limit}$ below,
a \textbf{limiting probability} $\mu$ which is the \textbf{binomial
law} $\mathscr{B}\left(m,\dfrac{1}{2}\right).$ That is, in the \textbf{stationary
regime}, the two colours occur in equal proportions, up to random
fluctuations.

Moreover, by Theorem $\ref{th:positive_recurrence_crit},$ the \textbf{mean
return time to $k,$} starting from $k,$ is $\dfrac{1}{\mu_{k}}=\dfrac{2^{m}}{\binom{m}{k}}.$ 

\end{solutionexample}

We now give a necessary and sufficient condition for the existence
of a \textbf{limiting probability}, that is, a probability $\mu$
such that, for every $x,y\in E,$ the sequence with general term $M^{n}\left(x,y\right)$
converges to limit $\mu\left(y\right),$ independent of $x.$ For
such a probability $\mu,$ we then have, for every $y\in E,$
\[
\lim_{n\to+\infty}P_{x}\left(X_{n}=y\right)=\mu\left(y\right),
\]
\textbf{independently} of the starting point $x$ of the chain at
the initial instant.

\begin{theorem}{Necessary and Sufficient Condition for the Existence of a Limiting Probability}{prob_limit}

A homogeneous Markov chain with transition matrix $M$ has a \textbf{limiting
probability\index{limiting probability}} if and only if it admits
a unique positive recurrent aperiodic class $C$ such that $P_{x}\left(T^{1}_{y}<+\infty\right)=1,$
for every $x\in E$ and every $y\in C.$

\end{theorem}

\begin{proof}{}{}

Assume that a limiting probability $\mu$ exists, By Proposition $\ref{pr:ex_unic_inv_prob},$
it is the unique invariant probability. Then Theorem $\ref{th:positive_recurrence_crit}$
ensures that there exists a unique positive recurrent class $C.$
This class must be aperiodic. Indeed, suppose that it is periodic
with period $d,$ and denote $C_{k},\,k\in\left\llbracket 0,d-1\right\rrbracket ,$
the cyclic classes of $C,$ indexed as in Proposition $\ref{pr:per_rec_state}.$
Then, for every $x\in C_{0}$ and $y\in C_{1},$ we would have
\[
\lim_{n\to+\infty}M^{nd+1}\left(x,y\right)=\dfrac{d}{\mathbb{E}_{y}\left(T^{1}_{y}\right)}>0
\]
whereas for every $n\in\mathbb{N}^{\ast},$ $M^{nd}\left(x,y\right)=0.$
which is in contradiction with the existence of a limiting probability.
Hence $C$ is aperiodic.

Finally, for every $x\in E$ and $y\in C,$
\[
\lim_{n\to+\infty}M^{n}\left(x,y\right)=\dfrac{P_{x}\left(T^{1}_{y}<+\infty\right)}{\mathbb{E}_{y}\left(T^{1}_{y}\right)}=\mu\left(y\right)>0,
\]
Thus the function $x\mapsto P_{x}\left(T^{1}_{y}<+\infty\right)$
is constant. Since $y$ is recurrent, 
\[
P_{y}\left(T^{1}_{y}<+\infty\right)=1.
\]
and therefore $P_{x}\left(T^{1}_{y}<+\infty\right)=1,$ for every
$x\in E$ and $y\in C.$

Conversely, suppose there exists a unique positive recurrent aperiodic
class $C$ such that $P_{x}\left(T^{1}_{y}<+\infty\right)=1,$ for
every $x\in E$ and $y\in C.$ 

Then, for such $x$ and $y,$
\[
\lim_{n\to+\infty}M^{n}\left(x,y\right)=\dfrac{1}{\mathbb{E}_{y}\left(T^{1}_{y}\right)}>0.
\]

Moreover, if $x\in E$ and $y\notin C,$ then $y$ is either null
recurrent or transient and 
\[
\lim_{n\to+\infty}M^{n}\left(x,y\right)=0.
\]
Hence, for every $x,y\in E,$ the sequence with general term $M^{n}\left(x,y\right)$
converges to a limit $\pi\left(y\right)$ independent of $x.$ Since
the measure $\pi$ thus defined is not identically zero, it follows
by Proposition $\ref{pr:ex_unic_inv_prob}$ that $\pi$ is the unique
invariant probability.

\end{proof}

We conclude this section by studying a genetic model.

\begin{example}{Genetic Model}{genetic_model}

A hereditary trait of an individual generally depends on the presence,
in its genetic makeup, of genes of two types $G$ and $g,$ called
\textbf{alleles}\footnotemark. These appear in pairs: $GG,gg,Gg$
and $gG,$ the last two being genetically identical. Thus, to characterize
the trait, we consider only the unordered pairs $GG,\,gg,\,Gg$ called
genotypes. The gene $G$ is often \textbf{predominant}; therefore,
the genotypes $Gg$ and $GG$ produce the same hereditary trait, called
the \textbf{phenotype}\footnotemark. Depending if an individual has
the genotype $GG,\,gg$ or $Gg$ it is said \textbf{dominant, recessive}
or \textbf{heterozygous}.

An individual receives one gene independently and at random from each
parent.
\begin{itemize}
\item If both parents are dominant---respectively recessive---then the
individual is also dominant---respectively recessive.
\item If one parent is dominant and the other recessive, then the individual
is heterozygous.
\item If one of the parent is dominant and the other heterozygous, the individual
receives the gene $G$ from the dominant parent and receives either
a gene $G$ or $g$ from the heterozygous parent, with equal probability.
Hence, the individual has same probability to be dominant or heterozygous.\\
Similarly, if one parent is recessive and the other heterozygous,
the individual has equal probability to be recessive or heterozygous.
\item If both parents are heterozygous, the individual independently receives
from each parent the gene $G$ or $g$ with same probability. Hence,
the individual will be dominant with the probability $\dfrac{1}{4},$
recessive also with probability $\dfrac{1}{4}$ and heterozygous with
the probability $\dfrac{1}{2}.$
\end{itemize}
Now, consider the following process: an individual with a given character
partners with an heterozygous individual and has children. One of
their children is chosen at random; this child partners an heterozygous
individual and has children, and so forth. Denote by $X_{n}$ the
genetic type of the observed descendant in the $n-$th generation.
The process $\left(X_{n}\right)_{n\in\mathbb{N}^{\ast}}$ is a Markov
chain with transition matrix
\[
\begin{array}{ccc}
 &  & \begin{array}{ccc}
\boldsymbol{GG} & \boldsymbol{Gg} & \boldsymbol{gg}\end{array}\\
M_{h}= & \begin{array}{c}
\boldsymbol{GG}\\
\boldsymbol{Gg}\\
\boldsymbol{gg}
\end{array} & \left(\begin{array}{ccc}
\dfrac{1}{2} & \dfrac{1}{2} & 0\\
\,\dfrac{1}{4}\, & \,\dfrac{1}{2}\, & \,\dfrac{1}{4}\,\\
0 & \dfrac{1}{2} & \dfrac{1}{2}
\end{array}\right)
\end{array}
\]

The associated graph to this Markov chain is

\begin{center}\begin{tikzpicture}[
  scale=0.75,
  transform shape,
  >=Stealth,
  every node/.style={
    circle,
    draw=black,
    minimum size=12mm,
    inner sep=1pt
  },
  edge/.style={
    ->,
    draw=black,
    line width=0.75pt
  }
]

% Nodes (horizontal layout)
\node (GG) at (0,0) {\(GG\)};
\node (Gg) at (4,0) {\(Gg\)};
\node (gg) at (2,-4) {\(gg\)};

% Self-loops
\draw[edge] (GG) to[out=210,in=150,looseness=8] (GG);
\draw[edge] (Gg) to[out=30,in=-30,looseness=8] (Gg);
\draw[edge] (gg) to[out=300,in=240,looseness=8] (gg);

% Transitions between states
\draw[edge] (GG) to[bend left=20] (Gg);

\draw[edge] (Gg) to[bend left=20] (GG);
\draw[edge] (Gg) to[bend left=20] (gg);

\draw[edge] (gg) to[bend left=20] (Gg);

\end{tikzpicture}
\end{center}

The chain is \textbf{irreducible, aperiodic, and positive recurrent}.
Hence, it admits a unique invariant probability law. Let us compute
it.

\end{example}

\addtocounter{footnote}{-1}

\footnotetext{From the Greek, allêlon, meaning ``one another'',
an \textbf{allele}\index{allele}, or allelomorph, designates a hereditary
variant that contrasts with another---for instance, pea seeds that
are smooth or wrinkled---or the gene responsible for this trait.}

\stepcounter{footnote}

\footnotetext{A \textbf{phenotype\index{phenotype}} denotes the
observable characteristics of a living organism, as opposed to the
genotype---its hereditary genetic constitution. Several different
genotypes may produce the same observable appearance of an individual,
although differences may reappear in later generations.}

\begin{solutionexample}{}{}

Let us first determine the invariant measures $v=\left(a,b,c\right).$
The vector $v$ is a left eigenvector of $M_{h}$ associated with
the eigenvalue 1. Hence it satisfies the equation $vM_{h}=v$ which
corresponds to the system
\[
\left\{ \begin{array}{l}
\dfrac{1}{2}a+\dfrac{1}{4}b\,\,\,\,\,\,\,\,\,\,=a\\
\dfrac{1}{2}a+\dfrac{1}{4}b+\dfrac{1}{2}c\,=b\\
\,\,\,\,\,\,\,\,\,\dfrac{1}{4}b+\dfrac{1}{2}c\,=c
\end{array}\right..
\]
This system has solutions of the form $\left(a,2a,a\right).$ Thus
the invariant measures are given by $\left(a,2a,a\right),$ with $a\geqslant0,$
arbitrary. Consequently, there exists a unique invariant probability
$v_{0}.$ It is determined by the normalization condition $a+2a+a=1,$
which yields $v_{0}=\left(\dfrac{1}{4},\dfrac{1}{2},\dfrac{1}{4}\right).$

Let $E=\left\{ GG,gg,Gg\right\} $ denote the space of states. Since
the chain is irreducible and aperiodic, it follows that, for every
$x,y\in E,$
\[
\lim_{n\to+\infty}M^{n}_{h}\left(x,y\right)=\dfrac{1}{E_{y}\left(T^{1}_{y}\right)}=v_{0}\left(y\right).
\]
In particular, we obtain the mean return times in a point:
\[
\mathbb{E}_{GG}\left(T^{1}_{GG}\right)=\mathbb{E}_{gg}\left(T^{1}_{gg}\right)=4\,\,\,\,\text{and}\,\,\,\,\mathbb{E}_{Gg}\left(T^{1}_{Gg}\right)=2.
\]
Now, instead to partner the selected descendant with an heterozygous
individual, suppose that the descendant partners with a dominant.
Then the process $\left(X_{n}\right)_{n\in\mathbb{N}^{\ast}}$ is
a Markov chain with transition matrix
\[
\begin{array}{ccc}
 &  & \begin{array}{ccc}
\boldsymbol{GG} & \boldsymbol{Gg} & \boldsymbol{gg}\end{array}\\
M_{h}= & \begin{array}{c}
\boldsymbol{GG}\\
\boldsymbol{Gg}\\
\boldsymbol{gg}
\end{array} & \left(\begin{array}{ccc}
1 & 0 & 0\\
\,\dfrac{1}{2}\, & \,\dfrac{1}{2}\, & \,\,0\,\,\\
0 & 1 & 0
\end{array}\right)
\end{array}
\]

The graph associated to this Markov chain is

\begin{center}\begin{tikzpicture}[
  scale=0.75,
  transform shape,
  >=Stealth,
  every node/.style={
    circle,
    draw=black,
    minimum size=12mm,
    inner sep=1pt
  },
  edge/.style={
    ->,
    draw=black,
    line width=0.75pt
  }
]

% Nodes (horizontal layout)
\node (GG) at (0,0) {\(GG\)};
\node (Gg) at (3,0) {\(Gg\)};
\node (gg) at (6,0) {\(gg\)};

% Self-loops
\draw[edge] (GG) to[out=120,in=60,looseness=8] (GG);
\draw[edge] (Gg) to[out=120,in=60,looseness=8] (Gg);

% Transitions
\draw[edge] (Gg) -- (GG);
\draw[edge] (gg) -- (Gg);

\end{tikzpicture}
\end{center}

Each states forms its own communication class. The state $gg$ is
\textbf{inessential}, the states $gg$ and $Gg$ are \textbf{transient},
and $GG$ is an \textbf{absorbing state}. 

\end{solutionexample}

\section{Strong Law of Large Numbers}\label{sec:Strong-Law-of}

We now state a theorem of strong law of large numbers for the case
of a homogeneous Markov chain. We then apply it, in the finite case,
to the estimation of its transition matrix.

\subsection{Strong Law Theorem}

\begin{theorem}{Chacon-Orstein Theorem}{chacon_orstein_th}

Let $X$ be a process such that, for every $x\in E,$ it is a homogeneous
Markov chain on the underlying filtered probabilized space $\left(\Omega,\mathscr{A},\left(\mathscr{A}_{n}\right)_{n\in\mathbb{N}},P_{x}\right)$
with initial law $\delta_{x}$ and transition matrix $M.$ Suppose
that $X$ admits a unique positive recurrent class $C$---recall
that in this case there exists a unique invariant probability $\pi.$
Moreover\footnotemark assume that there exists $y\in C$ such that,
for every $x\in E,$ $P_{x}\left(T^{1}_{y}<+\infty\right)=1.$

Let $f$ and $g$ be two $\pi-$integrable functions defined on $E.$
Suppose that $g$ does not take the value 0. Then, for every $x\in E,$
the sequence with general term
\[
\dfrac{\sum^{n}_{j=1}f\left(X_{j}\right)}{\sum^{n}_{j=1}g\left(X_{j}\right)}
\]
converges $P_{x}-$almost surely, and
\[
\lim_{n\to+\infty}\dfrac{\sum^{n}_{j=1}f\left(X_{j}\right)}{\sum^{n}_{j=1}g\left(X_{j}\right)}=\dfrac{\sum_{x\in E}f\left(x\right)\pi\left(x\right)}{\sum_{x\in E}g\left(x\right)\pi\left(x\right)}\,\,\,\,P_{x}-\text{almost surely.}
\]

\end{theorem}

\footnotetext{If $E$ is finite, this hypothesis is automatically
satisfied.}

\begin{proof}{}{}
\begin{itemize}
\item \textbf{Recurrence of $y$}\\
For every $x\in E,$ $P_{x}\left(R_{y}\right)=1.$ Indeed, a trivial
modification of the proof of Lemma $\ref{lm:rec_vis_and_proba}$ shows
that, for every $p\in\mathbb{N}^{\ast},$
\[
P_{x}\left(T^{p+1}_{y}<+\infty\right)=P_{x}\left(T^{1}_{y}<+\infty\right)\left[P_{y}\left(T^{1}_{y}<+\infty\right)\right]^{p}.
\]
By the hypothesis, it follows that
\[
P_{x}\left(R_{y}\right)=\lim_{n\to+\infty}\searrow P_{x}\left(T^{1}_{y}<+\infty\right)\left[P_{y}\left(T^{1}_{y}<+\infty\right)\right]^{p}=1.
\]
\item \textbf{Construction of an invariant measure}\\
Let $\mu$ be the measure on $E$ defined, for every $x\in E,$ by
\[
\mu\left(x\right)=\mathbb{E}_{y}\left[\sum^{T^{1}_{y}}_{n=1}\boldsymbol{1}_{\left(X_{n}=x\right)}\right],
\]
that is, the mean number of visits to $x$ before the first return
to $y$ by the chain starting from $y$ at the initial instant. We
show that $\mu$ is an invariant measure. Indeed,
\begin{align*}
\mu M\left(x\right) & =\sum_{z\in E}\mu\left(z\right)M\left(z,x\right)=\sum_{z\in E}\left[\mathbb{E}_{y}\left(\sum^{T^{1}_{y}}_{n=1}\boldsymbol{1}_{\left(X_{n}=z\right)}\right)\right]M\left(z,x\right)\\
 & =\mathbb{E}_{y}\left[\sum^{T^{1}_{y}}_{n=1}\left(\sum_{z\in E}\boldsymbol{1}_{\left(X_{n}=z\right)}M\left(z,x\right)\right)\right]\\
 & =\mathbb{E}_{y}\left[\sum^{T^{1}_{y}}_{n=1}M\left(X_{n},x\right)\right],
\end{align*}
By partitioning,
\[
\mu M\left(x\right)=\mathbb{E}_{y}\left[\boldsymbol{1}_{\left(T^{1}_{y}=1\right)}M\left(X_{1},x\right)\right]+\mathbb{E}_{y}\left[\boldsymbol{1}_{\left(T^{1}_{y}\geqslant2\right)}\left(\sum^{T^{1}_{y}}_{n=1}M\left(X_{n},x\right)+M\left(X_{T^{1}_{y}},x\right)\right)\right].
\]
Since $X_{T^{1}_{y}}=X_{0}=y$ $P_{y}-$almost surely, we obtain
\[
\mu M\left(x\right)=\mathbb{E}_{y}\left[\boldsymbol{1}_{\left(T^{1}_{y}=1\right)}M\left(X_{0},x\right)\right]+\mathbb{E}_{y}\left[\boldsymbol{1}_{\left(T^{1}_{y}\geqslant2\right)}\sum^{T^{1}_{y}}_{n=0}M\left(X_{n},x\right)\right].
\]
Thus
\[
\mu M\left(x\right)=\mathbb{E}_{y}\left[\sum^{T^{1}_{y}-1}_{n=0}M\left(X_{n},x\right)\right].
\]
Using the Markov property, this equality can be written as
\[
\mu M\left(x\right)=\sum^{+\infty}_{n=0}\mathbb{E}_{y}\left[\boldsymbol{1}_{\left(n<T^{1}_{y}\right)}P_{X_{n}}\left(X_{1}=x\right)\right]=\sum^{+\infty}_{n=0}\mathbb{E}_{y}\left[\boldsymbol{1}_{\left(n<T^{1}_{y}\right)}\mathbb{E}^{\mathscr{A}_{n}}_{y}\boldsymbol{1}_{\left(X_{n+1}=x\right)}\right].
\]
Since $\left(n<T^{1}_{y}\right)\in\mathscr{A}_{n},$ we obtain
\[
\mu M\left(x\right)=\sum^{+\infty}_{n=0}\mathbb{E}_{y}\left[\boldsymbol{1}_{\left(n<T^{1}_{y}\right)}\boldsymbol{1}_{\left(X_{n+1}=x\right)}\right]=\mathbb{E}_{y}\left[\sum^{T^{1}_{y}}_{n=1}\boldsymbol{1}_{\left(X_{n}=x\right)}\right]=\mu\left(x\right).
\]
This proves that $\mu$ is an invariant measure. The measure is bounded.
Since $E$ is countable and $y$ is positive recurrent, 
\[
\mu\left(E\right)=\sum_{x\in E}\mathbb{E}_{y}\left[\sum^{T^{1}_{y}}_{n=1}\boldsymbol{1}_{\left(X_{n}=x\right)}\right]=\mathbb{E}_{y}\left[\sum^{T^{1}_{y}}_{n=1}\left(\sum_{x\in E}\boldsymbol{1}_{\left(X_{n}=x\right)}\right)\right]=\mathbb{E}_{y}\left[T^{1}_{y}\right]<+\infty.
\]
Therefore $\dfrac{\mu}{\mathbb{E}_{y}\left[T^{1}_{y}\right]}$ is
an invariant probability. By uniqueness of the invariant probability,
we must have
\[
\pi=\dfrac{\mu}{\mathbb{E}_{y}\left[T^{1}_{y}\right]}.
\]
This provides an intuitive interpretation of the invariant probability.\\
If $f$ is a nonnegative function on $E,$ a direct computation gives
\begin{equation}
\intop f\text{d}\mu=\sum_{x\in E}f\left(x\right)\mu\left(x\right)=\mathbb{E}_{y}\left[\sum^{T^{1}_{y}}_{n=1}f\left(X_{n}\right)\right]=\mathbb{E}_{y}\left[\sum^{T^{1}_{y}-1}_{n=0}f\left(X_{n}\right)\right].\label{eq:eq17_82}
\end{equation}
If $f$ has arbitrary sign, one obtains classically an integrability
criterion. The formulas $\refpar{eq:eq17_82}$ remain valid for $\mu-$integrable
functions.
\item \textbf{Decomposition along successive visits to $y$}\\
The principle of the proof is to decompose the sums under consideration
according to successive return times to $y.$ The resulting pieces
are \textbf{independent and with same law}. We then apply the law
of large numbers to these independent random variables.\\
Let $f$ be a $\mu-$integrable function. For every $p\in\mathbb{N},$
define the random variable $\mathbb{Z}_{p}$ by
\[
Z_{p}=\begin{cases}
\sum^{T^{p+1}_{y}-1}_{n=T^{p}_{y}}f\left(X_{n}\right), & \text{on\,}\left(T^{p}_{y}<+\infty\right),\\
0, & \text{on\,}\left(T^{p}_{y}=+\infty\right),
\end{cases}
\]
---recalling that $T^{0}_{y}=0.$ \\
Since, for every $x\in E$ and for every $p\in\mathbb{N}^{\ast},$
\[
P_{x}\left(T^{p}_{y}<+\infty\right)=1,
\]
the \textbf{random variables $Z_{p}$ are $P_{x}-$almost surely finite}.
\\
We now prove that, for every $x\in E,$ the variables $Z_{p}$ are
\textbf{$P_{x}-$independent and with same law}\footnotemark. Indeed,
let, for each $p\in\mathbb{N}^{\ast},$ $A_{p}$ be an arbitrary Borel
subset of $\mathbb{R}.$ Since, $P_{x}-$almost surely,
\[
T^{p+1}_{y}=T^{p}_{y}+\tau^{1}_{y}\left[\theta_{T^{p}_{y}}\left(X\right)\right],
\]
the strong Markov property yields
\[
\mathbb{E}^{\mathscr{A}_{T^{p}_{y}}}_{x}\left[\boldsymbol{1}_{\text{\ensuremath{\left(Z_{p}\in A_{p}\right)}}}\right]=\mathbb{E}_{X_{T^{p}_{y}}}\left[\boldsymbol{1}_{\left(\sum^{T^{1}_{y}-1}_{n=0}f\left(X_{n}\right)\in A_{p}\right)}\right]=P_{y}\left(Z_{0}\in A_{p}\right).
\]
Consequenly, for every $N,$ a classical conditioning argument gives
\begin{align*}
\mathbb{E}_{x}\left[\prod^{N}_{p=1}\boldsymbol{1}_{\left(Z_{p}\in A_{p}\right)}\right] & =\mathbb{E}_{x}\left[\left(\prod^{N-1}_{p=1}\boldsymbol{1}_{\left(Z_{p}\in A_{p}\right)}\right)\mathbb{E}^{\mathscr{A}_{T^{N-1}_{y}}}_{x}\left[\boldsymbol{1}_{\text{\ensuremath{\left(Z_{N}\in A_{N}\right)}}}\right]\right]\\
 & =\mathbb{E}_{x}\left[\prod^{N-1}_{p=1}\boldsymbol{1}_{\left(Z_{p}\in A_{p}\right)}\right]P_{y}\left(Z_{0}\in A_{N}\right).
\end{align*}
By iterating this identity backward, we obtain
\[
\mathbb{E}_{x}\left[\prod^{N}_{p=1}\boldsymbol{1}_{\left(Z_{p}\in A_{p}\right)}\right]=\prod^{N}_{p=1}P_{y}\left(Z_{0}\in A_{p}\right).
\]
This shows that the random variables $Z_{p},\,p\in\mathbb{N}^{\ast},$
have the same law under $P_{x}$ as $Z_{0}$ under $P_{y}$ and that
the $Z_{p}$ are $P_{x}-$independent.
\item \textbf{Integrability of $Z_{1}$}\\
We now show that $Z_{1}$ is $P_{x}-$integrable. The same argument
will then apply to every $Z_{p}.$
\[
\mathbb{E}_{x}\left(\left|Z_{1}\right|\right)=\mathbb{E}_{x}\left[\left|\sum^{T^{2}_{y}-1}_{n=T^{1}_{y}}f\left(X_{n}\right)\right|\right]=\mathbb{E}_{x}\left[\mathbb{E}^{\mathscr{A}_{T^{1}_{y}}}_{x}\left(\left|\sum^{T^{2}_{y}-1}_{n=T^{1}_{y}}f\left(X_{n}\right)\right|\right)\right].
\]
Using the strong Markov property, we obtain
\[
\mathbb{E}_{x}\left(\left|Z_{1}\right|\right)=\mathbb{E}_{x}\left[\mathbb{E}_{X_{T^{1}_{y}}}\left(\left|\sum^{T^{1}_{y}-1}_{n=0}f\left(X_{n}\right)\right|\right)\right]=\mathbb{E}_{y}\left[\left|\sum^{T^{1}_{y}-1}_{n=0}f\left(X_{n}\right)\right|\right].
\]
It follows that
\[
\mathbb{E}_{x}\left(\left|Z_{1}\right|\right)\leqslant E_{y}\left[\sum^{T^{1}_{y}-1}_{n=0}\left|f\left(X_{n}\right)\right|\right]=\intop\left|f\right|\text{d}\mu<+\infty.
\]
A similar computation shows that $\mathbb{E}_{x}\left(Z_{1}\right)=\intop f\text{d}\mu.$
\item It follows, by the second strong law of large numbers applied to the
independent random variables, that
\begin{equation}
\dfrac{1}{n}\sum^{n-1}_{p=1}Z_{p}=\dfrac{1}{n}\sum^{T^{n}_{y}-1}_{k=T^{1}_{y}}f\left(X_{k}\right)\stackrel[n\to+\infty]{P-\text{a.s.}}{\longrightarrow}\intop f\text{d}\mu.\label{eq:eq17_83}
\end{equation}
Let 
\[
v\left(n\right)=\sum^{n}_{j=1}\boldsymbol{1}_{\left(X_{j}=y\right)}
\]
be the nondecreasing sequence of random integers, representing the
number of visits to $y$ up to time $n.$ \\
By hypothesis,
\[
P_{x}\left(\lim_{n\to+\infty}v\left(n\right)=+\infty\right)=P_{x}\left(R_{y}\right)=1.
\]
Moreover, by definition of $v\left(n\right),$
\[
T^{v\left(n\right)}_{y}\leqslant n<T^{v\left(n\right)+1}_{y}.
\]
If, in addition, $f$ is nonnegative, we obtain the inequalities
\[
\dfrac{\sum^{T^{v\left(n\right)}_{y}}_{k=0}f\left(X_{k}\right)}{v\left(n\right)}\leqslant\dfrac{\sum^{n}_{k=0}f\left(X_{k}\right)}{\sum^{n}_{j=1}\boldsymbol{1}_{\left(X_{j}=y\right)}}\leqslant\dfrac{\sum^{T^{v\left(n\right)+1}_{y}}_{k=0}f\left(X_{k}\right)}{v\left(n\right)}.
\]
By $\refpar{eq:eq17_83},$ the extreme terms converge $P_{x}-$almost
surely. Therefore, the middle term also converges $P_{x}-$almost
surely. In this case, the theorem follows immediately by recalling
that $\pi$ is proportional to $\mu.$ Finally, the theorem extends
to the general case where $f$ has arbitrary sign by decomposing $f$
in its positive and negative parts.
\end{itemize}
\end{proof}

\footnotetext{Using this remark, one can deduce a central limit theorem
for homogeneous Markov chains satisfying the hypotheses of the present
theorem from a central limit theorem for sequences of independent
random variables.

}

In particular, we obtain the traditional formulation of the strong
law of large numbers for homogeneous Markov chains.

\begin{corollary}{Strong Law of Large Numbers For Homogeneous Markov Chains}{strong_law_large_numbers_markov}

Under the hypotheses of the Chacon-Orstein theorem $\ref{th:chacon_orstein_th},$
for every $\pi-$integrable function $f,$\boxeq{
\[
\dfrac{1}{n}\sum^{n}_{j=1}f\left(X_{j}\right)\stackrel[n\to+\infty]{P_{x}-\text{a.s.}}{\longrightarrow}\intop f\text{d}\pi.
\]
}

\end{corollary}

\begin{proof}{}{}

It suffices to apply the Chacon-Orstein theorem by taking $g$ to
be constant function equal to 1.

\end{proof}

\begin{remark}{}{}

Under the same hypotheses, by taking $f$ to be the indicator function
of a singleton, we obtain that for every $x,y\in E,$\boxeq{
\[
\dfrac{1}{n}\sum^{n}_{j=1}\boldsymbol{1}_{\left(X_{j}=y\right)}\stackrel[n\to+\infty]{P-\text{a.s.}}{\longrightarrow}\pi\left(y\right).
\]
}The quotient 
\[
\dfrac{1}{n}\sum^{n}_{j=1}\boldsymbol{1}_{\left(X_{j}=y\right)}
\]
represents the \textbf{average time\index{average time}} spent by
a trajectory in the state $y$ between times $1$ and $n.$ This result
therefore provides a method for \textbf{estimating the invariant probability
}\mindex{estimation!invariant probability}.

\end{remark}

\subsection{Estimation of the Transition Matrix}

We now suppose that $E=\left\{ x_{1},\cdots,x_{L}\right\} $ is finite
and that $X$ is an irreducible homogeneous Markov chain with transition
matrix $M.$ Then the chain is positive recurrent and there exists
a unique invariant probability $\pi.$ 

Define for $i,j\in\left\llbracket 1,L\right\rrbracket $ and $n\in\mathbb{N}^{\ast},$
the random variables $N^{n}_{i}$ and $N^{n}_{i,j}$ by
\[
N^{n}_{i}=\sum^{n-1}_{l=0}\boldsymbol{1}_{\left(X_{l}=x_{i}\right)}\,\,\,\,\text{and}\,\,\,\,N^{n}_{i,j}=\sum^{n-1}_{l=0}\boldsymbol{1}_{\left(X_{l}=x_{i}\right)}\boldsymbol{1}_{\left(X_{l+1}=x_{j}\right)}.
\]
These quantities represent, respectively, the number of visits to
$x_{i},$ and the number of transitions from $x_{i}$ to $x_{j}$
up to the time $n.$ 

We have
\[
\sum^{L}_{j=1}N^{n}_{i,j}=\sum^{n-1}_{l=0}\sum^{L}_{j=1}\boldsymbol{1}_{\left(X_{l}=x_{i}\right)}\boldsymbol{1}_{\left(X_{l+1}=x_{j}\right)}=\sum^{n-1}_{l=0}\boldsymbol{1}_{\left(X_{l}=x_{i}\right)},
\]
hence 
\[
N^{n}_{i}=\sum^{L}_{j=1}N^{n}_{i,j}.
\]
 Define 
\[
\widehat{M^{n}_{i,j}}=\dfrac{N^{n}_{i,j}}{N^{n}_{i}}.
\]
We study for $x\in E,$ the $P_{x}-$almost sure convergence of the
sequence with general term $\widehat{M^{n}_{i,j}}.$

\begin{proposition}{}{}

With the previous notation and hypotheses, for $i,j\in\left\llbracket 1,L\right\rrbracket $
and for every $x\in E,$
\[
\widehat{M^{n}_{i,j}}\stackrel[n\to+\infty]{P_{x}-\text{a.s.}}{\longrightarrow}M\left(x_{i},x_{j}\right).
\]

\end{proposition}

\begin{proof}{}{}

By the previous remark, we already have, for $i\in\left\llbracket 1,L\right\rrbracket ,$
\begin{equation}
\dfrac{N^{n}_{i}}{n}\stackrel[n\to+\infty]{P_{x}-\text{a.s.}}{\longrightarrow}\pi\left(x_{i}\right).\label{eq:eq17_84}
\end{equation}

Let $\pi\otimes M$ be the probability on $E\times E$ defined, for
every subset $A$ of $E\times E$ by
\begin{align*}
\left(\pi\otimes M\right)\left(A\right) & =\intop\left[\sum^{L}_{i=1}\boldsymbol{1}_{A}\left(x,x_{i}\right)M\left(x,x_{i}\right)\right]\text{d}\pi\left(x\right)\\
 & =\sum_{x\in E}\left[\sum^{L}_{i=1}\boldsymbol{1}_{A}\left(x,x_{i}\right)M\left(x,x_{i}\right)\right]\pi\left(x\right).
\end{align*}
Of course,
\[
\intop f\text{d}\left(\pi\otimes M\right)=\sum_{x\in E}\left[\sum^{L}_{i=1}f\left(x,x_{i}\right)M\left(x,x_{i}\right)\right]\pi\left(x\right).
\]
By adapting the proof of Corollary $\ref{co:strong_law_large_numbers_markov},$
we can prove that for every function $f$ on $E\times E$---the integrability
being automatic here---
\begin{equation}
\dfrac{1}{n}\sum^{n-1}_{l=0}f\left(X_{l},X_{l+1}\right)\stackrel[n\to+\infty]{P_{x}-\text{a.s.}}{\longrightarrow}\intop f\text{d}\pi\otimes M.\label{eq:eq17_85}
\end{equation}
Indeed, for any $y\in E,$ define
\[
Z_{p}=\begin{cases}
\sum^{T^{p+1}_{y}-1}_{n=T^{p}_{y}}f\left(X_{n},X_{n+1}\right), & \text{on }\left(T^{p}_{y}<+\infty\right),\\
0, & \text{on }\left(T^{p}_{y}=+\infty\right).
\end{cases}
\]
By a computation similar to that in the proof of Theorem $\ref{th:chacon_orstein_th},$
we obtain for every Borel set $B\subset\mathbb{R},$
\[
\mathbb{E}^{\mathscr{A}_{T^{p}_{y}}}_{x}\left[\boldsymbol{1}_{\left(Z_{p}\in B\right)}\right]=P_{y}\left(Z_{0}\in B\right).
\]
This also proves that the $Z_{p},\,p\in\mathbb{N}^{\ast},$ have under
$P_{x}$ the same law as $Z_{0}$ under $P_{y}$ and that the $Z_{p}$
are $P_{x}-$independent.

It remains to compute $\mathbb{E}_{x}\left(Z_{1}\right).$ Using an
argument similar to that in the proof of Theorem $\ref{th:chacon_orstein_th},$
the strong Markov property yields
\[
\mathbb{E}_{x}\left(Z_{1}\right)=\mathbb{E}_{y}\left[\sum^{T^{1}_{y}-1}_{n=0}f\left(X_{n},X_{n+1}\right)\right].
\]
Noting that $\left(n<T^{1}_{y}\right)\in\mathscr{A}_{n},$ it follows
by the simple Markov proverty that
\[
\mathbb{E}_{x}\left(Z_{1}\right)=\sum^{+\infty}_{n=0}\mathbb{E}_{y}\left[\boldsymbol{1}_{\left(n<T^{1}_{y}\right)}\mathbb{E}^{\mathscr{A}_{n}}_{y}\left[f\left(X_{n},X_{n+1}\right)\right]\right]=\sum^{+\infty}_{n=0}\mathbb{E}_{y}\left[\boldsymbol{1}_{\left(n<T^{1}_{y}\right)}\mathbb{E}_{X_{n}}\left[f\left(X_{n},X_{n+1}\right)\right]\right].
\]
This can also be written
\[
\mathbb{E}_{x}\left(Z_{1}\right)=\sum^{+\infty}_{n=0}\mathbb{E}_{y}\left[\boldsymbol{1}_{\left(n<T^{1}_{y}\right)}\left(\sum^{L}_{i=1}f\left(X_{n},x_{i}\right)M\left(X_{n},x_{i}\right)\right)\right],
\]
or equivalently
\[
\mathbb{E}_{x}\left(Z_{1}\right)=\mathbb{E}_{y}\left[\sum^{T^{1}_{y}-1}_{n=0}\left(\sum^{L}_{i=1}f\left(X_{n},x_{i}\right)M\left(X_{n},x_{i}\right)\right)\right].
\]
Using formulas $\refpar{eq:eq17_82}$ and the identity $\mu=\mathbb{E}_{y}\left(T^{1}_{y}\right)\pi$---see
Corollary $\ref{co:strong_law_large_numbers_markov}$---we obtain
\[
\mathbb{E}_{x}\left(Z_{1}\right)=\sum_{x\in E}\mu\left(x\right)\left[\sum^{L}_{i=1}f\left(x,x_{i}\right)M\left(x,x_{i}\right)\right]=\mathbb{E}_{y}\left(T^{1}_{y}\right)\intop f\text{d}\left(\pi\otimes M\right).
\]
The convergence stated in $\refpar{eq:eq17_85}$ then follows as the
final part of the proof of Corollary $\ref{co:strong_law_large_numbers_markov}$
and the remark that follows it.

By taking for $f$ the function defined by
\[
f\left(x,y\right)=\boldsymbol{1}_{\left\{ x_{i}\right\} }\left(x\right)\boldsymbol{1}_{\left\{ x_{j}\right\} }\left(y\right),
\]
we obtain
\[
N^{n}_{i,j}\stackrel[n\to+\infty]{P_{x}-\text{a.s.}}{\longrightarrow}\intop f\text{d}\left(\pi\otimes M\right)=\pi\left(x_{i}\right)M\left(x_{i},x_{j}\right).
\]
Taking into account $\refpar{eq:eq17_84},$ the proposition follows.

\end{proof}

\section*{Exercises}

\addcontentsline{toc}{section}{Exercises}

Unless explicitly mentioned, the random variables are defined on a
probabilized space \preds, and the processes are defined on the adequate
underlying filtered probabilized space $\left(\Omega,\mathscr{A},\left(\mathscr{A}_{n}\right)_{n\in\mathbb{N}},P\right).$

\begin{exercise}{Stock Management}{exercise17.1}

A smartphone salesperson has noticed that the number $A_{t}$ of buyers
of a given smartphone during week $t$ is independent of the number
of buyers in previous weeks. Its law is given by
\begin{gather*}
P\left(A_{t}=0\right)=0.4,\,\,\,\,P\left(A_{t}=1\right)=0.4,\,\,\,\,P\left(A_{t}=2\right)=0.15,\\
P\left(A_{t}=3\right)=0.05,\,\,\,\,P\left(A_{t}>3\right)=0.
\end{gather*}
The salesperson places orders at the end of each week, but only if
there are no smartphone left in stock at that time. In this case,
they order two smartphones, which are received on the first day of
the following week. 

Denote by $X_{t}$ the number of smartphones in stock at the end of
week $t.$ 

Prove that $X=\left(X_{t}\right)_{t\in\mathbb{\mathbb{N}}}$ is a
homogeneous Markov chain---with respect to the natural filtration
$\left(\mathscr{A}_{t}\right)_{t\in\mathbb{N}}$---taking values
in $E=\left\{ 0,1,2\right\} ,$ with transition matrix $M$ given
by
\[
M=\left(\begin{array}{ccc}
0.2 & 0.4 & 0.4\\
0.6 & 0.4 & 0\\
0.2 & 0.4 & 0.4
\end{array}\right).
\]

\end{exercise}

\begin{exercise}{Reliability. End of life}{exercise17.2}

Time is measured in a discrete manner---for instance in seconds---and
is therefore indexed by $\mathbb{N}.$ 

A machine runs permanently and contains a critical components that
fails easily. As soon as this component fails, it is instantaneously
replaced by an identical component. 

Let $X_{n}$ denote the random time between time $n$ and the time
of next failure occurring after $n.$ The sequence satisfies the relation
\[
X_{n+1}=\begin{cases}
X_{n}-1, & \text{if }X_{n}\geqslant1,\\
Z_{n}-1, & \text{if }X_{n}=0,
\end{cases}
\]
where $Z_{n}$ is the lifetime of the component installed at time
$n.$ We assume that the variables $Z_{n},\,n\in\mathbb{N}^{\ast},$
are independent random variables following the same law $\mu$ on
$\mathbb{N}^{\ast}$ generated by the germ $\left(p_{k}\right)_{k\in\mathbb{N}^{\ast}}.$ 

Prove that the process $X=\left(X_{n}\right)_{n\in\mathbb{N}^{\ast}}$
is a homogeneous Markov chain---with respect to its natural filtration
$\left(\mathscr{A}_{n}\right)_{n\in\mathbb{N}^{\ast}}$---taking
values in $\mathbb{N},$ and determine its transition matrix $M.$

\end{exercise}

\begin{exercise}{Simple and Strong Markov Properties}{exercise17.3}

Let $X=\left(X_{n}\right)_{n\in\mathbb{N}}$ be a process taking values
in $E=\left\{ 1,2,3\right\} $ such that, for each $i\in E,$ $X$
is a homogeneous Markov chain on the underlying filtered probabilized
space 
\[
\left(\Omega,\mathscr{A},\left(\mathscr{A}_{n}\right)_{n\in\mathbb{N}},P_{i}\right),
\]
with transition matrix
\[
M=\left(\begin{array}{ccc}
0 & \dfrac{1}{3} & \dfrac{2}{3}\\
\dfrac{1}{2} & \dfrac{1}{4} & \dfrac{1}{4}\\
\dfrac{1}{2} & \dfrac{1}{2} & 0
\end{array}\right).
\]

1. Let $f$ be the functional on $E^{\mathbb{N}}$ defined, for every
$x\in E^{\mathbb{N}},$ by
\[
f\left(x\right)=\boldsymbol{1}_{\left(\sum^{3}_{j=1}\boldsymbol{1}_{\left\{ 1\right\} }\left(x_{j}\right)=1\right)}.
\]
That is, $f\left(x\right)$ is equal to 1 if, among the three components
of $x$ with indices 1 to 3, exactly one is equal to 1. Otherwise,
$f\left(x\right)$ is equal to 0. 

Compute, for every $i\in E,$ 
\[
\mathbb{E}^{\mathscr{A}_{n}}_{i}\left(f\left(\theta_{n}\left(X\right)\right)\right).
\]

2. Let $T=\inf\left\{ n\in\mathbb{N}^{\ast}:\,X_{n}=1\right\} ,$
with $\inf\emptyset=+\infty.$ Justify the fact that, for every $i\in E,$
$P_{i}\left(T<+\infty\right)=1.$ Compute, for $i\in E,$ the conditional
expectation $\mathbb{E}^{\mathscr{A}_{T}}_{i}\left(f\left(\theta_{T}\left(X\right)\right)\right).$

\end{exercise}

\begin{exercise}{Discrete Birth-and-Death Process. Random Walk on $\mathbb{N}$ With Elastic Barriers. Dirichlet Problem}{exercise17.4}

We seek to model the evolution of a population---of individuals,
physical particles, or any similar entities---whose size may increase
or decrease by one unit at each time $n,$ with probabilities depending
on the current population size $X_{n}.$

Let $X=\left(X_{n}\right)_{n\in\mathbb{N}}$ be a process taking values
in $\mathbb{N},$ such that for every $x\in\mathbb{N},$ $X$ is a
homogeneous Markov chain on the underlying filtered probabilized space
$\left(\Omega,\mathscr{A},\left(\mathscr{A}_{n}\right)_{n\in\mathbb{N}},P_{x}\right),$
with transition matrix $M$ given by, for every $x\in\mathbb{N},$
\[
M\left(x,x+1\right)=p_{x}\,\,\,\,M\left(x,x-1\right)=q_{x},\,\,\,\,M\left(x,x\right)=r_{x}
\]
where 
\[
p_{x},q_{x},r_{x}\in\left[0,1\right]\,\,\,\,\text{and}\,\,\,\,p_{x}+q_{x}+r_{x}=1.
\]

Let $a,b\in\mathbb{N}$ such that $0\leqslant a<b.$ Denote by $T_{x}$
the hitting time of $x,$ that is $T_{x}=\inf\left\{ n\in\mathbb{N}^{\ast}:\,X_{n}=x\right\} ,$
with $\inf\emptyset=+\infty.$

1. Prove that, for every $x,$
\begin{equation}
P_{x}\left(T_{a}<T_{b}\right)=M\left(x,a\right)+\sum_{z\neq a,b}P_{z}\left(T_{a}<T_{b}\right)M\left(x,z\right)\label{eq:P_xT_alessT_b}
\end{equation}

2. Suppose that $a=0$ and that $a$ and $b$ are \textbf{barriers}\index{barrier},
that is $q_{0}=0$ and $p_{b}=0.$ They are \textbf{reflecting}\index{reflecting},
in the sense that $r_{a}>0$ and $r_{b}>0.$ Thus, if at some time
the whole popupation has disapeared, one element may be reintroduced
at the next instant; conversely, if the population reaches size $b,$
one element may be removed at the next instant, but no new one may
be added. 

Assume moreover that $p_{x}>0$ for every $x\in\left]0,b\right[.$

For every $x\in\left]0,b\right[,$ define $f\left(x\right)=P_{x}\left(T_{0}<T_{b}\right)$
which is the probability that, starting from population size $x,$
the population becomes extinct before reaching size $b.$ 

Prove that $f$ is solution of a second-order recurrence equation
with boundary conditions---\textbf{\index{Dirichlet problem}}that
is, a \textbf{Dirichlet problem}. 

Then compute explicitely $P_{x}\left(T_{0}<T_{b}\right)$ for $x\in\left]0,b\right[$
as a function of the sequence with general term $a_{x}$ defined by
\[
a_{0}=1,\,\,\,\,\text{and for }x\in\left]0,b\right[,\,\,\,\,a_{x}=\dfrac{q_{1}q_{2}\cdots q_{x}}{p_{1}p_{2}\cdots p_{x}}.
\]

\end{exercise}

\begin{exercise}{Law of the First Hitting Time of a State}{exercise17.5}

Let $E=\left\{ 1,2,3\right\} .$ Let $X=\left(X_{n}\right)_{n\in\mathbb{N}}$
be a process such that, for every $x\in E,$ $X$ is a homogeneous
Markov chain on the underlying filtered probabilized space 
\[
\left(\Omega,\mathscr{A},\left(\mathscr{A}_{n}\right)_{n\in\mathbb{N}},P_{x}\right),
\]
with transition matrix $M$ given by
\[
M=\left(\begin{array}{ccc}
1 & 0 & 0\\
\dfrac{1}{2} & \dfrac{1}{6} & \dfrac{1}{3}\\
\dfrac{1}{3} & \dfrac{3}{5} & \dfrac{1}{15}
\end{array}\right).
\]
1. Determine the graph associated with this chain and specify its
communication classes.

2. Under the probability $P_{x},$ we study the law of the first hitting
time at 3, namely $T_{3}=\inf\left\{ n\in\mathbb{N}^{\ast}:\,X_{n}=3\right\} ,$
with $\inf\emptyset=+\infty.$ For this purpose, define, for $x\in E,$
\[
f_{k}\left(x\right)=P_{x}\left(T_{3}=k\right)\equiv F_{k}\left(x,3\right).
\]
Prove that the sequence $\left(f_{k}\right)_{k\in\mathbb{N}^{\ast}}$
of vectors in $\mathbb{R}^{3}$ satisfies the recurrence relation
\begin{equation}
f_{k}=Qf_{k-1},\label{eq:recurrence_rel_f_k}
\end{equation}
where $Q$ is a $3\times3$ matrix to be determined. Then compute
$f_{k},$ for every $k\in\mathbb{N}^{\ast}.$

\textit{Hint: the Cayley-Hamilton theorem may be used to advantage.}

3. Compute, for every $x\in E,$ the probability $P_{x}\left(T_{3}=+\infty\right).$

4. For every $y\in E,$ let 
\[
N_{y}=\sum_{j\in\mathbb{N}}\boldsymbol{1}_{\left(X_{j}=y\right)}
\]
be the number of visits to $y.$ 

Determine $\mathbb{E}_{3}\left(N_{3}\right)$ and $\mathbb{E}_{2}\left(N_{3}\right).$

\end{exercise}

\begin{exercise}{Invariant Probability. Mean Return Time to a State}{exercise17.6}

Let $E=\left\{ 1,2,3,4,5\right\} .$ Let $X=\left(X_{n}\right)_{n\in\mathbb{N}}$
be a process that, for every $x\in E,$ is a homogeneous Markov chain
on the underlying filtered probabilized space $\left(\Omega,\mathscr{A},\left(\mathscr{A}_{n}\right)_{n\in\mathbb{N}},P_{x}\right),$
with transition matrix $M$ given by
\[
M=\left(\begin{array}{ccccc}
0 & \dfrac{2}{3} & \dfrac{1}{3} & 0 & 0\\
1 & 0 & 0 & 0 & 0\\
1 & 0 & 0 & 0 & 0\\
\dfrac{1}{2} & 0 & 0 & 0 & \dfrac{1}{2}\\
0 & 0 & \dfrac{1}{2} & \dfrac{1}{2} & 0
\end{array}\right).
\]

1. Determine the graph associated with this chain and specify the
communication classes. Also specify the period and the nature of these
classes.

2. Justify, without computation, the existence of a unique invariant
probability $\nu.$ Compute it, and deduce from it, for every $x\in\left\{ 1,2,3\right\} ,$
the quantity $\mathbb{E}_{x}\left(T_{x}\right),$ where $T_{x}$ is
the first hitting time in $x,$ hence $T_{x}=\inf\left\{ n\in\mathbb{N}^{\ast}:X_{n}=x\right\} ,$
with $\inf\emptyset=+\infty.$

3. Determine the potential matrix $R.$

4. Compute the probabilities $P_{4}\left(T_{5}<+\infty\right)$ and
$P_{5}\left(T_{5}<+\infty\right).$

5. For every $y\in E,$ let $N_{y}=\sum_{j\in\mathbb{N}}\boldsymbol{1}_{\left(X_{j}=y\right)}$
be the number of visits to $y.$ Compute, for every $m\in\mathbb{N},$
the probabilities $P_{4}\left(N_{5}=m\right)$ and $P_{5}\left(N_{5}=m\right).$

6. Prove that, for every $x,y\in E,$
\begin{equation}
\mathbb{E}_{x}\left[T_{y}\boldsymbol{1}_{\left(T_{y}<+\infty\right)}\right]=M\left(x,y\right)+\sum_{z\neq y}M\left(x,z\right)\left[\mathbb{E}_{z}\left(T_{y}\boldsymbol{1}_{\left(T_{y}<+\infty\right)}\right)+P_{z}\left(T_{y}<+\infty\right)\right].\label{eq:E_x_T_y_1_t_yless_inf}
\end{equation}
Deduce numerically, for every $x\in E,$ $\mathbb{E}_{x}\left(T_{1}\right).$

\end{exercise}

\begin{exercise}{First Hitting Time of a Set and Manufacturing Process}{exercise17.7}

Let $X=\left(X_{n}\right)_{n\in\mathbb{N}}$ be a process taking values
in $E,$ which, for every $x\in E,$ is a homogeneous Markov chain
on the underlying filtered probabilized space $\left(\Omega,\mathscr{A},\left(\mathscr{A}_{n}\right)_{n\in\mathbb{N}},P_{x}\right),$
with transition matrix $M.$

We adopt the following notation. For every subset $C$ of $E,$ we
define
\[
T_{C}=\inf\left\{ n\in\mathbb{N}^{\ast}:\,X_{n}\in C\right\} ,
\]
the first hitting time of $C.$ The mapping $\tau_{C}$ is defined
on $E^{\mathbb{N}},$ for every $u\in E^{\mathbb{N}},$ by
\[
\tau_{C}\left(u\right)=\inf\left\{ n\in\mathbb{N}^{\ast}:\,u_{n}\in C\right\} .
\]
We set $\inf\emptyset=+\infty.$ 

For every $x,y\in E,$ denote 
\[
F\left(x,y\right)=P_{x}\left(T_{y}<+\infty\right).
\]
We also introduce the shift operators $\theta_{p}$ on $E^{\mathbb{N}},$
defined, for every $u\in E^{\mathbb{N}}$ and every $n\in\mathbb{N},$
by $\left[\theta_{p}\left(u\right)\right]_{n}=u_{n+p}.$

1. For every $x,y\in E,$ define
\[
G\left(x,y\right)=\mathbb{E}_{x}\left[T_{y}\boldsymbol{1}_{\left(T_{y}<+\infty\right)}\right].
\]

Prove that
\begin{equation}
G\left(x,y\right)=P_{x}\left(T_{y}<+\infty\right)+\sum_{z\neq y}M\left(x,z\right)G\left(z,y\right).\label{eq:G_x_y}
\end{equation}

2. Let $C$ and $D$ be two disjoint subsets of $E.$ Define
\[
\phi_{C,D}\left(x\right)=P_{x}\left[\left(T_{C}<+\infty\right)\cap\left(\tau_{D}\left(\theta_{T_{C}}\left(X\right)\right)=+\infty\right)\right],
\]
that is, the $P_{x}-$probability of the set of trajectories that
reach $C$ in finite time and do not pass through $D$ before having
passed through $C.$ 

By conditioning with respect to $\mathscr{A}_{1},$ prove that $\phi_{\left(C,1\right)}$
is a solution of the system of equations
\begin{equation}
\forall x\in E,\,\,\,\,\phi_{C,D}\left(x\right)=\sum_{y\in C}M\left(x,y\right)P_{y}\left(T_{D}=+\infty\right)+\sum_{y\notin C}M\left(x,y\right)\phi_{C,D}\left(y\right).\label{eq:system_phi_CD_x}
\end{equation}
\textbf{Application}

The manufacturing process of a component requires three successive
steps, denoted 1, 2, 3. After the step $i,$ the component is tested.
If it is satisfactory, which happens with the probability $r_{i},$
then it moves to step $i+1.$ Otherwise, either it is considered completely
defective, with probability $p_{i}$ and is discarded---this is state
5---or it is considered only slightly defective, which happens with
probability $q_{i},$ and is sent back to step $i.$ 

State 4 is the one of a good component that has successfully passed
all the manufacturing steps. We assume that $p_{i}+q_{i}+r_{i}=1,$
for every $i=1,\,2,\,3.$

We model this manufacturing process by a homogeneous Markov chain
$X$ with values in $E=\left\{ 1,\,2,\,3,\,4,\,5\right\} $ where
$X_{n}$ represents the state of the component at step $n,$ and with
transition matrix $M$ given by
\[
\begin{array}{rcl}
 &  & \begin{array}{ccccc}
\thinspace\,\,\boldsymbol{1} & \thinspace\boldsymbol{2} & \thinspace\,\boldsymbol{3} & \thinspace\boldsymbol{4} & \thinspace\boldsymbol{5}\end{array}\\
M= & \begin{array}{c}
\boldsymbol{1}\\
\boldsymbol{2}\\
\boldsymbol{3}\\
\boldsymbol{4}\\
\boldsymbol{5}
\end{array} & \left(\begin{array}{ccccc}
q_{1} & r_{1} & 0 & 0 & p_{1}\\
0 & q_{2} & r_{2} & 0 & p_{2}\\
0 & 0 & q_{3} & r_{3} & p_{3}\\
0 & 0 & 0 & 1 & 0\\
0 & 0 & 0 & 0 & 1
\end{array}\right).
\end{array}
\]
We keep the notation introduced in the previous questions.

3. Determine the graph associated with this chain and specify the
communication classes. Also specify the nature of these classes.

4. Justify the equality $F\left(5,4\right)=0.$ 

If $K$ is the vector in $\mathbb{R}^{3}$ whose components are $F\left(i,4\right)$
for $i=1,2,3,$ prove that $K$ is solution of an equation of the
form $K=b+TK,$ where $b\in\mathbb{R}^{3}$ and $T$ is a matrix $3\times3$
to be specified. 

Deduce the value of $K.$ Determine the probability $P_{1}\left(T_{4}<+\infty\right).$

5. Write the potential matrix $R$ of $X.$

6. For every $x\in E,$ let $H\left(x\right)=\mathbb{E}_{x}\left[T_{4}\boldsymbol{1}_{\left(T_{4}<+\infty\right)}\right].$ 

Prove that $H\left(5\right)=0.$ 

Compute $H\left(x\right)$ for $x=1,2,3$ and interpret $H\left(1\right).$

7. Set $C=\left\{ 2\right\} $ and $D=\left\{ 4\right\} .$ Compute,
for every $x\in E,$ $\phi\left(x\right)\equiv\phi_{C,D}\left(x\right).$
Interpret $\phi\left(1\right).$

\end{exercise}

\begin{exercise}{A Coin-Tossing Game and Markov Chains}{exercise17.8}

We consider a sequence of tosses of a coin, not necessarily fair,
and we are interested in the outcomes obtained over two consecutive
tosses. In particular, we study the random variable giving the number
of tosses needed, for instance to obtain two consecutive tails.

Let $\left(X_{n}\right)_{n\in\mathbb{N}}$ be a Bernoulli process,
that is, a sequence of random variables defined on a probabilized
space $\left(\Omega,\mathscr{A},P\right)$ taking values in $\left\{ 0,1\right\} ,$
independent, following the same Bernoulli law with parameter $p,$
such that $P\left(X_{n}=1\right)=p,$ $P\left(X_{n}=0\right)=q,$
where $p+q=1.$ 

For every $n\in\mathbb{N},$ define the random variable $Y_{n}=\left(X_{n},X_{n+1}\right)$
and the \salg $\mathscr{A}_{n}=\sigma\left(X_{j}:\,0\leqslant j\leqslant n+1\right).$ 

Let $E$ be the set of four points
\[
\alpha_{1}=\left(1,0\right),\,\,\,\,\alpha_{2}=\left(1,1\right),\,\,\,\,\alpha_{3}=\left(0,0\right)\,\,\,\,\text{and}\,\,\,\,\alpha_{4}=\left(0,1\right).
\]

1. Compute, for every real valued function $f$ on $E$ and every
$n\in\mathbb{N},$ the conditional expectation $\mathbb{E}^{\mathscr{A}_{n}}\left[f\left(Y_{n+1}\right)\right]$
and deduce that the process $Y=\left(Y_{n}\right)_{n\in\mathbb{N}}$
is a homogeneous Markov chain taking values in $E,$ with transition
matrix $M$ to be determined. Prove that, for every $n\geqslant2,$
$M^{n}=M^{2}.$

2. Determine the graph associated with this chain and specify the
communication classes. Determine the nature and period of the points
of $E.$

For every $x\in E,$ denote by $P_{x}$ the conditional probability
$P_{x}=P\left(\cdot\mid Y_{0}=x\right).$

3. Under the probability $P_{x},$ we study the law of the first hitting
time $T$ of the chain $Y$ in $\alpha_{2}.$ Let $T=\inf\left\{ n\in\mathbb{N}^{\ast}:\,Y_{n}=\alpha_{2}\right\} ,$
with the convention $\inf\emptyset=+\infty.$ For this purpose, denote,
for $x=1,2,3,4,$ 
\[
f_{k}\left(i\right)=P_{\alpha_{i}}\left(T=k\right)\equiv F_{k}\left(\alpha_{i},\alpha_{2}\right).
\]

a. Prove that the sequence $\left(f_{k}\right)_{k\in\mathbb{N}^{\ast}}$
of vectors in $\mathbb{R}^{4}$ satisfies a first order recurrence
relation. 

b. Deduce that the sequence of probabilities $\left(f_{k}\left(1\right)\right)_{k\in\mathbb{N}^{\ast}}$
satisfies the second-order recurrence equation
\begin{equation}
x_{k}=qx_{k-1}+pqx_{k-2}\label{eq:second_ord_rec_eq_x_k}
\end{equation}
 with initial conditions to be specified.

c. Then determine two particular solutions of the form $\left(\lambda^{k}\right)_{k\geqslant2}$
corresponding to the two values of $\lambda,$ denoted $\lambda_{1}$
and $\lambda_{2},$ expressed as functions of $p$ and $q.$ Deduce
the value of the probability $f_{k}\left(1\right)$ as a function
of $\lambda_{1}$ and $\lambda_{2}.$

4. Compute the mean time $\mathbb{E}_{\alpha_{1}}\left(T\right),$
expressing it only as a function of $p.$

5. Justify the existence of a unique invariant probability $\nu,$
and compute it. Deduce the value of the mean time $\mathbb{E}_{\alpha_{2}}\left(T\right)$
and compare it with $\mathbb{E}_{\alpha_{1}}\left(T\right).$

\end{exercise}

\begin{exercise}{Random Walk on the Integer Interval $\left\llbracket 0,N\right\rrbracket$ with Reflecting Barriers. Invariant Probability }{exercise17.9}

Let $X=\left(X_{n}\right)_{n\in\mathbb{N}}$ be a process taking values
in $E=\left\llbracket 0,N\right\rrbracket $ which, for every $x\in E,$
is a homogeneous Markov chain on the underlying filtered probabilized
space $\left(\Omega,\mathscr{A},\left(\mathscr{A}_{n}\right)_{n\in\mathbb{N}},P_{x}\right),$
with transition matrix $M$ given by
\[
\left\{ \begin{array}{lll}
M\left(x,x+1\right)=p\qquad & M\left(x,x-1\right)=q,\qquad & \text{if }1\leqslant x\leqslant N-1,\\
M\left(0,1\right)=1, & M\left(N,N-1\right)=1,
\end{array}\right.
\]
where $p,q\in\left]0,1\right[$ and $p+q=1.$ 

For every $x\in E,$ denote by $T_{x}$ the hitting time in $x,$
that is $T_{x}=\inf\left\{ n\in\mathbb{N}^{\ast}:X_{n}=x\right\} ,$
with $\inf\emptyset=+\infty.$

1. Specify the communication class or classes.

2. Justify the existence of a unique invariant probability $\nu$
and compute it as a function of $p,q,N.$ Deduce the value of the
mean time $E_{0}\left(T_{0}\right).$

\end{exercise}

\begin{exercise}{Random Walk on $\mathbb{N}$ with a Barrier of Arbitrary Type. Invariant Measure and Probability. Limiting Probability}{exercise17.10}

Let $X=\left(X_{n}\right)_{n\in\mathbb{N}}$ be a process taking values
in $\mathbb{N}$ which, for every $x\in\mathbb{N},$ is a homogeneous
Markov chain on the underlying filtered probabilized space $\left(\Omega,\mathscr{A},\left(\mathscr{A}_{n}\right)_{n\in\mathbb{N}},P_{x}\right),$
with transition matrix $M$ given by
\[
\left\{ \begin{array}{lll}
M\left(x,x+1\right)=p\qquad & M\left(x,x-1\right)=q,\qquad & \text{if }x\in\mathbb{N}^{\ast},\\
M\left(0,0\right)=\alpha, & M\left(0,1\right)=1-\alpha,
\end{array}\right.
\]
where $p,q\in\left]0,1\right[,$ $p+q=1,$ and $\alpha\in\left[0,1\right].$ 
\begin{itemize}
\item If $\alpha=0,$ the point $0$ is called \textbf{reflecting barrier}\index{reflecting barrier}\mindex{barrier!reflecting}.
\item If $\alpha\in\left]0,1\right[,$ the point 0 is called an \textbf{elastic
barrier}\index{elastic barrier}\mindex{barrier!elastic}.
\item If $\alpha=1,$ 0 is called \textbf{absorbing barrier}\index{absorbing barrier}\mindex{barrier!absorbing}.
\end{itemize}
For every $x\in\mathbb{N},$ denote $T_{x}$ the first hitting time
of $x,$ that is $T_{x}=\inf\left\{ n\in\mathbb{N}^{\ast}:\,X_{n}=x\right\} ,$
with $\inf\emptyset=+\infty.$

1. Specify the communication class or classes. Study the periodicity
of the states.

\textbf{We first consider the case where $\alpha\in\left]0,1\right[.$}

2. a. Prove, by computation, the existence of an invariant measure
$\nu.$ 

b. Study, depending on the respective values of $p$ and $q,$ the
existence and uniqueness of an invariant probability, and compute
it in the case where it exists and is unique.

c. Deduce that, in the case where $\alpha\in\left[0,1\right[$ and
$p<q,$ the nature of the points of $\mathbb{N},$ and give, for every
$x\in\mathbb{N},$ the value $\mathbb{E}_{x}\left(T_{x}\right)$ of
the mean return time to $x.$

3. In the case where $p\geqslant q,$ study the nature of the points
of $\mathbb{N,}$ and if $p>q,$ compute, for every $x\in\mathbb{N}^{\ast},$
the probability $P_{x}\left(T_{0}=+\infty\right).$

4. In the case where $p<q,$ justify the $P_{x}-$almost sure convergence
of the sequence with general term $\dfrac{1}{n}\sum^{n}_{j=1}\text{e}^{-aX_{j}},$
for any $a>0.$ 

\textbf{We now study the case where $\alpha=1.$}

5. Compute $P_{0}\left(T_{0}<+\infty\right)$ and $\mathbb{E}_{0}\left(T_{0}\right).$
Deduce the nature of the point 0.

6. Compute, for every $x\in\mathbb{N}^{\ast},$ the probability $P_{x}\left(T_{0}<+\infty\right).$
Determine the nature of the points of $\mathbb{N}^{\ast}.$

7. Study the convergence of the sequence with general term $M^{n}\left(x,y\right),$
and specify, when appropriate, its limit, first when $x\in\mathbb{N}$
and $y\in\mathbb{N}^{\ast},$ and then when $x\in\mathbb{N}^{\ast}$
and $y=0.$

\end{exercise}

\begin{exercise}{Galton-Watson Process and Martingales. Evolution of the Size of a Population}{exercise17.11}

We study the evolution of the sizes of successive generations in a
population where each individual gives birth to a random number of
descendants following the same probability law $\mu.$ In particular,
we aim to evaluate the probability of extinction of the population.
The model is defined as follows.

Consider a probability law $\mu$ on $\mathbb{N}$ such that $0<\mu\left(\left\{ 0\right\} \right)<1$
and $0<m<+\infty,$ where $m$ denotes the expectation of $\mu,$
defined by 
\[
\mathbb{E}\left(\mu\right)=\sum^{+\infty}_{n=0}n\mu\left(\left\{ n\right\} \right).
\]
Finally, denote by $g$ the generative function of $\mu,$ defined
on $\left[0,1\right],$ by
\[
g\left(u\right)=\sum^{+\infty}_{n=0}u^{n}\mu\left(\left\{ n\right\} \right).
\]

Consider a family, indexed on $\mathbb{N}\times\mathbb{N}^{\ast}$
of random variables $Y_{n,i},$ defined on the probabilized space
\preds taking values in $\mathbb{N},$ independent and following
the same law $\mu.$

The variables $Y_{n,i}$ represents the number of direct descendants
of the $i-$th individual of the $n-$th generation. Let $a\geqslant1$
be an integer representing the size of the initial population. 

The process $X,$ called a \textbf{branching process}\index{branching process}\mindex{process!branching}---individuals
may be identified with the \textbf{vertices} of a tree, in the graph-theoretic
sense, or with \textbf{nodes} of a genealogical tree---or \textbf{Galton-Watson
process}\index{Galton-Watson process}\mindex{process!Galton-Watson},
is defined by
\[
X_{0}=a
\]
and, for every $n\in\mathbb{N},$
\[
X_{n+1}=\boldsymbol{1}_{\left(X_{n}\geqslant1\right)}\sum^{X_{n}}_{j=1}Y_{n,j},
\]
with the convention 
\[
\sum^{0}_{j=1}Y_{n,j}=0.
\]
The variable $X_{n}$ represents the number of individuals in the
$n-$th generation.

The natural filtration $\left(\mathscr{A}_{n}\right)_{n\in\mathbb{N}}$
of the process $X$ will be the only filtration considered in what
follows.

1. Prove that $X$ is a homogeneous Markov chain, and determine its
transition matrix $M.$

2. Show that $X$ is a \textbf{martingale}, a \textbf{submartingale}
or \textbf{supermartingale}, depending on the values of $m.$

3. For every $n\in\mathbb{N},$ define the random variable $Y_{n}=\dfrac{X_{n}}{m^{n}}.$
Prove that $Y$ is a \textbf{nonnegative martingale}.

4. Assume that $m>1.$ We admit that there exists a unique real number
$s\in\left]0,1\right[$ such that $g\left(s\right)=s.$ For every
$n\in\mathbb{N},$ define the random variable $Z_{n}=s^{X_{n}}.$
Prove that $Z$ is a \textbf{equi-integrable submartingale}.

5. Prove that the sequence $\left(X_{n}\right)_{n\in\mathbb{N}}$
converges $P-$almost surely to a random variable $X_{\infty}.$ Study
separately the cases $0<m\leqslant1$ and $m>1.$ Identify the limit
$X_{\infty}$ in the case $0<m<1.$

6. Let $j\in\mathbb{N}^{\ast}.$ For every $k>N,$ compute the probability
\[
P\left[\bigcap^{k}_{n=N}\left(X_{n}=j\right)\right]
\]
as a function of $M\left(j,j\right)$ and $P\left(X_{N}=j\right).$
Deduce that 
\[
P\left[\liminf_{n\to+\infty}\left(X_{n}=j\right)\right]=0.
\]

7. Prove that, for every $j\in\mathbb{N}^{\ast},$ 
\[
P\left(X_{\infty}=j\right)=0.
\]
Consequently, $X_{\infty}\in\left\{ 0,+\infty\right\} $ $P-$almost
surely. Justify that every points of $\mathbb{N}^{\ast}$ is transient.
If $m>1,$ deduce from the fourth question that
\[
P\left(X_{\infty}=0\right)=s^{a}\,\,\,\,\text{and}\,\,\,\,P\left(X_{\infty}=+\infty\right)=1-s^{a}.
\]

8. Let $T$ denote the \textbf{extinction time\index{extinction time}}
of the process $X,$ that is, the first hitting time of 0, defined
by $T=\inf\left\{ n\in\mathbb{N}^{\ast}:\,X_{n}=0\right\} ,$ with
the convention $\inf\emptyset=+\infty.$ Verify that $P-$almost surely
\[
\left(X_{\infty}=+\infty\right)=\liminf_{n\to+\infty}\left(X_{n}\neq0\right)=\left(T=+\infty\right).
\]
Deduce the value of the probability $P\left(T<+\infty\right)$ for
the different values of $m.$ The stopping time $T$ is the\textbf{
extinction date\index{extinction date}} of the population.

\end{exercise}

\begin{exercise}{Contagious Diseases Diffusion Model of P\'olya---continued}{exercise17.12}

In this exercise we complete the study of the P\'olya model---see
Example $\ref{ex:polya_diff_mod}$ and $\ref{ex:polya_diff_mod_2}$---described
in terms of ball draws from a urn. 

More precisely, we aim to prove that the law of the random variable
$Y_{\infty},$ which is the $P-$almost sure limit of the sequence
$\left(Y_{n}\right)$ of proportions of blue balls contained in the
urn after the $n$-th draws---and after adding the drawn ball together
with $c$ balls of the same color\footnotemark---is the Beta law
$\beta\left(\dfrac{b}{c},\dfrac{r}{c}\right)$ of the first kind on
$\left[0,1\right].$ 

The method consists in computing the moments of all orders of $Y_{\infty}.$

We keep all the notations of Example $\ref{ex:polya_diff_mod_2}.$

1. Let $l\geqslant1$ be an arbitrary integer. Define the process
$Z=\left(Z_{n}\right)_{n\in\mathbb{N}^{\ast}}$ by
\[
Z_{n}=\prod^{l-1}_{j=0}Y_{n+j}=\dfrac{B_{n}\left(B_{n}+c\right)\cdots\left[B_{n}+\left(l-1\right)c\right]}{k_{n}k_{n+1}\cdots k_{n+\left(l-1\right)}}.
\]
Prove that $Z$ is a bounded martingale, and that the sequence $\left(Z_{n}\right)_{n\in\mathbb{N}^{\ast}}$
converges $P-$almost surely and in $\mathscr{L}^{1}.$

2. Deduce the values of $\mathbb{E}\left(Y^{l}_{\infty}\right).$
Express the result using the Gamma function $\Gamma.$

3. Let $U$ be a real-valued random variable following the Beta law
$\beta\left(a,b\right)$ of first kind on $\left[0,1\right].$ Compute
its moment of order $l.$

4. Using the results of questions 2 and 3, prove that the law of $Y_{\infty}$
is the Beta law $\beta\left(\dfrac{b}{c},\dfrac{r}{c}\right)$ of
the first kind on $\left[0,1\right].$

\end{exercise}

\footnotetext{Recall that the process $Y$ is a non-homogeneous Markov
chain and a martingale.}

\section*{Solutions of Exercises}

\addcontentsline{toc}{section}{Solutions of Exercises}

\begin{solution}{}{solexercise17.1}

Let $f$ the function defined on $E^{2}$ by
\[
f\left(x,y\right)=\boldsymbol{1}_{\left(x\geqslant1\right)\cap\left(x\geqslant y\right)}\times\left(x-y\right)+\boldsymbol{1}_{\left(x=0\right)\cap\left(2\geqslant y\right)}\times\left(2-y\right).
\]
We have\footnotemark
\[
X_{t+1}=f\left(X_{t},A_{t+1}\right).
\]
Hence, the random variable $\left(X_{0},X_{1},\cdots,X_{t}\right)$
is a measurable function of $\left(A_{1},A_{2},\cdots,A_{t}\right).$
Since the random variables $A_{t},\,t\in\mathbb{N}$ are independent,
the random variables $\left(X_{0},X_{1},\cdots,X_{t}\right)$ and
$A_{t+1}$ are also independent. Therefore, for every function $g$
defined on $E,$ and for every $\left(x_{0},x_{1},\cdots,x_{t}\right)\in E^{t+1},$
\begin{align*}
 & \mathbb{E}^{\left(X_{0},X_{1},\cdots,X_{t}\right)=\left(x_{0},x_{1},\cdots,x_{t}\right)}\left(g\left(X_{t+1}\right)\right)\\
 & \qquad\qquad=\mathbb{E}^{\left(X_{0},X_{1},\cdots,X_{t}\right)=\left(x_{0},x_{1},\cdots,x_{t}\right)}\left(g\circ f\left(x_{t},A_{t+1}\right)\right)\\
 & \qquad\qquad=\mathbb{E}\left(g\circ f\left(x_{t},A_{t+1}\right)\right).
\end{align*}
Denoting, for every $x\in E,$ 
\[
M\left(x,g\right)=\mathbb{E}\left(g\circ f\left(x,A_{t+1}\right)\right),
\]
we obtain
\[
\mathbb{E}^{\mathscr{A}_{t}}\left(g\left(X_{t+1}\right)\right)=M\left(X_{t},g\right).
\]
This proves that $X$ is a homogeneous Markov chain with transition
matrix, the matrix $M$ whose entries are given by
\[
M\left(x,y\right)=M\left(x,\boldsymbol{1}_{\left\{ y\right\} }\right)=P\left(f\left(x,A_{t+1}\right)=y\right).
\]
The law of $A_{t+1}$ leads to the announced matrix $M.$ Indeed,
we obtain succesively:
\begin{itemize}
\item \textbf{Case $x=0$}
\[
f\left(0,A_{t+1}\right)=\left(2-A_{t+1}\right)\boldsymbol{1}_{\left(A_{t+1}\leqslant2\right)}=2\boldsymbol{1}_{\left(A_{t+1}=0\right)}+\boldsymbol{1}_{\left(A_{t+1}=1\right)},
\]
which yields
\begin{align*}
M\left(0,0\right) & =P\left(A_{t+1}\geqslant2\right)=0.15+0.05=0.2,\\
M\left(0,1\right) & =P\left(A_{t+1}=1\right)=0.4,\\
M\left(0,2\right) & =P\left(A_{t+1}=0\right)=0.4.
\end{align*}
\item \textbf{Case $x=1$}
\[
f\left(1,A_{t+1}\right)=\left(1-A_{t+1}\right)\boldsymbol{1}_{\left(1\geqslant A_{t+1}\right)}=\boldsymbol{1}_{\left(A_{t+1}=0\right)},
\]
which yields
\begin{align*}
M\left(1,0\right) & =P\left(A_{t+1}\geqslant1\right)=0.4+0.15+0.05=0.6,\\
M\left(1,1\right) & =P\left(A_{t+1}=0\right)=0.4,\\
M\left(1,2\right) & =P\left(\emptyset\right)=0.
\end{align*}
\item \textbf{Case $x=2$}
\[
f\left(2,A_{t+1}\right)=\left(2-A_{t+1}\right)\boldsymbol{1}_{\left(2\geqslant A_{t+1}\right)}=2\boldsymbol{1}_{\left(A_{t+1}=0\right)}+\boldsymbol{1}_{\left(A_{t+1}=1\right)},
\]
which yields
\begin{align*}
M\left(2,0\right) & =P\left(A_{t+1}\geqslant2\right)=0.15+0.05=0.2,\\
M\left(2,1\right) & =P\left(A_{t+1}=1\right)=0.4,\\
M\left(2,2\right) & =P\left(A_{t+1}=0\right)=0.4.
\end{align*}
\end{itemize}
\end{solution}

\footnotetext{The process is therefore auto-regressive---see Example
$\ref{ex:auto-reg_proc}.$

}

\begin{solution}{}{solexercise17.2}

Let $f$ be the function defined on $\mathbb{N}^{2}$ by
\[
f\left(x,y\right)=\boldsymbol{1}_{\left(x\geqslant1\right)}\times\left(x-1\right)+\boldsymbol{1}_{\left(x=0\right)}\times\left(y-1\right).
\]
We have 
\[
X_{n+1}=f\left(X_{n},Z_{n}\right).
\]
In particular, since the random variables $Z_{n},\,n\in\mathbb{N},$
are independent, the random variables $\left(X_{1},\cdots,X_{n}\right)$
and $Z_{n}$ are also independent. 

Therefore, for every bounded function $g$ on $\mathbb{N},$ and for
every $\left(x_{1},\cdots,x_{n}\right)\in\mathbb{N}^{n},$
\begin{align*}
 & \mathbb{E}^{\left(X_{1},\cdots,X_{n}\right)=\left(x_{1},\cdots,x_{n}\right)}\left(g\left(X_{n+1}\right)\right)\\
 & \qquad\qquad=\mathbb{E}^{\left(X_{1},\cdots,X_{n}\right)=\left(x_{1},\cdots,x_{n}\right)}\left(g\circ f\left(x_{n},Z_{n}\right)\right)\\
 & \qquad\qquad=\mathbb{E}\left(g\circ f\left(x_{n},Z_{n}\right)\right).
\end{align*}

Denoting, for every $x\in\mathbb{N},$ 
\[
M\left(x,g\right)=\mathbb{E}\left(g\circ f\left(x_{n},Z_{n}\right)\right)=\intop_{\mathbb{R}}g\circ f\left(x,z\right)\text{d}\mu\left(z\right),
\]
we obtain
\[
\mathbb{E}^{\mathscr{A}_{n}}\left(g\left(X_{n+1}\right)\right)=M\left(X_{n},g\right).
\]
This proves that $X$ is a homogeneous Markov chain with transition
matrix $M,$ whose entries are given by
\[
M\left(x,y\right)=M\left(x,\boldsymbol{1}_{\left\{ y\right\} }\right)=P\left(f\left(x,Z_{n}\right)=y\right)=\mu\left(\left\{ z\in\mathbb{N}:\,f\left(x,z\right)=y\right\} \right).
\]

They can also be obtained as follows:
\begin{itemize}
\item If $x\in\mathbb{N}^{\ast}$ and $y\in\mathbb{N},$
\[
M\left(x,y\right)=P\left(X_{n+1}=y\mid X_{n}=x\right)=P\left(X_{n}-1=y\mid X_{n}=x\right).
\]
Thus,
\[
M\left(x,y\right)=P\left(X_{n}=y+1\mid X_{n}=x\right)=\boldsymbol{1}_{\left(y=x-1\right)}.
\]
\item If $x=0$ and $y\in\mathbb{N},$
\[
M\left(0,y\right)=P\left(X_{n+1}=y\mid X_{n}=0\right)=P\left(Z_{n}-1=y\mid X_{n}=0\right).
\]
and thus, since $Z_{n}$ and $X_{n}$ are independent,
\[
M\left(0,y\right)=P\left(Z_{n}=y+1\right)=p_{y+1}.
\]
\end{itemize}
Hence, the infinite transition matrix $M$ is of the form
\[
\left(\begin{array}{cccccc}
p_{1} & p_{2} & p_{3} & p_{4} & \cdots & \cdots\\
1 & 0 & 0 & 0 & \cdots & \cdots\\
0 & 1 & 0 & 0 & \cdots & \cdots\\
0 & 0 & 1 & 0 & 0 & \cdots\\
0 & 0 & 0 & 1 & 0 & \cdots\\
\cdots & \cdots & \cdots & \cdots & \cdots & \cdots
\end{array}\right)
\]

\end{solution}

\begin{solution}{}{solexercise17.3}

\textbf{1. Computation of $\mathbb{E}^{\mathscr{A}_{n}}_{i}\left[f\left(\theta_{n}\left(X\right)\right)\right],$
for every $i\in E$}

The simple Markov property ensures that
\[
\mathbb{E}^{\mathscr{A}_{n}}_{i}\left[f\left(\theta_{n}\left(X\right)\right)\right]=\mathbb{E}_{X_{n}}\left[f\left(X\right)\right].
\]

Thus, we compute, for every $i\in E,$ $\mathbb{E}_{i}\left[f\left(X\right)\right].$
We have
\[
\mathbb{E}_{i}\left[f\left(X\right)\right]=P_{i}\left(\sum^{3}_{j=1}\boldsymbol{1}_{\left\{ 1\right\} }\left(X_{j}\right)=1\right),
\]
and consequently
\begin{multline*}
\mathbb{E}_{i}\left[f\left(X\right)\right]=P_{i}\left(X_{1}=1,X_{2}\neq1,X_{3}\neq1\right)+P_{i}\left(X_{1}\neq1,X_{2}=1,X_{3}\neq1\right)\\
+P_{i}\left(X_{1}\neq1,X_{2}\neq1,X_{3}=1\right).
\end{multline*}
Now, for the first term,
\begin{align*}
P_{i}\left(X_{1}=1,X_{2}\neq1,X_{3}\neq1\right) & =\sum_{j,k=2,3}P_{i}\left(X_{1}=1,X_{2}=j,X_{3}=k\right)\\
 & =\sum_{j,k=2,3}M\left(i,1\right)M\left(1,j\right)M\left(j,k\right)\\
 & =M\left(i,1\right)\left[\sum_{j=2,3}M\left(1,j\right)\sum_{k=2,3}M\left(j,k\right)\right].
\end{align*}
Since $\sum_{k=2,3}M\left(j,k\right)=1-M\left(j,1\right),$
\begin{align*}
P_{i}\left(X_{1}=1,X_{2}\neq1,X_{3}\neq1\right) & =M\left(i,1\right)\left[\sum_{j=2,3}M\left(1,j\right)-\sum_{j=2,3}M\left(1,j\right)M\left(j,1\right)\right]\\
 & =M\left(i,1\right)\left[1-\dfrac{1}{3}\cdot\dfrac{1}{2}-\dfrac{2}{3}\cdot\dfrac{1}{2}\right],
\end{align*}
we obtain
\begin{equation}
P_{i}\left(X_{1}=1,X_{2}\neq1,X_{3}\neq1\right)=\dfrac{1}{2}M\left(i,1\right).\label{eq:P_i_X_1_1_X_2_neq_1_X_3_neq_1}
\end{equation}
Using the values of the transition matrix,
\begin{align*}
P_{i}\left(X_{1}\neq1,X_{2}=1,X_{3}\neq1\right) & =\sum_{j,k=2,3}P_{i}\left(X_{1}=j,X_{2}=1,X_{3}=k\right)\\
 & =\sum_{j,k=2,3}M\left(i,j\right)M\left(j,1\right)M\left(1,k\right)\\
 & =\left[\sum_{j=2,3}M\left(i,j\right)M\left(j,1\right)\right]\left[\sum_{k=2,3}M\left(1,k\right)\right],
\end{align*}
and, since $\sum_{k=2,3}M\left(1,k\right)=1,$
\begin{equation}
P_{i}\left(X_{1}\neq1,X_{2}=1,X_{3}\neq1\right)=\dfrac{1}{2}\left[M\left(i,2\right)+M\left(i,3\right)\right].\label{eq:P_iX_1_neq1X_2_1_X_3neq1}
\end{equation}

Similarly, for the second term,
\begin{align*}
P_{i}\left(X_{1}\neq1,X_{2}\neq1,X_{3}=1\right) & =\sum_{j,k=2,3}P_{i}\left(X_{1}=j,X_{2}=k,X_{3}=1\right)\\
 & =\sum_{j,k=2,3}M\left(i,j\right)M\left(j,k\right)M\left(k,1\right).
\end{align*}
Hence,
\begin{multline*}
P_{i}\left(X_{1}\neq1,X_{2}\neq1,X_{3}=1\right)=M\left(i,2\right)\left[\sum_{k=2,3}M\left(2,k\right)M\left(k,1\right)\right]\\
+M\left(i,3\right)\left[\sum_{k=2,3}M\left(3,k\right)M\left(k,1\right)\right].
\end{multline*}
Using the numerical values,
\[
P_{i}\left(X_{1}\neq1,X_{2}\neq1,X_{3}=1\right)=M\left(i,2\right)\left[\dfrac{1}{4}\cdot\dfrac{1}{2}+\dfrac{1}{4}\cdot\dfrac{1}{2}\right]+M\left(i,3\right)\left[\dfrac{1}{2}\cdot\dfrac{1}{2}+0\times\dfrac{1}{2}\right].
\]
Thus,
\begin{equation}
P_{i}\left(X_{1}\neq1,X_{2}\neq1,X_{3}=1\right)=\dfrac{1}{4}\left[M\left(i,2\right)+M\left(i,3\right)\right].\label{eq:P_i_X_1_neq_1_X_2_neq_2_X_3_1}
\end{equation}

By gathering $\refpar{eq:P_i_X_1_1_X_2_neq_1_X_3_neq_1},$ $\refpar{eq:P_iX_1_neq1X_2_1_X_3neq1}$
and $\refpar{eq:P_i_X_1_neq_1_X_2_neq_2_X_3_1},$ we obtain
\begin{align*}
\mathbb{E}_{i}\left[f\left(X\right)\right] & =\dfrac{1}{2}M\left(i,1\right)+\dfrac{3}{4}\left[M\left(i,2\right)+M\left(i,3\right)\right]\\
 & =\dfrac{1}{2}M\left(i,1\right)+\dfrac{3}{4}\left[1-M\left(i,1\right)\right]=\dfrac{3}{4}-\dfrac{1}{4}M\left(i,1\right).
\end{align*}
Thus\boxeq{
\[
\mathbb{E}^{\mathscr{A}_{n}}_{i}\left[f\left(\theta_{n}\left(X\right)\right)\right]=\dfrac{3}{4}-\dfrac{1}{4}M\left(X_{n},1\right).
\]
}

\textbf{2. Justification that $P_{i}\left(T<+\infty\right)=1.$ Computation
of $\mathbb{E}^{\mathscr{A}_{T}}_{i}\left(f\left(\theta_{T}\left(X\right)\right)\right)$}

The graph associated with this Markov chain is shown below.

\begin{center}\begin{tikzpicture}[
  >=Stealth,
  every node/.style={
    circle,
    draw=black,
    minimum size=10mm,
    inner sep=0pt
  },
  edge/.style={
    ->,
    draw=black,
    line width=0.75pt
  }
]

% Nodes (triangular layout)
\node (3) at (0,-3.2) {3};
\node (1) at (-1.8,0) {1};
\node (2) at (1.8,0) {2};

% Self-loop at 2
\draw[edge] (2) to[out=30,in=330,looseness=8] (2);

% Bidirectional edges (bent)
\draw[edge] (1) to[bend left=20] (2);
\draw[edge] (2) to[bend left=20] (1);

\draw[edge] (1) to[bend left=20] (3);
\draw[edge] (3) to[bend left=20] (1);

\draw[edge] (2) to[bend left=20] (3);
\draw[edge] (3) to[bend left=20] (2);

\end{tikzpicture}
\end{center}

This graph allows to ensure that all states communicate. Hence, the
chain is finite and irreducible, and therefore positive recurrent.
In particular, for every $i\in E,$ $P_{i}\left(T<+\infty\right)=1.$
By the strong Markov property and the result of the previous question,
\[
\mathbb{E}^{\mathscr{A}_{T}}_{i}\left[f\left(\theta_{T}\left(X\right)\right)\right]=\mathbb{E}_{X_{T}}\left[f\left(X\right)\right]=\dfrac{3}{4}-\dfrac{1}{4}M\left(X_{T},1\right).
\]
Since $X_{T}=1$ and $M\left(1,1\right)=0,$ we obtain\boxeq{
\[
\mathbb{E}^{\mathscr{A}_{T}}_{i}\left[f\left(\theta_{T}\left(X\right)\right)\right]=\dfrac{3}{4}.
\]
}

\end{solution}

\begin{solution}{}{solexercise17.4}

\textbf{1. Proof of $P_{x}\left(T_{a}<T_{b}\right)=M\left(x,a\right)+\sum_{z\neq a,b}P_{z}\left(T_{a}<T_{b}\right)M\left(x,z\right)$}

Let $\theta$ be the \textbf{shift operator} on $\mathbb{N}^{\mathbb{N}}.$
We write
\[
P_{x}\left(T_{a}<T_{b}\right)=\mathbb{E}_{x}\left[\boldsymbol{1}_{\left(X_{1}=a\right)}\right]+\mathbb{E}_{x}\left[\boldsymbol{1}_{\left(X_{1}\neq a,b\right)}\boldsymbol{1}_{\left(\tau_{a}\left(\theta\left(X\right)\right)<\tau_{b}\left(\theta\left(X\right)\right)\right)}\right].
\]
Conditioning with respect to $\mathscr{A}_{1},$ the simple Markov
property yields
\begin{align*}
P_{x}\left(T_{a}<T_{b}\right) & =M\left(x,a\right)+\mathbb{E}_{x}\left[\boldsymbol{1}_{\left(X_{1}\neq a,b\right)}\mathbb{E}^{\mathscr{A}_{1}}_{x}\left[\boldsymbol{1}_{\left(\tau_{a}\left(\theta\left(X\right)\right)<\tau_{b}\left(\theta\left(X\right)\right)\right)}\right]\right]\\
 & =M\left(x,a\right)+\mathbb{E}_{x}\left[\boldsymbol{1}_{\left(X_{1}\neq a,b\right)}\mathbb{E}_{X_{1}}\left[\boldsymbol{1}_{\left(\tau_{a}\left(X\right)<\tau_{b}\left(X\right)\right)}\right]\right].
\end{align*}
Hence,
\[
P_{x}\left(T_{a}<T_{b}\right)=M\left(x,a\right)+\mathbb{E}_{x}\left[\sum_{z\neq a,b}\boldsymbol{1}_{\left(X_{1}=z\right)}\mathbb{E}_{z}\left[\boldsymbol{1}_{\left(T_{a}<T_{b}\right)}\right]\right].
\]
Factoring the expectation gives
\[
P_{x}\left(T_{a}<T_{b}\right)=M\left(x,a\right)+\sum_{z\neq a,b}\left(\mathbb{E}_{x}\left[\boldsymbol{1}_{\left(X_{1}=z\right)}\right]\mathbb{E}_{z}\left[\boldsymbol{1}_{\left(T_{a}<T_{b}\right)}\right]\right),
\]
which yields relation $\refpar{eq:P_xT_alessT_b}.$

\textbf{2. Explicit computation of $P_{x}\left(T_{0}<T_{b}\right)$}

The function $f$ satisfies the system
\[
\left\{ \begin{array}{rl}
f\left(x\right) & =p_{x}f\left(x+1\right)+q_{x}f\left(x-1\right)+r_{x}f\left(x\right),\,\,\,\,\text{if }1<x<b,\\
f\left(1\right) & =p_{1}f\left(2\right)+q_{1}+r_{1}f\left(1\right)\\
f\left(b-1\right) & =q_{b-1}f\left(b-2\right)+r_{b-1}f\left(b-1\right).
\end{array}\right.
\]
This system can be rewritten by extending the function $f$ with the
conditions $f\left(0\right)=1$ and $f\left(b\right)=0,$ although
this extension should not be interpreted probabilistically.

We then obtain
\[
f\left(x+1\right)-f\left(x\right)=\dfrac{q_{x}}{p_{x}}\left[f\left(x\right)-f\left(x-1\right)\right].
\]
Iterating this relation yields
\[
f\left(x+1\right)-f\left(x\right)=\dfrac{q_{x}q_{x-1}\cdots q_{1}}{p_{x}p_{x-1}\cdots p_{1}}\left[f\left(1\right)-f\left(0\right)\right].
\]
Thus, for $0\leqslant x<b,$
\begin{equation}
f\left(x\right)-f\left(x+1\right)=a_{x}\left[f\left(0\right)-f\left(1\right)\right].\label{eq:fx-fxm1}
\end{equation}
Summing over $x$ from 0 to $b-1$, allow us to identify $f\left(0\right)-f\left(1\right),$
\[
\left[\sum^{b-1}_{x=0}a_{x}\right]\left[f\left(0\right)-f\left(1\right)\right]=\sum^{b-1}_{x=0}\left[f\left(x\right)-f\left(x+1\right)\right]=f\left(0\right)-f\left(b\right)=1.
\]
Substituting into relation $\refpar{eq:fx-fxm1}$ gives
\[
f\left(x\right)-f\left(x+1\right)=\dfrac{a_{x}}{\sum^{b-1}_{y=0}a_{y}}.
\]
Summing again the increments of $f,$ this time starting from $x,$
we obtain
\[
f\left(x\right)=\sum^{b-1}_{y=x}\left[f\left(x\right)-f\left(x+1\right)\right]=\dfrac{\sum^{b-1}_{y=x}a_{y}}{\sum^{b-1}_{y=0}a_{y}}.
\]
Hence, for every $x\in\left]0,b\right[,$\boxeq{
\[
P_{x}\left(T_{0}<T_{b}\right)=f\left(x\right)=\dfrac{\sum^{b-1}_{y=x}a_{y}}{\sum^{b-1}_{y=0}a_{y}}.
\]
}

\end{solution}

\begin{solution}{}{solexercise17.5}

\textbf{1. Graph associated with the chain. Communication classes}

The graph associated with this chain is

\begin{center}\begin{tikzpicture}[
  >=Stealth,
  every node/.style={
    circle,
    draw=black,
    minimum size=10mm,
    inner sep=0pt
  },
  edge/.style={
    ->,
    draw=black,
    line width=0.75pt
  }
]

% Nodes (triangular layout)
\node (3) at (0,-1.8) {3};
\node (1) at (-2,1.5) {1};
\node (2) at (2,1.5) {2};

% Self-loops
\draw[edge] (1) to[out=210,in=150,looseness=8] (1);
\draw[edge] (2) to[out=-30,in=30,looseness=8] (2);
\draw[edge] (3) to[out=240,in=300,looseness=8] (3);

% Bidirectional transitions (bent)
\draw[edge] (2) to[bend right=20] (1);
\draw[edge] (3) to[bend left=20] (1);
\draw[edge] (2) to[bend left=20] (3);
\draw[edge] (3) to[bend left=20] (2);

\end{tikzpicture}
\end{center}

This graph associated with this chain highlights two communication
classes $\left\{ 1\right\} $ and $\left\{ 2,3\right\} .$ The state
$1$ is recurrent, and in fact absorbing. The class $\left\{ 2,3\right\} $
is transient.

\textbf{2. $\left(f_{k}\right)_{k\in\mathbb{N}^{\ast}}$ solution
of $f_{k}=Qf_{k-1}.$ Computation of $f_{k}$}

By Proposition $\ref{pr:first_hitting_point}$ and more precisely,
the equality $\refpar{eq:F_k(x_y)},$
\[
\begin{cases}
F_{1}\left(x,3\right)=M\left(x,3\right),\\
F_{k}\left(x,3\right)=\sum_{z\in E\backslash\left\{ 3\right\} }M\left(x,z\right)F_{k-1}\left(z,3\right), & \text{if }k\geqslant2.
\end{cases}
\]
It follows that the sequence of vector $f_{k}$ is the solution of
$\refpar{eq:fx-fxm1},$ where $f_{1}$ is the last column of $M,$
and $Q$ is the matrix obtained from $M$ by replacing the last column
with zeros,
\[
f_{1}=\left(\begin{array}{c}
0\\
\dfrac{1}{3}\\
\dfrac{1}{15}
\end{array}\right)\,\,\,\,\text{and}\,\,\,\,Q=\left(\begin{array}{ccc}
1 & 0 & 0\\
\dfrac{1}{2} & \dfrac{1}{6} & 0\\
\dfrac{1}{3} & \dfrac{3}{5} & 0
\end{array}\right).
\]
Thus, 
\[
f_{k}=Q^{k-1}f_{1}.
\]
Rather than computing the powers of $Q$ directly, we compute $f_{k}$
using the Cayley-Hamilton theorem. Since $Q$ has the distinct simple
eigenvalues 1, $\dfrac{1}{6}$ and 0, which are the roots of its characteristic
polynomial,
\[
Q\left(Q-1\right)\left(Q-\dfrac{1}{6}\right)=0.
\]
Expanding this relations yields
\[
Q^{3}-Q^{2}=\dfrac{1}{6}\left(Q^{2}-Q\right).
\]
Hence, for every $n\geqslant3,$
\[
Q^{n}-Q^{n-1}=\dfrac{1}{6}\left(Q^{n-1}-Q^{n-2}\right),
\]
Summing these relations and applying them to $f_{1}$ gives
\[
\sum^{n}_{j=3}\left(Q^{j}-Q^{j-1}\right)f_{1}=\dfrac{1}{6}\sum^{n}_{j=3}\left(Q^{j-1}-Q^{j-2}\right)f_{1},
\]
which implies
\[
Q^{n}f_{1}-Q^{2}f_{1}=\dfrac{1}{6}\left(Q^{n-1}f_{1}-Qf_{1}\right).
\]
Thus, for every $n\geqslant3,$
\[
f_{n+1}=\dfrac{1}{6}f_{n}+\left(Q^{2}f_{1}-\dfrac{1}{6}Qf_{1}\right).
\]
Consequently,
\begin{align*}
f_{n+1} & =\dfrac{1}{6^{n-2}}f_{3}+\left(1+\dfrac{1}{6}+\dfrac{1}{6^{2}}+\cdots+\dfrac{1}{6^{n-3}}\right)\left(Q^{2}f_{1}-\dfrac{1}{6}Qf_{1}\right)\\
 & =\dfrac{1}{6^{n-2}}Q^{2}f_{1}+\dfrac{1-\dfrac{1}{6^{n-2}}}{1-\dfrac{1}{6}}\left(Q^{2}f_{1}-\dfrac{1}{6}Qf_{1}\right)\\
 & =\dfrac{6}{5}\left(1-\dfrac{1}{6^{n-1}}\right)Q^{2}f_{1}-\dfrac{1}{5}\left(1-\dfrac{1}{6^{n-2}}\right)Qf_{1}.
\end{align*}
Since
\[
f_{2}\equiv Qf_{1}=\left(\begin{array}{c}
0\\
\dfrac{1}{18}\\
\dfrac{1}{5}
\end{array}\right)\,\,\,\,\text{and}\,\,\,\,f_{3}=Q^{2}f_{1}=\dfrac{1}{6}\left(\begin{array}{c}
0\\
\dfrac{1}{18}\\
\dfrac{1}{5}
\end{array}\right)=\dfrac{1}{6}Qf_{1},
\]
we obtain, after simplification, for $n\geqslant3,$
\[
f_{n+1}=\dfrac{1}{6^{n-1}}\left(\begin{array}{c}
0\\
\dfrac{1}{18}\\
\dfrac{1}{5}
\end{array}\right).
\]
Note that this formula is in fact still valid for $n=1,2.$

\textbf{3. Computation of $P_{x}\left(T_{3}=+\infty\right)$}

We have
\[
P_{x}\left(T_{3}=+\infty\right)=1-P_{x}\left(T_{3}<+\infty\right)=1-\sum_{k\in\mathbb{N}^{\ast}}P_{x}\left(T_{3}=k\right)=1-\sum_{k\in\mathbb{N}^{\ast}}f_{k}\left(x\right).
\]
Clearly, $P_{1}\left(T_{3}=+\infty\right)=1.$ Moreover, for every
$n\in\mathbb{N}^{\ast},$
\[
f_{n}\left(2\right)=\dfrac{1}{3}\cdot\dfrac{1}{6^{n-1}},
\]
hence
\[
P_{2}\left(T_{3}=+\infty\right)=1-\sum_{k\in\mathbb{N}^{\ast}}\dfrac{1}{3}\cdot\dfrac{1}{6^{k-1}}=\dfrac{3}{5}.
\]
Finally, since for every $n\geqslant2,$ $f_{n}\left(3\right)=\dfrac{1}{5}\cdot\dfrac{1}{6^{n-2}},$
and $f_{1}\left(3\right)=\dfrac{1}{15},$ we obtain
\[
P_{3}\left(T_{3}=+\infty\right)=1-\left[\dfrac{1}{15}+\sum_{k\geqslant2}\dfrac{1}{5}\cdot\dfrac{1}{6^{k-2}}\right]=\dfrac{52}{75}.
\]

\textbf{4. Determination of $\mathbb{E}_{3}\left(N_{3}\right)$ and
$\mathbb{E}_{2}\left(N_{3}\right).$}

The mean number $R\left(x,y\right)$ of visits to state $y$ by the
chain starting from $x$ at time 0 is given by relations $\refpar{eq:potential_mat_pr}.$
Thus,
\[
\mathbb{E}_{3}\left(N_{3}\right)=R\left(3,3\right)=\dfrac{1}{1-P_{3}\left(T_{3}<+\infty\right)}=\dfrac{75}{52}.
\]
Similarly,
\[
\mathbb{E}_{2}\left(N_{3}\right)=R\left(2,3\right)=F\left(2,3\right)R\left(3,3\right)=\dfrac{15}{26}.
\]

\end{solution}

\begin{solution}{}{solexercise17.6}

\textbf{1. Graph associated with this chain. Communication classes.
Period and nature of the classes}

The graph associated with the chain is

\begin{center}\begin{tikzpicture}[
  >=Stealth,
  every node/.style={
    circle,
    draw=black,
    minimum size=10mm,
    inner sep=0pt
  },
  edge/.style={
    ->,
    draw=black,
    line width=0.75pt
  }
]

% Nodes (pentagon-like layout)
\node (1) at (0,0) {1};
\node (2) at (4,0) {2};
\node (3) at (-4,0) {3};
\node (4) at (0,-4) {4};
\node (5) at (-4,-4) {5};

% Transitions involving 1, 2, 3
\draw[edge] (1) to[bend left=15] (2);
\draw[edge] (2) to[bend left=15] (1);

\draw[edge] (1) to[bend left=15] (3);
\draw[edge] (3) to[bend left=15] (1);

% Transitions involving 4 and 5
\draw[edge] (4) -- (1);
\draw[edge] (4) to[bend left=15] (5);

\draw[edge] (5) to[bend left=15] (4);
\draw[edge] (5) -- (3);

\end{tikzpicture}
\end{center}

This graph highlights two communitation classes $C=\left\{ 1,2,3\right\} $
and $D=\left\{ 4,5\right\} .$

We have $M\left(4,4\right)=0$ and $M^{2}\left(4,4\right)\geqslant M\left(4,5\right)M\left(5,4\right)>0.$
A simple argument by induction shows that, for every $n\in\mathbb{N}^{\ast},$
$M^{2n}\left(4,4\right)>0$ and $M^{2n+1}\left(4,4\right)=0.$ It
follows that state $4,$ and hence also state 5, has period $2.$
This is the period of the class $D.$

The same argument shows that state $2,$ and hence also states 1 and
3, has period 2. This is the period of the class $C.$ 

State 5 leads to 3, but 3 does not lead to 5. Hence, 5 is not recurrent:
it is transient. Therefore, the class $D$ is transient. Since $E$
is finite and $C$ is a closed class, $C$ is positive recurrent.

\textbf{2. Existence of a unique invariant probability $\nu.$ Computation
of $\nu.$ Computation of $\mathbb{E}_{x}\left(T_{x}\right)$ for
every $x\in\left\{ 1,2,3\right\} $}

The existence of a unique invariant probability $\nu$ follows from
the existence of a unique positive recurrent class. Let $\nu=\left(a,b,c,d,e\right)$
be a row vector corresponding to the invariant measure on $E.$ Then
$\nu$ is a solution of the system
\[
\left(a,b,c,d,e\right)\left(\begin{array}{ccccc}
0 & \dfrac{2}{3} & \dfrac{1}{3} & 0 & 0\\
1 & 0 & 0 & 0 & 0\\
1 & 0 & 0 & 0 & 0\\
\dfrac{1}{2} & 0 & 0 & 0 & \dfrac{1}{2}\\
0 & 0 & \dfrac{1}{2} & \dfrac{1}{2} & 0
\end{array}\right)=\left(a,b,c,d,e\right).
\]
Hence,
\[
\left\{ \begin{array}{rl}
b+c+\dfrac{d}{2} & =a\\
\dfrac{2}{3} & =b\\
\dfrac{a}{3}+\dfrac{e}{2} & =c\\
\dfrac{e}{2} & =d\\
\dfrac{d}{2} & =e
\end{array}\right.\Longleftrightarrow\left\{ \begin{array}{l}
d=e=0\\
b=\dfrac{2}{3}a\\
c=\dfrac{a}{3}
\end{array}\right..
\]
The invariant probability is such that $a+b+c+d+e=1,$ which yields
$a=\dfrac{1}{2}.$ Thus, the invariant probability $\nu$ is 
\[
\nu=\left(\dfrac{1}{2},\dfrac{1}{3},\dfrac{1}{6},0,0\right).
\]
 We know that for every $x\in C,$ $\nu\left(x\right)=\dfrac{1}{\mathbb{E}_{x}\left(T_{x}\right)},$
which yields\boxeq{
\[
\mathbb{E}_{1}\left(T_{1}\right)=2\,\,\,\,\mathbb{E}_{2}\left(T_{2}\right)=3\,\,\,\,\mathbb{E}_{3}\left(T_{3}\right)=6.
\]
}

\textbf{3. Potential matrix $R$}

The potential matrix $R$ has the block structure
\[
\begin{array}{ccc}
 &  & \begin{array}{cc}
\,\,C\,\, & D\,\end{array}\\
R= & \begin{array}{c}
C\\
D
\end{array} & \left(\begin{array}{cc}
+\infty & \,\boldsymbol{0}\,\\
+\infty & \,S\,
\end{array}\right),
\end{array}
\]
where $Q=M_{\mid_{D\times D}}$ and $S=\sum^{+\infty}_{n=0}Q^{n},$
so that $S=\left(I-Q\right)^{-1}.$

We have
\[
Q=\left(\begin{array}{cc}
0 & \dfrac{1}{2}\\
\dfrac{1}{2} & 0
\end{array}\right)\,\,\,\,\text{hence\,\,\,\,}I-Q=\left(\begin{array}{cc}
1 & -\dfrac{1}{2}\\
-\dfrac{1}{2} & 1
\end{array}\right)
\]
Therefore, 
\[
S=\left(\begin{array}{cc}
\dfrac{4}{3} & \dfrac{2}{3}\\
\dfrac{2}{3} & \dfrac{4}{3}
\end{array}\right).
\]

\textbf{4. Computation of $P_{4}\left(T_{5}<+\infty\right)$ and $P_{5}\left(T_{5}<+\infty\right)$}

From the equalities $\refpar{eq:potential_mat_pr},$ we deduce that
\[
P_{4}\left(T_{5}<+\infty\right)=\dfrac{R\left(4,5\right)}{R\left(5,5\right)}=\dfrac{\frac{2}{3}}{\frac{4}{3}}=\dfrac{1}{2},
\]
and
\[
P_{5}\left(T_{5}<+\infty\right)=1-\dfrac{1}{R\left(5,5\right)}1-\dfrac{1}{\frac{4}{3}}=\dfrac{1}{4}.
\]

\textbf{5. Computation of $P_{4}\left(N_{5}=m\right)$ and $P_{5}\left(N_{5}=m\right)$
for every $m\in\mathbb{N}$ }

The law of $N_{5}$ under $P_{4}$ is given---see Proposition $\ref{pr:reformulation_m_visit_withF}$---by
\[
P_{4}\left(N_{5}=m\right)=\begin{cases}
P_{4}\left(T_{5}=+\infty\right), & \text{if\,}m=0,\\
P_{4}\left(T_{5}<+\infty\right)\left[P_{5}\left(T_{5}<+\infty\right)\right]^{m-1}P_{5}\left(T_{5}=+\infty\right), & \text{if\,}m\in\mathbb{N}^{\ast}.
\end{cases}
\]
Hence,
\[
P_{4}\left(N_{5}=m\right)=\begin{cases}
\dfrac{1}{2}, & \text{if\,}m=0,\\
\dfrac{3}{8}\left(\dfrac{1}{4}\right)^{m-1}, & \text{if}\,m\in\mathbb{N}^{\ast}.
\end{cases}
\]
Similarly, the law of $N_{5}$ under $P_{5}$ is given by, for every
$m\in\mathbb{N}^{\ast},$
\[
P_{5}\left(N_{5}=m\right)=\left[P_{5}\left(T_{5}<+\infty\right)\right]^{m-1}P_{5}\left(T_{5}=+\infty\right),
\]
hence
\[
P_{5}\left(N_{5}=m\right)=\dfrac{3}{4}\left(\dfrac{1}{4}\right)^{m-1}.
\]
This is the geometric law $\mathscr{G}_{\mathbb{N}^{\ast}}\left(\dfrac{3}{4}\right).$

\textbf{6. Proof of $\mathbb{E}_{x}\left[T_{y}\boldsymbol{1}_{\left(T_{y}<+\infty\right)}\right]=M\left(x,y\right)+\sum_{z\neq y}M\left(x,z\right)\left[\mathbb{E}_{z}\left(T_{y}\boldsymbol{1}_{\left(T_{y}<+\infty\right)}\right)+P_{z}\left(T_{y}<+\infty\right)\right].$
Numerical computation of $\mathbb{E}_{x}\left(T_{1}\right)$}

Let $F_{k}\left(x,y\right)=P_{x}\left(T_{y}=k\right).$ By Proposition
$\ref{pr:first_hitting_point},$ and more precisely the equality $\refpar{eq:F_k(x_y)},$
\[
\begin{cases}
F_{1}\left(x,y\right)=M\left(x,y\right),\\
F_{k}\left(x,y\right)=\sum_{z\neq y}M\left(x,z\right)F_{k-1}\left(z,y\right), & \text{if }k\geqslant2.
\end{cases}
\]
It follows that
\begin{align*}
 & \mathbb{E}_{x}\left[T_{y}\boldsymbol{1}_{\left(T_{y}<+\infty\right)}\right]=\sum_{k\in\mathbb{N}^{\ast}}kF_{k}\left(x,y\right)\\
 & \qquad\qquad=M\left(x,y\right)+\sum^{+\infty}_{k=2}k\left[\sum_{z\neq y}M\left(x,z\right)F_{k-1}\left(z,y\right)\right]\\
 & \qquad\qquad=M\left(x,y\right)+\sum_{z\neq y}M\left(x,z\right)\left[\sum^{+\infty}_{l=1}\left(l+1\right)F_{l}\left(z,y\right)\right]\\
 & \qquad\qquad=M\left(x,y\right)+\sum_{z\neq y}M\left(x,z\right)\left[\mathbb{E}_{z}\left(T_{y}\boldsymbol{1}_{\left(T_{y}<+\infty\right)}+P_{z}\left(T_{y}<+\infty\right)\right)\right].
\end{align*}
This is the equality $\refpar{eq:fx-fxm1}$---since all terms are
nonnegative, the interchange of sums is justified.

In particular, since for every $x\in E,$ $P_{x}\left(T_{1}<+\infty\right)=1,$
it follows that
\[
\mathbb{E}_{x}\left(T_{1}\right)=M\left(x,1\right)+\sum_{z\neq1}M\left(x,z\right)\mathbb{E}_{z}\left(T_{1}\right)+\sum_{z\neq1}M\left(x,z\right),
\]
hence
\[
\mathbb{E}_{x}\left(T_{1}\right)=1+\sum_{z\neq1}M\left(x,z\right)\mathbb{E}_{z}\left(T_{1}\right).
\]
Let $g$ be the vector with components $g\left(x\right)=\mathbb{E}_{x}\left(T_{1}\right),\,x\in E.$
Then $g$ satisfies
\[
g=\left(\begin{array}{c}
1\\
1\\
1\\
1\\
1
\end{array}\right)+\left(\begin{array}{ccccc}
0 & \dfrac{2}{3} & \dfrac{1}{3} & 0 & 0\\
0 & 0 & 0 & 0 & 0\\
0 & 0 & 0 & 0 & 0\\
0 & 0 & 0 & 0 & \dfrac{1}{2}\\
0 & 0 & \dfrac{1}{2} & \dfrac{1}{2} & 0
\end{array}\right)g,
\]
or, equivalently, its components satisfy the system
\[
\left\{ \begin{array}{rl}
g_{1} & =1+\dfrac{2}{3}g_{2}+\dfrac{1}{3}g_{3}\\
g_{2} & =1\\
g_{3} & =1\\
g_{4} & =1+\dfrac{1}{2}g_{5}\\
g_{5} & =1+\dfrac{1}{2}g_{3}+\dfrac{1}{2}g_{4}
\end{array}\right.\Longleftrightarrow\left\{ \begin{array}{rl}
g_{1} & =2\\
g_{2} & =1\\
g_{3} & =1\\
g_{4} & =\dfrac{7}{3}\\
g_{5} & =\dfrac{8}{3}.
\end{array}\right.
\]
Therefore,\boxeq{
\[
\mathbb{E}_{1}\left[T_{1}\right]=2\,\,\,\,\mathbb{E}_{2}\left[T_{1}\right]=1\,\,\,\,\mathbb{E}_{3}\left[T_{1}\right]=1\,\,\,\,\mathbb{E}_{4}\left[T_{1}\right]=\dfrac{7}{3}\,\,\,\,\mathbb{E}_{5}\left[T_{1}\right]=\dfrac{8}{3}.
\]
}

\end{solution}

\begin{solution}{}{solexercise17.7}

\textbf{1. Proof of $G\left(x,y\right)=P_{x}\left(T_{y}<+\infty\right)+\sum_{z\neq y}M\left(x,z\right)G\left(z,y\right).$}

This question is essentially the same as Question 6 of Exercise $\ref{exo:exercise17.6},$
up to one additional line of computation. We thus refer to those calculations,
which remain valid whether $E$ is finite or infinite. Equality $\refpar{eq:E_x_T_y_1_t_yless_inf}$
can be written, with the notation of the present exercise and after
regrouping the terms, as
\[
G\left(x,y\right)=M\left(x,y\right)+\sum_{z\neq y}M\left(x,z\right)P_{z}\left(T_{y}<+\infty\right)+\sum_{z\neq y}M\left(x,z\right)G\left(z,y\right).
\]
Proposition $\ref{pr:first_hitting_point}$ and more precisely equality
$\refpar{eq:F_as_fixed_point},$ then yields equality $\refpar{eq:G_x_y}.$

2. \textbf{Proof that $\phi_{\left(C,1\right)}$ is solution of $\forall x\in E,\,\,\,\,\phi_{C,D}\left(x\right)=\sum_{y\in C}M\left(x,y\right)P_{y}\left(T_{D}=+\infty\right)+\sum_{y\notin C}M\left(x,y\right)\phi_{C,D}\left(y\right)$}

We have
\begin{align*}
\phi_{C,D}\left(x\right) & =\mathbb{E}_{x}\left[\boldsymbol{1}_{\left(X_{1}\in C\right)}\boldsymbol{1}_{\left(\tau_{D}\left(\theta_{1}\left(X\right)\right)=+\infty\right)}\right]\\
 & \qquad+\mathbb{E}_{x}\left[\boldsymbol{1}_{\left(X_{1}\notin C\right)}\boldsymbol{1}_{\left(\tau_{C}\left(\theta_{1}\left(X\right)\right)<+\infty\right)}\boldsymbol{1}_{\left(\tau_{D}\left(\theta_{\tau_{C}\left(\theta_{1}\left(X\right)\right)}\left(X\right)\right)=+\infty\right)}\right].
\end{align*}
Since $X_{1}$ is $\mathscr{A}_{1}-$measurable,
\begin{align*}
\phi_{C,D}\left(x\right) & =\mathbb{E}_{x}\left[\boldsymbol{1}_{\left(X_{1}\in C\right)}\mathbb{E}^{\mathscr{A}_{1}}_{x}\left(\boldsymbol{1}_{\left(\tau_{D}\left(\theta_{1}\left(X\right)\right)=+\infty\right)}\right)\right]\\
 & \qquad+\mathbb{E}_{x}\left[\boldsymbol{1}_{\left(X_{1}\notin C\right)}\mathbb{E}^{\mathscr{A}_{1}}_{x}\left(\boldsymbol{1}_{\left(\tau_{C}\left(\theta_{1}\left(X\right)\right)<+\infty\right)}\boldsymbol{1}_{\left(\tau_{D}\left(\theta_{\tau_{C}\left(\theta_{1}\left(X\right)\right)}\left(X\right)\right)=+\infty\right)}\right)\right].
\end{align*}
By the simple Markov property,
\begin{align*}
\phi_{C,D}\left(x\right) & =\mathbb{E}_{x}\left[\boldsymbol{1}_{\left(X_{1}\in C\right)}\mathbb{E}_{X_{1}}\left(\boldsymbol{1}_{\left(T_{D}=+\infty\right)}\right)\right]\\
 & \qquad+\mathbb{E}_{x}\left[\boldsymbol{1}_{\left(X_{1}\notin C\right)}\mathbb{E}_{X_{1}}\left(\boldsymbol{1}_{\left(T_{C}<+\infty\right)}\boldsymbol{1}_{\left(\tau_{D}\left(\theta_{T_{C}}\left(X\right)\right)=+\infty\right)}\right)\right],
\end{align*}
Thus,
\begin{align*}
\phi_{C,D}\left(x\right) & =\sum_{y\in C}M\left(x,y\right)P_{y}\left(T_{D}=+\infty\right)\\
 & \qquad+\sum_{y\notin C}M\left(x,y\right)P_{y}\left[\left(T_{C}<+\infty\right)\cap\left(\tau_{D}\left(\theta_{T_{C}}\left(X\right)\right)=+\infty\right)\right],
\end{align*}
which is the equality $\refpar{eq:system_phi_CD_x}.$

\textbf{3. Graph associated with this chain. Communication classes.
Nature of these classes}

The graph associated with this chain is

\begin{center}\begin{tikzpicture}[
  >=Stealth,
  every node/.style={
    circle,
    draw=black,
    minimum size=10mm,
    inner sep=0pt
  },
  edge/.style={
    ->,
    draw=black,
    line width=0.75pt
  }
]

% Nodes (left-to-right chain with absorbing states 4 and 5)
\node (1) at (0,0) {1};
\node (2) at (2.5,0) {2};
\node (3) at (5,0) {3};
\node (4) at (5,-2.5) {4};
\node (5) at (2.5,-2.5) {5};

% Self-loops (q1,q2,q3 and absorbing 4,5)
\draw[edge] (1) to[out=150,in=210,looseness=8] (1);
\draw[edge] (2) to[out=120,in=60,looseness=8] (2);
\draw[edge] (3) to[out=30,in=-30,looseness=8] (3);
\draw[edge] (4) to[out=-120,in=-60,looseness=8] (4);
\draw[edge] (5) to[out=-120,in=-60,looseness=8] (5);

% Forward transitions along the chain (r1,r2,r3)
\draw[edge] (1) -- (2);
\draw[edge] (2) -- (3);
\draw[edge] (3) -- (4);

% Transitions to absorbing state 5 (p1,p2,p3)
\draw[edge] (1) -- (5);
\draw[edge] (2) -- (5);
\draw[edge] (3) -- (5);

\end{tikzpicture}\end{center}

The graph shows that no state communicates with any other. Hence,
there are five communication classes, namely the singletons $\left\{ i\right\} ,\,1\leqslant i\leqslant5.$
States 1, 2 and 3 are transient, whereas states 4 and 5 are absorbing---since
the corresponding classes are closed.

\textbf{4. Justification of $F\left(5,4\right)=0.$ Proof that $K$
is solution of the form $K=b+TK.$ Value of $K.$ Computation of $P_{1}\left(T_{4}<+\infty\right)$}

States $4$ and $5$ do not lead to one another. In particular, $F\left(5,4\right)=0.$
By Proposition $\ref{pr:first_hitting_point}$ and more precisely
the equality $\refpar{eq:F_as_fixed_point},$ we have, for $1\leqslant i\leqslant5,$
\[
F\left(i,4\right)=M\left(i,4\right)+\sum_{j\neq4}M\left(i,j\right)F\left(j,4\right).
\]
Hence, since $F\left(5,4\right)=0,$
\[
F\left(i,4\right)=M\left(i,4\right)+\sum^{3}_{j=1}M\left(i,j\right)F\left(j,4\right).
\]
Therefore,
\[
K=b+TK,\,\,\,\,\text{where}\,\,\,\,b=\left(\begin{array}{c}
0\\
0\\
r_{3}
\end{array}\right)\,\,\,\,\text{and}\,\,\,\,T=\left(\begin{array}{ccc}
q_{1} & r_{1} & 0\\
0 & q_{2} & r_{2}\\
0 & 0 & q_{3}
\end{array}\right),
\]
This equation is equivalent to
\[
\left(\begin{array}{ccc}
1-q_{1} & -r_{1} & 0\\
0 & 1-q_{2} & -r_{2}\\
0 & 0 & 1-q_{3}
\end{array}\right)K=\left(\begin{array}{c}
0\\
0\\
r_{3}
\end{array}\right)\,\,\,\,\text{or}\,\,\,\,\left\{ \begin{array}{rl}
\left(1-q_{1}\right)k_{1}-r_{1}k_{2} & =0\\
\left(1-q_{2}\right)k_{2}-r_{2}k_{3} & =0\\
\left(1-q_{3}\right)k_{3}-r_{3} & =0.
\end{array}\right.
\]
This yields,
\[
k_{3}=\dfrac{r_{3}}{1-q_{3}}\,\,\,\,k_{2}=\dfrac{r_{2}r_{3}}{\left(1-q_{2}\right)\left(1-q_{3}\right)}\,\,\,\,k_{1}=\dfrac{r_{1}r_{2}r_{3}}{\left(1-q_{1}\right)\left(1-q_{2}\right)\left(1-q_{3}\right)}.
\]
It follows that
\[
P_{1}\left(T_{4}<+\infty\right)=\dfrac{r_{1}r_{2}r_{3}}{\left(1-q_{1}\right)\left(1-q_{2}\right)\left(1-q_{3}\right)}.
\]

\textbf{5. Potential matrix $R$ of $X$}

Let $\text{Tr}=\left\{ 1,2,3\right\} $ be the ordered set of transient
states, and $A=\left\{ 4,5\right\} $ the ordered set of absorbing
points---hence positive recurrent, since $E$ is finite. The potential
matrix $R$ then has the block structure
\[
\begin{array}{rcl}
 &  & \begin{array}{rl}
\,\,\,\,\,\,\,\text{Tr}\,\,\,\,\,\,\,\,\, & \,\,\,A\end{array}\\
R= & \begin{array}{c}
\text{Tr}\\
A
\end{array} & \left(\begin{array}{rl}
\left(I-T\right)^{-1} & +\infty\\
\boldsymbol{0}\,\,\,\,\,\,\,\,\, & +\infty
\end{array}\right)
\end{array}
\]
where $T$ is the matrix introduced in the previous question. 

Since the matrix $I-T$ is triangular, its inverse can be obtained,
for instance, by solving the associated system
\[
\left\{ \begin{array}{rl}
\left(1-q_{1}\right)x_{1}-r_{1}x_{3} & =a\\
\left(1-q_{2}\right)x_{2}-r_{2}x_{3} & =b\\
\left(1-q_{3}\right)x_{3} & =c.
\end{array}\right.
\]
This gives
\begin{gather*}
x_{3}=\dfrac{1}{1-q_{3}}c\,\,\,\,x_{2}=\dfrac{1}{1-q_{2}}\left[b+\dfrac{r_{2}}{1-q_{3}}c\right]\\
x_{1}=\dfrac{1}{1-q_{1}}\left[a+\dfrac{r_{1}}{1-q_{2}}b+\dfrac{r_{1}r_{2}}{\left(1-q_{2}\right)\left(1-q_{3}\right)}c\right].
\end{gather*}
Hence,
\[
\left(I-T\right)^{-1}=\left(\begin{array}{ccc}
\dfrac{1}{1-q_{1}} & \dfrac{r_{1}}{\left(1-q_{1}\right)\left(1-q_{2}\right)} & \dfrac{r_{1}r_{2}}{\left(1-q_{1}\right)\left(1-q_{2}\right)\left(1-q_{3}\right)}\\
0 & \dfrac{1}{1-q_{2}} & \dfrac{r_{2}}{\left(1-q_{2}\right)\left(1-q_{3}\right)}\\
0 & 0 & \dfrac{1}{1-q_{3}}
\end{array}\right)
\]

\textbf{6. Proof that $H\left(5\right)=0.$ Computation of $H\left(x\right)$
for $x=1,2,3.$ Interpretation of $H\left(1\right)$}

Since $P_{5}\left(T_{4}<+\infty\right)=0,$ we have $H\left(5\right)=0.$
Moreover, with the notation introduced earlier, 
\[
H\left(x\right)=G\left(x,4\right).
\]
It then follows from $\refpar{eq:G_x_y}$ that, for $1\leqslant x\leqslant5,$
\[
H\left(x\right)=P_{x}\left(T_{4}<+\infty\right)+\sum_{z\neq4}M\left(x,z\right)H\left(z\right).
\]
Hence, under vector form, since $H\left(5\right)=0,$
\[
\widehat{H}=K+T\widehat{H},
\]
where $\widehat{H}$ is the vector in $\mathbb{R}^{3}$ with components
$H\left(i\right),\,i=1,2,3.$ 

Thus, 
\[
\left(I-T\right)^{-1}\widehat{H}=K.
\]
Using the results of the previous question, a matrix computation shows
that the components of $\widehat{H}$ are given by
\begin{align*}
\widehat{H}_{1} & =\dfrac{r_{1}r_{2}r_{3}}{\left(1-q_{1}\right)^{2}\left(1-q_{2}\right)\left(1-q_{3}\right)}+\dfrac{r_{1}r_{2}r_{3}}{\left(1-q_{1}\right)\left(1-q_{2}\right)^{2}\left(1-q_{3}\right)}+\dfrac{r_{1}r_{2}r_{3}}{\left(1-q_{1}\right)\left(1-q_{2}\right)\left(1-q_{3}\right)^{2}},\\
\widehat{H}_{2} & =\dfrac{r_{2}r_{3}}{\left(1-q_{2}\right)^{2}\left(1-q_{3}\right)}+\dfrac{r_{2}r_{3}}{\left(1-q_{2}\right)\left(1-q_{3}\right)^{3}},\\
\widehat{H}_{3} & =\dfrac{r_{3}}{\left(1-q_{3}\right)^{2}}.
\end{align*}
The average time needed to obtain a good component is
\[
H\left(1\right)=\mathbb{E}_{1}\left[T_{4}\boldsymbol{1}_{\left(T_{4}<+\infty\right)}\right]
\]
and therefore\boxeq{
\[
H\left(1\right)=\dfrac{r_{1}r_{2}r_{3}}{\left(1-q_{1}\right)\left(1-q_{2}\right)\left(1-q_{3}\right)}\left[\dfrac{1}{1-q_{1}}+\dfrac{1}{1-q_{2}}+\dfrac{1}{1-q_{3}}\right].
\]
}

\textbf{7. Computation of $\phi\left(x\right)\equiv\phi_{C,D}\left(x\right).$
Interpretation of $\phi\left(1\right)$}

If $C=\left\{ 2\right\} $ and $D=\left\{ 4\right\} ,$ then $\phi\left(1\right)$
is the probability that a component reaches the second stage of production
and is never good. From $\refpar{eq:system_phi_CD_x},$ it follows
that
\begin{equation}
\phi=P_{2}\left(T_{4}=+\infty\right)\left(\begin{array}{l}
r_{1}\\
q_{2}\\
0\\
0\\
0
\end{array}\right)+Q\phi,\label{eq:phi_eq}
\end{equation}
where $Q$ is the matrix $M$ with the second column replaced by zeros.

Since $P_{x}\left(T_{2}<+\infty\right)=0$ for $x=3,4,5,$ it follows
from the definition of $\phi$ that $\phi\left(x\right)=0$ if $x=3,4,5.$
The equation $\refpar{eq:phi_eq}$ therefore yields
\[
\left\{ \begin{array}{l}
\left(1-q_{1}\right)\phi\left(1\right)=r_{1}P_{2}\left(T_{4}=+\infty\right)\\
\phi\left(2\right)=q_{2}P_{2}\left(T_{4}=+\infty\right).
\end{array}\right.
\]
Now, from question 4,
\[
P_{2}\left(T_{4}=+\infty\right)=1-K\left(2\right)=1-\dfrac{r_{2}r_{3}}{\left(1-q_{2}\right)\left(1-q_{3}\right)}.
\]
Hence,\boxeq{
\[
\begin{array}{l}
\phi\left(1\right)=\dfrac{r_{1}}{1-q_{1}}\left[1-\dfrac{r_{2}r_{3}}{\left(1-q_{2}\right)\left(1-q_{3}\right)}\right],\\
\phi\left(2\right)=q_{2}\left[1-\dfrac{r_{2}r_{3}}{\left(1-q_{2}\right)\left(1-q_{3}\right)}\right].
\end{array}
\]
}

\end{solution}

\begin{solution}{}{solexercise17.8}

\textbf{1. Computation of $\mathbb{E}^{\mathscr{A}_{n}}\left[f\left(Y_{n+1}\right)\right].$
Proof that $Y=\left(Y_{n}\right)_{n\in\mathbb{N}}$ is a homogeneous
Markov chain taking values in $E.$ Determination of $M.$ Proof of
$M^{n}=M^{2}$}

For every function $f$ on $E,$ and every $\left(x_{0},\cdots,x_{n+1}\right)\in\left\{ 0,1\right\} ^{n+2},$
we have
\begin{align*}
\mathbb{E}^{\left(X_{0},\cdots,X_{n+1}\right)=\left(x_{0},\cdots,x_{n+1}\right)}\left[f\left(Y_{n+1}\right)\right] & =\mathbb{E}^{\left(X_{0},\cdots,X_{n+1}\right)=\left(x_{0},\cdots,x_{n+1}\right)}\left[f\left(X_{n+1},X_{n+2}\right)\right]\\
 & =\mathbb{E}^{\left(X_{0},\cdots,X_{n+1}\right)=\left(x_{0},\cdots,x_{n+1}\right)}\left[f\left(x_{n+1},X_{n+2}\right)\right].
\end{align*}
Hence, since the random variables $\left(X_{0},\cdots,X_{n+1}\right)$
and $f\left(x_{n+1},X_{n+2}\right)$ are independent
\begin{align*}
\mathbb{E}^{\left(X_{0},\cdots,X_{n+1}\right)=\left(x_{0},\cdots,x_{n+1}\right)}\left[f\left(Y_{n+1}\right)\right] & =\mathbb{E}\left[f\left(x_{n+1},X_{n+2}\right)\right]\\
 & =pf\left(x_{n+1},1\right)+qf\left(x_{n+1},0\right).
\end{align*}
It follows that
\[
\mathbb{E}^{\mathscr{A}_{n}}\left[f\left(Y_{n+1}\right)\right]=M\left(Y_{n},f\right),
\]
where, $M\left(\cdot,f\right)$ is defined, for every $y\in E,$ by
\[
M\left(y,f\right)=pf\left(\pi_{2}\left(y\right),1\right)+qf\left(\pi_{2}\left(y\right),0\right),
\]
and $\pi_{2}$ denotes the second projection of $E=\left\{ 0,1\right\} ^{2}$
onto $\left\{ 0,1\right\} $ defined by $\pi_{2}\left(y_{1},y_{2}\right)=y_{2}.$
This proves that $Y$ is a homogeneous Markov chain with transition
matrix $M,$ whose entries 
\[
M\left(\alpha_{i},\alpha_{j}\right)=M\left(\alpha_{i},\boldsymbol{1}_{\left\{ \alpha_{j}\right\} }\right).
\]
 Thus
\[
\begin{array}{ccc}
 &  & \begin{array}{cccc}
\alpha_{1} & \alpha_{2} & \alpha_{3} & \alpha_{4}\end{array}\\
M= & \begin{array}{c}
\alpha_{1}\\
\alpha_{2}\\
\alpha_{3}\\
\alpha_{4}
\end{array} & \left(\begin{array}{cccc}
\,0\, & \,0\, & \,p\, & \,p\,\\
\,q\, & p & 0 & 0\\
0 & 0 & q & p\\
q & p & 0 & 0
\end{array}\right).
\end{array}
\]
We have
\[
\begin{array}{cc}
M^{2}= & \left(\begin{array}{cccc}
pq & p^{2} & q^{2} & pq\\
pq & p^{2} & q^{2} & pq\\
pq & p^{2} & q^{2} & pq\\
pq & p^{2} & q^{2} & pq
\end{array}\right)\end{array}.
\]
If we suppose that, for $n\geqslant2,$ we have $M^{n}=M^{2},$ then
\[
M^{n+1}=\left(\begin{array}{cccc}
p^{2}q+pq^{2} & p^{3}+p^{2}q & pq^{2}+q^{3} & p^{2}q+pq^{2}\\
p^{2}q+pq^{2} & p^{3}+p^{2}q & pq^{2}+q^{3} & p^{2}q+pq^{2}\\
p^{2}q+pq^{2} & p^{3}+p^{2}q & pq^{2}+q^{3} & p^{2}q+pq^{2}\\
p^{2}q+pq^{2} & p^{3}+p^{2}q & pq^{2}+q^{3} & p^{2}q+pq^{2}
\end{array}\right),
\]
which, using the identity $p+q=1,$ proves that $M^{n+1}=M^{2}.$
Hence, for every $n\geqslant2,$ $M^{n}=M^{2}.$

\textbf{2. Associated graph. Nature and period of the states of $E$}

The graph associated with this chain is

\begin{center}\begin{tikzpicture}[
  >=Stealth,
  every node/.style={
    circle,
    draw=black,
    minimum size=12mm,
    inner sep=1pt
  },
  edge/.style={
    ->,
    draw=black,
    line width=0.75pt
  }
]

% Nodes (rectangular layout)
\node (a1) at (0,1.5) {\(\alpha_1\)};
\node (a2) at (0,-1.5) {\(\alpha_2\)};
\node (a3) at (4,-1.5) {\(\alpha_3\)};
\node (a4) at (4,1.5) {\(\alpha_4\)};

% Self-loops
\draw[edge] (a2) to[out=210,in=150,looseness=8] (a2);
\draw[edge] (a3) to[out=30,in=-30,looseness=8] (a3);

% Transitions
\draw[edge] (a1) -- (a3);
\draw[edge] (a1) to[bend right=20] (a4);

\draw[edge] (a2) -- (a1);
\draw[edge] (a4) to[bend right=20] (a1);

\draw[edge] (a3) -- (a4);

\draw[edge] (a4) -- (a2);

\end{tikzpicture}
\end{center}

This graph shows that all states communicate. Hence, there is only
one communication class, and the chain is \textbf{irreducible}. Since
$E$ is finite, the chain is positive recurrent. Moreover, since $M\left(\alpha_{2},\alpha_{2}\right)>0,$
the state $\alpha_{2}$ is aperiodic. Therefore the whole chain $y$
is \textbf{aperiodic}.

\textbf{3. a. Proof that the sequence $\left(f_{k}\right)_{k\in\mathbb{N}^{\ast}}$
of vectors in $\mathbb{R}^{4}$ satisfies a first-order recurrence
relation }

By Proposition $\ref{pr:first_hitting_point},$ and more precisely
the equality $\refpar{eq:F_k(x_y)},$ we have, for every $x\in E,$
\[
\begin{cases}
f_{1}\left(x\right)=M\left(x,\alpha_{2}\right),\\
f_{k}\left(x\right)=\sum_{z\in E\backslash\left\{ \alpha_{2}\right\} }M\left(x,z\right)f_{k-1}\left(z\right), & \text{if }k\geqslant2.
\end{cases}
\]
It follows that the sequence of vectors $f_{k}$ satisfies the recurrence
relation $f_{k}=Qf_{k-1},$ $k\geqslant2,$ where $f_{1}$ is the
second column of $M,$ and $Q$ is the matrix $M$ by zeroing its
second column, namely
\[
f_{1}=\left(\begin{array}{c}
0\\
p\\
0\\
p
\end{array}\right)\,\,\,\,\text{and}\,\,\,\,Q=\left(\begin{array}{cccc}
0 & 0 & q & p\\
q & 0 & 0 & 0\\
0 & 0 & q & p\\
q & 0 & 0 & 0
\end{array}\right).
\]

\textbf{b. Proof that $\left(f_{k}\left(1\right)\right)_{k\in\mathbb{N}^{\ast}}$
satisfies $x_{k}=qx_{k-1}+pqx_{k-2}$ }

That is, the components of the sequence $\left(f_{k}\right)_{k\in\mathbb{N}^{\ast}}$
satisfy the system
\[
\left\{ \begin{array}{l}
f_{k}\left(1\right)=qf_{k-1}\left(3\right)+pf_{k-1}\left(4\right)\\
f_{k}\left(2\right)=qf_{k-1}\left(1\right)\\
f_{k}\left(3\right)=qf_{k-1}\left(3\right)+pf_{k-1}\left(4\right)\\
f_{k}\left(4\right)=qf_{k-1}\left(1\right)
\end{array}\right..
\]
This system is equivalent to the system
\[
\left\{ \begin{array}{l}
f_{k}\left(1\right)=f_{k}\left(3\right)=qf_{k-1}\left(3\right)+pf_{k-1}\left(4\right)\\
f_{k}\left(2\right)=f_{k}\left(4\right)=qf_{k-1}\left(1\right)
\end{array}\right.,
\]
Hence, if $k\geqslant3,$
\[
f_{k}\left(1\right)=qf_{k-1}\left(1\right)+pqf_{k-2}\left(1\right).
\]
Therefore, the sequence of probabilities $\left(f_{k}\left(1\right)\right)_{k\in\mathbb{N}^{\ast}}$
satisfies the second-order recurrence relation, for $k\geqslant3$
\begin{equation}
x_{k}=qx_{k-1}+pqx_{k-2}\label{eq:x_k_order_two_eq_rec_f_k_1}
\end{equation}
with initial conditions\boxeq{
\begin{equation}
x_{1}=0,\,\,\,\,x_{2}=p^{2}.\label{eq:init_cond_x_k}
\end{equation}
}Indeed,
\[
x_{1}=f_{1}\left(1\right)=M\left(\alpha_{1},\alpha_{2}\right)=0,
\]
and, since the random variables $X_{n}$ are independent,
\begin{align*}
x_{2} & =f_{2}\left(1\right)=P_{\alpha_{1}}\left(Y_{1}\neq\alpha_{2},Y_{2}=\alpha_{2}\right)=\dfrac{P\left(Y_{0}=\alpha_{1},Y_{1}\neq\alpha_{2},Y_{2}=\alpha_{2}\right)}{P\left(Y_{0}=\alpha_{1}\right)}\\
 & =\dfrac{P\left(X_{0}=1,X_{1}=0,X_{2}=1,X_{3}=1\right)}{P\left(X_{0}=1,X_{1}=0\right)}=p^{2}.
\end{align*}

\textbf{c. Particular solution of the form $\left(\lambda^{k}\right)_{k\geqslant2}$ }

The characteristic equation associated with $\refpar{eq:x_k_order_two_eq_rec_f_k_1}$
is
\[
\lambda^{2}-q\lambda-pq=0.
\]
Its discriminant is $\Delta=q^{2}+4pq=q\left(1+3p\right)$ and it
has the two roots
\[
\lambda_{1}=\dfrac{q-\sqrt{q\left(1+3p\right)}}{2}\,\,\,\,\text{and}\,\,\,\,\lambda_{2}=\dfrac{q+\sqrt{q\left(1+3p\right)}}{2}.
\]
The general solution of $\refpar{eq:x_k_order_two_eq_rec_f_k_1}$
is therefore of the form $x_{k}=\beta_{1}\lambda^{k}_{1}+\beta_{2}\lambda^{k}_{2}.$
The sequence $\left(f_{k}\left(1\right)\right)_{k\in\mathbb{N}^{\ast}}$
is the solution corresponding to the initial conditions given by $\refpar{eq:init_cond_x_k}.$
The constants $\beta_{1}$ and $\beta_{2}$ satisfy
\[
\left\{ \begin{array}{l}
\beta_{1}\lambda_{1}+\beta_{2}\lambda_{2}=0\\
\beta_{1}\lambda^{2}_{1}+\beta_{2}\lambda^{2}_{2}=p^{2}
\end{array}\right.\Longleftrightarrow\left\{ \begin{array}{l}
\beta_{1}=\dfrac{p^{2}}{\lambda_{1}\left(\lambda_{1}-\lambda_{2}\right)}\\
\beta_{2}=\dfrac{p^{2}}{\lambda_{2}\left(\lambda_{2}-\lambda_{1}\right)}
\end{array}\right..
\]
Thus, after observing that the resulting formula is also valid for
$k=2,$\boxeq{
\[
P_{\alpha_{1}}\left(T=k\right)=f_{k}\left(1\right)=\begin{cases}
0, & \text{if }k=1,\\
\dfrac{p^{2}}{\lambda_{2}-\lambda_{1}}\left[\lambda^{k-1}_{2}-\lambda^{k-1}_{1}\right], & \text{if }k\geqslant2.
\end{cases}
\]
}

\textbf{4. Computation of $\mathbb{E}_{\alpha_{1}}\left(T\right)$}

We then have
\[
\mathbb{E}_{\alpha_{1}}\left(T\right)=\sum^{+\infty}_{k=2}kP_{\alpha_{1}}\left(T=k\right)=\dfrac{p^{2}}{\lambda_{2}-\lambda_{1}}\sum^{+\infty}_{k=2}k\left[\lambda^{k-1}_{2}-\lambda^{k-1}_{1}\right].
\]
By the classical differentiation formula for power series, if $\left|x\right|<1,$
\[
\sum^{+\infty}_{k=2}kx^{k-1}=\left(\sum^{+\infty}_{k=0}x^{k}\right)^{\prime}-1=\dfrac{1}{\left(1-x\right)^{2}}-1=\dfrac{x\left(2-x\right)}{\left(1-x\right)^{2}}.
\]
Taking into account the values of $\lambda_{1}$ and $\lambda_{2}$---and
noting that $\left|\lambda_{1}\right|<1$ and $\left|\lambda_{2}\right|<1$---a
straightforward computation then yields\boxeq{
\[
\mathbb{E}_{\alpha_{1}}\left(T\right)=\dfrac{2-q}{p^{2}}=\dfrac{1+p}{p^{2}}.
\]
}

\begin{remark}{}{}We give, as an indication, the table of numerical
values, for some values of $p,$

\begin{center}%
\begin{tabular}{|c|c|c|c|c|}
\hline 
$p$ &
0.1 &
0.5 &
0.8 &
0.9\tabularnewline
\hline 
$\mathbb{E}_{\alpha_{1}}\left(T\right)$ &
110 &
6 &
2.81 &
2.34\tabularnewline
\hline 
\end{tabular}\end{center}

\end{remark}

\textbf{5. Justification of the existence of a unique invariant probability
$\nu.$ Computation of $\nu.$ Value of $\mathbb{E}_{\alpha_{2}}\left(T\right)$
and comparison with $\mathbb{E}_{\alpha_{1}}\left(T\right)$}

Since the chain $Y$ is irreducible and positive recurrent, there
exists a unique invariant probability $\nu.$ If $\nu$ is an invariant
measure, then in particular $\nu M^{2}=\nu.$ Hence, for every $y\in E,$
\[
\sum^{4}_{j=1}v\left(\alpha_{j}\right)M^{2}\left(\alpha_{j},y\right)=v\left(y\right).
\]
But, here, $M^{2}\left(\alpha_{j},y\right)$ is independent of $j.$
It follows that, for every $y\in E,$
\[
M^{2}\left(\alpha_{j},y\right)=\dfrac{v\left(y\right)}{\sum^{4}_{j=1}v\left(\alpha_{j}\right)}.
\]
Hence, $v=\left(pq,p^{2},q^{2},pq\right)$ is the invariant probability.
It follows that\boxeq{
\[
\mathbb{E}_{\alpha_{2}}\left(T\right)=\dfrac{1}{v\left(\alpha_{2}\right)}=\dfrac{1}{p^{2}}.
\]
}It is therefore clear that $\mathbb{E}_{\alpha_{2}}\left(T\right)<\mathbb{E}_{\alpha_{1}}\left(T\right).$

\begin{remark}{}{}We give, as an indication, the table of numerical
values, for some values of $p,$

\begin{center}%
\begin{tabular}{|c|c|c|c|c|}
\hline 
$p$ &
0.1 &
0.5 &
0.8 &
0.9\tabularnewline
\hline 
$\mathbb{E}_{\alpha_{2}}\left(T\right)$ &
100 &
4 &
1.56 &
1.23\tabularnewline
\hline 
\end{tabular}\end{center}

\end{remark}

\end{solution}

\begin{solution}{}{solexercise17.9}

\textbf{1. Communication class(es) }

The associated graph associate with this chain is

\begin{center}\begin{tikzpicture}[
  scale=0.75,
  transform shape,
  >=Stealth,
  every node/.style={
    circle,
    draw=black,
    minimum size=8mm,
    inner sep=1pt
  },
  edge/.style={
    ->,
    draw=black,
    line width=0.75pt
  }
]

% Nodes (rectangular layout)
\node (a0) at (0,0) {0};
\node (a1) at (2.5,0) {1};
\node (a2) at (5,0) {2};
\node[draw=none] (a3) at (7.5,0) {};
\node[draw=none] at (8.0,0)  {\(\cdot\)};
\node[draw=none] at (8.5,0)  {\(\cdot\)};
\node[draw=none] at (9.0,0)  {\(\cdot\)};
\node[draw=none] at (9.5,0) {\(\cdot\)};
\node[draw=none] (aNm2) at (10,0) {};
\node (aNm1) at (12.5,0) {$N-1$};
\node (aN) at (15,0) {$N$};

% Transitions
\draw[edge] (a0) to[bend right=20] (a1);
\draw[edge] (a1) to[bend right=20] (a0);
\draw[edge] (a1) to[bend right=20] (a2);
\draw[edge] (a2) to[bend right=20] (a1);
\draw[edge] (a2) to[bend right=20] (a3);
\draw[edge] (a3) to[bend right=20] (a2);
\draw[edge] (aNm2) to[bend right=20] (aNm1);
\draw[edge] (aNm1) to[bend right=20] (aNm2);
\draw[edge] (aNm1) to[bend right=20] (aN);
\draw[edge] (aN) to[bend right=20] (aNm1);

\end{tikzpicture}
\end{center}

All states communicate. Hence, there is only one communication class,
and the chain is irreducible. Since $E$ is finite, this chain is
therefore positive recurrent.

\textbf{2. Existence and computation of a unique invariant probability
$\nu$. Value of the mean time $\mathbb{E}_{0}\left(T_{0}\right)$}

It follows that there exists a unique invariant probability $\nu.$
Let us write that $\nu$ is an invariant measure. Then, for every
$x\in E,$
\[
\nu\left(x\right)=\sum_{z\in E}\nu\left(z\right)M\left(z,x\right).
\]
Taking into account the values of $M,$ this is equivalent to the
system
\[
\begin{cases}
\nu\left(x\right)=p\nu\left(x-1\right)+q\nu\left(x+1\right), & \text{if}\,\,\,\,2\leqslant x\leqslant N-2,\\
\nu\left(1\right)=\nu\left(0\right)+q\nu\left(2\right)\\
\nu\left(0\right)=q\nu\left(1\right)\\
\nu\left(N\right)=p\nu\left(N-1\right)\\
\nu\left(N-1\right)=p\nu\left(N-2\right)+\nu\left(N\right).
\end{cases}
\]
Then, for $2\leqslant y\leqslant N-2,$
\[
\nu\left(y+1\right)-\nu\left(y\right)=\dfrac{p}{q}\left[\nu\left(y\right)-\nu\left(y-1\right)\right].
\]
Summing over $y,$ for $2\leqslant x\leqslant N-2$
\[
\nu\left(x+1\right)-\nu\left(2\right)=\dfrac{p}{q}\sum^{x}_{y=2}\left[\nu\left(y\right)-\nu\left(y-1\right)\right]=\dfrac{p}{q}\left[\nu\left(x\right)-\nu\left(1\right)\right].
\]
Hence,
\[
\nu\left(x+1\right)=\dfrac{1}{q}\left(\dfrac{1}{q}-1\right)\nu\left(0\right)+\dfrac{p}{q}\nu\left(x\right)-\dfrac{p}{q^{2}}\nu\left(0\right),
\]
which simplifies to
\begin{equation}
\nu\left(x+1\right)=\dfrac{p}{q}\nu\left(x\right)\,\,\,\,\text{if }2\leqslant x\leqslant N-2.\label{eq:nu_xp1}
\end{equation}
Moreover,
\[
q\nu\left(2\right)=\nu\left(1\right)-\nu\left(0\right)=\nu\left(1\right)-q\nu\left(1\right)=p\nu\left(1\right),
\]
which proves that the equality $\refpar{eq:nu_xp1}$ is also valid
for $x=1.$ It then follows from equality $\refpar{eq:nu_xp1}$ that---paying
attention to the index shift---for $2\leqslant x\leqslant N-1,$
\begin{equation}
\nu\left(x\right)=\left(\dfrac{p}{q}\right)^{x-1}\nu\left(1\right).\label{eq:nu_x_f_nu1}
\end{equation}

Let us now compute the mass of this measure. Two cases must be distinguished:
\begin{itemize}
\item \textbf{Case $p\neq q$}\\
We have
\begin{align*}
\sum^{N}_{x=0}\nu\left(x\right) & =\nu\left(0\right)+\dfrac{1}{q}\nu\left(0\right)+\sum^{N-1}_{x=2}\left[\left(\dfrac{p}{q}\right)^{x-1}\dfrac{1}{q}\nu\left(0\right)\right]+\left(\dfrac{p}{q}\right)^{N-1}\nu\left(0\right)\\
 & =\nu\left(0\right)\left[1+\dfrac{1}{q}+\dfrac{1}{q}\dfrac{p}{q}\dfrac{1-\left(\dfrac{p}{q}\right)^{N-2}}{1-\dfrac{p}{q}}+\left(\dfrac{p}{q}\right)^{N-1}\right]\\
 & =\nu\left(0\right)\left[1+\dfrac{1}{q}\left(1+\dfrac{p}{q-p}\right)+\left(\dfrac{p}{q}\right)^{N-1}\left(1+\dfrac{1}{q-p}\right)\right]\\
 & =\nu\left(0\right)\left[1+\dfrac{1}{q-p}+\dfrac{2p}{q-p}\left(\dfrac{p}{q}\right)^{N-1}\right].
\end{align*}
Thus,
\[
\sum^{N}_{x=0}\nu\left(x\right)=\nu\left(0\right)\dfrac{2q}{q-p}\left[1-\left(\dfrac{p}{q}\right)^{N}\right].
\]
To obtain the unique invariant probability, it is enough to normalize
this measure. Hence it is given, using the previous expressions, by
\[
\left\{ \begin{array}{l}
\nu\left(0\right)=\dfrac{1-\frac{p}{q}}{2\left[1-\left(\frac{p}{q}\right)^{N}\right]},\,\,\,\,\nu\left(1\right)=\dfrac{1-\frac{p}{q}}{2q\left[1-\left(\frac{p}{q}\right)^{N}\right]},\\
\nu\left(x\right)=\dfrac{1}{2q}\dfrac{1-\frac{p}{q}}{1-\left(\frac{p}{q}\right)^{N}}\left(\dfrac{p}{q}\right)^{x-1},\,\,\,\,\text{if }2\leqslant x\leqslant N-1,\\
\nu\left(N\right)=\dfrac{1}{2}\dfrac{1-\frac{p}{q}}{1-\left(\frac{p}{q}\right)^{N}}\left(\dfrac{p}{q}\right)^{N-1}.
\end{array}\right.
\]
The mean return time to 0, namely $\mathbb{E}_{0}\left(T_{0}\right),$
is therefore equal to $\dfrac{1}{\nu\left(0\right)},$ hence\boxeq{
\[
\mathbb{E}_{0}\left(T_{0}\right)=\dfrac{2\left[1-\left(\frac{p}{q}\right)^{N}\right]}{1-\frac{p}{q}}.
\]
}
\item Case $p=q=\dfrac{1}{2}$\\
Since the equality $\refpar{eq:nu_xp1}$ holds for every $1\leqslant x\leqslant N-2,$
we have
\[
\nu\left(x\right)=\nu\left(1\right),\,\,\,\,\text{if }2\leqslant x\leqslant N-1.
\]
Moreover, in this case,
\[
\nu\left(0\right)=\dfrac{1}{2}\nu\left(1\right)\,\,\,\,\text{and}\,\,\,\,\nu\left(N\right)=\dfrac{1}{2}\nu\left(N-1\right).
\]
Thus in fact
\[
\nu\left(0\right)=\dfrac{1}{2}\nu\left(1\right)\,\,\,\,\text{and}\,\,\,\,\nu\left(N\right)=\dfrac{1}{2}\nu\left(N-1\right)=\dfrac{1}{2}\nu\left(1\right).
\]
Hence the mass of $\nu$ is 
\[
\nu\left(E\right)=\left(N-1\right)\nu\left(1\right).
\]
It follows that the unique invariant probability is given by
\[
\nu\left(x\right)=\begin{cases}
\dfrac{1}{N-1}, & \text{if }2\leqslant x\leqslant N-1,\\
\dfrac{1}{2\left(N-1\right)}, & \text{if }x=0\text{ or }2.
\end{cases}
\]
The mean return time to $0,$ namely $\mathbb{E}_{0}\left(T_{0}\right),$
is therefore equal to $\dfrac{1}{\nu\left(0\right)},$ thus\boxeq{
\[
\mathbb{E}_{0}\left(T_{0}\right)=2\left(N-1\right).
\]
}
\end{itemize}
\end{solution}

\begin{solution}{}{solexercise17.10}

\textbf{1. Communication class(es). Periodicity of the states}

We study the graph associated with this chain.
\begin{itemize}
\item If $\alpha\in\left[0,1\right[,$ the graph associated with this chain
is\\
\begin{center}\begin{tikzpicture}[
  scale=0.75,
  transform shape,
  >=Stealth,
  every node/.style={
    circle,
    draw=black,
    minimum size=8mm,
    inner sep=1pt
  },
  edge/.style={
    ->,
    draw=black,
    line width=0.75pt
  }
]

% Nodes (rectangular layout)
\node (a0) at (0,0) {0};
\node (a1) at (2.5,0) {1};
\node (a2) at (5,0) {2};
\node[draw=none] (a3) at (7.5,0) {};
\node[draw=none] at (7.75,0)  {\(\cdot\)};
\node[draw=none] at (8,0)  {\(\cdot\)};
\node[draw=none] at (8.25,0)  {\(\cdot\)};
\node[draw=none] (aNm2) at (8.5,0) {};
\node (aNm1) at (11.0,0) {$N$};
\node[draw=none] (aN) at (13.5,0) {};
\node[draw=none] at (13.75,0)  {\(\cdot\)};
\node[draw=none] at (14,0)  {\(\cdot\)};
\node[draw=none] at (14.25,0)  {\(\cdot\)};

% Transitions
\draw[edge] (a0) to[bend right=20] (a1);
\draw[edge] (a1) to[bend right=20] (a0);
\draw[edge] (a1) to[bend right=20] (a2);
\draw[edge] (a2) to[bend right=20] (a1);
\draw[edge] (a2) to[bend right=20] (a3);
\draw[edge] (a3) to[bend right=20] (a2);
\draw[edge] (aNm2) to[bend right=20] (aNm1);
\draw[edge] (aNm1) to[bend right=20] (aNm2);
\draw[edge] (aNm1) to[bend right=20] (aN);
\draw[edge] (aN) to[bend right=20] (aNm1);

\end{tikzpicture}
\end{center}All points communicate. Hence, there is only one communication
class, and the chain is \textbf{irreducible}. Since $\mathbb{N}$
is infinite, we cannot determine a priori the nature of this chain.
This is precisely the aim of the problem.
\item If $\alpha=1,$ the graph associated with this chain is\\
\begin{center}\begin{tikzpicture}[
  scale=0.75,
  transform shape,
  >=Stealth,
  every node/.style={
    circle,
    draw=black,
    minimum size=8mm,
    inner sep=1pt
  },
  edge/.style={
    ->,
    draw=black,
    line width=0.75pt
  }
]

% Nodes (rectangular layout)
\node (a0) at (0,0) {0};
\node (a1) at (2.5,0) {1};
\node (a2) at (5,0) {2};
\node[draw=none] (a3) at (7.5,0) {};
\node[draw=none] at (7.75,0)  {\(\cdot\)};
\node[draw=none] at (8,0)  {\(\cdot\)};
\node[draw=none] at (8.25,0)  {\(\cdot\)};
\node[draw=none] (aNm2) at (8.5,0) {};
\node (aNm1) at (11,0) {$N$};
\node[draw=none] (aN) at (13.5,0) {};
\node[draw=none] at (13.75,0)  {\(\cdot\)};
\node[draw=none] at (14,0)  {\(\cdot\)};
\node[draw=none] at (14.25,0)  {\(\cdot\)};

% Self-loops
\draw[edge] (a0) to[out=120,in=60,looseness=6] (a0);

% Transitions
%\draw[edge] (a0) to[bend right=20] (a1);
\draw[edge] (a1) to[bend left=20] (a0);
\draw[edge] (a1) to[bend right=20] (a2);
\draw[edge] (a2) to[bend right=20] (a1);
\draw[edge] (a2) to[bend right=20] (a3);
\draw[edge] (a3) to[bend right=20] (a2);
\draw[edge] (aNm2) to[bend right=20] (aNm1);
\draw[edge] (aNm1) to[bend right=20] (aNm2);
\draw[edge] (aNm1) to[bend right=20] (aN);
\draw[edge] (aN) to[bend right=20] (aNm1);

\end{tikzpicture}
\end{center}All points of $\mathbb{N}^{\ast}$ communicate, but 0
and 1 do not. Indeed, $M\left(0,0\right)=1,$ and if we suppose that
$M^{n}\left(0,0\right)=1,$ then, since $M\left(z,0\right)=0$ as
soon as $z\geqslant2,$
\begin{align*}
M^{n+1}\left(0,0\right) & =\sum_{z\in\mathbb{N}}M^{n}\left(0,z\right)M\left(z,0\right)\\
 & =M^{n}\left(0,0\right)M\left(0,0\right)+M^{n}\left(0,1\right)M\left(1,0\right)=1,
\end{align*}
which proves that $M^{n}\left(0,0\right)=1$ for every $n\in\mathbb{N}^{\ast},$
and therefore that $0$ does not lead to 1. In this case, there are
two communication classes $\left\{ 0\right\} $ and $\mathbb{N}^{\ast}.$
\end{itemize}
As for periodicity,
\begin{itemize}
\item If $\alpha=0,$ we have $M\left(0,0\right)=0$ and, by a straightforward
induction, for every $n\in\mathbb{N}^{\ast},$ $M^{2n}\left(0,0\right)>0$
and $M^{2n+1}\left(0,0\right)=0.$ Thus, 0 has period 2, and so does
every point of $\mathbb{N}.$
\item If $0<\alpha<1,$ we have $M\left(0,0\right)>0,$ which implies that
0 is aperiodic. The same holds for every point of $\mathbb{N}.$
\item If $\alpha=1,$ we have $M\left(0,0\right)>0,$ which implies that
the class reduced to the single point 0 is aperiodic. Moreover, $M\left(1,1\right)=0$
and, by a straightforward induction, for every $n\in\mathbb{N}^{\ast},$
$M^{2n}\left(1,1\right)>0$ and $M^{2n+1}\left(1,1\right)=0.$ Thus,
1 has period 2, and so does every points of the class $\mathbb{N}^{\ast}.$
\end{itemize}
\textbf{We first study the case $\alpha\in\left]0,1\right[.$}

\textbf{2. a. Proof of an invariant measure $\nu.$ }

An invariant measure $\nu$ verifies, for every $x\in\mathbb{N},$
\[
\nu\left(x\right)=\sum_{z\in\mathbb{N}}\nu\left(z\right)M\left(z,x\right),
\]
which, taking into account the value of $M,$ is equivalent to the
system
\begin{equation}
\left\{ \begin{array}{l}
\nu\left(x\right)=p\nu\left(x-1\right)+q\nu\left(x+1\right),\,\,\,\,\text{if\,}\,2\leqslant x,\\
\nu\left(1\right)=\left(1-\alpha\right)\nu\left(0\right)+q\nu\left(2\right),\\
\nu\left(0\right)=\alpha\nu\left(0\right)+q\nu\left(1\right).
\end{array}\right.\label{eq:nu_systemex17.10}
\end{equation}
Then,
\[
\nu\left(1\right)=\dfrac{1-\alpha}{q}\nu\left(0\right),
\]
 and thus,
\[
q\nu\left(2\right)=\dfrac{1-\alpha}{q}\nu\left(0\right)-\left(1-\alpha\right)\nu\left(0\right),
\]
which gives, after simplification,
\[
\nu\left(2\right)=\left(1-\alpha\right)\dfrac{p}{q^{2}}\nu\left(0\right).
\]
It follows, by substitution into the system $\refpar{eq:nu_systemex17.10},$
that
\[
q\nu\left(3\right)=\left(1-\alpha\right)p\nu\left(0\right)\left[\dfrac{1}{q^{2}}-\dfrac{1}{q}\right]=\left(1-\alpha\right)\dfrac{p^{2}}{q^{2}}\nu\left(0\right),
\]
and thus
\[
\nu\left(3\right)=\left(1-\alpha\right)\dfrac{p^{2}}{q^{3}}\nu\left(0\right).
\]

Now, suppose that, for every $i$ such that $1\leqslant i\leqslant j,$
\begin{equation}
\nu\left(i\right)=\dfrac{1-\alpha}{q}\left(\dfrac{p}{q}\right)^{i-1}\nu\left(0\right).\label{eq:nu_i_general_form}
\end{equation}
By a similar computation,
\begin{align*}
q\nu\left(j+1\right) & =\dfrac{1-\alpha}{q}\left(\dfrac{p}{q}\right)^{j-1}\nu\left(0\right)-p\dfrac{1-\alpha}{q}\left(\dfrac{p}{q}\right)^{j-2}\nu\left(0\right)\\
 & =\dfrac{1-\alpha}{q}\left(\dfrac{p}{q}\right)^{j-2}p\left[\dfrac{1}{q}-1\right]\nu\left(0\right)\\
 & =\left(1-\alpha\right)\left(\dfrac{p}{q}\right)^{j}\nu\left(0\right),
\end{align*}
and therefore
\[
\nu\left(j+1\right)=\dfrac{1-\alpha}{q}\left(\dfrac{p}{q}\right)^{j}\nu\left(0\right).
\]
It follows that $\refpar{eq:nu_i_general_form}$ holds for every $j\in\mathbb{N}^{\ast}.$
Thus, invariant measures $\nu$ do exist. They are all proportional
and are determined up to a multiplicative constant by equality $\refpar{eq:nu_i_general_form}.$

\textbf{b. Existence and uniqueness of an invariant probability. Computation
when it exists}
\begin{itemize}
\item If $\alpha=1,$ the Dirac measure is the unique invariant probability.
\item If $\alpha\in\left[0,1\right[,$ \textbf{these invariant measures
have finite mass only if $p<q.$} In this case, their mass is given
by
\begin{align*}
\nu\left(\mathbb{N}\right) & =\nu\left(0\right)+\sum_{j\in\mathbb{N}^{\ast}}\nu\left(j\right)=\nu\left(0\right)+\sum_{j\in\mathbb{N}^{\ast}}\dfrac{1-\alpha}{q}\left(\dfrac{p}{q}\right)^{j-1}\nu\left(0\right)\\
 & =\left[1+\dfrac{1-\alpha}{q}\dfrac{1}{1-\frac{p}{q}}\right]\nu\left(0\right).
\end{align*}
Hence,
\[
\nu\left(\mathbb{N}\right)=\dfrac{2q-\alpha}{q-p}\nu\left(0\right).
\]
Therefore, the unique invariant probability $\nu$ is given by\boxeq{
\begin{equation}
\nu\left(j\right)=\begin{cases}
\dfrac{q-p}{2q-\alpha}, & \text{if }j=0,\\
\dfrac{q-p}{2q-\alpha}\cdot\dfrac{1-\alpha}{q}\cdot\left(\dfrac{p}{q}\right)^{j-1}, & \text{if }j\in\mathbb{N}^{\ast}.
\end{cases}\label{eq:nu_j_for_alpha_0_1_p_less_q}
\end{equation}
}
\end{itemize}
\textbf{c. When $\alpha\in\left[0,1\right[$ and $p<q,$ nature of
the points of $\mathbb{N}.$ Value of $\mathbb{E}_{x}\left(T_{x}\right)$}

Hence, if $\alpha\in\left[0,1\right[$ and $p<q,$ the chain is irreducible
and positive recurrent, and $\mathbb{E}_{j}\left(T_{j}\right)=\dfrac{1}{\nu\left(j\right)},$
thus\boxeq{
\[
\nu\left(j\right)=\begin{cases}
\dfrac{2q-\alpha}{q-p}, & \text{if }j=0,\\
\dfrac{2q-\alpha}{1-\alpha}\cdot\dfrac{q}{q-p}\left(\dfrac{q}{p}\right)^{j-1}, & \text{if }j\in\mathbb{N}^{\ast}.
\end{cases}
\]
}

\textbf{3. When $p\geqslant q,$ nature of the points of $\mathbb{N},$
and if $p>q,$ computation of $P_{x}\left(T_{0}=+\infty\right)$}

If $\alpha\in\left[0,1\right[$ and $p\geqslant q,$ there is no invariant
probability, and the chain is irreducible, hence either transient
or null recurrent. To determine its nature, we apply the analytical
recurrence criterion---Corollary $\ref{co:chain_recurrent}$: the
chain is recurrent if and only if the system
\begin{equation}
h=Qh\,\,\,\,0\leqslant h\leqslant1,\label{eq:cond_rec_chain}
\end{equation}
has the unique solution $h=0,$ where $Q$ is the restriction of $M$
to $\mathbb{N}^{\ast}\times\mathbb{N}^{\ast}.$ Now,
\[
Qh\left(x\right)=\sum_{y\in\mathbb{N}^{\ast}}Q\left(x,y\right)h\left(y\right)=\begin{cases}
qh\left(x-1\right)+ph\left(x+1\right), & \text{if }x\geqslant2,\\
ph\left(2\right), & \text{if }x=1.
\end{cases}
\]
Thus, the system $\refpar{eq:cond_rec_chain}$ can be written as
\[
\left\{ \begin{array}{l}
h\left(1\right)=ph\left(2\right),\\
h\left(x\right)=qh\left(x-1\right)+ph\left(x+1\right),\,\,\,\,\text{if }x\geqslant2,\\
0\leqslant h\left(x\right)\leqslant1.
\end{array}\right.
\]
The second equation yields, for $x\geqslant2,$
\[
q\left[h\left(x\right)-h\left(x-1\right)\right]=p\left[h\left(x+1\right)-h\left(x\right)\right].
\]
Hence, by backward iteration
\[
h\left(x+1\right)-h\left(x\right)=\left(\dfrac{q}{p}\right)^{x-1}\left[h\left(2\right)-h\left(1\right)\right]=\left(\dfrac{q}{p}\right)^{x}h\left(1\right).
\]
We note that this equality is still valid for $x=1.$ It follows,
by summation, that for every $x\in\mathbb{N}^{\ast},$
\[
h\left(x+1\right)-h\left(1\right)=\sum^{x}_{j=1}\left(\dfrac{q}{p}\right)^{j}h\left(1\right),
\]
or equivalently,
\begin{equation}
h\left(x+1\right)=h\left(1\right)\sum^{x}_{j=1}\left(\dfrac{q}{p}\right)^{j}.\label{eq:h_x_pp_1_form}
\end{equation}

\begin{itemize}
\item If $p=q,$ then $h\left(x\right)=xh\left(1\right),$ and the unique
solution to the system $\refpar{eq:cond_rec_chain}$ is $h\equiv0.$
In this case, \textbf{the chain is null recurrent}.
\item If $p>q,$ for every $x\geqslant2,$
\[
h\left(x\right)=h\left(1\right)\dfrac{1-\left(\dfrac{q}{p}\right)^{x}}{1-\dfrac{q}{p}}=h\left(1\right)\dfrac{p}{p-q}\left[1-\left(\dfrac{q}{p}\right)^{x}\right].
\]
Hence, if we choose $h\left(1\right)=\dfrac{p-q}{p}=1-\dfrac{q}{p}<1,$
the function $h$ defined, for every $x\in\mathbb{N}^{\ast},$ by
\[
h\left(x\right)=1-\left(\dfrac{q}{p}\right)^{x}
\]
is a nonzero solution of the system $\refpar{eq:cond_rec_chain}.$
Therefore, in this case, \textbf{the chain is transient}. The probability
\[
P_{x}\left(T_{0}=+\infty\right)=P_{x}\left[\bigcap_{j\in\mathbb{N}^{\ast}}\left(X_{j}\in\mathbb{N}^{\ast}\right)\right]
\]
is equal to $f\left(x\right),$ where $f$ is the maximal solution
of the system $\refpar{eq:cond_rec_chain}$---see Proposition $\ref{pr:max_sol_fix_point_eq}.$
This is precisely the function $h,$ since it corresponds to the largest
possible value of $h\left(1\right)$ under the constraint $0\leqslant h\leqslant1.$
Hence, for every $x\in\mathbb{N}^{\ast},$ \boxeq{
\[
P_{x}\left(T_{0}=+\infty\right)=1-\left(\dfrac{q}{p}\right)^{x}.
\]
}
\end{itemize}
\textbf{4. When $p<q,$ proof of the $P_{x}-$almost sure convergence
of the sequence with general term $\dfrac{1}{n}\sum^{n}_{j=1}\exp\left(-aX_{j}\right),$
for any $a>0$ }

If $\alpha\in\left[0,1\right[$ and $p<q,$ the chain is irreducible
and positive recurrent. Therefore, the strong law of large numbers
for homogeneous Markov chains applies to every function that is integrable
with respect to the invariant probability $\nu.$ This is the case
here, since the function under consideration is bounded. By $\refpar{eq:nu_j_for_alpha_0_1_p_less_q},$
\begin{align*}
\intop_{\mathbb{R}}\text{e}^{-ax}\text{d}\nu\left(x\right) & =\dfrac{q-p}{2q-\alpha}+\dfrac{q-p}{2q-\alpha}\cdot\dfrac{1-\alpha}{q}\sum_{j\in\mathbb{N}^{\ast}}\text{e}^{-aj}\left(\dfrac{p}{q}\right)^{j-1}\\
 & =\dfrac{q-p}{2q-\alpha}\left[1+\dfrac{1-\alpha}{q}\text{e}^{-a}\sum_{j\in\mathbb{N}^{\ast}}\left(\dfrac{p}{q}\text{e}^{-a}\right)^{j-1}\right]\\
 & =\dfrac{q-p}{2q-\alpha}\left[1+\dfrac{1-\alpha}{q}\dfrac{\text{e}^{-a}}{1-\frac{p}{q}\text{e}^{-a}}\right].
\end{align*}
Hence,
\[
\intop_{\mathbb{R}}\text{e}^{-ax}\text{d}\nu\left(x\right)=\dfrac{q-p}{2q-\alpha}\left[1+\left(1-\alpha\right)\dfrac{\text{e}^{-a}}{q-p\text{e}^{-a}}\right].
\]

The strong law of large numbers proves the $P_{x}-$almost sure convergence
of the sequence with general term $\dfrac{1}{n}\sum^{n}_{j=1}\text{e}^{-aX_{j}}$
and ensures that $P_{x}-$almost surely,
\[
\lim_{n\to+\infty}\dfrac{1}{n}\sum^{n}_{j=1}\text{e}^{-aX_{j}}=\dfrac{q-p}{2q-\alpha}\left[1+\left(1-\alpha\right)\dfrac{\text{e}^{-a}}{q-p\text{e}^{-a}}\right].
\]

\textbf{5. Case where $\alpha=1.$ Computation of $P_{0}\left(T_{0}<+\infty\right)$
and $\mathbb{E}_{0}\left(T_{0}\right).$ Nature of the point 0}

If $\alpha=1,$ we have seen that the Dirac measure is the unique
invariant probability. Hence, 0 is recurrent, and therefore 
\[
P_{0}\left(T_{0}<+\infty\right)=1.
\]
 A more elementary argument is to say that
\[
P_{0}\left(T_{0}<+\infty\right)\geqslant P_{0}\left(X_{1}=0\right)=M\left(0,0\right)=1.
\]
Thus, $T_{0}=1$ $P_{0}-$almost surely, and consequently $E_{0}\left(T_{0}\right)=1.$
Hence, the point $0$ is positive recurrent.

\textbf{6. Computation of $P_{x}\left(T_{0}<+\infty\right).$ Determination
of the nature of the points of $\mathbb{N}^{\ast}.$}

The probability 
\[
P_{x}\left(T_{0}<+\infty\right)=P_{x}\left[\bigcap_{j\in\mathbb{N}^{\ast}}\left(X_{j}\in\mathbb{N}^{\ast}\right)\right]
\]
is equal to $f\left(x\right),$ where $f$ is the maximal solution
of the system $\refpar{eq:cond_rec_chain}$---see Proposition $\ref{pr:max_sol_fix_point_eq}.$
Its general solution $h$ is still of the form $\refpar{eq:h_x_pp_1_form}:$
\begin{itemize}
\item If $p\leqslant q,$ and if $h\left(1\right)>0,$ then 
\[
\lim_{x\to+\infty}h\left(x\right)=+\infty.
\]
Therefore, the only solution of the system $\refpar{eq:cond_rec_chain}$
is the zero solution, which implies that, for every $x\in\mathbb{N}^{\ast},$
$P_{x}\left(T_{0}<+\infty\right)=1.$
\item If $p>q,$ the maximal solution of the system $\refpar{eq:cond_rec_chain}$
is the same as in the case $\alpha<1,$ and we still have, for every
$x\in\mathbb{N}^{\ast},$
\[
P_{x}\left(T_{0}=+\infty\right)=1-\left(\dfrac{q}{p}\right)^{x},
\]
 or equivalently,\boxeq{
\[
P_{x}\left(T_{0}<+\infty\right)=\left(\dfrac{q}{p}\right)^{x}.
\]
}
\end{itemize}
We have seen that, when $\alpha=1,$ the state 0 does not lead to
1. The state $0$ is recurrent, whereas 1 is transient, since otherwise
the points 0 and 1 would communicate. Hence, \textbf{when $\alpha=1,$
the class $\mathbb{N}^{\ast}$ is transient}.

\textbf{7. Convergence of the sequence with general term $M^{n}\left(x,y\right)$---limit,
if it exists---when $x\in\mathbb{N}$ and $y\in\mathbb{N}^{\ast},$
and when $x\in\mathbb{N}^{\ast}$ and $y=0$}

If $\alpha=1,$ then for $x\in\mathbb{N}$ and $y\in\mathbb{N}^{\ast},$
\[
\lim_{n\to+\infty}M^{n}\left(x,y\right)=0.
\]
Moreover, since 0 is recurrent and aperiodic, we have, for $x\in\mathbb{N}^{\ast},$
\[
\lim_{n\to+\infty}M^{n}\left(x,0\right)=\dfrac{P_{x}\left(T_{0}<+\infty\right)}{\mathbb{E}_{0}\left(T_{0}\right)}=P_{x}\left(T_{0}<+\infty\right).
\]
Hence,\boxeq{
\[
\lim_{n\to+\infty}M^{n}\left(x,0\right)=\begin{cases}
\left(\dfrac{q}{p}\right)^{x}, & \text{if }p>q,\\
1, & \text{if }p\leqslant q.
\end{cases}
\]
}

\end{solution}

\begin{solution}{}{solexercise17.11}

\textbf{1. $X$ is a homogeneous Markov chain with transition matrix
$M$ }

For every bounded function $f$ on $\mathbb{N},$ for, every $n\in\mathbb{N}^{\ast},$
and for every $\left(x_{1},\cdots,x_{n}\right)\in\mathbb{N}^{n},$
we have
\begin{align*}
 & \mathbb{E}^{\left(X_{0},X_{1},\cdots,X_{n}\right)=\left(a,x_{1},\cdots,x_{n}\right)}\left[f\left(X_{n+1}\right)\right]\\
 & \qquad\qquad=\mathbb{E}^{\left(X_{0},X_{1},\cdots,X_{n}\right)=\left(a,x_{1},\cdots,x_{n}\right)}\left[f\left(\boldsymbol{1}_{\left(x_{n}\geqslant1\right)}\sum^{x_{n}}_{j=1}Y_{n,j}\right)\right].
\end{align*}
Hence, since the random variables $\left(X_{0},X_{1},\cdots,X_{n}\right)$
and $Y_{n,j},\,j\in\mathbb{N}^{\ast},$ are independent,
\[
\mathbb{E}^{\left(X_{0},X_{1},\cdots,X_{n}\right)=\left(a,x_{1},\cdots,x_{n}\right)}\left[f\left(X_{n+1}\right)\right]=\mathbb{E}\left[f\left(\boldsymbol{1}_{\left(x_{n}\geqslant1\right)}\sum^{x_{n}}_{j=1}Y_{n,j}\right)\right].
\]
Since the random variables $Y_{n,j},\,j\in\mathbb{N}^{\ast},$ are
independent and following the same law $\mu,$ 
\[
\mathbb{E}^{\left(X_{0},X_{1},\cdots,X_{n}\right)=\left(a,x_{1},\cdots,x_{n}\right)}\left[f\left(X_{n+1}\right)\right]=\begin{cases}
f\left(0\right), & \text{if }x_{n}=0,\\
\intop f\text{d}\mu^{\ast x_{n}}, & \text{if }x_{n}\in\mathbb{N}^{\ast}.
\end{cases}
\]
It follows that
\begin{equation}
\mathbb{\mathbb{E}}^{\mathscr{A}_{n}}\left[f\left(X_{n+1}\right)\right]=M\left(X_{n},f\right),\label{eq:exp_f_X_npl1_know_A_n}
\end{equation}
where, $\delta_{0}$ denotes the Dirac measure at 0, and $M\left(\cdot,f\right)$
is defined, for every $x\in\mathbb{N},$ by
\[
M\left(x,f\right)=\boldsymbol{1}_{\left\{ 0\right\} }\left(x\right)\intop f\text{d}\delta_{0}+\boldsymbol{1}_{\mathbb{N}^{\ast}}\left(x\right)\intop f\text{d}\mu^{\ast x}.
\]
This proves that $X$ is a homogeneous Markov chain with transition
matrix $M,$ whose entries are 
\[
M\left(x,y\right)=M\left(x,\boldsymbol{1}_{\left\{ y\right\} }\right),
\]
hence\boxeq{
\[
M\left(x,y\right)=\begin{cases}
\boldsymbol{1}_{\left\{ 0\right\} }\left(y\right), & \text{if }x=0,\\
\mu^{\ast x}\left(\left\{ y\right\} \right), & \text{if }x\in\mathbb{N}^{\ast}.
\end{cases}
\]
}

\textbf{2. $X$ is a martingale, a submartingale, or a supermartingale,
depending on the value of $m$}

The process $X$ is adapted and nonnegative. The same computation
of the conditional expectation for nonnegative random variables---not
necessarily integrable---yields, for every $n\in\mathbb{N}^{\ast}$
and every $\left(x_{1},\cdots,x_{n}\right)\in\mathbb{N}^{n},$
\[
\mathbb{E}^{\left(X_{0},X_{1},\cdots,X_{n}\right)=\left(a,x_{1},\cdots,x_{n}\right)}\left[X_{n+1}\right]=\mathbb{E}^{\left(X_{0},X_{1},\cdots,X_{n}\right)=\left(a,x_{1},\cdots,x_{n}\right)}\left[\boldsymbol{1}_{\left(x_{n}\geqslant1\right)}\sum^{x_{n}}_{j=1}Y_{n,j}\right].
\]
Hence, since the random variables $\left(X_{0},X_{1},\cdots,X_{n}\right)$
and $Y_{n,j},\,j\in\mathbb{N}^{\ast},$ are independent,
\[
\mathbb{E}^{\left(X_{0},X_{1},\cdots,X_{n}\right)=\left(a,x_{1},\cdots,x_{n}\right)}\left[X_{n+1}\right]=\boldsymbol{1}_{\left(x_{n}\geqslant1\right)}\mathbb{E}\left[\sum^{x_{n}}_{j=1}Y_{n,j}\right].
\]
Since the random variables $Y_{n,j},\,j\in\mathbb{N}^{\ast}$ are
independent and following the same law $\mu,$ of expectation $m,$
\begin{align*}
\mathbb{E}^{\left(X_{0},X_{1},\cdots,X_{n}\right)=\left(a,x_{1},\cdots,x_{n}\right)}\left[X_{n+1}\right] & =\boldsymbol{1}_{\left(x_{n}\geqslant1\right)}\sum^{x_{n}}_{j=1}\mathbb{E}\left(Y_{n,j}\right)\\
 & =\boldsymbol{1}_{\mathbb{N}^{\ast}}\left(x_{n}\right)mx_{n}=mx_{n}.
\end{align*}
It follows that\boxeq{
\begin{equation}
\mathbb{\mathbb{E}}^{\mathscr{A}_{n}}\left(X_{n+1}\right)=mX_{n}.\label{eq:exp_f_X_npl1_know_A_n-1}
\end{equation}
}

The process $X$ is a martingale---integrable, since then $\mathbb{E}\left(X_{n+1}\right)=\mathbb{E}\left(X_{0}\right)=a$---if
$m=1,$ a nonnegative supermartingale if $m<1,$ and a nonnegative
submartingale if $m>1.$

\textbf{3. $Y$ is a nonnegative martingale}

Dividing both sides of the equality $\refpar{eq:exp_f_X_npl1_know_A_n-1}$
by $m^{n+1},$ we obtain
\[
\mathbb{\mathbb{E}}^{\mathscr{A}_{n}}\left(\dfrac{X_{n+1}}{m^{n+1}}\right)=\dfrac{X_{n}}{m^{n}},
\]
which proves that $Y$ is a positive martingale.

\textbf{4. $Z_{n}=s^{X_{n}}.$ Proof that $Z$ is an equi-integrable
submartingale}

Since $s\in\left]0,1\right[,$ the function $p:x\mapsto s^{x}$ is
bounded by 1 on $\mathbb{N}.$ It follows immediately that the sequence
of random variables $Z_{n}$ is equi-integrable. Moreover, the equality
$\refpar{eq:exp_f_X_npl1_know_A_n}$ yields, for this function
\[
\mathbb{\mathbb{E}}^{\mathscr{A}_{n}}\left(s^{X_{n+1}}\right)=M\left(X_{n},p\right)
\]
where, for every $x\in\mathbb{N},$
\[
M\left(x,p\right)=\boldsymbol{1}_{\left\{ 0\right\} }\left(x\right)s^{0}+\boldsymbol{1}_{\mathbb{N}^{\ast}}\left(x\right)\intop s^{y}\text{d}\mu^{\ast x}\left(y\right).
\]
Now---by a classical result on convolutions of measures---,
\[
\intop s^{y}\text{d}\mu^{\ast x}\left(y\right)=\left[\intop s^{y}\text{d}\mu\left(y\right)\right]^{x}=\left[g\left(s\right)\right]^{x}.
\]
Since $g\left(s\right)=s,$ it follows that
\[
M\left(x,p\right)=\boldsymbol{1}_{\left\{ 0\right\} }\left(x\right)s^{0}+\boldsymbol{1}_{\mathbb{N}^{\ast}}\left(x\right)s^{x}=s^{x},
\]
and thus
\[
\mathbb{\mathbb{E}}^{\mathscr{A}_{n}}\left(s^{X_{n+1}}\right)=s^{X_{n}}.
\]
Hence, $Z$ is an equi-integrable martingale.

\textbf{5. $\left(X_{n}\right)_{n\in\mathbb{N}}$ converges $P-$almost
surely to a random variable $X_{\infty}$---separate study of the
cases $0<m\leqslant1$ and $m>1.$ Limit $X_{\infty}$ when $0<m<1$}

There are three cases to distinguish.
\begin{itemize}
\item If $m\leqslant1,$ $X$ is a nonnegative supermartingale, so it converges,
hence $P-$almost surely, to a random variable $X_{\infty}$ such
that, for every $n\in\mathbb{N},$
\[
0\leqslant\mathbb{\mathbb{E}}^{\mathscr{A}_{n}}\left(X_{\infty}\right)\leqslant X_{n}.
\]
Moreover $X_{n}$ is integrable since 
\[
0\leqslant\mathbb{\mathbb{E}}^{\mathscr{A}_{0}}\left(X_{n}\right)\leqslant X_{0}=a.
\]
If follows that $X_{\infty}$ is integrable and that 
\[
0\leqslant\mathbb{E}\left(X_{\infty}\right)\leqslant\mathbb{E}\left(X_{n}\right).
\]
Moreover, by the equality $\refpar{eq:exp_f_X_npl1_know_A_n-1},$
\[
\mathbb{E}\left(X_{n+1}\right)=m\mathbb{E}\left(X_{n}\right)=m^{n+1}a.
\]

\begin{itemize}
\item If $m<1,$ then 
\[
\lim_{n\to+\infty}\mathbb{E}\left(X_{n}\right)=0
\]
and thus, by the Fatou lemma, 
\[
0\leqslant\mathbb{E}\left(X_{\infty}\right)\leqslant\liminf_{n\to+\infty}\mathbb{E}\left(X_{n}\right)=0.
\]
Since $X_{\infty}$ is nonnegative, $X_{\infty}=0,$ $P-$almost surely.
\end{itemize}
\begin{remark}{}{}

If $m=1,$ the last argument does not allow us, a priori, to concluded
anything about the finiteness of $X_{\infty}.$

\end{remark}
\item If $m=1,$ $X$ is a nonnegative martingale such that $\mathbb{E}\left(X_{n}\right)=a.$
That is, it is bounded in $\text{L}^{1}$ and therefore converges
$P-$almost surely to a nonnegative random variable $X_{\infty},$
which is $P-$almost surely finite.
\item If $m>1,$ $Z=s^{X}$ is an equi-integrable martingale---hence bounded
in $\text{L}^{1}.$ It converges $P-$almost surely and in $\text{L}^{1}$
to a nonnegative random variable $U_{\infty}$ $P-$almost surely
finite. It follows that the sequence $\left(X_{n}\right)_{n\in\mathbb{N}}$
converges $P-$almost surely to a nonnegative random variable $X_{\infty}.$
It satisfies
\begin{equation}
\left[\lim_{n\to+\infty}s^{X_{n}}\right]\boldsymbol{1}_{\left(X_{\infty}=+\infty\right)}=0\label{eq:eq_lim_sX_n}
\end{equation}
\end{itemize}
\textbf{6. Probability $P\left[\bigcap^{k}_{n=N}\left(X_{n}=j\right)\right]$
as a function of $M\left(j,j\right)$ and $P\left(X_{N}=j\right).$
Proof of $P\left[\liminf_{n\to+\infty}\left(X_{n}=j\right)\right]=0$}

Let $j\in\mathbb{N}^{\ast}.$ Since $X$ is a homogeneous Markov chain
with transition matrix $M,$ for every $k>N,$
\begin{align*}
 & P\left[\bigcap^{k}_{n=N}\left(X_{n}=j\right)\right]=P\left(X_{k}=j\mid X_{k-1}=j\right)P\left(X_{k-1}=j\mid X_{k-2}=j\right)\cdots P\left(X_{N}=j\right),
\end{align*}
thus,
\[
P\left[\bigcap^{k}_{n=N}\left(X_{n}=j\right)\right]=\left[M\left(j,j\right)\right]^{k-N}P\left(X_{N}=j\right).
\]
But, since $\mu\left(\left\{ 0\right\} \right)>0,$ we have $\mu^{\ast j}\left(\left\{ 0\right\} \right)>0,$
and thus $M\left(j,j\right)=\mu^{\ast}\left(\left\{ j\right\} \right)<1.$
It follows that, for every $N\in\mathbb{N}^{\ast},$
\[
\lim_{k\to+\infty}P\left[\bigcap^{k}_{n=N}\left(X_{n}=j\right)\right]=P\left[\bigcap_{n\geqslant N}\left(X_{n}=j\right)\right]=0,
\]
and, therefore immediately
\[
P\left[\liminf_{n\to+\infty}\left(X_{n}=j\right)\right]=0.
\]

\textbf{7. Proof that, for every $j\in\mathbb{N}^{\ast},$ $P\left(X_{\infty}=j\right)=0,$
and that $X_{\infty}\in\left\{ 0,+\infty\right\} $ $P-$almost surely.
Every point of $\mathbb{N}^{\ast}$ is transient. If $m>1,$ proof
that $P\left(X_{\infty}=0\right)=s^{a}\,\text{and}\,P\left(X_{\infty}=+\infty\right)=1-s^{a}$}

Since the sequence $\left(X_{n}\right)_{n\in\mathbb{N}}$ takes values
in $\mathbb{N}$ and converges $P-$almost surely in $\overline{\mathbb{N}},$
for every $j\in\mathbb{N}^{\ast},$
\[
\left(X_{\infty}=j\right)\subset\liminf_{n\to+\infty}\left(X_{n}=j\right),
\]
which, by the previous question, proves that for every $j\in\mathbb{N}^{\ast},$
\[
P\left(X_{\infty}=j\right)=0,
\]
and thus $X_{\infty}\in\left\{ 0,+\infty\right\} $ $P-$almost surely.
If $j\in\mathbb{N}^{\ast}$ were recurrent, then we would have 
\[
P\left[\limsup_{n\to+\infty}\left(X_{n}=j\right)\right]=1,
\]
and since the sequence $\left(X_{n}\right)_{n\in\mathbb{N}}$ converges
to $X_{\infty}$ on $\limsup_{n\to+\infty}\left(X_{n}=j\right),$
we would have $X_{\infty}=j$ $P-$almost surely, which is a contradiction.
Thus, every point of $\mathbb{N}^{\ast}$ is transient.

If $m>1,$ we have seen that $Z$ is an equi-integrable martingale
and that the sequence $\left(Z_{n}\right)_{n\in\mathbb{N}}$ converges
in $\text{L}^{1}.$ Thus,
\[
s^{a}=\mathbb{E}\left(s^{X_{0}}\right)=\lim_{n\to+\infty}\mathbb{E}\left(s^{X_{n}}\right)=\mathbb{E}\left(\lim_{n\to+\infty}s^{X_{n}}\right).
\]
Consequently, since $X_{\infty}\in\left\{ 0,+\infty\right\} $ $P-$almost
surely, it follows from the equality $\refpar{eq:eq_lim_sX_n}$ that
\[
s^{a}=\mathbb{E}\left(s^{X_{\infty}}\boldsymbol{1}_{\left(X_{\infty}=0\right)}\right)=P\left(X_{\infty}=0\right).
\]
Thus,\boxeq{
\[
P\left(X_{\infty}=0\right)=s^{a}\,\,\,\,\text{and}\,\,\,\,P\left(X_{\infty}=+\infty\right)=1-s^{a}.
\]
}

\textbf{8. Proof that $P-$almost surely $\left(X_{\infty}=+\infty\right)=\liminf_{n\to+\infty}\left(X_{n}\neq0\right)=\left(T=+\infty\right).$
Value of the probability $P\left(T<+\infty\right)$ for the different
values of $m$}

Since we have the implication
\[
X_{n}\left(\omega\right)=0\Longrightarrow\left(\forall p\geqslant n,\,X_{p}\left(\omega\right)=0\right),
\]
we have
\[
\liminf_{n\to+\infty}\left(X_{n}\neq0\right)\subset\bigcap_{n\in\mathbb{N}^{\ast}}\left(X_{n}\geqslant1\right)=\left(T=+\infty\right),
\]
the converse inclusion being immediate. Hence,
\[
\liminf_{n\to+\infty}\left(X_{n}\neq0\right)=\left(T=+\infty\right).
\]

Moreover, we proved that, $P-$almost surely,
\[
\liminf_{n\to+\infty}\left(X_{n}\neq0\right)=\left(X_{\infty}=+\infty\right).
\]
It follows that $P\left(T=+\infty\right)=P\left(X_{\infty}=+\infty\right),$
and thus $P\left(T<+\infty\right)=P\left(X_{\infty}=0\right).$
\begin{itemize}
\item If $m>1,$ then\boxeq{
\[
P\left(T<+\infty\right)=s^{a}\,\,\,\,\text{and}\,\,\,\,P\left(T=+\infty\right)=1-s^{a}.
\]
}
\item If $m<1,$ we have seen that $X_{\infty}=0$ $P-$almost surely. Thus
$P\left(T<+\infty\right)=1.$
\item If $m=1,$ $X$ is a martingale, and we saw that $P\left(X_{\infty}<+\infty\right)=1.$
Thus $P\left(X_{\infty}=0\right)=1,$ and consequently, $P\left(T<+\infty\right)=1.$
\end{itemize}
\end{solution}

\begin{solution}{}{solexercise17.12}

\textbf{1. Proof that $Z$ is a bounded martingale. Proof of the $P-$almost
sure convergence and $\mathscr{L}^{1}$ convergence of $\left(Z_{n}\right)_{n\in\mathbb{N}^{\ast}}$}

Note that
\[
B_{n+1}=B_{n}+cX_{n+1}.
\]

This allows us to make $Z_{n}$ explicit on each of the events $\left(X_{n+1}=0\right)$
and $\left(X_{n+1}=1\right).$
\begin{itemize}
\item On $\left(X_{n+1}=0\right),$ we have $B_{n+1}=B_{n}$ and thus
\[
Z_{n+1}=\dfrac{B_{n}\left(B_{n}+c\right)\cdots\left[B_{n}+\left(l-1\right)c\right]}{k_{n+1}k_{n+2}\cdots k_{n+\left(l-1\right)}k_{n+l}}=\dfrac{k_{n}}{k_{n+l}}Z_{n}.
\]
\item On $\left(X_{n+1}=1\right),$ we have $B_{n+1}=B_{n}+c$ and thus
\[
Z_{n+1}=\dfrac{\left(B_{n}+c\right)\left(B_{n}+2c\right)\cdots\left[B_{n}+c+\left(l-1\right)c\right]}{k_{n+1}k_{n+2}\cdots k_{n+\left(l-1\right)}k_{n+l}}=\dfrac{k_{n}}{k_{n+l}}Z_{n}\dfrac{B_{n}+lc}{B_{n}}.
\]
\end{itemize}
Hence, taking into account the adaptedness of the processes $Z$ and
$B,$
\begin{align*}
\mathbb{E}^{\mathscr{A}_{n}}\left(Z_{n+1}\right) & =\mathbb{E}^{\mathscr{A}_{n}}\left(\boldsymbol{1}_{\left(X_{n+1}=0\right)}Z_{n+1}\right)+\mathbb{E}^{\mathscr{A}_{n}}\left(\boldsymbol{1}_{\left(X_{n+1}=1\right)}Z_{n+1}\right)\\
 & =\dfrac{k_{n}}{k_{n+l}}Z_{n}\mathbb{E}^{\mathscr{A}_{n}}\left(\boldsymbol{1}_{\left(X_{n+1}=0\right)}\right)+\dfrac{k_{n}}{k_{n+l}}Z_{n}\dfrac{B_{n}+lc}{B_{n}}\mathbb{E}^{\mathscr{A}_{n}}\left(\boldsymbol{1}_{\left(X_{n+1}=1\right)}\right).
\end{align*}
Hence, by the equalities $\refpar{eq:cond_pr_X_nplus1_kn_gen_sig_alg_Y_1Y_n}$
of the Example $\ref{ex:polya_diff_mod_2},$ together with the equality
$B_{n}=k_{n}Y_{n},$
\[
\mathbb{E}^{\mathscr{A}_{n}}\left(Z_{n+1}\right)=\dfrac{k_{n}}{k_{n+l}}Z_{n}\left[\left(1-Y_{n}\right)+\dfrac{B_{n}+lc}{B_{n}}Y_{n}\right]=\dfrac{k_{n}}{k_{n+l}}Z_{n}\left[1+\dfrac{lc}{B_{n}}Y_{n}\right].
\]
Now, since $B_{n}=k_{n}Y_{n}$ and $k_{n+l}=k_{n}+lc,$
\[
\dfrac{k_{n}}{k_{n+l}}\left[1+\dfrac{lc}{k_{n}}\right]=1,
\]
and therefore\boxeq{
\[
\mathbb{E}^{\mathscr{A}_{n}}\left(Z_{n+1}\right)=Z_{n}.
\]
}

Hence, \textbf{$Z$ is a martingale}. Since $\left|Y_{n}\right|\leqslant1,$
for every $n\in\mathbb{N}^{\ast},$ we have also $\left|Z_{n}\right|\leqslant1$
for every $n\in\mathbb{N}^{\ast}.$ Thus, the martingale $Z$ is equi-integrable,
and therefore converges $P-$almost surely and in $\mathscr{L}^{1}$
to an integrable random variable $Z_{\infty}.$

\textbf{2. Values of $\mathbb{E}\left(Y^{l}_{\infty}\right)$ }

Moreover, since for every $n\in\mathbb{N}^{\ast},$
\[
Z_{n}=\prod^{l}_{j=0}Y_{n+j},
\]
and since the sequence $\left(Y_{n}\right)_{n\in\mathbb{N}^{\ast}}$
converges $P-$almost surely to $Y_{\infty},$ we have
\[
Z_{\infty}=Y^{l}_{\infty}\,\,\,\,P-\text{almost surely.}
\]
Since the convergence of the martingale $Z$ also holds in $\mathscr{L}^{1},$
\boxeq{
\begin{equation}
\mathbb{E}\left(Z_{1}\right)=\mathbb{E}\left(Z_{\infty}\right)=\mathbb{E}\left(Y^{l}_{\infty}\right).\label{eq:exp_Z1}
\end{equation}
}

It therefore remains to compute $\mathbb{E}\left(Z_{1}\right).$ Since
$B_{1}=b+cX_{1},$
\[
Z_{1}=\dfrac{\left[b+cX_{1}\right]\left[b+c\left(1+X_{1}\right)\right]\cdots\left[b+c\left(l-1+X_{1}\right)\right]}{k_{1}k_{2}\cdots k_{l}}.
\]
But the law of $X_{1}$ is
\[
P_{X_{1}}=\dfrac{b}{b+r}\delta_{1}+\dfrac{r}{b+r}\delta_{0},
\]
hence
\begin{align*}
 & \mathbb{E}\left[\left[b+cX_{1}\right]\left[b+c\left(1+X_{1}\right)\right]\cdots\left[b+c\left(l-1+X_{1}\right)\right]\right]\\
 & \qquad=\dfrac{b}{b+r}\left[\left(b+c\right)\left(b+2c\right)\cdots\left(b+lc\right)\right]+\dfrac{r}{b+r}\left[b\left(b+c\right)\left(b+2c\right)\cdots\left(b+\left(l-1\right)c\right)\right]\\
 & \qquad=\dfrac{b}{b+r}\left(b+c\right)\left(b+2c\right)\cdots\left(b+\left(l-1\right)c\right)\left[b+lc+r\right].
\end{align*}
Moreover,
\[
k_{1}k_{2}\cdots k_{l}=\left(b+r+c\right)\left(b+r+2c\right)\cdots\left(b+r+lc\right),
\]
which yields
\[
\mathbb{E}\left(Z_{1}\right)=\dfrac{b}{b+r}\dfrac{\left(b+c\right)\left(b+2c\right)\cdots\left(b+\left(l-1\right)c\right)\left[b+lc+r\right]}{\left(b+r+c\right)\left(b+r+2c\right)\cdots\left(b+r+lc\right)}.
\]
After simplications,\boxeq{
\[
\mathbb{E}\left(Y^{l}_{\infty}\right)=\mathbb{E}\left(Z_{1}\right)=\dfrac{\Gamma\left(l+\dfrac{b}{c}\right)\Gamma\left(\dfrac{b+r}{c}\right)}{\Gamma\left(l+\dfrac{b+r}{c}\right)\Gamma\left(\dfrac{b}{c}\right)}.
\]
}

\textbf{3. Computation of the moment of order $l$ for a random variable
$U$ following the Beta law $\beta\left(a,b\right)$ of the first
kind on $\left[0,1\right]$}

\[
\mathbb{E}\left(U^{l}\right)=\dfrac{\Gamma\left(a+b\right)}{\Gamma\left(a\right)\Gamma\left(b\right)}\intop^{1}_{0}x^{a+l-1}\left(1-x\right)^{b-1}\text{d}x=\dfrac{\Gamma\left(a+b\right)}{\Gamma\left(a\right)\Gamma\left(b\right)}\beta\left(a+l,b\right),
\]
and therefore\boxeq{
\[
\mathbb{E}\left(U^{l}\right)=\dfrac{\Gamma\left(l+a\right)\Gamma\left(a+b\right)}{\Gamma\left(l+a+b\right)\Gamma\left(a\right)}.
\]
}

\textbf{4. Proof that the law of $Y_{\infty}$ is the Beta law $\beta\left(\dfrac{b}{c},\dfrac{r}{c}\right)$
of the first kind on $\left[0,1\right]$}

We thus observe that $Y_{\infty}$ has the same moments as a random
variable $U$ following the beta law $\beta\left(\dfrac{b}{c},\dfrac{r}{c}\right)$
of the first kind on $\left[0,1\right].$ Since these random variables
are \textbf{bounded}, their \textbf{characteristic functions are analytic}
on $\mathbb{R}$ and coincide everywhere---see Property $\ref{pr:tayl_exp_pow_ser_ch_f}.$
The Levy theorem then implies that $Y_{\infty}$ \textbf{follows the
Beta law $\beta\left(\dfrac{b}{c},\dfrac{r}{c}\right)$ of the first
kind on $\left[0,1\right].$}

\begin{remark}{}{}

In particular, if $b=r=c,$ the law of $Y_{\infty}$ is the \textbf{uniform
law} on $\left[0,1\right].$

\end{remark}

\end{solution}

\bibliographystyle{plainnat}
\phantomsection\addcontentsline{toc}{chapter}{\bibname}\bibliography{97_tmp_book_jyo_biblio}

\printindex[idx][Glossary]

\printindex[fam][Famous Characters]

\pagenumbering{gobble}

\thispagestyle{empty}

\vspace*{\fill}

\pagebreak\thispagestyle{empty}

\textbf{This volume is the English translation of the two French books
Probabilités Tome 1 and Tome 2, authored by Professor Jean-Yves Ouvrard
and originally published in 1998 and 2000 by Éditions Cassini.}

\textbf{At the time of their first publication, these books were entirely
original in French in their approach and content. They have since
remained authoritative references in the teaching of probability theory,
covering material from the undergraduate level (first part) through
to the master's programs and the preparation of France's highly competitive
examinations for mathematics teaching (second part).}

\textbf{The first part introduces the fundamentals of probability
theory, including event algebras, random variables, independence,
conditional probabilities, moments of discrete and continuous random
variables, generating functions, and limit theorems.}

\textbf{The second part, following a brief introduction to measure
theory, develops more advanced topics: probability measures and their
complements, distributions and moments of random variables, modes
of convergence, laws of large numbers, conditional expectation, Fourier
transforms and characteristic functions, Gaussian random variables,
convergence of measures, convergence in distribution, discrete-time
stochastic processes, martingales, and Markov chains.}

\textbf{The reader’s work is greatly facilitated by the inclusion,
in every chapter, of numerous exercises, all accompanied by detailed
solutions that often provide substantial extensions to the theoretical
material.}

\vspace{0.5cm}

\textbf{About the author.} Professor Jean-Yves Ouvrard has served
on numerous juries for the Agrégation de mathématiques, France’s highest
competitive examination for mathematics teaching, and was a faculty
member at Université Joseph Fourier in Grenoble, France.

\textbf{About the translator.} Dr. Xavier Ouvrard currently in computer
science works at EPFL, and has a strong interdicsciplinary background
across several scientific fields.

\vspace{0.5cm}

This work is published under the Creative Commons Attribution--NonCommercial--ShareAlike
4.0 International License (\href{https://creativecommons.org/licenses/by-nc-sa/4.0/}{CC BY-NC-SA 4.0}).

The digital content may be copied, shared, and adapted under the terms
of this license, provided that it is distributed free of charge and
without any form of financial compensation.

\vspace*{\fill}

\begin{center}\includegraphics[scale=0.3]{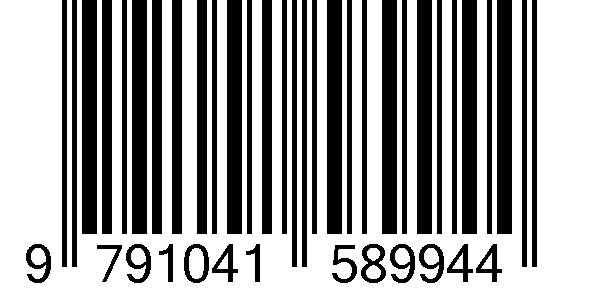}\end{center}

\end{document}